\documentclass{amsbook}

\usepackage{sty/preamble}

\usepackage{geometry}
\setamsbookgeometry 
\savegeometry{abg} 

\title{Multifunctorial Equivariant Algebraic $K$-Theory}

\authorinfoDY

\hypersetup{pdfauthor=\authors}

\date{\today}

\begin{document}
\frontmatter

\begin{abstract}
  A central question in equivariant algebraic $K$-theory asks whether there exists an equivariant $K$-theory machine from genuine symmetric monoidal $G$-categories to orthogonal $G$-spectra that preserves equivariant algebraic structures.  We answer this question positively by constructing an enriched multifunctor $\Kgo$ from the $G$-categorically enriched multicategory of $\Op$-pseudoalgebras to the symmetric monoidal category of orthogonal $G$-spectra, for a compact Lie group $G$ and a 1-connected pseudo-commutative $G$-categorical operad $\Op$.  As the main application of its enriched multifunctoriality, $\Kgo$ preserves all equivariant algebraic structures parametrized by multicategories enriched in either $G$-spaces or $G$-categories.  For example, for a finite group $G$ and the $G$-Barratt-Eccles operad $\GBE$, $\Kggbe$ transports equivariant $\Einf$-algebras, in the sense of Guillou-May or Blumberg-Hill, of genuine symmetric monoidal $G$-categories to equivariant $\Einf$-algebras of orthogonal $G$-spectra.

The construction $\Kgo$ is a composite of three enriched multifunctors.  The first step is a $G$-categorically enriched multifunctor from the $G$-categorically enriched multicategory of $\Op$-pseudoalgebras to the symmetric monoidal closed category of $\Gskg$-categories.  The second step is a $G$-topologically enriched symmetric monoidal functor from $\Gskg$-categories to $\Gskg$-spaces induced by the classifying space functor.  The last step is a $G$-topologically enriched unital symmetric monoidal functor from $\Gskg$-spaces to orthogonal $G$-spectra. 

\end{abstract}

\maketitle

\cleardoublepage
\thispagestyle{empty}
\vspace*{13.5pc}
\begin{center}
The author dedicates this work to Jacqueline.  
\end{center}
\cleardoublepage

\pdfbookmark{\contentsname}{Contents}
\tableofcontents


\chapter*{Introduction}
\newcommand{\sect}[1]{\section*{#1}}
\newcommand{\prefacepartNumName}[1]{\Cref{#1}: \nameref{#1}}
\newcommand{\prefacechapNumName}[1]{\medskip\begin{center}\Cref{#1}: \nameref{#1} \end{center}}
\newcommand{\tref}[1]{\text{\ref{#1}}}

A central question in equivariant algebraic $K$-theory asks whether there exists an equivariant $K$-theory machine from genuine symmetric monoidal $G$-categories to orthogonal $G$-spectra that preserves equivariant algebraic structures.  We answer this question positively by constructing an enriched multifunctor \pcref{thm:Kgo_multi}
\[\MultpsO \fto{\Kgo} \GSp\]
from the $G$-categorically enriched multicategory $\MultpsO$ of $\Op$-pseudoalgebras to the symmetric monoidal category $\GSp$ of orthogonal $G$-spectra, for a compact Lie group $G$ and a 1-connected pseudo-commutative $G$-categorical operad $\Op$.  Our enriched multifunctors are always assumed to preserve the symmetric group action.  As the main application of its enriched multifunctoriality, $\Kgo$ preserves all equivariant algebraic structures parametrized by multicategories enriched in either $G$-spaces or $G$-categories.  For example, for a finite group $G$ and the $G$-Barratt-Eccles operad $\GBE$, $\Kggbe$ transports equivariant $\Einf$-algebras, in the sense of Guillou-May or Blumberg-Hill, of genuine symmetric monoidal $G$-categories to equivariant $\Einf$-algebras of orthogonal $G$-spectra.  

The rest of this introduction consists of the following sections: \nameref{sec:litreview}, \nameref{sec:main_results}, \nameref{sec:chsummary}, and \nameref{sec:read_sug}.

\sect{Literature Review}\label{sec:litreview}

Before we describe our enriched multifunctor $\Kgo$ and its applications to equivariant algebraic structures in more detail, we briefly review some relevant literature on (non)equivariant algebraic $K$-theory.  This brief literature review is intended to provide motivation, background, and historical context.  The main text of this work does not use any of these (non)equivariant $K$-theory constructions.  Readers familiar with existing literature in this area can jump straight to the section titled \nameref{sec:main_results}.  

\subsection*{Nonequivariant Algebraic $K$-Theory}
We begin with a review of $K$-theory machines in the nonequivariant context.  One reason why spectra are important is that they represent \index{generalized cohomology theory}generalized cohomology theory.  As categories are generally easier to manipulate than spectra, an early motivating question in infinite loop space theory is this:
\qtn{Is there a functor from categories to spectra?}
There are several such $K$-theory machines due to Segal, May, and Elmendorf-Mandell, which we briefly review below.  

Segal $K$-theory\index{Segal K-theory@Segal $K$-theory}\index{K-theory@$K$-theory!Segal} \cite{segal} is a composite of three functors
\[\begin{tikzpicture}[vcenter]
\def\h{2.5} \def\u{.7}
\draw[0cell]
(0,0) node (a) {\PermCat}
(a)++(\h,0) node (b) {\Fskelcat}
(b)++(\h,0) node (c) {\phantom{\Fskelsset}}
(c)++(0,.045) node (c') {\Fskelsset}
(c)++(\h,0) node (d) {\Sp}
;
\draw[1cell=.9]
(a) edge node {\Jse} (b)
(b) edge node {\Ner} (c)
(c) edge node {\Klse} (d)
;
\draw[1cell=1]
(a) [rounded corners=3pt] |- ($(b)+(0,\u)$)
-- node {\Kse} ($(c)+(0,\u)$) -| (d)
;
\end{tikzpicture}\]
from the category $\PermCat$ of small permutative categories and strictly unital symmetric monoidal functors, to the category $\Sp$ of symmetric spectra.  In the intermediate diagram categories, $\Fskel$ is a small skeleton of the category of pointed finite sets and pointed functions.  The levelwise nerve functor $\Ner$ and the functor $\Klse$, which is also called Segal $K$-theory in the literature, are both symmetric monoidal functors.

To get a flavor of Segal $K$-theory, we first note that both $\Ner$, given by post-composing with the nerve functor, and $\Klse$, given objectwise by a coend, are somewhat formal categorical constructions.  On the other hand, Segal $J$-theory $\Jse$ involves a notion of gluing that is reminiscent of other $K$-theory constructions.  Recall that $K_0$\index{K0@$K_0$} of a ring is the Grothendieck group\index{Grothendieck group} of the commutative monoid of equivalence classes of finitely generated projective modules, under direct sums of modules.  Complex topological $K$-theory\index{K-theory@$K$-theory!topological} is the Grothendieck group of the commutative monoid of equivalence classes of complex vector bundles, under the Whitney sum, which means fiberwise direct sum.  Segal $J$-theory of a small permutative category $(\C,\oplus)$ is based on gluing morphisms 
\begin{equation}\label{segal_gluing}
x_s \oplus x_t \to x_{s \sqcup t}
\end{equation}
for disjoint subsets $s$ and $t$ of $\ufs{n} = \{1,2,\ldots,n\}$.  We will not use Segal $K$-theory in this work; the reader is referred to \cite[Ch.\! 8]{cerberusIII} for more discussion of Segal $K$-theory.  

By the work of Bousfield-Friedlander \cite{bousfield_friedlander}, Thomason \cite{thomason}, and Mandell \cite{mandell_inverseK}, each of the functors $\Jse$, $\Ner$, $\Klse$, and $\Kse$ is an \index{equivalences of homotopy theories}equivalence of homotopy theories.  Furthermore, Mandell's inverse $K$-theory\index{K-theory@$K$-theory!inverse} functor \cite{mandell_inverseK}, from $\Fskelcat$ to $\PermCat$, is a pseudo symmetric $\Cat$-multifunctor by a result in \cite[Ch.\! 10]{yau-multgro}; see also \cite{johnson-yau-invK}.

May's operadic infinite loop space machine\index{infinite loop space machine!May} \cite{may} produces connective spectra from small permutative categories via $\Einf$-spaces.  The first step of May's machine sends a small permutative category $\C$ to the classifying space $\cla\C$ as an algebra over the topological Barratt-Eccles $\Einf$-operad $\cla\BE$.  Using the product trick in \cite[3.10]{may}, $\Einf$-spaces are also given by algebras over the Steiner $\Einf$-operad $\Stein$.  The second step of May's machine sends a $\Stein$-algebra $X$ to a connective spectrum $\mathbbm{E}X$.  The latter is defined using a two-sided monadic bar construction and a natural $\Stein$-action on the 0th-space of a spectrum.  

At first glance, the two infinite loop space machines due to Segal and May look very different.  The uniqueness theorem of May and Thomason \cite{may-thomason} proves that the two machines are equivalent when applied to equivalent input data.  The reader is referred to \cite{may-groupcompletion,may77,may-construction,may-good,may-precisely} for more in-depth surveys of May's machine.

\subsection*{Structured Ring Spectra from Categorical Data}

Spectra of interest usually have extra multiplicative structures as algebras over some operad, such as an $E_n$-operad for some $n \in \{1, 2, \ldots,\infty\}$.  A fundamental question in multiplicative infinite loop space theory is this:
\qtn{What extra structure does a small permutative category need to have in order for its $K$-theory spectrum to be a (commutative) ring spectrum or, more generally, a $\cQ$-algebra for some operad $\cQ$?}
Depending on the intended applications, there are essentially two approaches to the question above: a monadic approach due to May and a multicategorical approach due to Elmendorf-Mandell.  We briefly review these approaches below.


Building on the additive infinite loop space machine in \cite{may}, May's multiplicative machine in \cite{may-multiplicative,may-construction,may-precisely} sends bipermutative categories \cite{may-einfinity} to $\Lr$-spectra, which are also called $\Einf$-ring spectra, for the linear isometries $\Einf$-operad $\Lr$.  We note that $\Einf$-ring spectra are not defined as algebras over some $\Einf$-operad.  The first main step of May's multiplicative machine sends a bipermutative category to a $(\Stein,\Lr)$-space, which is also called an $\Einf$-ring space, for the canonical $\Einf$-operad pair $(\Stein,\Lr)$.  This step is a composite of several steps, two of which involve two-sided monadic bar constructions.  To encode the relationship between the additive part and the multiplicative part, this step uses Beck's theory of monadic distributivity \cite{beck-distributive} in a crucial way.

We mention a few notable properties of the first main step of May's multiplicative machine.  First, the input data, bipermutative categories, can be obtained as strictifications of Laplaza's symmetric bimonoidal categories with invertible distributivity morphisms; see \cite[Ch.\! 5]{cerberusI} for a detailed proof.  Thus, examples of bipermutative categories abound.  Second, for the categorical Barratt-Eccles operad $\BE$, $\BE$-algebras in $\Cat$ are small permutative categories.  However, bipermutative categories are \emph{not} the $\BE$-algebras in the $\Cat$-enriched multicategory $\PermCat$.  Instead, $\BE$-algebras in $\PermCat$ are bipermutative categories in the redefined sense of Elmendorf-Mandell \cite{elmendorf-mandell}; see \cite[Section 11.5]{cerberusIII} for a detailed proof.  Third, May's multiplicative machine starts by sending a bipermutative category to a pseudofunctor $\Fsk \!\txint\! \Fsk \to \Cat$.  Then it uses Street's strictification \cite{street_lax} to strictify that pseudofunctor to a functor in the strict sense.

The second main step of May's multiplicative machine sends a $(\Stein,\Lr)$-space $X$ to a connective $\Lr$-spectrum $\mathbbm{E}X$.  The spectrum $\mathbbm{E}X$ is the same one from the additive infinite loop space machine discussed in the previous subsection.  Thus, it is once again given by a two-sided monadic bar construction.

The multicategorical approach to the question above begins with the observation that the category $\PermCat$ extends to a categorically enriched \index{multicategory!permutative category}\index{permutative category!multicategory}multicategory, using multilinear functors and multilinear transformations for small permutative categories.  The category $\Sp$ of symmetric spectra is symmetric monoidal under the smash product, so $\Sp$ is also a multicategory.  Since multifunctors are closed under composition, the question above has a natural positive answer \emph{if} the $K$-theory machine under consideration is a multifunctor that respects the categorical enrichment of $\PermCat$ and the simplicial set enrichment of $\Sp$.  As far as we know, May's machine \cite{may} is not a multifunctor from $\PermCat$ to $\Sp$, so it does not generally send a $\cQ$-algebra in $\PermCat$ to a $\cQ$-algebra in $\Sp$ for a categorical operad $\cQ$.  

Segal $K$-theory, $\Kse$, is also not a multifunctor.  As we mention above, $\Ner$ and $\Klse$ are symmetric monoidal functors, so they are also multifunctors.  Thus, the non-multifunctoriality of $\Kse$ is entirely due to the fact that its first functor, $\Jse$, is not compatible with the multicategory structures on $\PermCat$ and $\Fskelcat$.  The root cause of the non-multifunctoriality of $\Jse$ is that the indexing category $\Fsk$ is not refined enough.  More specifically, the smash product of two pointed finite sets, $S_1 \sma S_2$, may contain basepoint-free subsets that are not of the form $T_1 \times T_2$ for some subsets $T_1 \subseteq S_1$ and $T_2 \subseteq S_2$.  This issue of $\Jse$ is explained in \cite[page 182]{elmendorf-mandell} and \cite[10.8.6]{cerberusIII}.

\subsection*{Elmendorf-Mandell $K$-Theory}

In view of the non-multifunctoriality of the $K$-theory machines of Segal and May, a natural question is this:
\qtn{Is there an enriched multifunctorial $K$-theory machine that is objectwise equivalent to Segal $K$-theory?}
The work of Elmendorf-Mandell \cite{elmendorf-mandell} constructs such an enriched multifunctor, $\Kem$, as a composite of three multifunctors\index{Elmendorf-Mandell K-theory@Elmendorf-Mandell $K$-theory}\index{K-theory@$K$-theory!Elmendorf-Mandell} as follows.
\begin{equation}\label{Kem_composite}
\begin{tikzpicture}[vcenter]
\def\h{2.5} \def\u{.7}
\draw[0cell]
(0,0) node (a) {\PermCat}
(a)++(\h,0) node (b) {\Gstarcat}
(b)++(\h,0) node (c) {\phantom{\Gstarsset}}
(c)++(0,.045) node (c') {\Gstarsset}
(c)++(\h,0) node (d) {\Sp}
;
\draw[1cell=.9]
(a) edge node {\Jem} (b)
(b) edge node {\Ner} (c)
(c) edge node {\Klem} (d)
;
\draw[1cell=1]
(a) [rounded corners=3pt] |- ($(b)+(0,\u)$)
-- node {\Kem} ($(c)+(0,\u)$) -| (d)
;
\end{tikzpicture}
\end{equation}
Similar to Segal $K$-theory, the levelwise nerve functor $\Ner$ and $\Klem$ are symmetric monoidal functors, so they are also multifunctors.  Unlike Segal $J$-theory, Elmendorf-Mandell $J$-theory $\Jem$ is a categorically enriched multifunctor.  For each categorically enriched multicategory $\cQ$, by changing enrichment and composing multifunctors, $\Kem$ sends each $\cQ$-algebra in $\PermCat$ to a $\cQ$-algebra in $\Sp$.  For example, $\Kem$ sends each $E_n$-algebra of permutative categories to an $E_n$-algebra of symmetric spectra for any $n \in \{1,2,\ldots,\infty\}$.  Moreover, each of the multifunctors $\Jem$, $\Ner \circ\, \Jem$, and $\Kem$ is an equivalence of homotopy theories \cite{johnson-yau-multiK}.

The crucial difference between Segal $K$-theory, $\Kse$, and Elmendorf-Mandell $K$-theory, $\Kem$, is that the latter uses a more refined indexing category $\Gsk$ in place of $\Fsk$.  An object of $\Gsk$ is a finite tuple of pointed finite sets.  In place of the smash product of pointed finite sets, $\Gsk$ is equipped with a permutative structure $\oplus$ given by concatenation of finite tuples.  The indexing category $\Gsk$ avoids the problem mentioned above that $\Fsk$ is not refined enough.  It is the main reason why the first step, $\Jem$, is a categorically enriched multifunctor.  The multifunctor $\Jem$ also involves gluing morphisms similar to those for Segal $J$-theory \cref{segal_gluing}, but the indexing parameters involve finite tuples of finite sets.  A detailed discussion of Elmendorf-Mandell $K$-theory can be found in \cite[Part 2]{cerberusIII}.  

Each of the $K$-theory machines of May and Elmendorf-Mandell has its own advantages.  An advantage of May's multiplicative machine is that it produces $\Einf$-ring spectra from bipermutative categories via $\Einf$-ring spaces, which are used in some applications \cite{coher_lada_may,may-einfinity,may-good}.  By definition, $\Einf$-ring spaces are $(\Stein,\Lr)$-spaces, and they are not defined as algebras over some operad.  Thus, there are no direct way to obtain $(\Stein,\Lr)$-spaces using the Elmendorf-Mandell multifunctors.  On the other hand, Elmendorf-Mandell multifunctorial $K$-theory can be used with any categorically enriched multicategory $\cQ$.  For example, composition with the multifunctor $\Kem$ produces spectral $E_n$-algebras for any $n \in \{1,2,\ldots,\infty\}$, diagrams of $E_n$-algebras, and modules over $E_n$-algebras, from corresponding types of algebras in $\PermCat$.

\subsection*{Equivariant Spectra from Categorical Data}

In the $G$-equivariant context for a compact Lie group $G$, there is a symmetric monoidal category $\GSp$ of orthogonal $G$-spectra \cite{mandell_may}.  Just as it is the case nonequivariantly, orthogonal $G$-spectra of interest usually have extra equivariant algebraic structures parametrized by some operads in $G$-spaces or $G$-categories.  Thus, we ask a similar motivating question as above:
\qtn{Is there a $G$-equivariant algebraic $K$-theory functor from genuine symmetric monoidal $G$-categories to orthogonal $G$-spectra that preserves all equivariant algebraic structures, including equivariant $\Einf$-algebras in the sense of Guillou-May and $\Ninf$-algebras in the sense of Blumberg-Hill?}
There are several existing equivariant $K$-theory machines that construct $G$-spectra from naive or genuine permutative $G$-categories, from space-level data, or from presheaves of nonequivariant data.  We briefly review some of them below.  The question of preserving equivariant algebraic structures is much harder.  In fact, \emph{none} of the existing equivariant $K$-theory machines is an enriched multifunctor starting from genuine symmetric monoidal $G$-categories.  Thus, none of them preserves equivariant algebraic structures in general.

Partially based on unpublished work of Segal, for a finite group $G$, Shimakawa's equivariant $K$-theory\index{Shimakawa K-theory@Shimakawa $K$-theory}\index{K-theory@$K$-theory!Shimakawa}\index{infinite loop space machine!Shimakawa} machine $\Ksh$ \cite{shimakawa89} constructs Lewis-May $G$-spectra from naive permutative $G$-categories.  The functor $\Ksh$ is based on Segal $K$-theory and uses the indexing $G$-category $\FG$ of pointed finite $G$-sets and pointed functions.  The issue mentioned above of $\Fskel$ not being refined enough is still present in the indexing $G$-category $\FG$.  Thus, just like Segal's functors $\Jse$ and $\Kse$ in the nonequivariant setting, the first constituent functor of Shimakawa $K$-theory and $\Ksh$ itself are not compatible with multiplicative structures.  May's operadic infinite loop space machine \cite{may} is extended to the $G$-equivariant context in the work \cite{gm17} of \index{infinite loop space machine!Guillou-May}Guillou and May.  Like its nonequivariant counterpart, the Guillou-May equivariant machine is not known to be multiplicative. 

Starting from space-level data, there is an equivalence between the category of $\FG$-$G$-spaces and the category of $\Fsk$-$G$-spaces \cite{shimakawa91}.  Using Shimakawa's results \cite{shimakawa89,shimakawa91} and extending the work \cite{bousfield_friedlander} of Bousfield-Friedlander, the work \cite{ostermayr} of Ostermayr proves that there is an equivalence between the category of very special $\Fsk$-$G$-spaces and the category of connective $G$-symmetric spectra.  Extending earlier works \cite{may-thomason,shimakawa89}, the work \cite{mmo} of May, Merling, and Osorno\index{infinite loop space machine!May-Merling-Osorno} constructs orthogonal $G$-spectra from $\FG$-$G$-spaces and proves that it is equivalent to the Guillou-May equivariant machine in an appropriate sense.  The work \cite{gmmo19} of Guillou, May, Merling, and Osorno\index{infinite loop space machine!Guillou-May-Merling-Osorno} constructs a symmetric monoidal functor from the symmetric monoidal category of $\Fsk$-$G$-spaces to orthogonal $G$-spectra.  The work \cite{kong_may_zou} of Kong, May, and Zou extends and conceptualizes earlier works \cite{gm17,may,may-multiplicative,may-precisely,mmo}, producing $\Einf$-ring $G$-spectra from $\Einf$-ring $G$-spaces in specific examples.  

There are also equivariant $K$-theory machines based on presheaves of nonequivariant spaces or spectra.  The work \cite{costenoble_waner} of Costenoble and Waner constructs Lewis-May $G$-spectra from topological presheaves on the orbit category of $G$ equipped with an equivariant $\Einf$ structure.  Using presheaves of nonequivariant orthogonal spectra as input, the main result of Guillou and May in \cite{gm-presheaves} proves that there is a Quillen equivalence between the category of orthogonal $G$-spectra and the category of spectral Mackey functors.

The equivariant $K$-theory\index{K-theory@$K$-theory!Guillou-May-Merling-Osorno} machine of Guillou, May, Merling, and Osorno \cite[Theorem A]{gmmo23}, denoted $\Kgm$ there, constructs orthogonal $G$-spectra from $\Op$-algebras, where $\Op$ is any $G$-categorical $\Einf$-operad that is levelwise a translation category.  For example, for the $G$-Barratt-Eccles operad $\GBE$, $\GBE$-algebras are genuine permutative $G$-categories.  The equivariant machine $\Kgm$ also applies to $\Op$-pseudoalgebras when it is precomposed with the strictification functor in \cite{gmmo20} that first strictifies $\Op$-pseudoalgebras to $\Op$-algebras.  The construction $\Kgm$ is a \emph{nonsymmetric} multifunctor.  This means that $\Kgm$ preserves the multicategorical compositions and units of $\Op$-algebras and orthogonal $G$-spectra, but \emph{not} the symmetric group actions.  The non-symmetry of $\Kgm$ is caused by one of the intermediate categories $\mathbf{Mult}(\FG\mh\mathbf{PsAlg})$ and the strictification functor $\mathrm{St}$ in its construction, which are not compatible with symmetric group actions.  In particular, $\Kgm$ does not generally preserve equivariant algebraic structures, such as equivariant $\Einf$-algebras in the sense of Guillou-May \cite{gm17} and $\Ninf$-algebras in the sense of Blumberg-Hill \cite{blumberg_hill}.  

In summary, all existing $G$-equivariant algebraic $K$-theory functors in the literature, from naive or genuine permutative $G$-categories to $G$-spectra, are \emph{not} enriched multifunctors that preserve equivariant algebraic structures.  This finishes our brief review of relevant literature.

\sect{Main Results}
\label{sec:main_results}

With the above motivation, background, and context in mind, we describe the main results of this work.

\subsection*{Multifunctorial Equivariant Algebraic $K$-Theory}

A $\Tinf$-operad \pcref{as:OpA} is a 1-connected pseudo-commutative operad enriched in $G$-categories.  For each compact Lie group $G$ and $\Tinf$-operad $\Op$, the main result of this work \pcref{thm:Jgo_multifunctor,thm:ggcat_ggtop,thm:Kg_smgtop,thm:Kgo_multi} constructs a composite enriched multifunctor $\Kgo$ as follows.
\begin{equation}\label{Kgo_composite}
\begin{tikzpicture}[vcenter]
\def\h{2.5} \def\u{.7}
\draw[0cell]
(0,0) node (a) {\phantom{\MultpsO}}
(a)++(0,-.04) node (a') {\MultpsO}
(a)++(\h,0) node (b) {\GGCat}
(b)++(\h,0) node (c) {\GGTop}
(c)++(\h,0) node (d) {\GSp}
;
\draw[1cell=.9]
(a) edge node {\Jgo} (b)
(b) edge node {\clast} (c)
(c) edge node {\Kg} (d)
;
\draw[1cell=1]
(a') [rounded corners=3pt] |- ($(b)+(0,\u)$)
-- node {\Kgo} ($(c)+(0,\u)$) -| (d)
;
\end{tikzpicture}
\end{equation}
The first part, $\Jgo$, is a $\Gcat$-multifunctor, which means a $G$-categorically enriched multifunctor.  The other two parts, $\clast$ and $\Kg$, are symmetric monoidal $\Gtop$-functors, which mean $G$-topologically enriched symmetric monoidal functors.  We emphasize that each of the four multifunctors---$\Kgo$, $\Jgo$, $\clast$, and $\Kg$---is enriched in either $G$-categories or $G$-spaces, not just categories or spaces.  The enrichment in $\Gcat$ and $\Gtop$ is crucial in applications because equivariant $\Einf$-algebras and $\Ninf$-algebras are parametrized by operads enriched in $G$-categories or $G$-spaces.  Below we discuss $\Jgo$, $\clast$, $\Kg$, applications, and relation to earlier work.

\parhead{The first step: $\Jgo$}.
The domain $\MultpsO$ is the $G$-categorically enriched multicategory with $\Op$-pseudoalgebras as objects, $k$-lax $\Op$-morphisms as $k$-ary 1-cells, and $k$-ary $\Op$-transformations as $k$-ary 2-cells.  In the notation $\MultpsO$, the superscript $\mathsf{ps}$ refers to $\Op$-pseudoalgebras, and the subscript $\mathsf{lax}$ refers to $k$-lax $\Op$-morphisms.  For example, for the $G$-Barratt-Eccles operad $\GBE$, $\GBE$-pseudoalgebras are genuine symmetric monoidal $G$-categories.  In this case, our equivariant algebraic $K$-theory multifunctor
\[\MultpsGBE \fto{\Kggbe} \GSp\]
sends genuine symmetric monoidal $G$-categories to orthogonal $G$-spectra, in a way that respects the enriched multicategorical structures of its domain and codomain.  

The first step of $\Kgo$ is the $\Gcat$-multifunctor \pcref{thm:Jgo_multifunctor}
\[\MultpsO \fto{\Jgo} \GGCat\]
from the $\Gcat$-multicategory $\MultpsO$ of $\Op$-pseudoalgebras to the symmetric monoidal closed diagram category $\GGCat$ of $\Gskg$-categories.  The $\Gcat$-multifunctor $\Jgo$ becomes a $\Gtop$-multifunctor by changing enrichment along the classifying space functor $\cla \cn \Gcat \to \Gtop$.  The existence of the $\Gcat$-multicategory $\MultpsO$ is proved in detail in \cref{ch:multpso}.  The diagram category $\GGCat$ is discussed in \cref{ch:ggcat}.  Its indexing category $\Gsk$ is the same one used in Elmendorf-Mandell $K$-theory \cref{Kem_composite}, and $\Gcat$ is the symmetric monoidal closed category of small $G$-categories.

The $\Gcat$-multifunctor $\Jgo$, constructed in \cref{ch:jemg}, is the most nontrivial part of $\Kgo$.  To get a flavor of $\Jgo$, recall that Segal $J$-theory $\Jse$ \cref{segal_gluing} and Elmendorf-Mandell $J$-theory $\Jem$ \cref{Kem_composite} both involve gluing morphisms that glue two objects together at a time.  A key component of the construction of $\Jgo$ is the notion of a gluing morphism that glues $r$ objects together at a time \cref{gluing-morphism}:
\begin{equation}\label{Jgo_gluing}
\gaA_r\big(x\sscs \ang{a_{\ang{s} \compk\, s_{k,i}}}_{i \in \ufs{r}} \big) 
\fto{\glu_{x; \ang{s} \csp k, \ang{s_{k,i}}_{i \in \ufs{r}}}} a_{\ang{s}}.
\end{equation}
Its details are explained in \cref{def:nsystem}.  One of the parameters of this gluing morphism is an object $x \in \Op(r)$ for the $\Tinf$-operad $\Op$.  Thus, instead of gluing two objects at a time, for multifunctorial $J$-theory $\Jgo$, we have a gluing morphism that glues $r$ objects $\ang{a_{\ang{s} \compk\, s_{k,i}}}_{i \in \ufs{r}}$ together for each object $x \in \Op(r)$.  Essentially all the definitions, axioms, constructions, and proofs in \cref{ch:jemg} are designed to make this idea of gluing morphism work.  

A nontrivial feature of the $\Gcat$-multifunctor $\Jgo$ is the fact that, at the object level, it sends $\Op$-pseudoalgebras to $\Gskg$-categories, which are pointed functors $\Gsk \to \Gcatst$ in the usual 1-categorical sense.  We emphasize that $\Jgo$ does \emph{not} involve any strictification functor in its construction.  In particular, $\Jgo$ does not go through some categories of $\Op$-algebras or pseudofunctors $\Gsk \to \Gcatst$.  The associativity constraint $\phiA$ of an $\Op$-pseudoalgebra $\A$ is incorporated into the construction of the $\Gskg$-category $\Jgo\A$, in the associativity axiom \cref{system_associativity} and the commutativity axiom \cref{system_commutativity} of an $\angordn$-system in $\A$.  See \cref{rk:no_strictification} for related discussion.

\parhead{The second step: $\clast$}.
The second step of $\Kgo$, constructed in \cref{ch:ggtop}, is the symmetric monoidal functor 
\[\GGCat \fto{\clast} \GGTop\]
from the symmetric monoidal closed category $\GGCat$ of $\Gskg$-categories to the symmetric monoidal closed category $\GGTop$ of $\Gskg$-spaces.  It is induced by the classifying space functor $\cla \cn \Gcat \to \Gtop$.  In the codomain, $\Gtop$ is the symmetric monoidal closed category of $G$-spaces.   Via the self-enrichment of $\GGCat$ and $\GGTop$, $\clast$ becomes a symmetric monoidal functor in the $\Gtop$-enriched sense.  The symmetric monoidal $\Gtop$-functor $\clast$ is a somewhat formal categorical construction.  Its details are interesting, but nowhere near as intense as $\Jgo$.

\parhead{The last step: $\Kg$}.
The last step of $\Kgo$, constructed in \cref{ch:semg}, is the unital symmetric monoidal $\Gtop$-functor 
\[\GGTop \fto{\Kg} \GSp\]
from the symmetric monoidal $\Gtop$-category of $\Gskg$-spaces to the symmetric monoidal $\Gtop$-category $\GSp$ of orthogonal $G$-spectra.  The unitality of $\Kg$ means that its unit constraint is a $G$-equivariant $\gsp$-module isomorphism
\[\gsp \fto[\iso]{\Kgzero} \Kg\gu.\]
Thus, $\Kg$ has the $G$-sphere $\gsp$ in its image.  The proof that $\Kg$ is a unital symmetric monoidal $\Gtop$-functor relies on a detailed, point-set level understanding of its domain and codomain.  To facilitate this proof, \cref{ch:spectra} provides a self-contained review of orthogonal $G$-spectra, supplying lots of details that are not explicitly available in the literature.

\subsection*{Applications: $\Kgo$ Preserves All Equivariant Algebraic Structures}

Enriched multifunctoriality is the main new feature of our $G$-equivariant algebraic $K$-theory $\Kgo$ that distinguishes it from all other equivariant $K$-theory machines.  By simply composing enriched multifunctors, $\Kgo$ transports all equivariant algebraic structures parametrized by multicategories enriched in $\Gcat$ or $\Gtop$, from the $\Gcat$-multicategory $\MultpsO$ of $\Op$-pseudoalgebras to corresponding types of equivariant algebraic structures in orthogonal $G$-spectra.  In more detail, we consider any $\Gcat$-multicategory $\cQ$ and $\Gcat$-multifunctor $f \cn \cQ \to \MultpsO$, which parametrizes a $\cQ$-algebra of $\Op$-pseudoalgebras.  Then the composite 
\[\cQ \fto{f} \MultpsO \fto{\Kgo} \GSp\]
is a $\Gtop$-multifunctor \pcref{thm:Kgo_preservation}.  Thus, $\Kgo$ transports the purely categorical data of a $\cQ$-algebra $f$ in the $\Gcat$-multicategory $\MultpsO$ to the $\cQ$-algebra $\Kgo \circ f$ in the symmetric monoidal $\Gtop$-category $\GSp$ of orthogonal $G$-spectra.  

For example, we can take $\cQ$ to be the $G$-Barratt-Eccles operad $\GBE$ for a finite group $G$ \pcref{def:GBE}.  A $\Gcat$-multifunctor 
\[\GBE \fto{f} \MultpsO\]
parametrizes an equivariant $\Einf$-algebra in the sense of Guillou-May \cite{gm17}, and also an $\Ninf$-algebra in the sense of Blumberg-Hill \cite{blumberg_hill}, of $\Op$-pseudoalgebras.  The composite $\Gtop$-multifunctor 
\[\GBE \fto{f} \MultpsO \fto{\Kgo} \GSp\]
is an equivariant $\Einf$-algebra, and also an $\Ninf$-algebra, of orthogonal $G$-spectra.  We can further specialize the $\Tinf$-operad $\Op$ to $\GBE$.  A $\Gcat$-multifunctor 
\[\GBE \fto{f} \MultpsGBE\]
parametrizes an equivariant $\Einf$-algebra and an $\Ninf$-algebra of genuine symmetric monoidal $G$-categories.  The composite $\Gtop$-multifunctor 
\[\GBE \fto{f} \MultpsGBE \fto{\Kggbe} \GSp\]
is an equivariant $\Einf$-algebra and an $\Ninf$-algebra of orthogonal $G$-spectra.  See \cref{ex:Einf_GBE,ex:Jgo_preservation,ex:Kgo_preservation} for further discussion.

\subsection*{Relation to Earlier Work}

The construction of our $G$-equivariant algebraic $K$-theory multifunctor, $\Kgo$ \cref{Kgo_composite}, is explicit and entirely self-contained.  It does not depend on any other (non)equivariant $K$-theory constructions.  Conceptually, $\Kgo$ is a $G$-equivariant extension of Elmendorf-Mandell nonequivariant multifunctorial $K$-theory, $\Kem$ \cref{Kem_composite}.  Specifically, $\Kgo$ extends, and differs from, $\Kem$ in the following ways.

First, we work $G$-equivariantly, so the codomain of our $\Kgo$ is the symmetric monoidal category $\GSp$ of orthogonal $G$-spectra, instead of the symmetric monoidal category $\Sp$ of symmetric spectra.  In the two intermediate symmetric monoidal closed diagram categories, $\GGCat$ and $\GGTop$, we use the symmetric monoidal closed categories $\Gcat$ and $\Gtop$, instead of $\Cat$ and $\sSet$.  The last step of $\Kgo$ is the unital symmetric monoidal $\Gtop$-functor \pcref{ch:semg}
\[\GGTop \fto{\Kg} \GSp.\]  
Its construction is conceptually similar to the last step of $\Kem$, which is the symmetric monoidal functor 
\[\Gstarsset \fto{\Klem} \Sp.\]  
However, the detailed verification of the properties of $\Kg$ is more involved than its nonequivariant counterpart, $\Klem$.

There is a significant conceptual jump from the domain of $\Kem$ to the domain of $\Kgo$.  The domain of $\Kem$ is the $\Cat$-multicategory $\PermCat$.  Its objects are small permutative categories, which are, equivalently, $\BE$-algebras in $\Cat$ for the Barratt-Eccles operad $\BE$.   On the other hand, as we discuss in \cref{ch:multpso}, the domain of our $\Kgo$ is the $\Gcat$-multicategory $\MultpsO$, whose objects are $\Op$-pseudoalgebras for a $\Tinf$-operad $\Op$.  By definition, a $\Tinf$-operad is a 1-connected pseudo-commutative operad in $\Gcat$.  Thus, our framework allows more general $\Gcat$-operads than just $\BE$ and the $G$-Barratt-Eccles operad $\GBE$.  More importantly, instead of just $\Op$-algebras in the strict sense, our $\Kgo$ takes $\Op$-pseudoalgebras as input data.

Multifunctorial $J$-theory \pcref{ch:jemg}
\[\MultpsO \fto{\Jgo} \GGCat\] 
is the most crucial part of $\Kgo$, and it diverges substantially from Elmendorf-Mandell $J$-theory
\[\PermCat \fto{\Jem} \Gstarcat.\]  
Since the input for $\Jem$ is a small permutative category, by coherence properties of the categorical Barratt-Eccles operad $\BE$, it is sufficient to axiomatize gluing of two objects.  On the other hand, a general $\Tinf$-operad $\Op$, including the $G$-Barratt-Eccles operad $\GBE$ for any nontrivial group $G$, does not have the coherence properties of $\BE$.  This is the reason why, for multifunctorial $J$-theory $\Jgo$, we need to consider general gluing morphisms, as stated in \cref{Jgo_gluing}, that glue $r \geq 0$ objects together for each object $x \in \Op(r)$.  Furthermore, as we mention above, $\Jgo$ takes an $\Op$-pseudoalgebra as input and produces a $\Gskg$-category as output, \emph{without} using any strictification functor.

Another subtle difference between $\Kem$ and $\Kgo$ happens in their enrichment.  Each of the three constituent multifunctors of $\Kem$---$\Jem$, $\Ner$, and $\Klem$---is enriched in either small categories or simplicial sets.  On the other hand, enrichment in the $G$-equivariant setting requires more care.  The self-enrichment of $\Gcat$ involves not-necessarily equivariant functors, on which $G$ acts by conjugation, and similarly for the self-enrichment of $\Gtop$.  Our multifunctorial equivariant algebraic $K$-theory, 
\[\Kgo = \Kg \circ \clast \circ \Jgo,\]
respects the $G$-categorical or $G$-topological enrichment in each step.  To achieve that, we need to take extra care when we define objects, $k$-ary 1-cells, $k$-ary 2-cells, and the $G$-action in each step.

\sect{Chapter Summaries}\label{sec:chsummary}

This work consists of the following two parts and three appendices.
\begin{itemize}
\item \prefacepartNumName{part:jgo}
\item \prefacepartNumName{part:kg}
\end{itemize}
\cref{part:jgo} consists of \cref{ch:psalg,ch:multpso,ch:ggcat,ch:jemg}.  \cref{ch:psalg,ch:multpso,ch:ggcat} discuss pseudo-commutative operads, pseudoalgebras of an operad, the $\Gcat$-multicategory $\MultpsO$ of $\Op$-pseudoalgebras, and the symmetric monoidal closed category $\GGCat$ of $\Gskg$-categories.  \cref{ch:jemg} constructs the $\Gcat$-multifunctor $\Jgo$ from $\MultpsO$ to $\GGCat$, which is the first multifunctor in $\Kgo$.

\cref{part:kg} consists of \cref{ch:spectra,ch:ggtop,ch:semg}.  \cref{ch:spectra} reviews orthogonal $G$-spectra.  \cref{ch:ggtop} discusses the symmetric monoidal closed category $\GGTop$ of $\Gskg$-spaces and constructs the symmetric monoidal $\Gtop$-functor $\clast$ from $\GGCat$ to $\GGTop$, which is the second part of $\Kgo$.  \cref{ch:semg} constructs the unital symmetric monoidal $\Gtop$-functor $\Kg$ from $\GGTop$ to $\GSp$, which is the last part of $\Kgo$, and discusses applications of $\Kgo$ to equivariant algebraic structures in orthogonal $G$-spectra.  \cref{ch:prelim,ch:prelim_multicat,ch:nerve} review background material on (multi)categories and the classifying space functor.

A brief summary of each chapter is given below.  In the main text, each chapter has a detailed introduction that provides an overview of that chapter and describes how it connects with other chapters.  Each section also has a brief introduction and an outline.

\prefacechapNumName{ch:psalg}

This chapter discusses pseudo-commutative operads in $\Gcat$ and pseudoalgebras of a $\Gcat$-operad.  A reduced $\Gcat$-operad that is levelwise a translation category is a pseudo-commutative operad in a unique way.  For example, the Barratt-Eccles operad $\BE$, with 
\[\BE(n) = \ESigma_n,\] 
and the $G$-Barratt-Eccles operad $\GBE$, with 
\[\GBE(n) = \Catg(\EG,\ESigma_n),\] 
both admit unique pseudo-commutative structures.  Pseudo-commutativity is used in \cref{ch:multpso} to define $k$-lax $\Op$-morphisms; see the commutativity axiom \cref{laxf_com}.  Pseudoalgebras generalize algebras of an operad in the same way that monoidal categories generalize strict monoidal categories.  For an $\Op$-pseudoalgebra $\A$, the usual associativity axiom of an $\Op$-algebra is replaced by an associativity constraint $\phiA$ that satisfies several coherence axioms.  We prove in detail that pseudoalgebras of $\BE$ are, up to 2-equivalences, \emph{naive} symmetric monoidal $G$-categories \pcref{thm:BEpseudoalg}.  Pseudoalgebras of $\GBE$ are, by definition, \emph{genuine} symmetric monoidal $G$-categories.  

\prefacechapNumName{ch:multpso}

This chapter constructs the $\Gcat$-multicategory $\MultpsO$ for each pseudo-commutative operad $\Op$ in $\Gcat$ \pcref{thm:multpso}.  Its objects are $\Op$-pseudoalgebras, as discussed in \cref{ch:psalg}.  For example, for the $G$-Barratt-Eccles operad $\GBE$, $\MultpsGBE$ is a $\Gcat$-multicategory whose objects are genuine symmetric monoidal $G$-categories.  The $\Gcat$-enrichment of $\MultpsO$ uses the notions of $k$-lax $\Op$-morphisms and $k$-ary $\Op$-transformations.  A $k$-lax $\Op$-morphism generalizes an $\Op$-algebra morphism by (i) allowing the input to be a $k$-tuple of $\Op$-pseudoalgebras and (ii) relaxing the compatibility axiom with the $\Op$-action to action constraints that satisfy several coherence axioms.  The commutativity axiom \cref{laxf_com} uses the pseudo-commutative structure of $\Op$.  For a $k$-lax $\Op$-morphism, the underlying functor and the action constrains are \emph{not} required to be $G$-equivariant, so the group $G$ can act nontrivially by conjugation.  This is similar in spirit to the self-enrichment of $G$-spaces, where $G$ acts by conjugation on not-necessarily equivariant morphisms. 

\prefacechapNumName{ch:ggcat}

This chapter discusses the symmetric monoidal closed category $\GGCat$.  Its objects, called $\Gskg$-categories, are pointed functors $\Gsk \to \Gcatst$, where $\Gcatst$ is the symmetric monoidal closed category of small pointed $G$-categories.  The first half of this chapter discusses the Elmendorf-Mandell indexing category $\Gsk$.  The second half of this chapter completely unravels the symmetric monoidal closed structure on $\GGCat$ and its associated $\Gcat$-multicategory.  The latter is used in \cref{ch:jemg} in the construction of $\Jgo$.

\prefacechapNumName{ch:jemg}

This chapter constructs the first part of our $G$-equivariant algebraic $K$-theory multifunctor $\Kgo$:  the $\Gcat$-multifunctor \pcref{thm:Jgo_multifunctor}
\[\MultpsO \fto{\Jgo} \GGCat\]
from $\Op$-pseudoalgebras \pcref{ch:multpso} to $\Gskg$-categories \pcref{ch:ggcat} for a $\Tinf$-operad $\Op$.  By a $\Tinf$-operad, we mean a 1-connected pseudo-commutative operad in $\Gcat$ \pcref{as:OpA}.  For example, both $\BE$ and $\GBE$ are $\Tinf$-operads.  There is also a strong variant $\Gcat$-multifunctor
\[\MultpspsO \fto{\Jgosg} \GGCat\]
where the $k$-ary 1-cells in the domain, $\MultpspsO$, are $k$-ary $\Op$-pseudomorphisms.  This chapter is divided into three parts.  \cref{sec:jemg_objects,sec:jemg_morphisms,sec:jemg_morphisms_ii} construct the object assignment of $\Jgo$, which sends each $\Op$-pseudoalgebra $\A$ to a $\Gskg$-category $\Adash$. \cref{sec:jemg_zero,sec:jemg_pos_i,sec:jemg_pos_i_proof,sec:jemg_pos_ii} construct the assignment of $\Jgo$ on multimorphism $G$-categories.  \cref{sec:jemg_axioms,sec:Jgo_preserves} prove the $\Gcat$-multifunctor axioms for $\Jgo$ and discuss its preservation of equivariant algebraic structures.  \cref{sec:Jgo_preserves} includes self-contained reviews of equivariant $\Einf$-operads and $\Ninf$-operads, along with examples.  This finishes \cref{part:jgo}.

\prefacechapNumName{ch:spectra}

Assuming no prior knowledge of spectra, this chapter is an in-depth review of orthogonal $G$-spectra, with an emphasis on its symmetric monoidal $\Gtop$-category structure.  While nothing in this chapter is technically new, a lot of its details, even in the nonequivariant case, are not explicitly spelled out anywhere in the literature that we are aware of.  The proofs in \cref{ch:semg} require a nontrivial amount of point-set level details of orthogonal $G$-spectra.  For completeness and clarity, this chapter supplies these details.  This chapter is suitable for anyone who is new to orthogonal $G$-spectra.

\prefacechapNumName{ch:ggtop}

This chapter constructs the second part of our $G$-equivariant algebraic $K$-theory multifunctor $\Kgo$: the symmetric monoidal $\Gtop$-functor \pcref{thm:ggcat_ggtop}
\[\GGCat \fto{\clast} \GGTop\]
between symmetric monoidal closed categories, induced by the classifying space functor $\cla$.  This chapter includes a detailed description of the symmetric monoidal $\Gtop$-category $\GGTop$ of $\Gskg$-spaces.  This detailed description of $\GGTop$ is used in \cref{ch:semg} to construct $\Kg$.

\prefacechapNumName{ch:semg}

The last chapter constructs the unital symmetric monoidal $\Gtop$-functor \pcref{thm:Kg_smgtop}
\[\GGTop \fto{\Kg} \GSp\]
from $\Gskg$-spaces to orthogonal $G$-spectra, using details from \cref{ch:spectra,ch:ggtop}.  \cref{sec:Kgo_multi} is the culmination of this work, combining most of the main results from earlier chapters and sections to obtain the enriched multifunctor \pcref{thm:Kgo_multi}
\[\MultpsO \fto{\Kgo} \GSp\]
and its strong variant
\[\MultpspsO \fto{\Kgosg} \GSp.\]
The second half of this section discusses the fact that $\Kgo$ preserves all equivariant algebraic structures parametrized by multicategories enriched in $\Gcat$ or $\Gtop$.  In particular, for a finite group $G$, equivariant $\Einf$-algebras and $\Ninf$-algebras of $\Op$-pseudoalgebras are sent by $\Kgo$ to corresponding types of algebras in $\GSp$.  Further specializing $\Op$ to the $G$-Barratt-Eccles operad $\GBE$, $\Kggbe$ sends equivariant $\Einf$-algebras and $\Ninf$-algebras of genuine symmetric monoidal $G$-categories to corresponding types of algebras in $\GSp$.  This finishes \cref{part:kg}.

\prefacechapNumName{ch:prelim}

This appendix reviews enriched monoidal categories and 2-categories.  A more detailed introduction to 2-categories is \cite{johnson-yau}.  Enriched monoidal categories are discussed in detail in \cite[Part 1]{cerberusIII}.

\prefacechapNumName{ch:prelim_multicat}

This appendix reviews enriched multicategories, change of enrichment, and pointed diagram categories.  More detailed discussion of this material can be found in \cite{yau-operad,cerberusIII}.

\prefacechapNumName{ch:nerve}

This brief appendix reviews simplicial sets, the nerve functor, geometric realization, and the classifying space functor.

\sect{Reading Suggestions}\label{sec:read_sug}

The following reading suggestions may help the reader get to certain threads of discussion quicker.  The main text has extensive cross-references to make it easier to navigate back and forth.

\subsection*{Overview}

For a more detailed overview of this work, read all the chapter introductions and \cref{sec:Kgo_multi}.  Each section has a short introduction and an outline.  We suggest reading the section outline before diving into each section.  The following diagram shows the logical interdependence of the seven chapters.
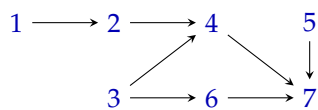
\[\begin{tikzpicture}
\def\h{1.3} \def\v{-1}
\draw[0cell]
(0,0) node (a1) {\tref{ch:psalg}}
(a1)++(\h,0) node (a2) {\tref{ch:multpso}}
(a2)++(0,\v) node (a3) {\tref{ch:ggcat}}
(a2)++(\h,0) node (a4) {\tref{ch:jemg}}
(a4)++(\h,0) node (a5) {\tref{ch:spectra}}
(a4)++(0,\v) node (a6) {\tref{ch:ggtop}}
(a5)++(0,\v) node (a7) {\tref{ch:semg}}
;
\draw[1cell]
(a1) edge (a2)
(a2) edge (a4)
(a3) edge (a4)
(a3) edge (a6)
(a4) edge (a7)
(a5) edge (a7)
(a6) edge (a7)
;
\end{tikzpicture}\]
For example, \cref{ch:jemg,ch:spectra,ch:ggtop} are largely independent of each other, and they are all used in \cref{ch:semg}.

\subsection*{Background Material}

For background material on enriched symmetric monoidal categories and enriched multicategories, read \cref{ch:prelim,ch:prelim_multicat}.  The classifying space functor, along with the nerve and geometric realization, are reviewed in \cref{ch:nerve}.  For orthogonal $G$-spectra, read \cref{ch:spectra}.

\subsection*{Related Literature}

For discussion about related literature, read the introduction of \cref{ch:spectra},  \cref{rk:pseudocom_generality,rk:gmmo-cg-notation,expl:pseudoalg_axioms,expl:BEpseudo_smcat,expl:GBEvsBE_pseudo,expl:shufni,expl:k_laxmorphism,ex:multBE,ex:multGBE,rk:no_strictification,expl:Einfinity_alg,ex:Einf_GBE,rk:gspace_gspectra}.

\subsection*{Symmetric Monoidal $G$-Categories}

For discussion about naive or genuine, permutative or symmetric monoidal $G$-categories, read the introduction of \cref{ch:psalg}, followed by \cref{sec:BE,sec:GBE,sec:pseudoalgebra,sec:naive_smc,sec:BEpseudoalg,sec:genuine_smgcat}.  Discussion that mentions these structured $G$-categories includes \cref{ex:multBE,ex:multGBE,ex:JgBE,ex:Einf_GBE,ex:Jgo_preservation,ex:Kgo_preservation}.

\subsection*{Equivariant $\Einf$-Operads and $\Ninf$-Operads}

For a review of these equivariant operads, read \cref{def:Einfty_operads,def:Ninfty_operads}.  For examples of these equivariant operads, read \cref{ex:Einf_GBE,ex:Einf_steiner,ex:Einf_linear}.  Discussion that uses the $G$-Barratt-Eccles operad $\GBE$ as a $G$-categorical $\Einf$-operad and $\Ninf$-operad includes \cref{ex:Jgo_preservation,ex:Kgo_preservation}.

\subsection*{Multifunctorial Equivariant Algebraic $K$-Theory}

Our $G$-equivariant algebraic $K$-theory multifunctor $\Kgo$ is a composite of three enriched multifunctors: $\Jgo$, $\clast$, and $\Kg$.  For the $\Gcat$-multifunctor $\Jgo$, read \cref{def:Jgo_multifunctor,thm:Jgo_multifunctor}.  For the symmetric monoidal $\Gtop$-functor $\clast$, read \cref{thm:ggcat_ggtop}.  For the unital symmetric monoidal $\Gtop$-functor $\Kg$, read \cref{thm:Kg_smgtop} and the definitions stated there.  For the $\Gtop$-multifunctor $\Kgo$, read \cref{thm:Kgo_multi,expl:Kgo_obj}.

\subsection*{Preservation of Equivariant Algebraic Structures}

For discussion about how multifunctorial $J$-theory $\Jgo$ sends equivariant algebraic structures in $\MultpsO$ to those in $\GGCat$, read \cref{thm:Jgo_preservation,Jgo_preserves_Einf,ex:Jgo_preservation}.  For discussion about how $\Kgo$ sends equivariant algebraic structures in $\MultpsO$ to those in $\GSp$, read \cref{thm:Kgo_preservation,ex:Kgo_preservation}.

\mainmatter

\part{From Operadic Pseudoalgebras to $\Gskg$-Categories}
\label{part:jgo}

\chapter{Pseudo-Commutative Operads and\\ Operadic Pseudoalgebras}
\label{ch:psalg}
This chapter discusses two concepts about $G$-categorical operads:
\begin{itemize}
\item pseudo-commutative operads (\cref{def:pseudocom_operad}) and 
\item pseudoalgebras of an operad (\cref{def:pseudoalgebra}). 
\end{itemize}
These notions are necessary to define the domain of our $G$-equivariant algebraic $K$-theory multifunctor $\Kgo$.  

Throughout this chapter, $G$ denotes an arbitrary group, and $\Gcat$ is the 2-category of small $G$-categories, $G$-functors, and $G$-natural transformations.   A pseudo-commutative structure on an operad is analogous to a symmetric monoidal structure on a monoidal category (\cref{expl:pseudocom_operad}).  Any reduced $\Gcat$-operad that is levelwise a translation category admits a unique pseudo-commutative structure (\cref{translation_pseudocom}).  Pseudoalgebras of an operad are equipped with an \emph{associativity constraint} that replaces the usual associativity axiom.  They generalize operadic algebras in the same way that monoidal categories generalize strict monoidal categories.

\subsection*{Symmetric Monoidal $G$-Categories}
One important example of a pseudo-commutative operad is the Barratt-Eccles $\Gcat$-operad $\BE$ (\cref{def:BE-Gcat}). Each level of $\BE$ is a translation category
\[\BE(n) = \ESigma_n\]
with $G$ acting trivially.  Here $\tn$ denotes the translation category construction (\cref{def:translation_cat}).  Algebras over $\BE$ are naive permutative $G$-categories (\cref{expl:naive_perm_Gcat}).  Moreover, $\BE$-pseudoalgebras correspond to naive symmetric monoidal $G$-categories (\cref{def:naive_smGcat}) via the 2-equivalences in \cref{thm:BEpseudoalg}.

Another motivating example of a pseudo-commutative operad is the $G$-Barratt-Eccles operad $\GBE$ (\cref{def:GBE}), with each level given by a translation category
\[\GBE(n) = \Catg(\EG,\ESigma_n) \iso \tn[G,\Sigma_n].\]
Its algebras and pseudoalgebras are, respectively, genuine permutative $G$-categories and genuine symmetric monoidal $G$-categories (\cref{def:GBE_algebra,def:GBE_pseudoalg}).  The following table summarizes the algebras and pseudoalgebras of $\BE$ and $\GBE$ in $\Gcat$.
\begin{center}
\resizebox{\columnwidth}{!}{%
{\renewcommand{\arraystretch}{1.3}%
{\setlength{\tabcolsep}{1em}
\begin{tabular}{c|cc}
& \multirow{2}{*}{\makecell{Barratt-Eccles operad \\ $\BE$ (\ref{def:BE-Gcat})}}
& \multirow{2}{*}{\makecell{$G$-Barratt-Eccles operad \\ $\GBE$ (\ref{def:GBE})}}\\
&&\\ \hline
\multirow{2}{*}{\makecell{Algebras}} 
& \multirow{2}{*}{\makecell{Naive permutative \\ $G$-categories (\ref{expl:naive_perm_Gcat})}} 
& \multirow{2}{*}{\makecell{Genuine permutative \\ $G$-categories (\ref{expl:GBE_algebra})}}\\
&&\\
\multirow{2}{*}{\makecell{Pseudoalgebras}} 
& \multirow{2}{*}{\makecell{Naive symmetric monoidal \\ $G$-categories (\ref{thm:BEpseudoalg})}} 
& \multirow{2}{*}{\makecell{Genuine symmetric monoidal \\ $G$-categories (\ref{expl:GBE_pseudoalg})}}\\
&&\\
\end{tabular}}}}
\end{center}

\connection
In \cref{ch:multpso}, we prove in detail that, for each pseudo-commutative operad $\Op$ in $\Gcat$, there is a $\Gcat$-multicategory $\MultpsO$ whose objects are $\Op$-pseudoalgebras; see \cref{thm:multpso}.  The pseudo-commutative structure of $\Op$ is used in the commutativity axiom \cref{laxf_com} of $k$-lax $\Op$-morphisms, which are the $k$-ary 1-cells of $\MultpsO$.  As examples of \cref{thm:multpso}, we obtain $\Gcat$-multicategories $\MultpsBE$ and $\MultpsGBE$, whose objects are, respectively, $\BE$-pseudoalgebras and $\GBE$-pseudoalgebras.  We emphasize that our enriched multicategories and multifunctors are always assumed to be symmetric; see \cref{def:enr-multicategory,def:enr-multicategory-functor}. 

As we discuss in \cref{ch:jemg}, our $G$-equivariant algebraic $K$-theory multifunctor
\[\MultpsO \fto{\Kgo} \GSp\] 
starts at the $\Gcat$-multicategory $\MultpsO$ of $\Op$-pseudoalgebras for a $\Tinf$-operad $\Op$ \pcref{as:OpA}, which means a 1-connected pseudo-commutative operad in $\Gcat$.  The first step of $\Kgo$ is a $\Gcat$-multifunctor (\cref{thm:Jgo_multifunctor})
\[\MultpsO \fto{\Jgo} \GGCat\] 
that sends each $\Op$-pseudoalgebra to a $\Gskg$-category, incorporating the associativity constraint in its construction.  We stress that $\Jgo$ does \emph{not} involve any strictification functors.  In particular, using $\Op = \BE$ and $\GBE$, $\Jgo$ applies directly to $\BE$-pseudoalgebras and $\GBE$-pseudoalgebras, without first strictifying them to, respectively, naive and genuine permutative $G$-categories.  See \cref{ex:JgBE}.

\organization
This chapter consists of the following sections.

\secname{sec:intrinsic_pairing}  This section discusses the intrinsic pairing of an operad, transpose permutations, and the 2-category $\Gcat$ of small $G$-categories, $G$-functors, and $G$-natural transformations.  The intrinsic pairing of the associative operad is discussed in \cref{ex:as_intrinsic}.  Each of these concepts is needed to define pseudo-commutative operads.  

\secname{sec:pseudocomoperads}  This section recalls the definition of a pseudo-commutative operad in $\Gcat$.  \cref{expl:pseudocom_operad} compares pseudo-commutative operads with symmetric monoidal categories.

\secname{sec:BE}  This section discusses translation categories and the Barratt-Eccles $\Gcat$-operad $\BE$.  A $\BE$-algebra is called a naive permutative $G$-category.  It consists of a small $G$-category and a $G$-equivariant permutative category structure.

\secname{sec:GBE}  This section introduces the $G$-Barratt-Eccles $\Gcat$-operad $\GBE$ from the work of Shimakawa \cite{shimakawa89}.  \cref{expl:GBE_algebra} unpacks the structure of a $\GBE$-algebra, which is called a genuine permutative $G$-category.  There is a functor $\Catg(\EG,-)$ from naive permutative $G$-categories to genuine permutative $G$-categories, which we describe in \cref{expl:naive_genuine_pGcat}.

\secname{sec:pseudoalgebra}  This section defines $\Op$-pseudoalgebras, lax $\Op$-morphisms, and $\Op$-transformations for a reduced $\Gcat$-operad $\Op$.  These structures are the objects, 1-cells, and 2-cells of a 2-category $\AlglaxO$.  There are also sub-2-categories $\AlgpspsO$ and $\AlgstO$ with the same objects and 2-cells as $\AlglaxO$, and with 1-cells given by, respectively, $\Op$-pseudomorphisms and strict $\Op$-morphisms.

\secname{sec:naive_smc}  This is the first of two sections that study in detail the 2-category $\AlglaxBE$ for the Barratt-Eccles $\Gcat$-operad $\BE$.  The main observation of this section, \cref{alglaxbe_smgcat}, states that there is a 2-functor
\[\AlglaxBE \fto{\Phi} \smgcat.\]
The codomain $\smgcat$ is the 2-category of naive symmetric monoidal $G$-categories, strictly unital symmetric monoidal $G$-functors, and monoidal $G$-natural transformations (\cref{def:smGcat_twocat}).  There are also variant 2-functors for $\AlgpspsBE$ and $\AlgstBE$, with suitably restricted codomains.  At the object level, $\Phi$ sends $\BE$-pseudoalgebras to naive symmetric monoidal $G$-categories (\cref{BEpseudo_smcat}).  An interesting part of this proof is the derivation of the pentagon axiom \cref{naive_pentagon} and the hexagon axiom \cref{naive_hexagon} from the axioms of a $\BE$-pseudoalgebra.

\secname{sec:BEpseudoalg}  The main result of this section, \cref{thm:BEpseudoalg}, states that each of the three variants of the 2-functor $\Phi$ in \cref{sec:naive_smc} is actually a 2-equivalence.  Thus, from a practical standpoint, the 2-categories $\AlglaxBE$ and $\smgcat$ are interchangeable via $\Phi$.  We emphasize that $\Phi$ is \emph{not} a 2-isomorphism because it is not a bijection on objects.  In other words, $\BE$-pseudoalgebras do \emph{not} correspond bijectively to naive symmetric monoidal $G$-categories.  \cref{expl:BEpseudo_smcat} discusses this subtlety in more detail.  Since we are not aware of any published proof or even a precise statement of \cref{thm:BEpseudoalg} in the literature, we include a detailed proof.

\secname{sec:genuine_smgcat}  This section explains the definitions of $\GBE$-pseudoalgebras, which are called genuine symmetric monoidal $G$-categories, lax $\GBE$-morphisms, and $\GBE$-transformations for the $G$-Barratt-Eccles $\Gcat$-operad $\GBE$.  When $G$ is nontrivial, $\GBE$-pseudoalgebras do not admit a description in terms of a finite list of generating operations, in contrast to \cref{thm:BEpseudoalg}.  The reason for this difference is that $\GBE$ does not have the same kind of coherence properties as $\BE$; see \cref{expl:GBEvsBE_pseudo}.

\section{Intrinsic Pairing and Equivariant Categories}
\label{sec:intrinsic_pairing}

This section discusses several preliminary concepts used in the definition of a pseudo-commutative operad, which is discussed in \cref{sec:pseudocomoperads}.  

\secoutline
\begin{itemize}
\item \cref{def:intrinsic_pairing} defines the intrinsic pairing of an operad. 
\item \cref{ex:as_intrinsic} illustrates the intrinsic pairing with the associative operad.
\item \cref{def:transpose_perm} defines transpose permutations.
\item \cref{ex:intrinsic_transpose} discusses the compatibility between the intrinsic pairing and the transpose permutation for the associative operad.
\item \cref{def:GCat,def:Catg,def:fixedpoint} define the Cartesian closed 2-category $\Gcat$ of small $G$-categories.
\end{itemize}

\subsection*{Intrinsic Pairing}

Pseudo-commutative operads involve the following concept of pairing from \cite[Definition 3.1]{gmmo23}.  Suppose $(\V,\times,\bone)$ is a Cartesian closed category (\cref{def:closedcat}) with a chosen terminal object $\bone$.  Recall from \cref{def:enr-multicategory} that a \emph{$\V$-operad} is a $\V$-multicategory with one object, usually denoted by $*$.  A $\V$-operad $\Op$ is \emph{reduced} if $\Op(0) = \bone$. 

\begin{definition}\label{def:intrinsic_pairing}
Suppose $(\Op,\ga,\opu)$ is a reduced $\V$-operad with $(\V,\times)$ Cartesian closed.  The \emph{intrinsic pairing}\index{intrinsic pairing}
\[\intr \cn (\Op,\Op) \to \Op\]
consists of the following family of composites in $\V$ for $j,k \geq 0$,
\begin{equation}\label{intr_jk}
\intr_{j,k} \cn \Op(j) \times \Op(k) \fto{1 \times \Delta^j} \Op(j) \times \Op(k)^j \fto{\ga} \Op(jk),
\end{equation}
where $\Delta^j$ is the $j$-fold diagonal.  We sometimes abbreviate $\intr_{j,k}$ to $\intr$.
\end{definition}

\begin{explanation}\label{expl:intrinsic_pairing}
Set-theoretically the intrinsic pairing of objects $a\in \Op(j)$ and $b \in \Op(k)$ is the composite
\[a \intr b = \ga\left(a; b,\ldots,b\right)\]
in $\Op(jk)$, where $(b,\ldots,b)$ has $j$ copies of $b$.  In particular, if $a = \opu \in \Op(1)$ is the operadic unit of $\Op$, then
\begin{equation}\label{opu_intr_b}
\opu \intr b = \ga\left(\opu; b\right) = b
\end{equation}
by the left unity axiom \cref{enr-multicategory-left-unity}.  Similarly, if $b = \opu$, then 
\begin{equation}\label{a_intr_opu}
a \intr \opu = \ga\left(a; \opu, \ldots, \opu\right) = a
\end{equation}
by the right unity axiom \cref{enr-multicategory-right-unity}.  
\end{explanation}

\begin{example}[Associative Operad]\label{ex:as_intrinsic}
To under the intrinsic pairing in \cref{def:intrinsic_pairing}, consider the \dindex{associative}{operad}associative operad $\As$ \cite[Section 14.2]{yau-operad}.  It is the $\Set$-operad with $\As(n) = \Sigma_n$, the symmetric group on $n$ letters, for each $n \geq 0$.  The right $\Sigma_n$-action on $\As(n)$ is given by the group multiplication of $\Sigma_n$.  The operadic composition
\[\Sigma_n \times \Sigma_{k_1} \times \cdots \times \Sigma_{k_n} \fto{\ga} \Sigma_{k_1 + \cdots + k_n}\]
is given by
\begin{equation}\label{as_gamma}
\ga\left(\phi; \phi_1,\ldots,\phi_n\right) 
= \phi\ang{k_1,\ldots,k_n} \circ \left(\phi_1 \times \cdots \times \phi_n\right)
\end{equation}
In this composite,
\begin{itemize}
\item $\phi_1 \times \cdots \times \phi_n$ is the block sum \cref{blocksum}, and
\item $\phi\ang{k_1,\ldots,k_n}$ is the block permutation \cref{blockpermutation} that permutes $n$ consecutive blocks of lengths $k_1, \ldots, k_n$, as $\phi$ permutes $\{1,\ldots,n\}$.
\end{itemize}
Algebras over $\As$ are precisely monoids (\cref{def:monoid}).  

The intrinsic pairing of permutations $\sigma \in \Sigma_j$ and $\tau \in \Sigma_k$ is the composite
\begin{equation}\label{as_intrinsic}
\sigma \intr \tau = \ga\left(\sigma; \tau,\ldots,\tau\right) = \sigma\ang{k,\ldots,k} \circ \tau^j \inspace \Sigma_{jk} 
\end{equation}
of
\begin{itemize}
\item the block sum $\tau^j = \tau \times \cdots \times \tau$ of $j$ copies of $\tau$ and
\item the block permutation $\sigma\ang{k,\ldots,k}$ induced by $\sigma$.
\end{itemize}
The intrinsic pairing in \cref{as_intrinsic} is compatible with the intrinsic pairing of $\Op$ in \cref{intr_jk}, in the sense that the following diagram in $\V$ commutes.
\begin{equation}\label{intrinsic_compatible}
\begin{tikzpicture}[xscale=3,yscale=1.3,vcenter]
\draw[0cell=.9]
(0,0) node (x11) {\Op(j) \times \Op(k)}
(x11)++(1,0) node (x12) {\Op(jk)}
(x11)++(0,-1) node (x21) {\Op(j) \times \Op(k)}
(x12)++(0,-1) node (x22) {\Op(jk)}
;
\draw[1cell=.9]  
(x11) edge node {\intr_{j,k}} (x12)
(x21) edge node {\intr_{j,k}} (x22)
(x11) edge node[swap] {\sigma \times \tau} (x21)
(x12) edge node {\sigma \intr \tau} (x22)
;
\end{tikzpicture}
\end{equation}
In this diagram, each of $\sigma$, $\tau$, and $\sigma \intr \tau$ denotes the symmetric group action on $\Op$ \cref{rightsigmaaction} for the indicated permutation.  To see that \cref{intrinsic_compatible} commutes, instead of drawing a large commutative diagram, we compute set-theoretically as follows for $x \in \Op(j)$ and $y \in \Op(k)$.  We write $?^j$ for the $j$-tuple consisting of $j$ copies of the entry $?$.
\[\begin{split}
x\sigma \intr_{j,k} y\tau 
&= \ga\big(x\sigma \sscs (y\tau)^j\big)\\
&= \ga\big(x\sigma \sscs y^j\big) \cdot \tau^j\\
&= \ga\big(x \sscs \sigma y^j\big) \cdot \big(\sigma\ang{k,\ldots,k} \circ \tau^j\big)\\
&= (x \intr_{j,k} y) \cdot (\sigma \intr \tau)
\end{split}\]
The second equality above holds by the bottom equivariance axiom \cref{enr-operadic-eq-2}.  The third equality holds by the symmetric group action axiom \cref{enr-multicategory-symmetry} and the top equivariance axiom \cref{enr-operadic-eq-1}.  The fourth equality holds by $\sigma y^j = y^j$ and \cref{as_intrinsic}.  In \cite[(3.6)]{gmmo23}, the intrinsic pairing $\sigma \intr \tau$ is denoted $\sigma \otimes \tau$.
\end{example}

\subsection*{Transpose Permutations}

Pseudo-commutative operads involve the following permutations.

\begin{definition}\label{def:transpose_perm}
For each $n \geq 0$, we denote by 
\begin{equation}\label{ufsn}
\ufs{n} = \begin{cases} \{1,2,\ldots,n\} & \text{if $n \geq 1$ and}\\
\emptyset & \text{if $n=0$}
\end{cases}
\end{equation}
the \dindex{unpointed}{finite set}unpointed finite set with $n$ elements, equipped with its natural ordering.
\begin{itemize}
\item A sum $\txsum_{j=1}^n$ is also denoted by $\txsum_{j \in \ufs{n}}$, and likewise for other operators that involve a running index.  
\item For $j,k \geq 0$, we define the \emph{$(j,k)$-transpose permutation}\dindex{transpose}{permutation} in $\Sigma_{jk}$
\begin{equation}\label{eq:transpose_perm}
\ufs{jk} \fto[\iso]{\twist_{j,k}} \ufs{kj}
\end{equation}
as the bijection given by
\[\twist_{j,k}\big(b + (a-1)k\big) = a + (b-1)j\]
for $a \in \ufs{j}$ and $b \in \ufs{k}$. 
\item A \emph{transpose permutation} is a $(j,k)$-transpose permutation for some $j,k \geq 0$.\defmark
\end{itemize}
\end{definition}

Our $\ufs{n}$ and $\twist_{j,k}$ are denoted by, respectively, $\ord{n}$ and $\tau_{j,k}$ in \cite{gmmo23}.

\begin{example}\label{ex:transpose_id}
If either $j \in \{0,1\}$ or $k \in \{0,1\}$, then $\twist_{j,k}$ is the identity permutation.  Moreover, the inverse of $\twist_{j,k}$ is $\twist_{k,j}$.
\end{example}

\begin{explanation}\label{expl:transpose_perm}
The \emph{lexicographic ordering}\index{lexicographic ordering} of the product $\ufs{j} \times \ufs{k}$ is the bijection 
\begin{equation}\label{lex_bijection}
\ufs{jk} \fto[\iso]{\lambda_{j,k}} \ufs{j} \times \ufs{k}
\end{equation}
given by 
\[\lambda_{j,k}\big(b + (a-1)k\big) = (a,b)\]
for $a \in \ufs{j}$ and $b \in \ufs{k}$.  We regard an element $(a,b) \in \ufs{j} \times \ufs{k}$ as the $(a,b)$-entry of a $j \times k$ matrix.  Using $\lambda_{j,k}$ and 
\[\ufs{k} \times \ufs{j} \fto[\iso]{\lambda_{k,j}^{-1}} \ufs{kj},\] 
the transpose permutation $\twist_{j,k}$ in \cref{eq:transpose_perm} takes the $(a,b)$-entry in the $j \times k$ matrix $\ufs{j} \times \ufs{k}$ to the $(b,a)$-entry in the $k \times j$ matrix $\ufs{k} \times \ufs{j}$.  This explains the name \emph{transpose permutation}.
\end{explanation}

\begin{example}[Intrinsic Pairing and Transpose Permutation]\label{ex:intrinsic_transpose}
The intrinsic pairing on $\As$ in \cref{as_intrinsic} is compatible with the transpose permutation in \cref{eq:transpose_perm}, in the sense that the following diagram commutes for permutations $\sigma \in \Sigma_j$ and $\tau \in \Sigma_k$.
\begin{equation}\label{intrinsic_transpose}
\begin{tikzpicture}[xscale=2.5,yscale=1.3,vcenter]
\draw[0cell=1]
(0,0) node (x11) {\ufs{jk}}
(x11)++(1,0) node (x12) {\ufs{jk}}
(x11)++(0,-1) node (x21) {\ufs{kj}}
(x12)++(0,-1) node (x22) {\ufs{kj}}
;
\draw[1cell=1]  
(x11) edge node {\sigma \intr \tau} (x12)
(x21) edge node {\tau \intr \sigma} (x22)
(x11) edge node[swap] {\twist_{j,k}} (x21)
(x12) edge node {\twist_{j,k}} (x22)
;
\end{tikzpicture}
\end{equation}
In fact, each of the two composites in \cref{intrinsic_transpose} is given by the assignment
\[\begin{tikzcd}
b + (a-1)k \rar[maps to] & \sigma(a) + \big(\tau(b) - 1\big)j
\end{tikzcd}\]
for $a \in \ufs{j}$ and $b \in \ufs{k}$.
\end{example}

\subsection*{Equivariant Categories, Functors, and Natural Transformations}

We denote by $\Cat$ the 2-category with small categories as objects, functors as 1-cells, and natural transformations as 2-cells (\cref{ex:catastwocategory}).  We regard each group $G$ also as a 2-category (\cref{def:twocategory})  with one object $*$, set of 1-cells $G$, only identity 2-cells, and horizontal composition of 1-cells given by the group multiplication.  Next, we define the 2-category of small $G$-categories and explain its Cartesian closed structure, which involves conjugation $G$-action.  This section ends with a discussion of $G$-fixed subcategories and how they relate to $G$-functors and $G$-natural transformations.

\begin{definition}\label{def:GCat}
For a group $G$, we define $\Gcat$ as the 2-category with the following data.
\begin{itemize}
\item Objects are 2-functors (\cref{def:twofunctor}) $G \to \Cat$.
\item 1-cells are 2-natural transformations (\cref{def:twonaturaltr}).
\item 2-cells are modifications (\cref{def:modification}).
\item Horizontal and vertical compositions are induced by the 2-category $\Cat$ (\cref{def:twonatcomposition,def:modcomposition}).
\end{itemize}
Objects, 1-cells, and 2-cells in $\Gcat$ are called, respectively, \index{G-category@$G$-category}\index{category!equivariant}\index{equivariant!category}\emph{small $G$-categories}, \index{G-functor@$G$-functor}\index{functor!equivariant}\index{equivariant!functor}\emph{$G$-functors}, and \index{G-natural transformation@$G$-natural transformation}\index{natural transformation!equivariant}\index{equivariant!natural transformation}\emph{$G$-natural transformations}.  We omit the word \emph{small} when there is no danger of confusion.  We use the same notation, $\Gcat$, for the underlying 1-category equipped with the symmetric monoidal structure given by Cartesian product and diagonal $G$-action.  The monoidal unit is the terminal $G$-category $\boldone$, with only one object $*$ and its identity morphism.
\end{definition}

\begin{example}\label{ex:CatasGcat}
If $G$ is the trivial group $\{e\}$, then $\Gcat$ is the 2-category $\Cat$.
\end{example}

\begin{explanation}[Unpacking $\Gcat$]\label{expl:GCat}
Unraveling \cref{def:twofunctor,def:twonaturaltr,def:modification} in the context of \cref{def:GCat}, we can describe the 2-category $\Gcat$ more explicitly as follows.  
\begin{description}
\item[Objects] A small $G$-category consists of a small category $\C$ and a \emph{$g$-action} isomorphism
\begin{equation}\label{gactioniso}
\C \fto[\iso]{g} \C \foreachspace g \in G
\end{equation}
such that the following two statements hold.
\begin{enumerate}
\item For the identity element $e \in G$, the $e$-action is the identity functor $1_\C$.
\item For $g, h \in G$, the $hg$-action is equal to the composition $h \circ g$.
\end{enumerate}
For $g \in G$ and an object or morphism $x \in \C$, we denote the $g$-action on $x$ by either $g \cdot x$ or $gx$.

\item[1-cells]
For small $G$-categories $\C$ and $\D$, a $G$-functor $F \cn \C \to \D$, which we also call a \emph{$G$-equivariant functor}, is a functor between the underlying categories such that the following diagram commutes for each $g \in G$.
\begin{equation}\label{Gfunctor}
\begin{tikzpicture}[xscale=2,yscale=1.3,vcenter]
\draw[0cell=.9]
(0,0) node (x11) {\C}
(x11)++(1,0) node (x12) {\D}
(x11)++(0,-1) node (x21) {\C}
(x12)++(0,-1) node (x22) {\D}
;
\draw[1cell=.9]  
(x11) edge node {F} (x12)
(x21) edge node {F} (x22)
(x11) edge node[swap] {g} (x21)
(x12) edge node {g} (x22)
;
\end{tikzpicture}
\end{equation}

\item[2-cells]
For $G$-functors $F,F' \cn \C \to \D$, a $G$-natural transformation $\theta \cn F \to F'$, which we also call a \emph{$G$-equivariant natural transformation}, is a natural transformation $F \to F'$ such that
\begin{equation}\label{Gnattr}
g \cdot \theta_c = \theta_{gc} \cn g \cdot Fc = F(gc) \to g \cdot F'c = F'(gc)
\end{equation}
for each $g \in G$ and object $c \in \C$.
\end{description}
We also use \cref{gactioniso,Gfunctor,Gnattr} as definitions of, respectively, $G$-categories, $G$-functors, and $G$-natural transformations, without assuming the underlying categories are small.
\end{explanation}

There is a closed structure on $\Gcat$ that uses the following $G$-categories.

\begin{definition}[Internal Hom for $\Gcat$]\label{def:Catg}
For small $G$-categories $\C$ and $\D$, we denote by \index{category!internal hom}$\Catg(\C,\D)$ the small $G$-category with
\begin{itemize}
\item functors $\C \to \D$ between the underlying categories as objects,
\item natural transformations as morphisms, 
\item identities and composition defined componentwise in $\D$, and
\item $G$ acting by conjugation.
\end{itemize}
The \emph{conjugation $G$-action}\index{conjugation G-action@conjugation $G$-action} means that, for functors $F,F' \cn \C \to \D$, a natural transformation $\theta \cn F \to F'$, and an element $g \in G$, the functor $g \cdot F$ and the natural transformation $g \cdot \theta$ are given by the following composite and whiskering.
\begin{equation}\label{conjugation-gaction}
\begin{tikzpicture}[xscale=1.5,yscale=1.5,vcenter]
\draw[0cell=.9]
(0,0) node (a) {\C}
(a)++(1,0) node (b) {\C}
(b)++(1.2,0) node (c) {\D}
(c)++(1,0) node (d) {\D}
;
\draw[1cell=.9]  
(a) edge node {g^{-1}} (b)
(b) edge[bend left] node {F} (c)
(b) edge[bend right] node[swap] {F'} (c)
(c) edge node {g} (d)
;
\draw[2cell]
node[between=b and c at .43, rotate=-90, 2label={above,\theta}] {\Rightarrow}
;
\end{tikzpicture}
\end{equation}
This finishes the definition of the small $G$-category $\Catg(\C,\D)$.
\end{definition}

We emphasize that the objects and morphisms in $\Catg(\C,\D)$ are not necessarily $G$-equivariant.  

\begin{explanation}[Cartesian Closed Structure on $\Gcat$]\label{expl:Gcat_closed}
Since $(\Cat,\times,\bone)$ is complete, cocomplete, and Cartesian closed (\cref{def:closedcat}), so is 
\[\left(\Gcat,\times,\bone\right),\] 
with limits and colimits computed in $\Cat$.  In particular, the product of small $G$-categories $\C$ and $\D$ is the product category $\C \times \D$ with diagonal $G$-action
\[G \fto{\Delta} G \times G \fto{\C \times \D} \Cat \times \Cat.\]
Products of $G$-functors or $G$-natural transformations are defined as products of underlying functors or natural transformations.  The closed structure on $\Gcat$ is given by the small $G$-categories $\Catg(\C,\D)$ in \cref{def:Catg}.  There is an isomorphism of categories
\begin{equation}\label{Catg_rightadj}
\Gcat(\A \times \C, \D) \iso \Gcat(\A, \Catg(\C,\D))
\end{equation}
that is natural in small $G$-categories $\A$, $\C$, and $\D$.
\end{explanation}

Next, we explain another way to relate $\Gcat$ and $\Catg$ via fixed points.

\begin{definition}[$G$-Fixed Subcategories]\label{def:fixedpoint}
Suppose $\C$ and $\D$ are $G$-categories.  
\begin{itemize}
\item An object or a morphism $c \in \C$ is \emph{$G$-fixed}\index{G-fixed@$G$-fixed} if 
\[gc = c \forallspace g \in G.\]
\item The \emph{$G$-fixed subcategory} of $\C$, denoted $\C^G$, is the subcategory consisting of $G$-fixed objects and $G$-fixed morphisms.
\item For a $G$-functor $F \cn \C \to \D$, its subfunctor between $G$-fixed subcategories, called the \emph{$G$-fixed subfunctor} of $F$, is denoted by
\[\C^G \fto{F^G} \D^G.\]
\item For $G$-functors $F, F' \cn \C \to \D$, the restriction of a $G$-natural transformation $\theta \cn F \to F'$ to $G$-fixed objects, called the \emph{$G$-fixed subnatural transformation} of $\theta$, is denoted by
\[F^G \fto{\theta^G} (F')^G.\]
\end{itemize}
This finishes the definition.
\end{definition}

\begin{explanation}[$\Gcat(\C,\D)$ as $G$-Fixed Subcategory]\label{expl:gcat_catg}
For small $G$-categories $\C$ and $\D$, the $G$-category $\Catg(\C,\D)$ has the conjugation $G$-action \cref{conjugation-gaction}.  A functor $F \cn \C \to \D$ is fixed by the conjugation $G$-action precisely when $F$ is $G$-equivariant, and the same is true for a natural transformation $\theta \cn F \to F'$.  Thus, there is an equality 
\[\Gcat(\C,\D) = \Catg(\C,\D)^G\]
between the category $\Gcat(\C,\D)$ of $G$-functors and $G$-natural transformations (\cref{expl:GCat}) and the $G$-fixed subcategory of $\Catg(\C,\D)$.
\end{explanation}

\section{Pseudo-Commutative Operads}
\label{sec:pseudocomoperads}

This section recalls the definition of a \emph{pseudo-commutative operad} from \cite{gmmo23}.  This definition is originally from \cite{corner-gurski,hyland-power}, which define \emph{pseudo-commutative 2-monads}.

\secoutline
\begin{itemize}
\item \cref{def:pseudocom_operad} defines pseudo-commutative operads in $\Gcat$.
\item \cref{expl:pseudocom_operad} further elaborates the axioms of a pseudo-commutative operad. 
\item \cref{rk:pseudocom_generality,rk:gmmo-cg-notation} compare the context of this chapter, \cite{gmmo23}, and \cite{corner-gurski}.
\end{itemize}

The following definition uses the Cartesian closed 2-category $\Gcat$ in \cref{def:GCat,expl:Gcat_closed}.  

\begin{definition}\label{def:pseudocom_operad}
Suppose $G$ is a group.  A \emph{pseudo-commutative operad}\dindex{pseudo-commutative}{operad} $(\Op,\ga,\opu,\pcom)$ in $\Gcat$ consists of the following data.
\begin{description}
\item[Reduced operad] $(\Op,\ga,\opu)$ is a reduced $\Gcat$-operad (\cref{def:enr-multicategory}) for the Cartesian closed category $\Gcat$.  The assumption of being reduced means $\Op(0) = \bone$, a terminal $G$-category.
\item[Pseudo-commutative structure] For $j,k \geq 0$, it is equipped with a $G$-natural isomorphism $\pcom_{j,k}$, called the \emph{$(j,k)$-pseudo-commutativity isomorphism}, as follows.
\begin{equation}\label{pseudocom_isos}
\begin{tikzpicture}[xscale=3,yscale=1.3,vcenter]
\draw[0cell=.9]
(0,0) node (x11) {\Op(j) \times \Op(k)}
(x11)++(1,0) node (x12) {\Op(jk)}
(x11)++(0,-1) node (x21) {\Op(k) \times \Op(j)}
(x12)++(0,-1) node (x22) {\Op(kj)}
;
\draw[1cell=.9]  
(x11) edge node {\intr_{j,k}} (x12)
(x21) edge node[swap] {\intr_{k,j}} (x22)
(x11) edge node[swap] {\twist} (x21)
(x12) edge node {\twist_{k,j}} (x22)
;
\draw[2cell]
node[between=x11 and x22 at .5, rotate=-90, 2label={above,\pcom_{j,k}}, 2label={below,\iso}] {\Rightarrow}
;
\end{tikzpicture}
\end{equation}
We call $\pcom = \{\pcom_{j,k}\}_{j,k \geq 0}$ the \emph{pseudo-commutative structure}.  In \cref{pseudocom_isos}, the left vertical arrow $\twist$ is the braiding for $\Gcat$, which swaps the two arguments.  The horizontal arrows are components of the intrinsic pairing in \cref{intr_jk}.  The right vertical arrow is the symmetric group action $G$-functor on $\Op$ for the $(k,j)$-transpose permutation $\twist_{k,j}$ \cref{eq:transpose_perm}.  For a pair of objects $(x,y) \in \Op(j) \times \Op(k)$, the $(x,y)$-component of $\pcom_{j,k}$ is an isomorphism
\[(x \intr_{j,k} y) \twist_{k,j} \fto[\iso]{\pcom_{j,k; x,y}} y \intr_{k,j} x \inspace \Op(kj).\]
The subscripts of $\pcom$ and $\intr$ are sometimes suppressed to simplify the typography.
\end{description}
The data $(\Op,\ga,\opu,\pcom)$ are required to satisfy the following four axioms whenever they are defined.
\begin{description}
\item[Unity] For the operadic unit $\opu \in \Op(1)$ and $y \in \Op(k)$, the $(\opu,y)$-component
\begin{equation}\label{pseudocom_unity}
(\opu \intr_{1,k} y) \twist_{k,1} \fto{\pcom_{1,k; \opu,y}} y \intr_{k,1} \opu \inspace \Op(k)
\end{equation}
is equal to the identity $1_y \cn y \to y$ in $\Op(k)$.  This axiom makes sense because
\[\opu \intr_{1,k} y = y = y \intr_{k,1} \opu\]
by \cref{opu_intr_b,a_intr_opu}, while $\twist_{k,1}$ is the identity permutation by \cref{ex:transpose_id}.

\item[Symmetry] The composite of the following pasting diagram is equal to the identity $G$-natural transformation of $\intr_{j,k}$.
\begin{equation}\label{pseudocom_sym}
\begin{tikzpicture}[xscale=3,yscale=1.3,vcenter]
\draw[0cell=.9]
(0,0) node (x11) {\Op(j) \times \Op(k)}
(x11)++(1,0) node (x12) {\Op(jk)}
(x11)++(0,-1) node (x21) {\Op(k) \times \Op(j)}
(x12)++(0,-1) node (x22) {\Op(kj)}
(x21)++(0,-1) node (x31) {\Op(j) \times \Op(k)}
(x22)++(0,-1) node (x32) {\Op(jk)}
;
\draw[1cell=.9]  
(x11) edge node {\intr_{j,k}} (x12)
(x21) edge node[swap] {\intr_{k,j}} (x22)
(x11) edge node[swap] {\twist} (x21)
(x12) edge node {\twist_{k,j}} (x22)
(x31) edge node[swap] {\intr_{j,k}} (x32)
(x21) edge node[swap] {\twist} (x31)
(x22) edge node {\twist_{j,k}} (x32)
;
\draw[2cell=.9]
node[between=x11 and x22 at .5, rotate=-90, 2label={above,\pcom_{j,k}}] {\Rightarrow}
node[between=x21 and x32 at .5, shift={(0,-.1)}, rotate=-90, 2label={above,\pcom_{k,j}}] {\Rightarrow}
;
\end{tikzpicture}
\end{equation}
This axiom makes sense because $\twist \circ \twist = 1$ and $\twist_{k,j} = \twist_{j,k}^{-1}$, so each of the left and right boundaries in \cref{pseudocom_sym} is the identity functor.

\item[Equivariance] For permutations $\sigma \in \Sigma_j$ and $\tau \in \Sigma_k$, the following pasting diagram equality holds.
\begin{equation}\label{pseudocom_equiv}
\begin{tikzpicture}[xscale=1,yscale=1,vcenter]
\def\h{2.5} \def\v{-1.3} \def\s{.9} 
\draw[0cell=\s]
(0,0) node (x11) {\Op(j) \times \Op(k)}
(x11)++(\h,0) node (x12) {\Op(jk)}
(x11)++(0,\v) node (x21) {\Op(k) \times \Op(j)}
(x12)++(0,\v) node (x22) {\Op(kj)}
(x21)++(0,\v) node (x31) {\Op(k) \times \Op(j)}
(x22)++(0,\v) node (x32) {\Op(kj)}
;
\draw[1cell=\s] 
(x11) edge node {\intr_{j,k}} (x12)
(x21) edge node[swap] {\intr_{k,j}} (x22)
(x31) edge node [swap]{\intr_{k,j}} (x32)
(x11) edge node[swap] {\twist} (x21)
(x21) edge node[swap] {\tau \times \sigma} (x31)
(x12) edge node {\twist_{k,j}} (x22)
(x22) edge node {\tau \intr \sigma} (x32)
;
\draw[2cell=\s]
node[between=x11 and x22 at .5, rotate=-90, 2label={above,\pcom_{j,k}}] {\Rightarrow}
;
\draw[0cell=1.5]
(x22)++(1.3,0) node (eq) {=}
;
\draw[0cell=\s]
(eq)++(1.6,-\v) node (y11) {\Op(j) \times \Op(k)}
(y11)++(\h,0) node (y12) {\Op(jk)}
(y11)++(0,\v) node (y21) {\Op(j) \times \Op(k)}
(y12)++(0,\v) node (y22) {\Op(jk)}
(y21)++(0,\v) node (y31) {\Op(k) \times \Op(j)}
(y22)++(0,\v) node (y32) {\Op(kj)}
;
\draw[1cell=\s] 
(y11) edge node {\intr_{j,k}} (y12)
(y21) edge node {\intr_{j,k}} (y22)
(y31) edge node[swap] {\intr_{k,j}} (y32)
(y11) edge node[swap] {\sigma \times \tau} (y21)
(y21) edge node[swap] {\twist} (y31)
(y12) edge node {\sigma \intr \tau} (y22)
(y22) edge node {\twist_{k,j}} (y32)
;
\draw[2cell=\s]
node[between=y21 and y32 at .5, rotate=-90, 2label={above,\pcom_{j,k}}] {\Rightarrow}
;
\end{tikzpicture}
\end{equation}
The two unlabeled regions are commutative by \cref{intrinsic_compatible}.  The left vertical boundaries of the two sides are equal by the naturality of the braiding $\twist$ for $\Gcat$.  The right vertical boundaries of the two sides are equal by \cref{intrinsic_transpose}.

\item[Operadic compatibility] Consider objects $x \in \Op(j)$, $y_i \in \Op(k_i)$ for $1 \leq i \leq j$, and $z \in \Op(\ell)$, along with the following notation for $r \in \{1,\ldots,j\}$.
\[\left\{\scalebox{.85}{$\begin{split}
k &= k_1 + \cdots + k_j \qquad y_r^\ell = (\overbracket[.5pt]{y_r,\ldots,y_r}^{\ell}) \in \Op(k_r)^\ell \\
y_\crdot^\ell &= \big(y_1^\ell, \ldots, y_j^\ell\big)\\
\twist_{\ell,k_\centerdot} &= \twist_{\ell,k_1} \times \cdots \times \twist_{\ell,k_j} \in \Sigma_{\ell k} \\
\Twist_{\ell,k_\centerdot} &= \big(\twist_{k_1,\ell} \times \cdots \times \twist_{k_j,\ell}\big) \twist_{\ell,k} \in \Sigma_{\ell k}\\
\boldy &= (y_1,\ldots,y_j) \in \Op(k_1) \times \cdots \times \Op(k_j)\\
\boldy^\ell &= (\overbracket[.5pt]{\boldy,\ldots,\boldy}^{\ell}) \in \left(\Op(k_1) \times \cdots \times \Op(k_j)\right)^\ell\\
z \intr \boldy & = \left(z \intr y_1, \ldots, z \intr y_j\right) \in \Op(\ell k_1) \times \cdots \times \Op(\ell k_j)\\ 
\boldy \intr z &= \left(y_1 \intr z, \ldots, y_j \intr z \right) \in \Op(k_1 \ell) \times \cdots \times \Op(k_j \ell)\\
(\bold y \intr z) \twist_{\ell,\bdot} &= \big((y_1 \intr z)\twist_{\ell,k_1}, \ldots, (y_j \intr z)\twist_{\ell,k_j}\big) \in \Op(\ell k_1) \times \cdots \times \Op(\ell k_j)\\
\pcom_{k_\centerdot,\ell} &= \big( \pcom_{k_1,\ell}, \ldots, \pcom_{k_j,\ell}\big)\\
\end{split}$}\right.\]
With the notation above, the following diagram in $\Op(\ell k)$ commutes.
\begin{equation}\label{pseudocom_compat}
\begin{tikzpicture}[xscale=1,yscale=1,vcenter]
\def\h{2.5} \def\f{.85} \def\u{-1.2} \def\v{-1.3}
\draw[0cell=.9]
(0,0) node (x) {\ga\big(x \sscs (\bold y \intr z) \twist_{\ell,\bdot}\big) \Twist_{\ell,k_\centerdot}}
(x)++(-\h,\u) node (x11) {\ga\left(x \sscs z \intr \boldy\right) \Twist_{\ell,k_\centerdot}}
(x11)++(0,\v) node (x21) {\ga\big(x \intr z \sscs y_\crdot^\ell\big) \Twist_{\ell,k_\centerdot}}
(x21)++(0,\v) node (x31) {\ga\big((x \intr z) \twist_{\ell,j} \sscs \boldy^\ell\big)}
(x31)++(\f,\v) node (x41) {\ga\big(z \intr x \sscs \boldy^\ell\big)}
(x)++(\h,\u) node (x12) {\ga\left(x \sscs \boldy \intr z\right) (\twist_{\ell, k_\centerdot} \Twist_{\ell,k_\centerdot})}
(x12)++(0,\v) node (x22) {\ga\left(x \sscs \boldy \intr z\right) \twist_{\ell,k}}
(x22)++(0,\v) node (x32) {\big(\ga\left(x \sscs \boldy\right) \intr z\big) \twist_{\ell,k}}
(x32)++(-\f,\v) node (x42) {z \intr \ga\left(x \sscs \boldy\right)}
;
\draw[1cell=.8]  
(x) edge[transform canvas={xshift={-1.3em}}] node[swap,pos=.5] {\ga\left(1 \sscs \pcom_{k_\centerdot,\ell}\right) \Twist_{\ell,k_\centerdot}} (x11)
(x11) edge[equal] node[swap] {\mathrm{(a)}} (x21)
(x21) edge[equal] node[swap] {\mathrm{(top)}} (x31)
(x31) edge[shorten <=-2pt, shorten >=-2pt] node[swap, pos=.2] {\ga( \pcom_{j,\ell} \sscs 1^{\ell j} )} (x41)
(x41) edge[equal] node {\mathrm{(a)}} (x42)
(x) edge[equal, transform canvas={xshift={1.3em}}] node {\mathrm{(bot)}} (x12)
(x12) edge[equal] node {\dagger} (x22)
(x22) edge[equal] node {\mathrm{(a)}} (x32)
(x32) edge[shorten <=-1pt, shorten >=-1pt] node[pos=.3] {\pcom_{k,\ell}} (x42)
;
\end{tikzpicture}
\end{equation}
The three equalities labeled $\mathrm{(a)}$ in \cref{pseudocom_compat} hold by the associativity axiom \cref{enr-multicategory-associativity}.  The equalities labeled $\mathrm{(top)}$ and $\mathrm{(bot)}$ hold by, respectively, the top equivariance axiom \cref{enr-operadic-eq-1} and the bottom equivariance axiom \cref{enr-operadic-eq-2}.  The equality labeled $\dagger$ follows from 
\[\twist_{\ell, k_\centerdot} \Twist_{\ell,k_\centerdot} = \twist_{\ell,k},\]
which holds because $\twist_{\ell,k_r} = \twist_{k_r,\ell}^{-1}$ for $r \in \{1,\ldots,j\}$.
\end{description}
This finishes the definition of a pseudo-commutative operad in $\Gcat$.  
\end{definition}

\begin{explanation}[Pseudo-Commutativity Axioms]\label{expl:pseudocom_operad}
To understand the axioms of a pseudo-commutative operad, it is helpful to compare them to the axioms of a symmetric monoidal category.  The following table provides a rough comparison of the two notions.
\begin{center}
\resizebox{\columnwidth}{!}{%
{\renewcommand{\arraystretch}{1.3}%
{\setlength{\tabcolsep}{1em}
\begin{tabular}{c|c} 
Pseudo-commutative operads (\ref{def:pseudocom_operad}) & Symmetric monoidal categories (\ref{def:symmoncat})\\ \hline 
intrinsic pairing \cref{intr_jk} & monoidal product (\ref{def:monoidalcategory}) \\ 
pseudo-commutative structure \cref{pseudocom_isos} & braiding \cref{braiding_bmc} \\
symmetry axiom \cref{pseudocom_sym} & symmetry axiom \cref{symmoncatsymhexagon} \\
equivariance axiom \cref{pseudocom_equiv} & no analogue \\
operadic compatibility \cref{pseudocom_compat} & right hexagon axiom \cref{hexagon-braided}\\ 
\end{tabular}}}}
\end{center}
\smallskip
We elaborate on this comparison below.
\begin{description}
\item[Symmetry] The symmetry axiom \cref{pseudocom_sym} means that the following diagram in $\Op(jk)$ commutes for objects $x \in \Op(j)$ and $y \in \Op(k)$.
\begin{equation}\label{pseudocom_sym_expl}
\begin{tikzpicture}[xscale=4,yscale=1.3,vcenter]
\draw[0cell=.9]
(0,0) node (x11) {x \intr_{j,k} y}
(x11)++(0,-1) node (x21) {(x \intr_{j,k} y) \twist_{k,j} \twist_{j,k}}
(x21)++(1,0) node (x22) {(y \intr_{k,j} x) \twist_{j,k}}
(x22)++(0,1) node (x12) {x \intr_{j,k} y}
;
\draw[1cell=.9]  
(x11) edge node {1} (x12)
(x11) edge[equal] (x21)
(x21) edge node {\pcom_{j,k} \twist_{j,k}} (x22)
(x22) edge node[swap] {\pcom_{k,j}} (x12)
;
\end{tikzpicture}
\end{equation}
This axiom is analogous to the symmetry axiom \cref{symmoncatsymhexagon} of a symmetric monoidal category 

\item[Equivariance]
The equivariance axiom \cref{pseudocom_equiv} means that the following diagram in $\Op(kj)$ commutes.
\begin{equation}\label{pseudocom_equiv_expl}
\begin{tikzpicture}[xscale=1,yscale=1,vcenter]
\def\h{4} \def\v{-1.3}
\draw[0cell=.9]
(0,0) node (x11) {(x \intr_{j,k} y) \twist_{k,j} (\tau \intr \sigma)}
(x11)++(0,\v) node (x21) {(y \intr_{k,j} x) (\tau \intr \sigma)}
(x11)++(\h,0) node (x12) {(x\sigma \intr_{j,k} y\tau) \twist_{k,j}}
(x12)++(0,\v) node (x22) {y\tau \intr_{k,j} x\sigma}
(x11)++(\h/2,1) node (x) {(x \intr_{j,k} y) (\sigma \intr \tau) \twist_{k,j}}
;
\draw[1cell=.9]  
(x11) edge node[swap] {\pcom_{j,k} (\tau \intr \sigma)} (x21)
(x21) edge[equal] (x22)
(x11) edge[equal, transform canvas={xshift={-1em}}] (x)
(x) edge[equal, transform canvas={xshift={1em}}] (x12)
(x12) edge node {\pcom_{j,k}} (x22)
;
\end{tikzpicture}
\end{equation}
This axiom encodes equivariant properties of $\pcom_{j,k}$ with respect to permutations of the form $\sigma \intr \tau$.  It does not have a symmetric monoidal category analogue.

\item[Operadic compatibility]
The operadic compatibility axiom \cref{pseudocom_compat} is analogous to the right hexagon axiom in \cref{hexagon-braided}, which also holds in a symmetric monoidal category (\cref{rk:smcat}).  It expresses the idea that the two ways of bringing $z$ pass $(x,\boldy)$ using the pseudo-commutative structure $\pcom$ are equal.  Along the left boundary of \cref{pseudocom_compat}, $z$ first passes $\boldy$ and then $x$.  Along the right boundary, $z$ passes $(x,\boldy)$ all at once.

When we apply the top equivariance axiom \cref{enr-operadic-eq-1} to obtain the equality labeled $\mathrm{(top)}$ in \cref{pseudocom_compat}, we are using the following two facts about the transpose permutation $\twist_{\ell,j}$ in \cref{eq:transpose_perm}.
\begin{enumerate}
\item $\twist_{\ell,j}$ sends $\boldy^\ell$ to $y^\ell_\crdot$:
\[\twist_{\ell,j} (\boldy^\ell) = y^\ell_\crdot \inspace \Op(k_1)^\ell \times \cdots \times \Op(k_j)^\ell.\]
\item There is an equality of permutations
\begin{equation}\label{block_transpose}
\Twist_{\ell,k_\centerdot} = \twist_{\ell,j} \big\langle \, \overbrace{k_1,\ldots,k_j}^{1} \, , \ldots, \, \overbrace{k_1,\ldots,k_j}^{\ell} \, \big\rangle \inspace \Sigma_{\ell k}.
\end{equation}
The right-hand side is the block permutation induced by $\twist_{\ell,j}$ that permutes $\ell j$ blocks of lengths $(k_1,\ldots,k_j), \ldots, (k_1,\ldots,k_j)$, where the $j$-tuple $(k_1,\ldots,k_j)$ is repeated $\ell$ times.
\end{enumerate}
Because of \cref{block_transpose}, we call $\Twist_{\ell,k_\centerdot}$ a \emph{block transpose}.\defmark
\end{description}
\end{explanation}

\begin{remark}[Generality of Pseudo-Commutative Operads]\label{rk:pseudocom_generality}
In \cref{def:pseudocom_operad}, we specify that $\Op$ is a reduced operad in $\Gcat$.  However, the same definition makes sense in greater generality, as explained in \cite[Section 2]{gmmo23}.  In more details, suppose $\cV$ is a complete and cocomplete Cartesian closed category.  Denote by $\Cat(\cV)$ the 2-category with
\begin{itemize}
\item internal $\cV$-categories\dindex{internal}{category} as objects,
\item internal $\cV$-functors as 1-cells, and
\item internal $\cV$-natural transformations as 2-cells.
\end{itemize}
Then \cref{def:pseudocom_operad} makes sense for any reduced operad $\Op$ in $\Cat(\cV)$.  We work with $\Gcat$ instead of $\Cat(\cV)$ in \cref{def:pseudocom_operad} because, when we discuss pseudo-commutative operads, we only use $\Gcat$ as our ambient Cartesian closed 2-category.
\end{remark}

\begin{remark}[Notational Comparison]\label{rk:gmmo-cg-notation}
For the reader's convenience, the table below summarizes the context, notation, and references regarding pseudo-commutative operads in this work, \cite{gmmo23}, and \cite{corner-gurski}.  References in the last two columns are those in \cite{gmmo23} and \cite{corner-gurski}, respectively. 
\begin{center}
\resizebox{\columnwidth}{!}{%
{\renewcommand{\arraystretch}{1.3}%
{\setlength{\tabcolsep}{1em}
\begin{tabular}{c|ccc}
& This chapter & \cite{gmmo23} & \cite{corner-gurski} \\ \hline
\multirow{2}{*}{\makecell{pseudo-commutative \\ operads}} & \multirow{2}{*}{\makecell{\cref{def:pseudocom_operad}}} & \multirow{2}{*}{\makecell{Definition 3.10}} & \multirow{2}{*}{\makecell{Theorem 4.4}}\\
&&&\\
context & $\Gcat$ (\ref{expl:Gcat_closed}) & $\Cat(\cV)$ (Section 2.1) & $\Cat$\\
basic structure & operads (\ref{def:enr-multicategory}) & operads & $\mathbf{G}$-operads (Def.\! 1.14)\\
operadic composition & $\ga$ & $\ga$ & $\mu$\\
intrinsic pairing & $\intr$ \cref{intr_jk} & $\intr$ (Def.\! 3.1) & $\mu(p; q,\ldots,q)$ or $\mu(p;\underline{q})$\\
\multirow{2}{*}{\makecell{pseudo-commutative \\ structure}} & \multirow{2}{*}{\makecell{$\pcom_{j,k}$ \cref{pseudocom_isos}}} & \multirow{2}{*}{\makecell{$\alpha_{j,k}$ (3.11)}} & \multirow{2}{*}{\makecell{$\lambda_{p,q}$}}\\
&&&\\
transpose permutation & $\twist_{j,k}$ \cref{eq:transpose_perm} & $\tau_{j,k}$ (Def.\! 3.8) & $\tau_{m,n}$ or $t_{m,n}$\\
block transpose & $\Twist_{\ell,k_\centerdot}$ \cref{block_transpose} & $D_{\ell, k_*}$ (Def.\! 11.1) & $\mu(t_{n,l}; \underline{e_{m_1}, \ldots, e_{m_l}})$\\
\end{tabular}}}}
\end{center}
\smallskip
In this work and \cite{gmmo23}, $G$ denotes a (finite) group.  On the other hand, in \cite[Definition 1.10]{corner-gurski}, $\mathbf{G}$ denotes an \dindex{action}{operad}action operad, which can be used to obtain a common generalization of (nonsymmetric) operads, braided operads, ribbon operads, and cactus operads.  We will not discuss $\mathbf{G}$-operads in this work.  The reader is referred to \cite{corner-gurski,yau-inf-operad} for further development of $\mathbf{G}$-operads.
\end{remark}

\section{Naive Permutative $G$-Categories}
\label{sec:BE}

This section discusses the categorical Barratt-Eccles operad $\BE$.  A $G$-equivariant analogue $\GBE$, called the $G$-Barratt-Eccles operad, is discussed in \cref{sec:GBE}.

\secoutline
\begin{itemize}
\item \cref{def:translation_cat} recalls translation categories.  
\item \cref{def:BE} defines the Barratt-Eccles operad $\BE$ in $\Cat$.
\item \cref{ex:barratt_eccles} discusses $\BE$-algebras in $\Cat$, which are small permutative categories.
\item \cref{def:BE-Gcat,expl:naive_perm_Gcat} discuss $\BE$-algebras in $\Gcat$, which are naive permutative $G$-categories.
\item \cref{translation_pseudocom} records the fact that each reduced $\Gcat$-operad that is levelwise a translation category has a unique pseudo-commutative structure (\cref{def:pseudocom_operad}).  In particular, both $\BE$ and $\GBE$ are pseudo-commutative operads in a unique way.   We regard $\BE$ and $\GBE$ as our primary examples of pseudo-commutative operads. 
\end{itemize}

\subsection*{Translation Categories}
The Barratt-Eccles operad is constructed using the following procedure.

\begin{definition}\label{def:translation_cat}
For a set $S$, the \emph{translation category}\dindex{translation}{category} $\tn S$ is the small category with
\begin{itemize}
\item object set $S$;
\item each hom set given by the one-element set $*$; and 
\item identity morphisms and composition uniquely determined by the fact that $*$ is a terminal object in the category $\Set$.
\end{itemize}  
For objects $a, b \in \tn S$, we denote the unique morphism from $a$ to $b$ by either \label{not:b_to_a}$[b,a] \cn a \to b$ or $! \cn a \to b$.
\end{definition}

\begin{example}[Translation Category of a $G$-Set]\label{ex:EG}
For a group $G$, a \emph{$G$-set}\index{G-set@$G$-set} is a pair $(S,\mu)$ consisting of a set $S$ and a group homomorphism 
\[\mu \cn G \to \Aut(S),\]
where $\Aut(S)$ is the group of bijections $S \fto{\iso} S$.  For a $G$-set $(S,\mu)$, the translation category $\tn S$ is a small $G$-category with $g$-action isomorphism \cref{gactioniso}
\[g \cn \tn S \fto{\iso} \tn S\]
given on objects by the $g$-action on $S$,
\[g \cdot a = \mu(g)(a) \forspace (g,a) \in G \times S.\]
The $g$-action on morphisms is given by the unique morphism
\[g \cdot [b,a] = [g \cdot b, g \cdot a] \cn g\cdot a \to g \cdot b \inspace \tn S\]
for objects $a, b \in \tn S$.

We equip $G$ with the structure of a $G$-set $(G,\mu)$ via the group multiplication $G \times G \to G$, so 
\[\mu(g)(h) = gh \forspace g,h \in G.\]
As a special case of the previous paragraph, the translation category $\EG$ is a small $G$-category.
\end{example}

Each of the four pseudo-commutativity axioms in \cref{def:pseudocom_operad} asks for the commutativity of some diagrams in some $\Op(n)$.  For translation categories, such diagrams are automatically commutative by the terminal property of $*$.  Therefore, we have the following observation from \cite[Corollary 4.9]{corner-gurski} and \cite[Lemma 3.12]{gmmo23}.

\begin{proposition}\label{translation_pseudocom}
Suppose $\Op$ is a reduced $\Gcat$-operad such that each $\Op(n)$ is a translation category.  Then $\Op$ admits a unique pseudo-commutative structure.
\end{proposition}

\subsection*{Barratt-Eccles Operad}
The following definition uses the Cartesian closed category $\Cat$ in \cref{ex:cat} and the associative operad $\As$ in \cref{ex:as_intrinsic}.

\begin{definition}\label{def:BE}
The \emph{Barratt-Eccles operad}\dindex{Barratt-Eccles}{operad} is the $\Cat$-operad $\BE$ given by the translation categories
\[\BE(n) = \tn\As(n) = \tn\Sigma_n \forspace n \geq 0.\]
The operad structure on objects is the one on $\As$.  On morphisms, the operad structure is uniquely determined by the terminal property of each hom set in $\tn\Sigma_n$.
\end{definition}

\begin{example}[(Bi)permutative Categories]\label{ex:barratt_eccles}
The Barratt-Eccles operad $\BE$ has the following algebras.  Recall from \cref{def:enr-multicategory-functor} that, for $\V$-multicategories $\M$ and $\N$, an $\M$-algebra in $\N$ is a $\V$-multifunctor $\M \to \N$.
\begin{itemize}
\item $\BE$-algebras in $\Cat$ are $\Cat$-multifunctors 
\[\BE \to \Cat.\]  
Coherence properties for $\BE$ imply that $\BE$-algebras in $\Cat$ are precisely small \dindex{permutative}{category}permutative categories.  See \cite[11.4.14 and 11.4.26]{cerberusIII} for details. 
\item Similarly, using a suitable $\Cat$-multicategory $\permcatsu$ of small permutative categories, $\BE$-algebras in $\permcatsu$ are precisely small \dindex{bipermutative}{category}bipermutative categories.  See \cite[11.5.5]{cerberusIII} for details.
\end{itemize}
The upshot is that $\BE$-algebras are the permutative objects in the target $\Cat$-multicategory.
\end{example}

\subsection*{Naive Permutative $G$-Categories}
We now consider the $G$-equivariant context.  The symmetric monoidal closed category \pcref{def:GCat,def:Catg} 
\[(\Gcat,\times,\boldone,\Catg)\] 
is regarded as a $\Gcat$-multicategory, where the $k$-ary multimorphism $G$-categories for $k \geq 0$ are given by the closed structure $\Catg$.  See \cref{theorem:v-closed-v-sm,proposition:monoidal-v-cat-v-multicat}.  The following definition is from \cite{gm17,gmmo20}.

\begin{definition}[Barratt-Eccles Operad in $\Gcat$]\label{def:BE-Gcat}
Given a group $G$, we regard the Barratt-Eccles operad $\BE$ in \cref{def:BE} as a reduced $\Gcat$-operad by letting $G$ act trivially on each translation category $\tn\Sigma_n$.  We refer to $\BE$-algebras in $\Gcat$ as \index{naive permutative G-category@naive permutative $G$-category}\index{G-category@$G$-category!naive permutative}\emph{naive permutative $G$-categories}.  
\end{definition}

Since the $\Gcat$-operad $\BE$ is reduced and levelwise a translation category, \cref{translation_pseudocom} applies to $\BE$.

\begin{corollary}\label{BE_pseudocom}
The $\Gcat$-operad $\BE$ admits a unique pseudo-commutative structure.
\end{corollary}

\begin{explanation}[Naive Permutative $G$-Categories]\label{expl:naive_perm_Gcat}
We unravel the structure of a naive permutative $G$-category, which is, by \cref{def:enr-multicategory-functor}, a $\Gcat$-multifunctor 
\[\upga \cn \BE \to \Gcat.\]  
The unique object $* \in \BE$ is sent to a small $G$-category $\upga(*) = \C$.  

When we interpret the component $G$-functors of $\upga$, we use the fact that $G$ acts trivially on $\BE$.  The $G$-functor at level 0
\[\upga \cn \BE(0) = \tn\Sigma_0 \iso \boldone \to \Catg(\boldone, \C)\]
corresponds to a $G$-fixed object $\pu \in \C$.  The $G$-functor at level 2
\[\upga \cn \BE(2) = \tn\Sigma_2 \to \Catg(\C \times \C, \C)\]
sends the object $\id_2 \in \Sigma_2$ to a $G$-functor \cref{Gfunctor}
\[\upga(\id_2) = \oplus \cn  \C \times \C \to \C.\]
For the non-identity permutation $\tau \in \Sigma_2$, there is a unique non-identity isomorphism 
\[ [\tau,\id_2] \cn \id_2 \fto{\iso} \tau \inspace \tn\Sigma_2.\]
Its image under $\upga$ is a $G$-natural isomorphism $\xi$ \cref{Gnattr} as follows, with $\twist$ swapping the two arguments.  
\[\begin{tikzpicture}[xscale=3, yscale=1.1]
\def\v{-1} \def\h{1} \def\m{1} \def\q{15}
\draw[0cell=1] 
(0,0) node (x11) {\C \times \C}
(x11)++(\h,0) node (x12) {\C}
(x11)++(\h/2,\v) node (x2) {\C \times \C}
;
\draw[1cell] 
(x11) edge node [pos=.45](i) {\oplus} (x12)
(x11) edge[bend right=\q] node[swap,pos=.5] (a) {\twist} (x2)
(x2) edge[bend right=\q] node[swap,pos=.5] {\oplus} (x12)
;
\draw[2cell] 
node[between=i and x2 at .5, shift={(0,0)}, rotate=-90, 2label={above,\xi}] {\Rightarrow}
;
\end{tikzpicture}\]
The data $(\C,\oplus,\pu,\xi)$ is a permutative category (\cref{def:symmoncat}) because the necessary axioms already hold in $\BE$.  By the coherence theorem for $\BE$ \cite[11.4.14]{cerberusIII}, the $\Gcat$-multifunctor $\upga$ is completely determined by the structure above.  

In summary, a naive permutative $G$-category is precisely a small $G$-category $\C$ together with a permutative category structure $(\oplus,\pu,\xi)$ that is $G$-equivariant.  The last part means that the monoidal product $\oplus$ and the braiding $\xi$ are $G$-equivariant (\cref{expl:GCat}), and the monoidal unit $\pu$ is a $G$-fixed object.
\end{explanation}

\section{Genuine Permutative $G$-Categories}
\label{sec:GBE}

This section discusses a $G$-equivariant analogue of the Barratt-Eccles operad, denoted $\GBE$, whose pseudoalgebras are genuine symmetric monoidal $G$-categories \pcref{sec:genuine_smgcat}.

\secoutline
\begin{itemize}
\item \cref{def:GBE} defines the $G$-Barratt-Eccles operad $\GBE$.  It is reduced and levelwise a translation category (\cref{GBEn}).  Thus, it is a pseudo-commutative operad in a unique way (\cref{GBE_pseudocom}).
\item \cref{expl:GBE} unpacks $\GBE$ in detail.
\item \cref{def:GBE_algebra,expl:GBE_algebra} discuss $\GBE$-algebras in $\Gcat$, which are called genuine permutative $G$-categories.  
\item \cref{naive_genuine_pGcat,expl:naive_genuine_pGcat} discuss the fact that each naive permutative $G$-category (\cref{expl:naive_perm_Gcat}) naturally yields a genuine permutative $G$-category via the functor $\Catg(\EG,-)$.
\end{itemize}

\subsection*{$G$-Barratt-Eccles Operad}

The following operad was introduced in \cite{shimakawa89} and studied further in \cite{gm17,gmm17}.

\begin{definition}\label{def:GBE}
Suppose $G$ is a group.  The \emph{$G$-Barratt-Eccles operad}\index{G-Barratt-Eccles operad@$G$-Barratt-Eccles operad}\index{operad!G-Barratt-Eccles@$G$-Barratt-Eccles} $\GBE$ is the reduced $\Gcat$-operad consisting of the small $G$-categories
\[\GBE(n) = \Catg(\EG, \BE(n)) = \Catg(\EG, \tn \Sigma_n) \forspace n \geq 0.\]
Here $\Catg(-,-)$ denotes the small $G$-category in \cref{def:Catg}, consisting of functors and natural transformations with the conjugation $G$-action.
\begin{itemize}
\item The translation $G$-category $\EG$ is the one in \cref{ex:EG}.  
\item The translation category $\tn \Sigma_n$ has the trivial $G$-action \pcref{def:BE-Gcat}. The right $\Sigma_n$-action on $\tn \Sigma_n$ is given by the group multiplication 
\[\Sigma_n \times \Sigma_n \to \Sigma_n\]
on objects, which extends uniquely to morphisms in $\tn \Sigma_n$. 
\item The $\Gcat$-operadic unit of $\GBE$ is given by the unique $G$-functor 
\[\opu \cn \EG \to \tn\Sigma_1 \iso \boldone.\] 
\item The right $\Sigma_n$-action on $\GBE(n)$ is induced by the one on $\tn\Sigma_n$.
\item The $\Gcat$-operadic multiplication of $\GBE$ is induced by the $\Cat$-operadic multiplication of the Barratt-Eccles operad $\BE$ (\cref{def:BE}) and the fact that $\Catg(\EG,-)$ preserves products, since it is a right adjoint \cref{Catg_rightadj}.
\end{itemize}
This finishes the definition of $\GBE$.
\end{definition}

To understand the $G$-Barratt-Eccles operad, we first observe that it is levelwise a translation category.  In \cref{expl:GBE} below, we unpack $\GBE$ further.

\begin{proposition}\label{GBEn}
For each $n \geq 0$, there is an isomorphism of $G$-categories
\[\GBE(n) \iso \tn[G,\Sigma_n]\]
where $[G,\Sigma_n]$ denotes the $G$-set of functions $G \to \Sigma_n$, with $G$ acting trivially on $\Sigma_n$.
\end{proposition}

\begin{proof}
To be precise, for $g \in G$ and a function $f \cn G \to \Sigma_n$, the function $g \cdot f$ is given by
\begin{equation}\label{gactionGSigman}
(g \cdot f)(k) = f\big(g^{-1} k\big) \forspace k \in G.
\end{equation}
This is the conjugation $G$-action on $[G,\Sigma_n$], since $G$ acts trivially on $\Sigma_n$.

Now we observe that there is an isomorphism of categories between $\GBE(n)$ and $\tn[G,\Sigma_n]$.  Since $\tn\Sigma_n$ is a translation category, a functor $\EG \to \tn\Sigma_n$ is uniquely determined by its object assignment, which is a function $G \to \Sigma_n$ that is not necessarily a group homomorphism.  Conversely, any function $G \to \Sigma_n$ is the object assignment of a unique functor $\EG \to \tn\Sigma_n$.  For the same reason, given any two functors 
\[F,H \cn \EG \to \tn\Sigma_n,\] 
there is a unique natural transformation $\theta \cn F \to H$, whose $g$-component must be the unique morphism 
\[\theta_g \cn Fg \to Hg \inspace \tn\Sigma_n\]
for each $g \in G$.  This shows that $\GBE(n)$ is isomorphic to the translation category $\tn[G,\Sigma_n]$.  

Moreover, the isomorphism of categories in the previous paragraph is a $G$-functor \cref{Gfunctor}.  Indeed, in both $\Catg(\EG,\tn\Sigma_n)$ and $\tn[G,\Sigma_n]$, the $G$-action comes from the conjugation $G$-action, with $G$ acting trivially on $\Sigma_n$ and $\tn\Sigma_n$.
\end{proof}

By \cref{GBEn}, the $G$-Barratt-Eccles operad $\GBE$ is reduced and levelwise a translation category.  Thus, \cref{translation_pseudocom} applies to $\GBE$.

\begin{corollary}\label{GBE_pseudocom}
The $\Gcat$-operad $\GBE$ admits a unique pseudo-commutative structure.
\end{corollary}

\begin{explanation}[$G$-Barratt-Eccles Operad]\label{expl:GBE}
We unravel the $\Gcat$-operad 
\[\GBE = \big\{\Catg(\EG,\tn\Sigma_n)\big\}_{n \geq 0}\]
in \cref{def:GBE}.
\begin{description}
\item[$G$-action] The $G$-action on a functor $F \cn \EG \to \tn\Sigma_n$ is given by
\begin{equation}\label{gFk}
(g \cdot F)(k) = F\big(g^{-1} k\big) \forspace g,k \in G, 
\end{equation}
as in \cref{gactionGSigman}.  Consider another functor $H \cn \EG \to \tn\Sigma_n$ and the unique natural transformation $\theta \cn F \to H$ in $\GBE(n)$.  The $g$-action on $\theta$ yields the natural transformation
\[g \cdot \theta \cn g \cdot F \to g \cdot H\]
whose $k$-component is the unique morphism 
\[(g \cdot \theta)_k = \big[H(g^{-1} k), F(g^{-1} k)\big] \cn F(g^{-1} k) \to H(g^{-1} k) \inspace \tn\Sigma_n.\]

\item[$\Sigma_n$-action] The right $\Sigma_n$-action on $\GBE(n)$ is induced by the one on $\tn\Sigma_n$.  For a functor $F \cn \EG \to \tn\Sigma_n$ and a permutation $\sigma \in \Sigma_n$, the right $\sigma$-action is given on objects by the group multiplication in $\Sigma_n$,
\begin{equation}\label{Fsigmak}
F^\sigma(k) = Fk \cdot \sigma \forspace k \in G.
\end{equation}
By \cref{GBEn}, this object assignment uniquely determines $F^\sigma$ on morphisms:
\[F^\sigma[k,g] = \big[Fk \cdot \sigma, Fg \cdot \sigma\big] 
\cn Fg \cdot \sigma \to Fk \cdot \sigma \inspace \tn\Sigma_n\]
for $k,g \in G$.  By \cref{gFk,Fsigmak}, the $G$-action and the right $\Sigma_n$-action on $\GBE(n)$ commute, since both $g \cdot F^\sigma$ and $(g \cdot F)^\sigma$ are given by
\[(k \in G) \mapsto \big(F(g^{-1} k)\big) \cdot \sigma.\]

For a natural transformation $\theta \cn F \to H$ in $\GBE(n)$, the right $\sigma$-action on $\theta$ yields the natural transformation 
\[\theta^\sigma \cn F^\sigma \to H^\sigma\]
whose $g$-component is the unique morphism
\[\theta^\sigma_g = \big[Hg \cdot \sigma, Fg \cdot \sigma\big] 
\cn Fg \cdot \sigma \to Hg \cdot \sigma \inspace \tn\Sigma_n.\]

\item[Operadic composition] We denote by $\ga$ the composition of $\BE$ and by $\ga^G$ the composition of $\GBE$.  For $n, p_1, \ldots, p_n \geq 0$, the composition $\ga^G$ is the following composite $G$-functor, where $p= \sum_{j=1}^n p_j$.
\begin{equation}\label{GBE_gamma}
\begin{tikzpicture}[vcenter]
\def\h{5.2} \def\v{-1.4}
\draw[0cell=.85]
(0,0) node (x11) {\GBE(n) \times \GBE(p_1) \times \cdots \times \GBE(p_n)}
(x11)++(\h,0) node (x12) {\GBE(p) = \Catg(\EG, \tn\Sigma_p)}
(x11)++(0,\v) node (x21) {\Catg(\EG,\tn\Sigma_n) \times \txprod_{j=1}^n \Catg(\EG,\tn\Sigma_{p_j})}
(x12)++(0,\v) node (x22) {\Catg\big(\EG, \tn\Sigma_n \times \txprod_{j=1}^n \tn\Sigma_{p_j} \big)}
;
\draw[1cell=.85]  
(x11) edge node {\ga^G} (x12)
(x11) edge[equal] (x21)
(x21) edge node {\iso} (x22)
(x22) edge[shorten <=-2pt] node[swap] {\Catg(1,\ga)} (x12)
;
\end{tikzpicture}
\end{equation}
More explicitly, consider functors 
\begin{equation}\label{functorsFHj}
\EG \fto{F} \tn\Sigma_n \andspace \EG \fto{H_j} \tn\Sigma_{p_j}
\end{equation}
for $1 \leq j \leq n$ and $g \in G$.  Then we have
\begin{equation}\label{gaGFHjg}
\begin{split}
&\ga^G\big(F \sscs H_1, \ldots, H_n\big)(g)\\
&= \ga\big(Fg \sscs H_1 g, \ldots, H_n g\big)\\
&= (Fg)\ang{p_1,\ldots,p_n} \circ \big(H_1 g \times \cdots \times H_n g\big) \inspace \Sigma_p.
\end{split}
\end{equation}
The last entry in \cref{gaGFHjg} is the composite permutation in \cref{as_gamma}.  Moreover, each of
\begin{itemize}
\item the value of $\ga^G\big(F \sscs H_1, \ldots, H_n\big)$ at a morphism in $\EG$ and
\item the value of $\ga^G$ at a morphism in $\GBE(n) \times \prod_{j=1}^n \GBE(p_j)$
\end{itemize}
is uniquely determined by the fact that $\tn\Sigma_p$ is a translation category.\defmark
\end{description}
\end{explanation}

Basically the same proof for \cref{GBEn} also proves the following result, which we record here for future reference.

\begin{proposition}\label{Catgceg}
Suppose $\C$ is a small $G$-category for a group $G$.  Then there is an isomorphism of $G$-categories
\[\Catg(\C,\EG) \iso \tn[\Obc,G],\]
where $[\Obc,G]$ denotes the $G$-set of functions $\Obc \to G$, with $G$ acting by conjugation.
\end{proposition}

\begin{proof}
Since $\EG$ is a translation category, a functor $F \cn \C \to \EG$ is uniquely determined by its object assignment $\Ob(F) \cn \Obc \to G$.  Conversely, each function $\Obc \to G$ is the object assignment of a unique functor $\C \to \EG$.  This bijection between the object sets of $\Catg(\C,\EG)$ and $\tn[\Obc,G]$ respects the $G$-action, since it is given by the conjugation $G$-action in each case.  This $G$-bijection for objects extends uniquely to morphisms because, given any two functors $F,H \cn \C \to \EG$, there is a unique natural transformation $F \to H$.
\end{proof}

\subsection*{Genuine Permutative $G$-Categories}

Recall from \cref{def:BE-Gcat} that naive permutative $G$-categories are algebras over the Barratt-Eccles operad $\BE$ in $\Gcat$.  The following genuine variant was introduced in \cite{shimakawa89} and given its current name in \cite[Definition 4.5]{gm17}.

\begin{definition}\label{def:GBE_algebra}
For the reduced $\Gcat$-operad $\GBE$ in \cref{def:GBE}, a $\GBE$-algebra in $\Gcat$ is called a \index{genuine permutative G-category@genuine permutative $G$-category}\index{G-category@$G$-category!genuine permutative}\emph{genuine permutative $G$-category}.
\end{definition}

Recall the notation $\ufs{n} = \{1,\ldots,n\}$ from \cref{ufsn}.  We use the notation
\begin{equation}\label{angaj_ufsn}
\ang{a_j}_{j \in \ufs{n}} = (a_1, a_2, \ldots , a_n)
\end{equation}
and $\ang{a_j}_{j=1}^n$ interchangeably.

\begin{explanation}[Genuine Permutative $G$-Categories]\label{expl:GBE_algebra}
Using \cref{GBEn}, we unpack the structure of a genuine permutative $G$-category, which is, by \cref{def:enr-multicategory-functor}, a $\Gcat$-multifunctor
\[\upga \cn \GBE \to \Gcat.\]
The unique object $* \in \GBE$ is sent to a small $G$-category $\upga(*) = \C$.

For each $n \geq 0$, the component $G$-functor
\[\upga_n \cn \GBE(n) = \Catg(\EG,\ESigma_n) \to \Catg(\C^n,\C)\]
is adjoint to a $G$-functor as follows, which we also denote by $\upga_n$ to simplify the notation.
\begin{equation}\label{gan_GBE_alg}
\upga_n \cn \GBE(n) \times \C^n \to \C
\end{equation}
\begin{itemize}
\item Since $\GBE(0) \iso \boldone$, $\upga_0$ is determined by a $G$-fixed object
\[\upga_0(*) \in \C,\]
which is called the \emph{basepoint} of $\C$.
\item Since $\GBE(1) \iso \boldone$, by the unity axiom \cref{enr-multifunctor-unit}, $\upga_1$ is determined by the identity functor of $\C$.
\end{itemize}

\begin{description}
\item[$G$-equivariance] For $n \geq 2$, an object in $\GBE(n)$ is a functor $F \cn \EG \to \tn\Sigma_n$, which is determined by its object assignment $F \cn G \to \Sigma_n$.  In the domain of $\upga_n$, $G$ acts diagonally.  Thus, $G$-equivariance of $\upga_n$ means the following equality in $\C$ for $g \in G$ and objects $\ang{a_j}_{j\in \ufs{n}} \in \C^n$.
\begin{equation}\label{mun_Gequivariant}
\ga_n\left(g \cdot F \sscs \ang{g \cdot a_j}_{j\in \ufs{n}}\right) 
= g \cdot \ga_n\left(F \sscs \ang{a_j}_{j\in \ufs{n}}\right)
\end{equation}
There is a similar equality for natural transformations $\theta \cn F \to H$ and morphisms in $\C^n$.  On the left-hand side of \cref{mun_Gequivariant}, $g \cdot F$ is the functor in \cref{gFk}.

\item[$\Sigma_n$-equivariance]
The symmetry axiom \cref{enr-multifunctor-equivariance} for $\upga$ means that, for each functor $F \cn \EG \to \tn\Sigma_n$, permutation $\sigma \in \Sigma_n$, and objects $\ang{a_j}_{j\in \ufs{n}} \in \C^n$, the following equality holds in $\C$.
\begin{equation}\label{mun_symmetry}
\upga_n\left(F^\sigma \sscs \ang{a_j}_{j\in \ufs{n}}\right) 
= \upga_n\big(F \sscs \ang{a_{\sigmainv(j)}}_{j\in \ufs{n}}\big)
\end{equation}
There is a similar equality for natural transformations $F \to H$ and morphisms in $\C^n$.  On the left-hand side of \cref{mun_symmetry}, $F^\sigma$ is the functor in \cref{Fsigmak}.

\item[Composition]
The composition axiom \cref{v-multifunctor-composition} for $\upga$ means that, for functors $F$ and $H_j$ as defined in \cref{functorsFHj} and objects $a_{j,i} \in \C$ for $1 \leq j \leq n$ and $1 \leq i \leq p_j$, the following equality holds in $\C$, where $p = \sum_{j=1}^n p_j$.
\begin{equation}\label{mun_composition}
\begin{split}
&\upga_p\big( \ga^G\big( F \sscs \ang{H_j}_{j\in \ufs{n}} \big) \sscs \ang{\ang{a_{j,i}}_{i \in \ufs{p}_j}}_{j\in \ufs{n}} \big)\\
&= \upga_n\big( F \sscs \big\langle\upga_{p_j} \big( H_j \sscs \ang{a_{j,i}}_{i\in \ufs{p}_j} \big) \big\rangle_{j\in \ufs{n}} \big)
\end{split}
\end{equation}
There is a similar equality for natural transformations and morphisms.  In \cref{mun_composition}, $\ga^G\big(F \sscs \ang{H_j}_{j\in \ufs{n}}\big)$ is the functor in \cref{gaGFHjg}.  If $p_j=0$ for any $j$, then the object $\upga_{p_j} \big( H_j \sscs \ang{a_{j,i}}_{i \in \ufs{p}_j} \big)$ in \cref{mun_composition} is interpreted as the $G$-fixed object $\upga_0(*) \in \C$.\defmark
\end{description}
\end{explanation}

\subsection*{Naive to Genuine Permutative $G$-Categories}

The functor $\Catg(\EG,-)$ on $\Gcat$ preserves small limits, including products, because it is a right adjoint \cref{Catg_rightadj}.  Since the $G$-Barratt-Eccles operad $\GBE$ (\cref{def:GBE}) is obtained from the Barratt-Eccles operad $\BE$ (\cref{def:BE}) by applying the functor $\Catg(\EG,-)$, we have the following observation from \cite[page 256]{shimakawa89} that relates naive and genuine permutative $G$-categories.  See also \cite[Proposition 4.6]{gm17}.

\begin{proposition}\label{naive_genuine_pGcat}
For each group $G$, each naive permutative $G$-category $\C$ yields a genuine permutative $G$-category $\Catg(\EG,\C)$.
\end{proposition}

\begin{explanation}[Naive to Genuine]\label{expl:naive_genuine_pGcat}
\cref{naive_genuine_pGcat} states that, if $\C$ is a $\BE$-algebra in $\Gcat$, then $\Catg(\EG,\C)$ is a $\GBE$-algebra.  Here we explain the details of this construction.  

The $\BE$-algebra structure on $\C$ is given by a $\Gcat$-multifunctor
\[\upga \cn \BE \to \Gcat.\]
As we discuss in \cref{expl:naive_perm_Gcat}, this means that $\C$ is a small $G$-category together with a $G$-equivariant permutative category structure $(\oplus,\pu,\xi)$.  By \cref{def:Catg}, $\Catg(\EG,\C)$ is the small category of functors $\EG \to \C$ and natural transformations, with $G$ acting by conjugation.

The $\GBE$-algebra structure on $\Catg(\EG,\C)$ is given by a $\Gcat$-multifunctor
\[\upgabar \cn \GBE \to \Gcat.\]
In the adjoint form \cref{gan_GBE_alg}, for each $n \geq 0$, the component $G$-functor $\upgabar_n$ is the following composite.
\begin{equation}\label{CatgEgC_GBE_alg}
\begin{tikzpicture}[xscale=5,yscale=1.4,vcenter]
\draw[0cell=.9]
(0,0) node (x11) {\GBE(n) \times \Catg(\EG,\C)^n}
(x11)++(0,-1) node (x21) {\Catg(\EG,\tn\Sigma_n) \times \Catg(\EG,\C)^n}
(x21)++(1,0) node (x22) {\Catg\big(\EG, \tn\Sigma_n \times \C^n\big)}
(x22)++(0,1) node (x12) {\Catg(\EG,\C)}
;
\draw[1cell=.9]  
(x11) edge node {\upgabar_n} (x12)
(x11) edge[equal] (x21)
(x21) edge node {\iso} (x22)
(x22) edge[transform canvas={xshift={-1em}}] node[swap] {\Catg(1,\upga_n)} (x12)
;
\end{tikzpicture}
\end{equation}
In the right vertical arrow, with $G$ acting trivially on $\ESigma_n$, the $G$-functor
\[\upga_n \cn \ESigma_n \times \C^n \to \C\]
is the adjoint of the component $G$-functor
\[\upga_n \cn \BE(n) = \ESigma_n \to \Catg(\C^n,\C).\]
More explicitly, for functors
\[\EG \fto{F} \tn\Sigma_n \andspace \EG \fto{H_j} \C\]
for $1 \leq j \leq n$ and $g \in G$, we have
\[\upgabar_n\left(F \sscs \ang{H_j}_{j\in \ufs{n}}\right)(g) 
= \upga_n\left(Fg \sscs \ang{H_j g}_{j\in \ufs{n}}\right)\]
in $\C$.
\end{explanation}

\begin{remark}\label{rk:naive_genuine}
As stated in \cite[page 3300]{gm17}, all known examples of genuine permutative $G$-categories are of the form $\Catg(\EG,\C)$ for some naive permutative $G$-categories $\C$.  The work of \cite[Theorem A]{lenz-genuine} claims that, in fact, $\Catg(\EG,-)$ is an equivalence of quasi-categories.  Since we will not use that result in this work, we refer the interested reader to \cite{lenz-genuine} for details.
\end{remark}

\section{Pseudoalgebras of an Operad}
\label{sec:pseudoalgebra}

This section discusses pseudoalgebras of an operad, along with appropriate notions of lax morphisms and transformations.

\secoutline
\begin{itemize}
\item \cref{def:pseudoalgebra} defines $\Op$-pseudoalgebras for a reduced $\Gcat$-operad $\Op$, with $G$ an arbitrary group.  Recall that \emph{reduced} means $\Op(0) = \boldone$, the terminal $G$-category.
\item As a recipe for generating new pseudoalgebras from given ones, \cref{catgego} shows that each $\Op$-pseudoalgebra $\A$ yields a $\Catg(\EG,\Op)$-pseudoalgebra $\Catg(\EG,\A)$.  In \cref{naive_genuine_smgcat}, we use \cref{catgego} to generate genuine symmetric monoidal $G$-categories from $\BE$-pseudoalgebras. 
\item \cref{phi_id} shows that the axioms of an $\Op$-pseudoalgebra implies additional unity properties.
\item \cref{def:laxmorphism} defines lax $\Op$-morphisms, along with their pseudo and strict variants.
\item \cref{def:algtwocells} defines $\Op$-transformations.
\item \cref{oalgps_twocat} shows that there is a 2-category $\AlglaxO$ of $\Op$-pseudoalgebras, lax $\Op$-morphisms, and $\Op$-transformations.  There are also variant 2-categories, denoted $\AlgpspsO$ and $\AlgstO$, with the same objects and 2-cells as $\AlglaxO$, and with 1-cells given by, respectively, $\Op$-pseudomorphisms and strict $\Op$-morphisms. 
\end{itemize}

\subsection*{Looking Ahead}
To understand the definitions and constructions in this section, in \cref{sec:naive_smc,sec:BEpseudoalg}, we discuss $\AlglaxBE$ for the Barratt-Eccles $\Gcat$-operad in \cref{def:BE-Gcat}.  We show in \cref{thm:BEpseudoalg} that $\AlglaxBE$ is 2-equivalent to a 2-category $\smgcat$ (\cref{def:smGcat_twocat}) with naive symmetric monoidal $G$-categories as objects, strictly unital symmetric monoidal $G$-functors as 1-cells, and monoidal $G$-natural transformations as 2-cells.

Moreover, \cref{thm:multpso} shows that, if $\Op$ is pseudo-commutative, then there is a $\Gcat$-multicategory $\MultpsO$ with $\Op$-pseudoalgebras as objects, lax $\Op$-morphisms as $G$-fixed 1-ary 1-cells, and $\Op$-transformations as $G$-fixed 1-ary 2-cells.  However, pseudo-commutativity is not needed in this section and \cref{sec:naive_smc,sec:BEpseudoalg}.

\subsection*{Pseudoalgebras}

Recall the Cartesian closed 2-category $\Gcat$ \pcref{expl:GCat,expl:Gcat_closed}.  The next definition is \cite[Definition 2.14]{gmmo20}.

\begin{definition}\label{def:pseudoalgebra}
Suppose $G$ is a group, and $(\Op,\ga,\opu)$ is a reduced $\Gcat$-operad.  An \emph{$\Op$-pseudoalgebra}\index{pseudoalgebra}\index{operad!pseudoalgebra} is a triple $(\A,\gaA,\phiA)$ consisting of the following data.
\begin{description}
\item[Underlying $G$-category] $\A$ is a small $G$-category.
\item[$\Op$-action] For each $n \geq 0$, it is equipped with a $G$-functor
\begin{equation}\label{gaAn}
\Op(n) \times \A^n \fto{\gaA_n} \A,
\end{equation}
called the \index{pseudoalgebra!action G-functor@action $G$-functor}\index{action G-functor@action $G$-functor}\emph{$n$-th $\Op$-action $G$-functor}.  When $n=0$, the $G$-functor
\[\Op(0) = \boldone \fto{\gaA_0} \A\]
specifies a $G$-fixed object (\cref{def:fixedpoint})
\begin{equation}\label{pseudoalg_zero}
\zero = \gaA_0(*) \in \A,
\end{equation} 
called the \emph{basepoint}\index{pseudoalgebra!basepoint}\index{basepoint} of $\A$.  We sometimes abbreviate $\gaA_n$ to $\gaA$.
\item[Associativity constraint] For each sequence 
\[\big(n; \ang{m_j}_{j \in \ufs{n}}\big) = \big(n; m_1,\ldots,m_n\big)\] 
with $n > 0$ and each $m_j \geq 0$, it is equipped with a $G$-natural isomorphism, called the \emph{associativity constraint}\index{pseudoalgebra!associativity constraint}\index{associativity constraint} of $\A$, as follows, where $m = \sum_{j=1}^n m_j$.
\begin{equation}\label{phiA}
\begin{tikzpicture}[xscale=1,yscale=1,vcenter]
\def\h{5.5} \def\g{.5} \def\v{-1.3}
\draw[0cell=.9]
(0,0) node (x11) {\Op(n) \times \txprod_{j=1}^n \big( \Op(m_j) \times \A^{m_j} \big)}
(x11)++(\h,0) node (x12) {\Op(n) \times \A^n}
(x12)++(\g,\v) node (x) {\A}
(x11)++(0,2*\v) node (x21) {\big(\Op(n) \times \txprod_{j=1}^n \Op(m_j) \big) \times \A^{m_1 + \cdots + m_n}}
(x21)++(\h,0) node (x22) {\Op(m) \times \A^m}
;
\draw[1cell=.9]  
(x11) edge node {1 \times \txprod_{j=1}^n\, \gaA_{m_j}} (x12)
(x12) edge node[pos=.7] {\gaA_n} (x)
(x11) edge node[swap] {\pi} (x21)
(x21) edge node {\ga \times 1} (x22)
(x22) edge node[swap,pos=.8] {\gaA_m} (x)
;
\draw[2cell]
node [between=x11 and x22 at .5, shift={(0,.2)}, rotate=-90, 2label={below,\iso}, 2label={above,\phiA_{(n;\, m_1,\ldots,m_n)}}] {\Rightarrow}
;
\end{tikzpicture}
\end{equation}
The arrow $\pi$ shuffles the $\Op(m_j)$ factors to the left, keeping their relative order unchanged.  We sometimes abbreviate $\phiA_{(n;\, m_1,\ldots,m_n)} = \phiA_{(n;\, \ang{m_j}_{j \in \ufs{n}})}$ to $\phiA$.  A typical component of $\phiA_{(n;\, \ang{m_j}_{j \in \ufs{n}})}$ is an isomorphism
\begin{equation}\label{phiA_component}
\begin{tikzpicture}[xscale=5.8,yscale=1,baseline={(a.base)}]
\draw[0cell=.8]
(0,0) node (a) {\gaA_n\Big(x \sscs \bang{ \gaA_{m_j}\big(x_j \sscs \ang{a_{j,i}}_{i\in \ufs{m}_j} \big) }_{j\in \ufs{n}}\Big)}
(a)++(1,0) node (b) {\gaA_m\Big( \ga\big( x \sscs \ang{x_j}_{j\in \ufs{n}} \big) \sscs \bang{\ang{a_{j,i}}_{i\in \ufs{m}_j} }_{j\in \ufs{n}} \Big)}
;
\draw[1cell=.8]
(a) edge node {\phiA_{(n;\, \ang{m_j}_{j \in \ufs{n}})}} node[swap] {\iso} (b)
;
\end{tikzpicture}
\end{equation}
in $\A$ for objects $x \in \Op(n)$, $x_j \in \Op(m_j)$, and $a_{j,i} \in \A$ for $1 \leq j \leq n$ and $1 \leq i \leq m_j$.  By convention, $\phiA_{(0;\ang{})}$ is $1_\zero$, the identity morphism of the basepoint $\zero \in \A$.
\end{description}
The data above are subject to the axioms  \cref{pseudoalg_action_sym,pseudoalg_action_unity,pseudoalg_basept_axiom,pseudoalg_topeq,pseudoalg_boteq,pseudoalg_unity,pseudoalg_comp_axiom} below.
\begin{description}
\item[Action equivariance] 
For each permutation $\sigma \in \Sigma_n$ with $n \geq 2$, the following diagram of $G$-functors commutes. 
\begin{equation}\label{pseudoalg_action_sym}
\begin{tikzpicture}[xscale=3,yscale=1.3,vcenter]
\draw[0cell=.9]
(0,0) node (x11) {\Op(n) \times \A^n}
(x11)++(1,0) node (x12) {\Op(n) \times \A^n}
(x11)++(0,-1) node (x21) {\Op(n) \times \A^n}
(x21)++(1,0) node (x22) {\A}
;
\draw[1cell=.9]  
(x11) edge node {1 \times \sigma} (x12)
(x12) edge node {\gaA_n} (x22)
(x11) edge node[swap] {\sigma \times 1} (x21)
(x21) edge node {\gaA_n} (x22)
;
\end{tikzpicture}
\end{equation}
In the top horizontal arrow, $\sigma$ permutes the $n$ copies of $\A$ from the left:
\[\sigma(a_1,a_2,\ldots,a_n) = \big(a_{\sigmainv(1)}, a_{\sigmainv(2)}, \ldots, a_{\sigmainv(n)}\big).\]
In the left vertical arrow, $\sigma$ is the right $\sigma$-action \cref{rightsigmaaction} on $\Op(n)$.
\item[Action unity] 
The following diagram commutes. 
\begin{equation}\label{pseudoalg_action_unity}
\begin{tikzpicture}[xscale=2.5,yscale=1.3,vcenter]
\draw[0cell=.9]
(0,0) node (x11) {\boldone \times \A}
(x11)++(1,0) node (x12) {\A}
(x11)++(0,-1) node (x21) {\Op(1) \times \A}
(x21)++(1,0) node (x22) {\A}
;
\draw[1cell=.9]  
(x11) edge node {\iso} (x12)
(x12) edge[equal] (x22)
(x11) edge node[swap] {\opu \times 1} (x21)
(x21) edge node {\gaA_1} (x22)
;
\end{tikzpicture}
\end{equation}
The remaining axioms involve the associativity constraint.
\item[Basepoint] 
Consider objects $x \in \Op(n)$ and 
\[\anga = \big(a_1,\ldots,a_{j-1}, a_{j+1}, \ldots, a_n) \in \A^{n-1}\] 
for any $j \in \{1,\ldots,n\}$.  We denote the objects\label{not:dyj}
\[\begin{split}
\dy_j x &= \ga\big(x \sscs \opu^{j-1}, *, \opu^{n-j}\big) \in \Op(n-1) \andspace\\
\dy_j \anga &= \big(a_1,\ldots,a_{j-1}, \zero, a_{j+1}, \ldots, a_n) \in \A^n
\end{split}\]
with $* \in \Op(0) = \boldone$ the unique object, $\opu \in \Op(1)$ the operadic unit, $\opu^k$ the $k$-tuple of copies of $\opu$, and $\zero = \gaA_0(*)$ the basepoint of $\A$ \cref{pseudoalg_zero}.
Then the following component of $\phiA$ is the identity morphism.
\begin{equation}\label{pseudoalg_basept_axiom} 
\gaA_n\big( x \sscs \dy_j\anga \big) \fto{\phiA_{(n;\, 1^{j-1},0,1^{n-j})} \,=\, 1} 
\gaA_{n-1}\big( \dy_j x \sscs \anga \big)
\end{equation}
In \cref{pseudoalg_basept_axiom}, we use the component of $\phiA$ \cref{phiA_component} with $m_j = 0$, $x_j = * \in \Op(0)$, $m_r = 1$ for $r \neq j$, and $x_r = \opu \in \Op(1)$.  In the domain of $\phiA$, we use the equality
\[\gaA_n\big( x \sscs \dy_j\anga \big) = 
\gaA_n\left( x \sscs \bang{\gaA_1(\opu; a_i)}_{i=1}^{j-1} \scs \gaA_0(*) \scs \bang{\gaA_1(\opu; a_i)}_{i=j+1}^n \right)\]
in $\A$, which holds by the action unity axiom \cref{pseudoalg_action_unity}.

\item[Top equivariance] 
For each permutation $\sigma \in \Sigma_n$ and objects $x \in \Op(n)$, $x_j \in \Op(m_j)$, and $a_{j,i} \in \A$ as in \cref{phiA_component}, the following diagram in $\A$ commutes.
\begin{equation}\label{pseudoalg_topeq}
\begin{tikzpicture}[xscale=1,yscale=1,vcenter]
\def\h{3} \def\g{0} \def\v{1.4} \def\u{-1}
\draw[0cell=.8]
(0,0) node (x11) {\gaA_n\Big(x\sigma \sscs \big\langle \gaA_{m_{\sigma(j)}}\big(x_{\sigma(j)} \sscs \ang{a_{\sigma(j), i}}_{i \in \ufs{m}_{\sigma(j)}} \big) \big\rangle_{j\in \ufs{n}}\Big)}
(x11)++(\h,\v) node (x12) {\gaA_m\left( \ga\big( x\sigma \sscs \ang{x_{\sigma(j)}}_{j\in \ufs{n}} \big) \sscs \bang{\ang{a_{\sigma(j),i}}_{i \in \ufs{m}_{\sigma(j)}} }_{j\in \ufs{n}} \right)}
(x11)++(\h+.5,\u) node (x) {\gaA_m\left( \ga\big( x \sscs \ang{x_j}_{j\in \ufs{n}} \big) \sigmabar \sscs \bang{\ang{a_{\sigma(j),i}}_{i \in \ufs{m}_{\sigma(j)}} }_{j\in \ufs{n}} \right)}
(x11)++(0,2*\u) node (x21) {\gaA_n\Big(x \sscs \big\langle \gaA_{m_j}\big(x_j \sscs \ang{a_{j,i}}_{i \in \ufs{m}_j} \big) \big\rangle_{j\in \ufs{n}}\Big)}
(x21)++(\h,-\v) node (x22) {\gaA_m\left( \ga\big( x \sscs \ang{x_j}_{j\in \ufs{n}} \big) \sscs \bang{\ang{a_{j,i}}_{i \in \ufs{m}_j} }_{j\in \ufs{n}} \right)}
;
\draw[1cell=.85]  
(x11) edge[transform canvas={xshift={-2em}}] node[pos=.2] {\phiA_{(n;\, m_{\sigma(1)}, \ldots, m_{\sigma(n)})}} (x12)
(x12) edge[equal] node[pos=.6] {\spadesuit} (x)
(x) edge[equal] (x22)
(x11) edge[equal, transform canvas={xshift={-1em}}](x21)
(x21) edge[transform canvas={xshift={-2em}}] node[swap,pos=.2] {\phiA_{(n;\, m_1,\ldots,m_n)}} (x22)
;
\end{tikzpicture}
\end{equation}
In \cref{pseudoalg_topeq}, the equality labeled $\spadesuit$ holds by the top equivariance axiom \cref{enr-operadic-eq-1} for $\Op$, where
\[\sigmabar = \sigma\bang{m_{\sigma(1)}, \ldots, m_{\sigma(n)}} \in \Sigma_{m_1+\cdots+m_n}\]
is the block permutation \cref{blockpermutation} induced by $\sigma$ that permutes blocks of lengths $m_{\sigma(1)}, \ldots, m_{\sigma(n)}$.  The other two equalities hold by the action equivariance axiom \cref{pseudoalg_action_sym}.

\item[Bottom equivariance] 
For permutations $\tau_j \in \Sigma_{m_j}$ for $1 \leq j \leq n$, the following diagram in $\A$ commutes.
\begin{equation}\label{pseudoalg_boteq}
\begin{tikzpicture}[xscale=1,yscale=1,vcenter]
\def\h{3} \def\g{0} \def\v{1.4} \def\u{-1}
\draw[0cell=.8]
(0,0) node (x11) {\gaA_n\Big(x \sscs \big\langle \gaA_{m_j}\big(x_j \tau_j \sscs \ang{a_{j,\tau_j(i)}}_{i \in \ufs{m}_j} \big) \big\rangle_{j\in \ufs{n}}\Big)}
(x11)++(\h,\v) node (x12) {\gaA_m\left( \ga\big( x \sscs \ang{x_j \tau_j}_{j\in \ufs{n}} \big) \sscs \bang{\ang{a_{j,\tau_j(i)}}_{i \in \ufs{m}_j} }_{j\in \ufs{n}} \right)}
(x11)++(\h+.5,\u) node (x) {\gaA_m\left( \ga\big( x \sscs \ang{x_j}_{j\in \ufs{n}} \big) \tau^\times \sscs \bang{\ang{a_{j,\tau_j(i)}}_{i \in \ufs{m}_j} }_{j\in \ufs{n}} \right)}
(x11)++(0,2*\u) node (x21) {\gaA_n\Big(x \sscs \big\langle \gaA_{m_j}\big(x_j \sscs \ang{a_{j,i}}_{i \in \ufs{m}_j} \big) \big\rangle_{j\in \ufs{n}}\Big)}
(x21)++(\h,-\v) node (x22) {\gaA_m\left( \ga\big( x \sscs \ang{x_j}_{j\in \ufs{n}} \big) \sscs \bang{\ang{a_{j,i}}_{i \in \ufs{m}_j} }_{j\in \ufs{n}} \right)}
;
\draw[1cell=.85]  
(x11) edge[transform canvas={xshift={-2em}}] node[pos=.3] {\phiA_{(n;\, m_1, \ldots, m_n)}} (x12)
(x12) edge[equal] node[pos=.6] {\clubsuit} (x)
(x) edge[equal] (x22)
(x11) edge[equal, transform canvas={xshift={-1em}}](x21)
(x21) edge[transform canvas={xshift={-2em}}] node[swap,pos=.2] {\phiA_{(n;\, m_1,\ldots,m_n)}} (x22)
;
\end{tikzpicture}
\end{equation}
In \cref{pseudoalg_boteq}, the equality labeled $\clubsuit$ holds by the bottom equivariance axiom \cref{enr-operadic-eq-2} for $\Op$, where
\[\tau^\times = \tau_1 \times \cdots \times \tau_n\in \Sigma_{m_1+\cdots+m_n}\]
is the block sum \cref{blocksum}.  The other two equalities hold by the action equivariance axiom \cref{pseudoalg_action_sym}.

\item[Unity] 
For objects $x \in \Op(n)$ and $\anga = \ang{a_j}_{j\in \ufs{n}} \in \A^n$, the following diagram in $\A$ commutes, where $1^n = \ang{1}_{j \in \ufs{n}}$.
\begin{equation}\label{pseudoalg_unity}
\begin{tikzpicture}[xscale=5,yscale=1.2,vcenter]
\draw[0cell=.9]
(0,0) node (x11) {\gaA_n\Big( x \sscs \bang{ \gaA_1(\opu \sscs a_j) }_{j\in \ufs{n}} \Big)}
(x11)++(1,0) node (x12) {\gaA_n\Big( \ga(x \sscs \opu^n) \sscs \anga \Big)}
(x11)++(0,-1) node (x21) {\gaA_n\big(x \sscs \anga \big)}
(x21)++(1,0) node (x22) {\gaA_n\big(x \sscs \anga \big)}
(x21)++(0,-1) node (x31) {\gaA_1\Big( \opu \sscs \gaA_n\big( x \sscs \anga \big) \Big)} 
(x31)++(1,0) node (x32) {\gaA_n\Big( \ga( \opu \sscs x) \sscs \anga\Big)}
;
\draw[1cell=.9]  
(x11) edge node {\phiA_{(n;\, 1^n)}} (x12)
(x21) edge node {1} (x22)
(x31) edge node {\phiA_{(1;\, n)}} (x32)
(x11) edge[equal] (x21)
(x12) edge[equal] (x22)
(x21) edge[equal] (x31)
(x22) edge[equal] (x32)
;
\end{tikzpicture}
\end{equation}
In \cref{pseudoalg_unity}, the two equalities along the left boundary hold by the action unity axiom \cref{pseudoalg_action_unity}.  The two equalities along the right boundary hold by the unity axioms \cref{enr-multicategory-right-unity,enr-multicategory-left-unity} for $\Op$.

\item[Composition] 
Consider objects 
\[x \in \Op(n), \quad x_j \in \Op(m_j), \quad x_{j,i} \in \Op(k_{j,i}), \andspace a_{j,i,h} \in \A\]
for $1 \leq j \leq n$, $1 \leq i \leq m_j$, and $1 \leq h \leq k_{j,i}$, along with the following notation, where we use $\bdot\,$ to denote a running index in a finite sequence.
\[\left\{\scalebox{.9}{$
\begin{aligned}
m_{\crdot} &= \ang{m_j}_{j\in \ufs{n}} & m &= \txsum_{j=1}^n m_j & k_j &= \txsum_{i=1}^{m_j} k_{j,i} & k &= \txsum_{j=1}^n k_j\\ 
k_{\crdot} &= \ang{k_j}_{j\in \ufs{n}} & k_{j,\crdot} &= \ang{k_{j,i}}_{i \in \ufs{m}_j} & k_{\crdot,\crdot} &= \ang{k_{j,\crdot} }_{j\in \ufs{n}} &&\\
x_{\crdot} &= \ang{x_j}_{j\in \ufs{n}} & x_{j,\crdot} &= \ang{x_{j,i}}_{i \in \ufs{m}_j} & x_{\crdot,\crdot} &= \ang{x_{j,\crdot} }_{j\in \ufs{n}} &&\\
a_{j,i,\crdot} &= \ang{a_{j,i,h}}_{h \in \ufs{k}_{j,i}} & a_{j,\crdot,\crdot} &= \ang{a_{j,i,\crdot}}_{i \in \ufs{m}_j}
& a_{\crdot,\crdot,\crdot} &= \ang{a_{j,\crdot,\crdot}}_{j\in \ufs{n}} &&
\end{aligned}
$}
\right.\]
Then the following diagram in $\A$ commutes.
\begin{equation}\label{pseudoalg_comp_axiom}
\begin{tikzpicture}[vcenter]
\def\h{4} \def\g{0} \def\v{1.4} \def\u{-1}
\draw[0cell=.8]
(0,0) node (x11) {\gaA_n\Big(x \sscs \bang{\gaA_{m_j} \big( x_j \sscs \ang{ \gaA_{k_{j,i}} ( x_{j,i} \sscs a_{j,i,\crdot} ) }_{i \in \ufs{m}_j} \big) }_{j\in \ufs{n}}\Big)}
(x11)++(\h,\v) node (x12) {\gaA_n \Big(x \sscs \bang{\gaA_{k_j}\big( \ga\left(x_j \sscs x_{j,\crdot} \right) \sscs a_{j,\crdot,\crdot} \big) }_{j\in \ufs{n}} \Big)}
(x11)++(\h+.5,\u) node (x) {\gaA_k\Big(\ga\big(x \sscs \bang{ \ga\left(x_j \sscs x_{j,\crdot} \right) }_{j\in \ufs{n}} \big) \sscs a_{\crdot,\crdot,\crdot} \Big)}
(x11)++(0,2*\u) node (x21) {\gaA_m\Big(\ga\left(x \sscs x_\crdot \right) \sscs  \bang{ \bang{\gaA_{k_{j,i}}\left( x_{j,i} \sscs a_{j,i,\crdot} \right) }_{i \in \ufs{m}_j} }_{j\in \ufs{n}} \Big)}
(x21)++(\h,-\v) node (x22) {\gaA_k \Big( \ga\big( \ga\left(x \sscs x_\crdot \right) \sscs x_{\crdot,\crdot} \big) \sscs a_{\crdot,\crdot,\crdot} \Big)}
;
\draw[1cell=.85]  
(x11) edge[transform canvas={xshift={-3em}}, shorten >=-3pt] node[pos=.4] {\gaA_n\big(1_x \sscs \ang{ \phiA_{(m_j;\, k_{j,\crdot})} }_{j\in \ufs{n}} \big)} (x12)
(x12) edge node[pos=.6] {\phiA_{(n;\, k_{\crdot})}} (x)
(x) edge[equal] (x22)
(x11) edge[transform canvas={xshift={-1em}}] node[swap] {\phiA_{(n;\, m_{\crdot})}} (x21)
(x21) edge[transform canvas={xshift={-3em}}, shorten >=-1em] node[swap,pos=.5] {\phiA_{(m;\, k_{\crdot,\crdot} )}} (x22)
;
\end{tikzpicture}
\end{equation}
In \cref{pseudoalg_comp_axiom}, the equality holds by the associativity axiom \cref{enr-multicategory-associativity} for $\Op$.
\end{description}
This finishes the definition of an $\Op$-pseudoalgebra.
\end{definition}

\begin{explanation}[Pseudoalgebras]\label{expl:pseudoalg_axioms}
Consider \cref{def:pseudoalgebra}.
\begin{description}
\item[Action] The action equivariance axiom \cref{pseudoalg_action_sym} and the action unity axiom \cref{pseudoalg_action_unity} are the usual ones for an algebra over an operad; see \cite[Chapter 13]{yau-operad} for an elementary discussion.  If each component of the associativity constraint $\phiA$ \cref{phiA} is the identity $G$-natural transformation, then the axioms \cref{pseudoalg_basept_axiom,pseudoalg_topeq,pseudoalg_boteq,pseudoalg_unity,pseudoalg_comp_axiom} hold automatically.  In this case, $(\A,\gaA)$ is an $\Op$-algebra in $\Gcat$ in the usual sense.  In other words, in an $\Op$-pseudoalgebra, only the associativity axiom of an $\Op$-algebra is relaxed.  The associativity axiom is replaced by the associativity constraint $\phiA$, which is a $G$-natural isomorphism that satisfies the coherence axioms \cref{pseudoalg_basept_axiom,pseudoalg_topeq,pseudoalg_boteq,pseudoalg_unity,pseudoalg_comp_axiom}.  Each of these coherence axioms is modeled after an operad axiom in \cref{def:enr-multicategory}.

\item[Basepoint] In the basepoint axiom \cref{pseudoalg_basept_axiom}, we think of $\dy_j$ as a degeneracy operator.  Thus, this axiom says that the associativity constraint $\phiA$ is the identity in degenerate cases.

\item[Generality] As we discuss in \cref{rk:pseudocom_generality}, \cref{def:pseudoalgebra} also makes sense if $\Gcat$ is replaced by the 2-category $\Cat(\cV)$ of internal $\cV$-categories, internal $\cV$-functors, and internal $\cV$-natural transformations, for a complete and cocomplete Cartesian closed category $\cV$.  This is the setting of \cite[Def.\! 2.14]{gmmo20}.  The table below summarizes the correspondence of axioms between this work and \cite{gmmo20}.
\begin{center}
\resizebox{.9\width}{!}{%
{\renewcommand{\arraystretch}{1.2}%
{\setlength{\tabcolsep}{1em}
\begin{tabular}{c|c}
\cref{def:pseudoalgebra} & \cite[Def.\! 2.14]{gmmo20} \\ \hline
Action equivariance \cref{pseudoalg_action_sym} & Axiom 2.16\\
Action unity \cref{pseudoalg_action_unity} & Axiom 2.18\\
Basepoint \cref{pseudoalg_basept_axiom} & Axiom 2.17\\
Equivariance \cref{pseudoalg_topeq,pseudoalg_boteq} & Axiom 2.19\\
Unity \cref{pseudoalg_unity} & Axiom 2.21\\
Composition \cref{pseudoalg_comp_axiom} & Axiom 2.22
\end{tabular}}}}
\end{center}

\item[Variant] \cite[Def.\! 2.14]{gmmo20} is adapted from \cite[Def.\! 2.2]{corner-gurski}, which defines pseudoalgebras over a $\mathbf{G}$-operad $\Op$ for an action operad $\mathbf{G}$.  However, in \cite{corner-gurski}, $\Op$ is not assumed to be reduced.  Moreover, in \cite{corner-gurski}, the operadic unit is not assumed to act strictly as the identity on its pseudoalgebras, in contrast to the action unity axiom \cref{pseudoalg_action_unity}.\defmark
\end{description}
\end{explanation}

\begin{explanation}[$\BE$-Pseudoalgebras]\label{expl:P_pseudoalgebras}
For the Barratt-Eccles reduced $\Gcat$-operad $\BE$ (\cref{def:BE-Gcat}), we show in \cref{BEpseudo_smcat,smcat_BEpseudo} that each $\BE$-pseudoalgebra yields a naive symmetric monoidal $G$-category, and vice versa.  One round trip composite---from a naive symmetric monoidal $G$-category to a $\BE$-pseudoalgebra and back---is the identity by \cref{Psi_Phi}.  The other composite is \emph{not} the identity (\cref{expl:BEpseudo_smcat}).  Nevertheless, these assignments constitute the object assignments of an inverse pair of 2-equivalences; see \cref{thm:BEpseudoalg}.
\end{explanation}

For a group $G$, recall the translation category $\EG$ in \cref{ex:EG}.  The functor $\Catg(\EG,-)$ preserves all small limits, including finite products, because it is a right adjoint by \cref{Catg_rightadj}.  Moreover, since $\Catg(-,-)$ is equipped with the conjugation $G$-action, $\Catg(\EG,-)$ preserves $G$-functors and $G$-natural transformations.  These facts yield the following $\Op$-pseudoalgebra analogue of \cref{naive_genuine_pGcat}.  See also \cref{naive_genuine_smgcat} and \cite[Proposition 3.7]{gmmo20}.

\begin{proposition}\label{catgego}
Suppose $\Op$ is a $\Gcat$-operad for a group $G$.  Then the following statements hold.
\begin{enumerate}
\item\label{ng_i} Applying $\Catg(\EG,-)$ levelwise to $\Op$ and to the $\Gcat$-operad structure of $\Op$ yields a $\Gcat$-operad $\Catg(\EG,\Op)$, which is, furthermore, reduced if $\Op$ is.
\item\label{ng_ii} Applying $\Catg(\EG,-)$ to an $\Op$-pseudoalgebra $(\A,\gaA,\phiA)$ yields a $\Catg(\EG,\Op)$-pseudoalgebra $\Catg(\EG,\A)$.
\end{enumerate}
\end{proposition}

\subsection*{Further Unity Properties}

The basepoint axiom \cref{pseudoalg_basept_axiom} and the unity axiom \cref{pseudoalg_unity} assert that some components of $\phiA$ are identities.  In \cref{def:pseudoalgebra}, $\Op$ is assumed to be reduced, so $\Op(0)$ is the terminal $G$-category with only one object $*$ and its identity morphism.  The following observation says that components of $\phiA$ involving $* \in \Op(0)$ are identities.

\begin{lemma}\label{phi_id}
For each $\Op$-pseudoalgebra $(\A,\gaA,\phiA)$ and object $x \in \Op(n)$ with $n>0$, the component of the associativity constraint 
\[\gaA_n\big(x ; \ang{\zero}_{j \in \ufs{n}}\big) = \gaA_n\big(x ; \ang{\gaA_0(*)}_{j \in \ufs{n}} \big) \fto{\phiA_{(n;\,0^n)}} \gaA_0(*) = \zero\]
is equal to $1_\zero$, where $0^n = \ang{0}_{j \in \ufs{n}}$.
\end{lemma}

\begin{proof}
We prove this assertion by induction.  If $n=1$, then $\phiA_{(1;\, 0)} = 1_\zero$ by the basepoint axiom \cref{pseudoalg_basept_axiom}.

Inductively, suppose $n>1$ and $\phiA_{(p;\, 0^p)} = 1_\zero$ for $p < n$.  We consider the composition axiom \cref{pseudoalg_comp_axiom} with the following parameters.
\begin{itemize}
\item For $1 \leq j \leq n-1$, we take $m_j = 1$, $x_j = \opu \in \Op(1)$, $k_{j,1} = 0$, and $x_{j,1} = * \in \Op(0)$.
\item We take $m_n = 0$ and $x_n = * \in \Op(0)$.
\end{itemize}
In the diagram \cref{pseudoalg_comp_axiom}, the upper right arrow, denoted $\phiA_{(n;\, k_\crdot)}$ there, is the morphism $\phiA_{(n;\,0^n)}$ under consideration.  Each of the other three arrows in \cref{pseudoalg_comp_axiom} is the identity for the following reasons.
\begin{enumerate}
\item The upper left arrow, denoted $\gaA_n\big(1_x; \ang{\phiA_{(m_j; \, k_{j,\crdot})}}_{j \in \ufs{n}} \big)$ there, is the identity by 
\begin{itemize}
\item the initial case---that is, $\phiA_{(1;\, 0)} = 1_\zero$---for $j < n$, 
\item the convention $\phiA_{(0; \ang{})} = 1_\zero$ for $j=n$, and
\item the functoriality of $\gaA_n$.
\end{itemize}
\item The left vertical arrow, denoted $\phiA_{(n;\, m_\crdot)}$ there, is the identity by the basepoint axiom \cref{pseudoalg_basept_axiom}, since $m_\crdot = (1, \ldots, 1, 0)$.
\item The lower left arrow, denoted $\phiA_{(m;\, k_{\crdot,\crdot})}$ there, is the identity by the induction hypothesis, since 
\[\ga(x; x_\crdot) = \ga\big(x; \opu, \ldots, \opu, * \big) \in \Op(n-1).\]
\end{enumerate}
Thus, the remaining arrow, $\phiA_{(n;\,0^n)}$, is also the identity.
\end{proof}

\subsection*{Lax Morphisms between Pseudoalgebras}

For an operad $\Op$, recall that an $\Op$-algebra morphism between $\Op$-algebras 
\[f \cn (\A,\gaA) \to (\B,\gaB)\]
is a morphism $f \cn \A \to \B$ that strictly preserves the $\Op$-algebra structure in the sense that the following diagram commutes for each $n \geq 0$.
\begin{equation}\label{algebra_morphism}
\begin{tikzpicture}[xscale=3.5,yscale=1.3,vcenter]
\draw[0cell=1]
(0,0) node (x11) {\Op(n) \times \A^n}
(x11)++(1,0) node (x12) {\Op(n) \times \B^n}
(x11)++(0,-1) node (x21) {\A}
(x21)++(1,0) node (x22) {\B}
;
\draw[1cell=1]  
(x11) edge node {1 \times f^n} (x12)
(x12) edge node {\gaB} (x22)
(x11) edge node[swap] {\gaA} (x21)
(x21) edge node{f} (x22)
;
\end{tikzpicture}
\end{equation}
See \cite[Def.\! 13.2.8]{yau-operad} for an elementary discussion.  The lax variant of an $\Op$-algebra morphism replaces the commutative diagram \cref{algebra_morphism} with a natural transformation that satisfies its own coherence axioms.  The next definition is a lax variant of \cite[Def.\! 2.23]{gmmo20}, which is, in turn, adapted from \cite[Def.\! 2.4]{corner-gurski}.

\begin{definition}\label{def:laxmorphism}
Suppose $(\Op,\ga,\opu)$ is a reduced $\Gcat$-operad for a group $G$.  Suppose $(\A,\gaA,\phiA)$ and $(\B,\gaB,\phiB)$ are $\Op$-pseudoalgebras \pcref{def:pseudoalgebra}.  A \emph{lax $\Op$-morphism}\index{lax morphism}\index{operad!lax morphism}
\[\big(\A,\gaA,\phiA\big) \fto{(f,\actf)} \big(\B,\gaB,\phiB\big)\]
consists of a $G$-functor $f \cn \A \to \B$ and a sequence of $G$-natural transformations $\actf_n$, one for each $n \geq 0$, as follows.
\begin{equation}\label{laxmorphism_constraint}
\begin{tikzpicture}[xscale=3.5,yscale=1.3,vcenter]
\draw[0cell=.9]
(0,0) node (x11) {\Op(n) \times \A^n}
(x11)++(1,0) node (x12) {\Op(n) \times \B^n}
(x11)++(0,-1) node (x21) {\A}
(x21)++(1,0) node (x22) {\B}
;
\draw[1cell=.9]  
(x11) edge node {1 \times f^n} (x12)
(x12) edge node {\gaB_n} (x22)
(x11) edge node[swap] {\gaA_n} (x21)
(x21) edge node[swap] {f} (x22)
;
\draw[2cell]
node[between=x11 and x22 at .5, shift={(0,0)}, rotate=-90, 2label={above,\actf_n}] {\Rightarrow}
;
\end{tikzpicture}
\end{equation}
We call $\actf_n$ an \emph{action constraint}\index{action constraint} and sometimes abbreviate it to $\actf$ or $\act$.  A typical component of $\actf_n$ is a morphism in $\B$ 
\begin{equation}\label{actf_component}
\gaB_n\big(x \sscs \bang{f(a_j)}_{j\in \ufs{n}} \big) \fto{\actf_n} 
f\big(\gaA_n (x \sscs \ang{a_j}_{j\in \ufs{n}} ) \big)
\end{equation}
for objects $x \in \Op(n)$ and $\ang{a_j}_{j\in \ufs{n}} \in \A^n$.  The pair $(f,\actf)$ is subject to the axioms \cref{laxmorphism_basepoint,laxmorphism_unity,laxmorphism_equiv,laxmorphism_associativity} below.
\begin{description}
\item[Basepoint] $\actf_0$ is the identity morphism:
\begin{equation}\label{laxmorphism_basepoint}
\zero^\B = \gaB_0(*) \fto{\actf_0 \,=\, 1} f\big(\gaA_0(*)\big) = f(\zero^\A).
\end{equation}

\item[Unity] The following diagram in $\B$ commutes for each object $a \in \A$, with $\opu \in \Op(1)$ the operadic unit.
\begin{equation}\label{laxmorphism_unity}
\begin{tikzpicture}[xscale=3.5,yscale=1,vcenter]
\draw[0cell=.9]
(0,0) node (x11) {\gaB_1\big( \opu \sscs f(a) \big)}
(x11)++(1,0) node (x12) {f\big(\gaA_1(\opu \sscs a) \big)}
(x11)++(0,-1) node (x21) {f(a)}
(x21)++(1,0) node (x22) {f(a)}
;
\draw[1cell=.9]  
(x11) edge node {\actf_1} (x12)
(x21) edge node {1} (x22)
(x11) edge[equal] (x21)
(x12) edge[equal] (x22)
;
\end{tikzpicture}
\end{equation}
The two equalities hold by the action unity axiom \cref{pseudoalg_action_unity}.

\item[Equivariance] For each permutation $\sigma \in \Sigma_n$ and objects $x \in \Op(n)$ and $\ang{a_j}_{j\in \ufs{n}} \in \A^n$, the following diagram in $\B$ commutes.
\begin{equation}\label{laxmorphism_equiv}
\begin{tikzpicture}[xscale=4.5,yscale=1.2,vcenter]
\draw[0cell=.9]
(0,0) node (x11) {\gaB_n\big(x\sigma \sscs \ang{f(a_{\sigma(j)})}_{j\in \ufs{n}} \big)}
(x11)++(1,0) node (x12) {f\big(\gaA_n \big(x\sigma \sscs \ang{a_{\sigma(j)}}_{j\in \ufs{n}}  \big) \big)}
(x11)++(0,-1) node (x21) {\gaB_n\big(x \sscs \bang{f(a_j)}_{j\in \ufs{n}} \big)}
(x21)++(1,0) node (x22) {f\big(\gaA_n (x \sscs \ang{a_j}_{j\in \ufs{n}} ) \big)}
;
\draw[1cell=.9]  
(x11) edge node {\actf_n} (x12)
(x21) edge node {\actf_n} (x22) 
(x11) edge[equal] (x21)
(x12) edge[equal] (x22)
;
\end{tikzpicture}
\end{equation}
The two equalities hold by the action equivariance axiom \cref{pseudoalg_action_sym}.

\item[Associativity] Consider objects 
\[x \in \Op(n),\quad x_j \in \Op(m_j), \andspace a_{j,i} \in \A\]
for $1 \leq j \leq n$ and $1 \leq i \leq m_j$, along with the following notation, where we use $\bdot\,$ to denote a running index in a finite sequence.
\[\left\{\scalebox{.9}{$
\begin{aligned}
m_{\crdot} &= \ang{m_j}_{j\in \ufs{n}} & a_{j,\crdot} &= \ang{a_{j,i}}_{i \in \ufs{m}_j} & a_{\crdot,\crdot} &= \ang{a_{j,\crdot}}_{j\in \ufs{n}} & x_{\crdot} &= \ang{x_j}_{j\in \ufs{n}}\\
\actf_{m_{\crdot}} &= \bang{\actf_{m_j}}_{j\in \ufs{n}} & fa_{j,\crdot} &= \bang{f(a_{j,i})}_{i \in \ufs{m}_j} & fa_{\crdot,\crdot} &= \bang{fa_{j,\crdot}}_{j\in \ufs{n}} & m &= \txsum_{j=1}^n m_j
\end{aligned}
$}
\right.\]
Then the following diagram in $\B$ commutes.
\begin{equation}\label{laxmorphism_associativity}
\begin{tikzpicture}[xscale=1,yscale=1,vcenter]
\def\h{3.3} \def\g{0} \def\v{1.2} \def\u{-1}
\draw[0cell=.8]
(0,0) node (x11) {\gaB_n\big(x \sscs \big\langle \gaB_{m_j} ( x_j \sscs fa_{j,\crdot} ) \big\rangle_{j\in \ufs{n}} \big)}
(x11)++(\h,\v) node (x12) {\gaB_m\big( \ga(x\sscs x_\crdot) \sscs fa_{\crdot,\crdot} \big)}
(x11)++(\h+.5,\u) node (x) {f \gaA_m\big( \ga(x \sscs x_\crdot) \sscs a_{\crdot,\crdot} \big)}
(x11)++(0,2*\u) node (x21) {\gaB_n\big(x \sscs \big\langle f\gaA_{m_j} (x_j \sscs a_{j,\crdot} ) \big\rangle_{j\in \ufs{n}} \big)}
(x21)++(\h,-\v) node (x22) {f\gaA_n\big( x \sscs \bang{\gaA_{m_j}(x_j \sscs a_{j,\crdot}) }_{j\in \ufs{n}} \big)}
;
\draw[1cell=.85]  
(x11) edge[transform canvas={xshift={-1em}}, shorten >=0ex] node[pos=.4] {\phiB_{(n;\, m_\crdot)}} (x12)
(x12) edge node[pos=.6] {\actf_m} (x)
(x11) edge[transform canvas={xshift={0em}}] node[swap] {\gaB_n\big( 1_x \sscs \actf_{m_{\crdot}} \big)} (x21)
(x21) edge[transform canvas={xshift={-1em}}, shorten >=0pt] node[swap,pos=.5] {\actf_n} (x22)
(x22) edge node[swap,pos=.6] {f\phiA_{(n;\, m_\crdot)}} (x)
;
\end{tikzpicture}
\end{equation}
\end{description}
This finishes the definition of a lax $\Op$-morphism.  We call $(f,\actf)$ an \emph{$\Op$-pseudomorphism}\index{pseudomorphism}\index{operad!pseudomorphism} if each $\actf_n$ is a $G$-natural isomorphism.  We call $(f,\actf)$ a \emph{strict $\Op$-morphism}\index{strict morphism}\index{operad!strict morphism} if each $\actf_n$ is the identity.
\end{definition}

\begin{explanation}[Lax Morphisms]\label{expl:laxmorphism}
We emphasize that for a lax $\Op$-morphism $(f,\actf)$ \pcref{def:laxmorphism}, the action constraints $\actf_n$ are not required to be invertible.  For the Barratt-Eccles reduced $\Gcat$-operad $\BE$, we show in \cref{laxmorphism_smf,smf_laxmorphism,Psi_Phi_smf,Phi_injective_onecells} that lax $\BE$-morphisms correspond bijectively to strictly unital symmetric monoidal $G$-functors.  Furthermore, when $\Op$ is a pseudo-commutative $\Gcat$-operad, \cref{thm:multpso} shows that there is a $\Gcat$-multicategory $\MultpsO$ whose $G$-fixed 1-ary 1-cells are precisely lax $\Op$-morphisms.  See \cref{expl:onelax_morphism}.
\end{explanation}

\subsection*{Transformations between Lax Morphisms}

The following definition is \cite[Def.\! 2.24]{gmmo20}, extended to lax morphisms.  The $\mathbf{G}$-operad variant is \cite[Def.\! 2.7]{corner-gurski}.

\begin{definition}\label{def:algtwocells}
In the context of \cref{def:laxmorphism}, suppose 
\[(f,\actf), (h,\acth) \cn \big(\A,\gaA,\phiA\big) \to \big(\B,\gaB,\phiB\big)\]
are lax $\Op$-morphisms between $\Op$-pseudoalgebras.  An \emph{$\Op$-transformation}\index{transformation}\index{operad!transformation}
\[(f,\actf) \fto{\omega} (h,\acth)\]
is a $G$-natural transformation $\omega \cn f \to h$ such that the diagram in $\B$
\begin{equation}\label{Otransformation_ax}
\begin{tikzpicture}[xscale=4,yscale=1.5,vcenter]
\draw[0cell=.9]
(0,0) node (x11) {\gaB_n\big(x \sscs \ang{f(a_j)}_{j\in \ufs{n}} \big)}
(x11)++(1,0) node (x12) {f\gaA_n \big(x \sscs \anga \big)}
(x11)++(0,-1) node (x21) {\gaB_n\big(x \sscs \bang{h(a_j)}_{j\in \ufs{n}} \big)}
(x21)++(1,0) node (x22) {h\gaA_n \big(x \sscs \anga \big)}
;
\draw[1cell=.9]  
(x11) edge node {\actf_n} (x12)
(x21) edge node {\acth_n} (x22) 
(x11) edge[transform canvas={xshift={2em}}] node[swap] {\gaB_n\big(1_x \sscs \ang{\omega_{a_j}}_{j\in \ufs{n}} \big)} (x21)
(x12) edge[transform canvas={xshift={-1em}}, shorten >=1pt] node {\omega_{\gaA_n(x;\, \anga)}} (x22)
;
\end{tikzpicture}
\end{equation}
commutes for objects $x \in \Op(n)$ and $\anga = \ang{a_j}_{j\in \ufs{n}} \in \A^n$.
\end{definition}

\begin{explanation}[Transformations]\label{expl:transformation}
We emphasize that an $\Op$-transformation is not required to be invertible.  Moreover, we note the following.
\begin{enumerate}
\item By the basepoint axiom \cref{laxmorphism_basepoint} for $f$ and $h$, the axiom \cref{Otransformation_ax} in the case $n=0$ is the equality
\begin{equation}\label{omega-zero}
\omega_{\zero^\A} = 1_{\zero^\B} \cn f(\zero^\A) = \zero^\B \to h(\zero^\A) = \zero^\B.
\end{equation}
\item For $n=1$ and $x = \opu \in \Op(1)$ the operadic unit, the diagram \cref{Otransformation_ax} is automatically commutative by the action unity axiom \cref{pseudoalg_action_unity} for $\A$ and $\B$ and the unity axiom \cref{laxmorphism_unity} for $f$ and $h$.
\item For the Barratt-Eccles reduced $\Gcat$-operad $\BE$ (\cref{def:BE-Gcat}), we show in \cref{BEtransformation,Phi_twosurject} that $\BE$-transformations correspond bijectively to monoidal $G$-natural transformations.
\item When $\Op$ is a pseudo-commutative $\Gcat$-operad, we show in \cref{thm:multpso} that there is a $\Gcat$-multicategory $\MultpsO$ whose $G$-fixed 1-ary 2-cells are precisely $\Op$-transformations.  See \cref{expl:ktransform_basept}.\defmark
\end{enumerate}
\end{explanation}

\subsection*{2-Category Structure}

Starting from a $\Gcat$-multicategory, taking the $G$-fixed subcategories of the multimorphism $G$-categories yields a $\Cat$-multicategory; see \cref{thm:multpso_gfixed}.  Further restricting to the 1-ary 1-cells and 1-ary 2-cells yields a 2-category.  Thus, the following observation is a consequence of the existence of the $\Gcat$-multicategory $\MultpsO$, which is proved in detail in \cref{thm:multpso}.  See \cref{ex:underlying_iicat}.  We state this simpler result here for future reference.

\begin{proposition}\label{oalgps_twocat}
Suppose $G$ is a group, and $(\Op,\ga,\opu)$ is a reduced $\Gcat$-operad.  Then there is a 2-category $\AlglaxO$\index{pseudoalgebra!2-category}\index{2-category!pseudoalgebra} with 
\begin{itemize}
\item $\Op$-pseudoalgebras (\cref{def:pseudoalgebra}) as objects,
\item lax $\Op$-morphisms (\cref{def:laxmorphism}) as 1-cells, and
\item $\Op$-transformations (\cref{def:algtwocells}) as 2-cells.
\end{itemize}
Moreover, there are variant 2-categories $\AlgpspsO$ and $\AlgstO$ with the same objects and 2-cells as $\AlglaxO$, and with 1-cells given by, respectively, $\Op$-pseudomorphisms and strict $\Op$-morphisms.
\end{proposition}

\begin{proof}
The horizontal composition of two composable lax $\Op$-morphisms
\[\begin{tikzpicture}[xscale=1,yscale=1,vcenter]
\def\h{3}
\draw[0cell=.9]
(0,0) node (a) {\big(\A,\gaA,\phiA\big)}
(a)++(\h,0) node (b) {\big(\B,\gaB,\phiB\big)}
(b)++(\h,0) node (c) {\big(\C,\gaC,\phiC\big)}
;
\draw[1cell=.9]
(a) edge node {(f,\actf)} (b)
(b) edge node {(h,\acth)} (c)
;
\end{tikzpicture}\]
is given by the composite $G$-functor
\begin{equation}\label{Omorphism_comp}
\A \fto{hf} \C
\end{equation}
and the action constraint $\acthf_n$ defined by the following pasting of $G$-natural transformations.
\begin{equation}\label{Omorphism_paste}
\begin{tikzpicture}[xscale=3,yscale=1.3,vcenter]
\draw[0cell=.9]
(0,0) node (x11) {\Op(n) \times \A^n}
(x11)++(1,0) node (x12) {\Op(n) \times \B^n}
(x12)++(1,0) node (x13) {\Op(n) \times \C^n}
(x11)++(0,-1) node (x21) {\A}
(x21)++(1,0) node (x22) {\B}
(x22)++(1,0) node (x23) {\C}
;
\draw[1cell=.9]  
(x11) edge node {1 \times f^n} (x12)
(x12) edge node {1 \times h^n} (x13)
(x11) edge node[swap] {\gaA_n} (x21)
(x12) edge node[swap] {\gaB_n} (x22)
(x13) edge node {\gaC_n} (x23)
(x21) edge node[swap] {f} (x22)
(x22) edge node[swap] {h} (x23)
;
\draw[2cell]
node[between=x11 and x22 at .5, shift={(-.1,0)}, rotate=-90, 2label={above,\actf_n}] {\Rightarrow}
node[between=x12 and x23 at .5, shift={(0,0)}, rotate=-90, 2label={above,\acth_n}] {\Rightarrow}
;
\end{tikzpicture}
\end{equation}
Note that if both $(f,\actf)$ and $(h,\acth)$ are $\Op$-pseudomorphisms, then the above pasting is invertible.  Thus, $(hf,\acthf)$ is also an $\Op$-pseudomorphism.  Similarly, if both $(f,\actf)$ and $(h,\acth)$ are strict $\Op$-morphisms, then the above pasting is the identity.  Thus, $(hf,\acthf)$ is also a strict $\Op$-morphism.

An $\Op$-transformation $\omega$ between lax $\Op$-morphisms is a $G$-natural transformation with an extra property \cref{Otransformation_ax}, but no extra structure.  Horizontal and vertical compositions of 2-cells are defined as those of $G$-natural transformations, which are, in turn, those of natural transformations.
\end{proof}

\section{Naive Symmetric Monoidal $G$-Categories}
\label{sec:naive_smc}

This section discusses a symmetric monoidal variant  of naive permutative $G$-categories, called naive symmetric monoidal $G$-categories.  Recall from \cref{def:BE-Gcat} that naive permutative $G$-categories are algebras over the Barratt-Eccles operad $\BE$ in $\Gcat$.  By \cref{expl:naive_perm_Gcat}, a naive permutative $G$-category is precisely a small $G$-category equipped with a $G$-equivariant permutative category structure.  Naive symmetric monoidal $G$-categories are closely related to $\BE$-pseudoalgebras in the sense of \cref{def:pseudoalgebra}.

\secoutline 
\begin{itemize}
\item \cref{def:naive_smGcat} defines naive symmetric monoidal $G$-categories.
\item \cref{BEpseudo_smcat} shows that each $\BE$-pseudoalgebra $\A$ in $\Gcat$ yields a naive symmetric monoidal $G$-category $\Phi\A$.
\item \cref{laxmorphism_smf} shows that each lax $\BE$-morphism $(f,\actf$) between $\BE$-pseudoalgebras in $\Gcat$ yields a strictly unital symmetric monoidal $G$-functor $\Phi(f,\actf)$.
\item \cref{BEtransformation} shows that each $\BE$-transformation $\omega$ yields a monoidal $G$-natural transformation $\Phi(\omega) = \omega$.
\item \cref{alglaxbe_smgcat} shows that these assignments together constitute a 2-functor
\[\Phi \cn \AlglaxBE \to \smgcat.\]
The domain $\AlglaxBE$ is the 2-category of $\BE$-pseudoalgebras, lax $\BE$-morphisms, and $\BE$-transformations \pcref{oalgps_twocat}.  The codomain $\smgcat$ is the 2-category of naive symmetric monoidal $G$-categories, strictly unital symmetric monoidal $G$-functors, and monoidal $G$-natural transformations (\cref{def:smGcat_twocat}).  
\end{itemize}

As we discuss in \cref{sec:BEpseudoalg}, even more is true.  \cref{thm:BEpseudoalg} shows that the 2-functor $\Phi$ is a 2-equivalence.  Thus, in practice the 2-categories $\AlglaxBE$ and $\smgcat$ are interchangeable.  Whether one uses $\AlglaxBE$ or $\smgcat$ depends on the situation.  On the one hand, $\AlglaxBE$ provides the operadic description of its objects, 1-cells, and 2-cells.  On the other hand, $\smgcat$ describes these notions in terms of shorter lists of generating structures and axioms.  We emphasize that $\Phi$ is \emph{not} bijective on objects.  In particular, it is not an isomorphism of 1-categories or 2-categories.  We discuss this subtlety further in \cref{expl:BEpseudo_smcat}.

\subsection*{Naive Symmetric Monoidal $G$-Categories from $\BE$-Pseudoalgebras}

Recall that a (symmetric) monoidal category is \emph{strictly unital} if the left and right unit isomorphisms are identities (\cref{def:monoidalcategory,def:symmoncat}).  

\begin{definition}\label{def:naive_smGcat}
For a group $G$, a \emph{naive symmetric monoidal $G$-category}\index{naive symmetric monoidal G-category@naive symmetric monoidal $G$-category}\index{G-category@$G$-category!naive symmetric monoidal} 
\[(\A,\otimes,\zero,\alpha,\lambda,\rho,\xi)\] 
consists of
\begin{itemize}
\item a small $G$-category $\A$ and
\item a strictly unital symmetric monoidal structure $(\otimes,\zero,\alpha,\lambda=1,\rho=1,\xi)$
\end{itemize}
such that the following two statements hold.
\begin{enumerate}
\item The monoidal unit $\zero \in \A$ is $G$-fixed.
\item The monoidal product $\otimes$, the associativity isomorphism $\alpha$, and the braiding $\xi$ are $G$-equivariant.\defmark
\end{enumerate}
\end{definition}

Regarding the Barratt-Eccles operad $\BE$ as a reduced $\Gcat$-operad with $G$ acting trivially (\cref{def:BE-Gcat}), we consider $\BE$-pseudoalgebras \pcref{def:pseudoalgebra}.  The next result provides a canonical passage from $\BE$-pseudoalgebras in $\Gcat$ to naive symmetric monoidal $G$-categories.

\begin{lemma}\label{BEpseudo_smcat}
Each $\BE$-pseudoalgebra in $\Gcat$ yields a naive symmetric monoidal $G$-category, as defined in \cref{naive_monunit,naive_monproduct,naive_braiding,naive_associativity,naive_unit_isos} below.
\end{lemma}

\begin{proof}
Suppose $(\A,\gaA,\phiA)$ is a $\BE$-pseudoalgebra in $\Gcat$.  We construct the associated naive symmetric monoidal $G$-category 
\begin{equation}\label{naive_smgcat}
\Phi(\A,\gaA,\phiA) = \big(\A,\otimes,\zero,\alpha, \lambda=1, \rho=1, \xi\big)
\end{equation}
as follows.

\parhead{Monoidal unit, monoidal product, and braiding}.    
By \cref{def:pseudoalgebra}, $\A$ is a small $G$-category.  The monoidal unit is defined as the $G$-fixed basepoint in \cref{pseudoalg_zero}:
\begin{equation}\label{naive_monunit}
\zero = \gaA_0(*).
\end{equation}
Applied to the identity permutation $\id_2 \in \Sigma_2$, the $\BE$-action $G$-functor
\[\gaA_2 \cn \BE(2) \times \A^2 = \tn\Sigma_2 \times \A^2 \to \A\]
yields the $G$-equivariant monoidal product
\begin{equation}\label{naive_monproduct}
\otimes = \gaA_2\big(\id_2 \sscs -,-\big) \cn \A^2 \to \A.
\end{equation}
By the action equivariance axiom \cref{pseudoalg_action_sym}, for the non-identity permutation $\tau \in \Sigma_2$, there is an equality of functors
\[\gaA_2\big(\tau \sscs -,-\big) = \otimes \circ \tau \cn \A^2 \to \A.\]
Applied to the non-identity isomorphism
\[[\tau,\id_2] \cn \id_2 \fto{\iso} \tau \inspace \tn\Sigma_2,\]
the $\BE$-action $G$-functor $\gaA_2$ yields the $G$-equivariant braiding
\begin{equation}\label{naive_braiding}
\xi = \gaA_2\big( [\tau, \id_2] \sscs -,- \big) \cn \otimes \to \otimes \circ \tau.
\end{equation}

\parhead{Associativity isomorphism}.  To construct the $G$-equivariant associativity isomorphism $\alpha$, we use the object equalities 
\[\ga\big(\id_2\sscs \id_2,\id_1\big) = \id_3 = \ga\big(\id_2 \sscs \id_1,\id_2 \big) \in \tn\Sigma_3.\]
The component of $\alpha$ at a triple of objects $a,b,c \in \A$ is defined as the following composite in $\A$.
\begin{equation}\label{naive_associativity}
\begin{tikzpicture}[xscale=1,yscale=1,vcenter]
\def\h{5}
\draw[0cell=.9]
(0,0) node (x11) {(a \otimes b) \otimes c}
(x11)++(\h,0) node (x12) {a \otimes (b \otimes c)}
;
\draw[0cell=.8]
(x11)++(0,-1) node (x21) {\gaA_2\big(\id_2 \sscs \gaA_2( \id_2 \sscs a,b), \gaA_1( \id_1 \sscs c)\big)}
(x21)++(\h,0) node (x22) {\gaA_2\big(\id_2 \sscs \gaA_1(\id_1 \sscs a), \gaA_2(\id_2 \sscs b,c )\big)}
;
\draw[0cell=.9]
(x21)++(0,-1.5) node (x31) {\gaA_3\big(\ga (\id_2 \sscs \id_2,\id_1) \sscs a,b,c \big)}
(x31)++(\h,0) node (x32) {\gaA_3\big(\ga(\id_2 \sscs \id_1,\id_2) \sscs a,b,c \big)}
(x31)++(\h/2,-1) node (x) {\gaA_3\big(\id_3 \sscs a,b,c \big)}
;
\draw[1cell=.9]  
(x11) edge node {\alpha_{a,b,c}} (x12)
(x11) edge[equal] (x21)
(x12) edge[equal] (x22)
(x21) edge node[swap] {\phiA_{(2;\, 2,1)}} (x31)
(x31) edge[equal, shorten <=-1ex, transform canvas={xshift={-1ex}}] (x)
(x) edge[equal, shorten >=-1ex, transform canvas={xshift={1ex}}] (x32)
(x32) edge node[swap] {\big(\phiA_{(2;\, 1,2)} \big)^{-1}} (x22)
;
\end{tikzpicture}
\end{equation}
We note that the invertibility of $\phiA$ is necessary to define the second constituent morphism, $\big(\phiA_{(2;\, 1,2)} \big)^{-1}$, in $\alpha_{a,b,c}$.  

\parhead{Unit isomorphisms}.  The left unit isomorphism $\lambda$ and the right unit isomorphism $\rho$ at an object $a \in \A$ are defined as follows.
\begin{equation}\label{naive_unit_isos}
\begin{tikzpicture}[xscale=1,yscale=1,vcenter]
\def\h{5}
\draw[0cell=.9]
(0,0) node (x11) {\zero \otimes a}
(x11)++(\h,0) node (x12) {a \otimes \zero}
(x11)++(0,-1) node (x21) {\gaA_2\big( \id_2 \sscs \gaA_0(*), \gaA_1(\id_1 \sscs a)\big)}
(x21)++(\h,0) node (x22) {\gaA_2\big( \id_2 \sscs \gaA_1(\id_1 \sscs a), \gaA_0(*)\big)}
(x21)++(\h/2,-1) node (x) {\gaA_1(\id_1 \sscs a) = a}
;
\draw[1cell=.9]  
(x11) edge[equal] (x21)
(x12) edge[equal] (x22)
(x21) edge[shorten <=-1ex, transform canvas={xshift={-1ex}}] node[swap,pos=.3] {\lambda_a \,=\, \phiA_{(2;\, 0,1)}} (x)
(x22) edge[shorten <=-1ex, transform canvas={xshift={1ex}}] node[pos=.3] {\rho_a \,=\, \phiA_{(2;\, 1,0)}} (x)
;
\end{tikzpicture}
\end{equation}
Each of the components 
\[\lambda_a = \phiA_{(2;\, 0,1)} \andspace \rho_a = \phiA_{(2;\, 1,,0)}\] 
is the identity morphism by the basepoint axiom \cref{pseudoalg_basept_axiom}.  The desired strictly unital symmetric monoidal category axioms are all derived from the $\BE$-pseudoalgebra axioms in \cref{def:pseudoalgebra}, as we discuss next.

\parhead{Pentagon axiom}.  For the pentagon axiom \cref{monoidalcataxioms}, we consider objects $a,b,c,d \in \A$ along with the following abbreviations.
\begin{equation}\label{notation_abcd}
\begin{aligned}
ab &= a \otimes b = \gaA_2\big(\id_2 \sscs a,b\big) & abc &= \gaA_3\big(\id_3 \sscs a,b,c \big)\\
abcd &= \gaA_4\big(\id_4 \sscs a,b,c,d\big) & &
\end{aligned}
\end{equation}
Iterating these notation, we have the following objects that appear in the diagram \cref{naive_pentagon} below.
\[\scalebox{.9}{$
\begin{aligned}
(ab)cd &= \gaA_3\big(\id_3 \sscs \gaA_2(\id_2 \sscs a,b), c,d\big)&
ab(cd) &= \gaA_3\big(\id_3 \sscs a,b, \gaA_2(\id_2 \sscs c,d)\big)\\
(abc)d &= \gaA_2\big(\id_2 \sscs \gaA_3(\id_3 \sscs a,b,c), d\big) &
a(bcd) &= \gaA_2\big(\id_2 \sscs a, \gaA_3(\id_3 \sscs b,c,d)\big)
\end{aligned}$}\]
\[\begin{split}
(a(bc))d &= \gaA_2\big(\id_2 \sscs \gaA_2\big(\id_2 \sscs a, \gaA_2(\id_2 \sscs b,c) \big), d\big)\\
a(bc)d &= \gaA_3\big(\id_3 \sscs a, \gaA_2(\id_2 \sscs b,c) , d \big)\\
a((bc)d) &= \gaA_2\big(\id_2 \sscs a, \gaA_2\big(\id_2 \sscs \gaA_2(\id_2 \sscs b,c), d \big) \big)
\end{split}\]
Moreover, we display each component of $\alpha$ \cref{naive_associativity} as a zigzag, noting that the second constituent arrow is invertible.  With the above notation and convention, the pentagon diagram \cref{monoidalcataxioms} for $(\A,\otimes,\alpha)$ is given by the boundary of the following diagram in $\A$.
\begin{equation}\label{naive_pentagon}
\begin{tikzpicture}[xscale=1,yscale=1,vcenter]
\def\h{7} \def\a{.7} \def\b{.4} \def\u{1.3} \def\v{1.5}
\draw[0cell=.85]
(0,0) node (x1) {((ab)c)d}
(x1)++(\h,0) node (x2) {a(b(cd))}
(x1)++(\h/2,\v) node (y) {(ab)(cd)}
node[between=x1 and y at .5] (y1) {(ab)cd}
node[between=x2 and y at .5] (y2) {ab(cd)}
(x1)++(\b,-\u) node (z1) {(abc)d}
(x2)++(-\b,-\u) node (z2) {a(bcd)}
(z1)++(\b,-\u) node (w1) {(a(bc))d}
(z2)++(-\b,-\u) node (w2) {a((bc)d)}
node[between=w1 and w2 at .5] (w) {a(bc)d}
node[between=y and w at .5] (z) {abcd} 
;
\draw[1cell=.8]  
(x1) edge node[pos=.7] {\phiA_{(2;\, 2,1)}} (y1) 
(y) edge node[swap,pos=.2] {\phiA_{(2;\, 1,2)}} (y1)
(y) edge node[pos=.2] {\phiA_{(2;\, 2,1)}} (y2)
(x2) edge node[swap,pos=.5] {\phiA_{(2;\, 1,2)}} (y2)
(x1) edge node[swap,pos=.1] {\phiA_{(2;\, 2,1)} 1_d} (z1) 
(w1) edge node[pos=.8] {\phiA_{(2;\, 1,2)} 1_d} (z1)
(w1) edge node {\phiA_{(2;\, 2,1)}} (w)
(w2) edge node[swap,pos=.4] {\phiA_{(2;\, 1,2)}} (w)
(w2) edge node[swap,pos=.8] {1_a \phiA_{(2;\, 2,1)}} (z2)
(x2) edge node[pos=.2] {1_a \phiA_{(2;\, 1,2)}} (z2)
(y1) edge node[swap,pos=.4,inner sep=0pt] {\phiA_{(3;\, 2,1,1)}} (z) 
(y) edge node[swap,pos=.4] {\phiA_{(2;\, 2,2)}} (z)
(y2) edge node[pos=.2] {\phiA_{(3;\, 1,1,2)}} (z)
(z1) edge node[pos=.5] {\phiA_{(2;\, 3,1)}} (z)
(z2) edge node[swap,pos=.4] {\phiA_{(2;\, 1,3)}} (z)
(w) edge node {\phiA_{(3;\, 1,2,1)}} (z)
;
\end{tikzpicture}
\end{equation}
Each of the six regions in \cref{naive_pentagon} commutes by the composition axiom \cref{pseudoalg_comp_axiom}.  For the two upper middle triangles, we also need to use the unity axiom \cref{pseudoalg_unity} and the fact that $\gaA$ preserves identity morphisms.  Since each component of $\phiA$ is an isomorphism, the commutative diagram \cref{naive_pentagon} implies the pentagon axiom \cref{monoidalcataxioms} for $(\A,\otimes,\alpha)$.

\parhead{Middle unity axiom}.  Next, we consider the middle unity axiom \cref{monoidalcataxioms} for $\A$.  With the notation in \cref{notation_abcd}, the middle unity diagram is the boundary of the following diagram.
\begin{equation}\label{naive_midunity}
\begin{tikzpicture}[xscale=1,yscale=1,vcenter]
\def\h{5} \def\v{1}
\draw[0cell=.9]
(0,0) node (x1) {(a\zero)b}
(x1)++(\h,0) node (x2) {a(\zero b)}
(x1)++(\h/2,\v) node (y) {a\zero b}
(x1)++(\h/2,-\v) node (z) {ab} 
;
\draw[1cell=.9]  
(x1) edge node[pos=.5] {\phiA_{(2;\, 2,1)}} (y)
(x2) edge node[swap,pos=.5] {\phiA_{(2;\, 1,2)}} (y)
(x1) edge node[swap,pos=.5] {\phiA_{(2;\, 1,0)} 1_b} (z)
(x2) edge node[pos=.5] {1_a \phiA_{(2;\, 0,1)}} (z)
(y) edge node[swap,pos=.5] {\phiA_{(3;\, 1,0,1)}} (z)
;
\end{tikzpicture}
\end{equation}
Each of the two triangles in \cref{naive_midunity} commutes by the unity axiom \cref{pseudoalg_unity} and the composition axiom \cref{pseudoalg_comp_axiom}.  Since each component of $\phiA$ is an isomorphism, the commutative diagram \cref{naive_midunity} implies the middle unity axiom \cref{monoidalcataxioms} for $\A$.

\parhead{Symmetry axiom}.  The $G$-equivariant braiding $\xi$ is defined above in \cref{naive_braiding}.  The symmetry axiom \cref{symmoncatsymhexagon} of a symmetric monoidal category holds for $(\A,\xi)$ because the composite
\[\begin{tikzpicture}[xscale=1,yscale=1,baseline={(x1.base)}]
\def\h{2.5} \def\v{1}
\draw[0cell=1]
(0,0) node (x1) {\id_2}
(x1)++(\h,0) node (x2) {\tau}
(x2)++(\h,0) node (x3) {\id_2}
;
\draw[1cell=1]  
(x1) edge node {[\tau, \id_2]} (x2)
(x2) edge node {[\id_2, \tau]} (x3)
;
\end{tikzpicture}\]
in $\tn\Sigma_2$ is the identity morphism of $\id_2 \in \Sigma_2$.

\parhead{Hexagon axiom}.  First, we extend the notation in \cref{notation_abcd} as follows for objects $a,b,c \in \A$ and objects or morphisms $\sigma \in \tn\Sigma_n$.
\[\begin{split}
(ab)_{\sigma} &= \gaA_2\big(\sigma \sscs a,b\big) \forspace \sigma \in \tn\Sigma_2\\
(abc)_{\sigma} &= \gaA_3\big(\sigma \sscs a,b,c \big) \forspace \sigma \in \tn\Sigma_3
\end{split}\]
For example, with $\tau \in \Sigma_2$ denoting the non-identity permutation, we have the following objects in $\A$ that appear in the diagram \cref{naive_hexagon} below.
\[\begin{split}
(ab)_{\tau} c &= \gaA_2\big(\id_2 \sscs \gaA_2(\tau \sscs a,b), c\big) = (ba)c\\
\big(a(bc)\big)_{\tau} &= \gaA_2\big(\tau \sscs a, \gaA_2(\id_2 \sscs b,c)\big) = (bc)a
\end{split}\]
The top and bottom composites in the hexagon diagram \cref{symmoncatsymhexagon} are, respectively, the top-right and left-bottom boundary composites in the following diagram, which we explain further below.
\begin{equation}\label{naive_hexagon}
\begin{tikzpicture}[xscale=1,yscale=1,vcenter]
\def\h{5} \def\g{2.5} \def\e{1.5} \def\v{-1.2} \def\u{-1}
\draw[0cell=.85]
(0,0) node (x11) {(ab)c}
(x11)++(\h,0) node (x12) {(ab)_\tau c}
(x12)++(\g,0) node (x13) {(ba)c}
(x11)++(0,\v) node (x21) {abc}
(x21)++(\h,0) node (x22) {(abc)_{\tau \times \id_1}}
(x22)++(\g,0) node (x23) {bac}
(x21)++(0,\v) node (x31) {a(bc)}
(x23)++(0,\v) node (x32) {b(ac)}
(x31)++(0,\v) node (x41) {(a(bc))_\tau}
(x32)++(0,\v) node (x42) {b(ac)_\tau}
(x41)++(0,\u) node (x51) {(bc)a}
(x42)++(0,\u) node (x53) {b(ca)}
node[between=x51 and x53 at .5] (x52) {bca}
node[between=x41 and x42 at .35] (y1) {(abc)_{\tau\ang{1,2}}}
node[between=x41 and x42 at .65] (y2) {(bac)_{\id_1 \times \tau}}
node[between=x11 and x22 at .5, shift={(0,.1)}] () {\mathbf{n}}
node[between=y2 and x32 at .7] () {\mathbf{n}}
node[between=x31 and y1 at .3] () {\mathbf{n}}
node[between=x12 and x23 at .5] () {\mathbf{b}}
node[between=y2 and x53 at .5] () {\mathbf{b}}
node[between=x51 and y1 at .5] () {\mathbf{t}}
node[between=x31 and x22 at .65] () {\mathbf{f}}
node[between=x22 and y2 at .6] () {\mathbf{a}}
node[between=y1 and y2 at .5, shift={(0,-.5)}] () {\mathbf{a}}
;
\draw[1cell=.8]  
(x11) edge node {(ab)_{[\tau, \id_2]} 1_c} (x12)
(x12) edge[equal] (x13)
(x11) edge node[swap] {\phiA_{(2;\, 2,1)}} (x21)
(x12) edge node[swap] {\phiA_{(2;\, 2,1)}} (x22)
(x13) edge node {\phiA_{(2;\, 2,1)}} (x23)
(x21) edge node {(abc)_{[\tau \times \id_1, \id_3]}} (x22)
(x22) edge[equal] (x23)
(x31) edge node {\phiA_{(2;\, 1,2)}} (x21)
(x32) edge node[swap] {\phiA_{(2;\, 1,2)}} (x23)
(x31) edge node [swap] {(a(bc))_{[\tau,\id_2]}} (x41)
(x32) edge node {1_b (ac)_{[\tau,\id_2]}} (x42)
(x41) edge[equal] (x51)
(x42) edge[equal] (x53)
(y1) edge[equal] (y2)
(y1) edge[equal] (x52)
(y2) edge[equal] (x52)
(x41) edge node {\phiA_{(2;\, 1,2)}} (y1)
(x42) edge node[swap] {\phiA_{(2;\, 1,2)}} (y2)
(x51) edge node[swap,pos=.4] {\phiA_{(2;\, 2,1)}} (x52)
(x53) edge node[pos=.4] {\phiA_{(2;\, 1,2)}} (x52)
(x21) edge node[pos=.3] {(abc)_{[\tau\ang{1,2}, \id_3]}} (y1)
(x22) edge node[swap,pos=.6] {(abc)_{\theta \cdot (\tau \times \id_1)}} (y1)
(x23) edge node[swap] {(bac)_\theta} (y2)
;
\end{tikzpicture}
\end{equation}
In the diagram \cref{naive_hexagon}, the notation $[y,x] \cn x \to y$ denotes the unique morphism from $x$ to $y$ in a translation category (\cref{def:translation_cat}).  Thus, each of the morphisms
\[\begin{aligned}
[\tau,\id_2] &\cn \id_2 \to \tau && \text{in $\tn\Sigma_2$,}\\
[\tau \times \id_1, \id_3]  &\cn \id_3 \to \tau \times \id_1 && \text{in $\tn\Sigma_3$,}\\
[\tau\ang{1,2}, \id_3] &\cn \id_3 \to \tau\ang{1,2} && \text{in $\tn\Sigma_3$, and}\\
\theta \defeq [\id_1 \times \tau, \id_3] &\cn \id_3 \to \id_1 \times \tau && \text{in $\tn\Sigma_3$}
\end{aligned}\]
is the unique morphism with the indicated domain and codomain, where $\tau\ang{1,2} \in \Sigma_3$ is the block permutation that sends $(a,b,c)$ to $(b,c,a)$.  The morphism $\theta$ appears in the two arrows bordering the parallelogram labeled $\mathbf{a}$, and it renders the following diagram in $\tn\Sigma_3$ commutative.
\begin{equation}\label{tauonetwo}
\begin{tikzpicture}[xscale=1,yscale=1,vcenter]
\def\h{3.5}
\draw[0cell]
(0,0) node (x1) {\id_3}
(x1)++(\h,0) node (x3) {\tau\ang{1,2}}
(x1)++(\h/2,-1) node (x2) {\tau \times \id_1}
(x3)++(2.6,0) node (x4) {(\id_1 \times \tau)(\tau \times \id_1)}
;
\draw[1cell]
(x1) edge node {[\tau\ang{1,2}, \id_3]} (x3)
(x3) edge[equal] (x4)
(x1) edge node[swap,pos=.4] {[\tau \times \id_1, \id_3]} (x2)
(x2) edge node[swap,pos=.7] {\theta \cdot (\tau \times \id_1)} (x3)
;
\end{tikzpicture}
\end{equation}

Now we explain why \cref{naive_hexagon} is commutative.  Each of the three regions labeled $\mathbf{n}$ commutes by the naturality of $\phiA$.  Each of the two regions labeled $\mathbf{b}$ commutes by the bottom equivariance axiom \cref{pseudoalg_boteq}.  The region labeled $\mathbf{t}$ commutes by the top equivariance axiom \cref{pseudoalg_topeq}.  Each of the two regions labeled $\mathbf{a}$ commutes by the action equivariance axiom \cref{pseudoalg_action_sym}.  The region labeled $\mathbf{f}$ commutes by the functoriality of $\gaA_3$ and the commutative diagram \cref{tauonetwo}.   Since each arrow in \cref{naive_hexagon} is an isomorphism, the commutativity of this diagram implies the hexagon axiom for $\A$.

This finishes the proof that, for a $\BE$-pseudoalgebra $(\A,\gaA,\phiA)$ in $\Gcat$, $(\A,\otimes,\zero,\alpha,\xi)$ in \cref{naive_smgcat} is a naive symmetric monoidal $G$-category.
\end{proof}

The construction $\Phi$ in \cref{BEpseudo_smcat} has a right inverse $\Psi$, which we discuss in \cref{smcat_BEpseudo,Psi_Phi}.

\subsection*{Symmetric Monoidal $G$-Functors from Lax $\BE$-Morphisms}

\begin{definition}\label{def:smGfunctor}
Suppose $G$ is a group.  For naive symmetric monoidal $G$-categories $(\A,\otimes,\zero,\alpha,\xi)$ and $(\B,\otimes,\zero,\alpha,\xi)$, a \emph{symmetric monoidal $G$-functor}\index{symmetric monoidal G-functor@symmetric monoidal $G$-functor}\index{G-functor@$G$-functor!symmetric monoidal}
\[(f,f^2,f^0) \cn \A \to \B\]
is a symmetric monoidal functor (\cref{def:monoidalfunctor}) such that the following three conditions hold.
\begin{enumerate}
\item $f$ is a $G$-functor.
\item The unit constraint $f^0 \cn \zero^\B \to f(\zero^\A)$ is $G$-fixed.
\item The monoidal constraint $f^2$ is a $G$-natural transformation.
\end{enumerate}
Moreover, we say that $(f,f^2,f^0)$ is
\begin{itemize}
\item \emph{strictly unital} if $f^0 = 1_{\zero^\B}$;
\item \emph{strictly unital strong} if $f^0 = 1_{\zero^\B}$ and $f^2$ is invertible; and
\item \emph{strict} if $f^0$ and $f^2$ are identities.\defmark
\end{itemize}
\end{definition}

We continue to use the Barratt-Eccles reduced $\Gcat$-operad $\BE$ in \cref{def:BE-Gcat}.   Next, we observe that lax $\BE$-morphisms in the sense of \cref{def:laxmorphism} yield strictly unital symmetric monoidal $G$-functors.  The following observation uses the naive symmetric monoidal $G$-category $\Phi(\A,\gaA,\phiA)$ in \cref{naive_smgcat}.

\begin{lemma}\label{laxmorphism_smf}
Each lax $\BE$-morphism between $\BE$-pseudoalgebras in $\Gcat$
\[(f,\actf) \cn (\A,\gaA,\phiA) \to (\B,\gaB,\phiB)\]
yields a strictly unital symmetric monoidal $G$-functor
\begin{equation}\label{Phi_factf}
\Phi(f,\actf) = \big(f,f^2 = \actf_2,f^0 = \actf_0\big) \cn \Phi(\A,\gaA,\phiA) \to \Phi(\B,\gaB,\phiB).
\end{equation}
Moreover, if $(f,\actf)$ is a $\BE$-pseudomorphism, then $\Phi(f,\actf)$ is strictly unital strong.  If $(f,\actf)$ is a strict $\BE$-morphism, then $\Phi(f,\actf)$ is strict.
\end{lemma}

\begin{proof}
In this proof, we refer to the naive symmetric monoidal $G$-category
\[\Phi(\A,\gaA,\phiA) = \big(\A,\otimes,\zero,\alpha, \lambda=1, \rho=1, \xi\big)\]
defined in \cref{naive_monunit,naive_monproduct,naive_braiding,naive_associativity,naive_unit_isos} for both $\A$ and $\B$.  We add superscripts---for example, $\zero^\A$---if we want to emphasize the category to which the structure belongs.  By \cref{def:laxmorphism}, $f \cn \A \to \B$ is a $G$-functor.  We construct its unit and monoidal constraints as follows.

\parhead{Unit constraint}.  We define the unit constraint $f^0$ to be the 0-th action constraint
\begin{equation}\label{laxm_smf_unit}
\zero^\B = \gaB_0(*) \fto{f^0 \,=\, \actf_0} f(\zero^\A) = f\gaA_0(*).
\end{equation}
This morphism is the identity morphism by the basepoint axiom \cref{laxmorphism_basepoint} of the lax $\BE$-morphism $(f,\actf)$.  Moreover, $f^0$ is $G$-fixed since the basepoint $\zero^\B$ is $G$-fixed.

\parhead{Monoidal constraint}.  We define the monoidal constraint $f^2$ to be the second action constraint 
\begin{equation}\label{laxm_smf_mconstraint}
\begin{tikzpicture}[xscale=4,yscale=1,vcenter]
\draw[0cell=.9]
(0,0) node (x11) {fa_1 \otimes fa_2}
(x11)++(1,0) node (x12) {f(a_1 \otimes a_2)}
(x11)++(0,-1) node (x21) {\gaB_2\big(\id_2 \sscs fa_1, fa_2\big)}
(x21)++(1,0) node (x22) {f \gaA_2\big(\id_2 \sscs a_1, a_2\big)}
;
\draw[1cell=.9]  
(x11) edge node {f^2_{a_1,a_2}} (x12)
(x21) edge node {\actf_2} (x22)
(x11) edge[equal] (x21)
(x12) edge[equal] (x22)
;
\end{tikzpicture}
\end{equation}
for objects $a_1,a_2 \in \A$.  Since $\actf_2$ is a $G$-natural transformation, so is $f^2$.  Next, we verify for $(f,f^2,f^0=1)$ the axioms of a symmetric monoidal functor in \cref{def:monoidalfunctor}.

\parhead{Unity axioms}.  The left unit isomorphism $\lambda$ is defined as the associativity constraint $\vphi_{(2;\, 0,1)}$ in \cref{naive_unit_isos}.  For an object $a \in \A$, the left unity diagram \cref{monoidalfunctorunity} for $(f,f^2,f^0)$ is the following diagram in $\B$.
\begin{equation}\label{laxm_smf_leftu}
\begin{tikzpicture}[xscale=1,yscale=1,vcenter]
\def\h{5} \def\g{.3} \def\v{1.3}
\draw[0cell=.85]
(0,0) node (x11) {\gaB_2\Big(\id_2 \sscs \gaB_0(*), \gaB_1(\id_1 \sscs fa) \Big)}
(x11)++(\h,0) node (x12) {\gaB_1\big(\id_1 \sscs fa\big)}
(x12)++(0,-\v) node (x) {f\gaA_1\big(\id_1 \sscs a\big)}
(x)++(1.5,0) node (x') {fa}
(x11)++(0,-2*\v) node (x21) {\gaB_2\Big(\id_2 \sscs f\gaA_0(*), f\gaA_1(\id_1 \sscs a)\Big)}
(x21)++(\h,0) node (x22) {f\gaA_2\Big(\id_2 \sscs \gaA_0(*), \gaA_1(\id_1 \sscs a)\Big)}
;
\draw[1cell=.9]  
(x11) edge node {\phiB_{(2;\, 0,1)}} (x12)
(x12) edge node[pos=.5] {\actf_1} (x)
(x) edge[equal] (x')
(x11) edge node[swap] {\gaB_2\big(1 \sscs \actf_0, \actf_1\big)} (x21)
(x21) edge node {\actf_2} (x22)
(x22) edge node[swap,pos=.5] {f\phiA_{(2;\, 0,1)}} (x)
;
\end{tikzpicture}
\end{equation}
In \cref{laxm_smf_leftu}, the arrow 
\[\actf_1 \cn \gaB_1(\id_1 \sscs fa) \to f \gaA_1(\id_1 \sscs a)\]
is the identity morphism by the unity axiom \cref{laxmorphism_unity} of the lax $\BE$-morphism $(f,\actf)$.  The diagram \cref{laxm_smf_leftu} is an instance of the associativity diagram \cref{laxmorphism_associativity}, which is, therefore, commutative.  The right unity diagram is proved similarly, by switching $\gaB_0(*)$ and $\gaB_1(\id_1 \sscs fa)$ in the upper left corner and using $\vphi_{(2;\, 1,0)}$ in place of $\vphi_{(2;\, 0,1)}$.  

\parhead{Associativity axiom}.  The associativity isomorphism $\alpha$ is defined in \cref{naive_associativity}.  In the diagram \cref{laxm_smf_as} below, we reuse the notation in \cref{notation_abcd} and the convention in \cref{naive_pentagon}.  For example, we have the following objects.
\[\begin{split}
\big((fa)(fb)\big) (fc) &= \gaB_2\big(\id_2 \sscs \gaB_2(\id_2 \sscs fa, fb), fc \big)\\
(fa)(fb)(fc) &= \gaB_3\big(\id_3 \sscs fa, fb, fc\big)\\
\big(f(ab) \big) (fc) &= \gaB_2\big(\id_2 \sscs f \gaA_2(\id_2 \sscs a,b), fc \big)
\end{split}\]
Moreover, each component of $\alpha$ is displayed as a zigzag of two components of the associativity constraint $\vphi$, noting that the second arrow is invertible.  For objects $a,b,c \in \A$, the associativity diagram \cref{monoidalfunctorassoc} for $(f,f^2)$ is the boundary of the following diagram in $\B$.
\begin{equation}\label{laxm_smf_as}
\begin{tikzpicture}[xscale=1,yscale=1,vcenter]
\def\h{3.5} \def\v{-1.3}
\draw[0cell=.85]
(0,0) node (x11) {\big((fa)(fb)\big) (fc)}
(x11)++(\h,0) node (x12) {(fa)(fb)(fc)}
(x12)++(\h,0) node (x13) {(fa) \big((fb)(fc)\big)}
(x11)++(0,\v) node (x21) {\big(f(ab) \big) (fc)}
(x13)++(0,\v) node (x23) {(fa) \big(f(bc)\big)}
(x21)++(0,\v) node (x31) {f\big((ab)c\big)}
(x31)++(\h,0) node (x32) {f(abc)}
(x32)++(\h,0) node (x33) {f\big(a(bc)\big)}
;
\draw[1cell=.85]  
(x11) edge node {\phiB_{(2;\, 2,1)}} (x12)
(x13) edge node[swap] {\phiB_{(2;\, 1,2)}} (x12)
(x11) edge node[swap] {\actf_2 \actf_1} (x21)
(x21) edge node[swap] {\actf_2} (x31)
(x12) edge node {\actf_3} (x32)
(x13) edge node {\actf_1 \actf_2} (x23)
(x23) edge node {\actf_2} (x33)
(x31) edge node {f\phiA_{(2;\, 2,1)}} (x32)
(x33) edge node[swap] {f\phiA_{(2;\, 1,2)}} (x32)
;
\end{tikzpicture}
\end{equation}
In \cref{laxm_smf_as}, each instance of $\actf_1$ is an identity morphism by the unity axiom \cref{laxmorphism_unity}.  Each of the two squares is an instance of the associativity diagram \cref{laxmorphism_associativity}, which is, therefore, commutative.  Since $\vphi$ is componentwise an isomorphism, this proves the associativity axiom \cref{monoidalfunctorassoc} for $(f,f^2)$.

\parhead{Symmetry axiom}.  The braiding $\xi$ is defined in \cref{naive_braiding}.  For objects $a,b \in \A$, the symmetry diagram \cref{monoidalfunctorbraiding} for $(f,f^2)$ is the boundary of the following diagram in $\B$, where $\tau \in \Sigma_2$ is the non-identity permutation.
\begin{equation}\label{laxm_smf_sym}
\begin{tikzpicture}[xscale=1,yscale=1,vcenter]
\def\h{4.8} \def\g{2.5} \def\v{-1.4}
\draw[0cell=.8]
(0,0) node (x11) {\gaB_2\big(\id_2 \sscs fa,fb\big)}
(x11)++(\h,0) node (x12) {\gaB_2\big(\tau \sscs fa,fb\big)}
(x12)++(\g,0) node (x13) {\gaB_2\big(\id_2 \sscs fb,fa\big)}
(x11)++(0,\v) node (x21) {f\gaA_2\big(\id_2 \sscs a,b\big)}
(x21)++(\h,0) node (x22) {f\gaA_2\big(\tau \sscs a,b\big)}
(x22)++(\g,0) node (x23) {f\gaA_2\big(\id_2 \sscs b,a\big)}
;
\draw[1cell=.85]  
(x11) edge node {\gaB_2\big([\tau,\id_2] \sscs 1,1\big)} (x12)
(x12) edge[equal] (x13)
(x21) edge node {f \gaA_2\big([\tau,\id_2] \sscs 1,1\big)} (x22)
(x22) edge[equal] (x23)
(x11) edge node[swap] {\actf_2} (x21)
(x12) edge node {\actf_2} (x22)
(x13) edge node {\actf_2} (x23)
;
\end{tikzpicture}
\end{equation}
In \cref{laxm_smf_sym}, the left region commutes by the naturality of the action constraint $\actf_2$ \cref{laxmorphism_constraint}.  The right region commutes by the equivariance axiom \cref{laxmorphism_equiv}.

This finishes the proof that $\Phi(f,\actf)$ is a strictly unital symmetric monoidal $G$-functor.  Moreover, if $(f,\actf)$ is a $\BE$-pseudomorphism, then each $\actf_n$ is invertible.  Thus, $f^2 = \actf_2$ is invertible, and $\Phi(f,\actf)$ is strictly unital strong.  If $(f,\actf)$ is a strict $\BE$-morphism, then each $\actf_n$ is the identity.  Thus, $f^2 = \actf_2$ is the identity, and $\Phi(f,\actf)$ is strict.
\end{proof}

\begin{remark}[Degenerate Action Constraints are Identities]\label{rk:actf_basepoint}
In the left unity diagram \cref{laxm_smf_leftu} above, each of the four arrows along the left, top, and right boundaries is an identity morphism.  Thus, the component of $\actf_2$ along the bottom boundary is also an identity morphism.
\end{remark}

The assignment $\Phi$ in \cref{laxmorphism_smf} has a right inverse $\Psi$, which we discuss in \cref{smf_laxmorphism,Psi_Phi_smf}.  Moreover, $\Phi$ is injective by \cref{Phi_injective_onecells}.

\subsection*{Monoidal $G$-Natural Transformations from $\BE$-Transformations}

Recall from \cref{Gnattr} that a $G$-natural transformation is a natural transformation that commutes with the $G$-action.   Also recall monoidal natural transformations from \cref{def:monoidalnattr}.

\begin{definition}\label{def:monGnat}
For a group $G$, suppose
\[(f,f^2,f^0 = 1), (h,h^2,h^0 = 1) \cn (\A, \otimes, \zero, \alpha,\xi) \to (\B, \otimes, \zero, \alpha,\xi)\]
are strictly unital symmetric monoidal $G$-functors (\cref{def:smGfunctor}) between naive symmetric monoidal $G$-categories (\cref{def:naive_smGcat}).  A \emph{monoidal $G$-natural transformation}\index{monoidal G-natural transformation@monoidal $G$-natural transformation}\index{G-natural transformation@$G$-natural transformation!monoidal}
\[\psi \cn (f,f^2,f^0 = 1) \to (h,h^2,h^0 = 1)\]
is a monoidal natural transformation that is also $G$-equivariant.
\end{definition}

For the Barratt-Eccles reduced $\Gcat$-operad $\BE$ (\cref{def:BE-Gcat}), we now show that $\BE$-transformations (\cref{def:algtwocells}) yield monoidal $G$-natural transformations.  The following observation uses the naive symmetric monoidal $G$-category $\Phi(\A,\gaA,\phiA)$ in \cref{naive_smgcat} and the strictly unital symmetric monoidal $G$-functor $\Phi(f,\actf)$ in \cref{laxmorphism_smf}.

\begin{lemma}\label{BEtransformation}
For lax $\BE$-morphisms between $\BE$-pseudoalgebras in $\Gcat$
\[(f,\actf), (h,\acth) \cn (\A,\gaA,\phiA) \to (\B,\gaB,\phiB),\]
each $\BE$-transformation
\[(f,\actf) \fto{\omega } (h,\acth)\]
yields a monoidal $G$-natural transformation
\[\Phi(f,\actf) \fto{\omega} \Phi(h,\acth).\]
\end{lemma}

\begin{proof}
By \cref{laxm_smf_unit,laxm_smf_mconstraint}, the strictly unital symmetric monoidal $G$-functor $\Phi(f,\actf)$ is given by
\[\Phi(f,\actf) = \big(f,f^2 = \actf_2,f^0 = \actf_0\big) \cn \Phi(\A,\gaA,\phiA) \to \Phi(\B,\gaB,\phiB)\]
and similarly for $\Phi(h,\acth)$.  By \cref{def:algtwocells}, $\omega \cn f \to h$ is a $G$-natural transformation.  It remains to verify the axioms of a monoidal natural transformation in \cref{def:monoidalnattr}.

\parhead{Unity}.  By the basepoint axiom \cref{laxmorphism_basepoint}, we have 
\[\actf_0 = \acth_0 = 1_{\zero^\B}.\] 
Thus, the desired unity axiom for a monoidal natural transformation \cref{monnattr} is equivalent to the equality
\[\omega_{\zero^\A} = 1_{\zero^\B}.\]
This equality holds by \cref{omega-zero}, which is the $n=0$ case of the compatibility axiom \cref{Otransformation_ax} of the $\BE$-transformation $\omega$.

\parhead{Compatibility}.  By construction \cref{laxm_smf_mconstraint}, the monoidal constraints of $\Phi(f,\actf)$ and $\Phi(h,\acth)$ are, respectively, the second action constraints $\actf_2$ and $\acth_2$ for the identity permutation $\id_2 \in \tn\Sigma_2$.  Thus, the desired compatibility axiom of $\omega$ with these monoidal constraints \cref{monnattr} is precisely the commutative diagram \cref{Otransformation_ax} for $x = \id_2 \in \tn\Sigma_2$.
\end{proof}

\cref{BEtransformation} has a converse; see \cref{Phi_twosurject}.

\subsection*{2-Categories of Naive Symmetric Monoidal $G$-Categories}

There is a 2-category of small symmetric monoidal categories, symmetric monoidal functors, and monoidal natural transformations.  The next definition defines the naive $G$-equivariant analogue.

\begin{definition}\label{def:smGcat_twocat}
For a group $G$, we define the 2-category $\smgcat$\index{naive symmetric monoidal G-category@naive symmetric monoidal $G$-category!2-category}\index{2-category!naive symmetric monoidal G-category@naive symmetric monoidal $G$-category} as follows.
\begin{itemize}
\item Objects are naive symmetric monoidal $G$-categories (\cref{def:naive_smGcat}).
\item 1-cells are strictly unital symmetric monoidal $G$-functors (\cref{def:smGfunctor}).
\item 2-cells are monoidal $G$-natural transformations (\cref{def:monGnat}).
\item Horizontal composition of 1-cells is composition of symmetric monoidal functors.  The three $G$-equivariance properties in \cref{def:smGfunctor} are closed under composition.
\item Horizontal and vertical compositions of 2-cells are those of monoidal natural transformations.  The property of being $G$-equivariant is closed under these compositions.
\end{itemize}
Moreover, we define the 2-categories $\smgcatsg$ and $\smgcatst$ with the same objects and 2-cells as $\smgcat$.
\begin{itemize}
\item 1-cells in $\smgcatsg$ are strictly unital strong symmetric monoidal $G$-functors.
\item 1-cells in $\smgcatst$ are strict symmetric monoidal $G$-functors.\defmark
\end{itemize}
\end{definition}

For the Barratt-Eccles reduced $\Gcat$-operad $\BE$ (\cref{def:BE-Gcat}), recall from \cref{oalgps_twocat} that there are 2-categories
\[\AlgstBE \bigsubset \AlgpspsBE \bigsubset \AlglaxBE.\]
In each of them, the objects are $\BE$-pseudoalgebras in $\Gcat$, and the 2-cells are $\BE$-transformations.  The subscripts specify the 1-cells:
\begin{itemize}
\item strict $\BE$-morphisms in $\AlgstBE$, 
\item $\BE$-pseudomorphisms in $\AlgpspsBE$, and
\item lax $\BE$-morphisms in $\AlglaxBE$.  
\end{itemize}
Recall from \cref{def:twofunctor} that a 2-functor between 2-categories consists of assignments on objects, 1-cells, and 2-cells that strictly preserve identity 1-cells, identity 2-cells, horizontal composition of both 1-cells and 2-cells, and vertical composition of 2-cells.  Next, we observe that there are 2-functors from $\BE$-pseudoalgebras to naive symmetric monoidal $G$-categories.

\begin{proposition}\label{alglaxbe_smgcat}
There is a 2-functor\index{2-functor}
\begin{equation}\label{Phi_twofunctor}
\AlglaxBE \fto{\Phi} \smgcat
\end{equation}
given by
\begin{itemize}
\item the object assignment 
\[(\A,\gaA,\phiA) \mapsto \Phi(\A,\gaA,\phiA)\]
in \cref{BEpseudo_smcat},
\item the 1-cell assignment 
\[(f,\actf) \mapsto \Phi(f,\actf) = \big(f, f^2 = \actf_2, f^0 = 1\big)\]
in \cref{laxmorphism_smf}, and
\item the 2-cell assignment 
\[\omega \mapsto \Phi(\omega) = \omega\]
in \cref{BEtransformation}.  
\end{itemize}
Moreover, restricting to sub-2-categories defines 2-functors
\begin{equation}\label{Phi_restrictions}
\begin{split}
& \AlgpspsBE \fto{\Phi} \smgcatsg \andspace\\
& \AlgstBE \fto{\Phi} \smgcatst.
\end{split}
\end{equation}
\end{proposition}

\begin{proof}
We check that $\Phi$ in \cref{Phi_twofunctor} preserves all the 2-categorical structures.

\parhead{2-cells}.  In both 2-categories $\AlglaxBE$ and $\smgcat$, the 2-cells are $G$-natural transformations with extra properties, but not extra structure.  In both cases, horizontal and vertical compositions of 2-cells are defined as those of $G$-natural transformations.  Since the 2-cell assignment is $\omega \mapsto \omega$, we conclude that $\Phi$ preserves identity 2-cells and both horizontal and vertical compositions of 2-cells.

\parhead{Identity 1-cells}.  If $(f,\actf)$ is an identity lax $\BE$-morphism, then $f$ is an identity functor and each $\actf_n$ is an identity natural transformation.  Thus, $(f,\actf_2,1)$ is an identity symmetric monoidal functor.  This shows that $\Phi$ preserves identity 1-cells.  

\parhead{Horizontal composition of 1-cells}.  As we explain in \cref{Omorphism_comp,Omorphism_paste}, the horizontal composition of lax $\BE$-morphisms 
\[\begin{tikzpicture}[xscale=1,yscale=1,vcenter]
\def\h{3}
\draw[0cell=.9]
(0,0) node (a) {\big(\A,\gaA,\phiA\big)}
(a)++(\h,0) node (b) {\big(\B,\gaB,\phiB\big)}
(b)++(\h,0) node (c) {\big(\C,\gaC,\phiC\big)}
;
\draw[1cell=.9]
(a) edge node {(f,\actf)} (b)
(b) edge node {(h,\acth)} (c)
;
\end{tikzpicture}\]
is given by composing the given $G$-functors and pasting the corresponding action constraints for each $n \geq 0$. 
\begin{itemize}
\item For $n=0$, the pasting \cref{Omorphism_paste} yields the action constraint
\[\zero^\C \fto{\acth_0} h(\zero^\B) \fto{h(\actf_0)} hf(\zero^\A).\]
Since $\actf_0 = 1$ and $\acth_0 = 1$, this composite is equal to $1_{\zero^\C}$.  The unit constraint of the composite of $\Phi(f,\actf)$ and $\Phi(h,\acth)$ is also $1_{\zero^\C}$. 
\item For $n=2$, the pasting \cref{Omorphism_paste} yields the following action constraint for objects $\sigma \in \BE(2) = \tn\Sigma_2$ and $a_1,a_2 \in \A$.  
\[\begin{tikzpicture}[xscale=1,yscale=1,vcenter]
\def\h{3} \def\v{1}
\draw[0cell=1]
(0,0) node (a) {\gaC\big(\sigma \sscs hf(a_1), hf(a_2)\big)}
(a)++(3.1,\v) node (b) {h\gaB\big(\sigma \sscs f(a_1), f(a_2)\big)}
(b)++(2.7,-\v) node (c) {hf\gaA\big(\sigma \sscs a_1,a_2\big)}
;
\draw[1cell=.9]
(a) edge[bend left=15] node {\acth_2} (b)
(b) edge[bend left=15] node {h(\actf_2)} (c)
;
\end{tikzpicture}\]
Specifying to $\sigma = \id_2$, this composite is equal to the monoidal constraint of the composite of $\Phi(f,\actf)$ and $\Phi(h,\acth)$. 
\end{itemize}
Thus, $\Phi$ preserves horizontal composition of 1-cells.  

This proves the first assertion about the existence of the 2-functor $\Phi$ in \cref{Phi_twofunctor}.  The restrictions in \cref{Phi_restrictions} are well-defined 2-functors by the second assertion in \cref{laxmorphism_smf}: $\Phi$ sends $\BE$-pseudomorphisms and strict $\BE$-morphisms to, respectively, strictly unital strong and strict symmetric monoidal $G$-functors.
\end{proof}

\section{Naive Symmetric Monoidal $G$-Categories as Pseudoalgebras}
\label{sec:BEpseudoalg}

This section proves that the 2-functor in \cref{alglaxbe_smgcat}
\[\AlglaxBE \fto{\Phi} \smgcat\]
and its two variants are 2-equivalences.  Therefore, the 2-categories $\AlglaxBE$ and $\smgcat$ can be used interchangeably via $\Phi$.  To prove that $\Phi$ is a 2-equivalence, we use the characterization of 2-equivalences in terms of essential surjectivity on objects and fully faithfulness on 1-cells and 2-cells. 

\secoutline

\begin{itemize}
\item \cref{smcat_BEpseudo} proves that each naive symmetric monoidal $G$-category $\A$ yields a $\BE$-pseudoalgebra $\Psi\A$, which has the same underlying $G$-category $\A$.
\item \cref{Psi_Phi} shows that $\Psi$ is a right inverse of $\Phi$, in the sense that $\Phi\Psi$ is the identity assignment on naive symmetric monoidal $G$-categories.  In particular, $\Phi$ is surjective on objects.  However, the other composite, $\Psi\Phi$, is \emph{not} the identity; see \cref{expl:BEpseudo_smcat}.
\item \cref{smf_laxmorphism} extends $\Psi$ to an assignment from strictly unital symmetric monoidal $G$-functors to lax $\BE$-morphisms.  Moreover, if $(f,f^2)$ is strictly unital strong, then $\Psi(f,f^2)$ is a $\BE$-pseudomorphism.  If $(f,f^2)$ is strict, then so is $\Psi(f,f^2)$.
\item \cref{Psi_Phi_smf} proves that $\Phi\Psi$ is the identity assignment on strictly unital symmetric monoidal $G$-functors.  In particular, $\Phi$ is surjective on 1-cells.
\item \cref{Phi_injective_onecells} proves that $\Phi$ is injective on lax $\BE$-morphisms.
\item \cref{Phi_twosurject} proves that $\Phi$ is surjective on 2-cells.
\item \cref{thm:BEpseudoalg} is the main result of this section.  It shows that all three variants of $\Phi$ are 2-equivalences.  We recall 2-equivalences in \cref{def:twoequivalence} and a local characterization in \cref{thm:twoequivalences}.
\end{itemize}

\subsection*{$\BE$-Pseudoalgebras from Naive Symmetric Monoidal $G$-Categories}

\cref{BEpseudo_smcat} shows that $\BE$-pseudoalgebras in $\Gcat$ (\cref{def:pseudoalgebra}) yield naive symmetric monoidal $G$-categories (\cref{def:naive_smGcat}).  The next result provides an assignment in the opposite direction.

\begin{lemma}\label{smcat_BEpseudo}
Each naive symmetric monoidal $G$-category yields a $\BE$-pseudoalgebra, as defined in \cref{Psi_basepoint,Psi_operadicunit,gamma_idn,gammaAn_sigma,gammaAn_morphism,phiA_coherence_iso,phiA_coherence_general} below.
\end{lemma}

\begin{proof}
Suppose we are given a naive symmetric monoidal $G$-category (\cref{def:naive_smGcat})
\[\big(\A,\otimes,\zero,\alpha,\lambda=1,\rho=1,\xi\big).\]
We construct a $\BE$-pseudoalgebra
\begin{equation}\label{BEpseudo_from_smcat}
\Psi(\A,\otimes,\zero,\alpha,\xi) = (\A,\gaA,\phiA)
\end{equation}
as follows.

\parhead{$\BE$-action}.
First, we construct the $\BE$-action $G$-functors 
\[\gaA_n \cn \BE(n) \times \A^n \to \A\]
in \cref{gaAn} as follows.  
\begin{itemize}
\item For $n=0$, the basepoint is defined as the $G$-fixed monoidal unit:
\begin{equation}\label{Psi_basepoint} 
\gaA_0(*) = \zero \in \A.
\end{equation}
\item For $n=1$, the operadic unit $\id_1 \in \BE(1) = \tn\Sigma_1$ acts as the identity:
\begin{equation}\label{Psi_operadicunit}
\gaA_1\big(\id_1 \sscs -\big) = 1_{\A} \cn \A \to \A.
\end{equation}
The action unity axiom \cref{pseudoalg_action_unity} holds by this definition.
\item For $n \geq 2$, we first define the $G$-functor
\begin{equation}\label{gamma_idn}
\gaA_n\big(\id_n \sscs -,\ldots,- \big) = \otimes^n \cn \A^n \to \A
\end{equation}
to be the $n$-fold monoidal product using the left normalized bracketing in \cref{expl:leftbracketing}.  For example, for $n=2$ and $n=3$, we have
\[\begin{split}
\gaA_2\big(\id_2 \sscs -,-\big) = \otimes &\cn \A^2 \to \A \andspace\\
\gaA_3\big(\id_3 \sscs -,-,-\big) = \otimes \circ (\otimes \times 1_{\A}) &\cn \A^3 \to \A.
\end{split}\]
For any permutation $\sigma \in \Sigma_n$, the action equivariance axiom \cref{pseudoalg_action_sym} forces the definition
\begin{equation}\label{gammaAn_sigma}
\gaA_n\big(\sigma \sscs -,\ldots,- \big) = \gaA_n\big(\id_n \sscs \sigma(-,\ldots,-)\big) = \otimes^n \circ \sigma,
\end{equation}
which we call the \emph{$\sigma$-permuted monoidal product}.
\item Consider the unique isomorphism 
\[[\phi,\sigma] \cn \sigma \to \phi \inspace \tn\Sigma_n\]
for permutations $\phi,\sigma \in \Sigma_n$.  We define the natural isomorphism
\begin{equation}\label{gammaAn_morphism}
\gaA_n\big([\phi,\sigma] \sscs -,\ldots,-\big) \cn \otimes^n \circ\, \sigma \fto{\iso} \otimes^n \circ\,\phi
\end{equation}
to be the unique coherence isomorphism in the symmetric monoidal category $(\A,\otimes,\alpha,\xi)$ \cite[XI.1 Theorem 1]{maclane}, from the $\sigma$-permuted monoidal product to the $\phi$-permuted monoidal product.  The coherence isomorphism \cref{gammaAn_morphism} is $G$-equivariant because it is composed of $\otimes$, $\alpha$, and $\xi$, which are all $G$-equivariant.
\end{itemize}
The action equivariance axiom \cref{pseudoalg_action_sym} holds by \cref{gammaAn_sigma,gammaAn_morphism}.

\parhead{Associativity constraint}.  The construction of the associativity constraint $\phiA$ in \cref{phiA} is similar to the construction of $\gaA$ above.  We start with identity permutations and then use equivariance axioms to infer the rest.  The detailed construction is given below.

With the notation in \cref{phiA_component} and $m_\crdot = (m_1,\ldots,m_n)$, we first define the component isomorphism
\[\begin{tikzpicture}[xscale=6,yscale=1,baseline={(a.base)}]
\draw[0cell=.9]
(0,0) node (a) {\gaA_n\big(\id_n \sscs \big\langle \gaA_{m_j}\big(\id_{m_j} \sscs \ang{a_{j,i}}_{i \in \ufs{m}_j} \big) \big\rangle_{j\in \ufs{n}}\big)}
(a)++(1,0) node (b) {\gaA_m\big( \id_m \sscs \ang{\ang{a_{j,i}}_{i \in \ufs{m}_j} }_{j\in \ufs{n}} \big)}
;
\draw[1cell=.9]
(a) edge node {\phiA_{(n;\, m_\crdot)}} node[swap] {\iso} (b)
;
\end{tikzpicture}\]
to be the unique coherence isomorphism 
\begin{equation}\label{phiA_coherence_iso}
\begin{tikzpicture}[xscale=3.5,yscale=1,baseline={(a.base)}]
\draw[0cell=.9]
(0,0) node (a) {\bigotimes_{j=1}^n \left(\bigotimes_{i=1}^{m_j} a_{j,i}\right)}
(a)++(1,0) node (b) {\bigotimes_{j,i} a_{j,i}\phantom{\bigotimes_{j=1}^n}}
;
\draw[1cell=.9]
(a) edge node {\phiA_{(n;\, m_\crdot)}} node[swap] {\iso} (b)
;
\end{tikzpicture} 
\end{equation}
in the strictly unital monoidal category $(\A,\otimes,\alpha)$ \cite[VII.2 Corollary]{maclane}.  In \cref{phiA_coherence_iso}, each of $\bigotimes_{j=1}^n$, $\bigotimes_{i=1}^{m_j}$, and $\bigotimes_{j,i}$ is left normalized (\cref{expl:leftbracketing}), where $j$ and $i$ run through, respectively, $\ufs{n} = \{1,\ldots,n\}$ and $\ufs{m}_j = \{1,\ldots,m_j\}$.  The coherence isomorphism $\phiA_{(n;\, m_\crdot)}$ moves parentheses using $\alpha^{-1}$ without permuting the objects $a_{j,i}$ and using the unit isomorphisms, which are identities, if $m_j = 0$.

Next, for permutations $\sigma \in \Sigma_n$ and $\tau_j \in \Sigma_{m_j}$ for $1 \leq j \leq n$, we define the component isomorphism
\[\begin{tikzpicture}[xscale=5.7,yscale=1,baseline={(a.base)}]
\draw[0cell=.85]
(0,0) node (a) {\gaA_n\big(\sigma \sscs \big\langle \gaA_{m_j}\big(\tau_j \sscs \ang{a_{j,i}}_{i \in \ufs{m}_j} \big) \big\rangle_{j\in \ufs{n}}\big)}
(a)++(1,0) node (b) {\gaA_m\big( \ga(\sigma \sscs \ang{\tau_j}_{j\in \ufs{n}} ) \sscs \ang{\ang{a_{j,i}}_{i \in \ufs{m}_j} }_{j\in \ufs{n}} \big)}
;
\draw[1cell=.85]
(a) edge node {\phiA_{(n;\, m_\crdot)}} node[swap] {\iso} (b)
;
\end{tikzpicture}\]
using \cref{phiA_coherence_iso}, the top equivariance axiom \cref{pseudoalg_topeq}, and the bottom equivariance axiom \cref{pseudoalg_boteq}.  More precisely, it is the following unique coherence isomorphism in the strictly unital monoidal category $(\A,\otimes,\alpha)$ that moves parentheses, where $m_{\sigmainv} = \bang{m_{\sigmainv(j)}}_{j\in \ufs{n}}$.  
\begin{equation}\label{phiA_coherence_general}
\begin{tikzpicture}[xscale=5.5,yscale=1,baseline={(a.base)}]
\draw[0cell=.9]
(0,0) node (a) {\bigotimes_{j=1}^n \left(\bigotimes_{i=1}^{m_{\sigmainv(j)}} a_{\sigmainv(j),\tauinv_{\sigmainv(j)}(i)}\right)}
(a)++(1,0) node (b) {\bigotimes_{j,i} a_{\sigmainv(j),\tauinv_{\sigmainv(j)}(i)}\phantom{\bigotimes_{j=1}^n}}
;
\draw[1cell=.9]
(a) edge node {\phiA_{(n;\, m_{\sigmainv})}} node[swap] {\iso} (b)
;
\end{tikzpicture} 
\end{equation}
In the codomain in \cref{phiA_coherence_general}, $j$ and $i$ run through, respectively, $\{1,\ldots,n\}$ and $\{1,\ldots,m_{\sigmainv(j)}\}$.  The coherence isomorphism \cref{phiA_coherence_general} is a $G$-natural isomorphism because it is composed of $\otimes$ and $\alpha^{-1}$, which are $G$-equivariant.  Moreover, the top and bottom equivariance axioms, \cref{pseudoalg_topeq,pseudoalg_boteq}, hold by definition \cref{phiA_coherence_general}.  Next, we verify the rest of the axioms of a $\BE$-pseudoalgebra in \cref{def:pseudoalgebra}.

\parhead{Basepoint axiom}.  
Suppose we are given a permutation $\sigma \in \Sigma_n$ and objects 
\[\anga = \big(a_1,\ldots,a_{j-1},a_{j+1},\ldots,a_n\big) \in \A^{n-1}.\]  
We define the objects 
\[\begin{split}
a_j &= \gaA_0(*) = \zero \in \A,\\
\dy_j\anga &= \big(a_1,\ldots,a_{j-1}, a_j, a_{j+1},\ldots,a_n\big) \in \A^n
\end{split}\]  
and the block permutation \cref{blockpermutation}
\[\dy_j\sigma = \ga\Big(\sigma \sscs \id_1^{j-1}, \id_0 \,, \id_1^{n-j} \Big)
= \sigma\bang{1^{j-1}, 0, 1^{n-j}} \in \Sigma_{n-1}.\]
Then the associativity constraint in the basepoint axiom \cref{pseudoalg_basept_axiom} for $(\A,\gaA,\phiA)$ is the unique coherence isomorphism that moves parentheses
\begin{equation}\label{phiA_nonezeroone}
\begin{tikzpicture}[xscale=1,yscale=1,baseline={(a.base)}]
\def\v{-1}
\draw[0cell=1]
(0,0) node (a) {\gaA_n\big(\sigma \sscs \dy_j\anga\big)}
(a)++(5.2,0) node (b) {\gaA_{n-1}\big(\dy_j\sigma \sscs \anga \big)}
(a)++(0,\v) node (a2) {\txotimes_{\ell=1}^n a_{\sigmainv(\ell)}}
(b)++(0,\v) node (b2) {\txotimes_{k=1}^{n-1} a_{(\dy_j\sigma)^{-1}(k) + \varepsilon}}
;
\draw[1cell=.9]
(a) edge node {\phiA_{(n;\, 1^{j-1}, 0, 1^{n-j})}} (b)
(a) edge[equal] (a2)
(b) edge[equal, shorten >=-.5ex] (b2)
;
\end{tikzpicture}
\end{equation}
where
\[\varepsilon = \begin{cases}
0 & \text{ if $1\leq (\dy_j\sigma)^{-1}(k) \leq j-1$,}\\
1 & \text{ if $j \leq (\dy_j\sigma)^{-1}(k) \leq n-1$.}
\end{cases}\]
The unique coherence isomorphism in \cref{phiA_nonezeroone} is the identity morphism for the following two reasons.
\begin{itemize}
\item $a_j = \zero$ is the strict monoidal unit of $(\A,\otimes,\zero)$.
\item The way $\sigma$ permutes $\dy_j\anga \in \A^n$ restricts to the way the block permutation $\dy_j\sigma$ permutes $\anga \in \A^{n-1}$.  Thus, the domain and codomain of $\phiA_{(n;\, 1^{j-1}, 0, 1^{n-j})}$ are the same object.
\end{itemize}

\parhead{Unity axiom}.  In \cref{phiA_coherence_general}, suppose either 
\begin{itemize}
\item each $m_j = 1$, which implies $\tau_j = \id_1 \in \Sigma_1$, or
\item $n=1$, which implies $\sigma = \id_1 \in \Sigma_1$.
\end{itemize}
Then the parenthesis-moving coherence isomorphism $\phiA_{(n;\, m_{\sigmainv})}$ is the identity because its domain and codomain are the same object.  This proves the unity axiom \cref{pseudoalg_unity}.

\parhead{Composition axiom}.  In the diagram \cref{pseudoalg_comp_axiom} for the composition axiom, each of the top and bottom composites is a parenthesis-moving coherence isomorphism in $(\A,\otimes,\alpha)$.  Therefore, they are equal by the uniqueness of coherence isomorphisms in monoidal categories \cite[VII.2 Corollary]{maclane}.  

This finishes the proof that each naive symmetric monoidal $G$-category $(\A,\otimes,\zero,\alpha,\xi)$ yields a $\BE$-pseudoalgebra $(\A,\gaA,\phiA)$.
\end{proof}

Next, we consider the composite of the constructions $\Psi$ in \cref{BEpseudo_from_smcat} and $\Phi$ in \cref{naive_smgcat}.

\begin{lemma}\label{Psi_Phi}
Suppose $(\A,\otimes,\zero,\alpha,\xi)$ is a naive symmetric monoidal $G$-category.  Then
\[\Phi \Psi(\A,\otimes,\zero,\alpha,\xi) = (\A,\otimes,\zero,\alpha,\xi)\]
as naive symmetric monoidal $G$-categories.
\end{lemma}

\begin{proof}
We abbreviate $\Phi\Psi(\A,\otimes,\zero,\alpha,\xi)$ to $\Phi\Psi\A$.
\begin{itemize}
\item By \cref{naive_smgcat,BEpseudo_from_smcat}, the underlying $G$-category of $\Phi\Psi\A$ is $\A$.
\item By \cref{naive_monunit,Psi_basepoint}, the monoidal unit of $\Phi\Psi\A$ is $\zero \in \A$.
\item By \cref{naive_monproduct,gamma_idn}, the monoidal product of $\Phi\Psi\A$ is $\otimes$.
\item By \cref{naive_unit_isos}, the left and right unit isomorphisms of $\Phi\Psi\A$ are identities.
\item By \cref{naive_braiding,gammaAn_morphism}, the braiding of $\Phi\Psi\A$ is the unique coherence isomorphism 
\[\otimes \to \otimes \circ \tau\]
in $(\A,\otimes,\xi)$.  This is equal to $\xi$ by uniqueness \cite[XI.1 Theorem 1]{maclane}.
\item By \cref{naive_associativity,phiA_coherence_iso} the associativity isomorphism of $\Phi\Psi\A$ is the unique parenthesis-moving coherence isomorphism 
\[\otimes \circ (\otimes \times 1) \to \otimes \circ (1 \times \otimes)\]
in $(\A,\otimes,\alpha)$.  This is equal to $\alpha$ by uniqueness \cite[VII.2 Corollary]{maclane}.
\end{itemize}
This proves that $\Phi\Psi\A$ is equal to $\A$ as $G$-categories and as symmetric monoidal categories.  Thus, they are equal as naive symmetric monoidal $G$-categories.
\end{proof}

\begin{explanation}[The Other Composite]\label{expl:BEpseudo_smcat}
\cref{Psi_Phi} shows that $\Phi\Psi$ is the identity assignment on naive symmetric monoidal $G$-categories.  Conversely, for each $\BE$-pseudoalgebra $(\A,\gaA,\phiA)$, \cref{BEpseudo_smcat,smcat_BEpseudo} show that there is a $\BE$-pseudoalgebra $\Psi\Phi(\A,\gaA,\phiA)$, which is \emph{not} equal to the original $(\A,\gaA,\phiA)$.  To see this, let us denote
\[\Psi\Phi(\A,\gaA,\phiA) = (\A,\ga',\vphi')\]
and observe that $\ga'$ is generally different from $\gaA$.  In more detail, by \cref{naive_monproduct,gamma_idn}, for objects $a,b,c \in \A$ we have the object
\begin{equation}\label{gaprime_three}
\ga'_3\big(\id_3 \sscs a,b,c\big) = \gaA_2\Big(\id_2 \sscs \gaA_2 \big(\id_2 \sscs a,b\big), c\Big).
\end{equation}
This object is isomorphic to $\gaA_3\big(\id_3 \sscs a,b,c\big)$ via the associativity constraint \cref{phiA_component}
\[\gaA_2\Big(\id_2 \sscs \gaA_2 \big(\id_2 \sscs a,b\big), c\Big) \fto[\iso]{\phiA_{2;\, 2,1}} 
\gaA_3\big(\id_3 \sscs a,b,c\big),\]
which is not the identity in general.  Thus, $\Psi\Phi\A$ is not equal to $\A$.  

The fact that $\Psi\Phi\A \neq \A$ contradicts a claim in \cite[page 227]{gmmo20}.  It is incorrectly claimed there that $\Phi$ and $\Psi$ constitute a bijective correspondence between small strictly unital symmetric monoidal categories and $\BE$-pseudoalgebras in $\Cat$.  In practice, however, it does not matter that $\Phi$ and $\Psi$ are not strictly inverse to each other.  As we show in \cref{thm:BEpseudoalg} below, the 2-functor 
\[\AlglaxBE \fto{\Phi} \smgcat\]
in \cref{Phi_twofunctor} is a 2-equivalence.  Thus, the 2-category $\AlglaxBE$, whose objects are $\BE$-pseudoalgebras, and the 2-category $\smgcat$, whose objects are naive symmetric monoidal $G$-categories, are interchangeable via $\Phi$.
\end{explanation}

\subsection*{Lax $\BE$-Morphisms from Symmetric Monoidal $G$-Functors}

\cref{laxmorphism_smf} shows that lax $\BE$-morphisms between $\BE$-pseudoalgebras in $\Gcat$ (\cref{def:laxmorphism}) yield strictly unital symmetric monoidal $G$-functors (\cref{def:smGfunctor}).  The next result provides an assignment in the opposite direction, using the $\BE$-pseudoalgebra $\Psi\A$ in \cref{smcat_BEpseudo}.

\begin{lemma}\label{smf_laxmorphism}
Each strictly unital symmetric monoidal $G$-functor between naive symmetric monoidal $G$-categories
\[\big(f,f^2,f^0 = 1\big) \cn (\A,\otimes,\zero,\alpha,\xi) \to (\B,\otimes,\zero,\alpha,\xi)\]
yields a lax $\BE$-morphism between $\BE$-pseudoalgebras in $\Gcat$
\begin{equation}\label{Psi_fftwo}
\Psi(f,f^2) = (f,\actf) \cn \Psi\A \to \Psi\B.
\end{equation}
Moreover, if $(f,f^2)$ is strictly unital strong, then $\Psi(f,f^2)$ is a $\BE$-pseudomorphism.  If $(f,f^2)$ is strict, then $\Psi(f,f^2)$ is a strict $\BE$-morphism.
\end{lemma}

\begin{proof}
We refer to the $\BE$-pseudoalgebra in \cref{BEpseudo_from_smcat} 
\[\Psi\A = (\A,\gaA,\phiA)\]
for both $\A$ and $\B$.  

\parhead{Action constraints}.  We construct the action constraint $G$-natural transformations in \cref{laxmorphism_constraint}
\begin{equation}\label{Psif_constraint}
\begin{tikzpicture}[xscale=4,yscale=1.3,vcenter]
\draw[0cell=.9]
(0,0) node (x11) {\BE(n) \times (\Psi\A)^n}
(x11)++(1,0) node (x12) {\BE(n) \times (\Psi\B)^n}
(x11)++(0,-1) node (x21) {\Psi\A}
(x21)++(1,0) node (x22) {\Psi\B}
;
\draw[1cell=.9]  
(x11) edge node {1 \times f^n} (x12)
(x12) edge node {\gaB_n} (x22)
(x11) edge node[swap] {\gaA_n} (x21)
(x21) edge node[swap] {f} (x22)
;
\draw[2cell]
node[between=x11 and x22 at .5, shift={(0,0)}, rotate=-90, 2label={above,\actf_n}] {\Rightarrow}
;
\end{tikzpicture}
\end{equation}
as follows.
\begin{itemize}
\item For $n=0$, we define
\[\actf_0 = f^0 = 1_{\zero^\B} \cn \zero^\B \to f(\zero^\A).\]
This is well defined by \cref{Psi_basepoint}.  The basepoint axiom \cref{laxmorphism_basepoint} holds by this definition.
\item For $n=1$, we define
\[\actf_1 = 1_{f(a)} \cn \gaB_1\big(\id_1 \sscs f(a)\big) \to f\big(\gaA_1(\id_1 \sscs a)\big)\]
for objects $a \in \A$.  This is well defined by \cref{Psi_operadicunit}.  The unity axiom \cref{laxmorphism_unity} holds by this definition.
\item For $n \geq 2$, we first define the following $\id_n$-component of $\actf_n$ for objects $a_1,\ldots,a_n \in \A$.
\begin{equation}\label{Psif_constraint_id}
\begin{tikzpicture}[xscale=4,yscale=1.3,vcenter]
\draw[0cell=.9]
(0,0) node (x11) {\gaB_n\Big( \id_n \sscs \ang{fa_j}_{j\in \ufs{n}} \Big)}
(x11)++(1,0) node (x12) {f\Big(\gaA_n\big(\id_n \sscs \ang{a_j}_{j\in \ufs{n}} \big)\Big)}
(x11)++(0,-1) node (x21) {\txotimes_{j=1}^n fa_j}
(x21)++(1,0) node (x22) {f\left(\txotimes_{j=1}^n a_j\right)}
;
\draw[1cell=.9]  
(x11) edge node {\actf_n} (x12)
(x21) edge node {f^2} (x22)
(x11) edge[equal] (x21)
(x12) edge[equal] (x22)
;
\end{tikzpicture}
\end{equation}
In \cref{Psif_constraint_id}, each $\txotimes_{j=1}^n$ is left normalized (\cref{expl:leftbracketing}), and $f^2$ means the iterate of the monoidal constraint $f^2$ that moves $f$ outside.  More precisely, $\actf_2$ is an instance of $f^2$, and $\actf_n$ for $n \geq 3$ are defined by induction.  For example, $\actf_3$ is the composite 
\[(fa_1 \otimes fa_2) \otimes fa_3 \fto{f^2 \otimes 1} f(a_1 \otimes a_2) \otimes fa_3 \fto{f^2} 
f\big((a_1 \otimes a_2) \otimes a_3\big).\]
\item For a permutation $\sigma \in \Sigma_n$, the $\id_n$-component \cref{Psif_constraint_id} and the desired equivariance axiom \cref{laxmorphism_equiv} force the following definition of $\actf_n$.
\begin{equation}\label{Psif_constraint_sigma}
\begin{tikzpicture}[xscale=4,yscale=1.3,vcenter]
\draw[0cell=.9]
(0,0) node (x11) {\gaB_n\Big( \sigma \sscs \ang{fa_j}_{j\in \ufs{n}} \Big)}
(x11)++(1,0) node (x12) {f\Big(\gaA_n\big(\sigma \sscs \ang{a_j}_{j\in \ufs{n}} \big)\Big)}
(x11)++(0,-1) node (x21) {\txotimes_{j=1}^n fa_{\sigmainv(j)}}
(x21)++(1,0) node (x22) {f\left(\txotimes_{j=1}^n a_{\sigmainv(j)}\right)}
;
\draw[1cell=.9]  
(x11) edge node {\actf_n} (x12)
(x21) edge node {f^2} (x22)
(x11) edge[equal] (x21)
(x12) edge[equal] (x22)
;
\end{tikzpicture}
\end{equation}
Since the monoidal constraint $f^2$ is a $G$-natural transformation, so is $\actf_n$.  Moreover, the equivariance axiom \cref{laxmorphism_equiv} holds by definition \cref{Psif_constraint_sigma}.
\end{itemize}

\parhead{Associativity axiom}.  It remains to verify the associativity axiom \cref{laxmorphism_associativity} for $(f,\actf)$.  By construction \cref{phiA_coherence_general}, each of $\phiA$ and $\phiB$ is a parenthesis-moving coherence isomorphism.  Each instance of $\actf$ is an iterate of the monoidal constraint $f^2$ that moves $f$ outside, as defined in \cref{Psif_constraint_sigma}.  Therefore, the two composites in the desired associativity diagram \cref{laxmorphism_associativity} are equal by Epstein's Coherence Theorem for monoidal functors \cite{epstein}.  This theorem generalizes the associativity axiom \cref{monoidalfunctorassoc} for monoidal functors.  For example, the case of \cref{laxmorphism_associativity} with 
\[x = x_2 = \id_2 \in \Sigma_2 \andspace x_1 = \id_1 \in \Sigma_1\]
holds by the associativity axiom \cref{monoidalfunctorassoc}.  This finishes the proof that 
\[\Psi(f,f^2) = (f,\actf) \cn \Psi\A \to \Psi\B\] 
is a lax $\BE$-morphism.  

\parhead{Other cases}.  For the other assertions, suppose $(f,f^2)$ is strictly unital strong, so $f^2$ is invertible.  Thus, $\actf_n$ in \cref{Psif_constraint_sigma}, which is an iterate of $f^2$, is also invertible.  This means that $\Psi(f,f^2)$ is a $\BE$-pseudomorphism.  Similarly, if $(f,f^2)$ is strict, then $f^2$ is the identity.  Thus, $\actf_n$ is the identity, and $\Psi(f,f^2)$ is a strict $\BE$-morphism.
\end{proof}

Next, we consider the composite of the constructions $\Psi$ in \cref{Psi_fftwo} and $\Phi$ in \cref{Phi_factf}.

\begin{lemma}\label{Psi_Phi_smf}
For each strictly unital symmetric monoidal $G$-functor between naive symmetric monoidal $G$-categories
\[\big(f,f^2,f^0 = 1\big) \cn (\A,\otimes,\zero,\alpha,\xi) \to (\B,\otimes,\zero,\alpha,\xi),\]
there is an equality
\[\Phi\Psi(f,f^2) = (f,f^2)\]
of strictly unital symmetric monoidal $G$-functors.
\end{lemma}

\begin{proof}
\cref{Psi_Phi} yields the equalities of naive symmetric monoidal $G$-categories
\[\Phi\Psi\A = \A \andspace \Phi\Psi\B = \B.\]
Thus, the two strictly unital symmetric monoidal $G$-functors in question, $\Phi\Psi(f,f^2)$ and $(f,f^2)$, have the same domain, $\A$, and the same codomain, $\B$.
\begin{itemize}
\item By \cref{Phi_factf,Psi_fftwo}, the underlying $G$-functor of $\Phi\Psi(f,f^2)$ is the same as that of $\Psi(f,f^2)$, which is $f$.  
\item By \cref{laxm_smf_mconstraint} the monoidal constraint of $\Phi\Psi(f,f^2)$ is the second action constraint of $\Psi(f,f^2)$ for $\id_2 \in \Sigma_2$, which is $f^2$ by \cref{Psif_constraint_id}.
\end{itemize}
Therefore, $\Phi\Psi(f,f^2)$ is equal to $(f,f^2)$.
\end{proof}

\cref{Psi_Phi_smf} implies that the assignment $\Phi$ in \cref{Phi_factf} is surjective onto strictly unital symmetric monoidal $G$-functors.  Next, we observe that $\Phi$ is also injective on 1-cells.

\begin{lemma}\label{Phi_injective_onecells}
Given two lax $\BE$-morphisms between $\BE$-pseudoalgebras in $\Gcat$
\[(f,\actf), (h,\acth) \cn (\A,\gaA,\phiA) \to (\B,\gaB,\phiB),\]
if there is an equality of strictly unital symmetric monoidal $G$-functors
\[\Phi(f,\actf) = \Phi(h,\acth) \cn \Phi\A \to \Phi\B,\]
then there is an equality
\[(f,\actf) = (h,\acth)\]
of lax $\BE$-morphisms.
\end{lemma}

\begin{proof}
First, we observe the following.
\begin{itemize}
\item By the basepoint axiom \cref{laxmorphism_basepoint}, $\actf_0$ is the identity morphism $1_{\zero^\B}$.
\item By the unity axiom \cref{laxmorphism_unity}, $\actf_1$ is the identity.
\item By \cref{Phi_factf,laxm_smf_mconstraint}, the underlying $G$-functor of $\Phi(f,\actf)$ is $f$, and its monoidal constraint is the component of $\actf_2$ for $\id_2 \in \Sigma_2$.
\end{itemize}
Thus, it suffices to show that a lax $\BE$-morphism $(f,\actf)$ is uniquely determined by the data in the above list and the associativity constraints of $\A$ and $\B$.  

Moreover, by the equivariance axiom \cref{laxmorphism_equiv}, for each $n \geq 2$, the action constraint $\actf_n$ is uniquely determined by its component for $\id_n \in \Sigma_n$.  Thus, it suffices to show that $\actf_n$ for each $n \geq 2$ and $\id_n \in \Sigma_n$ is determined by the list of data in the previous paragraph and the associativity constraints of $\A$ and $\B$.  We prove this assertion by induction, starting with the base case $n=2$.  For the induction step, we use the object equality
\[\id_n = \ga\big(\id_2 \sscs \id_{n-1}, \id_1\big) \in \BE(n) = \tn\Sigma_n\]
and consider the permutations
\[x = \id_2 \in \BE(2), \quad x_1 = \id_{n-1} \in \BE(n-1), \andspace x_2 = \id_1 \in \BE(1)\]
and objects $\anga = \ang{a_j}_{j\in \ufs{n}} \in \A^n$,
\[\anga' = \ang{a_j}_{j=1}^{n-1}, \quad f\anga = \ang{fa_j}_{j\in \ufs{n}}, \andspace f\anga' = \ang{fa_j}_{j=1}^{n-1}.\]
Using the invertibility of $\phiB$, the associativity axiom \cref{laxmorphism_associativity} yields the following factorization of $\actf_n$ for $\id_n \in \BE(n)$.
\begin{equation}\label{actfn_factorization}
\begin{tikzpicture}[xscale=1,yscale=1,vcenter]
\def\h{3} \def\g{0} \def\v{1.4} \def\u{-1}
\draw[0cell=.8]
(0,0) node (x11) {\gaB_2\Big(\id_2 \sscs \gaB_{n-1}\big( \id_{n-1} \sscs f\anga' \big), fa_n\Big)}
(x11)++(\h,\v) node (x12) {\gaB_n \big(\id_n \sscs f\anga \big)}
(x11)++(\h+.5,\u) node (x) {f \gaA_n \big(\id_n \sscs \anga \big)}
(x11)++(0,2*\u) node (x21) {\gaB_2\Big(\id_2 \sscs f\gaA_{n-1}\big( \id_{n-1} \sscs \anga'\big), fa_n\Big)}
(x21)++(\h,-\v) node (x22) {f\gaA_2\Big(\id_2 \sscs \gaA_{n-1}\big( \id_{n-1} \sscs \anga'\big), a_n\Big)}
;
\draw[1cell=.85]  
(x12) edge[transform canvas={xshift={0em}}, shorten >=0pt] node[swap] {(\phiB)^{-1}} (x11)
(x12) edge node[pos=.6] {\actf_n} (x)
(x11) edge[transform canvas={xshift={1em}}] node[swap] {\gaB_2\big( 1 \sscs \actf_{n-1}, \actf_1 \big)} (x21)
(x21) edge[transform canvas={xshift={0em}}, shorten >=0pt] node[swap,pos=.1] {\actf_2} (x22)
(x22) edge node[swap,pos=.6] {f\phiA} (x)
;
\end{tikzpicture}
\end{equation}
The induction hypothesis applies to $\actf_{n-1}$, finishing the induction.
\end{proof}

\subsection*{$\BE$-Transformations from Monoidal $G$-Natural Transformations}
\cref{BEtransformation} shows that each $\BE$-transformation is also a monoidal $G$-natural transformation via the construction $\Phi$ on objects and 1-cells.  The following observation establishes the converse and shows that $\Phi$ is surjective on 2-cells.

\begin{lemma}\label{Phi_twosurject}
Given lax $\BE$-morphisms between $\BE$-pseudoalgebras in $\Gcat$
\[(f,\actf), (h,\acth) \cn (\A,\gaA,\phiA) \to (\B,\gaB,\phiB),\]
each monoidal $G$-natural transformation
\[\Phi(f,\actf) = \big(f,\actf_2,\actf_0 = 1\big) \fto{\psi} \Phi(h,\acth) = \big(h,\acth_2,\acth_0=1\big)\]
yields a $\BE$-transformation
\[(f,\actf) \fto{\psi} (h,\acth).\]
\end{lemma}

\begin{proof}
By \cref{def:monGnat}, $\psi$ is a $G$-natural transformation.  Next, we check the compatibility axiom \cref{Otransformation_ax} for $n \geq 0$.  
\begin{itemize}
\item As we discuss in \cref{omega-zero}, the compatibility axiom \cref{Otransformation_ax} for $n=0$ means the equality
\[\psi_{\zero^\A} = 1_{\zero^\B}.\]
This equality holds by (i) the unity axiom \cref{monnattr} of the monoidal natural transformation $\psi$ and (ii) the fact that the unit constraints of $\Phi(f,\actf)$ and $\Phi(h,\acth)$ are identities.
\item Since $\BE(1) \iso \boldone$, the axiom \cref{Otransformation_ax} for $n=1$ holds by the action unity axiom \cref{pseudoalg_action_unity} of $\BE$-pseudoalgebras and the equalities
\[\actf_1 = 1 \andspace \acth_1 = 1.\] 
The last two equalities hold by the unity axiom \cref{laxmorphism_unity} of lax $\BE$-morphisms.
\end{itemize}

By the equivariance axiom \cref{laxmorphism_equiv} of lax $\BE$-morphisms, the axiom \cref{Otransformation_ax} for $n \geq 2$ holds if and only if they hold for $\id_n \in \Sigma_n$.  We prove the axiom \cref{Otransformation_ax} for $n \geq 2$ and $\id_n \in \Sigma_n$ by induction.  For $n=2$ and $\id_2 \in \Sigma_2$, the axiom \cref{Otransformation_ax} holds by the compatibility axiom \cref{monnattr} of the monoidal natural transformation $\psi$.

For the induction step, we use the factorization \cref{actfn_factorization} for both $\actf_n$ and $\acth_n$.  The diagram \cref{Otransformation_ax} for $\id_n$ is the boundary of the following diagram.
\begin{equation}\label{Ptwocell_axiom}
\begin{tikzpicture}[xscale=1,yscale=1,vcenter]
\def\v{-1.5} \def\a{15}
\draw[0cell=.8]
(0,0) node (x1) {\gaB_n\big(\id_n \sscs f\anga\big)}
(x1)++(0,\v) node (x2) {\gaB_2\Big(\id_2 \sscs \gaB_{n-1} \big( \id_{n-1} \sscs f\anga'\big), fa_n \Big)}
(x2)++(0,\v) node (x3) {\gaB_2\Big(\id_2 \sscs f\gaA_{n-1} \big( \id_{n-1} \sscs \anga'\big), fa_n \Big)}
(x3)++(0,\v) node (x4) {f\gaA_2\Big(\id_2 \sscs \gaA_{n-1} \big( \id_{n-1} \sscs \anga'\big), a_n \Big)}
(x4)++(0,\v) node (x5) {f\gaA_n\big( \id_n \sscs \anga\big)}
(x1)++(5,0) node (y1) {\gaB_n\big(\id_n \sscs h\anga\big)}
(y1)++(0,\v) node (y2) {\gaB_2\Big(\id_2 \sscs \gaB_{n-1} \big( \id_{n-1} \sscs h\anga'\big), ha_n \Big)}
(y2)++(0,\v) node (y3) {\gaB_2\Big(\id_2 \sscs h\gaA_{n-1} \big( \id_{n-1} \sscs \anga'\big), ha_n \Big)}
(y3)++(0,\v) node (y4) {h\gaA_2\Big(\id_2 \sscs \gaA_{n-1} \big( \id_{n-1} \sscs \anga'\big), a_n \Big)}
(y4)++(0,\v) node (y5) {h\gaA_n\big( \id_n \sscs \anga\big)}
;
\draw[1cell=.7]
(x1) edge node[swap] {(\phiB)^{-1}} (x2)
(x2) edge node[swap] {\gaB_2(1 \sscs \actf_{n-1}, \actf_1)} (x3)
(x3) edge node[swap] {\actf_2} (x4)
(x4) edge node[swap] {f\phiA} (x5)
(y1) edge node {(\phiB)^{-1}} (y2)
(y2) edge node {\gaB_2(1 \sscs \acth_{n-1}, \acth_1)} (y3)
(y3) edge node {\acth_2} (y4)
(y4) edge node {h\phiA} (y5)
(x1) edge node {\gaB_n\big(1 \sscs \psi\anga\big)} (y1)
(x2) edge[bend left=\a] node {\gaB_2\big(1 \sscs \gaB_{n-1}(1 \sscs \psi\anga') , \psi\big)} (y2)
(x3) edge[bend left=\a] node {\gaB_2\big(1 \sscs \psi,\psi\big)} (y3)
(x4) edge node {\psi} (y4)
(x5) edge node {\psi} (y5)
;
\end{tikzpicture}
\end{equation}
From top to bottom, the four regions in \cref{Ptwocell_axiom} commute for the following reasons.
\begin{itemize}
\item The top region commutes by the naturality of the associativity constraint $\phiB$.
\item The second region commutes by the functoriality of $\gaB_2$, the induction hypothesis for the case $n-1$, and the unity axiom \cref{laxmorphism_unity}.
\item The third region commutes by the compatibility axiom \cref{monnattr} of the monoidal natural transformation $\psi$.
\item The bottom region commutes by the naturality of $\psi$.
\end{itemize}
This finishes the induction.
\end{proof}

\subsection*{2-Equivalences}

\cref{alglaxbe_smgcat} shows that there is a 2-functor 
\[\AlglaxBE \fto{\Phi} \smgcat\]
from the 2-category $\AlglaxBE$ (\cref{oalgps_twocat}) ,with
\begin{itemize}
\item $\BE$-pseudoalgebras in $\Gcat$ as objects,
\item lax $\BE$-morphisms as 1-cells, and
\item $\BE$-transformations as 2-cells,
\end{itemize}
to the 2-category $\smgcat$ (\cref{def:smGcat_twocat}), with
\begin{itemize}
\item naive symmetric monoidal $G$-categories as objects,
\item strictly unital symmetric monoidal $G$-functors as 1-cells, and
\item monoidal $G$-natural transformations as 2-cells.
\end{itemize}
Now we observe that $\Phi$ is a 2-equivalence in the sense of \cref{def:twoequivalence}.  Moreover, the pseudo and strict variants (for 1-cells) are also true.  Recall that $G$ is an arbitrary group, and $\BE$ is the Barratt-Eccles reduced $\Gcat$-operad in \cref{def:BE-Gcat}. 

\begin{theorem}\label{thm:BEpseudoalg}
The 2-functors
\[\begin{split}
\Phi &\cn \AlglaxBE \to \smgcat,\\
\Phi &\cn \AlgpspsBE \to \smgcatsg, \andspace\\
\Phi &\cn \AlgstBE \to \smgcatst
\end{split}\]
in \cref{Phi_twofunctor,Phi_restrictions} are \index{2-equivalence}2-equivalences.
\end{theorem}

\begin{proof}
Suppose $\Phi$ is any one of the three variant 2-functors in the statement above. 
\begin{itemize}
\item \cref{Psi_Phi} shows that $\Phi$ is surjective, and hence essentially surjective, on objects.
\item $\Phi$ is surjective on 1-cells by \cref{smf_laxmorphism,Psi_Phi_smf}.
\item $\Phi$ is injective on 1-cells by \cref{Phi_injective_onecells}.
\item $\Phi$ is injective on 2-cells by \cref{BEtransformation}.
\item $\Phi$ is surjective on 2-cells by \cref{Phi_twosurject}.
\end{itemize}
Thus, \cref{thm:twoequivalences} implies that $\Phi$ is a 2-equivalence.
\end{proof}

\begin{remark}[Non-Symmetric Variants]\label{rk:Phi_nonsymmetric}
We emphasize that each 2-equivalence $\Phi$ in \cref{thm:BEpseudoalg} is not an isomorphism because it is not bijective on objects; see \cref{expl:BEpseudo_smcat}.  Nonsymmetric analogues of \cref{thm:BEpseudoalg} for the lax and pseudo cases are given in \cite[Corollary 3.2.5]{leinster}, whose proofs proceed somewhat differently from ours given in \cref{sec:naive_smc} and this section.
\end{remark}

\section{Genuine Symmetric Monoidal $G$-Categories}
\label{sec:genuine_smgcat}

This section recalls the definition of a genuine symmetric monoidal $G$-category and explains the associated notions of lax morphisms and transformations.

\secoutline
\begin{itemize}
\item \cref{def:GBE_pseudoalg} defines genuine symmetric monoidal $G$-categories.
\item \cref{expl:GBEvsBE_pseudo} discusses why there are no analogues of \cref{thm:BEpseudoalg} for genuine symmetric monoidal $G$-categories.
\item \cref{naive_genuine_smgcat} records the fact that each $\BE$-pseudoalgebra yields a genuine symmetric monoidal $G$-category by applying $\Catg(\EG,-)$.
\item \cref{expl:GBE_pseudoalg} unpacks the structure of a genuine symmetric monoidal $G$-category in detail.
\item \cref{expl:GBE_laxmorphism,expl:GBE_transformation} unpack the notions of lax morphisms and transformations for genuine symmetric monoidal $G$-categories.
\end{itemize}

\recollection
For an arbitrary group $G$, recall the $G$-Barratt-Eccles reduced $\Gcat$-operad $\GBE$ in \cref{def:GBE}.  Its component $G$-categories are translation categories (\cref{GBEn})
\[\GBE(n) = \Catg(\EG,\ESigma_n) \iso \tn[G,\Sigma_n] \forspace n \geq 0.\]
Its $\Gcat$-operad structure is induced by the one on the Barratt-Eccles operad $\BE$ (\cref{def:BE-Gcat}) via the product-preserving functor $\Catg(\EG,-)$.  Since $\GBE$ is reduced, which means $\GBE(0) = \boldone$, and levelwise a translation category, it admits a unique pseudo-commutative structure (\cref{GBE_pseudocom}).  Its algebras in $\Gcat$ are genuine permutative $G$-categories (\cref{def:GBE_algebra}).

Recall $\Op$-pseudoalgebras in \cref{def:pseudoalgebra}.  The following definition is \cite[Definition 0.3]{gmmo20}.

\begin{definition}\label{def:GBE_pseudoalg}
For the $G$-Barratt-Eccles $\Gcat$-operad $\GBE$ in \cref{def:GBE}, $\GBE$-pseudoalgebras in $\Gcat$ are called \index{G-Barratt-Eccles operad@$G$-Barratt-Eccles operad!pseudoalgebra}\index{pseudoalgebra!G-Barratt-Eccles operad@$G$-Barratt-Eccles operad}\index{G-category@$G$-category!genuine symmetric monoidal}\index{genuine symmetric monoidal G-category@genuine symmetric monoidal $G$-category}\emph{genuine symmetric monoidal $G$-categories}.
\end{definition}

\begin{explanation}[Comparison with $\BE$-Pseudoalgebras]\label{expl:GBEvsBE_pseudo}
We emphasize that $\GBE$-pseudoalgebras do \emph{not} admit 2-equivalences like the ones in \cref{thm:BEpseudoalg} for $\BE$-pseudoalgebras.  The main reason is that \cref{thm:BEpseudoalg} relies crucially on some coherence properties of the Barratt-Eccles operad $\BE$ that are not shared by $\GBE$ when $G$ is not the trivial group.  For example, consider $\alpha$ in \cref{naive_associativity}, $\gaA_n(\sigma\sscs \ldots)$ in \cref{gammaAn_sigma}, and $\phiA$ in \cref{phiA_coherence_general}.  In those constructions, we use the fact that the objects in $\BE$ are operadically generated by the three objects 
\[\id_n \in \ESigma_n \forspace n \in \{0,1,2\}.\]
There is no such coherence result for $\GBE$ when $G$ is nontrivial.  Thus, a $\GBE$-pseudoalgebra is not generally generated by a finite list of generating operations.  On the other hand, \cref{naive_genuine_smgcat} states that $\BE$-pseudoalgebras yield $\GBE$-pseudoalgebras via $\Catg(\EG,-)$.

There is, however, a general strictification left 2-adjoint that goes from $\Op$-pseudoalgebras to $\Op$-algebras; see \cite[Theorem 0.2]{gmmo20}.  That theorem applies to both $\BE$ and $\GBE$.  Furthermore, as we mention in \cref{rk:naive_genuine}, $\Catg(\EG,-)$ is an equivalence of quasi-categories from $\BE$-algebras to $\GBE$-algebras by \cite[Theorem A]{lenz-genuine}.  Since we do not use those results in this work, we refer the reader to those sources for details.
\end{explanation}

The following recipe for generating genuine symmetric monoidal $G$-categories is the special case of \cref{catgego} \eqref{ng_ii} for the Barratt-Eccles $\Gcat$-operad $\BE$ (\cref{def:BE-Gcat}).

\begin{proposition}\label{naive_genuine_smgcat}
For each group $G$, each $\BE$-pseudoalgebra $\A$ yields a $\GBE$-pseudoalgebra $\Catg(\EG,\A)$.
\end{proposition}

\begin{explanation}[$\GBE$-Pseudoalgebras]\label{expl:GBE_pseudoalg}
We unpack \cref{def:GBE_pseudoalg}.  Suppose $(\A,\gaA,\phiA)$ is a $\GBE$-pseudoalgebra in $\Gcat$.  Then $\A$ is a small $G$-category.  

\parhead{$\GBE$-action}.  Consider the $\GBE$-action $G$-functors \cref{gaAn}
\[\GBE(n) \times \A^n = \Catg(\EG,\ESigma_n) \times \A^n \fto{\gaA_n} \A.\]
For $n=0$, $\gaA_0$ is specified by a $G$-fixed object 
\[\zero = \gaA_0(*) \in \A,\] 
called the basepoint.

\parhead{Unity}. 
Since $\GBE(1) \iso \boldone$, the action unity axiom \cref{pseudoalg_action_unity} implies that $\gaA_1$ is specified by the identity functor
\[\gaA_1(\opu \sscs -) = 1_{\A} \cn \A \to \A.\]
Here 
\[\EG \fto{\opu} \ESigma_1 \iso \boldone\] 
is the unique functor to the terminal category; it is the operadic unit of $\GBE$.

\parhead{Equivariance}. 
For $n \geq 2$, the $G$-equivariance of the $G$-functor $\gaA_n$ means
\[\gaA_n\big( g\cdot F \sscs \ang{g\cdot a_j}_{j\in \ufs{n}}\big) 
= g \cdot \gaA_n\big(F \sscs \ang{a_j}_{j\in \ufs{n}} \big)\]
for $g \in G$, $F \in \GBE(n)$, and $\ang{a_j}_{j\in \ufs{n}} \in \A^n$, and similarly for morphisms in $\GBE(n)$ and $\A^n$.  Here $g \cdot F$ is the functor in \cref{gFk}.  The action equivariance axiom \cref{pseudoalg_action_sym} means
\[\gaA_n\big(F^\sigma \sscs \ang{a_j}_{j\in \ufs{n}}\big) 
= \gaA_n\big(F \sscs \ang{a_{\sigmainv(j)}}_{j\in \ufs{n}}\big)\]
for each permutation $\sigma \in \Sigma_n$, and similarly for morphisms in $\GBE(n) \times \A^n$.  Here $F^\sigma$ is the functor in \cref{Fsigmak}.

\parhead{Associativity constraints}.  A typical component of the $G$-natural isomorphism $\phiA$ \cref{phiA} is an isomorphism
\begin{equation}\label{phiA_GBE}
\begin{tikzpicture}[baseline={(a.base)}]
\draw[0cell=.85]
(0,0) node (a) {\gaA_n\big(F \sscs \bang{\gaA_{m_j} (H_j \sscs \ang{a_{j,i}}_{i \in \ufs{m}_j} ) }_{j\in \ufs{n}}\big)}
(a)++(5.5,0) node (b) {\gaA_m\big( \ga^G \big( F \sscs \ang{H_j}_{j\in \ufs{n}} \big) \sscs \ang{\ang{a_{j,i}}_{i \in \ufs{m}_j} }_{j\in \ufs{n}} \big)}
;
\draw[1cell=.8]
(a) edge node {\phiA} node[swap] {\iso} (b)
;
\end{tikzpicture}
\end{equation}
in $\A$ for functors 
\[\EG \fto{F} \ESigma_n \andspace \EG \fto{H_j} \ESigma_{m_j}\] 
and objects $a_{j,i} \in \A$.  In the codomain, $\ga^G$ denotes the operadic composition in $\GBE$, as we discuss in \cref{GBE_gamma,functorsFHj,gaGFHjg}.  The $G$-equivariance of $\phiA$ means that, for each $g \in G$, the component in \cref{phiA_GBE} for $g \cdot F$, $g \cdot H_j$, and $g \cdot a_{j,i}$ is equal to $g \cdot \phiA$.

The remaining axioms of a $\GBE$-pseudoalgebra,  \cref{pseudoalg_basept_axiom,pseudoalg_topeq,pseudoalg_boteq,pseudoalg_unity,pseudoalg_comp_axiom}, are interpreted using $\gaA$ and $\phiA$ above.
\end{explanation}

\begin{explanation}[Lax $\GBE$-Morphisms]\label{expl:GBE_laxmorphism}
Suppose 
\[(\A,\gaA,\phiA) \fto{(f,\actf)} (\B,\gaB,\phiB)\]
is a lax $\GBE$-morphism between $\GBE$-pseudoalgebras \pcref{def:laxmorphism}.  Then $f \cn \A \to \B$ is a $G$-functor between small $G$-categories. 

\parhead{Action constraints}.  For each $n \geq 0$, a typical component of the $G$-natural transformation $\actf_n$ \cref{laxmorphism_constraint} is a morphism in $\B$
\begin{equation}\label{actfn_GBE}
\gaB_n\big(F \sscs \ang{f(a_j)}_{j\in \ufs{n}} \big) \fto{\actf_n} 
f\big(\gaA_n (F \sscs \ang{a_j}_{j\in \ufs{n}} ) \big)
\end{equation}
for a functor $F \cn \EG \to \ESigma_n$ and objects $\ang{a_j}_{j\in \ufs{n}} \in \A^n$.  Its $G$-equivariance means that, for each $g \in G$, $\actf_n$ in \cref{actfn_GBE} for $g \cdot F$ and $g \cdot a_j$ is equal to $g \cdot \actf_n$.  The lax $\GBE$-morphism $(f,\actf)$ is a $\GBE$-pseudomorphism, respectively strict $\GBE$-morphism, if and only if each component $\actf_n$ \cref{actfn_GBE} is an isomorphism, respectively identity morphism.

\parhead{Axioms}.  By the basepoint axiom \cref{laxmorphism_basepoint} and the unity axiom \cref{laxmorphism_unity}, $\actf_0$ and $\actf_1$ are identity morphisms.  The equivariance axiom \cref{laxmorphism_equiv} and the associativity axiom \cref{laxmorphism_associativity} are interpreted using $\actf_n$ in \cref{actfn_GBE}.
\end{explanation}

\begin{explanation}[$\GBE$-Transformations]\label{expl:GBE_transformation}
Suppose 
\[(f,\actf), (h,\acth) \cn \big(\A,\gaA,\phiA\big) \to \big(\B,\gaB,\phiB\big)\]
are two lax $\GBE$-morphisms between $\GBE$-pseudoalgebras, and suppose 
\[(f,\actf) \fto{\omega} (h,\acth)\]
is a $\GBE$-transformation \pcref{def:algtwocells}.  Then $\omega \cn f \to h$ is a $G$-natural transformation between the underlying $G$-functors.  The compatibility axiom \cref{Otransformation_ax} requires the diagram in $\B$
\[\begin{tikzpicture}[xscale=3.6,yscale=1.5,vcenter]
\draw[0cell=.9]
(0,0) node (x11) {\gaB_n\big(F \sscs \ang{f(a_j)}_{j\in \ufs{n}} \big)}
(x11)++(1,0) node (x12) {f\gaA_n \big(F \sscs \anga \big)}
(x11)++(0,-1) node (x21) {\gaB_n\big(F \sscs \ang{h(a_j)}_{j\in \ufs{n}} \big)}
(x21)++(1,0) node (x22) {h\gaA_n \big(F \sscs \anga \big)}
;
\draw[1cell=.9]  
(x11) edge node {\actf_n} (x12)
(x21) edge node {\acth_n} (x22) 
(x11) edge[transform canvas={xshift={2em}}] node[swap] {\gaB_n\big(1_F \sscs \ang{\omega_{a_j}}_{j\in \ufs{n}} \big)} (x21)
(x12) edge[transform canvas={xshift={-1em}}, shorten >=1pt] node {\omega_{\gaA_n(F;\, \anga)}} (x22)
;
\end{tikzpicture}\]
to commute for each functor $F \cn \EG \to \ESigma_n$ and objects $\anga = \ang{a_j}_{j\in \ufs{n}} \in \A^n$.
\end{explanation}

\chapter{$\Gcat$-Multicategories of Operadic Pseudoalgebras}
\label{ch:multpso}
Throughout this chapter,  $(\Op,\ga,\opu,\pcom)$ denotes a pseudo-commutative operad in $\Gcat$ (\cref{def:pseudocom_operad}) for an arbitrary group $G$.  The main result of this chapter, \cref{thm:multpso}, proves in detail the existence of a $\Gcat$-multicategory $\MultpsO$ that has
\begin{itemize}
\item $\Op$-pseudoalgebras (\cref{def:pseudoalgebra}) as objects, 
\item $k$-lax $\Op$-morphisms (\cref{def:k_laxmorphism}) as $k$-ary 1-cells, and
\item $k$-ary $\Op$-transformations (\cref{def:kary_transformation}) as $k$-ary 2-cells.
\end{itemize}
We emphasize that $\MultpsO$ is enriched in $\Gcat$ (\cref{def:GCat}), not just $\Cat$.  Thus, $\MultpsO$ has multimorphism $G$-categories, and its multicategorical structure morphisms are $G$-functors.  For a $\Tinf$-operad $\Op$ \pcref{as:OpA}, the $\Gcat$-multicategory $\MultpsO$ is the domain of our $G$-equivariant algebraic $K$-theory $\Gtop$-multifunctor $\Kgo$, as we discuss in \cref{sec:Kgo_multi}.

\subsection*{Pseudo and Strict Variants}
There are variant $\Gcat$-multicategories, denoted $\MultpspsO$ and $\MultstO$, with the same objects and $k$-ary 2-cells as $\MultpsO$.  Their $k$-ary 1-cells are, respectively, $k$-ary $\Op$-pseudomorphisms and $k$-ary strict $\Op$-morphisms \pcref{def:k_laxmorphism}.  

\subsection*{$\Cat$-Enriched Variants}
For each variant $\va \in \{\sflax,\sfps,\sfst\}$, passing to $G$-fixed subcategories in the $\Gcat$-enrichment yields a $\Cat$-multicategory $\MultvO^G$ with
\begin{itemize}
\item the same objects as $\MultvO$, which are $\Op$-pseudoalgebras, 
\item $G$-equivariant $k$-ary 1-cells, and
\item $G$-equivariant $k$-ary 2-cells.
\end{itemize}  
See \cref{thm:multpso_gfixed}.  The $\Cat$-multicategory $\MultpsO^G$ extends---in objects, $k$-ary 1-cells, and $k$-ary 2-cells---an analogous multicategory in \cite{gmmo23}.  See \cref{expl:k_laxmorphism}.

\summary 
The following table summarizes the $\Gcat$-multicategories $\MultvO$ for $\va \in \{\sflax,\sfps,\sfst\}$.
\begin{center}
\resizebox{.9\width}{!}{%
{\renewcommand{\arraystretch}{1.3}%
{\setlength{\tabcolsep}{1em}
\begin{tabular}{c|cccc}
& $\MultpsO$ & $\MultpspsO$ & $\MultstO$ & \eqref{def:multicatO} \\ \hline
objects & \multicolumn{3}{c}{$\Op$-pseudoalgebras} & \eqref{def:pseudoalgebra} \\
$k$-ary 1-cells & lax & pseudo & strict & \eqref{def:k_laxmorphism} \\
& \multicolumn{3}{c}{$G$-action on $k$-ary 1-cells} & \eqref{def:k_laxmorphism_g} \\
$k$-ary 2-cells & \multicolumn{3}{c}{$k$-ary $\Op$-transformations} & \eqref{def:kary_transformation} \\
& \multicolumn{3}{c}{$G$-action on $k$-ary 2-cells} & \eqref{expl:O_tr_g} \\
& \multicolumn{3}{c}{symmetric group action} & \eqref{multpso_sym_functor} \\
& \multicolumn{3}{c}{multicategorical composition} & \eqref{def:gam_functor} \\
\end{tabular}}}}
\end{center}

\subsection*{Application to Symmetric Monoidal $G$-Categories}
For the Barratt-Eccles operad $\BE$ (\cref{def:BE-Gcat}), $\MultpsBE$ is a $\Gcat$-multicategory with $\BE$-pseudoalgebras as objects.  Using the 2-equivalences in \cref{thm:BEpseudoalg}, one may regard $\BE$-pseudoalgebras as naive symmetric monoidal $G$-categories.  On the other hand, for the $G$-Barratt-Eccles operad $\GBE$ (\cref{def:GBE}), $\MultpsGBE$ is a $\Gcat$-multicategory whose objects are $\GBE$-pseudoalgebras, which are \emph{genuine} symmetric monoidal $G$-categories.  These applications of \cref{thm:multpso} are discussed in \cref{ex:multBE,ex:multGBE}.

\connection
Restricting to 1-ary 1-cells and 1-ary 2-cells, the $\Cat$-multicategories $\MultvO^G$ yield the 2-categories $\AlgpsvO$ in \cref{oalgps_twocat}; see \cref{ex:underlying_iicat}.  The importance of the $\Gcat$-multicategory $\MultpsO$ is that it serves as the domain of our $G$-equivariant algebraic $K$-theory multifunctor \pcref{sec:Kgo_multi}
\[\MultpsO \fto{\Kgo} \GSp,\]
assuming $\Op$ is a $\Tinf$-operad \pcref{as:OpA}.  For example, the Barratt-Eccles operad $\BE$ and the $G$-Barratt-Eccles operad $\GBE$ are both $\Tinf$-operads.  The first step of $\Kgo$ is a $\Gcat$-multifunctor
\[\MultpsO \fto{\Jgo} \GGCat\] 
from the $\Gcat$-multicategory of $\Op$-pseudoalgebras to the $\Gcat$-multicategory of $\Gskg$-categories.  See \cref{thm:Jgo_multifunctor}.  The $\Gcat$-multifunctor $\Jgo$ takes into account the associativity constraints of $\Op$-pseudoalgebras, so it requires no strictification.  For example, the $\Gcat$-multifunctors $\Jgbe$ and $\Jggbe$ apply directly to, respectively, $\BE$-pseudoalgebras and $\GBE$-pseudoalgebras.  See \cref{ex:JgBE}.

\organization
This chapter consists of the following sections.

\secname{sec:higher_laxmorphisms}  The objects of $\MultpsO$ are $\Op$-pseudoalgebras, which we discuss in \cref{sec:pseudoalgebra}.  The purpose of this section is to define the $k$-ary multimorphism $G$-categories of $\MultpsO$, with $k$-lax $\Op$-morphisms as objects and $k$-ary $\Op$-transformations as morphisms.  The commutativity axiom \cref{laxf_com} of a $k$-lax $\Op$-morphism uses the pseudo-commutative structure of $\Op$.  We also define the $k$-ary multimorphism $G$-categories for the pseudo variant $\MultpspsO$ and the strict variant $\MultstO$, which restrict the objects to, respectively, $k$-ary $\Op$-pseudomorphisms and $k$-ary strict $\Op$-morphisms.

\secname{sec:multpso_sym}  This section constructs the symmetric group action $G$-functors on $\MultvO$ for each variant $\va \in \{\sflax,\sfps,\sfst\}$.  \cref{klax_sigma_welldef,multpso_sym_functor} check that the symmetric group action is well defined.

\secname{sec:multpso_comp}  This section constructs the multicategorical composition $G$-functors on $\MultvO$.  \cref{gam_basept_unity_eq,gam_associativity,gam_commutativity,gam_Otr_welldef} check in detail that the multicategorical composition is well defined.  \cref{multo_ga_g_obj,multo_ga_g_mor} check that composition is $G$-equivariant.  The most intense part of this chapter is \cref{gam_commutativity}, which checks the commutativity axiom for each composite.

\secname{sec:multpsodef}  This section proves the main result of this chapter, \cref{thm:multpso}, on the existence of the $\Gcat$-multicategory $\MultvO$ for each variant $\va \in \{\sflax,\sfps,\sfst\}$.  The symmetric group action and the composition are already shown to be well defined in previous sections.  Thus, the proof of the main theorem consists of checking the axioms of a $\Gcat$-multicategory.  The second half of this section discusses the $\Cat$-multicategory $\MultvO^G$ (\cref{thm:multpso_gfixed}), its underlying 2-category $\AlgpsvO$ (\cref{ex:underlying_iicat}), and applications of these results to the Barratt-Eccles operad $\BE$ (\cref{ex:multBE}) and the $G$-Barratt-Eccles operad $\GBE$ (\cref{ex:multGBE}).

\section{Higher Arity Lax Morphisms and Transformations}
\label{sec:higher_laxmorphisms}

This section begins the construction of the $\Gcat$-multicategory $\MultpsO$ by defining its $k$-ary multimorphism $G$-categories.  Recall from \cref{ufsn} the unpointed finite set 
\[\ufs{n} = \{1 < 2 < \cdots < n\}\] 
with its natural ordering.  We often denote $\txsum_{j=1}^n$ by $\txsum_{j \in \ufs{n}}$, and likewise for other operators that involve a running index.

\secoutline

\begin{itemize}
\item We begin with a discussion of \emph{strength}, denoted $\shuf_{n,i}$, which combines diagonal morphisms and permutations.  Strength is used in the definition of the $k$-ary 1-cells in $\MultpsO$.
\item For $\Op$-pseudoalgebras $\ang{\A_i}_{i \in \ufs{k}}$ and $\B$, the objects of the $k$-ary multimorphism $G$-category 
\[\MultpsO\scmap{\ang{\A_i}_{i \in \ufs{k}} ; \B}\] 
are $k$-lax $\Op$-morphisms (\cref{def:k_laxmorphism,def:k_laxmorphism_g}), which generalize lax $\Op$-morphisms (\cref{expl:onelax_morphism}).  The pseudo-commutative structure of $\Op$ is used in the commutativity axiom \cref{laxf_com}, which only happens if $k > 1$.  This axiom explains why a pseudo-commutative structure on $\Op$ is needed.
\item \cref{laxf_unity_properties} proves two unity properties of $k$-lax $\Op$-morphisms that are used in later discussion.
\item The morphisms in the $k$-ary multimorphism $G$-categories are $k$-ary $\Op$-transformations (\cref{def:kary_transformation}), which generalize $\Op$-transformations (\cref{expl:ktransform_basept}).
\item \cref{def:MultpsO_karycat} defines the $k$-ary multimorphism $G$-categories and the full sub-$G$-categories for $k$-ary strict $\Op$-morphisms and $k$-ary $\Op$-pseudomorphisms.
\item \cref{expl:k_laxmorphism} compares our context with the one in \cite{gmmo23}.
\end{itemize}

\subsection*{Strength}
\cref{expl:leftbracketing} of left normalized bracketing is in effect.  To improve readability, we generally omit parentheses for iterated Cartesian products and tacitly use the associativity isomorphisms wherever they are needed.  The \emph{strength} functor in the following definition is needed to define $k$-lax $\Op$-morphisms; see \cref{klax_constraint}.  Recall that a $\Gcat$-operad $\Op$ is \emph{reduced}\index{reduced}\index{operad!reduced} if $\Op(0) \iso \boldone$, a terminal $G$-category.

\begin{definition}\label{def:Apq}
Suppose $\A_1,\ldots,\A_k$ are $G$-categories (\cref{def:GCat}) for some $k \geq 1$, and $\Op$ is a reduced operad in $\Gcat$.  For $1 \leq p \leq q \leq k$ and $n \geq 0$, we define the following product $G$-categories.
\begin{equation}\label{Apq}
\begin{split}
\A_{[p,q]} &= \A_p \times \A_{p+1} \times \cdots \times \A_q\\
\big(\A_{[p,q]}\big)^n &= \left(\A_p \times \A_{p+1} \times \cdots \times \A_q\right)^n\\
(\A^n)_{[p,q]} &= \A_p^n \times \A_{p+1}^n \times \cdots \times \A_q^n
\end{split}
\end{equation}
By convention, an empty product and $\A_{[u,v]}$ with $u>v$ are the terminal $G$-category $\boldone$.
Using \cref{Apq}, for each $n \geq 0$ and $i \in \ufs{k} = \{1,\ldots,k\}$, we define the \index{strength}\emph{$(n,i)$-strength}, denoted $\shuf_{n,i}$, to be the following composite $G$-functor.
\begin{equation}\label{shufni}
\begin{tikzpicture}[xscale=1,yscale=1,vcenter]
\def\v{-1.5}
\draw[0cell=.9]
(0,0) node (x1) {\A_{[1,i-1]} \times \Op(n) \times \A_i^n \times \A_{[i+1,k]}}
(x1)++(5,0) node (x) {\Op(n) \times \big(\A_{[1,k]}\big)^n}
(x1)++(0,\v) node (x2) {(\A^n)_{[1,i-1]} \times \Op(n) \times \A_i^n \times (\A^n)_{[i+1,k]}}
(x)++(0,\v) node (x3) {\Op(n) \times (\A^n)_{[1,k]}}
;
\draw[1cell=.9]  
(x1) edge node {\shuf_{n,i}} (x)
(x1) edge[transform canvas={xshift=3em}] node[swap] {(\Delta^n)_{[1,i-1]} \times 1 \times (\Delta^n)_{[i+1,k]}} (x2)
(x2) edge node {\twist} node[swap] {\iso} (x3)
(x3) edge[transform canvas={xshift=-1em}] node {\iso} node[swap] {1 \times \twist_{k,n}} (x)
;
\end{tikzpicture}
\end{equation}
In \cref{shufni}, each arrow
\[\Delta^n \cn \A_p \to \A_p^n \forspace p \in \ufs{k} \setminus \{i\}\]
is the $n$-fold diagonal.  The arrow $\twist$ shuffles $\Op(n)$ to the left of $(\A^n)_{[1,i-1]}$.  The arrow $\twist_{k,n}$ permutes the $kn$ factors of $\A$'s as the $(k,n)$-transpose permutation in \cref{eq:transpose_perm}.
\end{definition}

\begin{explanation}[Strength]\label{expl:shufni}
For an explicit description of the $(n,i)$-strength, we consider a tuple 
\begin{equation}\label{strength_dom}
\begin{split}
\bz &= \big(\ang{a_p}_{p=1}^{i-1} \scs x, \ang{a_{i,j}}_{j \in \ufs{n}} \scs \ang{a_p}_{p=i+1}^k\big) \\
&\in \A_{[1,i-1]} \times \Op(n) \times \A_i^n \times \A_{[i+1,k]}
\end{split}
\end{equation}
of all objects or all morphisms.  Then the $(n,i)$-strength $\shuf_{n,i}$ in \cref{shufni} is given by
\begin{equation}\label{shufni_objects}
\shuf_{n,i}\bz = \big(x, \ang{a^j}_{j \in \ufs{n}}\big)
\end{equation}
where
\begin{equation}\label{aj}
a^j = \big(\ang{a_p}_{p=1}^{i-1} \scs a_{i,j}, \ang{a_p}_{p=i+1}^{k}\big) \in \A_{[1,k]}
\end{equation}
for each $j \in \ufs{n}$.  In the case $k=i=1$, the $(n,1)$-strength
\[\Op(n) \times \A^n \fto{\shuf_{n,1}} \Op(n) \times \A^n\]
is the identity $G$-functor.  In the context of 2-endofunctors, the concept of strength is originally introduced in \cite[Section 3.1]{hyland-power}, with further elaboration in \cite[Def.\! 4.1]{corner-gurski}.  It is adapted to the multicategorical context in \cite[3.13]{gmmo23}.
\end{explanation}

\subsection*{Higher Arity Lax Morphisms}
Recall from \cref{def:pseudocom_operad,def:pseudoalgebra} pseudo-commutative operads in $\Gcat$ and $\Op$-pseudoalgebras.  \cref{def:k_laxmorphism} below defines $k$-lax $\Op$-morphisms.  They are the $k$-ary 1-cells of the $\Gcat$-multicategory $\MultpsO$, which we define in \cref{sec:multpsodef}.  \cref{expl:k_laxmorphism} discusses the relationship between the following definition and \cite[Def.\! 3.14]{gmmo23}.  Further discussion of the axioms are given in \cref{expl:onelax_morphism,expl:klax_axioms}.  The $G$-action on a $k$-lax $\Op$-morphism is defined in \cref{def:k_laxmorphism_g}.

\begin{definition}\label{def:k_laxmorphism}
Suppose $(\Op,\ga,\opu,\pcom)$ is a pseudo-commutative operad in $\Gcat$, and suppose 
\[(\A_i,\gaAi,\phiAi) \andspace (\B,\gaB,\phiB)\]
are $\Op$-pseudoalgebras for $1 \leq i \leq k$ for some $k \geq 1$.  A \emph{$k$-lax $\Op$-morphism}\index{k-lax O-morphism@$k$-lax $\Op$-morphism}\index{lax morphism!k-ary@$k$-ary}
\begin{equation}\label{klax_Omorphism}
\bang{(\A_i,\gaAi,\phiAi)}_{i \in \ufs{k}} \fto{(f, \laxf)} (\B,\gaB,\phiB)
\end{equation}
consists of the following data.
\begin{description}
\item[Underlying functor] It is equipped with an \emph{underlying functor}\dindex{underlying}{functor}
\begin{equation}\label{underlying_functor_f}
\txprod_{i \in \ufs{k}}\, \A_i \fto{f} \B.
\end{equation}
\item[Action constraints] For each $n \geq 0$ and $i \in \ufs{k}$, it is equipped with a natural transformation $\laxf_{n,i}$, called the \index{action constraint}\emph{$(n,i)$-action constraint}, as follows, where $\shuf_{n,i}$ is the $(n,i)$-strength in \cref{shufni}.
\begin{equation}\label{klax_constraint}
\begin{tikzpicture}[xscale=1,yscale=1,vcenter]
\def\a{.5} \def\h{4.5} \def\v{1.3}
\draw[0cell=.9]
(0,0) node (x1) {\A_{[1,i-1]} \times \Op(n) \times \A_i^n \times \A_{[i+1,k]}}
(x1)++(\a,\v) node (x2) {\Op(n) \times (\A_{[1,k]})^n}
(x2)++(\h,0) node (x3) {\Op(n) \times \B^n}
(x3)++(0,-2*\v) node (x4) {\B}
(x1)++(0,-\v) node (x) {\A_{[1,i-1]} \times \A_i \times \A_{[i+1,k]}}
;
\draw[1cell=.9]  
(x1) edge node[pos=.2] {\shuf_{n,i}} (x2)
(x2) edge node {1 \times f^n} (x3)
(x3) edge node {\gaB_n} (x4)
(x1) edge node[swap,pos=.5] {1 \times \gaAi_n \times 1} (x)
(x) edge node {f} (x4)
;
\draw[2cell]
node[between=x2 and x4 at .5, shift={(.3,0)}, rotate=-90, 2label={above,\laxf_{n,i}}] {\Rightarrow}
;
\end{tikzpicture}
\end{equation}
\end{description}
With the notation in \cref{strength_dom,shufni_objects,aj}, the $\bz$-component of $\laxf_{n,i}$ is a morphism in $\B$ as follows.
\begin{equation}\label{shufni_component}
\begin{tikzpicture}[xscale=1,yscale=1,baseline={(a.base)}]
\draw[0cell=.9]
(0,0) node (a) {\gaB_n\big(x \sscs \ang{fa^j}_{j \in \ufs{n}}\big)}
(a)++(5.5,0) node (b) {f\Big(\ang{a_p}_{p=1}^{i-1}, \gaAi_n\big( x \sscs \ang{a_{i,j}}_{j \in \ufs{n}} \big), \ang{a_p}_{p=i+1}^{k}\Big)}
;
\draw[1cell=.9]
(a) edge node {\laxf_{n,i;\, \bz}} (b)
;
\end{tikzpicture}
\end{equation}
In \cref{shufni_component}, we call each $a_p \in \A_p$ for $p \in \ufs{k} \setminus \{i\}$ an \index{inactive variable}\emph{inactive variable}.  To simplify the presentation below, we sometimes denote 
\begin{equation}\label{ftilde}
\fgr = f\big( \ang{a_p}_{p=1}^{i-1}, - , \ang{a_p}_{p=i+1}^{k}\big) \cn \A_i \to \B,
\end{equation}
so the $\bz$-component of $\laxf_{n,i}$ becomes
\begin{equation}\label{shufni_simplified}
\begin{tikzpicture}[xscale=1,yscale=1,baseline={(a.base)}]
\draw[0cell=1]
(0,0) node (a) {\gaB_n\big(x \sscs \bang{\fgr a_{i,j}}_{j \in \ufs{n}}\big)}
(a)++(4.5,0) node (b) {\fgr \gaAi_n\big( x \sscs \ang{a_{i,j}}_{j \in \ufs{n}} \big).}
;
\draw[1cell=.9]
(a) edge node {\laxf_{n,i;\, \bz}} (b)
;
\end{tikzpicture}
\end{equation}
Unless otherwise specified, we use the notation in \cref{shufni_component,ftilde,shufni_simplified} and  abbreviate $\laxf_{n,i; \bz}$ to $\laxf_{n,i}$.  The data $(f,\laxf)$ are subject to the axioms \cref{laxf_basepoint,laxf_unity,laxf_eq,laxf_associativity,laxf_com} below whenever they are defined.
\begin{description}
\item[Basepoint]  
$\laxf_{0,i}$ is the identity morphism:
\begin{equation}\label{laxf_basepoint}
\begin{tikzpicture}[xscale=1,yscale=1,baseline={(a.base)}]
\draw[0cell=1]
(0,0) node (a) {\zero^\B}
(a)++(4.5,0) node (b) {f \big(\ang{a_p}_{p=1}^{i-1}, \zero^{\A_i}, \ang{a_p}_{p=i+1}^{k} \big).}
;
\draw[1cell=.9]
(a) edge node {\laxf_{0,i} \,=\, 1} (b)
;
\end{tikzpicture}
\end{equation}
Here $\zero^\B = \gaB_0(*)$ and $\zero^{\A_i} = \gaAi_0(*)$ are the basepoints \cref{pseudoalg_zero} of $\B$ and $\A_i$, respectively.

\item[Unity] 
For objects $\anga = \ang{a_p}_{p \in \ufs{k}} \in \A_{[1,k]}$ and the operadic unit $\opu \in \Op(1)$, $\laxf_{1,i}$ is the identity morphism:
\begin{equation}\label{laxf_unity}
\begin{tikzpicture}[xscale=1,yscale=1,baseline={(a.base)}]
\def\v{-1.4}
\draw[0cell=.9]
(0,0) node (a1) {f\anga = \gaB_1\big( \opu \sscs f\anga\big)}
(a1)++(5.2,0) node (a2) {f\big(\ang{a_p}_{p=1}^{i-1}, \gaAi_1(\opu \sscs a_i), \ang{a_p}_{p=i+1}^{k} \big).}
;
\draw[1cell=.9]
(a1) edge node {\laxf_{1,i} \,=\, 1} (a2)
;
\end{tikzpicture}
\end{equation}
This axiom makes sense by the action unity axiom \cref{pseudoalg_action_unity} for $\B$ and $\A_i$.

\item[Equivariance]
For each permutation $\sigma \in \Sigma_n$, the following diagram in $\B$ commutes, with $a^j$ as defined in \cref{aj}.
\begin{equation}\label{laxf_eq} 
\begin{tikzpicture}[xscale=1,yscale=1,vcenter]
\def\v{-1.4}
\draw[0cell=.85]
(0,0) node (a1) {\gaB_n\big(x \sscs \ang{f a^{\sigmainv(j)}}_{j \in \ufs{n}} \big)}
(a1)++(0,\v) node (a2) {\gaB_n \big( x\sigma \sscs \ang{fa^j}_{j \in \ufs{n}} \big)}
(a1)++(5.6,0) node (b1) {f\big(\ang{a_p}_{p=1}^{i-1} \scs \gaAi_n\big(x \sscs \ang{a_{i,\sigmainv(j)}}_{j \in \ufs{n}} \big), \ang{a_p}_{p=i+1}^k \big)}
(b1)++(0,\v) node (b2) {f\big(\ang{a_p}_{p=1}^{i-1} \scs \gaAi_n\big(x\sigma \sscs \ang{a_{i,j}}_{j \in \ufs{n}} \big), \ang{a_p}_{p=i+1}^k \big)}
;
\draw[1cell=.9]
(a1) edge node {\laxf_{n,i}} (b1)
(a2) edge node {\laxf_{n,i}} (b2)
(a1) edge[equal] (a2)
(b1) edge[equal] (b2)
;
\end{tikzpicture}
\end{equation}
In \cref{laxf_eq}, the two equalities hold by the action equivariance axiom \cref{pseudoalg_action_sym} for $\B$ and $\A_i$.

\item[Associativity] 
Given objects
\[x \in \Op(n) \cq x_j \in \Op(m_j) \cq a_p \in \A_p, \andspace a_{i,j.h} \in \A_i\]
for $p \in \ufs{k} \setminus \{i\}$, $j \in \ufs{n}$, and $h \in \ufs{m_j}$, we use the following notation, with $\fgr \cn \A_i \to \B$ as defined in \cref{ftilde}.
\[\left\{
\begin{gathered}
m = \txsum_{j \in \ufs{n}}\, m_j \qquad m_{\crdot} = \ang{m_j}_{j \in \ufs{n}} \qquad x_\crdot = \ang{x_j}_{j \in \ufs{n}}\\
\begin{aligned}
a_{i,j,\crdot} &= \bang{a_{i,j,h}}_{h \in \ufs{m_j}} \in \A_i^{m_j} & a_{i,\crdot,\crdot} &= \bang{a_{i,j,\crdot}}_{j \in \ufs{n}} \in \A_i^m\\
\fgr a_{i,j,\crdot} &= \bang{\fgr a_{i,j,h}}_{h \in \ufs{m_j}} \in \B^{m_j} \phantom{M} & \fgr a_{i,\crdot,\crdot} &= \bang{\fgr a_{i,j,\crdot}}_{j \in \ufs{n}} \in \B^m\\
\end{aligned}
\end{gathered}
\right.\]
Then the following diagram in $\B$ commutes.
\begin{equation}\label{laxf_associativity}
\begin{tikzpicture}[xscale=1,yscale=1,vcenter]
\def\h{3.5} \def\g{0} \def\v{1} \def\u{-.8}
\draw[0cell=.8]
(0,0) node (x11) {\gaB_n\big(x \sscs \big\langle \gaB_{m_j}( x_j \sscs \fgr a_{i,j,\crdot} ) \big\rangle_{j \in \ufs{n}} \big)}
(x11)++(\h,\v) node (x12) {\gaB_m\big( \ga(x\sscs x_\crdot) \sscs \fgr a_{i,\crdot,\crdot} \big)}
(x11)++(\h+.5,\u) node (x) {\fgr \gaAi_m\big( \ga(x \sscs x_\crdot) \sscs a_{i,\crdot,\crdot} \big)}
(x11)++(0,2*\u) node (x21) {\gaB_n\big(x \sscs \big\langle \fgr \gaAi_{m_j}(x_j \sscs a_{i,j,\crdot}) \big\rangle_{j \in \ufs{n}} \big)}
(x21)++(\h,-\v) node (x22) {\fgr \gaAi_n\big( x \sscs \bang{\gaAi_{m_j}(x_j \sscs a_{i,j,\crdot}) }_{j \in \ufs{n}} \big)}
;
\draw[1cell=.8]  
(x11) edge[transform canvas={xshift={-.7em}}, shorten >=-3pt] node[pos=.5] {\phiB_{(n;\, m_\crdot)}} (x12)
(x12) edge node[pos=.6] {\laxf_{m,i}} (x)
(x11) edge[transform canvas={xshift={1.4em}}] node[swap] {\gaB_n\big( 1_x \sscs \ang{\laxf_{m_j, i}}_{j \in \ufs{n}} \big)} (x21)
(x21) edge[transform canvas={xshift={-.7em}}, shorten >=0pt] node[swap,pos=.5] {\laxf_{n,i}} (x22)
(x22) edge node[swap,pos=.6] {\fgr \phiAi_{(n;\, m_\crdot)}} (x)
;
\end{tikzpicture}
\end{equation}

\item[Commutativity] 
For $1 \leq i < j \leq k$ and objects
\[(x,y) \in \Op(m) \times \Op(n) \cq a_p \in \A_p \cq a_{i,q} \in \A_i, \andspace a_{j,r} \in \A_j\]
for $p \in \ufs{k} \setminus \{i,j\}$, $q \in \ufs{m}$, and $r \in \ufs{n}$, we use the following notation.
\[\left\{
\begin{gathered}
\fhat = f\big(\ang{a_\ell}_{\ell=1}^{i-1} , - , \ang{a_\ell}_{\ell=i+1}^{j-1} , - , \ang{a_\ell}_{\ell=j+1}^{k} \big) \cn \A_i \times \A_j \to \B\\
a_{i,\crdot} = \ang{a_{i,q}}_{q \in \ufs{m}} \in \A_i^m \qquad
a_{j,\crdot} = \ang{a_{j,r}}_{r \in \ufs{n}} \in \A_j^n\\
\begin{aligned}
\fhat\big(a_{i,\crdot} \scs a_{j,r} \big) &= \bang{\fhat(a_{i,q} \scs a_{j,r})}_{q \in \ufs{m}} \in \B^m\\ 
\fhat\big(a_{i,q} \scs a_{j,\crdot} \big) &= \bang{\fhat(a_{i,q} \scs a_{j,r})}_{r \in \ufs{n}} \in \B^n
\end{aligned}
\end{gathered}
\right.\]
Then the following diagram in $\B$ commutes, with $\intr$ the intrinsic pairing of $\Op$ in \cref{intr_jk}, $\twist_{m,n} \in \Sigma_{mn}$ the $(m,n)$-transpose permutation in \cref{eq:transpose_perm}, and $\pcom_{m,n}$ the $(m,n)$-pseudo-commutativity isomorphism in \cref{pseudocom_isos}.
\begin{equation}\label{laxf_com}
\begin{tikzpicture}[xscale=1,yscale=1,vcenter]
\def\h{2.5} \def\v{1.5} \def\u{1.3} \def\a{1.5em}
\draw[0cell=.8]
(0,0) node (x1) {\gaB_m\big( x \sscs \bang{ \gaB_n\big(y \sscs \fhat(a_{i,q} \scs a_{j,\crdot} ) \big)}_{q \in \ufs{m}} \big)}
(x1)++(\h,-\u) node (x2) {\gaB_m\big(x \sscs \bang{\fhat\big(a_{i,q} \scs \gaAj_n(y \sscs a_{j,\crdot} ) \big) }_{q \in \ufs{m}} \big)}
(x2)++(\h,\u) node (x3) {\fhat\big(\gaAi_m (x \sscs a_{i,\crdot} ), \gaAj_n(y \sscs a_{j,\crdot} ) \big)}
(x1)++(0,\v) node (y1) {\gaB_{mn}\big(x \intr y \sscs \bang{\fhat(a_{i,q} \scs a_{j,\crdot} )}_{q \in \ufs{m}} \big)}
(y1)++(0,\v) node (y2) {\gaB_{mn}\big((y \intr x)\twist_{m,n} \sscs \bang{\fhat(a_{i,q} \scs a_{j,\crdot} )}_{q \in \ufs{m}}  \big)}
(y2)++(\h,\u) node (y3) {\gaB_{nm}\big(y \intr x \sscs \bang{\fhat(a_{i,\crdot} \scs a_{j,r} ) }_{r \in \ufs{n}} \big)}
(y3)++(\h,-\u) node (y4) {\gaB_n\big(y \sscs \bang{\gaB_m\big(x \sscs \fhat(a_{i,\crdot} \scs a_{j,r} ) \big) }_{r \in \ufs{n}} \big)}
(y4)++(0,-\v) node (y5) {\gaB_n\big(y \sscs \bang{\fhat\big(\gaAi_m(x \sscs a_{i,\crdot} ), a_{j,r} \big) }_{r \in \ufs{n}} \big)}
;
\draw[1cell=.8]
(x1) edge[transform canvas={xshift=-\a}] node[swap,pos=.2] {\gaB_m(1_x \sscs \ang{\laxf_{n,j}}_{q \in \ufs{m}} )} (x2)
(x2) edge[transform canvas={xshift=\a},shorten >=1ex] node[swap,pos=.8] {\laxf_{m,i}} (x3)
(x1) edge node {\phiB} node[swap] {\iso} (y1)
(y1) edge node {\gaB_{mn}(\pcom_{m,n}\twist_{m,n} \sscs 1^{mn} )} node[swap] {\iso} (y2)
(y2) edge[transform canvas={xshift=-\a},equal,shorten <=1ex] node[pos=.4] {\mathbf{eq}} (y3)
(y3) edge[transform canvas={xshift=\a},shorten >=1ex] node[pos=.6] {(\phiB)^{-1}} node[swap,pos=.2] {\iso} (y4)
(y4) edge node {\gaB_n(1_y \sscs \ang{\laxf_{m,i}}_{r \in \ufs{n}} )} (y5)
(y5) edge node {\laxf_{n,j}} (x3)
;
\end{tikzpicture}
\end{equation}
In \cref{laxf_com}, the equality labeled $\mathbf{eq}$ holds by the action equivariance axiom \cref{pseudoalg_action_sym} of $\B$.
\end{description}
This finishes the definition of a $k$-lax $\Op$-morphism.  Moreover, we call $(f,\laxf)$
\begin{itemize}
\item \emph{$G$-equivariant} if $f$ and $\laxf_{n,i}$ are all $G$-equivariant;
\item a \emph{$k$-ary $\Op$-pseudomorphism}\index{k-ary O-pseudomorphism@$k$-ary $\Op$-pseudomorphism}\index{pseudomorphism!k-ary@$k$-ary} if each $\laxf_{n,i}$ is invertible;
\item a \emph{$k$-ary strict $\Op$-morphism} if each $\laxf_{n,i}$ is the identity; and
\item an \emph{identity $\Op$-morphism} if $f$ is the identity functor with $k=1$ and each $\laxf_{n,1}$ the identity.\defmark
\end{itemize} 
\end{definition}

\begin{definition}[$G$-Action]\label{def:k_laxmorphism_g}
In the same context as \cref{def:k_laxmorphism}, suppose
\[\bang{(\A_i,\gaAi,\phiAi)}_{i \in \ufs{k}} \fto{(f, \laxf)} (\B,\gaB,\phiB)\]
is a $k$-lax $\Op$-morphism.  For each element $g \in G$, the $k$-lax $\Op$-morphism
\begin{equation}\label{klax_Omorphism_g}
\bang{(\A_i,\gaAi,\phiAi)}_{i \in \ufs{k}} \fto{g \cdot (f, \laxf) = (gf\ginv, \laxgf)} (\B,\gaB,\phiB)
\end{equation}
is defined by the conjugation $g$-action as follows.  

\parhead{Underlying functor}. 
The underlying functor of $g \cdot (f,\laxf)$, denoted $gf\ginv$, is defined as the composite
\begin{equation}\label{klax_mor_g}
\txprod_{i \in \ufs{k}}\, \A_i \fto{(\ginv)^k} \txprod_{i \in \ufs{k}}\, \A_i \fto{f} \B \fto{g} \B,
\end{equation}
where each $g$ denotes the $g$-action functor \cref{gactioniso} on $\A_i$ or $\B$.  In other words, $gf\ginv$ is the conjugation $g$-action on $f$.

\parhead{Action constraints}.
For each $n \geq 0$ and $i \in \ufs{k}$, to define the $(n,i)$-action constraint $\laxgf_{n,i}$ of $g \cdot (f,\laxf)$, we consider an arbitrary object $\bz$, as defined in \cref{strength_dom}, and its image under the diagonal $\ginv$-action: 
\begin{equation}\label{strength_dom_ginv}
\begin{split}
\ginv\bz &= \big(\ang{\ginv a_p}_{p=1}^{i-1} \scs \ginv x, \ang{\ginv a_{i,j}}_{j \in \ufs{n}} \scs \ang{\ginv a_p}_{p=i+1}^k\big) \\
&\in \A_{[1,i-1]} \times \Op(n) \times \A_i^n \times \A_{[i+1,k]}.
\end{split}
\end{equation}
The natural transformation
\begin{equation}\label{klax_constraint_g}
\begin{tikzpicture}[xscale=1,yscale=1,vcenter]
\def\a{.5} \def\h{5} \def\v{1.3}
\draw[0cell=.9]
(0,0) node (x1) {\A_{[1,i-1]} \times \Op(n) \times \A_i^n \times \A_{[i+1,k]}}
(x1)++(\a,\v) node (x2) {\Op(n) \times (\A_{[1,k]})^n}
(x2)++(\h,0) node (x3) {\Op(n) \times \B^n}
(x3)++(0,-2*\v) node (x4) {\B}
(x1)++(0,-\v) node (x) {\A_{[1,i-1]} \times \A_i \times \A_{[i+1,k]}}
;
\draw[1cell=.9]  
(x1) edge node[pos=.2] {\shuf_{n,i}} (x2)
(x2) edge node {1 \times (gf\ginv)^n} (x3)
(x3) edge node {\gaB_n} (x4)
(x1) edge node[swap,pos=.5] {1 \times \gaAi_n \times 1} (x)
(x) edge node {gf\ginv} (x4)
;
\draw[2cell]
node[between=x2 and x4 at .5, shift={(0,0)}, rotate=-90, 2label={above,\laxgf_{n,i}}] {\Rightarrow}
;
\end{tikzpicture}
\end{equation}
has $\bz$-component morphism defined by the following commutative diagram in $\B$.
\begin{equation}\label{laxgf_z}
\begin{tikzpicture}[vcenter]
\def\b{-2.5em} \def\g{4.3} \def\h{3} \def\v{-1} 
\draw[0cell=.8]
(0,0) node (a11) {\gaB_n\big(x; \bang{gf\big(\ang{\ginv a_p}_{p=1}^{i-1}, \ginv a_{i,j}, \ang{\ginv a_p}_{p=i+1}^k \big)}_{j \in \ufs{n}} \big)}
(a11)++(\h,\v) node (a12) {gf\big(\ang{\ginv a_p}_{p=1}^{i-1}, \ginv \gaAi_n(x; \ang{a_{i,j}}_{j \in \ufs{n}}), \ang{\ginv a_p}_{p=i+1}^k \big)}
(a11)++(\g,\v) node (a12') {\phantom{\ang{}_{p=i+1}^k}}
(a11)++(0,2*\v) node (a21) {g\gaB_n\big(\ginv x; \bang{f\big(\ang{\ginv a_p}_{p=1}^{i-1}, \ginv a_{i,j}, \ang{\ginv a_p}_{p=i+1}^k \big)}_{j \in \ufs{n}} \big)}
(a21)++(0,1.5*\v) node (a22) {\phantom{\ang{}_{p=i+1}^k}}
(a22)++(\g,0) node (a22') {\phantom{\ang{}_{p=i+1}^k}}
(a22)++(1.7,0) node (a22'') {gf\big(\ang{\ginv a_p}_{p=1}^{i-1}, \gaAi_n(\ginv x; \ang{\ginv a_{i,j}}_{j \in \ufs{n}}), \ang{\ginv a_p}_{p=i+1}^k \big)}
;
\draw[1cell=.8]
(a11) [rounded corners=2pt] -| node[pos=.3] {\laxgf_{n,i;\, \bz}} (a12')
;
\draw[1cell=.8]
(a12') edge[equal] (a22')
(a11) edge[equal, transform canvas={xshift=\b}] (a21)
(a21) edge[transform canvas={xshift=\b}] node[swap] {g\laxf_{n,i;\, \ginv\bz}} (a22)
;
\end{tikzpicture}
\end{equation}
The boundary of the diagram \cref{laxgf_z} is defined as follows.
\begin{itemize}
\item $\laxf_{n,i; \ginv\bz}$ is the $\ginv\bz$-component morphism of $\laxf_{n,i}$, and $g\laxf_{n,i; \ginv\bz}$ is its image under the $g$-action.
\item The left and right vertical equalities hold by the $G$-equivariance of, respectively, $\gaB_n$ and $\gaAi_n$ \cref{gaAn}.
\end{itemize}
In other words, $\laxgf_{n,i}$ is obtained from the conjugation $g$-action on $\laxf_{n,i}$.  
\begin{itemize}
\item The naturality of $\laxgf_{n,i}$ with respect to $\bz$ follows from
\begin{itemize}
\item the naturality of $\laxf_{n,i}$ and
\item the functoriality of the $g$-action on $\B$ and the $\ginv$-action on $\A_i$.
\end{itemize}
\item Each of the axioms \cref{laxf_basepoint,laxf_unity,laxf_eq,laxf_associativity,laxf_com} for $g \cdot (f,\laxf)$ holds by the same axiom for $(f,\laxf)$, the functoriality of the $g$-action, and the following facts about $\Op$, $\ang{\A_i}_{i \in \ufs{k}}$, and $\B$.
\begin{itemize}
\item The basepoint $* \in \Op(0) = \boldone$ is $G$-fixed.  This is needed for the basepoint axiom \cref{laxf_basepoint}.
\item The operadic unit, symmetric group action, and composition on the $\Gcat$-operad $\Op$ are $G$-functors.  
\begin{enumerate}
\item In particular, the unit $\opu \in \Op(1)$ is $G$-fixed.  This is used in the unity axiom \cref{laxf_unity}.  
\item For each permutation $\si \in \Si_n$, the right $\si$-action commutes with the $g$-action.  This is used in the equivariance axiom \cref{laxf_eq}.
\end{enumerate}
\item In each $\Op$-pseudoalgebra $\A$, such as $\B$ and $\A_i$, the $\Op$-action $\gaA_n$ \cref{gaAn} and the associativity constraints $\phiA$ \cref{phiA} are $G$-equivariant.  This is used in the associativity axiom \cref{laxf_associativity} and the commutativity axiom \cref{laxf_com}.
\item For the commutativity axiom \cref{laxf_com}, we also need the $G$-equivariance of (i) the intrinsic pairing $\intr_{m,n}$ \cref{intr_jk} and (ii) the pseudo-commutativity isomorphism $\pcom_{m,n}$ \cref{pseudocom_isos}.
\end{itemize}
\end{itemize}
This finishes the definition of the $k$-lax $\Op$-morphism $g \cdot (f,\laxf)$.  Moreover, we say that $(f,\laxf)$ is \emph{$G$-fixed}\index{G-fixed@$G$-fixed!k-lax O-morphism@$k$-lax $\Op$-morphism} if it is fixed by the $g$-action for each $g \in G$.
\end{definition}

\begin{explanation}[$G$-Action]\label{expl:k_laxmor_g}
The $g$-action \cref{klax_Omorphism_g}
\[(f,\laxf) \mapsto g \cdot (f, \laxf) = (gf\ginv, \laxgf)\]
preserves $k$-ary $\Op$-pseudomorphisms, $k$-ary strict $\Op$-morphisms, and identity $\Op$-morphisms.  Moreover, $(f,\laxf)$ is $G$-equivariant---which means that $f$ and $\laxf_{n,i}$ are all $G$-equivariant---if and only if $(f,\laxf)$ is $G$-fixed.
\end{explanation}

\begin{explanation}[1-Lax $\Op$-Morphisms]\label{expl:onelax_morphism}
A lax $\Op$-morphism, as defined in \cref{def:laxmorphism}, is precisely a $G$-fixed 1-lax $\Op$-morphism in \cref{def:k_laxmorphism_g} for the following reasons.
\begin{itemize}
\item In a lax $\Op$-morphism $(f,\actf)$, $f$ and $\actf$ are $G$-equivariant.  As we noted in \cref{expl:k_laxmor_g}, a $k$-lax $\Op$-morphism $(f,\laxf)$ is $G$-fixed if and only if $f$ and $\laxf_{n,i}$ are $G$-equivariant.  
\item If $k=i=1$ in \cref{def:k_laxmorphism}, then the strength $\shuf_{n,1}$ in \cref{shufni_objects} is the identity, and the $(n,1)$-action constraint $\laxf_{n,1}$ in \cref{klax_constraint} becomes the action constraint $\actf_n$ in \cref{laxmorphism_constraint}.
\item The 1-lax $\Op$-morphism axioms \cref{laxf_basepoint,laxf_unity,laxf_eq,laxf_associativity} become the lax $\Op$-morphism axioms \cref{laxmorphism_basepoint,laxmorphism_unity,laxmorphism_equiv,laxmorphism_associativity}.
\item The commutativity axiom \cref{laxf_com} is empty if $k=1$.
\end{itemize}
Moreover, a $G$-fixed 1-ary $\Op$-pseudomorphism is precisely an $\Op$-pseudomorphism in \cref{def:laxmorphism}.  A $G$-fixed 1-ary strict $\Op$-morphism is precisely a strict $\Op$-morphism.
\end{explanation}

\begin{explanation}[Axioms]\label{expl:klax_axioms}
The basepoint axiom \cref{laxf_basepoint} and the naturality of $\laxf_{0,i}$ imply the following equalities for objects $a_p$ and morphisms $h_p$ in $\A_p$ for $p \in \ufs{k} \setminus \{i\}$.
\begin{equation}\label{laxf_basept}
\begin{split}
f\big(\ang{a_p}_{p=1}^{i-1}, \zero^{\A_i}, \ang{a_p}_{p=i+1}^{k}\big) &= \zero^\B\\
f\big(\ang{h_p}_{p=1}^{i-1}, 1_{\zero^{\A_i}}, \ang{h_p}_{p=i+1}^{k} \big) &= 1_{\zero^\B}
\end{split}
\end{equation}

The commutativity axiom \cref{laxf_com} describes how the action constraints $\laxf_{m,i}$ and $\laxf_{n,j}$ commute in an appropriate sense.  This axiom involves the following instances of the intrinsic pairing \cref{intr_jk}.
\[\begin{split}
x \intr y &= \ga\big(x \sscs \ang{y}_{q \in \ufs{m}}\big) \in \Op(mn)\\
y \intr x &= \ga\big(y \sscs \ang{x}_{r \in \ufs{n}}\big) \in \Op(nm)
\end{split}\]
There is an $(m,n)$-pseudo-commutativity isomorphism \cref{pseudocom_isos}
\begin{equation}\label{pcom_laxf_com}
(x \intr y)\twist_{n,m} \fto[\iso]{\pcom_{m,n}} y \intr x \inspace \Op(nm).
\end{equation}
Applying the $(m,n)$-transpose permutation $\twist_{m,n} = (\twist_{n,m})^{-1}$ to \cref{pcom_laxf_com} yields the isomorphism
\[x \intr y \fto[\iso]{\pcom_{m,n} \twist_{m,n}} (y \intr x)\twist_{m,n} \inspace \Op(mn)\]
that appears in the commutativity axiom \cref{laxf_com}.
\end{explanation}

Before we discuss $\Op$-transformations in higher arity, we record the following unity properties of the action constraint $\laxf$.

\begin{lemma}\label{laxf_unity_properties}
Suppose $(f,\laxf)$ is a $k$-lax $\Op$-morphism as defined in \cref{klax_Omorphism}.  In the context of \cref{shufni_simplified}, a component of the $(n,i)$-action constraint
\[\begin{tikzpicture}[xscale=1,yscale=1,baseline={(a.base)}]
\draw[0cell=1]
(0,0) node (a) {\gaB_n\big(x \sscs \ang{\fgr a_{i,j}}_{j \in \ufs{n}}\big)}
(a)++(4.2,0) node (b) {\fgr \gaAi_n\big( x \sscs \ang{a_{i,j}}_{j \in \ufs{n}} \big)}
;
\draw[1cell=.9]
(a) edge node {\laxf_{n,i}} (b)
;
\end{tikzpicture}\]
is equal to the identity morphism $1_\zero$ in $\B$ if either one of the following two conditions holds:
\begin{enumerate}
\item\label{laxf_u_i} $a_{i,j} = \zero \in \A_i$ for each $j \in \ufs{n}$.
\item\label{laxf_u_ii} $a_\ell = \zero \in \A_\ell$ for some $\ell \in \ufs{k} \setminus \{i\}$.
\end{enumerate}
\end{lemma}

\begin{proof}
For case \eqref{laxf_u_i}, we use the associativity axiom \cref{laxf_associativity} with $m_j = 0$ for each $j \in \ufs{n}$.  This implies $m = \sum_{j\in \ufs{n}}\, m_j = 0$ and 
\[x_j = \ga(x \sscs x_\crdot) = * \in \Op(0).\]
We observe that four of the five morphisms in the associativity diagram \cref{laxf_associativity} are equal to $1_\zero$ in $\B$.
\begin{itemize}
\item The basepoint axiom \cref{laxf_basepoint} for $f$ implies 
\[\laxf_{m,i} = \laxf_{0,i} = 1_\zero.\]
\item The morphism 
\[\gaB_n\big(1_x \sscs \ang{\laxf_{m_j,i}}_{j \in \ufs{n}}\big)
= \gaB_n\big(1_x \sscs \ang{\laxf_{0,i}}_{j \in \ufs{n}}\big)\]
is equal to $1_\zero$ by the basepoint axiom \cref{laxf_basepoint} for $f$, the functoriality of $\gaB_n$, and the basepoint axiom \cref{pseudoalg_basept_axiom} for $\B$.
\item The morphism $\phiB_{(n;\, m_\crdot)} = 1_\zero$ by the basepoint axiom \cref{pseudoalg_basept_axiom} for $\B$.
\item The morphism $\phiAi_{(n;\, m_\crdot)} = 1_\zero$ in $\A_i$ by the basepoint axiom \cref{pseudoalg_basept_axiom} for $\A_i$.  Thus, $\fgr \phiAi_{(n;\, m_\crdot)} = 1_\zero$ in $\B$ by the basepoint condition \cref{laxf_basept} for $f$.
\end{itemize}
Thus, the remaining morphism in the associativity diagram \cref{laxf_associativity}, $\laxf_{n,i}$, is also equal to $1_\zero$.

For case \eqref{laxf_u_ii}, we use the commutativity axiom \cref{laxf_com} with $m=0$ and $x=* \in \Op(0)$.  This implies
\[x \intr y = y \intr x = * \in \Op(0).\]
We observe that six of the seven morphisms in the commutativity diagram \cref{laxf_com} are equal to $1_\zero$ in $\B$, going clockwise from the lower right.
\begin{itemize}
\item $\laxf_{m,i} = \laxf_{0,i} = 1_\zero$ by the basepoint axiom \cref{laxf_basepoint} for $f$.
\item The morphism
\[\gaB_m\big(1_x \sscs \ang{\laxf_{n,j}}_{q \in \ufs{m}}\big) = \gaB_0(1_*)\]
is equal to $1_\zero$ by the functoriality of $\gaB_0$.
\item The morphism labeled $\phiB$ is equal to $1_\zero$ by definition.
\item The morphism 
\[\gaB_{mn}(\pcom_{m,n} \twist_{m,n} \sscs 1^{mn}) = \gaB_0(1_*)\] 
is equal to $1_\zero$ by the functoriality of $\gaB_0$.
\item The morphism labeled $(\phiB)^{-1}$ is equal to $1_\zero$ by the basepoint axiom \cref{pseudoalg_basept_axiom} for $\B$.
\item The morphism 
\[\gaB_n\big(1_y \sscs \ang{\laxf_{m,i}}_{r \in \ufs{n}}\big) 
= \gaB_n\big(1_y \sscs \ang{1_\zero}_{r \in \ufs{n}}\big)\]
is equal to $1_\zero$ by the functoriality of $\gaB_n$ and the basepoint axiom \cref{pseudoalg_basept_axiom} for $\B$.
\end{itemize}
Thus, the remaining morphism in the commutativity diagram \cref{laxf_com}, $\laxf_{n,j}$, is also equal to $1_\zero$.
\end{proof}

\subsection*{Higher Arity Transformations}

Next, we define transformations between  $k$-lax $\Op$-morphisms.

\begin{definition}\label{def:kary_transformation}
Suppose $(\Op,\ga,\opu,\pcom)$ is a pseudo-commutative operad in $\Gcat$ \pcref{def:pseudocom_operad}, and $k \geq 1$.  Suppose we are given two $k$-lax $\Op$-morphisms (\cref{def:k_laxmorphism}) between $\Op$-pseudoalgebras as follows.
\[\big(f, \laxf \big), \big(h, \laxh \big) \cn \bang{(\A_i,\gaAi,\phiAi)}_{i \in \ufs{k}} \to (\B,\gaB,\phiB)\]
A \emph{$k$-ary $\Op$-transformation}\index{transformation!k-ary@$k$-ary}\index{k-ary O-transformation@$k$-ary $\Op$-transformation}
\begin{equation}\label{kary_Otransformation}
\begin{tikzpicture}[xscale=1,yscale=1,baseline={(a.base)}]
\draw[0cell=1]
(0,0) node (a) {\phantom{Z}}
(a)++(2,0) node (b) {\phantom{Z}}
(a)++(-1.3,0) node (a') {\bang{(\A_i,\gaAi,\phiAi)}_{i \in \ufs{k}}}
(b)++(.7,0) node (b') {(\B,\gaB,\phiB)}
;
\draw[1cell=.9]  
(a) edge[bend left=30] node {(f,\laxf)} (b)
(a) edge[bend right=30] node[swap] {(h,\laxh)} (b)
;
\draw[2cell]
node[between=a and b at .45, shift={(0,0)}, rotate=-90, 2label={above,\theta}] {\Rightarrow}
;
\end{tikzpicture}
\end{equation}
is a natural transformation 
\begin{equation}\label{ktr_natural}
\begin{tikzpicture}[baseline={(a.base)}]
\draw[0cell=1]
(0,0) node (a) {\phantom{Z}}
(a)++(2,0) node (b) {\B}
(a)++(-.4,-.05) node (a') {\txprod_{i \in \ufs{k}}\, \A_i}
;
\draw[1cell=.9]  
(a) edge[bend left=30] node {f} (b)
(a) edge[bend right=30] node[swap] {h} (b)
;
\draw[2cell]
node[between=a and b at .45, shift={(0,0)}, rotate=-90, 2label={above,\theta}] {\Rightarrow}
;
\end{tikzpicture}
\end{equation}
such that, using the notation in \cref{shufni_component}, the following diagram in $\B$ commutes.
\begin{equation}\label{ktransform_multilinearity}
\begin{tikzpicture}[xscale=1,yscale=1,vcenter]
\def\v{-1.5}
\draw[0cell=.9]
(0,0) node (a) {\gaB_n\big(x \sscs \ang{fa^j}_{j \in \ufs{n}}\big)}
(a)++(5,0) node (b) {f\big(\ang{a_p}_{p=1}^{i-1}, \gaAi_n\big( x \sscs \ang{a_{i,j}}_{j \in \ufs{n}} \big), \ang{a_p}_{p=i+1}^{k}\big)}
(a)++(0,\v) node (a2) {\gaB_n\big(x \sscs \ang{ha^j}_{j \in \ufs{n}}\big)}
(b)++(0,\v) node (b2) {h\big(\ang{a_p}_{p=1}^{i-1}, \gaAi_n\big( x \sscs \ang{a_{i,j}}_{j \in \ufs{n}} \big), \ang{a_p}_{p=i+1}^{k}\big)}
;
\draw[1cell=.9]
(a) edge node {\laxf_{n,i}} (b)
(a2) edge node {\laxh_{n,i}} (b2)
(a) edge[transform canvas={xshift=2em}] node[swap] {\gaB_n\big(1_x \sscs \ang{\theta}_{j \in \ufs{n}} \big)} (a2)
(b) edge[transform canvas={xshift=-6.5em}] node {\theta} (b2)
;
\end{tikzpicture}
\end{equation}
We call this condition the \emph{multilinearity axiom}.  We say that $\theta$ is \emph{$G$-equivariant} if it is $G$-equivariant as a natural transformation \cref{Gnattr}. 

\parhead{$G$-action}.  We define the $G$-action on a $k$-ary $\Op$-transformation $\theta$ by letting each element $g \in G$ act on $\theta$, denoted $g \cdot \theta$, componentwise by conjugation.  We say that $\theta$ is \emph{$G$-fixed}\index{G-fixed@$G$-fixed!k-ary O-transformation@$k$-ary $\Op$-transformation} if it is fixed by the $g$-action for each $g \in G$.
\end{definition}

\begin{explanation}[$G$-Action]\label{expl:O_tr_g}
For a $k$-ary $\Op$-transformation $\theta$ as defined in \cref{kary_Otransformation}, an element $g \in G$, and objects $a_i \in \A_i$ for $i \in \ufs{k}$, the $k$-ary $\Op$-transformation 
\[g \cdot (f,\laxf) = (gf\ginv, \laxgf) \fto{g \cdot \theta} g \cdot (h,\laxh) = (gh\ginv,\laxgh)\] 
has $\ang{a_i}_{i \in \ufs{k}}$-component morphism given by
\begin{equation}\label{O_tr_g}
gf\ang{\ginv a_i}_{i \in \ufs{k}}
\fto{(g \cdot \theta)_{\ang{a_i}_{i \in \ufs{k}}} = g \theta_{\ang{\ginv a_i}_{i \in \ufs{k}}}}
gh\ang{\ginv a_i}_{i \in \ufs{k}}.
\end{equation}
The naturality of $g \cdot \theta$ follows from
\begin{itemize}
\item the naturality of $\theta$ and 
\item the functoriality of the $g$-action on $\B$.
\end{itemize}
The multilinearity axiom \cref{ktransform_multilinearity} for $g \cdot \theta$ holds by
\begin{itemize}
\item the same axiom for $\theta$,
\item the functoriality of the $g$-action on $\B$, and
\item the $G$-equivariance of $\gaB_n$ and $\gaAi_n$ \cref{gaAn}.
\end{itemize}
Moreover, $\theta$ is $G$-equivariant if and only if it is $G$-fixed, which means $g \cdot \theta = \theta$.
\end{explanation}

\begin{explanation}[Special Cases]\label{expl:ktransform_basept}\
\begin{description}
\item[$\Op$-transformations] Recall from \cref{expl:onelax_morphism} that a $G$-fixed 1-lax $\Op$-morphism is precisely a lax $\Op$-morphism.  Similarly, a $G$-fixed 1-ary $\Op$-transformation is precisely an $\Op$-transformation in the sense of \cref{def:algtwocells}.  Indeed, the multilinearity axiom \cref{ktransform_multilinearity} for $k=1$ is precisely the commutative diagram \cref{Otransformation_ax} that defines an $\Op$-transformation.  
\item[The case $n=0$] The multilinearity axiom \cref{ktransform_multilinearity} for $n=0$ is equivalent to the following equality for each object $\anga = \ang{a_p}_{p \in \ufs{k}} \in \A_{[1,k]}$ with some $a_p = \zero^{\A_p} = \gaAp_0(*)$.
\begin{equation}\label{ktransform_basepoint}
\begin{tikzpicture}[xscale=1,yscale=1,baseline={(a.base)}]
\draw[0cell=1]
(0,0) node (a) {f\anga = \zero^\B}
(a)++(4,0) node (b) {h\anga = \zero^\B}
;
\draw[1cell=1]
(a) edge node {\theta_{\anga} = 1_{\zero^\B}} (b)
;
\end{tikzpicture}
\end{equation}
Indeed, if $n=0$, then both $\laxf_{0,i}$ and $\laxh_{0,i}$ are equal to $1_{\zero^\B}$ by the basepoint axiom \cref{laxf_basepoint}, and $\gaB_0(1_*) = 1_{\zero^\B}$ by functoriality.
\item[The case $n=1$] The multilinearity axiom \cref{ktransform_multilinearity} is automatically true for $n=1$ and the operadic unit $x = \opu \in \Op(1)$.  Indeed, both $\laxf_{1,i}$ and $\laxh_{1,i}$ are identity morphisms by the unity axiom \cref{laxf_unity}.  Moreover, 
\[\gaB_1\big(\opu \sscs -\big) = 1_\B \andspace \gaAi_1\big(\opu \sscs -\big) = 1_{\A_i}\]
by the action unity axiom \cref{pseudoalg_action_unity}.\defmark
\end{description}
\end{explanation}

\subsection*{Multimorphism $G$-Categories}
Next, we define the $k$-ary multimorphism $G$-categories of $\MultpsO$.

\begin{definition}\label{def:MultpsO_karycat}
Suppose $(\Op,\ga,\opu,\pcom)$ is a pseudo-commutative operad in $\Gcat$ (\cref{def:pseudocom_operad}) for an arbitrary group $G$, and suppose
\[(\A_i,\gaAi,\phiAi) \andspace (\B,\gaB,\phiB)\]
are $\Op$-pseudoalgebras (\cref{def:pseudoalgebra}) for $i \in \ufs{k} = \{1,\ldots,k\}$ for some $k \geq 1$.  
\begin{enumerate}
\item\label{multpso_k} We define the $G$-category
\[\begin{split}
& \MultpsO\lrscmap{\ang{\A_i}_{i \in \ufs{k}} ; \B} \\
&= \MultpsO\lrscmap{\bang{(\A_i,\gaAi,\phiAi)}_{i \in \ufs{k}} ; (\B,\gaB,\phiB)}
\end{split}\]
with the following data.
\begin{itemize}
\item Objects are $k$-lax $\Op$-morphisms $\ang{\A_i}_{i \in \ufs{k}} \to \B$ (\cref{def:k_laxmorphism}).
\item Morphisms are $k$-ary $\Op$-transformations (\cref{def:kary_transformation}).
\item Identity morphisms are identity natural transformations.
\item Composition is vertical composition of natural transformations.
\item The conjugation $G$-action is defined in \cref{def:k_laxmorphism_g,def:kary_transformation}.
\end{itemize}
This is a well-defined $G$-category because the multilinearity axiom \cref{ktransform_multilinearity} is closed under vertical composition and is satisfied by identity natural transformations.  The two conditions for a $G$-category, which are stated immediately after \cref{gactioniso}, hold because the $G$-action on $k$-lax $\Op$-morphisms and $k$-ary $\Op$-transformations are both defined by the conjugation $g$-action for $g \in G$, as displayed in \cref{klax_mor_g,laxgf_z,O_tr_g}.
\item\label{multpspso_k} We define the full sub-$G$-categories 
\[\MultstO\lrscmap{\ang{\A_i}_{i \in \ufs{k}} ; \B} 
\bigsubset \MultpspsO\lrscmap{\ang{\A_i}_{i \in \ufs{k}} ; \B}\]
of $\MultpsO\lrscmap{\ang{\A_i}_{i \in \ufs{k}} ; \B}$ with, respectively, $k$-ary strict $\Op$-morphisms and $k$-ary $\Op$-pseudomorphisms as objects.
\item\label{multpso_zero} For the case $k=0$, we define
\[\MultstO\scmap{\ang{};\B} = \MultpspsO\scmap{\ang{};\B} = \MultpsO\scmap{\ang{};\B} = \B,\]
which is a $G$-category.\defmark
\end{enumerate}
\end{definition}

\begin{explanation}[Literature Comparison]\label{expl:k_laxmorphism}
\cref{def:k_laxmorphism,def:MultpsO_karycat} are related to \cite[Def.\! 3.14]{gmmo23} in the following ways.  
\begin{enumerate}
\item In \cite{gmmo23}, it is assumed that each $\A_i$ and $\B$ are $\Op$-algebras, instead of $\Op$-pseudoalgebras.  However, in \cite[3.4]{gmmo23}, it is remarked that $\Op$-pseudoalgebras can also be used in their definition. 
\item In \cite{gmmo23}, it is assumed that each $\laxf_{n,i}$, denoted $\delta_i(n)$ there, is invertible.  Our action constraints $\laxf_{n,i}$ in \cref{klax_constraint} are not required to be invertible.
\item In \cite{gmmo23}, $f$ and $\delta_i(n)$ are assumed to be $G$-equivariant.  In contrast, a $k$-lax $\Op$-morphism is $G$-equivariant if and only if it is fixed by the $G$-action \pcref{expl:k_laxmor_g}.
\end{enumerate}
The following table summarizes our context and the one in \cite{gmmo23}.
\begin{center}
\resizebox{\columnwidth}{!}{%
{\renewcommand{\arraystretch}{1.3}%
{\setlength{\tabcolsep}{1em}
\begin{tabular}{c|ccc|c}
& \multicolumn{3}{c|}{\cref{def:MultpsO_karycat}} & \cite[Def.\! 3.14]{gmmo23}\\ 
& $\MultstO$ & $\MultpspsO$ & $\MultpsO$ & $\mathbf{Mult}(\Op)$\\ \hline
ambient category & \multicolumn{3}{c|}{$\Gcat$ (\ref{def:GCat})} & $\Cat(\cV)$ (2.1)\\
objects & \multicolumn{3}{c|}{$\Op$-pseudoalgebras (\ref{def:pseudoalgebra})} & $\Op$-algebras (2.10)\\
$k$-ary $\Op$-morphisms & strict & pseudo & lax (\ref{def:k_laxmorphism}) & $G$-equiv. pseudo (2.11)\\
$k$-ary 2-cells & \multicolumn{3}{c|}{$k$-ary $\Op$-transformations (\ref{def:kary_transformation})} & not specified\\
multimorphism & \multicolumn{3}{c|}{$G$-categories (\ref{def:MultpsO_karycat})} & sets\\
\end{tabular}}}}
\end{center}
Our proof in the rest of this chapter demonstrates that invertibility of the action constraints $\laxf_{n,i}$ is not necessary to obtain a $\Gcat$-multicategory $\MultpsO$.  

In \cite{gmmo23}, the existence of their multicategory $\mathbf{Mult}(\Op)$ is inferred from the following two results.
\begin{enumerate}
\item For a pseudo-commutative operad $\Op$ in $\Cat$, the associated 2-monad has a symmetric pseudo-commutative structure.  This statement is extracted from \cite[Theorems 4.4 and 4.6]{corner-gurski}; see also \cite[Remark 11.2]{gmmo23}.
\item For a symmetric pseudo-commutative 2-monad $T$, there is a $\Cat$-multicategory consisting of $T$-algebras, multilinear pseudomorphisms, and algebra 2-cells.  This assertion is \cite[Proposition 18]{hyland-power}.
\end{enumerate}
As far as we know, the cited results in \cite{corner-gurski} are never published.  Moreover, the cited result in \cite{hyland-power} has little explanation.  In particular, it does not explain the commutativity axiom, as stated in \cite[page 169]{hyland-power}, for a composite, which we consider to be the hardest part of the proof.  In our multicategorical context, the commutativity axiom \cref{laxf_com} for a composite is verified in \cref{gam_commutativity}.  The commutativity axiom \cref{laxf_com} for a $k$-lax $\Op$-morphism is the precise reason why we need a pseudo-commutative structure on $\Op$, as the pseudo-commutative structure is not used in any other axioms in \cref{def:pseudoalgebra,def:k_laxmorphism,def:kary_transformation}. 
\end{explanation}

\section{Symmetric Group Action}
\label{sec:multpso_sym}

Recall from \cref{def:enr-multicategory} that each enriched multicategory has a symmetric group action.  As the next step in the construction of the $\Gcat$-multicategory $\MultpsO$, this section defines its symmetric group action $G$-functors.

\secoutline

\begin{itemize}
\item \cref{def:multpso_sym} constructs the symmetric group action on $k$-lax $\Op$-morphisms, with further elaboration given in \cref{expl:laxfsi_components}.  \cref{klax_sigma_welldef} checks that this assignment is well defined.
\item \cref{def:ktransformation_sym} constructs the symmetric group action on $k$-ary $\Op$-transformations.
\item \cref{multpso_sym_functor} checks that the symmetric group action is a $G$-functor and satisfies the expected axioms.
\end{itemize}
Throughout this section, we assume $(\Op,\ga,\opu,\pcom)$ is a pseudo-commutative operad in $\Gcat$ (\cref{def:pseudocom_operad}), and 
\[(\A_i,\gaAi,\phiAi) \andspace (\B,\gaB,\phiB)\]
are $\Op$-pseudoalgebras (\cref{def:pseudoalgebra}) for $i \in \ufs{k} = \{1,\ldots,k\}$ for some $k \geq 1$. 

\subsection*{Permuted Lax Morphisms}
We first define the symmetric group action on $k$-lax $\Op$-morphisms (\cref{def:k_laxmorphism}).

\begin{definition}\label{def:multpso_sym}
Suppose
\[\bang{(\A_i,\gaAi,\phiAi)}_{i \in \ufs{k}} \fto{(f, \laxf)} (\B,\gaB,\phiB)\]
is a $k$-lax $\Op$-morphism, and $\sigma \in \Sigma_k$ is a permutation.  We define a $k$-lax $\Op$-morphism
\begin{equation}\label{klax_sigma}
\bang{\A_{\sigma(i)}, \gaAsii, \phiAsii}_{i \in \ufs{k}} \fto{(f\sigma, \laxfsi)} (\B,\gaB,\phiB)
\end{equation}
as follows.
\begin{description}
\item[Underlying functor] Its underlying functor $f\sigma$ is the following composite, in which the first arrow permutes the factors according to $\sigma$.
\begin{equation}\label{klax_sigma_functor}
\begin{tikzpicture}[baseline={(a.base)}]
\draw[0cell=.9]
(0,0) node (a) {\prod_{i=1}^k\, \A_{\sigma(i)}}
(a)++(2.5,0) node (b) {\prod_{i=1}^k\, \A_i}
(b)++(2,0) node (c) {\B}
;
\draw[1cell=.9]
(a) edge node {\sigma} node[swap] {\iso} (b)
(b) edge node {f} (c)
;
\end{tikzpicture}
\end{equation}
\item[Action constraints] For $n \geq 0$ and $i \in \ufs{k}$, we denote by 
\[\A' = \prod_{p=1}^{\sigmainv(i)-1} \A_{\sigma(p)} \andspace 
\A'' = \prod_{p=\sigmainv(i)+1}^k \A_{\sigma(p)}.\]
The $(n,\sigmainv(i))$-action constraint is defined as the whiskering 
\begin{equation}\label{laxfsi_whiskering}
\laxfsi_{n,\sigmainv(i)} = \laxf_{n,i} * \sigma
\end{equation}
in the following pasting diagram, where the two unlabeled regions commute.
\begin{equation}\label{fsigma_constraint}
\begin{tikzpicture}[baseline={(x1.base)}]
\def\a{.5} \def\h{2.7} \def\v{1.4}
\draw[0cell=.75]
(0,0) node (x1) {\A_{[1,i-1]} \times \Op(n) \times \A_i^n \times \A_{[i+1,k]}}
(x1)++(\a,\v) node (x2) {\Op(n) \times (\A_{[1,k]})^n}
(x2)++(\h,0) node (x3) {\Op(n) \times \B^n}
(x3)++(0,-2*\v) node (x4) {\B}
(x1)++(0,-\v) node (x) {\A_{[1,i-1]} \times \A_i \times \A_{[i+1,k]}}
(x1)++(-3.8,0) node (y1) {\A' \times \Op(n) \times \A_i^n \times \A''}
(y1)++(\a,\v) node (y2) {\Op(n) \times \big(\txprod_{p=1}^k \A_{\sigma(p)} \big)^n}
(y1)++(0,-\v) node (y) {\A' \times \A_i \times \A'' = \txprod_{p=1}^k \A_{\sigma(p)}}
;
\draw[1cell=.75]  
(x1) edge node[pos=.2] {\shuf_{n,i}} (x2)
(x2) edge node {1 \times f^n} (x3)
(x3) edge node {\gaB_n} (x4)
(x1) edge node[swap,pos=.5] {1 \times \gaAi_n \times 1} (x)
(x) edge node {f} (x4)
(y1) edge node[pos=.2] {\shuf_{n,\sigmainv(i)}} (y2)
(y2) edge node {1 \times \sigma^n} (x2)
(y1) edge node {\sigma} (x1)
(y1) edge node[swap] {1 \times \gaAi_n \times 1} (y)
(y) edge node {\sigma} (x)
;
\draw[2cell]
node[between=x2 and x4 at .5, shift={(0,.6)}, rotate=-90, 2label={above,\laxf_{n,i}}] {\Rightarrow}
;
\end{tikzpicture}
\end{equation}
The naturality of $\laxfsi_{n,\sigmainv(i)}$ follows from the functoriality of permutation by $\sigma$ and the naturality of $\laxf_{n,i}$.
\end{description}
This finishes the definition of the data $(f\sigma, \laxfsi)$.  \cref{klax_sigma_welldef} proves that it is a $k$-lax $\Op$-morphism.
\end{definition}

\begin{explanation}[Permuted Action Constraints]\label{expl:laxfsi_components}
Consider an object in the middle left entry of \cref{fsigma_constraint}:
\begin{equation}\label{laxfsi_domain_obj}
\begin{split}
z = &~ \big( \ang{a_{\sigma(p)}}_{p=1}^{\sigmainv(i)-1} \scs x, \ang{a_{i,j}}_{j \in \ufs{n}} \scs \ang{a_{\sigma(p)}}_{p=\sigmainv(i)+1}^k \big)\\ 
&\in \A' \times \Op(n) \times \A_i^n \times \A''.
\end{split}
\end{equation}
The middle horizontal permutation $\sigma$ in \cref{fsigma_constraint} regards $\Op(n) \times \A_i^n$ as the $\sigmainv(i)$-th object, so
\begin{equation}\label{laxfsi_si_z}
\si z = \big( \ang{a_p}_{p=1}^{i-1} \scs x, \ang{a_{i,j}}_{j \in \ufs{n}} \scs \ang{a_p}_{p=i+1}^k \big).
\end{equation}
The $z$-component of $\laxfsi_{n,\sigmainv(i)}$ is the $\si z$-component of $\laxf_{n,i}$ \cref{shufni_component}:
\begin{equation}\label{laxfsi_z}
\begin{tikzpicture}[xscale=1,yscale=1,baseline={(a.base)}]
\draw[0cell=.75]
(0,0) node (a) {\gaB_n\Big(x \sscs \bang{f\big(\ang{a_p}_{p=1}^{i-1}, a_{i,j}, \ang{a_p}_{p=i+1}^{k}\big)}_{j \in \ufs{n}}\Big)}
(a)++(5.6,0) node (b) {f\left(\ang{a_p}_{p=1}^{i-1}, \gaAi_n\big( x \sscs \ang{a_{i,j}}_{j \in \ufs{n}} \big), \ang{a_p}_{p=i+1}^{k}\right).}
;
\draw[1cell=.75]
(a) edge node {\laxf_{n,i}} (b)
;
\end{tikzpicture}
\end{equation}
The following computation shows that the upper left region in \cref{fsigma_constraint} commutes, where $z$ may denote either an object or a morphism.
\begin{equation}\label{shuf_sigma_z}
\begin{split}
&\shuf_{n,i}(\si z)\\
&= \shuf_{n,i} \big( \ang{a_p}_{p=1}^{i-1} \scs x, \ang{a_{i,j}}_{j \in \ufs{n}} \scs \ang{a_p}_{p=i+1}^k\big)\\
&= \big( x, \bang{ \ang{a_p}_{p=1}^{i-1} \scs a_{i,j} \scs \ang{a_p}_{p=i+1}^k }_{j \in \ufs{n}} \big)\\
&= \left(1 \times \sigma^n\right) \big(x, \bang{ \ang{a_{\sigma(p)}}_{p=1}^{\sigmainv(i)-1} \scs a_{i,j} \scs \ang{a_{\sigma(p)}}_{p=\sigmainv(i)+1}^k }_{j \in \ufs{n}} \big)\\
&= \left(1 \times \sigma^n\right) \shuf_{n,\sigmainv(i)}(z)
\end{split}
\end{equation}
The lower left rectangle in \cref{fsigma_constraint} commutes by the naturality of the braiding for the Cartesian product.  
\end{explanation}

Now we verify that $(f\sigma, \laxfsi)$ is well defined.

\begin{lemma}\label{klax_sigma_welldef}
The pair $(f\sigma, \laxfsi)$ defined in \cref{klax_sigma} is a $k$-lax $\Op$-morphism.
\end{lemma}

\begin{proof}
We need to verify the axioms \cref{laxf_basepoint,laxf_unity,laxf_eq,laxf_associativity,laxf_com} of a $k$-lax $\Op$-morphism.  By \cref{laxfsi_whiskering,laxfsi_z}, each component of the $(n,\sigmainv(i))$-action constraint $\laxfsi_{n,\sigmainv(i)}$ is the corresponding component of $\laxf_{n,i}$ precomposed with the permutation $\sigma \in \Sigma_k$. 

\parhead{Basepoint and Unity}.  The basepoint axiom \cref{laxf_basepoint} and the unity axiom \cref{laxf_unity} hold for $\laxfsi$ because $\laxf_{0,i}$ and $\laxf_{1,i}$ for the operadic unit $\opu \in \Op(1)$ are identities.

\parhead{Equivariance}.  We consider a permutation $\rho \in \Sigma_n$ and an object
\[\begin{split}
y = &~ \big( \ang{a_{\sigma(p)}}_{p=1}^{\sigmainv(i)-1} \scs x\rho, \ang{a_{i,j}}_{j \in \ufs{n}} \scs \ang{a_{\sigma(p)}}_{p=\sigmainv(i)+1}^k \big)\\ 
&\in \A' \times \Op(n) \times \A_i^n \times \A''
\end{split}\]
with $x \in \Op(n)$.  The object $y$ is obtained from the object $z$ in \cref{laxfsi_domain_obj} by replacing $x$ with $x\rho$.  The $y$-component of $\laxfsi_{n,\sigmainv(i)}$ is the top horizontal arrow in the following diagram in $\B$.
\[\begin{tikzpicture}[xscale=1,yscale=1,vcenter]
\def\v{-1.4}
\draw[0cell=.75]
(0,0) node (a1) {\gaB_n \Big( x\rho \sscs \bang{f\big(\ang{a_p}_{p=1}^{i-1}, a_{i,j}, \ang{a_p}_{p=i+1}^{k} \big)}_{j \in \ufs{n}} \Big)}
(a1)++(0,\v) node (a2) {\gaB_n\Big(x \sscs \bang{f \big(\ang{a_p}_{p=1}^{i-1}, a_{i,\rhoinv(j)}, \ang{a_p}_{p=i+1}^{k} \big)}_{j \in \ufs{n}} \Big)}
(a1)++(6,0) node (b1) {f\Big(\ang{a_p}_{p=1}^{i-1}, \gaAi_n\big(x\rho \sscs \ang{a_{i,j}}_{j \in \ufs{n}} \big) \scs \ang{a_p}_{p=i+1}^{k} \Big)}
(b1)++(0,\v) node (b2) {f\Big(\ang{a_p}_{p=1}^{i-1}, \gaAi_n\big(x \sscs \ang{a_{i,\rhoinv(j)}}_{j \in \ufs{n}} \big) \scs \ang{a_p}_{p=i+1}^{k} \Big)}
;
\draw[1cell=.75]
(a1) edge node {\laxf_{n,i}} (b1)
(a2) edge node {\laxf_{n,i}} (b2)
(a1) edge[equal] (a2)
(b1) edge[equal] (b2)
;
\end{tikzpicture}\]
This diagram commutes by the equivariance axiom \cref{laxf_eq} for $(f,\laxf)$.  The bottom horizontal arrow is the component of $\laxfsi_{n,\sigmainv(i)}$ at the object
\[\big( \ang{a_{\sigma(p)}}_{p=1}^{\sigmainv(i)-1} \scs x, \ang{a_{i,\rhoinv(j)}}_{j \in \ufs{n}} \scs \ang{a_{\sigma(p)}}_{p=\sigmainv(i)+1}^k \big).\]
This proves the equivariance axiom \cref{laxf_eq} for $\big(f\sigma, \laxfsi\big)$.

\parhead{Associativity}.  We consider objects
\[x \in \Op(n) \cq x_j \in \Op(m_j) \cq a_p \in \A_p, \andspace a_{i,j.h} \in \A_i\]
for $p \in \ufs{k} \setminus \{i\}$, $j \in \ufs{n}$, and $h \in \ufs{m_j}$, along with the notation defined just before \cref{laxf_associativity}.  The functor $\fgr$ in \cref{ftilde} is equal to
\[\fsigr = (f\sigma) \Big( \ang{a_{\sigma(\ell)}}_{\ell=1}^{\sigmainv(i)-1} \scs - \scs \ang{a_{\sigma(\ell)}}_{\ell=\sigmainv(i)+1}^{k} \Big) \cn \A_i \to \B,\]
where the argument, denoted $-$, appears in the $\sigmainv(i)$-th slot and is reserved for $\A_i$.  The desired associativity axiom \cref{laxf_associativity} for $\big( f\sigma, \laxfsi\big)$ asserts the commutativity of the following diagram in $\B$.
\[\begin{tikzpicture}[xscale=1,yscale=1,vcenter]
\def\h{3} \def\g{0} \def\v{1.2} \def\u{-.9}
\draw[0cell=.8]
(0,0) node (x11) {\gaB_n\big(x \sscs \ang{\gaB_{m_j} ( x_j \sscs \fsigr a_{i,j,\crdot} ) }_{j \in \ufs{n}} \big)}
(x11)++(\h,\v) node (x12) {\gaB_m\big( \ga(x\sscs x_\crdot) \sscs \fsigr a_{i,\crdot,\crdot} \big)}
(x11)++(\h+.5,\u) node (x) {\fsigr \gaAi_m\big( \ga(x \sscs x_\crdot) \sscs a_{i,\crdot,\crdot} \big)}
(x11)++(0,2*\u) node (x21) {\gaB_n\big(x \sscs \ang{ \fsigr \gaAi_{m_j} (x_j \sscs a_{i,j,\crdot}) }_{j \in \ufs{n}} \big)}
(x21)++(\h,-\v) node (x22) {\fsigr \gaAi_n\big( x \sscs \ang{\gaAi_{m_j}(x_j \sscs a_{i,j,\crdot}) }_{j \in \ufs{n}} \big)}
;
\draw[1cell=.7]  
(x11) edge[transform canvas={xshift={-.5em}}, shorten >=0pt] node[pos=.5] {\phiB_{(n;\, m_\crdot)}} (x12)
(x12) edge node[pos=.6] {\laxfsi_{m,\sigmainv(i)}} (x)
(x11) edge[transform canvas={xshift={1em}}] node[swap] {\gaB_n\big( 1_x \sscs \ang{\laxfsi_{m_j, \sigmainv(i)}}_{j \in \ufs{n}} \big)} (x21)
(x21) edge[transform canvas={xshift={-.5em}}, shorten >=0pt] node[swap,pos=.3,inner sep=-2pt] {\laxfsi_{n,\sigmainv(i)}} (x22)
(x22) edge node[swap,pos=.6] {\fsigr \phiAi_{(n;\, m_\crdot)}} (x)
;
\end{tikzpicture}\]
This diagram is equal to \cref{laxf_associativity}, which commutes by the associativity axiom of $(f,\laxf)$.

\parhead{Commutativity}.  For $1 \leq i < j \leq k$, we consider objects
\[(x,y) \in \Op(m) \times \Op(n) \cq a_p \in \A_p \cq a_{i,q} \in \A_i, \andspace a_{j,r} \in \A_j\]
for $p \in \ufs{k} \setminus \{i,j\}$, $q \in \ufs{m}$, and $r \in \ufs{n}$, along with the notation defined just before \cref{laxf_com}.  The functor $\fhat$ there is equal to
\[\fsihat = (f\sigma) \big( a_{\sigma(1)}, \ldots, -, \ldots, -, \ldots, a_{\sigma(k)} \big) \cn \A_i \times \A_j \to \B,\]
where the two arguments, each denoted $-$, appear in slots $\sigmainv(i)$ and $\sigmainv(j)$ and are reserved for $\A_i$ and $\A_j$.  The commutativity axiom for $\big(f\sigma, \laxfsi\big)$ asserts the commutativity of the following diagram in $\B$.
\[\begin{tikzpicture}[xscale=1,yscale=1,vcenter]
\def\h{2.5} \def\v{1.5} \def\u{1.3} \def\a{1em}
\draw[0cell=.8]
(0,0) node (x1) {\gaB_m\big( x \sscs \bang{ \gaB_n\big(y \sscs \fsihat(a_{i,q} \scs a_{j,\crdot} ) \big)}_{q \in \ufs{m}} \big)}
(x1)++(\h,-\u) node (x2) {\gaB_m\big(x \sscs \bang{\fsihat\big(a_{i,q} \scs \gaAj_n(y \sscs a_{j,\crdot} ) \big) }_{q \in \ufs{m}} \big)}
(x2)++(\h,\u) node (x3) {\fsihat\big(\gaAi_m (x \sscs a_{i,\crdot} ), \gaAj_n(y \sscs a_{j,\crdot} ) \big)}
(x1)++(0,\v) node (y1) {\gaB_{mn}\big(x \intr y \sscs \bang{\fsihat(a_{i,q} \scs a_{j,\crdot} )}_{q \in \ufs{m}} \big)}
(y1)++(0,\v) node (y2) {\gaB_{mn}\big((y \intr x)\twist_{m,n} \sscs \bang{\fsihat(a_{i,q} \scs a_{j,\crdot} )}_{q \in \ufs{m}}  \big)}
(y2)++(\h,\u) node (y3) {\gaB_{nm}\big(y \intr x \sscs \bang{\fsihat(a_{i,\crdot} \scs a_{j,r} ) }_{r \in \ufs{n}} \big)}
(y3)++(\h,-\u) node (y4) {\gaB_n\big(y \sscs \bang{\gaB_m\big(x \sscs \fsihat(a_{i,\crdot} \scs a_{j,r} ) \big) }_{r \in \ufs{n}} \big)}
(y4)++(0,-\v) node (y5) {\gaB_n\big(y \sscs \bang{\fsihat\big(\gaAi_m(x \sscs a_{i,\crdot} ), a_{j,r} \big) }_{r \in \ufs{n}} \big)}
;
\draw[1cell=.75]
(x1) edge[transform canvas={xshift=-\a}] node[swap,pos=.2] {\gaB_m\big(1_x \sscs \bang{\laxfsi_{n,\sigmainv(j)}}_{q \in \ufs{m}} \big)} (x2)
(x2) edge[transform canvas={xshift=\a},shorten >=1ex] node[swap,pos=.8] {\laxfsi_{m,\sigmainv(i)}} (x3)
(x1) edge node {\phiB} node[swap] {\iso} (y1)
(y1) edge node {\gaB_{mn}\big(\pcom_{m,n}\twist_{m,n} \sscs 1^{mn} \big)} node[swap] {\iso} (y2)
(y2) edge[transform canvas={xshift=-\a},equal,shorten <=1ex] node[pos=.4] {\mathbf{eq}} (y3)
(y3) edge[transform canvas={xshift=\a},shorten >=1ex] node[pos=.6] {(\phiB)^{-1}} node[swap,pos=.2] {\iso} (y4)
(y4) edge node {\gaB_n\big(1_y \sscs \bang{\laxfsi_{m,\sigmainv(i)}}_{r \in \ufs{n}} \big)} (y5)
(y5) edge node {\laxfsi_{n,\sigmainv(j)}} (x3)
;
\end{tikzpicture}\]
This diagram is equal to \cref{laxf_com}, which commutes by the commutativity axiom of $(f,\laxf)$.
\end{proof}

\subsection*{Permuted Transformations}

Next, we define the symmetric group action on $k$-ary $\Op$-transformations.

\begin{definition}\label{def:ktransformation_sym}
In the context of \cref{def:kary_transformation}, suppose
\[\begin{tikzpicture}[xscale=1,yscale=1,baseline={(a.base)}]
\draw[0cell=1]
(0,0) node (a) {\phantom{Z}}
(a)++(2,0) node (b) {\phantom{Z}}
(a)++(-1.3,0) node (a') {\bang{(\A_i,\gaAi,\phiAi)}_{i \in \ufs{k}}}
(b)++(.7,0) node (b') {(\B,\gaB,\phiB)}
;
\draw[1cell=.9]  
(a) edge[bend left=30] node {(f,\laxf)} (b)
(a) edge[bend right=30] node[swap] {(h,\laxh)} (b)
;
\draw[2cell]
node[between=a and b at .45, shift={(0,0)}, rotate=-90, 2label={above,\theta}] {\Rightarrow}
;
\end{tikzpicture}\]
is a $k$-ary $\Op$-transformation.  For each permutation $\sigma \in \Sigma_k$, we define a $k$-ary $\Op$-transformation
\[\begin{tikzpicture}[xscale=1,yscale=1,baseline={(a.base)}]
\draw[0cell=1]
(0,0) node (a) {\phantom{Z}}
(a)++(2,0) node (b) {\phantom{Z}}
(a)++(-1.8,0) node (a') {\bang{(\A_{\sigma(i)},\gaAsii,\phiAsii)}_{i \in \ufs{k}}}
(b)++(.7,0) node (b') {(\B,\gaB,\phiB)}
;
\draw[1cell=.9]  
(a) edge[bend left=30] node {(f\sigma,\laxfsi)} (b)
(a) edge[bend right=30] node[swap] {(h\sigma,\laxhsi)} (b)
;
\draw[2cell]
node[between=a and b at .4, shift={(0,0)}, rotate=-90, 2label={above,\theta^\sigma}] {\Rightarrow}
;
\end{tikzpicture}\]
as the natural transformation 
\begin{equation}\label{Otr_sigma}
\theta^\sigma = \theta * \sigma
\end{equation}
given by the following whiskering.
\[\begin{tikzpicture}[xscale=1,yscale=1,baseline={(a.base)}]
\draw[0cell=1]
(0,0) node (a) {\phantom{Z}}
(a)++(2,0) node (b) {\B}
(a)++(-.3,0) node (a') {\prod_{i=1}^k \A_i}
(a')++(-2.5,0) node (c) {\prod_{i=1}^k \A_{\sigma(i)}}
;
\draw[1cell=.9]  
(c) edge node {\sigma} (a')
(a) edge[bend left=30] node {f} (b)
(a) edge[bend right=30] node[swap] {h} (b)
;
\draw[2cell]
node[between=a and b at .45, shift={(0,0)}, rotate=-90, 2label={above,\theta}] {\Rightarrow}
;
\end{tikzpicture}\]
The naturality of $\theta^\sigma$ follows from the functoriality of $\sigma$ and the naturality of $\theta$.  The multilinearity axiom \cref{ktransform_multilinearity} for $\theta^\sigma$ follows from the same axiom for $\theta$, using \cref{klax_sigma_functor,laxfsi_whiskering,Otr_sigma}.
\end{definition}

\subsection*{$G$-Functoriality}
Denoting $\ang{\A} = \ang{\A_i}_{i \in \ufs{k}}$, recall from \cref{def:MultpsO_karycat} the $G$-categories\label{not:multpsvo_ab} 
\[\MultvO\scmap{\ang{\A}; \B}\]
with the subscript $\va \in \{\sflax,\sfps,\sfst\}$ specifying the objects:
\begin{itemize}
\item $k$-lax $\Op$-morphisms $\ang{\A} \to \B$ for $\va = \sflax$,
\item $k$-ary $\Op$-pseudomorphisms $\ang{\A} \to \B$ for $\va = \sfps$, and
\item $k$-ary strict $\Op$-morphisms $\ang{\A} \to \B$ for $\va = \sfst$.
\end{itemize}
The morphisms in $\MultvO\scmap{\ang{\A}; \B}$ are $k$-ary $\Op$-transformations (\cref{def:kary_transformation}).

\begin{lemma}\label{multpso_sym_functor}
For each permutation $\sigma \in \Sigma_k$ and each of the three variants $\va \in \{\sflax,\sfps,\sfst\}$, there is an isomorphism of $G$-categories
\begin{equation}\label{multpso_sigma_functor}
\MultvO\scmap{\ang{\A_i}_{i \in \ufs{k}};\B} \fto[\iso]{\sigma} 
\MultvO\scmap{\ang{\A_{\sigma(i)}}_{i \in \ufs{k}};\B}
\end{equation}
given by the assignments in \cref{def:multpso_sym,def:ktransformation_sym}:
\[\begin{split}
(f,\laxf) & \mapsto (f\sigma, \laxfsi) \andspace \\
\theta & \mapsto \theta * \sigma.
\end{split}\]
Moreover, this is the identity functor if $\sigma = \id_k$, and the symmetric group action axiom \cref{enr-multicategory-symmetry} holds.
\end{lemma}

\begin{proof}
The proofs for the cases $\va \in \{\sfps,\sfst\}$ are obtained from the proof for the case $\va = \sflax$ by restricting to $k$-ary $\Op$-pseudomorphisms and $k$-ary strict $\Op$-morphisms.  Thus, we concentrate on the case $\va = \sflax$ in the rest of this proof.

\parhead{Well defined}.  The object assignment is well defined by \cref{klax_sigma_welldef}, which proves that $(f\sigma, \laxfsi)$ is a $k$-lax $\Op$-morphism whenever $(f,\laxf)$ is so.  Moreover, if the action constraints $\laxf_{n,i}$ are all invertible, respectively identities, then the same is true for $\laxfsi_{n,\sigmainv(i)}$ by \cref{laxfsi_whiskering}.  The morphism assignment, $\theta \mapsto \theta * \sigma$, is well defined, as we explain at the end of \cref{def:ktransformation_sym}. 

\parhead{Functoriality}.  These object and morphism assignments define a functor because whiskering with $\sigma$ preserves identities and vertical composition of natural transformations.

\parhead{Axioms}.  The assertions about $\id_k$ and the symmetric group action axiom hold by \cref{klax_sigma_functor,laxfsi_whiskering,Otr_sigma}.

\parhead{$G$-equivariance}.  To show that the functor $\si$ is $G$-equivariant, recall that the $G$-action on $k$-lax $\Op$-morphisms and $k$-ary $\Op$-transformations are defined in, respectively, \cref{def:k_laxmorphism_g,def:kary_transformation}.  For a $k$-lax $\Op$-morphism $(f,\laxf)$ in the domain of $\si$ and an element $g \in G$, first applying the $g$-action and then $\si$ yields the $k$-lax $\Op$-morphism
\begin{equation}\label{multo_si_g_i}
\big((gf\ginv)\si, \lax^{(g \cdot f)\si} \big).
\end{equation}
On the other hand, first applying $\si$ and then the $g$-action yields
\begin{equation}\label{multo_si_g_ii}
\big(g(f\si)\ginv, \lax^{g \cdot (f\si)} \big).
\end{equation}
We need to show that the $k$-lax $\Op$-morphisms in \cref{multo_si_g_i,multo_si_g_ii} have the same underlying functors and the same action constraints.
\begin{description}
\item[Underlying functors]
By \cref{klax_mor_g,klax_sigma_functor}, the functors $(gf\ginv)\si$ and $g(f\si)\ginv$ are, respectively, the left-bottom composite and the top composite in the following diagram.
\begin{equation}\label{multo_si_gaction}
\begin{tikzpicture}[vcenter]
\def\v{-1.5}
\draw[0cell=.9]
(0,0) node (a11) {\txprod_{i \in \ufs{k}}\, \A_{\si(i)}}
(a11)++(3,0) node (a12) {\txprod_{i \in \ufs{k}}\, \A_{\si(i)}}
(a11)++(0,\v) node (a21) {\txprod_{i \in \ufs{k}}\, \A_i}
(a12)++(0,\v) node (a22) {\txprod_{i \in \ufs{k}}\, \A_i}
(a22)++(2,0) node (a23) {\B}
(a23)++(1.5,0) node (a24) {\B}
;
\draw[1cell=.9]
(a11) edge node {(\ginv)^k} (a12)
(a11) edge node[swap] {\si} (a21)
(a12) edge node {\si} (a22)
(a21) edge node {(\ginv)^k} (a22)
(a22) edge node {f} (a23) 
(a23) edge node {g} (a24)
;
\end{tikzpicture}
\end{equation}
The left square in \cref{multo_si_gaction} commutes by the naturality of the braiding for the Cartesian product.  Thus, the $k$-lax $\Op$-morphisms in \cref{multo_si_g_i,multo_si_g_ii} have the same underlying functor.
\item[Action constraints]
Using the naturality of the braiding for the Cartesian product,  \cref{strength_dom_ginv,laxgf_z,laxfsi_whiskering,laxfsi_si_z}, the following computation shows that the $k$-lax $\Op$-morphisms in \cref{multo_si_g_i,multo_si_g_ii} have the same $(n,i)$-action constraint for $n \geq 0$, $i \in \ufs{k}$, and an object $z$ as defined in \cref{laxfsi_domain_obj}.
\[\begin{split}
\lax^{(g \cdot f)\si}_{n,i;\, z} 
&= \lax^{g \cdot f}_{n, \si(i);\, \si(z)} \\
&= g\laxf_{n, \si(i);\, \ginv \si(z)} \\
&= g\laxf_{n,\si(i);\, \si(\ginv z)} \\
&= g\lax^{f\si}_{n,i;\, \ginv z} \\
&= \lax^{g \cdot (f\si)}_{n,i;\, z}
\end{split}\]
Thus, the $k$-lax $\Op$-morphisms in \cref{multo_si_g_i,multo_si_g_ii} are equal.
\end{description}

The proof that the functor $\si$ commutes with the $g$-action on each $k$-ary $\Op$-transformation $\theta$ reuses the diagram \cref{multo_si_gaction}, with $f$ replaced by $\theta$, along with \cref{O_tr_g,Otr_sigma}.
\end{proof}

\section{Multicategorical Composition}
\label{sec:multpso_comp}

As a brief recapitulation, we have constructed
\begin{itemize}
\item the $k$-ary multimorphism $G$-categories $\MultpsO\scmap{\ang{\A_i}_{i\in \ufs{k}}; \B}$, as well as their pseudo and strict variants, in \cref{def:MultpsO_karycat} and
\item the symmetric group action $G$-functors in \cref{multpso_sym_functor}.
\end{itemize}
Continuing our construction of the $\Gcat$-multicategory $\MultpsO$, this section constructs its multicategorical composition $G$-functor \cref{eq:enr-defn-gamma}.

\secoutline

\begin{itemize}
\item The composition on objects---which are $k$-lax $\Op$-morphisms for $k \geq 0$---is given in  \cref{def:gam_Omorphisms}.  This definition is an adaptation of the composition of monoidal functors; see \cref{expl:laxfh_component}. 
\item \cref{gam_basept_unity_eq,gam_associativity,gam_commutativity} verify the necessary axioms.  The most intense part of this section is \cref{gam_commutativity}, which verifies the commutativity axiom for a composite.
\item \cref{multo_ga_g_obj} proves that composition on objects is $G$-equivariant.
\item The composition on morphisms---which are $k$-ary $\Op$-transformations for $k \geq 0$---is given in \cref{def:gam_Otr}.  \cref{gam_Otr_welldef,multo_ga_g_mor} verify the necessary axioms and $G$-equivariance.
\item The composition $G$-functors of $\MultpsO$, $\MultpspsO$, and $\MultstO$ are given in \cref{def:gam_functor}.
\end{itemize}

\begin{assumption}\label{assum:multpso_comp}
Throughout this section, we assume $(\Op,\ga,\opu,\pcom)$ is a pseudo-commutative operad in $\Gcat$ (\cref{def:pseudocom_operad}), and 
\[(\Aij,\gaAij,\phiAij) \cq (\B_i,\gaBi,\phiBi), \andspace (\C,\gaC,\phiC)\]
are $\Op$-pseudoalgebras (\cref{def:pseudoalgebra}) for $i \in \ufs{k} = \{1,\ldots,k\}$ for some $k \geq 1$ and $j \in \ufs{d}_i = \{1,\ldots,d_i\}$ for some $d_i \geq 0$.  To simplify the presentation, we use the following notation, with $\crdot\,$ denoting a running index, $1 \leq p \leq q \leq k$, and $1 \leq r \leq t \leq d_i$.
\[\left\{
\begin{gathered}
\begin{aligned}
d &= \txsum_{i=1}^k d_i & d_\crdot &= \ang{d_i}_{i \in \ufs{k}} & \Bdot &= \ang{\B_i}_{i \in \ufs{k}}\\
\Aidot &= \ang{\Aij}_{j \in \ufs{d}_i} \phantom{M} & \Addot &= \ang{\Aidot}_{i \in \ufs{k}} \phantom{M} & \B_{[1,k]} &= \txprod_{i=1}^k\, \B_i\\
\A_i &= \txprod_{j=1}^{d_i}\, \Aij & \A_{[p,q]} &= \txprod_{i=p}^q\, \A_i & \A_{i[r,t]} &= \txprod_{j=r}^t\, \A_{i,j}\\
\end{aligned}\\
\begin{split}
\A_{<(i,j)} &= \A_{[1,i-1]} \times \A_{i[1,j-1]} \\
&= \big(\txprod_{u=1}^{i-1}\, \A_u\big) \times \big(\txprod_{v=1}^{j-1}\, \A_{i,v}\big)\\ 
\A_{>(i,j)} &= \A_{i[j+1,d_i]} \times \A_{[i+1,k]} \\
&= \big(\txprod_{v=j+1}^{d_i}\, \A_{i,v}\big) \times \big(\txprod_{u=i+1}^k\, \A_u\big)
\end{split}
\end{gathered}
\right.\]
By convention, $\A_{[p,q]} = \boldone$ if $p>q$, and $\A_{i[r,t]} = \boldone$ if $r > t$.
\end{assumption}

\subsection*{Multicategorical Composition of Lax $\Op$-Morphisms}

The multicategorical composition of $\MultpsO$ is a $G$-functor.  We first define its object assignment.  To avoid confusion between several composition, we use the symbol $\gam$ for the composition of $\MultpsO$.

\begin{definition}[$\gam$ on Objects]\label{def:gam_Omorphisms}
Under \cref{assum:multpso_comp}, suppose we are given a $k$-lax $\Op$-morphism (\cref{def:k_laxmorphism})
\[\B_\crdot \fto{(f,\laxf)} \C\]
and a $d_i$-lax $\Op$-morphism
\[\A_{i\crdot} \fto{(h_i, \laxhi)} \B_i \foreachspace i \in \ufs{k}.\]
If any $d_i = 0$, then $h_i$ is interpreted as an object in $\B_i$ by \cref{def:MultpsO_karycat} \eqref{multpso_zero}.  We define the \emph{composite} $d$-lax $\Op$-morphism
\begin{equation}\label{gam_fh}
\gam\big((f,\laxf) \sscs \ang{(h_i, \laxhi)}_{i \in \ufs{k}} \big) 
= (fh,\laxfh) \cn \A_{\crdots} \to \C
\end{equation}
as follows.
\begin{description}
\item[Underlying functor] Denoting by
\[h = \txprod_{i=1}^k h_i \cn \txprod_{i=1}^k \txprod_{j=1}^{d_i} \Aij = \A_{[1,k]} \to \txprod_{i=1}^k \B_i = \B_{[1,k]},\]
the underlying functor of \cref{gam_fh} is given by the following composite.
\begin{equation}\label{gam_fh_functor}
\begin{tikzpicture}[xscale=1,yscale=1,baseline={(a.base)}]
\draw[0cell=1]
(0,0) node (a) {\A_{[1,k]}}
(a)++(2.5,0) node (b) {\B_{[1,k]}}
(b)++(2,0) node (c) {\C}
;
\draw[1cell=.9]  
(a) edge node {h} (b)
(b) edge node {f} (c)
;
\end{tikzpicture}
\end{equation}
\item[Action constraints] For each $n \geq 0$, $(i,j) \in \ufs{k} \times \ufs{d}_i$, and 
\begin{equation}\label{elldj}
\ell = \big(\txsum_{t=1}^{i-1} d_t\big) + j,
\end{equation}
the $(n,\ell)$-action constraint of \cref{gam_fh} is a natural transformation as follows.
\[\begin{tikzpicture}[xscale=1,yscale=1,vcenter]
\def\h{3.2} \def\g{1} \def\v{1.5} \def\t{0} \def\u{1.3}
\draw[0cell=.8]
(0,0) node (x1) {\A_{<(i,j)} \times \Op(n) \times \Aij^n \times \A_{>(i,j)}}
(x1)++(0,\v) node (x2) {\Op(n) \times \big(\A_{[1,k]}\big)^n}
(x2)++(\h+.3,\t) node (x3) {\Op(n) \times \big(\B_{[1,k]}\big)^n}
(x3)++(\h,-\t) node (x4) {\Op(n) \times \C^n}
(x4)++(0,-\v) node (x5) {\C}
(x1)++(\g,-\u) node (y1) {\A_{<(i,j)} \times \Aij \times \A_{>(i,j)} = \A_{[1,k]}}
(x5)++(-\g,-\u) node (y2) {\B_{[1,k]}}
;
\draw[1cell=.8]  
(x1) edge node {\shuf_{n,\ell}} (x2)
(x2) edge node {1 \times h^n} (x3)
(x3) edge node {1 \times f^n} (x4)
(x4) edge node {\gaC_n} (x5)
(x1) edge node[swap,pos=.1] {1 \times \gaAij_n \times 1} (y1)
(y1) edge node {h} (y2)
(y2) edge node[swap,pos=.6] {f} (x5)
;
\draw[2cell]
node[between=x1 and x5 at .5, rotate=-90, 2label={above,\laxfh_{n,\ell}}] {\Rightarrow}
;
\end{tikzpicture}\]
To define its components, we consider objects
\begin{equation}\label{laxfh_objects}
x \in \Op(n) \cq a_{p,r} \in \A_{p,r}, \andspace a_{i,j,\crdot} = \ang{a_{i,j,m}}_{m \in \ufs{n}} \in \Aij^n
\end{equation}
for $(p,r) \in \ufs{k} \times \ufs{d}_p \setminus \{(i,j)\}$.  We use the following notation to simplify the inactive variables, which are those in $\A_{p,r}$ with $(p,r) \neq (i,j)$.
\begin{equation}\label{htilftil}
\left\{
\begin{aligned}
a_{p,\crdot} &= \ang{a_{p,r}}_{r \in \ufs{d}_p} \in \A_{p,1} \times \cdots \times \A_{p,d_p} = \A_p\\
\hgr_i &= h_i\big(\ang{a_{i,r}}_{r=1}^{j-1} \scs - \scs \ang{a_{i,r}}_{r=j+1}^{d_i} \big) \cn \Aij \to \B_i\\
\fgr &= f\big(\ang{h_p a_{p,\crdot}}_{p=1}^{i-1} \scs - \scs \ang{h_p a_{p,\crdot}}_{p=i+1}^k \big) \cn \B_i \to \C\\
\hgr_i a_{i,j,\crdot} &= \bang{\hgr_i a_{i,j,m}}_{m \in \ufs{n}} \in \B_i^n \\
\fgr \hgr_i a_{i,j,\crdot} &= \bang{\fgr \hgr_i a_{i,j,m}}_{m \in \ufs{n}} \in \C^n
\end{aligned}
\right.
\end{equation}
The component of $\laxfh_{n,\ell}$ at the objects in \cref{laxfh_objects} is defined as the following composite morphism in $\C$.
\begin{equation}\label{laxfh_component}
\begin{tikzpicture}[xscale=1,yscale=1,vcenter]
\def\h{2.5} \def\v{1.2} 
\draw[0cell=.9]
(0,0) node (x1) {\gaC_n\big(x \sscs \fgr \hgr_i a_{i,j,\crdot} \big)}
(x1)++(\h,-\v) node (x2) {\fgr \gaBi_n\big(x \sscs \hgr_i a_{i,j,\crdot} \big)}
(x1)++(2*\h,0) node (x3) {\fgr \hgr_i \gaAij_n\big(x \sscs a_{i,j,\crdot} \big)}
;
\draw[1cell=.9]
(x1) edge node {\laxfh_{n,\ell}} (x3)
(x1) edge[bend right=25] node[swap,pos=.2] {\laxf_{n,i}} (x2)
(x2) edge[bend right=25] node[swap,pos=.8] {\fgr \laxhi_{n,j}} (x3)
;
\end{tikzpicture}
\end{equation}
The naturality of $\laxfh_{n,\ell}$ follows from (i) the naturality of $\laxf_{n,i}$ and $\laxhi_{n,j}$ and (ii) the functoriality of $\fgr$.
\end{description}
\cref{gam_basept_unity_eq,gam_associativity,gam_commutativity} prove that the pair $(fh,\laxfh)$ in \cref{gam_fh} is a $d$-lax $\Op$-morphism in the sense of \cref{def:k_laxmorphism}.
\end{definition}

\begin{explanation}[Composite Action Constraints]\label{expl:laxfh_component}
The component of $\laxfh_{n,\ell}$ in \cref{laxfh_component} is a generalization of the monoidal constraint of a composite of two monoidal functors.  Indeed, for composable monoidal functors (\cref{def:monoidalfunctor})
\[\A \fto{(H,H^2,H^0)} \B \fto{(F,F^2,F^0)} \C,\]
the composite monoidal functor $FH$ has monoidal constraint given componentwise by the following composite for objects $a_1, a_2 \in \A$.
\[FHa_1 \otimes FHa_2 \fto{F^2} F(Ha_1 \otimes Ha_2) \fto{F(H^2)} FH(a_1 \otimes a_2)\]
The argument for \cref{gam_basept_unity_eq,gam_associativity,gam_commutativity} is an elaborate version of the proof that $FH$ is a monoidal functor.
\end{explanation}

\subsection*{Proofs}
\cref{gam_basept_unity_eq,gam_associativity,gam_commutativity,multo_ga_g_obj} verify that the composite $(fh,\laxfh)$ is a $d$-lax $\Op$-morphism (\cref{def:k_laxmorphism}), starting with the simpler axioms, and that $\gam$ is $G$-equivariant.

\begin{lemma}\label{gam_basept_unity_eq}
The pair $(fh,\laxfh)$ in \cref{gam_fh} satisfies the axioms \cref{laxf_basepoint,laxf_unity,laxf_eq}.
\end{lemma}

\begin{proof}
\parhead{Basepoint}.  The basepoint axiom \cref{laxf_basepoint} asks for the equality
\[\begin{tikzpicture}[xscale=1,yscale=1,baseline={(a.base)}]
\draw[0cell=1]
(0,0) node (a) {\zero^\C}
(a)++(3,0) node (b) {\fgr \hgr_i \zero^{\Aij}.}
;
\draw[1cell=.9]
(a) edge node {\laxfh_{0,\ell} \,=\, 1} (b)
;
\end{tikzpicture}\]
By the definition \cref{laxfh_component} of $\laxfh_{0,\ell}$, this axiom holds by (i) the basepoint axiom for $\laxf_{0,i}$ and $\laxhi_{0,j}$ and (ii) the functoriality of $\fgr$.

\parhead{Unity}.  The unity axiom \cref{laxf_unity} asks for the equality
\[\begin{tikzpicture}[xscale=1,yscale=1,baseline={(a.base)}]
\draw[0cell=1]
(0,0) node (a1) {\gaC_1\big( \opu \sscs \fgr \hgr_i a_{i,j,1}\big)}
(a1)++(5,0) node (a2) {\fgr \hgr_i \gaAij_1\big(\opu \sscs a_{i,j,1}\big).}
;
\draw[1cell=.9]
(a1) edge node {\laxfh_{1,\ell} \,=\, 1} (a2)
;
\end{tikzpicture}\]
By the definition \cref{laxfh_component} of $\laxfh_{1,\ell}$, this axiom holds by (i) the unity axiom for $\laxf_{1,i}$ and $\laxhi_{1,j}$ and (ii) the functoriality of $\fgr$.

\parhead{Equivariance}.  By the definition \cref{laxfh_component} of $\laxfh_{n,\ell}$, the equivariance axiom \cref{laxf_eq} asks for the commutativity of the following diagram for each permutation $\sigma \in \Sigma_n$.
\[\begin{tikzpicture}[baseline={(a.base)}]
\def\v{-1.4}
\draw[0cell=.9]
(0,0) node (a1) {\gaC_n\big( x\sigma \sscs \bang{\fgr \hgr_i a_{i,j,m}}_{m \in \ufs{n}}\big)}
(a1)++(4.2,0) node (a2) {\gaC_n\big( x \sscs \bang{\fgr \hgr_i a_{i,j,\sigmainv(m)}}_{m \in \ufs{n}}\big)}
(a1)++(0,\v) node (b1) {\fgr \gaBi_n\big( x\sigma \sscs \bang{\hgr_i a_{i,j,m}}_{m \in \ufs{n}}\big)}
(a2)++(0,\v) node (b2) {\fgr \gaBi_n\big( x \sscs \bang{\hgr_i a_{i,j,\sigmainv(m)}}_{m \in \ufs{n}}\big)}
(b1)++(0,\v) node (c1) {\fgr \hgr_i \gaAij_n\big( x\sigma \sscs \bang{a_{i,j,m}}_{m \in \ufs{n}}\big)} 
(b2)++(0,\v) node (c2) {\fgr \hgr_i \gaAij_n\big( x \sscs \bang{a_{i,j,\sigmainv(m)}}_{m \in \ufs{n}}\big)} 
;
\draw[1cell=.9]
(a1) edge node[swap] {\laxf_{n,i}} (b1)
(b1) edge node[swap] {\fgr \laxhi_{n,j}} (c1)
(a2) edge node {\laxf_{n,i}} (b2)
(b2) edge node {\fgr \laxhi_{n,j}} (c2)
(a1) edge[equal] (a2)
(b1) edge[equal] (b2)
(c1) edge[equal] (c2)
;
\end{tikzpicture}\]
In this diagram, the top region commutes by the equivariance axiom for $\laxf_{n,i}$.  The bottom region commutes by the equivariance axiom for $\laxhi_{n,j}$ and the functoriality of $\fgr$.
\end{proof}

\begin{lemma}\label{gam_associativity}
The pair $(fh,\laxfh)$ in \cref{gam_fh} satisfies the associativity axiom \cref{laxf_associativity}.
\end{lemma}

\begin{proof}
For $1 \leq i \leq k$ and $1 \leq j \leq d_i$, we consider objects
\[\big(x \sscs \ang{x_t}_{t\in \ufs{n}} \big) \in \Op(n) \times \txprod_{t=1}^n \Op(m_t) \cq a_{p,r} \in \A_{p,r}, \andspace \adot_{t,u} \in \Aij\]
for $(p,r) \in \ufs{k} \times \ufs{d}_p \setminus \{(i,j)\}$, $t \in \ufs{n}$, and $u \in \ufs{m}_t$, along with the notation in \cref{htilftil} and the following notation.
\[\left\{
\begin{gathered}
\begin{aligned}
m &= \txsum_{t \in \ufs{n}} m_t & \adot_{t,\crdot} &= \ang{\adot_{t,u}}_{u \in \ufs{m}_t} \in \Aij^{m_t}\\
\bx &= \ga\big(x \sscs \ang{x_t}_{t \in \ufs{n}}\big) \in \Op(m) \phantom{M}
& \adot_{\crdot,\crdot} &= \ang{\adot_{t,\crdot}}_{t \in \ufs{n}} \in \Aij^m\\
\end{aligned}
\end{gathered}
\right.\]
The dot in $\adot_{t,u}$ indicates that it is an object in $\Aij$.  We use this notation instead of the more cumbersome notation $a_{i,j,t,u}$.

By the definition \cref{laxfh_component} of $\laxfh$ and the functoriality of $\gaC_n$, the associativity axiom for $(fh,\laxfh)$ asserts the commutativity of the following boundary diagram in $\C$.
\[\begin{tikzpicture}[baseline={(a.base)}]
\def\a{10} \def\b{15} \def\g{1}\def\h{3} \def\e{3.5} \def\u{-1.2} \def\v{-1.5}
\draw[0cell=.7]
(0,0) node (x1) {\gaC_n\big( x \sscs \bang{\gaC_{m_t} \big(x_t \sscs \fgr \hgr_i \adot_{t,\crdot} \big)}_{t \in \ufs{n}}\big)}
(x1)++(\e,\u) node (x2) {\gaC_m\big( \bx \sscs \fgr \hgr_i \adot_{\crdot,\crdot} \big)}
(x2)++(0,\v) node (x3)  {\fgr \gaBi_m \big( \bx \sscs \hgr_i \adot_{\crdot,\crdot} \big)}
(x3)++(0,\v) node (x4) {\fgr \hgr_i \gaAij_m \big( \bx \sscs \adot_{\crdot,\crdot} \big)}
(x1)++(-\h,\u) node (y1) {\gaC_n\big( x \sscs \bang{\fgr \gaBi_{m_t} \big(x_t \sscs \hgr_i \adot_{t,\crdot} \big)}_{t \in \ufs{n}}\big)}
(y1)++(-\g,\v) node (y2) {\gaC_n\big( x \sscs \bang{\fgr \hgr_i \gaAij_{m_t} \big(x_t \sscs \adot_{t,\crdot} \big)}_{t \in \ufs{n}}\big)}
(y2)++(\g,\v) node (y3) {\fgr \gaBi_n\big( x \sscs \bang{\hgr_i \gaAij_{m_t} \big(x_t \sscs \adot_{t,\crdot} \big)}_{t \in \ufs{n}}\big)}
(y3)++(\h,\u) node (y4) {\fgr \hgr_i \gaAij_n\big( x \sscs \bang{\gaAij_{m_t} \big(x_t \sscs \adot_{t,\crdot} \big)}_{t \in \ufs{n}}\big)}
node[between=y2 and x3 at .5] (z) {\fgr \gaBi_n\big( x \sscs \bang{ \gaBi_{m_t} \big(x_t \sscs \hgr_i \adot_{t,\crdot} \big)}_{t \in \ufs{n}}\big)}
;
\draw[1cell=.75]
(x1) edge[bend left=\a] node {\phiC} (x2)
(x2) edge node {\laxf_{m,i}} (x3)
(x3) edge node {\fgr \laxhi_{m,j}} (x4)
(x1) edge[bend right=\a] node[swap,pos=.7] {\gaC_n\big(1_x \sscs \bang{\laxf_{m_t,i}}_{t \in \ufs{n}} \big)} (y1)
(y1) edge[bend right=10] node[swap,pos=.8] {\gaC_n\big(1_x \sscs \bang{\fgr \laxhi_{m_t,j}}_{t \in \ufs{n}} \big)} (y2)
(y2) edge[bend right=10] node[swap,pos=.1] {\laxf_{n,i}} (y3)
(y3) edge[bend right=\a] node[swap,pos=.3] {\fgr \laxhi_{n,j}} (y4)
(y4) edge[bend right=\a] node[swap,pos=.7] {\fgr \hgr_i \phiAij} (x4)
(y1) edge[bend left=\b] node {\laxf_{n,i}} (z)
(z) edge[bend left=\b] node[pos=.4] {\fgr \gaBi_n\big(1_x \sscs \bang{\laxhi_{m_t,j}}_{t \in \ufs{n}} \big)} (y3)
(z) edge node {\fgr \phiBi} (x3)
;
\end{tikzpicture}\]
In this diagram, the left quadrilateral commutes by the naturality of $\laxf_{n,i}$.  The top right region commutes by the associativity axiom \cref{laxf_associativity} for $(f,\laxf)$.  The bottom right region commutes by the associativity axiom for $(h_i, \laxhi)$ and the functoriality of $\fgr$.
\end{proof}

\begin{lemma}\label{gam_commutativity}
The pair $(fh,\laxfh)$ in \cref{gam_fh} satisfies the commutativity axiom \cref{laxf_com}.
\end{lemma}

\begin{proof}
The commutativity axiom is about two distinct indices in 
\[\ufs{d} = \{1,\ldots,d\} \iso \txcoprod_{v = 1}^k \ufs{d}_v\] 
with $d = \sum_{i \in \ufs{k}}\, d_i$.  There are two cases.  In the first case,  the two indices lie in two different $\ufs{d}_v$'s.  In the second case, the two indices lie in the same $\ufs{d}_v$.

\parhead{Case 1}.  For integers
\[1 \leq i < j \leq k \cq 1 \leq p \leq d_i, \andspace 1 \leq q \leq d_j,\]
we consider objects
\[\left\{
\begin{aligned}
(x,y) &\in \Op(m) \times \Op(n) \phantom{MM} & \adot_\crdot &= \ang{\adot_t}_{t \in \ufs{m}} \in \A_{i,p}^m\\
a_{r,s} &\in \A_{r,s} & \addot_\crdot &= \ang{\addot_u}_{u \in \ufs{n}} \in \A_{j,q}^n
\end{aligned}
\right.\]
for $(r,s) \in \ufs{k} \times \ufs{d}_r \setminus \{(i,p), (j,q)\}$.  The dot in $\adot_t$ indicates that it is an object in $\A_{i,p}$.  We use this notation instead of the more cumbersome notation $a_{i,p,t}$, and similarly for $\addot_u$.  To simplify the next diagram, we use the following notation.
\[\left\{
\begin{gathered}
\begin{aligned}
\hgr_i &= h_i\big(\ang{a_{i,s}}_{s=1}^{p-1} \scs - \scs \ang{a_{i,s}}_{s=p+1}^{d_i} \big) \cn \A_{i,p} \to \B_i\\
\hgr_j &= h_j\big(\ang{a_{j,s}}_{s=1}^{q-1} \scs - \scs \ang{a_{j,s}}_{s=q+1}^{d_j} \big) \cn \A_{j,q} \to \B_j
\end{aligned}\\
\begin{aligned}
\hgr_i \adot_\crdot &= \ang{\hgr_i \adot_t}_{t \in \ufs{m}} \in \B_i^m \phantom{MMM}&
\hgr_j \addot_\crdot &= \ang{\hgr_j \addot_u}_{u \in \ufs{n}} \in \B_j^n\\
b_i &= \gaBi_m\big(x \sscs \hgr_i \adot_\crdot \big) \in \B_i &
b_j &= \gaBj_n\big(y \sscs \hgr_j \addot_\crdot \big) \in \B_j\\ 
b_i' &= \hgr_i \gaAip_m \big(x \sscs \adot_\crdot\big) \in \B_i &
b_j' &= \hgr_j \gaAjq_n \big(y \sscs \addot_\crdot\big) \in \B_j
\end{aligned}
\end{gathered}
\right.\]
\[\left\{
\begin{gathered}
\begin{aligned}
a_{r,\crdot} &= \ang{a_{r,s}}_{s \in \ufs{d}_r} \in \A_{r,1} \times \cdots \times \A_{r,d_r} \forspace r \not\in\{i,j\}\\
\fhat &= f\scalebox{.9}{$\big(\ang{h_r a_{r,\crdot}}_{r=1}^{i-1} \scs - \scs \ang{h_r a_{r,\crdot}}_{r=i+1}^{j-1} \scs - \scs \ang{h_r a_{r,\crdot}}_{r=j+1}^k \big)$} \cn \B_i \times \B_j \to \C
\end{aligned}\\
\begin{aligned}
\fhat\big(\hgr_i \adot_\crdot, ?\big) &= \bang{\fhat\big(\hgr_i \adot_t, ?\big)}_{t \in \ufs{m}} \in \C^m\\
\fhat\big(?, \hgr_j \addot_\crdot\big) &= \bang{\fhat\big(?, \hgr_j \addot_u\big)}_{u \in \ufs{n}} \in \C^n
\end{aligned}
\end{gathered}
\right.\]
The action constraints
\[\begin{split}
b_i \fto{\laxhi_{m,p}} b_i' \andspace b_j \fto{\laxhj_{n,q}} b_j'
\end{split}\]
in, respectively, $\B_i$ and $\B_j$ appear in the next diagram.  By the definition \cref{laxfh_component} of $\laxfh$, the desired commutativity diagram \cref{laxf_com} for $(fh,\laxfh)$ is the boundary of the following diagram in $\C$.
\[\begin{tikzpicture}[baseline={(a.base)}]
\def\h{5.5} \def\u{-1.2} \def\v{-1.5}
\draw[0cell=.7]
(0,0) node (x11) {\gaC_{mn} \big(x \intr y \sscs \bang{\fhat \big(\hgr_i \adot_t, \hgr_j \addot_\crdot \big)}_{t\in \ufs{m}} \big)}
(x11)++(\h,0) node (x12) {\gaC_{mn} \big((y \intr x)\twist_{m,n} \sscs \bang{\fhat \big(\hgr_i \adot_t, \hgr_j \addot_\crdot \big)}_{t\in \ufs{m}} \big)}
(x11)++(0,\v) node (x21) {\gaC_m\big(x \sscs \bang{\gaC_n\big(y \sscs \fhat \big(\hgr_i \adot_t, \hgr_j \addot_\crdot \big) \big) }_{t \in \ufs{m}} \big)}
(x12)++(0,\v) node (x22) {\gaC_n\big(y \sscs \bang{\gaC_m \big(x \sscs \fhat \big(\hgr_i \adot_\crdot, \hgr_j \addot_u \big) \big) }_{u \in \ufs{n}} \big)}
(x21)++(0,\v) node (x31) {\gaC_m\big(x \sscs \fhat \big(\hgr_i \adot_\crdot, b_j \big) \big)}
(x22)++(0,\v) node (x32) {\gaC_n\big(y \sscs \fhat \big(b_i, \hgr_j \addot_\crdot \big) \big)}
(x31)++(0,\v) node (x41) {\gaC_m \big(x \sscs \fhat\big(\hgr_i \adot_\crdot, b_j' \big) \big)}
(x32)++(0,\v) node (x42) {\gaC_n\big(y \sscs \fhat \big(b_i', \hgr_j \addot_\crdot \big) \big)}
(x41)++(0,\v) node (x51) {\fhat\big(b_i, b_j' \big)}
(x42)++(0,\v) node (x52) {\fhat\big(b_i', b_j \big)}
(x51)++(\h/2,-1) node (x) {\fhat\big(b_i', b_j' \big)}
(x51)++(\h/2,1) node (y) {\fhat\big(b_i, b_j \big)}
;
\draw[0cell=.8]
node[between=x21 and x32 at .5] () {\mathbf{c}}
node[between=x51 and x52 at .5] () {\mathbf{f}}
;
\draw[1cell=.7]
(x21) edge node {\phiC} node[swap] {\iso} (x11)
(x11) edge[bend left=10] node {\gaC_{mn} \big(\pcom_{m,n} \twist_{m,n} \sscs 1^{mn}\big)} node[swap] {\iso} (x12)
(x12) edge node {\mathbf{eq} + (\phiC)^{-1}} node[swap] {\iso} (x22)
(x22) edge node {\gaC_n\big(1_y \sscs \bang{\laxf_{m,i} }_{u \in \ufs{n}} \big)} (x32)
(x32) edge node {\gaC_n\big(1_y \sscs \bang{\fhat \big(\laxhi_{m,p} \scs 1\big)}_{u \in \ufs{n}} \big)} (x42)
(x42) edge node {\laxf_{n,j}} (x52)
(x52) edge node {\fhat\big(1, \laxhj_{n,q} \big)} (x)
(x21) edge node[swap] {\gaC_m\big(1_x \sscs \bang{\laxf_{n,j} }_{t \in \ufs{m}} \big)} (x31)
(x31) edge node[swap] {\gaC_m\big(1_x \sscs \bang{\fhat\big(1, \laxhj_{n,q} \big)}_{t \in \ufs{m}} \big)} (x41)
(x41) edge node[swap] {\laxf_{m,i}} (x51)
(x51) edge node[swap] {\fhat\big(\laxhi_{m,p} \scs 1 \big)} (x)
(x31) edge node {\laxf_{m,i}} (y)
(y) edge node[swap,pos=.3] {\fhat\big(1, \laxhj_{n,q} \big)} (x51)
(x32) edge node[swap] {\laxf_{n,j}} (y)
(y) edge node[pos=.3] {\fhat\big(\laxhi_{m,p} \scs 1 \big)} (x52)
;
\end{tikzpicture}\]
In the diagram above, the arrow $\mathbf{eq} + (\phiC)^{-1}$ consists of an instance of the action equivariance axiom \cref{pseudoalg_action_sym} for $\C$, followed by $(\phiC)^{-1}$.  The region labeled $\mathbf{c}$ commutes by the commutativity axiom \cref{laxf_com} for $(f,\laxf)$.  The region labeled $\mathbf{f}$ commutes by the functoriality of $\fhat$.  The left and right triangles commute by the naturality of, respectively, $\laxf_{m,i}$ and $\laxf_{n,j}$.

\parhead{Case 2}.  For integers
\[1 \leq i \leq k \andspace 1 \leq p < q \leq d_i,\]
we consider objects
\[\left\{
\begin{aligned}
(x,y) &\in \Op(m) \times \Op(n) \phantom{MM} & \adot_\crdot &= \ang{\adot_t}_{t \in \ufs{m}} \in \A_{i,p}^m\\
a_{r,s} &\in \A_{r,s} & \addot_\crdot &= \ang{\addot_u}_{u \in \ufs{n}} \in \A_{i,q}^n
\end{aligned}
\right.\]
for $(r,s) \in \ufs{k} \times \ufs{d}_r \setminus \{(i,p), (i,q)\}$, along with the following notation.
\[\left\{
\begin{aligned}
\hhat_i &= h_i \scalebox{.9}{$\big(\bang{a_{i,s}}_{s=1}^{p-1} \scs - \scs  \bang{a_{i,s}}_{s=p+1}^{q-1} \scs - \scs  \bang{a_{i,s}}_{s=q+1}^{d_i} \big)$} \cn \A_{i,p} \times \A_{i,q} \to \B_i\\
a_{r,\crdot} &= \ang{a_{r,s}}_{s\in \ufs{d}_r} \in \A_{r,1} \times \cdots \times \A_{r,d_r} \forspace r \neq i\\
\fgr &= f\big(\ang{h_r a_{r,\crdot}}_{r=1}^{i-1} \scs - \scs \ang{h_r a_{r,\crdot}}_{r=i+1}^k \big) \cn \B_i \to \C\\
b_{t,u} &= \hhat_i\big( \adot_t \scs \addot_u \big) \in \B_i
\end{aligned}
\right.\]
\[\left\{
\begin{aligned}
b_{\crdot,u} &= \ang{b_{t,u}}_{t \in \ufs{m}} \in \B_i^m \phantom{MMM} & 
b_{\crdot,u}^\ga &= \gaBi_m\big(x \sscs b_{\crdot,u} \big) \in \B_i\\
b_{t,\crdot} &= \ang{b_{t,u}}_{u \in \ufs{n}} \in \B_i^n &
b_{t,\crdot}^\ga &= \gaBi_n\big(y \sscs b_{t,\crdot} \big) \in \B_i\\
\fgr b_{\crdot,u} &= \bang{\fgr b_{t,u}}_{t \in \ufs{m}} \in \C^m &
\adot &= \gaAip_m\big(x \sscs \adot_\crdot\big) \in \Aip\\
\fgr b_{t,\crdot} &= \bang{\fgr b_{t,u}}_{u \in \ufs{n}} \in \C^n &
\addot &= \gaAiq_n\big(y \sscs \addot_\crdot) \in \Aiq
\end{aligned}
\right.\]
\[\left\{
\begin{aligned}
\fgr_{x \intr y} &= \fgr \gaBi_{mn} \big(x \intr y \sscs \ang{b_{t,\crdot}}_{t \in \ufs{m}} \big) \in \C\\
\fgr_{x,y} &= \fgr \gaBi_m\big(x \sscs \ang{b_{t,\crdot}^\ga }_{t \in \ufs{m}} \big) \in \C\\
\fgr_{y \intr x} &= \fgr \gaBi_{mn} \big( (y \intr x) \twist_{m,n} \sscs \ang{b_{t,\crdot}}_{t \in \ufs{m}} \big) \in \C\\
\fgr_{y,x} &= \fgr \gaBi_n \big(y \sscs \ang{b_{\crdot,u}^\ga}_{u \in \ufs{n}} \big) \in \C
\end{aligned}
\right.\]
By the definition \cref{laxfh_component} of $\laxfh$, the desired commutativity diagram \cref{laxf_com} for $(fh,\laxfh)$ is the boundary of the following diagram in $\C$.
\[\begin{tikzpicture}[xscale=1,yscale=1,baseline={(a.base)}]
\def\a{2.5} \def\b{.75} \def\h{7} \def\g{.5} \def\u{-1.2} \def\v{-1.5}
\draw[0cell=.7]
(0,0) node (x11) {\gaC_{mn} \big(x \intr y \sscs \ang{\fgr b_{t,\crdot}}_{t \in \ufs{m}} \big)}
(x11)++(\h-2*\g,0) node (x12) {\gaC_{mn} \big((y \intr x)\twist_{m,n} \sscs \ang{\fgr b_{t,\crdot}}_{t \in \ufs{m}} \big)}
(x11)++(-\g,\v) node (x21) {\gaC_m\big(x \sscs \bang{\gaC_n(y \sscs \fgr b_{t,\crdot} )}_{t \in \ufs{m}} \big)}
(x12)++(\g,\v) node (x22) {\gaC_n\big(y \sscs \bang{\gaC_m(x \sscs \fgr b_{\crdot, u} ) }_{u \in \ufs{n}} \big)}
(x21)++(0,\v) node (x31) {\gaC_m\big(x \sscs \ang{\fgr b_{t,\crdot}^\ga }_{t \in \ufs{m}} \big)}
(x22)++(0,\v) node (x32) {\gaC_n\big(y \sscs \ang{\fgr b_{\crdot,u}^\ga }_{u \in \ufs{n}} \big)}
(x31)++(0,\v) node (x41) {\gaC_m\big(x \sscs \ang{\fgr \hhat_i(\adot_t, \addot) }_{t \in \ufs{m}} \big)}
(x32)++(0,\v) node (x42) {\gaC_n\big(y \sscs \ang{\fgr \hhat_i(\adot, \addot_u) }_{u \in \ufs{n}} \big)}
(x41)++(0,\v) node (x51) {\fgr \gaBi_m\big(x \sscs \ang{\hhat_i(\adot_t, \addot) }_{t \in \ufs{m}} \big)}
(x42)++(0,\v) node (x52) {\fgr \gaBi_n\big(y \sscs \ang{\hhat_i(\adot, \addot_u) }_{u \in \ufs{n}} \big)}
(x51)++(\h/2,\u) node (x) {\fgr \hhat_i (\adot, \addot)}
(x21)++(\a,-\b) node (y11) {\fgr_{x \intr y}}
(x22)++(-\a,-\b) node (y12) {\fgr_{y \intr x}}
(x41)++(\a,0) node (y21) {\fgr_{x,y}}
(x42)++(-\a,0) node (y22) {\fgr_{y,x}}
;
\draw[0cell=.8]
node[between=x51 and x52 at .5, shift={(0,-.3)}] () {\mathbf{c}}
node[between=x31 and y11 at .6] () {\mathbf{a}}
node[between=x32 and y12 at .6] () {\mathbf{eq}+\mathbf{a}}
;
\draw[1cell=.7]
(x21) edge node[pos=.3] {\phiC} node[swap,pos=.6] {\iso} (x11)
(x11) edge node {\gaC_{mn} \big(\pcom_{m,n} \twist_{m,n} \sscs 1^{mn} \big)} node[swap] {\iso} (x12)
(x12) edge node[pos=.75] {\mathbf{eq} + (\phiC)^{-1}} node[swap,pos=.4] {\iso} (x22)
(x22) edge node {\gaC_n\big(1_y \sscs \ang{\laxf_{m,i} }_{u \in \ufs{n}} \big)} (x32)
(x32) edge node {\gaC_n\big(1_y \sscs \ang{\fgr \laxhi_{m,p} }_{u \in \ufs{n}} \big)} (x42)
(x42) edge node {\laxf_{n,i}} (x52)
(x52) edge node {\fgr \laxhi_{n,q}} (x)
(x21) edge node[swap] {\gaC_m \big(1_x \sscs \ang{\laxf_{n,i}}_{t \in \ufs{m}} \big)} (x31)
(x31) edge node[swap] {\gaC_m \big(1_x \sscs \ang{\fgr \laxhi_{n,q} }_{t \in \ufs{m}} \big)} (x41)
(x41) edge node[swap] {\laxf_{m,i}} (x51)
(x51) edge node[swap] {\fgr \laxhi_{m,p}} (x)
(y11) edge node {} (y12)
(x11) edge node {\laxf_{mn,i}} (y11)
(y11) edge node {\fgr (\phiBi)^{-1}} (y21)
(y21) edge node[pos=.2] {\fgr \gaBi_m\big(1_x \sscs \ang{\laxhi_{n,q}}_{t \in \ufs{m}} \big)} (x51)
(x31) edge node {\laxf_{m,i}} (y21)
(x12) edge node[swap] {\laxf_{mn,i}} (y12)
(y12) edge node {} (y22)
(y22) edge node[swap,pos=.6] {\fgr \gaBi_n\big(1_y \sscs \ang{\laxhi_{m,p}}_{u \in \ufs{n}} \big)} (x52)
(x32) edge node[swap] {\laxf_{n,i}} (y22)
;
\end{tikzpicture}\]
In the diagram above, the arrow $\mathbf{eq} + (\phiC)^{-1}$ consists of an instance of the action equivariance axiom \cref{pseudoalg_action_sym} for $\C$, followed by $(\phiC)^{-1}$.  The two unlabeled arrows in the interior are given as follows, where $\mathbf{eq}$ is an instance of the action equivariance axiom \cref{pseudoalg_action_sym} for $\B_i$.
\[\begin{tikzpicture}[xscale=1,yscale=1,baseline={(a.base)}]
\draw[0cell=.9]
(0,0) node (y11) {\fgr_{x \intr y}}
(y11)++(5,0) node (y12) {\fgr_{y \intr x}}
(y12)++(0,-1.5) node (y22) {\fgr_{y,x}}
;
\draw[1cell=.9]
(y11) edge node {\fgr \gaBi_{mn}\big(\pcom_{m,n} \twist_{m,n} \sscs 1^{mn} \big)} (y12)
(y12) edge node {\fgr\big(\mathbf{eq} + (\phiBi)^{-1})} (y22)
;
\end{tikzpicture}\]
The region labeled $\mathbf{c}$ commutes by the commutativity axiom \cref{laxf_com} for $(h_i, \laxhi)$ and the functoriality of $\fgr$.  The region labeled $\mathbf{a}$ commutes by the associativity axiom \cref{laxf_associativity} for $(f,\laxf)$.  The region labeled $\mathbf{eq}+\mathbf{a}$ commutes by the equivariance axiom \cref{laxf_eq} and the associativity axiom \cref{laxf_associativity} for $(f,\laxf)$.  The three unlabeled regions, going clockwise from the top, commute by the naturality of, respectively, $\laxf_{mn,i}$, $\laxf_{n,i}$, and $\laxf_{m,i}$.
\end{proof}

\cref{gam_basept_unity_eq,gam_associativity,gam_commutativity} prove that the pair $(fh,\laxfh)$ defined in \cref{gam_fh} is a $d$-lax $\Op$-morphism in the sense of \cref{def:k_laxmorphism}.  The next result shows that $\gam$ preserves the $G$-action (\cref{def:k_laxmorphism_g}).

\begin{lemma}\label{multo_ga_g_obj}
In the context of \cref{def:gam_Omorphisms}, for each element $g \in G$, there is an equality of $d$-lax $\Op$-morphisms
\begin{equation}\label{multo_gagob}
\gam\big( g \cdot (f,\laxf) ; \ang{g \cdot (h_i, \laxhi) }_{i \in \ufs{k}} \big) 
= g \cdot \gam\big( (f,\laxf) ; \ang{(h_i, \laxhi) }_{i \in \ufs{k}} \big).
\end{equation}
\end{lemma}

\begin{proof}
We need to show that the two $d$-lax $\Op$-morphisms in \cref{multo_gagob} have (i) the same underlying functor and (ii) the same action constraints.

\parhead{Underlying functor}.  Using the notation in \cref{assum:multpso_comp,gam_fh_functor}, we consider the following diagram of functors.
\begin{equation}\label{multo_gagob_diag}
\begin{tikzpicture}[vcenter]
\def\v{1.4} \def\t{25}
\draw[0cell=1]
(0,0) node (a) {\A_{[1,k]}}
(a)++(3,0) node (b) {\B_{[1,k]}}
(b)++(3,0) node (c) {\C}
(a)++(0,\v) node (a1) {\A_{[1,k]}}
(b)++(0,\v) node (b1) {\B_{[1,k]}}
(c)++(0,\v) node (c1) {\C}
;
\draw[1cell=.9]  
(a) edge node {\txprod_{i=1}^k\, h_i} (b)
(b) edge node {f} (c)
(a1) edge node[swap] {(\ginv)^d} (a)
(b) edge[bend left=\t] node {g^k} (b1)
(b1) edge[bend left=\t] node {(\ginv)^k} (b)
(c) edge node[swap] {g} (c1)
;
\end{tikzpicture}
\end{equation}
By \cref{klax_mor_g}, the underlying functor of the left-hand side of \cref{multo_gagob} is the composite in \cref{multo_gagob_diag} including the loop in the middle.  The underlying functor of the right-hand side of \cref{multo_gagob} is the composite in \cref{multo_gagob_diag} excluding the loop in the middle.  Since the composite of the loop is the identity functor on $\B_{[1,k]}$, the two sides of \cref{multo_gagob} have the same underlying functor.

\parhead{Action constraints}.  We use the notation in \cref{elldj,laxfh_objects,htilftil}.  By \cref{laxgf_z,laxfh_component}, the component of $\lax^{(g \cdot f)(g \cdot h)}_{n,\ell}$---which is the $(n,\ell)$-action constraint of the left-hand side of \cref{multo_gagob}---at the objects in \cref{laxfh_objects} is given by the following composite morphism in $\C$.  In the following diagram, each inactive variable $a_{p,r}$, with $(p,r) \neq (i,j)$, is replaced by $\ginv a_{p,r}$ in the definitions of $\fgr$ and $\hgr_i$.
\begin{equation}\label{multo_ga_obj_const}
\begin{tikzpicture}[vcenter]
\def\v{-1} \def\u{-1.5}
\draw[0cell=.8]
(0,0) node (a11) {\gaC_n\big(x; g\fgr\hgr_i \ginv a_{i,j,\crdot} \big)}
(a11)++(4.6,0) node (a12) {g\fgr\hgr_i \ginv \gaAij_n(x; a_{i,j,\crdot})}
(a11)++(0,\v) node (a21) {g\gaC_n\big(\ginv x; \fgr\ginv g \hgr_i \ginv a_{i,j,\crdot} \big)}
(a12)++(0,\v) node (a22) {(g\fgr \ginv) \big(g\hgr_i \gaAij_n(\ginv x; \ginv a_{i,j,\crdot}) \big)}
(a21)++(0,\u) node (a31) {g\fgr \gaBi_n\big(\ginv x; \ginv g\hgr_i \ginv a_{i,j,\crdot} \big)}
(a22)++(0,\u) node (a32) {(g\fgr\ginv) \big(g\gaBi_n(\ginv x; \hgr_i \ginv a_{i,j,\crdot}) \big)}
;
\draw[1cell=.8]
(a11) edge node {\lax^{(g \cdot f)(g \cdot h)}_{n,\ell}} (a12)
(a11) edge[equal] (a21)
(a21) edge node[swap] {g\laxf_{n,i}} (a31)
(a31) edge[equal] (a32)
(a32) edge[shorten >=-.5ex] node[swap] {(g\fgr\ginv)(g \laxhi_{n,j})} (a22)
(a22) edge[equal] (a12)
;
\end{tikzpicture}
\end{equation}
On the other hand, by the functoriality of the $g$-action on $\C$, \cref{laxgf_z}, and \cref{laxfh_component}, the composite in \cref{multo_ga_obj_const} is also equal to the $(n,\ell)$-action constraint of the right-hand side of \cref{multo_gagob}, evaluated at the objects in \cref{laxfh_objects}.  Thus, the two sides of \cref{multo_gagob} have the same action constraints.
\end{proof}

\subsection*{Multicategorical Composition of $\Op$-Transformations}

\cref{def:gam_Omorphisms} specifies the object assignment of the multicategorical composition of $\MultpsO$.  The next definition specifies the morphism assignment.

\begin{definition}[$\gam$ on Morphisms]\label{def:gam_Otr}
Under \cref{assum:multpso_comp}, suppose we are given
\[\begin{tikzpicture}[baseline={(a.base)}]
\draw[0cell]
(0,0) node (a) {\phantom{\B_i}} 
(a)++(-.11,0) node () {\A_{i\crdot}}
(a)++(2,0) node (b) {\B_i}
(b)++(1.5,0) node () {\text{and}}
;
\draw[1cell]
(a) edge[bend left] node {(h_i,\laxhi)} (b)
(a) edge[bend right] node[swap] {(h_i',\laxhip)} (b)
;
\draw[2cell]
node[between=a and b at .45, rotate=-90, 2label={above,\theta_i}] {\Rightarrow}
;
\begin{scope}[shift={(5,0)}]
\draw[0cell]
(0,0) node (a) {\B_\crdot}
(a)++(2,0) node (b) {\C}
;
\draw[1cell]
(a) edge[bend left] node {(f,\laxf)} (b)
(a) edge[bend right] node[swap] {(f',\laxfp)} (b)
;
\draw[2cell]
node[between=a and b at .45, rotate=-90, 2label={above,\theta}] {\Rightarrow}
;
\end{scope}
\end{tikzpicture}\]
with
\begin{itemize}
\item $\theta$ a $k$-ary $\Op$-transformation (\cref{def:kary_transformation}) between $k$-lax $\Op$-morphisms and 
\item $\theta_i$ a $d_i$-ary $\Op$-transformation between $d_i$-lax $\Op$-morphisms for $i \in \ufs{k}$.
\end{itemize}  
We define the \emph{composite} $d$-ary $\Op$-transformation
\begin{equation}\label{gam_Otr}
\begin{tikzpicture}[baseline={(a.base)}]
\def\t{25}
\draw[0cell]
(0,0) node (a) {\A_{\crdots}}
(a)++(3.5,0) node (b) {\C}
;
\draw[1cell=.85]
(a) edge[bend left=\t] node {\gam\big((f,\laxf) \sscs \ang{(h_i, \laxhi)}_{i \in \ufs{k}} \big)} (b)
(a) edge[bend right=\t] node[swap] {\gam\big((f',\laxfp) \sscs \ang{(h_i', \laxhip)}_{i \in \ufs{k}} \big)} (b)
;
\draw[2cell=1]
node[between=a and b at .25, rotate=-90, 2label={above,\gam\big(\theta ; \ang{\theta_i}_{i \in \ufs{k}}\big)}] {\Rightarrow}
;
\end{tikzpicture}
\end{equation}
as the following horizontal composite natural transformation.
\begin{equation}\label{gam_thetas}
\begin{tikzpicture}[baseline={(a.base)}]
\def\t{25}
\draw[0cell]
(0,0) node (a) {\txprod_{i=1}^k \txprod_{j=1}^{d_i}\,  \Aij}
(a)++(.9,0) node (a') {\phantom{C}}
(a')++(2.8,0) node (b') {\phantom{C}}
(b')++(.45,0) node (b) {\txprod_{i=1}^k\, \B_i}
(b)++(.45,0) node (b'') {\phantom{C}}
(b'')++(2,0) node (c) {\C}
;
\draw[1cell=.85]
(a') edge[bend left=\t] node {\txprod_{i=1}^k h_i} (b')
(a') edge[bend right=\t] node[swap] {\txprod_{i=1}^k h_i'} (b')
(b'') edge[bend left] node {f} (c)
(b'') edge[bend right] node[swap] {f'} (c)
;
\draw[2cell=1]
node[between=a' and b' at .3, rotate=-90, 2label={above,\txprod_{i=1}^k \theta_i}] {\Rightarrow}
node[between=b'' and c at .45, rotate=-90, 2label={above,\theta}] {\Rightarrow}
;
\end{tikzpicture}
\end{equation}
\cref{gam_Otr_welldef} below proves that this definition is well defined.
\end{definition}

\begin{lemma}\label{gam_Otr_welldef}
$\gam\big(\theta ; \ang{\theta_i}_{i \in \ufs{k}}\big)$ in \cref{gam_Otr} is a $d$-ary $\Op$-transformation.
\end{lemma}

\begin{proof}
To prove the multilinearity axiom \cref{ktransform_multilinearity} for $\gam\big(\theta ; \ang{\theta_i}_{i \in \ufs{k}}\big)$, we use the notation in \cref{elldj,laxfh_objects,htilftil,laxfh_component} along with the following notation.
\[\left\{
\scalebox{.9}{$\begin{aligned}
b_i &= \hgr_i \gaAij_n\big(x \sscs a_{i,j,\crdot}\big) \phantom{M} & b_{<i} &= \bang{h_p a_{p,\crdot}}_{p=1}^{i-1} \phantom{M} & b_{>i} &= \bang{h_p a_{p,\crdot}}_{p=i+1}^k\\
b_i' &= \hgr_i' \gaAij_n\big(x \sscs a_{i,j,\crdot}\big) & b_{<i}' &= \bang{h_p' a_{p,\crdot}}_{p=1}^{i-1} & b_{>i}' &= \bang{h_p' a_{p,\crdot}}_{p=i+1}^k\\
&& \theta_{<i} &= \ang{\theta_p}_{p=1}^{i-1} & \theta_{>i} &= \ang{\theta_p}_{p=i+1}^k\\
\fch &= f(1^{i-1}, - , 1^{k-i}) & \fch' &= f'(1^{i-1}, - , 1^{k-i}) &&
\end{aligned}$}
\right.\]
The desired multilinearity diagram is the boundary of the following diagram in $\C$.
\[\begin{tikzpicture}
\def\t{7} \def\g{3} \def\h{3.5} \def\u{1.5} \def\v{\u/2} \def\w{1.5em} \def\x{-1em}
\draw[0cell=.7]
(0,0) node (a1) {\gaC_n\big(x \sscs \bang{f\big(b_{<i}  \scs \hgr_i a_{i,j,m} \scs b_{>i} \big) }_{m \in \ufs{n}} \big)}
(a1)++(\h,\v) node (a2) {f\big(b_{<i} \spc \gaBi_n\big(x \sscs \hgr_i a_{i,j,\crdot} \big), b_{>i} \big)}
(a2)++(\g,-\v) node (a3) {f\big(b_{<i} \spc b_i, b_{>i} \big)}
(a1)++(0,-\u) node (b1) {\gaC_n\big(x \sscs \bang{f\big(b_{<i}'  \scs \hgr_i' a_{i,j,m} \scs b_{>i}' \big) }_{m \in \ufs{n}} \big)}
(b1)++(\h,\v) node (b2) {f\big(b_{<i}' \spc \gaBi_n\big(x \sscs \hgr_i' a_{i,j,\crdot} \big), b_{>i}' \big)}
(b2)++(\g,-\v) node (b3) {f\big(b_{<i}' \spc b_i' , b_{>i}' \big)}
(b1)++(0,-\u) node (c1) {\gaC_n\big(x \sscs \bang{f' \big(b_{<i}'  \scs \hgr_i' a_{i,j,m} \scs b_{>i}' \big) }_{m \in \ufs{n}} \big)}
(c1)++(\h,\v) node (c2) {f' \big(b_{<i}' \spc \gaBi_n\big(x \sscs \hgr_i' a_{i,j,\crdot} \big), b_{>i}' \big)}
(c2)++(\g,-\v) node (c3) {f' \big(b_{<i}' \spc b_i' , b_{>i}' \big)}
;
\draw[1cell=.7]
(a1) edge[bend left=\t] node {\laxf_{n,i}} (a2)
(a2) edge[bend left=\t] node[pos=.2] {\fch \laxhi_{n,j}} (a3)
(b1) edge[bend left=\t] node {\laxf_{n,i}} (b2)
(b2) edge[bend left=\t] node[pos=.2] {\fch \laxhip_{n,j}} (b3)
(c1) edge[bend left=\t] node {\laxfp_{n,i}} (c2)
(c2) edge[bend left=\t] node[pos=.2] {\fch' \laxhip_{n,j}} (c3)
(a1) edge[transform canvas={xshift=\w}] node[swap] {\gaC_n\big(1_x \sscs \bang{f\big( \theta_{<i} \scs \theta_i \scs \theta_{>i} \big)}_{m \in \ufs{n}} \big)} (b1)
(b1) edge[transform canvas={xshift=\w}] node[swap] {\gaC_n\big(1_x \sscs \theta^n \big)} (c1)
(a2) edge[transform canvas={xshift=-3em}] node {f\big(\theta_{<i} \scs \gaBi_n\big(1_x \sscs \theta_i^n\big) \scs \theta_{>i} \big)} (b2)
(b2) edge[transform canvas={xshift=-3em}] node {\theta} (c2)
(a3) edge[transform canvas={xshift=\x}] node {f\big(\theta_{<i} \scs \theta_i \scs \theta_{>i} \big)} (b3)
(b3) edge[transform canvas={xshift=\x}] node {\theta} (c3)
;
\end{tikzpicture}\]
In this diagram, the upper left region commutes by the naturality of $\laxf_{n,i}$.  The lower right region commutes by the naturality of $\theta$.  The lower left region commutes by the multilinearity axiom \cref{ktransform_multilinearity} for $\theta$.  The upper right region commutes by the multilinearity axiom for $\theta_i$ and the functoriality of $f$.
\end{proof}

\cref{multo_ga_g_obj} proves that $\gam$ preserves the $G$-action on objects.  The following observation proves that $\gam$ preserves the $G$-action on morphisms (\cref{def:kary_transformation}).

\begin{lemma}\label{multo_ga_g_mor}
In the context of \cref{def:gam_Otr}, for each element $g \in G$, there is an equality of $d$-ary $\Op$-transformations
\begin{equation}\label{multo_gagmor}
\gam\big( g \cdot \theta ; \ang{g \cdot \theta_i}_{i \in \ufs{k}} \big) 
= g \cdot \gam( \theta ; \ang{\theta_i}_{i \in \ufs{k}}).
\end{equation}
\end{lemma}

\begin{proof}
By \cref{O_tr_g}, the $g$-action $g \cdot \theta$ is given by the conjugation $g$-action.  By \cref{gam_thetas}, the composite $\gam(\theta; \ang{\theta_i}_{i \in \ufs{k}})$ is given by the horizontal composition $\theta \ast \txprod_{i=1}^k\, \theta_i$.  The equality \cref{multo_gagmor} is proved by reusing the diagram \cref{multo_gagob_diag} and the paragraph immediately under that diagram, with $(f,\ang{h_i}_{i \in \ufs{k}})$ replaced by $(\theta, \ang{\theta_i}_{i \in \ufs{k}})$.
\end{proof}

\subsection*{Multicategorical Composition $G$-Functors}

Next, we define the multicategorical composition of $\MultpsO$.  Recall its $k$-ary multimorphism $G$-categories in \cref{def:MultpsO_karycat}. 

\begin{definition}\label{def:gam_functor}
Under \cref{assum:multpso_comp}, we define the $G$-functor
\[\MultpsO\big(\B_\crdot \sscs \C\big ) \times \txprod_{i \in \ufs{k}}\, \MultpsO\big(\A_{i\crdot} \sscs \B_i \big) 
\fto{\gam} \MultpsO\big( \A_{\crdots} \sscs \C\big)\]
with object and morphism assignments given by, respectively, \cref{def:gam_Omorphisms,def:gam_Otr}.
\begin{itemize}
\item The object assignment is well defined and $G$-equivariant by \cref{gam_basept_unity_eq,gam_associativity,gam_commutativity,multo_ga_g_obj}. 
\item The morphism assignment is well defined and $G$-equivariant by \cref{gam_Otr_welldef,multo_ga_g_mor}. 
\item The functoriality of $\gam$ follows from \cref{gam_thetas}, since horizontal composition of natural transformations preserves identities and vertical composition.
\end{itemize} 
Moreover, we define the pseudo and strict variants of the $G$-functor $\gam$ in the same way by replacing $\MultpsO$ with, respectively, $\MultpspsO$ and $\MultstO$ in \cref{def:MultpsO_karycat} \eqref{multpspso_k}.  The resulting $G$-functors $\gam$ are well defined because, if $\laxf_{n,i}$ and $\laxhi_{n,j}$ in the diagram \cref{laxfh_component} are isomorphisms, respectively, identities, then so is $\laxfh_{n,\ell}$.
\end{definition}

\section{$\Gcat$-Multicategory of Pseudoalgebras}
\label{sec:multpsodef}

This section constructs
\begin{itemize}
\item the $\Gcat$-multicategories $\MultvO$ for $\va \in \{\sflax,\sfps,\sfst\}$ and
\item the $\Cat$-multicategories $\MultvO^G$ induced by the $G$-fixed subcategory functor on the $\Gcat$-enrichment.
\end{itemize}
Each of these enriched multicategories has $\Op$-pseudoalgebras as objects.

\secoutline
\begin{itemize}
\item \cref{def:multicatO} defines the $\Gcat$-multicategory $\MultvO$.
\item \cref{thm:multpso} proves in detail that $\MultvO$ satisfies the axioms of a $\Gcat$-multicategory.
\item \cref{thm:multpso_gfixed} records the fact that $\MultvO$ yields a $\Cat$-multicategory $\MultvO^G$ with the same objects---namely, $\Op$-pseudoalgebras---by replacing each multimorphism $G$-category with its $G$-fixed subcategory.
\item \cref{ex:underlying_iicat} records the fact that the underlying 2-category of $\MultvO^G$ is the 2-category $\AlgpsvO$ in \cref{oalgps_twocat}.
\item \cref{ex:multBE,ex:multGBE} apply \cref{thm:multpso,thm:multpso_gfixed} to the Barratt-Eccles operad $\BE$ and the $G$-Barratt-Eccles operad $\GBE$.  The objects of the resulting $\Gcat$-multicategories and $\Cat$-multicategories are, respectively,
\begin{itemize}
\item $\BE$-pseudoalgebras, which are, up to a 2-equivalence, naive symmetric monoidal $G$-categories, and
\item $\GBE$-pseudoalgebras, which are genuine symmetric monoidal $G$-categories.  
\end{itemize}
\end{itemize}

\subsection*{$\Gcat$-Multicategory Structure}

For the next definition, recall enriched multicategories from \cref{def:enr-multicategory} and the Cartesian closed category $\Gcat$ of small $G$-categories (\cref{expl:Gcat_closed}).  

\begin{definition}\label{def:multicatO}
Suppose $(\Op,\ga,\opu,\pcom)$ is a pseudo-commutative operad in $\Gcat$ (\cref{def:pseudocom_operad}) for an arbitrary group $G$.  We define the $\Gcat$-multicategory $\MultpsO$\index{pseudoalgebra!multicategory}\index{multicategory!pseudoalgebra} as follows.
\begin{description}
\item[Objects] The objects of $\MultpsO$ are $\Op$-pseudoalgebras (\cref{def:pseudoalgebra}).
\item[Multimorphism $G$-categories] For $\Op$-pseudoalgebras $\ang{\A_i}_{i \in \ufs{k}}$ and $\B$, the $k$-ary multimorphism $G$-category
\[\MultpsO\scmap{\ang{\A_i}_{i \in \ufs{k}} ; \B}\]
is defined in \cref{def:MultpsO_karycat}. 
\begin{itemize}
\item Its objects are $k$-lax $\Op$-morphisms (\cref{def:k_laxmorphism}). 
\item Its morphisms are $k$-ary $\Op$-transformations (\cref{def:kary_transformation}).
\item The conjugation $G$-action on objects and morphisms are defined in, respectively, \cref{def:k_laxmorphism_g,def:kary_transformation}.
\end{itemize}
\item[Units] For an $\Op$-pseudoalgebra $\B$, the unit $G$-functor
\begin{equation}\label{Opsalg_unit}
\boldone \fto{\opu_\B} \MultpsO(\B;\B)
\end{equation}
is given by the identity $\Op$-morphism $1 \cn \B \to \B$, which consists of the identity functor on $\B$ and identity action constraints.  The unit $\opu_\B$ is $G$-equivariant because the identity $\Op$-morphism is fixed by the conjugation $G$-action.
\item[Symmetric group action] For each permutation $\sigma \in \Sigma_k$, the right $\sigma$-action 
\[\MultpsO\scmap{\ang{\A_i}_{i \in \ufs{k}};\B} \fto[\iso]{\sigma} \MultpsO\scmap{\ang{\A_{\sigma(i)}}_{i \in \ufs{k}};\B}\]
is the isomorphism of $G$-categories given in \cref{multpso_sym_functor}.
\item[Composition] The composition $G$-functor $\gam$ is defined in \cref{def:gam_functor}.
\end{description}
This finishes the definition of $\MultpsO$.  Moreover, we define the $\Gcat$-multicategories 
\[\MultpspsO \andspace \MultstO\]
with the same objects by replacing each multimorphism $G$-category of $\MultpsO$ by, respectively, its pseudo and strict variants in \cref{def:MultpsO_karycat} \eqref{multpspso_k}.
\end{definition}

\begin{theorem}\label{thm:multpso}
In the context of \cref{def:multicatO}, for each of the three variants $\va \in \{\sflax,\sfps,\sfst\}$, $\MultvO$ is a $\Gcat$-multicategory.
\end{theorem}

\begin{proof}
We verify the $\Gcat$-multicategory axioms in \cref{def:enr-multicategory} for $\MultpsO$.  The same argument applies to $\MultpspsO$ and $\MultstO$ by restricting to, respectively, $k$-ary $\Op$-pseudomorphisms and $k$-ary strict $\Op$-morphisms.

\parhead{Unity}.  The unit of each $\Op$-pseudoalgebra consists of an identity functor and identity action constraints.  Thus, the unity axioms \cref{enr-multicategory-right-unity,enr-multicategory-left-unity} follow from 
\begin{itemize}
\item the definitions \cref{gam_fh_functor,laxfh_component} of a composite $d$-lax $\Op$-morphism and
\item the definition \cref{gam_thetas} of a composite $d$-ary $\Op$-transformation.
\end{itemize}

\parhead{Symmetry}.  The symmetric group action axioms \cref{enr-multicategory-symmetry} are verified in \cref{multpso_sym_functor}.

\parhead{Associativity}.  In addition to 
\begin{itemize}
\item the $k$-lax $\Op$-morphism $f \cn \B_\crdot \to \C$ and
\item the $d_i$-lax $\Op$-morphisms $h_i \cn \A_{i\crdot} \to \B_i$ for $i \in \ufs{k}$
\end{itemize}
in \cref{def:gam_Omorphisms}, suppose we are given
\begin{itemize}
\item $\Op$-pseudoalgebras $\Z_{i,j,s}$ for $(i,j,s) \in \ufs{k} \times \ufs{d}_i \times \ufs{e}_{i,j}$ and
\item $e_{i,j}$-lax $\Op$-morphisms
\[\bang{\Z_{i,j,s}}_{s \in \ufs{e}_{i,j}} \fto{(y_{i,j} \scs \laxyij)} \Aij\]
for $(i,j) \in \ufs{k} \times \ufs{d}_i$.
\end{itemize}
The associativity axiom \cref{enr-multicategory-associativity} on objects asks for the following equality of $e$-lax $\Op$-morphisms, where $e = \sum_{i \in \ufs{k}} \sum_{j \in \ufs{d}_i} e_{i,j}$.
\begin{equation}\label{gam_fhy}
\gam\big(\gam\big(f \sscs \ang{h_i}_{i \in \ufs{k}}) \sscs \ang{\ang{y_{i,j}}_{j \in \ufs{d}_i}}_{i \in \ufs{k}} \big) 
= \gam\big(f \sscs \bang{\gam (h_i \sscs \ang{y_{i,j}}_{j \in \ufs{d}_i} ) }_{i \in \ufs{k}} \big)
\end{equation}
By \cref{gam_fh_functor}, the underlying functor of each side of \cref{gam_fhy} is the following composite.
\begin{equation}\label{fhiyij}
\begin{tikzpicture}[baseline={(a.base)}]
\draw[0cell=.85]
(0,0) node (a) {\prod_{i=1}^k \prod_{j=1}^{d_i} \prod_{s=1}^{e_{i,j}} \Z_{i,j,s}}
(a)++(3.2,0) node (b) {\prod_{i=1}^k \prod_{j=1}^{d_i} \Aij}
(b)++(2.2,0) node (c) {\prod_{i=1}^k \B_i}
(c)++(1.6,0) node (d) {\C}
;
\draw[1cell=.85]
(a) edge node {\txprod_i \txprod_j \, y_{i,j}} (b)
(b) edge node {\txprod_i \, h_i} (c)
(c) edge node {f} (d)
;
\end{tikzpicture}
\end{equation}
Moreover, since composition of $k$-ary $\Op$-transformations for $k \geq 0$ is defined by horizontal composition of natural transformations \cref{gam_thetas}, the diagram \cref{fhiyij} also shows that $\gam$ is associative on morphisms.

To prove that the two sides of \cref{gam_fhy} have the same action constraints, suppose $(i,j,s) \in \ufs{k} \times \ufs{d}_i \times \ufs{e}_{i,j}$ and
\[\ell = \big( \txsum_{p=1}^{i-1} \txsum_{q=1}^{d_p}\, e_{p,q}\big) + \big(\txsum_{q=1}^{j-1}\, e_{i,q} \big) + s.\]
For $n \geq 0$, we consider objects
\[x \in \Op(n) \cq z_{p,q,r} \in \Z_{p,q,r}, \andspace \zdot_\crdot = \ang{\zdot_m}_{m \in \ufs{n}} \in \Z_{i,j,s}^n\]
for $(p,q,r) \in \ufs{k} \times \ufs{d}_p \times \ufs{e}_{p,q} \setminus \{(i,j,s)\}$, along with the following notation.
\[\left\{
\begin{aligned}
\ygr_{i,j} &= y_{i,j} \big(\ang{z_{i,j,r}}_{r=1}^{s-1} \scs - \scs \ang{z_{i,j,r}}_{r=s+1}^{e_{i,j}} \big) \cn \Z_{i,j,s} \to \Aij\\
z_{p,q,\crdot} &= \ang{z_{p,q,r}}_{r \in \ufs{e}_{p,q}} \forspace (p,q) \neq (i,j)\\
\hgr_i &= h_i\big(\ang{y_{i,q} z_{i,q,\crdot}}_{q=1}^{j-1} \scs - \scs \ang{y_{i,q} z_{i,q,\crdot}}_{q=j+1}^{d_i} \big) \cn \Aij \to \B_i\\
b_p &= h_p \bang{y_{p,q} z_{p,q,\crdot} }_{q \in \ufs{d}_p} \in \B_p \forspace p \neq i\\
\fgr &= f\big(\ang{b_p}_{p=1}^{i-1} \scs - \scs \ang{b_p}_{p=i+1}^k \big) \cn \B_i \to \C\\
\ygr_{i,j} \zdot_\crdot &= \bang{\ygr_{i,j} \zdot_m}_{m \in \ufs{n}} \in \Aij^n\\
\hgr_i \ygr_{i,j} \zdot_\crdot &= \bang{\hgr_i \ygr_{i,j} \zdot_m}_{m \in \ufs{n}} \in \B_i^n\\
\fgr \hgr_i \ygr_{i,j} \zdot_\crdot &= \bang{\fgr \hgr_i \ygr_{i,j} \zdot_m}_{m \in \ufs{n}} \in \C^n
\end{aligned}
\right.\]
With the notation above, the following diagram in $\C$ commutes by the functoriality of $\fgr$.
\begin{equation}\label{laxfhy_assoc}
\begin{tikzpicture}[vcenter]
\def\t{15} \def\h{3} \def\v{1}
\draw[0cell=.9]
(0,0) node (a) {\gaC_n\big(x \sscs \fgr \hgr_i \ygr_{i,j} \zdot_\crdot \big)}
(a)++(0,-1.5) node (b) {\fgr \gaBi_n\big(x \sscs \hgr_i \ygr_{i,j} \zdot_\crdot \big)}
(b)++(\h,-\v) node (c) {\fgr \hgr_i \gaAij_n \big(x \sscs \ygr_{i,j} \zdot_\crdot \big)}
(c)++(\h,\v) node (d) {\fgr \hgr_i \ygr_{i,j} \gaZijs_n \big(x \sscs \zdot_\crdot\big)}
;
\draw[1cell=.9]
(a) edge node[swap] {\laxf_{n,i}} (b)
(b) edge[bend right=\t] node[swap,pos=.2] {\fgr \laxhi_{n,j}} (c)
(c) edge[bend right=\t] node[swap,pos=.7] {\fgr \hgr_i \laxyij_{n,s}} (d)
(b) edge node {\fgr \big( \hgr_i \laxyij_{n,s} \circ \laxhi_{n,j} \big)} (d)
;
\end{tikzpicture}
\end{equation}
By \cref{laxfh_component}, the bottom composite in \cref{laxfhy_assoc} is a component of the $(n,\ell)$-action constraint of the left-hand side of \cref{gam_fhy}.  The top composite in \cref{laxfhy_assoc} is the same component of the $(n,\ell)$-action constraint of the right-hand side of \cref{gam_fhy}.  This proves the associativity axiom \cref{enr-multicategory-associativity} for $\MultpsO$.

\parhead{Top equivariance}.  Using the notation in \cref{def:gam_Omorphisms}, the object part of the top equivariance axiom \cref{enr-operadic-eq-1} requires, for each permutation $\sigma \in \Sigma_k$, the equality 
\begin{equation}\label{multpso_topeq}
\gam\big(f\sigma \sscs \ang{h_{\sigma(i)}}_{i \in \ufs{k}} \big) 
= \gam\big(f \sscs \ang{h_i}_{i \in \ufs{k}} \big) \sigmabar 
\cn \bang{\ang{\A_{\sigma(i),j}}_{j \in \ufs{d}_{\sigma(i)}} }_{i \in \ufs{k}} \to \C
\end{equation}
of $d$-lax $\Op$-morphisms.  In \cref{multpso_topeq},
\[\sigmabar = \sigma\bang{d_{\sigma(1)}, \ldots, d_{\sigma(k)}} \in \Sigma_d\]
is the block permutation induced by $\sigma$ that permutes blocks of lengths $d_{\sigma(1)}, \ldots, d_{\sigma(k)}$.  To see that the two sides of \cref{multpso_topeq} have the same underlying functor, we consider the following diagram.
\begin{equation}\label{multpso_topeq_functor}
\begin{tikzpicture}[vcenter]
\def\v{-2}
\draw[0cell=.9]
(0,0) node (a1) {\prod_{i=1}^k \prod_{j=1}^{d_{\sigma(i)}} \A_{\sigma(i),j}}
(a1)++(4,0) node (a2) {\prod_{i=1}^k \B_{\sigma(i)}}
(a2)++(2.5,0) node (a3) {\C}
(a1)++(0,\v) node (b1) {\prod_{i=1}^k \prod_{j=1}^{d_i} \Aij}
(a2)++(0,\v) node (b2) {\prod_{i=1}^k \B_i}
(a3)++(0,\v) node (b3) {\C}
;
\draw[1cell=.9]
(a1) edge node {\txprod_i\, h_{\sigma(i)}} (a2)
(a2) edge node {f\sigma} (a3)
(b1) edge node {\txprod_i\, h_i} (b2)
(b2) edge node {f} (b3)
(a1) edge node[swap] {\sigmabar} (b1)
(a2) edge node {\sigma} (b2)
(a3) edge[equal] (b3)
;
\end{tikzpicture}
\end{equation}
In this diagram, the left and right regions commute by, respectively, the naturality of the braiding for the Cartesian product and the definition of $f\sigma$ in \cref{klax_sigma_functor}.  The composites along the top and left-bottom boundaries are the underlying functors of, respectively, the left-hand side and the right-hand side of \cref{multpso_topeq}.  Thus, the two sides of \cref{multpso_topeq} have the same underlying functor.  Moreover, by \cref{Otr_sigma,gam_thetas}, the diagram \cref{multpso_topeq_functor} also proves the morphism part of the top equivariance axiom.

To verify that the two sides of \cref{multpso_topeq} have the same action constraints, suppose 
\[\ell = \big(\txsum_{p=1}^{\sigmainv(i) - 1} d_{\sigma(p)} \big) + j \qquad  \text{for some $(i,j) \in \ufs{k} \times \ufs{d}_i$}.\]  
For $n \geq 0$, using \cref{laxfsi_whiskering,laxfh_component}, the following equalities prove that the two sides of \cref{multpso_topeq} have the same $(n,\ell)$-action constraint.  
\[\begin{split}
\laxof{{\gam(f\sigma \sscs \ang{h_{\sigma(i)}})}}_{n,\ell} 
&= (f\sigma) \big(h_{\sigma(1)} , \ldots , \laxhi_{n,j}, \ldots, h_{\sigma(k)} \big) \circ \laxfsi_{n,\sigmainv(i)}\\
&= f\big(h_1, \ldots , \laxhi_{n,j}, \ldots, h_k\big) \circ \big( \laxf_{n,i} * \sigma\big)\\
&= \big( f\big(h_1, \ldots , \laxhi_{n,j}, \ldots, h_k\big) \circ \laxf_{n,i} \big) * \sigmabar\\
&= \laxof{{ \gam(f \sscs \ang{h_i})}}_{n,\, d_1 + \cdots + d_{i-1} + j} * \sigmabar\\ 
&= \laxof{{ \gam(f \sscs \ang{h_i})}}_{n,\, \sigmabar(\ell)} * \sigmabar\\
&= \laxof{{ \gam(f \sscs \ang{h_i}) \sigmabar}}_{n,\ell} 
\end{split}\]
This finishes the proof of the top equivariance axiom \cref{enr-operadic-eq-1} for $\MultpsO$.

\parhead{Bottom equivariance}.  
The object part of the bottom equivariance axiom \cref{enr-operadic-eq-2} requires, for permutations $\tau_i \in \Sigma_{d_i}$ for $i \in \ufs{k}$, the equality
\begin{equation}\label{multpso_boteq}
\gam\big(f \sscs \ang{h_i \tau_i}_{i \in \ufs{k}} \big) 
= \gam\big(f \sscs \ang{h_i}_{i \in \ufs{k}}\big) \tautimes 
\cn \bang{ \ang{\A_{i,\tau_i(j)}}_{j \in \ufs{d}_i}}_{i \in \ufs{k}} \to \C
\end{equation}
of $d$-lax $\Op$-morphisms.  In \cref{multpso_boteq},
\[\tautimes = \tau_1 \times \cdots \times \tau_k \in \Sigma_d\]
is the block sum defined in \cref{blocksum}.  To see that the two sides of \cref{multpso_boteq} have the same underlying functor, we consider the following commutative diagram.
\begin{equation}\label{multpso_boteq_functor}
\begin{tikzpicture}[vcenter]
\def\t{15} \def\v{1.3}
\draw[0cell=.9]
(0,0) node (a) {\prod_{i=1}^k \prod_{j=1}^{d_i} \A_{i,\tau_i(j)}}
(a)++(2.5,-\v) node (b) {\prod_{i=1}^k \prod_{j=1}^{d_i} \Aij}
(b)++(2,\v) node (c) {\prod_{i=1}^k \B_i}
(c)++(2,0) node (d) {\C}
;
\draw[1cell=.9]
(a) edge[bend right=\t,shorten <=-1ex] node[swap,pos=.3] {\tautimes} (b)
(b) edge[bend right=\t] node[swap] {\txprod_i\, h_i} (c)
(c) edge node {f} (d)
(a) edge node {\txprod_i\, h_i\tau_i} (c)
;
\end{tikzpicture}
\end{equation}
The top and bottom composites in \cref{multpso_boteq_functor} are the underlying functors of, respectively, the left-hand side and right-hand side of \cref{multpso_boteq}.  Moreover, by \cref{Otr_sigma,gam_thetas}, the diagram \cref{multpso_boteq_functor} also proves the morphism part of the bottom equivariance axiom.

To verify that the two sides of \cref{multpso_boteq} have the same action constraints, suppose 
\[\ell = \big(\txsum_{p=1}^{i - 1} d_p \big) + \tauinv_i(j) \qquad  \text{for some $(i,j) \in \ufs{k}\times \ufs{d}_i$}.\]  
For $n \geq 0$, using \cref{laxfsi_whiskering,laxfh_component}, the following equalities prove that the two sides of \cref{multpso_boteq} have the same $(n,\ell)$-action constraint.  
\[\begin{split}
\laxof{{\gam(f \sscs \ang{h_i \tau_i})}}_{n,\ell} 
&= f\big(h_1 \tau_1, \ldots, \laxof{{h_i\tau_i}}_{n,\tauinv_i(j)}, \ldots, h_k \tau_k \big) \circ \laxf_{n,i}\\
&= f\big(h_1 \tau_1, \ldots, \laxhi_{n,j} * \tau_i, \ldots, h_k \tau_k \big) \circ \laxf_{n,i}\\
&\overset{\blacktriangle}{=} \big( f\big(h_1, \ldots , \laxhi_{n,j}, \ldots, h_k\big) \circ \laxf_{n,i} \big) * \tautimes\\
&= \laxof{{\gam(f \sscs \ang{h_i})}}_{n,\, d_1 + \cdots + d_{i-1} + j} * \tautimes\\
&= \laxof{{\gam(f \sscs \ang{h_i})}}_{n,\, \tautimes(\ell)} * \tautimes\\
&= \laxof{{\gam(f \sscs \ang{h_i}) \tautimes}}_{n,\ell}
\end{split}\]
The equality labeled $\blacktriangle$ above uses the following facts, in which the objects are defined in \cref{laxfh_objects}.
\[\begin{split}
\tau_p \ang{a_{p,\tau_p(r)}}_{r \in \ufs{d}_p} 
&= \ang{a_{p,r}}_{r \in \ufs{d}_p} \forspace p \in \ufs{k} \setminus \{i\} \\
\tau_i \big(a_{i,\tau_i(1)}, \ldots, a_{i,j,m}, \ldots, a_{i,\tau_i(d_i)} \big)
&= \big(a_{i,1}, \ldots, a_{i,j,m}, \ldots, a_{i,d_i} \big)
\end{split}\]
This proves the bottom equivariance axiom \cref{enr-operadic-eq-2} for $\MultpsO$.
\end{proof}

\subsection*{$\Cat$-Multicategories of Pseudoalgebras}
Taking $G$-fixed subcategories and restricting $G$-functors to $G$-fixed subcategories define a strict symmetric monoidal functor
\begin{equation}\label{Gfixed}
\Gcat \fto{\gfixed} \Cat.
\end{equation}
Thus, applying $\gfixed$ to the multimorphism $G$-categories of a $\Gcat$-multicategory yields a $\Cat$-multicategory.  Applying this procedure to the $\Gcat$-multicategory $\MultvO$ in \cref{thm:multpso}, along with \cref{expl:k_laxmor_g,expl:O_tr_g}, yields the following result.

\begin{theorem}\label{thm:multpso_gfixed}
For each group $G$ and each pseudo-commutative operad $\Op$ in $\Gcat$, there is a $\Cat$-multicategory
\[\MultpsO^G\]
defined by the following data.
\begin{itemize}
\item Its objects are $\Op$-pseudoalgebras (\cref{def:pseudoalgebra}).
\item For $k \geq 0$, each $k$-ary multimorphism category has the form
\[\MultpsO\scmap{\ang{\A_i}_{i \in \ufs{k}} ; \B}^G,\]
where the $G$-category $\MultpsO\scmap{\ang{\A_i}_{i \in \ufs{k}} ; \B}$ is defined in \cref{def:MultpsO_karycat}.
\begin{itemize}
\item Its objects are $G$-equivariant $k$-lax $\Op$-morphisms (\cref{def:k_laxmorphism}). 
\item Its morphisms are $G$-equivariant $k$-ary $\Op$-transformations (\cref{def:kary_transformation}).
\end{itemize}
\item The unit of an $\Op$-pseudoalgebra $\B$ is given by the identity $\Op$-morphism $1 \cn \B \to \B$.
\item The symmetric group action and the composition are the restrictions of the corresponding structures for $\MultpsO$ (\cref{def:multicatO}) to $G$-fixed subcategories.
\end{itemize}
Moreover, there are pseudo and strict variant $\Cat$-multicategories, denoted
\[\MultpspsO^G \andspace \MultstO^G,\]
that replace the $k$-ary 1-cells by, respectively, $G$-equivariant $k$-ary $\Op$-pseudomorphisms and $G$-equivariant $k$-ary strict $\Op$-morphisms.  
\end{theorem}

\begin{example}[Underlying 2-Categories]\label{ex:underlying_iicat}
By \cref{ex:unarycategory}, the $\Cat$-multicategories in \cref{thm:multpso_gfixed},
\[\MultstO^G \bigsubset \MultpspsO^G \bigsubset \MultpsO^G,\]
have underlying 2-categories
\begin{equation}\label{underlying_algvo}
\AlgstO \bigsubset \AlgpspsO \bigsubset \AlglaxO,
\end{equation}
as stated in \cref{oalgps_twocat}.  Note that $\Op$ is assumed to have a pseudo-commutative structure in \cref{thm:multpso,thm:multpso_gfixed}, but not in \cref{oalgps_twocat}.  In other words, the 2-categories in \cref{underlying_algvo} exist without assuming any pseudo-commutative structure on $\Op$.  The reason is that, among \cref{def:pseudoalgebra,def:k_laxmorphism,def:kary_transformation}, the pseudo-commutative structure on $\Op$ is only used in the commutativity axiom \cref{laxf_com}, which only happens in arity $k \geq 2$.  On the other hand, the underlying 2-category $\AlgpsvO$ only uses the 1-ary 1-cells and 1-ary 2-cells of $\MultvO^G$.  Thus, \cref{oalgps_twocat} is valid as stated.
\end{example}

\subsection*{$\Gcat$-Multicategories of Symmetric Monoidal $G$-Categories}

Next, we apply \cref{thm:multpso,thm:multpso_gfixed} to the Barratt-Eccles operad $\BE$ and the $G$-Barratt-Eccles operad $\GBE$.  See the introduction of \cref{ch:psalg} for a summary of their algebras and pseudoalgebras.

\begin{example}[Naive Variants]\label{ex:multBE}
By \cref{BE_pseudocom}, the Barratt-Eccles operad $\BE$\index{Barratt-Eccles operad}\index{pseudoalgebra}\index{naive symmetric monoidal G-category@naive symmetric monoidal $G$-category} in $\Gcat$ (\cref{def:BE-Gcat}) has a unique pseudo-commutative structure.  By \cref{thm:multpso}, there are $\Gcat$-multicategories
\begin{equation}\label{MultBE_gcat}
\MultstBE \bigsubset \MultpspsBE \bigsubset \MultpsBE,
\end{equation}
each with $\BE$-pseudoalgebras as objects and $k$-ary $\BE$-transformations as $k$-ary 2-cells.  The subscript $\va$ specifies the $k$-ary 1-cells:
\begin{itemize}
\item $k$-lax $\BE$-morphisms for $\va = \sflax$,
\item $k$-ary $\BE$-pseudomorphisms for $\va = \sfps$, and
\item $k$-ary strict $\BE$-morphisms for $\va = \sfst$.  
\end{itemize}

By \cref{thm:multpso_gfixed}, there are $\Cat$-multicategories
\[\MultstBE^G \bigsubset \MultpspsBE^G \bigsubset \MultpsBE^G,\]
each with $\BE$-pseudoalgebras as objects.  In each case, the $k$-ary 1-cells and $k$-ary 2-cells are the $G$-equivariant versions of those in \cref{MultBE_gcat}.  Moreover, by \cref{ex:underlying_iicat}, these $\Cat$-multicategories yield the 2-categories
\[\AlgstBE \bigsubset \AlgpspsBE \bigsubset \AlglaxBE,\]
which also have $\BE$-pseudoalgebras as objects.  By \cref{thm:BEpseudoalg}, these 2-categories are 2-equivalent via $\Phi$ to, respectively, 
\[\smgcatst \bigsubset \smgcatsg \bigsubset \smgcat\]
in \cref{def:smGcat_twocat}.  The objects in the latter three 2-categories are \emph{naive} symmetric monoidal $G$-categories (\cref{def:naive_smGcat}).  Thus, using the 2-equivalences $\Phi$, we may regard each $\Gcat$-multicategory $\MultvBE$, with $\va \in \{\sflax,\sfps,\sfst\}$, as having naive symmetric monoidal $G$-categories as objects.

\cref{expl:k_laxmorphism} compares the $\Gcat$-multicategories $\MultvBE$ and the multicategories in \cite{gmmo23}.  In particular, the objects of the multicategory $\mathbf{Mult}(\BE)$ in \cite{gmmo23} are $\BE$-algebras in $\Gcat$, which are, by \cref{def:BE-Gcat}, naive \emph{permutative} $G$-categories.  See \cref{expl:naive_perm_Gcat} for an unraveled description of naive permutative $G$-categories.
\end{example}

\begin{example}[Genuine Variants]\label{ex:multGBE}
The $G$-Barratt-Eccles operad $\GBE$\index{G-Barratt-Eccles operad@$G$-Barratt-Eccles operad}\index{pseudoalgebra}\index{genuine symmetric monoidal G-category@genuine symmetric monoidal $G$-category} in \cref{def:GBE} has a unique pseudo-commutative structure by \cref{GBE_pseudocom}.  By \cref{thm:multpso,thm:multpso_gfixed,ex:underlying_iicat}, for each of the three variants $\va \in \{\sflax,\sfps,\sfst\}$, there are
\begin{itemize}
\item a $\Gcat$-multicategory $\MultvGBE$,
\item a $\Cat$-multicategory $\MultvGBE^G$, and
\item a 2-category $\Algpsv(\GBE)$.
\end{itemize}
In each case, the objects are $\GBE$-pseudoalgebras, which are, by \cref{def:GBE_pseudoalg}, \emph{genuine} symmetric monoidal $G$-categories.  In contrast, the objects of the multicategory $\mathbf{Mult}(\GBE)$ in \cite{gmmo23} are $\GBE$-algebras in $\Gcat$, which are, by \cref{def:GBE_algebra}, genuine \emph{permutative} $G$-categories.  See \cref{expl:GBE_algebra} for an unraveled description of genuine permutative $G$-categories. 
\end{example}

\chapter{Symmetric Monoidal Closed Category of $\Gskg$-Categories}
\label{ch:ggcat}
For an arbitrary group $G$, this chapter constructs a symmetric monoidal closed category $\GGCat$ whose objects are pointed functors
\[(\Gsk,\vstar) \fto{f} (\Gcatst,\boldone),\]
called $\Gskg$-categories.  The $\Gcat$-multicategory associated to $\GGCat$ serves as an intermediate enriched multicategory for our $G$-equivariant algebraic $K$-theory multifunctor 
\[\MultpsO \fto{\Kgo} \GSp\] 
from $\Op$-pseudoalgebras to orthogonal $G$-spectra.  The first half of this chapter reviews in detail the indexing category $\Gsk$, which originated in the work of Elmendorf-Mandell \cite{elmendorf-mandell} on nonequivariant multifunctorial $K$-theory.

\summary
The following table summarizes corresponding concepts for the $\Gcat$-multicategories $\MultpsO$ and $\GGCat$.
\begin{center}
\resizebox{.9\columnwidth}{!}{%
{\renewcommand{\arraystretch}{1.3}%
{\setlength{\tabcolsep}{1ex}
\begin{tabular}{c|cr|cr}
& $\MultpsO$ & \eqref{thm:multpso} & $\GGCat$ & \eqref{expl:ggcat_gcatenr} \\ \hline
objects & $\Op$-pseudoalgebras & \eqref{def:pseudoalgebra} & $\Gskg$-categories & \eqref{expl:ggcat_obj} \\
$k$-ary 1-cells & $k$-lax $\Op$-morphisms & \eqref{def:k_laxmorphism} & \multicolumn{2}{c}{\eqref{ggcat_zero_obj_comp}, \eqref{ggcat_k1}} \\
$k$-ary 2-cells & $k$-ary $\Op$-transformations & \eqref{def:kary_transformation} & \multicolumn{2}{c}{\eqref{ggcat_zero_mor_comp}, \eqref{ggcat_k2}} \\
$G$-action & \multicolumn{2}{c|}{\eqref{def:k_laxmorphism_g}, \eqref{expl:O_tr_g}} & \multicolumn{2}{c}{\eqref{ggcat_zero_mor_g}, \eqref{ggcat_k_g}} \\
symmetry & \multicolumn{2}{c|}{\eqref{multpso_sym_functor}} & \multicolumn{2}{c}{\eqref{expl:ggcat_symmetry}} \\
composition & \multicolumn{2}{c|}{\eqref{def:gam_functor}} & \multicolumn{2}{c}{\eqref{expl:ggcat_composition}} \\
symmetric monoidal & \multicolumn{2}{c|}{no} & yes & \eqref{def:GGCat} \\
\end{tabular}}}}
\end{center}

\connection
The first step of our $G$-equivariant algebraic $K$-theory multifunctor $\Kgo$ is a $\Gcat$-multifunctor (\cref{thm:Jgo_multifunctor})
\[\MultpsO \fto{\Jgo} \GGCat,\]
with $\MultpsO$ the $\Gcat$-multicategory in \cref{thm:multpso}.  In \cref{ch:ggtop,ch:semg}, the symmetric monoidal closed category $\GGCat$ is the domain of the second part of $\Kgo$:
\[\GGCat \fto{\clast} \GGTop \fto{\Kg} \GSp.\]
This is a composite of two symmetric monoidal functors in the $\Gtop$-enriched sense.

\organization
This chapter consists of the following sections.

\secname{sec:Fsk}  This section reviews the permutative category $\Fsk$ of pointed finite sets and pointed functions, along with the reindexing functors between its smash powers.

\secname{sec:Gsk}  This section is a self-contained review of the permutative category $\Gsk$.  It is used as the indexing category of $\GGCat$ and $\GGTop$.

\secname{sec:ggcat_sm}  This section constructs the symmetric monoidal closed category $\GGCat$.  Its objects, called $\Gskg$-categories, are pointed functors
\[(\Gsk,\vstar) \fto{f} (\Gcatst,\boldone),\]
and its morphisms are pointed natural transformations.  \cref{expl:ggcat_unit,expl:ggcat_smag,expl:ggcat_brkst} unravel the symmetric monoidal closed structure on $\GGCat$.  Its internal hom involves the internal hom for $\Gcatst$, denoted $\Catgst$ \cref{catgst_cd}.  The latter consists of not-necessarily $G$-equivariant pointed functors and pointed natural transformations, with $G$ acting by conjugation.

\secname{sec:ggcat_multicat}  This section describes in detail the $\Gcat$-multicategory associated to the symmetric monoidal closed category $\GGCat$.  Along the way, this section provides an independent proof that $\GGCat$ is a $\Gcat$-multicategory.  The detailed description of the $\Gcat$-multicategory $\GGCat$ in this section is used in \cref{ch:jemg} to construct the $\Gcat$-multifunctor $\Jgo$.

\section{The Indexing Category $\Fsk$}
\label{sec:Fsk}

To prepare for the discussion of the indexing category $\Gsk$ in \cref{sec:Gsk}, this section reviews the permutative category $\Fsk$ of pointed finite sets and pointed functions, and reindexing functors between the smash powers of $\Fsk$.

\secoutline
\begin{itemize}
\item \cref{def:Fsk} defines the small pointed category $\Fsk$ of pointed finite sets and pointed functions.
\item \cref{def:Fsk_permutative,Fsk_permutative} specify a permutative category structure on $\Fsk$.
\item \cref{def:Fsk_smashpower,def:injections} define smash powers of $\Fsk$ and reindexing functors between them.
\end{itemize}

\subsection*{The Category $\Fsk$}

\begin{definition}\label{def:Fsk}\
\begin{enumerate}
\item A \emph{pointed set}\index{pointed set} $(S, \bp)$ is a set $S$ equipped with a distinguished element $\bp$, called the \emph{basepoint}.  A \emph{pointed finite set}\index{pointed finite set} is a pointed set whose underlying set is finite.
\item The pointed finite set
\begin{equation}\label{ordn}
\ord{n} = \{0 < 1 < \cdots < n\}
\end{equation}
is equipped with the basepoint $0$ and its natural ordering.  Unless otherwise specified, from now on, a \emph{pointed finite set} means one of the form $\ordn$.
\item A \emph{pointed function}\index{pointed function} between pointed sets is a basepoint-preserving function.
\item The small pointed category $\Fsk$ has objects $\ordn$ for $n \geq 0$, pointed functions as morphisms, and basepoint $\ord{0}$.\defmark
\end{enumerate}
\end{definition}

\begin{remark}\label{rk:Fsk}
The category $\Fsk$ is the opposite of the category $\Ga$ in \cite{segal}.  The object $\ord{0}$ is initial and terminal in $\Fsk$.
\end{remark}

Recall that a permutative category is a strict symmetric monoidal category (\cref{def:symmoncat}).

\begin{definition}\label{def:Fsk_permutative}
We define a permutative category structure $(\sma, \ord{1}, \xi)$ on the small pointed category $\Fsk$ as follows.
\begin{description}
\item[Monoidal product on objects] The monoidal product is the functor
\[\Fsk \times \Fsk \fto{\sma} \Fsk\]
with object assignment
\begin{equation}\label{m-sma-n}
\ord{m} \sma \ord{n} = \ord{mn}.
\end{equation}
In \cref{m-sma-n}, $\ord{mn}$ is identified with the smash product $\ord{m} \sma \ord{n}$ using the lexicographic ordering in \cref{lex_bijection} 
\[\ufs{mn} \fto[\iso]{\la_{m,n}} \ufs{m} \times \ufs{n}\]
away from the basepoint $0 \in \ord{mn}$.
\item[Monoidal product on morphisms] 
We extend $\sma$ to morphisms in the same way using the lexicographic ordering.  More explicitly, for pointed functions
\begin{equation}\label{fmp-hnq}
\ord{m} \fto{f} \ord{p} \andspace \ord{n} \fto{h} \ord{q},
\end{equation}
the pointed function
\[\ord{m} \sma \ord{n} \fto{f \sma h} \ord{p} \sma \ord{q}\]
is defined by
\begin{equation}\label{fsmah-ab}
\begin{split}
(f \sma h)(0) &= 0 \andspace\\
(f \sma h)(a,b) &= \big(f(a), h(b)\big)
\end{split}
\end{equation}
for $(a,b) \in \ufs{m} \times \ufs{n}$. 
\item[Monoidal unit] The monoidal unit is the object $\ord{1} = \{0 < 1\} \in \Fsk$.
\item[Braiding] The component of the braiding $\xi$ at a pair of objects $(\ord{m}, \ord{n}) \in \Fsk^2$ is the pointed bijection
\begin{equation}\label{Fsk_braiding}
\ord{m} \sma \ord{n} \fto[\iso]{\xi_{\ord{m}, \ord{n}}} \ord{n} \sma \ord{m}
\end{equation}
given by the $(m,n)$-transpose permutation $\twist_{m,n}$ in \cref{eq:transpose_perm} away from the basepoint $0 \in \ord{mn}$.  More explicitly, using the lexicographic ordering, we have
\[\xi_{\ord{m}, \ord{n}} (a,b) = (b,a).\]
\end{description}
This finishes the definition of $(\sma,\ord{1},\xi)$. \cref{Fsk_permutative} verifies that this triple is a  permutative category structure on $\Fsk$.
\end{definition}

\begin{lemma}\label{Fsk_permutative}
In \cref{def:Fsk_permutative}, the quadruple $(\Fsk,\sma,\ord{1},\xi)$ is a permutative category.
\end{lemma}

\begin{proof}
The assignments \cref{m-sma-n,fsmah-ab} of $\sma$ on, respectively, objects and morphisms define a functor because $f$ and $h$ in \cref{fsmah-ab} are applied entrywise.  Moreover, by those definitions, $\sma$ is associative and unital with respect to the object $\ord{1} \in \Fsk$.  Thus, $(\Fsk,\sma,\ord{1})$ is a strict monoidal category (\cref{def:monoidalcategory}) once we define the associativity and unit isomorphisms to be identities.   

The braiding $\xi$ defined in \cref{Fsk_braiding} is a natural isomorphism because
\[\begin{split}
(h \sma f) \xi_{\ord{m}, \ord{n}}  (a,b) 
&= \big(h(b), f(a)\big)\\ 
&= \xi_{\ord{m}, \ord{n}} (f \sma h)(a,b)
\end{split}\]
for $(a,b) \in \ufs{m} \times \ufs{n}$.  Moreover, the braiding $\xi$ satisfies the symmetry and hexagon axioms in \cref{symmoncatsymhexagon} because transpose permutations satisfy those axioms.  Thus,  the quadruple $\big(\Fsk,\sma,\ord{1},\xi\big)$ is a permutative category.
\end{proof}

\subsection*{Smash Powers of $\Fsk$}

To define the indexing category $\Gsk$ in \cref{sec:Gsk}, we need the following smash powers of the pointed category $\Fsk$.

\begin{definition}\label{def:Fsk_smashpower}
For $q \geq 0$, we define the small pointed category $\Fsk^{(q)}$ as follows.
\begin{description}
\item[$q=0$] The pointed category\label{not:Fsk_smash_zero}
\[\Fsk^{(0)} = \big\{\vstar \rightleftarrows \ang{}\big\}\]
has two objects: the empty tuple $\ang{}$ and the initial-terminal basepoint $\vstar$.  The only non-identity morphisms are the unique morphisms $\vstar \to \ang{}$ and $\ang{} \to \vstar$.
\item[$q>0$] In this case, we define the small pointed category
\[\Fsk^{(q)} = \Fsk^{\sma q}\]
as the $q$-fold smash power of $\Fsk$, as defined in \cref{eq:smash} for $\C = (\Cat, \times, \boldone)$. 
\begin{description}
\item[Objects]  We denote the basepoint of $\Fsk^{(q)}$ by $\vstar$, which is both initial and terminal.  A typical object in $\Fsk^{(q)}$ is denoted by
\begin{equation}\label{angordn}
\ang{\ord{n}} = \ang{\ord{n}_i}_{i \in \ufs{q}} = (\ord{n}_1, \ldots, \ord{n}_q)
\end{equation}
with each $\ord{n}_i$ a pointed finite set \cref{ordn}.  If any $\ord{n}_i = \ord{0}$, then $\ang{\ord{n}} = \vstar$. 
\item[Morphisms] A typical morphism in $\Fsk^{(q)}$ is denoted by
\begin{equation}\label{angpsi}
\ang{\psi} = \ang{\psi_i}_{i \in \ufs{q}} \cn \ang{\ord{m}_i}_{i \in \ufs{q}} \to \ang{\ord{n}_i}_{i \in \ufs{q}}
\end{equation}
with each
\[\ord{m}_i \fto{\psi_i} \ord{n}_i\]
a morphism in $\Fsk$, which means a pointed function.  If any $\psi_i$ factors through $\ord{0}$, then $\ang{\psi}$ is the \index{0-morphism}\emph{0-morphism}, which means that it factors through the basepoint $\vstar$.  A morphism that is not the 0-morphism is called a \index{nonzero morphism}\emph{nonzero morphism}.
\end{description}
\end{description}  
This finishes the definition of the small pointed category $\Fsk^{(q)}$.  A non-basepoint object or a nonzero morphism in $\Fsk^{(q)}$ is said to have \emph{length}\index{length} $q$.
\end{definition}

\subsection*{Reindexing}
There are morphisms in the category $\Gsk$ that are not of the form $\ang{\psi}$ as defined in \cref{angpsi}.  To define a general morphism in $\Gsk$, we use the unpointed finite sets $\ufs{n} = \{1, 2,\ldots , n\}$ defined in \cref{ufsn} and the reindexing device in the following definition.

\begin{definition}\label{def:injections}
We define the category $\Inj$ with
\begin{itemize}
\item objects $\ufs{n}$ for $n \geq 0$ and 
\item injections as morphisms.
\end{itemize}  
Given an injection $h \cn \ufs{q} \to \ufs{r}$, we define a pointed functor  
\[\Fsk^{(q)} \fto{h_*} \Fsk^{(r)}\]
as follows.
\begin{description}
\item[$q = r = 0$] 
In this case, $h$ is the identity function on $\ufs{0} = \emptyset$, and $h_*$ is defined as the identity functor on $\Fsk^{(0)} = \{\vstar \rightleftarrows \ang{}\}$.
\item[$q = 0 < r$] 
In this case, the pointed functor $h_*$ is determined by the following object assignment.
\begin{equation}\label{hangempty}
\left\{\begin{split}
h_* \vstar &= \vstar\\
h_*\ang{} &= \ang{\ord{1}}_{j \in \ufs{r}} = \left(\ord{1}, \ldots, \ord{1}\right) \in \Fsk^{(r)}
\end{split}\right.
\end{equation}
\item[$q>0$] 
Given an object $\ang{\ord{n}_i}_{i \in \ufs{q}}$ \cref{angordn} and a morphism $\ang{\psi_i}_{i \in \ufs{q}}$ \cref{angpsi} in $\Fsk^{(q)}$, we define the object and morphism
\begin{equation}\label{reindexing_functor}
\begin{split}
h_*\ang{\ord{n}_i}_{i \in \ufs{q}} &= \ang{\ord{n}_{\hinv(j)}}_{j \in \ufs{r}} \andspace\\
h_*\ang{\psi_i}_{i \in \ufs{q}} &= \ang{\psi_{\hinv(j)}}_{j \in \ufs{r}}
\end{split}
\end{equation}
in $\Fsk^{(r)}$, where, if $\hinv(j) = \emptyset$, then
\begin{equation}\label{ordn_empty}
\begin{split}
\ord{n}_{\emptyset} &= \ord{1} \andspace\\ 
\psi_{\emptyset} &= 1_{\ord{1}} \cn \ord{1} \to \ord{1}.
\end{split}
\end{equation}
\end{description}
We call $h$ a \emph{reindexing injection}\index{reindexing injection} and $h_*$ a \index{reindexing functor}\emph{reindexing functor}.
\end{definition}

\begin{explanation}\label{expl:reindexing}
The reindexing functor $h_*$ in \cref{def:injections} is well defined for the following reasons.
\begin{itemize}
\item Suppose $\ord{n}_i = \ord{0}$ for some $i \in \ufs{q}$. Since $h \cn \ufs{q} \to \ufs{r}$ is an injection, at least one entry in $h_*\ang{\ord{n}_i}_{i \in \ufs{q}}$ is $\ord{0}$.  Thus, we have that
\[\ang{\ord{n}_i}_{i \in \ufs{q}} = \vstar \in \Fsk^{(q)} \impliespace 
h_*\ang{\ord{n}_i}_{i \in \ufs{q}} = \vstar \in \Fsk^{(r)}.\]
\item Similarly, if some $\psi_i$ factors through $\ord{0}$, then at least one entry in $h_*\ang{\psi_i}_{i \in \ufs{q}}$ factors through $\ord{0}$.  Thus, if $\ang{\psi}$ is the 0-morphism, then so is $h_*\ang{\psi_i}_{i \in \ufs{q}}$.
\end{itemize}  
Furthermore, the assignment $h \mapsto h_*$ is functorial in the following sense.
\begin{itemize}
\item If $h \cn \ufs{q} \to \ufs{q}$ is the identity function, then $h_*$ is the identity functor on $\Fsk^{(q)}$.
\item For another injection $f \cn \ufs{p} \to \ufs{q}$, there is an equality of functors
\begin{equation}\label{hf_reindexing}
(hf)_* = h_* \circ f_* \cn \Fsk^{(p)} \to \Fsk^{(r)}.
\end{equation}
\end{itemize}
These properties follow from \cref{hangempty,reindexing_functor,ordn_empty}: $h_*$ permutes the entries according to the reindexing injection $h$, and inserts $\ord{1}$ or $1_{\ord{1}}$ for entries not in the image of $h$.
\end{explanation}

\section{The Indexing Category $\Gsk$}
\label{sec:Gsk}

This section reviews the permutative category $\Gsk$.  The indexing category $\Gsk$ is used in \cref{sec:ggcat_sm} to construct the symmetric monoidal closed category $\GGCat$.  The category $\Gsk$ originated in the work of Elmendorf-Mandell \cite{elmendorf-mandell}.  The reader is referred to \cite[Chapters 9--13]{cerberusIII} for a thorough discussion of $\Gsk$ in the context of nonequivariant multifunctorial $K$-theory.  

\secoutline
\begin{itemize}
\item \cref{def:Gsk} defines the small pointed category $\Gsk$, with further elaboration in \cref{expl:Gsk_composite}.
\item \cref{def:Gsk_permutative,Gsk_permutative} equip $\Gsk$ with the structure of a permutative category.
\item \cref{def:smashFskGsk,sma_symmon} construct a strict symmetric monoidal pointed functor
\[\big(\Gsk,\oplus,\ang{},\xi\big) \fto{\sma} \big(\Fsk,\sma,\ord{1},\xi\big).\]
\end{itemize}

\subsection*{The Category $\Gsk$}

Using the small pointed categories $\Fsk^{(q)}$ in \cref{def:Fsk_smashpower}, now we define the indexing category $\Gsk$.  See \cref{expl:Gsk} for its relationship with the original Elmendorf-Mandell category.

\begin{definition}\label{def:Gsk}
The small pointed category $\Gsk$ is defined as follows.
\begin{description}
\item[Objects] The object pointed set of $\Gsk$ is defined as the wedge
\begin{equation}\label{Gsk_objects}
\Ob(\Gsk) = \bigvee_{q \geq 0} \Ob(\Fsk^{(q)}),
\end{equation}
which identifies the basepoints $\vstar \in \Fsk^{(q)}$ for $q \geq 0$.  The identified object is the basepoint $\vstar \in \Gsk$, which is defined to be both initial and terminal in $\Gsk$.
\item[Morphisms]
Given an object $\angordm \in \Fsk^{(p)}$ and an object $\angordn \in \Fsk^{(q)}$, the morphism pointed set is defined as the following wedge.
\begin{equation}\label{Gsk_morphisms}
\begin{split}
\Gsk\big(\angordm, \angordn\big)
&= \bigvee_{f \in \Inj(\ufs{p},\, \ufs{q})}~ \Fsk^{(q)}\big(f_* \angordm, \angordn \big)\\
&= \bigvee_{f \in \Inj(\ufs{p},\, \ufs{q})}~ \bigwedge_{i \in \ufs{q}}~ \Fsk\big(\ord{m}_{\finv(i)} \scs \ord{n}_i \big)
\end{split}
\end{equation}
The basepoint of each morphism pointed set is the \index{0-morphism}\emph{0-morphism}, which is the unique morphism that factors through the initial-terminal basepoint $\vstar \in \Gsk$.  A morphism that is not the 0-morphism is called a \index{nonzero morphism}\emph{nonzero morphism}.  The set of nonzero morphisms $\angordm \to \angordn$ is denoted by $\Gskpunc(\angordm,\angordn)$.

In \cref{Gsk_morphisms}, for each reindexing injection $f \cn \ufs{p} \to \ufs{q}$, 
\[\Fsk^{(p)} \fto{f_*} \Fsk^{(q)}\]
is the reindexing functor defined in \cref{def:injections}. 
\begin{itemize}
\item If $q=0$ in \cref{Gsk_morphisms}, then $p=0$ and $f \cn \ufs{0} \to \ufs{0}$ is $1_\emptyset$.  In this case, the pointed set
\begin{equation}\label{Gsk_empty_mor}
\Gsk\big(\ang{}, \ang{}\big) = \Fsk^{(0)}\big(\ang{}, \ang{}\big)
\end{equation}
has two elements: the identity morphism of $\ang{}$ and the 0-morphism $\ang{} \to \vstar \to \ang{}$.
\item For $q>0$ in \cref{Gsk_morphisms}, a typical morphism in $\Gsk$ is usually denoted by 
\begin{equation}\label{fangpsi}
\angordm \fto{(f, \ang{\psi})} \angordn
\end{equation}
with
\begin{itemize}
\item a reindexing injection $f \cn \ufs{p} \to \ufs{q}$ and
\item a morphism 
\[\ang{\psi} = \ang{\psi_i}_{i \in \ufs{q}} \cn f_* \angordm \to \angordn
\inspace \Fsk^{(q)}.\]
\end{itemize}
If some pointed function
\[\ord{m}_{\finv(i)} \fto{\psi_i} \ord{n}_i\]
factors through $\ord{0} \in \Fsk$, then $(f, \ang{\psi})$ is the 0-morphism, factoring through the basepoint $\vstar$.  
\end{itemize}
When we do not need to specify the reindexing injection $f$ and the pointed functions $\psi_i$, we denote a typical morphism in $\Gsk$ by a generic symbol, such as $\upom$.
\item[Identities] The identity morphism of an object $\angordm$ of length $p \geq 0$ is 
\[\big(1_{\ufs{p}}, \ang{1_{\ordm_k}}_{k \in \ufs{p}}\big),\]
which consists of the identity function on $\ufs{p}$ and the identity function on $\ordm_k$ for each $k \in \ufs{p}$.
\item[Composition] Suppose we are given objects $\angordm \in \Fsk^{(p)}$, $\angordn \in \Fsk^{(q)}$, and $\angordl \in \Fsk^{(r)}$, and composable morphisms in $\Gsk$
\begin{equation}\label{Gsk_composable}
\angordm \fto{(f,\ang{\psi})} \angordn \fto{(h,\ang{\phi})} \angordl.
\end{equation}
\begin{itemize}
\item If $q=0$, then $p=0$, and $(f,\ang{\psi})$ is either the identity morphism $1_{\ang{}}$ or the 0-morphism.  In these cases, the composites with $(h,\ang{\phi})$ are, respectively, $(h,\ang{\phi})$ and the 0-morphism.
\item If $q>0$, then $r>0$, and the composite is defined as
\begin{equation}\label{Gsk_composite}
(h,\ang{\phi}) \circ (f,\ang{\psi}) = \big(hf, \ang{\phi} \circ h_*\ang{\psi} \big) \cn \angordm \to \angordl.
\end{equation}
\end{itemize}
\end{description}
This finishes the definition of $\Gsk$.  \cref{Gsk_Gcategory} proves that $\Gsk$ is a well-defined pointed category.
\end{definition}

\begin{explanation}[Composition in $\Gsk$]\label{expl:Gsk_composite}
On the right-hand side of \cref{Gsk_composite}, the first entry is the composite injection
\[\ufs{p} \fto{f} \ufs{q} \fto{h} \ufs{r}.\]
The second entry is the $r$-tuple
\[\begin{split}
&\ang{\phi} \circ h_*\ang{\psi} \\
&= \bang{\ord{m}_{(hf)^{-1}(j)} = \ordm_{\finv\hinv(j)} \fto{\psi_{\hinv(j)}} \ord{n}_{\hinv(j)} \fto{\phi_j} \ord{\ell}_j}_{j \in \ufs{r}}
\end{split}\]
of pointed functions between pointed finite sets.
\end{explanation}

\begin{lemma}\label{Gsk_Gcategory}
\cref{def:Gsk} defines a small pointed category $(\Gsk,\vstar)$.
\end{lemma}

\begin{proof}
The composition defined in \cref{Gsk_composite} is well defined because if some $\phi_j$ or some $\psi_i$ factors through $\ord{0}$, then at least one entry in $\ang{\phi} \circ h_*\ang{\psi}$ factors through $\ord{0}$.  Identity morphisms are two-sided units for composition \cref{Gsk_composite} because (i) $h_*$ is a functor and (ii) $1_*$ is an identity functor.  Composition is associative by \cref{hf_reindexing,expl:Gsk_composite}.
\end{proof}

\begin{explanation}\label{expl:Gsk}
The category $\Gsk$ in \cref{def:Gsk} agrees with $\Gsk$ in \cite[Def.\! 9.1.7]{cerberusIII}, but it is denoted by $\Gsk_*$ in \cite[page 190]{elmendorf-mandell}.  The category denoted by $\Gsk$ in \cite[Def.\! 5.1]{elmendorf-mandell} is different from our $\Gsk$ in the sense that the former does not have the initial-terminal basepoint $\vstar$ and the associated identification.
\end{explanation}

\subsection*{Permutative Structure}

Next, we define a permutative category structure (\cref{def:symmoncat}) on $\Gsk$ and compare it with the permutative category $\Fsk$ in \cref{def:Fsk_permutative}.

\begin{definition}\label{def:Gsk_permutative}
We define a permutative category structure $(\oplus, \ang{}, \xi)$ on the small pointed category $\Gsk$ as follows.
\begin{description}
\item[Monoidal product on objects] The monoidal product is the functor
\[\Gsk \times \Gsk \fto{\oplus} \Gsk\]
with object assignment given by concatenation\index{concatenation}
\begin{equation}\label{Gsk_oplus_obj}
\angordm \oplus \angordn = 
\big(\ord{m}_1, \ldots, \ord{m}_p, \ord{n}_1, \ldots, \ord{n}_q\big)
\end{equation}
for $\angordm \in \Fsk^{(p)}$ and $\angordn \in \Fsk^{(q)}$.  If any $\ord{m}_i$ or $\ord{n}_j$ is $\ord{0}$, then the right-hand side of \cref{Gsk_oplus_obj} has at least one entry of $\ord{0}$.  Thus, we have
\begin{equation}\label{Gsk_oplus_vstar}
\vstar \oplus \angordn = \vstar = \angordm \oplus \vstar,
\end{equation}
and $\vstar$ is a null object in the sense of \cref{definition:zero}.
\item[Monoidal product on morphisms]
By \cref{Gsk_oplus_vstar} and the fact that the basepoint $\vstar$ is initial and terminal in $\Gsk$, $\oplus$ is uniquely defined when one morphism has either domain or codomain $\vstar$.

To define $\oplus$ for other morphisms, we consider non-basepoint objects $\angordm \in \Fsk^{(p)}$, $\angordn \in \Fsk^{(q)}$, $\angordj \in \Fsk^{(r)}$, and $\angordl \in \Fsk^{(s)}$, and morphisms
\begin{equation}\label{fpsi_hphi}
\angordm \fto{(f,\ang{\psi})} \angordn \andspace \angordj \fto{(h,\ang{\phi})} \angordl
\end{equation}
in $\Gsk$.  Then we define the morphism
\begin{equation}\label{Gsk_oplus_morphism}
\begin{split}
&(f, \ang{\psi}) \oplus (h, \ang{\phi}) \\
&= \big(f \oplus h, \ang{\psi} \oplus \ang{\phi}\big) \cn
\angordm \oplus \angordj \to \angordn \oplus \angordl
\end{split}
\end{equation}
as follows.  The reindexing injection
\begin{equation}\label{f_oplus_h}
\ufs{p+r} \fto{f \oplus h} \ufs{q+s}
\end{equation}
is defined by
\[(f \oplus h)(i) = \begin{cases} f(i) & \text{if $1 \leq i \leq p$ and}\\
q + h(i-p) & \text{if $p+1 \leq i \leq p+r$}.
\end{cases}\]
The morphism $\ang{\psi} \oplus \ang{\phi}$ is the concatenation of the $q$-tuple $\ang{\psi}$ and the $s$-tuple $\ang{\phi}$, as displayed below.
\begin{equation}\label{psi_oplus_phi}
\begin{tikzpicture}[vcenter]
\draw[0cell=.9]
(0,0) node (a) {(f \oplus h)_* (\angordm \oplus \angordj)}
(a)++(0,-1) node (b) {f_* \angordm \oplus h_* \angordj}
(b)++(5,0) node (c) {\angordn \oplus \angordl}
(c)++(0,1) node (d) {\angordn \oplus \angordl}
;
\draw[1cell=.9]
(a) edge[equal] (b)
(c) edge[equal] (d)
(a) edge node {\ang{\psi} \oplus \ang{\phi}} (d)
(b) edge node {\big(\ang{\psi_i}_{i \in \ufs{q}} \, , \ang{\phi_k}_{k \in \ufs{s}} \big)} (c)
;
\end{tikzpicture}
\end{equation}
\item[Monoidal unit]
The monoidal unit is the empty tuple $\ang{} \in \Gsk$.
\item[Braiding] The component of the braiding $\xi$ at a pair of objects $(\angordm, \angordn) \in \Gsk^2$ is defined as $1_\vstar$ if either $\angordm$ or $\angordn$ is the basepoint $\vstar \in \GG$.  

For non-basepoint objects, it is defined as the isomorphism
\begin{equation}\label{Gsk_braiding}
\xi_{\angordm, \angordn} = (\tau_{p,q} \, , \ang{1})
\cn \angordm \oplus \angordn \fto{\iso} \angordn \oplus \angordm.
\end{equation}
The reindexing injection
\[\ufs{p+q} \fto[\iso]{\tau_{p,q}} \ufs{q+p}\]
is the block permutation that interchanges the first $p$ elements with the last $q$ elements:
\[\tau_{p,q}(i) = \begin{cases} 
q+i & \text{if $1 \leq i \leq p$ and}\\
i-p & \text{if $p+1 \leq i \leq p+q$}.
\end{cases}\]
The morphism $\ang{1}$ in \cref{Gsk_braiding} is the $(q+p)$-tuple with each entry given by an identity function of some $\ordn_i$ for $i \in \ufs{q}$ or some $\ordm_k$ for $k \in \ufs{p}$.
\end{description}
This finishes the definition of $(\oplus, \ang{}, \xi)$.  \cref{Gsk_permutative} verifies that this triple is a permutative category structure on $\Gsk$.
\end{definition}

\begin{lemma}\label{Gsk_permutative}
In \cref{def:Gsk_permutative}, the quadruple $(\Gsk,\oplus,\ang{},\xi)$ is a permutative category.
\end{lemma}

\begin{proof}
\parhead{Monoidal structure}. There is an equality of reindexing injections
\[f'f \oplus h'h = (f' \oplus h') \circ (f \oplus h)\]
whenever the left-hand side is defined.  Functoriality of $\oplus$ follows from the equality
\[(f' \oplus h')_* \big(\ang{\psi} \oplus \ang{\phi} \big) = f'_*\ang{\psi} \oplus h'_*\ang{\phi}\]
of tuples of pointed functions.  By definitions \cref{Gsk_oplus_obj,Gsk_oplus_vstar,Gsk_oplus_morphism}, $\oplus$ is associative with the empty tuple $\ang{}$ as a two-sided unit.  Thus, once we define the associativity and unit isomorphisms to be identities, $(\Gsk,\oplus,\ang{})$ is a strict monoidal category.

\parhead{Permutative structure}. The braiding component $\xi_{\angordm,\angordn}$ defined in \cref{Gsk_braiding} is a well-defined morphism in $\Gsk$ because
\[(\tau_{p,q})_* (\angordm \oplus \angordn) = \angordn \oplus \angordm.\]
Naturality of $\xi$ means the commutativity of the following diagram in $\Gsk$.
\begin{equation}\label{Gsk_xi_natural}
\begin{tikzpicture}[vcenter]
\def\v{-1.4}
\draw[0cell=1]
(0,0) node (a1) {\angordm \oplus \angordj}
(a1)++(3.8,0) node (a2) {\angordj \oplus \angordm}
(a1)++(0,\v) node (b1) {\angordn \oplus \angordl}
(a2)++(0,\v) node (b2) {\angordl \oplus \angordn}
;
\draw[1cell=.9]
(a1) edge node {(\tau_{p,r} \spc \ang{1})} (a2)
(b1) edge node {(\tau_{q,s} \spc \ang{1})} (b2)
(a1) edge[transform canvas={xshift=1.5em}] node[swap] {(f,\ang{\psi}) \oplus (h,\ang{\phi})} (b1)
(a2) edge[transform canvas={xshift=-1.5em}] node {(h,\ang{\phi}) \oplus (f,\ang{\psi})} (b2)
;
\end{tikzpicture}
\end{equation}
The reindexing injections of the two composites in \cref{Gsk_xi_natural} are the following two composites.
\[\begin{tikzpicture}
\def\v{-1.3}
\draw[0cell=1]
(0,0) node (a1) {\ufs{p+r}}
(a1)++(3,0) node (a2) {\ufs{r+p}}
(a1)++(0,\v) node (b1) {\ufs{q+s}}
(a2)++(0,\v) node (b2) {\ufs{s+q}}
;
\draw[1cell=.9]
(a1) edge node {\tau_{p,r}} (a2)
(b1) edge node {\tau_{q,s}} (b2)
(a1) edge node[swap] {f \oplus h} (b1)
(a2) edge node {h \oplus f} (b2)
;
\end{tikzpicture}\]
This diagram commutes because the block permutation $\tau_{?,?}$ is natural with respect to functions of unpointed finite sets.  In the second entry, the two composites in \cref{Gsk_xi_natural} are equal because
\[(\tau_{q,s})_* \big(\ang{\psi} \oplus \ang{\phi}\big) = \ang{\phi} \oplus \ang{\psi}.\]
This shows that the braiding $\xi$ is a natural isomorphism.

By definition \cref{Gsk_braiding}, the second entry of each component of $\xi$ is a tuple of identity morphisms.  Thus, the symmetry and hexagon axioms in \cref{symmoncatsymhexagon} follow from the following equalities of permutations.
\[\begin{split}
\tau_{q,p} \circ \tau_{p,q} &= 1_{\ufs{p+q}}\\
\big(1_{\ufs{q}} \oplus \tau_{p,r}\big)\circ \big(\tau_{p,q} \oplus 1_{\ufs{r}}\big) &= \tau_{p,q+r}
\end{split}\]
This proves that $\big(\Gsk,\oplus,\ang{},\xi\big)$ is a permutative category.
\end{proof}

\subsection*{Comparing $\Gsk$ and $\Fsk$}

Recall from \cref{Fsk_permutative} that $(\Fsk, \sma, \ord{1}, \xi)$ is a permutative category.  Next, we compare the permutative categories $\Gsk$ and $\Fsk$ via the functor in the following definition.

\begin{definition}\label{def:smashFskGsk}
We define a functor
\[\Gsk \fto{\sma} \Fsk\]
as follows.
\begin{description}
\item[Objects] The object assignment is defined as follows for $\ang{\ord{m}_k}_{k \in \ufs{p}} \in \Gsk \setminus \{\vstar, \ang{}\}$.  
\begin{equation}\label{smash_Gskobjects}
\left\{
\begin{split}
\sma \vstar &= \ord{0}\\
\sma \ang{} &= \ord{1}\\
\sma \ang{\ord{m}_k}_{k \in \ufs{p}} &= \sma_{k \in \ufs{p}} \,\ord{m}_k
= \ord{\txprod_{k \in \ufs{p}}\, m_k}
\end{split}\right.
\end{equation}
In \cref{smash_Gskobjects}, $\sma_{k \in \ufs{p}}$ is the $p$-fold iterate of the monoidal product of $\Fsk$ defined in \cref{m-sma-n}.  This is well defined because if some $\ord{m}_k = \ord{0}$, then $\ord{\txprod_{k \in \ufs{p}}\, m_k} = \ord{0}$.
\item[Morphisms]
The two morphisms in $\Gsk\big(\ang{}, \ang{}\big)$ are the identity morphism and the 0-morphism \cref{Gsk_empty_mor}.  They are sent by $\sma$ to, respectively, the identity morphism and the 0-morphism in $\Fsk\big(\ord{1}, \ord{1}\big)$.

For a morphism in $\Gsk$ as defined in \cref{fangpsi}, 
\[\angordm \fto{(f, \ang{\psi})} \angordn,\]
the morphism
\[\sma \angordm \fto{\sma (f, \ang{\psi})} \!\smam\angordn \inspace \Fsk\]
is defined as the following composite pointed function.
\begin{equation}\label{smash_fpsi}
\begin{tikzpicture}[vcenter]
\def\t{.6ex}
\draw[0cell=.9]
(0,0) node (a) {\bigwedge_{k \in \ufs{p}} \ord{m}_k}
(a)++(2.8,0) node (b) {\bigwedge_{\finv(i) \neq \emptyset} \ord{m}_{\finv(i)}}
(b)++(3,0) node (c) {\bigwedge_{i \in \ufs{q}} \ord{m}_{\finv(i)}}
(c)++(2.5,0) node (d) {\bigwedge_{i \in \ufs{q}} \ord{n}_i}
;
\draw[1cell=.9]
(a) edge[transform canvas={yshift=\t}] node {f_*} node[swap] {\iso} (b)
(b) edge[transform canvas={yshift=\t}] node {\iso} (c)
(c) edge[transform canvas={yshift=\t}] node {\sma_i\, \psi_i} (d)
;
\end{tikzpicture}
\end{equation}
The three pointed functions in \cref{smash_fpsi} are defined as follows. 
\begin{itemize}
\item The pointed bijection $f_*$ permutes the $p$ entries according to the reindexing injection $f \cn \ufs{p} \to \ufs{q}$.  The indexing set in the codomain is given by $\big\{ i \in \ufs{q} \cn \finv(i) \neq \emptyset\}$.
\item Using the definition \cref{ordn_empty}, the middle pointed bijection in \cref{smash_fpsi} inserts a copy of the smash unit 
\[\ord{1} = \ord{m}_{\emptyset}\]
for each index $i \in \ufs{q}$ satisfying $\finv(i) = \emptyset$.
\item The arrow $\sma_i\, \psi_i$ in \cref{smash_fpsi} is the smash product of the pointed functions $\psi_i$ for $i \in \ufs{q}$.
\end{itemize}
This yields a well-defined morphism $\sma(f,\ang{\psi})$ because if any $\psi_i$ factors through $\ord{0}$, then the composite in \cref{smash_fpsi} also factors through $\ord{0}$.
\end{description}
This finishes the definition of the functor $\sma$.
\end{definition}

\begin{explanation}[Smash of Morphisms]\label{expl:smash_fpsi}
We explain the morphism $\sma(f,\ang{\psi})$ defined in \cref{smash_fpsi} more explicitly.  The domain $\sma_{k \in \ufs{p}}\, \ord{m}_k$ is the smash product defined in \cref{m-sma-n}.  It uses the lexicographic ordering \cref{lex_bijection} to make the identification
\[\ord{\txprod_{k \in \ufs{p}} \, m_k} = \txsma_{k \in \ufs{p}} \ord{m}_k.\]
Under this identification, the morphism $\sma(f,\ang{\psi})$ sends an element $\ang{a_k}_{k \in \ufs{p}} \in \sma_{k \in \ufs{p}}\, \ord{m}_k$ to
\begin{equation}\label{smafangpsiangs}
\sma(f,\ang{\psi}) \ang{a_k}_{k \in \ufs{p}} = \bang{\psi_i a_{\finv(i)}}_{i \in \ufs{q}} \in \txsma_{i \in \ufs{q}} \ord{n}_i,
\end{equation}
where
\[a_{\emptyset} = 1 \in \ord{1} = \ord{m}_\emptyset\]
for each index $i \in \ufs{q}$ with $\finv(i) = \emptyset$.
\end{explanation}

\begin{lemma}\label{sma_symmon}
The assignments in \cref{def:smashFskGsk} define a strict symmetric monoidal pointed functor
\[\big(\Gsk,\oplus,\ang{},\xi\big) \fto{\sma} \big(\Fsk,\sma,\ord{1},\xi\big).\]
\end{lemma}

\begin{proof}
\parhead{Pointed functoriality}.  To see that $\sma$ is a pointed functor, we observe the following.
\begin{itemize}
\item By the first line of the definition \cref{smash_Gskobjects}, $\sma$ sends the basepoint $\vstar \in \Gsk$ to the basepoint $\ord{0} \in \Fsk$. 
\item The morphism assignment of $\sma$ preserves identity morphisms $(1_{\ufs{p}} , \ang{1})$ because, in this case, each of the three arrows in \cref{smash_fpsi} is an identity function. 
\item To see that $\sma$ preserves composition of morphisms in $\Gsk$ \cref{Gsk_composite}, we consider composable morphisms in $\Gsk$ \cref{Gsk_composable}
\[\angordm \fto{(f,\ang{\psi})} \angordn \fto{(h,\ang{\phi})} \angordl,\]
and the following two morphisms in $\Fsk$.
\begin{equation}\label{sma_hphi_fpsi}
\begin{split}
\big(\! \smam(h,\angphi)\big) \circ \big(\! \smam(f,\angpsi)\big) & \cn \smam\angordm \to \smam\angordl \\
\sma\big((h,\angphi) \circ (f,\angpsi)\big) & \cn \smam\angordm \to \smam\angordl
\end{split}
\end{equation}
If either $(f,\angpsi)$ or $(h,\angphi)$ is the 0-morphism in $\Gsk$, then the two morphisms in \cref{sma_hphi_fpsi} are both given by the 0-morphism in $\Fsk$.  Thus, we assume that $(f,\angpsi)$ and $(h,\angphi)$ are nonzero morphisms.  By \cref{expl:Gsk_composite,expl:smash_fpsi}, the following computation proves that the two morphisms in \cref{sma_hphi_fpsi} are equal on each element $\ang{a_k}_{k \in \ufs{p}} \in \sma_{k \in \ufs{p}}\, \ordm_k$.  
\[\begin{split}
& \big(\! \smam(h,\angphi)\big) \circ \big(\! \smam(f,\angpsi)\big) \ang{a_k}_{k \in \ufs{p}} \\
&= \big(\! \smam(h,\angphi)\big) \bang{\psi_i a_{\finv (i)}}_{i \in \ufs{q}} \\
&= \bang{\phi_j \psi_{\hinv(j)} a_{\finv \hinv(j)}}_{j \in \ufs{r}} \\
&= \bang{\phi_j \psi_{\hinv(j)} a_{(hf)^{-1}(j)}}_{j \in \ufs{r}} \\
&= \big(\! \smam\big(hf, \ang{\phi} \circ h_*\ang{\psi} \big) \big) \ang{a_k}_{k \in \ufs{p}} \\
&= \big(\! \smam\big((h,\angphi) \circ (f,\angpsi)\big)\big) \ang{a_k}_{k \in \ufs{p}} 
\end{split}\]
\end{itemize} 
Thus, $\sma \cn \GG \to \FG$ is a pointed functor.

\parhead{Strict symmetric monoidality}.  Next, we check that $\sma \cn \Gsk \to \Fsk$ is a strict symmetric monoidal functor.  We need to check that $\sma$ strictly preserves the monoidal units, the monoidal products, and the braidings.  The functor $\sma$ sends the monoidal unit $\ang{} \in \Gsk$ to the monoidal unit $\ord{1} \in \Fsk$ by the second line of the definition \cref{smash_Gskobjects}.  

The functor $\sma \cn \Gsk \to \Fsk$ preserves the monoidal products of its domain and codomain for the following reasons.
\begin{itemize}
\item For objects, this follows from the definition \cref{Gsk_oplus_obj} of $\oplus$ on objects, which is given by concatenation, and the definition \cref{smash_Gskobjects} of $\sma$ on objects, which is given by smash products of pointed finite sets.
\item For morphisms, we consider two morphisms in $\Gsk$ \cref{fpsi_hphi}
\[\angordm \fto{(f,\ang{\psi})} \angordn \andspace 
\angordj \fto{(h,\ang{\phi})} \angordl.\]
By \cref{Gsk_oplus_morphism,smafangpsiangs}, both morphisms
\[\begin{split}
\sma\big(f \oplus h, \angpsi \oplus \angphi \big) & \cn \smam(\angordm \oplus \angordj) \to \smam(\angordn \oplus \angordl) \andspace\\
\big(\!\smam(f,\angpsi)\big) \sma \big(\!\smam(h,\angphi)\big) & \cn (\sma\angordm) \sma (\sma\angordj) \to (\sma\angordn) \sma (\sma\angordl)
\end{split}\]
in $\Fsk$ send an element
\[\big(\ang{a_k}_{k \in \ufs{p}} \spc \ang{b_t}_{t \in \ufs{r}}\big) 
\in (\sma\angordm) \sma (\sma\angordj)\]
to the element
\[\big(\ang{\psi_i a_{\finv(i)}}_{i \in \ufs{q}} \spc \ang{\phi_u b_{\hinv(u)}}_{u \in \ufs{s}} \big) 
\in (\sma\angordn) \sma (\sma\angordl).\]
\end{itemize}
This proves that $\sma$ preserves the monoidal products.

Using the notation in \cref{Fsk_braiding,Gsk_braiding}, the functor $\sma \cn \Gsk \to \Fsk$ preserves the braidings of its domain and codomain if and only if the following diagram in $\Fsk$ commutes.
\begin{equation}\label{sma_braiding}
\begin{tikzpicture}[vcenter]
\def\v{-1.2}
\draw[0cell=.9]
(0,0) node (a1) {(\sma \angordm) \sma (\sma \angordn)}
(a1)++(5.3,0) node (a2) {(\sma \angordn) \sma (\sma \angordm)}
(a1)++(0,\v) node (b1) {\sma\big(\angordm \oplus \angordn \big)}
(a2)++(0,\v) node (b2) {\sma\big(\angordn \oplus \angordm\big)}
;
\draw[1cell=.85]
(a1) edge[equal] (b1)
(a2) edge[equal] (b2)
(a1) edge node {\xi_{\sma \angordm ,\, \sma \angordn}} (a2)
(b1) edge node {\sma\xi_{\angordm,\angordn} = \sma \big(\tau_{p,q} \spc \ang{1} \big)} (b2)
;
\end{tikzpicture}
\end{equation}
The top and bottom horizontal arrows in \cref{sma_braiding} are given as follows.
\begin{itemize}
\item By definition \cref{Fsk_braiding}, we have
\[\xi_{\sma \angordm ,\, \sma \angordn} = (\tau_{p,q})_*,\]
which block permutes a $p$-tuple with a $q$-tuple.
\item By definition \cref{smash_fpsi}, $\sma (\tau_{p,q} , \ang{1})$ is the composite of $(\tau_{p,q})_*$ with two identity morphisms. 
\end{itemize} 
This proves that the diagram \cref{sma_braiding} commutes and that $\sma$ preserves the braidings.  This finishes the proof that $\sma$ is a strict symmetric monoidal pointed functor.
\end{proof}

\section{Symmetric Monoidal Closed Category Structure}
\label{sec:ggcat_sm}

For an arbitrary group $G$, this section discusses a symmetric monoidal closed category $\GGCat$, in the sense of \cref{def:closedcat}.  Its objects are pointed functors
\[(\Gsk,\vstar) \fto{f} (\Gcat,\boldone),\]
called \emph{$\Gskg$-categories}, with $\Gsk$ the small pointed category defined in \cref{def:Gsk}.  The $\Gcat$-multicategory associated to $\GGCat$ is discussed in \cref{sec:ggcat_multicat}.

\secoutline
\begin{itemize}
\item The symmetric monoidal closed structure on $\GGCat$ uses a pointed variant of $\Gcat$.  \cref{def:ptGcat,def:gcatst} define the 2-category $\Gcatst$ of small pointed $G$-categories, pointed $G$-functors, and pointed $G$-natural transformations.
\item \cref{expl:Gcatst} discusses the symmetric monoidal closed structure on $\Gcatst$.
\item \cref{def:GGCat} defines the symmetric monoidal closed category $\GGCat$.  It is further elaborated in \cref{expl:ggcat_obj,expl:ggcat_mor,expl:ggcat_unit,expl:ggcat_smag,expl:ggcat_brkst}.  In particular, by \cref{ggcat_obj}, a $\Gskg$-category is equivalent to a pointed functor $\Gsk \to \Gcatst$.
\end{itemize}

\subsection*{Pointed $G$-Categories}
Recall from \cref{expl:GCat} that a \emph{$G$-category} is a category equipped with a $G$-action.  \cref{def:pointed-category} is adapted to the $G$-equivariant setting as follows.

\begin{definition}\label{def:ptGcat}
We define the following for an arbitrary group $G$.
\begin{itemize}
\item A \emph{pointed $G$-category}\index{pointed G-category@pointed $G$-category}\index{G-category@$G$-category!pointed} is a $G$-category equipped with a $G$-fixed \emph{basepoint}.
\item A \emph{pointed $G$-functor}\index{pointed G-functor@pointed $G$-functor}\index{G-functor@$G$-functor!pointed} between pointed $G$-categories is a $G$-equivariant basepoint-preserving functor. 
\item A \emph{pointed $G$-natural transformation}\index{pointed G-natural transformation@pointed $G$-natural transformation}\index{G-natural transformation@$G$-natural transformation!pointed} between pointed $G$-functors is a $G$-equivariant pointed natural transformation.  Pointedness means that the basepoint-component is the identity morphism at the basepoint.
\item A \emph{pointed $G$-modification}\index{pointed G-modification@pointed $G$-modification}\index{G-modification@$G$-modification!pointed} between pointed $G$-natural transformations is a $G$-equivariant pointed modification.  Pointedness means that the basepoint-component is the identity 2-cell of the identity 1-cell at the basepoint.\defmark
\end{itemize}
\end{definition}

Recall from \cref{def:GCat} the 2-category $\Gcat$ of small $G$-categories, $G$-functors, and $G$-natural transformations.  Its underlying 1-category is complete, cocomplete, and Cartesian closed (\cref{expl:Gcat_closed}), with terminal object $\boldone$ and internal hom $\Catg$ (\cref{def:Catg}).  The next definition is the pointed analogue of $\Gcat$.

\begin{definition}\label{def:gcatst}
For an arbitrary group $G$, the 2-category\index{pointed G-category@pointed $G$-category!2-category}\index{2-category!pointed G-category@pointed $G$-category}
\begin{equation}\label{Gcatst}
\Gcatst
\end{equation}
is defined by the following data.
\begin{itemize}
\item Objects are small pointed $G$-categories.
\item 1-cells are pointed $G$-functors.
\item 2-cells are pointed $G$-natural transformations.
\item Identity 1-cells, identity 2-cells, vertical composition of 2-cells, and horizontal composition of 1-cells and 2-cells are given by the corresponding structures for functors and natural transformations.
\end{itemize}
The underlying 1-category of $\Gcatst$ is denoted by the same notation.  For the trivial group $G$, $\Gcatst$ is denoted by $\Catst$.
\end{definition}

\begin{explanation}[Symmetric Monoidal Closed Structure]\label{expl:Gcatst}
There is a symmetric monoidal closed structure on $\Gcatst$ given by the smash product.  In more detail, for a $G$-category $\C$, a $G$-fixed object $c \in \C$ is uniquely determined by a $G$-functor 
\[f \cn \boldone \to \C \stspace f(*) = c.\]  
The 1-category $\Gcatst$, with small pointed $G$-categories as objects and pointed $G$-functors as morphisms, is precisely the category of pointed objects in $\Gcat$ equipped with the terminal object $\boldone$, as defined in \cref{def:pointed-objects}.  Thus, the constructions in \cref{def:wedge-smash-phom}---wedge, smash product, smash unit, and pointed hom---apply to small pointed $G$-categories.

Applying \cref{theorem:pC-sm-closed} to the complete and cocomplete Cartesian closed category \pcref{expl:Gcat_closed}
\[(\Gcat, \times, \boldone, \Catg)\]
with terminal object $\boldone$ yields the complete and cocomplete symmetric monoidal closed category
\begin{equation}\label{Gcatst_smc}
(\Gcatst, \sma, \bonep, \Catgst).
\end{equation}
\begin{itemize}
\item For small pointed $G$-categories $\C$ and $\D$, by \cref{eq:smash}, the small pointed $G$-category $\C \sma \D$ is defined in $\Catst$, with $G$ acting diagonally, which means
\[g (c,d) = (gc,gd)\]
for $g \in G$, $c \in \C$, and $d \in \D$. 
\item The monoidal unit \cref{smash-unit-object} of $\Gcatst$ is the discrete category with two objects, 
\begin{equation}\label{gcatst_unit}
\bonep = \boldone \bincoprod \boldone,
\end{equation} 
on which $G$ acts trivially. 
\item The internal hom \cref{eq:pHom-pullback} for $\Gcatst$ is given by the small pointed $G$-category 
\begin{equation}\label{catgst_cd}
\Catgst(\C,\D)
\end{equation}
with
\begin{itemize}
\item pointed functors $\C \to \D$ as objects;
\item basepoint given by the constant functor at the basepoint in $\D$;
\item pointed natural transformations as morphisms;
\item identities and composition defined componentwise in $\D$; and
\item $G$ acting by conjugation, in the sense of \cref{conjugation-gaction}.
\end{itemize}
\end{itemize}
Taking the $G$-fixed subcategory \cref{Gfixed} yields the category
\begin{equation}\label{catgst_gcatst}
\Catgst(\C,\D)^G = \Gcatst(\C,\D),
\end{equation}
which has
\begin{itemize}
\item pointed $G$-functors $\C \to \D$ as objects and
\item pointed $G$-natural transformations as morphisms.\defmark
\end{itemize}
\end{explanation}

\subsection*{$\Gskg$-Categories: Symmetric Monoidal Closed Structure}

\cref{def:GGCat} below uses
\begin{itemize}
\item the complete and cocomplete Cartesian closed category
\[(\Gcat, \times, \boldone, \Catg),\]
with terminal object $\boldone$ and internal hom $\Catg$ (\cref{expl:Gcat_closed}), and
\item the small permutative category 
\[(\Gsk,\oplus,\ang{},\xi)\]
with null object $\vstar$ (\cref{def:Gsk,def:Gsk_permutative}).
\end{itemize} 

\begin{definition}[$\Gskg$-Categories]\label{def:GGCat}
For an arbitrary group $G$, applying \cref{thm:Dgm-pv-convolution-hom} with $\V = \Gcat$ and $\Dgm = \Gsk$, we define the complete and cocomplete symmetric monoidal closed category
\begin{equation}\label{GGCat_smc}
\big(\GGCat, \smag, \gu, \brkst\big).
\end{equation}
The objects of $\GGCat$ are called \index{G-G-category@$\Gskg$-category}\emph{$\Gskg$-categories}. 
\end{definition}

\cref{expl:ggcat_obj,expl:ggcat_mor,expl:ggcat_unit,expl:ggcat_smag,expl:ggcat_brkst} unravel the symmetric monoidal closed category $\GGCat$, starting with its objects and morphisms, followed by the monoidal unit $\gu$, the monoidal product $\smag$, and the internal hom $\brkst$.

\begin{explanation}[$\Gskg$-Categories]\label{expl:ggcat_obj}
By \cref{def:pointed-category}, a $\Gskg$-category is a pointed functor
\begin{equation}\label{ggcat_objects}
(\Gsk, \vstar) \fto{f} (\Gcat,\boldone).
\end{equation}
\begin{itemize}
\item $f$ sends each object $\angordm \in \Gsk$, as defined in \cref{Gsk_objects}, to a small $G$-category $f\angordm$ such that $f\vstar = \boldone$. 
\item $f$ sends each morphism $\upom \cn \angordm \to \angordn$ in $\Gsk$, as defined in \cref{Gsk_morphisms}, to a $G$-functor
\begin{equation}\label{f_upom}
f\angordm \fto{f\upom} f\angordn
\end{equation}
such that $f$ preserves identities and composition of morphisms.  
\item Each $G$-category $f\angordm$ is regarded as a pointed $G$-category with the \emph{canonical basepoint}\dindex{canonical}{basepoint} given by the $G$-functor
\begin{equation}\label{fm_pointed}
f(\vstar \to \angordm) \cn f\vstar = \boldone \to f\angordm,
\end{equation}
where $\vstar \to \angordm$ is the unique morphism in $\Gsk(\vstar,\angordm)$.  Since $\vstar \in \Gsk$ is an initial object, the functoriality of $f$ implies that each $G$-functor $f\upom$ is pointed.
\end{itemize}
Thus, a $\Gskg$-category is the same thing as a pointed functor
\begin{equation}\label{ggcat_obj}
(\Gsk, \vstar) \fto{f} (\Gcatst,\boldone).
\end{equation}
We use the description \cref{ggcat_obj} of a $\Gskg$-category when we want to emphasize that the $G$-categories $f\angordm$ are pointed.
\end{explanation}

\begin{explanation}[Morphisms of $\Gskg$-Categories]\label{expl:ggcat_mor}
By \cref{def:pointed-category}, for $\Gskg$-categories $f$ and $f'$, a morphism $\theta \cn f \to f'$ in $\GGCat$ is a pointed natural transformation as displayed below.
\begin{equation}\label{ggcat_morphisms}
\begin{tikzpicture}[vcenter]
\def\t{28}
\draw[0cell]
(0,0) node (a1) {\Gsk}
(a1)++(1.8,0) node (a2) {\phantom{\Gskel}}
(a2)++(.25,0) node (a2') {\Gcat}
;
\draw[1cell=.9]
(a1) edge[bend left=\t] node {f} (a2)
(a1) edge[bend right=\t] node[swap] {f'} (a2)
;
\draw[2cell]
node[between=a1 and a2 at .45, rotate=-90, 2label={above,\theta}] {\Rightarrow}
;
\end{tikzpicture}
\end{equation}
In more detail, $\theta$ consists of, for each object $\angordm \in \Gsk$ \cref{Gsk_objects}, an $\angordm$-component $G$-functor
\begin{equation}\label{ggcat_mor_component}
f\angordm \fto{\theta_{\angordm}} f'\angordm
\end{equation}
such that, for each morphism $\upom \cn \angordm \to \angordn$ in $\Gsk$ \cref{Gsk_morphisms}, the following naturality diagram of $G$-functors commutes.
\begin{equation}\label{ggcat_mor_naturality}
\begin{tikzpicture}[vcenter]
\def\v{-1.4}
\draw[0cell]
(0,0) node (a11) {f\angordm}
(a11)++(2.5,0) node (a12) {f'\angordm}
(a11)++(0,\v) node (a21) {f\angordn}
(a12)++(0,\v) node (a22) {f'\angordn}
;
\draw[1cell=.9]
(a11) edge node {\theta_{\angordm}} (a12)
(a12) edge node {f'\upom} (a22)
(a11) edge node[swap] {f\upom} (a21)
(a21) edge node {\theta_{\angordn}} (a22)
;
\end{tikzpicture}
\end{equation}
Identity morphisms and composition in $\GGCat$ are defined componentwise using the component $G$-functors in \cref{ggcat_mor_component}.

As a natural transformation, $\theta$ is automatically pointed because its $\vstar$-component,
\begin{equation}\label{ggcat_mor_vstar}
f\vstar = \boldone \fto{\theta_\vstar} f'\vstar = \boldone,
\end{equation}
is necessarily the identity functor on the terminal $G$-category $\boldone$. 

In \cref{ggcat_morphisms}, we can equivalently use $\Gcatst$ in place of $\Gcat$.  By \cref{fm_pointed}, each $G$-category $f\angordm$ is canonically pointed.  Each component $G$-functor $\theta_{\angordm}$ is also pointed by the naturality diagram \cref{ggcat_mor_naturality}, applied to the unique morphism $\vstar \to \angordm$.  Thus, a morphism $\theta \cn f \to f'$ in $\GGCat$ is the same thing as a pointed natural transformation as displayed below.
\begin{equation}\label{ggcat_mor}
\begin{tikzpicture}[vcenter]
\def\t{28}
\draw[0cell]
(0,0) node (a1) {\Gsk}
(a1)++(1.8,0) node (a2) {\phantom{\Gskel}}
(a2)++(.3,0) node (a2') {\Gcatst}
;
\draw[1cell=.9]
(a1) edge[bend left=\t] node {f} (a2)
(a1) edge[bend right=\t] node[swap] {f'} (a2)
;
\draw[2cell]
node[between=a1 and a2 at .45, rotate=-90, 2label={above,\theta}] {\Rightarrow}
;
\end{tikzpicture}
\end{equation}
We use the description \cref{ggcat_mor} of a morphism of $\Gskg$-categories when we want to emphasize that the component $G$-functors $\theta_{\angordm}$ are pointed.
\end{explanation}

\begin{explanation}[Monoidal Unit of $\GGCat$]\label{expl:ggcat_unit}
By \cref{reindexing_functor,ordn_empty,fangpsi,smash_Gskobjects,ptdayunit}, the monoidal unit of $\GGCat$ is the pointed functor
\begin{equation}\label{GGCat_unit}
(\Gsk, \vstar) \fto{\gu} (\Gcatst,\boldone)
\end{equation}
given on objects by
\begin{equation}\label{gu_angordm}
\begin{split}
\gu\angordm &= \txwedge_{\Gskpunc(\ang{}, \angordm)} \, \bonep \\
&\iso \txsma_{i \in \ufs{p}}\, \ordm_i = \sma\angordm
\end{split}
\end{equation}
for $\angordm = \ang{\ordm_i}_{i \in \ufs{p}} \in \Gsk$.
\begin{itemize}
\item In \cref{gu_angordm}, $\Gskpunc$ denotes the set of nonzero morphisms in $\Gsk$, and $\bonep$ is the smash unit in $\Gcatst$ \cref{gcatst_unit}.  An empty wedge is defined to be $\bone$.
\item The isomorphism in \cref{gu_angordm} follows from the fact that each nonzero morphism 
\[\ang{} \to \angordm \inspace \Gsk\] 
is uniquely determined by an element in the unpointed finite set $\ufs{m}_i$, one for each $i \in \ufs{p}$.
\end{itemize}
The pointed finite set $\sma\angordm$ is regarded as a discrete pointed $G$-category with the trivial $G$-action.  On morphisms, $\gu$ is given by the assignment $\upom \mapsto \!\smam\upom$ defined in \cref{smash_fpsi}.
\end{explanation}

\begin{explanation}[Monoidal Product of $\GGCat$]\label{expl:ggcat_smag}
By \cref{ptdayconv}, for two $\Gskg$-categories $f$ and $f'$, their monoidal product is the $\Gskg$-category given by the coend 
\begin{equation}\label{ggcat_ptday}
f \smag f' = \ecint^{(\angordm, \angordmp) \in \Gsk^2} 
\mquad \bigvee_{\Gskpunc(\angordm \oplus \angordmp, -)} \mquad (f\angordm \sma f'\angordmp)
\end{equation}
taken in $\Gcatst$.  An empty wedge and $(f \smag f')(\vstar)$ are both defined to be $\boldone$.  By  \cref{expl:pointedday}, given another $\Gskg$-category $f''$, a pointed natural transformation
\[\begin{tikzcd}[column sep=normal]
f \smag f' \ar{r}{\theta} & f''
\end{tikzcd}\]
consists of component pointed $G$-functors
\begin{equation}\label{ggcat_sma_ext}
f\angordm \sma f'\angordmp \fto{\theta_{\angordm,\angordmp}} f''(\angordm \oplus \angordmp)
\end{equation}
for $(\angordm, \angordmp) \in \Gsk^2$ such that, for each pair of morphisms \cref{Gsk_morphisms}
\[\angordm \fto{\upom} \angordn \andspace 
\angordmp \fto{\upom'} \angordnp\]
in $\Gsk$, the following diagram of pointed $G$-functors commutes.
\begin{equation}\label{ggcat_sma_nat}
\begin{tikzpicture}[vcenter]
\def\v{-1.4}
\draw[0cell=1]
(0,0) node (x11) {f\angordm \sma f'\angordmp}
(x11)++(4.5,0) node (x12) {f''(\angordm \oplus \angordmp)}
(x11)++(0,\v) node (x21) {f\angordn \sma f'\angordnp}
(x12)++(0,\v) node (x22) {f''(\angordn \oplus \angordnp)}
;
\draw[1cell=.9] 
(x11) edge node {\theta_{\angordm,\angordmp}} (x12)
(x21) edge node {\theta_{\angordn,\angordnp}} (x22)
(x11) edge[transform canvas={xshift=1em}] node[swap] {f\upom \sma f'\upom'} (x21)
(x12) edge[transform canvas={xshift=-1.5em}] node {f''(\upom \oplus \upom')} (x22)
;
\end{tikzpicture}
\end{equation}
This description of $f \smag f'$ is called the \emph{external characterization}.
\end{explanation}

\begin{explanation}[Internal Hom of $\GGCat$]\label{expl:ggcat_brkst}
By \cref{ptdayhom}, for $\Gskg$-categories $f$ and $f'$, the internal hom $\Gskg$-category is given by the end 
\begin{equation}\label{ggcat_inthom}
\brkstg{f}{f'} = \ecint_{\angordn \in \Gsk} \Catgst\big(f\angordn, f'(- \oplus \angordn)\big)
\end{equation}
taken in $\Gcatst$, where $\Catgst$ is the internal hom for $\Gcatst$ \cref{catgst_cd}.  More details are given below.

\parhead{Pointed $G$-categories: objects}.  For each object $\angordm \in \Gsk$ \cref{Gsk_objects}, an object $\theta$ in the pointed $G$-category
\begin{equation}\label{ggcat_inthom_comp}
\brkstg{f}{f'}\angordm = \ecint_{\angordn \in \Gsk} \Catgst\big(f\angordn, f'(\angordm \oplus \angordn)\big)
\end{equation}
consists of $\angordn$-component pointed functors
\begin{equation}\label{ggcat_inthom_theta}
f\angordn \fto{\theta_{\angordn}} f'(\angordm \oplus \angordn) \forspace \angordn \in \Gsk
\end{equation}
such that, for each morphism $\upom \cn \angordn \to \angordl$ in $\Gsk$ \cref{Gsk_morphisms}, the diagram 
\begin{equation}\label{ggcat_inthom_theta_nat}
\begin{tikzpicture}[vcenter]
\def\v{-1.4}
\draw[0cell]
(0,0) node (a11) {f\angordn}
(a11)++(3,0) node (a12) {f'(\angordm \oplus \angordn)}
(a11)++(0,\v) node (a21) {f\angordl}
(a12)++(0,\v) node (a22) {f'(\angordm \oplus \angordl)}
;
\draw[1cell=.9]
(a11) edge node {\theta_{\angordn}} (a12)
(a12) edge[transform canvas={xshift=-2em}] node {f'(1_{\angordm} \oplus \upom)} (a22)
(a11) edge node[swap] {f\upom} (a21)
(a21) edge node {\theta_{\angordl}} (a22)
;
\end{tikzpicture}
\end{equation}
of pointed functors commutes.  We emphasize that the components $\theta_{\angordn}$ are \emph{not} required to be $G$-equivariant, so \cref{ggcat_inthom_theta_nat} is not generally a diagram in $\Gcat$.  In particular, $\theta$ is not a natural transformation $f \to f'(\angordm \oplus -)$.

The basepoint of $\brkstg{f}{f'}\angordm$ has each component functor given by the constant functor at the basepoint of the codomain.  

\parhead{Pointed $G$-categories: morphisms}.
For objects $\theta,\ups \in \brkstg{f}{f'}\angordm$, a morphism $\Theta \cn \theta \to \ups$ consists of $\angordn$-component pointed natural transformations
\begin{equation}\label{ggcat_inthom_Theta}
\begin{tikzpicture}[vcenter]
\def\t{25}
\draw[0cell]
(0,0) node (a1) {\phantom{f'}}
(a1)++(2.5,0) node (a2) {\phantom{f'}}
(a1)++(-.2,0) node (a1') {f\angordn}
(a2)++(.9,0) node (a2') {f'(\angordm \oplus \angordn)}
;
\draw[1cell=.85]
(a1) edge[bend left=\t] node {\theta_{\angordn}} (a2)
(a1) edge[bend right=\t] node[swap] {\ups_{\angordn}} (a2)
;
\draw[2cell]
node[between=a1 and a2 at .37, rotate=-90, 2label={above, \Theta_{\angordn}}] {\Rightarrow}
;
\end{tikzpicture}
\end{equation}
for $\angordn \in \Gsk$ such that, for each morphism $\upom \cn \angordn \to \angordl$ in $\Gsk$, the following two whiskered pointed natural transformations are equal.
\begin{equation}\label{ggcat_inthom_Theta_modax}
\begin{tikzpicture}[vcenter]
\def\t{25} \def\v{-1.6}
\draw[0cell]
(0,0) node (a1) {\phantom{f'}}
(a1)++(2.5,0) node (a2) {\phantom{f'}}
(a1)++(-.2,0) node (a1') {f\angordn}
(a2)++(.9,0) node (a2') {f'(\angordm \oplus \angordn)}
(a1)++(0,\v) node (b1) {\phantom{f'}}
(a2)++(0,\v) node (b2) {\phantom{f'}}
(b1)++(-.2,0) node (b1') {f\angordl}
(b2)++(.9,0) node (b2') {f'(\angordm \oplus \angordl)}
;
\draw[1cell=.8]
(a1) edge[bend left=\t] node[pos=.4] {\theta_{\angordn}} (a2)
(a1) edge[bend right=\t] node[swap,pos=.6] {\ups_{\angordn}} (a2)
(b1) edge[bend left=\t] node[pos=.4] {\theta_{\angordl}} (b2)
(b1) edge[bend right=\t] node[swap,pos=.6] {\ups_{\angordl}} (b2)
(a1) edge node[swap] {f\upom} (b1)
(a2) edge node {f'(1_{\angordm} \oplus \upom)} (b2)
;
\draw[2cell]
node[between=a1 and a2 at .37, rotate=-90, 2label={above, \Theta_{\angordn}}] {\Rightarrow}
node[between=b1 and b2 at .37, rotate=-90, 2label={above, \Theta_{\angordl}}] {\Rightarrow}
;
\end{tikzpicture}
\end{equation}
The natural transformations $\Theta_{\angordn}$ are \emph{not} required to be $G$-equivariant.  Identities and composition of morphisms are given componentwise using \cref{ggcat_inthom_Theta}.

\parhead{$G$-action}.  The $G$-action on the pointed $G$-category $\brkstg{f}{f'}\angordm$ is given componentwise by the conjugation $G$-action.  In other words, for a morphism $\Theta \cn \theta \to \ups$, $g \in G$, and $\angordn \in \Gsk$, the $\angordn$-components of $g \cdot \theta$, $g \cdot \ups$, and $g \cdot \Theta$ are defined as the following composites and whiskering.
\begin{equation}\label{ggcat_inthom_theta_g}
\begin{tikzpicture}[vcenter]
\def\v{-1.5} \def\t{18}
\draw[0cell]
(0,0) node (a11) {\phantom{f'}}
(a11)++(3.5,0) node (a12) {\phantom{f'}}
(a11)++(0,\v) node (a21) {\phantom{f'}}
(a12)++(0,\v) node (a22) {\phantom{f'}}
(a11)++(-.2,0) node (a11') {f\angordn}
(a12)++(.9,0) node (a12') {f'(\angordm \oplus \angordn)}
(a21)++(-.2,0) node (a21') {f\angordn}
(a22)++(.9,0) node (a22') {f'(\angordm \oplus \angordn)}
;
\draw[1cell=.8]
(a11) edge node[swap] {\ginv} (a21)
(a22) edge node[swap] {g} (a12)
(a11) edge[bend left=\t] node[pos=.42] {(g \cdot \theta)_{\angordn}} (a12)
(a21) edge[bend left=\t] node[pos=.4] {\theta_{\angordn}} (a22)
(a11) edge[bend right=\t] node[swap,pos=.58] {(g \cdot \ups)_{\angordn}} (a12)
(a21) edge[bend right=\t] node[swap,pos=.6] {\ups_{\angordn}} (a22)
;
\draw[2cell=.9]
node[between=a11 and a12 at .33, rotate=-90, 2label={above, (g \cdot \Theta)_{\angordn}}] {\Rightarrow}
node[between=a21 and a22 at .4, rotate=-90, 2label={above, \Theta_{\angordn}}] {\Rightarrow}
;
\end{tikzpicture}
\end{equation}
The components of $g \cdot \theta$ satisfy the property \cref{ggcat_inthom_theta_nat} by 
\begin{itemize}
\item the same property for $\theta$ and
\item the fact \cref{f_upom} that $f\upom$ and $f'(1_{\angordm} \oplus \upom)$ are pointed $G$-functors.
\end{itemize}
Similarly, the components of $g \cdot \Theta$ satisfies the property \cref{ggcat_inthom_Theta_modax} by
\begin{itemize}
\item the same property for $\Theta$ and
\item the $G$-equivariance of $f\upom$ and $f'(1_{\angordm} \oplus \upom)$.
\end{itemize}
This finishes the description of the pointed $G$-category $\brkstg{f}{f'}\angordm$ in \cref{ggcat_inthom_comp}.

\parhead{Pointed $G$-functors}.  The internal hom $\Gskg$-category $\brkstg{f}{f'}$ sends a  morphism $\uprho \cn \angordm \to \angordmp$ in $\Gsk$ to the pointed $G$-functor defined by the following commutative diagram.
\begin{equation}\label{ggcat_inthom_gfun}
\begin{tikzpicture}[vcenter]
\def\v{-1.5}
\draw[0cell=.9]
(0,0) node (a11) {\brkstg{f}{f'}\angordm}
(a11)++(0,\v) node (a21) {\brkstg{f}{f'}\angordmp}
(a11)++(4,0) node (a12) {\txint_{\angordn \in \Gsk} \Catgst\big(f\angordn, f'(\angordm \oplus \angordn)\big)}
(a12)++(0,\v) node (a22) {\txint_{\angordn \in \Gsk} \Catgst\big(f\angordn, f'(\angordmp \oplus \angordn)\big)}
;
\draw[1cell=.9]
(a11) edge[equal] (a12)
(a21) edge[equal] (a22)
(a11) edge[transform canvas={xshift=1em}] node[swap] {\brkstg{f}{f'} \uprho} (a21)
(a12) edge[transform canvas={xshift=-4.5em}] node {\txint_{\angordn} \Catgst(1_{f\angordn}, f'(\uprho \oplus 1_{\angordn}))} (a22)
;
\end{tikzpicture}
\end{equation}
In terms of the components in \cref{ggcat_inthom_theta,ggcat_inthom_Theta}, the pointed $G$-functor $\brkstg{f}{f'} \uprho$ is given by post-composing (for objects) and post-whiskering (for morphisms) with the pointed $G$-functor $f'(\uprho \oplus 1_{\angordn})$,
\begin{equation}\label{ggcat_inthom_gfun_comp}
\begin{tikzpicture}[baseline={(a1.base)}]
\def\t{28}
\draw[0cell]
(0,0) node (a1) {\phantom{f'}}
(a1)++(2,0) node (a2) {\phantom{f'}}
(a1)++(-.2,0) node (a1') {f\angordn}
(a2)++(.9,0) node (a2') {f'(\angordm \oplus \angordn)}
(a2')++(.9,0) node (a2'') {\phantom{f'}}
(a2'')++(2.5,0) node (a3) {\phantom{f'}}
(a3)++(.95,0) node (a3') {f'(\angordmp \oplus \angordn)}
;
\draw[1cell=.85]
(a1) edge[bend left=\t] node {\theta_{\angordn}} (a2)
(a1) edge[bend right=\t] node[swap] {\ups_{\angordn}} (a2)
(a2'') edge node {f'(\uprho \oplus 1_{\angordn})} (a3)
;
\draw[2cell=.9]
node[between=a1 and a2 at .35, rotate=-90, 2label={above, \Theta_{\angordn}}] {\Rightarrow}
;
\end{tikzpicture}
\end{equation}
for $\angordn \in \Gsk$. 
\begin{itemize}
\item The properties \cref{ggcat_inthom_theta_nat,ggcat_inthom_Theta_modax} hold for the images of $\brkstg{f}{f'} \uprho$ by the functoriality of $f'$ and $\oplus$ (\cref{Gsk_permutative}).
\item The functoriality of $\brkstg{f}{f'} \uprho$ with respect to $\Theta$ follows from the fact that post-whiskering with $f'(\uprho \oplus 1_{\angordn})$ preserves identities and vertical composition of natural transformations.
\item The $G$-equivariance of $\brkstg{f}{f'} \uprho$ follows from \cref{ggcat_inthom_theta_g} and the $G$-equivariance \cref{f_upom} of $f'(\uprho \oplus 1_{\angordn})$.
\item The functoriality of $\brkstg{f}{f'}$ with respect to $\uprho$ follows from the functoriality of $f'$ and $\oplus$.
\end{itemize}
This finishes the description of the internal hom $\Gskg$-category $\brkstg{f}{f'}$.
\end{explanation}

\section{$\Gcat$-Multicategory Structure}
\label{sec:ggcat_multicat}

For an arbitrary group $G$, this section describes the $\Gcat$-multicategory structure on the symmetric monoidal closed category $\GGCat$ in \cref{def:GGCat}.

\secoutline
\begin{itemize}
\item \cref{expl:ggcat_gcatenr} discusses how the symmetric monoidal closed category $\GGCat$ acquires the structure of a $\Gcat$-multicategory via general enriched category theory.  The rest of this section explains this $\Gcat$-multicategory in detail.
\item \cref{expl:ggcat_zero,expl:ggcat_positive} describe the multimorphism $G$-categories of $\GGCat$ in, respectively, arity 0 and positive arity.
\item \cref{expl:ggcat_symmetry} discusses the symmetric group action $G$-functors on $\GGCat$ and the symmetric group action axioms \cref{enr-multicategory-symmetry}.
\item \cref{expl:ggcat_composition} discusses the composition $G$-functors on $\GGCat$, ending with a justification of the unity and associativity axioms, \cref{enr-multicategory-right-unity,enr-multicategory-left-unity,enr-multicategory-associativity}.
\item \cref{expl:ggcat_topeq,expl:ggcat_boteq} prove, respectively, the top and bottom equivariance axioms, \cref{enr-operadic-eq-1,enr-operadic-eq-2}, for the $\Gcat$-multicategory $\GGCat$.
\end{itemize}

\begin{explanation}[$\Gcat$-Multicategory Structure]\label{expl:ggcat_gcatenr}
The symmetric monoidal closed category $\GGCat$ (\cref{def:GGCat}) yields a \index{G-G-category@$\Gskg$-category!symmetric monoidal category}\index{symmetric monoidal category!G-G-category@$\Gskg$-category}\index{G-G-category@$\Gskg$-category!multicategory}\index{multicategory!G-G-category@$\Gskg$-category}$\Gcat$-multicategory (\cref{def:enr-multicategory}) as follows, where $\Gcat$ is the Cartesian closed category in \cref{expl:Gcat_closed}. 
\begin{itemize}
\item First, \cref{theorem:v-closed-v-sm} implies that $\GGCat$ has the structure of a symmetric monoidal $(\GGCat)$-category 
\begin{equation}\label{ggcat_smggcat}
\underline{\GGCat}
\end{equation}
in the sense of \cref{definition:monoidal-vcat,definition:braided-monoidal-vcat,definition:symm-monoidal-vcat}. 
\item Evaluation at the monoidal unit $\ang{} \in \Gsk$ (\cref{def:Fsk_smashpower,def:Gsk}) yields a symmetric monoidal functor \cref{evtu_unpt}
\begin{equation}\label{evang_gcat}
\GGCat \fto{\ev_{\ang{}}} \Gcat.
\end{equation}
By \cref{thm:change-enrichment} \eqref{change-enr-i}, changing enrichment along $\ev_{\ang{}}$ makes $\underline{\GGCat}$ into a symmetric monoidal $\Gcat$-category 
\begin{equation}\label{ggcat_smgcat}
\underline{\GGCat}_{\ev_{\ang{}}}
\end{equation}
with the same objects as $\GGCat$. 
\item By \cref{proposition:monoidal-v-cat-v-multicat}, the endomorphism construction in \cref{definition:EndK} makes $\underline{\GGCat}_{\ev_{\ang{}}}$ into a $\Gcat$-multicategory
\begin{equation}\label{ggcat_gcatmulti}
\End\big(\underline{\GGCat}_{\ev_{\ang{}}}\big).
\end{equation}
\end{itemize}
To simplify the notation, we denote each of \cref{ggcat_smggcat,ggcat_smgcat,ggcat_gcatmulti} by $\GGCat$.  In each case, the objects are $\Gskg$-categories \cref{ggcat_objects}, which are pointed functors
\[(\Gsk,\vstar) \fto{f} (\Gcat,\boldone).\]
Equipped with the canonical basepoints \cref{fm_pointed}, a $\Gskg$-category is equivalently a pointed functor \cref{ggcat_obj}
\[(\Gsk,\vstar) \fto{f} (\Gcatst,\boldone).\]
The rest of this section unravels the $\Gcat$-multicategory $\GGCat$, using \cref{def:Gsk,def:Gsk_permutative,definition:canonical-v-enrichment,def:enr-multicategory,definition:EndK,definition:Dgm-pV-convolution-hom}.
\end{explanation}

\begin{explanation}[Arity 0 Multimorphism $G$-Categories]\label{expl:ggcat_zero}
For a $\Gskg$-category $f \cn \Gsk \to \Gcatst$, we unravel the 0-ary multimorphism $G$-category \cref{multimorphism_object}
\begin{equation}\label{ggcat_zero_gcat}
\begin{split}
\GGCat(\ang{}; f) &= \brkstg{\gu}{f} \ang{} \\
&= \ecint_{\angordn \in \Gsk} \Catgst(\gu\angordn, f\angordn).
\end{split}
\end{equation}
In \cref{ggcat_zero_gcat}, $\gu$ is the monoidal unit of $\GGCat$ \cref{GGCat_unit}, and $\brkstg{\gu}{f}$ is the internal hom $\Gskg$-category \cref{ggcat_inthom}.

\parhead{Objects}.  By \cref{expl:ggcat_unit,expl:ggcat_brkst}, an object $\theta \in \GGCat(\ang{}; f)$ consists of $\angordn$-component pointed functors
\begin{equation}\label{ggcat_zero_obj_comp}
\sma\angordn = \gu\angordn \fto{\theta_{\angordn}} f\angordn \forspace \angordn \in \Gsk
\end{equation}
such that, for each morphism $\upom \cn \angordn \to \angordl$ in $\Gsk$, the following diagram of pointed functors commutes.
\begin{equation}\label{ggcat_zero_obj_nat}
\begin{tikzpicture}[vcenter]
\def\v{-1.4} \def\h{-1.4}
\draw[0cell]
(0,0) node (a11) {\gu\angordn}
(a11)++(2.3,0) node (a12) {f\angordn}
(a11)++(0,\v) node (a21) {\gu\angordl}
(a12)++(0,\v) node (a22) {f\angordl}
(a11)++(\h,0) node (a11') {\sma\angordn}
(a21)++(\h,0) node (a21') {\sma\angordl}
;
\draw[1cell=.9]
(a11') edge[equal] (a11)
(a21') edge[equal] (a21)
(a11) edge node {\theta_{\angordn}} (a12)
(a12) edge node {f\upom} (a22)
(a11) edge node[swap] {\sma\upom = \gu\upom} (a21)
(a21) edge node {\theta_{\angordl}} (a22)
;
\end{tikzpicture}
\end{equation}
Since the category $\gu\angordn$ is discrete, with only identity morphisms, $\theta_{\angordn}$ is determined by its object assignment.  We emphasize that the pointed functors $\theta_{\angordn}$ are \emph{not} required to be $G$-equivariant.

\parhead{Morphisms}.  For objects $\theta,\ups \in \GGCat(\ang{}; f)$, a morphism $\Theta \cn \theta \to \ups$ consists of $\angordn$-component pointed natural transformations
\begin{equation}\label{ggcat_zero_mor_comp}
\begin{tikzpicture}[vcenter]
\def\t{25}
\draw[0cell]
(0,0) node (a1) {\phantom{f'}}
(a1)++(2.5,0) node (a2) {\phantom{f'}}
(a1)++(-.2,0) node (a1') {\gu\angordn}
(a2)++(.2,0) node (a2') {f\angordn}
;
\draw[1cell=.85]
(a1) edge[bend left=\t] node {\theta_{\angordn}} (a2)
(a1) edge[bend right=\t] node[swap] {\ups_{\angordn}} (a2)
;
\draw[2cell]
node[between=a1 and a2 at .37, rotate=-90, 2label={above, \Theta_{\angordn}}] {\Rightarrow}
;
\end{tikzpicture}
\end{equation}
for $\angordn \in \Gsk$ such that, for each morphism $\upom \cn \angordn \to \angordl$ in $\Gsk$, the following two whiskered pointed natural transformations are equal.
\begin{equation}\label{ggcat_zero_mor_ax}
\begin{tikzpicture}[vcenter]
\def\t{25} \def\v{-1.6}
\draw[0cell]
(0,0) node (a1) {\phantom{f'}}
(a1)++(2.5,0) node (a2) {\phantom{f'}}
(a1)++(-.2,0) node (a1') {\gu\angordn}
(a2)++(.2,0) node (a2') {f\angordn}
(a1)++(0,\v) node (b1) {\phantom{f'}}
(a2)++(0,\v) node (b2) {\phantom{f'}}
(b1)++(-.2,0) node (b1') {\gu\angordl}
(b2)++(.2,0) node (b2') {f\angordl}
;
\draw[1cell=.8]
(a1) edge[bend left=\t] node[pos=.4] {\theta_{\angordn}} (a2)
(a1) edge[bend right=\t] node[swap,pos=.6] {\ups_{\angordn}} (a2)
(b1) edge[bend left=\t] node[pos=.4] {\theta_{\angordl}} (b2)
(b1) edge[bend right=\t] node[swap,pos=.6] {\ups_{\angordl}} (b2)
(a1) edge node[swap] {\gu\upom} (b1)
(a2) edge node {f\upom} (b2)
;
\draw[2cell]
node[between=a1 and a2 at .37, rotate=-90, 2label={above, \Theta_{\angordn}}] {\Rightarrow}
node[between=b1 and b2 at .37, rotate=-90, 2label={above, \Theta_{\angordl}}] {\Rightarrow}
;
\end{tikzpicture}
\end{equation}
The naturality of $\Theta_{\angordn}$ is automatic because the category $\gu\angordn$ is discrete.  The pointed natural transformations $\Theta_{\angordn}$ are \emph{not} required to be $G$-equivariant.  Identities and composition of morphisms are given componentwise using \cref{ggcat_zero_mor_comp}.

\parhead{$G$-action functors}.  The $G$-action on the $G$-category $\GGCat(\ang{}; f)$ is given componentwise by the conjugation $G$-action.  In other words, for a morphism $\Theta \cn \theta \to \ups$, an element $g \in G$, and an object $\angordn \in \Gsk$, the $\angordn$-components of $g \cdot \theta$, $g \cdot \ups$, and $g \cdot \Theta$ are given by the following composites and whiskering.
\begin{equation}\label{ggcat_zero_mor_g}
\begin{tikzpicture}[vcenter]
\def\v{-1.5} \def\t{18}
\draw[0cell]
(0,0) node (a11) {\phantom{f'}}
(a11)++(3.5,0) node (a12) {\phantom{f'}}
(a11)++(0,\v) node (a21) {\phantom{f'}}
(a12)++(0,\v) node (a22) {\phantom{f'}}
(a11)++(-.2,0) node (a11') {\gu\angordn}
(a12)++(.2,0) node (a12') {f\angordn}
(a21)++(-.2,0) node (a21') {\gu\angordn}
(a22)++(.2,0) node (a22') {f\angordn}
;
\draw[1cell=.8]
(a11) edge node[swap] {1 = \ginv} (a21)
(a22) edge node[swap] {g} (a12)
(a11) edge[bend left=\t] node[pos=.42] {(g \cdot \theta)_{\angordn}} (a12)
(a21) edge[bend left=\t] node[pos=.4] {\theta_{\angordn}} (a22)
(a11) edge[bend right=\t] node[swap,pos=.58] {(g \cdot \ups)_{\angordn}} (a12)
(a21) edge[bend right=\t] node[swap,pos=.6] {\ups_{\angordn}} (a22)
;
\draw[2cell=.9]
node[between=a11 and a12 at .33, rotate=-90, 2label={above, (g \cdot \Theta)_{\angordn}}] {\Rightarrow}
node[between=a21 and a22 at .4, rotate=-90, 2label={above, \Theta_{\angordn}}] {\Rightarrow}
;
\end{tikzpicture}
\end{equation}
Since the $G$-action on the discrete category $\gu\angordn = \sma\angordn$ is trivial, the $G$-action on $\GGCat(\ang{}; f)$ is given componentwise by post-composing (for objects) or post-whiskering (for morphisms) with the $G$-action on $f\angordn$.  

The properties \cref{ggcat_zero_obj_nat,ggcat_zero_mor_ax} are satisfied by, respectively, $g \cdot \theta$ and $g \cdot \Theta$ by
\begin{itemize}
\item the same properties for $\theta$ and $\Theta$, and
\item the $G$-equivariance \cref{f_upom} of $f\upom$. 
\end{itemize}
The functoriality of the $g$-action on $\GGCat(\ang{}; f)$ follows from the fact that, in \cref{ggcat_zero_mor_g}, whiskering with $g$ preserves identities and vertical composition of natural transformations.  This finishes the description of the 0-ary multimorphism $G$-category $\GGCat(\ang{}; f)$.
\end{explanation}

\begin{explanation}[Positive Arity Multimorphism $G$-Categories]\label{expl:ggcat_positive}
For $k > 0$ and pointed functors
\[\ang{f_i}_{i \in \ufs{k}} \spc f \cn (\Gsk,\vstar) \to (\Gcatst,\boldone),\]
we unravel the $k$-ary multimorphism $G$-category \cref{multimorphism_object}
\begin{equation}\label{ggcat_kary_gcat}
\begin{split}
& \GGCat\big(\ang{f_i}_{i \in \ufs{k}}; f \big) 
= \Brkstg{\!\txsmag_{i \in \ufs{k}}\, f_i}{f} \ang{} \\
&= \ecint_{\angordn \in \Gsk} \Catgst\Big( \big(\!\txsmag_{i \in \ufs{k}}\, f_i\big) \angordn, f\angordn \Big).
\end{split}
\end{equation}
To simplify the notation, we let the subscript $i$ run through $\ufs{k}$, so $\smag_i$, $\sma_i$,  $\ang{\cdots}_i$, and $\oplus_i$ mean, respectively, $\smag_{i \in \ufs{k}}$, $\sma_{i \in \ufs{k}}$, $\ang{\cdots}_{i \in \ufs{k}}$, and $\oplus_{i \in \ufs{k}}$.  In \cref{ggcat_kary_gcat}, $\brkst$ is the internal hom $\Gskg$-category \cref{ggcat_inthom}, and $\smag_i$ is an iterate of the monoidal product in \cref{ggcat_ptday}.  We have the pointed $G$-category
\begin{equation}\label{smag_fin}
(\txsmag_{i} f_i) \angordn 
= \dint^{\ang{\angordmi}_{i} \,\in\, \Gsk^k} \mquad \bigvee_{\Gskpunc(\oplus_{i} \angordmi, \angordn)} \mquad \big(\!\txsma_i f_i\angordmi \big).
\end{equation}
In \cref{smag_fin}, each $\angordmi \in \Gsk$ is an object.  We use the notation
\begin{equation}\label{angordmdot}
\angordmdot = \ang{\angordmi}_{i \in \ufs{k}} \in \Gsk^k.
\end{equation}

\parhead{Objects}.  By \cref{ggcat_kary_gcat}, \cref{smag_fin}, and the universal properties of ends, coends, and wedges, an object $\theta \in \GGCat(\ang{f_i}_{i}; f)$, which is also denoted by
\begin{equation}\label{ggcat_k1}
\ang{f_i}_{i \in \ufs{k}} \fto{\theta} f,
\end{equation}
consists of $\angordmdot$-component pointed functors
\begin{equation}\label{ggcat_k1_comp}
\txsma_{i \in \ufs{k}}\, f_i\angordmi \fto{\theta_{\angordmdot}} 
f\big(\!\txoplus_{i \in \ufs{k}}\, \angordmi \big)
\forspace \angordmdot \in \Gsk^k
\end{equation}
such that, given morphisms
\begin{equation}\label{Gk_mor}
\bang{\angordmi \fto{\upom^i} \angordni \inspace \Gsk}_{i \in \ufs{k}} \spc
\end{equation}
the following diagram of pointed functors commutes.
\begin{equation}\label{ggcat_k1_nat}
\begin{tikzpicture}[vcenter]
\def\v{-1.4}
\draw[0cell=1]
(0,0) node (x11) {\txsma_{i \in \ufs{k}}\, f_i\angordmi}
(x11)++(4,0) node (x12) {f\big(\!\txoplus_{i \in \ufs{k}}\, \angordmi \big)}
(x11)++(0,\v) node (x21) {\txsma_{i \in \ufs{k}}\, f_i\angordni}
(x12)++(0,\v) node (x22) {f\big(\!\txoplus_{i \in \ufs{k}}\, \angordni \big)}
;
\draw[1cell=.9] 
(x11) edge node {\theta_{\angordmdot}} (x12)
(x21) edge node {\theta_{\angordndot}} (x22)
(x11) edge[transform canvas={xshift=1.3em}] node[swap] {\txsma_i f_i \upom^i} (x21)
(x12) edge[transform canvas={xshift=-2em}] node {f(\oplus_i\, \upom^i)} (x22)
;
\end{tikzpicture}
\end{equation}
We note the following about the components of $\theta$.
\begin{itemize}
\item While $\txsma_i f_i\upom_i$ and $f(\oplus_i \upom_i)$ are pointed $G$-functors \cref{f_upom}, the pointed functors $\theta_{\angordmdot}$ are \emph{not} required to be $G$-equivariant.  Thus, \cref{ggcat_k1_nat} is not generally a diagram in $\Gcat$.
\item In the $k$-tuple $\angordmdot \in \Gsk^k$, if any $\angordmi = \vstar \in \Gsk$, then
\[\begin{split}
f_i \angordmi &= f_i \vstar = \boldone,\\
f\big(\!\txoplus_{i \in \ufs{k}}\, \angordmi \big) &= f\vstar = \boldone, \andspace \\
\theta_{\angordmdot} &= 1_{\boldone}.
\end{split}\]
\end{itemize}

\parhead{Morphisms}.  For objects $\theta, \ups \in \GGCat(\ang{f_i}_{i}; f)$, a morphism $\Theta \cn \theta \to \ups$, which is also denoted by
\begin{equation}\label{ggcat_k2}
\begin{tikzpicture}[baseline={(a1.base)}]
\def\t{26}
\draw[0cell]
(0,0) node (a1) {\phantom{f}}
(a1)++(2,0) node (a2) {f,}
(a1)++(-.4,0) node (a1') {\ang{f_i}_{i \in \ufs{k}}}
;
\draw[1cell=.9]
(a1) edge[bend left=\t] node {\theta} (a2)
(a1) edge[bend right=\t] node[swap] {\ups} (a2)
;
\draw[2cell]
node[between=a1 and a2 at .4, rotate=-90, 2label={above,\Theta}] {\Rightarrow}
;
\end{tikzpicture}
\end{equation}
consists of $\angordmdot$-component pointed natural transformations
\begin{equation}\label{ggcat_k2_comp}
\begin{tikzpicture}[vcenter]
\def\t{28}
\draw[0cell]
(0,0) node (a1) {\phantom{X}} 
(a1)++(2.5,0) node (a2) {\phantom{X}}
(a1)++(-.7,0) node (a1') {\txsma_{i \in \ufs{k}}\, f_i\angordmi}
(a2)++(.9,0) node (a2') {f\big(\!\txoplus_{i \in \ufs{k}}\, \angordmi \big)}
;
\draw[1cell=.8]
(a1) edge[bend left=\t] node {\theta_{\angordmdot}} (a2)
(a1) edge[bend right=\t] node[swap] {\ups_{\angordmdot}} (a2)
;
\draw[2cell]
node[between=a1 and a2 at .33, rotate=-90, 2label={above,\Theta_{\angordmdot}}] {\Rightarrow}
;
\end{tikzpicture}
\end{equation}
for $\angordmdot \in \Gsk^k$ such that, given morphisms $\ang{\upom^i}_{i \in \ufs{k}}$ as defined in \cref{Gk_mor},  the following two whiskered pointed natural transformations are equal.
\begin{equation}\label{ggcat_k2_modax}
\begin{tikzpicture}[vcenter]
\def\t{25} \def\h{2.5} \def\v{-1.6}
\draw[0cell=.9]
(0,0) node (a1) {\phantom{X}} 
(a1)++(\h,0) node (a2) {\phantom{X}}
(a1)++(-.65,0) node (a1') {\txsma_{i \in \ufs{k}}\, f_i\angordmi}
(a2)++(.8,0) node (a2') {f\big(\!\txoplus_{i \in \ufs{k}}\, \angordmi \big)}
(a1)++(0,\v) node (b1){\phantom{X}} 
(a2)++(0,\v) node (b2) {\phantom{X}}
(b1)++(-.6,0) node (b1') {\txsma_{i \in \ufs{k}}\, f_i\angordni}
(b2)++(.8,0) node (b2') {f\big(\!\txoplus_{i \in \ufs{k}}\, \angordni \big)}
;
\draw[1cell=.8]
(a1) edge[bend left=\t] node[pos=.4] {\theta_{\angordmdot}} (a2)
(a1) edge[bend right=\t] node[swap,pos=.6] {\ups_{\angordmdot}} (a2)
(b1) edge[bend left=\t] node[pos=.4] {\theta_{\angordndot}} (b2)
(b1) edge[bend right=\t] node[swap,pos=.6] {\ups_{\angordndot}} (b2)
(a1) edge[shorten <=.5ex, shorten >=.5ex] node[swap] {\txsma_i f_i \upom^i} (b1)
(a2) edge[shorten <=.5ex, shorten >=.5ex] node {f(\oplus_i\, \upom^i \big)} (b2)
;
\draw[2cell=.9]
node[between=a1 and a2 at .37, rotate=-90, 2label={above,\Theta_{\angordmdot}}] {\Rightarrow}
node[between=b1 and b2 at .37, rotate=-90, 2label={above,\Theta_{\angordndot}}] {\Rightarrow}
;
\end{tikzpicture}
\end{equation}
We note the following about the components of $\Theta$.
\begin{itemize}
\item The pointed natural transformations $\Theta_{\angordmdot}$ are \emph{not} required to be $G$-equivariant.
\item If any $\angordmi = \vstar \in \Gsk$, then $\Theta_{\angordmdot}$ is the identity natural transformation of the identity functor $1_{\boldone}$.
\end{itemize}
Identities and composition of morphisms are given componentwise using \cref{ggcat_k2_comp}.

\parhead{$G$-action functors}.  The $G$-action on the $G$-category $\GGCat(\ang{f_i}_i; f)$ is given componentwise by the conjugation $G$-action.  In other words, for a morphism $\Theta \cn \theta \to \ups$, an element $g \in G$, and an object $\angordmdot \in \Gsk^k$, the $\angordmdot$-components of $g \cdot \theta$, $g \cdot \ups$, and $g \cdot \Theta$ are given by the following composites and whiskering.
\begin{equation}\label{ggcat_k_g}
\begin{tikzpicture}[vcenter]
\def\t{18} \def\h{3.5} \def\v{-1.6}
\draw[0cell=.9]
(0,0) node (a1) {\phantom{X}} 
(a1)++(\h,0) node (a2) {\phantom{X}}
(a1)++(-.65,0) node (a1') {\txsma_{i \in \ufs{k}}\, f_i\angordmi}
(a2)++(.8,0) node (a2') {f\big(\!\txoplus_{i \in \ufs{k}}\, \angordmi \big)}
(a1)++(0,\v) node (b1){\phantom{X}} 
(a2)++(0,\v) node (b2) {\phantom{X}}
(b1)++(-.65,0) node (b1') {\txsma_{i \in \ufs{k}}\, f_i\angordmi}
(b2)++(.8,0) node (b2') {f\big(\!\txoplus_{i \in \ufs{k}}\, \angordmi \big)}
;
\draw[1cell=.75]
(a1) edge[bend left=\t] node[pos=.42] {(g \cdot \theta)_{\angordmdot}} (a2)
(a1) edge[bend right=\t] node[swap,pos=.58] {(g \cdot \ups)_{\angordmdot}} (a2)
(b1) edge[bend left=\t] node[pos=.4] {\theta_{\angordmdot}} (b2)
(b1) edge[bend right=\t] node[swap,pos=.6] {\ups_{\angordmdot}} (b2)
(a1) edge[shorten <=.5ex, shorten >=.5ex] node[swap] {\txsma_i \ginv} (b1)
(b2) edge[shorten <=.5ex, shorten >=.5ex] node[swap] {g} (a2)
;
\draw[2cell=.9]
node[between=a1 and a2 at .3, rotate=-90, 2label={above,(g \cdot \Theta)_{\angordmdot}}] {\Rightarrow}
node[between=b1 and b2 at .4, rotate=-90, 2label={above,\Theta_{\angordmdot}}] {\Rightarrow}
;
\end{tikzpicture}
\end{equation}
The properties \cref{ggcat_k1_nat,ggcat_k2_modax} are satisfied by, respectively, $g \cdot \theta$ and $g \cdot \Theta$ by
\begin{itemize}
\item the same properties for $\theta$ and $\Theta$,
\item the functoriality of $\sma$, and
\item the $G$-equivariance \cref{f_upom} of $f_i\upom_i$ and $f(\oplus_i\, \upom^i)$. 
\end{itemize}
The functoriality of the $g$-action on $\GGCat(\ang{f_i}_{i}; f)$ follows from the fact that, in \cref{ggcat_k_g}, whiskering with $\sma_i\, \ginv$ and $g$ preserves identities and vertical composition of natural transformations.  This finishes the description of the $k$-ary multimorphism $G$-category $\GGCat(\ang{f_i}_{i}; f)$.

\parhead{Units}.  For an object $f \in \GGCat$, the $f$-colored unit \cref{ccoloredunit} is the $G$-functor
\begin{equation}\label{ggcat_funit}
\boldone \fto{1_f} \GGCat(f;f)
\end{equation}
determined by the morphism $1_f \cn f \to f$ whose $\angordm$-component \cref{ggcat_k1_comp} is the identity functor on $f\angordm$ for each $\angordm \in \Gsk$.  Note that $1_f$ is fixed by the $G$-action \cref{ggcat_k_g}.
\end{explanation}

\begin{explanation}[Symmetric Group Action]\label{expl:ggcat_symmetry}
To unravel the symmetric group action on the $\Gcat$-multicategory $\GGCat$, we consider the context of \cref{expl:ggcat_positive}, with objects $\angordmdot \in \Gsk^k$, 
\[\theta,\ups \cn \ang{f_i}_{i \in \ufs{k}} \to f,\]
and a morphism $\Theta \cn \theta \to \ups$, as defined in \cref{angordmdot,ggcat_k1,ggcat_k2} .  For each permutation $\si \in \Si_k$, the images of $\theta$, $\ups$, and $\Theta$ under the right $\si$-action $G$-functor \cref{rightsigmaaction}
\begin{equation}\label{ggcat_sigma}
\GGCat\big(\ang{f_i}_{i \in \ufs{k}}; f \big) \fto[\iso]{\si} \GGCat\big(\ang{f_{\si(i)}}_{i \in \ufs{k}}; f \big)
\end{equation}
have $\angordmdot$-components given by  the following composites and whiskering in $\Catst$.
\begin{equation}\label{ggcat_si_action}
\begin{tikzpicture}[vcenter]
\def\s{25} \def\v{-2} \def\t{.45} \def\u{.55}
\draw[0cell=.9]
(0,0) node (a) {\txsma_{i \in \ufs{k}}\, f_{\si(i)} \angordmi}
(a)++(.8,0) node (a') {\phantom{X}}
(a')++(2.8,0) node (b') {\phantom{X}}
(b')++(.8,0) node (b) {f \big(\! \txoplus_{i \in \ufs{k}}\, \angordmi \big)}
(a')++(0,\v) node (c') {\phantom{X}} 
(c')++(-.93,0) node (c) {\txsma_{i \in \ufs{k}}\, f_i \ang{\ordm^{\sigmainv(i)}}}
(b')++(0,\v) node (d') {\phantom{X}}
(d')++(1.1,0) node (d) {f \big(\! \txoplus_{i \in \ufs{k}}\, \ang{\ordm^{\sigmainv(i)}} \big)}
;
\draw[1cell=.9]
(a') edge[bend left=\s] node[pos=\t] {\theta^\si_{\angordmdot}} (b')
(a') edge[bend right=\s] node[swap,pos=\u] {\ups^\si_{\angordmdot}} (b')
(c') edge[bend left=\s] node[pos=\t] {\theta_{\si\angordmdot}} (d')
(c') edge[bend right=\s] node[swap,pos=\u] {\ups_{\si\angordmdot}} (d')
(a') edge[shorten <=1ex, shorten >=1ex] node[swap] {\si} (c')
(d') edge[shorten <=1ex, shorten >=1ex] node[swap] {f(\sigmainv)} (b')
;
\draw[2cell=.9]
node[between=a' and b' at .34, rotate=-90, 2label={above,\Theta^\si_{\angordmdot}}] {\Rightarrow}
node[between=c' and d' at .34, rotate=-90, 2label={above,\Theta_{\si\angordmdot}}] {\Rightarrow}
;
\end{tikzpicture}
\end{equation}
The arrows in \cref{ggcat_si_action} are given as follows.
\begin{itemize}
\item The left vertical arrow denoted by $\si$ uses the braiding for $(\Catst,\sma)$ to permute the $k$ smash factors according to the permutation $\si \in \Si_k$.
\item Along the bottom, $\theta_{\si\angordmdot}$, $\ups_{\si\angordmdot}$, and $\Theta_{\si\angordmdot}$ are the components of, respectively, $\theta$, $\ups$, and $\Theta$ at the $k$-tuple of objects in $\Gsk$ given by
\begin{equation}\label{si_angordmdot}
\begin{split}
\si\angordmdot &= \bang{\ang{\ordm^{\sigmainv(i)}}}_{i \in \ufs{k}} \\
&= \big(\ang{\ordm^{\sigmainv(1)}}, \ldots, \ang{\ordm^{\sigmainv(k)}} \big).
\end{split}
\end{equation}
\item In the right vertical arrow, $\sigmainv$ is an iterate of the braiding in $\Gsk$ \cref{Gsk_braiding} that permutes the $k$ $\oplus$-factors according to the permutation $\sigmainv \in \Si_k$:
\begin{equation}\label{sigmainv_Gsk}
\txoplus_{i \in \ufs{k}}\, \ang{\ordm^{\sigmainv(i)}} \fto[\iso]{\sigmainv} 
\txoplus_{i \in \ufs{k}}\, \ang{\ordm^i}.
\end{equation}
Since this is an isomorphism in $\Gsk$, the right vertical arrow $f(\sigmainv)$ is a pointed $G$-isomorphism.  To make the iterated braiding $\sigmainv$ more explicit, suppose the object $\angordmi \in \Gsk$ has length $r_i$ for each $i \in \ufs{k}$.  Then $\sigmainv$ in \cref{sigmainv_Gsk} is the isomorphism 
\begin{equation}\label{sigmainv_Gsk_exp}
\big(\sigmainv\ang{r_{\sigmainv(i)}}_{i \in \ufs{k}} \spc \ang{1}\big) \cn 
\txoplus_{i \in \ufs{k}}\, \ang{\ordm^{\sigmainv(i)}} \fiso
\txoplus_{i \in \ufs{k}}\, \ang{\ordm^i}
\end{equation}
in $\Gsk$ \cref{Gsk_morphisms}.
\begin{itemize}
\item Its reindexing injection 
\[\ufs{\txsum_{i \in \ufs{k}}\, r_{\sigmainv(i)}} \fiso \ufs{\txsum_{i \in \ufs{k}}\, r_i}\]
is the block permutation
\[\sigmainv\ang{r_{\sigmainv(i)}}_{i \in \ufs{k}} \in 
\Sigma_{r_{\sigmainv(1)}+\cdots+r_{\sigmainv(k)}} = \Sigma_{r_1 + \cdots + r_k}\]
that permutes $k$ blocks of lengths $r_{\sigmainv(1)}, r_{\sigmainv(2)}, \ldots, r_{\sigmainv(k)}$ according to $\sigmainv \in \Sigma_k$. 
\item Each entry in $\ang{1}$ is an identity function of some pointed finite set in some $\angordmi$.
\end{itemize}
\end{itemize}

\parhead{Functoriality}.  The images $\theta^\si$ and $\Theta^\si$ satisfy, respectively, the properties \cref{ggcat_k1_nat,ggcat_k2_modax}, by 
\begin{itemize}
\item the corresponding properties for $\theta$ and $\Theta$,
\item the naturality of the braidings for $(\Catst, \sma)$ and $(\Gsk,\oplus)$ \cref{Gsk_xi_natural}, and
\item the functoriality of $f$.
\end{itemize}
The functoriality of $\si$ \cref{ggcat_sigma} with respect to $\Theta$ follows from the fact that, in \cref{ggcat_si_action}, whiskering with $\si$ and $f(\sigmainv)$ preserves identities and vertical composition of natural transformations.  

\parhead{$G$-equivariance}.  The following commutative diagram of pointed functors shows that, for each $g \in G$ and $\angordmdot \in \Gsk^k$,  $(g \cdot \theta)^\si$ and $g \cdot \theta^\si$ have the same $\angordmdot$-components \cref{ggcat_k1_comp}.
\begin{equation}\label{ggcat_si_g}
\begin{tikzpicture}[vcenter]
\def\h{2.5} \def\u{-1} \def\v{-1.4}
\draw[0cell=.9]
(0,0) node (a1) {\txsma_i f_{\si(i)} \angordmi}
(a1)++(-\h,\u) node (a21) {\txsma_i f_i \angordmsiinvi}
(a1)++(\h,\u) node (a22) {\txsma_i f_{\si(i)} \angordmi}
(a1)++(0,2*\u) node (a3) {\txsma_i f_i \angordmsiinvi}
(a3)++(0,\v) node (a4) {f\big(\! \txoplus_i \angordmsiinvi\big)}
(a4)++(-\h,\u) node (a51) {f\big(\! \txoplus_i \angordmsiinvi\big)}
(a4)++(\h,\u) node (a52) {f\big(\! \txoplus_i \angordmi\big)}
(a4)++(0,2*\u) node (a6) {f\big(\! \txoplus_i \angordmi\big)}
;
\draw[1cell=.8]
(a1) edge node[swap] {\si} (a21)
(a21) edge node[pos=.6] {\sma_i\, \ginv} (a3)
(a1) edge node {\sma_i\, \ginv} (a22)
(a22) edge node[swap,pos=.5] {\si} (a3)
(a3) edge node {\theta_{\si\angordmdot}} (a4)
(a4) edge node[swap,pos=.5] {g} (a51)
(a51) edge node[swap] {f(\sigmainv)} (a6)
(a4) edge node [pos=.5]{f(\sigmainv)} (a52)
(a52) edge node {g} (a6)
(a21) edge node[swap] {(g \cdot \theta)_{\si\angordmdot}} (a51)
(a22) edge node {\theta^\si_{\angordmdot}} (a52)
;
\end{tikzpicture}
\end{equation}
The following statements hold for the diagram \cref{ggcat_si_g}.
\begin{itemize}
\item By \cref{ggcat_k_g,ggcat_si_action}, the composites along the left and right boundaries in \cref{ggcat_si_g} are, respectively, $(g \cdot \theta)^\si_{\angordmdot}$ and $(g \cdot \theta^\si)_{\angordmdot}$.
\item The top quadrilateral commutes by the naturality of the braiding for $(\Catst, \sma)$.
\item The bottom quadrilateral commutes by the $G$-equivariance \cref{f_upom} of the pointed functor $f(\sigmainv)$.
\item The middle left and right trapezoids are, respectively, \cref{ggcat_k_g,ggcat_si_action}.
\end{itemize}
Thus, $(g \cdot \theta)^\si$ and $g \cdot \theta^\si$ are equal, which proves that the functor $\si$ is $G$-equivariant on objects.  The $G$-equivariance of $\si$ on morphisms is proved by reusing the diagram \cref{ggcat_si_g} and replacing $\theta$ by a morphism $\Theta$ \cref{ggcat_k2}.

\parhead{Axioms}.  The symmetric group action axioms \cref{enr-multicategory-symmetry} hold by \cref{ggcat_si_action} and the functoriality of $f$.

This finishes the description of the right $\sigma$-action $G$-functor on $\GGCat$, along with the symmetric group action axioms.
\end{explanation}

\begin{explanation}[Multicategorical Composition]\label{expl:ggcat_composition}
For $k>0$, $r_i \geq 0$, and pointed functors
\[f_{i\crdot} = \ang{f_{i,j}}_{j \in \ufs{r}_i} \spc f_\crdot = \ang{f_i}_{i \in \ufs{k}} \spc f \cn (\Gsk,\vstar) \to (\Gcatst,\boldone),\]
we unravel the composition $G$-functor \cref{eq:enr-defn-gamma}
\begin{equation}\label{ggcat_gamma}
\GGCat(f_\crdot ; f) \times 
\txprod_{i \in \ufs{k}}\, \GGCat(f_{i\crdot} ; f_i) \fto{\ga} \GGCat(f_{\crdots} ; f)
\end{equation}
where $f_{\crdots} = \ang{f_{i\crdot}}_{i \in \ufs{k}}$.  In the discussion below, if any $r_i = 0$, then $f_{i\crdot}$ is interpreted as the monoidal unit $\gu$, as described in \cref{expl:ggcat_zero}.

\parhead{Objects}.  The $G$-functor $\ga$ sends objects \cref{ggcat_k1}
\begin{equation}\label{theta_thetai}
\big(f_\crdot \fto{\theta} f \spsc \ang{f_{i\crdot} \fto{\theta_i} f_i}_{i \in \ufs{k}} \big)
\end{equation}
to the object
\[f_{\crdots} \fto{\ga(\theta; \theta_\crdot)} f,\]
where $\theta_\crdot = \ang{\theta_i}_{i \in \ufs{k}}$.  To describe its components, we consider the following objects, where $r = \sum_{i \in \ufs{k}}\, r_i$.
\begin{equation}\label{angordmij}
\left\{\begin{aligned}
\angordmij & \in \Gsk & \fm^i &= \txoplus_{j \in \ufs{r}_i}\, \angordmij \in \Gsk\\
\angordmidot &= \ang{\angordmij}_{j \in \ufs{r}_i} \in \Gsk^{r_i} \phantom{M} &
\fmdot &= \ang{\fm^i}_{i \in \ufs{k}} \in \Gsk^k\\
\angordmddot &= \ang{\angordmidot}_{i \in \ufs{k}} \in \Gsk^r &
\fm &= \txoplus_{i \in \ufs{k}}\, \fm^i \in \Gsk
\end{aligned}\right.
\end{equation}
The $\angordmddot$-component pointed functor of $\ga(\theta; \theta_\crdot)$, in the sense of \cref{ggcat_k1_comp}, is the composite
\begin{equation}\label{ggcat_gamma_comp}
\begin{tikzpicture}[baseline={(a.base)}]
\def\u{.7}
\draw[0cell=.9]
(0,0) node (a) {\txsma_{i \in \ufs{k}} \txsma_{j \in \ufs{r}_i} f_{i,j} \angordmij}
(a)++(4.2,0) node (b) {\txsma_{i \in \ufs{k}} f_i \fm^i}
(b)++(2.2,0) node (c) {f \fm}
;
\draw[1cell=.9]
(a) edge node {\sma_{i \in \ufs{k}} \, \theta_{i,\angordmidot}} (b)
(b) edge node {\theta_{\fmdot}} (c)
;
\draw[1cell=.9]
(a) [rounded corners=2pt] |- ($(b)+(-1,\u)$)
-- node[pos=0] {\ga(\theta; \theta_\crdot)_{\angordmddot}} ($(b)+(1,\u)$) -| (c)
;
\end{tikzpicture}
\end{equation}
where $\theta_{i,\angordmidot}$ is the $\angordmidot$-component of $\theta_i$.  The property \cref{ggcat_k1_nat} holds for $\ga(\theta; \theta_\crdot)$ by
\begin{itemize}
\item the same property for $\theta$ and $\theta_i$, and
\item the functoriality of the smash product.
\end{itemize}

If some $r_i = 0$, then 
\[\angordmidot = \fm^i = \ang{} \in \Gsk,\] 
which is the empty sequence, and 
\[\ord{1} = \sma\ang{} = \gu\ang{} \fto{\theta_{i,\ang{}}} f_i\ang{}\]
is the $\ang{}$-component pointed functor \cref{ggcat_zero_obj_comp} of $\theta_i \in \GGCat(\ang{};f_i)$.

\parhead{Morphisms}.
The $G$-functor $\ga$ sends morphisms \cref{ggcat_k2} 
\begin{equation}\label{Theta_Thetai}
\begin{split}
&\big(\theta \fto{\Theta} \ups \spsc \ang{\theta_i \fto{\Theta_i} \ups_i}_{i \in \ufs{k}} \big) \\
& \in  \GGCat(f_\crdot ; f) \times \txprod_{i \in \ufs{k}}\, \GGCat(f_{i\crdot} ; f_i) 
\end{split}
\end{equation}
to the morphism
\[\begin{tikzpicture}[vcenter]
\def\t{27}
\draw[0cell]
(0,0) node (a) {\phantom{f}}
(a)++(3,0) node (b) {f}
(a)++(-.1,0) node (a') {f_{\crdots}}
;
\draw[1cell=.8]
(a) edge[bend left=\t] node {\ga(\theta;\theta_\crdot)} (b)
(a) edge[bend right=\t] node[swap] {\ga(\ups;\ups_\crdot)} (b)
;
\draw[2cell]
node[between=a and b at .28, rotate=-90, 2label={above,\ga(\Theta;\Theta_\crdot)}] {\Rightarrow}
;
\end{tikzpicture}\]
where $\Theta_\crdot = \ang{\Theta_i}_{i \in \ufs{k}}$.  The $\angordmddot$-component pointed natural transformation of $\ga(\Theta; \Theta_\crdot)$, in the sense of \cref{ggcat_k2_comp}, is the horizontal composite
\begin{equation}\label{ggcat_2cell_comp}
\begin{tikzpicture}[baseline={(a.base)}]
\def\u{.7} \def\s{25} \def\t{35}
\draw[0cell=.85]
(0,0) node (a1) {\txsma_{i \in \ufs{k}} \txsma_{j \in \ufs{r}_i} f_{i,j} \angordmij}
(a1)++(1.1,0) node (a2) {\phantom{A}}
(a2)++(3.2,0) node (a3) {\phantom{A}}
(a3)++(.5,0) node (a4) {\txsma_{i \in \ufs{k}} f_i \fm^i}
(a4)++(.5,0) node (a5) {\phantom{A}}
(a5)++(1.8,0) node (a6) {\phantom{A}}
(a6)++(.1,0) node (a7) {f \fm}
;
\draw[1cell=.8]
(a2) edge[bend left=\s] node {\sma_{i \in \ufs{k}} \, \theta_{i,\angordmidot}} (a3)
(a2) edge[bend right=\s] node[swap] {\sma_{i \in \ufs{k}} \, \ups_{i,\angordmidot}} (a3)
(a5) edge[bend left=\t] node {\theta_{\fmdot}} (a6)
(a5) edge[bend right=\t] node[swap] {\ups_{\fmdot}} (a6)
;
\draw[2cell=.9]
node[between=a2 and a3 at .25, rotate=-90, 2label={above, \sma_{i \in \ufs{k}}\, \Theta_{i,\angordmidot}}] {\Rightarrow}
node[between=a5 and a6 at .35, rotate=-90, 2label={above,\Theta_{\fmdot}}] {\Rightarrow}
;
\end{tikzpicture}
\end{equation}
where $\Theta_{i,\angordmidot}$ is the $\angordmidot$-component of $\Theta_i$.  The property \cref{ggcat_k2_modax} holds for $\ga(\Theta; \Theta_\crdot)$ by
\begin{itemize}
\item the same property for $\Theta$ and $\Theta_i$, and
\item the 2-functoriality of the smash product.
\end{itemize}

\parhead{$G$-Functoriality}.  The functoriality of $\ga$ with respect to $(\Theta; \Theta_\crdot)$ follows from \cref{ggcat_2cell_comp}, the 2-functoriality of $\sma$, and the fact that $\Catst$ is a 2-category.  The functor $\ga$ is $G$-equivariant by
\begin{itemize}
\item \cref{ggcat_k_g}, which defines the componentwise conjugation $G$-action, 
\item the definitions \cref{ggcat_gamma_comp,ggcat_2cell_comp} of $\ga$, and
\item the functoriality of the smash product.
\end{itemize}

\parhead{Unity and associativity}.  The axioms \cref{enr-multicategory-right-unity,enr-multicategory-left-unity,enr-multicategory-associativity} hold for $\GGCat$ by
\begin{itemize}
\item the definitions \cref{ggcat_funit,ggcat_gamma_comp,ggcat_2cell_comp} of the units and $\ga$;
\item the unity and associativity of functor composition and horizontal composition of natural transformations; and
\item the 2-functoriality of the smash product.
\end{itemize}
This finishes the description of the composition $G$-functor $\ga$ on $\GGCat$ \cref{ggcat_gamma}, along with the unity and associativity axioms.
\end{explanation}

\begin{explanation}[Top Equivariance]\label{expl:ggcat_topeq}
We provide a self-contained proof of the top equivariance axiom \cref{enr-operadic-eq-1} for the $\Gcat$-multicategory $\GGCat$.

Using the notation in \cref{expl:ggcat_positive,expl:ggcat_symmetry,expl:ggcat_composition}, the axiom \cref{enr-operadic-eq-1} asks for the commutativity of the following diagram of $G$-functors for each permutation $\si \in \Sigma_k$, where $\sigmabar \in \Si_r$ denotes the block permutation induced by $\si$ that permutes blocks of lengths $r_{\si(1)} , r_{\si(2)}, \ldots, r_{\si(k)}$.
\begin{equation}\label{ggcat_topeq}
\begin{tikzpicture}[vcenter]
\def\v{-1.5}
\draw[0cell=.85]
(0,0) node (a1) {\GGCat(f_\crdot; f) \times \txprod_{i \in \ufs{k}} \,\GGCat(f_{i\crdot}; f_i)}
(a1)++(5,0) node (a2) {\GGCat(f_{\crdots}; f)}
(a1)++(0,\v) node (b1) {\GGCat(f_\crdot \si; f) \times \txprod_{i\in \ufs{k}} \,\GGCat(f_{\si(i)\crdot} ; f_{\si(i)})}
(a2)++(0,\v) node (b2) {\GGCat(f_{\crdots} \sigmabar ; f)}
;
\draw[1cell=.85]
(a1) edge node {\ga} (a2)
(b1) edge node {\ga} (b2)
(a1) edge[transform canvas={xshift=1em}] node[swap] {\si \times \sigmainv} (b1)
(a2) edge node {\sigmabar} (b2)
;
\end{tikzpicture}
\end{equation}
To check the commutativity of \cref{ggcat_topeq} on objects, we consider objects
\[\begin{split}
& \big(f_\crdot \fto{\theta} f \spsc \ang{f_{i\crdot} \fto{\theta_i} f_i}_{i \in \ufs{k}} \big) \andspace\\
& \angordmij \in \Gsk \forspace (i,j) \in \ufs{k} \times \ufs{r}_i,
\end{split}\]
as defined in \cref{theta_thetai,angordmij}.  By \cref{ggcat_si_action,ggcat_gamma_comp}, the left-bottom and top-right composites in the diagram \cref{ggcat_topeq} yield, respectively, the top-left-bottom and right-bottom boundary composites of the following diagram of pointed functors.
\begin{equation}\label{ggcat_topeq_comp}
\begin{tikzpicture}[vcenter]
\def\h{5.3} \def\v{-1.5}
\draw[0cell=.8]
(0,0) node (a1) {\txsma_{i \in \ufs{k}}\, f_{\si(i)} \fm^{\si(i)}}
(a1)++(\h,0) node (a2) {\txsma_{i \in \ufs{k}}\, \txsma_{j \in \ufs{r}_{\si(i)}} f_{\si(i),j} \angordmsiij}
(a1)++(0,\v) node (b1) {\txsma_{i \in \ufs{k}}\, f_i \fm^i}
(a2)++(0,\v) node (b2) {\txsma_{i \in \ufs{k}}\, \txsma_{j \in \ufs{r}_i} f_{i,j} \angordmij}
(b1)++(\h/3,\v) node (c1) {f\fm}
(b2)++(0,\v) node (c2) {f\big(\! \txoplus_{i \in \ufs{k}}\, \fm^{\si(i)}\big)}
;
\draw[1cell=.8]
(a2) edge node[swap] {\sma_{i \in \ufs{k}}\, \theta_{\si(i), \angordmsiidot}} (a1)
(a1) edge node[swap] {\si} (b1)
(b1) edge node[swap,pos=.4] {\theta_{\fmdot}} (c1)
(c1) edge node[pos=.7] {f(\sigmabar^{-1})} (c2)
(a2) edge node {\sigmabar} (b2)
(b2) edge node[swap] {\sma_{i \in \ufs{k}}\, \theta_{i,\angordmidot}} (b1)
(b2) edge node[swap,pos=.5] {\ga(\theta; \theta_\crdot)_{\angordmddot}} (c1)
;
\end{tikzpicture}
\end{equation}
In \cref{ggcat_topeq_comp}, the top rectangle commutes by the naturality of the braiding for $(\Catst,\sma)$.  The lower triangle commutes by the definition \cref{ggcat_gamma_comp} of $\ga$.  This proves the commutativity of the diagram \cref{ggcat_topeq} on objects.  

The commutativity of the diagram \cref{ggcat_topeq} on morphisms is proved by 
\begin{itemize}
\item reusing the diagram \cref{ggcat_topeq_comp},
\item replacing the objects $(\theta; \theta_\crdot)$ by the morphisms $(\Theta; \Theta_\crdot)$ \cref{Theta_Thetai}, and
\item using \cref{ggcat_2cell_comp} and the 2-naturality of the braiding for $(\Catst,\sma)$.
\end{itemize}
Thus, the top equivariance diagram \cref{ggcat_topeq} for $\GGCat$ commutes.
\end{explanation}

\begin{explanation}[Bottom Equivariance]\label{expl:ggcat_boteq}
We provide a self-contained proof of the bottom equivariance axiom \cref{enr-operadic-eq-2} for the $\Gcat$-multicategory $\GGCat$.

Using the notation in \cref{expl:ggcat_positive,expl:ggcat_symmetry,expl:ggcat_composition}, the axiom \cref{enr-operadic-eq-2} asks for the commutativity of the following diagram of $G$-functors for permutations $\ang{\tau_i \in \Si_{r_i}}_{i \in \ufs{k}}$, where $\tautimes$ denotes the block sum $\tau_1 \times \cdots \times \tau_k \in \Si_r$.
\begin{equation}\label{ggcat_boteq}
\begin{tikzpicture}[vcenter]
\def\v{-1.5}
\draw[0cell=.85]
(0,0) node (a1) {\GGCat(f_\crdot; f) \times \txprod_{i \in \ufs{k}}\,\GGCat(f_{i\crdot}; f_i)}
(a1)++(5,0) node (a2) {\GGCat(f_{\crdots}; f)}
(a1)++(0,\v) node (b1) {\GGCat(f_\crdot; f) \times \txprod_{i \in \ufs{k}}\, \GGCat(f_{i\crdot}\tau_i ; f_{i})}
(a2)++(0,\v) node (b2) {\GGCat(f_{\crdots} \tautimes ; f)}
;
\draw[1cell=.85]
(a1) edge node {\ga} (a2)
(b1) edge node {\ga} (b2)
(a1) edge[transform canvas={xshift=2em}] node[swap] {1 \times \txprod_{i \in \ufs{k}} \, \tau_i} (b1)
(a2) edge node {\tautimes} (b2)
;
\end{tikzpicture}
\end{equation}
To check the commutativity of \cref{ggcat_boteq} on objects, we use the objects $(\theta; \theta_\crdot)$ and $\angordmij$ defined in \cref{theta_thetai,angordmij}, along with the following notation.  
\[\begin{aligned}
\tau_i\angordmidot &= \bang{\angordmitauiinvj}_{j \in \ufs{r}_i} \in \Gsk^{r_i} \phantom{M}
& \tau_i \fm^i &= \txoplus_{j \in \ufs{r}_i}\, \angordmitauiinvj \in \Gsk \\
\tau\fmdot &= \ang{\tau_i \fm^i}_{i \in \ufs{k}} \in \Gsk^k 
& \tau\fm &= \txoplus_{i \in \ufs{k}}\, \tau_i \fm^i \in \Gsk
\end{aligned}\]
By \cref{ggcat_si_action,ggcat_gamma_comp}, the left-bottom and top-right composites in the diagram \cref{ggcat_boteq} yield, respectively, the bottom and left-top-right boundary composites of the following diagram of pointed functors.
\begin{equation}\label{ggcat_boteq_comp}
\begin{tikzpicture}[vcenter]
\def\h{3} \def\v{0} \def\u{1.5} \def\t{0}
\draw[0cell=.8]
(0,0) node (a1) {\txsma_{i \in \ufs{k}} \txsma_{j \in \ufs{r}_i} f_{i,\tau_i(j)} \angordmij}
(a1)++(4.3,\v) node (a2) {\txsma_{i \in \ufs{k}}\, f_i \fm^i}
(a2)++(2.5,-\v) node (a3) {f\fm}
(a1)++(0,\u) node (b1) {\txsma_{i \in \ufs{k}} \txsma_{j \in \ufs{r}_i} f_{i,j} \angordmitauiinvj}
(a2)++(0,\u) node (b2) {\txsma_{i \in \ufs{k}}\, f_i (\tau_i \fm^i)}
(a3)++(0,\u) node (b3) {f(\tau\fm)}
;
\draw[1cell=.8]
(a1) edge node[pos=.5] {\sma_i\, \theta_{i,\angordmidot}^{\tau_i}} (a2)
(a2) edge node[pos=.5] {\theta_{\fmdot}} (a3)
(a1) edge[transform canvas={xshift=0em}] node {\sma_i \tau_i} (b1)
(b1) edge[bend left=\t] node {\sma_i\, \theta_{i,\tau_i\angordmidot}} (b2)
(b2) edge[bend left=\t] node {\theta_{\tau\fmdot}} (b3)
(b3) edge[transform canvas={xshift=-0em}] node[swap,pos=.4] {f(\oplus_i \tauinv_i)} (a3)
(b2) edge node[swap,pos=.4] {\sma_i\, f_i(\tauinv_i)} (a2)
;
\end{tikzpicture}
\end{equation}
In \cref{ggcat_boteq_comp}, the left region commutes by
\begin{itemize}
\item the functoriality of the smash product on $\Catst$ and
\item the definition \cref{ggcat_si_action} of $\theta_i^{\tau_i}$.
\end{itemize}  
The right region commutes by the property \cref{ggcat_k1_nat} of $\theta$.  This proves the commutativity of the diagram \cref{ggcat_boteq} on objects.  

The commutativity of the diagram \cref{ggcat_boteq} on morphisms is proved by
\begin{itemize}
\item reusing the diagram \cref{ggcat_boteq_comp},
\item replacing the objects $(\theta; \theta_\crdot)$ by the morphisms $(\Theta; \Theta_\crdot)$ \cref{Theta_Thetai}, and
\item using the 2-functoriality of $\sma$, \cref{ggcat_k2_modax,ggcat_si_action,ggcat_2cell_comp} for morphisms.
\end{itemize}
Thus, the bottom equivariance diagram \cref{ggcat_boteq} for $\GGCat$ commutes.
\end{explanation}

\chapter{$\Gcat$-Multifunctors from Operadic Pseudoalgebras\\ to $\Gskg$-Categories}
\label{ch:jemg}
This chapter constructs a $\Gcat$-multifunctor (\cref{def:enr-multicategory-functor})
\[\MultpsO \fto{\Jgo} \GGCat,\]
called \emph{multifunctorial $J$-theory}, from the $\Gcat$-multicategory $\MultpsO$ of $\Op$-pseudoalgebras \pcref{thm:multpso} to the $\Gcat$-multicategory $\GGCat$ of $\Gskg$-categories \pcref{expl:ggcat_gcatenr}.  It is the first, and the most nontrivial, step of our $G$-equivariant algebraic $K$-theory multifunctor from $\Op$-pseudoalgebras to orthogonal $G$-spectra.  There is also a strong variant
\[\MultpspsO \fto{\Jgosg} \GGCat,\]
called \emph{multifunctorial strong $J$-theory}.  See \cref{thm:Jgo_multifunctor}.  At the object level, each of $\Jgo$ and $\Jgosg$ sends $\Op$-pseudoalgebras \pcref{def:pseudoalgebra} to $\Gskg$-categories \pcref{expl:ggcat_obj}.

While the $\Gcat$-multicategory $\MultpsO$ exists whenever $\Op$ is a pseudo-commutative operad in $\Gcat$ \pcref{def:multicatO}, the $\Gcat$-multifunctor $\Jgo$ also requires $\Op(1)$ to be a terminal $G$-category.  \cref{as:OpA} below is in effect throughout this chapter.

\begin{assumption}[$\Tinf$-Operads]\label{as:OpA}
We assume that $(\Op,\ga,\opu,\pcom)$ is a pseudo-commutative operad in $\Gcat$ for an arbitrary group $G$ \pcref{def:GCat,def:pseudocom_operad,def:enr-multicategory} such that $\Op(1)$ is a terminal $G$-category.  Such a $\Gcat$-operad is called a \index{T-infinity operad@$\Tinf$-operad}\index{operad!T-infinity@$\Tinf$}\emph{$\Tinf$-operad}.  
\end{assumption}

A pseudo-commutative operad $\Op$ is already assumed to be reduced, which means that $\Op(0)$ is a terminal $G$-category.  Thus, a $\Tinf$-operad is precisely a pseudo-commutative operad in $\Gcat$ that is also 1-connected in the sense that
\begin{equation}\label{as:OpA_i}
\Op(0) = \boldone = \Op(1),
\end{equation}
where $\boldone$ is a terminal $G$-category with only one object and its identity morphism.  For a $\Tinf$-operad $\Op$, \cref{thm:multpso} yields a $\Gcat$-multicategory $\MultpsO$, which is the domain of $\Jgo$.  \cref{expl:nsystem_assumptions} discusses how \cref{as:OpA} is used in the construction of $\Jgo$.  

\subsection*{Application to Symmetric Monoidal $G$-Categories}
By \cref{BE_pseudocom,GBE_pseudocom}, each of the Barratt-Eccles $\Gcat$-operad $\BE$ (\cref{def:BE-Gcat}) and the $G$-Barratt-Eccles operad $\GBE$ (\cref{def:GBE}) admits a unique pseudo-commutative structure.  Moreover, the $G$-categories
\[\BE(1) = \ESigma_1 = \boldone \andspace \GBE(1) = \Catg(\EG,\ESigma_1)\]
are both terminal.  Thus, $\BE$ and $\GBE$ are $\Tinf$-operads.  By \cref{thm:Jgo_multifunctor}, there are $\Gcat$-multifunctors
\[\begin{split}
\MultpsBE \fto{\Jgbe} \GGCat,\\
\MultpsGBE \fto{\Jggbe} \GGCat,
\end{split}\]
and their strong variants.  At the object level, the $\Gcat$-multifunctors $\Jgbe$ and $\Jggbe$ send, respectively, naive and genuine symmetric monoidal $G$-categories to $\Gskg$-categories.  See \cref{ex:JgBE} for more details.

\subsection*{Multifunctorial $J$-Theory Preserves Equivariant $\Einf$-Algebras}
Since $\Jgo$ is a $\Gcat$-multifunctor (\cref{def:enr-multicategory-functor}), it preserves any algebraic structure parametrized by a $\Gcat$-multifunctor.  This means that, given any $\Gcat$-multifunctor $f \cn \cQ \to \MultpsO$ from a $\Gcat$-multicategory $\cQ$, the composite
\[\cQ \fto{f} \MultpsO \fto{\Jgo} \GGCat\]
is also a $\Gcat$-multifunctor.  In particular, for any $\Tinf$-operad $\Op$, $\Jgo$ preserves (i) equivariant $\Einf$-algebras in the sense of Guillou-May \cite{gm17} and (ii) $\Ninf$-algebras in the sense of Blumberg-Hill \cite{blumberg_hill}, where $G$ is assumed to be a finite group.  See \cref{thm:Jgo_preservation,Jgo_preserves_Einf}.
\begin{description}
\item[$\Einf$] For a $G$-categorical $\Einf$-operad $\cQ$ \pcref{def:Einfty_operads}, a $\Gcat$-multifunctor $f \cn \cQ \to \MultpsO$ is an equivariant $\Einf$-algebra in the $\Gcat$-multicategory $\MultpsO$ of $\Op$-pseudoalgebras.  The composite $\Jgo \circ f$ is an equivariant $\Einf$-algebra in the $\Gcat$-multicategory $\GGCat$ of $\Gskg$-categories.   
\item[$\Ninf$] For a $G$-categorical $\Ninf$-operad $\cQ$ \pcref{def:Ninfty_operads}, a $\Gcat$-multifunctor $f \cn \cQ \to \MultpsO$ is an $\Ninf$-algebra in $\MultpsO$.  The composite $\Jgo \circ f$ is an $\Ninf$-algebra in $\GGCat$.   
\end{description}
When these examples are applied to the $\Tinf$-operad $\Op = \GBE$ \pcref{as:OpA}, we conclude that $\Jggbe$ sends equivariant $\Einf$-algebras in the $\Gcat$-multicategory $\MultpsGBE$ of genuine symmetric monoidal $G$-categories to equivariant $\Einf$-algebras in $\GGCat$, and similarly for $\Ninf$-algebras.  See \cref{ex:Jgo_preservation} for more details.

\subsection*{Pseudoalgebras to Functors without Strictification}  An object in the domain $\MultpsO$ of $\Jgo$ is an $\Op$-\emph{pseudoalgebra} $(\A,\gaA,\phiA)$ \pcref{def:pseudoalgebra}, whose compatibility with the operadic composition of $\Op$ is controlled by the associativity constraint $\phiA$ \cref{phiA}.  On the other hand, an object in the codomain $\GGCat$ of $\Jgo$ is a $\Gskg$-category \pcref{expl:ggcat_obj}, which means a pointed functor 
\[(\Gsk,\vstar) \fto{f} (\Gcatst, \boldone)\]
in the usual 1-categorical sense.  We emphasize that $\Jgo$ does \emph{not} involve any strictification functor from $\Op$-pseudoalgebras to strict $\Op$-algebras, or from pseudofunctors $\Gsk \to \Gcatst$ to functors.  Instead, the key to $\Jgo$ turning $\Op$-pseudoalgebras into $\Gskg$-categories is our concept of \emph{$\angordn$-systems in $\A$} \pcref{def:nsystem}.  See \cref{rk:no_strictification} for more discussion.

\summary
The following table summarizes the $\Gcat$-multifunctor $\Jgo$, along with its domain and codomain $\Gcat$-multicategories.
\begin{center}
\resizebox{\columnwidth}{!}{%
{\renewcommand{\arraystretch}{1.4}%
{\setlength{\tabcolsep}{1ex}
\begin{tabular}{c|cr|cr|cr}
& $\MultpsO$ & \eqref{thm:multpso} & $\GGCat$ & \eqref{expl:ggcat_gcatenr} & $\Jgo$ & \eqref{thm:Jgo_multifunctor} \\ \hline
objects & $\Op$-pseudoalgebras & \eqref{def:pseudoalgebra} & $\Gskg$-categories & \eqref{expl:ggcat_obj} & $\Jgo \A = \Adash$ & \eqref{def:Aangordn_gcat} \\
$k$-ary 1-cells & $k$-lax $\Op$-morphisms & \eqref{def:k_laxmorphism} & \multicolumn{2}{c|}{\eqref{ggcat_zero_obj_comp}, \eqref{ggcat_k1}} & \multicolumn{2}{c}{\eqref{Jgoa}, \eqref{Jgo_f}} \\
$k$-ary 2-cells & $k$-ary $\Op$-transformations & \eqref{def:kary_transformation} & \multicolumn{2}{c|}{\eqref{ggcat_zero_mor_comp}, \eqref{ggcat_k2}} & \multicolumn{2}{c}{\eqref{Jgotha}, \eqref{Jgotheta}} \\
$G$-action & \multicolumn{2}{c|}{\eqref{def:k_laxmorphism_g}, \eqref{expl:O_tr_g}} & \multicolumn{2}{c|}{\eqref{ggcat_zero_mor_g}, \eqref{ggcat_k_g}} & \multicolumn{2}{c}{\eqref{def:Jgo_zero}, \eqref{Jgo_k_Geq}, \eqref{def:Jgo_pos_mor}} \\
units & \multicolumn{2}{c|}{\eqref{Opsalg_unit}} & \multicolumn{2}{c|}{\eqref{ggcat_funit}} & \multicolumn{2}{c}{\eqref{Jgo_unity}} \\
symmetry & \multicolumn{2}{c|}{\eqref{multpso_sym_functor}} &  \multicolumn{2}{c|}{\eqref{expl:ggcat_symmetry}} & \multicolumn{2}{c}{\eqref{Jgo_sigma}} \\
composition & \multicolumn{2}{c|}{\eqref{def:gam_functor}} & \multicolumn{2}{c|}{\eqref{expl:ggcat_composition}}  & \multicolumn{2}{c}{\eqref{Jgo_gamma}}\\
\end{tabular}}}}
\end{center}

\connection
The domain and codomain of $\Jgo$ are, respectively, the $\Gcat$-multicategories $\MultpsO$ and $\GGCat$ constructed in \cref{ch:multpso,ch:ggcat}.  Two key examples of $\Tinf$-operads, the Barratt-Eccles $\Gcat$-operad $\BE$ and the $G$-Barratt-Eccles operad $\GBE$, are discussed in \cref{ch:psalg}.  The $\Gcat$-multifunctor $\Jgo$ is used in \cref{sec:Kgo_multi} to construct our $G$-equivariant algebraic $K$-theory multifunctor $\Kgo$.

\organization
To construct the $\Gcat$-multifunctor $\Jgo$, we need to
\begin{itemize}
\item define its assignments on objects, $k$-ary 1-cells, and $k$-ary 2-cells for $k \geq 0$, and 
\item prove that these assignments are well defined and satisfy the unity, composition, and symmetric group action axioms in \cref{def:enr-multicategory-functor}.
\end{itemize}  
This chapter first discusses the object assignment (\cref{sec:jemg_objects,sec:jemg_morphisms,sec:jemg_morphisms_ii}), followed by the multimorphism assignment (\cref{sec:jemg_zero,sec:jemg_pos_i,sec:jemg_pos_i_proof,sec:jemg_pos_ii}), the $\Gcat$-multifunctor axioms (\cref{sec:jemg_axioms}), and applications \pcref{sec:Jgo_preserves}.  The main definitions deal with multifunctorial $J$-theory, $\Jgo$, and point out the necessary restrictions for the strong variant, $\Jgosg$.  This chapter consists of the following sections.

\secname{sec:jemg_objects}  The object assignment of $\Jgo$ sends each $\Op$-pseudoalgebra $\A$ to a $\Gskg$-category, which means a pointed functor
\[(\Gsk,\vstar) \fto{\Jgo\A = \Adash} (\Gcatst,\boldone).\]
This section defines the object assignment of the functor $\Adash$, which sends each object $\angordn \in \Gsk$ to a small pointed $G$-category $\Aangordn$.  Its objects are called \emph{$\angordn$-systems} in $\A$ (\cref{def:nsystem,def:nsystem_morphism,def:Aangordn_system,def:Aangordn_gcat}).  Each $\angordn$-system consists of
\begin{itemize}
\item a collection of component objects in $\A$, parametrized by tuples of subsets of $\ufs{n}_j$, and
\item \emph{gluing morphisms} among component objects, parametrized by objects of $\Op$ and suitable partitions of subsets.
\end{itemize}
A \emph{strong $\angordn$-system} is an $\angordn$-system with invertible gluing morphisms.  

\secname{sec:jemg_morphisms}  This section constructs the morphism assignment of the functor $\Jgo\A = \Adash$, whose object assignment is given in \cref{sec:jemg_objects}.  This morphism assignment sends each morphism $\upom \cn \angordm \to \angordn$ in $\Gsk$ to a pointed $G$-functor (\cref{def:ftil_functor,def:psitil_functor,def:Afangpsi})
\[\Aangordm \fto{\sys{(\Jgo\A)}{(\upom)} = \Aupom} \Aangordn.\]
On objects, $\Aupom$ sends an $\angordm$-system to an $\angordn$-system.  

\secname{sec:jemg_morphisms_ii}  Having defined the object and morphism assignments of $\Jgo\A = \Adash$, this section proves that $\Adash$ is actually a $\Gskg$-category, that is, a pointed functor.  See \cref{A_ptfunctor}.  This finishes the construction of the object assignment, $\A \mapsto \Jgo\A$, of the $\Gcat$-multifunctor $\Jgo$.

\secname{sec:jemg_zero}  This section defines the multimorphism $G$-functor in arity 0 (\cref{def:Jgo_zero})
\[\MultpsO(\ang{} \sscs \A) = \A \fto{\Jgo} \GGCat(\ang{} \sscs \Adash).\]
The remainder of this section proves that this 0-ary multimorphism $G$-functor is well defined.  See \cref{Jgoa_nbetad_welldef,Jgoa_nat,Jgotha_gmod}.

\secname{sec:jemg_pos_i}  This section constructs the $G$-equivariant object assignment of the $k$-ary multimorphism $G$-functor 
\[\MultpsO(\ang{\A_i}_{i \in \ufs{k}} \sscs \B) \fto{\Jgo} 
\GGCat(\ang{\Aidash}_{i \in \ufs{k}} \sscs \Bdash)\]
for $k \geq 1$.  See \cref{def:Jgo_pos_obj}.  Several statements used in this definition are proved in the following section.

\secname{sec:jemg_pos_i_proof}  This section proves \cref{Jgof_obj_welldef,Jgof_natural,Jgo_k_Geq}, which are used in \cref{def:Jgo_pos_obj}.  The most intense part of this section is \cref{Jgof_obj_welldef}, which proves that each component of $\Jgo f$ has a well-defined object assignment.  The verification of the associativity axiom and the commutativity axiom involves some nontrivial diagrams.

\secname{sec:jemg_pos_ii}  This section constructs the $G$-equivariant morphism assignment for the $k$-ary multimorphism $G$-functor of $\Jgo$, whose object assignment is constructed in \cref{sec:jemg_pos_i,sec:jemg_pos_i_proof}.  See \cref{def:Jgo_pos_mor}.  The second half of this section proves \cref{Jgotheta_system_mor,Jgotheta_m_natural,Jgotheta_modification}, which are used in \cref{def:Jgo_pos_mor}.  This completes the construction of the $k$-ary multimorphism $G$-functors of $\Jgo$ for $k \geq 0$.

\secname{sec:jemg_axioms}  This section defines the $J$-theory $\Gcat$-multifunctor $\Jgo$ (\cref{def:Jgo_multifunctor}) and verifies the $\Gcat$-multifunctor axioms.  The symmetric group action axiom, the composition axiom, and the unity axiom are proved in, respectively, \cref{Jgo_sigma}, \cref{Jgo_gamma}, and \cref{thm:Jgo_multifunctor}.  \cref{ex:JgBE} is an application of \cref{thm:Jgo_multifunctor} to the Barratt-Eccles $\Gcat$-operad $\BE$ and the $G$-Barratt-Eccles operad $\GBE$.  For example, $\Jggbe$ sends each genuine symmetric monoidal $G$-category to a $\Gskg$-category, in a way that respects the $\Gcat$-multicategory structures of $\MultpsGBE$ and $\GGCat$.

\secname{sec:Jgo_preserves}  The main result of this section is \cref{thm:Jgo_preservation}, which states that $\Jgo$ preserves equivariant algebraic structures parametrized by $\Gcat$-multifunctors from $\Gcat$-multicategories.  In particular, $\Jgo$ preserves equivariant $\Einf$-algebras and $\Ninf$-algebras \pcref{Jgo_preserves_Einf}.  \cref{ex:Jgo_preservation} applies \cref{thm:Jgo_preservation} to $\GBE$ and $\BE$.  One conclusion is that multifunctorial $J$-theory $\Jggbe$ for $\GBE$ sends equivariant $\Einf$-algebras and $\Ninf$-algebras in $\MultpsGBE$, whose objects are genuine symmetric monoidal $G$-categories, to, respectively, equivariant $\Einf$-algebras and $\Ninf$-algebras in $\GGCat$.

\section{$\Gskg$-Categories from Operadic Pseudoalgebras: Objects}
\label{sec:jemg_objects}

Throughout this chapter, we let $(\A,\gaA,\phiA)$ denote an arbitrary $\Op$-pseudoalgebra (\cref{def:pseudoalgebra}).  Under \cref{as:OpA}, this section begins the construction of the $\Gcat$-multifunctor
\[\MultpsO \fto{\Jgo} \GGCat\]
by constructing the object assignment of the $\Gskg$-category (\cref{expl:ggcat_obj})
\[(\Gsk,\vstar) \fto{\Jgo\A} (\Gcatst,\boldone).\]
In other words, for each object $\angordn \in \Gsk$, this section constructs a small pointed $G$-category
\[\sys{(\Jgo\A)}{\angordn} = \Aangordn\]
such that $\sys{\A}{\vstar} = \boldone$.  \cref{sec:jemg_morphisms,sec:jemg_morphisms_ii} construct the morphism assignment of the $\Gskg$-category $\Jgo\A$.

\secoutline

The small pointed $G$-category $\Aangordn$ is constructed in several steps.
\begin{itemize}
\item After establishing some notation, \cref{def:nsystem} defines \emph{$\angordn$-systems in $\A$} for each object $\angordn \in \Gsk$.  An $\angordn$-system in $\A$ consists of a collection of objects in $\A$ indexed by subsets of $\ufs{n}_j$, together with gluing morphisms relating these objects.  
\item \cref{expl:nsystem_assumptions,expl:nsystem_base,expl:nsystem} further elaborate \cref{def:nsystem} of $\angordn$-systems.
\item \cref{def:nsystem_morphism} defines morphisms and the category of $\angordn$-systems $\Aangordn$.  There is also a strong variant, denoted $\syssg{\A}{\ordtu{n}}$, in which the gluing morphisms are isomorphisms.
\item \cref{def:Aangordn_system} defines, for each object $\angordn \in \Gsk$, the pointed categories $\Aangordn$ and $\Aangordnsg$.  \cref{expl:Aangordn_system,ex:nsystem_ones} further elaborate this definition.
\item \cref{def:Aangordn_gcat} extends the pointed categories $\Aangordn$ and $\Aangordnsg$ to pointed $G$-categories with nontrivial $G$-actions.
\end{itemize}

\subsection*{Categories of Systems}

Recall the following notation.
\begin{itemize}
\item $\ufs{m} = \{1 < 2 <\cdots <m\}$ denotes an unpointed finite set for $m \geq 0$ \cref{ufsn}.
\item $\ord{m} = \{0 < 1 < \cdots < m\}$ denotes a pointed finite set \cref{ordn}.  Thus, there is an equality
\[\ufs{m} = \ord{m} \setminus \{0\}.\]
\item $\Fsk$ denotes the small pointed category of pointed finite sets $\ordn$ for $n \geq 0$ and pointed functions (\cref{def:Fsk}).  
\end{itemize}
We use the following notation to denote substitution of entries in a tuple and partitions of sets.

\begin{notation}[Substitution and Partition]\label{not:compk}
For $m > 0$, suppose $\ang{s}  = \ang{s_i}_{i \in \ufs{m}}$ is an $m$-tuple of symbols. 
\begin{itemize}
\item For any $k \in \ufs{m}$ and any symbol $t$, we denote by\index{substitution}
\begin{equation}\label{compk}
\ang{s} \compk t = \big(s_1, \ldots, s_{k-1}, t, s_{k+1}, \ldots, s_m\big)
\end{equation}
the $m$-tuple obtained from $\ang{s}$ by replacing its $k$-th entry, $s_k$, by $t$. 
\item For $k \neq \ell \in \ufs{m}$ and a pair of symbols $(t,u)$, we denote by
\begin{equation}\label{compkl}
\ang{s} \compk t \compell u = \big(\ang{s} \compk t\big) \compell u = \big(\ang{s} \compell u\big) \compk t
\end{equation}
the $m$-tuple obtained from $\ang{s}$ by replacing (i) its $k$-th entry, $s_k$, by $t$ and (ii) its $\ell$-th entry, $s_\ell$, by $u$.
\end{itemize}
For a set $S$, a \emph{partition}\index{partition} of $S$ means a $p$-tuple $\ang{S_j}_{j \in \ufs{p}}$ for some $p \geq 0$ such that the following three statements hold.
\begin{itemize}
\item Each $S_j$ is a possibly-empty subset of $S$.
\item These subsets are pairwise disjoint, which means that 
\[S_i \cap S_j = \emptyset \forspace i \neq j \in \ufs{p}.\]
\item The disjoint union of $\ang{S_j}_{j \in \ufs{p}}$ is $S$.  
\end{itemize}
Such a partition of $S$ is usually denoted by
\begin{equation}\label{partition}
S = \coprod_{j \in \ufs{p}}\, S_j.
\end{equation}
The case $p=0$ can only happen if $S$ is empty.
\end{notation}

Now we define a category $\Aangordn$ associated to each $\Op$-pseudoalgebra $\A$ (\cref{def:pseudoalgebra}) and each object $\angordn \in \Gskel$ with length $q > 0$ (\cref{def:Gsk}).  \cref{def:nsystem} defines the objects of this category; its morphisms are defined in \cref{def:nsystem_morphism}.

\begin{definition}[Systems]\label{def:nsystem}
Suppose $(\Op,\ga,\opu,\pcom)$ is a $\Tinf$-operad \pcref{as:OpA}, and $(\A,\gaA,\phiA)$ is an $\Op$-pseudoalgebra.  Suppose 
\[\angordn = \ang{\ord{n}_j}_{j \in \ufs{q}} = \big(\ord{n}_1, \ldots, \ord{n}_q\big) \in \Gskel\]
is an object of length $q > 0$, with each $\ord{n}_j \in \Fskel$ a pointed finite set.  An  \emph{$\angordn$-system in $\A$}\index{system} is defined to be a pair
\begin{equation}\label{nsystem}
(a,\glu)
\end{equation}
consisting of the following data.
\begin{description}
\item[Component objects]
For each $q$-tuple of subsets
\begin{equation}\label{marker}
\ang{s} = \ang{s_j}_{j \in \ufs{q}} = \ang{s_j \subseteq \ufs{n}_j}_{j \in \ufs{q}}
\end{equation}
with $\ufs{n}_j = \ord{n}_j \setminus \{0\}$, $(a,\glu)$ consists of an \emph{$\ang{s}$-component object}\index{component object}
\begin{equation}\label{a_angs}
a_{\ang{s}} \in \A.
\end{equation}
The $q$-tuple $\ang{s}$ is called the \emph{marker}\index{marker} of the object $a_{\ang{s}}$.
\item[Gluing]
For each $r \geq 0$, object $x \in \Op(r)$, marker $\ang{s}$, index $k \in \ufs{q}$, and partition of $s_k$ into $r$ subsets
\begin{equation}\label{subset-partition}
s_k = \coprod_{i \in \ufs{r}} \, s_{k,i} \subseteq \ufs{n}_k,
\end{equation}
$(a,\glu)$ consists of a \emph{gluing morphism}\index{gluing morphism} at $\big(x; \angs, k, \ang{s_{k,i}}_{i \in \ufs{r}}\big)$:
\begin{equation}\label{gluing-morphism}
\gaA_r\big(x\sscs \ang{a_{\ang{s} \compk\, s_{k,i}}}_{i \in \ufs{r}} \big) 
\fto{\glu_{x; \ang{s} \csp k, \ang{s_{k,i}}_{i \in \ufs{r}}}} a_{\ang{s}} \inspace \A.
\end{equation}
In \cref{gluing-morphism}, each marker 
\begin{equation}\label{angs_compk_ski}
\ang{s} \compk s_{k,i} = \big(s_1, \ldots, s_{k-1}, s_{k,i} \scs s_{k+1}, \ldots, s_q \big)
\end{equation}
is obtained from $\ang{s}$ by replacing the $k$-th subset, $s_k$, by $s_{k,i} \subseteq \ufs{n}_k$, as defined in \cref{compk}.  The gluing morphism in \cref{gluing-morphism} is also denoted by $(a,\glu)_{x; \ang{s} \csp k, \ang{s_{k,i}}_{i \in \ufs{r}}}$.
\end{description}
The pair $(a,\glu)$ is required to satisfy the axioms \cref{system_obj_unity,system_naturality,system_unity_i,system_unity_iii,system_equivariance,system_associativity,system_commutativity} below whenever they are defined.
\begin{description}
\item[Object unity]
If $s_j = \emptyset \subseteq \ufs{n}_j$ for some $j \in \ufs{q}$, then
\begin{equation}\label{system_obj_unity}
a_{\ang{s}} = \zero = \gaA_0(*) \in \A,
\end{equation}
the basepoint of $\A$ as defined in \cref{pseudoalg_zero}, where $* \in \Op(0)$ is the unique object.
\item[Naturality]
For each morphism $h \cn x \to y$ in the $G$-category $\Op(r)$ with $r \geq 0$, the following diagram in $\A$ commutes.
\begin{equation}\label{system_naturality}
\begin{tikzpicture}[vcenter]
\def\v{-1.5} \def\u{-.08}
\draw[0cell]
(0,0) node (a1) {\gaA_r\big(x \sscs \ang{a_{\ang{s} \compk\, s_{k,i}}}_{i \in \ufs{r}} \big)}
(a1)++(0,\v) node (b1) {\gaA_r\big(y \sscs \ang{a_{\ang{s} \compk\, s_{k,i}}}_{i \in \ufs{r}} \big)}
(a1)++(5,0) node (a2) {\phantom{a_{\ang{s}}}}
(a2)++(0,\u) node (a2') {a_{\ang{s}}}
(a2)++(0,\v) node (b2) {\phantom{a_{\ang{s}}}}
(b2)++(0,\u) node (b2') {a_{\ang{s}}}
;
\draw[1cell=1]
(a1) edge node {\glu_{x; \ang{s} \csp k, \ang{s_{k,i}}_{i \in \ufs{r}}}} (a2)
(b1) edge node {\glu_{y; \ang{s} \csp k, \ang{s_{k,i}}_{i \in \ufs{r}}}} (b2)
(a1) edge[transform canvas={xshift=2.5em}] node {\iso} node[swap] {\gaA_r(h \sscs \ang{1}_{i \in \ufs{r}} )} (b1)
(a2') edge[equal] (b2')
;
\end{tikzpicture}
\end{equation}
\item[Unity]
The gluing morphism $\glu_{x; \ang{s} \csp k, \ang{s_{k,i}}_{i \in \ufs{r}}}$ in \cref{gluing-morphism} is the identity morphism in each of the following two cases.
\begin{itemize}
\item If $s_j = \emptyset$ for some $j \in \ufs{q}$, then the following diagram in $\A$ commutes.
\begin{equation}\label{system_unity_i}
\begin{tikzpicture}[vcenter]
\def\v{-1.2} \def\u{-.08}
\draw[0cell]
(0,0) node (a1) {\gaA_r\big(x \sscs \ang{a_{\ang{s} \compk\, s_{k,i}}}_{i \in \ufs{r}} \big)}
(a1)++(0,\v) node (b1) {\gaA_r\big(x \sscs \ang{\zero}_{i \in \ufs{r}} \big)}
(b1)++(2.6,0) node (b) {\gaA_0(*) = \zero}
(a1)++(5,0) node (a2) {\phantom{a_{\ang{s}}}}
(a2)++(0,\u) node (a2') {a_{\ang{s}}}
(a2)++(0,\v) node (b2) {\zero}
;
\draw[1cell]
(a1) edge node {\glu_{x;\, \ang{s},\, k, \ang{s_{k,i}}_{i \in \ufs{r}}}} (a2)
(b) edge node {1} (b2)
(a1) edge[equal] (b1)
(a2') edge[equal] (b2)
(b1) edge[equal] node {\mathbf{b}} (b)
;
\end{tikzpicture}
\end{equation}
\begin{itemize}
\item The vertical equalities in \cref{system_unity_i} follow from the object unity axiom \cref{system_obj_unity} and, if $j=k$, the fact that $s_{k,i} \subseteq s_k$.
\item The equality labeled $\mathbf{b}$ follows from \cref{phi_id}.
\end{itemize}
\item If $r=1$---which implies $x = \opu \in \Op(1)$, the operadic unit---then the following diagram in $\A$ commutes.
\begin{equation}\label{system_unity_iii}
\begin{tikzpicture}[vcenter]
\def\v{-1.2} \def\u{-.08} \def\t{.6}
\draw[0cell]
(0,0) node (a1) {\gaA_1\big(\opu \sscs a_{\ang{s}} \big)}
(a1)++(2,0) node (a2) {\phantom{a_{\ang{s}}}}
(a2)++(0,\u) node (a2') {a_{\ang{s}}}
(a2)++(2,0) node (a3) {\phantom{a_{\ang{s}}}}
(a3)++(0,\u) node (a3') {a_{\ang{s}}}
;
\draw[1cell]
(a1) edge[equal] (a2)
(a2) edge node {1} (a3)
;
\draw[1cell=1]
(a1) [rounded corners=2pt, shorten <=-.5ex] |- ($(a2)+(-1,\t)$)
-- node {\glu_{\opu;\, \ang{s},\, k, \{s_k\}}} ($(a2)+(1,\t)$) -| (a3')
;
\end{tikzpicture}
\end{equation}
The equality in \cref{system_unity_iii} follows from the action unity axiom \cref{pseudoalg_action_unity} for $\A$.
\end{itemize}
\item[Equivariance]
For each permutation $\si \in \Sigma_r$, the following diagram in $\A$ commutes.
\begin{equation}\label{system_equivariance}
\begin{tikzpicture}[vcenter]
\def\v{-1.5} \def\u{-.08}
\draw[0cell]
(0,0) node (a1) {\gaA_r\big(x\si \sscs \ang{a_{\ang{s} \compk\, s_{k,i}}}_{i \in \ufs{r}} \big)}
(a1)++(0,\v) node (b1) {\gaA_r\big(x \sscs \ang{a_{\ang{s} \compk\, s_{k,\sigmainv(i)}}}_{i \in \ufs{r}} \big)}
(a1)++(6,0) node (a2) {\phantom{a_{\ang{s}}}}
(a2)++(0,\u) node (a2') {a_{\ang{s}}}
(a2)++(0,\v) node (b2) {\phantom{a_{\ang{s}}}}
(b2)++(0,\u) node (b2') {a_{\ang{s}}}
;
\draw[1cell=1]
(a1) edge node {\glu_{x\si ; \ang{s} \csp k, \ang{s_{k,i}}_{i \in \ufs{r}}}} (a2)
(b1) edge node {\glu_{x; \ang{s} \csp k, \ang{s_{k,\sigmainv(i)}}_{i \in \ufs{r}}}} (b2)
(a1) edge[equal] (b1)
(a2') edge[equal] (b2')
;
\end{tikzpicture}
\end{equation}
In \cref{system_equivariance}, the left vertical equality follows from the action equivariance axiom \cref{pseudoalg_action_sym} for $\A$.  The bottom horizontal gluing morphism uses the partition
\[s_k = \coprod_{i \in \ufs{r}}\, s_{k,\sigmainv(i)}.\]
\item[Associativity]
Suppose we are given objects
\[\begin{split}
\left(x \sscs \ang{x_i}_{i \in \ufs{r}} \right) &\in \Op(r) \times \txprod_{i \in \ufs{r}}\, \Op(t_i) \andspace\\
\boldx = \ga\big(x\sscs \ang{x_i}_{i \in \ufs{r}} \big) &\in \Op(t)
\end{split}\]
with $t = \sum_{i \in \ufs{r}} t_i$, a marker $\ang{s} = \ang{s_j \subseteq \ufs{n}_j}_{j \in \ufs{q}}$ as defined in \cref{marker}, an index $k \in \ufs{q}$, and partitions
\[s_k = \coprod_{i \in \ufs{r}} s_{k,i} \andspace s_{k,i} = \coprod_{\ell \in \ufs{t}_i} s_{k,i,\ell}\]
of, respectively, $s_k$ and $s_{k,i}$ for each $i \in \ufs{r}$.  Then the following diagram in $\A$ commutes.
\begin{equation}\label{system_associativity}
\begin{tikzpicture}[vcenter]
\def\h{4} \def\v{-2} \def\t{2em}
\draw[0cell=.85]
(0,0) node (a) {\gaA_r\Big(x \sscs \bang{\gaA_{t_i}\big(x_i \sscs \ang{a_{\ang{s} \compk\, s_{k,i,\ell}}}_{\ell \in \ufs{t}_i}  \big)}_{i \in \ufs{r}} \Big)}
(a)++(-\h,\v/2) node (b) {\gaA_r\big(x \sscs \ang{a_{\ang{s} \compk\, s_{k,i}}}_{i \in \ufs{r}} \big)}
(a)++(0,\v) node (c) {\gaA_t\big(\boldx \sscs \ang{\ang{a_{\ang{s} \compk\, s_{k,i,\ell}}}_{\ell \in \ufs{t}_i}}_{i \in \ufs{r}} \big)}
(c)++(-.5em,0) node (c') {\phantom{\gaA_t\Big(\boldx \sscs \bang{\bang{a_{\ang{s} \compk\, s_{k,i,\ell}}}_{\ell \in \ufs{t}_i}}_{i \in \ufs{r}} \Big)}}
(b)++(0,\v) node (d) {\phantom{a_{\ang{s}}}}
(d)++(\t,0) node (d') {a_{\ang{s}}}
;
\draw[1cell=.85]
(a) edge[bend right=8, shorten >=-0ex] node[swap,pos=.6,inner sep=0pt] {\gaA_r\big(x \sscs \ang{\glu_{x_i;\, \ang{s} \compk\, s_{k,i} \scs k, \ang{s_{k,i,\ell}}_{\ell \in \ufs{t}_i}}}_{i \in \ufs{r}} \big)} (b)
(b) edge[bend right=10,transform canvas={xshift=\t}] node[swap] {\glu_{x; \ang{s} \csp k, \ang{s_{k,i}}_{i \in \ufs{r}}}} (d)
(a) edge[bend left=10,transform canvas={xshift=-\t}] node {\phiA_{(r;\, t_1,\ldots,t_r)}} node[swap] {\iso} (c)
(c') edge[bend left=6] node {\glu_{\boldx; \ang{s} \csp k, \ang{\ang{s_{k,i,\ell}}_{\ell \in \ufs{t}_i}}_{i \in \ufs{r}}}} (d')
;
\end{tikzpicture}
\end{equation}
In \cref{system_associativity}, $\phiA_{(r;\, t_1,\ldots,t_r)}$ is the component of the associativity constraint \cref{phiA} of $\A$ at the objects $x$, $\ang{x_i}_{i \in \ufs{r}}$, and $\ang{\ang{a_{\ang{s} \compk\, s_{k,i,\ell}}}_{\ell \in \ufs{t}_i}}_{i \in \ufs{r}}$.  The bottom gluing morphism uses the partition
\[s_k = \coprod_{i \in \ufs{r}} \coprod_{\ell \in \ufs{t}_i}\, s_{k,i,\ell}.\] 
\item[Commutativity]
Suppose we are given a pair of objects
\[(x,y) \in \Op(r) \times \Op(t),\]
a marker $\ang{s} = \ang{s_j \subseteq \ufs{n}_j}_{j \in \ufs{q}}$, two distinct indices $k,\ell \in \ufs{q}$, and partitions
\[s_k = \coprod_{i \in \ufs{r}} s_{k,i} \subseteq \ufs{n}_k \andspace 
s_\ell = \coprod_{p \in \ufs{t}} s_{\ell,p} \subseteq \ufs{n}_\ell\]
of, respectively, $s_k$ and $s_\ell$.
 Then the following diagram in $\A$ commutes.
\begin{equation}\label{system_commutativity}
\begin{tikzpicture}[vcenter]
\def\h{4} \def\v{-1.7} \def\u{8} \def\t{2.5em} 
\draw[0cell=.8]
(0,0) node (a1) {\gaA_{tr} \big((x \intr y)\twist_{t,r} \sscs \ang{\ang{a_{\ang{s} \compk\, s_{k,i}\, \compell\, s_{\ell,p}}}_{i \in \ufs{r}}}_{p \in \ufs{t}}\big)} 
(a1)++(-\h,\v/2) node (a2) {\gaA_{tr} \big(y \intr x \sscs \ang{\ang{a_{\ang{s} \compk\, s_{k,i}\, \compell\, s_{\ell,p}}}_{i \in \ufs{r}}}_{p \in \ufs{t}}\big)}
(a2)++(0,\v) node (a3) {\gaA_t \big(y \sscs \bang{\gaA_r \big(x \sscs \ang{a_{\ang{s} \compk\, s_{k,i}\, \compell\, s_{\ell,p}}}_{i \in \ufs{r}}\big)}_{p \in \ufs{t}}\big)} 
(a3)++(0,\v) node (a4) {\gaA_t\big(y \sscs \ang{a_{\ang{s} \compell\, s_{\ell,p}}}_{p \in \ufs{t}}\big)}
(a4)++(\t,\v) node (a5) {a_{\ang{s}}}
(a4)++(0,\v) node (a5') {\phantom{a_{\ang{s}}}}
(a1)++(0,\v) node (b2) {{\gaA_{rt} \big(x \intr y \sscs \ang{\ang{a_{\ang{s} \compk\, s_{k,i}\, \compell\, s_{\ell,p}}}_{p \in \ufs{t}}}_{i \in \ufs{r}}\big)}}
(b2)++(0,\v) node (b3) {\gaA_r\big(x \sscs \bang{\gaA_t \big(y \sscs \ang{a_{\ang{s} \compk\, s_{k,i}\, \compell\, s_{\ell,p}}}_{p \in \ufs{t}}\big)}_{i \in \ufs{r}}\big)} 
(b3)++(0,\v) node (b4) {\gaA_r\big(x \sscs \ang{a_{\ang{s} \compk\, s_{k,i}}}_{i \in \ufs{r}}\big)}
(b4)++(-1em,0) node (b4') {\phantom{\gaA_r\big(x \sscs \bang{a_{\ang{s} \compk\, s_{k,i}}}_{i \in \ufs{r}}\big)}}
;
\draw[1cell=.8]
(a1) edge[bend right=\u] node[swap] {\gaA_{tr}\big(\pcom_{r,t} \sscs \ang{1}_{i \in \ufs{r},\, p \in \ufs{t}} \big)} (a2)
(a2) edge[transform canvas={xshift=\t}] node[swap] {\big(\phiA_{(t;\, r,\ldots,r)}\big)^{-1}} (a3)
(a3) edge[transform canvas={xshift=\t}] node[swap] {\gaA_t\big(y \sscs \ang{\glu_{x;\, \ang{s} \compell\, s_{\ell,p} \scs k, \ang{s_{k,i}}_{i \in \ufs{r}}}}_{p \in \ufs{t}} \big)} (a4)
(a4) edge[transform canvas={xshift=\t}] node[swap] {\glu_{y; \ang{s} \csp \ell, \ang{s_{\ell,p}}_{p \in \ufs{t}}}} (a5')
(a1) edge[transform canvas={xshift=-\t},equal] node {\mathbf{eq}} (b2)
(b2) edge[transform canvas={xshift=-\t}] node {\big(\phiA_{(r;\, t,\ldots,t)}\big)^{-1}} (b3)
(b3) edge[transform canvas={xshift=-\t}] node {\gaA_r\big(x \sscs \ang{\glu_{y; \ang{s} \compk\, s_{k,i} \scs \ell, \ang{s_{\ell,p}}_{p \in \ufs{t}}}}_{i \in \ufs{r}}\big)} (b4)
(b4') edge[bend left=\u] node[pos=.5] {\glu_{x; \ang{s} \csp k, \ang{s_{k,i}}_{i \in \ufs{r}}}} (a5)
;
\end{tikzpicture}
\end{equation}
\begin{itemize}
\item In \cref{system_commutativity}, the marker
\[\ang{s} \compk\, s_{k,i}\, \compell\, s_{\ell,p}\]
is obtained from $\ang{s}$ by replacing
\begin{itemize}
\item its $k$-th subset by $s_{k,i}$ and
\item its $\ell$-th subset by $s_{\ell,p}$,
\end{itemize} 
as defined in \cref{compkl}.
\item The equality labeled $\mathbf{eq}$ follows from the action equivariance axiom \cref{pseudoalg_action_sym} for $\A$, applied to the transpose permutation $\twist_{t,r} \in \Sigma_{tr}$ defined in \cref{eq:transpose_perm}.  The transpose permutation $\twist_{t,r}$ switches the order of indexing from $\ang{\ang{\cdots}_{i \in \ufs{r}}}_{p \in \ufs{t}}$ to $\ang{\ang{\cdots}_{p \in \ufs{t}}}_{i \in \ufs{r}}$.
\item Each $\intr$ is an instance of the intrinsic pairing of $\Op$ (\cref{def:intrinsic_pairing}), with
\[\begin{split}
x \intr y &= \ga\big(x \sscs \ang{y}_{i \in \ufs{r}}\big) \in \Op(rt) \andspace\\
y \intr x &= \ga\big(y \sscs \ang{x}_{p \in \ufs{t}}\big) \in \Op(tr).
\end{split}\]
The top boundary arrow in \cref{system_commutativity} involves the isomorphism
\[(x \intr y) \twist_{t,r} \fto[\iso]{\pcom_{r,t}} y \intr x,\]
which is the $(x,y)$-component of the $(r,t)$-pseudo-commutativity isomorphism \cref{pseudocom_isos} of $\Op$. 
\item Both $\phiA_{(t;\, r,\ldots,r)}$ and $\phiA_{(r;\, t,\ldots,t)}$ are instances of the associativity constraint \cref{phiA} of $\A$.  
\end{itemize}
\end{description}
This finishes the definition of an $\angordn$-system $(a,\glu)$.  

Moreover, we define the following.
\begin{itemize}
\item A \emph{strong $\angordn$-system}\dindex{strong}{system} is an $\angordn$-system such that each gluing morphism $\glu_{x; \ang{s} \csp k, \ang{s_{k,i}}_{i \in \ufs{r}}}$ in \cref{gluing-morphism} is an isomorphism.
\item The \emph{base $\angordn$-system}\dindex{base}{system} is the $\angordn$-system $(\zero,1_\zero)$\label{not:basesystem} with
\begin{itemize}
\item each component object given by the basepoint $\zero \in \A$ and
\item each gluing morphism given by the identity morphism $1_\zero$ of $\zero$.
\end{itemize}
\cref{expl:nsystem_base} shows that $(\zero,1_\zero)$ is well defined.
\defmark
\end{itemize}
\end{definition}

\begin{remark}\label{rk:angordn_system}
\cref{def:nsystem} of an $\angordn$-system is an operadic generalization of Elmendorf-Mandell $J$-theory \cite[Construction 4.4]{elmendorf-mandell}, which may be regarded as the case where $\Op$ is the Barratt-Eccles operad $\BE$.  Elmendorf-Mandell $J$-theory itself is a refinement of Segal $J$-theory.  The reader is referred to \cite[Ch.\! 8 and 10]{cerberusIII} for a detailed discussion of Segal $J$-theory and Elmendorf-Mandell $J$-theory.
\end{remark}

\begin{example}[Barratt-Eccles Operads]\label{ex:n_systems}
\cref{def:nsystem} requires a $\Tinf$-operad \pcref{as:OpA}.  By \cref{BE_pseudocom,GBE_pseudocom}, the Barratt-Eccles $\Gcat$-operad $\BE$ (\cref{def:BE-Gcat}), with 
\[\BE(n) = \ESigma_n,\] and the $G$-Barratt-Eccles operad $\GBE$ (\cref{def:GBE}), with 
\[\GBE(n) = \Catg(\EG,\ESigma_n) \iso \tn[G,\Sigma_n],\]
are both $\Tinf$-operads.
By \cref{thm:BEpseudoalg}, $\BE$-pseudoalgebras correspond to \emph{naive} symmetric monoidal $G$-categories (\cref{def:naive_smGcat}) via each of the 2-equivalences $\Phi$.  By \cref{def:GBE_pseudoalg}, $\GBE$-pseudoalgebras are \emph{genuine} symmetric monoidal $G$-categories.
\end{example}

\begin{explanation}[Assumptions about $\Op$]\label{expl:nsystem_assumptions}
\cref{as:OpA} is used in \cref{def:nsystem} as follows. 
\begin{enumerate}
\item The reduced assumption, $\Op(0) = \{*\}$, is used in the unity axiom \cref{system_unity_i} to ensure the equality
\[\gaA_r\big(x \sscs \ang{\zero}_{i \in \ufs{r}} \big) = \gaA_0(*).\]
\item The assumption $\Op(1) = \{\opu\}$ is used in the unity axiom \cref{system_unity_iii} to ensure $x = \opu$ when $r=1$.
\item The pseudo-commutative structure $\pcom$ of $\Op$ is used in the commutativity axiom \cref{system_commutativity}.\defmark
\end{enumerate}
\end{explanation}

\begin{explanation}[Base $\angordn$-System]\label{expl:nsystem_base}
To see that the base $\angordn$-system $(\zero,1_\zero)$ is well defined, we observe that it satisfies all the axioms of an $\angordn$-system as follows.
\begin{description}
\item[Unity and equivariance] The object unity axiom \cref{system_obj_unity}, the unity axioms \cref{system_unity_i,system_unity_iii}, and the equivariance axiom \cref{system_equivariance} hold because each gluing morphism is the identity morphism $1_\zero \cn \zero \to \zero$.
\item[Naturality] The naturality axiom \cref{system_naturality} holds because the left vertical arrow in that diagram,
\[\gaA_r\big(x \sscs \ang{\zero}_{i \in \ufs{r}}\big) \fto{\gaA_r(h \sscs \ang{1_\zero}_{i \in \ufs{r}})} 
\gaA_r\big(y \sscs \ang{\zero}_{i \in \ufs{r}}\big),\]
is equal to $1_\zero$.  This is true by
\begin{itemize}
\item $r$ applications of the basepoint axiom \cref{pseudoalg_basept_axiom} for $\A$ and 
\item the naturality of the associativity constraint $\phiA$ \cref{phiA}.
\end{itemize}
\item[Associativity] The associativity axiom \cref{system_associativity} holds because the isomorphism in that diagram,
\[\gaA_r\big(x \sscs \ang{\gaA_{t_i} (x_i \sscs \ang{\zero}_{\ell \in \ufs{t}_i} ) }_{i \in \ufs{r}} \big) \fto{\phiA_{(r;\, t_1,\ldots,t_r)}}
\gaA_t\big(\boldx \sscs \ang{\ang{\zero}_{\ell \in \ufs{t}_i}}_{i \in \ufs{r}} \big),\]
is equal to $1_\zero$ by the composition axiom \cref{pseudoalg_comp_axiom} for $\A$ and \cref{phi_id}.
\item[Commutativity] The commutativity axiom \cref{system_commutativity} holds for the base $\angordn$-system for the following reasons.
\begin{enumerate}
\item Each of $\phiA_{(t;\, r,\ldots,r)}$ and $\phiA_{(r;\, t,\ldots,t)}$ in \cref{system_commutativity} is equal to $1_\zero$ by the composition axiom \cref{pseudoalg_comp_axiom} for $\A$ and \cref{phi_id}.
\item The arrow $\gaA_{tr}\big(\pcom_{r,t} \sscs \ang{1_\zero}_{i \in \ufs{r},\, p \in \ufs{t}}\big)$ in \cref{system_commutativity} is equal to $1_\zero$ by
\begin{itemize}
\item $tr$ applications of the basepoint axiom \cref{pseudoalg_basept_axiom} for $\A$ and
\item the naturality of the associativity constraint $\phiA$ \cref{phiA}.
\end{itemize}
\end{enumerate}
\end{description}
Thus, the base $\angordn$-system $(\zero,1_\zero)$ is well defined.
\end{explanation}

\begin{explanation}[$\angordn$-Systems]\label{expl:nsystem}
Consider \cref{def:nsystem} of $\angordn$-systems.
\begin{enumerate}
\item\label{expl:nsystem_zero} If $r=0$ in the gluing morphism $\glu_{x; \ang{s} \csp k, \ang{s_{k,i}}_{i \in \ufs{r}}}$---which implies $x = * \in \Op(0)$ and $s_k = \emptyset$---then the unity axiom \cref{system_unity_i} reduces to the following commutative diagram in $\A$.
\begin{equation}\label{system_unity_ii}
\begin{tikzpicture}[vcenter]
\def\v{-1.2} \def\u{-.08}
\draw[0cell]
(0,0) node (a1) {\gaA_0(*)}
(a1)++(0,\v) node (b1) {\zero}
(a1)++(3.5,0) node (a2) {\phantom{a_{\ang{s}}}}
(a2)++(0,\u) node (a2') {a_{\ang{s}}}
(a2)++(0,\v) node (b2) {\zero}
;
\draw[1cell]
(a1) edge node {\glu_{*; \ang{s} \csp k, \ang{}}} (a2)
(b1) edge node {1} (b2)
(a1) edge[equal] (b1)
(a2') edge[equal] (b2)
;
\end{tikzpicture}
\end{equation}
\item\label{expl:nsystem_i} Suppose $\ord{n}_j = \ord{0}$ for some $j \in \ufs{q}$.  Then 
\[\angordn = \vstar,\]
the initial-terminal basepoint of $\Gskel$ \cref{Gsk_objects}.  In this case, the only $\angordn$-system is given by the base $\angordn$-system $(\zero,1_\zero)$.  In more detail, the object unity axiom \cref{system_obj_unity} implies that, for each marker $\ang{s}$ as defined in \cref{marker}, 
\[a_{\ang{s}} = \zero \in \A\]
because $s_j \subseteq \ufs{0} = \emptyset$ is the empty set.  Moreover, the unity axiom \cref{system_unity_i} implies that each gluing morphism $\glu$ is the identity morphism 
\[1_\zero \cn \zero \to \zero\]
of the basepoint $\zero$.
\item\label{expl:nsystem_ii} Among the axioms of an $\angordn$-system, only the associativity axiom \cref{system_associativity} and the commutativity axiom \cref{system_commutativity} involve the associativity constraint $\phiA$ of the $\Op$-pseudoalgebra $\A$.  The invertibility of $\phiA$ is necessary in the commutativity axiom.
\item\label{expl:nsystem_iii} The commutativity axiom \cref{system_commutativity} requires two distinct indices $k, \ell \in \ufs{q}$, so it only happens when $\angordn$ has length $q \geq 2$.\defmark
\item\label{expl:nsystem_iv} In the partition \cref{subset-partition} of $s_k$, suppose that $r>0$ and there exists an index $\jm \in \ufs{r}$ such that 
\[s_{k,i} = \begin{cases} s_k & \text{if $i = \jm$, and}\\
\emptyset & \text{if $i  \in \ufs{r} \setminus \{\jm\}$}.
\end{cases}\]
By the object unity axiom \cref{system_obj_unity}, the gluing morphism in \cref{gluing-morphism} has the form
\begin{equation}\label{skj_glu}
\gaA_r\big(x \sscs \zero^{\jm-1}, a_{\ang{s}}, \zero^{r-\jm} \big) \fto{\glu_{x; \angs \csp k, \ang{s_{k,\ell}}_{\ell \in \ufs{r}}}} a_{\angs},
\end{equation}
where $\zero^t$ denotes a $t$-tuple of copies of $\zero$.  This gluing morphism is the identity morphism of $a_{\angs}$ by the axioms \cref{pseudoalg_basept_axiom,system_unity_i,system_unity_iii,system_associativity}.
\end{enumerate}
\end{explanation}

Next, we define morphisms and the category of $\angordn$-systems in an $\Op$-pseudoalgebra $\A$.

\begin{definition}[Morphisms of Systems]\label{def:nsystem_morphism}
Under \cref{as:OpA}, suppose $(a,\glu^a)$ and $(b,\glu^b)$ are $\angordn$-systems in $\A$ (\cref{def:nsystem}).  A \emph{morphism}\index{system!morphism} of $\angordn$-systems
\[(a,\glu^a) \fto{\theta} (b,\glu^b)\]
consists of, for each marker $\ang{s} = \ang{s_j \subseteq \ufs{n}_j}_{j \in \ufs{q}}$, an \emph{$\ang{s}$-component morphism}\index{component morphism}
\begin{equation}\label{theta_angs}
a_{\ang{s}} \fto{\theta_{\ang{s}}} b_{\ang{s}} \inspace \A
\end{equation}
such that the following two axioms are satisfied.
\begin{description}
\item[Unity] If $s_j = \emptyset \subseteq \ufs{n}_j$ for some $j \in \ufs{q}$, then there is an equality of morphisms
\begin{equation}\label{nsystem_mor_unity}
a_{\ang{s}} = \zero \fto{\theta_{\ang{s}} \,=\, 1_\zero} b_{\ang{s}} = \zero. 
\end{equation}
The object equalities $a_{\ang{s}} = b_{\ang{s}} = \zero$ follow from the object unity axiom \cref{system_obj_unity}.
\item[Compatibility] For each $r \geq 0$, object $x \in \Op(r)$, marker $\ang{s}$, index $k \in \ufs{q}$, and partition
\[s_k = \coprod_{i \in \ufs{r}} \, s_{k,i} \subseteq \ufs{n}_k,\]
the following diagram in $\A$ commutes.
\begin{equation}\label{nsystem_mor_compat}
\begin{tikzpicture}[vcenter]
\def\v{-1.5} \def\u{-.08}
\draw[0cell]
(0,0) node (a1) {\gaA_r\big(x \sscs \ang{a_{\ang{s} \compk\, s_{k,i}}}_{i \in \ufs{r}} \big)}
(a1)++(0,\v) node (b1) {\gaA_r\big(x \sscs \ang{b_{\ang{s} \compk\, s_{k,i}}}_{i \in \ufs{r}} \big)}
(a1)++(5,0) node (a2) {\phantom{a_{\ang{s}}}}
(a2)++(0,\u) node (a2') {a_{\ang{s}}}
(a2)++(0,\v) node (b2) {\phantom{b_{\ang{s}}}}
(b2)++(0,\u) node (b2') {b_{\ang{s}}}
;
\draw[1cell=.9]
(a1) edge node {\glu^a_{x; \ang{s} \csp k, \ang{s_{k,i}}_{i \in \ufs{r}}}} (a2)
(b1) edge node {\glu^b_{x; \ang{s} \csp k, \ang{s_{k,i}}_{i \in \ufs{r}}}} (b2)
(a1) edge[transform canvas={xshift=3em}] node[swap] {\gaA_r\big(1_x \sscs \ang{\theta_{\ang{s} \compk\, s_{k,i}}}_{i \in \ufs{r}} \big)} (b1)
(a2') edge node {\theta_{\ang{s}}} (b2')
;
\end{tikzpicture}
\end{equation}
\end{description}
This finishes the definition of a morphism of $\angordn$-systems.  Moreover, we define the following.
\begin{itemize}
\item Composition and identity morphisms of $\angordn$-systems in $\A$ are defined componentwise. 
\item We denote by $\Aangordn$\label{not:Aangordn} the category of $\angordn$-systems in $\A$ and their morphisms.
\item We denote by $\Aangordnsg$ the full subcategory of $\Aangordn$ with strong $\angordn$-systems in $\A$ as objects.\defmark
\end{itemize}
\end{definition}

\begin{example}[Base Case]\label{ex:nsystem_basecase}
Suppose $\ord{n}_j = \ord{0}$ for some $j \in \ufs{q}$ in \cref{def:nsystem}.  As we discuss in \cref{expl:nsystem} \pcref{expl:nsystem_i}, in this case, the only $\angordn$-system in $\A$ is the base $\angordn$-system $(\zero,1_\zero)$.  Moreover, by the unity axiom \cref{nsystem_mor_unity}, the only morphism on $(\zero,1_\zero)$ is the identity morphism.  Thus, in this case, there is an isomorphism
\[\Aangordn = \Aangordnsg \iso \boldone,\] 
a terminal category. 
\end{example}

\begin{example}[Compatibility for $r=0$]\label{ex:nsystem_mor_zero}
Consider the compatibility axiom \cref{nsystem_mor_compat} with $r=0$.  In this case, we have that 
\[x = * \in \Op(0) \cq s_k = \emptyset, \andspace \gaA_0(1_*) = 1_\zero.\]  
Moreover, by \cref{system_unity_ii}, 
\[\glu^a_{*; \angs \csp k, \ang{}} = \glu^b_{*; \angs \csp k, \ang{}} = 1_\zero.\]
Thus, the compatibility axiom \cref{nsystem_mor_compat} with $r=0$ states that 
\[\theta_{\angs} = 1_\zero\]
whenever some entry in $\angs$ is $\emptyset$.  This is precisely the unity axiom \cref{nsystem_mor_unity}.  When we check that something is a morphism of systems, we need to check separately the $r=0$ case of the compatibility axiom, which is the unity axiom.  This is the reason why \cref{def:nsystem_morphism} states the unity axiom separately.
\end{example}

\subsection*{Pointed Categories of Systems}

A \emph{pointed category} is a category equipped with a chosen object, called the basepoint (\cref{def:pointed-category}).  Next, we define the pointed categories of systems for an $\Op$-pseudoalgebra $(\A,\gaA,\phiA)$ (\cref{def:pseudoalgebra}) and objects in $\Gsk$ \cref{Gsk_objects}.

\begin{definition}\label{def:Aangordn_system}
Under \cref{as:OpA}, suppose $\angordn \in \Gsk$ is an object.  We define the pointed category of \emph{$\angordn$-systems in $\A$}\index{system!pointed category} as follows.
\begin{equation}\label{A_angordn}
\Aangordn = \begin{cases}
\boldone & \text{if $\angordn = \vstar$},\\
(\A,\zero) & \text{if $\angordn = \ang{}$, and}\\
\big(\Aangordn, (\zero,1_\zero)\big) & \text{if $\angordn \neq \vstar, \ang{}$.}
\end{cases}
\end{equation}
Moreover, we define the pointed category $\Aangordnsg$ of \emph{strong $\angordn$-systems in $\A$} using \cref{A_angordn} and replacing the category $\Aangordn$ with $\Aangordnsg$ in the last case.
\end{definition}

\begin{explanation}\label{expl:Aangordn_system}
Let us consider \cref{A_angordn}.
\begin{enumerate}
\item\label{expl:nbeta_system_i} In the first case, $\boldone$ is the terminal $G$-category.
\item\label{expl:nbeta_system_ii} In the second case, $\A$ is the underlying small $G$-category of the $\Op$-pseudoalgebra $(\A,\gaA,\phiA)$.  It is equipped with the $G$-fixed basepoint $\zero = \gaA_0(*)$, as defined in \cref{pseudoalg_zero}.
\item\label{expl:nbeta_system_iii} In the last case, $\Aangordn$ is the category of $\angordn$-systems in $\A$ (\cref{def:nsystem_morphism}).  Its basepoint is given by the base $\angordn$-system $(\zero,1_\zero)$ in \cref{expl:nsystem_base}.\defmark
\end{enumerate}
\end{explanation}

\begin{example}[Tuple of $\ord{1}$]\label{ex:nsystem_ones}
We consider \cref{def:Aangordn_system} with
\[\angordn = \ordtu{1}_{j \in \ufs{q}} = \big(\ord{1},\ldots,\ord{1}\big) \in \Gsk\]
consisting of $q$ copies of $\ord{1}$.  The only subsets of $\ufs{1} = \{1\}$ are $\emptyset$ and $\{1\}$.
\begin{itemize}
\item By the object unity axiom \cref{system_obj_unity}, a $\ordtu{1}$-system $(a,\glu)$ has $\ang{s}$-component object given by
\[a_{\ang{s}} = \begin{cases}
a_{\ang{\{1\}}_{j \in \ufs{q}}} & \text{if $\ang{s} = \ang{\{1\}}_{j \in \ufs{q}}$ and},\\
\zero & \text{otherwise.}
\end{cases}\]
\item By the unity axioms \cref{system_unity_i,system_unity_iii}, the associativity axiom \cref{system_associativity}, and the basepoint axiom \cref{pseudoalg_basept_axiom} of $\A$, each gluing morphism $\glu$ is an identity morphism.  
\item By the unity axiom \cref{nsystem_mor_unity}, a morphism $\theta$ of $\ordtu{1}$-systems has $\ang{s}$-component morphism given by
\[\theta_{\ang{s}} = \begin{cases}
\theta_{\ang{\{1\}}_{j \in \ufs{q}}} & \text{if $\ang{s} = \ang{\{1\}}_{j \in \ufs{q}}$ and},\\
1_\zero & \text{otherwise.}
\end{cases}\]
Moreover, the compatibility axiom \cref{nsystem_mor_compat} for $r \geq 1$ is a consequence of the naturality of the associativity constraint $\phiA$ \cref{phiA}.  
\end{itemize}
Thus, there is a canonical isomorphism of pointed categories
\begin{equation}\label{Aordtuone}
(\A,\zero) \iso \big(\A\ordtu{1}, (\zero,1_\zero)\big) 
= \big(\syssg{\A}{\ordtu{1}}, (\zero,1_\zero)\big).
\end{equation}
Under this pointed isomorphism, an object $b \in \A$ corresponds to the $\ordtu{1}$-system with
\begin{itemize}
\item $\ang{\{1\}}_{j \in \ufs{q}}$-component object given by $b$, 
\item all other component objects given by the basepoint $\zero \in \A$, and 
\item all gluing morphisms given by identities.
\end{itemize}
Such a $\ordtu{1}$-system is necessarily strong because each gluing morphism is invertible. 
\end{example}

\subsection*{Pointed $G$-Categories of Systems}

A \emph{pointed $G$-category} is a $G$-category equipped with a $G$-fixed basepoint (\cref{def:GCat,def:ptGcat}).  The next definition equips each of the pointed categories $\Aangordn$ and $\Aangordnsg$ with a $G$-action.

\begin{definition}\label{def:Aangordn_gcat}
Suppose $\Op$ is a $\Tinf$-operad \pcref{as:OpA}, and $\A$ is an $\Op$-pseudoalgebra \pcref{def:pseudoalgebra}.  The pointed category $\Aangordn$ (\cref{def:nsystem,def:nsystem_morphism,def:Aangordn_system}) is extended to a pointed $G$-category\index{system!pointed G-category@pointed $G$-category} as follows.  

If $\angordn = \vstar$ or $\ang{}$, then, by \cref{expl:Aangordn_system} \pcref{expl:nbeta_system_i} and \pcref{expl:nbeta_system_ii},
\[\sys{\A}{\vstar} = \boldone \andspace \sys{\A}{\ang{}} = (\A,\zero)\] 
are already pointed $G$-categories.  Suppose $\angordn \in \Gsk \setminus \{\vstar, \ang{}\}$ is an object of length $q > 0$. 

\parhead{$G$-action on systems}.  Suppose $g \in G$ and $(a,\glu) \in \Aangordn$ is an $\angordn$-system in $\A$ (\cref{def:nsystem}).  The $\angordn$-system in $\A$ 
\begin{equation}\label{nsystem_gaction}
g \cdot (a,\glu) = \big(ga, g\glu \big)
\end{equation}
is defined as follows.
\begin{description}
\item[Component objects] For each marker $\ang{s} = \ang{s_j \subseteq \ufs{n}_j}_{j \in \ufs{q}}$, the $\ang{s}$-component object of $(ga,g\glu)$ is defined as
\begin{equation}\label{ga_scomponent}
(ga)_{\ang{s}} = g a_{\ang{s}} \in \A,
\end{equation}
the image of $a_{\angs}$ under the $g$-action on $\A$.
\item[Gluing] Suppose we are given $r \geq 0$, an object $x \in \Op(r)$, a marker $\ang{s}$, an index $k \in \ufs{q}$, and a partition
\[s_k = \coprod_{i \in \ufs{r}}\, s_{k,i} \subseteq \ufs{n}_k\]
of $s_k$ into $r$ subsets.  The gluing morphism of $(ga,g\glu)$ at $(x; \angs, k, \ang{s_{k,i}}_{i \in \ufs{r}})$ is defined by the following commutative diagram.
\begin{equation}\label{ga_gluing}
\begin{tikzpicture}[vcenter]
\def\v{-1}
\draw[0cell=.8]
(0,0) node (a1) {\gaA_r\big(x \sscs \ang{(ga)_{\ang{s} \compk\, s_{k,i}} }_{i \in \ufs{r}} \big)}
(a1)++(0,\v) node (a2) {\gaA_r\big(x \sscs \ang{g a_{\ang{s} \compk\, s_{k,i}}}_{i \in \ufs{r}}\big)}
(a2)++(0,\v) node (a3) {g \gaA_r \big(\ginv x \sscs \ang{a_{\ang{s} \compk\, s_{k,i}}}_{i \in \ufs{r}} \big)}
(a1)++(5,0) node (b1) {(ga)_{\ang{s}}}
(b1)++(0,2*\v) node (b3) {g a_{\ang{s}}}
;
\draw[1cell=.9]
(a1) edge node {(g\glu)_{x; \ang{s},\, k, \ang{s_{k,i}}_{i \in \ufs{r}}}} (b1)
(a3) edge node {g \glu_{\ginv x; \ang{s},\, k, \ang{s_{k,i}}_{i \in \ufs{r}}}} (b3)
(a1) edge[equal,shorten >=-.5ex] (a2)
(a2) edge[equal,shorten >=-.5ex] node[swap,inner sep=4pt] {\mathbf{f}} (a3)
(b1) edge[equal] (b3)
;
\end{tikzpicture}
\end{equation}
\begin{itemize}
\item The two unlabeled equalities in \cref{ga_gluing} follow from the definition \cref{ga_scomponent} of $(ga)_{\ang{s}}$.
\item The equality labeled $\mathbf{f}$ follows from the $G$-functoriality of $\gaA_r$ \cref{gaAn}.
\item In the bottom horizontal arrow in \cref{ga_gluing}, 
\[\gaA_r \big(\ginv x \sscs \ang{a_{\ang{s} \compk\, s_{k,i}}}_{i \in \ufs{r}} \big)
\fto{\glu_{\ginv x; \ang{s},\, k, \ang{s_{k,i}}_{i \in \ufs{r}}}} a_{\ang{s}}\]
is the gluing morphism of $(a,\glu)$ at $(\ginv x; \ang{s}, k, \ang{s_{k,i}}_{i \in \ufs{r}})$, and $g\glu_{\cdots}$ is its image under the $g$-action on $\A$.
\end{itemize}
\end{description}

\parhead{Axioms}.  Each of the axioms of an $\angordn$-system,  \cref{system_obj_unity,system_naturality,system_unity_i,system_unity_iii,system_equivariance,system_associativity,system_commutativity}, for $(ga,g\glu)$ follows from the corresponding axiom for $(a,\glu)$ and the following facts.
\begin{itemize}
\item The $g$-action on $\A$ is a functor.
\item The basepoint $\zero \in \A$, its identity morphism $1_\zero$, and the base $\angordn$-system $(\zero,1_\zero)$ are $G$-fixed. 
\item The object $* \in \Op(0)$ and the operadic unit $\opu \in \Op(1)$ are $G$-fixed.
\item The following structures are $G$-equivariant: the right symmetric group action \cref{rightsigmaaction} and the composition \cref{eq:enr-defn-gamma} on the $\Gcat$-operad $\Op$, the pseudo-commutative structure $\pcom$ on $\Op$ \cref{pseudocom_isos}, the $\Op$-action $\gaA_r$ \cref{gaAn}, and the associativity constraint $\phiA$ \cref{phiA}.
\end{itemize}
For example, to prove the commutativity axiom \cref{system_commutativity} for $(ga,g\glu)$, we start with the commutative diagram \cref{system_commutativity} with $(x,y)$ replaced by $(\ginv x, \ginv y)$.  Then we apply the $g$-action on $\A$ to the latter commutative diagram.  This finishes the definition of the $\angordn$-system $g \cdot (a,\glu) = (ga,g\glu)$.

\parhead{$G$-action on morphisms}.  For a morphism of $\angordn$-systems in $\A$ (\cref{def:nsystem_morphism})
\[(a,\glu^a) \fto{\theta} (b,\glu^b),\]
the morphism of $\ordtu{n}$-systems
\begin{equation}\label{gtheta}
(ga,g\glu^a) \fto{g\theta} (gb, g\glu^b)
\end{equation}
is defined by, for each marker $\ang{s} = \ang{s_j \subseteq \ufs{n}_j}_{j \in \ufs{q}}$, the $\ang{s}$-component morphism
\begin{equation}\label{gtheta_angs}
(ga)_{\ang{s}} = ga_{\ang{s}}
\fto{(g\theta)_{\ang{s}} \,=\, g\theta_{\ang{s}}}
(gb)_{\ang{s}} = gb_{\ang{s}}.
\end{equation}
\begin{itemize}
\item The unity axiom \cref{nsystem_mor_unity} holds for $g\theta$ because the identity morphism $1_\zero$ is $G$-fixed. 
\item The compatibility axiom \cref{nsystem_mor_compat} holds for $g\theta$ by the compatibility axiom for $\theta$, the functoriality of the $g$-action on $\A$, and the $G$-equivariance of $\gaA_r$.
\end{itemize}
This finishes the definition of the morphism $g\theta$ of $\angordn$-systems.

\parhead{Functoriality}.  The functoriality of the $g$-action on $\Aangordn$ follows from
\begin{itemize}
\item the definition \cref{gtheta_angs} of $(g\theta)_{\ang{s}}$ and
\item the functoriality of the $g$-action on $\A$.
\end{itemize}
The two defining conditions for a $G$-action---which are stated immediately after \cref{gactioniso}---are satisfied by $\Aangordn$ by
\begin{itemize}
\item the definitions \cref{ga_scomponent,ga_gluing,gtheta_angs}, and 
\item the corresponding conditions for the $G$-actions on $\A$ and $\Op(r)$.
\end{itemize}  
This finishes the definition of the pointed $G$-category $\Aangordn$.

\parhead{Strong variant}.  Moreover, we extend the pointed category 
\begin{equation}\label{sgAordnbe}
\Aangordnsg
\end{equation}
of strong $\angordn$-systems in $\A$ to a pointed $G$-category by restricting the $G$-action on $\Aangordn$ defined above to the full subcategory $\Aangordnsg$.  This is well defined because, if each gluing morphism $\glu_{\cdots}$ is an isomorphism in \cref{ga_gluing}, then so is its image $g\glu_{\cdots}$ under the $g$-action.
\end{definition}

\section{$\Gskg$-Categories from Operadic Pseudoalgebras: Morphisms}
\label{sec:jemg_morphisms}

For a $\Tinf$-operad $\Op$ \pcref{as:OpA}, \cref{sec:jemg_objects} constructs the object assignment of the $\Gskg$-category (\cref{expl:ggcat_obj})
\[(\Gsk,\vstar) \fto{\Jgo\A = \Adash} (\Gcatst,\boldone)\]
for an $\Op$-pseudoalgebra $(\A,\gaA,\phiA)$ (\cref{def:pseudoalgebra}).  This section constructs the morphism assignment of the $\Gskg$-category $\Jgo\A$.  In other words, for each morphism \cref{fangpsi}
\[\angordm \fto{\upom} \angordn \inspace \Gsk,\]
this section constructs a pointed $G$-functor
\[\Aangordm \fto{\sys{(\Jgo\A)}{(\upom)} = \Aupom} \Aangordn\]
from the pointed $G$-category of $\angordm$-systems in $\A$ to the pointed $G$-category of $\angordn$-systems in $\A$ (\cref{def:Aangordn_gcat}).  The proof that the object and morphism assignments define a $\Gskg$-category, that is, a pointed functor
\[(\Gsk,\vstar) \fto{\Jgo\A = \Adash} (\Gcatst,\boldone),\] 
is given in \cref{sec:jemg_morphisms_ii}.  For the rest of this section, $(\A,\gaA,\phiA)$ denotes an arbitrary $\Op$-pseudoalgebra.

\secoutline

The pointed $G$-functor $\Aupom$ is constructed in several steps.  When the domain object  $\angordm$ has positive length, the morphism $\upom$ has the form $(f,\angpsi)$.  The pointed $G$-functor $\Afpsi$ is the composite of two pointed $G$-functors, denoted $\ftil$ and $\psitil$. 
\begin{itemize}
\item \cref{def:ftil_functor} constructs the pointed $G$-functor
\[\Aangordm \fto{\ftil} \Afangordm\]
from $\angordm$-systems to $f_*\angordm$-systems.
\item \cref{def:psitil_functor} constructs the pointed $G$-functor
\[\Afangordm \fto{\psitil} \Aangordn\]
from $f_*\angordm$-systems to $\angordn$-systems.
\item \cref{def:Afangpsi} defines the pointed $G$-functor $\Aupom$.  There are five possible cases, depending on $\angordm$ and $\angordn$.  If 
\[\angordm = \vstar \cq \angordn = \vstar, \orspace \angordm = \ang{} = \angordn,\] 
then $\Aupom$ is either the constant functor at the basepoint or the identity functor.  These cases are defined in \cref{A_vstar_n,A_m_vstar,A_empty_empty}.
\item The two nontrivial cases of $\Aupom$ are defined in \cref{A_fangpsi,AF_empty_n}.  \cref{expl:A_fangpsi,expl:AF_empty_n} completely unpack $\Aupom$ in these two cases.
\end{itemize}

\subsection*{$G$-Functors between Categories of Systems}

In what follows, we sometimes restrict attention to morphisms in $\Gsk$ of the form stated in \cref{as:fpsi} below.  We will state it explicitly whenever this assumption is in effect.

\begin{assumption}\label{as:fpsi}
We assume
\[\angordm = \ang{\ordm_i}_{i \in \ufs{p}} \fto{(f, \ang{\psi})} 
\angordn = \ang{\ordn_j}_{j \in \ufs{q}}\]
is a morphism in $\Gsk$ consisting of a reindexing injection 
\[\ufs{p} \fto{f} \ufs{q} \withspace p>0\] 
and a morphism 
\[f_*\angordm =  \bang{\ordm_{\finv(j)}}_{j \in \ufs{q}} 
\fto{\ang{\psi} = \ang{\psi_j}_{j \in \ufs{q}}} \angordn \inspace \Fsk^{(q)},\]
as defined in \cref{reindexing_functor,fangpsi}.
\end{assumption}

The following definition constructs the first pointed $G$-functor that comprises $\sys{\A}{(f,\angpsi)}$.
It uses the pointed $G$-category $\Aangordm$ of $\angordm$-systems in an $\Op$-pseudoalgebra $(\A,\gaA,\phiA)$, with the basepoint given by the base $\angordm$-system $(\zero,1_\zero)$; see \cref{def:nsystem,def:nsystem_morphism,def:Aangordn_system,def:Aangordn_gcat}.

\begin{definition}[The $G$-Functor $\ftil$]\label{def:ftil_functor}
Under \cref{as:OpA,as:fpsi},  the pointed $G$-functor
\begin{equation}\label{ftil_functor}
\Aangordm \fto{\ftil} \Afangordm
\end{equation}
is defined as follows.
\begin{description}
\item[Objects] Suppose $(a,\glu) \in \Aangordm$ is an $\angordm$-system in $\A$.  The $f_*\angordm$-system
\begin{equation}\label{ftil_aglu}
\ftil(a,\glu) = \left(\atil, \glutil\right) \in \Afangordm
\end{equation}
is defined as follows.
\begin{description}
\item[Component objects] Given a marker
\begin{equation}\label{angs_finv}
\ang{s} = \bang{s_j \subseteq \ufs{m}_{\finv(j)}}_{j \in \ufs{q}},
\end{equation}
we first define the marker
\begin{equation}\label{ftil_angs}
\ftil_*\ang{s} = \bang{s_{f(i)} \subseteq \ufs{m}_i}_{i \in \ufs{p}}
\end{equation}
obtained from $\ang{s}$ by
\begin{itemize}
\item removing those $s_j$ with $\finv(j) = \emptyset$ and
\item permuting the resulting $p$-tuple according to $\finv$.
\end{itemize}  
We define the $\ang{s}$-component object of $(\atil,\glutil)$ as follows.
\begin{equation}\label{atil_component}
\atil_{\ang{s}} = \begin{cases}
\zero & \text{if $s_j = \emptyset$ for some $j \in \ufs{q}$, and}\\
a_{\ftil_*\ang{s}} & \text{if $s_j \neq \emptyset$ for each $j \in \ufs{q}$.}
\end{cases}
\end{equation}
\begin{itemize}
\item In \cref{atil_component}, $\zero = \gaA_0(*)$ is the basepoint of $\A$, and the first case is forced by the object unity axiom \cref{system_obj_unity}.
\item In the second case, $a_{\ftil_*\ang{s}}$ is the $\ftil_*\ang{s}$-component object of $(a,\glu)$.  Note that for an index $j \in \ufs{q}$ such that $\finv(j) = \emptyset$---which implies $\ufs{m}_{\finv(j)} = \{1\}$---the condition $s_j \neq \emptyset$ means $s_j = \{1\}$.
\end{itemize}
\item[Gluing] Suppose we are given $r \geq 0$, an object $x \in \Op(r)$, a marker $\ang{s}$ as defined in \cref{angs_finv}, an index $k \in \ufs{q}$, and a partition
\[s_k = \coprod_{\ell \in \ufs{r}}\, s_{k,\ell} \subseteq \ufs{m}_{\finv(k)}.\]
There are three possible cases of the corresponding gluing morphism of $(\atil,\glutil)$ as follows.
\begin{itemize}
\item If $s_j \neq \emptyset$ for each $j \in \ufs{q}$ and $\finv(k) \neq \emptyset$, then we define the corresponding gluing morphism of $(\atil,\glutil)$ via the following commutative diagram.
\begin{equation}\label{glutil_component}
\begin{tikzpicture}[vcenter]
\def\v{-1.2}
\draw[0cell=.9]
(0,0) node (a1) {\gaA_r\big(x \sscs \ang{\atil_{\ang{s} \compk\, s_{k,\ell}}}_{\ell \in \ufs{r}}\big)}
(a1)++(5.5,0) node (a2) {\atil_{\ang{s}}}
(a1)++(0,\v) node (b1) {\gaA_r\big(x \sscs \ang{a_{\ftil_*\ang{s} \,\comp_{\finv(k)}\, s_{k,\ell}}}_{\ell \in \ufs{r}}\big)}
(a2)++(0,\v) node (b2) {a_{\ftil_*\ang{s}}}
;
\draw[1cell=.9]
(a1) edge node {\glutil_{x;\, \ang{s},\, k, \ang{s_{k,\ell}}_{\ell \in \ufs{r}}}} (a2)
(b1) edge node {\glu_{x;\, \ftil_*\ang{s},\, \finv(k), \ang{s_{k,\ell}}_{\ell \in \ufs{r}}}} (b2)
(a1) edge[equal,shorten >=-1ex] (b1)
(a2) edge[equal] (b2)
;
\end{tikzpicture}
\end{equation}
The bottom horizontal arrow $\glu_{\cdots}$ in \cref{glutil_component} is the indicated gluing morphism of $(a,\glu)$.
\item If $s_j = \emptyset$ for some $j \in \ufs{q}$, then we define
\begin{equation}\label{glutil_component_ii}
\glutil_{x;\, \ang{s},\, k, \ang{s_{k,\ell}}_{\ell \in \ufs{r}}} = 1_\zero \cn \zero \to \zero.
\end{equation}
This definition is forced by the unity axiom \cref{system_unity_i}.
\item If $s_j \neq \emptyset$ for each $j \in \ufs{q}$ and $\finv(k) = \emptyset$, then we define
\begin{equation}\label{glutil_component_iii}
\glutil_{x;\, \ang{s},\, k, \ang{s_{k,\ell}}_{\ell \in \ufs{r}}} 
= 1_{a_{\ftil_*\ang{s}}} \cn a_{\ftil_*\ang{s}} \to a_{\ftil_*\ang{s}}.
\end{equation}
\cref{expl:ftil_functor} discusses this gluing morphism in more detail.
\end{itemize}
\item[Axioms] The object unity axiom \cref{system_obj_unity} and the unity axioms \cref{system_unity_i,system_unity_iii} of an $f_*\angordm$-system are already incorporated into the definitions \cref{atil_component,glutil_component_ii,glutil_component_iii}.  Using \cref{atil_component,glutil_component,glutil_component_ii,glutil_component_iii}, each of the other axioms of an $f_*\angordm$-system for $(\atil,\glutil)$---namely, naturality \cref{system_naturality}, equivariance \cref{system_equivariance}, associativity \cref{system_associativity}, and commutativity \cref{system_commutativity}---is an instance of the corresponding axiom for $(a,\glu)$.
\end{description}
\item[Morphisms] Suppose we are given a morphism of $\angordm$-systems
\[(a,\glu^a) \fto{\theta} (b,\glu^b).\]
The morphism of $f_*\angordm$-systems
\[(\atil,\glutil^a) \fto{\ftil(\theta) = \thatil} (\btil,\glutil^b)\]
has, for each marker $\ang{s}$ as defined in \cref{angs_finv}, $\ang{s}$-component morphism defined as follows.
\begin{equation}\label{thatil_component}
\begin{split}
&\big(\atil_{\ang{s}} \fto{\thatil_{\ang{s}}} \btil_{\ang{s}}\big)\\ 
&= \begin{cases}
\zero \fto{1_\zero} \zero & \text{if $s_j = \emptyset$ for some $j \in \ufs{q}$, and}\\
a_{\ftil_*\ang{s}} \fto{\theta_{\ftil_*\ang{s}}} b_{\ftil_*\ang{s}} & \text{if $s_j \neq \emptyset$ for each $j \in \ufs{q}$.}
\end{cases}
\end{split}
\end{equation}
In \cref{thatil_component}, $\theta_{\ftil_*\ang{s}}$ is the $\ftil_*\ang{s}$-component morphism of $\theta$.  The first case of \cref{thatil_component} is forced by the unity axiom \cref{nsystem_mor_unity}.  The compatibility axiom \cref{nsystem_mor_compat} for $\thatil$ follows from the same axiom for $\theta$ and the definitions \cref{atil_component,glutil_component,glutil_component_ii,glutil_component_iii,thatil_component}.
\item[Pointed $G$-functoriality]  The functoriality of $\ftil$ follows from the definition \cref{thatil_component} of $\thatil_{\ang{s}}$. The functor $\ftil$ is pointed---that is, it sends the base $\angordm$-system $(\zero,1_\zero)$ to the base $f_*\angordm$-system---by \cref{atil_component,glutil_component,glutil_component_ii,glutil_component_iii}.  The pointed functor $\ftil$ is $G$-equivariant
\begin{itemize} 
\item on component objects by \cref{ga_scomponent,atil_component},
\item on gluing morphisms by \cref{ga_gluing,glutil_component,glutil_component_ii,glutil_component_iii}, and
\item on morphisms by \cref{gtheta_angs,thatil_component},
\end{itemize}
along with the functoriality of the $g$-action on $\A$ and the fact that $\zero \in \A$ is $G$-fixed.
\end{description}
This finishes the construction of the pointed $G$-functor $\ftil$ in \cref{ftil_functor}.
\end{definition}

\begin{explanation}\label{expl:ftil_functor}
The gluing morphism $\glutil$ in \cref{glutil_component_iii} is well defined for the following reasons.  This concerns the case where $s_j \neq \emptyset$ for each $j \in \ufs{q}$---which implies $r > 0$---and $\finv(k) = \emptyset$.  In the partition
\[\emptyset \neq s_k = \coprod_{\ell \in \ufs{r}}\, s_{k,\ell} \subseteq \ufs{m}_{\finv(k)} = \ufs{m}_{\emptyset} = \{1\},\]
it must be the case that 
\begin{equation}\label{skj}
s_{k,\ell} = \begin{cases}
s_k = \{1\} & \text{for one index $\jm \in \ufs{r}$, and}\\
\emptyset & \text{if $\ell \in \ufs{r} \setminus \{\jm\}$.}
\end{cases}
\end{equation}
By the definition \cref{atil_component} of $\atil_{\ang{s}}$, it follows that
\[\atil_{\ang{s} \compk\, s_{k,\ell}} = \begin{cases}
a_{\ftil_*\ang{s}} & \text{if $\ell = \jm$, and}\\
\zero & \text{if $\ell \in \ufs{r} \setminus \{\jm\}$}.
\end{cases}\]
Thus, the gluing morphism under consideration has the following domain and codomain.
\[\begin{tikzpicture}[vcenter]
\draw[0cell=.9]
(0,0) node (a1) {\gaA_r\big(x \sscs \zero, \ldots, \zero, a_{\ftil_*\ang{s}}, \zero, \ldots, \zero \big)}
(a1)++(5,0) node (a2) {\atil_{\ang{s}}}
(a1)++(0,-1) node (b1) {a_{\ftil_*\ang{s}}}
(a2)++(0,-1) node (b2) {a_{\ftil_*\ang{s}}}
;
\draw[1cell=.9]
(a1) edge node {\glutil_{x; \ang{s},\, k, \ang{s_{k,\ell}}_{\ell \in \ufs{r}}}} (a2)
(a1) edge[equal,shorten <=-.5ex] (b1)
(a2) edge[equal] (b2)
;
\end{tikzpicture}\]
The left equality follows from $r-1$ applications of the basepoint axiom \cref{pseudoalg_basept_axiom} for $\A$, the action unity axiom \cref{pseudoalg_action_unity}, and $\Op(1) = \{\opu\}$.  Thus, the definition in \cref{glutil_component_iii},
\[\glutil_{x;\, \ang{s},\, k, \ang{s_{k,\ell}}_{\ell \in \ufs{r}}} = 1_{a_{\ftil_*\ang{s}}} \cn 
a_{\ftil_*\ang{s}} \to a_{\ftil_*\ang{s}},\]
is well defined.  In fact, this gluing morphism is forced to be the identity morphism by \cref{expl:nsystem} \pcref{expl:nsystem_iv}.
\end{explanation}

Next, we define the second pointed $G$-functor that comprises $\Afpsi$.

\begin{definition}[The $G$-Functor $\psitil$]\label{def:psitil_functor}
Under \cref{as:OpA,as:fpsi}, the pointed $G$-functor
\begin{equation}\label{psitil_functor}
\Afangordm \fto{\psitil} \Aangordn
\end{equation}
is defined as follows.
\begin{description}
\item[Objects] Suppose $(a,\glu) \in \Afangordm$ is an $f_*\angordm$-system.  The $\angordn$-system
\[\psitil(a,\glu) = (a^{\psitil}, \glu^{\psitil}) \in \Aangordn\]
is defined as follows.
\begin{description}
\item[Component objects] Given a marker
\begin{equation}\label{angs_ordtun}
\ang{s} = \ang{s_j \subseteq \ufs{n}_j}_{j \in \ufs{q}},
\end{equation}
recalling that 
\[\ord{m}_{\finv(j)} \fto{\psi_j} \ord{n}_j\]
is a pointed function for each $j \in \ufs{q}$, we define the marker
\[\psiinv\ang{s} = \bang{\psiinv_j s_j \subseteq \ufs{m}_{\finv(j)}}_{j \in \ufs{q}}.\]
We define the $\ang{s}$-component object of $(a^{\psitil}, \glu^{\psitil})$ as
\begin{equation}\label{apsitil_angs}
a^{\psitil}_{\ang{s}} = a_{\psiinv\ang{s}},
\end{equation}
which is the $\psiinv\ang{s}$-component object of $(a,\glu)$.
\item[Gluing] Suppose we are given $r \geq 0$, an object $x \in \Op(r)$, a marker $\ang{s}$ as defined in \cref{angs_ordtun}, an index $k \in \ufs{q}$, and a partition
\begin{equation}\label{skpartition}
s_k = \coprod_{\ell \in \ufs{r}}\, s_{k,\ell} \subseteq \ufs{n}_k.
\end{equation}
We define the corresponding gluing morphism of $(a^{\psitil}, \glu^{\psitil})$ via the following commutative diagram.
\begin{equation}\label{glupsitil}
\begin{tikzpicture}[vcenter]
\def\v{-1.2}
\draw[0cell=.9]
(0,0) node (a1) {\gaA_r\big(x \sscs \ang{a^{\psitil}_{\ang{s} \compk\, s_{k,\ell}}}_{\ell \in \ufs{r}}\big)}
(a1)++(5.5,0) node (a2) {a^{\psitil}_{\ang{s}}}
(a1)++(0,\v) node (b1) {\gaA_r\big(x \sscs \ang{a_{\psiinv\ang{s} \,\comp_{k}\, \psiinv_k s_{k,\ell}}}_{\ell \in \ufs{r}}\big)}
(a2)++(0,\v) node (b2) {a_{\psiinv\ang{s}}}
;
\draw[1cell=.9]
(a1) edge node {\glu^{\psitil}_{x;\, \ang{s},\, k, \ang{s_{k,\ell}}_{\ell \in \ufs{r}}}} (a2)
(b1) edge node {\glu_{x;\, \psiinv\ang{s},\, k, \ang{\psiinv_k s_{k,\ell}}_{\ell \in \ufs{r}}}} (b2)
(a1) edge[equal,shorten >=-1ex] (b1)
(a2) edge[equal] (b2)
;
\end{tikzpicture}
\end{equation}
\item[Axioms]  Each of the axioms of an $\angordn$-system, \cref{system_obj_unity,system_naturality,system_unity_i,system_unity_iii,system_equivariance,system_associativity,system_commutativity}, holds for $(a^{\psitil}, \glu^{\psitil})$ by the corresponding axiom for $(a,\glu)$ and the fact that taking pre-images preserves the empty set and partitions.
\end{description}
\item[Morphisms] Suppose we are given a morphism of $f_*\angordm$-systems
\[(a,\glu^a) \fto{\theta} (b,\glu^b).\]
The morphism of $\angordn$-systems
\[(a^{\psitil}, \glu^{a, \psitil}) \fto{\psitil(\theta) = \tha^{\psitil}} (b^{\psitil}, \glu^{b, \psitil})\]
has, for each marker $\ang{s}$ as defined in \cref{angs_ordtun}, $\ang{s}$-component morphism defined as the $\psiinv\ang{s}$-component morphism of $\theta$:
\begin{equation}\label{thapsitil_component}
a^{\psitil}_{\ang{s}} = a_{\psiinv\ang{s}} \fto{\tha^{\psitil}_{\ang{s}} \,=\, \theta_{\psiinv\ang{s}}} 
b^{\psitil}_{\ang{s}} = b_{\psiinv\ang{s}}.
\end{equation}
Each of the unity axiom \cref{nsystem_mor_unity} and the compatibility axiom \cref{nsystem_mor_compat} of a morphism of $\angordn$-systems holds for $\tha^{\psitil}$ by the corresponding axiom for $\theta$ and the fact that taking pre-images preserves the empty set and partitions.
\item[Pointed $G$-functoriality]
The functoriality of $\psitil$ follows from the definition \cref{thapsitil_component} of $\tha^{\psitil}_{\ang{s}}$.  The functor $\psitil$ is pointed---that is, it sends the base $f_*\angordm$-system $(\zero,1_\zero)$ to the base $\angordn$-system---by \cref{apsitil_angs,glupsitil}.  The pointed functor $\psitil$ is $G$-equivariant
\begin{itemize} 
\item on component objects by \cref{ga_scomponent,apsitil_angs},
\item on gluing morphisms by \cref{ga_gluing,glupsitil}, and
\item on morphisms by \cref{gtheta_angs,thapsitil_component}.
\end{itemize}
\end{description}
This finishes the construction of the pointed $G$-functor $\psitil$.
\end{definition}

Next, we define the pointed $G$-functor $\Aupom$ associated to a general morphism $\upom$ in $\Gsk$ \cref{Gsk_morphisms}.

\begin{definition}[The $G$-Functor $\Aupom$]\label{def:Afangpsi}
For a $\Tinf$-operad $\Op$ \pcref{as:OpA}, an $\Op$-pseudoalgebra $\A$ \pcref{def:pseudoalgebra}, and a morphism 
\[\angordm \fto{\upom} \angordn \inspace \Gsk,\]
the pointed $G$-functor
\begin{equation}\label{AF}
\Aangordm \fto{\Aupom} \Aangordn
\end{equation}
is defined as follows.  First, suppose $\upom$ has the form
\[\angordm = \ang{\ordm_i}_{i \in \ufs{p}} \fto{\upom = (f,\angpsi)} 
\angordn = \ang{\ordn_j}_{j \in \ufs{q}}\]
with $\angordm, \angordn \in \Gsk \setminus \{\vstar, \ang{}\}$, as defined in \cref{as:fpsi}.  We define $\Aupom$ as the composite pointed $G$-functor
\begin{equation}\label{A_fangpsi}
\begin{tikzpicture}[vcenter]
\def\h{2.5} \def\t{.7}
\draw[0cell]
(0,0) node (a1) {\Aangordm}
(a1)++(\h,0) node (a2) {\Afangordm}
(a2)++(\h,0) node (a3) {\Aangordn}
;
\draw[1cell=.9]
(a1) edge node{\ftil} (a2)
(a2) edge node {\psitil} (a3)
;
\draw[1cell=.9]
(a1) [rounded corners=2pt, shorten <=0ex] |- ($(a2)+(-1,\t)$)
-- node {\Afpsi} ($(a2)+(1,\t)$) -| (a3)
;
\end{tikzpicture}
\end{equation}
where $\ftil$ and $\psitil$ are the pointed $G$-functors in, respectively, \cref{def:ftil_functor,def:psitil_functor}.  \cref{expl:A_fangpsi} discusses the functor $\Afpsi$ in more detail.  

\parhead{Marginal cases}.  The remaining cases of $\Aupom$ are defined as follows.
\begin{itemize}
\item If $\angordm = \vstar \in \Gsk$, then 
\begin{equation}\label{A_vstar_n}
\sys{\A}{\vstar} = \boldone \fto{\Aupom} \Aangordn
\end{equation}
sends the unique object of $\boldone$ to the base $\angordn$-system $(\zero,1_\zero)$.  This defines a pointed $G$-functor because $(\zero,1_\zero)$ is $G$-fixed.
\item If $\angordn = \vstar \in \Gsk$, then 
\begin{equation}\label{A_m_vstar}
\Aangordm \fto{\Aupom} \sys{\A}{\vstar} = \boldone
\end{equation}
is the unique pointed $G$-functor to the terminal pointed $G$-category.
\item Suppose $\angordm = \angordn = \ang{}$, the empty sequence.  By \cref{Gsk_empty_mor}, $\upom \cn \ang{} \to \ang{}$ is either
\begin{itemize}
\item the identity morphism, $1_{\ang{}}$, or
\item the 0-morphism, $\ang{} \to \vstar \to \ang{}$.
\end{itemize}  
In these two cases, we define 
\begin{equation}\label{A_empty_empty}
\sys{\A}{\ang{}} = \A \fto{\Aupom} \sys{\A}{\ang{}} = \A
\end{equation}
as, respectively, the identity functor and the constant functor at the $G$-fixed basepoint $\zero \in \A$.
\item Suppose $\angordm = \ang{}$ and $\angordn \in \Gsk \setminus \{\vstar, \ang{}\}$ has length $q>0$.  Then the morphism $\upom$ factors as the following composite in $\Gsk$ \cref{Gsk_composite}.
\begin{equation}\label{Fiqangpsi}
\begin{tikzpicture}[vcenter]
\def\t{.7} \def\h{3}
\draw[0cell]
(0,0) node (a1) {\ang{}}
(a1)++(\h,0) node (a2) {\ordtu{1}_{j \in \ufs{q}}}
(a2)++(\h,0) node (a3) {\angordn}
;
\draw[1cell=.9]
(a1) edge node {(\im_q,\ang{1_{\ord{1}}}_{j \in \ufs{q}})} (a2)
(a2) edge node {(1_{\ufs{q}}, \ang{\psi_j}_{j \in \ufs{q}})} (a3)
;
\draw[1cell=1]
(a1) [rounded corners=2pt, shorten <=0ex] |- ($(a2)+(-1,\t)$)
-- node {\upom \,=\, (\im_q, \ang{\psi_j}_{j \in \ufs{q}})} ($(a2)+(1,\t)$) -| (a3)
;
\end{tikzpicture}
\end{equation}
In this composite,
\begin{itemize}
\item $\im_q \cn \emptyset \to \ufs{q}$ is the unique reindexing injection, and
\item $\psi_j \cn \ord{1} \to \ord{n}_j$ is a pointed function for each $j \in \ufs{q}$.
\end{itemize}
We define $\Aupom$ as the following composite pointed $G$-functor.
\begin{equation}\label{AF_empty_n}
\begin{tikzpicture}[vcenter]
\def\t{.7}
\draw[0cell]
(0,0) node (a1) {\sys{\A}{\ang{}}}
(a1)++(.9,0) node (a1') {\phantom{\A}}
(a1')++(0,.02) node (a1'') {\A}
(a1')++(2,0) node (a2) {\phantom{\sys{\A}{\ordtu{1}_{j \in \ufs{q}}}}}
(a2)++(0,-.03) node (a2') {\sys{\A}{\ordtu{1}_{j \in \ufs{q}}}}
(a2)++(3.5,0) node (a3) {\Aangordn}
;
\draw[1cell=.9]
(a1) edge[equal] (a1')
(a1') edge node {\iso} (a2)
(a2) edge node {\sys{\A}{(1_{\ufs{q}}, \ang{\psi_j}_{j \in \ufs{q}})}} (a3)
;
\draw[1cell=1]
(a1) [rounded corners=2pt, shorten <=0ex] |- ($(a2)+(0,\t)$)
-- node[pos=.3] {\Aupom} ($(a2)+(1,\t)$) -| (a3)
;
\end{tikzpicture}
\end{equation}
In \cref{AF_empty_n},
\begin{itemize}
\item the functor labeled $\iso$ is the pointed $G$-isomorphism in \cref{Aordtuone}, and
\item $\sys{\A}{(1_{\ufs{q}}, \ang{\psi_j}_{j \in \ufs{q}})}$ is the pointed $G$-functor in \cref{A_fangpsi} with $f = 1_{\ufs{q}}$.
\end{itemize}
\cref{expl:AF_empty_n} discusses the pointed $G$-functor in \cref{AF_empty_n} in more detail.
\end{itemize}
This finishes the construction of the pointed $G$-functor $\Aupom$.  

\parhead{Strong variant}.  We define the strong variant
\begin{equation}\label{AF_sg}
\Aangordmsg \fto{\Aupomsg} \Aangordnsg
\end{equation}
by restricting the definitions above to strong $\angordm$-systems.  This is well defined because both $\ftil$ and $\psitil$ preserve invertibility of gluing morphisms by, respectively, \cref{glutil_component,glutil_component_ii,glutil_component_iii} and \cref{glupsitil}.
\end{definition}

\subsection*{Unpacking $\Aupom$}

The rest of this section unpacks the pointed functor $\Aupom$ in the two nontrivial cases, \cref{A_fangpsi,AF_empty_n}.

\begin{explanation}\label{expl:A_fangpsi}
We unpack the pointed $G$-functor in \cref{A_fangpsi}
\[\Afpsi \cn \Aangordm \fto{\ftil} \Afangordm \fto{\psitil} \Aangordn\]
using \cref{def:ftil_functor,def:psitil_functor}.

\parhead{Objects}.  For an $\angordm$-system $(a,\glu) \in \Aangordm$, the $\angordn$-system
\[\Afpsi(a,\glu) = \psitil \ftil(a,\glu) = (\atil^{\psitil}, \glutil^{\psitil}) \in \Aangordn\]
has, for each marker $\ang{s} = \ang{s_j \subseteq \ufs{n}_j}_{j \in \ufs{q}}$, the following $\ang{s}$-component object.
\begin{equation}\label{A_fangpsi_obj_component}
\begin{split}
&\left(\Afpsi(a,\glu)\right)_{\ang{s}} \\
&= \atil^{\psitil}_{\ang{s}} = \atil_{\psiinv\ang{s}} = \atil_{\ang{\psiinv_j s_j}_{j \in \ufs{q}}}\\
&= \begin{cases} 
\zero & \text{if $\psiinv_j s_j = \emptyset$ for some $j \in \ufs{q}$, and}\\
a_{\ang{\psiinv_{f(i)} s_{f(i)}}_{i \in \ufs{p}}} & \text{if $\psiinv_j s_j \neq \emptyset$ for each $j \in \ufs{q}$.}
\end{cases}
\end{split}
\end{equation}
In the context of \cref{skpartition}, if $\psiinv_j s_j \neq \emptyset$ for each $j \in \ufs{q}$ and $\finv(k) \neq \emptyset$, then the gluing morphism of $\Afpsi(a,\glu)$ is given as follows.
\begin{equation}\label{A_fangpsi_gluing}
\begin{tikzpicture}[vcenter]
\def\v{-1} \def\u{-1.1} \def\h{5}
\draw[0cell=.85]
(0,0) node (a1) {\gaA_r\big(x \sscs \ang{\atil^{\psitil}_{\ang{s} \compk\, s_{k,\ell}}}_{\ell \in \ufs{r}}\big)}
(a1)++(\h,0) node (a2) {\atil^{\psitil}_{\ang{s}}}
(a1)++(0,\v) node (b1) {\gaA_r\big(x \sscs \ang{a_{\ang{\psiinv_{f(i)} s_{f(i)}}_{i \in \ufs{p}} \,\comp_{\finv(k)}\, \psiinv_k s_{k,\ell}}}_{\ell \in \ufs{r}}\big)}
(a2)++(0,\v) node (b2) {a_{\ang{\psiinv_{f(i)} s_{f(i)}}_{i \in \ufs{p}}}}
;
\draw[1cell=.9]
(a1) edge node {\glutil^{\psitil}_{x;\, \ang{s},\, k, \ang{s_{k,\ell}}_{\ell \in \ufs{r}}}} (a2)
(a1) edge[equal,shorten >=-1ex] (b1)
(a2) edge[equal] (b2)
;
\draw[1cell=.9]
(b1) [rounded corners=2pt, shorten <=0ex] |- ($(b1)+(.3*\h,\u)$)
-- node {\glu_{x;\, \ang{\psiinv_{f(i)} s_{f(i)}}_{i \in \ufs{p}},\, \finv(k), \ang{\psiinv_k s_{k,\ell}}_{\ell \in \ufs{r}}}} ($(b1)+(.7*\h,\u)$) -| (b2)
;
\end{tikzpicture}
\end{equation}
Other gluing morphisms of the $\angordn$-system $\Afpsi(a,\glu)$ are identity morphisms, as stated in \cref{glutil_component_ii,glutil_component_iii}.

\parhead{Morphisms}.  For a morphism of $\angordm$-systems
\[(a,\glu^a) \fto{\theta} (b,\glu^b),\]
the morphism of $\angordn$-systems
\[\Afpsi(\theta) = \thatil^{\psitil} \cn (\atil^{\psitil}, \glutil^{a,\psitil}) 
\to (\btil^{\psitil}, \glutil^{b,\psitil})\]
has the following $\ang{s}$-component morphism.
\begin{equation}\label{A_fangpsi_mor_component}
\begin{split}
&\left(\Afpsi(\theta)\right)_{\ang{s}} \\
&= \thatil^{\psitil}_{\ang{s}} = \thatil_{\psiinv\ang{s}} = \thatil_{\ang{\psiinv_j s_j}_{j \in \ufs{q}}}\\
&= \begin{cases} 
1_\zero & \text{if $\psiinv_j s_j = \emptyset$ for some $j \in \ufs{q}$, and}\\
\theta_{\ang{\psiinv_{f(i)} s_{f(i)}}_{i \in \ufs{p}}} & \text{if $\psiinv_j s_j \neq \emptyset$ for each $j \in \ufs{q}$.}
\end{cases}
\end{split}
\end{equation}
This finishes our unpacking of the pointed $G$-functor $\Afpsi$ in \cref{A_fangpsi}. 
\end{explanation}

\begin{explanation}\label{expl:AF_empty_n}
We unpack the pointed $G$-functor in \cref{AF_empty_n}
\[\Aupom \cn \sys{\A}{\ang{}} = \A \fiso \sys{\A}{\ordtu{1}}_{j \in \ufs{q}} \fto{\sys{\A}{(1_{\ufs{q}}, \ang{\psi_j}_{j \in \ufs{q}})}} \Aangordn.\]
Recall from \cref{Fiqangpsi} that $\upom$ has the form
\[\upom = (\im_q, \ang{\psi_j}_{j \in \ufs{q}}) \cn \ang{} \to \angordn,\]
with 
\begin{itemize}
\item $\angordm = \ang{}$,
\item $\angordn \in \Gsk \setminus \{\vstar, \ang{}\}$ of length $q>0$,
\item $\im_q \cn \emptyset \to \ufs{q}$ the unique reindexing injection, and
\item $\psi_j \cn \ord{1} \to \ord{n}_j$ a pointed function for each $j \in \ufs{q}$.
\end{itemize}  

\parhead{Objects}.  For an object $a \in \A$ and a marker $\ang{s} = \ang{s_j \subseteq \ufs{n}_j}_{j \in \ufs{q}}$, the $\ang{s}$-component object of the $\angordn$-system $(\Aupom)(a)$ is given as follows.
\begin{equation}\label{AFa_angs}
\big((\Aupom)(a)\big)_{\ang{s}} = \begin{cases}
\zero & \text{if $\psi_j(1) \not\in s_j$ for some $j \in \ufs{q}$, and}\\
a & \text{if $\psi_j(1) \in s_j$ for each $j \in \ufs{q}$.}
\end{cases}
\end{equation}
Each gluing morphism of the $\angordn$-system $(\Aupom)(a)$ is the identity morphism.  We regard \cref{AFa_angs} as a special case of \cref{A_fangpsi_obj_component}.

\parhead{Morphisms}.  For a morphism $\theta \cn a \to b$ in $\A$, the $\ang{s}$-component of the morphism 
\[(\Aupom)(a) \fto{(\Aupom)(\theta)} (\Aupom)(b) \inspace \Aangordn\]
is given as follows.
\begin{equation}\label{AF_theta_angs}
\big((\Aupom)(\theta)\big)_{\ang{s}} = \begin{cases}
1_\zero & \text{if $\psi_j(1) \not\in s_j$ for some $j \in \ufs{q}$, and}\\
\theta & \text{if $\psi_j(1) \in s_j$ for each $j \in \ufs{q}$.}
\end{cases}
\end{equation}
We regard \cref{AF_theta_angs} as a special case of \cref{A_fangpsi_mor_component}.  This finishes our unpacking of the pointed $G$-functor $\Aupom$ in \cref{AF_empty_n}.
\end{explanation}

\section{$\Gskg$-Categories from Operadic Pseudoalgebras: Functoriality}
\label{sec:jemg_morphisms_ii}

As a reminder, \cref{as:OpA} is in effect throughout this chapter, and $(\A,\gaA,\phiA)$ is an arbitrary $\Op$-pseudoalgebra  (\cref{def:pseudoalgebra}).  We are in the process of constructing the object assignment of the $\Gcat$-multifunctor
\[\MultpsO \fto{\Jgo} \GGCat.\]
\cref{sec:jemg_objects,sec:jemg_morphisms} construct, respectively, the object and morphism assignments of 
\[(\Gsk,\vstar) \fto{\Jgo\A = \Adash} (\Gcatst,\boldone).\]
This section proves that $\Adash$ is a $\Gskg$-category (\cref{expl:ggcat_obj})---that is, a pointed functor---finishing the construction of the object assignment, $\A \mapsto \Jgo\A$, of $\Jgo$.  The strong variant $\Asgdash$, which involves strong systems, is also a $\Gskg$-category. 

\begin{lemma}\label{A_ptfunctor}
Under \cref{as:OpA}, the object and morphism assignments 
\[\angordn \mapsto \Aangordn \andspace \upom \mapsto \Aupom\]
in, respectively,  \cref{def:Aangordn_gcat,def:Afangpsi} define a pointed functor
\[(\Gsk,\vstar) \fto{\Adash} (\Gcatst,\boldone).\]
Moreover, the strong variant $\Asgdash$, defined in \cref{sgAordnbe,AF_sg}, is also a pointed functor.
\end{lemma}

\begin{proof}
We prove that $\Adash$ is a pointed functor.  The same proof applies to the strong case, $\Asgdash$, by restricting to strong systems, which are those with invertible gluing morphisms.

\parhead{Pointedness}.  The assignment $\Adash$ is pointed because $\sys{\A}{\vstar} = \boldone$ by definition \cref{A_angordn}.  Next, we prove that $\Adash$ is a functor.

\parhead{Preservation of identities}. If $\upom$ is the identity morphism of either $\vstar$ or $\ang{} \in \Gsk$, then, by \cref{A_vstar_n,A_empty_empty}, $\Aupom$ is the identity functor of, respectively, $\sys{\A}{\vstar} = \boldone$ or $\sys{\A}{\ang{}} = \A$.  

For the identity morphism of an object $\angordm \in \Gsk \setminus \{\vstar,\ang{}\}$, we observe the following.
\begin{itemize}
\item Suppose, in \cref{def:ftil_functor}, that the reindexing injection is $f = 1_{\ufs{p}} \cn \ufs{p} \to \ufs{p}$ .  Then
\[\Aangordm \fto{\ftil} \Aangordm\]
is the identity functor because 
\[\ftil_*\ang{s} = \ang{s_{f(i)}}_{i \in \ufs{p}} = \ang{s}\]
for each marker $\ang{s} = \ang{s_i \subseteq \ufs{m}_i}_{i \in \ufs{p}}$.
\item Suppose, in \cref{def:psitil_functor}, that the morphism is the identity
\[(f,\angpsi) = \left(1_{\ufs{p}}, \ang{1_{\ord{m}_i}}_{i \in \ufs{p}}\right) \cn \angordm \to \angordm\]
Then
\[\Aangordm \fto{\psitil} \Aangordm\]
is the identity functor because 
\[\psiinv\ang{s} = \ang{\psiinv_i s_i}_{i \in \ufs{p}} = \ang{s}\]
for each marker $\ang{s}$.
\end{itemize}
Thus, if $(f,\angpsi)$ is an identity morphism in $\Gsk$, then $\Afpsi = \psitil\ftil$, as defined in \cref{A_fangpsi}, is the composite of two identity functors.

\parhead{Preservation of composition}.  To show that $\Adash$ preserves composition, we consider composable morphisms 
\[\angordm \fto{\upom} \angordn \fto{\upom'} \angordl \inspace \Gsk\]
in all possible cases as follows.

\parhead{Case 1}.  Suppose $\angordm$, $\angordn$, or $\angordl$ is the initial-terminal basepoint $\vstar \in \Gsk$.  Then, by \cref{A_vstar_n,A_m_vstar}, each of $\sys{\A}{(\upom' \upom)}$ and $(\Aupomp)(\Aupom)$ is the constant functor at the base $\angordl$-system $(\zero,1_\zero) \in \Aangordl$.  Thus, we assume that $\angordm, \angordn, \angordl \neq \vstar$ for the rest of this proof.

\parhead{Case 2}.  Suppose $\angordl = \ang{}$.  Then, by \cref{Gsk_empty_mor},
\[\angordm = \angordn = \ang{}.\] 
Each of $\upom$ and $\upom'$ is either (i) the identity morphism, $1_{\ang{}}$, or (ii) the 0-morphism, $\ang{} \to \vstar \to \ang{}$.
\begin{itemize}
\item If either $\upom$ or $\upom'$ is the 0-morphism, then so is the composite $\upom'\upom$.   By \cref{A_empty_empty}, each of $\sys{\A}{(\upom' \upom)}$ and $(\Aupomp) (\Aupom)$ is the constant functor at the basepoint $\zero \in \A = \A\ang{}$.
\item If both $\upom$ and $\upom'$ are equal to $1_{\ang{}}$, then so is $\upom'\upom$.  In this case, each of $\sys{\A}{(\upom' \upom)}$ and $(\Aupomp) (\Aupom)$ is the identity functor.
\end{itemize}
For the rest of this proof, we further assume that $\angordl \in \Gsk \setminus \{\vstar, \ang{}\}$ has length $t > 0$.

\parhead{Case 3}.  Suppose $\angordn = \ang{}$.  Then $\angordm = \ang{}$, and $\upom$ is either $1_{\ang{}}$ or the 0-morphism.
\begin{itemize}
\item If $\upom = 1_{\ang{}}$, then $\Aupom = 1$ and
\[\sys{\A}{(\upom'\upom)} = \Aupomp = (\Aupomp) (\Aupom).\]
\item If $\upom$ is the 0-morphism, then so is $\upom'\upom$.  Each of $\sys{\A}{(\upom' \upom)}$ and $(\Aupomp) (\Aupom)$ is the constant functor at the base $\angordl$-system $(\zero,1_\zero) \in \Aangordl$.
\end{itemize}
For the rest of this proof, we further assume that $\angordn \in \Gsk \setminus \{\vstar,\ang{}\}$ has length $q > 0$.  Thus, $\upom'$ has the form
\begin{equation}\label{morphism_Fggprime}
\angordn = \ang{\ordn_j}_{j \in \ufs{q}} \fto{\upom' = (f',\ang{\psi'_c}_{c \in \ufs{t}})} 
\angordl = \ang{\ordl_c}_{c \in \ufs{t}}.
\end{equation}
It consists of
\begin{itemize}
\item a reindexing injection $f' \cn \ufs{q} \to \ufs{t}$ and
\item a pointed function 
\[\ord{n}_{(f')^{-1}c} \fto{\psi'_c} \ord{\ell}_c\] 
for each $c \in \ufs{t}$,
\end{itemize}
where, by \cref{ordn_empty},
\[\ord{n}_{\emptyset} = \ord{1} \ifspace \fpinv c = \emptyset.\]

\parhead{Case 4}.  Suppose $\angordm = \ang{}$.  As we discuss in \cref{Fiqangpsi}, $\upom$ has the form 
\begin{equation}\label{morphism_Fgg}
\ang{} \fto{\upom = (\im_q, \ang{\psi_j}_{j \in \ufs{q}})} \angordn.
\end{equation}
It consists of
\begin{itemize}
\item the unique reindexing injection $\im_q \cn \emptyset \to \ufs{q}$ and
\item a pointed function $\psi_j \cn \ord{1} \to \ord{n}_j$ for each $j \in \ufs{q}$.
\end{itemize}
Combining \cref{morphism_Fggprime,morphism_Fgg}, the composite $\upom'\upom$ has the following form.
\begin{equation}\label{morphism_FggprimeFgg}
\begin{tikzpicture}[vcenter]
\def\t{.7}
\draw[0cell]
(0,0) node (a1) {\ang{}}
(a1)++(3.5,0) node (a2) {\angordn}
(a2)++(3.5,0) node (a3) {\angordl}
;
\draw[1cell=.9]
(a1) edge node {\upom = (\im_q, \ang{\psi_j}_{j \in \ufs{q}})} (a2)
(a2) edge node {\upom' = (f', \ang{\psi'_c}_{c\in \ufs{t}})} (a3)
;
\draw[1cell=1]
(a1) [rounded corners=2pt, shorten <=0ex] |- ($(a2)+(-1,\t)$)
-- node {\upom' \upom = (\im_t, \ang{\psi'_c \psi_{\fpinv c}}_{c \in \ufs{t}})} ($(a2)+(1,\t)$) -| (a3)
;
\end{tikzpicture}
\end{equation}
It consists of
\begin{itemize}
\item the unique reindexing injection $\im_t \cn \emptyset \to \ufs{t}$ and
\item the composite pointed function
\[\ord{1} \fto{\psi_{\fpinv c}} \ord{n}_{\fpinv c} \fto{\psi'_c} \ord{\ell}_c\]
for each $c \in \ufs{t}$,
where, by \cref{ordn_empty},
\[\ord{1} \fto{\psi_{\emptyset} = 1_{\ord{1}}} \ord{n}_{\emptyset} = \ord{1}
\ifspace \fpinv c = \emptyset.\]
\end{itemize}

We want to show that the following diagram of pointed $G$-functors commutes.
\begin{equation}\label{A_FggprimeFgg}
\begin{tikzpicture}[vcenter]
\def\t{.7}
\draw[0cell]
(0,0) node (a1) {\sys{\A}{\ang{}} = \A}
(a1)++(3.5,0) node (a2)  {\Aangordn}
(a2)++(3.5,0) node (a3) {\Aangordl}
;
\draw[1cell=.9]
(a1) edge node {\sys{\A}{(\im_q, \ang{\psi_j}_{j \in \ufs{q}})}} (a2)
(a2) edge node {\sys{\A}{(f', \ang{\psi'_c}_{c\in \ufs{t}})}} (a3)
;
\draw[1cell=1]
(a1) [rounded corners=3pt, shorten <=0ex] |- ($(a2)+(-1,\t)$)
-- node {\sys{\A}{(\upom'\upom)} \,=\, \sys{\A}{(\im_t, \ang{\psi'_c \psi_{\fpinv c}}_{c \in \ufs{t}})}} ($(a2)+(1,\t)$) -| (a3)
;
\end{tikzpicture}
\end{equation}
To check the commutativity of \cref{A_FggprimeFgg} on objects, suppose we are given an object $a \in \A$ and a marker $\ang{s} = \ang{s_c \subseteq \ufs{\ell}_c}_{c \in \ufs{t}}$.
Using \cref{A_fangpsi_obj_component,AFa_angs}, we compute the $\ang{s}$-component object of the $\angordl$-system
\[(\Aupomp)(\Aupom)(a) \in \Aangordl\]
as follows.  
\begin{itemize}
\item If $(\psi'_c)^{-1} s_c = \emptyset$ for some $c \in \ufs{t}$, then
\begin{equation}\label{AF'AFa_angs_i}
\left((\Aupomp)(\Aupom)(a)\right)_{\ang{s}} = \zero.
\end{equation}
\item If $(\psi'_c)^{-1} s_c \neq \emptyset$ for each $c \in \ufs{t}$, then
\begin{equation}\label{AF'AFa_angs_ii}
\begin{split}
&\left((\Aupomp)(\Aupom)(a)\right)_{\ang{s}} \\
&= \left((\Aupom)(a)\right)_{\ang{(\psi'_{f'(j)})^{-1} s_{f'(j)}}_{j \in \ufs{q}}}\\
&= \begin{cases}
\zero & \text{if $\psi_j(1) \not\in (\psi'_{f'(j)})^{-1} s_{f'(j)}$ for some $j \in \ufs{q}$, and}\\
a & \text{if $\psi_j(1) \in (\psi'_{f'(j)})^{-1} s_{f'(j)}$ for each $j \in \ufs{q}$.}
\end{cases}
\end{split}
\end{equation}
\end{itemize}
For an index $c \in \ufs{t}$ such that $\fpinv c = \emptyset$---which implies $\psi_{\fpinv c} = 1_{\ord{1}}$---the condition $(\psi'_c)^{-1} s_c \neq \emptyset$ means
\[(\psi'_c)^{-1} s_c = \ufs{n}_{\emptyset} = \{1\},\]
which is equivalent to 
\begin{equation}\label{psiprimec_one}
\psi'_c(1) = \psi'_c \psi_{\fpinv c} (1) \in s_c.
\end{equation}  
Thus, we can combine \cref{AF'AFa_angs_i,AF'AFa_angs_ii} as follows.
\begin{equation}\label{AF'AFa_angs}
\begin{split}
&\left((\Aupomp)(\Aupom)(a)\right)_{\ang{s}} \\
&= \begin{cases}
\zero & \text{if $\psi'_c \psi_{\fpinv c}(1) \not\in s_c$ for some $c \in \ufs{t}$, and}\\
a & \text{if $\psi'_c \psi_{\fpinv c}(1) \in s_c$ for each $c \in \ufs{t}$.}
\end{cases}
\end{split}
\end{equation}
By \cref{AFa_angs,A_FggprimeFgg}, the right-hand side of \cref{AF'AFa_angs} is also equal to $\sys{\A}{(\upom'\upom)}(a)_{\ang{s}}$.
\begin{itemize}
\item By \cref{expl:A_fangpsi,expl:AF_empty_n}, each gluing morphism of each of the $\angordl$-systems $(\Aupomp)(\Aupom)(a)$ and $\sys{\A}{(\upom'\upom)}(a)$ is an identity morphism.  Thus, these two $\angordl$-systems are equal, which proves that the diagram \cref{A_FggprimeFgg} commutes on objects.
\item The same proof as above, using \cref{A_fangpsi_mor_component,AF_theta_angs} instead of \cref{A_fangpsi_obj_component,AFa_angs}, shows that the diagram \cref{A_FggprimeFgg} commutes on morphisms.
\end{itemize}
This finishes the proof that $\Adash$ preserves composition in the present case.  

\parhead{Case 5}.  Suppose $\angordm \in \Gsk \setminus \{\vstar, \ang{}\}$ has length $p>0$.  Thus, $\upom$ has the form
\begin{equation}\label{Fgg_general}
\angordm = \ang{\ordm_i}_{i \in \ufs{p}} \fto{\upom = (f,\ang{\psi_j}_{j \in \ufs{q}})} 
\angordn = \ang{\ordn_j}_{j \in \ufs{q}},
\end{equation}
as stated in \cref{as:fpsi}.  Combining \cref{Gsk_composite,morphism_Fggprime,Fgg_general}, the composite $\upom'\upom$ has the following form.
\begin{equation}\label{FggprimeFgg_general}
\begin{tikzpicture}[vcenter]
\def\t{.7}
\draw[0cell]
(0,0) node (a1) {\angordm}
(a1)++(3.5,0) node (a2) {\angordn}
(a2)++(3.5,0) node (a3) {\angordl}
;
\draw[1cell=.9]
(a1) edge node {\upom = (f, \ang{\psi_j}_{j \in \ufs{q}})} (a2)
(a2) edge node {\upom' = (f', \ang{\psi'_c}_{c\in \ufs{t}})} (a3)
;
\draw[1cell=1]
(a1) [rounded corners=2pt, shorten <=0ex] |- ($(a2)+(-1,\t)$)
-- node {\upom' \upom = (f'f, \ang{\psi'_c \psi_{\fpinv c}}_{c \in \ufs{t}})} ($(a2)+(1,\t)$) -| (a3)
;
\end{tikzpicture}
\end{equation}
We want to show that the following diagram of pointed $G$-functors commutes.
\begin{equation}\label{A_FggprimeFgg_general}
\begin{tikzpicture}[vcenter]
\def\t{.8}
\draw[0cell]
(0,0) node (a1) {\Aangordm}
(a1)++(4.3,0) node (a2) {\Aangordn}
(a2)++(4.3,0) node (a3) {\Aangordl}
;
\draw[1cell=.9]
(a1) edge node {\Aupom = \sys{\A}{(f, \ang{\psi_j}_{j \in \ufs{q}})}} (a2)
(a2) edge node {\Aupomp = \sys{\A}{(f', \ang{\psi'_c}_{c\in \ufs{t}})}} (a3)
;
\draw[1cell=1]
(a1) [rounded corners=2pt, shorten <=0ex] |- ($(a2)+(-1,\t)$)
-- node {\sys{\A}{(\upom' \upom)} = \sys{\A}{(f'f, \ang{\psi'_c \psi_{\fpinv c}}_{c \in \ufs{t}})}} ($(a2)+(1,\t)$) -| (a3)
;
\end{tikzpicture}
\end{equation}

\parhead{Component objects}.  To check the commutativity of \cref{A_FggprimeFgg_general} on objects, suppose we are given an $\angordm$-system $(a,\glu) \in \Aangordm$ and a marker $\ang{s} = \ang{s_c \subseteq \ufs{\ell}_c}_{c \in \ufs{t}}$.  Using \cref{A_fangpsi_obj_component}, we compute the $\ang{s}$-component object of the $\angordl$-system 
\[(\Aupomp)(\Aupom)(a,\glu) \in \Aangordl\]
as follows.  
\begin{itemize}
\item If $(\psi'_c)^{-1} s_c = \emptyset$ for some $c \in \ufs{t}$, then
\begin{equation}\label{AF'AF_angs_i}
\left((\Aupomp)(\Aupom)(a,\glu)\right)_{\ang{s}} = \zero.
\end{equation}
\item If $(\psi'_c)^{-1} s_c \neq \emptyset$ for each $c \in \ufs{t}$, then
\begin{equation}\label{AF'AF_angs_ii}
\begin{split}
&\left((\Aupomp)(\Aupom)(a,\glu)\right)_{\ang{s}} \\
&= \left((\Aupom)(a,\glu)\right)_{\ang{(\psi'_{f'(j)})^{-1} s_{f'(j)}}_{j \in \ufs{q}}}\\
&= \scalebox{.85}{$\begin{cases}
\zero & \text{if $\psiinv_j (\psi'_{f'(j)})^{-1} s_{f'(j)} = \emptyset$ for some $j \in \ufs{q}$, and}\\
a_{\ang{\psiinv_{f(i)} (\psi'_{f'f(i)})^{-1} s_{f'f(i)}}_{i \in \ufs{p}}} & \text{if $\psi_j^{-1} (\psi'_{f'(j)})^{-1} s_{f'(j)} \neq \emptyset$ for each $j \in \ufs{q}$.}
\end{cases}$}
\end{split}
\end{equation}
\end{itemize}
Since 
\[\psi_{\fpinv c} = 1_{\ord{1}} \ifspace \fpinv c = \emptyset,\]
by \cref{psiprimec_one}, we can combine \cref{AF'AF_angs_i,AF'AF_angs_ii} as follows.
\begin{equation}\label{AF'AF_angs}
\begin{split}
&\left((\Aupomp)(\Aupom)(a,\glu)\right)_{\ang{s}}\\
&= \scalebox{.85}{$\begin{cases}
\zero & \text{if $(\psi'_c \psi_{\fpinv c})^{-1} s_c = \emptyset$ for some $c \in \ufs{t}$, and}\\
a_{\ang{(\psi'_{f'f(i)} \psi_{f(i)})^{-1} s_{f'f(i)}}_{i \in \ufs{p}}} & \text{if $(\psi'_c \psi_{\fpinv c})^{-1} s_c \neq \emptyset$ for each $c \in \ufs{t}$.}
\end{cases}$}
\end{split}
\end{equation}
By \cref{A_fangpsi_obj_component,A_FggprimeFgg_general}, the $\ang{s}$-component object $\big(\sys{\A}{(\upom'\upom)} (a,\glu)\big)_{\ang{s}}$ is also given by \cref{AF'AF_angs}.  Moreover, the argument in this paragraph so far also applies to morphisms of $\angordm$-systems using \cref{A_fangpsi_mor_component}, which proves that the diagram \cref{A_FggprimeFgg_general} commutes on morphisms.

\parhead{Gluing}.  To see that the $\angordl$-systems $(\Aupomp)(\Aupom)(a,\glu)$ and $\sys{\A}{(\upom'\upom)}(a,\glu)$ have the same gluing morphisms, suppose we are given $r \geq 0$, an object $x \in \Op(r)$, a marker $\ang{s} = \ang{s_c \subseteq \ufs{\ell}_c}_{c \in \ufs{t}}$, an index $k \in \ufs{t}$, and a partition
\[s_k = \coprod_{u \in \ufs{r}}\, s_{k,u} \subseteq \ufs{\ell}_k.\]
\begin{itemize}
\item Suppose 
\[\begin{split}
(\psi'_c \psi_{\fpinv c} )^{-1} s_c &\neq \emptyset \quad \text{for each $c \in \ufs{t}$ and}\\
(f'f)^{-1}(k) &\neq \emptyset.
\end{split}\]  
Then \cref{A_fangpsi_gluing}---applied to each of $\Aupom$, $\Aupomp$, and $\sys{\A}{(\upom'\upom)}$---implies that the corresponding gluing morphism of each of $(\Aupomp)(\Aupom)(a,\glu)$ and $\sys{\A}{(\upom'\upom)}(a,\glu)$ is given by the following gluing morphism of $(a,\glu)$.
\[\glu_{x;\, \ang{(\psi'_{f'f(i)} \psi_{f(i)})^{-1} s_{f'f(i)}}_{i \in \ufs{p}} \scs (f'f)^{-1}(k) \scs \ang{(\psi'_k \psi_{\fpinv k} )^{-1} s_{k,u}}_{u \in \ufs{r}}}\]
\item All other gluing morphisms of $(\Aupomp)(\Aupom)(a,\glu)$ and $\sys{\A}{(\upom'\upom)}(a,\glu)$ are identity morphisms.
\end{itemize}
This proves that the diagram \cref{A_FggprimeFgg_general} commutes on objects.  This completes the proof that $\Adash$ preserves composition of morphisms.
\end{proof}

\section{Multimorphism $G$-Functors in Arity 0}
\label{sec:jemg_zero}

For a $\Tinf$-operad $\Op$ \pcref{as:OpA}, \cref{sec:jemg_objects,sec:jemg_morphisms,sec:jemg_morphisms_ii} construct the object assignment of
\[\MultpsO \fto{\Jgo} \GGCat,\]
where $\MultpsO$ and $\GGCat$ are the $\Gcat$-multicategories in, respectively, \cref{thm:multpso,expl:ggcat_gcatenr}.   The object assignment of $\Jgo$ sends each $\Op$-pseudoalgebra $(\A,\gaA,\phiA)$ (\cref{def:pseudoalgebra}) to a $\Gskg$-category (\cref{expl:ggcat_obj}), that is, a pointed functor
\[(\Gsk,\vstar) \fto{\Jgo\A = \Adash} (\Gcatst,\boldone);\]
see \cref{A_ptfunctor}.  Continuing the construction of the $\Gcat$-multifunctor $\Jgo$, this section constructs the 0-ary multimorphism $G$-functor \cref{multifunctor_component}
\[\MultpsO(\ang{} \sscs \A) \fto{(\Jgo)_{(\ang{};\, \A)}} \GGCat(\ang{} \sscs \Adash)\]
for each $\Op$-pseudoalgebra $\A$.  The positive arity multimorphism $G$-functors of $\Jgo$ are constructed in subsequent sections.  To simplify the notation, we usually abbreviate $(\Jgo)_{(\ang{};\, \A)}$ to $\Jgo$.

\secoutline
\begin{itemize}
\item \cref{def:Jgo_zero} constructs the 0-ary multimorphism $G$-functor 
\[\MultpsO(\ang{} \sscs \A) \fto{\Jgo} \GGCat(\ang{} \sscs \Adash)\] 
and its strong variant, denoted $\Jgosg$, where $\Adash$ in the codomain is replaced by $\Asgdash$.
\item The second half of this section proves \cref{Jgoa_nbetad_welldef,Jgoa_nat,Jgotha_gmod}, which are used in \cref{def:Jgo_zero} to make sure that the 0-ary multimorphism $G$-functor $\Jgo$ is well defined.
\end{itemize}

\recollection

For \cref{def:Jgo_zero} below, we recall the following.
\begin{itemize}
\item By \cref{def:MultpsO_karycat} \pcref{multpso_zero}, 
\[\MultpsO(\ang{} \sscs \A)\] 
is the $G$-category $\A$.
\item The $G$-category 
\[\GGCat(\ang{} \sscs \Adash)\]
is described in \cref{expl:ggcat_zero}, applied to the $\Gskg$-category $\Adash$ in \cref{A_ptfunctor}.
\item The monoidal unit $\Gskg$-category $\gu$ is described in \cref{expl:ggcat_unit}.
\item The pointed $G$-category of $\angordn$-systems in $\A$, $\Aangordn$, is constructed in \cref{def:nsystem,def:nsystem_morphism,def:Aangordn_system,def:Aangordn_gcat}.
\item For each morphism $\upom$ in $\Gsk$, the pointed $G$-functor $\Aupom$ is defined in \cref{def:ftil_functor,def:psitil_functor,def:Afangpsi}.
\end{itemize}

\begin{definition}\label{def:Jgo_zero}
For a $\Tinf$-operad $\Op$ \pcref{as:OpA} and an $\Op$-pseudoalgebra $(\A,\gaA,\phiA)$ \pcref{def:pseudoalgebra}, the $G$-functor
\begin{equation}\label{Jgo_zero}
\MultpsO(\ang{} \sscs \A) = \A \fto{\Jgo} \GGCat(\ang{} \sscs \Adash)
\end{equation}
is defined as follows.  A strong variant $\Jgosg$ is defined in \cref{Jgosg_zero}.
\begin{description}
\item[Objects] For an object $a \in \A$, the component pointed functors of the object
\begin{equation}\label{Jgoa}
\Jgo a \in \GGCat(\ang{} \sscs \Adash),
\end{equation}
in the sense of \cref{ggcat_zero_obj_comp}, are defined as follows for objects $\angordn \in \Gsk$ \cref{Gsk_objects}.
\begin{description}
\item[Basepoint] For the basepoint $\vstar \in \Gsk$, the $\vstar$-component of $\Jgo a$ is the identity functor
\begin{equation}\label{Jgoa_vstar}
\gu\vstar = \ord{0} \fto{(\Jgo a)_\vstar = 1} \sys{\A}{\vstar} = \boldone.
\end{equation}
\item[Empty tuple] For the empty tuple $\ang{} \in \Gsk$, the pointed functor
\begin{equation}\label{Jgoa_empty}
\gu\ang{} = \ord{1} \fto{(\Jgo a)_{\ang{}}} \sys{\A}{\ang{}} = \A
\end{equation}
is defined by the following object assignment.
\begin{equation}\label{Jgoa_c}
(\Jgo a)_{\ang{}} (c) = \begin{cases}
\zero & \text{if $c = 0 \in \ord{1}$, and}\\
a & \text{if $c = 1 \in \ord{1}$.}
\end{cases}
\end{equation}
\item[Nonempty tuples] For an object $\angordn \in \Gsk \setminus \{\vstar, \ang{}\}$ of length $q>0$, the pointed functor
\begin{equation}\label{Jgoa_nbeta}
\gu\angordn = \sma_{j \in \ufs{q}}\, \ord{n}_j \fto{(\Jgo a)_{\angordn}} \Aangordn
\end{equation}
sends an object 
\begin{equation}\label{angdj}
\angd = \ang{d_j}_{j \in \ufs{q}} \in \sma_{j \in \ufs{q}}\, \ord{n}_j
\end{equation}
to the $\angordn$-system
\begin{equation}\label{Jgoa_nbetad}
(\Jgo a)_{\angordn} \angd \in \Aangordn.
\end{equation}
For each marker $\ang{s} = \ang{s_j \subseteq \ufs{n}_j}_{j \in \ufs{q}}$, its $\ang{s}$-component object \cref{a_angs} is defined as follows.
\begin{equation}\label{Jgoa_nbetads}
\big((\Jgo a)_{\angordn} \angd\big)_{\ang{s}} 
= \begin{cases}
\zero & \text{if $d_j \not\in s_j$ for some $j \in \ufs{q}$, and}\\
a & \text{if $d_j \in s_j$ for each $j \in \ufs{q}$.}
\end{cases}
\end{equation}
Each gluing morphism \cref{gluing-morphism} of $(\Jgo a)_{\angordn} \angd$ is defined as an identity morphism.  \cref{Jgoa_nbetad_welldef} shows that $(\Jgo a)_{\angordn} \angd$ is an $\angordn$-system.  Note that $(\Jgo a)_{\angordn} \angd$ is a \emph{strong} $\angordn$-system in $\A$, since each gluing morphism is the identity.
\item[Pointed functoriality]  $(\Jgo a)_{\angordn}$ is a pointed functor for the following reasons.
\begin{itemize}
\item Its domain, $\gu\angordn$, is a discrete category.
\item The basepoint $\ang{0} \in \gu\angordn$ is sent to the base $\angordn$-system $(\zero,1_\zero)$ (\cref{expl:nsystem_base}) by the first case in \cref{Jgoa_nbetads}.
\end{itemize}
\end{description}
\cref{Jgoa_nat} proves that $\Jgo a$ is an object of $\GGCat(\ang{} \sscs \Adash)$.
\item[Morphisms] For a morphism $\theta \cn a \to b$ in $\A$, the component pointed natural transformations of the morphism
\begin{equation}\label{Jgotha}
\Jgo a \fto{\Jgo\theta} \Jgo b,
\end{equation}
in the sense of \cref{ggcat_zero_mor_comp}, are defined as follows for objects $\angordn \in \Gsk$.
\begin{description}
\item[Basepoint] For the basepoint $\vstar \in \Gsk$, the $\vstar$-component of $\Jgo \theta$ is the identity natural transformation of the identity functor $1_{\boldone}$:
\begin{equation}\label{Jgotha_vstar}
(\Jgo a)_\vstar = 1_{\boldone} 
\fto{(\Jgo\theta)_{\vstar} = 1_{1_{\boldone}}} (\Jgo b)_\vstar = 1_{\boldone}.
\end{equation}
\item[Empty tuple] For the empty tuple $\ang{} \in \Gsk$, the pointed natural transformation
\begin{equation}\label{Jgotha_empty}
\begin{tikzpicture}[baseline={(a.base)}]
\draw[0cell]
(0,0) node (a) {\gu\ang{}}
(a)++(-.9,0) node (a') {\phantom{\ord{1}}}
(a')++(0,-.01) node (a'') {\ord{1}}
(a)++(3,0) node (b) {\sys{\A}{\ang{}}}
(b)++(.9,0) node (b') {\phantom{\A}}
(b')++(0,.03) node (b'') {\A}
;
\draw[1cell=.8]
(a') edge[equal] (a)
(b) edge[equal] (b')
(a) edge[bend left=20] node {(\Jgo a)_{\ang{}}} (b)
(a) edge[bend right=20] node[swap] {(\Jgo b)_{\ang{}}} (b)
;
\draw[2cell=.9]
node[between=a and b at .33, rotate=-90, 2label={above, (\Jgo\theta)_{\ang{}}}] {\Rightarrow}
;
\end{tikzpicture}
\end{equation}
is defined by the following component morphisms in $\A$.
\begin{equation}\label{Jgotha_empty_c}
(\Jgo\theta)_{\ang{}} (c) 
= \begin{cases}
\zero \fto{1_{\zero}} \zero & \text{if $c = 0 \in \ord{1}$, and}\\
a \fto{\theta} b & \text{if $c = 1 \in \ord{1}$.}
\end{cases}
\end{equation}
\item[Nonempty tuples] For an object $\angordn \in \Gsk \setminus \{\vstar, \ang{}\}$ of length $q>0$, the pointed natural transformation
\begin{equation}\label{Jgotha_nbeta}
\begin{tikzpicture}[baseline={(a.base)}]
\draw[0cell]
(0,0) node (a) {\gu\angordn}
(a)++(-1.4,0) node (a') {\sma\angordn}
(a)++(3.2,0) node (b) {\Aangordn}
;
\draw[1cell=.8]
(a') edge[equal] (a)
(a) edge[bend left=20] node {(\Jgo a)_{\angordn}} (b)
(a) edge[bend right=20] node[swap] {(\Jgo b)_{\angordn}} (b)
;
\draw[2cell=.9]
node[between=a and b at .33, rotate=-90, 2label={above, (\Jgo\theta)_{\angordn}}] {\Rightarrow}
;
\end{tikzpicture}
\end{equation}
is defined by the following $\ang{s}$-component morphism \cref{theta_angs} for each object $\angd \in \sma_{j\in \ufs{q}}\, \ord{n}_j$, as defined in \cref{angdj}, and each marker $\ang{s} = \ang{s_j \subseteq \ufs{n}_j}_{j \in \ufs{q}}$.
\begin{equation}\label{Jgotha_nbetads}
\begin{split}
&\big((\Jgo \theta)_{\angordn} \angd\big)_{\ang{s}} \\
&= \begin{cases}
\zero \fto{1_\zero} \zero & \text{if $d_j \not\in s_j$ for some $j \in \ufs{q}$, and}\\
a \fto{\theta} b & \text{if $d_j \in s_j$ for each $j \in \ufs{q}$.}
\end{cases}
\end{split}
\end{equation}
\item[Axioms]  The assignment 
\[(\Jgo a)_{\angordn} \angd \fto{(\Jgo \theta)_{\angordn} \angd} (\Jgo b)_{\angordn} \angd\]
satisfies the unity axiom \cref{nsystem_mor_unity} by the first case in  \cref{Jgotha_nbetads}.  The compatibility axiom \cref{nsystem_mor_compat} holds for $(\Jgo \theta)_{\angordn} \angd$ by 
\begin{itemize}
\item the fact that each gluing morphism of $(\Jgo a)_{\angordn} \angd$ and $(\Jgo b)_{\angordn} \angd$ is an identity morphism;
\item the naturality of the associativity constraint $\phiA$ \cref{phiA}; and
\item the action unity axiom \cref{pseudoalg_action_unity} and the basepoint axiom \cref{pseudoalg_basept_axiom}.
\end{itemize}
Thus, $(\Jgo \theta)_{\angordn} \angd$ is a morphism of $\angordn$-systems. 
\item[Pointed naturality] The assignment 
\begin{equation}\label{Jgotheta_n_d}
\angd \mapsto (\Jgo \theta)_{\angordn} \angd
\end{equation}
is a natural transformation because $\gu\angordn$ is a discrete category, with only identity morphisms.  The natural transformation $(\Jgo \theta)_{\angordn}$ is pointed because, for the basepoint $\ang{0} \in \sma_{j \in \ufs{q}}\, \ord{n}_j$, each component of $(\Jgo \theta)_{\angordn} \ang{0}$ is $1_\zero$ by the first case in \cref{Jgotha_nbetads}.
\end{description}
\cref{Jgotha_gmod} proves that 
\[\Jgo a \fto{\Jgo\theta} \Jgo b\] 
is a morphism in $\GGCat(\ang{}; \Adash)$.
\item[Functoriality]
The functoriality of $\Jgo$ in \cref{Jgo_zero} with respect to morphisms in $\A$ follows from \cref{Jgotha_vstar}, \cref{Jgotha_empty_c}, \cref{Jgotha_nbetads}, and the following facts.
\begin{itemize}
\item Identity morphisms and composition in $\GGCat(\ang{}; \Adash)$ are defined componentwise using the component pointed natural transformations \cref{ggcat_zero_mor_comp}.
\item Identity morphisms and composition in $\Aangordn$ are defined componentwise in $\A$ (\cref{def:nsystem_morphism}).
\end{itemize}
\item[$G$-equivariance] The functor $\Jgo$ in \cref{Jgo_zero} is $G$-equivariant by
\begin{itemize}
\item \cref{ga_scomponent,ga_gluing,gtheta_angs,Jgoa_c,Jgoa_nbetads,Jgotha_empty_c,Jgotha_nbetads};
\item the fact \cref{ggcat_zero_mor_g} that the $G$-action on $\GGCat(\ang{}; \Adash)$ is given componentwise by post-composing or post-whiskering with the $G$-action on $\Aangordn$; 
\item the fact that $\zero \in \A$ is $G$-fixed \cref{pseudoalg_zero}; and
\item the functoriality of the $g$-action on $\A$ for each $g \in G$.
\end{itemize}
\end{description}
This finishes the definition of the $G$-functor \cref{Jgo_zero} 
\[\MultpsO(\ang{} \sscs \A) = \A \fto{\Jgo} \GGCat(\ang{} \sscs \Adash).\]

\parhead{Strong variant}.  The $G$-functor
\begin{equation}\label{Jgosg_zero}
\MultpspsO(\ang{} \sscs \A) = \A \fto{\Jgosg} \GGCat(\ang{} \sscs \Asgdash)
\end{equation}
is defined by modifying the construction of $\Jgo$ given above as follows.
\begin{itemize}
\item The domain of $\Jgosg$ is the $G$-category $\MultpspsO(\ang{} \sscs \A)$, which is again given by the $G$-category $\A$ (\cref{def:MultpsO_karycat} \pcref{multpso_zero}).
\item In the codomain of $\Jgosg$, the $\Gskg$-category $\Asgdash$ is defined in \cref{sgAordnbe,AF_sg} and verified in \cref{A_ptfunctor}.
\end{itemize}
After these changes, the definition of $\Jgosg$ is the same as that of $\Jgo$.  
\begin{description}
\item[Objects] More precisely, for each object $a \in \A$, the object
\begin{equation}\label{Jgosga}
\Jgosg a \in \GGCat(\ang{} \sscs \Asgdash)
\end{equation}
has each $\angordn$-component pointed functor
\begin{equation}\label{Jgosga_nbeta}
\gu\angordn \fto{(\Jgosg a)_{\angordn}} \Aangordnsg
\end{equation}
defined in the same way as $(\Jgo a)_{\angordn}$ in \cref{Jgoa_vstar,Jgoa_empty,Jgoa_nbeta}.  This is well defined because, as we point out right after \cref{Jgoa_nbetads}, for each object $\angd \in \sma_{j \in \ufs{q}}\, \ord{n}_j$, the $\angordn$-system $(\Jgo a)_{\angordn} \angd$ is strong, with each gluing morphism given by an identity.
\item[Morphisms]
For a morphism $\theta \cn a \to b$ in $\A$, the morphism
\begin{equation}\label{Jgosgtha}
\Jgosg a \fto{\Jgosg\theta} \Jgosg b
\end{equation}
is defined in the same way as $\Jgo\theta$ in \cref{Jgotha}, where $\Adash$, $\Jgo a$, and $\Jgo b$ are replaced by, respectively, $\Asgdash$, $\Jgosg a$, and $\Jgosg b$.  This is well defined because $\Aangordnsg$ is the \emph{full} subcategory of $\Aangordn$, with strong $\angordn$-systems as objects.\defmark
\end{description}
\end{definition}

\subsection*{Proofs}

The rest of this section proves \cref{Jgoa_nbetad_welldef,Jgoa_nat,Jgotha_gmod}, which are used in \cref{def:Jgo_zero}.

\begin{lemma}\label{Jgoa_nbetad_welldef}
The data $(\Jgo a)_{\angordn} \angd$ defined in \cref{Jgoa_nbetad} form an $\angordn$-system in $\A$.  
\end{lemma}

\begin{proof}
In this proof, we abbreviate $(\Jgo a)_{\angordn} \angd$ to $(z,\glu)$.  

\parhead{Gluing}.  First, we observe that the gluing morphisms $\glu$, which are defined as identities, are well defined.  Given $r \geq 0$, an object $x \in \Op(r)$, a marker $\ang{s} = \ang{s_j \subseteq \ufs{n}_j}_{j \in \ufs{q}}$, an index $k \in \ufs{q}$, and a partition
\[s_k = \coprod_{\ell \in \ufs{r}}\, s_{k,\ell} \subseteq \ufs{n}_k,\]
the corresponding gluing morphism has the form
\begin{equation}\label{gluz}
\glu_{x;\, \ang{s},\, k, \ang{s_{k,\ell}}_{\ell \in \ufs{r}}} = 1_{z_{\ang{s}}}\cn 
\gaA_r\big(x \sscs \ang{z_{\ang{s} \,\compk\, s_{k,\ell}}}_{\ell \in \ufs{r}} \big) \to z_{\ang{s}}.
\end{equation}
We want to show that the domain and codomain of this gluing morphism are equal.  There are three cases to consider.

\parhead{Case 1}.  If $r = 0$, then $x = * \in \Op(0)$ and $s_k = \emptyset$.  The domain of the gluing morphism in \cref{gluz} is 
\[\gaA_0(*) = \zero \in \A.\]  
Since $d_k \not\in s_k = \emptyset$, by \cref{Jgoa_nbetads}, the codomain of the gluing morphism in \cref{gluz} is
\[z_{\ang{s}} = \zero.\]
Thus, the gluing morphism
\begin{equation}\label{gluz_i}
\glu_{*;\, \ang{s},\, k, \ang{}} = 1_\zero \cn \zero \to \zero
\end{equation}
is well defined in this case.

\parhead{Case 2}.  Suppose $r>0$ and $d_j \in s_j$ for each $j \in \ufs{q}$.  By \cref{Jgoa_nbetads}, the codomain of the gluing morphism in \cref{gluz} is
\[z_{\ang{s}} = a.\]
On the other hand, $d_k \in s_k$ implies that
\[\begin{cases}
d_k \in s_{k,\jm} & \text{for a unique index $\jm \in \ufs{r}$, and}\\
d_k \not\in s_{k,\ell} & \text{if $\ell \in \ufs{r} \setminus \{\jm\}$.}
\end{cases}\]
By \cref{Jgoa_nbetads}, this implies that
\[z_{\ang{s} \compk\, s_{k,\ell}} = \begin{cases}
a & \text{if $\ell = \jm$, and}\\
\zero & \text{if $\ell \in \ufs{r} \setminus \{\jm\}$.}
\end{cases}\]
Thus, the domain of the gluing morphism in \cref{gluz} is the object
\[\gaA_r\big(x \sscs \ang{\zero}_{\ell \in \ufs{r}} \,\comp_{\jm}\, a\big) = \gaA_1(\opu \sscs a) = a.\]
\begin{itemize}
\item The first equality follows from 
\begin{itemize}
\item $(r-1)$ applications of the basepoint axiom \cref{pseudoalg_basept_axiom} of $\A$ and
\item $\Op(1) = \{\opu\}$.
\end{itemize}
\item The second equality follows from the action unity axiom \cref{pseudoalg_action_unity} of $\A$.
\end{itemize}  
Thus, the gluing morphism
\begin{equation}\label{gluz_ii}
\glu_{x;\, \ang{s},\, k, \ang{s_{k,\ell}}_{\ell \in \ufs{r}}} = 1_a \cn a \to a
\end{equation}
is well defined in this case.

\parhead{Case 3}.  Suppose $r>0$ and $d_j \not\in s_j$ for some $j \in \ufs{q}$.  By \cref{Jgoa_nbetads},
\[z_{\ang{s}} = \zero = z_{\ang{s} \compk\, s_{k,\ell}}\]
for each $\ell \in \ufs{r}$.  Thus, the domain of the gluing morphism in \cref{gluz} is the object
\[\gaA_r\big(x \sscs \ang{\zero}_{\ell \in \ufs{r}}\big) = \gaA_0(*) = \zero.\]
The first equality follows from
\begin{itemize}
\item $r$ applications of the basepoint axiom \cref{pseudoalg_basept_axiom} of $\A$ and 
\item $\Op(0) = \{*\}$.
\end{itemize}  
Thus, the gluing morphism
\begin{equation}\label{gluz_iii}
\glu_{x;\, \ang{s},\, k, \ang{s_{k,\ell}}_{\ell \in \ufs{r}}} = 1_\zero \cn \zero \to \zero
\end{equation}
is well defined in this case.  We have shown that each gluing morphism of $(z,\glu)$ is well defined.

\parhead{Axioms}.  To see that the pair $(z,\glu)$ satisfies the axioms \cref{system_obj_unity,system_naturality,system_unity_i,system_unity_iii,system_equivariance,system_associativity,system_commutativity} for an $\angordn$-system in $\A$, we reuse the argument in \cref{expl:nsystem_base}, which proves that the base $\angordn$-system $(\zero,1_\zero)$ is well defined.  The key point is that each gluing morphism of $(z,\glu)$, as discussed in \cref{gluz_i,gluz_ii,gluz_iii}, is the identity morphism of either $\zero$ or $a$.
\end{proof}


\begin{lemma}\label{Jgoa_nat}
For each object $a \in \A$, the object specified in \cref{Jgoa} 
\[\Jgo a \in \GGCat(\ang{}; \Adash)\] 
is well defined.
\end{lemma}

\begin{proof}
By \cref{Jgoa_nbetad_welldef} and the discussion between \cref{Jgoa_nbetads,Jgotha}, each $(\Jgo a)_{\angordn}$ is a pointed functor.  By \cref{expl:ggcat_zero}, we need to show that, for each morphism $\upom \cn \angordm \to \angordn$ in $\Gsk$, the following diagram of pointed functors commutes.
\begin{equation}\label{Jgoa_natural}
\begin{tikzpicture}[vcenter]
\def\v{-1.5}
\draw[0cell]
(0,0) node (a1) {\gu\angordm}
(a1)++(3,0) node (a2) {\Aangordm}
(a1)++(0,\v) node (b1) {\gu\angordn}
(a2)++(0,\v) node (b2) {\Aangordn}
;
\draw[1cell=.9]
(a1) edge node[swap] {\gu\upom} (b1)
(a2) edge node {\Aupom} (b2)
(a1) edge node {(\Jgo a)_{\angordm}} (a2)
(b1) edge node {(\Jgo a)_{\angordn}} (b2)
;
\end{tikzpicture}
\end{equation}
There are four cases to consider, depending on what $\angordm$ and $\angordn$ are.

\parhead{Case 1: basepoint}.  Suppose either $\angordm$ or $\angordn$ is the basepoint $\vstar \in \Gsk$.  The diagram \cref{Jgoa_natural} commutes because
\begin{itemize}
\item $\gu\vstar = \sys{\A}{\vstar} = \boldone$ and
\item each of the four arrows is a pointed functor.
\end{itemize}  
Thus, we assume that $\angordm, \angordn \in \Gsk \setminus \{\vstar\}$.

\parhead{Case 2: empty codomain}. Suppose $\angordn = \ang{}$, the empty tuple.  Then $\angordm = \ang{}$ by \cref{Gsk_empty_mor}, and $\upom$ is either (i) the identity morphism, $1_{\ang{}}$, or (ii) the 0-morphism $\ang{} \to \vstar \to \ang{}$.
\begin{itemize}
\item If $\upom = 1_{\ang{}}$, then each composite in \cref{Jgoa_natural} is equal to $(\Jgo a)_{\angordm}$.
\item If $\upom$ is the 0-morphism, then each composite in \cref{Jgoa_natural} is the constant functor at the base $\angordn$-system $(\zero,1_\zero) \in \Aangordn$.
\end{itemize}
Thus, we assume that $\angordn \in \Gsk \setminus \{\vstar, \ang{}\}$ has length $q>0$.

\parhead{Case 3: empty domain}.  Suppose $\angordm = \ang{}$.  Then $\upom$ has the following form \cref{Fiqangpsi}.
\[\begin{tikzpicture}[vcenter]
\def\t{.7} \def\h{3}
\draw[0cell]
(0,0) node (a1) {\ang{}}
(a1)++(\h,0) node (a2) {\ordtu{1}_{j \in \ufs{q}}}
(a2)++(\h,0) node (a3) {\angordn}
;
\draw[1cell=.9]
(a1) edge node {(\im_q,\ang{1_{\ord{1}}}_{j \in \ufs{q}})} (a2)
(a2) edge node {(1_{\ufs{q}}, \ang{\psi_j}_{j \in \ufs{q}})} (a3)
;
\draw[1cell=1]
(a1) [rounded corners=2pt, shorten <=0ex] |- ($(a2)+(-1,\t)$)
-- node {\upom \,=\, (\im_q, \ang{\psi_j}_{j \in \ufs{q}})} ($(a2)+(1,\t)$) -| (a3)
;
\end{tikzpicture}\]
The diagram \cref{Jgoa_natural} takes the following form.
\begin{equation}\label{Jgoa_natural_emptynbeta}
\begin{tikzpicture}[vcenter]
\def\v{-1.5}
\draw[0cell]
(0,0) node (a1) {\gu\ang{}}
(a1)++(3,0) node (a2) {\sys{\A}{\ang{}}}
(a1)++(-.9,0) node (a1') {\ord{1}}
(a2)++(.9,0) node (a2') {\phantom{\A}}
(a2')++(0,.03) node (a2'') {\A}
(a1)++(0,\v) node (b1) {\gu\angordn}
(a2)++(0,\v) node (b2) {\Aangordn}
(b1)++(-1.5,0) node (b1') {\phantom{\sma_{j\in \ufs{q}}\, \ord{n}_j}}
(b1')++(0,-.06) node (b1'') {\sma_{j\in \ufs{q}}\, \ord{n}_j}
;
\draw[1cell=.9]
(a1) edge node[swap] {\gu\upom} (b1)
(a2) edge node {\Aupom} (b2)
(a1) edge node {(\Jgo a)_{\ang{}}} (a2)
(b1) edge node {(\Jgo a)_{\angordn}} (b2)
(a1') edge[equal] (a1)
(a2) edge[equal] (a2')
(b1') edge[equal] (b1)
;
\end{tikzpicture}
\end{equation}
Since $\ord{1} = \{0,1\}$, it suffices to show that the diagram \cref{Jgoa_natural_emptynbeta} commutes at the non-basepoint $1 \in \ord{1}$.
\begin{itemize}
\item Along the left-bottom composite of \cref{Jgoa_natural_emptynbeta}, the functor $\gu\upom$ sends $1 \in \ord{1}$ to the object
\[(\gu\upom)(1) = \ang{\psi_j(1)}_{j \in \ufs{q}} \in \gu\angordn.\]
By \cref{Jgoa_nbetads}, its image under the functor $(\Jgo a)_{\angordn}$ is the $\angordn$-system in $\A$ with, for each marker $\ang{s} = \ang{s_j \in \ufs{n}_j}_{j \in \ufs{q}}$, the following $\ang{s}$-component object.
\begin{equation}\label{Jgoa_nbeta_psione}
\begin{split}
&\big((\Jgo a)_{\angordn} \ang{\psi_j(1)}_{j \in \ufs{q}} \big)_{\ang{s}} \\
&= \begin{cases}
\zero & \text{if $\psi_j(1) \not\in s_j$ for some $j \in \ufs{q}$, and}\\
a & \text{if $\psi_j(1) \in s_j$ for each $j \in \ufs{q}$.}
\end{cases}
\end{split}
\end{equation}
\item By \cref{Jgoa_c}, the top-right composite in \cref{Jgoa_natural_emptynbeta} yields the $\ang{s}$-component object
\[\big((\Aupom) (\Jgo a)_{\ang{}} (1)\big)_{\ang{s}} = \big((\Aupom)(a)\big)_{\ang{s}}.\]
By \cref{AFa_angs}, the object $\big((\Aupom)(a)\big)_{\ang{s}}$ is also given by the right-hand side of \cref{Jgoa_nbeta_psione}.
\end{itemize}  
Thus, the $\angordn$-systems 
\[(\Jgo a)_{\angordn} (\gu\upom)(1) \andspace (\Aupom) (\Jgo a)_{\ang{}}(1)\]
have the same $\ang{s}$-component object for each marker $\ang{s}$.  Moreover, each of their gluing morphisms is the identity morphism, so they are equal as $\angordn$-systems.  This proves that the diagram \cref{Jgoa_natural_emptynbeta} commutes.

\parhead{Case 4: nonempty domain and codomain}.  Suppose $\angordm \in \Gsk \setminus \{\vstar,\ang{}\}$ has length $p>0$.  The morphism $\upom$ has the form
\[\angordm = \ang{\ordm_i}_{i \in \ufs{p}} \fto{\upom = (f, \ang{\psi})} 
\angordn = \ang{\ordn_j}_{j \in \ufs{q}}\]
as stated in \cref{as:fpsi}.  To prove the commutativity of the diagram \cref{Jgoa_natural}, we first observe the following.
\begin{itemize}
\item The pointed function
\[\sma_{i \in \ufs{p}}\, \ord{m}_i \fto{\gu\upom = \sma(f,\angpsi)} \sma_{j \in \ufs{q}}\, \ord{n}_j\]
is defined in \cref{smash_fpsi} and described explicitly in \cref{expl:smash_fpsi}. 
\item The pointed functor 
\[\Aangordm \fto{\Aupom = \psitil\ftil} \Aangordn\] 
is defined in \cref{A_fangpsi} and described explicitly in \cref{expl:A_fangpsi}. 
\item The pointed functors 
\[\begin{split}
\gu\angordm & \fto{(\Jgo a)_{\angordm}} \Aangordm \andspace\\ 
\gu\angordn & \fto{(\Jgo a)_{\angordn}} \Aangordn
\end{split}\]
are defined in \cref{Jgoa_nbeta}. 
\item For each object
\[\angc = \ang{c_i}_{i \in \ufs{p}} \in \sma_{i \in \ufs{p}}\, \ord{m}_i = \gu\angordm,\]
each of the two composites in \cref{Jgoa_natural} yields an $\angordn$-system in $\A$ with identity gluing morphisms.
\end{itemize}
Thus, it suffices to show that the two composites in \cref{Jgoa_natural} yield $\angordn$-systems with the same $\ang{s}$-component object for each object $\angc \in \gu\angordm$ and each marker $\ang{s} = \ang{s_j \subseteq \ufs{n}_j}_{j \in \ufs{q}}$.

At an object $\angc \in \gu\angordm$, the left-bottom composite in \cref{Jgoa_natural} yields the following $\ang{s}$-component object, where $c_\emptyset = 1$ for each index $j \in \ufs{q}$ with $\finv(j) = \emptyset$.
\begin{equation}\label{Jgoa_natc_i}
\begin{split}
&\big((\Jgo a)_{\angordn} (\gu\upom) \angc\big)_{\angs} \\
&= \big((\Jgo a)_{\angordn} \ang{\psi_j c_{\finv(j)}}_{j \in \ufs{q}}\big)_{\angs}\\
&= \begin{cases}
\zero & \text{if $\psi_j c_{\finv(j)} \not\in s_j$ for some $j \in \ufs{q}$, and}\\
a & \text{if $\psi_j c_{\finv(j)} \in s_j$ for each $j \in \ufs{q}$.}
\end{cases}
\end{split}
\end{equation}
The top-right composite in \cref{Jgoa_natural} yields the following $\ang{s}$-component object.
\begin{equation}\label{Jgoa_natc_ii}
\begin{split}
&\big((\Aupom) (\Jgo a)_{\angordm} \angc\big)_{\angs}\\
&= \scalebox{.9}{$\begin{cases}
\zero & \text{if $\psiinv_j s_j = \emptyset$ for some $j \in \ufs{q}$, and}\\
\big((\Jgo a)_{\angordm} \angc\big)_{\ang{\psiinv_{f(i)} s_{f(i)}}_{i \in \ufs{p}}} & \text{if $\psiinv_j s_j \neq \emptyset$ for each $j \in \ufs{q}$.}
\end{cases}$}
\end{split}
\end{equation}
The object in the second case in \cref{Jgoa_natc_ii} is given by
\[\begin{split}
&\big((\Jgo a)_{\angordm} \angc\big)_{\ang{\psiinv_{f(i)} s_{f(i)}}_{i \in \ufs{p}}}\\
&= \scalebox{.9}{$\begin{cases}
\zero & \text{if $c_i \not\in \psiinv_{f(i)} s_{f(i)}$ (i.e., $\psi_{f(i)} c_i \not\in s_{f(i)}$) for some $i \in \ufs{p}$, and}\\
a & \text{if $c_i \in \psiinv_{f(i)} s_{f(i)}$ (i.e., $\psi_{f(i)} c_i \in s_{f(i)}$) for each $i \in \ufs{p}$.}
\end{cases}$}
\end{split}\]
Thus, \cref{Jgoa_natc_i,Jgoa_natc_ii} describe the same $\angs$-component object.  This proves that the diagram \cref{Jgoa_natural} commutes in this case.
\end{proof}

The following lemma is used in \cref{def:Jgo_zero} to make sure that the morphism assignment of $\Jgo$ in \cref{Jgo_zero} is well defined.

\begin{lemma}\label{Jgotha_gmod} 
For each morphism $\theta \cn a \to b$ in $\A$, the data
\[\Jgo a \fto{\Jgo\theta} \Jgo b\] 
specified in \cref{Jgotha} form a morphism in $\GGCat(\ang{};\Adash)$.
\end{lemma}

\begin{proof}
It is explained in \cref{Jgotheta_n_d} that each component $(\Jgo \theta)_{\angordn}$ is a pointed natural transformation.  By \cref{ggcat_zero_mor_ax}, we need to check that, for each morphism $\upom \cn \angordm \to \angordn$ in $\Gsk$, the following two whiskered pointed natural transformations are equal.
\begin{equation}\label{Jgotha_modaxiom}
\begin{tikzpicture}[vcenter]
\def\v{-2} \def\s{.4} \def\t{.6} \def\d{20}
\draw[0cell]
(0,0) node (a1) {\gu\angordm}
(a1)++(3.5,0) node (a2) {\Aangordm}
(a1)++(0,\v) node (b1) {\gu\angordn}
(a2)++(0,\v) node (b2) {\Aangordn}
;
\draw[1cell=.8]
(a1) edge node[swap] {\gu\upom} (b1)
(a2) edge node {\Aupom} (b2)
(a1) edge[bend left=\d] node[pos=\s] {(\Jgo a)_{\angordm}} (a2)
(a1) edge[bend right=\d] node[swap,pos=\t] {(\Jgo b)_{\angordm}} (a2)
(b1) edge[bend left=\d] node[pos=\s] {(\Jgo a)_{\angordn}} (b2)
(b1) edge[bend right=\d] node[swap,pos=\t] {(\Jgo b)_{\angordn}} (b2)
;
\draw[2cell=.9]
node[between=a1 and a2 at .35, rotate=-90, 2label={above,(\Jgo\theta)_{\angordm}}] {\Rightarrow}
node[between=b1 and b2 at .35, rotate=-90, 2label={above,(\Jgo\theta)_{\angordn}}] {\Rightarrow}
;
\end{tikzpicture}
\end{equation}
To prove that the two whiskered natural transformations in \cref{Jgotha_modaxiom} are equal, we reuse the proof of the commutativity of the diagram \cref{Jgoa_natural} in \cref{Jgoa_nat}.  Instead of the expressions \cref{A_fangpsi_obj_component,AFa_angs,Jgoa_c,Jgoa_nbetads} for objects, here we use, respectively, \cref{A_fangpsi_mor_component,AF_theta_angs,Jgotha_empty_c,Jgotha_nbetads} for morphisms.
\end{proof}

\section{Multimorphism $G$-Functors in Positive Arity: Objects}
\label{sec:jemg_pos_i}

For a $\Tinf$-operad $\Op$ \pcref{as:OpA}, so far in this chapter we have constructed the following parts of the desired $\Gcat$-multifunctor 
\[\MultpsO \fto{\Jgo} \GGCat.\]
\begin{itemize}
\item \cref{sec:jemg_objects,sec:jemg_morphisms,sec:jemg_morphisms_ii} construct the object assignment of $\Jgo$, which sends each $\Op$-pseudoalgebra $\A$ to a $\Gskg$-category (\cref{expl:ggcat_obj})
\[(\Gsk,\vstar) \fto{\Jgo\A = \Adash} (\Gcatst,\boldone).\]
\item \cref{sec:jemg_zero} constructs the 0-ary multimorphism $G$-functors
\[\MultpsO(\ang{} \sscs \A) = \A \fto{\Jgo} \GGCat(\ang{} \sscs \Adash).\]
\end{itemize}

\begin{assumption}\label{as:BAi}
We assume that 
\[(\B,\gaB,\phiB) \andspace \bang{(\A_i,\gaAi,\phiAi)}_{i \in \ufs{k}}\]
are $\Op$-pseudoalgebras (\cref{def:pseudoalgebra}) for some $k \geq 1$.  
\end{assumption}

Continuing the construction of the $\Gcat$-multifunctor $\Jgo$, this section constructs the $G$-equivariant object assignment of the $k$-ary multimorphism $G$-functor \cref{multifunctor_component}
\[\MultpsO(\ang{\A_i}_{i \in \ufs{k}} \sscs \B ) \fto{(\Jgo)_{(\ang{\A_i}_{i \in \ufs{k}} ;\, \B)}} 
\GGCat(\ang{\Aidash}_{i \in \ufs{k}} \sscs \Bdash).\]
See \cref{def:Jgo_pos_obj}.  To simplify the notation, we abbreviate this $G$-functor to $\Jgo$.  There is also a strong variant of $\Jgo$, denoted $\Jgosg$, that involves $\MultpspsO$ in the domain---with $k$-ary $\Op$-pseudomorphisms instead of $k$-lax $\Op$-morphisms---and strong systems in the codomain.  To improve readability, several statements used in \cref{def:Jgo_pos_obj} are proved in \cref{sec:jemg_pos_i_proof}.  The morphism assignment of $\Jgo$ is constructed in \cref{sec:jemg_pos_ii}.

\recollection
For \cref{def:Jgo_pos_obj} below, we briefly recall the following.
\begin{itemize}
\item By \cref{def:multicatO}, an object in the $G$-category 
\[\MultpsO(\ang{\A_i}_{i \in \ufs{k}} \sscs \B)\] 
is a $k$-lax $\Op$-morphism (\cref{def:k_laxmorphism})
\[\bang{(\A_i,\gaAi,\phiAi)}_{i \in \ufs{k}} \fto{(f, \laxf)} (\B,\gaB,\phiB).\]
\item For $\Gskg$-categories $\ang{f_i}_{i \in \ufs{k}}$ and $f$, the $G$-category 
\[\GGCat(\ang{f_i}_{i \in \ufs{k}} \sscs f)\] 
is described in \cref{expl:ggcat_positive}, in particular \cref{ggcat_k1,ggcat_k1_comp,Gk_mor,ggcat_k1_nat}.  \cref{def:Jgo_pos_obj} below involves the $\Gskg$-categories \pcref{A_ptfunctor}
\[f_i = \Aidash \andspace f = \Bdash \cn (\Gsk,\vstar) \to (\Gcatst,\boldone).\]
\end{itemize}

\begin{definition}\label{def:Jgo_pos_obj}
Under \cref{as:OpA,as:BAi}, the $G$-equivariant object assignment
\begin{equation}\label{Jgo_pos_obj_assignment}
\MultpsO(\ang{\A_i}_{i \in \ufs{k}} \sscs \B) \fto{\Jgo} 
\GGCat(\ang{\Aidash}_{i \in \ufs{k}} \sscs \Bdash)
\end{equation}
is defined as follows.  A strong variant $\Jgosg$ is defined in \cref{Jgosg_pos_obj_assignment}.  Given a $k$-lax $\Op$-morphism 
\[\bang{(\A_i,\gaAi,\phiAi)}_{i \in \ufs{k}} \fto{(f, \laxf)} (\B,\gaB,\phiB),\]
to define the object
\begin{equation}\label{Jgo_f}
\ang{\Aidash}_{i \in \ufs{k}} \fto{\Jgo f} \Bdash
\end{equation}
in the sense of \cref{ggcat_k1}, we consider the following objects for $i \in \ufs{k}$ and $j \in \ufs{r}_i$ for each $i$.
\begin{equation}\label{angordmdot_bm}
\left\{\begin{aligned}
\ordmij & \in \Fsk & \angordmi &= \angordmij_{j \in \ufs{r}_i} \in \Gsk\\
\angordmdot &= \ang{\angordmi}_{i \in \ufs{k}} \in \Gsk^k \phantom{M} &
\bm &= \txoplus_{i \in \ufs{k}}\, \angordmi \in \Gsk
\end{aligned}\right.
\end{equation}
The $\angordmdot$-component pointed functor
\begin{equation}\label{Jgof_component}
\txsma_{i \in \ufs{k}}\, \Aiangordmi \fto{(\Jgo f)_{\angordmdot}} \Bboldm,
\end{equation}
in the sense of \cref{ggcat_k1_comp}, is defined as follows.  
\begin{description}
\item[Base case]  If $\angordmi = \vstar$ for some $i \in \ufs{k}$, then $\bm = \vstar$ by \cref{Gsk_oplus_vstar}.  Thus, by the first case in \cref{A_angordn}, $(\Jgo f)_{\angordmdot}$ is the unique isomorphism between terminal $G$-categories.  In the rest of this definition, we assume $\angordmi \in \Gsk \setminus \{\vstar\}$ for each $i \in \ufs{k}$, which means $\ordmij \neq \ord{0}$ for any $j \in \ufs{r}_i$. 
\item[Component objects] For the object assignment of $(\Jgo f)_{\angordmdot}$, we consider an $\angordmi$-system in $\A_i$ (\cref{def:nsystem})
\begin{equation}\label{ai_gluai}
(a_i,\glu^{a_i}) \in \Aiangordmi \foreachspace i \in \ufs{k}
\end{equation}
and the object
\begin{equation}\label{aglua_sma}
(a,\glu^a) = \ang{(a_i, \glu^{a_i})}_{i \in \ufs{k}} \in \txsma_{i \in \ufs{k}}\, \Aiangordmi.
\end{equation}

\parhead{Base systems}.  Since we want $(\Jgo f)_{\angordmdot}$ to be a pointed functor, if any $(a_i,\glu^{a_i})$ is the base $\angordmi$-system $(\zero, 1_\zero)$ in $\A_i$ (\cref{expl:nsystem_base}), then we must define
\begin{equation}\label{Jgof_m_basesystem}
\big((\Jgo f)_{\angordmdot} (a,\glu^a), \glu\big) = (\zero,1_\zero) \in \Bboldm,
\end{equation}
the base $\bm$-system in $\B$.  

\parhead{Non-base systems}.  Suppose $(a_i,\glu^{a_i})$ is not the base $\angordmi$-system for any $i \in \ufs{k}$.  For markers
\begin{equation}\label{marker_sij}
\angs = \bang{\bang{s_{i,j} \subseteq \ufsmij}_{j \in \ufs{r}_i}}_{i \in \ufs{k}} \andspace
s_{i\crdot} = \ang{s_{i,j}}_{j \in \ufs{r}_i},
\end{equation}
using the underlying functor \cref{underlying_functor_f} of $(f,\laxf)$, the $\bm$-system in $\B$
\begin{equation}\label{Jgof_m_objects}
\big((\Jgo f)_{\angordmdot} (a,\glu^a), \glu\big) \in \Bboldm
\end{equation}
has $\angs$-component object defined as
\begin{equation}\label{Jgof_m_obj_comp}
\big((\Jgo f)_{\angordmdot} (a,\glu^a)\big)_{\angs} 
= f\ang{a_{i,\, s_{i\crdot}}}_{i \in \ufs{k}} \in \B,
\end{equation}
where
\[a_{i,\, s_{i\crdot}} = (a_i)_{s_{i\crdot}} \in \A_i\]
denotes the $s_{i\crdot}$-component object \cref{a_angs} of $(a_i,\glu^{a_i})$.   If $r_i = 0$ for some index $i \in \ufs{k}$, then
\[\angordmi = \ang{} \andspace a_i \in \sys{\A_i}{\ang{}} = \A_i.\]
For such an index $i$, the object $a_{i,\, s_{i\crdot}}$ in \cref{Jgof_m_obj_comp} is interpreted as $a_i$ .

\item[Gluing] To define the gluing morphisms $\glu$ of the $\bm$-system in \cref{Jgof_m_objects}, suppose we are given an object $x \in \Op(t)$ for some $t \geq 0$, a marker $\angs$ as defined in \cref{marker_sij}, a pair of indices $(p,q) \in \ufs{k} \times \ufs{r}_p$, and a partition
\[s_{p,q} = \coprod_{\ell \in \ufs{t}}\, s_{p,q,\ell} \subseteq \ufs{m}_{p,q}.\]
To simplify the presentation below, we first use the $\comp_?$ notation in \cref{compk} to define the pointed functor
\begin{equation}\label{ftil_ApB}
\fgr = f\big(\ang{a_{i,\, s_{i\crdot}}}_{i \in \ufs{k}}\, \compp - \big) \cn (\A_p,\zero) \to (\B,\zero).
\end{equation}
In terms of $\fgr$, the $\angs$-component object in \cref{Jgof_m_obj_comp} is 
\[\big((\Jgo f)_{\angordmdot} (a,\glu^a)\big)_{\angs} = \fgr a_{p,\, s_{p\crdot}}.\]
Denoting by $u = \big(\!\sum_{i=1}^{p-1}\, r_i\big) + q$, we note that the markers
\[\angs \compu s_{p,q,\ell} \andspace s_{p\crdot} \compq s_{p,q,\ell}\]
are obtained from, respectively, $\angs$ and $s_{p\crdot}$ by replacing $s_{p,q}$ with $s_{p,q,\ell}$.  The corresponding gluing morphism $\glu$ in \cref{Jgof_m_objects} is defined as the following composite in $\B$.
\begin{equation}\label{Jgof_m_gluing}
\begin{tikzpicture}[vcenter]
\def\u{-1} \def\h{4.5} \def\a{10} \def\b{.7}
\draw[0cell=.85]
(0,0) node (a1) {\gaB_t\big(x \sscs \ang{((\Jgo f)_{\angordmdot} (a,\glu^a))_{\angs \,\compu\, s_{p,q,\ell}}}_{\ell \in \ufs{t}}\big)}
(a1)++(\h,0) node (a2) {\big((\Jgo f)_{\angordmdot} (a,\glu^a)\big)_{\angs}}
(a1)++(0,\u) node (b1) {\gaB_t\big(x \sscs \ang{\fgr a_{p,\, (s_{p\crdot} \,\compq\, s_{p,q,\ell})}}_{\ell \in \ufs{t}}\big)}
(a2)++(0,\u) node (b2) {\fgr a_{p,\, s_{p\crdot}}}
(b1)++(\h/2,-1.2) node (c) {\fgr \gaAp_t \big(x \sscs \ang{a_{p,\, (s_{p\crdot} \,\compq\, s_{p,q,\ell})}}_{\ell \in \ufs{t}} \big)}
;
\draw[1cell=.9]
(a1) edge[equal,shorten <=-.5ex,shorten >=-.5ex] (b1)
(a2) edge[equal] (b2)
;
\draw[1cell=.9]
(a1) [rounded corners=2pt, shorten <=-.2ex] |- ($(a1)+(1,\b)$)
-- node {\glu_{x;\, \angs,\, u, \ang{s_{p,q,\ell}}_{\ell \in \ufs{t}}}} ($(a2)+(-1,\b)$) -| (a2)
;
\draw[1cell=.9]
(b1) [rounded corners=2pt, shorten <=-.2ex] |- node[pos=.2,swap] {\laxf_{t,p;\, \bolda}} (c);
\draw[1cell=.9]
(c) [rounded corners=2pt, shorten <=-0ex] -| node[pos=.8] {\fgr \glu^{a_p}_{x;\, s_{p\crdot} \scs q, \ang{s_{p,q,\ell}}_{\ell \in \ufs{t}}}} (b2);
\end{tikzpicture}
\end{equation}
The lower-left and lower-right morphisms in \cref{Jgof_m_gluing} are given as follows.
\begin{itemize}
\item $\laxf_{t,p;\, \bolda}$ is the component of the $(t,p)$-action constraint \cref{klax_constraint} of the $k$-lax $\Op$-morphism $(f,\laxf)$ at the $(k+t)$-tuple of objects
\begin{equation}\label{laxf_bolda}
\bolda = \big(\ang{a_{i,\, s_{i\crdot}}}_{i=1}^{p-1}, x, \ang{a_{p,\, (s_{p\crdot} \,\compq\, s_{p,q,\ell})}}_{\ell \in \ufs{t}} \spc \ang{a_{i,\, s_{i\crdot}}}_{i=p+1}^k \big).
\end{equation}
When there is no danger of confusion, we abbreviate $\laxf_{t,p;\, \bolda}$ to $\laxf_{t,p}$.
\item The morphism in $\A_p$
\[\gaAp_t \big(x \sscs \ang{a_{p,\, (s_{p\crdot} \,\compq\, s_{p,q,\ell})}}_{\ell \in \ufs{t}} \big) 
\fto{\glu^{a_p}_{x;\, s_{p\crdot} \scs q, \ang{s_{p,q,\ell}}_{\ell \in \ufs{t}}}} a_{p,\, s_{p\crdot}}\]
is a component of the gluing morphism \cref{gluing-morphism} of the $\angordmpe$-system $(a_p, \glu^{a_p})$ in \cref{ai_gluai}.
\end{itemize}
\cref{Jgof_obj_welldef} proves that the data $\big((\Jgo f)_{\angordmdot} (a,\glu^a), \glu\big)$---which consist of the $\angs$-component objects in \cref{Jgof_m_obj_comp} and the gluing morphisms in \cref{Jgof_m_gluing}---form an $\bm$-system in $\B$.  Thus, the object assignment
\[(a,\glu^a) \mapsto \big((\Jgo f)_{\angordmdot} (a,\glu^a), \glu\big)\] 
is well defined.
\item[Morphisms] For the morphism assignment of $(\Jgo f)_{\angordmdot}$, we consider a morphism of $\angordmi$-systems in $\A_i$ (\cref{def:nsystem_morphism})
\begin{equation}\label{thetai_aaprime}
(a_i, \glu^{a_i}) \fto{\theta_i} (a_i', \glu^{a_i'}) \foreachspace i \in \ufs{k}
\end{equation}
and the morphism
\begin{equation}\label{theta_angthetai}
(a,\glu^a) = \ang{(a_i, \glu^{a_i})}_{i \in \ufs{k}} \fto{\theta = \ang{\theta_i}_{i \in \ufs{k}}} 
(a',\glu^{a'}) = \ang{(a_i', \glu^{a_i'})}_{i \in \ufs{k}}
\end{equation}
in $\txsma_{i \in \ufs{k}}\, \Aiangordmi$.

\parhead{Base systems}.  If, for some index $i \in \ufs{k}$,
\[(\zero,1_\zero) \fto{\theta_i = 1} (\zero, 1_\zero),\] 
the identity morphism of the base $\angordmi$-system in $\A_i$, then we define
\begin{equation}\label{Jgof_m_base}
(\zero,1_\zero) \fto{(\Jgo f)_{\angordmdot} \theta = 1} (\zero,1_\zero),
\end{equation}
the identity morphism of the base $\bm$-system in $\B$.

\parhead{Non-base systems}.  Suppose $\theta_i \neq 1_{(\zero,1_\zero)}$ for any $i \in \ufs{k}$.  For each marker \cref{marker_sij}
\[\angs = \ang{\ang{s_{i,j} \subseteq \ufsmij}_{j \in \ufs{r}_i}}_{i \in \ufs{k}},\] 
the $\angs$-component of the morphism of $\bm$-systems in $\B$
\begin{equation}\label{Jgof_m_theta}
(\Jgo f)_{\angordmdot} (a,\glu^a) \fto{(\Jgo f)_{\angordmdot} \theta} 
(\Jgo f)_{\angordmdot} (a',\glu^{a'})
\end{equation}
is defined as
\begin{equation}\label{Jgof_m_theta_comp}
f\ang{a_{i,\, s_{i\crdot}}}_{i \in \ufs{k}} 
\fto{((\Jgo f)_{\angordmdot} \theta)_{\angs} = f\ang{\theta_{i,\, s_{i\crdot}}}_{i \in \ufs{k}}} 
f\ang{a'_{i,\, s_{i\crdot}}}_{i \in \ufs{k}},
\end{equation}
where
\[a_{i,\, s_{i\crdot}} \fto{\theta_{i,\, s_{i\crdot}}} a'_{i,\, s_{i\crdot}}\]
is the $s_{i\crdot}$-component morphism \cref{theta_angs} of $\theta_i$.  For each index $i \in \ufs{k}$ such that $r_i = 0$, we interpret $\theta_{i,\, s_{i\crdot}}$ as the morphism
\[a_i \fto{\theta_i} a_i' \inspace \sys{\A_i}{\ang{}} = \A_i.\]
Note that, by the basepoint conditions \cref{laxf_basept} for $(f,\laxf)$, the case \cref{Jgof_m_base} may be regarded as a sub-case of \cref{Jgof_m_theta}.
\item[Unity and compatibility]  $(\Jgo f)_{\angordmdot} \theta$ in \cref{Jgof_m_theta} is a morphism of $\bm$-systems in $\B$ for the following reasons.
\begin{itemize}
\item The unity axiom \cref{nsystem_mor_unity} holds for $(\Jgo f)_{\angordmdot} \theta$ by
\begin{itemize}
\item the unity axiom for each $\theta_i$ and 
\item the basepoint condition \cref{laxf_basept} for the $k$-lax $\Op$-morphism $f$.
\end{itemize}
\item The compatibility axiom \cref{nsystem_mor_compat} holds for $(\Jgo f)_{\angordmdot} \theta$ by
\begin{itemize}
\item the definition \cref{Jgof_m_gluing} of the gluing morphism $\glu$, 
\item the naturality of $\laxf_{t,p}$,
\item the functoriality of $f \cn \txprod_{i \in \ufs{k}}\, \A_i \to \B$, and
\item the compatibility axiom for $\theta_p$.
\end{itemize}
\end{itemize}
Thus, the morphism assignment of $(\Jgo f)_{\angordmdot}$ is well defined.
\item[Pointed functoriality]
The assignment 
\begin{equation}\label{Jgof_m_ptfunctor}
\theta \mapsto (\Jgo f)_{\angordmdot} \theta
\end{equation} 
given by \cref{Jgof_m_base,Jgof_m_theta}
preserves identity morphisms and composition by 
\begin{itemize}
\item the definition \cref{Jgof_m_theta_comp} of the $\angs$-component morphism and 
\item the functoriality of $f$.
\end{itemize}
Thus,
\begin{itemize}
\item the object assignment in \cref{Jgof_m_basesystem,Jgof_m_objects}, and
\item the morphism assignment in \cref{Jgof_m_base,Jgof_m_theta}
\end{itemize}
define a pointed functor $(\Jgo f)_{\angordmdot}$, as stated in \cref{Jgof_component}.
\end{description}

\cref{Jgof_natural} proves that, for each $k$-lax $\Op$-morphism $(f,\laxf)$, $\Jgo f$ \cref{Jgo_f} is a well-defined object.  This finishes the construction of the object assignment
\[\MultpsO(\ang{\A_i}_{i \in \ufs{k}} \sscs \B ) \fto{\Jgo} 
\GGCat(\ang{\Aidash}_{i \in \ufs{k}} \sscs \Bdash),\]
sending $(f,\laxf)$ to $\Jgo f$, as stated in \cref{Jgo_pos_obj_assignment}.  \cref{Jgo_k_Geq} proves that this object assignment is $G$-equivariant.

\parhead{Strong variant}.  The $G$-equivariant object assignment
\begin{equation}\label{Jgosg_pos_obj_assignment}
\MultpspsO(\ang{\A_i}_{i \in \ufs{k}} \sscs \B) \fto{\Jgosg} 
\GGCat(\ang{\Aisgdash}_{i \in \ufs{k}} \sscs \Bsgdash)
\end{equation}
is defined by modifying the construction of $\Jgo$ given in \cref{Jgo_pos_obj_assignment} as follows.
\begin{itemize}
\item The domain of $\Jgosg$ is the $G$-category in \cref{def:MultpsO_karycat} \eqref{multpspso_k}, whose objects are $k$-ary $\Op$-pseudomorphisms.  Recall from \cref{def:k_laxmorphism} that a $k$-ary $\Op$-pseudomorphism is a $k$-lax $\Op$-morphism $(f,\laxf)$ with invertible action constraint $\laxf$. 
\item In the codomain, the $\Gskg$-categories $\Aisgdash$ and $\Bsgdash$, which involve strong systems \pcref{def:nsystem}, are given by \cref{A_ptfunctor}.
\item For each $k$-ary $\Op$-pseudomorphism 
\[\ang{(\A_i,\gaAi,\phiAi)}_{i \in \ufs{k}} \fto{(f,\laxf)} (\B,\gaB,\phiB),\]
the $\angordmdot$-component pointed functor of $\Jgosg f$,
\begin{equation}\label{Jgosgf_component}
\txsma_{i \in \ufs{k}}\, \Aisgangordmi \fto{(\Jgosg f)_{\angordmdot}} \Bsgboldm,
\end{equation}
is the restriction of $(\Jgo f)_{\angordmdot}$, as defined in \cref{Jgof_component}, to the full subcategory $\Aisgangordmi \subseteq \Aiangordmi$ for $i \in \ufs{k}$.
\end{itemize}
To see that $\Jgosg$ is well defined, suppose that $(a_i,\glu^{a_i})$ in \cref{ai_gluai} is a strong $\angordmi$-system in $\A_i$ for each $i \in \ufs{k}$.  Then each gluing morphism defined in \cref{Jgof_m_gluing} 
\[\glu_{x;\, \angs,\, u,\, \ang{s_{p,q,\ell}}_{\ell \in \ufs{t}}} 
= \fgr \glu^{a_p}_{x;\, s_{p\crdot} \scs q,\, \ang{s_{p,q,\ell}}_{\ell \in \ufs{t}}} \circ \laxf_{t,p}\] 
is an isomorphism, since $\laxf_{t,p}$ and $\glu^{a_p}_{\cdots}$ are now isomorphisms.  Thus, the pair
\[\big((\Jgosg f)_{\angordmdot} (a,\glu^a), \glu\big)\]
is an object in $\Bsgboldm$, that is, a strong $\bm$-system in $\B$.
\end{definition}

\section{Proofs}
\label{sec:jemg_pos_i_proof}

This section proves \cref{Jgof_obj_welldef,Jgof_natural,Jgo_k_Geq}, which are used in \cref{def:Jgo_pos_obj}.  The setting is the same as \cref{sec:jemg_pos_i}, so \cref{as:OpA,as:BAi} are in effect throughout this section.

\secoutline
\begin{itemize}
\item \cref{Jgof_obj_welldef} proves that each $k$-lax $\Op$-morphism for $\Op$-pseudoalgebras 
\[\ang{\A_i}_{i \in \ufs{k}} \fto{(f,\laxf)} \B\] 
yields a well-defined object assignment
\[\txsma_{i \in \ufs{k}}\, \Aiangordmi \ni (a,\glu^a) 
\mapsto \big((\Jgo f)_{\angordmdot} (a,\glu^a), \glu\big) \in \Bboldm.\] 
\item \cref{Jgof_natural} proves that $\Jgo f$, as defined in \cref{Jgo_f}, is an object in the $k$-ary multimorphism $G$-category 
\[\GGCat(\ang{\Aidash}_{i \in \ufs{k}}; \Bdash).\]
\item \cref{Jgo_k_Geq} proves that the assignment 
\[(f,\laxf) \mapsto \Jgo f\] 
is $G$-equivariant.
\end{itemize}

\begin{lemma}\label{Jgof_obj_welldef}
The pair defined in \cref{Jgof_m_objects},
\[\big((\Jgo f)_{\angordmdot} (a,\glu^a), \glu\big),\]
is an $\bm$-system in $\B$.
\end{lemma}

\begin{proof}
In this proof, we abbreviate $\big((\Jgo f)_{\angordmdot} (a,\glu^a), \glu\big)$ to $(z,\glu)$.  We need to check that $(z,\glu)$ satisfies the axioms \cref{system_obj_unity,system_naturality,system_unity_i,system_unity_iii,system_equivariance,system_associativity,system_commutativity} in \cref{def:nsystem} for an $\bm$-system in $\B$.  Suppose
\[\angs = \ang{\ang{s_{i,j} \subseteq \ufs{m}_{i,j}}_{j \in \ufs{r}_i}}_{i \in \ufs{k}}\]
is a marker, with $s_{i\crdot} = \ang{s_{i,j}}_{j \in \ufs{r}_i}$, as defined in \cref{marker_sij}.

\parhead{Object unity}.  To check the object unity axiom \cref{system_obj_unity}, suppose $s_{\im,\jm} = \emptyset$ for some $(\im,\jm) \in \ufs{k} \times \ufs{r}_\im$.  Then 
\[z_{\angs} = f\ang{a_{i,\, s_{i\crdot}}}_{i \in \ufs{k}} = \zero \in \B\]
by 
\begin{itemize}
\item the definition \cref{Jgof_m_obj_comp}, 
\item the object unity axiom 
\[a_{\im,\, s_{\im\crdot}} = \zero \in \A_\im\]
for the $\angordmim$-system $(a_\im, \glu^{a_\im})$ in $\A_\im$, and 
\item the basepoint condition \cref{laxf_basept} for $f$.
\end{itemize}

\parhead{Naturality}.  For the naturality axiom \cref{system_naturality}, we assume the context of \cref{Jgof_m_gluing} and that $h \cn x \to y$ is a morphism in $\Op(t)$.  The naturality diagram \cref{system_naturality} for $(z,\glu)$ is the boundary diagram below.
\[\begin{tikzpicture}[vcenter]
\def\h{3.5} \def\v{-1.3}
\draw[0cell=.85]
(0,0) node (a1) {\gaB_t\big(x \sscs \ang{\fgr a_{p,\, (s_{p\crdot} \,\compq\, s_{p,q,\ell})}}_{\ell \in \ufs{t}}\big)}
(a1)++(\h,\v) node (a2) {\fgr \gaAp_t \big(x \sscs \ang{a_{p,\, (s_{p\crdot} \,\compq\, s_{p,q,\ell})}}_{\ell \in \ufs{t}} \big)}
(a2)++(\h,\v) node (a3) {\fgr a_{p,\, s_{p\crdot}}}
(a1)++(0,2*\v) node (b1) {\gaB_t\big(y \sscs \ang{\fgr a_{p,\, (s_{p\crdot} \,\compq\, s_{p,q,\ell})}}_{\ell \in \ufs{t}}\big)}
(a2)++(0,2*\v) node (b2) {\fgr \gaAp_t \big(y \sscs \ang{a_{p,\, (s_{p\crdot} \,\compq\, s_{p,q,\ell})}}_{\ell \in \ufs{t}} \big)}
;
\draw[1cell=.8]
(a1) edge node[pos=.75] {\laxf_{t,p}} (a2)
(a2) edge node[pos=.6] {\fgr \glu^{a_p}_{x;\, s_{p\crdot} \scs q, \ang{s_{p,q,\ell}}_{\ell \in \ufs{t}}}} (a3)
(a1) edge[transform canvas={xshift=1ex}] node[swap] {\gaB_t(h \sscs \ang{1}_{\ell \in \ufs{t}})} (b1)
(b1) edge node[swap,pos=.25] {\laxf_{t,p}} (b2)
(b2) edge node[swap,pos=.6] {\fgr \glu^{a_p}_{y;\, s_{p\crdot} \scs q, \ang{s_{p,q,\ell}}_{\ell \in \ufs{t}}}} (a3)
(a2) edge[transform canvas={xshift=-1em}] node {\fgr \gaAp_t(h \sscs \ang{1}_{\ell \in \ufs{t}})} (b2)
;
\end{tikzpicture}\]
In the diagram above, the left region commutes by the naturality of the $(t,p)$-action constraint $\laxf_{t,p}$ \cref{klax_constraint}.  The right region commutes by 
\begin{itemize}
\item the functoriality of $\fgr$ and 
\item the naturality axiom for the $\angordmpe$-system $(a_p, \glu^{a_p})$ in $\A_p$.
\end{itemize}

\parhead{Unity}. To prove the unity axiom \cref{system_unity_i} for $(z,\glu)$, we show that, in the component of $\glu$ in \cref{Jgof_m_gluing}, each of its two constituent morphisms---$\laxf_{t,p}$ and $\fgr \glu^{a_p}_{\cdots}$---is equal to $1_\zero$ in $\B$.  There are two cases.
\begin{enumerate}
\item Suppose that $s_{i,j} = \emptyset$ for some $i \in \ufs{k} \setminus \{p\}$ and some $j \in \ufs{r}_i$.  The object unity axiom \cref{system_obj_unity} for $a_i$ implies 
\[a_{i,\, s_{i\crdot}} = \zero \in \A_i.\] 
\begin{itemize}
\item $\laxf_{t,p} = 1_\zero$ by \cref{laxf_unity_properties} \eqref{laxf_u_ii}. 
\item $\fgr \glu^{a_p}_{\cdots} = 1_\zero$ by the basepoint condition \cref{laxf_basept} for $f$.
\end{itemize} 
\item Suppose $s_{p,j} = \emptyset$ for some $j \in \ufs{r}_p$.  The object unity axiom \cref{system_obj_unity} for $a_p$ implies
\[a_{p,\, (s_{p\crdot} \,\compq\, s_{p,q,\ell})} = a_{p,\, s_{p\crdot}} = \zero \in \A_p.\]
\begin{itemize}
\item $\laxf_{t,p} = 1_\zero$ by \cref{laxf_unity_properties} \eqref{laxf_u_i}.
\item The unity axiom \cref{system_unity_i} for $(a_p,\glu^{a_p})$ implies
\[\glu^{a_p}_{x;\, s_{p\crdot} \scs q, \ang{s_{p,q,\ell}}_{\ell \in \ufs{t}}} = 1_\zero \inspace \A_p.\]
Thus, the basepoint condition \cref{laxf_basept} for $f$ implies 
\[\fgr \glu^{a_p}_{x;\, s_{p\crdot} \scs q, \ang{s_{p,q,\ell}}_{\ell \in \ufs{t}}} = 1_\zero \inspace \B.\]
\end{itemize}
\end{enumerate}
In either case, the gluing morphism $\glu$ in \cref{Jgof_m_gluing} is equal to $1_\zero$.  

To prove the second unity axiom \cref{system_unity_iii}, suppose in \cref{Jgof_m_gluing} that $t = 1$, which implies $x = \opu \in \Op(1)$, the operadic unit.
\begin{itemize}
\item $\laxf_{1,p} = 1_{f\ang{a_{i,\, s_{i\crdot}}}_{i \in \ufs{k}}}$ in $\B$ by the unity axiom \cref{laxf_unity} for $f$.
\item $\glu^{a_p}_{\opu;\, s_{p\crdot} \scs q, \{s_{p,q}\}} = 1_{a_{p,\, s_{p\crdot}}}$ in $\A_p$ by the second unity axiom \cref{system_unity_iii} for $(a_p,\glu^{a_p})$.  The functoriality of $\fgr$ implies 
\[\fgr \glu^{a_p}_{\opu;\, s_{p\crdot} \scs q, \{s_{p,q}\}} = 1_{f\ang{a_{i,\, s_{i\crdot}}}_{i \in \ufs{k}}} \inspace \B.\]
\end{itemize}
Thus, the gluing morphism $\glu$ in \cref{Jgof_m_gluing} is the identity morphism.  

\parhead{Equivariance}.  In the context of \cref{Jgof_m_gluing}, suppose $\sigma \in \Sigma_t$ is a permutation.  The equivariance diagram \cref{system_equivariance} for $(z,\glu)$ is the boundary of the following diagram.
\[\begin{tikzpicture}[vcenter]
\def\h{3.5} \def\v{-1.3}
\draw[0cell=.85]
(0,0) node (a1) {\gaB_t\big(x\sigma \sscs \ang{\fgr a_{p,\, (s_{p\crdot} \,\compq\, s_{p,q,\ell})}}_{\ell \in \ufs{t}}\big)}
(a1)++(\h,\v) node (a2) {\fgr \gaAp_t \big(x\sigma \sscs \ang{a_{p,\, (s_{p\crdot} \,\compq\, s_{p,q,\ell})}}_{\ell \in \ufs{t}} \big)}
(a2)++(\h,\v) node (a3) {\fgr a_{p,\, s_{p\crdot}}}
(a1)++(0,2*\v) node (b1) {\gaB_t\big(x \sscs \ang{\fgr a_{p,\, (s_{p\crdot} \,\compq\, s_{p,q,\sigmainv(\ell)})}}_{\ell \in \ufs{t}}\big)}
(a2)++(0,2*\v) node (b2) {\fgr \gaAp_t \big(x \sscs \ang{a_{p,\, (s_{p\crdot} \,\compq\, s_{p,q,\sigmainv(\ell)})}}_{\ell \in \ufs{t}} \big)}
;
\draw[1cell=.8]
(a1) edge node[pos=.75] {\laxf_{t,p}} (a2)
(a2) edge node[pos=.5] {\fgr \glu^{a_p}_{x\sigma ;\, s_{p\crdot} \scs q, \ang{s_{p,q,\ell}}_{\ell \in \ufs{t}}}} (a3)
(a1) edge[equal] (b1)
(b1) edge node[swap,pos=.25] {\laxf_{t,p}} (b2)
(b2) edge node[swap,pos=.5] {\fgr \glu^{a_p}_{x;\, s_{p\crdot} \scs q, \ang{s_{p,q,\sigmainv(\ell)}}_{\ell \in \ufs{t}}}} (a3)
(a2) edge[equal] (b2)
;
\end{tikzpicture}\]
In the diagram above, the left region commutes by the equivariance axiom \cref{laxf_eq} for $(f,\laxf)$.  The right region commutes by
\begin{itemize}
\item the functoriality of $\fgr$ and
\item the equivariance axiom \cref{system_equivariance} for $(a_p,\glu^{a_p})$.
\end{itemize}

\parhead{Associativity}.  To prove the associativity axiom \cref{system_associativity} for $(z,\glu)$, suppose we are given objects
\[\begin{split}
(x \sscs \ang{x_\ell}_{\ell \in \ufs{t}}) &\in \Op(t) \times \txprod_{\ell \in \ufs{t}}\, \Op(v_\ell) \andspace\\
\bx = \ga(x \sscs \ang{x_\ell}_{\ell \in \ufs{t}}) &\in \Op(v)
\end{split}\]
with $v = \sum_{\ell \in \ufs{t}}\, v_\ell$, a pair of indices $(p,q) \in \ufs{k} \times \ufs{r}_p$, and partitions
\[s_{p,q} = \coprod_{\ell \in \ufs{t}}\, s_{p,q,\ell} \andspace 
s_{p,q,\ell} = \coprod_{c \in \ufs{v}_\ell}\, s_{p,q,\ell,c} \subseteq \ufs{m}_{p,q}.\]
We use the following notation.
\[\left\{\begin{gathered}
\begin{aligned}
a_p^{q\ell} &= a_{p,\, (s_{p\crdot} \,\compq\, s_{p,q,\ell})} \phantom{M} &
a_p^{q\ell c} &= a_{p,\, (s_{p\crdot} \,\compq\, s_{p,q,\ell,c})}
\end{aligned}\\
\scalebox{.9}{$\begin{aligned}
\glu^{a_p}_{x_\ell} = \glu^{a_p}_{x_\ell;\, (s_{p\crdot} \,\compq\, s_{p,q,\ell}) \scs q,\, \ang{s_{p,q,\ell,c}}_{c \in \ufs{v}_{\ell}}} & \cn \gaAp_{v_\ell}\big(x_\ell \sscs \ang{a_p^{q\ell c}}_{c \in \ufs{v}_\ell} \big) \to a_p^{q\ell}\\
\glu^{a_p}_x = \glu^{a_p}_{x;\, s_{p\crdot} \scs q,\, \ang{s_{p,q,\ell}}_{\ell \in \ufs{t}}} & \cn \gaAp_t\big(x \sscs \ang{a_p^{q\ell}}_{\ell \in \ufs{t}}\big) \to a_{p,\, s_{p\crdot}}\\
\glu^{a_p}_{\bx} = \glu^{a_p}_{\bx;\, s_{p\crdot} \scs q,\, \ang{\ang{s_{p,q,\ell,c}}_{c\in \ufs{v}_\ell}}_{\ell \in \ufs{t}}} & \cn \gaAp_v\big(\bx \sscs \ang{\ang{a_p^{q\ell c}}_{c\in \ufs{v}_\ell}}_{\ell \in \ufs{t}} \big) \to a_{p,\, s_{p\crdot}}
\end{aligned}$}
\end{gathered}\right.\]
With $\fgr$ defined as in \cref{ftil_ApB}, the associativity diagram \cref{system_associativity} for $(z,\glu)$ is the boundary diagram below.
\[\begin{tikzpicture}[vcenter]
\def\h{7} \def\g{1} \def\v{-1.4} \def\u{-1}
\draw[0cell=.75]
(0,0) node (a1) {\gaB_t\big(x \sscs \bang{\gaB_{v_\ell} \big(x_\ell \sscs \ang{\fgr a_p^{q\ell c}}_{c \in \ufs{v}_\ell}\big) }_{\ell \in \ufs{t}} \big)}
(a1)++(\h-2*\g,0) node (a2) {\gaB_v\big(\bx \sscs \ang{\ang{\fgr a_p^{q\ell c}}_{c \in \ufs{v}_\ell}}_{\ell \in \ufs{t}}\big)}
(a2)++(\g,\v) node (a3) {\fgr \gaAp_v\big(\bx \sscs \ang{\ang{a_p^{q\ell c}}_{c \in \ufs{v}_\ell}}_{\ell \in \ufs{t}} \big)}
(a3)++(0,2*\u) node (a4) {\fgr a_{p,\, s_{p\crdot}}}
(a1)++(-\g,\v) node (b1) {\gaB_t\big(x \sscs \bang{\fgr \gaAp_{v_\ell}\big(x_\ell \sscs \ang{a_p^{q\ell c}}_{c \in \ufs{v}_\ell} \big) }_{\ell \in \ufs{t}}\big)}
(b1)++(0,2*\u) node (b2) {\gaB_t\big(x \sscs \bang{\fgr a_p^{q\ell}}_{\ell \in \ufs{t}}\big)}
(b2)++(\h/2,\u) node (b3) {\fgr \gaAp_t\big(x \sscs \ang{a_p^{q\ell}}_{\ell \in \ufs{t}}\big)}
(b1)++(\h/2,\u) node (b) {\fgr \gaAp_t \big(x \sscs \bang{\gaAp_{v_\ell} \big(x_\ell \sscs \ang{a_p^{q\ell c}}_{c \in \ufs{v}_\ell} \big)}_{\ell \in \ufs{t}}\big)}
;
\draw[1cell=.75]
(a1) edge node {\phiB} (a2)
(a2) edge node[pos=.8] {\laxf_{v,p}} (a3)
(a3) edge node {\fgr \glu^{a_p}_{\bx}} (a4)
(a1) edge node[swap,pos=.8] {\gaB_t(1_x \sscs \ang{\laxf_{v_\ell,p}}_{\ell \in \ufs{t}})} (b1)
(b1) edge node[swap] {\gaB_t(1_x \sscs \ang{\fgr \glu^{a_p}_{x_\ell}}_{\ell \in \ufs{t}})} (b2)
(b2) edge node[swap] {\laxf_{t,p}} (b3)
(b3) edge node[swap] {\fgr \glu^{a_p}_{x}} (a4)
(b1) edge node[pos=.8] {\laxf_{t,p}} (b)
(b) edge node[pos=.4] {\fgr \phiAp} (a3)
(b) edge node[pos=.35] {\fgr \gaAp_t(1_x \sscs \ang{\glu^{a_p}_{x_\ell}}_{\ell \in \ufs{t}})} (b3)
;
\end{tikzpicture}\]
In the diagram above, $\phiB$ and $\phiAp$ are the associativity constraints \cref{phiA} of, respectively, $\B$ and $\A_p$.  The three regions commute for the following reasons.
\begin{itemize}
\item The top pentagon commutes by the associativity axiom \cref{laxf_associativity} for $f$. 
\item The lower left quadrilateral commutes by the naturality of the $(t,p)$-action constraint $\laxf_{t,p}$.  
\item The lower right quadrilateral commutes by 
\begin{itemize}
\item the functoriality of $\fgr$ and 
\item the associativity axiom \cref{system_associativity} for $(a_p,\glu^{a_p})$.
\end{itemize}
\end{itemize}

\parhead{Commutativity}.  To prove the commutativity axiom \cref{system_commutativity} for $(z,\glu)$, suppose we are given a pair of objects 
\[(x,y) \in \Op(t) \times \Op(v).\]
The proof of the commutativity axiom splits into two cases.

\parhead{Case 1}.  We consider two pairs of indices
\[(p,q) \in \ufs{k} \times \ufs{r}_p \andspace (c,d) \in \ufs{k} \times \ufs{r}_c\]
with $p<c$, and partitions
\[s_{p,q} = \coprod_{\ell \in \ufs{t}}\, s_{p,q,\ell} \subseteq \ufs{m}_{p,q} \andspace
s_{c,d} = \coprod_{h \in \ufs{v}}\, s_{c,d,h} \subseteq \ufs{m}_{c,d}.\]
Using the $\comp_?$ notation in \cref{compkl}, we define the pointed functor
\[\fhat = f\big(\ang{a_{i,\, s_{i\crdot}}}_{i \in \ufs{k}} \compp - \compc - \big) \cn \A_p \times \A_c \to \B.\]
In terms of $\fhat$, the $\angs$-component object of $(z,\glu)$ in \cref{Jgof_m_obj_comp} is 
\[z_{\angs} = \fhat\big(a_{p,\, s_{p\crdot}} \scs a_{c,\, s_{c\crdot}}\big).\]
We use the following notation. 
\[\left\{\begin{gathered}
\begin{aligned}
\ang{\cdots}_\ell &= \ang{\cdots}_{\ell \in \ufs{t}} & 
\ang{\cdots}_h &= \ang{\cdots}_{h \in \ufs{v}} \\
\ang{\cdots}_{\ell,h} &= \ang{\ang{\cdots}_{\ell \in \ufs{t}}}_{h \in \ufs{v}} &
\ang{\cdots}_{h,\ell} &= \ang{\ang{\cdots}_{h \in \ufs{v}}}_{\ell \in \ufs{t}} 
\end{aligned}\\
\begin{aligned}
a_p^{q\ell} &= a_{p,\, (s_{p\crdot} \,\compq\, s_{p,q,\ell})} \\
a_c^{dh} &= a_{c,\, (s_{c\crdot} \,\compd\, s_{c,d,h})} \\
a_{p,x}^{q\ell} &= \gaAp_t\big(x \sscs \ang{a_p^{q\ell}}_{\ell} \big) 
\fto{\glu^{a_p}_x = \glu^{a_p}_{x;\, s_{p\crdot},\, q,\, \ang{s_{p,q,\ell}}_{\ell}}} 
a_{p,\, s_{p\crdot}}\\
a_{c,y}^{dh} &= \gaAc_v\big(y \sscs \ang{a_c^{dh}}_{h} \big)
\fto{\glu^{a_c}_y = \glu^{a_c}_{y;\, s_{c\crdot},\, d,\, \ang{s_{c,d,h}}_{h}}} a_{c,\, s_{c\crdot}}
\end{aligned}
\end{gathered}\right.\]
The commutativity diagram \cref{system_commutativity} for $(z,\glu)$ is the boundary diagram below.
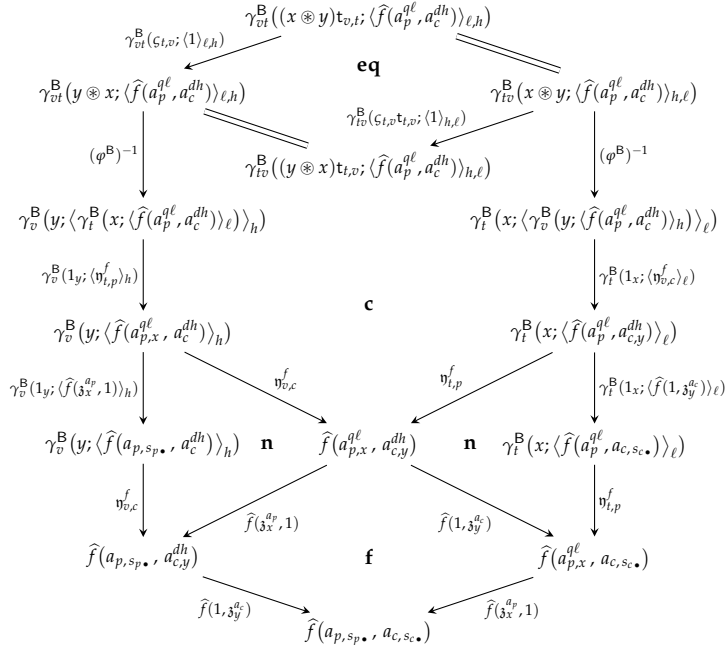
\[\begin{tikzpicture}[vcenter]
\def\g{3} \def\u{-1} \def\v{-1.7} \def\w{-1.5} \def\e{1em}
\draw[0cell=.7]
(0,0) node (a1) {\gaB_{vt}\big((x \intr y)\twist_{v,t} \sscs \ang{\fhat(a_p^{q\ell}, a_c^{dh})}_{\ell,h}\big)}
(a1)++(0,2*\u) node (a') {\gaB_{tv}\big((y \intr x)\twist_{t,v} \sscs \ang{\fhat(a_p^{q\ell}, a_c^{dh})}_{h,\ell} \big)}
(a1)++(-\g,\u) node (a21) {\gaB_{vt}\big(y \intr x \sscs \ang{\fhat(a_p^{q\ell}, a_c^{dh})}_{\ell,h}\big)}
(a1)++(\g,\u) node (a22) {\gaB_{tv}\big(x \intr y \sscs \ang{\fhat(a_p^{q\ell}, a_c^{dh})}_{h,\ell}\big)}
(a21)++(0,\v) node (a31) {\gaB_v\big(y \sscs \bang{\gaB_t\big(x \sscs \ang{\fhat(a_p^{q\ell}, a_c^{dh})}_\ell \big)}_h \big)}
(a22)++(0,\v) node (a32) {\gaB_t\big(x \sscs \bang{\gaB_v\big(y \sscs \ang{\fhat(a_p^{q\ell}, a_c^{dh})}_h \big)}_\ell \big)}
(a31)++(0,\w) node (a41) {\gaB_v\big(y \sscs \bang{\fhat(a_{p,x}^{q\ell} \scs a_c^{dh})}_h \big)}
(a32)++(0,\w) node (a42) {\gaB_t\big(x \sscs \bang{\fhat(a_p^{q\ell}, a_{c,y}^{dh})}_\ell \big)}
(a41)++(\g,\w) node (a'') {\fhat\big(a_{p,x}^{q\ell} \scs a_{c,y}^{dh}\big)}
(a41)++(0,\w) node (a51) {\gaB_v\big(y \sscs \bang{\fhat\big(a_{p,\, s_{p\crdot}} \scs a_c^{dh}\big)}_h \big)}
(a42)++(0,\w) node (a52) {\gaB_t\big(x \sscs \bang{\fhat\big(a_p^{q\ell}, a_{c,\, s_{c\crdot}}\big)}_\ell \big)}
(a51)++(0,\w) node (a61) {\fhat\big(a_{p,\, s_{p\crdot}} \scs a_{c,y}^{dh} \big)}
(a52)++(0,\w) node (a62) {\fhat\big(a_{p,x}^{q\ell} \scs a_{c,\, s_{c\crdot}}\big)}
(a61)++(\g,\u) node (a7) {\fhat\big(a_{p,\, s_{p\crdot}} \scs a_{c,\, s_{c\crdot}}\big)}
;
\draw[0cell=.8]
node[between=a1 and a' at .35] {\mathbf{eq}}
node[between=a' and a'' at .5] {\mathbf{c}}
node[between=a'' and a7 at .6] {\mathbf{f}}
node[between=a'' and a51 at .45, shift={(0,0)}] {\mathbf{n}}
node[between=a'' and a52 at .45, shift={(0,0)}] {\mathbf{n}}
;
\draw[1cell=.65]
(a1) edge[transform canvas={xshift=-\e}, shorten >=.5ex] node[swap] {\gaB_{vt}(\pcom_{t,v} \sscs \ang{1}_{\ell,h})} (a21)
(a21) edge node[swap] {(\phiB)^{-1}} (a31)
(a31) edge node[swap] {\gaB_v(1_y \sscs \ang{\laxf_{t,p}}_h)} (a41)
(a41) edge node[swap] {\gaB_v(1_y \sscs \ang{\fhat(\glu^{a_p}_x , 1)}_h)} (a51)
(a51) edge node[swap] {\laxf_{v,c}} (a61)
(a61) edge node[swap] {\fhat(1, \glu^{a_c}_y)} (a7)
(a1) edge[equal,transform canvas={xshift=\e}] (a22)
(a22) edge node {(\phiB)^{-1}} (a32)
(a32) edge node {\gaB_t(1_x \sscs \ang{\laxf_{v,c}}_\ell)} (a42)
(a42) edge node {\gaB_t(1_x \sscs \ang{\fhat(1, \glu^{a_c}_y)}_\ell)} (a52)
(a52) edge node {\laxf_{t,p}} (a62)
(a62) edge node {\fhat(\glu^{a_p}_x, 1)} (a7)
(a21) edge[equal,transform canvas={xshift=0}] (a')
(a22) edge[transform canvas={xshift=0}] node[swap,pos=.6] {\gaB_{tv}(\pcom_{t,v}\twist_{t,v} \sscs \ang{1}_{h,\ell})} (a')
(a41) edge node[pos=.6] {\laxf_{v,c}} (a'')
(a42) edge node[swap,pos=.6] {\laxf_{t,p}} (a'')
(a'') edge node[pos=.6] {\fhat(\glu^{a_p}_x , 1)} (a61)
(a'') edge node[swap,pos=.6] {\fhat(1,\glu^{a_c}_y)} (a62)
;
\end{tikzpicture}\]
The following statements hold for the diagram above.
\begin{itemize}
\item The top quadrilateral labeled $\mathbf{eq}$ commutes by the action equivariance axiom \cref{pseudoalg_action_sym} for $\B$, applied to the transpose permutation $\twist_{t,v}$ \cref{eq:transpose_perm}.
\item The middle region labeled $\mathbf{c}$ commutes by the commutativity axiom \cref{laxf_com} for $(f,\laxf)$.
\item The bottom quadrilateral labeled $\mathbf{f}$ commutes by the functoriality of $\fhat$.
\item The left and right triangles labeled $\mathbf{n}$ commute by the naturality of, respectively, $\laxf_{v,c}$ and $\laxf_{t,p}$.
\end{itemize}
This proves the commutativity axiom \cref{system_commutativity} for $(z,\glu)$ in the present case.  Note that this case does not use the commutativity axiom for $(a_i, \glu^{a_i}) \in \Aiangordmi$ for any $i \in \ufs{k}$.

\parhead{Case 2}.  We consider two pairs of indices
\[(p,q) \in \ufs{k} \times \ufs{r}_p \andspace (p,d) \in \ufs{k} \times \ufs{r}_p\]
with $q<d$, and partitions
\[s_{p,q} = \coprod_{\ell \in \ufs{t}}\, s_{p,q,\ell} \subseteq \ufs{m}_{p,q} \andspace
s_{p,d} = \coprod_{h \in \ufs{v}}\, s_{p,d,h} \subseteq \ufs{m}_{p,d}.\]
We use the following notation.
\[\left\{\begin{gathered}
\scalebox{.9}{$
\ang{\cdots}_\ell = \ang{\cdots}_{\ell \in \ufs{t}} \qquad
\ang{\cdots}_h = \ang{\cdots}_{h \in \ufs{v}} \qquad
\ang{\cdots}_{\ell,h} = \ang{\ang{\cdots}_{\ell \in \ufs{t}}}_{h \in \ufs{v}}$}
\\
a_p^{q\ell,dh} = a_{p,\,(s_{p\crdot} \,\compq\, s_{p,q,\ell} \,\compd\, s_{p,d,h})} \qquad
\ang{\cdots}_{h,\ell} = \ang{\ang{\cdots}_{h \in \ufs{v}}}_{\ell \in \ufs{t}}
\\
\scalebox{.8}{$\begin{aligned}
a_{p,x}^{q\ell,dh} &= \gaAp_t\big(x \sscs \ang{a_p^{q\ell,dh}}_\ell\big) 
\fto{\glu^{a_p}_{x,\, dh} = \glu^{a_p}_{x;\, (s_{p\crdot} \compd\, s_{p,d,h}),\, q,\, \ang{s_{p,q,\ell}}_\ell}}
a_p^{dh} = a_{p,\, (s_{p\crdot} \,\compd\, s_{p,d,h})}
\\
a_{p,y}^{q\ell,dh} &= \gaAp_v\big(y \sscs \ang{a_p^{q\ell,dh}}_h\big)
\fto{\glu^{a_p}_{y,\, q\ell} = \glu^{a_p}_{y;\, (s_{p\crdot} \compq\, s_{p,q,\ell}),\, d,\, \ang{s_{p,d,h}}_h}}
a_p^{q\ell} = a_{p,\, (s_{p\crdot} \,\compq\, s_{p,q,\ell})}
\\
a_{p,y}^{dh} &= \gaAp_v\big(y \sscs \ang{a_p^{dh}}_{h} \big) 
\fto{\glu^{a_p}_y = \glu^{a_p}_{y;\, s_{p\crdot},\, d,\, \ang{s_{p,d,h}}_{h}}} a_{p,\, s_{p\crdot}}
\\
a_{p,x}^{q\ell} &= \gaAp_t\big(x \sscs \ang{a_p^{q\ell}}_{\ell} \big)
\fto{\glu^{a_p}_x = \glu^{a_p}_{x;\, s_{p\crdot},\, q,\, \ang{s_{p,q,\ell}}_{\ell}}} a_{p,\, s_{p\crdot}}
\end{aligned}$}
\end{gathered}\right.\]
With $\fgr$ defined as in \cref{ftil_ApB}, the commutativity diagram \cref{system_commutativity} for $(z,\glu)$ is the boundary diagram below.
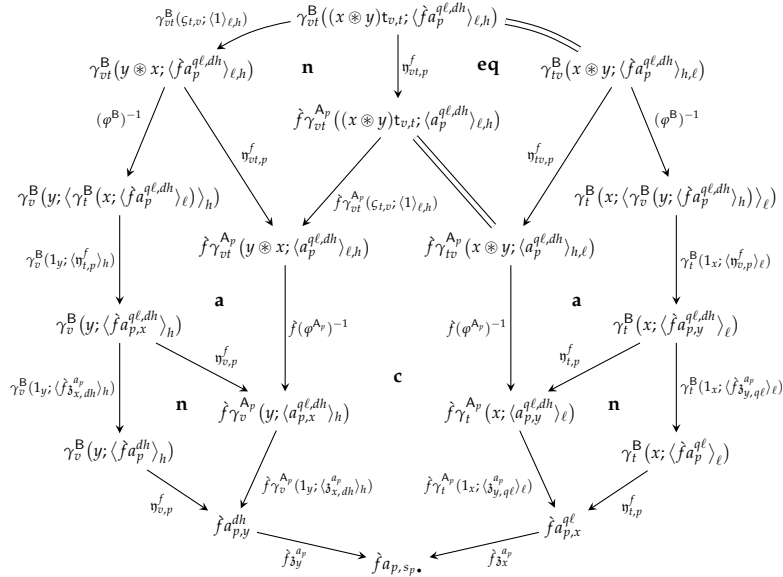
\[\begin{tikzpicture}[vcenter]
\def\f{3} \def\h{.7} \def\g{1.5} \def\t{-.5} \def\u{-.7} \def\v{-1.7} \def\w{-1} \def\x{-2.2} \def\e{1em} \def\d{1.5} \def\s{10}
\draw[0cell=.65]
(0,0) node (a1) {\gaB_{vt}\big((x \intr y)\twist_{v,t} \sscs \ang{\fgr a_p^{q\ell,dh}}_{\ell,h}\big)}
(a1)++(-\f,\u) node (a21) {\gaB_{vt}\big(y \intr x \sscs \ang{\fgr a_p^{q\ell,dh}}_{\ell,h}\big)}
(a1)++(\f,\u) node (a22) {\gaB_{tv}\big(x \intr y \sscs \ang{\fgr a_p^{q\ell,dh}}_{h,\ell}\big)}
(a21)++(-\h,\v) node (a31) {\gaB_v\big(y \sscs \bang{\gaB_t\big(x \sscs \ang{\fgr a_p^{q\ell,dh}}_\ell \big)}_h \big)}
(a22)++(\h,\v) node (a32) {\gaB_t\big(x \sscs \bang{\gaB_v\big(y \sscs \ang{\fgr a_p^{q\ell,dh}}_h \big)}_\ell \big)}
(a31)++(0,\v) node (a41) {\gaB_v\big(y \sscs \bang{\fgr a_{p,x}^{q\ell, dh}}_h \big)}
(a32)++(0,\v) node (a42) {\gaB_t\big(x \sscs \bang{\fgr a_{p,y}^{q\ell, dh}}_\ell \big)}
(a41)++(0,\v) node (a51) {\gaB_v\big(y \sscs \bang{\fgr a_p^{dh}}_h \big)}
(a42)++(0,\v) node (a52) {\gaB_t\big(x \sscs \bang{\fgr a_p^{q\ell}}_\ell \big)}
(a51)++(\g,\w) node (a61) {\fgr a_{p,y}^{dh}}
(a52)++(-\g,\w) node (a62) {\fgr a_{p,x}^{q\ell}}
(a61)++(\f+\h-\g,\t) node (a7) {\fgr a_{p,\, s_{p\crdot}}}
;
\draw[0cell=.65]
(a1)++(0,.8*\v) node (b1) {\fgr \gaAp_{vt}\big((x \intr y)\twist_{v,t} \sscs \ang{a_p^{q\ell,dh}}_{\ell,h}\big)}
(b1)++(-\d,\v) node (b21) {\fgr \gaAp_{vt}\big(y \intr x \sscs \ang{a_p^{q\ell,dh}}_{\ell,h}\big)}
(b1)++(\d,\v) node (b22) {\fgr \gaAp_{tv}\big(x \intr y \sscs \ang{a_p^{q\ell,dh}}_{h,\ell}\big)}
(b21)++(0,\x) node (b31) {\fgr \gaAp_v\big(y \sscs \ang{a_{p,x}^{q\ell,dh}}_h \big)}
(b22)++(0,\x) node (b32) {\fgr \gaAp_t\big(x \sscs \ang{a_{p,y}^{q\ell,dh}}_\ell \big)}
;
\draw[0cell=.75]
node[between=a21 and a22 at .3] {\mathbf{n}}
node[between=a21 and a22 at .7] {\mathbf{eq}}
node[between=a41 and b21 at .5, shift={(.3,-.3)}] {\mathbf{a}}
node[between=a42 and b22 at .5, shift={(-.3,-.3)}] {\mathbf{a}}
node[between=a41 and a61 at .4, shift={(.3,0)}] {\mathbf{n}}
node[between=a42 and a62 at .4, shift={(-.3,0)}] {\mathbf{n}}
node[between=b31 and b32 at .5, shift={(0,.6)}] {\mathbf{c}}
;
\draw[1cell=.6]
(a1) edge[bend right=\s] node[swap] {\gaB_{vt}(\pcom_{t,v} \sscs \ang{1}_{\ell,h})} (a21)
(a21) edge node[swap] {(\phiB)^{-1}} (a31)
(a31) edge node[swap] {\gaB_v(1_y \sscs \ang{\laxf_{t,p}}_h)} (a41)
(a41) edge node[swap] {\gaB_v(1_y \sscs \ang{\fgr \glu^{a_p}_{x,\,dh}}_h)} (a51)
(a51) edge node[swap] {\laxf_{v,p}} (a61)
(a61) edge node[swap] {\fgr \glu^{a_p}_y} (a7)
(a1) edge[equal,bend left=\s] (a22)
(a22) edge node {(\phiB)^{-1}} (a32)
(a32) edge node {\gaB_t(1_x \sscs \ang{\laxf_{v,p}}_\ell)} (a42)
(a42) edge node {\gaB_t(1_x \sscs \ang{\fgr \glu^{a_p}_{y,\,q\ell}}_\ell)} (a52)
(a52) edge node {\laxf_{t,p}} (a62)
(a62) edge node {\fgr \glu^{a_p}_x} (a7)
;
\draw[1cell=.6]
(a1) edge node {\laxf_{vt,p}} (b1)
(a21) edge node[pos=.6] {\laxf_{vt,p}} (b21)
(a22) edge node[swap,pos=.6] {\laxf_{tv,p}} (b22)
(b1) edge node[pos=.6,inner sep=0pt] {\fgr \gaAp_{vt}(\pcom_{t,v} \sscs \ang{1}_{\ell,h})} (b21)
(b1) edge[equal] (b22)
(b21) edge node {\fgr (\phiAp)^{-1}} (b31)
(b22) edge node[swap] {\fgr (\phiAp)^{-1}} (b32)
(b31) edge node {\fgr \gaAp_v(1_y \sscs \ang{\glu^{a_p}_{x,\,dh}}_h)} (a61)
(b32) edge node[swap] {\fgr \gaAp_t(1_x \sscs \ang{\glu^{a_p}_{y,\,q\ell}}_\ell)} (a62)
(a41) edge node[pos=.6] {\laxf_{v,p}} (b31)
(a42) edge node[swap,pos=.6] {\laxf_{t,p}} (b32)
;
\end{tikzpicture}\]
The following statements hold for the diagram above.
\begin{itemize}
\item The quadrilateral labeled $\mathbf{eq}$ commutes by the action equivariance axiom \cref{laxf_eq} for $(f,\laxf)$.
\item The three quadrilaterals labeled $\mathbf{n}$ commute by the naturality of, respectively, $\laxf_{vt,p}$, $\laxf_{v,p}$, and $\laxf_{t,p}$.
\item The two pentagons labeled $\mathbf{a}$ commute by the associativity axiom \cref{laxf_associativity} for $(f,\laxf)$.
\item The middle region labeled $\mathbf{c}$ commutes by 
\begin{itemize}
\item the functoriality of $\fgr$ and 
\item the commutativity axiom \cref{system_commutativity} for $(a^p, \glu^{a_p})$.
\end{itemize}
\end{itemize}
This proves that $(z,\glu)$ satisfies the commutativity axiom \cref{system_commutativity} and that it is an $\bm$-system in $\B$.
\end{proof}

\begin{lemma}\label{Jgof_natural}
For each $k$-lax $\Op$-morphism $(f, \laxf) \cn \ang{\A_i}_{i \in \ufs{k}} \to \B$, the data defined in \cref{Jgo_f},
\[\ang{\Aidash}_{i \in \ufs{k}} \fto{\Jgo f} \Bdash,\]
form an object in $\GGCat(\ang{\Aidash}_{i \in \ufs{k}} \sscs \Bdash)$.
\end{lemma}

\begin{proof}
By \cref{Jgof_m_ptfunctor,Jgof_obj_welldef}, each component $(\Jgo f)_{\angordmdot}$, as stated in \cref{Jgof_component}, is a well-defined pointed functor.  To prove that $\Jgo f$ is a well-defined object, by \cref{ggcat_k1_nat}, we need to show that, given morphisms
\begin{equation}\label{upomdot}
\begin{aligned}
\angordmdot = \ang{\angordmi}_{i \in \ufs{k}} 
& \fto{\upomdot = \ang{\upom^i}_{i \in \ufs{k}}} 
\angordndot = \ang{\angordni}_{i \in \ufs{k}} && \text{in $\Gsk^k$ and}\\
\bm = \txoplus_{i \in \ufs{k}}\, \angordmi 
& \fto{\upom = \txoplus_{i \in \ufs{k}}\, \upom^i}
\bn = \txoplus_{i \in \ufs{k}}\, \angordni && \text{in $\Gsk$},
\end{aligned}
\end{equation}
the following diagram of pointed functors commutes.
\begin{equation}\label{Jgof_nat_diagram}
\begin{tikzpicture}[vcenter]
\def\v{-1.4}
\draw[0cell=1]
(0,0) node (a1) {\txsma_{i \in \ufs{k}}\, \Aiangordmi}
(a1)++(3.5,0) node (a2) {\Bboldm}
(a1)++(0,\v) node (b1) {\txsma_{i \in \ufs{k}}\, \Aiangordni}
(a2)++(0,\v) node (b2) {\Bboldn}
;
\draw[1cell=.9]
(a1) edge node {(\Jgo f)_{\angordmdot}} (a2)
(b1) edge node {(\Jgo f)_{\angordndot}} (b2)
(a1) edge[transform canvas={xshift=1.3em}] node[swap] {\txsma_{i \in \ufs{k}}\, \Aiupomi} (b1)
(a2) edge node {\Bupom} (b2)
;
\end{tikzpicture}
\end{equation}

\parhead{Base case}.  First, suppose that for some $i \in \ufs{k}$,
\[\angordmi \fto{\upom^i} \angordni\]
is the 0-morphism, factoring through the basepoint $\vstar \in \Gsk$.  Then $\upom$ is also the 0-morphism by \cref{Gsk_oplus_vstar}.  By \cref{A_vstar_n,A_m_vstar}, each of the two composites in \cref{Jgof_nat_diagram} is the constant functor at the basepoint in $\Bboldn$, which is given by the base $\bn$-system $(\zero,1_\zero)$ in $\B$ (\cref{expl:nsystem_base}).  

\parhead{Non-base case}.  Suppose that, for each $i \in \ufs{k}$, the morphism
\begin{equation}\label{Fihipsii}
\begin{split}
\upom_i &= \big(h_i, \psi_i = \ang{\psi_{i,j}}_{j \in \ufs{t}_i}\big)\\ 
&= \Big(\ufs{r}_i \fto{h_i} \ufs{t}_i \scs \bang{\ordm^{i,\hinv_i(j)} 
\fto{\psi_{i,j}} \ordn^{i,j}}_{j \in \ufs{t}_i} \Big) 
\cn \angordmi \to \angordni
\end{split}
\end{equation}
is not a 0-morphism in $\Gsk$ \cref{fangpsi}. 
\begin{itemize}
\item In \cref{Fihipsii}, $r_i$ and $t_i$ are the lengths of, respectively, the objects 
\[\angordmi = \ang{\ordm^{i,d}}_{d \in \ufs{r}_i} \andspace 
\angordni = \ang{\ordn^{i,j}}_{j \in \ufs{t}_i} \in \Gsk.\] 
\item $h_i$ is a reindexing injection (\cref{def:injections}). 
\item Each $\psi_{i,j}$ is a pointed function.
\item With the notation in \cref{f_oplus_h,psi_oplus_phi,Fihipsii}, the morphism $\upom \cn \bm \to \bn$ is given by
\[\upom = \txoplus_{i \in \ufs{k}}\, (h_i, \psi_i) 
= \big(\! \txoplus_{i \in \ufs{k}}\, h_i, \txoplus_{i \in \ufs{k}}\, \psi_i \big).\]
\end{itemize}
With ${h_i}_*$ denoting the reindexing functor of $h_i$, as defined in \cref{reindexing_functor,ordn_empty}, we define the following objects.
\[\begin{split}
{h_i}_* \angordmi &= \bang{\ordm^{i,\hinv_i(j)}}_{j \in \ufs{t}_i} \in \Gsk \\
h\angordmdot &= \bang{{h_i}_* \angordmi}_{i \in \ufs{k}} \in \Gsk^k \\
h\bm &= \txoplus_{i \in \ufs{k}}\, {h_i}_* \angordmi \in \Gsk
\end{split}\]
By (i) the functoriality of the smash product and (ii) \cref{A_fangpsi} applied to $\Aiupomi$ and $\Bupom$, the diagram \cref{Jgof_nat_diagram} is the boundary diagram below, where $\oplus_i$ means $\oplus_{i \in \ufs{k}}$.
\begin{equation}\label{Jgof_nat_factor}
\begin{tikzpicture}[vcenter]
\def\v{-1.5}
\draw[0cell=.9]
(0,0) node (a1) {\txsma_{i \in \ufs{k}}\, \Aiangordmi}
(a1)++(4,0) node (a2) {\Bboldm}
(a1)++(0,\v) node (b1) {\txsma_{i \in \ufs{k}}\, \Aihiangordmi}
(a2)++(0,\v) node (b2) {\Bhboldm}
(b1)++(0,\v) node (c1) {\txsma_{i \in \ufs{k}}\, \Aiangordni}
(b2)++(0,\v) node (c2) {\Bboldn}
;
\draw[1cell=.85]
(a1) edge node {(\Jgo f)_{\angordmdot}} (a2)
(b1) edge node {(\Jgo f)_{h\angordmdot}} (b2)
(c1) edge node {(\Jgo f)_{\angordndot}} (c2)
(a1) edge[transform canvas={xshift=1em}] node[swap] {\txsma_{i \in \ufs{k}}\, \htil_i} (b1)
(a2) edge node {\wt{\oplus_i\, h_i}} (b2)
(b1) edge[transform canvas={xshift=1em}] node[swap] {\txsma_{i \in \ufs{k}}\, \psitil_i} (c1)
(b2) edge node {\wt{\oplus_i\, \psi_i}} (c2)
;
\end{tikzpicture}
\end{equation}
It suffices to show that the top and bottom halves in \cref{Jgof_nat_factor} commute.

\parhead{Top half: component objects}.  To show that the top half of \cref{Jgof_nat_factor} commutes, we consider an object
\[(a,\glu^a) = \ang{(a_i, \glu^{a_i})}_{i \in \ufs{k}} \in \txsma_{i \in \ufs{k}}\, \Aiangordmi\]
and markers
\[\angs = \ang{\ang{s_{i,j} \subseteq \ufsm^{i,\hinv_i(j)}}_{j \in \ufs{t}_i}}_{i \in \ufs{k}} 
\andspace s_{i\crdot} = \ang{s_{i,j}}_{j \in \ufs{t}_i}.\]
If any $s_{i,j} = \emptyset$, then, by the object unity axiom \cref{system_obj_unity} for $h\bm$-systems, each of the two composites in the top half of \cref{Jgof_nat_factor} yields the $\angs$-component object $\zero \in \B$.  Thus, we assume that $s_{i,j} \neq \emptyset$ for all $i \in \ufs{k}$ and $j \in \ufs{t}_i$.

Using \cref{ftil_angs,atil_component,Jgof_m_obj_comp}, the following computation shows that the two composites in the top half of \cref{Jgof_nat_factor} yield the same $\angs$-component object.
\begin{equation}\label{Jgof_nat_computation}
\begin{split}
& \big((\Jgo f)_{h\angordmdot} (\txsma_{i \in \ufs{k}}\, \htil_i) (a,\glu^a)\big)_{\angs}\\
&= f\bang{\big(\htil_i (a_i,\glu^{a_i})\big)_{s_{i\crdot}}}_{i \in \ufs{k}}\\
&= f\bang{a_{i,\, \htil_{i*} s_{i\crdot}}}_{i \in \ufs{k}}\\
&= \big((\Jgo f)_{\angordmdot} (a,\glu^a)\big)_{\ang{\htil_{i*} s_{i\crdot}}_{i \in \ufs{k}}}\\
&= \big((\Jgo f)_{\angordmdot} (a,\glu^a)\big)_{{\wt{\oplus_i\, h_i}}_* \angs}\\
&= \big(\big(\wt{\oplus_i\, h_i}\big) (\Jgo f)_{\angordmdot} (a,\glu^a)\big)_{\angs}
\end{split}
\end{equation}
The middle equalities in \cref{Jgof_nat_computation} involve the following markers.
\[\begin{split}
\htil_{i*} s_{i\crdot} &= \ang{s_{i, h_i(d)} \subseteq \ufs{m}^{i,d}}_{d \in \ufs{r}_i}\\
{\wt{\oplus_i\, h_i}}_* \angs &= \ang{\ang{s_{i,h_i(d)}}_{d \in \ufs{r}_i}}_{i \in \ufs{k}} 
= \ang{\htil_{i*} s_{i\crdot}}_{i \in \ufs{k}}
\end{split}\]
Moreover, using \cref{thatil_component,Jgof_m_theta_comp} instead of \cref{atil_component,Jgof_m_obj_comp}, the computation in \cref{Jgof_nat_computation} also proves that the two composites in the top half of \cref{Jgof_nat_factor} are equal on morphisms.

\parhead{Top half: gluing}.  Next, we show that the two composites in the top half of \cref{Jgof_nat_factor} yield the same gluing morphism.  Suppose we are given an object $x \in \Op(v)$ for some $v \geq 0$, a pair of indices $(p,q) \in \ufs{k} \times \ufs{t}_p$, and a partition
\[s_{p,q} = \coprod_{\ell \in \ufs{v}}\, s_{p,q,\ell} \subseteq \ufs{m}^{p,\hinv_p(q)}.\]
We define the index $u = \big(\sum_{i=1}^{p-1} t_i\big) + q$ and the tuple
\begin{equation}\label{xsupql}
\big(x \sscs \angs, u, \ang{s_{p,q,\ell}}_{\ell \in \ufs{v}}\big).
\end{equation}
If $v = 0$ or $v=1$, then the unity axioms \cref{system_unity_iii,system_unity_ii} imply that the gluing morphism at \cref{xsupql} is the identity.  Thus, we assume $v > 1$.  There are three cases, depending on $s_{i,j}$ and $\hinv_p(q)$.
\begin{description}
\item[Case 1] If any $s_{i,j} = \emptyset$, then, by the unity axiom \cref{system_unity_i} for $h\bm$-systems, each composite in the top half of \cref{Jgof_nat_factor} has a gluing morphism at \cref{xsupql} given by $1_\zero$ in $\B$.
\item[Case 2] Suppose all $s_{i,j} \neq \emptyset$ and $\hinv_p(q) = \emptyset$. 
\begin{itemize}
\item By \cref{glutil_component_iii}, the clockwise composite in the top half of \cref{Jgof_nat_factor} has an identity gluing morphism at \cref{xsupql}.  The next bullet point shows that the counterclockwise composite also yields an identity gluing morphism at \cref{xsupql}.
\item Since
\[\emptyset \neq s_{p,q} \subseteq \ufsm^{p,\emptyset} = \ufs{1} = \{1\},\]
it follows that
\begin{equation}\label{s_pqell}
s_{p,q,\ell} = \begin{cases}
s_{p,q} = \{1\} & \text{for one index $\ell = c \in \ufs{v}$ and}\\
\emptyset & \text{for $\ell \in \ufs{v} \setminus \{c\}$.}
\end{cases}
\end{equation}
Defining the pointed functor
\begin{equation}\label{fgr_Ap_B}
\fgr = f\big(\ang{a_{i,\, \htil_{i*} s_{i\crdot}}}_{i \in \ufs{k}} \,\compp - \big) \cn \A_p \to \B,
\end{equation}
by \cref{Jgof_m_gluing}, the counterclockwise composite in the top half of \cref{Jgof_nat_factor} has a gluing morphism at \cref{xsupql} given by the composite 
\begin{equation}\label{fgr_hpap}
\fgr \glu^{\htil_{p} a_p}_{x;\, s_{p\crdot} \csp q, \ang{s_{p,q,\ell}}_{\ell \in \ufs{v}}}\circ \laxf_{v,p}.
\end{equation}
Each of these two morphisms is an identity morphism for the following reasons.
\begin{enumerate}
\item By \cref{atil_component,s_pqell}, the first morphism in \cref{fgr_hpap} is the following component of $\laxf_{v,p}$ \cref{klax_constraint}.
\[\gaB_v\big(x; (\fgr\zero)^{c-1}, \fgr a_{p,\, \htil_{p*} s_{p\crdot}}, (\fgr\zero)^{v-c} \big) \fto{\laxf_{v,p}} 
\fgr\gaAp_v\big(x; \zero^{c-1}, a_{p,\, \htil_{p*} s_{p\crdot}}, \zero^{v-c} \big)\]
This morphism is an identity morphism by the associativity axiom \cref{laxf_associativity} for $f$, in the case where the bottom arrow in the diagram \cref{laxf_associativity} is given by $\laxf_{v,p}$ above and
\begin{equation}\label{xell_starone}
x_\ell = \begin{cases}
* \in \Op(0) & \text{if $\ell \in \ufs{v} \setminus \{c\}$, and}\\
\opu \in \Op(1) & \text{if $\ell = c$.}
\end{cases}
\end{equation}
The other four arrows in the diagram \cref{laxf_associativity} under consideration are identity morphisms for the following reasons.
\begin{itemize}
\item By \cref{xell_starone} and the basepoint axiom \cref{pseudoalg_basept_axiom} for each of $\B$ and $\A_p$, both $\phiAp$ and $\phiB$ in that diagram are identity morphisms.
\item By $\Op(1) = \{\opu\}$, \cref{xell_starone}, and the unity axiom \cref{laxf_unity} for $f$, the upper right arrow in that diagram is the identity morphism.
\item By \cref{xell_starone}, the basepoint axiom \cref{laxf_basepoint} for $f$, and the unity axiom \cref{laxf_unity} for $f$, the left vertical arrow in that diagram is the identity morphism.
\end{itemize}
\item Since all $s_{p,j} \neq \emptyset$ and $\hinv_p(q) = \emptyset$, in the second morphism in \cref{fgr_hpap}, $\glu^{\htil_p a_p}_{\cdots}$ is an identity morphism by \cref{glutil_component_iii}.  Thus, $\fgr \glu^{\htil_p a_p}_{\cdots}$ is an identity morphism by the functoriality of $\fgr$.
\end{enumerate}
This proves that the composite in \cref{fgr_hpap} is an identity morphism.
\end{itemize}
\item[Case 3] Suppose all $s_{i,j} \neq \emptyset$ and $\hinv_p(q) \neq \emptyset$.
\begin{itemize}
\item Using \cref{glutil_component,Jgof_m_gluing,fgr_Ap_B}, the counterclockwise composite in the top half of \cref{Jgof_nat_factor} has a gluing morphism at the component \cref{xsupql} given by the following composite in $\B$.
\begin{equation}\label{Jgof_nat_topglu}
\begin{tikzpicture}[vcenter]
\draw[0cell=.9]
(0,0) node (a) {\gaB_v\big(x \sscs \bang{\fgr a_{p,\, (\htil_{p*} s_{p\crdot} \,\comp_{\hinv_p(q)}\, s_{p,q,\ell})}}_{\ell \in \ufs{v}} \big)}
(a)++(1.5,-1.5) node (b) {\fgr \gaAp_v\big(x \sscs \bang{a_{p,\, (\htil_{p*} s_{p\crdot} \,\comp_{\hinv_p(q)}\, s_{p,q,\ell})}}_{\ell \in \ufs{v}} \big)} 
(a)++(4.5,0) node (c) {\fgr a_{p,\, \htil_{p*} s_{p\crdot}}}
(a)++(-1.5,0) node (a') {\phantom{\fgr_{\hinv}}}
;
\draw[1cell=.85]
(a') [rounded corners=2pt, shorten <=-0ex] |- node[pos=.2,swap] {\laxf_{v,p}} (b)
;
\draw[1cell=.85]
(b) [rounded corners=2pt, shorten <=-0ex] -| node[pos=.8] {\fgr \glu^{a_p}_{x;\, \htil_{p*} s_{p\crdot} \scs \hinv_p(q), \ang{s_{p,q,\ell}}_{\ell \in \ufs{v}}}} (c)
;
\end{tikzpicture}
\end{equation}
\item By \cref{glutil_component}, the clockwise composite in the top half of \cref{Jgof_nat_factor} has a gluing morphism at the component \cref{xsupql} given by the gluing morphism of $(\Jgo f)_{\angordmdot} (a,\glu^a)$ at the following component.
\[\begin{split}
&\big(x \sscs \wt{\oplus_i\, h_i}_* \angs, (\oplus_i\, h_i)^{-1} (u), \ang{s_{p,q,\ell}}_{\ell \in \ufs{v}}\big)\\
&= \big(x \sscs \ang{\htil_{i*} s_{i\crdot}}_{i \in \ufs{k}} \scs \big(\!\txsum_{i=1}^{p-1} r_i\big) + \hinv_p(q), \ang{s_{p,q,\ell}}_{\ell \in \ufs{v}}\big)
\end{split}\]
By \cref{Jgof_m_gluing}, this gluing morphism of $(\Jgo f)_{\angordmdot} (a,\glu^a)$ is also given by the composite in \cref{Jgof_nat_topglu}.  
\end{itemize}
\end{description}
This proves that the top half of \cref{Jgof_nat_factor} commutes.

\parhead{Bottom half: component objects}.  To show that the bottom half of \cref{Jgof_nat_factor} commutes, we consider an object
\[(a,\glu^a) = \ang{(a_i, \glu^{a_i})}_{i \in \ufs{k}} \in \txsma_{i \in \ufs{k}}\, \Aihiangordmi\]
and markers
\[\angs = \ang{\ang{s_{i,j} \subseteq \ufs{n}^{i,j}}_{j \in \ufs{t}_i}}_{i \in \ufs{k}} \andspace
s_{i\crdot} = \ang{s_{i,j}}_{j \in \ufs{t}_i}.\]
If any $s_{i,j} = \emptyset$, then, by the object unity axiom \cref{system_obj_unity} of $\bn$-systems, each of the two composites in the bottom half of \cref{Jgof_nat_factor} yields the $\angs$-component object $\zero \in \B$.  Thus, we assume all $s_{i,j} \neq \emptyset$.

Using \cref{apsitil_angs,Jgof_m_obj_comp}, the following computation shows that the two composites in the bottom half of \cref{Jgof_nat_factor} yield the same $\angs$-component object.
\begin{equation}\label{Jgof_natbot_computation}
\begin{split}
& \big((\Jgo f)_{\angordndot} (\txsma_{i \in \ufs{k}}\, \psitil_i) (a,\glu^a)\big)_{\angs}\\
&= f\bang{\big(\psitil_i (a_i,\glu^{a_i})\big)_{s_{i\crdot}}}_{i \in \ufs{k}}\\
&= f\bang{a_{i,\, \psiinv_i s_{i\crdot}}}_{i \in \ufs{k}}\\
&= \big((\Jgo f)_{h\angordmdot} (a,\glu^a) \big)_{\ang{\psiinv_i s_{i\crdot}}_{i \in \ufs{k}}}\\
&= \big((\Jgo f)_{h\angordmdot} (a,\glu^a) \big)_{(\oplus_i\, \psi_i)^{-1} \angs}\\
&= \big(\big(\wt{\oplus_i\, \psi_i}\big) (\Jgo f)_{h\angordmdot} (a,\glu^a)\big)_{\angs}
\end{split}
\end{equation}
The middle equalities in \cref{Jgof_natbot_computation} involve the following markers.
\[\begin{split}
\psiinv_i s_{i\crdot} &= \bang{\psiinv_{i,j} s_{i,j} \subseteq \ufs{m}^{i, \hinv_i(j)}}_{j \in \ufs{t}_i}\\
(\oplus_i\, \psi_i)^{-1} \angs &= \ang{\ang{\psiinv_{i,j} s_{i,j}}_{j \in \ufs{t}_i}}_{i \in \ufs{k}} 
= \ang{\psiinv_i s_{i\crdot}}_{i \in \ufs{k}}
\end{split}\]
Moreover, using \cref{thapsitil_component,Jgof_m_theta_comp} instead of \cref{apsitil_angs,Jgof_m_obj_comp}, the computation in \cref{Jgof_natbot_computation} also proves that the two composites in the bottom half of \cref{Jgof_nat_factor} are equal on morphisms.

\parhead{Bottom half: gluing}.  Next, we show that the two composites in the bottom half of \cref{Jgof_nat_factor} yield the same gluing morphism.  Suppose we are given an object $x \in \Op(v)$ for some $v \geq 0$, a pair of indices $(p,q) \in \ufs{k} \times \ufs{t}_p$, and a partition
\[s_{p,q} = \coprod_{\ell \in \ufs{v}}\, s_{p,q,\ell} \subseteq \ufs{n}^{p,q}.\]
We define the index $u = \big(\sum_{i=1}^{p-1} t_i\big) + q$, the tuple
\begin{equation}\label{xangsupql}
\big(x \sscs \angs, u, \ang{s_{p,q,\ell}}_{\ell \in \ufs{v}}\big),
\end{equation}
and the pointed functor
\[\fgr = f\big(\ang{a_{i,\, \psiinv_i s_{i\crdot}}}_{i \in \ufs{k}} \,\compp - \big) \cn \A_p \to \B.\]
\begin{itemize}
\item By \cref{glupsitil,Jgof_m_gluing}, the counterclockwise composite in the bottom half of \cref{Jgof_nat_factor} has a gluing morphism at the component \cref{xangsupql} given by the following composite in $\B$.
\begin{equation}\label{Jgof_nat_botglu}
\begin{tikzpicture}[vcenter]
\draw[0cell=.9]
(0,0) node (a) {\gaB_v\big(x \sscs \bang{\fgr a_{p,\, (\psiinv_p s_{p\crdot} \,\compq\, \psiinv_{p,q}\, s_{p,q,\ell})}}_{\ell \in \ufs{v}} \big)}
(a)++(1.5,-1.5) node (b) {\fgr \gaAp_v\big(x \sscs \bang{a_{p,\, (\psiinv_p s_{p\crdot} \,\compq\, \psiinv_{p,q}\, s_{p,q,\ell})}}_{\ell \in \ufs{v}} \big)} 
(a)++(4.5,0) node (c) {\fgr a_{p,\, \psiinv_p s_{p\crdot}}}
(a)++(-1.5,0) node (a') {\phantom{\fgr_{\hinv}}}
;
\draw[1cell=.85]
(a') [rounded corners=2pt, shorten <=-0ex] |- node[pos=.2,swap] {\laxf_{v,p}} (b)
;
\draw[1cell=.85]
(b) [rounded corners=2pt, shorten <=-0ex] -| node[pos=.8] {\fgr \glu^{a_p}_{x;\, \psiinv_p s_{p\crdot} \scs q,\, \ang{\psiinv_{p,q}\, s_{p,q,\ell}}_{\ell \in \ufs{v}}}} (c)
;
\end{tikzpicture}
\end{equation}
\item By \cref{glupsitil}, the clockwise composite in the bottom half of \cref{Jgof_nat_factor} has a gluing morphism at the component \cref{xangsupql} given by the gluing morphism of $(\Jgo f)_{h\angordmdot} (a,\glu^a)$ at the following component.
\[\begin{split}
&\big(x \sscs (\oplus_i\, \psi_i)^{-1}\angs, u, \ang{\psiinv_{p,q}\, s_{p,q,\ell}}_{\ell \in \ufs{v}}\big)\\
&= \big(x \sscs \ang{\psiinv_i s_{i\crdot}}_{i \in \ufs{k}} \spc u, \ang{\psiinv_{p,q}\, s_{p,q,\ell}}_{\ell \in \ufs{v}}\big)
\end{split}\]
By \cref{Jgof_m_gluing}, this gluing morphism of $(\Jgo f)_{h\angordmdot} (a,\glu^a)$ is also given by the composite in \cref{Jgof_nat_botglu}.  
\end{itemize}
This proves that the bottom half of \cref{Jgof_nat_factor} commutes.
\end{proof}

\begin{lemma}\label{Jgo_k_Geq}
The object assignment in \cref{Jgo_pos_obj_assignment},
\[\MultpsO(\ang{\A_i}_{i \in \ufs{k}} \sscs \B ) \fto{\Jgo} 
\GGCat(\ang{\Aidash}_{i \in \ufs{k}} \sscs \Bdash),\]
is $G$-equivariant.
\end{lemma}

\begin{proof}
\cref{Jgof_natural} shows that the object assignment $f \mapsto \Jgo f$ in \cref{Jgo_pos_obj_assignment}, for $k$-lax $\Op$-morphisms $(f,\laxf) \cn \ang{\A_i}_{i \in \ufs{k}} \to \B$, is well defined.  Recall that the $G$-actions on the domain of $\Jgo$, the codomain of $\Jgo$, and $\angordn$-systems are defined in, respectively, \cref{klax_Omorphism_g,ggcat_k_g,nsystem_gaction}.  The $G$-equivariance of the object assignment $f \mapsto \Jgo f$ means that, for each $k$-lax $\Op$-morphism
\[\bang{(\A_i,\gaAi,\phiAi)}_{i \in \ufs{k}} \fto{(f, \laxf)} (\B,\gaB,\phiB),\]
each element $g \in G$, and each object $\angordmdot \in \Gsk^k$ \cref{angordmdot_bm}, the following two $\angordmdot$-component pointed functors \cref{Jgof_component} are equal.
\begin{equation}\label{Jgo_k1_Geq}
\begin{tikzpicture}[baseline={(a1.base)}]
\draw[0cell]
(0,0) node (a1) {\txsma_{i \in \ufs{k}}\, \Aiangordmi}
(a1)++(4,0) node (a2) {\Bboldm}
;
\draw[1cell=.9]
(a1) edge[transform canvas={yshift=.5ex}] node {(g \cdot \Jgo f)_{\angordmdot}} (a2)
(a1) edge[transform canvas={yshift=-.5ex}] node[swap] {(\Jgo (g \cdot f))_{\angordmdot}} (a2) 
;
\end{tikzpicture}
\end{equation}

\parhead{Component objects}.  To check that the two pointed functors in \cref{Jgo_k1_Geq} are equal on objects, we consider an object \cref{aglua_sma}
\[(a,\glu^a) = \ang{(a_i, \glu^{a_i})}_{i \in \ufs{k}} \in \txsma_{i \in \ufs{k}}\, \Aiangordmi.\]
Since each functor in \cref{Jgo_k1_Geq} is a pointed functor, we may assume that $(a_i,\glu^{a_i})$ is not the base $\angordmi$-system in $\A_i$ for any $i \in \ufs{k}$.  For markers $\angs$ and $s_{i\crdot}$ as defined in \cref{marker_sij}, the following computation proves that the two pointed functors in \cref{Jgo_k1_Geq} yield the same $\angs$-component objects \cref{a_angs}.
\begin{equation}\label{Jgo_k1_Geq_obj}
\begin{aligned}
& \big( (g \cdot \Jgo f)_{\angordmdot} (a,\glu) \big)_{\angs} && \\
&= \big(g \circ (\Jgo f)_{\angordmdot} (\txsma_{i \in \ufs{k}}\, \ginv) (a,\glu) \big)_{\angs} && \text{by \cref{ggcat_k_g}} \\
&= \big(g \circ (\Jgo f)_{\angordmdot} \ang{(\ginv a_i, \ginv \glu^{a_i})}_{i \in \ufs{k}} \big)_{\angs} && \text{by \cref{nsystem_gaction}} \\
&= g \big( (\Jgo f)_{\angordmdot} \ang{(\ginv a_i, \ginv \glu^{a_i})}_{i \in \ufs{k}} \big)_{\angs} && \text{by \cref{ga_scomponent}} \\
&= gf\ang{(\ginv a_i)_{s_{i\crdot}}}_{i \in \ufs{k}} && \text{by \cref{Jgof_m_obj_comp}} \\ 
&= gf\ang{\ginv a_{i,\, s_{i\crdot}}}_{i \in \ufs{k}} && \text{by \cref{ga_scomponent}} \\
&= (g \cdot f) \ang{a_{i,\, s_{i\crdot}}}_{i \in \ufs{k}} && \text{by \cref{klax_mor_g}} \\
&= \big((\Jgo (g \cdot f))_{\angordmdot} (a,\glu) \big)_{\angs} && \text{by \cref{Jgof_m_obj_comp}}
\end{aligned}
\end{equation}
Moreover, the computation \cref{Jgo_k1_Geq_obj} also proves that the two pointed functors in \cref{Jgo_k1_Geq} are equal on morphisms in $\txsma_{i \in \ufs{k}}\, \Aiangordmi$, where we use \cref{gtheta_angs,Jgof_m_theta_comp} for morphisms instead of \cref{ga_scomponent,Jgof_m_obj_comp} for objects.

\parhead{Gluing}.  To show that the two pointed functors in \cref{Jgo_k1_Geq} yield the same gluing morphism, suppose we are given an object $x \in \Op(t)$ for some $t \geq 0$, a marker $\angs$ as defined in \cref{marker_sij}, a pair of indices $(p,q) \in \ufs{k} \times \ufs{r}_p$, and a partition
\[s_{p,q} = \coprod_{\ell \in \ufs{t}}\, s_{p,q,\ell} \subseteq \ufs{m}_{p,q}.\]
We consider
\begin{itemize}
\item the pointed functor
\[\fgr = f\big(\ang{\ginv a_{i,\, s_{i\crdot}}}_{i \in \ufs{k}}\, \compp - \big) \cn (\A_p,\zero) \to (\B,\zero),\]
\item the index $u = \big(\!\sum_{i=1}^{p-1} r_i\big) + q$,
\item the $(k+t)$-tuple of objects $\bolda$ in \cref{laxf_bolda}, and
\item its image under the diagonal $\ginv$-action $\ginv\bolda$, as defined in \cref{strength_dom_ginv}.
\end{itemize}
The following computation proves that the two pointed functors in \cref{Jgo_k1_Geq} yield the same gluing morphism \cref{gluing-morphism} at the tuple
\[(x; \angs, u, \ang{s_{p,q,\ell}}_{\ell \in \ufs{t}}),\]
where $\ang{\cdots}_i = \ang{\cdots}_{i \in \ufs{k}}$ and $\ang{\cdots}_\ell = \ang{\cdots}_{\ell \in \ufs{t}}$.
\[\scalebox{.9}{$\begin{aligned}
& \big( (g \cdot \Jgo f)_{\angordmdot} (a,\glu) \big)_{x; \angs \csp u, \ang{s_{p,q,\ell}}_{\ell}} && \\
&= \big(g \circ (\Jgo f)_{\angordmdot} \ang{(\ginv a_i, \ginv \glu^{a_i})}_{i} \big)_{x; \angs \csp u, \ang{s_{p,q,\ell}}_{\ell}} && \text{by \cref{ggcat_k_g,nsystem_gaction}} \\
&= g\big((\Jgo f)_{\angordmdot} \ang{(\ginv a_i, \ginv \glu^{a_i})}_{i} \big)_{\ginv x; \angs \csp u, \ang{s_{p,q,\ell}}_{\ell}} && \text{by \cref{ga_gluing}} \\
&= g\big(\fgr (\ginv \glu^{a_p})_{\ginv x;\, s_{p\crdot} \csp q, \ang{s_{p,q,\ell}}_\ell} \circ \laxf_{t,p;\, \ginv\bolda} \big) && \text{by \cref{Jgof_m_gluing}} \\
&= g\big(\fgr \ginv \glu^{a_p}_{x;\, s_{p\crdot} \csp q, \ang{s_{p,q,\ell}}_\ell} \circ \laxf_{t,p;\, \ginv\bolda} \big) && \text{by \cref{ga_gluing}} \\
&= \big(g \fgr \ginv \glu^{a_p}_{x;\, s_{p\crdot} \csp q, \ang{s_{p,q,\ell}}_\ell} \big) \circ \big(g\laxf_{t,p;\, \ginv\bolda} \big) && \text{by functoriality} \\
&= \big(g \fgr \ginv \glu^{a_p}_{x;\, s_{p\crdot} \csp q, \ang{s_{p,q,\ell}}_\ell} \big) \circ \laxgf_{t,p;\, \bolda} && \text{by \cref{laxgf_z}}\\
&= \big((\Jgo (g \cdot f))_{\angordmdot} (a,\glu) \big)_{x; \angs \csp u, \ang{s_{p,q,\ell}}_{\ell}} && \text{by \cref{klax_mor_g,Jgof_m_gluing}} 
\end{aligned}$}
\]
This proves that the two pointed functors in \cref{Jgo_k1_Geq} are equal on objects.
\end{proof}

\section{Multimorphism $G$-Functors in Positive Arity: Morphisms}
\label{sec:jemg_pos_ii}

Under \cref{as:OpA,as:BAi}, \cref{sec:jemg_pos_i,sec:jemg_pos_i_proof} construct the $G$-equivariant object assignment of the $k$-ary multimorphism $G$-functor
\[\MultpsO(\ang{\A_i}_{i \in \ufs{k}} \sscs \B ) \fto{\Jgo} 
\GGCat(\ang{\Aidash}_{i \in \ufs{k}} \sscs \Bdash).\]
This section finishes the construction of this $G$-functor by defining its $G$-equivariant morphism assignment.  There is also a strong variant
\[\MultpspsO(\ang{\A_i}_{i \in \ufs{k}} \sscs \B ) \fto{\Jgosg} 
\GGCat(\ang{\Aisgdash}_{i \in \ufs{k}} \sscs \Bsgdash)\]
that involves (i) $k$-ary $\Op$-pseudomorphisms in the domain and (ii) strong systems in the codomain.

\secoutline
\begin{itemize}
\item Using the object assignments in \cref{def:Jgo_pos_obj}, \cref{def:Jgo_pos_mor} constructs the morphism assignments of $\Jgo$ and $\Jgosg$.
\item Several statements are used in \cref{def:Jgo_pos_mor} to make sure that the morphism assignment of $\Jgo$ is well defined.  These statements are proved in  \cref{Jgotheta_system_mor,Jgotheta_m_natural,Jgotheta_modification}.
\end{itemize}

\recollection
For \cref{def:Jgo_pos_mor} below, we briefly recall the following.

\begin{itemize}
\item The domain of $\Jgo$ is the $G$-category (\cref{def:MultpsO_karycat})
\[\MultpsO(\ang{\A_i}_{i \in \ufs{k}} \sscs \B)\]
with
\begin{itemize}
\item $k$-lax $\Op$-morphisms (\cref{def:k_laxmorphism}) as objects,
\item $k$-ary $\Op$-transformations (\cref{def:kary_transformation}) as morphisms, and
\item $G$ acting by conjugation \pcref{def:k_laxmorphism_g,def:kary_transformation}.
\end{itemize}
The domain of $\Jgosg$ is the full sub-$G$-category whose objects are $k$-ary $\Op$-pseudomorphisms, which have invertible action constraints.
\item Using the $\Gskg$-categories $\Aidash$ and $\Bdash$ in \cref{A_ptfunctor}, the codomain of $\Jgo$ is the $G$-category 
\[\GGCat(\ang{\Aidash}_{i \in \ufs{k}} \sscs \Bdash)\]
described in \cref{expl:ggcat_positive}, in particular \cref{ggcat_k2,ggcat_k2_comp,ggcat_k2_modax,ggcat_k_g} for morphisms.
The codomain of $\Jgosg$ is defined in the same way using the $\Gskg$-categories $\Aisgdash$ and $\Bsgdash$ \pcref{A_ptfunctor}.
\end{itemize}

The following definition uses the notation in \cref{def:Jgo_pos_obj}, in particular \cref{angordmdot_bm,Jgof_component,ai_gluai,aglua_sma,Jgof_m_basesystem,marker_sij,Jgof_m_objects,Jgof_m_obj_comp}.

\begin{definition}\label{def:Jgo_pos_mor}
Under \cref{as:OpA,as:BAi}, using the $G$-equivariant object assignment given in \cref{def:Jgo_pos_obj}, the $G$-functor
\begin{equation}\label{Jgo_kary_functor}
\MultpsO(\ang{\A_i}_{i \in \ufs{k}} \sscs \B ) \fto{\Jgo} 
\GGCat(\ang{\Aidash}_{i \in \ufs{k}} \sscs \Bdash)
\end{equation}
is defined as follows.  A strong variant $\Jgosg$ is defined in \cref{Jgosg_pos_mor_assignment}.  Given a $k$-ary $\Op$-transformation 
\begin{equation}\label{theta_fd}
\begin{tikzpicture}[baseline={(a.base)}]
\def\t{27}
\draw[0cell=1]
(0,0) node (a) {\phantom{Z}}
(a)++(2,0) node (b) {\phantom{Z}}
(a)++(-1.25,0) node (a') {\ang{(\A_i,\gaAi,\phiAi)}_{i \in \ufs{k}}}
(b)++(.7,0) node (b') {(\B,\gaB,\phiB)}
;
\draw[1cell=.85]  
(a) edge[bend left=\t] node {(f,\laxf)} (b)
(a) edge[bend right=\t] node[swap] {(d,\laxd)} (b)
;
\draw[2cell]
node[between=a and b at .42, rotate=-90, 2label={above,\theta}] {\Rightarrow}
;
\end{tikzpicture}
\end{equation}
between $k$-lax $\Op$-morphisms \cref{kary_Otransformation}, the morphism
\begin{equation}\label{Jgotheta}
\begin{tikzpicture}[baseline={(a.base)}]
\def\t{25}
\draw[0cell=1]
(0,0) node (a) {\phantom{Z}}
(a)++(2.3,0) node (b) {\phantom{Z}}
(a)++(-.7,0) node (a') {\ang{\Aidash}_{i \in \ufs{k}}}
(b)++(.3,0) node (b') {\Bdash,}
;
\draw[1cell=.85]  
(a) edge[bend left=\t] node {\Jgo f} (b)
(a) edge[bend right=\t] node[swap] {\Jgo d} (b)
;
\draw[2cell]
node[between=a and b at .37, rotate=-90, 2label={above,\Jgo\theta}] {\Rightarrow}
;
\end{tikzpicture}
\end{equation}
in the sense of \cref{ggcat_k2}, has its $\angordmdot$-component pointed natural transformation 
\begin{equation}\label{Jgotheta_m}
\begin{tikzpicture}[baseline={(a.base)}]
\def\t{25}
\draw[0cell=1]
(0,0) node (a) {\phantom{Z}}
(a)++(3,0) node (b) {\phantom{Z}}
(a)++(-.75,0) node (a') {\txsma_{i \in \ufs{k}}\, \Aiangordmi}
(b)++(.17,0) node (b') {\Bboldm,}
;
\draw[1cell=.8]  
(a) edge[bend left=\t] node {(\Jgo f)_{\angordmdot}} (b)
(a) edge[bend right=\t] node[swap] {(\Jgo d)_{\angordmdot}} (b)
;
\draw[2cell=.9]
node[between=a and b at .3, rotate=-90, 2label={above,(\Jgo\theta)_{\angordmdot}}] {\Rightarrow}
;
\end{tikzpicture}
\end{equation}
in the sense of \cref{ggcat_k2_comp}, defined as follows.
\begin{description}
\item[Base case]  If $\angordmi = \vstar$ for some $i \in \ufs{k}$, then both $(\Jgo f)_{\angordmdot}$ and $(\Jgo d)_{\angordmdot}$ are the unique isomorphism between terminal categories.  Thus, $(\Jgo \theta)_{\angordmdot}$ must be the identity natural transformation on $1_{\boldone}$.  In the rest of this definition, we assume $\angordmi \in \Gsk \setminus \{\vstar\}$ for each $i \in \ufs{k}$. 
\item[Non-base case]  
Given an object \cref{aglua_sma}
\[(a,\glu^a) = \ang{(a_i, \glu^{a_i})}_{i \in \ufs{k}} \in \txsma_{i \in \ufs{k}}\, \Aiangordmi,\]
we denote the $(a,\glu^a)$-component morphism of $(\Jgo\theta)_{\angordmdot}$ in $\Bboldm$ by
\begin{equation}\label{Jgotheta_m_aglua}
(\Jgo f)_{\angordmdot} (a,\glu^a) \fto{(\Jgo\theta)_{\angordmdot, (a,\glu^a)}} 
(\Jgo d)_{\angordmdot} (a,\glu^a).
\end{equation}

\parhead{Base systems}.  Since we want $(\Jgo \theta)_{\angordmdot}$ to be a pointed natural transformation, if any $(a_i,\glu^{a_i})$ is the base $\angordmi$-system $(\zero, 1_\zero)$ in $\A_i$ (\cref{expl:nsystem_base}), then we must define
\begin{equation}\label{Jgotheta_m_basesystem}
(\Jgo \theta)_{\angordmdot, (a,\glu^a)} = 1_{(\zero,1_\zero)} \inspace \Bboldm,
\end{equation}
which is the identity morphism of the base $\bm$-system in $\B$.  This is well defined by \cref{Jgof_m_basesystem}.

\parhead{Non-base systems}.  Suppose $(a_i,\glu^{a_i})$ is not the base $\angordmi$-system for any $i \in \ufs{k}$.  The $\angs$-component morphism \cref{theta_angs} of $(\Jgo\theta)_{\angordmdot, (a,\glu^a)}$ is defined to be the $\ang{a_{i,\, s_{i\crdot}}}_{i \in \ufs{k}}$-component of $\theta$ \cref{ktr_natural}, as displayed below, where the two equalities are from \cref{Jgof_m_obj_comp}.
\begin{equation}\label{Jgotheta_m_angs}
\begin{tikzpicture}[vcenter]
\def\u{-1.2}
\draw[0cell=.9]
(0,0) node (a1) {\big((\Jgo f)_{\angordmdot} (a,\glu^a)\big)_{\angs}}
(a1)++(5.7,0) node (a2) {\big((\Jgo d)_{\angordmdot} (a,\glu^a)\big)_{\angs}}
(a1)++(0,\u) node (b1) {f\ang{a_{i,\, s_{i\crdot}}}_{i \in \ufs{k}}}
(a2)++(0,\u) node (b2) {d\ang{a_{i,\, s_{i\crdot}}}_{i \in \ufs{k}}}
;
\draw[1cell=.9]
(a1) edge[equal] (b1)
(a2) edge[equal] (b2)
(a1) edge node {(\Jgo\theta)_{\angordmdot, (a,\glu^a), \angs}} (a2)
(b1) edge node {\theta_{\ang{a_{i,\, s_{i\crdot}}}_{i \in \ufs{k}}}} (b2)
;
\end{tikzpicture}
\end{equation}
\item[Well defined]  \cref{Jgotheta_system_mor} proves that $(\Jgo\theta)_{\angordmdot, (a,\glu^a)}$ in \cref{Jgotheta_m_aglua} is a morphism of $\bm$-systems in $\B$.  \cref{Jgotheta_m_natural} proves that $(\Jgo\theta)_{\angordmdot}$ in \cref{Jgotheta_m} is a pointed natural transformation.  \cref{Jgotheta_modification} proves that $\Jgo\theta$ in \cref{Jgotheta} is a well-defined morphism.  
\item[Functoriality]  
The assignment defined in \cref{Jgotheta},
\[\theta \mapsto \Jgo\theta,\]
preserves identity morphisms and composition because these two notions are defined componentwise in 
\begin{itemize}
\item the domain $\MultpsO(\ang{\A_i}_{i \in \ufs{k}} \sscs \B)$ of $\Jgo$ \pcref{def:MultpsO_karycat}, 
\item the codomain $\GGCat(\ang{\Aidash}_{i \in \ufs{k}} \sscs \Bdash)$ of $\Jgo$ \pcref{expl:ggcat_positive}, and 
\item the category of $\angordn$-systems (\cref{def:nsystem_morphism}).
\end{itemize}
\item[$G$-equivariance]
The $G$-equivariance of the assignment $\theta \mapsto \Jgo\theta$ is proved by
\begin{itemize}
\item reusing the computation \cref{Jgo_k1_Geq_obj},
\item replacing the $k$-lax $\Op$-morphism $f$ by the $k$-ary $\Op$-transformation $\theta$, and
\item using \cref{gtheta_angs,Jgotheta_m_angs,O_tr_g}.
\end{itemize}
Thus,
\begin{itemize}
\item the $G$-equivariant object assignment $f \mapsto \Jgo f$ in \cref{def:Jgo_pos_obj} and
\item the $G$-equivariant morphism assignment $\theta \mapsto \Jgo \theta$ in the current definition
\end{itemize}
together define the $G$-functor \cref{Jgo_kary_functor}
\[\MultpsO(\ang{\A_i}_{i \in \ufs{k}} \sscs \B ) \fto{\Jgo} 
\GGCat(\ang{\Aidash}_{i \in \ufs{k}} \sscs \Bdash).\]
\end{description}

\parhead{Strong variant}.  The $G$-functor
\begin{equation}\label{Jgosg_pos_mor_assignment}
\MultpspsO(\ang{\A_i}_{i \in \ufs{k}} \sscs \B) \fto{\Jgosg} 
\GGCat(\ang{\Aisgdash}_{i \in \ufs{k}} \sscs \Bsgdash)
\end{equation}
is defined by
\begin{itemize}
\item the $G$-equivariant object assignment in \cref{Jgosg_pos_obj_assignment} and 
\item the $G$-equivariant morphism assignment in the current definition.
\end{itemize} 
This yields a $G$-functor $\Jgosg$ because $\Bsgboldm$ is a full subcategory of $\Bboldm$ by \cref{def:nsystem_morphism}.
\end{definition}

\subsection*{Proofs}
The rest of this section proves \cref{Jgotheta_system_mor,Jgotheta_m_natural,Jgotheta_modification}, which are used in \cref{def:Jgo_pos_mor}.  

\begin{lemma}\label{Jgotheta_system_mor}
The data in \cref{Jgotheta_m_aglua},
\[(\Jgo f)_{\angordmdot} (a,\glu^a) \fto{(\Jgo\theta)_{\angordmdot, (a,\glu^a)}}
(\Jgo d)_{\angordmdot} (a,\glu^a),\]
form a morphism of $\bm$-systems in $\B$.
\end{lemma}

\begin{proof}
We need to show that the $\angs$-component morphisms $(\Jgo\theta)_{\angordmdot, (a,\glu^a), \angs}$ defined in \cref{Jgotheta_m_angs} satisfy the two axioms in \cref{def:nsystem_morphism} for a morphism of $\bm$-systems in $\B$.

\parhead{Unity}.  To prove the unity axiom \cref{nsystem_mor_unity}, suppose $s_{p,q} = \emptyset$ for some $(p,q) \in \ufs{k} \times \ufs{r}_p$.  The object unity axiom \cref{system_obj_unity} for the $\angordmpe$-system $(a_p,\glu^{a_p})$ in $\A_p$ implies
\[a_{p,\, s_{p\crdot}} = \zero \in \A_p.\]
The basepoint condition \cref{ktransform_basepoint} for $\theta$ implies
\[(\Jgo\theta)_{\angordmdot, (a,\glu^a), \angs} = 
\theta_{\ang{a_{i,\, s_{i\crdot}}}_{i \in \ufs{k}}} = 1_\zero \inspace \B.\]

\parhead{Compatibility}.  To prove the compatibility axiom \cref{nsystem_mor_compat}, suppose we are given an object $x \in \Op(t)$ for some $t \geq 0$, a pair of indices $(p,q) \in \ufs{k} \times \ufs{r}_p$, and a partition
\[s_{p,q} = \coprod_{\ell \in \ufs{t}}\, s_{p,q,\ell} \subseteq \ufs{m}^{p,q}.\]
Using the notation in \cref{ftil_ApB,Jgof_m_gluing}, the compatibility diagram \cref{nsystem_mor_compat} for $(\Jgo\theta)_{\angordmdot, (a,\glu^a)}$ is the boundary diagram below.
\[\begin{tikzpicture}[vcenter]
\def\u{-1.2} \def\v{-2} \def\h{3.5} \def\a{10}
\draw[0cell=.8]
(0,0) node (a1) {\gaB_t\big(x \sscs \ang{\fgr a_{p,\, (s_{p\crdot} \,\compq\, s_{p,q,\ell})}}_{\ell \in \ufs{t}}\big)}
(a1)++(\h,\u) node (a2) {\fgr \gaAp_t \big(x \sscs \ang{a_{p,\, (s_{p\crdot} \,\compq\, s_{p,q,\ell})}}_{\ell \in \ufs{t}} \big)}
(a1)++(2*\h,0) node (a3) {\fgr a_{p,\, s_{p\crdot}}}
(a1)++(0,\v) node (b1)  {\gaB_t\big(x \sscs \ang{\dgr a_{p,\, (s_{p\crdot} \,\compq\, s_{p,q,\ell})}}_{\ell \in \ufs{t}}\big)}
(a2)++(0,\v) node (b2) {\dgr \gaAp_t \big(x \sscs \ang{a_{p,\, (s_{p\crdot} \,\compq\, s_{p,q,\ell})}}_{\ell \in \ufs{t}} \big)}
(a3)++(0,\v) node (b3) {\dgr a_{p,\, s_{p\crdot}}}
;
\draw[1cell=.8]
(a1) edge node[pos=.7] {\laxf_{t,p;\, \bolda}} (a2)
(a2) edge node[pos=.5] {\fgr \glu^{a_p}} (a3)
(b1) edge node[swap,pos=.3] {\laxd_{t,p;\, \bolda}} (b2)
(b2) edge node[swap,pos=.5] {\dgr \glu^{a_p}} (b3)
(a1) edge node[swap] {\gaB_t(1_x \sscs \ang{\theta}_{\ell \in \ufs{t}})} (b1)
(a2) edge node [swap] {\theta} (b2)
(a3) edge node {\theta} (b3)
;
\end{tikzpicture}\]
In the diagram above, each instance of $\glu^{a_p}$ is the gluing morphism
\[\glu^{a_p}_{x;\, s_{p\crdot} \scs q, \ang{s_{p,q,\ell}}_{\ell \in \ufs{t}}} 
\cn \gaAp_t \big(x \sscs \ang{a_{p,\, (s_{p\crdot} \,\compq\, s_{p,q,\ell})}}_{\ell \in \ufs{t}} \big) 
\to a_{p,\, s_{p\crdot}}.\]
The left quadrilateral commutes by the multilinearity axiom \cref{ktransform_multilinearity} for $\theta$.  The right quadrilateral commutes by the naturality of $\theta$ \cref{ktr_natural}.
\end{proof}

\begin{lemma}\label{Jgotheta_m_natural}
The data in \cref{Jgotheta_m},
\[\begin{tikzpicture}[baseline={(a.base)}]
\def\t{25}
\draw[0cell=1]
(0,0) node (a) {\phantom{Z}}
(a)++(3,0) node (b) {\phantom{Z}}
(a)++(-.75,0) node (a') {\txsma_{i \in \ufs{k}}\, \Aiangordmi}
(b)++(.17,0) node (b') {\Bboldm,}
;
\draw[1cell=.8]  
(a) edge[bend left=\t] node {(\Jgo f)_{\angordmdot}} (b)
(a) edge[bend right=\t] node[swap] {(\Jgo d)_{\angordmdot}} (b)
;
\draw[2cell=.9]
node[between=a and b at .3, rotate=-90, 2label={above,(\Jgo\theta)_{\angordmdot}}] {\Rightarrow}
;
\end{tikzpicture}\]
form a pointed natural transformation.
\end{lemma}

\begin{proof}
\cref{Jgotheta_system_mor} shows that each component of $(\Jgo\theta)_{\angordmdot}$ is a morphism of $\bm$-systems in $\B$.  Moreover, $(\Jgo\theta)_{\angordmdot}$ is pointed by \cref{Jgotheta_m_basesystem}.  The naturality of $(\Jgo\theta)_{\angordmdot}$ with respect to morphisms in $\txsma_{i \in \ufs{k}}\, \Aiangordmi$ follows from
\begin{itemize}
\item the definitions \cref{Jgof_m_theta_comp,Jgotheta_m_angs}, 
\item the fact that composition of morphisms of $\bm$-systems is defined componentwise (\cref{def:nsystem_morphism}), and 
\item the naturality of $\theta$ \cref{ktr_natural}.
\end{itemize}
This proves that $(\Jgo\theta)_{\angordmdot}$ is a pointed natural transformation.
\end{proof}

\begin{lemma}\label{Jgotheta_modification}
The data in \cref{Jgotheta},
\[\begin{tikzpicture}[baseline={(a.base)}]
\def\t{25}
\draw[0cell=1]
(0,0) node (a) {\phantom{Z}}
(a)++(2.3,0) node (b) {\phantom{Z}}
(a)++(-.7,0) node (a') {\ang{\Aidash}_{i \in \ufs{k}}}
(b)++(.3,0) node (b') {\Bdash,}
;
\draw[1cell=.85]  
(a) edge[bend left=\t] node {\Jgo f} (b)
(a) edge[bend right=\t] node[swap] {\Jgo d} (b)
;
\draw[2cell]
node[between=a and b at .37, rotate=-90, 2label={above,\Jgo\theta}] {\Rightarrow}
;
\end{tikzpicture}\]
form a morphism in $\GGCat(\ang{\Aidash}_{i \in \ufs{k}}; \Bdash)$.
\end{lemma}

\begin{proof}
By \cref{ggcat_k2_modax}, we need to check that, given morphisms \cref{upomdot}
\[\begin{aligned}
\angordmdot = \ang{\angordmi}_{i \in \ufs{k}} 
& \fto{\upomdot = \ang{\upom^i}_{i \in \ufs{k}}} 
\angordndot = \ang{\angordni}_{i \in \ufs{k}} && \text{in $\Gsk^k$ and}\\
\bm = \txoplus_{i \in \ufs{k}}\, \angordmi 
& \fto{\upom = \txoplus_{i \in \ufs{k}}\, \upom^i}
\bn = \txoplus_{i \in \ufs{k}}\, \angordni && \text{in $\Gsk$},
\end{aligned}\]
the following two whiskered pointed natural transformations are equal.
\begin{equation}\label{Jgotheta_modax_diagram}
\begin{tikzpicture}[vcenter]
\def\v{-1.7} \def\a{20} \def\b{.6} \def\c{.1}
\draw[0cell=.8]
(0,0) node (a1) {\phantom{X}}
(a1)++(3.5,0) node (a2) {\phantom{X}}
(a1)++(0,\v) node (b1) {\phantom{X}}
(a2)++(0,\v) node (b2) {\phantom{X}}
(a1)++(-\b,0) node (a1') {\txsma_{i \in \ufs{k}}\, \Aiangordmi}
(a2)++(\c,0) node (a2') {\Bboldm}
(b1)++(-\b,0) node (b1') {\txsma_{i \in \ufs{k}}\, \Aiangordni}
(b2)++(\c-.02,0) node (b2') {\Bboldn}
;
\draw[1cell=.8]
(a1) edge[bend left=\a] node[pos=.44] {(\Jgo f)_{\angordmdot}} (a2)
(a1) edge[bend right=\a] node[pos=.6,swap] {(\Jgo d)_{\angordmdot}} (a2)
(b1) edge[bend left=\a] node[pos=.44] {(\Jgo f)_{\angordndot}} (b2)
(b1) edge[bend right=\a] node[pos=.6,swap] {(\Jgo d)_{\angordndot}} (b2)
(a1') edge[transform canvas={shift={(\b,0)}}] node[swap] {\txsma_{i \in \ufs{k}}\, \Aiupomi} (b1')
(a2) edge[transform canvas={shift={(0,0)}}] node {\Bupom} (b2)
;
\draw[2cell=.9]
node[between=a1 and a2 at .35, rotate=-90, 2label={above,(\Jgo\theta)_{\angordmdot}}] {\Rightarrow}
node[between=b1 and b2 at .35, rotate=-90, 2label={above,(\Jgo\theta)_{\angordndot}}] {\Rightarrow}
;
\end{tikzpicture}
\end{equation}
The rest of this proof follows the same pattern as the proof of \cref{Jgof_natural}, which proves the commutativity of the diagram \cref{Jgof_nat_diagram}.  

In more detail, first suppose $\upom_i$ is the 0-morphism for some $i \in \ufs{k}$.  Since $(\Jgo\theta)_{\angordndot}$ is pointed by \cref{Jgotheta_m_basesystem}, each whiskered pointed natural transformation in \cref{Jgotheta_modax_diagram} is constant at the identity morphism $1_{(\zero,1_\zero)}$ in $\Bboldn$.  

Thus, we assume that each morphism \cref{Fihipsii}
\[\angordmi \fto{\upom^i = (h_i, \psi_i)} \angordni \inspace \Gsk\]
is not a 0-morphism.  Using the notation in \cref{Jgof_nat_factor}, the diagram \cref{Jgotheta_modax_diagram} is the boundary diagram below.
\begin{equation}\label{Jgotheta_modax_factor}
\begin{tikzpicture}[vcenter]
\def\v{-1.7} \def\a{20} \def\b{.6} \def\c{.1}
\draw[0cell=.8]
(0,0) node (a1) {\phantom{X}}
(a1)++(3.5,0) node (a2) {\phantom{X}}
(a1)++(0,\v) node (b1) {\phantom{X}}
(a2)++(0,\v) node (b2) {\phantom{X}}
(b1)++(0,\v) node (c1) {\phantom{X}}
(b2)++(0,\v) node (c2) {\phantom{X}}
(a1)++(-\b,0) node (a1') {\txsma_{i \in \ufs{k}}\, \Aiangordmi}
(a2)++(\c,0) node (a2') {\Bboldm}
(b1)++(-.95,0) node (b1') {\txsma_{i \in \ufs{k}}\, \Aihiangordmi}
(b2)++(.3,0) node (b2') {\Bhboldm}
(c1)++(-\b,0) node (c1') {\txsma_{i \in \ufs{k}}\, \Aiangordni}
(c2)++(\c-.02,0) node (c2') {\Bboldn}
;
\draw[1cell=.8]
(a1) edge[bend left=\a] node[pos=.44] {(\Jgo f)_{\angordmdot}} (a2)
(a1) edge[bend right=\a] node[pos=.6,swap] {(\Jgo d)_{\angordmdot}} (a2)
(b1) edge[bend left=\a] node[pos=.44] {(\Jgo f)_{h \angordmdot}} (b2)
(b1) edge[bend right=\a] node[pos=.6,swap] {(\Jgo d)_{h \angordmdot}} (b2)
(c1) edge[bend left=\a] node[pos=.44] {(\Jgo f)_{\angordndot}} (c2)
(c1) edge[bend right=\a] node[pos=.6,swap] {(\Jgo d)_{\angordndot}} (c2)
(a1) edge[shorten <=.5ex, shorten >=1ex] node[swap] {\txsma_{i \in \ufs{k}}\, \htil_i} (b1)
(b1) edge[shorten <=.5ex, shorten >=1ex] node[swap] {\txsma_{i \in \ufs{k}}\, \psitil_i} (c1)
(a2) edge[shorten >=.5ex] node {\wt{\oplus_i\, h_i}} (b2)
(b2) edge[shorten >=.5ex] node {\wt{\oplus_i\, \psi_i}} (c2)
;
\draw[2cell=.9]
node[between=a1 and a2 at .35, rotate=-90, 2label={above,(\Jgo\theta)_{\angordmdot}}] {\Rightarrow}
node[between=b1 and b2 at .35, rotate=-90, 2label={above,(\Jgo\theta)_{h \angordmdot}}] {\Rightarrow}
node[between=c1 and c2 at .35, rotate=-90, 2label={above,(\Jgo\theta)_{\angordndot}}] {\Rightarrow}
;
\end{tikzpicture}
\end{equation}
To prove that the two whiskered pointed natural transformations in the top half of \cref{Jgotheta_modax_factor} are equal, we reuse the computation in \cref{Jgof_nat_computation} by
\begin{itemize}
\item replacing $f$ with $\theta$ and
\item using \cref{thatil_component,Jgotheta_m_angs}.
\end{itemize}
To prove that the two whiskered pointed natural transformations in the bottom half of \cref{Jgotheta_modax_factor} are equal, we reuse the computation in \cref{Jgof_natbot_computation} by replacing $f$ with $\theta$ and using \cref{thapsitil_component}.
\end{proof}

\section{$J$-Theory $\Gcat$-Multifunctors}
\label{sec:jemg_axioms}

For a $\Tinf$-operad $\Op$ \pcref{as:OpA}, using the object assignment and multimorphism $G$-functors defined in earlier sections, this section constructs the $J$-theory $\Gcat$-multifunctor
\[\MultpsO \fto{\Jgo} \GGCat\]
and the strong variant
\[\MultpspsO \fto{\Jgosg} \GGCat.\]
As a reminder, our enriched multifunctors (\cref{def:enr-multicategory-functor}) are always assumed to preserve the symmetric group action of its domain and codomain.  The $\Gcat$-multifunctor $\Jgo$ is the first step of our $G$-equivariant algebraic $K$-theory multifunctor, from $\Op$-pseudoalgebras to orthogonal $G$-spectra via $\Gskg$-categories.  

\secoutline
\begin{itemize}
\item \cref{def:Jgo_multifunctor} defines the $\Gcat$-multifunctors $\Jgo$ and $\Jgosg$.
\item \cref{Jgo_sigma} proves that $\Jgo$ and $\Jgosg$ preserve the symmetric group action.
\item \cref{Jgo_gamma} proves that $\Jgo$ and $\Jgosg$ preserve the multicategorical composition.
\item \cref{thm:Jgo_multifunctor} proves that $\Jgo$ and $\Jgosg$ are $\Gcat$-multifunctors.
\item \cref{ex:JgBE} applies \cref{thm:Jgo_multifunctor} to the Barratt-Eccles $\Gcat$-operad $\BE$ and the $G$-Barratt-Eccles operad $\GBE$.
\end{itemize}

\recollection
For \cref{def:Jgo_multifunctor} below, we briefly recall the following.
\begin{itemize}
\item $\MultpsO$ is the $\Gcat$-multicategory in \cref{thm:multpso}, with
\begin{itemize}
\item $\Op$-pseudoalgebras (\cref{def:pseudoalgebra}) as objects,
\item $k$-lax $\Op$-morphisms (\cref{def:k_laxmorphism}) as $k$-ary 1-cells,
\item $k$-ary $\Op$-transformations (\cref{def:kary_transformation}) as $k$-ary 2-cells, and
\item $G$ acting by conjugation on $k$-ary 1-cells and $k$-ary 2-cells \pcref{def:k_laxmorphism_g,expl:O_tr_g}.
\end{itemize}
For $\MultpspsO$, the $k$-ary 1-cells are $k$-ary $\Op$-pseudomorphisms, which have invertible action constraints.
\item $\GGCat$ is the $\Gcat$-multicategory in \cref{expl:ggcat_gcatenr}, with
\begin{itemize} 
\item $\Gskg$-categories \pcref{expl:ggcat_obj} as objects,
\item 0-ary multimorphism $G$-categories described in \cref{expl:ggcat_zero}, and
\item positive arity multimorphism $G$-categories described in \cref{expl:ggcat_positive}.
\end{itemize}
\end{itemize}

The next definition uses \cref{def:enr-multicategory-functor} for the Cartesian closed category $(\Gcat, \times, \boldone, \Catg)$ \pcref{expl:Gcat_closed}.

\begin{definition}[Multifunctorial $J$-Theory]\label{def:Jgo_multifunctor}
For a $\Tinf$-operad $\Op$ \pcref{as:OpA}, we define the $\Gcat$-multifunctors\index{multifunctorial J-theory@multifunctorial $J$-theory}
\[\begin{split}
\MultpsO & \fto{\Jgo} \GGCat \andspace\\
\MultpspsO & \fto{\Jgosg} \GGCat
\end{split}\]
as follows.
\begin{description}
\item[Objects]  $\Jgo$ and $\Jgosg$ send each $\Op$-pseudoalgebra $\A$ to, respectively, the $\Gskg$-categories
\[\begin{split}
(\Gsk,\vstar) & \fto{\Jgo\A = \Adash} (\Gcatst,\boldone) \andspace\\ 
(\Gsk,\vstar) & \fto{\Jgosg\A = \Asgdash} (\Gcatst,\boldone)
\end{split}\]
in \cref{A_ptfunctor}.  Their object and morphism assignments are given in, respectively, \cref{def:Aangordn_gcat,def:Afangpsi}.
\item[Multimorphisms in arity 0] For an $\Op$-pseudoalgebra $\A$, the 0-ary multimorphism $G$-functors 
\[\begin{split}
\MultpsO(\ang{} \sscs \A) = \A & \fto{\Jgo} \GGCat(\ang{} \sscs \Adash) \andspace \\
\MultpspsO(\ang{} \sscs \A) = \A & \fto{\Jgosg} \GGCat(\ang{} \sscs \Asgdash)
\end{split}\]
are defined in \cref{def:Jgo_zero}.
\item[Multimorphisms in positive arity]  Suppose $\B$ and $\ang{\A_i}_{i \in \ufs{k}}$ are $\Op$-pseudoalgebras with $k \geq 1$.  The $k$-ary multimorphism $G$-functors
\[\begin{split}
\MultpsO(\ang{\A_i}_{i \in \ufs{k}} \sscs \B) & \fto{\Jgo} 
\GGCat(\ang{\Aidash}_{i \in \ufs{k}} \sscs \Bdash) \andspace \\
\MultpspsO(\ang{\A_i}_{i \in \ufs{k}} \sscs \B) & \fto{\Jgosg} 
\GGCat(\ang{\Aisgdash}_{i \in \ufs{k}} \sscs \Bsgdash)
\end{split}\]
have object and morphism assignments given in, respectively, \cref{def:Jgo_pos_obj,def:Jgo_pos_mor}.
\end{description}
The $\Gcat$-multifunctor axioms for $\Jgo$ and $\Jgosg$ are verified in \cref{Jgo_sigma,Jgo_gamma,thm:Jgo_multifunctor}.  We refer to $\Jgo$ and $\Jgosg$ as, respectively, \emph{multifunctorial $J$-theory} and \emph{multifunctorial strong $J$-theory}\index{multifunctorial J-theory@multifunctorial $J$-theory!strong} for $\Op$.
\end{definition}

\subsection*{$\Gcat$-Multifunctor Axioms}

The rest of this section proves that each of $\Jgo$ and $\Jgosg$ satisfies the axioms in \cref{def:enr-multicategory-functor} for a $\Gcat$-multifunctor.  \cref{Jgo_sigma} below verifies the symmetric group action axiom.  The symmetric group action on $\MultpsO$ and $\GGCat$ are given in \cref{def:multpso_sym,def:ktransformation_sym,expl:ggcat_symmetry}.

\begin{lemma}\label{Jgo_sigma}
In \cref{def:Jgo_multifunctor}, each of $\Jgo$ and $\Jgosg$ satisfies the symmetric group action axiom \cref{enr-multifunctor-equivariance}.
\end{lemma}

\begin{proof}
We consider $\Jgo$.  The proof for the strong variant $\Jgosg$ is the same after restricting to $k$-ary $\Op$-pseudomorphisms in the domain and strong systems in the codomain.  The symmetric group action axiom \cref{enr-multifunctor-equivariance} for $\Jgo$ states that the following diagram of $G$-functors commutes for each permutation $\si \in \Si_k$ with $k \geq 2$.
\begin{equation}\label{Jgo_symmetry_diag}
\begin{tikzpicture}[vcenter]
\def\v{-1.5}
\draw[0cell=.85]
(0,0) node (a1) {\MultpsO(\ang{\A_i}_{i \in \ufs{k}} \sscs \B)}
(a1)++(5,0) node (a2) {\GGCat(\ang{\Aidash}_{i \in \ufs{k}} \sscs \Bdash)}
(a1)++(0,\v) node (b1) {\MultpsO(\ang{\A_{\si(i)}}_{i \in \ufs{k}} \sscs \B)}
(a2)++(0,\v) node (b2) {\GGCat(\ang{\Asiidash}_{i \in \ufs{k}} \sscs \Bdash)}
;
\draw[1cell=.9]
(a1) edge node {\Jgo} (a2)
(b1) edge node {\Jgo} (b2)
(a1) edge[transform canvas={xshift=2em}] node[swap] {\si} (b1)
(a2) edge[transform canvas={xshift=-3em}] node {\si} (b2)
;
\end{tikzpicture}
\end{equation}
It suffices to show that the diagram \cref{Jgo_symmetry_diag} commutes on
\begin{itemize}
\item objects, which are $k$-lax $\Op$-morphisms \pcref{def:k_laxmorphism}, and 
\item morphisms, which are $k$-ary $\Op$-transformations \pcref{def:kary_transformation}.
\end{itemize}

\parhead{Component objects}.  Using the notation in \cref{def:Jgo_pos_obj}, the right $\si$-action on $\angordmdot$, $\bm$, $(a,\glu^a)$, and $\angs$ are given as follows.
\begin{equation}\label{angordmdotsi}
\begin{split}
\angordmdot\si &= \ang{\angordmsii}_{i \in \ufs{k}} \in \Gsk^{k}\\
\bm\si &= \txoplus_{i \in \ufs{k}}\, \angordmsii \in \Gsk\\
(a,\glu^a)\si &= \ang{(a_{\si(i)}, \glu^{a_{\si(i)}})}_{i \in \ufs{k}} \in \txsma_{i \in \ufs{k}}\, \Asiiangordmsii\\
\angs\si &= \ang{s_{\si(i)\crdot}}_{i \in \ufs{k}} = \ang{\ang{s_{\si(i),j}}_{j \in \ufs{r}_{\si(i)}}}_{i \in \ufs{k}}
\end{split}
\end{equation}
We show that the two composites in \cref{Jgo_symmetry_diag} are equal when they apply to a $k$-lax $\Op$-morphism $(f,\laxf)$, the objects $\angordmdot\si$ and $(a,\glu^a)\si$, and the marker $\angs\si$.  By the object unity axiom \cref{system_obj_unity} for $\bm\si$-systems in $\B$, we may assume that all $s_{i,j} \neq \emptyset$.  

By \cref{ggcat_si_action}, the image of $(f,\laxf)$ under the clockwise composite $\si \circ \Jgo$ in \cref{Jgo_symmetry_diag} has $\angordmdot\si$-component pointed functor given by the following composite.
\begin{equation}\label{Jgofsi_malsi}
\begin{tikzpicture}[vcenter]
\def\v{-1.5}
\draw[0cell]
(0,0) node (a1) {\txsma_{i \in \ufs{k}}\, \Asiiangordmsii}
(a1)++(4.5,0) node (a2) {\Bboldmsi}
(a1)++(0,\v) node (b1) {\txsma_{i \in \ufs{k}}\, \Aiangordmi}
(a2)++(0,\v) node (b2) {\Bboldm}
;
\draw[1cell=.9]
(a1) edge node {(\Jgo f)^{\si}_{\angordmdot\si}} (a2)
(a1) edge[transform canvas={xshift=0em}] node[swap] {\si} (b1)
(b1) edge node {(\Jgo f)_{\angordmdot}} (b2)
(b2) edge node[swap] {\sys{\B}(\sigmainv)} (a2)
;
\end{tikzpicture}
\end{equation}
By \cref{sigmainv_Gsk_exp}, in the right vertical arrow in \cref{Jgofsi_malsi}, $\sigmainv$ is the following braiding isomorphism in $\Gsk$.
\[\bm = \txoplus_{i \in \ufs{k}}\, \angordmi \fto[\iso]{(\sigmainv\ang{r_i}_{i \in \ufs{k}} \spc \ang{1})} 
\txoplus_{i \in \ufs{k}}\, \angordmsii = \bm\si\]
It consists of
\begin{itemize}
\item the reindexing injection $\sigmainv\ang{r_i}_{i \in \ufs{k}}$, which is the block permutation induced by $\sigmainv \in \Si_k$ that permutes blocks of lengths $(r_1,\ldots,r_k)$, and
\item identity morphisms of $\ordm^{\si(i),j} \in \Fsk$ for $i \in \ufs{k}$ and $j \in \ufs{r}_{\si(i)}$.
\end{itemize}  
By \cref{def:ftil_functor,def:Afangpsi}, applying $\sys{\B}{(-)}$ to $\sigmainv$ yields the pointed $G$-functor 
\begin{equation}\label{Bsigmainv}
\sys{\B}(\sigmainv) = \wt{\sigmainv}.
\end{equation}
Using the notation in \cref{ftil_angs}, we have the marker
\begin{equation}\label{siangssi}
\wt{\sigmainv}_* (\angs\si) = \big(\sigmainv \ang{r_i}_{i \in \ufs{k}}\big)^{-1} (\angs\si) = \angs.
\end{equation}
The following computation shows that the two composites in \cref{Jgo_symmetry_diag} yield the same $\angs\si$-component object in $\B$.
\begin{equation}\label{Jgo_sym_computation}
\begin{aligned}
& \big((\Jgo f)^\si_{\angordmdot\si} ((a,\glu^a)\si)\big)_{\angs\si} &&\\
&= \big(\big(\sys{\B}(\sigmainv) \circ (\Jgo f)_{\angordmdot} \circ \si \big) \big((a, \glu^a)\si\big) \big)_{\angs\si} && \text{by \cref{Jgofsi_malsi}}\\
&= \big(\big(\wt{\sigmainv} \circ (\Jgo f)_{\angordmdot} \circ \si \big) \big((a, \glu^a)\si\big) \big)_{\angs\si} && \text{by \cref{Bsigmainv}}\\
&= \big(\wt{\sigmainv} (\Jgo f)_{\angordmdot} (a, \glu^a)\big)_{\angs\si} && \text{by \cref{angordmdotsi}}\\
&= \big((\Jgo f)_{\angordmdot} (a, \glu^a)\big)_{\wt{\sigmainv}_* (\angs\si)} && \text{by \cref{atil_component}}\\
&= \big((\Jgo f)_{\angordmdot} (a, \glu^a)\big)_{\angs} && \text{by \cref{siangssi}}\\
&= f\ang{a_{i,\, s_{i\crdot}}}_{i \in \ufs{k}} && \text{by \cref{Jgof_m_obj_comp}}\\
&= (f\si) \ang{a_{\si(i),\, s_{\si(i)\crdot}}}_{i \in \ufs{k}} && \text{by \cref{klax_sigma_functor}}\\
&= \big((\Jgo (f\si))_{\angordmdot\si}  \ang{(a_{\si(i)}, \glu^{a_{\si(i)}})}_{i \in \ufs{k}} \big)_{\ang{s_{\si(i)\crdot}}_{i \in \ufs{k}}} && \text{by \cref{Jgof_m_obj_comp}}\\
&= \big((\Jgo (f\si))_{\angordmdot\si} ((a,\glu^a)\si) \big)_{\angs\si} && \text{by \cref{angordmdotsi}}
\end{aligned}
\end{equation}

\parhead{Gluing}.  To show that the two composites in the diagram \cref{Jgo_symmetry_diag} yield the same gluing morphisms, we consider an object $x \in \Op(t)$ for some $t \geq 0$, a pair of indices $(p,q) \in \ufs{k} \times \ufs{r}_p$, a partition
\[s_{p,q} = \coprod_{\ell \in \ufs{t}}\, s_{p,q,\ell} \subseteq \ufs{m}^{p,q},\]
and the indices
\[\begin{split}
u &= \big(\txsum_{i=1}^{p-1} r_i\big) + q \andspace \\
u' &= \big(\txsum_{i=1}^{\sigmainv(p)-1} r_{\si(i)} \big) + q.
\end{split}\]
We observe the following.
\begin{itemize}
\item $u'$ is the image of $u$ under the block permutation $\sigmainv\ang{r_i}_{i \in \ufs{k}}$.
\item $s_{p,q}$ is the $u$-th entry in $\angs = \ang{\ang{s_{i,j}}_{j \in \ufs{r}_i}}_{i \in \ufs{k}}$.
\item $s_{p,q}$ is the $u'$-th entry in $\angs\si = \ang{\ang{s_{\si(i),j}}_{j \in \ufs{r}_{\si(i)}}}_{i \in \ufs{k}}$.
\end{itemize}  
By \cref{glutil_component}, \cref{siangssi}, and the first three equalities in \cref{Jgo_sym_computation}, the gluing morphism of the $\bm\si$-system
\[\begin{split}
& (\Jgo f)^\si_{\angordmdot\si} \big((a,\glu^a)\si\big) \\
&= \wt{\sigmainv} (\Jgo f)_{\angordmdot} (a, \glu^a) \in \Bboldmsi
\end{split}\]
at the component
\begin{equation}\label{xangssiuprime}
\big(x \sscs \angs\si, u', \ang{s_{p,q,\ell}}_{\ell \in \ufs{t}}\big)
\end{equation}
is given by the gluing morphism of the $\bm$-system $(\Jgo f)_{\angordmdot} (a,\glu^a)$ at the component
\[\begin{split}
&\big(x \sscs \wt{\sigmainv}_* (\angs\si), (\sigmainv \ang{r_i}_{i \in \ufs{k}})^{-1}(u'),\, \ang{s_{p,q,\ell}}_{\ell \in \ufs{t}} \big)\\
&= \big(x \sscs \angs, u, \ang{s_{p,q,\ell}}_{\ell \in \ufs{t}}\big).
\end{split}\]
This gluing morphism is given by the composite 
\begin{equation}\label{ftil_gluap_laxftp}
\fgr \glu^{a_p}_{x;\, s_{p\crdot} \scs q,\, \ang{s_{p,q,\ell}}_{\ell \in \ufs{t}}} \circ \laxf_{t,p;\, \bolda}
\end{equation}
defined in \cref{Jgof_m_gluing}.  

On the other hand, by \cref{laxfsi_whiskering,laxfsi_domain_obj,Jgof_m_gluing}, the gluing morphism of the $\bm\si$-system 
\[\begin{split}
& (\Jgo (f\si))_{\angordmdot\si} \big((a,\glu^a)\si\big) \\
&= (\Jgo (f\si))_{\angordmdot\si}  \ang{(a_{\si(i)}, \glu^{a_{\si(i)}})}_{i \in \ufs{k}} \in \Bboldmsi
\end{split}\]
at the component \cref{xangssiuprime} is given by the composite
\begin{equation}\label{fsi_boldasi_laxfsi}
\begin{split}
& \big\lbrace (f\si)\big(\ang{a_{\si(i),\, s_{\si(i)\crdot}}}_{i \in \ufs{k}} \,\comp_{\sigmainv(p)} \big(\glu^{a_p}_{x;\, s_{p\crdot} \scs q,\, \ang{s_{p,q,\ell}}_{\ell \in \ufs{t}}} \big)\big) \big\rbrace \circ \laxfsi_{t,\, \sigmainv(p);\, \bolda\si}\\
&= \big\lbrace f\big(\ang{a_{i,\, s_{i\crdot}}}_{i \in \ufs{k}} \,\compp \big(\glu^{a_p}_{x;\, s_{p\crdot} \scs q,\, \ang{s_{p,q,\ell}}_{\ell \in \ufs{t}}} \big)\big) \big\rbrace \circ \big(\laxf_{t,p} * \si\big)_{\bolda\si}\\
\end{split}
\end{equation}
where
\[\bolda\si =
\scalebox{.9}{$\big(\ang{a_{\si(i),\, s_{\si(i)\crdot}}}_{i=1}^{\sigmainv(p)-1}, x, \ang{a_{p,\, (s_{p\crdot} \,\compq\, s_{p,q,\ell})}}_{\ell \in \ufs{t}} \spc \ang{a_{\si(i),\, s_{\si(i)\crdot}}}_{i=\sigmainv(p)+1}^k \big)$}.\]
By \cref{laxfsi_si_z,laxfsi_z,ftil_ApB}, the composite morphisms in \cref{ftil_gluap_laxftp,fsi_boldasi_laxfsi} are equal.  This proves that the $\angordmdot\si$-component pointed functors 
\begin{equation}\label{Jgof_si_angordmdotsi}
\begin{tikzpicture}[baseline={(a1.base)}]
\draw[0cell]
(0,0) node (a1) {\txsma_{i \in \ufs{k}}\, \Asiiangordmsii}
(a1)++(4.5,0) node (a2) {\Bboldmsi}
;
\draw[1cell=.9]
(a1) edge[transform canvas={yshift=.5ex}] node {(\Jgo f)^{\si}_{\angordmdot\si}} (a2)
(a1) edge[transform canvas={yshift=-.5ex}] node[swap] {(\Jgo (f\si))_{\angordmdot\si}} (a2)
;
\end{tikzpicture}
\end{equation}
are equal on objects in $\txsma_{i \in \ufs{k}}\, \Asiiangordmsii$.  

To show that the two $\angordmdot\si$-component pointed functors in \cref{Jgof_si_angordmdotsi} are equal on morphisms in $\txsma_{i \in \ufs{k}}\, \Asiiangordmsii$, we reuse the computation in \cref{Jgo_sym_computation} with each $\angordmi$-system $(a_i,\glu^{a_i})$ in $\A_i$ replaced by a morphism in $\Aiangordmi$, along with \cref{thatil_component,Jgof_m_theta_comp}.  Thus, the diagram \cref{Jgo_symmetry_diag} commutes on objects.

\parhead{Morphisms}.  To show that the diagram \cref{Jgo_symmetry_diag} commutes on morphisms, we reuse the computation in \cref{Jgo_sym_computation} with the $k$-lax $\Op$-morphism $f$ replaced by a $k$-ary $\Op$-transformation (\cref{def:kary_transformation}), along with \cref{Otr_sigma,ggcat_si_action,thatil_component,Jgotheta_m_angs}.
\end{proof}

\cref{Jgo_gamma} below verifies that $\Jgo$ and $\Jgosg$ preserve the multicategorical composition on $\MultpsO$ and $\GGCat$, which are described in \cref{def:gam_functor,expl:ggcat_composition}.

\begin{lemma}\label{Jgo_gamma}
In \cref{def:Jgo_multifunctor}, each of $\Jgo$ and $\Jgosg$ satisfies the composition axiom \cref{v-multifunctor-composition}.
\end{lemma}

\begin{proof}
We consider $\Jgo$.  The proof for the strong variant $\Jgosg$ is the same after restricting to $k$-ary $\Op$-pseudomorphisms in the domain and strong systems in the codomain.

The composition axiom \cref{v-multifunctor-composition} for $\Jgo$ states that, for $k \geq 1$, $d_i \geq 0$, $\Op$-pseudoalgebras \pcref{def:pseudoalgebra}
\[\begin{split}
& (\C,\gaC,\phiC) \cq \Bdot = \ang{(\B_i,\gaBi,\phiBi)}_{i \in \ufs{k}} \spc\\ 
& \Aidot = \ang{(\Aij,\gaAij,\phiAij)}_{j \in \ufs{d}_i} \spc \andspace \Addot = \ang{\Aidot}_{i \in \ufs{k}},
\end{split}\]
and $\Gskg$-categories \pcref{A_ptfunctor}
\[\begin{split}
\Bdotdash &= \ang{\Bidash}_{i \in \ufs{k}} \spc\\ 
\Aidotdash &= \ang{\Aijdash}_{j \in \ufs{d}_i} \spc\andspace\\
\Addotdash &= \ang{\Aidotdash}_{i \in \ufs{k}} \, ,
\end{split}\]
the following diagram of $G$-functors commutes.
\begin{equation}\label{Jgo_comp_diagram}
\begin{tikzpicture}[vcenter]
\def\v{-2} \def\h{4}
\draw[0cell=.8]
(0,0) node (a1) {\MultpsO(\Bdot \sscs \C) \times \txprod_{i \in \ufs{k}}\, \MultpsO(\Aidot \sscs \B_i)}
(a1)++(\h,-.5*\v) node(a2) {\MultpsO(\Addot \sscs \C)}
(a1)++(0,\v) node (b1) {\GGCat(\Bdotdash \sscs \Cdash) \times \txprod_{i \in \ufs{k}}\, \GGCat(\Aidotdash \sscs \Bidash)}
(a2)++(0,\v) node (b2) {\GGCat(\Addotdash \sscs \Cdash)}  
;
\draw[1cell=.8]
(a1) edge[transform canvas={xshift=0em}, shorten <=0em] node[pos=.4] {\gam} (a2)
(b1) edge[transform canvas={xshift=0em}, shorten <=0em] node[pos=.4] {\ga} (b2)
(a2) edge[transform canvas={xshift=-2em}, shorten >=0ex] node {\Jgo} (b2)
(a1) edge[transform canvas={xshift=1em}] node[swap] {\Jgo \times \txprod_{i \in \ufs{k}}\, \Jgo} (b1)
;
\end{tikzpicture}
\end{equation}
It suffices to show that the diagram \cref{Jgo_comp_diagram} commutes on objects and morphisms \pcref{def:k_laxmorphism,def:kary_transformation}.

\parhead{Component objects}.  To show that \cref{Jgo_comp_diagram} commutes on objects, suppose we are given a $k$-lax $\Op$-morphism (\cref{def:k_laxmorphism})
\[\Bdot \fto{(f,\laxf)} \C,\]
a $d_i$-lax $\Op$-morphism
\[\Aidot \fto{(h_i, \laxhi)} \B_i \foreachspace i \in \ufs{k},\]
and the following objects \pcref{def:Gsk}, markers \cref{marker}, and systems \pcref{def:nsystem,def:nsystem_morphism}, where $d = \sum_{i \in \ufs{k}}\, d_i$ and $r_{i,j} \geq 0$.
\[\begin{gathered}
\begin{aligned}
\angordmijdot &= \ang{\ordmijl \in \Fsk}_{\ell \in \ufs{r}_{i,j}} \in \Gsk \phantom{M}
& s_{i,j\crdot} &= \ang{s_{i,j,\ell} \subseteq \ufs{m}^{i,j,\ell}}_{\ell \in \ufs{r}_{i,j}} \\
\angordmiddot &= \ang{\angordmijdot}_{j \in \ufs{d}_i} \in \Gsk^{d_i} 
& s_{i\crdots} &= \ang{s_{i,j\crdot}}_{j \in \ufs{d}_i} \\
\angordmdddot &= \ang{\angordmiddot}_{i \in \ufs{k}} \in \Gsk^d 
& \angs &= \ang{s_{i\crdots}}_{i \in \ufs{k}} \\
\fm^i &= \txoplus_{j \in \ufs{d}_i}\,\angordmijdot \in \Gsk & 
\fm &= \txoplus_{i \in \ufs{k}}\, \fm^i \in \Gsk \\
\fm^\crdot &= \ang{\fm^i}_{i \in \ufs{k}} \in \Gsk^k &&
\end{aligned}
\\
\begin{split}
(a_{i,j}, \glu^{a_{i,j}}) & \in \Aijangordmijdot \\
(a_{i\crdot}, \glu^{a_{i\crdot}}) &= \ang{(a_{i,j}, \glu^{a_{i,j}})}_{j \in \ufs{d}_i} \in \txsma_{j \in \ufs{d}_i} \Aijangordmijdot \\
(a,\glu^a) &= \ang{(a_{i\crdot}, \glu^{a_{i\crdot}})}_{i \in \ufs{k}} \in \txsma_{i \in \ufs{k}} \txsma_{j \in \ufs{d}_i} \Aijangordmijdot
\end{split}
\end{gathered}\]
Applying the two composites in \cref{Jgo_comp_diagram} to the objects
\[\big((f,\laxf) \sscs \ang{(h_i,\laxhi)}_{i \in \ufs{k}}\big),\]
the two resulting objects \cref{ggcat_k1}
\[\begin{tikzpicture}[vcenter]
\draw[0cell]
(0,0) node (a1) {\Addotdash}
(a1)++(4,0) node (a2) {\Cdash}
;
\draw[1cell=.9]
(a1) edge[transform canvas={yshift=.5ex}] node {\Jgo\gam(f \sscs \ang{h_i}_{i \in \ufs{k}})} (a2)
(a1) edge[transform canvas={yshift=-.5ex}] node[swap] {\ga(\Jgo f \sscs \ang{\Jgo h_i}_{i \in \ufs{k}})} (a2)
;
\end{tikzpicture}\]
have $\angordmdddot$-component pointed functors \cref{ggcat_k1_comp} given as follows, where, by  \cref{ggcat_gamma_comp}, the bottom composite is $\ga(\Jgo f \sscs \ang{\Jgo h_i}_{i \in \ufs{k}})_{\angordmdddot}$.
\begin{equation}\label{Jgo_comp_mal}
\begin{tikzpicture}[vcenter]
\def\h{5.5}
\draw[0cell=.9]
(0,0) node (a1) {\txsma_{i \in \ufs{k}} \txsma_{j \in \ufs{d}_i} \Aijangordmijdot}
(a1)++(1,0) node (a1') {\phantom{\txsma_{j \in \ufs{d}_i}}}
(a1)++(.6*\h,-1.2) node (a2) {\txsma_{i \in \ufs{k}} \Bifmi}
(a1)++(\h,0) node (a3) {\Cfm}
;
\draw[1cell=.85]
(a1) edge node {\Jgo\gam\big(f \sscs \ang{h_i}_{i \in \ufs{k}}\big)_{\angordmdddot}} (a3);
\draw[1cell=.85]
(a1') [rounded corners=2pt, shorten <=0ex] |- node[pos=.25,swap] {\txsma_{i \in \ufs{k}}\, (\Jgo h_i)_{\angordmiddot}} (a2);
\draw[1cell=.85]
(a2) [rounded corners=2pt, shorten <=0ex] -| node[pos=.75] {(\Jgo f)_{\fm^\crdot}} (a3);
\end{tikzpicture}
\end{equation}
To show that the diagram \cref{Jgo_comp_diagram} commutes on objects, it suffices to show that the diagram \cref{Jgo_comp_mal} commutes.

First, we deal with the following two marginal cases.
\begin{itemize}
\item Suppose some $\angordmijdot = \vstar$, the basepoint in $\Gsk$ \cref{Gsk_objects}.  By \cref{Gsk_oplus_vstar}, $\fm = \vstar$.  By the first case in \cref{A_angordn}, $\Cfm = \sys{\C}{\vstar}$ is the terminal category $\boldone$, so the diagram \cref{Jgo_comp_mal} commutes.  Thus, we assume that $\angordmijdot \in \Gsk \setminus \{\vstar\}$ for all $i \in \ufs{k}$ and $j \in \ufs{d}_i$. 
\item Suppose some $(a_{i,j}, \glu^{a_{i,j}}) = (\zero, 1_{\zero})$, the base $\angordmijdot$-system in $\Aij$ \pcref{expl:nsystem_base}.  By \cref{Jgof_m_basesystem}, each composite in the diagram \cref{Jgo_comp_mal} yields the base $\fm$-system $(\zero,1_{\zero})$ in $\C$.  Thus, we assume that $(a_{i,j}, \glu^{a_{i,j}})$ is not the base $\angordmijdot$-system for all $i \in \ufs{k}$ and $j \in \ufs{d}_i$. 
\end{itemize}

Using \cref{Jgof_m_obj_comp}---for each of $\gam\big(f \sscs \ang{h_i}_{i \in \ufs{k}}\big)$, $h_i$, and $f$---and \cref{gam_fh_functor}, the following computation proves that the two composites in \cref{Jgo_comp_mal} yield the same $\angs$-component object \cref{a_angs} in $\C$.
\begin{equation}\label{Jgo_comp_computation}
\begin{split}
&\big(\Jgo\gam(f \sscs \ang{h_i}_{i \in \ufs{k}})_{\angordmdddot} (a,\glu^a)\big)_{\angs}\\
&= \gam\big(f \sscs \ang{h_i}_{i \in \ufs{k}}\big) \ang{\ang{a_{i,j,\, s_{i,j\crdot}}}_{j \in \ufs{d}_i}}_{i \in \ufs{k}}\\
&= f\bang{h_i \ang{a_{i,j,\, s_{i,j\crdot}}}_{j \in \ufs{d}_i} }_{i \in \ufs{k}}\\
&= f\bang{\big( (\Jgo h_i)_{\angordmiddot} (a_{i\crdot}, \glu^{a_{i\crdot}}) \big)_{s_{i\crdot\crdot}}}_{i \in \ufs{k}}\\
&= \big((\Jgo f)_{\fm^\crdot} \ang{(\Jgo h_i)_{\angordmiddot} (a_{i\crdot}, \glu^{a_{i\crdot}}) }_{i \in \ufs{k}}\big)_{\angs}\\
&= \big(\ga(\Jgo f \sscs \ang{\Jgo h_i}_{i \in \ufs{k}})_{\angordmdddot} (a,\glu^a) \big)_{\angs}
\end{split}
\end{equation}
 
\parhead{Gluing}.  To show that the two composites in \cref{Jgo_comp_mal} yield the same gluing morphism \cref{gluing-morphism}, suppose we are given the following object, indices, and partition. 
\[\begin{split}
x \in \Op(t) \qquad (p,q,w) & \in \ufs{k} \times \ufs{d}_p \times \ufs{r}_{p,q} \\
s_{p,q,w} &= \coprod_{v \in \ufs{t}}\, s^v \subseteq \ufs{m}^{p,q,w}
\end{split}\]
We define the following indices and functors.
\[\begin{split}
u &= \big(\txsum_{i=1}^{p-1} \txsum_{j=1}^{d_i} r_{i,j}\big) + \big(\txsum_{j=1}^{q-1} r_{p,j}\big) + w\\
u' &= \big(\txsum_{i=1}^{p-1} d_i\big) + q\\
\fgr &= f\big(\ang{h_i\ang{a_{i,j,\, s_{i,j\crdot}}}_{j \in \ufs{d}_i}}_{i \in \ufs{k}} \,\compp\, - \big) \cn \B_p \to \C\\
\hgr_p &= h_p\big( \ang{a_{p,j,\, s_{p,j\crdot}}}_{j \in \ufs{d}_p} \,\compq\, - \big) \cn \A_{p,q} \to \B_p
\end{split}\]
Note that
\begin{itemize}
\item $s_{p,q,w}$ is the $u$-th entry in $\angs = \ang{\ang{\ang{s_{i,j,\ell}}_{\ell \in \ufs{r}_{i,j}}}_{j \in \ufs{d}_i}}_{i \in \ufs{k}}$ and
\item $s_{p,q\crdot}$ is the $u'$-th entry in $\ang{\ang{s_{i,j\crdot}}_{j \in \ufs{d}_i}}_{i \in \ufs{k}}$.
\end{itemize}
Using 
\begin{itemize}
\item the bottom composite in \cref{Jgo_comp_mal}, 
\item the functoriality of $\fgr$, and 
\item the definition \cref{Jgof_m_gluing} for $f$ and $h_p$, 
\end{itemize}
the gluing morphism of the $\fm$-system 
\[\begin{split}
& \ga\big(\Jgo f \sscs \ang{\Jgo h_i}_{i \in \ufs{k}}\big)_{\angordmdddot} (a,\glu^a) \\
&= (\Jgo f)_{\fm^\crdot} \ang{(\Jgo h_i)_{\angordmiddot} (a_{i\crdot}, \glu^{a_{i\crdot}}) }_{i \in \ufs{k}} \in \Cfm
\end{split}\]
at the component
\begin{equation}\label{xangsusv}
\big(x \sscs \angs, u, \ang{s^v}_{v \in \ufs{t}}\big)
\end{equation}
is the following composite in $\C$.
\begin{equation}\label{Jgo_comp_gluing}
\begin{tikzpicture}[vcenter]
\def\h{2.2} \def\v{-1.4}
\draw[0cell=.85]
(0,0) node (a1) {\gaC_t\big(x; \ang{\fgr \hgr_p a_{p,q, (s_{p,q\crdot} \,\compw\, s^v)}}_{v \in \ufs{t}}\big)}
(a1)++(0,\v) node (a2) {\fgr \gaBp_t\big(x; \ang{\hgr_p a_{p,q, (s_{p,q\crdot} \,\compw\, s^v)}}_{v \in \ufs{t}}\big)}
(a2)++(\h,\v) node (a3) {\fgr \hgr_p \gaApq_t\big(x; \ang{a_{p,q, (s_{p,q\crdot} \,\compw\, s^v)}}_{v \in \ufs{t}} \big)}
(a3)++(\h,-\v) node (a4) {\fgr \hgr_p a_{p,q,\, s_{p,q\crdot}}}
(a4)++(0,-\v) node (a5) {f\bang{h_i \ang{a_{i,j,\, s_{i,j\crdot}}}_{j \in \ufs{d}_i} }_{i \in \ufs{k}}}
;
\draw[1cell=.85]
(a1) edge node[swap] {\laxf_{t,p}} (a2)
(a2) edge[transform canvas={xshift=-1em}] node[swap,pos=.2] {\fgr \laxhp_{t,q}} (a3)
(a3) edge[transform canvas={xshift=1em}] node[swap,pos=.7] {\fgr \hgr_p \glu^{a_{p,q}}} (a4)
(a4) edge[equal] (a5)
;
\end{tikzpicture}
\end{equation}
In \cref{Jgo_comp_gluing}, $\glu^{a_{p,q}}$ is the following gluing morphism. 
\[\gaApq_t\big(x; \ang{a_{p,q, (s_{p,q\crdot} \,\compw\, s^v)}}_{v \in \ufs{t}} \big) 
\fto{\glu^{a_{p,q}}_{x;\, s_{p,q\crdot},\, w, \ang{s^v}_{v \in \ufs{t}}}} a_{p,q,\, s_{p,q\crdot}}\]
By \cref{laxfh_component}, the composite of the first two morphisms in \cref{Jgo_comp_gluing} is
\begin{equation}\label{Jgo_comp_gluing_first}
\fgr \laxhp_{t,q} \circ \laxf_{t,p} = \laxfh_{t,u'}
\end{equation}
where, with the notation in \cref{gam_fh},
\[fh = \gam\big(f \sscs \ang{h_i}_{i \in \ufs{k}}\big).\]
By \cref{Jgof_m_gluing} for $\gam\big(f \sscs \ang{h_i}_{i \in \ufs{k}}\big)$, the gluing morphism of the $\fm$-system
\[\Jgo \gam(f \sscs \ang{h_i}_{i \in \ufs{k}})_{\angordmdddot} (a,\glu^a) \in \Cfm\]
at the component \cref{xangsusv} is the composite
\begin{equation}\label{Jgo_comp_gluing_ii}
\fgr \hgr_p \glu^{a_{p,q}} \circ \laxfh_{t,u'}.
\end{equation}
By \cref{Jgo_comp_gluing_first}, the composites in \cref{Jgo_comp_gluing,Jgo_comp_gluing_ii} are equal.  This proves that the diagram \cref{Jgo_comp_mal} commutes on objects.

To prove that the diagram \cref{Jgo_comp_mal} commutes on morphisms, we reuse the computation in \cref{Jgo_comp_computation} with each $\angordmijdot$-system $(a_{i,j}, \glu^{a_{i,j}})$ replaced by a morphism in $\Aijangordmijdot$, along with \cref{Jgof_m_theta_comp}.  This proves that the diagram \cref{Jgo_comp_diagram} commutes on objects.

\parhead{Morphisms}.  To show that the diagram \cref{Jgo_comp_diagram} commutes on morphisms, we reuse the computation in \cref{Jgo_comp_computation} with
\begin{itemize}
\item $f$ replaced by a $k$-ary $\Op$-transformation \pcref{def:kary_transformation} and
\item each $h_i$ replaced by a $d_i$-ary $\Op$-transformation,
\end{itemize} 
along with \cref{gam_thetas,ggcat_2cell_comp,Jgof_m_theta_comp,Jgotheta_m_angs}.
\end{proof}

The following observation is the main result of this chapter.

\begin{theorem}\label{thm:Jgo_multifunctor}
For each $\Tinf$-operad $\Op$ \pcref{as:OpA}, the assignments in \cref{def:Jgo_multifunctor},
\[\begin{split}
\MultpsO & \fto{\Jgo} \GGCat \andspace\\ 
\MultpspsO & \fto{\Jgosg} \GGCat,
\end{split}\]
are \index{multifunctorial J-theory@multifunctorial $J$-theory!multifunctor}\index{multifunctor!multifunctorial J-theory@multifunctorial $J$-theory!}$\Gcat$-multifunctors.
\end{theorem}

\begin{proof}
We need to verify the three axioms in \cref{def:enr-multicategory-functor} for the Cartesian closed category $\Gcat$ \pcref{expl:Gcat_closed}.  \cref{Jgo_sigma,Jgo_gamma} verify, respectively, the symmetric group action axiom \cref{enr-multifunctor-equivariance} and the composition axiom \cref{v-multifunctor-composition}.  

The unity axiom \cref{enr-multifunctor-unit} for $\Jgo$ means that, for each $\Op$-pseudoalgebra $\B$ with identity $\Op$-morphism $1_\B \cn \B \to \B$ (\cref{def:k_laxmorphism}), there is an equality 
\begin{equation}\label{Jgo_unity}
\Jgo 1_\B = 1_{\Bdash} \cn \Bdash \to \Bdash.
\end{equation}
The identity $\Op$-morphism $1_\B$ consists of
\begin{itemize}
\item the identity functor on $\B$ and 
\item identity action constraints.
\end{itemize}  
Thus, the desired equality \cref{Jgo_unity} follows from the definitions \cref{Jgof_m_obj_comp,Jgof_m_gluing,Jgof_m_theta_comp} of $\Jgo f$ for $f = 1_\B$.  The unity axiom for $\Jgosg$ is the equality
\[\Jgosg 1_\B = 1_{\Bsgdash} \cn \Bsgdash \to \Bsgdash.\]
This equality holds for the same reasons as \cref{Jgo_unity}.
\end{proof}

\begin{remark}[Lax versus Strong]\label{rk:JgoJgosg}
$\MultpspsO$ is the sub-$\Gcat$-multicategory of $\MultpsO$ with
\begin{itemize}
\item the same objects, given by $\Op$-pseudoalgebras (\cref{def:pseudoalgebra});
\item the same $k$-ary 2-cells, given by $k$-ary $\Op$-transformations (\cref{def:kary_transformation}); and 
\item $k$-ary 1-cells given by $k$-ary $\Op$-pseudomorphisms (\cref{def:k_laxmorphism}).
\end{itemize}  
However, the $\Gcat$-multifunctor $\Jgosg$ in \cref{thm:Jgo_multifunctor} is \emph{not} the restriction of $\Jgo$ to $\MultpspsO$.  The reason is that $\Jgo$ is constructed using $\angordn$-systems (\cref{def:nsystem,def:nsystem_morphism,def:Aangordn_system,def:Aangordn_gcat}), while $\Jgosg$ involves strong $\angordn$-systems \cref{sgAordnbe}, with invertible gluing morphisms.
\end{remark}

\begin{remark}[No Strictification Necessary]\label{rk:no_strictification}
The $\Gcat$-multifunctors $\Jgo$ and $\Jgosg$ in \cref{thm:Jgo_multifunctor}, from $\Op$-pseudoalgebras to $\Gskg$-categories, do \emph{not} involve any strictification functor.  The associativity constraint $\phiA$ \cref{phiA} of an $\Op$-pseudoalgebra $\A$ is incorporated into the construction of the $\Gskg$-categories \pcref{A_ptfunctor}
\[\Jgo\A = \Adash \andspace \Jgosg\A = \Asgdash\] 
in the associativity axiom \cref{system_associativity} and the commutativity axiom \cref{system_commutativity}.  In contrast, the $G$-equivariant $K$-theory constructed in \cite{gmmo23} involves a strictification functor, denoted by $\mathrm{St}$ there.  This is one reason why that $K$-theory construction is a \index{multifunctor!nonsymmetric}\emph{nonsymmetric} multifunctor, as explained in \cite[Section 8.3]{gmmo23}.  See also \cref{rk:gmmo-cg-notation,expl:k_laxmorphism}.
\end{remark}

\begin{example}[$\Gskg$-Categories from Symmetric Monoidal $G$-Categories]\label{ex:JgBE}
Both
\begin{itemize}
\item the Barratt-Eccles $\Gcat$-operad $\BE$ (\cref{def:BE-Gcat}) and
\item the $G$-Barratt-Eccles operad $\GBE$ (\cref{def:GBE})
\end{itemize} 
are $\Tinf$-operads \pcref{ex:n_systems}.  By \cref{thm:Jgo_multifunctor}, $\BE$ and $\GBE$ yield four $J$-theory $\Gcat$-multifunctors as follows.
\medskip
\begin{center}
\resizebox{.9\width}{!}{%
{\setlength{\tabcolsep}{1em}
\begin{tabular}{c|cc}
$k$-ary 1-cells & $k$-lax $\Op$-morphisms & $k$-ary $\Op$-pseudomorphisms\\ \hline
naive & $\MultpsBE \fto{\Jgbe} \GGCat$ & $\MultpspsBE \fto{\Jsgbe} \GGCat$\\
genuine & $\MultpsGBE \fto{\Jggbe} \GGCat$ & $\MultpspsGBE \fto{\Jsggbe} \GGCat$
\end{tabular}}}
\end{center}
\medskip
At the object level, these $\Gcat$-multifunctors are given as follows.
\begin{itemize}
\item $\Jgbe$ and $\Jsgbe$ send $\BE$-pseudoalgebras---which correspond to \emph{naive} symmetric monoidal $G$-categories (\cref{def:naive_smGcat,expl:BEpseudo_smcat}) under the 2-equivalences in \cref{thm:BEpseudoalg}---to $\Gskg$-categories \pcref{A_ptfunctor}.
\item $\Jggbe$ and $\Jsggbe$ send $\GBE$-pseudoalgebras---which are \emph{genuine} symmetric monoidal $G$-categories (\cref{def:GBE_pseudoalg})---to $\Gskg$-categories \pcref{A_ptfunctor}.
\end{itemize}
Applications of these $\Gcat$-multifunctors are discussed in \cref{ex:Jgo_preservation,ex:Kgo_preservation}.
\end{example}

\section{$J$-Theory Preserves Equivariant Algebraic Structures}
\label{sec:Jgo_preserves}

This section discusses the fact that the $\Gcat$-multifunctors \pcref{thm:Jgo_multifunctor}
\[\begin{split}
\MultpsO & \fto{\Jgo} \GGCat \andspace\\
\MultpspsO & \fto{\Jgosg} \GGCat
\end{split}\] 
transport algebraic structures in $\MultpsO$ and $\MultpspsO$, parametrized by $\Gcat$-multifunctors, to the same types of algebraic structures in $\GGCat$.  This preservation property of $\Jgo$ is illustrated with equivariant $\Einf$-algebras and $\Ninf$-algebras in the $\Gcat$-multicategories $\MultpsGBE$ and $\MultpsBE$ of genuine and naive symmetric monoidal $G$-categories.

\secoutline
\begin{itemize}
\item \cref{thm:Jgo_preservation} records the fact that $\Jgo$ and $\Jgosg$ remain $\Gcat$-multifunctors after pre-composition with any $\Gcat$-multifunctor from a $\Gcat$-multicategory.
\item \cref{cla_multi} recalls the fact that the classifying space functor $\cla$ induces a change-of-enrichment functor from $\Gcat$-multicategories to $\Gtop$-multicategories.
\item \cref{def:Einfty_operads} recalls $G$-topological and $G$-categorical $\Einf$-operads in the sense of Guillou-May.
\item \cref{def:Ninfty_operads} recalls $\Ninf$-operads in the sense of Blumberg-Hill and extends them to $G$-categorical $\Ninf$-operads.
\item \cref{Jgo_preserves_Einf} records the fact that $\Jgo$ and $\Jgosg$ preserve equivariant $\Einf$-algebras and $\Ninf$-algebras.
\item \cref{ex:Einf_GBE,ex:Einf_steiner,ex:Einf_linear} discuss examples of equivariant $\Einf$-operads and $\Ninf$-operads.
\item \cref{ex:Jgo_preservation} applies \cref{thm:Jgo_preservation} to the $G$-Barratt-Eccles operad $\GBE$ and the Barratt-Eccles $\Gcat$-operad $\BE$.  It uses $\Jggbe$, $\Jsggbe$, $\Jgbe$, and $\Jsgbe$ to transport equivariant algebraic structures, including equivariant $\Einf$-algebras and $\Ninf$-algebras, of genuine or naive symmetric monoidal $G$-categories to the same types of equivariant algebraic structures of $\Gskg$-categories.
\end{itemize}

\begin{theorem}\label{thm:Jgo_preservation}
Suppose $\Op$ is a $\Tinf$-operad \pcref{as:OpA}, and $\cQ$ is a $\Gcat$-multicategory.
\begin{enumerate}
\item\label{Jgo_preservation_i} If $f \cn \cQ \to \MultpsO$ is a $\Gcat$-multifunctor, then the composite
\[\cQ \fto{f} \MultpsO \fto{\Jgo} \GGCat\]
is also a $\Gcat$-multifunctor.
\item\label{Jgo_preservation_ii} If $f \cn \cQ \to \MultpspsO$ is a $\Gcat$-multifunctor, then the composite
\[\cQ \fto{f} \MultpspsO \fto{\Jgosg} \GGCat\]
is also a $\Gcat$-multifunctor. 
\end{enumerate}
\end{theorem}

\begin{proof}
The assertions follow from \cref{thm:Jgo_multifunctor} because multifunctors enriched in a given symmetric monoidal category, including $\Gcat$, are closed under composition \cref{multifunctors_compose}.
\end{proof}

\begin{explanation}\label{expl:Jgo_preservation}
\cref{thm:Jgo_preservation} says that multifunctorial $J$-theory, $\Jgo$ \pcref{thm:Jgo_multifunctor}, sends each $\cQ$-algebra in the $\Gcat$-multicategory $\MultpsO$ \pcref{thm:multpso}, whose objects are $\Op$-pseudoalgebras \pcref{def:pseudoalgebra}, to a $\cQ$-algebra in the $\Gcat$-multicategory $\GGCat$ \pcref{expl:ggcat_gcatenr}, whose objects are $\Gskg$-categories \pcref{expl:ggcat_obj}.  Similarly, multifunctorial strong $J$-theory, $\Jgosg$, sends each $\cQ$-algebra in $\MultpspsO$ to a $\cQ$-algebra in $\GGCat$.  

Moreover, we note the following about the hypotheses of \cref{thm:Jgo_preservation}.
\begin{itemize}
\item There are no restrictions on the group $G$.  
\item There are no restrictions on the $\Gcat$-multicategory $\cQ$.  
\item The only restriction on the $\Gcat$-operad $\Op$ is that it is a $\Tinf$-operad \pcref{as:OpA}.  This means that $\Op$ is a pseudo-commutative operad \pcref{def:pseudocom_operad} with $\Op(1)$ given by a terminal $G$-category.
\end{itemize}
On the other hand, to apply \cref{thm:Jgo_preservation} when $\cQ$ is a $G$-categorical $\Einf$-operad or $\Ninf$-operad, the group $G$ will be assumed to be finite.
\end{explanation}

\subsection*{$G$-Spaces}
To apply \cref{thm:Jgo_preservation} to $G$-categorical $\Einf$-operads and $\Ninf$-operads $\cQ$, as we discuss in \cref{ex:Jgo_preservation} below, we first recall some relevant notions about $G$-spaces and the classifying space functor.

\begin{definition}[$G$-Spaces]\label{def:Gtop}
$\Top$ denotes the complete and cocomplete category of
\begin{itemize}
\item \index{space}\emph{spaces}, which mean compactly generated weak Hausdorff spaces, and
\item \index{morphism}\emph{morphisms}, which mean continuous maps between spaces.
\end{itemize}
Suppose $G$ is a group, which is also regarded as a category with one object $*$ and morphism set $G$, with composition and identity given by, respectively, the group multiplication and the group unit.
\begin{itemize}
\item A \emph{$G$-space}\index{G-space@$G$-space} is a functor $G \to \Top$.  Equivalently, a $G$-space is a space $X$ equipped with a $g$-action homeomorphism $g \cn X \fiso X$ for each $g \in G$ that satisfies the following two conditions.
\begin{enumerate}
\item For the identity element $e \in G$, the $e$-action is equal to $1_X$.
\item For $g,h \in G$, the $hg$-action is equal to the composite $h \circ g$.
\end{enumerate}
\item A \emph{$G$-morphism}\index{G-morphism@$G$-morphism} $f \cn X \to Y$ between $G$-spaces is a natural transformation between functors $G \to \Top$.  Equivalently, $f$ is a morphism of spaces that is \index{G-equivariant@$G$-equivariant}\emph{$G$-equivariant}, in the sense that
\[f(gx) = g(fx) \forspace (g,x) \in G \times X.\]
\item $\Gtop$ denotes the complete and cocomplete category of $G$-spaces and $G$-morphisms, with composition and identities defined in $\Top$.
\item For $G$-spaces $X$ and $Y$, $\Topg(X,Y)$ denotes the $G$-space of all morphisms $X \to Y$, with the compact-open topology and conjugation $G$-action.  For $g \in G$ and a morphism $h \cn X \to Y$, the conjugation $G$-action $g \cdot h$ is the composite
\begin{equation}\label{ginv_h_g}
X \fto{\ginv} X \fto{h} Y \fto{g} Y.
\end{equation}
\end{itemize}
Similar to $\Gcat$ \pcref{expl:Gcat_closed}, the quadruple
\begin{equation}\label{gtop_smclosed}
(\Gtop, \times, *, \Topg)
\end{equation} 
is a Cartesian closed category.
\end{definition}

Using the symmetric monoidal structure on $\Gtop$, we consider $\Gtop$-multicategories and $\Gtop$-multifunctors \pcref{def:enr-multicategory,def:enr-multicategory-functor}.

\subsection*{Change of Enrichment via Classifying Spaces}
The \index{classifying space}classifying space functor $\cla \cn \Cat \to \Top$ \cref{classifying_space} preserves finite products, so it yields a strong symmetric monoidal functor \pcref{def:monoidalfunctor}
\begin{equation}\label{cla_gcat_gtop}
(\Gcat,\times,\boldone) \fto{(\cla,\clatwo,\clazero)} (\Gtop, \times, *).
\end{equation}
The proof of \cref{cla_multi} below for symmetric monoidal functors, including $\cla$, can be found in \cite[Th.\! 11.5.1]{yau-operad}, \cite[Cor.\! 11.16 and Th.\! 12.11 (1)]{bluemonster}, and \cite[6.2.9]{cerberusIII}.  See \cref{mult_change_enr}.

\begin{lemma}\label{cla_multi}
Suppose $G$ is a group.\index{change of enrichment}
\begin{enumerate}
\item\label{cla_multi_i} $\cla$ sends each $\Gcat$-multicategory $\M$ to a $\Gtop$-multicategory $\cla\M$.
\item\label{cla_multi_ii} $\cla$ sends each $\Gcat$-multifunctor $f \cn \M \to \N$ to a $\Gtop$-multifunctor $\cla f \cn \cla\M \to \cla\N$.
\end{enumerate}
\end{lemma}

\begin{explanation}[Change of Enrichment]\label{expl:cla_multi}
For a $\Gcat$-multicategory $\M$, the $\Gtop$-multicategory $\cla\M$ is given as follows
\begin{itemize}
\item $\cla\M$ has the same objects as $\M$.
\item For each $(k+1)$-tuple of objects $(\ang{x_i}_{i \in \ufs{k}}; y)$ with $k \geq 0$, the $k$-ary multimorphism $G$-space of $\cla\M$ is given by
\begin{equation}\label{claMxy}
(\cla\M)(\ang{x_i}_{i \in \ufs{k}}; y) = \cla\big(\M(\ang{x_i}_{i \in \ufs{k}}; y)\big).
\end{equation}
\item The symmetric group action on $\cla\M$ is given by $\cla$ applied to the symmetric group action on $\M$.
\item The unit $1_x^{\cla\M}$ of each object $x \in \cla\M$ is given by composing the unit constraint $\clazero$ with $\cla$ applied to the unit $1_x$ of $x \in \M$, as displayed below.  
\begin{equation}\label{clazero_unit}
\begin{tikzpicture}[baseline={(a1.base)}]
\def\u{.6}
\draw[0cell]
(0,0) node (a1) {*}
(a1)++(1.3,0) node (a2) {\cla\boldone}
(a2)++(2.3,0) node (a3) {\cla\M(x;x)}
;
\draw[1cell=.9]
(a1) edge[equal] node {\clazero} (a2)
(a2) edge node {\cla 1_x} (a3)
(a1) [rounded corners=2pt, shorten >=0ex] |- ($(a2)+(0,\u)$)
-- node[pos=1] {1_x^{\cla\M}} ($(a2)+(.5,\u)$) -| (a3)
;
\end{tikzpicture}
\end{equation}
\item The composition $\ga^{\cla\M}$ of $\cla\M$ is given by composing an iterate of the monoidal constraint $\clatwo$ with $\cla$ applied to the composition $\ga$ of $\M$ \cref{eq:enr-defn-gamma}, as displayed below.
\begin{equation}\label{cla_gamma}
\begin{tikzpicture}[vcenter]
\draw[0cell=.9]
(0,0) node (a11) {\cla\M(\ang{x'}; x'') \times \txprod_{j=1}^n\, \cla\M(\ang{x_j}; x_j')}
(a11)++(5,0) node (a12) {\cla\M(\angx; x'')}
(a11)++(0,-1.5) node (a2) {\cla\big(\M(\ang{x'}; x'') \times \txprod_{j=1}^n\,\M(\ang{x_j}; x_j') \big)}
(a12)++(-.5,0) node (a12') {\phantom{()}}
;
\draw[1cell=.9]
(a11) edge node {\ga^{\cla\M}} (a12)
(a11) edge[transform canvas={xshift=-1em}] node {\iso} node[swap,pos=.4] {\clatwo} (a2)
(a2) [rounded corners=2pt, shorten >=0ex] -| node[pos=.2] {\cla\ga} (a12')
;
\end{tikzpicture}
\end{equation}
\end{itemize}
\cref{cla_multi} \pcref{cla_multi_i} is used below in the definitions of
\begin{itemize}
\item $G$-categorical $\Einf$-operads (\cref{def:Einfty_operads} \pcref{def:Einfty_ii}) and
\item $G$-categorical $\Ninf$-operads (\cref{def:Ninfty_operads} \pcref{def:Ninfty_iv}).
\end{itemize} 
For a $\Gcat$-multifunctor $f \cn \M \to \N$, the $\Gtop$-multifunctor $\cla f \cn \cla\M \to \cla\N$ has the same object assignment as $f$.  Each component of $\cla f$ is given by $\cla$ applied to the corresponding component of $f$, using \cref{claMxy}.

To simplify the notation, we sometimes denote $\cla\M$ and $\cla f$ by, respectively, $\M$ and $f$.
\end{explanation}

\subsection*{Equivariant $\Einf$-Operads and $\Ninf$-Operads}
Next, we recall equivariant $\Einf$-operads in the sense of Guillou-May \cite{gm17} and $\Ninf$-operads in the sense of Blumberg-Hill \cite{blumberg_hill}, in both $\Gtop$ and $\Gcat$, followed by some examples.  Both equivariant $\Einf$-operads and $\Ninf$-operads extend nonequivariant $\Einf$-operads to the $G$-equivariant context.  Since \cref{thm:Jgo_preservation} applies to all $\Gcat$-multicategories $\cQ$, it is applicable to all $G$-categorical $\Einf$-operads and $\Ninf$-operads.  For a symmetric monoidal category $\V$, recall from \cref{def:enr-multicategory-functor} that a $\V$-multifunctor $F \cn \M \to \N$ is also called an \emph{$\M$-algebra}\index{algebra}\index{multicategory!algebra} in $\N$.

\begin{definition}[Equivariant $\Einf$-Operads]\label{def:Einfty_operads}
Suppose $G$ is a finite group.
\begin{enumerate}
\item\label{def:Einfty_i} A \emph{$G$-topological $\Einf$-operad}\index{G-topological E-infinity-operad@$G$-topological $\Einf$-operad}\index{operad!G-topological E-infinity@$G$-topological $\Einf$} \cite[Def.\! 2.1]{gm17} is a $\Gtop$-operad $\cQ$ that satisfies the following three conditions.
\begin{itemize}
\item $\cQ(0)$ is \index{G-contractible@$G$-contractible}$G$-contractible.  This means that, for each subgroup $H \subseteq G$, the $H$-fixed point space $\cQ(0)^H$ is contractible.
\item Each $G \times \Si_n$-space $\cQ(n)$ is \index{Sigma-free@$\Si_n$-free}$\Si_n$-free.  This means that
\[x \si = x \impliespace \si = \id_n\]
for any pair $(x,\si) \in \cQ(n) \times \Si_n$.
\item Each $\cQ(n)$ is a \index{universal principal bundle}universal principal $(G, \Si_n)$-bundle.  This means that, for each subgroup $\Lambda \subseteq G \times \Si_n$, the $\Lambda$-fixed point space satisfies
\[\cQ(n)^\Lambda \simeq 
\begin{cases}
* & \text{if $\Lambda \cap \Si_n = \{\id_n\}$, and}\\
\emptyset & \text{otherwise}.
\end{cases}\]
\end{itemize}
\item\label{def:Einfty_i'} For a $G$-topological $\Einf$-operad $\cQ$ and a $\Gtop$-multicategory $\M$, a $\cQ$-algebra in $\M$ is also called an \index{equivariant E-infinity-algebra@equivariant $\Einf$-algebra}\emph{equivariant $\Einf$-algebra in $\M$}.
\item\label{def:Einfty_ii} A \emph{$G$-categorical $\Einf$-operad}\index{G-categorical E-infinity-operad@$G$-categorical $\Einf$-operad}\index{operad!G-categorical E-infinity@$G$-categorical $\Einf$} \cite[Def.\! 3.11]{gm17} is a $\Gcat$-operad $\cQ$ such that the $\Gtop$-operad $\cla\cQ$ is a $G$-topological $\Einf$-operad.
\item\label{def:Einfty_ii'} For a $G$-categorical $\Einf$-operad $\cQ$ and a $\Gcat$-multicategory $\M$, a $\cQ$-algebra in $\M$ is also called an \emph{equivariant $\Einf$-algebra in $\M$}.\defmark
\end{enumerate}
\end{definition}

\begin{explanation}[Equivariant $\Einf$-Algebras]\label{expl:Einfinity_alg}
Suppose $\M$ is a $\Gcat$-multicategory.
\begin{itemize}
\item An equivariant $\Einf$-algebra in $\M$ is given by a $\Gcat$-multifunctor $f \cn \cQ \to \M$, where $\cQ$ is a $G$-categorical $\Einf$-operad.
\item Changing enrichment along the classifying space functor $\cla$ \pcref{cla_multi}, the $\Gtop$-multifunctor
\[\cla\cQ \fto{\cla f} \cla\M\]
specifies an equivariant $\Einf$-algebra in the $\Gtop$-multicategory $\cla\M$.  
\end{itemize}
The following table summarizes the relationship between the terminology in \cite{gm17} and \cref{def:Einfty_operads}.
\begin{center}
\resizebox{.9\width}{!}{%
{\renewcommand{\arraystretch}{1.2}%
{\setlength{\tabcolsep}{1em}
\begin{tabular}{c|c}
\cite{gm17} & \cref{def:Einfty_operads} \\ \hline
$\Einf$-operad of $G$-spaces (Def.\! 2.1) & $G$-topological $\Einf$-operad \\
$\Einf$-operad of $G$-categories (Def.\! 3.11) & $G$-categorical $\Einf$-operad \\
$\Einfg$-category (Def.\! 4.10) & equivariant $\Einf$-algebra in $\Gcat$ \\
\end{tabular}}}}
\end{center}
In the bottom right entry of the previous table, $\Gcat$ is regarded as a $\Gcat$-multicategory using its symmetric monoidal closed structure \pcref{expl:Gcat_closed}, along with \cref{theorem:v-closed-v-sm}---which gives $\Gcat$ the structure of a symmetric monoidal $\Gcat$-category---and \cref{proposition:monoidal-v-cat-v-multicat}.  
\end{explanation}

Another way to extend $\Einf$-operads to the equivariant context is given in the following definition.

\begin{definition}[$\Ninf$-Operads]\label{def:Ninfty_operads}
Suppose $G$ is a finite group.
\begin{enumerate}
\item\label{def:Ninfty_i} A \emph{family of subgroups}\index{family of subgroups} $\Fam$ \cite[Def.\! 3.3]{blumberg_hill} is a collection of subgroups of $G$ closed under passage to subgroups and under conjugacy.
\item\label{def:Ninfty_ii} Suppose $\Fam$ is a family of subgroups of $G$.  A \index{universal space}\emph{universal space for $\Fam$} is a $G$-space $X$ such that, for each subgroup $H \subseteq G$, the $H$-fixed point space satisfies
\[X^H \simeq \begin{cases}
* & \text{if $H \in \Fam$, and}\\
\emptyset & \text{if $H \not\in \Fam$.}
\end{cases}\]
\item\label{def:Ninfty_iii} An \emph{$\Ninf$-operad}\index{N-infinity-operad@$\Ninf$-operad}\index{operad!N-infinity@$\Ninf$} \cite[Def.\! 3.7]{blumberg_hill} is a $\Gtop$-operad $\cQ$ that satisfies the following three conditions.
\begin{itemize}
\item $\cQ(0)$ is $G$-contractible.
\item Each $G \times \Si_n$-space $\cQ(n)$ is $\Si_n$-free.
\item Each $\cQ(n)$ is a universal space for some family $\Fam_n(\cQ)$ of subgroups of $G \times \Si_n$ that contains all the subgroups of the form $H \times \{\id_n\}$ with $H$ a subgroup of $G$.
\end{itemize}
\item\label{def:Ninfty_iii'} For an $\Ninf$-operad $\cQ$ and a $\Gtop$-multicategory $\M$, a $\cQ$-algebra in $\M$ is also called an \index{N-infinity-algebra@$\Ninf$-algebra}\emph{$\Ninf$-algebra in $\M$}.
\item\label{def:Ninfty_iv} A \emph{$G$-categorical $\Ninf$-operad}\index{G-categorical N-infinity-operad@$G$-categorical $\Ninf$-operad}\index{operad!G-categorical N-infinity@$G$-categorical $\Ninf$} is a $\Gcat$-operad $\cQ$ such that the $\Gtop$-operad $\cla\cQ$ is an $\Ninf$-operad.
\item\label{def:Ninfty_iv'} For a $G$-categorical $\Ninf$-operad $\cQ$ and a $\Gcat$-multicategory $\M$, a $\cQ$-algebra in $\M$ is also called an \emph{$\Ninf$-algebra in $\M$}.\defmark
\end{enumerate}
\end{definition}

\begin{explanation}[Equivariant $\Einf$ is $\Ninf$]\label{expl:Einf_Ninf}
Each $G$-topological $\Einf$-operad $\cQ$ (\cref{def:Einfty_operads} \pcref{def:Einfty_i}) is also an $\Ninf$-operad (\cref{def:Ninfty_operads} \pcref{def:Ninfty_iii}), where the family $\Fam_n(\cQ)$ of subgroups of $G \times \Si_n$ is chosen to be
\[\Fam_n(\cQ) = \big\{\Lambda \subseteq G \times \Si_n \text{ such that } \Lambda \cap \Si_n = \{\id_n\} \big\}.\]
Thus, each $G$-categorical $\Einf$-operad (\cref{def:Einfty_operads} \pcref{def:Einfty_ii}) is also a $G$-categorical $\Ninf$-operad (\cref{def:Ninfty_operads} \pcref{def:Ninfty_iv}).
\end{explanation}

\begin{explanation}[$\Ninf$-Algebras]\label{expl:Ninfinity_alg}
Suppose $\M$ is a $\Gcat$-multicategory.
\begin{itemize}
\item An $\Ninf$-algebra in $\M$ is given by a $\Gcat$-multifunctor $f \cn \cQ \to \M$, where $\cQ$ is a $G$-categorical $\Ninf$-operad. 
\item Changing enrichment along the classifying space functor $\cla$ \pcref{cla_multi}, the $\Gtop$-multifunctor
\[\cla\cQ \fto{\cla f} \cla\M\]
specifies an $\Ninf$-algebra in the $\Gtop$-multicategory $\cla\M$.
\end{itemize}
The notion of a $G$-categorical $\Ninf$-operad in \cref{def:Ninfty_operads} \pcref{def:Ninfty_iv} is not from \cite{blumberg_hill}.  We extend the notion of an $\Ninf$-operad from $\Gtop$ to $\Gcat$ by mimicking \cref{def:Einfty_operads} \pcref{def:Einfty_ii}.  
\end{explanation}

The following observation follows from \cref{thm:Jgo_preservation}, \cref{def:Einfty_operads} \pcref{def:Einfty_ii'}, and \cref{def:Ninfty_operads} \pcref{def:Ninfty_iv'}.

\begin{corollary}\label{Jgo_preserves_Einf}
Suppose $G$ is a finite group, and $\Op$ is a $\Tinf$-operad \pcref{as:OpA}.  Then each of the $\Gcat$-multifunctors
\[\begin{split}
\MultpsO & \fto{\Jgo} \GGCat \andspace \\
\MultpspsO & \fto{\Jgosg} \GGCat
\end{split}\]
preserves equivariant $\Einf$-algebras and $\Ninf$-algebras.
\end{corollary}

\subsection*{Examples of Equivariant $\Einf$-Operads and $\Ninf$-Operads}
\cref{ex:Einf_GBE,ex:Einf_steiner,ex:Einf_linear} discuss examples of equivariant $\Einf$-operads and $\Ninf$-operads, where $G$ is assumed to be a finite group.   For further examples, the reader is referred to \cite{gm17,blumberg_hill,gutierrez_white}.

\begin{example}[$G$-Barratt-Eccles Operad]\label{ex:Einf_GBE}
By either \cite[Theorem 3.11]{gmm17} or \cite[Theorem 3.10]{gm17}, for a finite group $G$, the $G$-Barratt-Eccles operad $\GBE$ \pcref{def:GBE} is
\begin{itemize}
\item a $G$-categorical $\Einf$-operad (\cref{def:Einfty_operads} \pcref{def:Einfty_ii}) and
\item a $G$-categorical $\Ninf$-operad (\cref{def:Ninfty_operads} \pcref{def:Ninfty_iv}).
\end{itemize}
Thus, a $\GBE$-algebra 
\[\GBE \fto{f} \M\]
in a $\Gcat$-multicategory $\M$ is both an equivariant $\Einf$-algebra and an $\Ninf$-algebra in $\M$. Changing enrichment along the classifying space functor $\cla$ \pcref{cla_multi}, $\cla\GBE$ is a $G$-topological $\Einf$-operad and an $\Ninf$-operad.  Moreover, 
\[\cla\GBE \fto{\cla f} \cla\M\]
is an equivariant $\Einf$-algebra and an $\Ninf$-algebra in the $\Gtop$-multicategory $\cla\M$.  This example is used in \cref{ex:Jgo_preservation} below.
\end{example}

The following definition is used in the subsequent two examples and \cref{ch:spectra}.

\begin{definition}\label{def:g_universe}\
\begin{enumerate}
\item A \emph{real $G$-inner product space}\index{G-inner product space@$G$-inner product space} is a pair $(X,\mu)$ consisting of
\begin{itemize}
\item a real inner product space $X$ and
\item a $G$-action $\mu \cn G \times X \to X$
\end{itemize}
such that, for each $g \in G$, $\mu(g,-)$ is a linear isometric isomorphism. 
\item A \emph{complete $G$-universe}\index{complete $G$-universe} $\univ$ is a real $G$-inner product space that contains countably many copies of each irreducible $G$-representation.  One choice of $\univ$ is the direct sum of countably many copies of the regular representation of $G$.\defmark
\end{enumerate}
\end{definition}

\begin{example}[Steiner Operads]\label{ex:Einf_steiner}
For each finite group $G$ and each finite dimensional real $G$-inner product space $V$, denote by $R_V$ the $G$-space of distance-reducing embeddings $V \to V$, with $G$ acting by conjugation.  A \emph{Steiner path} is a continuous morphism 
\[[0,1] \fto{h} R_V \stspace h(1) = 1_V.\]  
Denote by $P_V$ the $G$-space of Steiner paths, with the $G$-action induced by the one on $R_V$.  The function $\pi \cn P_V \to R_V$ is defined by $\pi(h) = h(0)$.

The \dindex{Steiner}{operad}\emph{$V$-Steiner operad}, denoted $\Stein_V$, is the $\Gtop$-operad where $\Stein_V(n)$ is given by the $G \times \Si_n$-space of $n$-tuples $\ang{h_j}_{j \in \ufs{n}}$ of Steiner paths such that $\ang{\pi(h_j)}_{j \in \ufs{n}}$ have disjoint images.
\begin{itemize}
\item The group $G$ acts diagonally on each $n$-tuple.
\item The right $\Si_n$-action on $\Stein_V(n)$ is given by right permutations of $n$-tuples.
\item The operadic unit in $\Stein_V(1)$ is given by the constant Steiner path at $1_V$.
\item The operadic composition is defined pointwise by disjoint union and composition.
\end{itemize}
By crossing with the identity of the orthogonal complement, an inclusion $V \to W$ of finite dimensional real $G$-inner product spaces induces a $\Gtop$-multifunctor $\Stein_V \to \Stein_W$.  For a complete $G$-universe $\univ$, the \emph{infinite Steiner operad} $\KU$ is defined as the union of the $\Gtop$-operads $\Stein_V$ for finite dimensional real $G$-inner product subspaces $V \subset \univ$.  By \cite[Example 2.2]{gm17}, $\KU$ is a $G$-topological $\Einf$-operad.  By \cite[Corollary 3.14]{blumberg_hill}, $\KU$ is an $\Ninf$-operad.
\end{example}

\begin{example}[Linear Isometries Operad]\label{ex:Einf_linear}
For a finite group $G$ and a complete $G$-universe $\univ$, the \dindex{linear isometries}{operad}linear isometries operad $\LU$ is the $\Gtop$-operad with $\LU(n)$ given by the $G \times \Si_n$-space of linear isometries $\txcoprod_{j=1}^n \univ \to \univ$. 
\begin{itemize}
\item The group $G$ acts on $\LU(n)$ by conjugation.
\item $\Si_n$ acts on $\LU(n)$ by permuting copies of $\univ$ in $\txcoprod_{j=1}^n \univ$ from the right. 
\item The operadic unit in $\LU(1)$ is given by $1_\univ$. 
\item The operad composition is given by disjoint union and composition.
\end{itemize}
By \cite[Example 2.3]{gm17}, $\LU$ is a $G$-topological $\Einf$-operad.  By \cite[Lemma 3.15]{blumberg_hill}, $\LU$ is an $\Ninf$-operad.
\end{example}

\subsection*{Applications to Symmetric Monoidal $G$-Categories}
The next example applies \cref{thm:Jgo_preservation} to obtain preservation results for equivariant algebraic structures, including equivariant $\Einf$-algebras and $\Ninf$-algebras, of genuine or naive symmetric monoidal $G$-categories.

\begin{example}\label{ex:Jgo_preservation}
Specializing \cref{thm:Jgo_preservation} \pcref{Jgo_preservation_i} to the case $\Op = \GBE$, the $G$-Barratt-Eccles operad \pcref{def:GBE}, we conclude that, for each $\Gcat$-multicategory $\cQ$ and $\Gcat$-multifunctor $f \cn \cQ \to \MultpsGBE$, the composite
\begin{equation}\label{Q_f_Jggbe}
\cQ \fto{f} \MultpsGBE \fto{\Jggbe} \GGCat
\end{equation}
is also a $\Gcat$-multifunctor.  Thus, multifunctorial $J$-theory $\Jggbe$ for $\GBE$ sends each $\cQ$-algebra in the $\Gcat$-multicategory $\MultpsGBE$---whose objects are $\GBE$-pseudoalgebras, which mean genuine symmetric monoidal $G$-categories \pcref{def:GBE_pseudoalg}---to a $\cQ$-algebra in the $\Gcat$-multicategory $\GGCat$, whose objects are $\Gskg$-categories \pcref{expl:ggcat_obj}.  

Specializing \cref{Q_f_Jggbe} to the case $\cQ = \GBE$, we conclude that, for each $\Gcat$-multifunctor $f \cn \GBE \to \MultpsGBE$, the composite
\begin{equation}\label{f_Jggbe}
\GBE \fto{f} \MultpsGBE \fto{\Jggbe} \GGCat
\end{equation}
is also a $\Gcat$-multifunctor.  We observe that \cref{f_Jggbe} uses the $G$-Barratt-Eccles operad $\GBE$ twice.
\begin{itemize}
\item First, as a $\Tinf$-operad \pcref{as:OpA}, $\Op = \GBE$ is used to form the $\Gcat$-multicategory $\MultpsGBE$ \pcref{thm:multpso}, whose objects are $\GBE$-pseudoalgebras \pcref{def:GBE_pseudoalg}.
\item Next, $\cQ = \GBE$ is used as the domain of $f$ to parametrize a $\GBE$-algebra in $\MultpsGBE$.
\end{itemize}
So far in this example, there are no restrictions on the group $G$.

Next, suppose $G$ is a finite group, so $\GBE$ is a $G$-categorical $\Einf$-operad and also a $G$-categorical $\Ninf$-operad \pcref{ex:Einf_GBE}.   As displayed in \cref{Q_f_Jggbe,f_Jggbe}, multifunctorial $J$-theory $\Jggbe$ for $\GBE$ sends equivariant $\Einf$-algebras and $\Ninf$-algebras in the $\Gcat$-multicategory $\MultpsGBE$ to, respectively, equivariant $\Einf$-algebras and $\Ninf$-algebras in the $\Gcat$-multicategory $\GGCat$.  

\parhead{Strong or naive variants}.  There are three variants of this entire example, as discussed below.
\begin{itemize} 
\item There is a strong variant that uses \cref{thm:Jgo_preservation} \pcref{Jgo_preservation_ii} and multifunctorial strong $J$-theory 
\[\MultpspsGBE \fto{\Jsggbe} \GGCat.\]
The $\Gcat$-multicategory $\MultpspsGBE$ \pcref{thm:multpso} has the same objects and $k$-ary 2-cells as $\MultpsGBE$.  The $k$-ary 1-cells of $\MultpspsGBE$ are $k$-ary $\GBE$-pseudomorphisms \pcref{def:k_laxmorphism}.
\item There are two variants that use the $\Tinf$-operad $\Op = \BE$, the Barratt-Eccles $\Gcat$-operad \pcref{def:BE-Gcat}, along with the $\Gcat$-multifunctor
\[\begin{split}
\MultpsBE & \fto{\Jgbe} \GGCat \orspace\\
\MultpspsBE & \fto{\Jsgbe} \GGCat.
\end{split}\]
The objects of the $\Gcat$-multicategories $\MultpsBE$ and $\MultpspsBE$ are $\BE$-pseudoalgebras.  They correspond to naive symmetric monoidal $G$-categories (\cref{def:naive_smGcat,expl:BEpseudo_smcat}) under the 2-equivalences in \cref{thm:BEpseudoalg}.
\end{itemize}
Further ramifications of this example in the context of orthogonal $G$-spectra are discussed in \cref{ex:Kgo_preservation}.
\end{example}

\part{From $\Gskg$-Categories to Orthogonal $G$-Spectra}
\label{part:kg}

\chapter{Equivariant Orthogonal Spectra}
\label{ch:spectra}
To prepare for our $G$-equivariant $K$-theory symmetric monoidal $\Gtop$-functor $\Kg$ in \cref{ch:semg}, this chapter provides a self-contained review of the category $\GSp$ of orthogonal $G$-spectra.  The $\Gtop$-enriched symmetric monoidal structure on $\GSp$ is discussed in detail.  This chapter assumes no prior knowledge of spectra, and it contains an abundance of details that are not explicitly available in the literature.  Throughout this chapter, $G$ denotes a compact Lie group.

\summary
Orthogonal $G$-spectra are first constructed in \cite{mandell_may}, following mostly the nonequivariant case in \cite{mmss}.  The following table summarizes corresponding notation and terminology in this work and \cite{mandell_may}.
\begin{center}
\resizebox{\columnwidth}{!}{%
{\renewcommand{\arraystretch}{1.3}%
{\setlength{\tabcolsep}{1ex}
\begin{tabular}{cr|cr}
\multicolumn{2}{c|}{This work} & \multicolumn{2}{c}{\cite[Chapter II]{mandell_may}} \\ \hline
$\Gtopst$-category & \eqref{def:enriched-category} & topological $G$-category & \\
$\Gtopst$-functor & \eqref{def:enriched-functor} & continuous $G$-functor & \\
$\Gtopst$-natural transformation & \eqref{def:enriched-natural-transformation} & natural $G$-map & \\
$\Gtopst$ & \eqref{Gtopst_smc} & $G\mathscr{T}$ & \\
internal hom pointed $G$-space $\Topgst(X,Y)$ & \eqref{Gtopst_smc} & function $G$-space $F(X,Y)$ & \\
$\Gtopst$-category $\Topgst$ & \eqref{topgst_gtopst_enr} & $\mathscr{T}_G$ & \\
collection $\IU$ & \eqref{def:indexing_gspace} & $\mathscr{V}$, $\mathscr{V}(U)$ & Def.\! 1.1 \\
$\Gtop$-category $\IU$ & \eqref{def:IU_spaces} & $\mathscr{I}^{\mathscr{V}}_G$, $\mathscr{I}_G$ & Def.\! 2.1 \\
small skeleton $\IUsk$ & \eqref{def:IU_spaces} & $sk\mathscr{I}_G$ & \\
$\IU$-space $\IU \to \Topgst$ & \eqref{def:iu_space} & $\mathscr{I}_G$-space $\mathscr{I}_G \to \mathscr{T}_G$ & Def.\! 2.3 \\
$\Gtop$-category $\IUT$ & \eqref{def:iu_morphism}, \eqref{def:iut_gtop_enr} & $\mathscr{I}_G\mathscr{T}$ & Def.\! 2.3 \\
$\Top$-category $\GIUT$ & \eqref{def:iu_morphism}, \eqref{def:giut_top_enr} & $G\mathscr{I}\mathscr{T}$ & Def.\! 2.3 \\
smash product of $\IU$-spaces $\smau$ & \eqref{x_smau_y} & $\sma$ & (3.6) \\
$G$-sphere $\gsp$ & \eqref{def:g_sphere} & $S^{\mathscr{V}}_G$, $S_G$ & Def.\! 2.1 \\
symmetric monoidal $\Gtop$-category $\GSp$ & \eqref{def:gsp_module}, \eqref{def:gsp_smgtop} & $\mathscr{I}_G\mathscr{S}$ & Def.\! 2.6 \\
smash product of $\gsp$-modules $\smasg$ & \eqref{def:gsp_sma} & $\sma_{S_G}$ & Def.\! 3.9 \\
\end{tabular}}}}
\end{center}
For further discussion of equivariant (stable) homotopy theory, the reader is referred to \cite[Ch.\! 0]{johnson-yau-mackey} and \cite{hill_hopkins_ravenel,lewis_may_steinberger,may_cohomology,tomdieck}.

\connection
In \cref{ch:semg}, the symmetric monoidal $\Gtop$-category $\GSp$ is the codomain of our $G$-equivariant $K$-theory unital symmetric monoidal $\Gtop$-functor
\[\GGTop \fto{\Kg} \GSp.\]
A lot of the detailed discussion in this chapter is used in \cref{ch:semg}.  The symmetric monoidal $\Gtop$-category $\GGTop$ of $\Gskg$-spaces is discussed in \cref{sec:ggtop}.  All the enriched monoidal category theory used in this chapter are reviewed in \cref{ch:prelim,ch:prelim_multicat}.

\organization
This chapter consists of the following sections.

\secname{sec:iu_space}  This section reviews the category $\IUT$ of $\IU$-spaces and its enrichment in $\Gtop$, where $\univ$ is a fixed complete $G$-universe.  The objects of $\IUT$, called $\IU$-spaces, are $\Gtopst$-functors 
\[\IU \fto{X} \Topgst.\]  
Thus, $X$ sends each object $U \in \IU$ to a pointed $G$-space $X_U$, and each isomorphism $f \in \IU$ to a pointed homeomorphism, which is \emph{not} generally $G$-equivariant.  The $G$-equivariance of an $\IU$-space $X$ takes the form \cref{x_gfginv}
\[X_{gf\ginv} = g \circ X_f \circ \ginv\]
for $g \in G$ and $f \in \IU$.

\secname{sec:iu_sm}  This section extends the $\Gtop$-category $\IUT$ to a symmetric monoidal $\Gtop$-category with the smash product $\smau$.  A nontrivial aspect of $\smau$ is that the $\IU$-space $X \smau Y$ is defined as a componentwise coend in $\Topst$ \cref{smau_obj}, not $\Gtopst$, because the component morphisms of $X$ and $Y$ are not $G$-equivariant in general.  The $G$-action on each component of $X \smau Y$ is defined on representatives of the coend \cref{smau_gaction}.

\secname{sec:gmonoid}  This section reviews commutative $G$-monoids in $\IUT$, modules over a commutative $G$-monoid, and their morphisms.  For a commutative $G$-monoid $A$, the category of $A$-modules and morphisms is a $\Gtop$-category, where the $\Gtop$-enrichment is inherited from $\IUT$.

\secname{sec:gspectra_gtop}  This section defines the $\Gtop$-category $\GSp$ of orthogonal $G$-spectra as the category of $\gsp$-modules, where $\gsp$ is the $G$-sphere.  An important subtlety to keep in mind is that, in order for $\GSp$ to be enriched in $\Gtop$, instead of just $\Top$, morphisms of orthogonal $G$-spectra are \emph{not} required to be $G$-equivariant componentwise.

\secname{sec:gspectra_smash}  This section uses the smash product of $\IU$-spaces to construct the smash product of orthogonal $G$-spectra, denoted $\smasg$.  Once again, there is a certain subtlety in the construction of $\smasg$.  In contrast to the construction of $\smau$, the smash product $X \smasg Y$ for two orthogonal $G$-spectra \cref{gsp_sma_coequal} is defined as a componentwise coequalizer in $\Gtopst$.  On the other hand, for morphisms of orthogonal $G$-spectra, the smash product involves coequalizers in $\Topst$ \cref{smasg_mor_coequal}.

\secname{sec:gspectra}  This section defines the symmetric monoidal structure on $\GSp$ in the $\Gtop$-enriched sense, with $\smasg$ as the enriched monoidal product.  The monoidal associator $\asg$ \cref{gsp_associator} and the braiding $\bsg$ \cref{gsp_braiding} for $\GSp$ are induced by the corresponding structures for $\IUT$ discussed in \cref{sec:iu_sm}.  The monoidal identity \cref{gsp_mon_id}, the right monoidal unitor $\rsg$ \cref{gsp_right_unitor}, and the left monoidal unitor $\ellsg$ \cref{gsp_left_unitor} for $\GSp$ involve the $G$-sphere $\gsp$ and the $\gsp$-action for an $\gsp$-module.

\section{The $\Gtop$-Category of $\IU$-Spaces}
\label{sec:iu_space}

Orthogonal $G$-spectra are $\IU$-spaces with extra structure.  This section reviews $\IU$-spaces, $\IU$-morphisms, and the $\Gtop$-category that they form.

\secoutline
\begin{itemize}
\item \cref{def:gtopst} defines the symmetric monoidal closed category $\Gtopst$ of pointed $G$-spaces.
\item \cref{def:indexing_gspace,def:IU_spaces} define the category $\IU$.
\item \cref{def:iu_space} defines $\IU$-spaces, with further discussion given in \cref{expl:iu_space}.
\item \cref{def:iu_morphism} defines ($G$-equivariant) $\IU$-morphisms, with further discussion given in \cref{expl:iu_morphism,expl:eqiu_morphism}.
\item \cref{def:iut_gtop_enr} defines the $\Gtop$-enrichment of the category $\IUT$ of $\IU$-spaces and $\IU$-morphisms, with further discussion given in \cref{expl:iut_gtop_enr}.
\item \cref{def:giut_top_enr} defines the $\Top$-enrichment of the category $\GIUT$ of $\IU$-spaces and $G$-equivariant $\IU$-morphisms.
\end{itemize}

\subsection*{Pointed $G$-Spaces}
Recall the complete and cocomplete symmetric monoidal closed category \pcref{def:Gtop} 
\[(\Gtop, \times, *, \Topg)\]
of $G$-spaces and $G$-morphisms.  The following definition is the topological analogue of \cref{def:gcatst,expl:Gcatst}.

\begin{definition}[Pointed $G$-Spaces]\label{def:gtopst}
Applying \cref{theorem:pC-sm-closed} to the complete and cocomplete Cartesian closed category
\[(\Gtop, \times, *, \Topg)\]
with terminal object $*$ \cref{gtop_smclosed}, we define the complete and cocomplete symmetric monoidal closed category
\begin{equation}\label{Gtopst_smc}
(\Gtopst, \sma, \stplus, \Topgst).
\end{equation}
If $G$ is the trivial group, then $\Gtopst$ is denoted by $\Topst$, whose objects and morphisms are called, respectively, \index{pointed space}\emph{pointed spaces} and \index{pointed morphism}\emph{pointed morphisms}.
\begin{itemize}
\item An object in $\Gtopst$ is called a \index{pointed G-space@pointed $G$-space}\emph{pointed $G$-space}.  It consists of a $G$-space $X$ \pcref{def:Gtop} equipped with a distinguished $G$-fixed object, called the \index{basepoint}\emph{basepoint}.
\item A morphism in $\Gtopst$ is called a \index{pointed G-morphism@pointed $G$-morphism}\emph{pointed $G$-morphism}.  It consists of a $G$-morphism between $G$-spaces that preserves the basepoints. 
\item For pointed $G$-spaces $X$ and $Y$, $\Gtopst(X,Y)$ also denotes the pointed space of pointed $G$-morphisms $X \to Y$, with the compact-open topology.  Its basepoint is the constant morphism at the basepoint of $Y$.
\item Composition and identities in $\Gtopst$ are defined in $\Gtop$.
\item $\sma$ is the smash product defined in \cref{eq:smash} with $\C = \Gtop$.
\item The \index{smash unit}\emph{smash unit} $\stplus$, as defined in \cref{smash-unit-object}, consists of two points, with $G$ acting trivially.
\item For pointed $G$-spaces $X$ and $Y$, the internal hom $\Topgst(X,Y)$, as defined in  \cref{eq:pHom-pullback}, is the $G$-subspace of $\Topg(X,Y)$ consisting of pointed morphisms $X \to Y$.
\begin{itemize}
\item The basepoint of $\Topgst(X,Y)$ is the constant morphism at the basepoint of $Y$.
\item $G$ acts on $\Topgst(X,Y)$ by conjugation \cref{ginv_h_g}.
\end{itemize}
The $G$-fixed point subspace of the pointed $G$-space $\Topgst(X,Y)$ is the pointed space
\[\Gtopst(X,Y) = \Topgst(X,Y)^G\]
of pointed $G$-morphisms.
\end{itemize}
The notation 
\begin{equation}\label{topgst_gtopst_enr}
\Topgst
\end{equation}
also denotes 
\begin{enumerate}
\item the category with 
\begin{itemize}
\item pointed $G$-spaces as objects and
\item pointed morphisms between them,
\end{itemize}
and
\item the $\Gtopst$-category \pcref{def:enriched-category} with
\begin{itemize}
\item pointed $G$-spaces as objects,
\item hom $\Gtopst$-object given by the pointed $G$-space $\Topgst(X,Y)$ for any pair $(X,Y)$ of pointed $G$-spaces, and
\item composition given by that of pointed morphisms.\defmark
\end{itemize}
\end{enumerate}
\end{definition}

\begin{convention}[Disjoint Basepoint]\label{conv:disjoint_gbasept}
For each $G$-space $X$, denote by $X_\splus$ the pointed $G$-space obtained from $X$ by adjoining a disjoint $G$-fixed basepoint\dindex{disjoint}{basepoint} $*$.  Applying this procedure to each hom $\Gtop$-object, a $\Gtop$-category becomes a $\Gtopst$-category.
\end{convention}

\subsection*{The Category $\IU$}
Fix a complete $G$-universe $\univ$ \pcref{def:g_universe} for the rest of this chapter.

\begin{definition}\label{def:indexing_gspace}\
\begin{itemize}
\item An \emph{indexing $G$-space}\index{indexing G-space@indexing $G$-space} in $\univ$ is a finite dimensional real $G$-inner product subspace $V \subset \univ$.
\item $\IU$ denotes the collection of all real $G$-inner product spaces that are isomorphic to some indexing $G$-spaces in $\univ$, via $G$-equivariant linear isometric isomorphisms.
\item For $V \in \IU$, the \emph{$V$-sphere}\index{sphere} $S^V$ is the pointed $G$-space given by the one-point compactification $V \sqcup \{\infty\}$ of $V$, with $G$-fixed basepoint $\infty$.\defmark
\end{itemize}
\end{definition}

\begin{definition}[Indexing Category $\IU$]\label{def:IU_spaces}
Extend the collection $\IU$ \pcref{def:indexing_gspace} into a category whose morphisms are linear isometric isomorphisms between real $G$-inner product spaces.  Note that each morphism in $\IU$ is invertible, but it is \emph{not} required to be $G$-equivariant.

Furthermore, $\IU$ is equipped with the following structures.
\begin{enumerate}
\item\label{def:iu_spaces_i} $(\IU,\oplus,0,\xi)$ is a symmetric monoidal category.
\begin{itemize}
\item $\oplus$ is the direct sum for real $G$-inner product spaces and linear isometric isomorphisms.
\item The monoidal unit is the one-point space 0.
\item The braiding 
\[V \oplus W \fto[\iso]{\xi_{V,W}} W \oplus V\]
swaps the two arguments.
\end{itemize}
\item\label{def:iu_spaces_ii} For objects $V,W \in \IU$, the set $\IU(V,W)$ is topologized as a $G$-subspace of $\Topg(V,W)$, so $G$ acts on $\IU(V,W)$ by conjugation \cref{ginv_h_g}.  With these hom $\Gtop$-objects, $\IU$ becomes a $\Gtop$-category \pcref{def:enriched-category}.
\item\label{def:iu_spaces_iii} By \cref{conv:disjoint_gbasept}, $\IU$ is also regarded as a $\Gtopst$-category.
\item\label{def:iu_spaces_iv} Denote by $\IUsk$ the small full subcategory of $\IU$ whose objects are the indexing $G$-spaces in $\univ$.  Note that $\IUsk$ is a small skeleton of $\IU$, since the inclusion $\IUsk \to \IU$ is an equivalence of categories. 
\item\label{def:iu_spaces_v} The small skeleton $\IUsk$ inherits a symmetric monoidal structure from $(\IU,\oplus)$ such that the inclusion functor 
\[(\IUsk,\oplus) \to (\IU,\oplus)\]
is strong symmetric monoidal.
\item\label{def:iu_spaces_vi} Using the Axiom of Choice, for each object $V \in \IU$, we choose a $G$-linear isometric isomorphism
\begin{equation}\label{VV'}
V \fto[\iso]{\upphi_V} V'
\end{equation}
with $V' \in \IUsk$, subject to the condition that 
\begin{equation}\label{upphi_one}
\upphi_V = 1_V \ifspace V \in \IUsk.
\end{equation}
\end{enumerate}
This finishes the definition.
\end{definition}

\subsection*{$\IU$-Spaces}

In the following definition, $(\pSet, \sma, \stplus)$ denotes the symmetric monoidal closed category of pointed sets and pointed functions, with the monoidal product given by the smash product $\sma$ \cref{eq:smash}.  The functor
\begin{equation}\label{gtopst_pset}
\Gtopst \fto{\und} \pSet
\end{equation}
that forgets the $G$-space structure is a strict symmetric monoidal faithful functor.

\begin{definition}\label{def:iu_space}
An \emph{$\IU$-space}\index{IU-space@$\IU$-space} is a $\Gtopst$-functor \pcref{def:enriched-functor}
\[\IU \fto{X} \Topgst\]
from the $\Gtopst$-category $\IU$ \pcref{def:IU_spaces} to the $\Gtopst$-category $\Topgst$ \cref{topgst_gtopst_enr}.  Moreover, changing enrichment along $\und$ \cref{gtopst_pset} for an $\IU$-space $X$, the resulting $\pSet$-functor is called the \emph{underlying functor}\index{underlying functor} of $X$ and is denoted by the same notation.
\end{definition}

\begin{explanation}[Unraveling $\IU$-Spaces]\label{expl:iu_space}
Interpreting \cref{def:enriched-functor} with $\V = \Gtopst$, an $\IU$-space $X \cn \IU \to \Topgst$ is determined by
\begin{itemize}
\item a pointed $G$-space $X_V$ for each object $V \in \IU$ and
\item a component pointed $G$-morphism between pointed $G$-spaces
\begin{equation}\label{iu_space_comp_mor}
\IU(V,W)_\splus \fto{X} \Topgst(X_V, X_W)
\end{equation}
for each pair $(V,W)$ of objects in $\IU$.
\end{itemize}
It is required that $X$ preserves composition and identities, in the sense that the following diagrams in $\Gtopst$ commute for objects $U,V,W \in \IU$, where $\mcomp$ denotes composition.
\begin{equation}\label{iu_space_axioms}
\begin{tikzpicture}[vcenter]
\def\v{-1.4}
\draw[0cell=.8]
(0,0) node (a11) {\IU(V,W)_\splus \sma \IU(U,V)_\splus}
(a11)++(3.7,0) node (a12) {\IU(U,W)_\splus}
(a11)++(0,\v) node (a21) {\Topgst(X_V,X_W) \sma \Topgst(X_U,X_V)}
(a12)++(0,\v) node (a22) {\Topgst(X_U,X_W)}
(a12)++(1.7,0) node (b11) {\stplus}
(b11)++(1.6,0) node (b12) {\IU(V,V)_\splus}
(b11)++(.9,\v) node (b22) {\Topgst(X_V,X_V)}
;
\draw[1cell=.8]
(a11) edge node {\mcomp} (a12)
(a12) edge node {X} (a22)
(a11) edge[transform canvas={xshift=1.2em}] node[swap] {X \sma X} (a21)
(a21) edge node {\mcomp} (a22)
(b11) edge node {1_V} (b12)
(b12) edge node[pos=.3] {X} (b22)
(b11) edge node[swap,pos=.3] {1_{X_V}} (b22)
;
\end{tikzpicture}
\end{equation}
The pointed $G$-space $X_V$ is also denoted by $X(V)$.

\parhead{Components}. The component pointed $G$-morphism in \cref{iu_space_comp_mor} is equivalent to a $G$-morphism between $G$-spaces
\begin{equation}\label{iu_space_comp'}
\IU(V,W) \fto{X} \Topgst(X_V, X_W).
\end{equation}
At the point-set level, the $G$-morphism in \cref{iu_space_comp'} sends each linear isometric isomorphism $f \cn V \fiso W$ in $\IU$ to a pointed homeomorphism between pointed $G$-spaces
\begin{equation}\label{iu_space_xf}
X_V \fto[\iso]{X_f} X_W.
\end{equation}
We emphasize that, as a point in $\Topgst(X_V, X_W)$, $X_f$ is \emph{not} required to be $G$-equivariant.

\parhead{Equivariance}. The group $G$ acts on the domain and the codomain of \cref{iu_space_comp_mor} by conjugation \cref{ginv_h_g}.  Thus, the $G$-equivariance of \cref{iu_space_comp_mor}---or, equivalently, of \cref{iu_space_comp'}---means that, for each element $g \in G$ and each linear isometric isomorphism $f \cn V \fiso W$ in $\IU$, the following diagram of pointed homeomorphisms commutes.
\begin{equation}\label{x_gfginv}
\begin{tikzpicture}[vcenter]
\def\v{-1.4}
\draw[0cell]
(0,0) node (a11) {X_V}
(a11)++(2.5,0) node (a12) {X_W}
(a11)++(0,\v) node (a21) {X_V}
(a12)++(0,\v) node (a22) {X_W}
;
\draw[1cell=.9]
(a11) edge node {X_{gf\ginv}} (a12)
(a11) edge node[swap] {\ginv} (a21)
(a21) edge node {X_f} (a22)
(a22) edge node[swap] {g} (a12)
;
\end{tikzpicture}
\end{equation}
In particular, if $f \cn V \fiso W$ is $G$-equivariant, then the commutative diagram \cref{x_gfginv} means that $X_f = X_{gf\ginv}$ is also $G$-equivariant.  Thus, a part of the $G$-equivariance of $X$ is that it preserves $G$-equivariant morphisms.

\parhead{Underlying functor}.  The underlying functor of $X$ is the $\pSet$-functor with
\begin{itemize}
\item object assignment 
\[(V \in \IU) \mapsto X_V,\] 
which is a pointed $G$-space, and
\item morphism assignment 
\[(V \fto{f} W) \mapsto X_f,\] 
which is a pointed homeomorphism.\defmark
\end{itemize}
\end{explanation}

\subsection*{$\IU$-Morphisms}
Next, we define morphisms between $\IU$-spaces.  The next definition is followed by detailed discussion.

\begin{definition}\label{def:iu_morphism}\
\begin{enumerate}
\item\label{def:iu_morphism_i}
For $\IU$-spaces $X,Y \cn \IU \to \Topgst$ \pcref{def:iu_space}, an \emph{$\IU$-morphism}\index{IU-morphism@$\IU$-morphism} $\theta \cn X \to Y$ is a $\pSet$-natural transformation \pcref{def:enriched-natural-transformation}
\begin{equation}\label{iu_mor}
\begin{tikzpicture}[vcenter]
\def\t{28}
\draw[0cell]
(0,0) node (a1) {\phantom{A}}
(a1)++(1.8,0) node (a2) {\phantom{A}}
(a1)++(-.1,0) node (a1') {\IU}
(a2)++(.21,-.04) node (a2') {\Topgst}
;
\draw[1cell=.9]
(a1) edge[bend left=\t] node {X} (a2)
(a1) edge[bend right=\t] node[swap] {Y} (a2)
;
\draw[2cell]
node[between=a1 and a2 at .45, rotate=-90, 2label={above,\theta}] {\Rightarrow}
;
\end{tikzpicture}
\end{equation}
between the underlying functors of $X$ and $Y$.
\item\label{def:iu_morphism_ii}  
The category
\begin{equation}\label{IUT}
\IUT
\end{equation}
has
\begin{itemize}
\item $\IU$-spaces \pcref{def:iu_space} as objects,
\item $\IU$-morphisms \cref{iu_mor} as morphisms, and
\item identities and composition given by those of $\pSet$-natural transformations.
\end{itemize}
\item\label{def:iu_morphism_iii}  
A \emph{$G$-equivariant $\IU$-morphism}\index{G-equivariant IU-morphism@$G$-equivariant $\IU$-morphism}\index{IU-morphism@$\IU$-morphism!G-equivariant@$G$-equivariant} $\psi \cn X \to Y$ is a $\Gtopst$-natural transformation \pcref{def:enriched-natural-transformation}.
\item\label{def:iu_morphism_iv}
The category
\begin{equation}\label{GIUT}
\GIUT
\end{equation}
has
\begin{itemize}
\item $\IU$-spaces \pcref{def:iu_space} as objects,
\item $G$-equivariant $\IU$-morphisms as morphisms, and
\item identities and composition given by those of $\Gtopst$-natural transformations.\defmark
\end{itemize}
\end{enumerate}
\end{definition}

\begin{explanation}[$\IU$-Morphisms]\label{expl:iu_morphism}
An $\IU$-morphism $\theta \cn X \to Y$ \cref{iu_mor} is determined by a pointed function
\begin{equation}\label{iu_morphism_comp}
\stplus \fto{\theta_V} \Topgst(X_V,Y_V) \foreachspace V \in \IU,
\end{equation}
which is given by a pointed morphism between pointed $G$-spaces
\begin{equation}\label{iu_mor_component}
X_V \fto{\theta_V} Y_V,
\end{equation}
such that the following naturality diagram commutes for each linear isometric isomorphism $f \cn V \fiso W$ in $\IU$.
\begin{equation}\label{iu_mor_natural}
\begin{tikzpicture}[vcenter]
\def\v{-1.4}
\draw[0cell]
(0,0) node (a11) {X_V}
(a11)++(2.3,0) node (a12) {Y_V}
(a11)++(0,\v) node (a21) {X_W}
(a12)++(0,\v) node (a22) {Y_W}
;
\draw[1cell=.9]
(a11) edge node {\theta_V} (a12)
(a12) edge node {Y_f} node[swap] {\iso} (a22)
(a11) edge node[swap] {X_f} node {\iso} (a21)
(a21) edge node {\theta_W} (a22)
;
\end{tikzpicture}
\end{equation}
The components $\theta_V$ \cref{iu_mor_component} are \emph{not} required to be $G$-equivariant.  Identities and composition are defined componentwise using \cref{iu_mor_component}.

\parhead{Small hom sets}.  For each object $V \in \IU$, the naturality diagram \cref{iu_mor_natural} implies that $\theta_V$ is equal to the following composite pointed morphism, where $\upphi_V \cn V \fiso V'$ is the chosen $G$-linear isometric isomorphism in \cref{VV'} with $V' \in \IUsk$.
\begin{equation}\label{iu_mor_extend}
\begin{tikzpicture}[vcenter]
\def\v{-1.4}
\draw[0cell]
(0,0) node (a11) {X_V}
(a11)++(2.3,0) node (a12) {Y_V}
(a11)++(0,\v) node (a21) {X_{V'}}
(a12)++(0,\v) node (a22) {Y_{V'}}
;
\draw[1cell=.9]
(a11) edge node {\theta_V} (a12)
(a11) edge node[swap] {X_{\upphi_V}} node {\iso} (a21)
(a21) edge node {\theta_{V'}} (a22)
(a22) edge node[swap,pos=.45] {Y_{\upphi_V^{-1}}} node {\iso} (a12)
;
\end{tikzpicture}
\end{equation}
Thus, each $\IU$-morphism $\theta \cn X \to Y$ is determined by the \emph{set} of components
\begin{equation}\label{theta_iusk}
\Big\{X_V \fto{\theta_{V}} Y_V \cn V \in \IUsk\Big\}.
\end{equation}
In summary, given $\IU$-spaces $X$ and $Y$, the collection of $\IU$-morphisms, $\IUT(X,Y)$, is a set.  This ensures that the category $\IUT$ \cref{IUT} is well defined.

\parhead{Determination on $\IUsk$}.  Conversely, to define an $\IU$-morphism $\theta \cn X \to Y$, it suffices to define the set of components in \cref{theta_iusk} such that the naturality diagram \cref{iu_mor_natural} commutes for each $f \cn V \fiso W$ in $\IUsk$.  Indeed, one must extend such a partially-defined $\theta$ to all of $\IU$ using \cref{iu_mor_extend}.  For a general isomorphism $f \cn V \fiso W$ in $\IU$, the naturality diagram \cref{iu_mor_natural} is the following diagram, where $V', W' \in \IUsk$.
\begin{equation}\label{iusk_mor}
\begin{tikzpicture}[vcenter]
\def\h{2} \def\v{-1.4}
\draw[0cell]
(0,0) node (a11) {X_V}
(a11)++(\h,0) node (a12) {X_{V'}}
(a12)++(\h,0) node (a13) {Y_{V'}}
(a13)++(\h,0) node (a14) {Y_V}
(a11)++(0,\v) node (a21) {X_W}
(a12)++(0,\v) node (a22) {X_{W'}}
(a13)++(0,\v) node (a23) {Y_{W'}}
(a14)++(0,\v) node (a24) {Y_W}
;
\draw[1cell=.9]
(a11) edge node {X_{\upphi_V}} node[swap] {\iso} (a12)
(a12) edge node {\theta_{V'}} (a13)
(a13) edge node {Y_{\upphi_V^{-1}}} node[swap] {\iso} (a14)
(a21) edge node {X_{\upphi_W}} node[swap] {\iso} (a22)
(a22) edge node {\theta_{W'}} (a23)
(a23) edge node {Y_{\upphi_W^{-1}}} node[swap] {\iso} (a24)
(a11) edge node[swap] {X_f} (a21)
(a14) edge node {Y_f} (a24)
;
\end{tikzpicture}
\end{equation}
The diagram \cref{iusk_mor} commutes because it is the naturality diagram \cref{iu_mor_natural} of $\theta$ for the composite
\[V' \fto[\iso]{\upphi_V^{-1}} V \fto{f} W \fto[\iso]{\upphi_W} W' \inspace \IUsk,\]
using the functoriality of $X$ and $Y$ \cref{iu_space_axioms}.  In summary, specifying an $\IU$-morphism is equivalent to specifying the set of components in \cref{theta_iusk} such that the naturality diagram \cref{iu_mor_natural} commutes for each morphism $f$ in $\IUsk$.
\end{explanation}

\begin{explanation}[Equivariant $\IU$-Morphisms]\label{expl:eqiu_morphism}
A $G$-equivariant $\IU$-morphism $\psi \cn X \to Y$ is determined by a pointed $G$-morphism
\begin{equation}\label{eqiu_morphism_comp}
\stplus \fto{\psi_V} \Topgst(X_V,Y_V) \foreachspace V \in \IU,
\end{equation}
which is given by a pointed $G$-morphism between pointed $G$-spaces
\begin{equation}\label{eqiu_mor_component}
X_V \fto{\psi_V} Y_V,
\end{equation}
such that the naturality diagram \cref{iu_mor_natural} commutes.  Thus, the diagram \cref{iu_mor_extend} ensures that the category $\GIUT$ \cref{GIUT} is well defined, in the sense that it has small hom sets.  Identities and composition are defined componentwise using \cref{eqiu_mor_component}.  

The discussion in the paragraphs of \cref{theta_iusk,iusk_mor} about determination on $\IUsk$ also applies to $G$-equivariant $\IU$-morphisms.  Thus, specifying a $G$-equivariant $\IU$-morphism $\psi$ is equivalent to specifying the \emph{set} of component pointed $G$-morphisms
\begin{equation}\label{psi_iusk}
\Big\{X_V \fto{\psi_{V}} Y_V \cn V \in \IUsk\Big\}
\end{equation}
such that the naturality diagram \cref{iu_mor_natural} commutes for each morphism $f$ in $\IUsk$.  The other components of $\psi$, for $V \in \IU$, are given by the diagrams \cref{iu_mor_extend} using the chosen $G$-linear isometric isomorphisms $\upphi_V$ in \cref{VV'}.
\end{explanation}

\subsection*{Topological Enrichment}
The rest of this section discusses enrichment of the categories $\IUT$ and $\GIUT$ in, respectively, $\Gtop$ and $\Top$.

\begin{definition}[$\Gtop$-Enrichment of $\IUT$]\label{def:iut_gtop_enr}\
\begin{itemize}
\item For any pair of $\IU$-spaces $X,Y \cn \IU \to \Topgst$ \pcref{def:iu_space}, the morphism set
\begin{equation}\label{IUTXY_topology}
\IUT(X,Y) \bigsubset \prod_{V \in \IUsk} \Topgst(X_V,Y_V)
\end{equation}
is given the $G$-subspace topology.
\item Equipped with these hom $G$-spaces, the category $\IUT$ in \cref{IUT} becomes a $\Gtop$-category \pcref{def:enriched-category}.\defmark
\end{itemize}
\end{definition}

\begin{explanation}[Conjugation $G$-Action on $\IU$-Morphisms]\label{expl:iut_gtop_enr}
\cref{def:iut_gtop_enr} uses \cref{iu_mor_extend,theta_iusk,iusk_mor} as follows.  Consider an $\IU$-morphism $\theta \cn X \to Y$ in the $G$-space $\IUT(X,Y)$ in \cref{IUTXY_topology} and an element $g \in G$.  The $g$-action on $\theta$ yields the $\IU$-morphism 
\[X \fto{g \cdot \theta} Y\]
whose $V$-component pointed morphism, for $V \in \IUsk$, is given by conjugating $\theta_V$ as follows.
\begin{equation}\label{iu_mor_gaction}
\begin{tikzpicture}[vcenter]
\def\v{-1.4}
\draw[0cell]
(0,0) node (a11) {X_V}
(a11)++(2.5,0) node (a12) {Y_V}
(a11)++(0,\v) node (a21) {X_{V}}
(a12)++(0,\v) node (a22) {Y_{V}}
;
\draw[1cell=.9]
(a11) edge node {(g \cdot \theta)_V} (a12)
(a11) edge node[swap] {\ginv} node {\iso} (a21)
(a21) edge node {\theta_{V}} (a22)
(a22) edge node[swap,pos=.5] {g} node {\iso} (a12)
;
\end{tikzpicture}
\end{equation}

\parhead{Naturality}. 
For an isomorphism $f \cn V \fiso W$ in $\IUsk$, the naturality diagram \cref{iu_mor_natural} for $g \cdot \theta$ is the following diagram.
\begin{equation}\label{iusk_mor_gaction}
\begin{tikzpicture}[vcenter]
\def\h{2} \def\v{-1.4}
\draw[0cell]
(0,0) node (a11) {X_V}
(a11)++(\h,0) node (a12) {X_{V}}
(a12)++(\h,0) node (a13) {Y_{V}}
(a13)++(\h,0) node (a14) {Y_V}
(a11)++(0,\v) node (a21) {X_W}
(a12)++(0,\v) node (a22) {X_{W}}
(a13)++(0,\v) node (a23) {Y_{W}}
(a14)++(0,\v) node (a24) {Y_W}
;
\draw[1cell=.9]
(a11) edge node {\ginv} node[swap] {\iso} (a12)
(a12) edge node {\theta_{V}} (a13)
(a13) edge node {g} node[swap] {\iso} (a14)
(a21) edge node {\ginv} node[swap] {\iso} (a22)
(a22) edge node {\theta_{W}} (a23)
(a23) edge node {g} node[swap] {\iso} (a24)
(a11) edge node[swap] {X_f} (a21)
(a14) edge node {Y_f} (a24)
;
\end{tikzpicture}
\end{equation}
The diagram \cref{iusk_mor_gaction} commutes because it is the naturality diagram for $\theta$ at the composite 
\[V \fto{g} V \fto{f} W \fto{\ginv} W,\]
since there are equalities
\[\begin{split}
X_{\ginv f g} &= \ginv \circ X_f \circ g \andspace\\
Y_{\ginv f g} &= \ginv \circ Y_f \circ g 
\end{split}\]
by the $G$-equivariance of $X$ and $Y$ \cref{x_gfginv}.

\parhead{Extension to $\IU$}. 
For each object $V \in \IU$, the $V$-component pointed morphism of $g \cdot \theta$ is given by the boundary of the following diagram, where  $\upphi_V \cn V \fiso V'$ is the chosen $G$-linear isometric isomorphism in \cref{VV'} with $V' \in \IUsk$.
\begin{equation}\label{iu_mor_gextend}
\begin{tikzpicture}[vcenter]
\def\h{2} \def\v{-1.4} \def\u{-.7}
\draw[0cell=.9]
(0,0) node (a11) {X_V}
(a11)++(\h,0) node (a12) {X_V}
(a12)++(\h,0) node (a13) {Y_V}
(a13)++(\h,0) node (a14) {Y_V}
(a11)++(0,\v) node (a21) {X_{V'}}
(a21)++(\h,0) node (a22) {X_{V'}}
(a22)++(\h,0) node (a23) {Y_{V'}}
(a23)++(\h,0) node (a24) {Y_{V'}}
;
\draw[1cell=.8]
(a11) [rounded corners=2pt] |- ($(a12)+(0,-\u)$) -- node {(g \cdot \theta)_{V}} ($(a13)+(0,-\u)$) -| (a14)
;
\draw[1cell=.8]
(a11) edge node {\ginv} (a12)
(a12) edge node {\theta_V} (a13)
(a13) edge node {g} (a14)
(a21) edge node {\ginv} (a22)
(a22) edge node {\theta_{V'}} (a23)
(a23) edge node {g} (a24)
(a11) edge node[swap] {X_{\upphi_V}} (a21)
(a12) edge node {X_{\upphi_V}} (a22)
(a23) edge node[swap,pos=.45] {Y_{\upphi_V^{-1}}} (a13)
(a24) edge node[swap,pos=.45] {Y_{\upphi_V^{-1}}} (a14)
(a21) [rounded corners=2pt] |- ($(a22)+(0,\u)$) -- node {(g \cdot \theta)_{V'}} ($(a23)+(0,\u)$) -| (a24)
;
\end{tikzpicture}
\end{equation}
\begin{itemize}
\item In the diagram \cref{iu_mor_gextend}, the middle left and right squares commute by
\begin{itemize}
\item the $G$-equivariance \cref{x_gfginv} of $X$ and $Y$, and
\item the fact that $\upphi_V$ and $\upphi_W$ are $G$-equivariant.
\end{itemize}
\item The middle square is the commutative diagram \cref{iu_mor_extend}.
\item The bottom rectangle is the definition \cref{iu_mor_gaction} of $(g \cdot \theta)_{V'}$.
\end{itemize}
Thus, the top commutative rectangle in \cref{iu_mor_gextend} implies that
\begin{equation}\label{gtheta_v}
(g \cdot \theta)_V = g \circ \theta_V \circ \ginv \forspace (g,V) \in G \times \IU.
\end{equation}
The discussion in the paragraph containing \cref{iusk_mor} shows that \cref{gtheta_v} yields a well-defined $\IU$-morphism $g \cdot \theta \cn X \to Y$.  In summary, the $G$-action on the hom $G$-space $\IUT(X,Y)$ is given componentwise by conjugation.

\parhead{Equivariant $\IU$-morphisms}. 
Recall that $G$-equivariant $\IU$-morphisms are componentwise pointed $G$-morphisms \cref{eqiu_mor_component}.  Passing from the hom $G$-space $\IUT(X,Y)$ to its $G$-fixed point space yields the equality
\begin{equation}\label{giutxy_fixedpt}
\GIUT(X,Y) = \IUT(X,Y)^G,
\end{equation}
where $\GIUT(X,Y)$ is the hom set in the category $\GIUT$ in \cref{GIUT}.
\end{explanation}

\begin{definition}[$\Top$-Enrichment of $\GIUT$]\label{def:giut_top_enr}\
\begin{itemize}
\item For any pair of $\IU$-spaces $X,Y \cn \IU \to \Topgst$ \pcref{def:iu_space}, the morphism set $\GIUT(X,Y)$ is given the subspace topology using \cref{IUTXY_topology,giutxy_fixedpt}.
\item Equipped with these hom spaces, the category $\GIUT$ in \cref{GIUT} becomes a $\Top$-category \pcref{def:enriched-category}.\defmark
\end{itemize}
\end{definition}

\section{The Symmetric Monoidal $\Gtop$-Category of $\IU$-Spaces}
\label{sec:iu_sm}

This section reviews the symmetric monoidal $\Gtop$-category structure on the $\Gtop$-category $\IUT$ of $\IU$-spaces and $\IU$-morphisms \pcref{def:iu_space,def:iu_morphism,def:iut_gtop_enr} for a compact Lie group $G$.

\secoutline
\begin{itemize}
\item \cref{def:iuspace_day,smau_welldef} construct the smash product $X \smau Y$ of two $\IU$-spaces.  This smash product is defined using a coend in $\Topst$ and \emph{not} in $\Gtopst$.  Its $G$-action is defined on representatives in \cref{smau_gaction}.  The computation \cref{smau_gact_welldef} proves that this $G$-action is well defined.
\item \cref{def:iumor_day} defines the smash product of two $\IU$-morphisms.
\item \cref{def:IU_smgtop,iut_smgtop} construct the symmetric monoidal $\Gtop$-category structure on $\IUT$.
\item \cref{def:giut_smtop} defines the symmetric monoidal $\Top$-category structure on the $\Top$-category $\GIUT$.
\end{itemize}

\subsection*{Smash Product of $\IU$-Spaces}
Recall that $(\IUsk,\oplus)$ is the small skeleton of $(\IU,\oplus)$ consisting of indexing $G$-spaces in the complete $G$-universe $\univ$ \pcref{def:g_universe,def:IU_spaces}.

\begin{definition}\label{def:iuspace_day}
Given $\IU$-spaces $X,Y \cn \IU \to \Topgst$ \pcref{def:iu_space}, the \emph{smash product}\index{smash product!IU-space@$\IU$-space}\index{IU-space@$\IU$-space!smash product}
\begin{equation}\label{x_smau_y}
\IU \fto{X \smau Y} \Topgst
\end{equation}
is the $\IU$-space whose value at each object $U \in \IU$ is defined as the coend
\begin{equation}\label{smau_obj}
(X \smau Y)_U = \int^{(V,W) \in (\IUsk)^2} \IU(V \oplus W, U)_\splus \sma (X_V \sma Y_W)
\end{equation}
taken in $\Topst$, with $G$ acting diagonally on representatives.  For each linear isometric isomorphism $k \cn U \fiso U'$ in $\IU$, the pointed homeomorphism
\begin{equation}\label{smau_mor}
(X \smau Y)_U \fto[\iso]{(X \smau Y)_k} (X \smau Y)_{U'}
\end{equation}
is induced by post-composition with $k$,
\[\IU(V \oplus W, U) \fto[\iso]{k \comp -} \IU(V \oplus W, U'),\]
on representatives.  A more detailed description of $X \smau Y$ is given in the proof of \cref{smau_welldef}.
\end{definition}

We emphasize that, while the $\IU$-spaces $X$ and $Y$ are, by definition, $\Gtopst$-functors, the coend in \cref{smau_obj} is \emph{not} taken in $\Gtopst$.  The reason is that, for a typical isomorphism $f$ in $\IU$, the pointed homeomorphism $X_f$ \cref{iu_space_xf} is not required to be $G$-equivariant.  The computation \cref{smau_gact_welldef} below explains that the diagonal $G$-action on representatives gives a well-defined $G$-action on $(X \smau Y)_U$.

\begin{lemma}\label{smau_welldef}
The smash product $X \smau Y$ in \cref{x_smau_y} is an $\IU$-space.
\end{lemma}

\begin{proof}
Recall the description of an $\IU$-space in \cref{expl:iu_space}.

\parhead{Component objects}.  
The pointed space $(X \smau Y)_U$ in \cref{smau_obj} is a quotient of the wedge
\begin{equation}\label{smau_wedge}
\bigvee_{(V,W) \in (\IUsk)^2} \IU(V \oplus W, U)_\splus \sma (X_V \sma Y_W).
\end{equation}
A typical point in $(X \smau Y)_U$ is represented by a triple
\begin{equation}\label{x_smau_y_rep}
z = \big(V \oplus W \fto{e} U \in \IU ; x \in X_V ; y \in Y_W \big).
\end{equation}
The coend \cref{smau_obj} identifies, for each quintuple
\begin{equation}\label{smau_coend}
\bz = \scalebox{.9}{$\Big(V \fto{f} V' ; W \fto{h} W' ; V' \oplus W' \fto{e} U ; x \in X_V ; y \in Y_W \Big)$}
\end{equation}
with $f,h \in \IUsk$ and $e \in \IU$, the following two triples in the wedge \cref{smau_wedge}.
\begin{equation}\label{smau_iden}
\begin{split}
\dzero \bz &= \big(V \oplus W \fto{f \oplus h} V' \oplus W' \fto{e} U ; x ; y \big) \\
\done \bz &= \big(V' \oplus W' \fto{e} U ; X_f x \in X_{V'} ; Y_h y \in Y_{W'}  \big)
\end{split}
\end{equation}

\parhead{$G$-action}.
Each element $g \in G$ acts diagonally on each representative $z$ \cref{x_smau_y_rep}:
\begin{equation}\label{smau_gaction}
gz = \scalebox{.9}{$\Big(V \oplus W \fto{\ginv \oplus \ginv} V \oplus W \fto{e} U \fto{g} U ; gx ; gy \Big)$}.
\end{equation}
We let $g$ act diagonally on $\bz$ \cref{smau_coend}, so
\begin{equation}\label{g_bz}
g\bz = \big(gf\ginv ; gh\ginv ; ge(\ginv \oplus \ginv) ; gx ; gy \big).
\end{equation}
Using the functoriality of $\oplus$ on $\IU$, \cref{x_gfginv}, and \cref{smau_iden,smau_gaction,g_bz}, the following computation proves that the $G$-action on the pointed space $(X \smau Y)_U$ is well defined.
\begin{equation}\label{smau_gact_welldef}
\begin{split}
& g(\dzero\bz) \\
&= \big(ge(f \oplus h)(\ginv \oplus \ginv) ; gx ; gy \big) \\
&= \big(ge(\ginv \oplus \ginv) (gf\ginv \oplus gh\ginv) ; gx ; gy \big) \\
&= \dzero (g\bz) \\
&= \done (g\bz) \\
&= \big(ge(\ginv \oplus \ginv) ; X_{gf\ginv}(gx) ; Y_{gh\ginv}(gy) \big) \\
&= \big(ge(\ginv \oplus \ginv) ; (g X_f \ginv)(gx) ; (g Y_h \ginv)(gy) \big)\\
&= \big(ge(\ginv \oplus \ginv) ; g (X_f x) ; g (Y_h y) \big) \\
&= g(\done \bz)
\end{split}
\end{equation}

\parhead{Component morphisms}.
For each isomorphism $k \cn U \fiso U'$ in $\IU$, the pointed homeomorphism \cref{smau_mor} 
\begin{equation}\label{xsmauy_mor}
(X \smau Y)_U \fto[\iso]{(X \smau Y)_k} (X \smau Y)_{U'}
\end{equation}
sends each representative $z = (e; x; y)$ in \cref{x_smau_y_rep} to the triple
\begin{equation}\label{smau_iumor}
(X \smau Y)_k (z) = \big(V \oplus W \fto{e} U \fto{k} U' ; x ; y \big).
\end{equation}
Using \cref{smau_iden,smau_iumor}, the following computation shows that $(X \smau Y)_k$ is well defined.
\[\begin{split}
(X \smau Y)_k (\dzero \bz)
&= \big(ke(f \oplus h) ; x ; y \big)\\
&= \dzero (f; h ; ke ; x ; y) \\
&= \done (f; h ; ke ; x ; y) \\
&= \big( ke ; X_f x ; Y_h y) \\
&= (X \smau Y)_k (\done \bz)\\
\end{split}\]
Moreover, the definition \cref{smau_iumor} implies that $X \smau Y$ preserves composition and identity morphisms in the sense of \cref{iu_space_axioms}.

\parhead{Equivariance}.  For each isomorphism $k \cn U \fiso U'$ in $\IU$ and $g \in G$, the following computation, using \cref{x_smau_y_rep,smau_gaction,smau_iumor}, proves that $X \smau Y$ is $G$-equivariant in the sense of \cref{x_gfginv}.
\[\begin{split}
& \big(g (X \smau Y)_k \ginv \big)(e; x; y)\\
&= \big(g (X \smau Y)_k \big) \big(\ginv e (g \oplus g) ; \ginv x ; \ginv y \big) \\
&= g \big(k\ginv e (g \oplus g) ; \ginv x ; \ginv y \big) \\
&= \big(gk \ginv e (g \oplus g) (\ginv \oplus \ginv) ; g(\ginv x) ; g(\ginv y)\big) \\
&= (gk \ginv e ; x ; y) \\
&= (X \smau Y)_{gk\ginv}(e; x; y)
\end{split}\]
This finishes the proof that $X \smau Y$ is an $\IU$-space.
\end{proof}

\subsection*{Smash Product of $\IU$-Morphisms}
Next, we extend the smash product $\smau$ from $\IU$-spaces to $\IU$-morphisms.

\begin{definition}\label{def:iumor_day}
Consider two $\IU$-morphisms \cref{iu_mor} $\theta$ and $\theta'$ as follows.
\[\begin{tikzpicture}[vcenter]
\def\t{28}
\draw[0cell]
(0,0) node (a1) {\phantom{A}}
(a1)++(1.8,0) node (a2) {\phantom{A}}
(a1)++(-.1,0) node (a1') {\IU}
(a2)++(.21,-.04) node (a2') {\Topgst}
;
\draw[1cell=.9]
(a1) edge[bend left=\t] node {X} (a2)
(a1) edge[bend right=\t] node[swap] {Y} (a2)
;
\draw[2cell]
node[between=a1 and a2 at .45, rotate=-90, 2label={above,\theta}] {\Rightarrow}
;
\begin{scope}[shift={(3.5,0)}]
\draw[0cell]
(0,0) node (a1) {\phantom{A}}
(a1)++(1.8,0) node (a2) {\phantom{A}}
(a1)++(-.1,0) node (a1') {\IU}
(a2)++(.21,-.04) node (a2') {\Topgst}
;
\draw[1cell=.9]
(a1) edge[bend left=\t] node {X'} (a2)
(a1) edge[bend right=\t] node[swap] {Y'} (a2)
;
\draw[2cell]
node[between=a1 and a2 at .43, rotate=-90, 2label={above,\theta'}] {\Rightarrow}
;
\end{scope}
\end{tikzpicture}\]
The $\IU$-morphism\index{smash product!IU-morphism@$\IU$-morphism}\index{IU-morphism@$\IU$-morphism!smash product}
\begin{equation}\label{tha_smau_thap}
\begin{tikzpicture}[vcenter]
\def\t{20}
\draw[0cell]
(0,0) node (a1) {\phantom{A}}
(a1)++(3,0) node (a2) {\phantom{A}}
(a1)++(-.1,0) node (a1') {\IU}
(a2)++(.21,-.04) node (a2') {\Topgst}
;
\draw[1cell=.8]
(a1) edge[bend left=\t] node {X \smau X'} (a2)
(a1) edge[bend right=\t] node[swap] {Y \smau Y'} (a2)
;
\draw[2cell=.9]
node[between=a1 and a2 at .33, rotate=-90, 2label={above,\theta \smau \theta'}] {\Rightarrow}
;
\end{tikzpicture}
\end{equation}
is defined by, for each object $U \in \IU$, the $U$-component pointed morphism given by the following commutative diagram, where $\int = \int^{(V,W) \in (\IUsk)^2}$.
\begin{equation}\label{tha_smau_thap_U}
\begin{tikzpicture}[vcenter]
\def\v{-1.5}
\draw[0cell=.9]
(0,0) node (a11) {(X \smau X')_U}
(a11)++(4,0) node (a12) {\txint \IU(V \oplus W, U)_\splus \sma (X_V \sma X'_W)}
(a11)++(0,\v) node (a21) {(Y \smau Y')_U}
(a12)++(0,\v) node (a22) {\txint \IU(V \oplus W, U)_\splus \sma (Y_V \sma Y'_W)}
;
\draw[1cell=.8]
(a11) edge[equal] (a12)
(a21) edge[equal] (a22)
(a11) edge[transform canvas={xshift=1em}] node[swap] {(\theta \smau \theta')_U} (a21)
(a12) edge[transform canvas={xshift=-4.5em}] node {\txint \IU(V \oplus W, U)_\splus \sma (\theta_V \sma \theta'_W)} (a22)
;
\end{tikzpicture}
\end{equation}
\begin{itemize}
\item The naturality diagram \cref{iu_mor_natural} commutes for $\theta \smau \theta'$ by \cref{smau_iumor}, \cref{tha_smau_thap_U}, the functoriality of $\sma$, and the universal property of coends.
\item The functoriality of 
\begin{equation}\label{smau_functor}
\IUT \times \IUT \fto{\smau} \IUT
\end{equation}
follows from the definition \cref{tha_smau_thap_U}, the functoriality of $\sma$, and the universal property of coends.
\end{itemize}
This finishes the definition.  
\end{definition}

\begin{explanation}\label{expl:iumor_day}
The $U$-component $(\theta \smau \theta')_U$ in \cref{tha_smau_thap_U} sends a triple \cref{x_smau_y_rep}
\[z = \big(V \oplus W \fto{e} U \in \IU ; x \in X_V ; x' \in X'_W \big),\]
which represents a point in $(X \smau X')_U$, to the triple
\begin{equation}\label{tha_smau_thap_uz}
(\theta \smau \theta')_U(z) = \big(e; \theta_V x \in Y_V; \theta'_W x' \in Y'_W \big),
\end{equation}
which represents a point in $(Y \smau Y')_U$.
\end{explanation}

\subsection*{$\IU$-Spaces Form a Symmetric Monoidal $\Gtop$-Category}
The following definition uses the notion of a symmetric monoidal $\V$-category \pcref{definition:monoidal-vcat,definition:braided-monoidal-vcat,definition:symm-monoidal-vcat} for the Cartesian closed category $(\Gtop,\times,*)$ \pcref{def:Gtop}.

\begin{definition}\label{def:IU_smgtop}
Define a symmetric monoidal $\Gtop$-category\index{IU-space@$\IU$-space!symmetric monoidal category}\index{symmetric monoidal category!IU-space@$\IU$-space}
\[(\IUT, \smau, \iu, \au, \ellu, \ru, \beu)\]
as follows.
\begin{description}
\item[Base $\Gtop$-category] It is the $\Gtop$-category $\IUT$ in \cref{def:iut_gtop_enr}.
\begin{itemize}
\item Its objects are $\IU$-spaces \pcref{def:iu_space}.
\item Each hom $G$-space $\IUT(X,Y)$ has $\IU$-morphisms $X \to Y$ \pcref{def:iu_morphism} as its points, with $G$ acting componentwise by conjugation \cref{gtheta_v}.
\end{itemize}
\item[Monoidal composition] It is the $\Gtop$-functor
\begin{equation}\label{IU_mon_composition}
\IUT \otimes \IUT \fto{\smau} \IUT
\end{equation}
whose underlying functor is $\smau$ in \cref{smau_functor}.  Its object and underlying morphism assignments are given in, respectively, \cref{def:iuspace_day,def:iumor_day}.
\item[Monoidal identity] 
With $\vtensorunit$ denoting the unit $\Gtop$-category \pcref{definition:unit-vcat}, the monoidal identity is the $\Gtop$-functor
\begin{equation}\label{IU_mon_identity}
\vtensorunit \fto{\iu} \IUT
\end{equation}
determined by the $\IU$-space $\iu \cn \IU \to \Topgst$ that is defined on objects by
\begin{equation}\label{iu_mon_id}
\iu_V = \begin{cases}
\stplus & \text{if $V=0$, and}\\
* & \text{if $V \in \IU \setminus \{0\}$.}
\end{cases}
\end{equation}
\item[Left monoidal unitor] 
It is the $\Gtop$-natural isomorphism
\begin{equation}\label{IU_ellu}
\begin{tikzpicture}[vcenter]
\def\v{1.4} \def\h{.15}
\draw[0cell=.9]
(0,0) node (a11) {\IUT}
(a11)++(4,0) node (a12) {\IUT}
(a11)++(\h,\v) node (a01) {\vtensorunit \otimes \IUT}
(a12)++(-\h,\v) node (a02) {(\IUT)^{\otimes 2}}
;
\draw[1cell=.9]
(a11) edge node[swap] {1} (a12)
(a11) edge node[pos=.2] {(\ell^\otimes)^{-1}} (a01)
(a01) edge node {\iu \otimes 1} (a02)
(a02) edge node[pos=.7] {\smau} (a12)
;
\draw[2cell]
node[between=a11 and a12 at .47, shift={(0,.5*\v)}, rotate=-90, 2label={above,\ellu}] {\Rightarrow}
;
\end{tikzpicture}
\end{equation}
whose component at an $\IU$-space $X$ is a $G$-morphism
\[* \fto{\ellu_X} \IUT(\iu \smau X, X),\]
which means a $G$-equivariant $\IU$-isomorphism (\cref{def:iu_morphism} \pcref{def:iu_morphism_iii})
\begin{equation}\label{ellux}
\iu \smau X \fto[\iso]{\ellu_X} X.
\end{equation}
By \cref{expl:eqiu_morphism}, it suffices to define the $U$-component pointed $G$-homeomorphism $\ellu_{X; U}$ for each object $U \in \IUsk$.  This component is defined by the following commutative diagram, where $\int = \int^{(V,W) \in (\IUsk)^2}$.
\begin{equation}\label{ellu_comp}
\begin{tikzpicture}[vcenter]
\def\v{-1.4} \def\h{4.3}
\draw[0cell=.9]
(0,0) node (a1) {\IU(0 \oplus U,U)_\splus \sma (\stplus \sma X_U)}
(a1)++(\h,0) node (a0) {\IU(U,U)_\splus \sma X_U}
(a1)++(0,\v) node (a2) {\txint \IU(V \oplus W,U)_\splus \sma (\iu_V \sma X_W)}
(a2)++(0,.7*\v) node (a3) {(\iu \smau X)_U}
(a3)++(\h,0) node (a4) {\phantom{X_U}}
(a4)++(0,-.02) node (a4') {X_U} 
;
\draw[1cell=.8]
(a1) edge node[swap] {\iso} node {\ell} (a0)
(a1) edge node[swap] {\inc_{(0,U)}} (a2)
(a2) edge[equal] (a3)
(a3) edge node {\ellu_{X; U}} (a4)
(a0) edge node {X} (a4')
;
\end{tikzpicture}
\end{equation}
The arrows in \cref{ellu_comp} are defined as follows.
\begin{itemize}
\item The isomorphism $\ell$ uses the (inverse) left unit isomorphisms
\[U \iso 0 \oplus U \andspace \stplus \sma\, X_U \iso X_U\]
for the symmetric monoidal categories $\IUsk$ and $\Gtopst$ \pcref{def:gtopst,def:IU_spaces}.
\item $\inc_{(0,U)}$ is part of the definition of the coend at the pair $(0,U) \in (\IUsk)^2$.
\item The equality comes from the definition \cref{smau_obj} of the smash product.
\item The arrow labeled $X$ sends a pair
\[\big(U \fto[\iso]{f} U \in \IU; x \in X_U\big) \tospace X_f(x) \in X_U,\]
where $X_f$ is the $f$-component pointed homeomorphism of $X$ \cref{iu_space_xf}.
\end{itemize}
\item[Right monoidal unitor] The definition of its $X$-component $G$-equivariant $\IU$-isomorphism 
\begin{equation}\label{ru_x}
X \smau \iu \fto[\iso]{\ru_X} X
\end{equation}
is analogous to the left monoidal unitor $\ellu_X$ defined in \cref{ellu_comp}, using the right unit isomorphisms for the symmetric monoidal categories $\IUsk$ and $\Gtopst$.
\item[Monoidal associator] 
It is the $\Gtop$-natural isomorphism
\begin{equation}\label{IU_au}
\begin{tikzpicture}[vcenter]
\def\h{2.5} \def\u{1} \def\v{-1.4}
\draw[0cell=.9]
(0,0) node (a11) {(\IUT)^{\otimes 2} \otimes \IUT}
(a11)++(\h,\u) node (a12) {(\IUT)^{\otimes 2}}
(a12)++(\h,-\u) node (a13) {\IUT}
(a11)++(0,\v) node (a21) {\IUT \otimes (\IUT)^{\otimes 2}}
(a13)++(0,\v) node (a22) {(\IUT)^{\otimes 2}}
;
\draw[1cell=.8]
(a11) edge node {\smau \otimes 1} (a12)
(a12) edge node {\smau} (a13)
(a11) edge node[swap] {a^\otimes} (a21)
(a21) edge node {1 \otimes \smau} (a22)
(a22) edge node[swap] {\smau} (a13)
;
\draw[2cell]
node[between=a11 and a13 at .46, shift={(0,.3*\v)}, rotate=-90, 2label={above,\au}] {\Rightarrow}
;
\end{tikzpicture}
\end{equation}
whose component at a triple $(X,Y,Z)$ of $\IU$-spaces is a $G$-equivariant $\IU$-isomorphism 
\begin{equation}\label{au_xyz}
(X \smau Y) \smau Z \fto[\iso]{\au_{X,Y,Z}} X \smau (Y \smau Z).
\end{equation}
For each object $U \in \IU$, its $U$-component pointed $G$-homeomorphism 
\begin{equation}\label{au_xyzu}
\big((X \smau Y) \smau Z\big)_U \fto[\iso]{\au_{X,Y,Z; U}} \big(X \smau (Y \smau Z)\big)_U
\end{equation}
is defined via the pointed $G$-homeomorphisms
\begin{equation}\label{au_xyzu_dom}
\begin{split}
& \big((X \smau Y) \smau Z\big)_U \\
&= \txint^{(V,V_3) \in (\IUsk)^2} \IU(V \oplus V_3, U)_\splus \sma \big((X \smau Y)_V \sma Z_{V_3} \big) \\
&= \scalebox{.8}{$\txint^{(V,V_3)} \IU(V \oplus V_3, U)_\splus \sma \big[\big( \txint^{(V_1,V_2)} \IU(V_1 \oplus V_2, V)_\splus \sma (X_{V_1} \sma Y_{V_2}) \big) \sma Z_{V_3} \big]$}\\
&\iso \scalebox{.8}{$\txint^{(V_1,V_2,V_3)} \big[\txint^V \IU(V \oplus V_3, U)_\splus \sma \IU(V_1 \oplus V_2, V)_\splus \big] \sma \big( (X_{V_1} \sma Y_{V_2}) \sma Z_{V_3} \big)$} \\
&\iso \scalebox{.9}{$\txint^{(V_1,V_2,V_3) \in (\IUsk)^3} \IU\big((V_1 \oplus V_2) \oplus V_3, U\big)_\splus \sma \big((X_{V_1} \sma Y_{V_2}) \sma Z_{V_3} \big)$}
\end{split}
\end{equation}
and
\begin{equation}\label{au_xyzu_cod}
\begin{split}
& \big(X \smau (Y \smau Z)\big)_U \\
&\iso \scalebox{.9}{$\txint^{(V_1,V_2,V_3) \in (\IUsk)^3} \IU\big(V_1 \oplus (V_2 \oplus V_3), U\big)_\splus \sma \big(X_{V_1} \sma (Y_{V_2} \sma Z_{V_3}) \big)$}.
\end{split}
\end{equation}
The first pointed $G$-homeomorphism in \cref{au_xyzu_dom} uses
\begin{itemize}
\item the fact that $\IU(V \oplus V_3, U)_\splus \sma -$ commutes with coends,
\item the associativity isomorphism for $\sma$, and
\item the \index{Fubini Theorem}Fubini Theorem for coends \cite[1.3.1]{loregian}.
\end{itemize}
The second pointed $G$-homeomorphism in \cref{au_xyzu_dom} uses the \index{Yoneda Density Theorem}Yoneda Density Theorem for coends; see \cite[A.5.7]{loregian} or \cite[3.7.15]{cerberusIII}.  The pointed $G$-homeomorphism \cref{au_xyzu_cod} is proved in the same way.

Using \cref{au_xyzu_dom,au_xyzu_cod}, the pointed $G$-homeomorphism $\au_{X,Y,Z;U}$ in \cref{au_xyzu} is defined by the universal property of coends and the (inverse) associativity isomorphisms
\begin{equation}\label{au_xyz_def}
\begin{split}
V_1 \oplus (V_2 \oplus V_3) &\iso (V_1 \oplus V_2) \oplus V_3 \andspace\\
(X_{V_1} \sma Y_{V_2}) \sma Z_{V_3} &\iso X_{V_1} \sma (Y_{V_2} \sma Z_{V_3})
\end{split}
\end{equation}
for $\IUsk$ and $\Gtopst$.
\item[Braiding] It is the $\Gtop$-natural isomorphism
\begin{equation}\label{IU_braiding}
\begin{tikzpicture}[vcenter]
\def\h{4}
\draw[0cell]
(0,0) node (a1) {(\IUT)^{\otimes 2}}
(a1)++(.5*\h,-1) node (a2) {(\IUT)^{\otimes 2}}
(a1)++(\h,0) node (a3) {\IUT}
;
\draw[1cell=.9]
(a1) edge node {\smau} (a3)
(a1) [rounded corners=2pt] |- node[swap,pos=.25] {\beta^\otimes} (a2)
;
\draw[1cell=.9]
(a2) [rounded corners=2pt] -| node[swap,pos=.75] {\smau} (a3)
;
\draw[2cell]
node[between=a1 and a3 at .5, shift={(0,-.45)}, rotate=-90, 2label={above,\beu}] {\Rightarrow}
;
\end{tikzpicture}
\end{equation}
whose component at a pair $(X,Y)$ of $\IU$-spaces is a $G$-equivariant $\IU$-isomorphism  
\begin{equation}\label{beu_xy}
X \smau Y \fto[\iso]{\beu_{X,Y}} Y \smau X.
\end{equation}
For each object $U \in \IU$, its $U$-component pointed $G$-homeomorphism 
\begin{equation}\label{beu_xyu}
(X \smau Y)_U \fto[\iso]{\beu_{X,Y; U}} (Y \smau X)_U
\end{equation}
sends a representative triple $z = (e; x ; y)$ in \cref{x_smau_y_rep} to the triple
\begin{equation}\label{beu_XYU_z}
\beu_{X,Y; U}(z) = \big(W \oplus V \fto[\iso]{\xi_{W,V}} V \oplus W \fto{e} U ; y ; x \big),
\end{equation}
which represents a point in the coend
\[(Y \smau X)_U = \int^{(W,V) \in (\IUsk)^2} \IU(W \oplus V, U)_\splus \sma (Y_W \sma X_V).\]
In \cref{beu_XYU_z}, $\xi_{W,V}$ denotes the $(W,V)$-component of the braiding for the symmetric monoidal category $(\IUsk,\oplus)$.  In the last two arguments, going from 
\[(x;y) \in X_V \sma Y_W \tospace (y;x) \in Y_W \sma X_V\]
uses the braiding for the symmetric monoidal category $(\Gtopst,\sma)$.
\end{description}
This finishes the definition of the symmetric monoidal $\Gtop$-category $\IUT$.  The underlying symmetric monoidal category of $\IUT$ is denoted by the same notation.
\end{definition}

\begin{proposition}\label{iut_smgtop}
$\IUT$ in \cref{def:IU_smgtop} is a symmetric monoidal $\Gtop$-category.
\end{proposition}

\begin{proof}
We first show that the data
\[(\smau, \iu, \au, \ellu, \ru, \beu)\]
are well defined.  Then we show that the symmetric monoidal $\Gtop$-category axioms are satisfied.

\parhead{Monoidal unit}. The monoidal unit $\iu$ in \cref{IU_mon_identity} is a $\Gtop$-functor because the componentwise conjugation $G$-action \cref{gtheta_v} fixes each identity $\IU$-morphism.  In particular, the identity $\IU$-morphism of the $\IU$-space $\iu$ defined in \cref{iu_mon_id} is $G$-fixed.

\parhead{Monoidal composition}. 
To prove that the monoidal composition $\smau$ in \cref{IU_mon_composition} is enriched in $\Gtop$, and not just $\Top$, we need to show that its components are $G$-equivariant.  In other words, given two $\IU$-morphisms \pcref{def:iumor_day}
\[X \fto{\theta} Y \andspace X' \fto{\theta'} Y'\]
and an element $g \in G$, we need to show that the following two $\IU$-morphisms are equal.
\begin{equation}\label{smau_g_equivariant}
\begin{tikzpicture}[baseline={(a1.base)}]
\draw[0cell]
(0,0) node (a1) {X \smau X'}
(a1)++(4.2,0) node (a2) {Y \smau Y'}
;
\draw[1cell=.9]
(a1) edge[transform canvas={yshift=.4ex}] node {(g \cdot \theta) \smau (g \cdot \theta')} (a2)
(a1) edge[transform canvas={yshift=-.5ex}] node[swap] {g \cdot (\theta \smau \theta')} (a2)
;
\end{tikzpicture}
\end{equation}
Using \cref{gtheta_v,smau_gaction,tha_smau_thap_uz}, the following computation proves that the two $\IU$-morphisms in \cref{smau_g_equivariant} are equal at each object $U \in \IU$ and each triple \cref{x_smau_y_rep}
\[\big(V \oplus W \fto{e} U \in \IU;  x \in X_V ; x' \in X'_W \big),\]
which represents a point in $(X \smau X')_U$.
\begin{equation}\label{smau_g_computation}
\begin{split}
& \big((g \cdot \theta) \smau (g \cdot \theta') \big)_U (e; x; x') \\
&= \big(e ; (g \cdot \theta)_V(x) ; (g \cdot \theta')_W(x') \big) \\
&= \big(e ; g\theta_V(\ginv x) ; g\theta'_W(\ginv x')  \big) \\
&= \big(g \ginv e (g \oplus g) (\ginv \oplus \ginv) ; g\theta_V(\ginv x) ; g\theta'_W(\ginv x') \big) \\
&= g\big(\ginv e (g \oplus g) ; \theta_V(\ginv x) ; \theta'_W(\ginv x') \big) \\
&= g\big((\theta \smau \theta')_U \big(\ginv e (g \oplus g) ; \ginv x; \ginv x' \big) \big) \\
&= g\big((\theta \smau \theta')_U \big(\ginv(e; x; x')\big) \big) \\
&= \big(g \cdot (\theta \smau \theta')\big)_U (e; x; x')
\end{split}
\end{equation}

\parhead{Braiding, unitors, and associator}. 
The assignment $\beu$, as defined in \cref{IU_braiding}, is a $\Gtop$-natural isomorphism \pcref{def:enriched-natural-transformation} for the following reasons, along with the universal property of coends.
\begin{itemize}
\item Using the naturality of the braiding $\xi$ for $\IUsk$, \cref{smau_coend,smau_iden,beu_XYU_z}, the following computation proves that each component $\beu_{X,Y; U}$ is well defined.
\[\begin{split}
& \beu_{X,Y;U} (\dzero \bz) \\
&= \beu_{X,Y;U} \big(e(f \oplus h) ; x ; y\big) \\
&= \big(e(f \oplus h)\xi_{W,V} ; y ; x\big) \\
&= \big(e \xi_{W',V'} (h \oplus f); y ; x\big) \\
&= \dzero\big(h; f; e\xi_{W',V'} ; y ; x \big) \\
&= \done\big(h; f; e\xi_{W',V'} ; y ; x \big) \\
&= \big(e \xi_{W',V'} ; Y_h y ; X_f x \big) \\
&= \beu_{X,Y;U} (e ; X_f x ; Y_h y) \\
&= \beu_{X,Y;U} (\done \bz)
\end{split}\]
\item The naturality of $\beu_{X,Y; U}$ in the variable $U \in \IU$, in the sense of \cref{iu_mor_natural}, follows from \cref{smau_iumor,beu_XYU_z}.
\item Under the componentwise conjugation $G$-action \cref{gtheta_v}, the condition that the point
\[\beu_{X,Y} \in \IUT(X \smau Y, Y \smau X)\]
is $G$-fixed means that each component $\beu_{X,Y; U}$ is $G$-equivariant.  Using \cref{smau_gaction}, \cref{beu_XYU_z}, and the naturality of the braidings for $\IUsk$ and $\Gtopst$, the following computation proves that $\beu_{X,Y; U}$ is $G$-equivariant.
\begin{equation}\label{beu_xyu_geq}
\begin{split}
& \beu_{X,Y;U} \big(g(e; x; y)\big) \\
&= \beu_{X,Y;U} \big(ge(\ginv \oplus \ginv) ; gx ; gy \big) \\
&= \big(ge(\ginv \oplus \ginv) \xi_{W,V} ; gy ; gx \big) \\ 
&= \big(ge \xi_{W,V} (\ginv \oplus \ginv) ; gy ; gx \big) \\
&= g\big(e \xi_{W,V} ; y ; x \big) \\
&= g (\beu_{X,Y;U}(e; x; y))
\end{split}
\end{equation}
\item The naturality of $\beu_{X,Y}$ in the variables $(X,Y) \in (\IUT)^{\otimes 2}$, in the sense of \cref{enr_naturality}, follows from \cref{tha_smau_thap_uz}, \cref{beu_XYU_z}, and the naturality of the braiding for $\Gtopst$.
\end{itemize}
This proves that $\beu$ is a $\Gtop$-natural isomorphism.  Analogously, each of $\ellu$, $\ru$, and $\au$ is a $\Gtop$-natural isomorphism.

\parhead{Axioms}. 
Each symmetric monoidal $\Gtop$-category axiom in \cref{definition:monoidal-vcat,definition:braided-monoidal-vcat,definition:symm-monoidal-vcat} holds for $\IUT$ by the universal property of coends and the symmetric monoidal category axioms for $(\IUsk,\oplus)$ and $(\Gtopst,\sma)$.  

As an example, consider the symmetry axiom \cref{vmonoidal-symmetry}, which states that the following diagram of $\IU$-isomorphisms commutes for each pair $(X,Y)$ of $\IU$-spaces.
\[\begin{tikzpicture}[vcenter]
\def\h{2.5} \def\u{.65}
\draw[0cell]
(0,0) node (a1) {X \smau Y}
(a1)++(\h,0) node (a2) {Y \smau X}
(a2)++(\h,0) node (a3) {X \smau Y}
;
\draw[1cell=.9]
(a1) edge node {\beu_{X,Y}} (a2)
(a2) edge node {\beu_{Y,X}} (a3)
(a1)  [rounded corners=2pt] |- ($(a2)+(-1,\u)$) -- node {1} ($(a2)+(1,\u)$) -| (a3)
;
\end{tikzpicture}\]
Each $\IU$-morphism is determined by its component pointed morphisms \cref{iu_mor_component}.  Thus, by the universal property of coends \cref{smau_obj}, it suffices to prove the commutativity of the previous diagram at the $U$-component for an object $U \in \IU$ and a representative triple $z = (e; x ; y)$ in \cref{x_smau_y_rep}.  This is proved using \cref{beu_XYU_z} as follows.
\begin{equation}\label{iut_symmetry_axiom}
\begin{split}
& \beu_{Y,X;U} \beu_{X,Y;U} (e; x; y)\\
&= \beu_{Y,X;U} \big(e\xi_{W,V} ; y ; x \big) \\
&= \big(e\xi_{W,V} \xi_{V,W} ; x ; y \big) \\
&= (e; x; y)
\end{split}
\end{equation}
The previous computation uses the symmetry axiom \cref{symmoncatsymhexagon} for the symmetric monoidal categories
\begin{itemize}
\item $(\Gtopst,\sma)$, applied to the pair $(x,y) \in X_V \sma Y_W$, and
\item $(\IUsk,\oplus)$, which yields $\xi_{W,V} \xi_{V,W} = 1_{V \oplus W}$.
\end{itemize}
This proves the symmetry axiom \cref{vmonoidal-symmetry} for $\IUT$.  

Every other symmetric monoidal $\Gtop$-category axiom \pcref{definition:monoidal-vcat,definition:braided-monoidal-vcat} for $\IUT$ is proved in an analogous way, by reducing it to the corresponding axioms for the symmetric monoidal categories $\Gtopst$ and $\IUsk$.
\end{proof}

\subsection*{Symmetric Monoidal $\Top$-Category}

Recall the category $\GIUT$ \cref{GIUT} of $\IU$-spaces and $G$-equivariant $\IU$-morphisms.  It has the structure of a $\Top$-category \pcref{def:giut_top_enr}, which is obtained from the $\Gtop$-category $\IUT$ by passing to $G$-fixed points of hom $G$-spaces \cref{giutxy_fixedpt}.

\begin{definition}\label{def:giut_smtop}
Passing to $G$-fixed points of each hom $G$-space, the symmetric monoidal $\Gtop$-category $\IUT$ \pcref{def:IU_smgtop} yields the symmetric monoidal $\Top$-category \pcref{definition:symm-monoidal-vcat}
\[(\GIUT, \smau, \iu, \au, \ellu, \ru, \beu).\]
Its underlying symmetric monoidal category is denoted by the same notation.
\end{definition}

\begin{explanation}[Unraveling $\GIUT$]\label{expl:giut_smtop}
The symmetric monoidal $\Top$-category $\GIUT$ is specified by the following data.

\parhead{Monoidal identity}.  It is the $\Top$-functor
\begin{equation}\label{giut_mon_id}
\vtensorunit \fto{\iu} \GIUT
\end{equation}
specified by the $\IU$-space $\iu$ in \cref{iu_mon_id}.  

\parhead{Monoidal composition}.  It is the $\Top$-functor
\begin{equation}\label{giut_mon_comp}
\GIUT \otimes \GIUT \fto{\smau} \GIUT
\end{equation}
whose object assignment is given by the smash product
\[(X,Y) \mapsto X \smau Y\] 
in \cref{def:iuspace_day}.  Its underlying morphism assignment is given by the smash product
\[(\theta, \theta') \mapsto \theta \smau \theta'\] 
in \cref{def:iumor_day}, applied to $G$-equivariant $\IU$-morphisms.  The computation \cref{smau_g_computation} shows that, if $\theta$ and $\theta'$ are $G$-equivariant $\IU$-morphisms, then so is $\theta \smau \theta'$.

\parhead{Monoidal associator, unitors, and braiding}.  
For each of $\au$, $\ellu$, $\ru$, and $\beu$, the components are the same as those for $\IUT$ in \cref{def:IU_smgtop}.  For example, the left monoidal unitor for $\GIUT$ is the $\Top$-natural isomorphism
\begin{equation}\label{giut_ellu}
\begin{tikzpicture}[vcenter]
\def\v{1.4} \def\h{.15}
\draw[0cell=.9]
(0,0) node (a11) {\GIUT}
(a11)++(4,0) node (a12) {\GIUT}
(a11)++(\h,\v) node (a01) {\vtensorunit \otimes \GIUT}
(a12)++(-\h,\v) node (a02) {(\GIUT)^{\otimes 2}}
;
\draw[1cell=.9]
(a11) edge node[swap] {1} (a12)
(a11) edge node[pos=.2] {(\ell^\otimes)^{-1}} (a01)
(a01) edge node {\iu \otimes 1} (a02)
(a02) edge node[pos=.7] {\smau} (a12)
;
\draw[2cell]
node[between=a11 and a12 at .47, shift={(0,.5*\v)}, rotate=-90, 2label={above,\ellu}] {\Rightarrow}
;
\end{tikzpicture}
\end{equation}
whose component at an $\IU$-space $X$,
\[* \fto{\ellu_X} \GIUT(\iu \smau X, X) = \big(\IUT(\iu \smau X, X) \big)^G,\]
is given by the $G$-equivariant $\IU$-isomorphism  
\[\iu \smau X \fto[\iso]{\ellu_X} X\]
in \cref{ellux}.

\parhead{Axioms}.  Each symmetric monoidal $\Top$-category axiom  \pcref{definition:monoidal-vcat,definition:braided-monoidal-vcat,definition:symm-monoidal-vcat} holds for $\GIUT$ by the same axiom for the symmetric monoidal $\Gtop$-category $\IUT$.  For example, the computation \cref{iut_symmetry_axiom} also proves the symmetry axiom \cref{vmonoidal-symmetry} for $\GIUT$.
\end{explanation}

\section{Commutative $G$-Monoids and Modules}
\label{sec:gmonoid}

This section reviews commutative $G$-monoids, modules, and their morphisms in the context of $\IU$-spaces.  As before, $G$ denotes a compact Lie group, and $\univ$ is a complete $G$-universe \pcref{def:g_universe}.

\secoutline
\begin{itemize}
\item \cref{def:com_g_monoid} defines (commutative) $G$-monoids in $\IUT$, along with their modules, morphisms, and $G$-action on morphisms.
\item \cref{expl:com_g_monoid,expl:gmonoid_module,expl:gmonoid_mod_mor} unpack the definitions of, respectively, a (commutative) $G$-monoid in $\IUT$, a module over a commutative $G$-monoid, and a morphism of modules. 
\end{itemize}

Recall (commutative) monoids and modules in a (symmetric) monoidal category \pcref{def:monoid,def:modules}.

\begin{definition}\label{def:com_g_monoid}
Consider the underlying symmetric monoidal category of $\GIUT$ in \cref{def:giut_smtop}.
\begin{enumerate}
\item\label{def:com_g_monoid_i} A \emph{$G$-monoid}\index{G-monoid@$G$-monoid} in $\IUT$ means a monoid in $\GIUT$.
\item\label{def:com_g_monoid_ii} A \emph{commutative $G$-monoid}\index{commutative G-monoid@commutative $G$-monoid} in $\IUT$ means a commutative monoid in $\GIUT$. 
\item\label{def:com_g_monoid_iii} For a commutative $G$-monoid $A$ in $\IUT$, an \emph{$A$-module}\index{module} means a right $A$-module in $\GIUT$.
\item\label{def:com_g_monoid_iv} For a commutative $G$-monoid $A$ in $\IUT$ and $A$-modules $(X,\umu^X)$ and $(Y,\umu^Y)$, a \emph{morphism of $A$-modules}\index{module!morphism}
\[(X,\umu^X)\fto{\theta} (Y,\umu^Y)\]
is an $\IU$-morphism $\theta \cn X \to Y$ \cref{iu_mor} such that the following diagram in $\IUT$ commutes.
\begin{equation}\label{amodule_mor}
\begin{tikzpicture}[vcenter]
\def\v{-1.3}
\draw[0cell]
(0,0) node (a11) {X \smau A}
(a11)++(2,0) node (a12) {X}
(a11)++(0,\v) node (a21) {Y \smau A}
(a12)++(0,\v) node (a22) {Y}
;
\draw[1cell=.9]
(a11) edge node {\umu^X} (a12)
(a12) edge node {\theta} (a22)
(a11) edge node[swap] {\theta \smau 1_A} (a21)
(a21) edge node {\umu^Y} (a22)
;
\end{tikzpicture}
\end{equation}
Note that $\theta$ is \emph{not} required to be $G$-equivariant.  The category of $A$-modules and morphisms of $A$-modules is denoted by $\AMod$.
\item\label{def:com_g_monoid_v} In the context of \pcref{def:com_g_monoid_iv}, the set of $A$-module morphisms $(X,\umu^X) \to (Y,\umu^Y)$ is regarded as a $G$-subspace of the $G$-space $\IUT(X,Y)$ in \cref{IUTXY_topology}.  Equipped with these hom $G$-spaces, $\AMod$ becomes a $\Gtop$-category.
\end{enumerate}
This finishes the definition.
\end{definition}

\cref{def:com_g_monoid} is discussed further in \cref{expl:com_g_monoid,expl:gmonoid_module,expl:gmonoid_mod_mor} below.

\begin{explanation}[Commutative $G$-Monoids]\label{expl:com_g_monoid}
Interpreting \cref{def:monoid} for the symmetric monoidal category 
\pcref{def:giut_smtop}
\[(\GIUT, \smau, \iu, \au, \ellu, \ru, \beu),\]
a \emph{$G$-monoid} in $\IUT$ (\cref{def:com_g_monoid} \pcref{def:com_g_monoid_i}) is a triple $(A,\mu,\eta)$ consisting of
\begin{itemize}
\item an $\IU$-space $A \cn \IU \to \Topgst$ \pcref{def:iu_space},
\item a $G$-equivariant $\IU$-morphism (\cref{def:iu_morphism} \pcref{def:iu_morphism_iii})
\begin{equation}\label{gmonoid_unit}
\iu \fto{\eta} A,
\end{equation}
called the \emph{unit}, and
\item a $G$-equivariant $\IU$-morphism 
\begin{equation}\label{gmonoid_mult}
A \smau A \fto{\mu} A,
\end{equation}
called the \emph{multiplication}.
\end{itemize}
These data are required to make the following associativity and unity diagrams in $\GIUT$ commute.
\begin{equation}\label{gmonoid_axioms}
\begin{tikzcd}[column sep=normal]
(A \smau A) \smau A \arrow{dd}[swap]{\mu \smau 1} \arrow[r, "\au", "\iso"']
& A \smau (A \smau A) \dar{1 \smau \mu}\\ 
& A \smau A \dar{\mu}\\  
A \smau A \arrow{r}{\mu} & A
\end{tikzcd}\qquad
\begin{tikzcd}[column sep=normal]
\iu \smau A \ar{d}[swap]{\eta \smau 1} \arrow[r, "\ellu", "\iso"'] & A \ar[equal]{d}\\ 
A \smau A \ar{r}{\mu} & A \ar[equal]{d}\\
A \smau \iu \ar{u}{1 \smau \eta} \arrow[r, "\ru", "\iso"'] & A
\end{tikzcd}
\end{equation}
Moreover, $(A,\mu,\eta)$ is a \emph{commutative $G$-monoid} in $\IUT$ (\cref{def:com_g_monoid} \pcref{def:com_g_monoid_ii}) if the following symmetry diagram in $\GIUT$ commutes.
\begin{equation}\label{gmonoid_symmetry}
\begin{tikzpicture}[vcenter]
\def\h{2.5}
\draw[0cell]
(0,0) node (a1) {A \smau A}
(a1)++(\h,0) node (a2) {A \smau A}
(a1)++(.5*\h,-1) node (a3) {A}
;
\draw[1cell=.9]
(a1) edge node {\beu} node[swap] {\iso} (a2)
(a2) edge node[pos=.4] {\mu} (a3)
(a1) edge node[swap,pos=.4] {\mu} (a3)
;
\end{tikzpicture}
\end{equation}
These data and axioms are further unpacked below.

\parhead{Unit}.  By \cref{expl:eqiu_morphism,iu_mon_id}, the unit $\eta \cn \iu \to A$ in \cref{gmonoid_unit} is given by a pointed $G$-morphism between pointed $G$-spaces
\begin{equation}\label{gmonoid_unit_comp}
\iu_0 = \stplus \fto{\eta_0} A_0,
\end{equation}
which means a $G$-fixed point in $A_0$.

\parhead{Multiplication}.  By \cref{smau_obj}, the domain of the multiplication $\mu \cn A \smau A \to A$ in \cref{gmonoid_mult} has, for each object $U \in \IU$, $U$-component pointed space given by the coend
\begin{equation}\label{gmonoid_mult_dom}
(A \smau A)_U = \int^{(V,W) \in (\IUsk)^2} \IU(V \oplus W, U)_\splus \sma (A_V \sma A_W)
\end{equation}
taken in $\Topst$, with $G$ acting diagonally on representatives \cref{smau_gaction}.  By the universal property of coends, the multiplication $\mu$ is determined by $(V,W)$-component pointed $G$-morphisms
\begin{equation}\label{gmonoid_mult_vw}
A_V \sma A_W \fto{\mu_{V,W}} A_{V \oplus W} \forspace (V,W) \in (\IUsk)^2
\end{equation}
such that, for each pair of linear isometric isomorphisms
\[V \fto[\iso]{f} V' \andspace W \fto[\iso]{h} W' \inspace \IUsk,\]
the following diagram of pointed morphisms commutes.
\begin{equation}\label{gmonoid_mult_nat}
\begin{tikzpicture}[vcenter]
\def\v{-1.4}
\draw[0cell]
(0,0) node (a11) {A_V \sma A_W}
(a11)++(3,0) node (a12) {A_{V \oplus W}}
(a11)++(0,\v) node (a21) {A_{V'} \sma A_{W'}}
(a12)++(0,\v) node (a22) {A_{V' \oplus W'}}
;
\draw[1cell=.9]
(a11) edge node {\mu_{V,W}} (a12)
(a12) edge node {A_{f \oplus h}} node[swap] {\iso} (a22)
(a11) edge node[swap] {A_f \sma A_h} node {\iso} (a21)
(a21) edge node {\mu_{V',W'}} (a22)
;
\end{tikzpicture}
\end{equation}
The arrows in the diagram \cref{gmonoid_mult_nat} have the following properties.
\begin{itemize}
\item $\mu_{V,W}$ and $\mu_{V',W'}$ are pointed $G$-morphisms between pointed $G$-spaces.
\item $A_f$, $A_h$, and $A_{f \oplus h}$ are pointed homeomorphisms, but \emph{not} generally $G$-equivariant \cref{iu_space_xf}.  Instead, since $A \cn \IU \to \Topgst$ is an $\IU$-space, each of $A_f$, $A_h$, and $A_{f \oplus h}$ satisfies the equivariance property in \cref{x_gfginv}.
\end{itemize}

\parhead{Unity}.  By the universal property of coends, \cref{smau_obj,iu_mon_id,ellu_comp,gmonoid_unit_comp,gmonoid_mult_vw}, the unity diagram in \cref{gmonoid_axioms} commutes if and only if the following diagram in $\Gtopst$ commutes for each object $U \in \IUsk$.
\begin{equation}\label{gmonoid_unity}
\begin{tikzpicture}[vcenter]
\def\h{3} \def\v{-1.2}
\draw[0cell]
(0,0) node (a11) {\stplus \sma A_U}
(a11)++(\h,0) node (a12) {A_U \sma \stplus}
(a11)++(-\h/2,\v) node (a21) {A_0 \sma A_U}
(a21)++(\h,0) node (a22) {A_U}
(a22)++(\h,0) node (a23) {A_U \sma A_0}
(a21)++(\h/2,\v) node (a31) {A_{0 \oplus U}}
(a31)++(\h,0) node (a32) {A_{U \oplus 0}}
;
\draw[1cell=.9]
(a11) edge node {\ell} node[swap] {\iso} (a22)
(a12) edge node[swap] {\rho} node {\iso} (a22)
(a11) edge node[swap] {\eta_0 \sma 1} (a21)
(a21) edge node[swap] {\mu_{0,U}} (a31)
(a31) edge node {\iso} node[swap] {A_\ell} (a22)
(a12) edge node {1 \sma \eta_0} (a23)
(a23) edge node {\mu_{U,0}} (a32)
(a32) edge node {A_\rho} node[swap] {\iso} (a22)
;
\end{tikzpicture}
\end{equation}
\begin{itemize}
\item In \cref{gmonoid_unity}, the top arrows $\ell$ and $\rho$ are, respectively, the left and right unit isomorphisms for the symmetric monoidal category $(\Gtopst,\sma,\stplus)$ \cref{Gtopst_smc}.
\item Along the bottom of \cref{gmonoid_unity}, the arrows 
\[\begin{tikzpicture}
\def\h{2}
\draw[0cell]
(0,0) node (a1) {0 \oplus U}
(a1)++(\h,0) node (a2) {U}
(a2)++(\h,0) node (a3) {U \oplus 0}
;
\draw[1cell=.9]
(a1) edge node {\ell} node[swap] {\iso} (a2)
(a3) edge node {\iso} node[swap] {\rho} (a2)
;
\end{tikzpicture}\]
are, respectively, the left and right unit isomorphisms for the symmetric monoidal category $(\IUsk,\oplus,0)$ \pcref{def:IU_spaces}, which are $G$-equivariant.  Their images under $A$ \cref{iu_space_xf} are the $G$-morphisms $A_\ell$ and $A_\rho$.
\end{itemize}
\parhead{Associativity}.  By the universal property of coends, \cref{smau_obj,au_xyzu,gmonoid_mult_vw}, the associativity diagram in \cref{gmonoid_axioms} commutes if and only if the following diagram in $\Gtopst$ commutes for each triple of objects $(U,V,W) \in (\IUsk)^3$.
\begin{equation}\label{gmonoid_assoc}
\begin{tikzpicture}[vcenter]
\def\h{4} \def\v{-1.4}
\draw[0cell]
(0,0) node (a11) {(A_U \sma A_V) \sma A_W}
(a11)++(\h,0) node (a12) {A_U \sma (A_V \sma A_W)}
(a11)++(0,\v) node (a21) {A_{U \oplus V} \sma A_W}
(a12)++(0,\v) node (a22) {A_U \sma A_{V \oplus W}}
(a21)++(0,\v) node (a31) {A_{(U \oplus V) \oplus W}}
(a31)++(\h,0) node (a32) {A_{U \oplus (V \oplus W)}}
;
\draw[1cell=.9]
(a11) edge node {\al} node[swap] {\iso} (a12)
(a12) edge node {1 \sma \mu_{V,W}} (a22)
(a22) edge node {\mu_{U,V \oplus W}} (a32)
(a11) edge node[swap] {\mu_{U,V} \sma 1} (a21)
(a21) edge node[swap] {\mu_{U \oplus V,W}} (a31)
(a31) edge node {A_\al} node[swap] {\iso} (a32)
;
\end{tikzpicture}
\end{equation}
In the diagram \cref{gmonoid_assoc}, the top horizontal arrow $\al$ is the associativity isomorphism for $(\Gtopst,\sma)$.  In the bottom horizontal arrow, 
\[(U \oplus V) \oplus W \fto[\iso]{\al} U \oplus (V \oplus W)\]
is the associativity isomorphism for $(\IUsk,\oplus)$, which is $G$-equivariant, and $A_\al$ is its image under $A$ \cref{iu_space_xf}.

\parhead{Symmetry}.  By the universal property of coends, \cref{smau_obj,beu_xyu,gmonoid_mult_vw}, the symmetry diagram \cref{gmonoid_symmetry} commutes if and only if the following diagram in $\Gtopst$ commutes for each pair of objects $(V,W) \in (\IUsk)^2$.
\begin{equation}\label{gmonoid_sym}
\begin{tikzpicture}[vcenter]
\def\v{1.4} 
\draw[0cell=1]
(0,0) node (a11) {A_V \sma A_W}
(a11)++(3.3,0) node (a12) {A_{V \oplus W}}
(a11)++(0,\v) node (a01) {A_W \sma A_V}
(a12)++(0,\v) node (a02) {A_{W \oplus V}}
;
\draw[1cell=.9]
(a11) edge node {\mu_{V,W}} (a12)
(a11) edge node {\xi} node[swap] {\iso} (a01)
(a01) edge node {\mu_{W,V}} (a02)
(a02) edge node {A_\xi} node[swap] {\iso} (a12)
;
\end{tikzpicture}
\end{equation}
\begin{itemize}
\item The left arrow $\xi$ in \cref{gmonoid_sym} is the braiding for the symmetric monoidal category $\Gtopst$.
\item In the right arrow in \cref{gmonoid_sym}, 
\[W \oplus V \fto[\iso]{\xi} V \oplus W\]
is the braiding for $\IUsk$, which is $G$-equivariant, and $A_\xi$ is its image under $A$ \cref{iu_space_xf}.
\end{itemize}
This finishes the unraveling of a (commutative) $G$-monoid in $\IUT$.
\end{explanation}

\begin{explanation}[Modules]\label{expl:gmonoid_module}
Interpreting \cref{def:modules} for the symmetric monoidal category 
\pcref{def:giut_smtop}
\[(\GIUT, \smau, \iu, \au, \ellu, \ru, \beu),\]
for a commutative $G$-monoid $(A,\mu,\eta)$ in $\IUT$, an \emph{$A$-module} (\cref{def:com_g_monoid} \pcref{def:com_g_monoid_iii}) is a pair $(X,\umu)$ consisting of
\begin{itemize}
\item an $\IU$-space $X \cn \IU \to \Topgst$ \pcref{def:iu_space} and
\item a $G$-equivariant $\IU$-morphism (\cref{def:iu_morphism} \pcref{def:iu_morphism_iii})
\begin{equation}\label{gmonoid_module_action}
X \smau A \fto{\umu} X,
\end{equation}
called the \emph{right $A$-action},
\end{itemize}
such that the following associativity and unity diagrams in $\GIUT$ commute.
\begin{equation}\label{gmonoid_module_axioms}
\begin{tikzcd}[column sep=normal]
(X \smau A) \smau A \arrow{dd}[swap]{\umu \smau 1} \arrow[r, "\au", "\iso"']
& X \smau (A \smau A) \dar{1 \smau \mu}\\ 
& X \smau A \dar{\umu}\\  
X \smau A \arrow{r}{\umu} & X
\end{tikzcd}\qquad
\begin{tikzcd}[column sep=normal]
X \smau A \ar{r}{\umu} & X \ar[equal]{d}\\
X \smau \iu \ar{u}{1 \smau \eta} \arrow[r, "\ru", "\iso"'] & X
\end{tikzcd}
\end{equation}
These data and axioms are further unpacked below.

\parhead{Right $A$-action}.  Similar to \cref{gmonoid_mult_vw}, the right $A$-action $\umu \cn X \smau A \to X$ in \cref{gmonoid_module_action} is determined by $(V,W)$-component pointed $G$-morphisms
\begin{equation}\label{gmonoid_mod_vw}
X_V \sma A_W \fto{\umu_{V,W}} X_{V \oplus W} \forspace (V,W) \in (\IUsk)^2
\end{equation}
such that, for each pair of linear isometric isomorphisms
\[V \fto[\iso]{f} V' \andspace W \fto[\iso]{h} W' \inspace \IUsk,\]
the following diagram of pointed morphisms commutes.
\begin{equation}\label{gmonoid_mod_nat}
\begin{tikzpicture}[vcenter]
\def\v{-1.4}
\draw[0cell]
(0,0) node (a11) {X_V \sma A_W}
(a11)++(3,0) node (a12) {X_{V \oplus W}}
(a11)++(0,\v) node (a21) {X_{V'} \sma A_{W'}}
(a12)++(0,\v) node (a22) {X_{V' \oplus W'}}
;
\draw[1cell=.9]
(a11) edge node {\umu_{V,W}} (a12)
(a12) edge node {X_{f \oplus h}} node[swap] {\iso} (a22)
(a11) edge node[swap] {X_f \sma A_h} node {\iso} (a21)
(a21) edge node {\umu_{V',W'}} (a22)
;
\end{tikzpicture}
\end{equation}
\begin{itemize}
\item $\umu_{V,W}$ and $\umu_{V',W'}$ are pointed $G$-morphisms between pointed $G$-spaces.
\item $X_f$, $A_h$, and $X_{f \oplus h}$ are pointed homeomorphisms that satisfy the equivariance property in \cref{x_gfginv}.
\end{itemize}

\parhead{Unity}.  Similar to the right half of \cref{gmonoid_unity}, the unity diagram in \cref{gmonoid_module_axioms} commutes if and only if the following diagram in $\Gtopst$ commutes for each object $U \in \IUsk$.
\begin{equation}\label{gmonoid_mod_unity}
\begin{tikzpicture}[vcenter]
\def\v{-1.4}
\draw[0cell]
(0,0) node (a11) {X_U \sma \stplus}
(a11)++(3,0) node (a12) {X_U}
(a11)++(0,\v) node (a21) {X_U \sma A_0}
(a12)++(0,\v) node (a22) {X_{U \oplus 0}}
;
\draw[1cell=.9]
(a11) edge node {\rho} node[swap] {\iso} (a12)
(a11) edge node[swap] {1 \sma \eta_0} (a21)
(a21) edge node {\umu_{U,0}} (a22)
(a22) edge node[swap,pos=.45] {X_{\rho}} node {\iso} (a12)
;
\end{tikzpicture}
\end{equation}

\parhead{Associativity}.  Similar to \cref{gmonoid_assoc}, the associativity diagram in \cref{gmonoid_module_axioms} commutes if and only if the following diagram in $\Gtopst$ commutes for each triple of objects $(U,V,W) \in (\IUsk)^3$.
\begin{equation}\label{gmonoid_mod_assoc}
\begin{tikzpicture}[vcenter]
\def\h{4} \def\v{-1.4}
\draw[0cell]
(0,0) node (a11) {(X_U \sma A_V) \sma A_W}
(a11)++(\h,0) node (a12) {X_U \sma (A_V \sma A_W)}
(a11)++(0,\v) node (a21) {X_{U \oplus V} \sma A_W}
(a12)++(0,\v) node (a22) {X_U \sma A_{V \oplus W}}
(a21)++(0,\v) node (a31) {X_{(U \oplus V) \oplus W}}
(a31)++(\h,0) node (a32) {X_{U \oplus (V \oplus W)}}
;
\draw[1cell=.9]
(a11) edge node {\al} node[swap] {\iso} (a12)
(a12) edge node {1 \sma \mu_{V,W}} (a22)
(a22) edge node {\umu_{U,V \oplus W}} (a32)
(a11) edge node[swap] {\umu_{U,V} \sma 1} (a21)
(a21) edge node[swap] {\umu_{U \oplus V,W}} (a31)
(a31) edge node {X_\al} node[swap] {\iso} (a32)
;
\end{tikzpicture}
\end{equation}
This finishes the unraveling of an $A$-module for a commutative $G$-monoid $(A,\mu,\eta)$ in $\IUT$.
\end{explanation}

\begin{explanation}[Module Morphisms]\label{expl:gmonoid_mod_mor}
Suppose $(X,\umu^X)$ and $(Y,\umu^Y)$ are $A$-modules for a commutative $G$-monoid $(A,\mu,\eta)$ in $\IUT$.  By \cref{expl:iu_morphism}, a \emph{morphism of $A$-modules} (\cref{def:com_g_monoid} \pcref{def:com_g_monoid_iv}) 
\[(X,\umu^X) \fto{\theta} (Y,\umu^Y)\]
is determined by $V$-component pointed morphisms
\begin{equation}\label{mod_mor_comp}
X_V \fto{\theta_V} Y_V \forspace V \in \IUsk
\end{equation}
such that the following two statements hold.
\begin{description}
\item[Naturality] The naturality diagram \cref{iu_mor_natural} commutes for each isomorphism $f \in \IUsk$.
\item[Compatibility] The following diagram of pointed morphisms commutes for each pair of objects $(V,W) \in (\IUsk)^2$.
\begin{equation}\label{gmonoid_mod_mor_axiom}
\begin{tikzpicture}[vcenter]
\def\v{-1.4}
\draw[0cell]
(0,0) node (a11) {X_V \sma A_W}
(a11)++(3,0) node (a12) {X_{V \oplus W}}
(a11)++(0,\v) node (a21) {Y_{V} \sma A_{W}}
(a12)++(0,\v) node (a22) {Y_{V \oplus W}}
;
\draw[1cell=.9]
(a11) edge node {\umu^X_{V,W}} (a12)
(a12) edge node {\theta_{V \oplus W}} (a22)
(a11) edge node[swap] {\theta_V \sma 1} (a21)
(a21) edge node {\umu^Y_{V,W}} (a22)
;
\end{tikzpicture}
\end{equation} 
\end{description}
While $\umu^X_{V,W}$ and $\umu^Y_{V,W}$ are pointed $G$-morphisms \cref{gmonoid_mod_vw}, the components of $\theta$ are not generally $G$-equivariant.

\parhead{Conjugation $G$-action}.  For each morphism $\theta$ of $A$-modules as above, the $G$-action (\cref{def:com_g_monoid} \pcref{def:com_g_monoid_v}) is given by the $G$-action on $\IUT(X,Y)$, which means componentwise conjugation \cref{gtheta_v}.  Using the $G$-equivariance of components of $\umu^X$ and $\umu^Y$, the following computation shows that the diagram \cref{gmonoid_mod_mor_axiom} remains commutative if each component of $\theta$ is conjugated by an element $g \in G$.
\begin{equation}\label{module_gaction_welldef}
\begin{split}
& \umu^Y_{V,W} \circ (g\theta_V \ginv \sma 1) \\
&= \umu^Y_{V,W} \circ (g \sma g) \circ (\theta_V \sma 1) \circ (\ginv \sma \ginv) \\
&= g \circ \umu^Y_{V,W} \circ (\theta_V \sma 1) \circ (\ginv \sma \ginv) \\
&= g \circ \theta_{V \oplus W} \circ \umu^X_{V,W} \circ (\ginv \sma \ginv) \\
&= g \circ \theta_{V \oplus W} \circ \ginv \circ \umu^X_{V,W}
\end{split}
\end{equation}
Thus, the componentwise conjugation $G$-action is well defined on morphisms of $A$-modules.
\end{explanation}

\section{The $\Gtop$-Category of Orthogonal $G$-Spectra}
\label{sec:gspectra_gtop}

Applying the general theory of commutative $G$-monoids and modules in \cref{sec:gmonoid}, this section reviews the $\Gtop$-category $\GSp$ of orthogonal $G$-spectra for a compact Lie group $G$ and a complete $G$-universe $\univ$ \pcref{def:g_universe}.

\secoutline
\begin{itemize}
\item \cref{def:g_sphere} defines the $G$-sphere $\gsp$, which is a commutative $G$-monoid in $\IUT$.
\item \cref{def:gsp_module} defines the $\Gtop$-category of orthogonal $G$-spectra as the $\Gtop$-category of $\gsp$-modules.
\item \cref{expl:gspectra,expl:gsp_morphism} unravel the structures of orthogonal $G$-spectra and morphisms between them.
\end{itemize}

The following definition uses the description of a commutative $G$-monoid in $\IUT$ given in \cref{expl:com_g_monoid}.

\begin{definition}[Equivariant Sphere]\label{def:g_sphere}
The \emph{$G$-sphere}\index{G-sphere@$G$-sphere}\index{sphere} $(\gsp,\mu,\eta)$ is the commutative $G$-monoid in $\IUT$ (\cref{def:com_g_monoid} \pcref{def:com_g_monoid_ii}) defined as follows. 
\begin{description}
\item[$\IU$-space] 
The underlying $\IU$-space \pcref{def:iu_space}
\begin{equation}\label{g_sphere}
\IU \fto{\gsp} \Topgst
\end{equation}
is defined on objects by
\begin{equation}\label{gsp_v}
\gsp(V) = S^V \forspace V \in \IU,
\end{equation}
where $S^V = V \sqcup \{\infty\}$ is the $V$-sphere \pcref{def:indexing_gspace}.  For $V,W \in \IU$, the component pointed $G$-morphism \cref{iu_space_comp_mor}
\begin{equation}\label{gsphere_mor}
\IU(V,W)_\splus \fto{\gsp} \Topgst(S^V, S^W)
\end{equation}
sends a linear isometric isomorphism $f \cn V \fiso W$ in $\IU$ to the pointed homeomorphism 
\[S^V \fto[\iso]{\gsp f} S^W\]
defined by
\begin{equation}\label{gsp_morphism}
(\gsp f)(x) = \begin{cases}
\infty & \text{if $x = \infty \in S^V$, and}\\
f(x) & \text{if $x \in V$.}
\end{cases}
\end{equation}
\begin{itemize}
\item The composition and identity diagrams \cref{iu_space_axioms} commute for $\gsp$ by \cref{gsp_morphism}. 
\item The equivariance diagram \cref{x_gfginv} commutes for $\gsp$ because each basepoint $\infty$ is $G$-fixed.
\end{itemize}
To simplify the notation, we also denote $\gsp f$ by $f$, with $f(\infty) = \infty$.
\item[Unit] 
The unit \cref{gmonoid_unit} 
\[\iu \fto{\eta} \gsp\]
is determined by the identity pointed morphism \cref{gmonoid_unit_comp}  
\begin{equation}\label{gsp_unit}
\iu_0 = \stplus \fto{\eta_0 = 1} \gsp(0) = S^0
\end{equation}
of the 0-sphere.
\item[Multiplication] 
The multiplication \cref{gmonoid_mult}
\begin{equation}\label{gsphere_mu}
\gsp \smau \gsp \fto{\mu} \gsp
\end{equation}
has, for each pair of objects $(V,W) \in (\IUsk)^2$, $(V,W)$-component pointed $G$-homeomorphism \cref{gmonoid_mult_vw} 
\begin{equation}\label{gsphere_multiplication}
S^V \sma S^W \fto[\iso]{\mu_{V,W}} S^{V \oplus W}
\end{equation}
defined by
\begin{equation}\label{gsp_mult}
\mu_{V,W}(x; y) = \begin{cases}
\infty & \text{if $x = \infty \in S^V$ or $y = \infty \in S^W$, and}\\
x \oplus y & \text{if $(x;y) \in V \sma W$.}
\end{cases}
\end{equation}
\item[Axioms] 
The commutative $G$-monoid axioms---namely, the commutative diagrams \cref{gmonoid_mult_nat,gmonoid_unity,gmonoid_assoc,gmonoid_sym}---hold for $\gsp$ by \cref{gsp_morphism,gsp_unit,gsp_mult}.  For example, the symmetry diagram \cref{gmonoid_sym} commutes by the following equalities for $(x;y) \in V \sma W$.
\[\begin{split}
(\gsp \xi) (\mu_{W,V}) \xi (x;y) 
&= (\gsp \xi) (\mu_{W,V}) (y;x) \\
&= \xi (y \oplus x) \\
&= x \oplus y \\
&= \mu_{V,W}(x;y)
\end{split}\]
\end{description}
This finishes the definition of the commutative $G$-monoid $(\gsp,\mu,\eta)$.
\end{definition}

\begin{explanation}[$G$-Sphere]\label{expl:gsphere_sm}
The triple
\[(\IU,\oplus,0) \fto{(\gsp,\mu,\eta_0)} (\Topgst,\sma,\stplus)\]
is a strong symmetric monoidal functor, with
\begin{itemize}
\item unit constraint given by $\eta_0$ in \cref{gsp_unit} and
\item monoidal constraint given by the components of $\mu$ in \cref{gsphere_multiplication}.
\end{itemize}
However, the multiplication $\mu$ in \cref{gsphere_mu} is \emph{not} an $\IU$-isomorphism.
\end{explanation}

\begin{definition}[Orthogonal $G$-Spectra]\label{def:gsp_module}
Using the commutative $G$-monoid $(\gsp,\mu,\eta)$ in $\IUT$ \pcref{def:g_sphere}, we define the $\Gtop$-category
\[\GSp = \gspmod\]
of $\gsp$-modules and morphisms of $\gsp$-modules (\cref{def:com_g_monoid} \pcref{def:com_g_monoid_v}).  An object in $\GSp$ is called an \emph{orthogonal $G$-spectrum}\index{orthogonal G-spectrum@orthogonal $G$-spectrum} or an \index{module!SG@$\gsp$}\index{SG-module@$\gsp$-module}\emph{$\gsp$-module}.
\end{definition}

\begin{explanation}[Orthogonal $G$-Spectra]\label{expl:gspectra}
By \cref{expl:gmonoid_module}, an orthogonal $G$-spectrum, which means an $\gsp$-module, is a pair $(X,\umu)$ consisting of
\begin{itemize}
\item an $\IU$-space $X \cn \IU \to \Topgst$ \pcref{def:iu_space} and
\item a $G$-equivariant $\IU$-morphism (\cref{def:iu_morphism} \pcref{def:iu_morphism_iii})
\begin{equation}\label{gspectra_action}
X \smau \gsp \fto{\umu} X,
\end{equation}
called the \index{right SG-action@right $\gsp$-action}\emph{right $\gsp$-action},
\end{itemize}
such that the associativity and unity diagrams in \cref{gmonoid_module_axioms} commute.

\parhead{Right $\gsp$-action}.  The right $\gsp$-action $\umu \cn X \smau \gsp \to X$ in \cref{gspectra_action} is determined by $(V,W)$-component pointed $G$-morphisms
\begin{equation}\label{gsp_action_vw}
X_V \sma S^W \fto{\umu_{V,W}} X_{V \oplus W} \forspace (V,W) \in (\IUsk)^2
\end{equation}
such that, for each pair of linear isometric isomorphisms
\[V \fto[\iso]{f} V' \andspace W \fto[\iso]{h} W' \inspace \IUsk,\]
the following naturality diagram of pointed morphisms commutes.
\begin{equation}\label{gsp_action_nat}
\begin{tikzpicture}[vcenter]
\def\v{-1.4}
\draw[0cell]
(0,0) node (a11) {X_V \sma S^W}
(a11)++(3,0) node (a12) {\phantom{X_{V \oplus W}}}
(a12)++(0,-.06) node (a12') {X_{V \oplus W}}
(a11)++(0,\v) node (a21) {X_{V'} \sma S^{W'}}
(a12)++(0,\v) node (a22) {\phantom{X_{V' \oplus W'}}}
(a22)++(0,-.07) node (a22') {X_{V' \oplus W'}}
;
\draw[1cell=.9]
(a11) edge[transform canvas={yshift=-.2ex}] node {\umu_{V,W}} (a12)
(a12') edge node {X_{f \oplus h}} node[swap] {\iso} (a22')
(a11) edge[shorten >=-.7ex] node[swap,pos=.6] {X_f \sma \gsp h} node[pos=.55] {\iso} (a21)
(a21) edge[transform canvas={yshift=-.2ex}] node {\umu_{V',W'}} (a22)
;
\end{tikzpicture}
\end{equation}

\parhead{Unity}.  The unity diagram in \cref{gmonoid_module_axioms} is equivalent to the following commutative diagram in $\Gtopst$ for each object $U \in \IUsk$.
\begin{equation}\label{gsp_unity}
\begin{tikzpicture}[vcenter]
\def\v{-1.4}
\draw[0cell]
(0,0) node (a11) {X_U \sma \stplus}
(a11)++(3,0) node (a12) {X_U}
(a11)++(0,\v) node (a21) {X_U \sma S^0}
(a12)++(0,\v) node (a22) {\phantom{X_{U \oplus 0}}}
(a22)++(0,-.05) node (a22') {X_{U \oplus 0}}
;
\draw[1cell=.9]
(a11) edge node {\rho} node[swap] {\iso} (a12)
(a11) edge[equal] node[swap] {1 \sma \eta_0} (a21)
(a21) edge node {\umu_{U,0}} (a22)
(a22') edge node[swap,pos=.45] {X_{\rho}} node {\iso} (a12)
;
\end{tikzpicture}
\end{equation}

\parhead{Associativity}.  The associativity diagram in \cref{gmonoid_module_axioms} is equivalent to the following commutative diagram in $\Gtopst$ for each triple of objects $(U,V,W) \in (\IUsk)^3$.
\begin{equation}\label{gsp_assoc}
\begin{tikzpicture}[vcenter]
\def\h{4} \def\v{-1.4}
\draw[0cell]
(0,0) node (a11) {(X_U \sma S^V) \sma S^W}
(a11)++(\h,0) node (a12) {X_U \sma (S^V \sma S^W)}
(a11)++(0,\v) node (a21) {X_{U \oplus V} \sma S^W}
(a12)++(0,\v) node (a22) {X_U \sma S^{V \oplus W}}
(a21)++(0,\v) node (a31) {X_{(U \oplus V) \oplus W}}
(a31)++(\h,0) node (a32) {X_{U \oplus (V \oplus W)}}
;
\draw[1cell=.9]
(a11) edge node {\al} node[swap] {\iso} (a12)
(a12) edge node {1 \sma \mu_{V,W}} node[swap] {\iso} (a22)
(a22) edge node {\umu_{U,V \oplus W}} (a32)
(a11) edge node[swap] {\umu_{U,V} \sma 1} (a21)
(a21) edge node[swap] {\umu_{U \oplus V,W}} (a31)
(a31) edge[transform canvas={yshift=.3ex}] node {X_\al} node[swap] {\iso} (a32)
;
\end{tikzpicture}
\end{equation}
This finishes the unraveling of an orthogonal $G$-spectrum $(X,\umu)$.
\end{explanation}

\begin{example}[$G$-Sphere]\label{ex:gsphere_spectrum}
The $G$-sphere $\gsp$ in \cref{def:g_sphere} is an orthogonal $G$-spectrum.  Its right $\gsp$-action $\umu_{V,W}$ \cref{gsp_action_vw} is given by the multiplication $\mu_{V,W}$ \cref{gsp_mult}.  The diagrams \cref{gsp_action_nat,gsp_unity,gsp_assoc} commute for $\gsp$ by \cref{gsp_morphism,gsp_unit,gsp_mult}.
\end{example}

\begin{explanation}[$\gsp$-Module Morphisms]\label{expl:gsp_morphism}
By \cref{expl:gmonoid_mod_mor}, a morphism of orthogonal $G$-spectra, which means a morphism of $\gsp$-modules (\cref{def:com_g_monoid} \pcref{def:com_g_monoid_iv}),
\[(X,\umu^X) \fto{\theta} (Y,\umu^Y)\]
is determined by $V$-component pointed morphisms between pointed $G$-spaces
\begin{equation}\label{gsp_mor_comp}
X_V \fto{\theta_V} Y_V \forspace V \in \IUsk
\end{equation}
such that the following two statements hold.
\begin{description}
\item[Naturality] The naturality diagram \cref{iu_mor_natural} commutes for each isomorphism $f \in \IUsk$.
\item[Compatibility] The following diagram of pointed morphisms commutes for each pair of objects $(V,W) \in (\IUsk)^2$.
\begin{equation}\label{gsp_mor_axiom}
\begin{tikzpicture}[vcenter]
\def\v{-1.4}
\draw[0cell]
(0,0) node (a11) {X_V \sma S^W}
(a11)++(3,0) node (a12) {\phantom{X_{V \oplus W}}}
(a12)++(0,-.06) node (a12') {X_{V \oplus W}}
(a11)++(0,\v) node (a21) {Y_{V} \sma S^{W}}
(a12)++(0,\v) node (a22) {\phantom{Y_{V \oplus W}}}
(a22)++(0,-.06) node (a22') {Y_{V \oplus W}}
;
\draw[1cell=.9]
(a11) edge node {\umu^X_{V,W}} (a12)
(a12') edge node {\theta_{V \oplus W}} (a22')
(a11) edge node[swap] {\theta_V \sma 1} (a21)
(a21) edge node {\umu^Y_{V,W}} (a22)
;
\end{tikzpicture}
\end{equation} 
\end{description}
The components of $\theta$ are \emph{not} $G$-equivariant in general.  Composition and identities for $\gsp$-module morphisms are defined componentwise using \cref{gsp_mor_comp}.

\parhead{$\Gtop$-enrichment}. The $\Gtop$-category $\GSp$ of orthogonal $G$-spectra has hom $G$-spaces 
\begin{equation}\label{gsp_gtop_enr}
\begin{split}
\GSp\big((X,\umu^X), (Y,\umu^Y)\big) 
&\bigsubset \IUT(X,Y) \\
&\bigsubset \prod_{V \in \IUsk} \Topgst(X_V,Y_V).
\end{split}
\end{equation}
\begin{itemize}
\item The first inclusion in \cref{gsp_gtop_enr} comes from the fact that each morphism of $\gsp$-modules is an $\IU$-morphism with extra properties \cref{gsp_mor_axiom}, but not extra structure.
\item The second inclusion in \cref{gsp_gtop_enr} is the one in \cref{IUTXY_topology}, which defines the $\Gtop$-enrichment of $\IUT$.
\end{itemize}
The group $G$ acts on morphisms of orthogonal $G$-spectra by componentwise conjugation \cref{gtheta_v}.  The computation \cref{module_gaction_welldef} shows that this $G$-action is well defined.
\end{explanation}

\section{Smash Product of Orthogonal $G$-Spectra}
\label{sec:gspectra_smash}

This section reviews the smash product of orthogonal $G$-spectra and their morphisms for a compact Lie group $G$ and a complete $G$-universe $\univ$ \pcref{def:g_universe}.  This smash product is used in \cref{sec:gspectra} to construct the symmetric monoidal $\Gtop$-category $\GSp$ of orthogonal $G$-spectra.

\secoutline
\begin{itemize}
\item \cref{def:gsp_sma} defines the smash product $X \smasg Y$ of two orthogonal $G$-spectra.
\item \cref{gsp_sma_welldef} proves that this construction is well defined.
\item \cref{def:gsp_mor_sma} defines the smash product $\theta \smasg \theta'$ of two morphisms of orthogonal $G$-spectra.
\item \cref{smasg_mor_welldef} proves that this construction is well defined.
\end{itemize}

\subsection*{Smash Product of $\gsp$-Modules}

Recall the symmetric monoidal $\Gtop$-category $(\IUT,\smau)$ and the $G$-sphere $\gsp$ \pcref{def:IU_smgtop,def:g_sphere}.

\begin{definition}\label{def:gsp_sma}
Given two orthogonal $G$-spectra $(X,\umu^X)$ and $(Y,\umu^Y)$ \pcref{def:gsp_module}, the \emph{smash product}\index{smash product!orthogonal G-spectrum@orthogonal $G$-spectrum}\index{orthogonal G-spectrum@orthogonal $G$-spectrum!smash product}
\[(X \smasg Y, \umu)\]
is the orthogonal $G$-spectrum whose underlying $\IU$-space is defined as the coequalizer
\begin{equation}\label{gsp_sma_coequal}
\begin{tikzpicture}[baseline={(a1.base)}]
\draw[0cell=.9]
(0,0) node (a1) {\phantom{(X \smau \gsp) \smau Y}}
(a1)++(0,-.04) node (a1') {(X \smau \gsp) \smau Y}
(a1)++(4,0) node (a2) {X \smau Y}
(a2)++(2,0) node (a3) {X \smasg Y}
;
\draw[1cell=.8]
(a1) edge[transform canvas={yshift=.4ex}] node {\umu^X \smau 1_Y} (a2)
(a1) edge[transform canvas={yshift=-.5ex}] node[swap] {(1_X \smau \umu^Y) \upbe} (a2)
(a2) edge node {\psma} (a3)
;
\end{tikzpicture}
\end{equation}
taken in $\Gtopst$ for each object $U \in \IU$.  In \cref{gsp_sma_coequal}, $\upbe$ denotes the $G$-equivariant $\IU$-isomorphism
\begin{equation}\label{gsp_sma_beta}
\begin{tikzpicture}[baseline={(a1.base)}]
\def\u{.65}
\draw[0cell=.9]
(0,0) node (a1) {(X \smau \gsp) \smau Y}
(a1)++(3.2,0) node (a2) {X \smau (\gsp \smau Y)}
(a2)++(3.8,0) node (a3) {X \smau (Y \smau \gsp)}
;
\draw[1cell=.8]
(a1) edge node {\au} node[swap] {\iso} (a2)
(a2) edge node {1_X \smau \beu} node[swap] {\iso} (a3)
(a1) [rounded corners=2pt] |- ($(a2)+(-1,\u)$) -- node {\upbe} ($(a2)+(1,\u)$) -| (a3) 
;
\end{tikzpicture}
\end{equation}
where $\au$ \cref{IU_au} and $\beu$ \cref{IU_braiding} are components of, respectively, the monoidal associator and braiding for $\IUT$.  The right $\gsp$-action
\begin{equation}\label{smasg_umu}
(X \smasg Y) \smau \gsp \fto{\umu} X \smasg Y
\end{equation}
is induced by $\umu^Y$, or equivalently $\umu^X$.
\end{definition}

\begin{lemma}\label{gsp_sma_welldef}
In \cref{def:gsp_sma}, $(X \smasg Y, \umu)$ is an orthogonal $G$-spectrum.
\end{lemma}

\begin{proof}
\parhead{Well-defined $\IU$-space}. 
To see that the $\IU$-space $X \smasg Y$ is well defined, recall that the right $\gsp$-actions 
\[X \smau \gsp \fto{\umu^X} X \andspace Y \smau \gsp \fto{\umu^Y} Y\]
are $G$-equivariant $\IU$-morphisms \cref{gspectra_action}.  By \cref{expl:giut_smtop}, each of the two arrows $\umu^X \smau 1_Y$ and $(1_X \smau \umu^Y) \upbe$ in \cref{gsp_sma_coequal} is a $G$-equivariant $\IU$-morphism.  Thus, their componentwise coequalizer in $\Gtopst$ is well defined.  Equivalently, one can take their componentwise coequalizer in $\Topst$ and the $G$-action inherited from each component of $X \smau Y$.

For each linear isometric isomorphism $f \cn V \fiso W$ in $\IU$, the pointed homeomorphism \cref{iu_space_xf}
\[(X \smasg Y)_V \fto[\iso]{(X \smasg Y)_f} (X \smasg Y)_W\]
is induced by $(X \smau Y)_f$ in \cref{xsmauy_mor}, using 
\begin{itemize}
\item the universal property of coequalizers in $\Topst$ and
\item the naturality \cref{iu_mor_natural} of the $\IU$-morphisms $\umu^X \smau 1_Y$ and $(1_X \smau \umu^Y) \upbe$.
\end{itemize} 
The composition, identity, and equivariance diagrams in \cref{iu_space_axioms,x_gfginv} commute for $X \smasg Y$ by
\begin{itemize}
\item the universal property of coequalizers and
\item the corresponding commutative diagrams for the $\IU$-space $X \smau Y$ \pcref{smau_welldef}.
\end{itemize}  
Thus, the $\IU$-space $X \smasg Y$ is well defined.  Its object assignment is explained further below.

\parhead{Component objects}.  For each object $U \in \IU$, \cref{smau_wedge,x_smau_y_rep,smau_coend,smau_iden} describe the pointed space $(X \smau Y)_U$ in terms of representatives.  There is an analogous description of the pointed space $(X \smasg Y)_U$ as follows.  Using the pointed $G$-homeomorphism \cref{au_xyzu_dom}
\[\begin{split}
& \big((X \smau \gsp) \smau Y\big)_U \\
&\iso \scalebox{.9}{$\txint^{(V_1,V_2,V_3) \in (\IUsk)^3} \IU\big((V_1 \oplus V_2) \oplus V_3, U\big)_\splus \sma \big((X_{V_1} \sma S^{V_2}) \sma Y_{V_3} \big)$},
\end{split}\]
consider a representative quadruple
\begin{equation}\label{xsgyu_rep}
\bv = \scalebox{.9}{$\big((V_1 \oplus V_2) \oplus V_3 \fto[\iso]{e} U\in \IU ; x \in X_{V_1} ; a \in S^{V_2} ; y \in Y_{V_3} \big)$}.
\end{equation}
The first two morphisms in \cref{gsp_sma_coequal} send $\bv$ to the following two triples, which represent two points in $(X \smau Y)_U$.
\begin{equation}\label{dzeroone_bv}
\begin{split}
(\umu^X \smau 1_Y)_U (\bv)
&= \big(e ; \umu^X_{V_1,V_2}(x,a) \in X_{V_1 \oplus V_2} ; y \big) \\
\big((1_X \smau \umu^Y) \upbe\big)_U (\bv) 
&= \big(e \al^{-1} (1 \oplus \xi) ; x ; \umu^Y_{V_3,V_2}(y,a) \in Y_{V_3 \oplus V_2}\big)
\end{split}
\end{equation}
In \cref{dzeroone_bv}, $e \al^{-1} (1 \oplus \xi)$ is the following composite isomorphism in $\IU$, where $\xi$ and $\al$ are the braiding and associativity isomorphism for the symmetric monoidal category $(\IUsk,\oplus)$ \pcref{def:IU_spaces}.
\[\begin{tikzpicture}[vcenter]
\def\v{-1.4}
\draw[0cell]
(0,0) node (a1) {V_1 \oplus (V_3 \oplus V_2)}
(a11)++(0,\v) node (a2) {V_1 \oplus (V_2 \oplus V_3)}
(a2)++(4,0) node (a3) {(V_1 \oplus V_2) \oplus V_3}
(a3)++(0,-\v) node (a4) {U}
;
\draw[1cell=.9]
(a1) edge node[swap] {1_{V_1} \oplus \xi} (a2)
(a2) edge node {\al^{-1}} (a3)
(a3) edge node[swap] {e} (a4)
;
\end{tikzpicture}\]
The pointed space $(X \smasg Y)_U$ is the quotient of the pointed space $(X \smau Y)_U$ that identifies the two points in \cref{dzeroone_bv} for each representative $\bv$ in \cref{xsgyu_rep}.

To understand \cref{dzeroone_bv} with simpler notation, let us abbreviate $\umu^X_{V_1,V_2}(x,a)$ to $xa$ and likewise for $\umu^Y$.  Then \cref{dzeroone_bv} identifies the two points
\begin{equation}\label{xsmasgy_iden}
(e; xa ; y) = (e' ; x ; ya),
\end{equation}
where $e'$ means $e$ pre-composed with a unique coherence isomorphism in $(\IUsk,\oplus)$ that adjusts for the change of domain.

\parhead{Right $\gsp$-action}.
By \cref{gsp_action_vw}, for each pair of objects $(U_1,U_2) \in (\IUsk)^2$, the $(U_1,U_2)$-component of the right $\gsp$-action $\umu$ \cref{smasg_umu} is the pointed $G$-morphism
\begin{equation}\label{smasg_sphere_action}
(X \smasg Y)_{U_1} \sma S^{U_2} \fto{\umu_{U_1,U_2}} (X \smasg Y)_{U_1 \oplus U_2}
\end{equation}
induced by $\umu^Y$, or equivalently $\umu^X$, as follows.  Consider a point $a \in S^{U_2}$ and a triple
\begin{equation}\label{x_smasg_y_rep}
\big(V \oplus W \fto[\iso]{e} U_1 \in \IUsk ; x \in X_V ; y \in Y_W\big),
\end{equation}
which represents a point in $(X \smasg Y)_{U_1}$.  
\begin{itemize}
\item The right $\gsp$-action on $Y$ yields the triple
\begin{equation}\label{smasg_action_y}
\big((e \oplus 1)\al^{-1} ; x ; \umu^Y_{W,U_2}(y,a) \in Y_{W \oplus U_2} \big),
\end{equation}
which represents a point in $(X \smasg Y)_{U_1 \oplus U_2}$.  The first entry in \cref{smasg_action_y} is the following composite isomorphism in $\IUsk$.
\[V \oplus (W \oplus U_2) \fto{\al^{-1}} (V \oplus W) \oplus U_2 \fto{e \oplus 1} U_1 \oplus U_2\]
\item The right $\gsp$-action on $X$ yields the triple
\begin{equation}\label{smasg_action_x}
\big((e \oplus 1)\al^{-1} (1 \oplus \xi) \al ; \umu^X_{V,U_2}(x, a) \in X_{V \oplus U_2} ; y \in Y_W \big),
\end{equation}
which represents a point in $(X \smasg Y)_{U_1 \oplus U_2}$.  The first entry in \cref{smasg_action_x} is the following composite isomorphism in $\IUsk$.
\[\begin{tikzpicture}[vcenter]
\def\h{3.5}\def\v{-1.4}
\draw[0cell=.9]
(0,0) node (a1) {(V \oplus U_2) \oplus W}
(a11)++(0,\v) node (a2) {V \oplus (U_2 \oplus W)}
(a2)++(\h,0) node (a3) {V \oplus (W \oplus U_2)}
(a3)++(\h,0) node (a4) {(V \oplus W) \oplus U_2}
(a4)++(0,-\v) node (a5) {U_1 \oplus U_2}
;
\draw[1cell=.85]
(a1) edge node[swap] {\al} (a2)
(a2) edge node {1 \oplus \xi} (a3)
(a3) edge node {\al^{-1}} (a4)
(a4) edge node[swap] {e \oplus 1} (a5)
;
\end{tikzpicture}\]
\end{itemize}
By the identification in \cref{dzeroone_bv}, the two triples in \cref{smasg_action_y,smasg_action_x} represent the same point in $(X \smasg Y)_{U_1 \oplus U_2}$; their common value is $\umu_{U_1,U_2} (e; x; y; a)$.

Using the notation in \cref{xsmasgy_iden,smasg_action_y,smasg_action_x}, the right $\gsp$-action $\umu_{U_1,U_2}$ is given by
\begin{equation}\label{smasg_exya}
\begin{split}
\umu_{U_1,U_2} (e; x; y; a) &= \big((e \oplus 1)' ; x ; ya \big) \\
&= \big((e \oplus 1)' ; xa ; y \big).
\end{split}
\end{equation}
In \cref{smasg_exya}, each $(e \oplus 1)'$ means $e \oplus 1$ pre-composed with a unique coherence isomorphism in $(\IUsk,\oplus)$ that adjusts for the change of domain.

\parhead{Well-defined $\gsp$-action}.  To see that $\umu$ in \cref{smasg_exya} is well defined, we need to show that, if the triple $(e; x; y)$ is given by either one of the two triples in \cref{xsmasgy_iden} with $\bv = (e; x; a; y)$ as defined in \cref{xsgyu_rep}, then the action $\umu$ yields the same result.  Using the uniqueness of coherence isomorphisms in $(\IUsk,\oplus)$ and the functoriality of $\oplus$, the following computation proves that $\umu$ is well defined, where $b \in S^{U'}$ and $(a,b) \in S^{V_2} \sma S^{U'} \iso S^{V_2 \oplus U'}$.
\[\begin{aligned}
&\umu_{U,U'} (e; xa; y; b) &&\\
&= \big( (e \oplus 1)' ; xa ; yb\big) && \text{by \cref{smasg_exya}} \\
&= \big( (e \oplus 1)' ; (xa)b ; y\big) && \text{by \cref{xsmasgy_iden}} \\
&= \big( (e \oplus 1)' ; x(a,b) ; y\big) && \text{by \cref{gsp_assoc}} \\
&= \big( (e \oplus 1)' ; x ; y(a,b)\big) && \text{by \cref{xsmasgy_iden}} \\
&= \big( (e \oplus 1)' ; x ; (ya)b\big) && \text{by \cref{gsp_assoc}} \\
&= \umu_{U,U'} (e'; x; ya; b) && \text{by \cref{smasg_exya}} 
\end{aligned}\]

\parhead{Unity and associativity}.
The pair $(X \smasg Y, \umu)$ satisfies the axioms \cref{gsp_unity,gsp_assoc} by 
\begin{itemize}
\item the description \cref{smasg_exya} of $\umu$ in terms of $\umu^Y$,
\item the corresponding unity and associativity axioms for $(Y,\umu^Y)$, and
\item the coherence of the monoidal category $(\IUsk,\oplus)$.
\end{itemize}

\parhead{Naturality}.  The naturality axiom \cref{gsp_action_nat} for $(X \smasg Y, \umu)$ states that, for each pair of linear isometric isomorphisms
\[U_1 \fto[\iso]{f} U_1' \andspace U_2 \fto[\iso]{h} U_2' \inspace \IUsk,\]
the following diagram of pointed morphisms commutes.
\begin{equation}\label{smasg_action_nat}
\begin{tikzpicture}[vcenter]
\def\v{-1.4}
\draw[0cell=.9]
(0,0) node (a11) {(X \smasg Y)_{U_1} \sma S^{U_2}}
(a11)++(4.2,0) node (a12) {\phantom{(X \smasg Y)_{U_1 \oplus U_2}}}
(a12)++(0,-.02) node (a12') {(X \smasg Y)_{U_1 \oplus U_2}}
(a11)++(0,\v) node (a21) {(X \smasg Y)_{U_1'} \sma S^{U_2'}}
(a12)++(0,\v) node (a22) {\phantom{(X \smasg Y)_{U_1' \oplus U_2'}}}
(a22)++(0,-.03) node (a22') {(X \smasg Y)_{U_1' \oplus U_2'}}
;
\draw[1cell=.8]
(a11) edge[transform canvas={yshift=-.1ex}] node {\umu_{U_1,U_2}} (a12)
(a12') edge[transform canvas={xshift=-2.5ex}] node[pos=.55] {(X \smasg Y)_{f \oplus h}} node[swap] {\iso} (a22')
(a11) edge[transform canvas={xshift=3ex}, shorten >=-.5ex] node[swap,pos=.6] {(X \smasg Y)_f \sma \gsp h} node[pos=.55] {\iso} (a21)
(a21) edge[transform canvas={yshift=0ex}] node {\umu_{U_1',U_2'}} (a22)
;
\end{tikzpicture}
\end{equation}
Using the notation in \cref{x_smasg_y_rep}, the following computation proves that the diagram \cref{smasg_action_nat} commutes.
\[\begin{aligned}
& (X \smasg Y)_{f \oplus h} (\umu_{U_1,U_2}) (e; x; y; a) && \\
&= (X \smasg Y)_{f \oplus h} \big((e \oplus 1_{U_2})\al^{-1}_{V,W,U_2} ; x ; \umu^Y_{W,U_2}(y,a) \big) && \text{by \cref{smasg_action_y}} \\
&= \big( (f \oplus h) (e \oplus 1_{U_2}) \al^{-1}_{V,W,U_2} ; x ; \umu^Y_{W,U_2}(y,a) \big) && \text{by \cref{smau_iumor}} \\
&= \big( (fe \oplus 1_{U_2'}) (1_{V \oplus W} \oplus h) \al^{-1}_{V,W,U_2}  ; x ; \umu^Y_{W,U_2}(y,a) \big) && \text{\scalebox{.9}{by functoriality of $\oplus$}}\\
&= \scalebox{.9}{$\big( (fe \oplus 1_{U_2'}) \al^{-1}_{V,W,U_2'} (1_V \oplus (1_W \oplus h) ); x ; \umu^Y_{W,U_2}(y,a) \big)$} && \text{\scalebox{.9}{by naturality of $\al$}} \\
&= \big( (fe \oplus 1_{U_2'}) \al^{-1}_{V,W,U_2'} ; x ; Y_{1_W \oplus h} \umu^Y_{W,U_2}(y,a) \big) && \text{by \cref{smau_iden}} \\
&= \big( (fe \oplus 1_{U_2'}) \al^{-1}_{V,W,U_2'} ; x ; \umu^Y_{W,U_2'}(y, (\gsp h)a) \big) && \text{by \cref{gsp_action_nat}} \\
&= \umu_{U_1',U_2'} \big(fe ; x; y; (\gsp h)a \big) && \text{by \cref{smasg_action_y}} \\
&= \umu_{U_1',U_2'} \big( ((X \smasg Y)_f \sma \gsp h) (e; x; y; a) \big) && \text{by \cref{smau_iumor}}\\
\end{aligned}\]
This finishes the proof that the smash product $(X \smasg Y, \umu)$ is an orthogonal $G$-spectrum.
\end{proof}

\subsection*{Smash Product of $\gsp$-Module Morphisms}

As we explain in the first paragraph of the proof of \cref{gsp_sma_welldef}, the objectwise coequalizer \cref{gsp_sma_coequal} in $\Gtopst$, which defines the smash product $X \smasg Y$, can also be taken in $\Topst$ along with the inherited $G$-action from $X \smau Y$.  This observation is used in the coequalizer \cref{smasg_mor_coequal} in the next definition.

\begin{definition}\label{def:gsp_mor_sma}
Given two morphisms of orthogonal $G$-spectra \pcref{expl:gsp_morphism}
\[(X,\umu^X) \fto{\theta} (X',\umu^{X'}) \andspace (Y,\umu^Y) \fto{\theta'} (Y',\umu^{Y'}),\]
the \emph{smash product}\index{smash product!orthogonal G-spectrum@orthogonal $G$-spectrum}\index{orthogonal G-spectrum@orthogonal $G$-spectrum!smash product}
\[(X \smasg Y,\umu) \fto{\theta \smasg \theta'} (X' \smasg Y',\umu')\]
is the morphism of orthogonal $G$-spectra induced by the $\IU$-morphism $\theta \smau \theta'$ \pcref{def:iumor_day} and the universal property of coequalizers, as displayed in the following diagram taken objectwise in $\Topst$.
 \begin{equation}\label{smasg_mor_coequal}
\begin{tikzpicture}[vcenter]
\def\h{4.4} \def\j{2.5}
\draw[0cell=.9]
(0,0) node (a1) {\phantom{(X \smau \gsp) \smau Y}}
(a1)++(0,-.04) node (a1') {(X \smau \gsp) \smau Y}
(a1)++(\h,0) node (a2) {X \smau Y}
(a2)++(\j,0) node (a3) {X \smasg Y}
;
\draw[1cell=.8]
(a1) edge[transform canvas={yshift=.4ex}] node {\umu^X \smau 1_Y} (a2)
(a1) edge[transform canvas={yshift=-.5ex}] node[swap] {(1_X \smau \umu^Y) \upbe} (a2)
(a2) edge node {\psma} (a3)
;
\draw[0cell=.9]
(a1)++(0,-1.5) node (b1) {\phantom{(X' \smau \gsp) \smau Y'}}
(b1)++(0,-.04) node (b1') {(X' \smau \gsp) \smau Y'}
(b1)++(\h,0) node (b2) {X' \smau Y'}
(b2)++(\j,0) node (b3) {X' \smasg Y'}
;
\draw[1cell=.8]
(b1) edge[transform canvas={yshift=.4ex}] node {\umu^{X'} \smau 1_{Y'}} (b2)
(b1) edge[transform canvas={yshift=-.5ex}] node[swap] {(1_{X'} \smau \umu^{Y'}) \upbe} (b2)
(b2) edge node {\psma} (b3)
;
\draw[1cell=.8]
(a1) edge[transform canvas={xshift=2em}] node[swap] {(\theta \smau 1_{\gsp}) \smau \theta'} (b1)
(a2) edge node {\theta \smau \theta'} (b2)
(a3) edge[dashed] node {\theta \smasg \theta'} (b3)
;
\end{tikzpicture}
\end{equation}
This finishes the definition of $\theta \smasg \theta'$.
\end{definition}

\begin{lemma}\label{smasg_mor_welldef}
In \cref{def:gsp_mor_sma}, $\theta \smasg \theta'$ is a morphism of orthogonal $G$-spectra.
\end{lemma}

\begin{proof}
\parhead{$\IU$-morphism}.  To see that $\theta \smasg \theta'$ is a well-defined $\IU$-morphism, consider an object $V \in \IU$.  The $V$-component pointed morphism \cref{iu_mor_component}
\[(X \smasg Y)_V \fto{(\theta \smasg \theta')_V} (X' \smasg Y')_V\]
is defined by applying the diagram \cref{smasg_mor_coequal} to $V$.  In the left half of \cref{smasg_mor_coequal}, the two rectangles are commutative by
\begin{itemize}
\item the naturality of $\upbe$ \cref{gsp_sma_beta} and
\item the compatibility diagram \cref{gsp_mor_axiom} for each of $\theta$ and $\theta'$.
\end{itemize}
The universal property of the coequalizer in the top row of \cref{smasg_mor_coequal} implies the unique existence of $(\theta \smasg \theta')_V$ that makes the right square commutes at $V$.  The naturality \cref{iu_mor_natural} for $\theta \smasg \theta'$ follows from the naturality of $\theta \smau \theta'$ and of coequalizers.  Thus, $\theta \smasg \theta'$ is an $\IU$-morphism \pcref{def:iu_morphism}.

\parhead{Compatibility with $\gsp$-action}.
The compatibility diagram \cref{gsp_mor_axiom} for $\theta \smasg \theta'$ is the following diagram of pointed morphisms for each pair of objects $(U_1,U_2) \in (\IUsk)^2$.
\begin{equation}\label{smasg_mor_axiom}
\begin{tikzpicture}[vcenter]
\def\v{-1.4}
\draw[0cell=1]
(0,0) node (a11) {(X \smasg Y)_{U_1} \sma S^{U_2}}
(a11)++(4.5,0) node (a12) {\phantom{(X \smasg Y)_{U_1 \oplus U_2}}}
(a12)++(0,-.02) node (a12') {(X \smasg Y)_{U_1 \oplus U_2}}
(a11)++(0,\v) node (a21) {(X' \smasg Y')_{U_1} \sma S^{U_2}}
(a12)++(0,\v) node (a22) {\phantom{(X' \smasg Y')_{U_1 \oplus U_2}}}
(a22)++(0,-.02) node (a22') {(X' \smasg Y')_{U_1 \oplus U_2}}
;
\draw[1cell=.9]
(a11) edge node {\umu_{U_1,U_2}} (a12)
(a12') edge[transform canvas={xshift=-2em}] node {(\theta \smasg \theta')_{U_1 \oplus U_2}} (a22')
(a11) edge[transform canvas={xshift=2em}] node[swap] {(\theta \smasg \theta')_{U_1} \sma 1} (a21)
(a21) edge node {\umu'_{U_1,U_2}} (a22)
;
\end{tikzpicture}
\end{equation} 
Using the notation in \cref{x_smasg_y_rep,smasg_exya}, the following computation proves that the diagram \cref{smasg_mor_axiom} commutes.
\[\begin{aligned}
& (\theta \smasg \theta')_{U_1 \oplus U_2} (\umu_{U_1,U_2}) (e; x; y; a) && \\
&= (\theta \smasg \theta')_{U_1 \oplus U_2} \big((e \oplus 1)' ; x ; ya \big) && \text{by \cref{smasg_exya}} \\
&= \big((e \oplus 1)' ; \theta_V x ; \theta'_{W \oplus U_2} (ya) \big) && \text{by \cref{tha_smau_thap_uz}} \\
&= \big((e \oplus 1)' ; \theta_V x ; (\theta'_{W} y)a \big) && \text{by \cref{gsp_mor_axiom}} \\
&= \umu'_{U_1,U_2} (e; \theta_V x; \theta'_W y ; a) && \text{by \cref{smasg_exya}} \\
&= \umu'_{U_1,U_2} \big( (\theta \smasg \theta')_{U_1} \sma 1 \big) (e; x; y; a) && \text{by \cref{tha_smau_thap_uz}}
\end{aligned}\]
This finishes the proof that $\theta \smasg \theta'$ is a morphism of orthogonal $G$-spectra.
\end{proof}

\section{The Symmetric Monoidal $\Gtop$-Category of\\ Orthogonal $G$-Spectra}
\label{sec:gspectra}

Using the smash product $\smasg$ constructed in \cref{sec:gspectra_smash}, this section shows that the $\Gtop$-category $\GSp$ of orthogonal $G$-spectra \pcref{def:gsp_module} extends to a symmetric monoidal $\Gtop$-category.  We continue to assume that $G$ is a compact Lie group and $\univ$ is a complete $G$-universe \pcref{def:g_universe}.

\secoutline
\begin{itemize}
\item \cref{def:gsp_smgtop} defines the symmetric monoidal $\Gtop$-category data for $\GSp$.
\item \cref{gspectra_smgtop} proves that these data constitute a symmetric monoidal $\Gtop$-category.
\end{itemize}

For \cref{def:gsp_smgtop} below, recall symmetric monoidal $\V$-category \pcref{definition:monoidal-vcat,definition:braided-monoidal-vcat,definition:symm-monoidal-vcat} and the Cartesian closed category $(\Gtop,\times,*)$ \pcref{def:Gtop}. 

\begin{definition}\label{def:gsp_smgtop}
Define a symmetric monoidal $\Gtop$-category\index{orthogonal G-spectrum@orthogonal $G$-spectrum!symmetric monoidal category}\index{symmetric monoidal category!orthogonal G-spectrum@orthogonal $G$-spectrum}
\[(\GSp, \smasg, \gsp, \asg, \ellsg, \rsg, \bsg)\]
as follows.
\begin{description}
\item[Base $\Gtop$-category]
It is the $\Gtop$-category $\GSp$ in \cref{def:gsp_module}.
\begin{itemize}
\item Its objects are $\gsp$-modules \pcref{expl:gspectra}.
\item In each hom $G$-space $\GSp(X,Y)$, the points are morphisms of $\gsp$-modules $X \to Y$ \pcref{expl:gsp_morphism}, with $G$ acting componentwise by conjugation \cref{gtheta_v}.
\end{itemize}
\item[Monoidal composition]
It is the $\Gtop$-functor
\begin{equation}\label{gsp_mon_comp}
\GSp \otimes \GSp \fto{\smasg} \GSp
\end{equation}
that has
\begin{itemize}
\item object assignment 
\[(X,Y) \mapsto X \smasg Y\]
given in \cref{def:gsp_sma}, and
\item underlying morphism assignment
\[(\theta, \theta') \mapsto \theta \smasg \theta'\]
given in \cref{def:gsp_mor_sma}.
\end{itemize}
\item[Monoidal identity]
With $\vtensorunit$ denoting the unit $\Gtop$-category \pcref{definition:unit-vcat}, the monoidal identity is the $\Gtop$-functor
\begin{equation}\label{gsp_mon_id}
\vtensorunit \fto{\gsp} \GSp
\end{equation}
determined by the $G$-sphere $\gsp$ \pcref{def:g_sphere}, regarded as an $\gsp$-module via its multiplication \pcref{ex:gsphere_spectrum}
\[\gsp \smau \gsp \fto{\mu} \gsp.\]
\item[Right monoidal unitor]
It is the $\Gtop$-natural isomorphism
\begin{equation}\label{gsp_right_unitor}
\begin{tikzpicture}[vcenter]
\def\v{1.4} \def\h{.15}
\draw[0cell=.9]
(0,0) node (a11) {\GSp}
(a11)++(3.5,0) node (a12) {\GSp}
(a11)++(\h,\v) node (a01) {\GSp \otimes \vtensorunit}
(a12)++(-\h,\v) node (a02) {\GSp^{\otimes 2}}
;
\draw[1cell=.9]
(a11) edge node[swap] {1} (a12)
(a11) edge node[pos=.2] {(r^\otimes)^{-1}} (a01)
(a01) edge node {1 \otimes \gsp} (a02)
(a02) edge node[pos=.7] {\smasg} (a12)
;
\draw[2cell]
node[between=a11 and a12 at .47, shift={(0,.5*\v)}, rotate=-90, 2label={above,\rsg}] {\Rightarrow}
;
\end{tikzpicture}
\end{equation}
whose component at an $\gsp$-module $(X,\umu^X)$ is a $G$-morphism
\[* \fto{\rsg_X} \GSp(X \smasg \gsp, X),\]
which means a $G$-equivariant $\gsp$-module isomorphism.  In other words, it is a $G$-equivariant $\IU$-isomorphism (\cref{def:iu_morphism} \pcref{def:iu_morphism_iii})
\begin{equation}\label{rsg_x}
X \smasg \gsp \fto[\iso]{\rsg_X} X
\end{equation}
that respects the $\gsp$-actions on its domain and codomain \cref{gsp_mor_axiom}.  The morphism $\rsg_X$ is induced by the right $\gsp$-action 
\[X \smau \gsp \fto{\umu^X} X\] 
via the coequalizer \cref{gsp_sma_coequal}.
\item[Left monoidal unitor]
It is the $\Gtop$-natural isomorphism
\begin{equation}\label{gsp_left_unitor}
\begin{tikzpicture}[vcenter]
\def\v{1.4} \def\h{.15}
\draw[0cell=.9]
(0,0) node (a11) {\GSp}
(a11)++(3.5,0) node (a12) {\GSp}
(a11)++(\h,\v) node (a01) {\vtensorunit \otimes \GSp}
(a12)++(-\h,\v) node (a02) {\GSp^{\otimes 2}}
;
\draw[1cell=.9]
(a11) edge node[swap] {1} (a12)
(a11) edge node[pos=.2] {(\ell^\otimes)^{-1}} (a01)
(a01) edge node {\gsp \otimes 1} (a02)
(a02) edge node[pos=.7] {\smasg} (a12)
;
\draw[2cell]
node[between=a11 and a12 at .47, shift={(0,.5*\v)}, rotate=-90, 2label={above,\ellsg}] {\Rightarrow}
;
\end{tikzpicture}
\end{equation}
whose component at an $\gsp$-module $(X,\umu^X)$ is the $G$-equivariant $\gsp$-module isomorphism
\begin{equation}\label{ellsg_x}
\gsp \smasg X \fto[\iso]{\ellsg_X} X
\end{equation}
induced by the composite $G$-equivariant $\IU$-morphism
\begin{equation}\label{ellsg_x_smau}
\gsp \smau X \fto[\iso]{\beu_{\gsp,X}} X \smau \gsp \fto{\umu^X} X
\end{equation}
via the coequalizer \cref{gsp_sma_coequal}.  In \cref{ellsg_x_smau}, $\beu_{\gsp,X}$ is a component of the braiding for $\IUT$ \cref{beu_xy}.
\item[Monoidal associator]
It is the $\Gtop$-natural isomorphism
\begin{equation}\label{gsp_associator}
\begin{tikzpicture}[vcenter]
\def\h{2.5} \def\u{1} \def\v{-1.4}
\draw[0cell=.9]
(0,0) node (a11) {(\GSp)^{\otimes 2} \otimes \GSp}
(a11)++(\h,\u) node (a12) {(\GSp)^{\otimes 2}}
(a12)++(\h,-\u) node (a13) {\GSp}
(a11)++(0,\v) node (a21) {\GSp \otimes (\GSp)^{\otimes 2}}
(a13)++(0,\v) node (a22) {(\GSp)^{\otimes 2}}
;
\draw[1cell=.8]
(a11) edge node {\smasg \otimes 1} (a12)
(a12) edge node {\smasg} (a13)
(a11) edge node[swap] {a^\otimes} (a21)
(a21) edge node {1 \otimes \smasg} (a22)
(a22) edge node[swap] {\smasg} (a13)
;
\draw[2cell]
node[between=a11 and a13 at .46, shift={(0,.2*\v)}, rotate=-90, 2label={above,\asg}] {\Rightarrow}
;
\end{tikzpicture}
\end{equation}
whose component at a triple $(X,Y,Z)$ of $\gsp$-modules is the $G$-equivariant $\gsp$-module isomorphism 
\begin{equation}\label{asg_xyz}
(X \smasg Y) \smasg Z \fto[\iso]{\asg_{X,Y,Z}} X \smasg (Y \smasg Z)
\end{equation}
induced by the $G$-equivariant $\IU$-isomorphism \cref{au_xyz}
\[(X \smau Y) \smau Z \fto[\iso]{\au_{X,Y,Z}} X \smau (Y \smau Z).\]
\item[Braiding]
It is the $\Gtop$-natural isomorphism
\begin{equation}\label{gsp_braiding}
\begin{tikzpicture}[vcenter]
\def\h{4}
\draw[0cell]
(0,0) node (a1) {(\GSp)^{\otimes 2}}
(a1)++(.5*\h,-1) node (a2) {(\GSp)^{\otimes 2}}
(a1)++(\h,0) node (a3) {\phantom{\GSp}}
(a3)++(0,-.04) node (a3') {\GSp}
;
\draw[1cell=.9]
(a1) edge[transform canvas={yshift=.2ex}] node {\smasg} (a3)
(a1) [rounded corners=2pt] |- node[swap,pos=.25] {\beta^\otimes} (a2)
;
\draw[1cell=.9]
(a2) [rounded corners=2pt] -| node[swap,pos=.75] {\smasg} (a3')
;
\draw[2cell]
node[between=a1 and a3 at .46, shift={(0,-.45)}, rotate=-90, 2label={above,\bsg}] {\Rightarrow}
;
\end{tikzpicture}
\end{equation}
whose component at a pair $(X,Y)$ of $\gsp$-modules is the $G$-equivariant $\gsp$-module isomorphism  
\begin{equation}\label{bsg_xy}
X \smasg Y \fto[\iso]{\bsg_{X,Y}} Y \smasg X
\end{equation}
induced by the $G$-equivariant $\IU$-isomorphism \cref{beu_xy}
\[X \smau Y \fto[\iso]{\beu_{X,Y}} Y \smau X\]
via the coequalizer \cref{gsp_sma_coequal}.
\end{description}
This finishes the definition of the symmetric monoidal $\Gtop$-category $\GSp$; more details are provided in the proof of \cref{gspectra_smgtop} below.

Moreover, $\GSp$ also denotes
\begin{itemize}
\item its underlying symmetric monoidal category and
\item the associated $\Gtop$-multicategory \pcref{proposition:monoidal-v-cat-v-multicat}.\defmark
\end{itemize} 
\end{definition}

\begin{theorem}\label{gspectra_smgtop}
$\GSp$ in \cref{def:gsp_smgtop} is a symmetric monoidal $\Gtop$-category.
\end{theorem}

\begin{proof}
We first observe that the data
\[(\smasg, \gsp, \asg, \ellsg, \rsg, \bsg)\]
are well defined.  Then we observe that the symmetric monoidal $\Gtop$-category axioms are satisfied.

\parhead{Monoidal unit}. The monoidal unit $\gsp$ in \cref{gsp_mon_id} is a $\Gtop$-functor because the componentwise conjugation $G$-action \cref{gtheta_v} on hom $G$-spaces of $\GSp$ fixes each identity $\gsp$-module morphism.  In particular, the identity morphism of the $\gsp$-module $(\gsp,\mu)$ is $G$-fixed.

\parhead{Monoidal composition}.  The $\Gtop$-functoriality of $\smasg$ in \cref{gsp_mon_comp} follows from
\begin{itemize}
\item the $\Gtop$-functoriality of
\[(\IUT)^{\otimes 2} \fto{\smau} \IUT\]
in \cref{IU_mon_composition} and
\item the universal property of the coequalizers in \cref{smasg_mor_coequal}.
\end{itemize}

\parhead{Braiding}.  The component $\bsg_{X,Y}$ \cref{bsg_xy} is induced by $\beu_{X,Y}$ \cref{beu_xy} via the following diagram in $\IUT$, in which each row is an instance of the coequalizer in \cref{gsp_sma_coequal}.
\begin{equation}\label{gsp_braiding_coequal}
\begin{tikzpicture}[vcenter]
\def\h{4.5} \def\v{-1.5}
\draw[0cell=.9]
(0,0) node (a1) {\phantom{(X \smau \gsp) \smau Y}}
(a1)++(0,-.04) node (a1') {(X \smau \gsp) \smau Y}
(a1)++(\h,0) node (a2) {X \smau Y}
(a2)++(.5*\h,0) node (a3) {X \smasg Y}
(a1)++(0,\v) node (b1) {\phantom{(Y \smau \gsp) \smau X}}
(b1)++(0,-.04) node (b1') {(Y \smau \gsp) \smau X}
(a2)++(0,\v) node (b2) {Y \smau X}
(a3)++(0,\v) node (b3) {Y \smasg X}
;
\draw[1cell=.8]
(a1) edge[transform canvas={yshift=.4ex}] node {\umu^X \smau 1_Y} (a2)
(a1) edge[transform canvas={yshift=-.5ex}] node[swap] {(1_X \smau \umu^Y) \upbe} (a2)
(a2) edge node {\psma} (a3)
(b1) edge[transform canvas={yshift=.4ex}] node {(1_Y \smau \umu^X) \upbe} (b2)
(b1) edge[transform canvas={yshift=-.5ex}] node[swap] {\umu^Y \smau 1_X} (b2)
(b2) edge node {\psma} (b3)
(a1) edge node[swap] {\iso} (b1)
(a2) edge node {\beu_{X,Y}} node[swap] {\iso} (b2)
(a3) edge[dashed] node {\bsg_{X,Y}} (b3)
;
\end{tikzpicture}
\end{equation}
\begin{itemize}
\item The left vertical arrow in \cref{gsp_braiding_coequal} is the unique coherence isomorphism in the symmetric monoidal $\Gtop$-category $\IUT$ that interchanges the variables $X$ and $Y$.
\item The two rectangles in the left half of \cref{gsp_braiding_coequal} commute by the naturality and uniqueness of coherence isomorphisms in $\IUT$.
\end{itemize}
The universal property of coequalizers implies the unique existence of the $\IU$-morphism $\bsg_{X,Y}$ such that the right square in \cref{gsp_braiding_coequal} commutes.

The assertion that $\bsg_{X,Y}$ is a $G$-equivariant $\IU$-isomorphism that is natural in $(X,Y)$ follows from
\begin{itemize}
\item the same fact for $\beu_{X,Y}$ and
\item the universal property of coequalizers.
\end{itemize}

To show that $\bsg_{X,Y}$ is compatible with the right $\gsp$-actions of its domain and codomain \cref{gsp_mor_axiom}, we consider the following diagram for each pair of objects $(U_1,U_2) \in (\IUsk)^2$.
\begin{equation}\label{bsg_compat}
\begin{tikzpicture}[vcenter]
\def\v{-1.4}
\draw[0cell=1]
(0,0) node (a11) {(X \smasg Y)_{U_1} \sma S^{U_2}}
(a11)++(4.5,0) node (a12) {\phantom{(X \smasg Y)_{U_1 \oplus U_2}}}
(a12)++(0,-.02) node (a12') {(X \smasg Y)_{U_1 \oplus U_2}}
(a11)++(0,\v) node (a21) {(Y \smasg X)_{U_1} \sma S^{U_2}}
(a12)++(0,\v) node (a22) {\phantom{(Y \smasg X)_{U_1 \oplus U_2}}}
(a22)++(0,-.02) node (a22') {(Y \smasg X)_{U_1 \oplus U_2}}
;
\draw[1cell=.8]
(a11) edge node {\umu^{X \smasg Y}_{U_1,U_2}} (a12)
(a12') edge[transform canvas={xshift=-2em}] node {\bsg_{X,Y; U_1 \oplus U_2}} (a22')
(a11) edge[transform canvas={xshift=2em}] node[swap] {\bsg_{X,Y; U_1} \sma 1} (a21)
(a21) edge node {\umu^{Y \smasg X}_{U_1,U_2}} (a22)
;
\end{tikzpicture}
\end{equation} 
Using the notation in \cref{xsmasgy_iden,x_smasg_y_rep,smasg_exya}, the following computation proves that the diagram \cref{bsg_compat} commutes.
\[\begin{aligned}
& (\bsg_{X,Y; U_1 \oplus U_2}) (\umu^{X \smasg Y}_{U_1,U_2}) (e; x; y; a) \\
&= \bsg_{X,Y; U_1 \oplus U_2} \big((e \oplus 1)' ; x; ya\big) && \text{by \cref{smasg_exya}} \\
&= \big((e \oplus 1)' ; ya ; x\big) && \text{by \cref{beu_XYU_z}} \\
&= \umu^{Y \smasg X}_{U_1,U_2} (e'; y; x; a) && \text{by \cref{smasg_exya}} \\
&= (\umu^{Y \smasg X}_{U_1,U_2}) (\bsg_{X,Y; U_1} \sma 1) (e; x; y; a) && \text{by \cref{beu_XYU_z}}
\end{aligned}\]
This proves that the braiding $\bsg$ \cref{gsp_braiding} is a $\Gtop$-natural isomorphism.

\parhead{Monoidal associator}.  The component $\asg_{X,Y,Z}$ \cref{asg_xyz} is induced by $\au_{X,Y,Z}$ \cref{au_xyz} via the following diagram in $\IUT$, in which each row is a colimit.
\begin{equation}\label{asg_colimit}
\begin{tikzpicture}[vcenter]
\def\h{3.5} \def\v{-1.5} \def\s{40} \def\t{1.15} \def\u{-.65}
\draw[0cell=.8]
(0,0) node (a1) {((X \smau \gsp) \smau Y) \smau Z}
(a1)++(\t,0) node (a1') {\phantom{\smau}}
(a1)++(\h,0) node (a2) {(X \smau Y) \smau Z}
(a2)++(\u,0) node (a2') {\phantom{\smau}}
(a2)++(.8*\h,0) node (a3) {(X \smasg Y) \smasg Z}
(a1)++(0,\v) node (b1) {(X \smau \gsp) \smau (Y \smau Z)}
(b1)++(\t,0) node (b1') {\phantom{\smau}}
(a2)++(0,\v) node (b2) {X \smau (Y \smau Z)}
(b2)++(\u,0) node (b2') {\phantom{\smau}}
(a3)++(0,\v) node (b3) {X \smasg (Y \smasg Z)}
;
\draw[1cell=.75]
(a1') edge[bend left=\s] node[pos=.3,inner sep=0pt] {\umu^X} (a2')
(a1) edge node[pos=.45] {\umu^Y} (a2)
(a1') edge[bend right=\s] node[pos=.7,swap,inner sep=1pt] {\umu^Z} (a2')
(a2) edge node {\psma} (a3)
(b1') edge[bend left=\s] node[pos=.3,inner sep=0pt] {\umu^X} (b2')
(b1) edge node[pos=.45] {\umu^Y} (b2)
(b1') edge[bend right=\s] node[pos=.7,swap,inner sep=1pt] {\umu^Z} (b2')
(b2) edge node {\psma} (b3)
(a1) edge node[swap] {\au_{X \smau \gsp,Y,Z}} node {\iso} (b1)
(a2) edge node {\au_{X,Y,Z}} node[swap] {\iso} (b2)
(a3) edge[dashed] node {\asg_{X,Y,Z}} (b3)
;
\end{tikzpicture}
\end{equation}
\begin{itemize}
\item In \cref{asg_colimit}, identity morphisms and unique coherence isomorphisms in $\IUT$ that permute $\gsp$ to the appropriate position are suppressed.  For example, the top $\umu^X$ means $(\umu^X \smau 1_Y) \smau 1_Z$, and the top $\umu^Y$ is the following composite.
\[\begin{tikzpicture}[vcenter]
\def\h{6} \def\v{-1.4}
\draw[0cell=.8]
(0,0) node (a1) {((X \smau \gsp) \smau Y) \smau Z}
(a1)++(0,\v) node (a2) {(X \smau (\gsp \smau Y)) \smau Z}
(a2)++(\h,0) node (a3) {(X \smau (Y \smau \gsp)) \smau Z}
(a1)++(\h,0) node (a4) {\phantom{(X \smau Y) \smau Z}}
(a4)++(-.45,0) node (a4') {(X \smau Y) \smau Z}
;
\draw[1cell=.8]
(a1) edge[transform canvas={xshift=2em}]  node{\iso} node[swap] {\au_{X,\gsp,Y} \smau 1_Z} (a2)
(a2) edge node[swap] {\iso} node {(1_X \smau \beu_{\gsp,Y}) \smau 1_Z} (a3)
(a3) edge[transform canvas={xshift=-2em}] node[swap] {(1_X \smau \umu^Y) \smau 1_Z} (a4)
;
\end{tikzpicture}\]
\item The left half of \cref{asg_colimit} commutes by the naturality and uniqueness of coherence isomorphisms in $\IUT$.
\end{itemize}
The universal property of colimits implies the unique existence of the $\IU$-morphism $\asg_{X,Y,Z}$ such that the right square in \cref{asg_colimit} commutes.

The assertion that $\asg_{X,Y,Z}$ is a $G$-equivariant $\IU$-isomorphism that is natural in $(X,Y,Z)$ follows from
\begin{itemize}
\item the same fact for $\au_{X,Y,Z}$ and
\item the universal property of colimits.
\end{itemize}

To show that $\asg_{X,Y,Z}$ is compatible with the right $\gsp$-actions of its domain and codomain \cref{gsp_mor_axiom}, consider the following diagram for each pair of objects $(U_1,U_2) \in (\IUsk)^2$.
\begin{equation}\label{asg_compat}
\begin{tikzpicture}[vcenter]
\def\v{-1.4}
\draw[0cell=.9]
(0,0) node (a11) {((X \smasg Y) \smasg Z)_{U_1} \sma S^{U_2}}
(a11)++(5.3,0) node (a12) {\phantom{((X \smasg Y) \smasg Z)_{U_1 \oplus U_2}}}
(a12)++(0,-.02) node (a12') {((X \smasg Y) \smasg Z)_{U_1 \oplus U_2}}
(a11)++(0,\v) node (a21) {(X \smasg (Y \smasg Z))_{U_1} \sma S^{U_2}}
(a12)++(0,\v) node (a22) {\phantom{(X \smasg (Y \smasg Z))_{U_1 \oplus U_2}}}
(a22)++(0,-.02) node (a22') {(X \smasg (Y \smasg Z))_{U_1 \oplus U_2}}
;
\draw[1cell=.8]
(a11) edge node {\umu^{(X \smasg Y) \smasg Z}_{U_1,U_2}} (a12)
(a12') edge[transform canvas={xshift=-2em}] node {\asg_{X,Y,Z; U_1 \oplus U_2}} (a22')
(a11) edge[transform canvas={xshift=2em}] node[swap] {\asg_{X,Y,Z; U_1} \sma 1} (a21)
(a21) edge node {\umu^{X \smasg (Y \smasg Z)}_{U_1,U_2}} (a22)
;
\end{tikzpicture}
\end{equation} 
The diagram \cref{asg_compat} commutes for the following two reasons.
\begin{itemize}
\item By \cref{smasg_exya}, both horizontal arrows in \cref{asg_compat} are induced by the right $\gsp$-action $\umu^Z \cn Z \smau \gsp \to Z$.
\item By \cref{au_xyz_def}, both vertical arrows in \cref{asg_compat} are induced by associativity isomorphisms for $(\IUsk,\oplus)$ and $(\Gtopst,\sma)$.
\end{itemize}
This proves that the monoidal associator $\asg$ \cref{gsp_associator} is a $\Gtop$-natural isomorphism.

\parhead{Monoidal unitors}.  The right monoidal unitor $\rsg_X$ \cref{rsg_x} is induced by the right $\gsp$-action $\umu^X$ via the following diagram in $\IUT$, in which the top row is an instance of the coequalizer \cref{gsp_sma_coequal}.
\begin{equation}\label{rsg_coequal}
\begin{tikzpicture}[vcenter]
\draw[0cell=.9]
(0,0) node (a1) {\phantom{(X \smau \gsp) \smau \gsp}}
(a1)++(0,-.04) node (a1') {(X \smau \gsp) \smau \gsp}
(a1)++(4,0) node (a2) {X \smau \gsp}
(a2)++(2.5,0) node (a3) {X \smasg \gsp}
(a2)++(0,-1.4) node (a4) {X}
;
\draw[1cell=.8]
(a1) edge[transform canvas={yshift=.4ex}] node {\umu^X \smau 1_{\gsp}} (a2)
(a1) edge[transform canvas={yshift=-.5ex}] node[swap] {(1_X \smau \mu) \upbe} (a2)
(a2) edge node {\psma} (a3)
(a2) edge node[swap] {\umu^X} (a4)
(a3) edge[dashed] node[pos=.4] {\rsg_X} (a4)
;
\end{tikzpicture}
\end{equation}
By 
\begin{itemize}
\item the associativity of $(X,\umu^X)$ \cref{gsp_assoc} and
\item the fact that $(\gsp,\mu,\eta)$ is a commutative $G$-monoid \pcref{def:g_sphere}, 
\end{itemize}
there is an equality of $\IU$-morphisms
\[\umu^X (\umu^X \smau 1_{\gsp}) = \umu^X \big((1_X \smau \mu) \upbe\big).\]
The universal property of coequalizers implies the unique existence of the $\IU$-morphism $\rsg_X$ such that the triangle in \cref{rsg_coequal} commutes.

Moreover, $\rsg$ \cref{gsp_right_unitor} is a $\Gtop$-natural isomorphism by  the universal property of coequalizers and statements \pcref{rsgx_i} through \pcref{rsgx_iv} below.
\begin{enumerate}
\item\label{rsgx_i} $\rsg_X$ is $G$-equivariant because $\umu^X$ is so \cref{gspectra_action}.
\item\label{rsgx_ii} $\rsg_X$ is natural in the $\gsp$-module $X$ because morphisms of $\gsp$-modules preserve the right $\gsp$-actions \cref{gsp_mor_axiom}.
\item\label{rsgx_iii} $\rsg_X$ is compatible with the right $\gsp$-actions of its domain and codomain \cref{gsp_mor_axiom} by the associativity axiom \cref{gsp_assoc} of $X$.
\item\label{rsgx_iv} The inverse of $\rsg_X$ is given by the following composite $\IU$-morphism, where $\ru_X$ is the $X$-component of $\ru$ \cref{ru_x}.
\begin{equation}\label{rsg_inverse}
\begin{tikzpicture}[vcenter]
\def\v{-1.4}
\draw[0cell]
(0,0) node (a11) {X}
(a11)++(3.5,0) node (a12) {X \smasg \gsp}
(a11)++(0,\v) node (a21) {X \smau \iu}
(a12)++(0,\v) node (a22) {X \smau \gsp}
;
\draw[1cell=.9]
(a11) edge node {(\rsg_X)^{-1}} (a12)
(a11) edge node[swap] {(\ru_X)^{-1}} (a21)
(a21) edge node {1_X \smau \eta} (a22)
(a22) edge node[swap] {\psma} (a12)
;
\end{tikzpicture}
\end{equation}
Indeed, the composite in \cref{rsg_inverse} followed by $\rsg_X$ is equal to $1_X$ by the following commutative diagram of $\IU$-morphisms.
\begin{equation}\label{rsg_inv_i}
\begin{tikzpicture}[vcenter]
\def\v{-1.4} \def\h{3.3}
\draw[0cell]
(0,0) node (a11) {X}
(a11)++(\h,0) node (a12) {X}
(a12)++(.7*\h,0) node (a22) {X \smasg \gsp}
(a11)++(0,\v) node (a31) {X \smau \iu}
(a12)++(0,\v) node (a32) {X \smau \gsp}
;
\draw[1cell=.8]
(a11) edge node[swap] {(\ru_X)^{-1}} (a31)
(a31) edge node[pos=.6] {1_X \smau \eta} (a32)
(a22) edge node[swap,pos=.5] {\rsg_X} (a12)
(a32) edge node {\umu_X} (a12)
(a31) edge node {\ru_X} (a12)
(a32) [rounded corners=2pt] -| node[pos=.2] {\psma} (a22)
;
\end{tikzpicture}
\end{equation}
In \cref{rsg_inv_i}, the right region commutes by \cref{rsg_coequal}.  The middle triangle commutes by the unity axiom of $X$ \cref{gsp_unity}.

Conversely, to prove that $\rsg_X$ followed by the composite in \cref{rsg_inverse} is equal to the identity on $X \smasg \gsp$, it suffices to show that the left-bottom composite in the diagram \cref{rsg_inv_ii} below is equal to the universal morphism $\psma$ in \cref{rsg_coequal}.  To save space in the following diagram, we abbreviate $\umu^X$ to $\umu$ and $\smau$ to concatenation.  For example, $X\gsp$ means $X \smau \gsp$.
\begin{equation}\label{rsg_inv_ii}
\begin{tikzpicture}[vcenter]
\def\h{2.5} \def\v{-1.4}
\draw[0cell=.8]
(0,0) node (a02) {X(\gsp \iu)}
(a02)++(-\h,\v) node (a11) {X\gsp}
(a11)++(\h,0) node (a12) {(X\gsp)\iu}
(a12)++(\h,0) node (a13) {(X\gsp)\gsp}
(a13)++(\h,0) node (a03) {X(\gsp \gsp)}
(a11)++(0,\v) node (a21) {X}
(a21)++(\h,0) node (a22) {X\iu}
(a22)++(\h,0) node (a23) {X\gsp}
(a23)++(\h,0) node (a24) {X\smasg \gsp}
;
\draw[1cell=.7]
(a02) [rounded corners=2pt] -| node[swap,pos=.75] {1(1\eta)} (a03)
;
\draw[1cell=.7]
(a11) edge node[swap,pos=.4] {\umu} (a21)
(a21) edge node {(\ru)^{-1}} (a22)
(a22) edge node {1\eta} (a23)
(a23) edge node[pos=.6] {\psma} (a24)
(a11) edge node[pos=.6] {(\ru)^{-1}} (a12)
(a12) edge node {(11)\eta} (a13)
(a11) edge node {1(\ru)^{-1}} (a02)
(a12) edge node[swap,pos=.6] {\au} (a02)
(a12) edge node[pos=.4] {\umu 1} (a22)
(a13) edge node{\au} (a03)
(a13) edge[transform canvas={xshift=-.8ex}] node[swap,pos=.4] {\umu 1} (a23)
(a13) edge[transform canvas={xshift=0ex}] node[pos=.4] {(1\mu)\upbe} (a23)
(a03) edge node[swap,pos=.4,inner sep=0pt] {1\mu} (a23)
;
\end{tikzpicture}
\end{equation}
The following statements hold for the diagram \cref{rsg_inv_ii}.
\begin{itemize}
\item Along the top-right boundary, the composite of the first three arrows is the identity morphism of $X\gsp$ by the right unity \cref{gmonoid_axioms} of $(\gsp,\mu,\eta)$ \pcref{def:g_sphere}.  Thus, the top-right boundary composite is equal to $\psma$.
\item The top left triangle commutes by the right unity property \cref{moncat-other-unit-axioms} of the symmetric monoidal category \pcref{def:IU_smgtop}
\[(\IUT,\smau,\iu,\au,\ellu,\ru,\beu).\]
\item The top right rectangle commutes by the naturality of $\au$. 
\item The lower left square commutes by the naturality of $\ru$.
\item The lower middle square commutes by the functoriality of $\smau$.
\item The lower right triangle commutes by the fact that $(\gsp,\mu)$ is a commutative $G$-monoid.
\item By \cref{rsg_coequal}, the parallel arrows 
\[\umu 1 \andspace (1\mu)\upbe = (1_X \smau \mu)\upbe\]
are coequalized by $\psma$.
\end{itemize}
Thus, the left-bottom composite in the diagram \cref{rsg_inv_ii} is equal to $\psma$.
\end{enumerate}
This proves that the right monoidal unitor $\rsg$ is a $\Gtop$-natural isomorphism.  

An analogous argument proves that the left monoidal unitor $\ellsg$ \cref{gsp_left_unitor} is a $\Gtop$-natural isomorphism.  This finishes the proof that the data
\[(\smasg, \gsp, \asg, \ellsg, \rsg, \bsg)\]
are well defined.

\parhead{Symmetric monoidal $\Gtop$-category axioms}.  The braiding $\bsg$ \cref{gsp_braiding} and the monoidal associator $\asg$ \cref{gsp_associator} for $\GSp$ are induced by, respectively, $\beu$ \cref{IU_braiding} and $\au$ \cref{IU_au} for $\IUT$, as explained in \cref{gsp_braiding_coequal,asg_colimit} above.  By the universal property of coequalizers and colimits, $\GSp$ satisfies each of
\begin{itemize}
\item the pentagon axiom \cref{vmonoidal-pentagon-axiom},
\item the left hexagon axiom \cref{hexagon-bvmL},
\item the right hexagon axiom \cref{hexagon-bvmR}, and
\item the symmetry axiom \cref{vmonoidal-symmetry}
\end{itemize}
because $\IUT$ satisfies the corresponding axiom \pcref{iut_smgtop}. 

The unity axiom \cref{vmonoidal-unit} for $\GSp$ states that, for $\gsp$-modules $(X,\umu^X)$ and $(Y,\umu^Y)$, the following diagram in $\GSp$ commutes. 
\begin{equation}\label{gspectra_midunity}
\begin{tikzpicture}[vcenter]
\def\h{4}
\draw[0cell]
(0,0) node (a1) {(X \smasg \gsp) \smasg Y}
(a1)++(.5*\h,1) node (a2) {X \smasg (\gsp \smasg Y)}
(a1)++(\h,0) node (a3) {X \smasg Y}
;
\draw[1cell=.9]
(a1) edge node {\rsg_X \smasg 1_Y} (a3)
(a1) [rounded corners=2pt] |- node[pos=.3] {\asg} (a2)
;
\draw[1cell=.9]
(a2) [rounded corners=2pt] -| node[pos=.7] {1_X \smasg \ellsg_Y} (a3)
;
\end{tikzpicture}
\end{equation}
To prove that the diagram \cref{gspectra_midunity} commutes, recall the following facts.
\begin{itemize}
\item The right monoidal unitor $\rsg_X$ \cref{rsg_x} is induced by the right $\gsp$-action $\umu^X$.
\item The left monoidal unitor $\ellsg_Y$ \cref{ellsg_x} is induced by the composite $\umu^Y \beu_{\gsp,Y}$.
\item A morphism of $\gsp$-modules is determined by its underlying $\IU$-morphism.
\end{itemize}
Using these facts and the universal property of colimits, to prove that \cref{gspectra_midunity} commutes, it suffices to prove that the following diagram in $\IUT$ commutes.
\[\begin{tikzpicture}
\def\v{1.4}
\draw[0cell=.9]
(0,0) node (a11) {(X \smau \gsp) \smau Y}
(a11)++(3.8,0) node (a12) {X \smau Y}
(a12)++(3.2,0) node (a13) {X \smasg Y}
(a11)++(0,\v) node (a01) {X \smau (\gsp \smau Y)}
(a12)++(0,\v) node (a02) {X \smau (Y \smau \gsp)}
(a13)++(0,\v) node (a03) {X \smau Y}
;
\draw[1cell=.8]
(a11) edge node {\umu^X \smau 1_Y} (a12)
(a12) edge node {\psma} (a13)
(a11) edge node {\au} (a01)
(a01) edge node {1_X \smau \beu} (a02)
(a02) edge node {1_X \smau \umu^Y} (a03)
(a03) edge node {\psma} (a13)
;
\end{tikzpicture}\]
This diagram commutes by the coequalizer \cref{gsp_sma_coequal}.  This finishes the proof that $\GSp$ is a symmetric monoidal $\Gtop$-category.
\end{proof}

\begin{remark}[Right Monoidal Unitor]\label{rk:rsg_coequal}
The diagrams \cref{rsg_inv_i,rsg_inv_ii} prove that the right monoidal unitor 
\[X \smasg \gsp \fto{\rsg_X} X,\]
which is induced by the right $\gsp$-action $\umu^X \cn X \smau \gsp \to X$ via the coequalizer in \cref{rsg_coequal}, is invertible with inverse given by the composite in \cref{rsg_inverse}.  This fact is a special case of a more general fact about monads and their algebras \cite[Lemma 4.3.3]{borceux2}: 
\begin{quote}
Given a monad $(T,\mu)$ on a category $\C$ and a $T$-algebra $(x,\umu)$, the diagram
\[\begin{tikzpicture}[vcenter]
\draw[0cell=1]
(0,0) node (a1) {(TTx,\mu_{Tx})}
(a1)++(3.3,0) node (a2) {(Tx,\mu_x)}
(a2)++(2.5,0) node (a3) {(x,\umu)}
;
\draw[1cell=.9]
(a1) edge[transform canvas={yshift=.4ex}] node {\mu_x} (a2)
(a1) edge[transform canvas={yshift=-.5ex}] node[swap] {T\umu} (a2)
(a2) edge node {\umu} (a3)
;
\end{tikzpicture}\]
is a coequalizer in the category of $T$-algebras.  Moreover, this coequalizer forgets to an absolute coequalizer in $\C$.
\end{quote}
We gave an explicit proof above for $\gsp$-modules for the reader's convenience.
\end{remark}

\chapter{From $\Gskg$-Categories to $\Gskg$-Spaces}
\label{ch:ggtop}
The first step of our $G$-equivariant algebraic $K$-theory $\Kgo$ is the $\Gcat$-multifunctor (\cref{thm:Jgo_multifunctor})
\[\MultpsO \fto{\Jgo} \GGCat,\]
where $\MultpsO$ and $\GGCat$ are the $\Gcat$-multicategories constructed in, respectively, \cref{thm:multpso,expl:ggcat_gcatenr}.  For an arbitrary group $G$, this chapter constructs the second part of $\Kgo$.  It is a symmetric monoidal functor between symmetric monoidal closed categories
\[\GGCat \fto{\clast} \GGTop\]
induced by the classifying space functor $\cla$.  By change of enrichment and the closed structures of its domain and codomain, $\clast$ induces a symmetric monoidal $\Gtop$-functor between symmetric monoidal $\Gtop$-categories.  Furthermore, via the endomorphism construction, $\clast$ yields a $\Gtop$-multifunctor between $\Gtop$-multicategories.  See \cref{thm:ggcat_ggtop}.

\summary
The following table summarizes the symmetric monoidal functor $\clast$, along with its domain and codomain symmetric monoidal closed categories.
\begin{center}
\resizebox{.9\width}{!}{%
{\renewcommand{\arraystretch}{1.3}%
{\setlength{\tabcolsep}{1ex}
\begin{tabular}{c|cr|cr|cr}
& $\GGCat$ & \eqref{def:GGCat} & $\GGTop$ & \eqref{def:ggtop_smc} & $\clast$ & \eqref{thm:ggcat_ggtop} \\ \hline
objects & $\Gskg$-categories & \eqref{expl:ggcat_obj} & $\Gskg$-spaces & \eqref{expl:ggtop_obj} & \multicolumn{2}{c}{\eqref{clast_obj}} \\
morphisms & \multicolumn{2}{c|}{\eqref{expl:ggcat_mor}} & \multicolumn{2}{c|}{\eqref{expl:ggtop_mor}} & \multicolumn{2}{c}{\eqref{clast_mor}} \\
monoidal unit & \multicolumn{2}{c|}{\eqref{expl:ggcat_unit}} & \multicolumn{2}{c|}{\eqref{expl:ggtop_unit}} & \multicolumn{2}{c}{\eqref{expl:clast_unit}} \\
monoidal product & \multicolumn{2}{c|}{\eqref{expl:ggcat_smag}} & \multicolumn{2}{c|}{\eqref{expl:ggtop_smag}} & \multicolumn{2}{c}{\eqref{expl:clast_monoidal}} \\
internal hom & \multicolumn{2}{c|}{\eqref{expl:ggcat_brkst}} & \multicolumn{2}{c|}{\eqref{expl:ggtop_brkst}} &&\\
\end{tabular}}}}
\end{center}

\connection
The domain of $\clast$ is the symmetric monoidal closed category $\GGCat$ discussed in \cref{sec:ggcat_sm}.  In \cref{ch:semg}, the symmetric monoidal $\Gtop$-category $\GGTop$ is the domain of the last part of $\Kgo$, which is a unital symmetric monoidal $\Gtop$-functor
\[\GGTop \fto{\Kg} \GSp\]
that lands in the symmetric monoidal $\Gtop$-category of orthogonal $G$-spectra \pcref{gspectra_smgtop}.  

\organization
This chapter consists of the following sections.

\secname{sec:ggtop}  This section defines and unravels the symmetric monoidal closed category $\GGTop$, whose objects are called \emph{$\Gskg$-spaces}.  It is the topological analogue of $\GGCat$ \pcref{def:GGCat}.  An important point to keep in mind is that, as a 1-category, the morphisms in $\GGTop$ have components that are $G$-equivariant \pcref{expl:ggtop_mor}.  However, the component objects of the internal hom $\Gskg$-spaces of the closed structure of $\GGTop$, denoted $\brkstg{f}{f'}\angordm$ in \cref{ggtop_inthom_comp}, involve pointed morphisms that are \emph{not} necessarily $G$-equivariant.

\secname{sec:ggcat_ggtop}  This section first constructs the symmetric monoidal functor
\[\GGCat \fto{\clast} \GGTop\]
induced by the classifying space functor $\cla$, along with its induced symmetric monoidal $\Gtop$-functor and $\Gtop$-multifunctor \pcref{thm:ggcat_ggtop}.  The rest of this section unravels the functor $\clast$, its unit constraint, monoidal constraint, and symmetric monoidal functor axioms.

\secname{sec:ggtop_smgtop}  This section discusses in detail the symmetric monoidal $\Gtop$-category associated to the symmetric monoidal closed category $\GGTop$.  The $\Gtop$-enrichment of $\GGTop$ uses the internal hom $\Gskg$-spaces evaluated at the empty tuple $\ang{} \in \Gsk$; see \cref{ggtop_gtop_enr}.  Thus, these hom $G$-spaces involve pointed morphisms \cref{ggtop_gtopenr_thetan} that are not generally $G$-equivariant.  The detailed description of the symmetric monoidal $\Gtop$-category $\GGTop$ in this section is used in \cref{ch:semg} to construct the last part of $\Kgo$.

\secname{sec:clast_ev_commute}  This section proves that the various symmetric monoidal $\Gtop$-categories and $\Gtop$-multicategories in \cref{thm:ggcat_ggtop} can also be constructed by changing the order of change of enrichment and the endomorphism construction.  These observations are not used in \cref{ch:semg}.

\section{The Symmetric Monoidal Closed Category of $\Gskg$-Spaces}
\label{sec:ggtop}

For an arbitrary group $G$, this section defines the symmetric monoidal closed category $\GGTop$ of $\Gskg$-spaces.  This category is the topological analogue of the symmetric monoidal closed category $\GGCat$ of $\Gskg$-categories \pcref{def:GGCat}.  In subsequent sections, the category $\GGTop$ is used as the next intermediate category of our $G$-equivariant algebraic $K$-theory symmetric monoidal $\Gtop$-functor.

\secoutline
\begin{itemize}
\item \cref{def:ggtop_smc} defines the symmetric monoidal closed category $\GGTop$, whose objects, called \emph{$\Gskg$-spaces}, are pointed functors $\Gsk \to \Gtopst$.
\item \cref{expl:ggtop_obj,expl:ggtop_mor,expl:ggtop_unit,expl:ggtop_smag,expl:ggtop_brkst} unravel, respectively, the objects, morphisms, monoidal unit, monoidal product, and internal hom of $\GGTop$.
\end{itemize}

The following definition is the topological analogue of \cref{def:GGCat}.

\begin{definition}[$\Gskg$-Spaces]\label{def:ggtop_smc}
Applying \cref{thm:Dgm-pv-convolution-hom} to
\begin{itemize}
\item the complete and cocomplete Cartesian closed category 
\[(\Gtop, \times, *, \Topg)\] 
with terminal object $*$ \cref{gtop_smclosed}, and
\item the small permutative category 
\[(\Gsk,\oplus,\ang{},\xi)\]
with null object $\vstar$ \pcref{def:Gsk,def:Gsk_permutative},
\end{itemize}
we define the complete and cocomplete symmetric monoidal closed category
\begin{equation}\label{ggtop_smc}
\big(\GGTop, \smag, \gu, \brkst\big).
\end{equation}
The objects of $\GGTop$ are called \index{G-G-space@$\Gskg$-space}\emph{$\Gskg$-spaces}. 
\end{definition}

The rest of this section unravels $\GGTop$, starting with its objects and morphisms, followed by the monoidal unit $\gu$, the monoidal product $\smag$, and the internal hom $\brkst$.  Recall the complete and cocomplete symmetric monoidal closed category $\Gtopst$ of pointed $G$-spaces and pointed $G$-morphisms \pcref{def:gtopst}.

\begin{explanation}[$\Gskg$-Spaces]\label{expl:ggtop_obj}
By \cref{def:pointed-category}, a $\Gskg$-space is a pointed functor
\begin{equation}\label{ggtop_objects}
(\Gsk, \vstar) \fto{f} (\Gtop,*).
\end{equation}
\begin{itemize}
\item $f$ sends each object $\angordm \in \Gsk$ \cref{Gsk_objects} to a $G$-space $f\angordm$ such that $f\vstar = *$. 
\item $f$ sends each morphism $\upom \cn \angordm \to \angordn$ in $\Gsk$ \cref{Gsk_morphisms} to a $G$-morphism between $G$-spaces
\begin{equation}\label{f_upom_ggtop}
f\angordm \fto{f\upom} f\angordn
\end{equation}
such that $f$ preserves identities and composition of morphisms.  
\item Each $G$-space $f\angordm$ is regarded as a pointed $G$-space with the \emph{canonical basepoint} given by the $G$-morphism
\begin{equation}\label{fm_pointed_ggtop}
f(\vstar \to \angordm) \cn f\vstar = * \to f\angordm.
\end{equation}
The functoriality of $f$ implies that each $f\upom$ is a pointed $G$-morphism.
\end{itemize}
Thus, a $\Gskg$-space is the same thing as a pointed functor
\begin{equation}\label{ggtop_obj}
(\Gsk, \vstar) \fto{f} (\Gtopst,*).
\end{equation}
This is the topological analogue of the description \cref{ggcat_obj} of a $\Gskg$-category.
\end{explanation}

\begin{explanation}[Morphisms of $\Gskg$-Spaces]\label{expl:ggtop_mor}
By \cref{def:pointed-category}, for $\Gskg$-spaces $f$ and $f'$, a morphism $\theta \cn f \to f'$ in $\GGTop$ is a pointed natural transformation as displayed below.
\begin{equation}\label{ggtop_morphisms}
\begin{tikzpicture}[vcenter]
\def\t{28}
\draw[0cell]
(0,0) node (a1) {\Gsk}
(a1)++(1.8,0) node (a2) {\phantom{\Gskel}}
(a2)++(.25,0) node (a2') {\Gtop}
;
\draw[1cell=.9]
(a1) edge[bend left=\t] node {f} (a2)
(a1) edge[bend right=\t] node[swap] {f'} (a2)
;
\draw[2cell]
node[between=a1 and a2 at .45, rotate=-90, 2label={above,\theta}] {\Rightarrow}
;
\end{tikzpicture}
\end{equation}
Such a morphism $\theta$ consists of, for each object \cref{Gsk_objects} $\angordm \in \Gsk$, an $\angordm$-component $G$-morphism between $G$-spaces
\begin{equation}\label{ggtop_mor_component}
f\angordm \fto{\theta_{\angordm}} f'\angordm
\end{equation}
such that, for each morphism \cref{Gsk_morphisms} $\upom \cn \angordm \to \angordn$ in $\Gsk$, the following naturality diagram of $G$-morphisms commutes.
\begin{equation}\label{ggtop_mor_naturality}
\begin{tikzpicture}[vcenter]
\def\v{-1.4}
\draw[0cell]
(0,0) node (a11) {f\angordm}
(a11)++(2.5,0) node (a12) {f'\angordm}
(a11)++(0,\v) node (a21) {f\angordn}
(a12)++(0,\v) node (a22) {f'\angordn}
;
\draw[1cell=.9]
(a11) edge node {\theta_{\angordm}} (a12)
(a12) edge node {f'\upom} (a22)
(a11) edge node[swap] {f\upom} (a21)
(a21) edge node {\theta_{\angordn}} (a22)
;
\end{tikzpicture}
\end{equation}
Identity morphisms and composition in $\GGTop$ are defined componentwise in $\Gtop$ using the component $G$-morphisms in \cref{ggtop_mor_component}.  As a natural transformation, $\theta$ is automatically pointed because $* \in \Gtop$ is terminal.

In \cref{ggtop_morphisms}, $\Gtop$ can be replaced by $\Gtopst$.  By \cref{fm_pointed_ggtop}, each $G$-space $f\angordm$ is canonically pointed.  Each component $G$-morphism $\theta_{\angordm}$ is also pointed by the naturality diagram \cref{ggtop_mor_naturality}, applied to the unique morphism $\vstar \to \angordm$.  Thus, a morphism $\theta \cn f \to f'$ in $\GGTop$ is the same thing as a pointed natural transformation as displayed below.
\begin{equation}\label{ggtop_mor}
\begin{tikzpicture}[vcenter]
\def\t{28}
\draw[0cell]
(0,0) node (a1) {\Gsk}
(a1)++(1.8,0) node (a2) {\phantom{\Gskel}}
(a2)++(.3,0) node (a2') {\Gtopst}
;
\draw[1cell=.9]
(a1) edge[bend left=\t] node {f} (a2)
(a1) edge[bend right=\t] node[swap] {f'} (a2)
;
\draw[2cell]
node[between=a1 and a2 at .45, rotate=-90, 2label={above,\theta}] {\Rightarrow}
;
\end{tikzpicture}
\end{equation}
This is the topological analogue of the description \cref{ggcat_mor} of a morphism of $\Gskg$-categories.
\end{explanation}

\begin{explanation}[Monoidal Unit of $\GGTop$]\label{expl:ggtop_unit}
By \cref{reindexing_functor,ordn_empty,fangpsi,smash_Gskobjects,ptdayunit}, the monoidal unit of $\GGTop$ is the pointed functor
\begin{equation}\label{ggtop_unit}
(\Gsk, \vstar) \fto{\gu} (\Gtopst,*)
\end{equation}
given on objects by
\begin{equation}\label{gu_angordm_ggtop}
\begin{split}
\gu\angordm &= \txwedge_{\Gskpunc(\ang{}, \angordm)} \, \stplus \\
&\iso \txsma_{i \in \ufs{p}}\, \ordm_i = \sma\angordm
\end{split}
\end{equation}
for $\angordm = \ang{\ordm_i}_{i \in \ufs{p}} \in \Gsk$.
\begin{itemize}
\item In \cref{gu_angordm_ggtop}, $\Gskpunc$ denotes the set of nonzero morphisms in $\Gsk$, and $\stplus$ is the smash unit in $\Gtopst$.  An empty wedge is defined to be $*$.
\item The isomorphism in \cref{gu_angordm_ggtop} follows from the fact that each nonzero morphism 
\[\ang{} \to \angordm \inspace \Gsk\] 
is uniquely determined by an element in the unpointed finite set $\ufs{m}_i$, one for each $i \in \ufs{p}$.
\end{itemize}
The pointed finite set $\sma\angordm$ is regarded as a discrete pointed $G$-space with the trivial $G$-action.  On morphisms, $\gu$ is given by the assignment $\upom \mapsto \!\smam\upom$ defined in \cref{smash_fpsi}.
\end{explanation}

\begin{explanation}[Monoidal Product of $\GGTop$]\label{expl:ggtop_smag}
By \cref{ptdayconv}, for two $\Gskg$-spaces $f$ and $f'$ \cref{ggtop_obj}, their monoidal product is the $\Gskg$-space given by the coend 
\begin{equation}\label{ggtop_ptday}
f \smag f' = \ecint^{(\angordm, \angordmp) \in \Gsk^2} 
\mquad \bigvee_{\Gskpunc(\angordm \oplus \angordmp, -)} \mquad (f\angordm \sma f'\angordmp)
\end{equation}
taken objectwise in $\Gtopst$.  An empty wedge and $(f \smag f')(\vstar)$ are both defined to be a single point.  

\parhead{Components}. 
For each object $\angordn \in \Gsk$, the pointed $G$-space
\begin{equation}\label{smag_ptspace}
(f \smag f')\angordn = \ecint^{(\angordm, \angordmp) \in \Gsk^2} 
\mquad \bigvee_{\Gskpunc(\angordm \oplus \angordmp, \angordn)} \mquad (f\angordm \sma f'\angordmp)
\end{equation}
is a quotient of the wedge
\begin{equation}\label{smag_ptspace_wedge}
\bigvee_{(\angordm, \angordmp) \in \Gsk^2} \,\bigvee_{\Gsk(\angordm \oplus \angordmp, \angordn)} \mquad (f\angordm \sma f'\angordmp).
\end{equation}
Each point of $(f \smag f')\angordn$ is represented by a triple
\begin{equation}\label{smag_ptspace_rep}
\by = \big(\angordm \oplus \angordmp \fto{\upom} \angordn \in \Gsk; x \in f\angordm ; x' \in f'\angordmp \big)
\end{equation}
for some pair of objects $(\angordm, \angordmp) \in \Gsk^2$.  The triple $\by$ represents the basepoint of $(f \smag f')\angordn$ if 
\begin{itemize}
\item $\upom$ is the 0-morphism, factoring through $\vstar \in \Gsk$, or
\item either $x$ or $x'$ is the basepoint. 
\end{itemize} 
The coend \cref{smag_ptspace} identifies, for each quintuple
\[\bw = 
\scalebox{.8}{$\big(\angordl \fto{\upla} \angordm ; \angordlp \fto{\upla'} \angordmp ; \angordm \oplus \angordmp \fto{\upom} \angordn ; x \in f\angordl ; x' \in f'\angordlp \big)$}\]
with $\upla, \upla', \upom \in \Gsk$, the following two triples in the wedge \cref{smag_ptspace_wedge}.
\begin{equation}\label{smag_ptspace_dzerone}
\begin{split}
\dzero\bw &= \big(\upom \circ (\upla \oplus \upla') ; x; x'\big)\\
\done\bw &= \big(\upom ; (f\upla)(x) \in f\angordm ; (f'\upla')(x') \in f'\angordmp)
\end{split}
\end{equation}

\parhead{$G$-action}. 
The group $G$ acts on the pointed $G$-space $(f \smag f')\angordn$ via the diagonal $G$-action on $f\angordm \sma f'\angordmp$.  In other words, for $g \in G$ and a representative $\by$ \cref{smag_ptspace_rep}, the $g$-action is given by
\begin{equation}\label{smag_rep_gaction}
g\by = (\upom ; gx ; gx').
\end{equation}
The identification \cref{smag_ptspace_dzerone} is invariant under the $G$-action \cref{smag_rep_gaction}, which means
\[g(\dzero\bw) = g(\done\bw),\]
because $f\upla$ and $f'\upla'$ are $G$-equivariant \cref{f_upom_ggtop}.

\parhead{External characterization}. 
By  \cref{expl:pointedday}, given another $\Gskg$-space $f''$, a morphism
\[f \smag f' \fto{\theta} f'' \inspace \GGTop\]
is a pointed natural transformation with component pointed $G$-morphisms between pointed $G$-spaces
\begin{equation}\label{ggtop_sma_ext}
f\angordm \sma f'\angordmp \fto{\theta_{\angordm,\angordmp}} f''(\angordm \oplus \angordmp)
\end{equation}
for $(\angordm, \angordmp) \in \Gsk^2$ such that, for each pair of morphisms \cref{Gsk_morphisms}
\[\angordm \fto{\upom} \angordn \andspace 
\angordmp \fto{\upom'} \angordnp\]
in $\Gsk$, the following diagram of pointed $G$-morphisms commutes.
\begin{equation}\label{ggtop_sma_nat}
\begin{tikzpicture}[vcenter]
\def\v{-1.4}
\draw[0cell=1]
(0,0) node (x11) {f\angordm \sma f'\angordmp}
(x11)++(4.5,0) node (x12) {f''(\angordm \oplus \angordmp)}
(x11)++(0,\v) node (x21) {f\angordn \sma f'\angordnp}
(x12)++(0,\v) node (x22) {f''(\angordn \oplus \angordnp)}
;
\draw[1cell=.9] 
(x11) edge node {\theta_{\angordm,\angordmp}} (x12)
(x21) edge node {\theta_{\angordn,\angordnp}} (x22)
(x11) edge[transform canvas={xshift=1em}] node[swap] {f\upom \sma f'\upom'} (x21)
(x12) edge[transform canvas={xshift=-1.5em}] node {f''(\upom \oplus \upom')} (x22)
;
\end{tikzpicture}
\end{equation}
This description of $f \smag f'$ is called the \emph{external characterization}.
\end{explanation}

\begin{explanation}[Internal Hom of $\GGTop$]\label{expl:ggtop_brkst}
By \cref{ptdayhom}, for $\Gskg$-spaces $f$ and $f'$, the internal hom $\Gskg$-space is given by the end 
\begin{equation}\label{ggtop_inthom}
\brkstg{f}{f'} = \ecint_{\angordn \in \Gsk} \Topgst\big(f\angordn, f'(- \oplus \angordn)\big)
\end{equation}
taken in $\Gtopst$, where $\Topgst$ is the internal hom for $\Gtopst$ \cref{Gtopst_smc}.

\parhead{Objects}.  In more details, for each object \cref{Gsk_objects} $\angordm \in \Gsk$, the pointed $G$-space
\begin{equation}\label{ggtop_inthom_comp}
\brkstg{f}{f'}\angordm 
= \ecint_{\angordn \in \Gsk} \Topgst\big(f\angordn, f'(\angordm \oplus \angordn)\big)
\end{equation}
is a $G$-subspace of the product
\begin{equation}\label{ggtop_inthom_topology}
\prod_{\angordn \in \Gsk} \Topgst\big(f\angordn, f'(\angordm \oplus \angordn) \big).
\end{equation}
A point $\theta \in \brkstg{f}{f'}\angordm$ consists of $\angordn$-component pointed morphisms between pointed spaces
\begin{equation}\label{ggtop_inthom_theta}
f\angordn \fto{\theta_{\angordn}} f'(\angordm \oplus \angordn) \forspace \angordn \in \Gsk
\end{equation}
such that, for each morphism \cref{Gsk_morphisms} $\upom \cn \angordn \to \angordl$ in $\Gsk$, the diagram 
\begin{equation}\label{ggtop_inthom_theta_nat}
\begin{tikzpicture}[vcenter]
\def\v{-1.4}
\draw[0cell]
(0,0) node (a11) {f\angordn}
(a11)++(3,0) node (a12) {f'(\angordm \oplus \angordn)}
(a11)++(0,\v) node (a21) {f\angordl}
(a12)++(0,\v) node (a22) {f'(\angordm \oplus \angordl)}
;
\draw[1cell=.9]
(a11) edge node {\theta_{\angordn}} (a12)
(a12) edge[transform canvas={xshift=-2em}] node {f'(1_{\angordm} \oplus \upom)} (a22)
(a11) edge node[swap] {f\upom} (a21)
(a21) edge node {\theta_{\angordl}} (a22)
;
\end{tikzpicture}
\end{equation}
of pointed morphisms commutes.  Since the components $\theta_{\angordn}$ are \emph{not} generally $G$-equivariant, \cref{ggtop_inthom_theta_nat} is not generally a diagram in $\Gtop$.  In particular, $\theta$ is not generally a natural transformation $f \to f'(\angordm \oplus -)$ for functors $\Gsk \to \Gtop$.  The basepoint of $\brkstg{f}{f'}\angordm$ has each component pointed morphism given by the constant morphism at the basepoint of the codomain.  

\parhead{$G$-action}.  The $G$-action on the pointed $G$-space $\brkstg{f}{f'}\angordm$ is given componentwise by the conjugation $G$-action.  In other words, for $g \in G$ and $\angordn \in \Gsk$, the $\angordn$-component of $g \cdot \theta$ is given by the following composite pointed morphism.
\begin{equation}\label{ggtop_inthom_theta_g}
\begin{tikzpicture}[vcenter]
\def\v{-1.5} \def\t{18}
\draw[0cell]
(0,0) node (a11) {\phantom{f'}}
(a11)++(2.5,0) node (a12) {\phantom{f'}}
(a11)++(0,\v) node (a21) {\phantom{f'}}
(a12)++(0,\v) node (a22) {\phantom{f'}}
(a11)++(-.2,0) node (a11') {f\angordn}
(a12)++(.9,0) node (a12') {f'(\angordm \oplus \angordn)}
(a21)++(-.2,0) node (a21') {f\angordn}
(a22)++(.9,0) node (a22') {f'(\angordm \oplus \angordn)}
;
\draw[1cell=.9]
(a11) edge node {(g \cdot \theta)_{\angordn}} (a12)
(a11) edge node[swap] {\ginv} (a21)
(a21) edge node {\theta_{\angordn}} (a22)
(a22) edge node[swap] {g} (a12)
;
\end{tikzpicture}
\end{equation}
The components of $g \cdot \theta$ satisfy the property \cref{ggtop_inthom_theta_nat} by 
\begin{itemize}
\item the same property for $\theta$ and
\item the $G$-equivariance \cref{f_upom_ggtop} of $f\upom$ and $f'(1_{\angordm} \oplus \upom)$.
\end{itemize}
This finishes the description of the pointed $G$-space $\brkstg{f}{f'}\angordm$ in \cref{ggtop_inthom_comp}.

\parhead{Morphisms}.  The internal hom $\Gskg$-space $\brkstg{f}{f'}$ sends a  morphism $\uprho \cn \angordm \to \angordmp$ in $\Gsk$ to the pointed $G$-morphism between pointed $G$-spaces defined by the following commutative diagram.
\begin{equation}\label{ggtop_inthom_gfun}
\begin{tikzpicture}[vcenter]
\def\v{-1.5}
\draw[0cell=.9]
(0,0) node (a11) {\brkstg{f}{f'}\angordm}
(a11)++(0,\v) node (a21) {\brkstg{f}{f'}\angordmp}
(a11)++(4,0) node (a12) {\txint_{\angordn \in \Gsk} \Topgst\big(f\angordn, f'(\angordm \oplus \angordn)\big)}
(a12)++(0,\v) node (a22) {\txint_{\angordn \in \Gsk} \Topgst\big(f\angordn, f'(\angordmp \oplus \angordn)\big)}
;
\draw[1cell=.9]
(a11) edge[equal] (a12)
(a21) edge[equal] (a22)
(a11) edge[transform canvas={xshift=1em}] node[swap] {\brkstg{f}{f'} \uprho} (a21)
(a12) edge[transform canvas={xshift=-4.5em}] node {\txint_{\angordn} \Topgst(1_{f\angordn}, f'(\uprho \oplus 1_{\angordn}))} (a22)
;
\end{tikzpicture}
\end{equation}
For a point $\theta \in \brkstg{f}{f'}\angordm$, with components $\theta_{\angordn}$ \cref{ggtop_inthom_theta} for $\angordn \in \Gsk$, the point
\[(\brkstg{f}{f'} \uprho)\theta \in \brkstg{f}{f'}\angordmp\]
has $\angordn$-component pointed morphism given by the composite
\begin{equation}\label{ggtop_inthom_gfun_comp}
\begin{tikzpicture}[baseline={(a1.base)}]
\def\t{28}
\draw[0cell]
(0,0) node (a1) {\phantom{f'}}
(a1)++(1.8,0) node (a2) {\phantom{f'}}
(a1)++(-.2,0) node (a1') {f\angordn}
(a2)++(.9,0) node (a2') {f'(\angordm \oplus \angordn)}
(a2')++(.9,0) node (a2'') {\phantom{f'}}
(a2'')++(2.7,0) node (a3) {\phantom{f'}}
(a3)++(.95,0) node (a3') {f'(\angordmp \oplus \angordn).}
;
\draw[1cell=.85]
(a1) edge node {\theta_{\angordn}} (a2)
(a2'') edge node {f'(\uprho \oplus 1_{\angordn})} (a3)
;
\end{tikzpicture}
\end{equation}
\begin{itemize}
\item The property \cref{ggtop_inthom_theta_nat} holds for $(\brkstg{f}{f'} \uprho)\theta$ by
\begin{itemize}
\item the same property for $\theta$ and
\item the functoriality of $f'$ and $\oplus$ (\cref{Gsk_permutative}).
\end{itemize}
\item The $G$-equivariance of $\brkstg{f}{f'} \uprho$ means
\begin{equation}\label{ffrho_geq}
(\brkstg{f}{f'} \uprho)(g \cdot \theta) = g \cdot \big((\brkstg{f}{f'} \uprho)\theta \big)
\end{equation}
for $g \in G$.  By \cref{ggtop_inthom_theta_g,ggtop_inthom_gfun_comp}, the equality \cref{ffrho_geq} means
\[\begin{split}
& f'(\uprho \oplus 1_{\angordn}) \circ (g \circ \theta_{\angordn} \circ \ginv) \\
&= g \circ \big( f'(\uprho \oplus 1_{\angordn}) \circ \theta_{\angordn}\big) \circ \ginv
\end{split}\]
for $\angordn \in \Gsk$.  This equality holds by the $G$-equivariance \cref{f_upom_ggtop} of the morphism $f'(\uprho \oplus 1_{\angordn})$.
\item The functoriality of $\brkstg{f}{f'}$ with respect to morphisms $\uprho \in \Gsk$ follows from 
\begin{itemize}
\item the universal property of ends and
\item the functoriality of $f'$ and $\oplus$.
\end{itemize}
\end{itemize}
This finishes the description of the internal hom $\Gskg$-space $\brkstg{f}{f'}$.
\end{explanation}

\section{The Symmetric Monoidal Functor from $\Gskg$-Categories\\ to $\Gskg$-Spaces}
\label{sec:ggcat_ggtop}

This section discusses the change-of-base symmetric monoidal $\Gtop$-functor 
\[\GGCat \fto{\clast} \GGTop\]
induced by the classifying space functor \cref{cla_gcat_gtop}
\[(\Gcat,\times,\boldone) \fto{(\cla, \clatwo,\clazero)} (\Gtop, \times, *).\]
The symmetric monoidal $\Gtop$-functor $\clast$ is the middle step of our $G$-equivariant algebraic $K$-theory multifunctor.

\secoutline
\begin{itemize}
\item \cref{thm:ggcat_ggtop} records the symmetric monoidal functor $\clast$ from $\GGCat$ to $\GGTop$---which are symmetric monoidal closed categories by, respectively, \cref{def:GGCat,def:ggtop_smc}---along with the induced
\begin{itemize}
\item symmetric monoidal $\Gtop$-functor $(\clast')_{\evang}$ and
\item $\Gtop$-multifunctor $\End((\clast')_{\evang})$.
\end{itemize}
\item \cref{expl:clast_functor,expl:clast_unit,expl:clast_monoidal,expl:clast_symmetry} unravel, respectively, the functor $\clast$, the unit constraint of $\clast$, the monoidal constraint of $\clast$, and the symmetric monoidal functor axioms for $\clast$.
\end{itemize}

The notation in \cref{clastpev} below is explained in the proof that follows; see also \cref{not:clast'}.

\begin{theorem}\label{thm:ggcat_ggtop}
The following statements hold for an arbitrary group $G$.  
\begin{enumerate}
\item\label{ggcat_ggtop_i} Post-composing and post-whiskering with the classifying space functor $\cla$ in \cref{cla_gcat_gtop} induce a symmetric monoidal functor\index{classifying space}\index{G-G-category@$\Gskg$-category}
\begin{equation}\label{clast}
\GGCat \fto{\clast} \GGTop
\end{equation}
between the symmetric monoidal closed categories in \cref{def:GGCat,def:ggtop_smc}.
\item\label{ggcat_ggtop_ii} By changing enrichment, $\clast$ induces a symmetric monoidal $\Gtop$-functor
\begin{equation}\label{clastpev}
\big(\uof{\GGCat}_{\clast}\big)_{\evang} \fto{(\clast')_{\evang}} (\uof{\GGTop})_{\evang}
\end{equation}
with the same object assignment as $\clast$.
\item\label{ggcat_ggtop_iii} Applying the endomorphism construction \pcref{def:EndF} to $(\clast')_{\evang}$ yields a $\Gtop$-multifunctor
\begin{equation}\label{End_clastpev}
\End \big(\big(\uof{\GGCat}_{\clast}\big)_{\evang}\big) \fto{\End ((\clast')_{\evang})} 
\End \big((\uof{\GGTop})_{\evang}\big)
\end{equation}
with the same object assignment as $\clast$.
\end{enumerate}
\end{theorem}

\begin{proof}
\parhead{Assertion \pcref{ggcat_ggtop_i}}.  This is an instance of \cref{thm:Dgm_f} in the following context.
\begin{itemize}
\item $\V$ is the complete and cocomplete Cartesian closed category 
\[(\Gcat, \times, \boldone, \Catg)\]
with terminal object $\boldone$ in \cref{expl:Gcat_closed}.
\item $\W$ is the complete and cocomplete Cartesian closed category 
\[(\Gtop, \times, *, \Topg)\] 
with terminal object $*$ in \cref{def:Gtop}. 
\item $\Dgm$ is the small permutative category 
\[(\Gsk, \oplus, \ang{}, \xi)\]
with null object $\vstar$ in \cref{def:Gsk_permutative}.
\item $f$ is the classifying space functor \cref{cla_gcat_gtop}
\[(\Gcat,\times,\boldone) \fto{(\cla, \clatwo,\clazero)} (\Gtop, \times, *),\]
which is strong symmetric monoidal and sends $\boldone$ to $*$.
\end{itemize}

\parhead{Assertion \pcref{ggcat_ggtop_ii}}.  The desired symmetric monoidal $\Gtop$-functor $(\clast')_{\evang}$ is obtained from the symmetric monoidal functor $\clast$ in two steps as follows.

First, by \cref{thm:change-enrichment} \pcref{change-enr-ii}, $\clast$ yields a symmetric monoidal $(\GGTop)$-functor
\begin{equation}\label{clastp}
\uof{\GGCat}_{\clast} \fto{\clast'} \uof{\GGTop}
\end{equation}
with the domain and codomain given as follows.
\begin{itemize}
\item The codomain $\uof{\GGTop}$ is the symmetric monoidal $(\GGTop)$-category obtained from the symmetric monoidal closed category $\GGTop$ \pcref{def:ggtop_smc} by applying \cref{theorem:v-closed-v-sm}.
\item In the domain of $\clast'$, $\uof{\GGCat}$ is the symmetric monoidal $(\GGCat)$-category obtained from the symmetric monoidal closed category $\GGCat$ \pcref{def:GGCat} by applying \cref{theorem:v-closed-v-sm}.
\item The domain of $\clast'$, $\uof{\GGCat}_{\clast}$, is the symmetric monoidal $(\GGTop)$-category obtained from $\uof{\GGCat}$ by changing enrichment along the symmetric monoidal functor $\clast$ in \cref{clast}, using \cref{thm:change-enrichment} \pcref{change-enr-i}.
\end{itemize}
Next, evaluation at the monoidal unit $\ang{} \in \Gsk$ defines a symmetric monoidal functor \cref{evtu_unpt}
\begin{equation}\label{evang_gtop}
\GGTop \fto{\evang} \Gtop.
\end{equation}
By \cref{thm:change-enrichment} \pcref{change-enr-i'}, changing enrichment of $\clast'$ along $\evang$ yields the desired symmetric monoidal $\Gtop$-functor
\[\big(\uof{\GGCat}_{\clast}\big)_{\evang} \fto{(\clast')_{\evang}} (\uof{\GGTop})_{\evang}.\]

\parhead{Assertion \pcref{ggcat_ggtop_iii}}.  This is an instance of \cref{EndF_multi}, applied to the symmetric monoidal $\Gtop$-functor $(\clast')_{\evang}$.
\end{proof}

\begin{notation}\label{not:clast'}
Using the notational convention for change of enrichment in \cref{expl:ggcat_gcatenr}, each of
\begin{itemize} 
\item the symmetric monoidal $(\GGTop)$-functor $\clast'$ in \cref{clastp},
\item the symmetric monoidal $\Gtop$-functor $(\clast')_{\evang}$ in \cref{clastpev}, and
\item the $\Gtop$-multifunctor $\End ((\clast')_{\evang})$ in \cref{End_clastpev}
\end{itemize}
is abbreviated to $\clast \cn \GGCat \to \GGTop$.  
\end{notation}

\cref{expl:clast_functor,expl:clast_unit,expl:clast_monoidal,expl:clast_symmetry} below unravel the symmetric monoidal functor 
\[\GGCat \fto{(\clast, \clasttwo, \clastzero)} \GGTop\]
in \cref{clast}, where $\clastzero$ and $\clasttwo$ denote, respectively, the unit constraint \cref{monoidal_unit} and the monoidal constraint \cref{monoidal_constraint} of $\clast$.

\begin{explanation}[The Functor $\clast$]\label{expl:clast_functor}
The object and morphism assignments of $\clast$ \cref{clast} are given as follows.

\parhead{Objects}.  An object in $\GGCat$ is a $\Gskg$-category \pcref{expl:ggcat_obj}, which means a pointed functor
\[(\Gsk,\vstar) \fto{f} (\Gcatst,\boldone).\]
The functor $\clast$ sends $f$ to the $\Gskg$-space $\clast f$ \pcref{expl:ggtop_obj} given by the composite pointed functor
\begin{equation}\label{clast_obj}
\begin{tikzpicture}[baseline={(a1.base)}]
\def\b{.6}
\draw[0cell]
(0,0) node (a1) {(\Gsk,\vstar)}
(a1)++(2.4,0) node (a2) {(\Gcatst,\boldone)}
(a2)++(2.9,0) node (a3) {(\Gtopst,*).}
;
\draw[1cell=.9]
(a1) edge node {f} (a2)
(a2) edge node {\cla} (a3)
;
\draw[1cell=.9]
(a1) [rounded corners=2pt, shorten <=0ex] |- ($(a2)+(-1,\b)$)
-- node {\clast f} ($(a2)+(1.25,\b)$) -| (a3)
;
\end{tikzpicture}
\end{equation}
This is well defined because $\cla\boldone = *$.

\parhead{Morphisms}.  A morphism in $\GGCat$ \pcref{expl:ggcat_mor} is a pointed natural transformation
\[\begin{tikzpicture}[vcenter]
\def\t{28}
\draw[0cell]
(0,0) node (a1) {\Gsk}
(a1)++(1.8,0) node (a2) {\phantom{\Gskel}}
(a2)++(.3,0) node (a2') {\Gcatst}
;
\draw[1cell=.9]
(a1) edge[bend left=\t] node {f} (a2)
(a1) edge[bend right=\t] node[swap] {f'} (a2)
;
\draw[2cell]
node[between=a1 and a2 at .43, rotate=-90, 2label={above,\theta}] {\Rightarrow}
;
\end{tikzpicture}\]
where pointedness is automatic, since $\boldone \in \Gcatst$ is terminal.  The functor $\clast$ sends $\theta$ to the morphism 
\[\clast f \fto{\clast\theta} \clast f'\] 
in $\GGTop$ given by the whiskering
\begin{equation}\label{clast_mor}
\begin{tikzpicture}[baseline={(a1.base)}]
\def\t{28}
\draw[0cell]
(0,0) node (a1) {\Gsk}
(a1)++(1.8,0) node (a2) {\phantom{\Gsk}}
(a2)++(.3,0) node (a2') {\Gcatst}
(a2')++(.2,0) node (a2'') {\phantom{Gsk}}
(a2'')++(2.3,0) node (a3) {\Gtopst.}
;
\draw[1cell=.9]
(a1) edge[bend left=\t] node {f} (a2)
(a1) edge[bend right=\t] node[swap] {f'} (a2)
(a2'') edge node {\cla} (a3)
;
\draw[2cell]
node[between=a1 and a2 at .43, rotate=-90, 2label={above,\theta}] {\Rightarrow}
;
\end{tikzpicture}
\end{equation}
Thus, for each object $\angordm \in \Gsk$ \cref{Gsk_objects}, the $\angordm$-component of $\clast\theta$ is the pointed $G$-morphism between pointed $G$-spaces
\[\cla f\angordm \fto{\cla\theta_{\angordm}} \cla f'\angordm\]
obtained from $\theta_{\angordm}$ \cref{ggcat_mor_component} by applying $\cla$.
\end{explanation}

\begin{explanation}[Unit Constraint of $\clast$]\label{expl:clast_unit}
By \cref{expl:ggcat_unit,expl:ggtop_unit}, the monoidal units
\[\begin{split}
(\Gsk,\vstar) & \fto{\gu} (\Gcatst,\boldone) \andspace\\
(\Gsk,\vstar) & \fto{\gu} (\Gtopst,*)
\end{split}\]
in, respectively, $\GGCat$ and $\GGTop$ are given on objects $\angordm \in \Gsk$ \cref{Gsk_objects} by
\[\begin{split}
\gu\angordm &= \txwedge_{\Gskpunc(\ang{}, \angordm)} \bonep \in \Gcatst \andspace\\
\gu\angordm &= \txwedge_{\Gskpunc(\ang{}, \angordm)} \stplus \in \Gtopst. \\
\end{split}\]
The unit constraint \cref{monoidal_unit} of the symmetric monoidal functor $\clast$ \cref{clast} is the pointed natural transformation
\begin{equation}\label{clastzero}
\begin{tikzpicture}[baseline={(a1.base)}]
\def\t{30}
\draw[0cell]
(0,0) node (a1) {\Gsk}
(a1)++(2,0) node (a2) {\phantom{\Gsk}}
(a2)++(.3,0) node (a2') {\Gtopst}
;
\draw[1cell=.9]
(a1) edge[bend left=\t] node {\gu} (a2)
(a1) edge[bend right=\t] node[swap] {\clast\gu} (a2)
;
\draw[2cell]
node[between=a1 and a2 at .38, rotate=-90, 2label={above,\clastzero}] {\Rightarrow}
;
\end{tikzpicture}
\end{equation}
whose $\angordm$-component is determined by the following commutative diagrams in $\Gtopst$ for nonzero morphisms $\upom \in \Gskpunc(\ang{}, \angordm)$.
\begin{equation}\label{clastzero_angordm}
\begin{tikzpicture}[vcenter]
\def\u{-1} \def\v{-1.4}
\draw[0cell]
(0,0) node (a11) {\gu\angordm}
(a11)++(4,0) node (a12) {(\clast\gu)\angordm}
(a11)++(0,\u) node (a21) {\txwedge_{\Gskpunc(\ang{}, \angordm)} \stplus}
(a12)++(0,\u) node (a22) {\cla(\txwedge_{\Gskpunc(\ang{}, \angordm)} \bonep)}
(a21)++(0,\v) node (a31) {\stplus}
(a22)++(0,\v) node (a32) {\cla(\bonep)}
;
\draw[1cell=.9]
(a11) edge node {(\clastzero)_{\angordm}} (a12)
(a11) edge[equal] (a21)
(a12) edge[equal] (a22)
(a31) edge node {\inc_\upom} (a21)
(a31) edge node {\jm} (a32)
(a32) edge node[swap] {\cla \inc_\upom} (a22)
;
\end{tikzpicture}
\end{equation}
In \cref{clastzero_angordm}, each $\inc_\upom$ is the wedge summand inclusion corresponding to a given nonzero morphism $\upom \cn \ang{} \to \angordm$ in $\Gsk$ \cref{fangpsi}.  The bottom horizontal morphism $\jm$ is the pointed $G$-morphism whose value at the non-basepoint of $\stplus$ is given by
\[* = \cla\boldone \to \cla(\bonep),\]
which is the image under $\cla$ of the non-basepoint inclusion $\boldone \to \bonep$.
\end{explanation}

\begin{explanation}[Monoidal Constraint of $\clast$]\label{expl:clast_monoidal}
By \cref{expl:ggcat_smag}, for two $\Gskg$-categories $f$ and $f'$, their monoidal product in $\GGCat$ is given by the coend 
\begin{equation}\label{ggcat_smag}
f \smag f' = \ecint^{(\angordm, \angordmp) \in \Gsk^2} 
\mquad \bigvee_{\Gskpunc(\angordm \oplus \angordmp, -)} \mquad (f\angordm \sma f'\angordmp)
\end{equation}
taken in $\Gcatst$.  The monoidal product for $\Gskg$-spaces is given by the same formula \cref{ggtop_ptday} taken in $\Gtopst$.  The $(f,f')$-component of the monoidal constraint \cref{monoidal_constraint} of the symmetric monoidal functor $\clast$ \cref{clast} is the pointed natural transformation
\begin{equation}\label{clasttwo}
\begin{tikzpicture}[baseline={(a1.base)}]
\def\t{25}
\draw[0cell]
(0,0) node (a1) {\Gsk}
(a1)++(2.75,0) node (a2) {\phantom{\Gsk}}
(a2)++(.35,0) node (a2') {\Gtopst}
;
\draw[1cell=.85]
(a1) edge[bend left=\t] node {\clast f \smag \clast f'} (a2)
(a1) edge[bend right=\t] node[swap] {\clast(f \smag f')} (a2)
;
\draw[2cell=.9]
node[between=a1 and a2 at .3, rotate=-90, 2label={above,(\clasttwo)_{f,f'}}] {\Rightarrow}
;
\end{tikzpicture}
\end{equation}
whose $\angordn$-component, for $\angordn \in \Gsk$, is determined by the following commutative diagrams in $\Gtopst$ for nonzero morphisms 
\[\angordm \oplus \angordmp \fto{\upom} \angordn \inspace \Gsk.\]  
Each wedge in the following diagram is indexed by the set $\Gskpunc(\angordm \oplus \angordmp, \angordn)$ of nonzero morphisms.
\begin{equation}\label{clasttwo_component}
\begin{tikzpicture}[vcenter]
\def\c{1em} \def\d{.4} \def\e{4.6} \def\u{1.4} \def\t{0}
\draw[0cell=.8]
(0,0) node (a11) {\cla f\angordm \sma \cla f'\angordmp}
(a11)++(\d,\u) node (a12) {\bigvee \cla f\angordm \sma \cla f' \angordmp}
(a12)++(\e,0) node (a13) {\ecint^{(\angordm, \angordmp) \in \Gsk^2} 
\bigvee \cla f\angordm \sma \cla f' \angordmp}
(a13)++(\d,0) node (a13') {\phantom{\bigvee_{()}}}
(a13)++(\d,-\u) node (a14) {(\clast f \smag \clast f')\angordn}
(a11)++(0,-1.4) node (a21) {\cla(f\angordm \sma f' \angordmp)}
(a21)++(\d,-\u) node (a22) {\cla\big(\bigvee f\angordm \sma f' \angordmp\big)}
(a22)++(\e,0) node (a23) {\cla\Big(\ecint^{(\angordm, \angordmp) \in \Gsk^2} 
\bigvee f\angordm \sma f' \angordmp \Big)}
(a23)++(\d,0) node (a23') {\phantom{\bigvee_{()}}}
(a23)++(\d,\u) node (a24) {(\clast(f \smag f'))\angordn}
;
\draw[1cell=.8]
(a11) edge[transform canvas={xshift=\c}] node[pos=.3] {\inc_\upom} (a12)
(a12) edge[transform canvas={yshift=\t}] node {\iota} (a13)
(a13') edge[equal] (a14)
(a14) edge node {(\clasttwo)_{f,f'; \angordn}} (a24)
(a11) edge[transform canvas={xshift=\c}] node[swap] {\clatwobar} (a21)
(a21) edge[transform canvas={xshift=\c}] node[swap,pos=.3] {\cla\inc_{\upom}} (a22)
(a22) edge[transform canvas={yshift=\t}] node {\cla\iota} (a23)
(a23') edge[equal] (a24)
;
\end{tikzpicture}
\end{equation}
The arrows in \cref{clasttwo_component} are given as follows.
\begin{itemize}
\item $(\clasttwo)_{f,f'; \angordn}$ denotes the $\angordn$-component of $(\clasttwo)_{f,f'}$.
\item Each $\iota$ is part of the definition of the coend in its codomain.
\item Each $\inc_\upom$ is the wedge summand inclusion for $\upom \in \Gskpunc(\angordm \oplus \angordmp, \angordn)$.
\item $\clatwobar$ is induced by the monoidal constraint $\clatwo$ of the strong symmetric monoidal functor \cref{cla_gcat_gtop}
\[(\Gcat,\times,\boldone) \fto{(\cla,\clatwo,\clazero)} (\Gtop,\times,*)\]
as displayed in the following diagram.
\begin{equation}\label{clasttwo_one}
\begin{tikzpicture}[vcenter]
\def\v{-1.4}
\draw[0cell=.9]
(0,0) node (a11) {\cla f \angordm \times \cla f'\angordmp}
(a11)++(5,0) node (a12) {\cla f \angordm \sma \cla f'\angordmp}
(a11)++(0,\v) node (a21) {\cla(f \angordm \times f' \angordmp)}
(a12)++(0,\v) node (a22) {\cla(f \angordm \sma f' \angordmp)}
;
\draw[1cell=.9]
(a11) edge node {\pjt_{\cla f\angordm, \cla f'\angordmp}} (a12)
(a12) edge[dashed] node {\clatwobar} (a22)
(a11) edge node {\iso} node[swap,pos=.45] {\clatwo} (a21)
(a21) edge node {\cla\pjt_{f\angordm, f'\angordmp}} (a22)
;
\end{tikzpicture}
\end{equation}
The basepoints of $\cla f \angordm$ and $\cla f'\angordmp$ are both sent by the left-bottom composite in \cref{clasttwo_one} to the basepoint of $\cla(f \angordm \sma f' \angordmp)$.  The universal property of the smash product \cref{eq:smash} implies the unique existence of $\clatwobar$ such that the diagram \cref{clasttwo_one} commutes.
\end{itemize}
This finishes our description of the monoidal constraint $\clasttwo$ of the symmetric monoidal functor $\clast$ \cref{clast}.
\end{explanation}

\begin{explanation}[Axioms]\label{expl:clast_symmetry}
The symmetric monoidal functor axioms \pcref{def:monoidalfunctor} for $(\clast, \clasttwo, \clastzero)$ \cref{clast} follow from
\begin{itemize}
\item the symmetric monoidal functor axioms for $(\cla,\clatwo,\clazero)$ \cref{cla_gcat_gtop} and
\item the universal properties of pushouts and coends.
\end{itemize}  
As an illustration, we explain the symmetry axiom \cref{monoidalfunctorbraiding}, which involves the monoidal constraint $\clasttwo$ \cref{clasttwo} and the braidings $\brg$ for both $(\GGCat,\smag)$ and $(\GGTop,\smag)$ \pcref{thm:Dgm-pv-convolution-hom}.

\parhead{Braiding for pointed Day convolution}. 
First, we recall that the braiding $\brg$ for $(\GGCat,\smag)$ involves the braidings $\xi$ for $(\Gcatst,\sma)$ and $(\Gsk,\oplus)$.  More explicitly, given two $\Gskg$-categories $f$ and $f'$, the $(f,f')$-component of $\brg$ is the pointed natural isomorphism
\[\begin{tikzpicture}[baseline={(a1.base)}]
\def\t{28}
\draw[0cell]
(0,0) node (a1) {\Gsk}
(a1)++(2.,0) node (a2) {\phantom{\Gsk}}
(a2)++(.3,0) node (a2') {\Gcatst}
;
\draw[1cell=.85]
(a1) edge[bend left=\t] node {f \smag f'} (a2)
(a1) edge[bend right=\t] node[swap] {f' \smag f} (a2)
;
\draw[2cell=.9]
node[between=a1 and a2 at .35, rotate=-90, 2label={above,\brg_{f,f'}}] {\Rightarrow}
;
\end{tikzpicture}\]
whose $\angordn$-component, for $\angordn \in \Gsk$, is determined by the following commutative diagrams in $\Gcatst$ for nonzero morphisms 
\[\angordm \oplus \angordmp \fto{\upom} \angordn \inspace \Gsk.\]
\begin{equation}\label{brg_components}
\begin{tikzpicture}[vcenter]
\def\c{1em} \def\d{.4} \def\e{4.6} \def\u{1.4} \def\t{.7ex}
\draw[0cell=.8]
(0,0) node (a11) {f\angordm \sma f'\angordmp}
(a11)++(\d,\u) node (a12) {\bigvee_{\Gskpunc(\angordm \oplus \angordmp, \angordn)} \mqq f\angordm \sma f'\angordmp}
(a12)++(\e,0) node (a13) {\ecint^{(\angordm, \angordmp) \in \Gsk^2} 
\mqq \bigvee_{\Gskpunc(\angordm \oplus \angordmp, \angordn)} \mqq f\angordm \sma f'\angordmp}
(a13)++(\d,0) node (a13') {\phantom{\bigvee_{()}}}
(a13)++(\d,-\u) node (a14) {(f \smag f')\angordn}
(a11)++(0,-1.4) node (a21) {f'\angordmp \sma f\angordm}
(a21)++(\d,-\u) node (a22) {\bigvee_{\Gskpunc(\angordmp \oplus \angordm, \angordn)} \mqq f'\angordmp \sma f\angordm}
(a22)++(\e,0) node (a23) {\ecint^{(\angordmp, \angordm) \in \Gsk^2} 
\mqq \bigvee_{\Gskpunc(\angordmp \oplus \angordm, \angordn)} \mqq f'\angordmp \sma f\angordm}
(a23)++(\d,0) node (a23') {\phantom{\bigvee_{()}}}
(a23)++(\d,\u) node (a24) {(f' \smag f)\angordn}
;
\draw[1cell=.8]
(a11) edge[transform canvas={xshift=\c}] node[pos=.3] {\inc_\upom} (a12)
(a12) edge[transform canvas={yshift=\t}] node {\iota} (a13)
(a13') edge[equal, shorten <=.7ex] (a14)
(a14) edge node {\brg_{f,f'; \angordn}} (a24)
(a11) edge[transform canvas={xshift=\c}] node[swap] {\xi_{f\angordm, f'\angordmp}} (a21)
(a21) edge[transform canvas={xshift=\c}] node[swap,pos=.3] {\inc_{\upom\xi}} (a22)
(a22) edge[transform canvas={yshift=\t}] node {\iota} (a23)
(a23') edge[equal] (a24)
;
\end{tikzpicture}
\end{equation}
The pointed $G$-functors in \cref{brg_components} are defined as follows.
\begin{itemize}
\item $\brg_{f,f'; \angordn}$ is the $\angordn$-component of $\brg_{f,f'}$.
\item Each $\iota$ is part of the definition of the coend in its codomain.
\item $\xi_{f\angordm, f'\angordmp}$ is the $(f\angordm, f'\angordmp)$-component of the braiding for $(\Gcatst,\sma)$.
\item The upper left arrow $\inc_{\upom}$ is the wedge summand inclusion corresponding to the nonzero morphism $\upom \cn \angordm \oplus \angordmp \to \angordn$.
\item The lower left arrow $\inc_{\upom\xi}$ is the wedge summand inclusion corresponding to the composite 
\[\angordmp \oplus \angordm \fto{\xi} \angordm \oplus \angordmp \fto{\upom} \angordn \inspace \Gsk,\]
where $\xi$ denotes the $(\angordmp,\angordm)$-component of the braiding for $(\Gsk,\oplus)$ \cref{Gsk_braiding}.
\end{itemize}
The braiding $\brg$ for $(\GGTop,\smag)$ admits the same description as \cref{brg_components}, with $\Gtopst$ in place of $\Gcatst$.

\parhead{Symmetry axiom}.
The axiom \cref{monoidalfunctorbraiding} for $(\clast,\clasttwo)$ requires the commutativity of the following diagram in $\GGTop$ for each pair $(f,f')$ of $\Gskg$-categories.
\begin{equation}\label{clast_sym}
\begin{tikzpicture}[vcenter]
\def\c{1ex} \def\v{-1.4}
\draw[0cell]
(0,0) node (a11) {\clast f \smag \clast f'}
(a11)++(4,0) node (a12) {\clast f' \smag \clast f}
(a11)++(0,\v) node (a21) {\clast(f \smag f')}
(a12)++(0,\v) node (a22) {\clast(f' \smag f)}
;
\draw[1cell=.9]
(a11) edge node {\brg_{\clast f, \clast f'}} (a12)
(a12) edge[transform canvas={xshift=-\c}] node {(\clasttwo)_{f',f}} (a22)
(a11) edge[transform canvas={xshift=\c}] node[swap] {(\clasttwo)_{f,f'}} (a21)
(a21) edge node {\clast \brg_{f,f'}} (a22)
;
\end{tikzpicture}
\end{equation}
By \cref{clasttwo_component}, \cref{brg_components}, and the universal property of coends, the symmetry diagram \cref{clast_sym} commutes if and only if the following diagram in $\Gtopst$ commutes for each nonzero morphism $\upom \cn \angordm \oplus \angordmp \to \angordn$ in $\Gsk$, where $\txwedge'$ means $\txwedge_{\Gskpunc(\angordmp \oplus \angordm, \angordn)}$.
\[\begin{tikzpicture}[vcenter]
\def\h{2} \def\u{-1} \def\v{-1.4} \def\c{1em}
\draw[0cell=.8]
(0,0) node (a1) {\cla f\angordm \sma \cla f'\angordmp}
(a1)++(-\h,\u) node (a21) {\cla(f\angordm \sma f'\angordmp)}
(a1)++(\h,\u) node (a22) {\cla f'\angordmp \sma \cla f\angordm}
(a21)++(\h,\u) node (a31) {\cla(f'\angordmp \sma f\angordm)}
(a31)++(0,\v) node (a41) {\cla\big(\txwedge' f'\angordmp \sma f\angordm \big)}
(a41)++(0,\v) node (a5) {\cla\big(\txint^{(\angordmp, \angordm)}\txwedge' f'\angordmp \sma f\angordm \big) = \big(\clast(f' \smag f)\big) \angordn}
;
\draw[1cell=.8]
(a1) edge[transform canvas={xshift=-\c}] node[swap] {\clatwobar} (a21)
(a1) edge[transform canvas={xshift=\c}] node {\xi} (a22)
(a21) edge[transform canvas={xshift=-\c}] node[swap] {\cla\xi} (a31)
(a22) edge[transform canvas={xshift=\c}]  node[pos=.35] {\clatwobar} (a31)
(a31) edge node[swap] {\cla\inc_{\upom\xi}} (a41)
(a41) edge node[swap] {\cla\iota} (a5)
;
\end{tikzpicture}\]
In the diagram above, the top quadrilateral commutes by the symmetry axiom \cref{monoidalfunctorbraiding} for the symmetric monoidal functor $(\cla,\clatwo)$ \cref{cla_gcat_gtop}, the universal property of the smash product, and \cref{clasttwo_one}.
\end{explanation}

\section{The Symmetric Monoidal $\Gtop$-Category of $\Gskg$-Spaces}
\label{sec:ggtop_smgtop}

This section unravels the symmetric monoidal $\Gtop$-category $\GGTop$.  It is the codomain of the symmetric monoidal $\Gtop$-functor $(\clast')_{\evang}$ and is denoted by $(\uof{\GGTop})_{\evang}$ in \cref{clastpev}.  This explicit description of $\GGTop$ is used in \cref{ch:semg} to construct a symmetric monoidal $\Gtop$-functor from $\GGTop$ to $\GSp$.  Recall a symmetric monoidal $\V$-category \pcref{definition:monoidal-vcat,definition:braided-monoidal-vcat,definition:symm-monoidal-vcat} and the Cartesian closed category $(\Gtop,\times,*)$ \pcref{def:Gtop}.  

\secoutline
\begin{itemize}
\item The hom $G$-spaces of $\GGTop$ are discussed in \cref{ggtop_gtop_enr,ggtop_gtopenr_theta_g}.
\item Its monoidal identity and monoidal composition are discussed in \cref{ggtop_mon_identity,ggtop_mon_composition}.
\item Its monoidal unitors are discussed in \cref{ggtop_ellg,rg_f}.
\item Its monoidal associator and braiding are discussed in \cref{ggtop_ag,ggtop_braiding}.
\end{itemize}

\begin{explanation}[$\Gtop$-Enrichment of $\GGTop$]\label{expl:ggtop_gtopenr}
As we explain in \cref{clastp,evang_gtop}, the symmetric monoidal $\Gtop$-category\index{G-G-space@$\Gskg$-space!symmetric monoidal category}\index{symmetric monoidal category!G-G-space@$\Gskg$-space} 
\begin{equation}\label{ggtop_smgtop}
\big(\GGTop, \smag, \gu, \ag, \ellg, \rg, \beg\big)
\end{equation}
is obtained from the symmetric monoidal $(\GGTop)$-category \pcref{expl:ggtop_unit,expl:ggtop_smag,expl:ggtop_brkst}
\[\uof{\GGTop}\]
by changing enrichment along the evaluation functor
\[\GGTop \fto{\evang} \Gtop\]
at the empty tuple $\ang{} \in \Gsk$ \cref{Gsk_objects}.  

\parhead{Base $\Gtop$-category}. 
As a $\Gtop$-category, the objects of $\GGTop$ are $\Gskg$-spaces \cref{ggtop_obj}, which mean pointed functors $\Gsk \to \Gtopst$.  

The hom $G$-spaces of $\GGTop$ are given by the internal hom $\Gskg$-spaces \cref{ggtop_inthom} evaluated at $\ang{} \in \Gsk$.  More explicitly, for $\Gskg$-spaces $f$ and $f'$, the hom $G$-space is given by the following end taken in $\Gtopst$.
\begin{equation}\label{ggtop_gtop_enr}
\begin{split}
\brkstg{f}{f'}\ang{} 
&= \ecint_{\angordn \in \Gsk} \Topgst(f\angordn, f'\angordn) \\
&\bigsubset \prod_{\angordn \in \Gsk} \Topgst(f\angordn, f'\angordn)
\end{split}
\end{equation}
A point $\theta \in \brkstg{f}{f'}\ang{}$ consists of $\angordn$-component pointed morphisms between pointed spaces
\begin{equation}\label{ggtop_gtopenr_thetan}
f\angordn \fto{\theta_{\angordn}} f'\angordn \forspace \angordn \in \Gsk
\end{equation}
such that, for each morphism \cref{Gsk_morphisms} $\upom \cn \angordn \to \angordl$ in $\Gsk$, the following diagram of pointed morphisms commutes.  
\begin{equation}\label{ggtop_gtopenr_theta_nat}
\begin{tikzpicture}[vcenter]
\def\v{-1.4}
\draw[0cell]
(0,0) node (a11) {f\angordn}
(a11)++(2.5,0) node (a12) {f'\angordn}
(a11)++(0,\v) node (a21) {f\angordl}
(a12)++(0,\v) node (a22) {f'\angordl}
;
\draw[1cell=.9]
(a11) edge node {\theta_{\angordn}} (a12)
(a12) edge[transform canvas={xshift=0em}] node {f'\upom} (a22)
(a11) edge node[swap] {f\upom} (a21)
(a21) edge node {\theta_{\angordl}} (a22)
;
\end{tikzpicture}
\end{equation}
While $f\upom$ and $f'\upom$ are $G$-equivariant \cref{f_upom_ggtop}, the components of $\theta$ are \emph{not} required to be $G$-equivariant.  In particular, $\theta$ is not generally a morphism of $\Gskg$-spaces \pcref{expl:ggtop_mor}.  

Composition \cref{mcomp_abc} and identities \cref{cone_a}, in the $\Gtop$-enriched sense, are defined using the components \cref{ggtop_gtopenr_thetan}.

\parhead{$G$-action}. 
The $G$-action on the hom $G$-space $\brkstg{f}{f'}\ang{}$ is given componentwise by conjugation \cref{ggtop_inthom_theta_g}.  In other words, for a point $\theta \in \brkstg{f}{f'}\ang{}$, $g \in G$, and $\angordn \in \Gsk$, the $\angordn$-component of $g \cdot \theta \in \brkstg{f}{f'}\ang{}$ is the following composite of pointed morphisms.
\begin{equation}\label{ggtop_gtopenr_theta_g}
\begin{tikzpicture}[vcenter]
\def\v{-1.5} \def\t{18}
\draw[0cell]
(0,0) node (a11) {f\angordn}
(a11)++(2.5,0) node (a12) {f'\angordn}
(a11)++(0,\v) node (a21) {f\angordn}
(a12)++(0,\v) node (a22) {f'\angordn}
;
\draw[1cell=.9]
(a11) edge node {(g \cdot \theta)_{\angordn}} (a12)
(a11) edge node[swap] {\ginv} (a21)
(a21) edge node {\theta_{\angordn}} (a22)
(a22) edge node[swap] {g} (a12)
;
\end{tikzpicture}
\end{equation}
A point $\theta \in \brkstg{f}{f'}\ang{}$ is a morphism $\theta \cn f \to f'$ of $\Gskg$-spaces \pcref{expl:ggtop_mor} if and only if it is $G$-fixed, which means that each component $\theta_{\angordn}$ is $G$-equivariant.

\parhead{Monoidal identity}. 
With $\vtensorunit$ denoting the unit $\Gtop$-category \pcref{definition:unit-vcat}, the monoidal identity \cref{monvcat_monid} is the $\Gtop$-functor
\begin{equation}\label{ggtop_mon_identity}
\vtensorunit \fto{\gu} \GGTop
\end{equation}
that sends the unique object of $\vtensorunit$ to the monoidal unit $\Gskg$-space $\gu$ in \cref{ggtop_unit}.  This assignment defines a $\Gtop$-functor \pcref{def:enriched-functor} because the identity morphism of each $\Gskg$-space, including $\gu$, is fixed by the conjugation $G$-action \cref{ggtop_gtopenr_theta_g}.

\parhead{Monoidal composition}.
The monoidal composition \cref{monvcat_moncomp} is the $\Gtop$-functor \pcref{def:enriched-functor,definition:vtensor-0}
\begin{equation}\label{ggtop_mon_composition}
\GGTop \otimes \GGTop \fto{\smag} \GGTop
\end{equation}
whose object assignment
\[(f; f') \mapsto f \smag f'\]
is given by the monoidal product of $\Gskg$-spaces \pcref{expl:ggtop_smag}. 

For $\Gskg$-spaces $f,f',h$, and $h'$, the $G$-morphism between hom $G$-spaces
\begin{equation}\label{ggtop_moncomp_mor}
\brkstg{f}{h}\ang{} \times \brkstg{f'}{h'}\ang{} \fto{\smag} \brkstg{f \smag f'}{h \smag h'}\ang{}  
\end{equation}
sends a pair of points \cref{ggtop_gtop_enr}
\[(\theta; \theta') \in \big(\txint_{\angordn \in \Gsk} \Topgst(f\angordn, h\angordn) \big) \times \big( \txint_{\angordn \in \Gsk} \Topgst(f'\angordn, h'\angordn) \big)\]
to the point
\begin{equation}\label{tha_smag_thap}
\theta \smag \theta' \in \txint_{\angordn \in \Gsk} \Topgst\big( (f \smag f')\angordn, (h \smag h')\angordn\big).
\end{equation}
For each object $\angordn \in \Gsk$, its $\angordn$-component pointed morphism is defined by the following commutative diagram using \cref{smag_ptspace}, where
\[\txint = \txint^{(\angordm, \angordmp) \in \Gsk^2} \andspace 
\txwedge = \txwedge_{\Gskpunc(\angordm \oplus \angordmp, \angordn)}.\]
\begin{equation}\label{ggtop_moncomp_morn}
\begin{tikzpicture}[vcenter]
\def\v{-1.4}
\draw[0cell]
(0,0) node (a11) {(f \smag f')\angordn}
(a11)++(3.3,0) node (a12) {\txint \txwedge (f\angordm \sma f'\angordmp)}
(a11)++(0,\v) node (a21) {(h \smag h')\angordn}
(a12)++(0,\v) node (a22) {\txint \txwedge (h\angordm \sma h'\angordmp)} 
;
\draw[1cell=.9]
(a11) edge[equal] (a12)
(a21) edge[equal] (a22)
(a11) edge[transform canvas={xshift=1.5em}] node[swap] {(\theta \smag \theta')_{\angordn}} (a21)
(a12) edge[transform canvas={xshift=-3em}] node {\txint \txwedge (\theta_{\angordm} \sma \theta'_{\angordmp})} (a22)
;
\end{tikzpicture}
\end{equation}
Using \cref{smag_rep_gaction,ggtop_gtopenr_theta_g}, the following computation proves the $G$-equivariance of the assignment $(\theta; \theta') \mapsto \theta \smag \theta'$ in \cref{ggtop_moncomp_mor}.
\[\begin{split}
& \big((g \cdot \theta) \smag (g \cdot \theta') \big)_{\angordn} \\
&= \txint \txwedge \big( (g \cdot \theta)_{\angordm} \sma (g \cdot \theta')_{\angordmp} \big) \\
&= \txint \txwedge \big( (g \theta_{\angordm} \ginv) \sma (g \theta'_{\angordmp} \ginv) \big) \\
&= \txint \txwedge \big((g \sma g) (\theta_{\angordm} \sma \theta'_{\angordmp}) (\ginv \sma \ginv) \big) \\
&= g \big(\txint \txwedge (\theta_{\angordm} \sma \theta'_{\angordmp})\big) \ginv \\
&= \big(g \cdot (\theta \smag \theta')\big)_{\angordn}
\end{split}\]
The $\Gtop$-functoriality of $\smag$ \cref{ggtop_moncomp_mor}, in the sense of \cref{def:enriched-functor}, follows from \cref{ggtop_moncomp_morn} and the fact that composition and identities are defined componentwise \cref{ggtop_gtopenr_thetan}.

\parhead{Monoidal unitors}. 
The left monoidal unitor \cref{eq:monoidal-unitors} is the $\Gtop$-natural isomorphism
\begin{equation}\label{ggtop_ellg}
\begin{tikzpicture}[vcenter]
\def\v{1.4} \def\h{.15}
\draw[0cell=.9]
(0,0) node (a11) {\GGTop}
(a11)++(4,0) node (a12) {\GGTop}
(a11)++(\h,\v) node (a01) {\vtensorunit \otimes \GGTop}
(a12)++(-\h,\v) node (a02) {(\GGTop)^{\otimes 2}}
;
\draw[1cell=.9]
(a11) edge node[swap] {1} (a12)
(a11) edge node[pos=.2] {(\ell^\otimes)^{-1}} (a01)
(a01) edge node {\gu \otimes 1} (a02)
(a02) edge node[pos=.7] {\smag} (a12)
;
\draw[2cell]
node[between=a11 and a12 at .47, shift={(0,.5*\v)}, rotate=-90, 2label={above,\ellg}] {\Rightarrow}
;
\end{tikzpicture}
\end{equation}
whose component at a $\Gskg$-space $f$ is the $G$-morphism
\[* \fto{\ellg_f} \brkstg{\gu \smag f}{f}\ang{}\]
given by the left unit isomorphism \pcref{expl:ggtop_mor,expl:ggtop_unit}
\begin{equation}\label{ellg_f}
\gu \smag f \fto[\iso]{\ellg_f} f \inspace \GGTop.
\end{equation}
More explicitly, for each object $\angordn \in \Gsk$, the $\angordn$-component of $\ellg_f$ is given by the following composite of pointed $G$-homeomorphisms.
\begin{equation}\label{ellg_f_n}
\begin{split}
& (\gu \smag f) \angordn \\
&= \int^{(\angordm, \angordmp) \in \Gsk^2} \mquad \bigvee_{\Gskpunc(\angordm \oplus \angordmp, \angordn)} \mquad (\gu\angordm \sma f\angordmp) \\
&\iso \int^{(\angordm, \angordmp) \in \Gsk^2} \mqq \bigvee_{\Gskpunc(\angordm \oplus \angordmp, \angordn) \times \Gskpunc(\ang{}, \angordm)} \mqq (\stplus \sma f\angordmp) \\
&\iso \int^{\angordmp \in \Gsk} \bigvee_{\Gskpunc(\angordmp, \angordn)} f\angordmp \iso f\angordn
\end{split}
\end{equation}
The details of \cref{ellg_f_n} are explained below.
\begin{itemize}
\item The first equality uses the definition \cref{smag_ptspace} of $\smag$. 
\item The first pointed $G$-homeomorphism uses the definition \cref{gu_angordm_ggtop} of $\gu\angordm$ and the fact that $- \sma f\angordmp$ commutes with wedges.
\item The second pointed $G$-homeomorphism uses
\begin{itemize}
\item the left unit isomorphism 
\[\stplus \sma f\angordmp \iso f\angordmp\] for $(\Gtopst,\sma)$,
\item the left unit identity
\[\ang{} \oplus \angordmp = \angordmp\]
for $(\Gsk,\oplus)$ \pcref{def:Gsk_permutative}, and
\item the Yoneda Density Theorem for coends (\cite[A.5.7]{loregian} or \cite[3.7.15]{cerberusIII}).
\end{itemize}
\item The last pointed $G$-homeomorphism also uses the Yoneda Density Theorem.
\end{itemize}
The naturality \cref{enr_naturality} of $\ellg$ follows from the fact that each pointed $G$-homeomorphism in \cref{ellg_f_n} is natural in the $\Gskg$-space $f$.

The right monoidal unitor $\rg$ is defined analogously to \cref{ellg_f_n} using the right unit isomorphism
\begin{equation}\label{rg_f}
f \smag \gu \fto[\iso]{\rg_f} f \inspace \GGTop
\end{equation}
for each $\Gskg$-space $f$.

\parhead{Monoidal associator}.  It is the $\Gtop$-natural isomorphism \cref{eq:monoidal-assoc}
\begin{equation}\label{ggtop_ag}
\begin{tikzpicture}[vcenter]
\def\h{2.5} \def\u{1} \def\v{-1.4}
\draw[0cell=.9]
(0,0) node (a11) {(\GGTop)^{\otimes 2} \otimes \GGTop}
(a11)++(\h,\u) node (a12) {(\GGTop)^{\otimes 2}}
(a12)++(\h,-\u) node (a13) {\GGTop}
(a11)++(0,\v) node (a21) {\GGTop \otimes (\GGTop)^{\otimes 2}}
(a13)++(0,\v) node (a22) {(\GGTop)^{\otimes 2}}
;
\draw[1cell=.8]
(a11) edge node {\smag \otimes 1} (a12)
(a12) edge node {\smag} (a13)
(a11) edge node[swap] {a^\otimes} (a21)
(a21) edge node {1 \otimes \smag} (a22)
(a22) edge node[swap] {\smag} (a13)
;
\draw[2cell]
node[between=a11 and a13 at .48, shift={(0,.3*\v)}, rotate=-90, 2label={above,\ag}] {\Rightarrow}
;
\end{tikzpicture}
\end{equation}
whose component at a triple $(f,f',f'')$ of $\Gskg$-spaces is given by the associativity isomorphism
\begin{equation}\label{ggtop_ag_fff}
(f \smag f') \smag f'' \fto[\iso]{\ag_{f,f',f''}} f \smag (f' \smag f'') \inspace \GGTop.
\end{equation}
For each object $\angordn \in \Gsk$, the $\angordn$-component of $\ag_{f,f',f''}$ is given by the following commutative diagram of pointed $G$-homeomorphisms, where
\[\txint = \txint^{(\angordm, \angordmp, \angordmpp) \in \Gsk^3} \andspace 
\txwedge = \txwedge_{\Gskpunc(\angordm \oplus \angordmp \oplus \angordmpp, \angordn)}.\]
\begin{equation}\label{ggtop_ag_comp}
\begin{tikzpicture}[vcenter]
\def\v{-1.4} 
\draw[0cell=.85]
(0,0) node (a11) {((f \smag f') \smag f'')\angordn}
(a11)++(4.3,0) node (a12) {\txint \txwedge [(f\angordm \sma f'\angordmp) \sma f''\angordmpp]}
(a11)++(0,\v) node (a21) {(f \smag (f' \smag f''))\angordn}
(a12)++(0,\v) node (a22) {\txint \txwedge [f\angordm \sma (f'\angordmp \sma f''\angordmpp)]} 
;
\draw[1cell=.8]
(a11) edge node {\iso} (a12)
(a21) edge node {\iso} (a22)
(a11) edge[transform canvas={xshift=1.5em}] node[swap] {\ag_{f,f',f''; \angordn}} (a21)
(a12) edge[transform canvas={xshift=-3em}] node[swap] {\iso} node {\txint \txwedge \al} (a22)
;
\end{tikzpicture}
\end{equation}
\begin{itemize}
\item In \cref{ggtop_ag_comp}, $\al$ is the associativity isomorphism for $(\Gtopst, \sma)$. 
\item The top and bottom pointed $G$-homeomorphisms in \cref{ggtop_ag_comp} are proved by the computation in
\cref{au_xyzu_dom,au_xyzu_cod}, up to a change of notation.  
\end{itemize}
The naturality \cref{enr_naturality} of $\ag$ follows from the fact that each of the right, top, and bottom pointed $G$-homeomorphisms in \cref{ggtop_ag_comp} is natural in the $\Gskg$-spaces $f$, $f'$, and $f''$.

\parhead{Braiding}.  It is the $\Gtop$-natural isomorphism
\begin{equation}\label{ggtop_braiding}
\begin{tikzpicture}[vcenter]
\def\h{4}
\draw[0cell]
(0,0) node (a1) {(\GGTop)^{\otimes 2}}
(a1)++(.5*\h,-1) node (a2) {(\GGTop)^{\otimes 2}}
(a1)++(\h,0) node (a3) {\GGTop}
;
\draw[1cell=.9]
(a1) edge node {\smag} (a3)
(a1) [rounded corners=2pt] |- node[swap,pos=.25] {\beta^\otimes} (a2)
;
\draw[1cell=.9]
(a2) [rounded corners=2pt] -| node[swap,pos=.75] {\smag} (a3)
;
\draw[2cell]
node[between=a1 and a3 at .5, shift={(0,-.45)}, rotate=-90, 2label={above,\beg}] {\Rightarrow}
;
\end{tikzpicture}
\end{equation}
whose component at a pair $(f,f')$ of $\Gskg$-spaces is given by the braiding
\begin{equation}\label{beg_ff}
f \smag f' \fto[\iso]{\beg_{f,f'}} f' \smag f \inspace \GGTop.
\end{equation}
For each object $\angordn \in \Gsk$, its $\angordn$-component is the pointed $G$-homeomorphism 
\[(f \smag f')\angordn \fto[\iso]{\beg_{f,f';\angordn}} (f' \smag f)\angordn\]
defined in \cref{brg_components}.

\parhead{Axioms}.  
Each of
\begin{itemize}
\item the monoidal identity $\gu$ \cref{ggtop_mon_identity},
\item the left unitor $\ellg$ \cref{ellg_f},
\item the right unitor $\rg$ \cref{rg_f},
\item the associator $\ag$ \cref{ggtop_ag_fff}, and
\item the braiding $\beg$ \cref{beg_ff}
\end{itemize} 
is defined componentwise by the corresponding structure for the symmetric monoidal category $\GGTop$ \cref{ggtop_smc}.  Each of the symmetric monoidal $\Gtop$-category axioms \pcref{definition:monoidal-vcat,definition:braided-monoidal-vcat,definition:symm-monoidal-vcat} for $\GGTop$ is an equality of $\Gtop$-natural transformations, which are determined by their components \pcref{def:enriched-natural-transformation}.  Thus, each of these symmetric monoidal $\Gtop$-category axioms holds by the corresponding symmetric monoidal category axiom for $\GGTop$ \pcref{def:monoidalcategory,def:braidedmoncat,def:symmoncat}.
\end{explanation}

\section{Commutation of Classifying Space and Evaluation}
\label{sec:clast_ev_commute}

The section explains that the symmetric monoidal $\Gtop$-category $(\uof{\GGCat}_{\clast})_{\evang}$ in \cref{clastpev} and the $\Gtop$-multifunctor $\End((\clast')_{\evang})$ in \cref{End_clastpev} can also be constructed by switching the order of change of enrichment and $\End$.  These observations are not used in \cref{ch:semg}; they are included for completeness.

\secoutline
\begin{itemize}
\item \cref{clast_ev_commute} proves that the symmetric monoidal functors $\clast$ and $\evang$ commute.  As a result, starting from the symmetric monoidal $(\GGCat)$-category $\uof{\GGCat}$, the symmetric monoidal $\Gtop$-categories obtained by changing enrichment along either $\evang \circ \clast$ or $\cla \circ \evang$ are equal \pcref{ggcat_clastpev}.
\item \cref{expl:End_clastpev} discusses the fact that the six $\Gtop$-multicategories associated to the symmetric monoidal $(\GGCat)$-category $\uof{\GGCat}$---by changing enrichment and applying the endomorphism construction in different orders---are the same. 
\end{itemize}


\begin{proposition}\label{clast_ev_commute}
For each group $G$, the diagram 
\begin{equation}\label{clast_ev}
\begin{tikzpicture}[vcenter]
\def\v{-1.4}
\draw[0cell]
(0,0) node (a11) {\GGCat}
(a11)++(3,0) node (a12) {\GGTop}
(a11)++(0,\v) node (a21) {\Gcat}
(a12)++(0,\v) node (a22) {\Gtop}
;
\draw[1cell=.9]
(a11) edge node {\clast} (a12)
(a12) edge node {\evang} (a22)
(a11) edge node[swap] {\evang} (a21)
(a21) edge node {\cla} (a22)
;
\end{tikzpicture}
\end{equation}
of symmetric monoidal functors commutes.
\end{proposition}

\begin{proof}
The diagram \cref{clast_ev} of symmetric monoidal functors commutes for the following reasons.
\begin{description}
\item[Functors] The diagram \cref{clast_ev} of functors commutes because, by \cref{expl:clast_functor}, each of the two composites sends a $\Gskg$-category $f$ to the $G$-space $\cla f\ang{}$, and likewise for morphisms in $\GGCat$.
\item[Unit constraints] By \cref{expl:clast_unit,evtu_unpt}, the unit constraint \cref{monoidal_unit} for each of the two composites in \cref{clast_ev} is the $G$-morphism
\[* \fto{} \cla\gu\ang{} = \cla(\bonep)\]
given by the non-basepoint of $\bonep$.
\item[Monoidal constraints] By \cref{expl:clast_monoidal,evtu_unpt}, the monoidal constraints \cref{monoidal_constraint} of the two composites in \cref{clast_ev} have, for any two $\Gskg$-categories $f$ and $f'$, $(f,f')$-components given by the following composites in $\Gtop$.
\[\begin{tikzpicture}[vcenter]
\def\v{-1.4}
\draw[0cell=.9]
(0,0) node (a11) {\cla f \ang{} \times \cla f' \ang{}}
(a11)++(4,0) node (a12) {\cla f \ang{} \sma \cla f' \ang{}}
(a11)++(0,\v) node (a21) {\cla(f \ang{} \times f' \ang{})}
(a12)++(0,\v) node (a22) {\cla(f \ang{} \sma f' \ang{})}
(a21)++(0,\v) node (a31) {\cla\Big(\int^{(\angordm, \angordmp)} \mqq \bigvee_{\Gskpunc(\angordm \oplus \angordmp, \ang{})} \mqq f\angordm \sma f'\angordmp \Big)}
(a22)++(0,\v) node (a32) {\cla\big((f \smag f') \ang{}\big)}
;
\draw[1cell=.85]
(a11) edge node {\pjt} (a12)
(a21) edge node {\cla\pjt} (a22)
(a31) edge[equal] (a32)
(a11) edge node {\iso} node[swap] {\clatwo} (a21)
(a12) edge node {\clatwobar} (a22)
(a22) edge node {\cla i} (a32)
;
\end{tikzpicture}\]
In the previous diagram, the top rectangle commutes by \cref{clasttwo_one}.  The $G$-functor 
\[f \ang{} \sma f' \ang{} \fto{i} (f \smag f') \ang{}\]
is given by
\begin{itemize}
\item the coend structure morphism for $\angordm = \ang{} = \angordmp$ and
\item the fact that
\[\Gskpunc(\ang{} \oplus \ang{}, \ang{}) = \Gskpunc(\ang{}, \ang{}) = \{1_{\ang{}}\}.\]
\end{itemize}
\end{description}
This proves that the two composites in \cref{clast_ev} are equal as symmetric monoidal functors.  
\end{proof}

By \cref{clast_ev_commute}, changing enrichment from $\GGCat$ to $\Gtop$ via either one of the two composites in \cref{clast_ev} yields the same symmetric monoidal $\Gtop$-category, thereby proving the following result.

\begin{corollary}\label{ggcat_clastpev}
For each group $G$, there is an equality 
\begin{equation}\label{clast_ev_change_enr}
(\uof{\GGCat}_{\clast})_{\evang} = (\uof{\GGCat}_{\evang})_{\cla}
\end{equation}
of symmetric monoidal $\Gtop$-categories.
\end{corollary}

\begin{explanation}\label{expl:clastpev}\
\begin{itemize}
\item The left-hand side of \cref{clast_ev_change_enr} is the symmetric monoidal $\Gtop$-category in \cref{clastpev}.  It is obtained from the symmetric monoidal $(\GGCat)$-category $\uof{\GGCat}$ in \cref{ggcat_smggcat} by changing enrichment along $\clast$ and then $\evang$. 
\item The right-hand side of \cref{clast_ev_change_enr} is obtained from the symmetric monoidal $\Gcat$-category $\uof{\GGCat}_{\evang}$ in \cref{ggcat_smgcat} by changing enrichment along the classifying space functor $\cla$ \cref{cla_gcat_gtop}.\defmark
\end{itemize}
\end{explanation}

\begin{explanation}[Change of Enrichment Commutes with $\End$]\label{expl:End_clastpev}
By \cref{End_dashf}, change of enrichment for enriched symmetric monoidal categories and for enriched multicategories are compatible via the endomorphism construction \pcref{definition:EndK,def:EndF}.  Thus, the commutative diagram \cref{clast_ev} of symmetric monoidal functors yields the following equalities of $\Gtop$-multicategories, in which the first entry is the domain in \cref{End_clastpev}. 
\[\begin{split}
\End \big(\big(\uof{\GGCat}_{\clast}\big)_{\evang}\big) 
&= \big(\End \big(\uof{\GGCat}_{\clast}\big)\big)_{\evang} \\
&= \big(\big(\End\, \uof{\GGCat}\big)_{\clast} \big)_{\evang} \\
&= \big(\big(\End\, \uof{\GGCat}\big)_{\evang} \big)_{\cla} \\
&= \big(\End \big(\uof{\GGCat}_{\evang}\big) \big)_{\cla} \\
&= \End \big(\big(\uof{\GGCat}_{\evang}\big)_{\cla} \big) \\
\end{split}\]
Similarly, for the codomain in \cref{End_clastpev}, there is an equality of $\Gtop$-multicategories
\[\End \big((\uof{\GGTop})_{\evang}\big) = \big(\End\, \uof{\GGTop}\big)_{\evang}.\]
Moreover, there is an equality of $\Gtop$-multifunctors
\[\End ((\clast')_{\evang}) = (\End\, \clast')_{\evang}\]
in which $\clast'$ is the symmetric monoidal $(\GGTop)$-functor in \cref{clastp}.
\end{explanation}

\chapter{Symmetric $\Gtop$-Monoidal Equivariant $K$-Theory}
\label{ch:semg}
This chapter completes the construction of our $G$-equivariant algebraic $K$-theory multifunctor by constructing the last step, which is a unital symmetric monoidal $\Gtop$-functor \pcref{thm:Kg_smgtop}
\[\GGTop \fto{(\Kg,\Kgtwo,\Kgzero)} \GSp.\]
The unit constraint $\Kgzero$ is an isomorphism, so the $G$-sphere is in the image of $\Kg$.  Throughout this chapter, $G$ is a compact Lie group, and $\univ$ is a complete $G$-universe \pcref{def:g_universe}.  
\begin{itemize}
\item The domain of $\Kg$ is the symmetric monoidal $\Gtop$-category $\GGTop$ \pcref{expl:ggtop_gtopenr} of $\Gskg$-spaces \cref{ggtop_obj}.
\item The codomain of $\Kg$ is the symmetric monoidal $\Gtop$-category $\GSp$ of orthogonal $G$-spectra \pcref{gspectra_smgtop}.  
\end{itemize}  
For a $\Tinf$-operad $\Op$ \pcref{as:OpA}, our $G$-equivariant algebraic $K$-theory $\Kgo$ is the following composite $\Gtop$-multifunctor \pcref{thm:Kgo_multi}.
\[\begin{tikzpicture}
\def\v{-1.4} 
\draw[0cell]
(0,0) node (a11) {\MultpsO}
(a11)++(3.2,0) node (a12) {\phantom{\GSp}}
(a12)++(0,.03) node (a12') {\GSp}
(a11)++(0,\v) node (a21) {\GGCat}
(a12)++(0,\v) node (a22) {\phantom{\GGTop}}
(a22)++(.3,0) node (a22') {\GGTop}
;
\draw[1cell=.9]
(a11) edge node {\Kgo} (a12)
(a11) edge[shorten <=-.5ex] node[swap] {\Jgo} (a21)
(a21) edge node {\clast} (a22')
(a22) edge node[swap] {\Kg} (a12')
;
\end{tikzpicture}\]
\begin{itemize}
\item The first step $\Jgo$ is the $\Gcat$-multifunctor in \cref{thm:Jgo_multifunctor}.  It is regarded as a $\Gtop$-multifunctor by changing enrichment along the classifying space functor $\cla \cn \Gcat \to \Gtop$ \pcref{cla_multi}.  The objects in $\MultpsO$ are $\Op$-pseudoalgebras \pcref{def:pseudoalgebra}.
\item The second step $\clast$ is the $\Gtop$-multifunctor in \cref{End_clastpev} induced by a symmetric monoidal $\Gtop$-functor.  
\end{itemize}
There is also a strong variant of $\Kgo$, denoted $\Kgosg$, whose first step is the $\Gcat$-multifunctor \pcref{thm:Jgo_multifunctor}
\[\MultpspsO \fto{\Jgosg} \GGCat.\]
All the desirable properties of $\Kgo$ are also true for $\Kgosg$.

\subsection*{Preservation of Equivariant Algebraic Structures}
Enriched multifunctoriality is the most important property of our $G$-equivariant algebraic $K$-theory machine $\Kgo$.  It preserves all algebraic structures parametrized by multicategories enriched in $\Gtop$ and $\Gcat$, by a change of enrichment along the classifying space functor $\cla$.  In other words, for each $\Gcat$-multicategory $\cQ$ and $\Gcat$-multifunctor $f \cn \cQ \to \MultpsO$, after changing enrichment along $\cla$, the composite
\[\cQ \fto{f} \MultpsO \fto{\Kgo} \GSp\]
is a $\Gtop$-multifunctor \pcref{thm:Kgo_preservation}.  Thus, $\Kgo$ sends each $\cQ$-algebra $f$ in the $\Gcat$-multicategory $\MultpsO$ to the $\cQ$-algebra $\Kgo \circ f$ in the symmetric monoidal $\Gtop$-category $\GSp$.  Here are two examples \pcref{ex:Kgo_preservation}.
\begin{itemize}
\item Taking $\cQ$ to be the $G$-Barratt-Eccles operad $\GBE$ for a finite group $G$ \pcref{def:GBE}, a $\Gcat$-multifunctor 
\[\GBE \fto{f} \MultpsO\]
is both an equivariant $\Einf$-algebra and an $\Ninf$-algebra of $\Op$-pseudoalgebras.  The composite $\Gtop$-multifunctor 
\[\GBE \fto{f} \MultpsO \fto{\Kgo} \GSp\]
is both an equivariant $\Einf$-algebra and an $\Ninf$-algebra of orthogonal $G$-spectra. 
\item Further specializing $\Op$ to $\GBE$, a $\Gcat$-multifunctor 
\[\GBE \fto{f} \MultpsGBE\]
is both an equivariant $\Einf$-algebra and an $\Ninf$-algebra of genuine symmetric monoidal $G$-categories \pcref{def:GBE_pseudoalg}.  The composite $\Gtop$-multifunctor 
\[\GBE \fto{f} \MultpsGBE \fto{\Kggbe} \GSp\]
is both an equivariant $\Einf$-algebra and an $\Ninf$-algebra of orthogonal $G$-spectra. 
\end{itemize} 

\summary
The following table summarizes $\Kgo = \Kg \circ \clast \circ \Jgo$, where sm means symmetric monoidal.
\begin{center}
\resizebox{\textwidth}{!}{%
{\renewcommand{\arraystretch}{1.3}%
{\setlength{\tabcolsep}{1ex}
\begin{tabular}{c|cr|cr}
&&& \multicolumn{2}{c}{Objects} \\ \hline
$\MultpsO$ & $\Gcat$-multicategory & \eqref{thm:multpso} & $\Op$-pseudoalgebras & \eqref{def:pseudoalgebra} \\
$\GGCat$ & symmetric monoidal closed category & \eqref{def:GGCat} & $\Gskg$-categories & \eqref{expl:ggcat_obj} \\
$\GGTop$ & symmetric monoidal closed category & \eqref{def:ggtop_smc} & $\Gskg$-spaces & \eqref{expl:ggtop_obj} \\
$\GSp$ & symmetric monoidal $\Gtop$-category & \eqref{gspectra_smgtop} & orthogonal $G$-spectra & \eqref{expl:gspectra} \\ \hline
$\Jgo$ & $\Gcat$-multifunctor $\MultpsO \to \GGCat$ & \eqref{thm:Jgo_multifunctor} &  \multicolumn{2}{c}{\eqref{def:Jgo_multifunctor}} \\
$\clast$ & sm $\Gtop$-functor $\GGCat \to \GGTop$ & \eqref{thm:ggcat_ggtop} & \multicolumn{2}{c}{\eqref{expl:clast_functor}, \eqref{expl:clast_unit}, \eqref{expl:clast_monoidal}, \eqref{expl:clast_symmetry}} \\
$\Kg$ & unital sm $\Gtop$-functor $\GGTop \to \GSp$ & \eqref{thm:Kg_smgtop} & \multicolumn{2}{c}{\eqref{def:ggspace_gspectra}, \eqref{def:ggtop_gsp_mor}, \eqref{def:Kg_zero}, \eqref{def:Kgtwo}} \\ \hline
$\Kgo$ & \multicolumn{3}{c}{$\Gtop$-multifunctor $\MultpsO \to \GGCat \to \GGTop \to \GSp$} & \eqref{thm:Kgo_multi}
\end{tabular}}}}
\end{center}

\connection
This chapter uses most major constructions and results from previous chapters.  In addition to the table above, this chapter includes extensive references back to previous chapters to help the reader navigate back and forth.

\organization
This chapter consists of the following sections.

\secname{sec:Kg_obj}  The unital symmetric monoidal $\Gtop$-functor 
\[\GGTop \fto{(\Kg,\Kgtwo,\Kgzero)} \GSp\]
is constructed in several steps.  This section constructs the object assignment of $\Kg$, which sends each $\Gskg$-space $X$ to an orthogonal $G$-spectrum $\Kg X$.

\secname{sec:Kg_gtop_functor}  This section constructs the $\Gtop$-functor $\Kg$ by defining its assignment on hom $G$-spaces and checking all the axioms of a $\Gtop$-functor.

\secname{sec:Kgzero}  This section constructs the unit constraint $\Kgzero$ for $\Kg$.  The invertibility of the $\Gtop$-natural transformation $\Kgzero$ implies that $\Kg$ includes the $G$-sphere $\gsp$ in its image.

\secname{sec:Kgtwo}  This section constructs the monoidal constraint $\Kgtwo$ for $\Kg$.  A substantial amount of detailed discussion in \cref{ch:spectra} about $\IU$-spaces and orthogonal $G$-spectra is used in this section and the next section.

\secname{sec:Kg}  This section proves that the triple $(\Kg,\Kgtwo,\Kgzero)$ satisfies all the axioms of a unital symmetric monoidal $\Gtop$-functor.  \cref{rk:gspace_gspectra} briefly discusses other equivariant infinite loop space machines in the literature.  Each such machine builds $G$-spectra from space-level data.

\secname{sec:Kgo_multi}  This section records our $G$-equivariant algebraic $K$-theory $\Gtop$-multifunctor
\[\MultpsO \fto{\Kgo} \GSp\]
and the fact that it preserves equivariant algebraic structures.  For example, $\Kgo$ constructs equivariant $\Einf$-algebras and $\Ninf$-algebras of orthogonal $G$-spectra from corresponding types of algebras in the $\Gcat$-multicategory $\MultpsO$ of $\Op$-pseudoalgebras.

\section{From $\Gskg$-Spaces to $G$-Spectra}
\label{sec:Kg_obj}

This section constructs the object assignment of 
\[\GGTop \fto{\Kg} \GSp,\]
which sends each $\Gskg$-space $X$ to an orthogonal $G$-spectrum $\Kg X$ \pcref{expl:ggtop_obj,expl:gspectra}.

\secoutline
\begin{itemize}
\item \cref{def:ggspace_gspectra} defines the data $\Kg X$ associated to each $\Gskg$-space $X$.
\item \cref{Kgx_welldef} proves that $\Kg X$ is a well-defined orthogonal $G$-spectrum.
\end{itemize}

\begin{definition}\label{def:ggspace_gspectra}
Given a $\Gskg$-space \cref{ggtop_obj}
\[(\Gsk,\vstar) \fto{X} (\Gtopst,*),\]
define the \index{orthogonal G-spectrum@orthogonal $G$-spectrum!multifunctorial K-theory@multifunctorial $K$-theory}\index{multifunctorial K-theory@multifunctorial $K$-theory!orthogonal G-spectrum@orthogonal $G$-spectrum}orthogonal $G$-spectrum \pcref{def:gsp_module}
\begin{equation}\label{Kg_object}
(\Kg X, \umu) \in \GSp
\end{equation}
as follows.
\begin{description}
\item[Object assignment of $\IU$-space]
The $\IU$-space \pcref{def:iu_space}
\begin{equation}\label{Kgx_iuspace}
\IU \fto{\Kg X} \Topgst
\end{equation}
sends each object $V \in \IU$ to the coend
\begin{equation}\label{Kgxv}
(\Kg X)_V = \int^{\angordn \in \Gsk} \Topgst(\sma\angordn, S^V) \sma X\angordn
\end{equation}
taken in $\Gtopst$.
\begin{itemize}
\item In \cref{Kgxv}, $\sma\angordn \in \Fsk$ uses the functor $\sma \cn \Gsk \to \Fsk$ \pcref{def:smashFskGsk}.  The pointed finite set $\sma\angordn \in \Fsk$ is regarded as a discrete pointed $G$-space with the trivial $G$-action.
\item $S^V$ is the $V$-sphere, with $G$-fixed basepoint $\infty$ \pcref{def:indexing_gspace}.
\item $\Topgst$ is the internal hom of $\Gtopst$ \cref{Gtopst_smc}, so $\Topgst(\sma\angordn, S^V)$ is the pointed $G$-space of pointed morphisms 
\[\sma\angordn \to S^V.\] 
\begin{itemize}
\item Its basepoint is given by the constant morphism at $\infty \in S^V$. 
\item The group $G$ acts on this pointed space by conjugation \cref{ginv_h_g}.  This means post-composition with the $G$-action on $S^V$ because $G$ acts trivially on $\sma\angordn$.
\end{itemize}
We also use the notation
\begin{equation}\label{SVsman}
(S^V)^{\sma\angordn} = \Topgst(\sma\angordn, S^V).
\end{equation}
\item Note that the coend \cref{Kgxv} can also be taken in $\Topst$, with $G$ acting diagonally on each pointed space $(S^V)^{\sma\angordn} \sma X\angordn$.
\end{itemize}
\item[Morphism assignment of $\IU$-space]
For a linear isometric isomorphism $f \cn V \fiso W$ in $\IU$, the pointed homeomorphism \cref{iu_space_xf}
\begin{equation}\label{Kgxf}
(\Kg X)_V \fto[\iso]{(\Kg X)_f} (\Kg X)_W
\end{equation}
is induced by the pointed homeomorphisms 
\[(S^V)^{\sma\angordn} \fto[\iso]{(\gsp f) \circ -} (S^W)^{\sma\angordn} \forspace \angordn \in \Gsk.\]
The pointed homeomorphism $\gsp f \cn S^V \fiso S^W$, which is also denoted by $f$, is defined in \cref{gsp_morphism}.
\item[$\gsp$-action]
The right $\gsp$-action \cref{gspectra_action}
\begin{equation}\label{Kgx_sphere_action}
\Kg X \smau \gsp \fto{\umu} \Kg X
\end{equation}
on $\Kg X$ has, for each pair of objects $(V,W) \in (\IUsk)^2$, component pointed $G$-morphism \cref{gsp_action_vw} defined by the following commutative diagram in $\Gtopst$.
\begin{equation}\label{Kgx_action_vw}
\begin{tikzpicture}[vcenter]
\def\h{5} \def\u{-1} \def\v{-1.4}
\draw[0cell=.9]
(0,0) node (a11) {(\Kg X)_V \sma S^W}
(a11)++(\h,0) node (a12) {(\Kg X)_{V \oplus W}}
(a11)++(0,\u) node (a21) {\big( \txint^{\angordn \in \Gsk} (S^V)^{\sma\angordn} \sma X\angordn \big) \sma S^W}
(a12)++(0,\u) node (a22) {\txint^{\angordn \in \Gsk} (S^{V \oplus W})^{\sma\angordn} \sma X\angordn}
(a21)++(0,\v) node (a3) {\txint^{\angordn \in \Gsk} \big( (S^V)^{\sma\angordn} \sma S^W \big) \sma X\angordn}
;
\draw[1cell=.8]
(a11) edge node {\umu_{V,W}} (a12)
(a11) edge[equal] (a21)
(a12) edge[equal] (a22)
(a21) edge node[swap] {\iso} (a3)
(a3) [rounded corners=2pt] -| node[pos=.25] {\assm = \txint\! \assm_{\angordn} \sma 1} (a22)
;
\end{tikzpicture}
\end{equation}
\begin{itemize}
\item In \cref{Kgx_action_vw}, the pointed $G$-homeomorphism denoted by $\iso$ uses
\begin{itemize}
\item the fact that $- \sma S^W$ commutes with coends, and
\item the associativity isomorphism and braiding for the symmetric monoidal category $(\Gtopst,\sma)$ \cref{Gtopst_smc} to move $S^W$ to the left of $X\angordn$.
\end{itemize}
\item The pointed $G$-morphism $\assm$ is induced by the pointed $G$-morphisms
\begin{equation}\label{assembly_sph}
(S^V)^{\sma\angordn} \sma S^W \fto{\assm_{\angordn}} (S^{V \oplus W})^{\sma\angordn}
\end{equation}
for $\angordn \in \Gsk$ defined by the assignment
\[\begin{split}
& \big(\mq \smam\angordn \fto{t} S^V ; w \big) \in (S^V)^{\sma\angordn} \sma S^W\\
&\mapsto \big(\mq \smam\angordn \fto{t} S^V \fto{- \oplus w} S^{V \oplus W} \big) \in (S^{V \oplus W})^{\sma\angordn}.
\end{split}\]
\end{itemize}
\end{description}
This finishes the definition of $(\Kg X, \umu)$.
\end{definition}

\begin{lemma}\label{Kgx_welldef}
In \cref{def:ggspace_gspectra}, $(\Kg X, \umu)$ is an orthogonal $G$-spectrum.
\end{lemma}

\begin{proof}
We need to prove the following two statements.
\begin{enumerate}
\item\label{Kgx_welldef_i} $\Kg X$, defined in \cref{Kgxv,Kgxf}, is an $\IU$-space. 
\item\label{Kgx_welldef_ii} The right $\gsp$-action $\umu$ \cref{Kgx_sphere_action} makes $\Kg X$ into an orthogonal $G$-spectrum.
\end{enumerate}

\parhead{Statement \pcref{Kgx_welldef_i}: component objects}. 
To show that $\Kg X$ is an $\IU$-space, we first use \cref{expl:ggtop_obj} to describe the pointed $G$-space $(\Kg X)_V$ in \cref{Kgxv} more explicitly.  For an object $\angordn \in \Gsk$, recall that $\sma\angordn$ \cref{smash_Gskobjects} is a discrete trivial $G$-space with basepoint $0 \in \sma\angordn$.  There is a pointed $G$-homeomorphism
\begin{equation}\label{gtopst_sman_sv}
\begin{split}
(S^V)^{\sma\angordn} &= \Topgst(\sma\angordn, S^V) \\
&\iso \prod_{(\sma\angordn) \setminus \{0\}} S^V.
\end{split}
\end{equation}
In the last product in \cref{gtopst_sman_sv}, $G$ acts diagonally.  A typical point in it is denoted by
\begin{equation}\label{sv_sman_product}
t = \ang{t_{\anga} \in S^V}_{\anga \in (\sma\angordn) \setminus \{0\}}.
\end{equation}
The $G$-fixed basepoint of $(S^V)^{\sma\angordn}$ is the tuple $\ang{\infty}$, where each entry is the basepoint $\infty \in S^V$.

The pointed $G$-space $(\Kg X)_V$ in \cref{Kgxv} is a quotient of the wedge
\begin{equation}\label{Kgxv_wedge}
\bigvee_{\angordn \in \Gsk} (S^V)^{\sma\angordn} \sma X\angordn.
\end{equation}
A typical point in $(\Kg X)_V$ is represented by a pair
\begin{equation}\label{Kgxv_rep}
(t \in (S^V)^{\sma\angordn} ; x \in X\angordn) \forsomespace \angordn \in \Gsk.
\end{equation}
The pair $(t;x)$ represents the basepoint of $(\Kg X)_V$ if either $t$ or $x$ is the basepoint.  The coend $(\Kg X)_V$ in \cref{Kgxv} identifies, for each triple
\begin{equation}\label{t_upom_x}
\bx = \big(\mq \smam\angordn \fto{t} S^V ; \angordm \fto{\upom} \angordn \in \Gsk ; x \in X\angordm\big),
\end{equation}
the following two pairs in the wedge \cref{Kgxv_wedge}.
\begin{equation}\label{Kgxv_dzeroone}
\begin{split}
\dzero\bx &= \big(\mq \smam\angordm \fto{\sma\upom} \smam\angordn \fto{t} S^V ; x \in X\angordm\big) \\
\done\bx &= \big(\mq \smam\angordn \fto{t} S^V ; (X\upom)(x) \in X\angordn \big)
\end{split}
\end{equation}

\parhead{$G$-action}.  
The group $G$ acts diagonally on representative pairs \cref{Kgxv_rep}.  This means that
\begin{equation}\label{g_tx}
g(t;x) = (gt; gx) \forspace g \in G,
\end{equation}
where $gt$ is the composite pointed morphism
\[\sma\angordn \fto{t} S^V \fto[\iso]{g} S^V.\]
The identification \cref{Kgxv_dzeroone} is invariant under the $G$-action \cref{g_tx}, which means
\[g(\dzero\bx) = g(\done\bx),\]
because
\[X\angordm \fto{X\upom} X \angordn\]
is a pointed $G$-morphism \cref{f_upom_ggtop}.

\parhead{Component morphisms}. 
For each linear isometric isomorphism $f \cn V \fiso W$ in $\IU$, the pointed homeomorphism \cref{Kgxf}
\[(\Kg X)_V \fto[\iso]{(\Kg X)_f} (\Kg X)_W\]
is defined on representatives \cref{Kgxv_rep} by
\begin{equation}\label{Kgxf_rep}
(\Kg X)_f (t;x) = \big(\mq \smam\angordn \fto{t} S^V \fto[\iso]{\gsp f} S^W ; x \big)
\end{equation}
with $\gsp f$ defined in \cref{gsp_morphism}.  We stress that the pointed homeomorphism $(\Kg X)_f$ is \emph{not} required to be $G$-equivariant.  The following discussion proves that $\Kg X$ is an $\IU$-space, as described in \cref{expl:iu_space}.
\begin{description}
\item[Functoriality] $\Kg X$ preserves composition and identities \cref{iu_space_axioms} by \cref{Kgxf_rep} and the functoriality of $\gsp$ \cref{g_sphere}.
\item[Equivariance] $\Kg X$ satisfies the equivariance property \cref{x_gfginv} by the following computation for a non-basepoint $(t;x)$, using \cref{g_tx,Kgxf_rep}.
\begin{equation}\label{Kgx_equivariance}
\begin{split}
(\Kg X)_{gf\ginv} (t;x)
&= (gf\ginv t ; x) \\
&= g(f \ginv t ; \ginv x) \\
&= \big(g \circ (\Kg X)_f\big) (\ginv t ; \ginv x) \\ 
&= \big(g \circ (\Kg X)_f \circ \ginv\big) (t;x)
\end{split}
\end{equation}
\item[Continuity] For each pair of objects $(V,W) \in (\IU)^2$, to see that the pointed function \cref{iu_space_comp_mor}
\begin{equation}\label{Kgx_continuous}
\IU(V,W)_\splus \fto{\Kg X} \Topgst\big( (\Kg X)_V, (\Kg X)_W \big)
\end{equation}
is continuous, observe that the pointed function
\[\IU(V,W)_\splus \sma (S^V)^{\sma\angordn} \to (S^W)^{\sma\angordn},\]
defined by composition, is continuous for each object $\angordn \in \Gsk$.  By 
\begin{itemize}
\item the universal property of coends and
\item the commutation of $\IU(V,W)_\splus \sma -$ with coends,
\end{itemize} 
the pointed function
\[\IU(V,W)_\splus \sma (\Kg X)_V \to (\Kg X)_W\]
is continuous.  Its adjoint, which is continuous, is the pointed function in \cref{Kgx_continuous}.
\end{description}
This finishes the proof that $\Kg X$ is an $\IU$-space.

\parhead{Statement \pcref{Kgx_welldef_ii}: $\gsp$-module}.  
Next, we prove that the pointed $G$-morphisms \cref{Kgx_action_vw}
\[(\Kg X)_V \sma S^W \fto{\umu_{V,W}} (\Kg X)_{V \oplus W} \forspace (V,W) \in (\IUsk)^2\]
equip $\Kg X$ with the structure of an $\gsp$-module.  For a pair \cref{Kgxv_rep} 
\[\big(t \in (S^V)^{\sma\angordn} ; x \in X\angordn \big)\]
that represents a point in $(\Kg X)_V$ and a point $w \in S^W$, $\umu_{V,W}$ is given by
\begin{equation}\label{Kgxv_action_rep}
\umu_{V,W}(t; x; w) = \big(t \oplus w \in (S^{V \oplus W})^{\sma\angordn} ; x \big),
\end{equation}
where $t \oplus w$ denotes the composite
\begin{equation}\label{t_dash_w}
\sma\angordn \fto{t} S^V \fto{- \oplus w} S^{V \oplus W}.
\end{equation}
More explicitly, if $t = \ang{t_{\anga} \in S^V}$, where $\anga$ runs through $(\sma\angordn) \setminus \{0\}$ \cref{sv_sman_product}, then 
\begin{equation}\label{t_oplus_w}
t \oplus w = \ang{t_{\anga} \oplus w \in S^{V \oplus W}}_{\anga \in (\sma\angordn) \setminus \{0\}}.
\end{equation}
The following discussion proves that $(\Kg X,\umu)$ is an $\gsp$-module, as described in \cref{expl:gspectra}.
\begin{description}
\item[Well-defined $\gsp$-action]
First, we prove that $\umu_{V,W}$ \cref{Kgx_action_vw} is a pointed $G$-morphism between pointed $G$-spaces.  Using \cref{t_upom_x,Kgxv_dzeroone,Kgxv_action_rep,t_dash_w}, the following computation proves that the assignment $\umu_{V,W}$ is well defined.  To simplify the notation, $\umu_{V,W}$ is abbreviated to $\umu$.
\[\begin{split}
\umu(\dzero\bx ; w) 
&= \umu(t \circ \sma\upom ; x ; w) \\
&= \big((- \oplus w) \circ t \circ \sma\upom ; x \big) \\
&= \dzero\big((- \oplus w) \circ t ; \upom ; x \big) \\
&= \done\big((- \oplus w) \circ t ; \upom ; x \big) \\
&= \big((- \oplus w) \circ t ; (X\upom)(x) \big) \\
&= \umu\big(t ; (X\upom)(x) ; w\big) \\
&= \umu(\done\bx ; w) 
\end{split}\]
Next, using \cref{g_tx,Kgxv_action_rep,t_oplus_w}, the following computation proves that the pointed function $\umu_{V,W}$ is $G$-equivariant.
\[\begin{split}
\umu\big( g(t; x; w)\big) 
&= \umu(gt; gx; gw) \\
&= (gt \oplus gw; gx) \\
&= \big(g(t \oplus w) ; gx\big) \\
&= g(t \oplus w; x) \\
&= g\umu(t; x; w)
\end{split}\]
The continuity of the pointed function \cref{Kgx_action_vw}
\[(\Kg X)_V \sma S^W \fto{\umu_{V,W}} (\Kg X)_{V \oplus W}\]
follows from 
\begin{itemize}
\item the universal property of coends,
\item the description \cref{gtopst_sman_sv} of $(S^V)^{\sma\angordn}$ as a product of copies of $S^V$, and
\item the continuity of $\assm_{\angordn}$ \cref{assembly_sph} for each object $\angordn \in \Gsk$.
\end{itemize}

\item[Naturality] The naturality \cref{gsp_action_nat} of $\umu_{V,W}$ in $V$ and $W$ states that, for each pair of linear isometric isomorphisms
\[V \fto[\iso]{f} V' \andspace W \fto[\iso]{h} W' \inspace \IUsk,\]
the following diagram commutes.
\begin{equation}\label{Kgx_action_nat}
\begin{tikzpicture}[vcenter]
\def\v{-1.4}
\draw[0cell]
(0,0) node (a11) {(\Kg X)_V \sma S^W}
(a11)++(4,0) node (a12) {(\Kg X)_{V \oplus W}}
(a11)++(0,\v) node (a21) {(\Kg X)_{V'} \sma S^{W'}}
(a12)++(0,\v) node (a22) {(\Kg X)_{V' \oplus W'}}
;
\draw[1cell=.9]
(a11) edge[transform canvas={yshift=-.2ex}] node {\umu_{V,W}} (a12)
(a12) edge[transform canvas={xshift=-1.5em}] node[pos=.55] {(\Kg X)_{f \oplus h}} node[swap] {\iso} (a22)
(a11) edge[transform canvas={xshift=2em}] node[swap,pos=.6] {(\Kg X)_f \sma h} node[pos=.55] {\iso} (a21)
(a21) edge[transform canvas={yshift=-.2ex}] node {\umu_{V',W'}} (a22)
;
\end{tikzpicture}
\end{equation}
Using \cref{Kgxf_rep,Kgxv_action_rep,t_oplus_w}, each of the two composites in \cref{Kgx_action_nat} sends a representative triple $(t; x; w) $ to the pair
\[\big( \ang{f t_{\anga} \oplus hw \in S^{V' \oplus W'}}_{\anga \in (\sma\angordn) \setminus \{0\}} ; x \big),\]
which represents a point in $(\Kg X)_{V' \oplus W'}$.
\item[Unity] 
Using
\begin{itemize}
\item \cref{Kgxv_action_rep} with $W = 0$ and
\item \cref{Kgxf_rep} for the right unit isomorphism $\rho \cn V \oplus 0 \iso V$ in $(\IUsk,\oplus)$,
\end{itemize}
the following equalities prove the unity axiom \cref{gsp_unity} for $(\Kg X,\umu)$.
\[\begin{split}
&(\Kg X)_\rho (\umu_{V,0}) (t ; x ; 0) \\
&= (\Kg X)_\rho \big(t \oplus 0 \in (S^{V \oplus 0})^{\sma\angordn} ; x \big) \\
&= (t ; x)
\end{split}\]
\item[Associativity]
The associativity axiom \cref{gsp_assoc} for $(\Kg X,\umu)$ states that,  for each triple of objects $(U,V,W) \in (\IUsk)^3$, the following diagram commutes.
\begin{equation}\label{Kgx_assoc}
\begin{tikzpicture}[vcenter]
\def\h{4} \def\v{-1.4}
\draw[0cell=.9]
(0,0) node (a11) {((\Kg X)_U \sma S^V) \sma S^W}
(a11)++(\h,0) node (a12) {(\Kg X)_U \sma (S^V \sma S^W)}
(a11)++(0,\v) node (a21) {(\Kg X)_{U \oplus V} \sma S^W}
(a12)++(0,\v) node (a22) {(\Kg X)_U \sma S^{V \oplus W}}
(a21)++(0,\v) node (a31) {(\Kg X)_{(U \oplus V) \oplus W}}
(a31)++(\h,0) node (a32) {(\Kg X)_{U \oplus (V \oplus W)}}
;
\draw[1cell=.8]
(a11) edge node {\al} node[swap] {\iso} (a12)
(a12) edge node {1 \sma \mu_{V,W}} node[swap] {\iso} (a22)
(a22) edge node {\umu_{U,V \oplus W}} (a32)
(a11) edge node[swap] {\umu_{U,V} \sma 1} (a21)
(a21) edge node[swap] {\umu_{U \oplus V,W}} (a31)
(a31) edge[transform canvas={yshift=.3ex}] node {(\Kg X)_\al} node[swap] {\iso} (a32)
;
\end{tikzpicture}
\end{equation}
To see that the diagram \cref{Kgx_assoc} commutes, consider $v \in S^V$, $w \in S^W$, and a pair \cref{Kgxv_rep}
\[\big(t \in (S^U)^{\sma\angordn} ; x \in X\angordn \big)\]
that represents a point in $(\Kg X)_U$.  By \cref{gsp_mult,Kgxf_rep,Kgxv_action_rep}, each of the two composites in \cref{Kgx_assoc} sends the quadruple
\[\big(((t ; x) ; v) ; w \big) \inspace ((\Kg X)_U \sma S^V) \sma S^W\]
to the pair
\[\big(t \oplus (v \oplus w) \in (S^{U \oplus (V \oplus W)})^{\sma\angordn} ; x \big),\]
which represents a point in $(\Kg X)_{U \oplus (V \oplus W)}$.
\end{description}
This finishes the proof that $(\Kg X,\umu)$ is an orthogonal $G$-spectrum.
\end{proof}

\section{A $\Gtop$-Functor from $\Gskg$-Spaces to $G$-Spectra}
\label{sec:Kg_gtop_functor}

Using the object assignment in \cref{def:ggspace_gspectra}, this section constructs the $\Gtop$-functor
\[\GGTop \fto{\Kg} \GSp\]
from the $\Gtop$-category $\GGTop$ of $\Gskg$-spaces \pcref{expl:ggtop_gtopenr} to the $\Gtop$-category $\GSp$ of orthogonal $G$-spectra \pcref{expl:gspectra,expl:gsp_morphism}.

\secoutline
\begin{itemize}
\item \cref{def:ggtop_gsp_mor} defines the data for $\Kg$.
\item \cref{Kg_gtop_functor} proves that $\Kg$ is a well-defined $\Gtop$-functor.
\end{itemize}

\begin{definition}\label{def:ggtop_gsp_mor}
The $\Gtop$-functor\index{G-G-space@$\Gskg$-space!to orthogonal G-spectra@to orthogonal $G$-spectra}\index{orthogonal G-spectrum@orthogonal $G$-spectrum!from G-G-space@from $\Gskg$-space} \pcref{def:enriched-functor}
\[\GGTop \fto{\Kg} \GSp\]
is defined as follows.
\begin{description}
\item[Objects] The object assignment 
\[X \mapsto (\Kg X,\umu),\] 
sending each $\Gskg$-space $X$ to an orthogonal $G$-spectrum $(\Kg X,\umu)$, is given in \cref{def:ggspace_gspectra}.
\item[Hom $G$-spaces]
Suppose $X, X' \cn \Gsk \to \Gtopst$ are $\Gskg$-spaces \cref{ggtop_obj}.  Using \cref{gsp_gtop_enr,ggtop_gtop_enr}, the $G$-morphism between hom $G$-spaces
\begin{equation}\label{Kg_mor_gspaces}
\brk{X}{X'} \fto{\Kg} \GSp\big((\Kg X,\umu), (\Kg X',\umu')\big)
\end{equation}
sends a point
\[\theta = \big\{X\angordn \fto{\theta_{\angordn}} X'\angordn \big\}_{\angordn \in \Gsk} \in \brk{X}{X'}\]
to the point
\begin{equation}\label{Kg_theta}
\Kg\theta = \big\{(\Kg\theta)_V \big\}_{V \in \IUsk} \in \GSp\big((\Kg X,\umu), (\Kg X',\umu')\big).
\end{equation}
For each object $V \in \IUsk$, the $V$-component pointed morphism of $\Kg\theta$ \cref{gsp_mor_comp} is defined by the following commutative diagram in $\Topst$.
\begin{equation}\label{Kg_theta_v}
\begin{tikzpicture}[vcenter]
\def\v{-1.4}
\draw[0cell=.9]
(0,0) node (a11) {(\Kg X)_V}
(a11)++(3,0) node (a12) {\txint^{\angordn \in \Gsk} (S^V)^{\sma\angordn} \sma X\angordn}
(a11)++(0,\v) node (a21) {(\Kg X')_V}
(a12)++(0,\v) node (a22) {\txint^{\angordn \in \Gsk} (S^V)^{\sma\angordn} \sma X'\angordn}
;
\draw[1cell=.8]
(a11) edge[equal] (a12)
(a21) edge[equal] (a22)
(a11) edge[transform canvas={xshift=1em}] node[swap] {(\Kg\theta)_V} (a21)
(a12) edge[transform canvas={xshift=-2em}] node {\txint^{\angordn} 1 \sma \theta_{\angordn}} (a22)
;
\end{tikzpicture}
\end{equation}
\end{description}
This finishes the definition of $\Kg$.
\end{definition}

\begin{explanation}\label{expl:Kg_theta_v}
$(\Kg\theta)_V$ sends a pair \cref{Kgxv_rep}
\[(t \in (S^V)^{\sma\angordn} ; x \in X\angordn),\]
which represents a point in $(\Kg X)_V$, to the pair
\begin{equation}\label{Kg_theta_vtx}
(\Kg\theta)_V (t; x) = (t; \theta_{\angordn} x \in X'\angordn),
\end{equation}
which represents a point in $(\Kg X')_V$.
\end{explanation}

\begin{lemma}\label{Kg_gtop_functor}
In \cref{def:ggtop_gsp_mor},
\[\GGTop \fto{\Kg} \GSp\] 
is a $\Gtop$-functor.
\end{lemma}

\begin{proof}
The object assignment of $\Kg$ is well defined by \cref{Kgx_welldef}.

\parhead{Component morphisms}. 
To show that the component $G$-morphism $\Kg$ in \cref{Kg_mor_gspaces} is well defined, we first observe that 
\[(\Kg X,\umu) \fto{\Kg\theta}(\Kg X',\umu'),\]
whose components are defined in \cref{Kg_theta_v}, is a morphism of $\gsp$-modules, using \cref{expl:gsp_morphism}.
\begin{description}
\item[Naturality] For each isomorphism $f \cn V \fiso W$ in $\IUsk$, the naturality diagram for $\Kg\theta$ is the following diagram of pointed morphisms.
\begin{equation}\label{Kg_theta_natural}
\begin{tikzpicture}[vcenter]
\def\v{-1.4}
\draw[0cell]
(0,0) node (a11) {(\Kg X)_V}
(a11)++(3.5,0) node (a12) {(\Kg X')_V}
(a11)++(0,\v) node (a21) {(\Kg X)_W}
(a12)++(0,\v) node (a22) {(\Kg X')_W}
;
\draw[1cell=.9]
(a11) edge node {(\Kg\theta)_V} (a12)
(a12) edge[transform canvas={xshift=-1em}] node {(\Kg X')_f} (a22)
(a11) edge[transform canvas={xshift=1em}] node[swap] {(\Kg X)_f} (a21)
(a21) edge node {(\Kg\theta)_W} (a22)
;
\end{tikzpicture}
\end{equation}
The diagram \cref{Kg_theta_natural} commutes by the universal property of coends, \cref{Kgxf}, \cref{Kg_theta_v}, and the following diagrams for $\angordn \in \Gsk$, which commute by the functoriality of the smash product for pointed spaces.
\[\begin{tikzpicture}[vcenter]
\def\v{-1.4}
\draw[0cell=.9]
(0,0) node (a11) {(S^V)^{\sma\angordn} \sma X\angordn}
(a11)++(4.2,0) node (a12) {(S^V)^{\sma\angordn} \sma X'\angordn}
(a11)++(0,\v) node (a21) {(S^W)^{\sma\angordn} \sma X\angordn}
(a12)++(0,\v) node (a22) {(S^W)^{\sma\angordn} \sma X'\angordn}
;
\draw[1cell=.8]
(a11) edge node {1 \sma \theta_{\angordn}} (a12)
(a12) edge[transform canvas={xshift=-2em}] node {(f \comp -) \sma 1} (a22)
(a11) edge[transform canvas={xshift=2em}] node[swap] {(f \comp -) \sma 1} (a21)
(a21) edge node {1 \sma \theta_{\angordn}} (a22)
;
\end{tikzpicture}\]
\item[Compatibility]
The compatibility diagram \cref{gsp_mor_axiom} for $\Kg\theta$ is the following diagram of pointed morphisms for objects $(V,W) \in (\IUsk)^2$.
\begin{equation}\label{Kg_theta_compat}
\begin{tikzpicture}[vcenter]
\def\v{-1.4}
\draw[0cell=.9]
(0,0) node (a11) {(\Kg X)_V \sma S^W}
(a11)++(3.6,0) node (a12) {(\Kg X)_{V \oplus W}}
(a11)++(0,\v) node (a21) {(\Kg X')_V \sma S^W}
(a12)++(0,\v) node (a22) {(\Kg X')_{V \oplus W}}
;
\draw[1cell=.8]
(a11) edge node {\umu_{V,W}} (a12)
(a12) edge[transform canvas={xshift=-1.5em}] node {(\Kg \theta)_{V \oplus W}} (a22)
(a11) edge[transform canvas={xshift=1em}] node[swap] {(\Kg\theta)_V \sma 1} (a21)
(a21) edge node {\umu'_{V,W}} (a22)
;
\end{tikzpicture}
\end{equation}
The diagram \cref{Kg_theta_compat} commutes by the universal property of coends, \cref{Kgx_action_vw}, \cref{Kg_theta_v}, and the following diagrams for $\angordn \in \Gsk$, which commute by the functoriality of $\sma$.
\[\begin{tikzpicture}[vcenter]
\def\v{-1.4}
\draw[0cell=.9]
(0,0) node (a11) {((S^V)^{\sma\angordn} \sma S^W) \sma X\angordn}
(a11)++(5,0) node (a12) {(S^{V \oplus W})^{\sma\angordn} \sma X\angordn}
(a11)++(0,\v) node (a21) {((S^V)^{\sma\angordn} \sma S^W) \sma X'\angordn}
(a12)++(0,\v) node (a22) {(S^{V \oplus W})^{\sma\angordn} \sma X'\angordn}
;
\draw[1cell=.8]
(a11) edge node {\assm_{\angordn} \sma 1} (a12)
(a12) edge[transform canvas={xshift=-2em}] node {1 \sma \theta_{\angordn}} (a22)
(a11) edge[transform canvas={xshift=2.5em}] node[swap] {1 \sma \theta_{\angordn}} (a21)
(a21) edge node {\assm_{\angordn} \sma 1} (a22)
;
\end{tikzpicture}\]
\end{description}
This finishes the proof that $\Kg\theta$ in \cref{Kg_theta} is a morphism of $\gsp$-modules.

\parhead{Equivariance}. 
The $G$-equivariance of the function $\theta \mapsto \Kg\theta$ in \cref{Kg_mor_gspaces} means the equality of $\gsp$-module morphisms
\begin{equation}\label{Kg_theta_gequiv}
\Kg(g \cdot \theta) = g \cdot (\Kg\theta)
\end{equation}
for each point $\theta \in \brk{X}{X'}$ and $g \in G$.  Using \cref{gsp_gtop_enr,ggtop_gtopenr_theta_g,g_tx}, the following computation for $V \in \IUsk$ proves the equality \cref{Kg_theta_gequiv}.
\[\begin{split}
\big( \Kg(g \cdot \theta)\big)_V
&= \txint^{\angordn} 1 \sma (g \theta_{\angordn} \ginv) \\
& = \txint^{\angordn} \big[(g \sma g) (1 \sma \theta_{\angordn}) (\ginv \sma \ginv) \big] \\
&= g \big(\txint^{\angordn} 1 \sma \theta_{\angordn} \big) \ginv \\
&= \big(g \cdot (\Kg\theta)\big)_V
\end{split}\]

\parhead{Continuity}. 
By \cref{ggtop_gtop_enr} and \cref{gsp_gtop_enr}, the domain and codomain of the function $\theta \mapsto \Kg\theta$ in \cref{Kg_mor_gspaces} are subspaces of, respectively, 
\[\prod_{\angordn \in \Gsk} \Topgst(X\angordn, X'\angordn) \andspace 
\prod_{V \in \IUsk} \Topgst\big((\Kg X)_V,(\Kg X')_V\big).\]
Thus, the function $\Kg$ is continuous by the definition \cref{Kg_theta_v} of $(\Kg \theta)_V$ for each $V \in \IUsk$.  Thus, $\Kg$ in \cref{Kg_mor_gspaces} is a $G$-morphism between hom $G$-spaces.

\parhead{$\Gtop$-functoriality}. 
$\Kg$ preserves composition and identities, in the $\Gtop$-enriched sense \cref{eq:enriched-composition}, by
\begin{itemize}
\item the fact that, in both its domain $\GGTop$ and codomain $\GSp$, composition and identities are defined componentwise;
\item the definition \cref{Kg_theta_v} of $(\Kg \theta)_V$; 
\item the universal property of coends; and
\item the functoriality of $\sma$.
\end{itemize}
This finishes the proof that $\Kg$ is a $\Gtop$-functor.
\end{proof}

\section{Unit Constraint}
\label{sec:Kgzero}

Having constructed the $\Gtop$-functor \pcref{def:ggtop_gsp_mor}
\[\GGTop \fto{\Kg} \GSp,\]
this section begins the extension of $\Kg$ to a symmetric monoidal $\Gtop$-functor \pcref{definition:monoidal-V-fun,definition:braided-monoidal-vfunctor} between the symmetric monoidal $\Gtop$-categories $\GGTop$ \pcref{expl:ggtop_gtopenr} and $\GSp$ \pcref{def:gsp_module,def:gsp_smgtop}.  This section constructs the $\Gtop$-enriched unit constraint, in the sense of \cref{enr_constraints}, for $\Kg$.  The monoidal constraint for $\Kg$ is constructed in \cref{sec:Kgtwo}.

\secoutline
\begin{itemize}
\item \cref{def:Kg_zero} constructs the unit constraint $\Kgzero$ for $\Kg$.
\item \cref{Kgzero_welldef} proves that $\Kgzero$ is a $\Gtop$-natural isomorphism.
\item \cref{expl:Kgzero} describes $\Kgzero$ and its inverse at the point-set level.
\end{itemize}

The following definition uses the notion of a $\V$-natural transformation for the Cartesian closed category $\Gtop$ \pcref{def:Gtop,def:enriched-natural-transformation}.

\begin{definition}\label{def:Kg_zero}
The \emph{unit constraint} for the $\Gtop$-functor \pcref{def:ggtop_gsp_mor}
\[\GGTop \fto{\Kg} \GSp\] 
is the $\Gtop$-natural isomorphism
\begin{equation}\label{Kgzero}
\begin{tikzpicture}[vcenter]
\def\h{3}
\draw[0cell]
(0,0) node (a1) {\vtensorunit}
(a1)++(.5*\h,-1) node (a2) {\GGTop}
(a1)++(\h,0) node (a3) {\GSp}
;
\draw[1cell=.9]
(a1) edge[transform canvas={yshift=.3ex}] node {\gsp} (a3)
(a1) [rounded corners=2pt] |- node[swap,pos=.25] {\gu} (a2)
;
\draw[1cell=.9]
(a2) [rounded corners=2pt] -| node[swap,pos=.75] {\Kg} (a3)
;
\draw[2cell]
node[between=a1 and a3 at .4, shift={(0,-.45)}, rotate=-90, 2label={above,\Kgzero}, 2label={below,\iso}] {\Rightarrow}
;
\end{tikzpicture}
\end{equation}
defined as follows.
\begin{itemize}
\item $\vtensorunit$ is the unit $\Gtop$-category \pcref{definition:unit-vcat}, with one object $*$ and unique hom $G$-space $\vtensorunit(*,*)$ given by the monoidal unit $* \in \Gtop$.
\item $\gsp \cn \vtensorunit \to \GSp$ is the monoidal identity of the symmetric monoidal $\Gtop$-category $\GSp$ \cref{gsp_mon_id}.  It sends the unique object $* \in \vtensorunit$ to the $G$-sphere $\gsp$ \pcref{def:g_sphere}.
\item $\gu \cn \vtensorunit \to \GGTop$ is the monoidal identity of the symmetric monoidal $\Gtop$-category $\GGTop$ \cref{ggtop_mon_identity}.  It sends the unique object $* \in \vtensorunit$ to the monoidal unit $\Gskg$-space $\gu$ in \cref{ggtop_unit}.
\end{itemize}
The unique component of $\Kgzero$ is the $G$-morphism
\[* \fto{\Kgzero} \GSp(\gsp,\Kg\gu)\]
determined by the $G$-equivariant $\gsp$-module morphism \pcref{expl:gsp_morphism}
\begin{equation}\label{Kgzero_gspmod_mor}
(\gsp,\mu) \fto{\Kgzero} (\Kg\gu,\umu).
\end{equation}
Using \cref{smash_Gskobjects,gsp_v,gu_angordm_ggtop,Kgxv}, for each object $V \in \IUsk$, the $V$-component pointed $G$-morphism \cref{gsp_mor_comp} of $\Kgzero$ is defined as the following composite in $\Gtopst$.
\begin{equation}\label{Kgzero_v}
\begin{tikzpicture}[vcenter]
\def\u{-1} \def\v{-1.3}
\draw[0cell=.9]
(0,0) node (a11) {\gsp(V)}
(a11)++(3.5,0) node (a12) {(\Kg\gu)_V}
(a11)++(0,\u) node (a21) {S^V}
(a12)++(0,\u) node (a22) {\txint^{\angordn \in \Gsk} (S^V)^{\sma \angordn} \sma \gu\angordn}
(a21)++(0,\v) node (a31) {\Topgst(\ord{1}, S^V) \sma \ord{1}}
(a22)++(0,\v) node (a32) {(S^V)^{\sma\ang{}} \sma \gu\ang{}}
;
\draw[1cell=.9]
(a11) edge node {\Kgzero_V} (a12)
(a11) edge[equal] (a21)
(a12) edge[equal] (a22)
(a21) edge node[swap] {\rhoinv} node {\iso} (a31)
(a31) edge[equal] (a32)
(a32) edge node[swap,pos=.4] {\inc_{\ang{}}} (a22)
;
\end{tikzpicture}
\end{equation}
The details of \cref{Kgzero_v} are explained below.
\begin{itemize}
\item The pointed $G$-homeomorphism $\rho$ is the right unit isomorphism for $(\Gtopst,\sma)$.
\item The bottom equality uses the equalities
\[\sma\ang{} = \ord{1} = \stplus = \txwedge_{\Gskpunc(\ang{}, \ang{})} \stplus = \gu\ang{}.\]
\item The $G$-morphism $\inc_{\ang{}}$ is part of the definition of the coend for the empty tuple $\ang{} \in \Gsk$.
\end{itemize}
This finishes the definition of $\Kgzero$.
\end{definition}

\begin{lemma}\label{Kgzero_welldef}
In \cref{def:Kg_zero}, $\Kgzero$ is a $\Gtop$-natural isomorphism.
\end{lemma}

\begin{proof}
\parhead{Equivariant component}.
We first observe that $\Kgzero$ in \cref{Kgzero_gspmod_mor} is a $G$-equivariant $\gsp$-module morphism \pcref{expl:gsp_morphism}.
\begin{description}
\item[Naturality] The naturality of $\Kgzero_V$ \cref{Kgzero_v} in the variable $V \in \IUsk$ \cref{iu_mor_natural} follows from 
\begin{itemize}
\item the definition \cref{gsp_morphism} of $\gsp$ on morphisms,
\item the naturality of the right unit isomorphism $\rho$ for $(\Gtopst,\sma)$, and
\item the naturality of the $G$-morphism $\inc_{\ang{}}$, using \cref{Kgxf}.
\end{itemize}
\item[Compatibility] The compatibility of $\Kgzero$ with the right $\gsp$-actions \cref{gsp_mor_axiom} follows from
\begin{itemize}
\item the definition \cref{gsphere_multiplication} of the multiplication $\mu_{V,W}$ for $\gsp$ and
\item the definition \cref{Kgx_action_vw} of the $\gsp$-action $\umu_{V,W}$ for $\Kg\gu$.
\end{itemize}
\item[$G$-equivariance] By \cref{gsp_gtop_enr}, the $G$-equivariance of $\Kgzero$ means that each component $\Kgzero_V$ \cref{Kgzero_v} is $G$-equivariant, which is true because it is a composite of two $G$-morphisms.
\end{description}
This finishes the proof that $\Kgzero$ in \cref{Kgzero_gspmod_mor} is a $G$-equivariant $\gsp$-module morphism.

\parhead{$\Gtop$-naturality}. 
The naturality diagram \cref{enr_naturality} for $\Kgzero$ \cref{Kgzero} commutes because $\vtensorunit$ is the unit $\Gtop$-category.  Thus, $\Kgzero$ is a $\Gtop$-natural transformation.

\parhead{Invertibility}. 
To show that $\Kgzero$ \cref{Kgzero} is a $\Gtop$-natural isomorphism, we observe that, for each object $V \in \IUsk$, its $V$-component $\Kgzero_V$ \cref{Kgzero_v} is also given by the following composite of pointed $G$-homeomorphisms.
\begin{equation}\label{Kgzero_iso}
\begin{split}
S^V &\iso (S^V)^{\sma \ang{}} \\
&\iso \txint^{\angordn \in \Gsk} \txwedge_{\Gskpunc(\ang{}, \angordn)} (S^V)^{\sma\angordn} \\
&\iso \txint^{\angordn \in \Gsk} \txwedge_{\Gskpunc(\ang{}, \angordn)} [(S^V)^{\sma\angordn} \sma \stplus] \\
&\iso \txint^{\angordn \in \Gsk} [(S^V)^{\sma\angordn} \sma (\txwedge_{\Gskpunc(\ang{}, \angordn)} \stplus)] \\ 
&= (\Kg\gu)_V
\end{split}
\end{equation}
The details of \cref{Kgzero_iso} are explained below.
\begin{itemize}
\item The first $\iso$ follows from $\sma\ang{} = \ord{1}$ \cref{smash_Gskobjects}.
\item The second $\iso$ uses the Yoneda Density Theorem for coends (\cite[A.5.7]{loregian} or \cite[3.7.15]{cerberusIII}).
\item The third $\iso$ uses the right unit isomorphism $\rho$ for $(\Gtopst, \sma, \stplus)$.
\item The fourth $\iso$ uses the commutation of $(S^V)^{\sma\angordn} \sma -$ with wedges.
\item The last equality uses the definitions \cref{gu_angordm_ggtop,Kgxv} of $\gu\angordn$ and $(\Kg\gu)_V$.
\end{itemize}
This finishes the proof that $\Kgzero$ is a $\Gtop$-natural isomorphism.
\end{proof}

\begin{explanation}[Unit Constraint for $\Kg$]\label{expl:Kgzero}
We unravel $\Kgzero$ and its inverse in terms of representatives \cref{Kgxv_rep}.  For each object $V \in \IUsk$, the pointed $G$-homeomorphism \cref{Kgzero_v}
\[S^V \fto[\iso]{\Kgzero_V} (\Kg\gu)_V = \txint^{\angordn \in \Gsk} (S^V)^{\sma\angordn} \sma \gu\angordn\]
sends each point $v \in V \subset S^V$ to the representative
\[\Kgzero_V(v) = \big(\mq\smam\ang{} = \ord{1} \fto{v} S^V ; 1 \in \ord{1} = \gu\ang{} \big).\]
The pointed morphism $\sma\ang{} \fto{v} S^V$ sends the non-basepoint $1 \in \sma\ang{} = \{0,1\}$ to the given point $v \in V$.  

\parhead{Inverse}.  The inverse 
\[(\Kg\gu)_V \fto{(\Kgzero_V)^{-1}} S^V\]
sends a representative pair \cref{Kgxv_rep}
\begin{equation}\label{Kgjv_rep}
\big(\mq \smam\angordn \fto{t} S^V ; \anga \in \gu\angordn \iso \sma\angordn\big)
\end{equation}
to the image of $\anga$ under $t$:
\[(\Kgzero_V)^{-1}(t ; \anga) = t_{\anga} \in S^V.\]
The pointed morphism $(\Kgzero_V)^{-1}$ is well defined because, with $X = \gu$ and $x = \anga \in \sma\angordm$ in \cref{Kgxv_dzeroone}, both $\dzero\bx$ and $\done\bx$ are sent by $(\Kgzero_V)^{-1}$ to $t_{(\sma\upom) \anga}$. 

\parhead{Composites}.  The composite
\[S^V \fto{\Kgzero_V} (\Kg\gu)_V \fto{(\Kgzero_V)^{-1}} S^V\]
is the identity morphism because $\ord{1} \fto{v} S^V$ sends $1 \in \ord{1}$ to $v \in V$.  

To see that the other composite
\[(\Kg\gu)_V \fto{(\Kgzero_V)^{-1}} S^V \fto{\Kgzero_V} (\Kg\gu)_V\]
is also equal to the identity, consider a pair $(t;\anga)$ \cref{Kgjv_rep} and the morphism \cref{Gsk_morphisms}
\[\ang{} \fto{\upom} \angordn \inspace \Gsk\]
determined by the point \cref{smash_fpsi}
\[(\sma\upom)(1) = \anga \in \sma\angordn.\]  
Using the identification \cref{Kgxv_dzeroone}, the following equalities in $(\Kg\gu)_V$ prove that $\Kgzero_V(\Kgzero_V)^{-1}$ is equal to the identity.
\[\begin{split}
& \Kgzero_V(\Kgzero_V)^{-1}(t; \anga) \\
&= \Kgzero_V(t_{\anga}) \\
&= \big(\mq \smam\ang{} \fto{t_{\anga}} S^V ; 1 \in \gu\ang{} \big) \\
&= \dzero(t ; \upom ; 1) \\
&= \done(t ; \upom ; 1) \\
&= (t; (\sma\upom)(1)) = (t; \anga)
\end{split}\]
Thus, $\Kgzero_V$ and $(\Kgzero_V)^{-1}$ are actually inverses of each other.
\end{explanation}

\section{Monoidal Constraint}
\label{sec:Kgtwo}

This section constructs the $\Gtop$-enriched monoidal constraint \cref{enr_constraints} for the $\Gtop$-functor \pcref{def:ggtop_gsp_mor}
\[\GGTop \fto{\Kg} \GSp\]
between the symmetric monoidal $\Gtop$-categories $\GGTop$ \pcref{expl:ggtop_gtopenr} and $\GSp$  \pcref{def:gsp_module,def:gsp_smgtop}.  

\secoutline
\begin{itemize}
\item \cref{def:Kgtwo} defines the monoidal constraint for $\Kg$, which is a $\Gtop$-natural transformation denoted $\Kgtwo$.
\item \cref{expl:ktwo_u_dom,expl:ktwo_u_cod,expl:ktwo_u} unravel \cref{def:Kgtwo} in detail.
\item \cref{ktwo_u_welldef,kgtwo_xy_welldef,kgtwo_welldef} prove that $\Kgtwo$ is a $\Gtop$-natural transformation.
\end{itemize}

\begin{definition}\label{def:Kgtwo}
The \emph{monoidal constraint} for the $\Gtop$-functor \pcref{def:ggtop_gsp_mor}
\[\GGTop \fto{\Kg} \GSp\] 
is the $\Gtop$-natural transformation \pcref{def:enriched-natural-transformation}
\begin{equation}\label{Kgtwo}
\begin{tikzpicture}[vcenter]
\def\v{-1.4}
\draw[0cell=.9]
(0,0) node (a11) {(\GGTop)^{\otimes 2}}
(a11)++(3.5,0) node (a12) {(\GSp)^{\otimes 2}}
(a11)++(0,\v) node (a21) {\GGTop}
(a12)++(0,\v) node (a22) {\GSp}
;
\draw[1cell=.8]
(a11) edge node {(\Kg)^{\otimes 2}} (a12)
(a12) edge node {\smasg} (a22)
(a11) edge node[swap] {\smag} (a21)
(a21) edge node {\Kg} (a22)
;
\draw[2cell]
node[between=a11 and a22 at .55, rotate=225, shift={(0,.1*\v)}, 2labelw={below,\Kgtwo,-1pt}] {\Rightarrow}
;
\end{tikzpicture}
\end{equation}
defined as follows.
\begin{itemize}
\item $\smag$ is the monoidal composition \cref{ggtop_mon_composition} for the symmetric monoidal $\Gtop$-category $\GGTop$.  It is defined
\begin{itemize}
\item on objects by the monoidal product of $\Gskg$-spaces \pcref{expl:ggtop_smag} and
\item on hom $G$-spaces in \cref{ggtop_moncomp_mor}.
\end{itemize}
\item $\smasg$ is the monoidal composition \cref{gsp_mon_comp} for the symmetric monoidal $\Gtop$-category $\GSp$.  It is defined on $\gsp$-modules and their morphisms in \cref{def:gsp_sma,def:gsp_mor_sma}.  These are induced by the monoidal composition $\smau$ for $\IUT$ \cref{IU_mon_composition}, which is defined on objects and underlying morphisms in \cref{def:iuspace_day,def:iumor_day}.
\end{itemize}
The component of $\Kgtwo$ at a pair $(X,Y)$ of $\Gskg$-spaces \cref{ggtop_obj} is the $G$-morphism
\[* \fto{\Kgtwo_{X,Y}} \GSp\big(\Kg X \smasg \Kg Y, \Kg(X \smag Y) \big)\]
determined by the $G$-equivariant $\gsp$-module morphism \pcref{expl:gsp_morphism}
\begin{equation}\label{Kgtwo_component}
\Kg X \smasg \Kg Y \fto{\Kgtwo_{X,Y}} \Kg(X \smag Y).
\end{equation}
The underlying $G$-equivariant $\IU$-morphism (\cref{def:iu_morphism} \pcref{def:iu_morphism_iii}) of $\Kgtwo_{X,Y}$ is induced by the $G$-equivariant $\IU$-morphism $\ktwo$ in the following diagram in $\IUT$ \pcref{def:IU_smgtop}.  The top row is an instance of the coequalizer in \cref{gsp_sma_coequal}, where $\umu$ and $\umu'$ are the right $\gsp$-actions on $\Kg X$ and $\Kg Y$ \cref{Kgx_sphere_action}.
\begin{equation}\label{smau_kgx}
\begin{tikzpicture}[vcenter]
\draw[0cell=.8]
(0,0) node (a1) {\phantom{(\Kg X \smau \gsp) \smau \Kg Y}}
(a1)++(0,-.04) node (a1') {(\Kg X \smau \gsp) \smau \Kg Y}
(a1)++(4,0) node (a2) {\Kg X \smau \Kg Y}
(a2)++(2.7,0) node (a3) {\Kg X \smasg \Kg Y}
(a2)++(0,-1.4) node (b) {\Kg(X \smag Y)}
;
\draw[1cell=.75]
(a1) edge[transform canvas={yshift=.4ex}] node {\umu \smau 1} (a2)
(a1) edge[transform canvas={yshift=-.5ex}] node[swap] {(1 \smau \umu') \upbe} (a2)
(a2) edge node {\psma} (a3)
(a2) edge node[swap] {\ktwo} (b)
(a3) edge[dashed] node {\Kgtwo_{X,Y}} (b)
;
\end{tikzpicture}
\end{equation}
Recall from \cref{expl:iu_morphism,expl:eqiu_morphism} that a $G$-equivariant  $\IU$-morphism is determined by its set of component pointed $G$-morphisms for objects in $\IUsk$ \pcref{def:IU_spaces}.

For each object $U \in \IUsk$, the $U$-component of $\ktwo$ is the pointed $G$-morphism $\ktwo_U$ defined by the following commutative diagrams in $\Gtopst$ for $(V,W) \in (\IUsk)^2$ and $(\angordm, \angordmp) \in \Gsk^2$, with further explanation given after the diagram.
\begin{equation}\label{ktwo_u}
\begin{tikzpicture}[vcenter]
\def\h{4.4} \def\u{-1} \def\v{-1.4} \def\w{-1.5} \def\x{-1.7}
\draw[0cell=.7]
(0,0) node (a11) {(\Kg X \smau \Kg Y)_U}
(a11)++(\h,0) node (a12) {\big(\Kg(X \smag Y)\big)_U}
(a11)++(0,\u) node[align=left] (a21) {$\txint^{V,W} \IU(V\oplus W, U)_\splus \sma$\\ 
$\phantom{\txint^{V,W}} [(\Kg X)_V \sma (\Kg Y)_W]$}
(a12)++(0,\u) node (a22) {\txint^{\angordn} (S^U)^{\sma\angordn} \sma (X \smag Y)\angordn}
(a21)++(0,\w) node[align=left] (a31) {$\txint^{V,W} \IU(V\oplus W, U)_\splus \sma$\\
$\phantom{\txint^{V,W}} \big[\big(\txint^{\angordm} (S^V)^{\sma\angordm} \sma X\angordm\big) 
\sma$\\
$\phantom{\txint^{V,W}} \big(\txint^{\angordmp} (S^W)^{\sma\angordmp} \sma Y\angordmp\big)\big]$}
(a22)++(0,\w) node[align=left] (a32) {$\txint^{\angordn} (S^U)^{\sma\angordn} \sma$\\
$\phantom{\txint^{\angordn}} \big[\txint^{\angordm, \angordmp} \txwedge (X \angordm \sma Y\angordmp)\big]$}
(a31)++(0,\x) node[align=left] (a41) {$\txint^{V,W,\angordm,\angordmp} \IU(V\oplus W, U)_\splus \sma$\\
$\phantom{\txint} \big[((S^V)^{\sma\angordm} \sma (S^W)^{\sma\angordmp}) \sma$\\
$\phantom{\txint\big[} (X\angordm \sma Y\angordmp)\big]$}
(a32)++(0,\x) node[align=left] (a42) {$\txint^{\angordn, \angordm, \angordmp} \txwedge \big[(S^U)^{\sma\angordn} \sma$\\
$\phantom{\txint^{\angordn, \angordm, \angordmp} \txwedge } (X\angordm \sma Y\angordmp) \big]$}
(a41)++(0,\x) node[align=left] (a51) {$\big[\IU(V\oplus W, U)_\splus \sma$\\
$\phantom{\big[} ((S^V)^{\sma\angordm} \sma (S^W)^{\sma\angordmp})\big] \sma$\\
$\phantom{\big[} (X\angordm \sma Y\angordmp)$}
(a42)++(0,\x) node[align=left] (a52) {$(S^U)^{\sma(\angordm \oplus \angordmp)} \sma$\\ 
$(X\angordm \sma Y\angordmp)$}
;
\draw[1cell=.7]
(a11) edge node {\ktwo_U} (a12)
(a11) edge[equal] node[swap] {(\mathrm{u})} (a21)
(a21) edge[equal] node[swap] {(\mathrm{k})} (a31)
(a41) edge node[pos=.4] {\iso} (a31)
(a51) edge node[pos=.4] {\inc} (a41)
(a51) edge node {\smas \sma 1} (a52)
(a52) edge node[swap] {\inc'} (a42)
(a42) edge node[swap] {\iso} (a32)
(a32) edge[equal] node[swap] {(\mathrm{g})} (a22)
(a22) edge[equal] node[swap] {(\mathrm{k})} (a12)
;
\end{tikzpicture}
\end{equation}
The details of \cref{ktwo_u} are explained below.
\begin{itemize}
\item The upper left equality labeled $(\mathrm{u})$ uses the definition \cref{smau_obj} of $\smau$, with $(V,W) \in (\IUsk)^2$.  The coend $\txint^{V,W}$ is taken in $\Topst$, with $G$ acting diagonally on representatives \cref{smau_gaction}.
\item The two equalities labeled $(\mathrm{k})$ use the definition \cref{Kgxv} of $\Kg$ on a $\Gskg$-space, with each of $\angordm$, $\angordmp$, and $\angordn$ running through $\Gsk$.
\item The equality labeled $(\mathrm{g})$ uses the definition \cref{smag_ptspace} of $\smag$, with the wedge $\txwedge$ indexed by the set $\Gskpunc(\angordm \oplus \angordmp, \angordn)$ of nonzero morphisms.
\item The two pointed $G$-homeomorphisms labeled $\iso$ use
\begin{itemize}
\item the commutation of $\sma$ with coends and wedges,
\item the Fubini Theorem for coends \cite[1.3.1]{loregian}, and
\item coherence isomorphisms for $(\Gtopst,\sma)$ to switch $X\angordm$ and $(S^W)^{\sma\angordmp}$.
\end{itemize}
\item The lower-left pointed $G$-morphism $\inc$ is the associativity isomorphism for $(\Gtopst,\sma)$ followed by the structure morphism of the coend in its codomain for each quadruple of objects 
\begin{equation}\label{vwmm}
\big((V,W) \in (\IUsk)^2 ; (\angordm, \angordmp) \in \Gsk^2\big).
\end{equation}
\item The lower-right pointed $G$-morphism $\inc'$ is the structure morphism of the coend in its codomain for the objects
\[\big(\angordn = \angordm \oplus \angordmp ; \angordm ; \angordmp \big) \in \Gsk^3,\]
where $(\angordm, \angordmp)$ is the pair in \cref{vwmm}, and the identity morphism
\begin{equation}\label{one_mm}
1_{\angordm \oplus \angordmp} \in \Gskpunc(\angordm \oplus \angordmp, \angordm \oplus \angordmp)
\end{equation}
for the wedge index.
\item In the bottom horizontal arrow, $1$ denotes the identity morphism of the pointed $G$-space $X\angordm \sma Y\angordmp$.  The pointed $G$-morphism $\smas$ is defined as the following composite.
\begin{equation}\label{smash_s}
\begin{tikzpicture}[vcenter]
\def\v{-1.4} \def\u{-1}
\draw[0cell=.8]
(0,0) node (a11) {\IU(V\oplus W, U)_\splus \sma [(S^V)^{\sma\angordm} \sma (S^W)^{\sma\angordmp}]}
(a11)++(4,.4*\v) node (a12) {(S^U)^{\sma(\angordm \oplus \angordmp)}}
(a11)++(0,\v) node (a21) {\IU(V\oplus W, U)_\splus \sma \Topgst\big((\sma\angordm) \sma (\sma\angordmp), S^V \sma S^W \big)}
(a21)++(0,\u) node (a31) {\IU(V\oplus W, U)_\splus \sma \Topgst(\sma(\angordm \oplus \angordmp), S^{V \oplus W})}
;
\draw[1cell=.8]
(a11) [rounded corners=2pt] -| node[pos=.28] {\smas} (a12)
;
\draw[1cell=.8]
(a11) edge node[swap] {1 \sma \sma} (a21)
(a21) edge node[swap] {1\, \smam \iso} (a31)
(a31) [rounded corners=2pt] -| node[pos=.2] {\mcomp} (a12)
;
\end{tikzpicture}
\end{equation}
The arrows in \cref{smash_s} are defined as follows.
\begin{itemize}
\item Each copy of $1$ denotes the identity morphism of the pointed $G$-space $\IU(V\oplus W, U)_\splus$ \pcref{def:IU_spaces}.
\item In the upper-left arrow, the pointed $G$-morphism
\[(S^V)^{\sma\angordm} \sma (S^W)^{\sma\angordmp} \fto{\sma} \Topgst\big((\sma\angordm) \sma (\sma\angordmp), S^V \sma S^W \big)\]
is defined by smashing a pointed morphism $\sma\angordm \to S^V$ with a pointed morphism $\sma\angordmp \to S^W$.
\item In the lower-left arrow, the pointed $G$-homeomorphism denoted by $\iso$ is given by post-composing with the pointed $G$-homeomorphism \cref{gsphere_multiplication}
\[S^V \sma S^W \fto[\iso]{\mu_{V,W}} S^{V \oplus W}\]
and using the object equality \pcref{sma_symmon}
\[(\sma\angordm) \sma (\sma\angordmp) = \sma(\angordm \oplus \angordmp) \inspace \Fsk\]
for the domain.
\item The pointed $G$-morphism $\mcomp$ is given by composing pointed morphisms
\[\sma(\angordm \oplus \angordmp) \to S^{V \oplus W} \to S^U,\]
using \cref{gsphere_mor} to first send a linear isometric isomorphism $V \oplus W \fiso U$ to a pointed $G$-homeomorphism $S^{V \oplus W} \fiso S^U$.
\end{itemize}
\end{itemize}
This finishes the definition of the pointed $G$-morphism $\ktwo_U$ in \cref{ktwo_u}, hence also $\ktwo$ and $\Kgtwo_{X,Y}$ in \cref{smau_kgx}.  \cref{ktwo_u_welldef} proves that $\ktwo$ is a $G$-equivariant $\IU$-morphism.  \cref{kgtwo_xy_welldef} proves that $\Kgtwo_{X,Y}$ is a $G$-equivariant $\IU$-morphism.  \cref{kgtwo_welldef} proves that $\Kgtwo$ is a $\Gtop$-natural transformation.
\end{definition}

\subsection*{Unraveling $\ktwo$}
Before we prove that the monoidal constraint $\Kgtwo$ in \cref{def:Kgtwo} is well defined, we unravel $\ktwo_U$ in \cref{ktwo_u}, starting with its domain and codomain.

\begin{explanation}[Domain of $\ktwo_U$]\label{expl:ktwo_u_dom}
Using the left column of the diagram \cref{ktwo_u}, the domain of $\ktwo_U$ is given by the coend
\[\begin{split}
& (\Kg X \smau \Kg Y)_U\\
&\iso \scalebox{.8}{$\dint^{V,W, \angordm,\angordmp} \hsp{-.5ex} \IU(V \oplus W,U)_\splus \sma \big[((S^V)^{\sma\angordm} \sma (S^W)^{\sma\angordmp}) \sma (X\angordm \sma Y\angordmp) \big]$}
\end{split}\]
taken in $\Topst$, with $(V,W) \in (\IUsk)^2$ and $(\angordm,\angordmp) \in \Gsk^2$.  Each of its points is represented by a quintuple
\begin{equation}\label{kxukyu_rep}
\begin{split}
\bolde = \big(e \in \IU(V \oplus W,U)_\splus ; & \,v \in (S^V)^{\sma\angordm} ; w \in (S^W)^{\sma\angordmp} ;\\
&\, x \in X\angordm ; y \in Y\angordmp \big),
\end{split}
\end{equation}
which is the basepoint if any one of its five entries is the basepoint. 

\parhead{$G$-action}.  By \cref{smau_gaction,g_tx}, the group $G$ acts diagonally on representatives, which means
\begin{equation}\label{ktwo_u_dom_gaction}
g\bolde = \big(ge(\ginv \oplus \ginv) ; gv ; gw ; gx ; gy \big)
\end{equation}
for $g \in G$.

\parhead{Identification}.  In terms of representatives \cref{kxukyu_rep}, there are two kinds of identifications in $(\Kg X \smau \Kg Y)_U$, as explained in \cref{kxukyu_iden_i,kxukyu_iden_ii} below.
\begin{itemize}
\item By \cref{Kgxv_wedge,Kgxv_rep,t_upom_x,Kgxv_dzeroone}, with $(e, x, y)$ as given in \cref{kxukyu_rep}, for each quadruple
\[\big(v \in (S^V)^{\sma\angordl} ; \angordm \fto{\upom} \angordl \in \Gsk ; w \in (S^W)^{\sma\angordlp} ; \angordmp \fto{\upom'} \angordlp \big),\]
the following two quintuples are identified.
\begin{equation}\label{kxukyu_iden_i}
\begin{split}
&\big(e ; v(\sma\upom) \in (S^V)^{\sma\angordm} ; w(\sma\upom') \in (S^W)^{\sma\angordmp} ; x ; y \big)\\
&\big(e ; v ; w ; (X\upom)(x) \in X\angordl ; (Y\upom')(y) \in Y\angordlp \big)
\end{split}
\end{equation}
\item By \cref{smau_wedge,x_smau_y_rep,smau_coend,smau_iden,Kgxf}, with $(v,w,x,y)$ as given in \cref{kxukyu_rep}, for each triple of isomorphisms in $\IUsk$
\[\big(V \fto{f} V' ; W \fto{h} W' ; V' \oplus W' \fto{e} U\big),\]
the following two quintuples are identified.
\begin{equation}\label{kxukyu_iden_ii}
\begin{split}
&\big(e(f \oplus h) \in \IU(V \oplus W, U) ; v; w; x; y \big) \\
&\big(e ; fv \in (S^{V'})^{\sma\angordm} ; hw \in (S^{W'})^{\sma\angordmp} ; x; y \big)
\end{split}
\end{equation}
\end{itemize}
This finishes the description of the domain of $\ktwo_U$ in \cref{ktwo_u}.
\end{explanation}

\begin{explanation}[Codomain of $\ktwo_U$]\label{expl:ktwo_u_cod}
Using the right column of the diagram \cref{ktwo_u}, the codomain of $\ktwo_U$ is given by the coend
\[\big(\Kg(X \smag Y) \big)_U
\iso \scalebox{.9}{$\dint^{\angordn, \angordm,\angordmp} \hsp{-2em} \bigvee_{\Gskpunc(\angordm \oplus \angordmp, \angordn)} \hsp{-2em} [(S^U)^{\sma\angordn} \sma (X\angordm \sma Y\angordmp) ]$}\]
taken in $\Topst$, with $(\angordn, \angordm,\angordmp) \in \Gsk^3$.  Each point in $(\Kg(X \smag Y))_U$ is represented by a quadruple
\begin{equation}\label{kxyu_rep}
\begin{split}
\bu = \big(u \in  (S^U)^{\sma\angordn} ; &\, \angordm \oplus \angordmp \fto{\upom} \angordn \in \Gsk;\\
&\, x \in X\angordm ; y \in Y\angordmp\big).
\end{split}
\end{equation}
It is the basepoint if either $\upom$ is the 0-morphism in $\Gsk$, or if any one of $u$, $x$, or $y$ is the basepoint. 

\parhead{$G$-action}.  By \cref{smag_rep_gaction,g_tx}, the group $G$ acts diagonally on representatives and trivially on $\Gsk$, which means
\begin{equation}\label{ktwo_u_cod_gaction}
g\bu = \big(gu ; \upom ; gx; gy \big)
\end{equation}
for $g \in G$.

\parhead{Identification}.  In terms of representatives \cref{kxyu_rep}, there are two kinds of identifications in $(\Kg(X \smag Y))_U$, as explained in \cref{kxyu_iden_i,kxyu_iden_ii} below.
\begin{itemize}
\item By \cref{ggtop_ptday,Kgxv_wedge,Kgxv_rep,t_upom_x,Kgxv_dzeroone}, with $(u,x,y)$ as given in \cref{kxyu_rep}, for each pair of morphisms in $\Gsk$
\[\big(\angordr \fto{\upde} \angordn ; \angordm \oplus \angordmp \fto{\upom} \angordr\big),\]
the following two quadruples are identified.
\begin{equation}\label{kxyu_iden_i}
\begin{split}
& \big(u(\sma\upde) \in (S^U)^{\sma\angordr} ; \upom ; x; y \big) \\
& \big(u ; \angordm \oplus \angordmp \fto{\upde\upom} \angordn ; x; y \big)
\end{split}
\end{equation}
\item By \cref{smag_ptspace_wedge,smag_ptspace_rep,smag_ptspace_dzerone}, with $(u,\upom)$ as given in \cref{kxyu_rep}, for each quadruple
\[\big(\angordl \fto{\upla} \angordm \in \Gsk; \angordlp \fto{\upla'} \angordmp \in \Gsk ; x \in X\angordl ; y \in Y\angordlp \big),\]
the following two quadruples are identified.
\begin{equation}\label{kxyu_iden_ii}
\begin{split}
&\big(u ; \angordl \oplus \angordlp \fto{\upom(\upla \oplus \upla')} \angordn ; x; y \big) \\
&\big(u ; \upom ; (X\upla)(x) \in X\angordm ; (Y\upla')(y) \in Y\angordmp\big)
\end{split}
\end{equation}
\end{itemize}
This finishes the description of the codomain of $\ktwo_U$ in \cref{ktwo_u}.
\end{explanation}

\begin{explanation}[Unraveling $\ktwo_U$]\label{expl:ktwo_u}
Using \cref{smash_s,one_mm,kxukyu_rep,kxyu_rep}, the pointed morphism \cref{ktwo_u}
\[(\Kg X \smau \Kg Y)_U \fto{\ktwo_U} \big(\Kg(X \smag Y)\big)_U\]
is defined on representatives by
\begin{equation}\label{ktwo_u_e}
\ktwo_U(e; v; w; x; y) = 
\big(e (\mu_{V,W}) (v \sma w) ;  1_{\angordm \oplus \angordmp} ; x; y\big),
\end{equation}
in which the first entry is the following composite.
\[\begin{tikzpicture}[vcenter]
\def\v{-1.4} 
\draw[0cell]
(0,0) node (a1) {\sma(\angordm \oplus \angordmp)}
(a1)++(0,\v) node (a2) {(\sma\angordm) \sma (\sma\angordmp)}
(a2)++(3.8,0) node (a3) {S^V \sma S^W}
(a3)++(2.8,0) node (a4) {S^{V \oplus W}}
(a4)++(0,-\v) node (a5) {S^U}
;
\draw[1cell=.9]
(a1) edge[equal] (a2)
(a2) edge node {v \sma w} (a3)
(a3) edge node {\mu_{V,W}} node[swap] {\iso} (a4)
(a4) edge node[swap] {e} node {\iso} (a5)
;
\end{tikzpicture}\]
The pointed $G$-homeomorphism $\mu_{V,W}$ is defined in \cref{gsphere_multiplication}.  The previous composite sends a pair of points
\[(t \in \sma\angordm ; t' \in \sma\angordmp)\]
to the point
\begin{equation}\label{ktwo_points}
e\big( v(t) \oplus w(t')\big)
\end{equation}
in $S^U$.
\end{explanation}

\subsection*{Proofs}
The rest of this section proves that $\Kgtwo$ is a well-defined $\Gtop$-natural transformation.  The following result proves that $\ktwo$ is well defined.

\begin{lemma}\label{ktwo_u_welldef}
In \cref{smau_kgx}, 
\[\Kg X \smau \Kg Y \fto{\ktwo} \Kg(X \smag Y)\]
is a $G$-equivariant $\IU$-morphism.
\end{lemma}

\begin{proof}
We first prove that, for each object $U \in \IUsk$, the $U$-component $\ktwo_U$ in \cref{ktwo_u} is a pointed $G$-morphism.

\parhead{Pointed morphism}.  To show that $\ktwo_U$ is a well-defined pointed morphism, we need to show that its definition \cref{ktwo_u_e} is independent of the choice of a representative in its domain.  Since the domain of $\ktwo_U$ has two kinds of identifications, \cref{kxukyu_iden_i,kxukyu_iden_ii}, we need to consider each case separately.  

The following equalities prove that the two quintuples in \cref{kxukyu_iden_i} are sent by $\ktwo_U$ to the same point in $(\Kg(X \smag Y))_U$.  We abbreviate $\mu_{V,W}$ to $\mu$.
\[\begin{aligned}
& \ktwo_U \big(e ; v(\sma\upom) ; w(\sma\upom') ; x ; y \big) && \\
&= \big(e\mu [v(\sma\upom) \sma w(\sma\upom')] ; 1_{\angordm \oplus \angordmp} ; x; y \big) && \text{by \cref{ktwo_u_e}}\\
&= \big(e\mu (v \sma w) [(\sma\upom) \sma (\sma\upom')] ; 1_{\angordm \oplus \angordmp} ; x; y \big) && \text{\scalebox{.85}{by functoriality of $\sma$}} \\
&= \big(e\mu (v \sma w) [\sma(\upom \oplus \upom')] ; 1_{\angordm \oplus \angordmp} ; x; y \big) && \text{\scalebox{.85}{by \cref{sma_symmon}}} \\
&= \big(e\mu (v \sma w) ; \upom \oplus \upom' ; x; y \big) && \text{by \cref{kxyu_iden_i}} \\
&= \big(e\mu (v \sma w) ; 1_{\angordl \oplus \angordlp} ; (X\upom)(x); (Y\upom')(y) \big) && \text{by \cref{kxyu_iden_ii}} \\
&= \ktwo_U \big(e; v; w; (X\upom)(x); (Y\upom')(y) \big) && \text{by \cref{ktwo_u_e}}
\end{aligned}\]
The following equalities prove that the two quintuples in \cref{kxukyu_iden_ii} are sent by $\ktwo_U$ to the same point in $(\Kg(X \smag Y))_U$.  We abbreviate $1_{\angordm \oplus \angordmp}$ to $1$.
\[\begin{aligned}
& \ktwo_U \big(e(f \oplus h); v; w; x; y \big) \\
&= \big(e(f \oplus h) (\mu_{V,W}) (v \sma w) ; 1; x; y \big) && \text{by \cref{ktwo_u_e}} \\
&= \big(e (\mu_{V',W'}) (f \sma h) (v \sma w) ; 1; x; y \big) && \text{by \cref{gsp_mult}} \\
&= \big(e (\mu_{V',W'}) (fv \sma hw) ; 1; x; y \big) && \text{by functoriality of $\sma$} \\
&= \ktwo_U(e; fv; hw; x; y) && \text{by \cref{ktwo_u_e}} 
\end{aligned}\]
This proves that $\ktwo_U$ is a pointed morphism.

\parhead{$G$-equivariance}.
The following computation for $g \in G$ proves that $\ktwo_U$ is $G$-equivariant.   We abbreviate $\mu_{V,W}$ to $\mu$ and $1_{\angordm \oplus \angordmp}$ to $1$.
\begin{equation}\label{ktwou_geq}
\begin{aligned}
& \ktwo_U \big(g(e; v; w; x; y)\big) && \\
&= \ktwo_U \big(ge(\ginv \oplus \ginv); gv; gw; gx; gy \big) && \text{by \cref{ktwo_u_dom_gaction}} \\
&= \big(ge(\ginv \oplus \ginv) \mu (gv \sma gw) ; 1; gx; gy \big) && \text{by \cref{ktwo_u_e}} \\
&= \big(ge\mu (\ginv \sma \ginv) (gv \sma gw) ; 1 ; gx; gy \big) && \text{by \cref{gsp_mult}} \\
&= \big(ge\mu (v \sma w) ; 1 ; gx; gy \big) && \text{\scalebox{.8}{by functoriality of $\sma$}} \\
&= g\big(e\mu (v \sma w) ; 1 ; x; y  \big) && \text{by \cref{ktwo_u_cod_gaction}} \\
&= g \big(\ktwo_U (e; v; w; x; y)\big) && \text{by \cref{ktwo_u_e}} 
\end{aligned}
\end{equation}
This proves that, for each object $U \in \IUsk$, $\ktwo_U$ is a pointed $G$-morphism.

\parhead{Naturality}.  
The naturality property \cref{iu_mor_natural} for $\ktwo$ states that, for each linear isometric isomorphism $f \cn U \fiso U'$ in $\IUsk$, the following diagram of pointed morphisms commutes.
\[\begin{tikzpicture}[vcenter]
\def\v{-1.4}
\draw[0cell]
(0,0) node (a11) {(\Kg X \smau \Kg Y)_U}
(a11)++(4.5,0) node (a12) {\big(\Kg(X \smag Y)\big)_U}
(a11)++(0,\v) node (a21) {(\Kg X \smau \Kg Y)_{U'}}
(a12)++(0,\v) node (a22) {\big(\Kg(X \smag Y)\big)_{U'}}
;
\draw[1cell=.8]
(a11) edge node {\ktwo_U} (a12)
(a12) edge[transform canvas={xshift=-2.5em}] node {(\Kg(X \smag Y))_f} (a22)
(a11) edge[transform canvas={xshift=2.5em}] node[swap] {(\Kg X \smau \Kg Y)_f} (a21)
(a21) edge node {\ktwo_{U'}} (a22)
;
\end{tikzpicture}\]
This diagram commutes because, by \cref{smau_mor,Kgxf_rep,ktwo_u_e}, each of the two composites sends a representative quintuple $\bolde$ in \cref{kxukyu_rep} to the representative quadruple
\[\big(fe(\mu_{V,W}) (v \sma w); 1_{\angordm \oplus \angordmp}; x; y\big)\]
in $(\Kg(X \smag Y))_{U'}$.  By \cref{expl:eqiu_morphism}, we have proved that $\ktwo$ is a $G$-equivariant $\IU$-morphism.
\end{proof}

The following result is the first step in showing that the components of $\Kgtwo$ are well defined.

\begin{lemma}\label{kgtwo_xy_welldef}
In \cref{smau_kgx}, 
\[\Kg X \smasg \Kg Y \fto{\Kgtwo_{X,Y}} \Kg(X \smag Y)\]
is a $G$-equivariant $\IU$-morphism.
\end{lemma}

\begin{proof}
The top row of the diagram \cref{smau_kgx} is a coequalizer.  To show that $\Kgtwo_{X,Y}$ is a well-defined $\IU$-morphism, we need to show that the following two composites are equal as $\IU$-morphisms.
\begin{equation}\label{ktwo_coequal}
\begin{tikzpicture}[vcenter]
\draw[0cell=.9]
(0,0) node (a1) {\phantom{(\Kg X \smau \gsp) \smau \Kg Y}}
(a1)++(0,-.04) node (a1') {(\Kg X \smau \gsp) \smau \Kg Y}
(a1)++(4.5,0) node (a2) {\Kg X \smau \Kg Y}
(a2)++(0,-1.4) node (b) {\Kg(X \smag Y)}
;
\draw[1cell=.8]
(a1) edge[transform canvas={yshift=.4ex}] node {\umu \smau 1} (a2)
(a1) edge[transform canvas={yshift=-.5ex}] node[swap] {(1 \smau \umu') \upbe} (a2)
(a2) edge node[swap] {\ktwo} (b)
;
\end{tikzpicture}
\end{equation}
Once this assertion is proved, the universal property of coequalizers implies the unique existence of the $\IU$-morphism $\Kgtwo_{X,Y}$ such that the triangle in the diagram \cref{smau_kgx} commutes.  Moreover, the $G$-equivariance of $\Kgtwo_{X,Y}$ follows from the universal property of coequalizers and the $G$-equivariance of $\ktwo$, which is proved in \cref{ktwou_geq} above.

To show that the two composites in \cref{ktwo_coequal} are equal at each object $U \in \IUsk$, we first use \cref{au_xyzu_dom,Kgxv} to express the $U$-component of its domain as the following coend in $\Topst$.
\begin{equation}\label{kxsky_u_coend}
\begin{split}
& [(\Kg X \smau \gsp) \smau \Kg Y]_U \\
&\iso \int^{(V_1,V_2,V_3) \in (\IUsk)^3,\, (\angordm, \angordmp) \in {\Gsk^2}} \IU\big((V_1 \oplus V_2) \oplus V_3, U\big)_\splus \sma\\
&\phantom{\iso\iso} \big[\big( (S^{V_1})^{\sma\angordm} \sma S^{V_2} \sma (S^{V_3})^{\sma\angordmp}\big)
\sma (X\angordm \sma Y\angordmp) \big]
\end{split}
\end{equation}
In terms of the coend \cref{kxsky_u_coend}, each point of $[(\Kg X \smau \gsp) \smau \Kg Y]_U$ is represented by a sextuple
\begin{equation}\label{kxsky_u_rep}
\begin{split}
\bolde &= \big( e \in \IU((V_1 \oplus V_2) \oplus V_3, U)_\splus ; v_1 \in (S^{V_1})^{\sma\angordm} ; \\
&\phantom{==} a \in S^{V_2} ; v_3 \in (S^{V_3})^{\sma\angordmp} ; x \in X\angordm ; y \in Y\angordmp\big).
\end{split}
\end{equation}
It represents the basepoint if any one of its six entries is the basepoint. 

The following equalities in $(\Kg(X\smag Y))_U$ prove that the two composites in \cref{ktwo_coequal} are equal on each representative $\bolde$ in \cref{kxsky_u_rep}.  The identity morphism $1_{\angordm \oplus \angordmp}$ is abbreviated to $1$.
\[\begin{split}
& (\ktwo_U) (\umu \smau 1)_U (\bolde) \\
&= (\ktwo_U) (e; v_1 \oplus a; v_3; x; y) \\
&= \big(e(\mu_{V_1 \oplus V_2, V_3}) [(v_1 \oplus a) \sma v_3] ; 1 ; x; y \big) \\
&= \big(e\al^{-1} (1_{V_1} \oplus \xi_{V_3,V_2}) (\mu_{V_1, V_3 \oplus V_2}) [v_1 \sma (v_3 \oplus a)] ; 1; x; y \big) \\
&= (\ktwo_U) \big(e\al^{-1} (1_{V_1} \oplus \xi_{V_3,V_2}) ; v_1 ; v_3 \oplus a ; x; y \big) \\
&= (\ktwo_U) \big((1 \smau \umu') \upbe\big)_U (\bolde) 
\end{split}\]
In the previous computation, the first and last equalities follow from
\begin{itemize}
\item the definitions \cref{dzeroone_bv} of $(\umu \smau 1)_U$ and $\big((1 \smau \umu') \upbe\big)_U$, and
\item the definition \cref{Kgxv_action_rep} of the right $\gsp$-actions on $\Kg X$ and $\Kg Y$.  
\end{itemize}
The second and fourth equalities follow from the definition of $\ktwo_U$ \cref{ktwo_u_e}.  The third equality follows from the following commutative diagram in $\Topst$.
\begin{equation}\label{kgtwo_xy_diag}
\begin{tikzpicture}[vcenter]
\def\v{-1.4} \def\u{.7} \def\h{3}
\draw[0cell=.9]
(0,0) node (a11) {(\sma\angordm) \sma (\sma\angordmp)}
(a11)++(\h,\u) node (a12) {S^{V_1 \oplus V_2} \sma S^{V_3}}
(a12)++(\h,-\u) node (a13) {S^{(V_1 \oplus V_2) \oplus V_3}}
(a11)++(0,\v) node (a21) {S^{V_1} \sma S^{V_3 \oplus V_2}}
(a21)++(\h,-\u) node (a22) {S^{V_1 \oplus (V_3 \oplus V_2)}}
(a22)++(\h,\u) node (a23) {S^{V_1 \oplus (V_2 \oplus V_3)}}
;
\draw[1cell=.8]
(a11) edge node[swap] {v_1 \sma (v_3 \oplus a)} (a21)
(a23) edge node[swap] {\al^{-1}} node[pos=.43] {\iso} (a13)
;
\draw[1cell=.8]
(a11) [rounded corners=2pt] |- node[pos=.75] {(v_1 \oplus a) \sma v_3} (a12);
\draw[1cell=.8]
(a12) [rounded corners=2pt] -| node[pos=.2] {\mu_{V_1 \oplus V_2, V_3}} (a13);
\draw[1cell=.8]
(a21) [rounded corners=2pt] |- node[pos=.8] {\mu_{V_1, V_3 \oplus V_2}} (a22);
\draw[1cell=.8]
(a22) [rounded corners=2pt] -| node[pos=.2] {1_{V_1} \oplus \xi_{V_3,V_2}} (a23);
\end{tikzpicture}
\end{equation}
More explicitly, the diagram \cref{kgtwo_xy_diag} commutes because each of the two composites sends a pair of points 
\[\big(t \in \sma\angordm ; t' \in \sma\angordmp \big)\]
to the point
\[\big( [v_1(t) \oplus a] \oplus v_3(t') \big) \in S^{(V_1 \oplus V_2) \oplus V_3}.\]
This finishes the proof that the two composites in \cref{ktwo_coequal} are equal.
\end{proof}

The last observation of this section proves that $\Kgtwo$ is well defined.

\begin{lemma}\label{kgtwo_welldef}
In \cref{def:Kgtwo}, $\Kgtwo$ is a $\Gtop$-natural transformation.
\end{lemma}

\begin{proof}
We first prove that the components of $\Kgtwo$ are $G$-equivariant $\gsp$-module morphisms.

\parhead{Well-defined components}.  \cref{kgtwo_xy_welldef} proves that, for each pair $(X,Y)$ of $\Gskg$-spaces \cref{ggtop_obj}, the $(X,Y)$-component \cref{Kgtwo_component}
\[\Kg X \smasg \Kg Y \fto{\Kgtwo_{X,Y}} \Kg(X \smag Y)\]
is a $G$-equivariant $\IU$-morphism.  To show that it is an $\gsp$-module morphism \cref{gsp_mor_axiom}, we need to show that, for each pair of objects $(U_1,U_2) \in (\IUsk)^2$, the following diagram of pointed morphisms commutes.
\begin{equation}\label{kgtwo_gsp_action}
\begin{tikzpicture}[vcenter]
\def\v{-1.4}
\draw[0cell=.9]
(0,0) node (a11) {(\Kg X \smasg \Kg Y)_{U_1} \sma S^{U_2}}
(a11)++(4.7,0) node (a12) {(\Kg X \smasg \Kg Y)_{U_1 \oplus U_2}}
(a11)++(0,\v) node (a21) {(\Kg(X \smag Y))_{U_1} \sma S^{U_2}}
(a12)++(0,\v) node (a22) {(\Kg(X \smag Y))_{U_1 \oplus U_2}}
;
\draw[1cell=.8]
(a11) edge node {\umu_{U_1,U_2}} (a12)
(a12) edge node {\Kgtwo_{X,Y; U_1 \oplus U_2}} (a22)
(a11) edge node[swap] {\Kgtwo_{X,Y; U_1} \sma 1} (a21)
(a21) edge node {\umu_{U_1,U_2}} (a22)
;
\end{tikzpicture}
\end{equation} 

Since $\Kg X \smasg \Kg Y$ is a coequalizer \cref{smau_kgx}, its points are represented by those in $\Kg X \smau \Kg Y$ \cref{kxukyu_rep}.  Thus, each point in $(\Kg X \smasg \Kg Y)_{U_1}$ is represented by a quintuple
\begin{equation}\label{kxkyuone_rep}
\begin{split}
\bolde = \big(e \in \IU(V \oplus W,U_1)_\splus ; & \,v \in (S^V)^{\sma\angordm} ; w \in (S^W)^{\sma\angordmp} ;\\
&\, x \in X\angordm ; y \in Y\angordmp \big).
\end{split}
\end{equation}
The following equalities, for $a \in S^{U_2}$, prove that the diagram \cref{kgtwo_gsp_action} commutes.  We abbreviate $1_{\angordm \oplus \angordmp}$ to 1 and $\umu_{U_1,U_2}$ to $\umu$.
\[\begin{split}
& \umu (\Kgtwo_{X,Y; U_1} \sma 1) (\bolde; a) \\
&= \umu \big(e(\mu_{V,W}) (v \sma w) ; 1; x; y; a \big) \\
&= \big([e(\mu_{V,W}) (v \sma w)] \oplus a; 1; x; y \big) \\
&= \big((e \oplus 1_{U_2}) (\al^{-1}) (\mu_{V,W \oplus U_2}) [v \sma (w \oplus a)] ; 1; x; y \big) \\
&= (\Kgtwo_{X,Y; U_1 \oplus U_2}) \big((e \oplus 1_{U_2}) \al^{-1} ; v; w \oplus a ; x; y \big) \\
&= (\Kgtwo_{X,Y; U_1 \oplus U_2}) (\umu)(\bolde; a)
\end{split}\]
The equalities above hold for the following reasons.
\begin{itemize}
\item The first and fourth equalities follow from the definition \cref{smau_kgx} of $\Kgtwo_{X,Y}$ in terms of $\ktwo$ \cref{ktwo_u_e}. 
\item The second equality holds by the definition \cref{Kgxv_action_rep} of the right $\gsp$-action on $\Kg (X \smag Y)$.
\item The last equality holds by the definitions of the right $\gsp$-actions on $\Kg X \smasg \Kg Y$ \cref{smasg_action_y} and $\Kg Y$ \cref{Kgxv_action_rep}.
\item The third equality follows from the following commutative diagram in $\Topst$.
\[\begin{tikzpicture}[vcenter]
\def\h{4.5} \def\v{-1.4}
\draw[0cell=.9]
(0,0) node (a11) {(\sma\angordm) \sma (\sma\angordmp)}
(a11)++(\h,0) node (a12) {S^V \sma S^{W \oplus U_2}}
(a12)++(\h/3,\v/2) node (a13) {S^{V \oplus (W \oplus U_2)}}
(a11)++(0,\v) node (a21) {S^V \sma S^W}
(a13)++(0,\v) node (a23) {S^{(V \oplus W) \oplus U_2}}
(a21)++(0,\v) node (a31) {S^{V \oplus W}}
(a31)++(\h/2,0) node (a32) {S^{U_1}}
(a32)++(\h/2,0) node (a33) {S^{U_1 \oplus U_2}}
;
\draw[1cell=.8]
(a11) edge node[swap] {v \sma w} (a21)
(a21) edge node[swap] {\mu_{V,W}} (a31)
(a31) edge node {e} (a32)
(a32) edge node {-\oplus a} (a33)
(a11) edge node {v \sma (w \oplus a)} (a12)
(a13) edge node {\al^{-1}} node[swap] {\iso} (a23)
(a12) [rounded corners=2pt] -| node[pos=.7] {\mu_{V,W \oplus U_2}} (a13)
;
\draw[1cell=.8]
(a23) [rounded corners=2pt] |- node[pos=.3] {e \oplus 1_{U_2}} (a33)

;
\end{tikzpicture}\]
More explicitly, each of the two composites in the previous diagram sends a pair of points
\[(t \in \sma\angordm; t' \in \sma\angordmp)\]
to the point
\[[e(v(t) \oplus w(t'))] \oplus a \in S^{U_1 \oplus U_2}.\]
\end{itemize}
This proves that each component of $\Kgtwo$ is a $G$-equivariant $\gsp$-module morphism.

\parhead{Naturality}.  The $\Gtop$-naturality \cref{enr_naturality} of $\Kgtwo$ means that, for each pair of points \cref{ggtop_gtop_enr}
\[\big(\theta^X \in \brkstg{X}{X'}\ang{} ; \theta^Y \in \brkstg{Y}{Y'}\ang{} \big),\]
the following diagram of $\gsp$-module morphisms commutes.
\begin{equation}\label{kgtwo_gtop_nat}
\begin{tikzpicture}[vcenter]
\def\v{-1.4}
\draw[0cell=.9]
(0,0) node (a11) {\Kg X \smasg \Kg Y}
(a11)++(3.5,0) node (a12) {\phantom{\Kg (X \smag Y)}}
(a12)++(0,-.04) node (a12') {\Kg (X \smag Y)}
(a11)++(0,\v) node (a21) {\Kg X' \smasg \Kg Y'}
(a12)++(0,\v) node (a22) {\phantom{\Kg (X' \smag Y')}}
(a22)++(0,-.04) node (a22') {\Kg (X' \smag Y')}
;
\draw[1cell=.8]
(a11) edge node {\Kgtwo_{X,Y}} (a12)
(a12') edge[transform canvas={xshift=-2em}, shorten <=-.5ex] node {\Kg(\theta^X \smag \theta^Y)} (a22')
(a11) edge[transform canvas={xshift=2em}] node[swap] {\Kg \theta^X \smasg \Kg \theta^Y} (a21)
(a21) edge node {\Kgtwo_{X',Y'}} (a22)
;
\end{tikzpicture}
\end{equation}
Since an $\gsp$-module morphism is determined by its components \cref{gsp_mor_comp}, it suffices to prove the commutativity of the diagram \cref{kgtwo_gtop_nat} at each representative quintuple $\bolde$ in \cref{kxkyuone_rep}.  The following equalities prove that the two composites in \cref{kgtwo_gtop_nat} are equal on $\bolde$.  We abbreviate $\mu_{V,W}$ to $\mu$ and $1_{\angordm \oplus \angordmp}$ to 1, and omit the subscript $U_1$.
\[\begin{split}
& (\Kgtwo_{X',Y'}) (\Kg \theta^X \smasg \Kg \theta^Y) (\bolde) \\
&= (\Kgtwo_{X',Y'}) \big(e; v; w; \theta^X_{\angordm}(x) ; \theta^Y_{\angordmp}(y)\big) \\
&= \big(e\mu(v \sma w) ; 1; \theta^X_{\angordm}(x) ; \theta^Y_{\angordmp}(y) \big) \\
&= \big(\Kg(\theta^X \smag \theta^Y)\big) \big(e\mu(v \sma w); 1; x; y\big) \\
&= \big(\Kg(\theta^X \smag \theta^Y)\big) (\Kgtwo_{X,Y}) (\bolde)
\end{split}\]
The equalities above hold for the following reasons.
\begin{itemize}
\item The first equality follows from \cref{tha_smau_thap_uz,smasg_mor_coequal,Kg_theta_vtx}. 
\item The second and last equalities follow from \cref{smau_kgx,ktwo_u_e}. 
\item The third equality follows from \cref{ggtop_moncomp_morn,Kg_theta_vtx}.
\end{itemize}
This finishes the proof that $\Kgtwo$ is a $\Gtop$-natural transformation.
\end{proof}

\section{Symmetric Monoidal $\Gtop$-Functoriality}
\label{sec:Kg}

We remind the reader that $G$ is a compact Lie group in this chapter.  The main result of this section constructs the last step of our multifunctorial $G$-equivariant algebraic $K$-theory.  This result states that the triple $(\Kg,\Kgtwo,\Kgzero)$ consisting of
\begin{itemize}
\item the $\Gtop$-functor \pcref{def:ggspace_gspectra,def:ggtop_gsp_mor}
\[\GGTop \fto{\Kg} \GSp\]
and
\item the $\Gtop$-natural transformations \pcref{def:Kg_zero,def:Kgtwo}
\[\begin{tikzpicture}[vcenter]
\def\v{-1.4} \def\h{2.5}
\draw[0cell=.9]
(0,0) node (a11) {(\GGTop)^{\otimes 2}}
(a11)++(3.2,0) node (a12) {(\GSp)^{\otimes 2}}
(a11)++(0,\v) node (a21) {\GGTop}
(a12)++(0,\v) node (a22) {\GSp}
;
\draw[1cell=.8]
(a11) edge node {(\Kg)^{\otimes 2}} (a12)
(a12) edge node {\smasg} (a22)
(a11) edge node[swap] {\smag} (a21)
(a21) edge node {\Kg} (a22)
;
\draw[2cell]
node[between=a11 and a22 at .55, rotate=225, shift={(0,.1*\v)}, 2labelw={below,\Kgtwo,-1pt}] {\Rightarrow}
;
\begin{scope}[shift={(4.8,0)}]
\draw[0cell=.9]
(0,0) node (a1) {\vtensorunit}
(a1)++(.5*\h,-1) node (a2) {\GGTop}
(a1)++(\h,0) node (a3) {\GSp}
;
\draw[1cell=.8]
(a1) edge[transform canvas={yshift=.3ex}] node {\gsp} (a3)
(a1) [rounded corners=2pt] |- node[swap,pos=.25] {\gu} (a2)
;
\draw[1cell=.8]
(a2) [rounded corners=2pt] -| node[swap,pos=.75] {\Kg} (a3)
;
\draw[2cell]
node[between=a1 and a3 at .4, shift={(0,-.45)}, rotate=-90, 2label={above,\Kgzero}, 2label={below,\iso}] {\Rightarrow}
;
\end{scope}
\end{tikzpicture}\]
\end{itemize}
is a unital symmetric monoidal $\Gtop$-functor between the symmetric monoidal $\Gtop$-categories
\begin{itemize}
\item $(\GGTop, \smag, \gu)$ of $\Gskg$-spaces \pcref{expl:ggtop_gtopenr} and
\item $(\GSp,\smasg,\gsp)$ of $\gsp$-modules \pcref{def:gsp_module,def:gsp_smgtop}.
\end{itemize}
Note that $\Kg$ has the $G$-sphere built into its definition: \cref{Kgzero_welldef} proves that the unit constraint \cref{Kgzero_gspmod_mor} 
\[\gsp \fto{\Kgzero} \Kg\gu\]
is a $G$-equivariant $\gsp$-module isomorphism.

\secoutline
\begin{itemize}
\item \cref{Kg_associativity_axiom,Kg_unity_axiom,Kg_braid_axiom} prove, respectively, the associativity axiom, the unity axioms, and the braid axiom for $(\Kg,\Kgtwo,\Kgzero)$.
\item \cref{thm:Kg_smgtop} is the main result of this section, establishing the existence of the symmetric monoidal $\Gtop$-functor $(\Kg,\Kgtwo,\Kgzero)$ and the induced $\Gtop$-multifunctor.
\item \cref{rk:gspace_gspectra} discusses other constructions in the literature that build $G$-spectra from space-level data.
\end{itemize}

\recollection
\cref{Kg_gtop_functor} proves that $\Kg$ is a $\Gtop$-functor.  \cref{Kgzero_welldef} proves that $\Kgzero$ is a $\Gtop$-natural isomorphism.  \cref{kgtwo_welldef} proves that $\Kgtwo$ is a $\Gtop$-natural transformation.  Thus, it remains to verify the axioms of a symmetric monoidal $\Gtop$-functor in \cref{definition:monoidal-V-fun,definition:braided-monoidal-vfunctor}.  To verify these axioms for $(\Kg,\Kgtwo,\Kgzero)$, we first recall the following facts.
\begin{itemize}
\item A $\Gtop$-natural transformation \pcref{def:enriched-natural-transformation} is determined by its components.
\item In each hom $G$-space in $\GSp$ \cref{gsp_gtop_enr}, the points are $\gsp$-module morphisms, which are determined by their underlying $\IU$-morphisms.
\item Each $\IU$-morphism is determined by its component pointed morphisms for objects in $\IUsk$ \cref{theta_iusk}.
\end{itemize}
Thus, it suffices to verify each axiom below for appropriate components.

\begin{lemma}\label{Kg_associativity_axiom}
The triple \pcref{def:ggspace_gspectra,def:ggtop_gsp_mor,def:Kg_zero,def:Kgtwo}
\[\GGTop \fto{(\Kg,\Kgtwo,\Kgzero)} \GSp\]
satisfies the associativity axiom \cref{enrmonfunctor-ass}.
\end{lemma}

\begin{proof}
The associativity axiom means the commutativity of the following diagram of $\IU$-morphisms for each triple $(X,Y,Z)$ of $\Gskg$-spaces \cref{ggtop_obj}, where $\asg$ \cref{asg_xyz} and $\ag$ \cref{ggtop_ag_fff} are the monoidal associators of, respectively, $\GSp$ and $\GGTop$.
\begin{equation}\label{Kg_associativity}
\begin{tikzpicture}[vcenter]
\def\v{-1.4}
\draw[0cell=.9]
(0,0) node (a11) {(\Kg X \smasg \Kg Y) \smasg \Kg Z}
(a11)++(4.7,0) node (a12) {\Kg X \smasg (\Kg Y \smasg \Kg Z)}
(a11)++(0,\v) node (a21) {\Kg(X \smag Y) \smasg \Kg Z}
(a12)++(0,\v) node (a22) {\Kg X \smasg \Kg(Y \smag Z)}
(a21)++(0,\v) node (a31) {\Kg((X \smag Y) \smag Z)}
(a22)++(0,\v) node (a32) {\Kg(X \smag (Y \smag Z))}
;
\draw[1cell=.8]
(a11) edge node {\asg} node[swap] {\iso} (a12)
(a12) edge node {1 \smasg \Kgtwo_{Y,Z}} (a22)
(a22) edge node {\Kgtwo_{X, Y \smag Z}} (a32)
(a11) edge node[swap] {\Kgtwo_{X,Y} \smasg 1} (a21)
(a21) edge node[swap] {\Kgtwo_{X \smag Y,Z}} (a31)
(a31) edge node {\Kg\ag} node[swap] {\iso} (a32)
;
\end{tikzpicture}
\end{equation}
The monoidal associator $\asg$ is induced by the monoidal associator $\au$ \cref{au_xyz} via the colimit in \cref{asg_colimit}.  Moreover, each component of $\Kgtwo$ is induced by $\ktwo$ via the coequalizer in \cref{smau_kgx}.  By the universal property of colimits, to prove the commutativity of \cref{Kg_associativity}, it suffices to prove that the following diagram of $\IU$-morphisms commutes.
\begin{equation}\label{Kg_ass_smau}
\begin{tikzpicture}[vcenter]
\def\v{-1.4}
\draw[0cell=.9]
(0,0) node (a11) {(\Kg X \smau \Kg Y) \smau \Kg Z}
(a11)++(4.7,0) node (a12) {\Kg X \smau (\Kg Y \smau \Kg Z)}
(a11)++(0,\v) node (a21) {\Kg(X \smag Y) \smau \Kg Z}
(a12)++(0,\v) node (a22) {\Kg X \smau \Kg(Y \smag Z)}
(a21)++(0,\v) node (a31) {\Kg((X \smag Y) \smag Z)}
(a22)++(0,\v) node (a32) {\Kg(X \smag (Y \smag Z))}
;
\draw[1cell=.8]
(a11) edge node {\au} node[swap] {\iso} (a12)
(a12) edge node {1 \smau \ktwo} (a22)
(a22) edge node {\ktwo} (a32)
(a11) edge node[swap] {\ktwo \smau 1} (a21)
(a21) edge node[swap] {\ktwo} (a31)
(a31) edge node {\Kg\ag} node[swap] {\iso} (a32)
;
\end{tikzpicture}
\end{equation}
To show that the diagram \cref{Kg_ass_smau} commutes, we first unravel its domain and codomain evaluated at each object $U \in \IUsk$.

\parhead{Domain}.  Using \cref{au_xyzu_dom,Kgxv,ktwo_u}, the $U$-component of the upper-left corner of \cref{Kg_ass_smau} is given by the coend
\begin{equation}\label{kgxyz_u}
\begin{split}
& [(\Kg X \smau \Kg Y) \smau \Kg Z]_U \\
&\iso \int \IU((V_1 \oplus V_2) \oplus V_3, U)_\splus \\
&\phantom{\iso \int } \sma \big[\big( (S^{V_1})^{\sma\angordm} \sma (S^{V_2})^{\sma\angordmp} \big) \sma (S^{V_3})^{\sma\angordmpp} \big]\\
&\phantom{\iso \int} \sma\big[(X\angordm \sma Y\angordmp) \sma Z\angordmpp \big]
\end{split}
\end{equation}
in $\Topst$ indexed by
\[(V_1,V_2,V_3) \in (\IUsk)^3 \andspace (\angordm, \angordmp, \angordmpp) \in \Gsk^3.\]
Each point of the pointed space in \cref{kgxyz_u} is represented by a septuple
\begin{equation}\label{kgxyz_u_rep}
\bolde = (e; v_1; v_2; v_3; x; y; z)
\end{equation}
with
\begin{itemize}
\item $e \in \IU((V_1 \oplus V_2) \oplus V_3, U)_\splus$,
\item $v_1 \in (S^{V_1})^{\sma\angordm}$, $v_2 \in (S^{V_2})^{\sma\angordmp}$, $v_3 \in (S^{V_3})^{\sma\angordmpp}$,
\item $x \in X\angordm$, $y \in Y\angordmp$, and $z \in Z\angordmpp$.
\end{itemize}

\parhead{Codomain}.  Using \cref{ggtop_ag_comp,Kgxv,ktwo_u}, the $U$-component of the lower-right corner of \cref{Kg_ass_smau} is given by the coend
\begin{equation}\label{kgxyz_u_cod}
\begin{split}
& [\Kg(X \smag (Y \smag Z))]_U \\
&\iso \int\bigvee \big[(S^U)^{\sma\angordn} \sma (X\angordm \sma (Y\angordmp \sma Z\angordmpp)) \big]
\end{split}
\end{equation}
in $\Topst$ indexed by 
\[(\angordn, \angordm, \angordmp, \angordmpp) \in \Gsk^4,\]
with the wedge indexed by the set of nonzero morphisms
\[\Gskpunc(\angordm \oplus \angordmp \oplus \angordmpp, \angordn).\]
Each point of the pointed space in \cref{kgxyz_u_cod} is represented by a quintuple
\begin{equation}\label{kgxyz_u_cod_rep}
\bu = (u; \upom, x; y; z)
\end{equation}
with 
\begin{itemize}
\item $u \in (S^U)^{\sma\angordn}$,
\item $\upom \in \Gsk(\angordm \oplus \angordmp \oplus \angordmpp, \angordn)$, 
\item $x \in X\angordm$, $y \in Y\angordmp$, and $z \in Z\angordmpp$.
\end{itemize}

\parhead{The composites}.  By \cref{au_xyz_def,gsphere_multiplication,ggtop_ag_comp,Kg_theta_v,ktwo_u_e}, the $U$-components of the left-bottom and the top-right composites in \cref{Kg_ass_smau} send a representative $\bolde$ in \cref{kgxyz_u_rep} to, respectively, the quintuples
\begin{equation}\label{kg_assoc_rep}
\begin{split}
& \big(e\mu (\mu \sma 1) [(v_1 \sma v_2) \sma v_3] ; 1_{\angordm \oplus \angordmp \oplus \angordmpp} ; x; y; z\big) \andspace\\
& \big(e\al^{-1} \mu (1 \sma \mu) [v_1 \sma (v_2 \sma v_3)] ; 1_{\angordm \oplus \angordmp \oplus \angordmpp} ; x; y; z\big),
\end{split}
\end{equation}
which represent two points in \cref{kgxyz_u_cod} with
\[\angordn = \angordm \oplus \angordmp \oplus \angordmpp \in \Gsk.\]
The two quintuples in \cref{kg_assoc_rep} are equal because the following diagram commutes.
\[\begin{tikzpicture}
\def\h{4.5} \def\u{-1} \def\v{-1.4}
\draw[0cell]
(0,0) node (a11) {\sma(\angordm \oplus \angordmp \oplus \angordmpp)}
(a11)++(0,\u) node (a21) {(\sma\angordm) \sma (\sma\angordmp) \sma (\sma\angordmpp)}
(a21)++(0,\v) node (a31) {(S^{V_1} \sma S^{V_2}) \sma S^{V_3}}
(a31)++(0,\v) node (a41) {(S^{V_1 \oplus V_2}) \sma S^{V_3}}
(a41)++(0,\v) node (a51) {S^{(V_1 \oplus V_2) \oplus V_3}}
(a31)++(\h,0) node (a32) {S^{V_1} \sma (S^{V_2} \sma S^{V_3})}
(a32)++(0,\v) node (a42) {S^{V_1} \sma S^{V_2 \oplus V_3}}
(a42)++(0,\v) node (a52) {S^{V_1 \oplus (V_2 \oplus V_3)}}
;
\draw[1cell=.9]
(a11) edge[equal] (a21)
(a21) edge node[swap] {(v_1 \sma v_2) \sma v_3} (a31)
(a31) edge node[swap] {\mu \sma 1} (a41)
(a41) edge node[swap] {\mu} (a51)
(a32) edge node {1 \sma \mu} (a42)
(a42) edge node {\mu} (a52)
(a52) edge node[swap] {\al^{-1}} node {\iso} (a51)
;
\draw[1cell=.9]
(a21) [rounded corners=2pt] -| node[pos=.25] {v_1 \sma (v_2 \sma v_3)} (a32);
\end{tikzpicture}\]
This finishes the proof of the associativity axiom \cref{enrmonfunctor-ass} for $(\Kg,\Kgtwo,\Kgzero)$.
\end{proof}

\begin{lemma}\label{Kg_unity_axiom}
The triple \pcref{def:ggspace_gspectra,def:ggtop_gsp_mor,def:Kg_zero,def:Kgtwo}
\[\GGTop \fto{(\Kg,\Kgtwo,\Kgzero)} \GSp\]
satisfies the unity axioms \cref{enrmonfunctor-lunity,enrmonfunctor-runity}.
\end{lemma}

\begin{proof}
We prove the right unity axiom \cref{enrmonfunctor-runity}; the proof of the left unity axiom \cref{enrmonfunctor-lunity} is similar.

The right unity axiom asserts the commutativity of the following diagram in $\IUT$ for each $\Gskg$-space $X$ \cref{ggtop_obj}, where $\rsg$ \cref{rsg_x} and $\rg$ \cref{rg_f} are the right monoidal unitors for, respectively, $\GSp$ and $\GGTop$.
\begin{equation}\label{Kg_runity_diagram}
\begin{tikzpicture}[vcenter]
\def\v{-1.4}
\draw[0cell]
(0,0) node (a11) {\Kg X \smasg \gsp}
(a11)++(3.8,0) node (a12) {\Kg X}
(a11)++(0,\v) node (a21) {\Kg X \smasg \Kg\gu}
(a12)++(0,\v) node (a22) {\phantom{\Kg(X \smag \gu)}} 
(a22)++(0,-.04) node (a22') {\Kg(X \smag \gu)}
;
\draw[1cell=.9]
(a11) edge node {\rsg} (a12)
(a11) edge[transform canvas={xshift=1.5em}] node[swap] {1 \smasg \Kgzero} (a21)
(a21) edge node {\Kgtwo_{X,\gu}} (a22)
(a22') edge node[swap] {\Kg\rg} (a12)
;
\end{tikzpicture}
\end{equation}
By the universal property of the coequalizer \cref{gsp_sma_coequal} that defines $\smasg$ in terms of $\smau$, to prove that the diagram \cref{Kg_runity_diagram} commutes, it suffices to prove that the boundary of the following diagram in $\IUT$ commutes.
\begin{equation}\label{Kg_runity_diag}
\begin{tikzpicture}[vcenter]
\def\v{-1.4} \def\h{-1.5}
\draw[0cell=.9]
(0,0) node (a11) {\Kg X \smasg \gsp}
(a11)++(3.8,0) node (a12) {\Kg X}
(a11)++(0,\v) node (a21) {\Kg X \smasg \Kg\gu}
(a12)++(0,\v) node (a22) {\phantom{\Kg(X \smag \gu)}} 
(a22)++(0,-.04) node (a22') {\Kg(X \smag \gu)}
(a11)++(0,-\v) node (a01) {\Kg X \smau \gsp}
(a21)++(0,\v) node (a31) {\Kg X \smau \Kg\gu}
;
\draw[1cell=.8]
(a11) edge node {\rsg} (a12)
(a11) edge node[pos=.5] {1 \smasg \Kgzero} (a21)
(a21) edge node[pos=.6] {\Kgtwo_{X,\gu}} (a22)
(a22') edge node[swap] {\Kg\rg} (a12)
(a01) edge node {\psma} (a11)
(a31) edge node[swap] {\psma} (a21)
(a01) [rounded corners=2pt] -| node[pos=.2] {\umu} (a12)
;
\draw[1cell=.8]
(a01) [rounded corners=2pt] -| ($(a11)+(\h,0)$) -- node {1 \smau \Kgzero} ($(a21)+(\h,0)$) |- (a31)
;
\draw[1cell=.8]
(a31) [rounded corners=2pt] -| node[pos=.2] {\ktwo} (a22')
;
\end{tikzpicture}
\end{equation}
In \cref{Kg_runity_diag}, the middle right region is \cref{Kg_runity_diagram}.  The left, top right, and bottom right regions commute by, respectively, the naturality of coequalizers, the definition of $\rsg$ \cref{rsg_x} in terms of $\umu$ \cref{Kgx_sphere_action}, and the definition of $\Kgtwo$ in terms of $\ktwo$ \cref{smau_kgx}.

The boundary diagram of \cref{Kg_runity_diag} is determined by the diagrams \cref{Kg_runity_vw} in $\Topst$ for pairs of objects $(V,W) \in (\IUsk)^2$, where $\ktwo_{V \oplus W}$ \cref{ktwo_u} denotes its restriction to the component with 
\[e = 1_{V \oplus W} \in \IU(V \oplus W, V \oplus W).\]
\begin{equation}\label{Kg_runity_vw}
\begin{tikzpicture}[vcenter]
\def\v{-1.4}
\draw[0cell]
(0,0) node (a11) {(\Kg X)_V \sma S^W}
(a11)++(4.6,0) node (a12) {(\Kg X)_{V \oplus W}}
(a11)++(0,\v) node (a21) {(\Kg X)_V \sma (\Kg\gu)_W}
(a12)++(0,\v) node (a22) {\phantom{(\Kg(X \smag \gu))_{V \oplus W}}} 
(a22)++(0,-.01) node (a22') {(\Kg(X \smag \gu))_{V \oplus W}}
;
\draw[1cell=.9]
(a11) edge node {\umu_{V,W}} (a12)
(a11) edge[transform canvas={xshift=1.5em}] node[swap] {1 \sma \Kgzero_W} (a21)
(a21) edge node {\ktwo_{V \oplus W}} (a22)
(a22') edge[transform canvas={xshift=-1em}] node[swap] {(\Kg\rg)_{V \oplus W}} (a12)
;
\end{tikzpicture}
\end{equation}
The following equalities in $(\Kg X)_{V \oplus W}$ prove that the diagram \cref{Kg_runity_vw} commutes on each representative \cref{Kgxv_rep}
\[\big(v \in (S^V)^{\sma\angordn} ; x \in X\angordn \big) \inspace (\Kg X)_V\]
and each point $w \in S^W = (S^W)^{\sma\ang{}}$.
\[\begin{aligned}
& ((\Kg\rg)_{V \oplus W}) (\ktwo_{V \oplus W}) (1 \sma \Kgzero_W) (v; x; w) \\
&= ((\Kg\rg)_{V \oplus W}) (\ktwo_{V \oplus W}) (v; x; w ; 1 \in \gu\ang{}) && \text{by \cref{Kgzero_v}} \\
&= ((\Kg\rg)_{V \oplus W}) \big(1_{V \oplus W} (\mu_{V,W}) (v \sma w) ; 1_{\angordn \oplus \ang{}} ; x ; 1\big) && \text{by \cref{ktwo_u_e}} \\
&= ((\Kg\rg)_{V \oplus W}) \big(v \oplus w \in (S^{V \oplus W})^{\sma\angordn} ; 1_{\angordn} ; x ; 1\big) && \text{by \cref{gsp_mult}} \\
&\overset{\vartriangle}{=} (v \oplus w; x) && \\
&= \umu_{V,W} (v; x; w) && \text{by \cref{Kgxv_action_rep}}
\end{aligned}\]
The equality $\vartriangle$ above uses the definition \cref{Kg_theta_v} of $(\Kg\rg)_{V \oplus W}$, along with the right unit version of \cref{ellg_f_n}.  This proves the commutativity of \cref{Kg_runity_vw}, hence also \cref{Kg_runity_diag} and \cref{Kg_runity_diagram}.
\end{proof}

\begin{lemma}\label{Kg_braid_axiom}
The triple \pcref{def:ggspace_gspectra,def:ggtop_gsp_mor,def:Kg_zero,def:Kgtwo}
\[\GGTop \fto{(\Kg,\Kgtwo,\Kgzero)} \GSp\]
satisfies the braid axiom \cref{enrmonfunctor-braided}.
\end{lemma}

\begin{proof}
The braid axiom asserts the commutativity of the following diagram in $\IUT$ for each pair $(X,Y)$ of $\Gskg$-spaces \cref{ggtop_obj}, where $\bsg$ \cref{bsg_xy} and $\beg$ \cref{beg_ff} are the braidings for, respectively, $\GSp$ and $\GGTop$.
\begin{equation}\label{Kg_braid_diagram}
\begin{tikzpicture}[vcenter]
\def\v{-1.4}
\draw[0cell]
(0,0) node (a11) {\Kg X \smasg \Kg Y}
(a11)++(4,0) node (a12) {\Kg Y \smasg \Kg X}
(a11)++(0,\v) node (a21) {\Kg(X \smag Y)}
(a12)++(0,\v) node (a22) {\Kg(Y \smag X)}
;
\draw[1cell=.9]
(a11) edge node {\bsg} (a12)
(a12) edge[transform canvas={xshift=-1em}] node {\Kgtwo_{Y,X}} (a22)
(a11) edge[transform canvas={xshift=1em}] node[swap] {\Kgtwo_{X,Y}} (a21)
(a21) edge node {\Kg\brg} (a22)
;
\end{tikzpicture}
\end{equation}
By the universal property of the coequalizer \cref{gsp_sma_coequal} that defines $\smasg$ in terms of $\smau$, to prove that the diagram \cref{Kg_braid_diagram} commutes, it suffices to prove that the boundary of the following diagram in $\IUT$ commutes.
\begin{equation}\label{Kg_braid_diag}
\begin{tikzpicture}[vcenter]
\def\v{-1.4} \def\h{-1.5}
\draw[0cell=.9]
(0,0) node (a11) {\Kg X \smasg \Kg Y}
(a11)++(3.5,0) node (a12) {\Kg Y \smasg \Kg X}
(a11)++(0,\v) node (a21) {\Kg(X \smag Y)}
(a12)++(0,\v) node (a22) {\Kg(Y \smag X)}
(a11)++(0,-\v) node (a01) {\Kg X \smau \Kg Y}
(a12)++(0,-\v) node (a02) {\Kg Y \smau \Kg X}
;
\draw[1cell=.8]
(a11) edge node {\bsg} (a12)
(a12) edge node[swap] {\Kgtwo_{Y,X}} (a22)
(a11) edge node {\Kgtwo_{X,Y}} (a21)
(a21) edge node {\Kg\brg} (a22)
(a01) edge node {\beu} (a02)
(a01) edge node {\psma} (a11)
(a02) edge node[swap] {\psma} (a12)
(a01) [rounded corners=2pt] -| node[pos=1] {\ktwo} ($(a11)+(\h,-.5*\v)$) -- ($(a11)+(\h,.5*\v)$) |- (a21)
;
\draw[1cell=.8]
(a02) [rounded corners=2pt] -| node[pos=1,swap] {\ktwo} ($(a12)+(-\h,-.5*\v)$) -- ($(a12)+(-\h,.5*\v)$) |- (a22)
;
\end{tikzpicture}
\end{equation}
In the diagram \cref{Kg_braid_diag}, the lower middle region is \cref{Kg_braid_diagram}.  The upper middle region commutes by the definition \cref{gsp_braiding_coequal} of $\bsg$ in terms of $\beu$.  The left and right regions commute by the definition of $\Kgtwo$ in terms of $\ktwo$ \cref{smau_kgx}.

For each object $U \in \IUsk$, the following equalities in $(\Kg(Y \smag X))_U$ prove that the boundary diagram in \cref{Kg_braid_diag} commutes at each representative quintuple \cref{kxukyu_rep}
\[\begin{split}
\bolde = \big(e \in \IU(V \oplus W,U)_\splus ; & \,v \in (S^V)^{\sma\angordm} ; w \in (S^W)^{\sma\angordmp} ;\\
&\, x \in X\angordm ; y \in Y\angordmp \big)
\end{split}\]
in $(\Kg X \smau \Kg Y)_U$.
\[\begin{aligned}
& (\Kg\brg)_U (\ktwo_U) (\bolde) \\
&= (\Kg\beg)_U \big(e\mu_{V,W} (v \sma w) ; 1_{\angordm \oplus \angordmp} ; x; y\big) && \text{by \cref{ktwo_u_e}} \\
&= \big(e\mu_{V,W} (v \sma w) ; \brg (1_{\angordm \oplus \angordmp} ; x; y) \big) && \text{by \cref{Kg_theta_v}} \\
&= \big(e\mu_{V,W} (v \sma w) ; 1_{\angordm \oplus \angordmp} \xi_{\angordmp,\angordm} ; y; x\big) && \text{by \cref{brg_components}} \\
&= \big(e\mu_{V,W} (v \sma w) ; \xi_{\angordmp,\angordm} 1_{\angordmp \oplus \angordm} ; y; x\big) && \\
&\overset{\vartriangle}{=} \scalebox{.8}{$\big(e \xi_{W,V} \mu_{W,V} (w \sma v) \xi_{\sma\angordm,\sma\angordmp} ; (Y \smag X)(\xi_{\angordmp,\angordm}) (1_{\angordmp \oplus \angordm} ; y; x) \big)$} && \text{by \cref{ggtop_ptday}}\\
&= \scalebox{.8}{$\big(e \xi_{W,V} \mu_{W,V} (w \sma v) (\xi_{\sma\angordm,\sma\angordmp}) (\sma\xi_{\angordmp,\angordm}) ; 1_{\angordmp \oplus \angordm} ; y; x \big)$} && \text{by \cref{Kgxv_dzeroone}} \\
&= \scalebox{.8}{$\big(e \xi_{W,V} \mu_{W,V} (w \sma v) (\xi_{\sma\angordm,\sma\angordmp}) (\xi_{\sma\angordmp,\sma\angordm}) ; 1_{\angordmp \oplus \angordm} ; y; x \big)$} && \text{by \cref{sma_braiding}} \\
&= \big(e \xi_{W,V} \mu_{W,V} (w \sma v) ; 1_{\angordmp \oplus \angordm} ; y; x \big) && \text{by \cref{symmoncatsymhexagon}} \\
&= (\ktwo_U) (e\xi_{W,V} ; w; v; y; x) && \text{by \cref{ktwo_u_e}} \\
&= (\ktwo_U) (\beu_U) (\bolde) && \text{by \cref{beu_XYU_z}} 
\end{aligned}\]
In the equality $\vartriangle$ above, the respective first entries of the two sides are equal by the following commutative diagram.
\[\begin{tikzpicture}[vcenter]
\def\v{-1.4}
\draw[0cell]
(0,0) node (a11) {(\sma\angordm) \sma (\sma\angordmp)}
(a11)++(3.7,0) node (a12) {S^V \sma S^W}
(a12)++(2.7,0) node (a13) {S^{V \oplus W}}
(a11)++(0,\v) node (a21) {(\sma\angordmp) \sma (\sma\angordm)}
(a12)++(0,\v) node (a22) {S^W \sma S^V}
(a13)++(0,\v) node (a23) {S^{W \oplus V}}
;
\draw[1cell=.9]
(a11) edge node {v \sma w} (a12)
(a12) edge node {\mu_{V,W}} (a13)
(a21) edge node {w \sma v} (a22)
(a22) edge node {\mu_{W,V}} (a23)
(a11) edge node[swap] {\xi_{\sma\angordm,\sma\angordmp}} (a21)
(a23) edge node[swap] {\xi_{W,V}} (a13)
;
\end{tikzpicture}\]
In the previous diagram, each of the two composites sends a pair of points
\[\big(t \in \sma\angordm ; t' \in \sma\angordmp \big)\]
to the point $v(t) \oplus w(t')$ in $S^{V \oplus W}$.  This finishes the proof of the commutativity of \cref{Kg_braid_diag}, hence also \cref{Kg_braid_diagram}.
\end{proof}

The next observation is the main result of this section.  It establishes the last step of our multifunctorial $G$-equivariant algebraic $K$-theory.  Recall the symmetric monoidal $\Gtop$-categories
\begin{itemize}
\item $\GGTop$ of $\Gskg$-spaces \pcref{expl:ggtop_gtopenr} and
\item $\GSp$ of $\gsp$-modules \pcref{def:gsp_module,def:gsp_smgtop}.
\end{itemize}

\begin{theorem}\label{thm:Kg_smgtop}
Suppose $G$ is a compact Lie group. 
\begin{enumerate}
\item\label{Kg_smgtop_i}The triple \pcref{def:ggspace_gspectra,def:ggtop_gsp_mor,def:Kg_zero,def:Kgtwo}
\begin{equation}\label{Kg_smfunctor}
\GGTop \fto{(\Kg,\Kgtwo,\Kgzero)} \GSp
\end{equation}
is a unital symmetric monoidal $\Gtop$-functor.\index{G-G-space@$\Gskg$-space!to orthogonal G-spectra@to orthogonal $G$-spectra}\index{orthogonal G-spectrum@orthogonal $G$-spectrum!from G-G-space@from $\Gskg$-space}
\item\label{Kg_smgtop_ii} Applying the endomorphism construction \pcref{def:EndF} to $\Kg$ yields a $\Gtop$-multifunctor 
\begin{equation}\label{EndKg}
\End(\GGTop) \fto{\End(\Kg)} \End(\GSp)
\end{equation}
with the same object assignment as $\Kg$.
\end{enumerate}
\end{theorem}

\begin{proof}
For assertion \pcref{Kg_smgtop_i}, the data $(\Kg,\Kgtwo,\Kgzero)$ are well defined by \cref{Kg_gtop_functor,Kgzero_welldef,kgtwo_welldef}.  The associativity axiom, the unity axioms, and the braid axiom of a symmetric monoidal $\Gtop$-functor \pcref{definition:monoidal-V-fun,definition:braided-monoidal-vfunctor} are proved in 
\cref{Kg_associativity_axiom,Kg_unity_axiom,Kg_braid_axiom}.  Assertion \pcref{Kg_smgtop_ii} is an instance of \cref{EndF_multi}, applied to the symmetric monoidal $\Gtop$-functor $\Kg$ in assertion \pcref{Kg_smgtop_i}.
\end{proof}

\begin{notation}\label{not:EndKg}
Continuing the notational convention of \cref{expl:ggcat_gcatenr,not:clast'}, the $\Gtop$-multifunctor $\End(\Kg)$ in \cref{thm:Kg_smgtop} \pcref{Kg_smgtop_ii} is also denoted by $\Kg \cn \GGTop \to \GSp$.
\end{notation}

\begin{remark}[Equivariant Infinite Loop Space Machines]\label{rk:gspace_gspectra}
Our construction $\Kg$ builds orthogonal $G$-spectra from $\Gskg$-spaces, and it is a unital symmetric monoidal functor in the $\Gtop$-enriched sense \pcref{thm:Kg_smgtop}.  In this Remark, we briefly discuss other constructions in the literature that build $G$-spectra from space-level data.  Since we do not use these other constructions in this work, we refer the interested reader to the cited work for more discussion.

The nonequivariant operadic infinite loop space machine of May \cite{may} is extended to the equivariant context in the work \cite{gm17} of Guillou and May.  The work \cite{mmo} of May, Merling, and Osorno modernizes and extends earlier work of Shimakawa \cite{shimakawa89}, which was written before the invention of modern categories of spectra.  The paper \cite{mmo} constructs orthogonal $G$-spectra from $\FG$-$G$-spaces.  The indexing category $\FG$ is the $G$-category of pointed finite $G$-sets and pointed functions, on which $G$ acts by conjugation.  There is a detailed comparison between this Segal-style equivariant machine and the operadic equivariant machine from \cite{gm17}.

Shimakawa proved in \cite{shimakawa91} that there is an equivalence between the category of $\FG$-$G$-spaces and the category of $\Fsk$-$G$-spaces.  Using Shimakawa's results \cite{shimakawa89,shimakawa91}, the work \cite{ostermayr} of Ostermayr proves that there is an equivalence between the category of very special $\Fsk$-$G$-spaces and the category of connective $G$-symmetric spectra, along the lines of \cite{bousfield_friedlander}.  The work \cite{gmmo19} of Guillou, May, Merling, and Osorno constructs a symmetric monoidal functor from the symmetric monoidal category of $\Fsk$-$G$-spaces to the symmetric monoidal category of orthogonal $G$-spectra.  The main constructions in \cite{gmmo19} are based on monadic bar constructions.

The work \cite{kong_may_zou} of Kong, May, and Zou conceptualizes and generalizes some of \cite{gm17,may,may-precisely,mmo}, proving homotopical variants of Beck's monadicity theorem.  In specific examples, this framework constructs $\Einf$ ring $G$-spectra from space-level data.  The work \cite{costenoble_waner} of Costenoble and Waner constructs Lewis-May $G$-spectra from topological presheaves on the orbit category of $G$ equipped with an equivariant $\Einf$ structure.  The main result of Guillou and May in \cite{gm-presheaves} proves that there is a Quillen equivalence between the category of orthogonal $G$-spectra and the category of spectral Mackey functors.
\end{remark}

\section{Multifunctorial Equivariant $K$-Theory}
\label{sec:Kgo_multi}

This section constructs our $G$-equivariant algebraic $K$-theory multifunctor
\[\MultpsO \fto{\Kgo} \GSp.\] 
Its domain is the $\Gcat$-multicategory $\MultpsO$ of $\Op$-pseudoalgebras.  Its codomain is the symmetric monoidal $\Gtop$-category $\GSp$ of orthogonal $G$-spectra.  The main application of the multifunctor $\Kgo$ is that it preserves all algebraic structures parametrized by multicategories enriched in $\Gcat$ or $\Gtop$.  Thus, $\Kgo$ sends the purely algebraic data of an equivariant $\Einf$-algebra or an $\Ninf$-algebra in $\MultpsO$ to an equivariant $\Einf$-algebra or an $\Ninf$-algebra in $\GSp$.

\secoutline
\begin{itemize}
\item \cref{thm:Kgo_multi} records the $\Gtop$-multifunctor $\Kgo$ and its strong variant $\Kgosg$.
\item \cref{expl:Kgo_obj} unravels $\Kgo$ and $\Kgosg$ on objects.
\item \cref{thm:Kgo_preservation} records the fact that $\Kgo$ and $\Kgosg$ preserve equivariant algebraic structures, including equivariant $\Einf$-algebras and $\Ninf$-algebras. 
\item \cref{ex:Kgo_preservation} illustrates \cref{thm:Kgo_preservation} with the $G$-Barratt-Eccles operad $\GBE$, whose pseudoalgebras are genuine symmetric monoidal $G$-categories.   In particular, $\Kgo$ and $\Kgosg$ send each equivariant $\Einf$-algebra or $\Ninf$-algebra of genuine symmetric monoidal $G$-categories to an equivariant $\Einf$-algebra or $\Ninf$-algebra of orthogonal $G$-spectra.
\end{itemize}

The following result is the main result of this chapter.  It combines \cref{thm:Kg_smgtop} with results from previous chapters, as specified in the statements, to construct $\Kgo$ and its strong variant.  It also uses the endomorphism construction \pcref{definition:EndK} to turn a symmetric monoidal $\Gtop$-category into a $\Gtop$-multicategory.  

\begin{theorem}\label{thm:Kgo_multi}
Suppose $G$ is a compact Lie group, and $\Op$ is a $\Tinf$-operad \pcref{as:OpA}.  
\begin{enumerate}
\item\label{thm:Kgo_multi_i} There is a $\Gtop$-multifunctor $\Kgo$ defined as the composite\index{multifunctorial K-theory@multifunctorial $K$-theory}\index{pseudoalgebra!to orthogonal G-spectra@to orthogonal $G$-spectra}\index{orthogonal G-spectrum@orthogonal $G$-spectrum!from pseudoalgebra}
\[\begin{tikzpicture}
\def\v{-1.4} 
\draw[0cell]
(0,0) node (a11) {\MultpsO}
(a11)++(3.2,0) node (a12) {\phantom{\GSp}}
(a12)++(0,.03) node (a12') {\GSp}
(a11)++(0,\v) node (a21) {\GGCat}
(a12)++(0,\v) node (a22) {\phantom{\GGTop}}
(a22)++(.3,0) node (a22') {\GGTop}
;
\draw[1cell=.9]
(a11) edge node {\Kgo} (a12)
(a11) edge[shorten <=-.5ex] node[swap] {\Jgo} (a21)
(a21) edge node {\clast} (a22')
(a22) edge node[swap] {\Kg} (a12')
;
\end{tikzpicture}\]
given by the following data.
\begin{itemize}
\item $\Jgo$ is the $\Gcat$-multifunctor in \cref{thm:Jgo_multifunctor}, with domain $\MultpsO$ \pcref{thm:multpso} and codomain $\GGCat$ \pcref{expl:ggcat_gcatenr}.  $\Jgo$ is regarded as a $\Gtop$-multifunctor by changing enrichment along the classifying space functor $\cla$ \pcref{cla_multi}.
\item $\clast$ is the $\Gtop$-multifunctor in \cref{End_clastpev}.  Its codomain is the $\Gtop$-multicategory associated to the symmetric monoidal $\Gtop$-category $\GGTop$ \pcref{expl:ggtop_gtopenr}.
\item $\Kg$ is the $\Gtop$-multifunctor in \cref{EndKg}.  Its codomain is the $\Gtop$-multicategory associated to the symmetric monoidal $\Gtop$-category $\GSp$ \pcref{gspectra_smgtop}.
\end{itemize}
\item\label{thm:Kgo_multi_ii} 
There is a strong variant of $\Kgo$ given by the composite $\Gtop$-multifunctor
\[\begin{tikzpicture}
\def\v{-1.4} 
\draw[0cell]
(0,0) node (a11) {\MultpspsO}
(a11)++(3.2,0) node (a12) {\phantom{\GSp}}
(a12)++(0,.03) node (a12') {\GSp}
(a11)++(0,\v) node (a21) {\GGCat}
(a12)++(0,\v) node (a22) {\phantom{\GGTop}}
(a22)++(.3,0) node (a22') {\GGTop}
;
\draw[1cell=.9]
(a11) edge node {\Kgosg} (a12)
(a11) edge[shorten <=-.5ex] node[swap] {\Jgosg} (a21)
(a21) edge node {\clast} (a22')
(a22) edge node[swap] {\Kg} (a12')
;
\end{tikzpicture}\]
in which $\Jgosg$ is the $\Gcat$-multifunctor in \cref{thm:Jgo_multifunctor}, regarded as a $\Gtop$-multifunctor by changing enrichment along $\cla$.
\end{enumerate}
\end{theorem}

\begin{proof}
In each of \pcref{thm:Kgo_multi_i} and \pcref{thm:Kgo_multi_ii}, the three constituent arrows are $\Gtop$-multifunctors.  Thus, their composite is also a $\Gtop$-multifunctor \cref{multifunctors_compose}.
\end{proof}

\begin{remark}\label{rk:Kgo_multi}
Consider \cref{thm:Kgo_multi}.
\begin{enumerate}
\item\label{rk:Kgo_multi_i} The assumption that $G$ be a compact Lie group is only needed for $\Kg$, whose codomain $\GSp$ \pcref{ch:spectra} requires this assumption. 
\item\label{rk:Kgo_multi_ii} The $\Tinf$-operad $\Op$ \pcref{as:OpA} is only used for the first step, which means $\Jgo$ or $\Jgosg$, since neither $\clast$ nor $\Kg$ involves $\Op$.  
\item\label{rk:Kgo_multi_iii} $\clast$ and $\Kg$ are not merely $\Gtop$-multifunctors; they are actually symmetric monoidal $\Gtop$-functors by \cref{clastpev,Kg_smfunctor}.
\item\label{rk:Kgo_multi_iv} By \cref{def:MultpsO_karycat,def:multicatO}, $\MultpsO$ and $\MultpspsO$ only differ in their $k$-ary 1-cells.  In $\MultpsO$, the $k$-ary 1-cells are $k$-lax $\Op$-morphisms \pcref{def:k_laxmorphism}.  In $\MultpspsO$, the $k$-ary 1-cells are $k$-ary $\Op$-pseudomorphisms, which have invertible action constraints.  For further discussion of $\Jgo$ and $\Jgosg$, see \cref{rk:JgoJgosg,rk:no_strictification,ex:JgBE}.\defmark
\end{enumerate}
\end{remark}

\begin{explanation}[$\Kgo$ on Objects]\label{expl:Kgo_obj}
The $\Gcat$-multifunctor\index{multifunctorial K-theory@multifunctorial $K$-theory!object assignment} \pcref{thm:Jgo_multifunctor}
\[\MultpsO \fto{\Jgo} \GGCat\]
sends an $\Op$-pseudoalgebra $\A$ to the $\Gskg$-category $\Jgo\A = \Adash$ \pcref{A_ptfunctor}.  The symmetric monoidal $\Gtop$-functor \pcref{thm:ggcat_ggtop}
\[\GGCat \fto{\clast} \GGTop\]
sends a $\Gskg$-category to a $\Gskg$-space by post-composing with the classifying space functor $\cla$ \cref{clast_obj}.  The symmetric monoidal $\Gtop$-functor \pcref{thm:Kg_smgtop}
\[\GGTop \fto{\Kg} \GSp\]
sends a $\Gskg$-space $X$ to an orthogonal $G$-spectrum $\Kg X$ \pcref{def:ggspace_gspectra}.  Thus, for each $\Op$-pseudoalgebra $\A$ \pcref{def:pseudoalgebra}, the orthogonal $G$-spectrum $\Kgo\A$ has, for each object $V \in \IU$, $V$-component pointed $G$-space given by the coend
\[(\Kgo\A)_V = \int^{\angordn \in \Gsk} (S^V)^{\sma\angordn} \sma \cla\Aangordn,\]
where $\Aangordn$ is the pointed $G$-category of $\angordn$-systems \pcref{def:nsystem,def:nsystem_morphism,def:Aangordn_system,def:Aangordn_gcat}.  The group $G$ acts on representatives of $(\Kgo\A)_V$ diagonally \cref{g_tx}, with the $G$-action on the category $\Aangordn$ given in \cref{def:Aangordn_gcat}.  The right $\gsp$-action is defined in \cref{Kgx_action_vw}, and it only uses the factor $(S^V)^{\sma\angordn}$.

The strong variant $\Kgosg$ is given objectwise by
\[(\Kgosg\A)_V = \int^{\angordn \in \Gsk} (S^V)^{\sma\angordn} \sma \cla\Aangordnsg,\]
where $\Aangordnsg$ is the full sub-$G$-category of $\Aangordn$ consisting of strong $\angordn$-systems \pcref{def:nsystem}.
\end{explanation}

\subsection*{Preservation of Equivariant Algebraic Structures}
As the main application of \cref{thm:Kgo_multi}, we observe that the $\Gtop$-multifunctors $\Kgo$ and $\Kgosg$ preserve all algebraic structures parametrized by $\Gtop$-multicategories, including equivariant $\Einf$-algebras and $\Ninf$-algebras \pcref{def:Einfty_operads,def:Ninfty_operads}.  As we discuss in \cref{ex:Kgo_preservation} below, by changing enrichment along the classifying space functor $\cla$, $\Kgo$ and $\Kgosg$ also preserve all algebraic structures parametrized by $\Gcat$-multicategories.

\begin{theorem}\label{thm:Kgo_preservation}
Suppose $\Op$ is a $\Tinf$-operad \pcref{as:OpA} with $G$ a compact Lie group, and $\cQ$ is a $\Gtop$-multicategory.
\begin{enumerate}
\item\label{Kgo_preservation_i} Suppose $f \cn \cQ \to \MultpsO$ is a $\Gtop$-multifunctor.  Then the composite
\[\cQ \fto{f} \MultpsO \fto{\Kgo} \GSp\]
is a $\Gtop$-multifunctor.
\item\label{Kgo_preservation_ii} Suppose $f \cn \cQ \to \MultpspsO$ is a $\Gtop$-multifunctor.   Then the composite
\[\cQ \fto{f} \MultpspsO \fto{\Kgosg} \GSp\]
is a $\Gtop$-multifunctor. 
\item\label{Kgo_preservation_iii} For a finite group $G$, each of the $\Gtop$-multifunctors $\Kgo$ and $\Kgosg$ preserves equivariant $\Einf$-algebras and $\Ninf$-algebras \pcref{def:Einfty_operads,def:Ninfty_operads}.
\end{enumerate}
\end{theorem}

\begin{proof}
Each assertion follows from \cref{thm:Kgo_multi} because multifunctors enriched in a given symmetric monoidal category, including $\Gtop$, are closed under composition \cref{multifunctors_compose}.
\end{proof}

The following example illustrates \cref{thm:Kgo_preservation} and also continues the thread of ideas discussed in \cref{ex:Einf_GBE,ex:Jgo_preservation}.

\begin{example}[$\Einf$ and $\Ninf$ $G$-Spectra from Categorical Data]\label{ex:Kgo_preservation}
In this example, we focus on $\Kgo$; the same discussion applies to the strong variant $\Kgosg$.

Suppose $\cQ$ is a $\Gcat$-multicategory, and $f \cn \cQ \to \MultpsO$ is a $\Gcat$-multifunctor.  By  changing enrichment along the classifying space functor $\cla$ \pcref{cla_multi}, $f$ becomes a $\Gtop$-multifunctor, which we denote by the same notation.  By \cref{thm:Kgo_preservation}, the $\cQ$-algebra $f$ is sent by $\Kgo$ to the $\cQ$-algebra 
\[\cQ \fto{f} \MultpsO \fto{\Kgo} \GSp\]
in the $\Gtop$-multicategory $\GSp$ of orthogonal $G$-spectra.

\parhead{Specializing to $\cQ = \GBE$}.  As we discuss in \cref{ex:Einf_GBE}, for a finite group $G$, the $G$-Barratt-Eccles operad $\GBE$ \pcref{def:GBE} is both a $G$-categorical $\Einf$-operad and a $G$-categorical $\Ninf$-operad \pcref{def:Einfty_operads,def:Ninfty_operads}.  A $\Gcat$-multifunctor
\[\GBE \fto{f} \MultpsO\]
is both an equivariant $\Einf$-algebra and an $\Ninf$-algebra in $\MultpsO$ \pcref{def:multicatO}, whose objects are $\Op$-pseudoalgebras \pcref{def:pseudoalgebra}.  After changing enrichment along the classifying space functor $\cla$, the composite $\Gtop$-multifunctor
\[\GBE \fto{f} \MultpsO \fto{\Kgo} \GSp\]
is both an equivariant $\Einf$-algebra and an $\Ninf$-algebra in the $\Gtop$-multicategory $\GSp$ of orthogonal $G$-spectra.

\parhead{Specializing to genuine symmetric monoidal $G$-categories}.  Similar to \cref{ex:Jgo_preservation}, we can further specialize the discussion in the previous paragraph to the case $\Op = \GBE$.  By \cref{def:GBE_pseudoalg}, $\GBE$-pseudoalgebras are genuine symmetric monoidal $G$-categories.  Thus, a $\Gcat$-multifunctor
\[\GBE \fto{f} \MultpsGBE\]
is both an equivariant $\Einf$-algebra and an $\Ninf$-algebra of genuine symmetric monoidal $G$-categories.  After changing enrichment along $\cla$, the composite $\Gtop$-multifunctor
\[\GBE \fto{f} \MultpsGBE \fto{\Kggbe} \GSp\]
is both an equivariant $\Einf$-algebra and an $\Ninf$-algebra of orthogonal $G$-spectra.
\end{example}

\part*{Appendices}
\appendix

\chapter{Symmetric Monoidal Categories}
\label{ch:prelim}
This appendix reviews 2-categories and enriched symmetric monoidal categories.  References are provided at the beginning of each section.  Below we list the sections in this appendix, along with the main concepts in each section.

\appsecname{sec:monoidalcat}

\begin{itemize}
\item \appentry{Monoidal categories}{def:monoidalcategory}
\item \appentry{Braided monoidal, symmetric monoidal, permutative, closed, and diagram categories}{def:braidedmoncat,def:symmoncat,def:closedcat,def:diagramcat}
\item \appentry{Monoids and modules}{def:monoid,def:modules}
\item \appentry{Monoidal functors and monoidal natural transformations}{def:monoidalfunctor,def:monoidalnattr}
\end{itemize}

\appsecname{sec:enrichedcat}
\begin{itemize}
\item \appentry{$\V$-categories, $\V$-functors, and $\V$-natural transformations}{def:enriched-category,def:enriched-functor,def:enriched-natural-transformation}
\end{itemize}

\appsecname{sec:twocategories}
\begin{itemize}
\item \appentry{2-categories, 2-functors, and 2-natural transformations}{def:twocategory,def:twofunctor,def:twonaturaltr}
\item \appentry{Modifications}{def:modification,def:modcomposition}
\item \appentry{2-equivalences}{def:twoequivalence,thm:twoequivalences}
\end{itemize}

\appsecname{sec:enrmonoidalcat}
\begin{itemize}
\item \appentry{Tensor product of $\V$-categories}{definition:vtensor-0}
\item \appentry{Symmetric monoidal category of small $\V$-categories}{definition:unit-vcat,definition:vtensor-unitors,definition:vtensor-assoc,definition:vtensor-beta,theorem:vcat-mon}
\item \appentry{Monoidal $\V$-categories}{definition:monoidal-vcat}
\item \appentry{Braided and symmetric monoidal $\V$-categories}{definition:braided-monoidal-vcat,definition:symm-monoidal-vcat}
\item \appentry{Symmtric monoidal $\Cat$-category of small $\V$-categories}{theorem:vcat-cat-mon}
\item \appentry{Monoidal $\V$-functors}{definition:monoidal-V-fun}
\item \appentry{Braided and symmetric monoidal $\V$-functors}{definition:braided-monoidal-vfunctor}
\item \appentry{Self-enrichment of symmetric monoidal closed categories}{definition:canonical-v-enrichment,theorem:v-closed-v-sm}
\item \appentry{Change of enrichment}{thm:change-enrichment}
\end{itemize}

\section{Monoidal Categories}
\label{sec:monoidalcat}

This section is a brief review of monoidal category theory.  More detailed references include \cite{joyal-street,maclane,cerberusI,cerberusII}.  Throughout this work, we assume Grothendieck's axiom of universes to deal with set-theoretic size issues.  More discussion of universes can be found in \cite[Section 1.1]{johnson-yau}, \cite[I.6]{maclane}, and \cite{agv,maclane-foundation}.

\begin{definition}\label{def:universe}
A \emph{Grothendieck universe}\index{universe} is a set\label{notation:universe} $\zU$ that satisfies the following four axioms.
\begin{enumerate}
\item\label{univ1} If $a \in \zU$ and $b \in a$, then $b\in\zU$.
\item\label{univ2} If $a \in \zU$, then the set of subsets of $a$ also belongs to $\zU$.
\item\label{univ3} If $a \in \zU$ and $b_x \in \zU$ for each $x \in a$, then $\bigcup_{x\in a} b_x$ belongs to $\zU$.
\item\label{univ4} The set $\bN$ of finite ordinals is an element in $\zU$.\defmark\end{enumerate}
\end{definition}

\begin{convention}[Axiom of Universes]\label{conv:universe}\index{Grothendieck Universe}\index{universe}\index{Axiom of Universes}
Throughout this work, we assume that every set belongs to some universe.  This is called the \emph{Axiom of Universes}.  We fix a universe $\zU$\label{not:universe}, in which an element is called a \emph{set}.  A subset of $\zU$ is called a \emph{class}.  A categorical structure is \emph{small} if it has a set of objects.  We silently replace $\zU$ by a larger universe $\zV$ in which $\zU$ is a set wherever necessary.
\end{convention}

\subsection*{Symmetric Monoidal Categories}

\begin{definition}\label{def:monoidalcategory}\index{monoidal category}\index{category!monoidal}
A \emph{monoidal category}
\[(\C,\otimes,\tu,\alpha,\lambda,\rho)\]
consists of the following data.
\begin{itemize}
\item $\C$ is a category.
\item \label{notation:monoidal-product}$\otimes \cn \C \times \C \to \C$ is a functor, called the \index{monoidal category!monoidal product}\emph{monoidal product}.
\item \label{not:monoidalunit}$\tu \in \C$ is an object, called the \index{monoidal category!monoidal unit}\emph{monoidal unit}.
\item \label{not:associativityiso}\label{not:unitisos}$\alpha$, $\lambda$, and $\rho$ are natural isomorphisms for objects $a,b,c \in \C$ as follows.  
\[\begin{tikzcd}[column sep=huge,row sep=tiny]
(a \otimes b) \otimes c \ar{r}{\alpha_{a,b,c}}[swap]{\iso} & a \otimes (b \otimes c)
\end{tikzcd}\]\vspace{-1em}
\[\begin{tikzcd}[column sep=large]
\tu \otimes a \ar{r}{\lambda_a}[swap]{\iso} & a & a \otimes \tu \ar{l}{\iso}[swap]{\rho_a}
\end{tikzcd}\]
They are called, respectively, the \index{associativity isomorphism}\index{monoidal category!associativity isomorphism}\emph{associativity isomorphism}, the \index{left unit isomorphism}\index{monoidal category!left unit isomorphism}\emph{left unit isomorphism}, and the \index{right unit isomorphism}\index{monoidal category!right unit isomorphism}\emph{right unit isomorphism}.
\end{itemize}
The following middle unity\index{monoidal category!middle unity axiom} and pentagon\index{monoidal category!pentagon axiom}\index{pentagon axiom} diagrams are required to commute for objects $a,b,c,d \in \C$, where $a \otimes b$ is abbreviated to $ab$ to save space.
\begin{equation}\label{monoidalcataxioms}
\begin{tikzpicture}[xscale=1,yscale=1,vcenter]
\tikzset{0cell/.append style={nodes={scale=.85}}}
\tikzset{1cell/.append style={nodes={scale=.85}}}
\def\h{1.6} \def\g{1.5} \def\v{.9} \def\u{2}
\draw[0cell]
(0,0) node (a) {(a \tu) b}
(a)++(0,-2) node (b) {a (\tu b)}
(a)++(1.5,-1) node (c) {ab}
;
\draw[1cell]  
(a) edge node[pos=.4] {\rho_a 1_b} (c)
(a) edge node[swap] {\alpha_{a,\tu,b}} (b)
(b) edge node[swap,pos=.4] {1_a \lambda_b} (c)
;
\begin{scope}[shift={(6,0)}]
\draw[0cell]
(0,0) node (x0) {(ab)(cd)}
(x0)++(-\h,-\v) node (x11) {((ab)c)d}
(x0)++(\h,-\v) node (x12) {a(b(cd))}
(x0)++(-\g,-\u) node (x21) {(a(bc))d}
(x0)++(\g,-\u) node (x22) {a((bc)d)}
;
\draw[1cell]
(x11) edge node[pos=.2] {\al_{ab,c,d}} (x0)
(x0) edge node[pos=.8] {\al_{a,b,cd}} (x12)
(x11) edge node[swap,pos=.2] {\al_{a,b,c} 1_d} (x21)
(x21) edge node {\al_{a,bc,d}} (x22)
(x22) edge node[swap,pos=.8] {1_a \al_{b,c,d}} (x12)
;
\end{scope}
\end{tikzpicture}
\end{equation}
This finishes the definition of a monoidal category.  We usually abbreviate a monoidal category to $(\C, \otimes, \tu)$ or $\C$.  Moreover, we define the following.
\begin{itemize}
\item A monoidal category is \emph{strictly unital}\index{monoidal category!strictly unital}\index{strictly unital monoidal category} if $\lambda$ and $\rho$ are identities.
\item A monoidal category is \emph{strict}\index{monoidal category!strict}\index{strict monoidal category} if $\alpha$, $\lambda$, and $\rho$ are identities.
\item A (monoidal) category is also called a \index{monoidal 1-category}\index{1-category}\emph{(monoidal) 1-category}.\defmark
\end{itemize}
\end{definition}

\begin{remark}[Unity Properties]
In each monoidal category, the equality\index{monoidal category!unity properties} 
\begin{equation}\label{lambda=rho}
\lambda_{\tu} = \rho_{\tu} \cn \tu \otimes \tu \to \tu
\end{equation}
holds, and the following unity diagrams commute.  
\begin{equation}\label{moncat-other-unit-axioms}
\begin{tikzcd}[column sep=normal]
(\tu \otimes a) \otimes b \dar[swap]{\lambda_a \otimes 1_b} \rar{\alpha_{\tu,a,b}}
& \tu \otimes (a \otimes b) \dar{\lambda_{a \otimes b}}\\ 
a \otimes b \rar[equal] & a \otimes b
\end{tikzcd}\qquad
\begin{tikzcd}[column sep=normal]
(a \otimes b) \otimes \tu \dar[swap]{\rho_{a \otimes b}} \rar{\alpha_{a,b,\tu}}
& a \otimes (b \otimes \tu) \dar{1_a \otimes \rho_b}\\ 
a \otimes b \rar[equal] & a \otimes b
\end{tikzcd}
\end{equation} 
Proofs of these properties can be found in \cite{kelly-coherence} and \cite[Section 2.2]{johnson-yau}.
\end{remark}

\begin{convention}[Left Normalized Bracketing]\label{expl:leftbracketing}
In the absence of clarifying parentheses, an iterated monoidal product is assumed \index{left normalized}\emph{left normalized}.  This means that the left half of each pair of parentheses is situated at the left end of the expression.
\end{convention}  

\begin{example}
As an example of our left normalized brancketing \cref{expl:leftbracketing}, we denote
\[a \otimes b \otimes c \otimes d = \big((a \otimes b) \otimes c\big) \otimes d.\]
In some diagrams, we omit the symbol $\otimes$ to save space, so the iterated monoidal product above is sometimes abbreviated to $abcd$.
\end{example}

\begin{definition}\label{def:braidedmoncat}\index{braided monoidal category}\index{monoidal category!braided}\index{category!braided monoidal}
A \emph{braided monoidal category} is a pair $(\C,\xi)$ consisting of
\begin{itemize}
\item a monoidal category $\C = (\C,\otimes,\tu,\al,\lambda,\rho)$ and
\item a natural isomorphism\label{notation:symmetry-iso}, called the \index{braiding}\emph{braiding},
\begin{equation}\label{braiding_bmc}
a \otimes b \fto[\iso]{\xi_{a,b}} b \otimes a
\end{equation}
for $a,b \in \C$.
\end{itemize}
The following hexagon diagrams\index{hexagon diagram} are required to commute for objects $a,b,c \in \C$. 
\begin{equation}\label{hexagon-braided}
\def\sb{\scalebox{.7}}
\begin{tikzpicture}[scale=.75,commutative diagrams/every diagram]
			\node (P0) at (0:2cm) {\sb{$b \otimes (c \otimes a)$}};
			\node (P1) at (60:2cm) {\makebox[3ex][l]{\sb{$b \otimes (a \otimes c)$}}};
			\node (P2) at (120:2cm) {\makebox[3ex][r]{\sb{$(b \otimes a) \otimes c$}}};
			\node (P3) at (180:2cm) {\sb{$(a \otimes b) \otimes c$}};
			\node (P4) at (240:2cm) {\makebox[3ex][r]{\sb{$a \otimes (b \otimes c)$}}};
			\node (P5) at (300:2cm) {\makebox[3ex][l]{\sb{$(b \otimes c) \otimes a$}}};
			\path[commutative diagrams/.cd, every arrow, every label]
			(P3) edge node[pos=.25] {\sb{$\xi_{a,b}\otimes 1_c$}} (P2)
			(P2) edge node {\sb{$\alpha$}} (P1)
			(P1) edge node[pos=.75] {\sb{$1_b \otimes \xi_{a,c}$}} (P0)
			(P3) edge node[swap,pos=.4] {\sb{$\alpha$}} (P4)
			(P4) edge node {\sb{$\xi_{a,b \otimes c}$}} (P5)
			(P5) edge node[swap,pos=.6] {\sb{$\alpha$}} (P0);
\end{tikzpicture}\qquad
\begin{tikzpicture}[scale=.75, commutative diagrams/every diagram]
			\node (P0) at (0:2cm) {\sb{$(c \otimes a) \otimes b$}};
			\node (P1) at (60:2cm) {\makebox[3ex][l]{\sb{$(a \otimes c) \otimes b$}}};
			\node (P2) at (120:2cm) {\makebox[3ex][r]{\sb{$a \otimes (c \otimes b)$}}};
			\node (P3) at (180:2cm) {\sb{$a \otimes (b \otimes c)$}};
			\node (P4) at (240:2cm) {\makebox[3ex][r]{\sb{$(a \otimes b) \otimes c$}}};
			\node (P5) at (300:2cm) {\makebox[3ex][l]{\sb{$c \otimes (a \otimes b)$}}};
			\path[commutative diagrams/.cd, every arrow, every label]
			(P3) edge node[pos=.25] {\sb{$1_a \otimes \xi_{b,c}$}} (P2)
			(P2) edge node {\sb{$\alpha^\inv$}} (P1)
			(P1) edge node[pos=.75] {\sb{$\xi_{a,c}\otimes 1_b$}} (P0)
			(P3) edge node[swap,pos=.4] {\sb{$\alpha^\inv$}} (P4)
			(P4) edge node {\sb{$\xi_{a \otimes b, c}$}} (P5)
			(P5) edge node[swap,pos=.6] {\sb{$\alpha^\inv$}} (P0);
\end{tikzpicture}
\end{equation}
This finishes the definition of a braided monoidal category.
\end{definition}

\begin{remark}[Additional Unity Properties]
In each braided monoidal category, the following unity diagram\index{braided monoidal category!unity properties} commutes for each object $a$.
\begin{equation}\label{braidedunity}
\begin{tikzcd}[column sep=large]
a \otimes \tu \dar[swap]{\rho_a} \rar{\xi_{a,\tu}} & \tu \otimes a \dar{\lambda_a} \rar{\xi_{\tu,a}} & a \otimes \tu \dar{\rho_a}\\ 
a \rar[equal] & a \rar[equal] & a
\end{tikzcd}
\end{equation}
See \cite[1.3.21]{cerberusII} for the proof. 
\end{remark}

\begin{definition}\label{def:symmoncat}\index{symmetric monoidal category}\index{monoidal category!symmetric}\index{category!symmetric monoidal}
A \emph{symmetric monoidal category} is a pair $(\C,\xi)$ consisting of
\begin{itemize}
\item a monoidal category $\C = (\C,\otimes,\tu,\al,\lambda,\rho)$ and 
\item a natural isomorphism, called the \emph{symmetry isomorphism} or the \emph{braiding}, 
\[a \otimes b \fto[\iso]{\xi_{a,b}} b \otimes a\]
for $a,b \in \C$.
\end{itemize}
The following symmetry\index{symmetry axiom} and \index{hexagon diagram}hexagon diagrams are required to commute for objects $a,b,c \in \C$. 
\begin{equation}\label{symmoncatsymhexagon}
\begin{tikzpicture}[xscale=3,yscale=1.2,vcenter]
\def\h{.1}
\draw[0cell=.8] 
(0,0) node (a) {a \otimes b}
(a)++(.7,0) node (c) {a \otimes b}
(a)++(.35,-1) node (b) {b \otimes a}
;
\draw[1cell=.8] 
(a) edge node {1_{a \otimes b}} (c)
(a) edge node [swap,pos=.3] {\xi_{a,b}} (b)
(b) edge node [swap,pos=.7] {\xi_{b,a}} (c)
;
\begin{scope}[shift={(1.5,.5)}]
\draw[0cell=.8] 
(0,0) node (x11) {(b \otimes a) \otimes c}
(x11)++(.9,0) node (x12) {b \otimes (a \otimes c)}
(x11)++(-\h,-1) node (x21) {(a \otimes b) \otimes c}
(x12)++(\h,-1) node (x22) {b \otimes (c \otimes a)}
(x11)++(0,-2) node (x31) {a \otimes (b \otimes c)}
(x12)++(0,-2) node (x32) {(b \otimes c) \otimes a}
;
\draw[1cell=.8]
(x21) edge node[pos=.25] {\xi_{a,b} \otimes 1_c} (x11)
(x11) edge node {\alpha} (x12)
(x12) edge node[pos=.75] {1_b \otimes \xi_{a,c}} (x22)
(x21) edge node[swap,pos=.25] {\alpha} (x31)
(x31) edge node {\xi_{a,b \otimes c}} (x32)
(x32) edge node[swap,pos=.75] {\alpha} (x22)
;
\end{scope}
\end{tikzpicture}
\end{equation}
A \emph{permutative category}\index{permutative category}\index{category!permutative} is a strict symmetric monoidal category, that is, a symmetric monoidal category with $\alpha$, $\lambda$, and $\rho$ identities.  For a permutative category, we sometimes denote its monoidal product by \label{not:opluse}$\oplus$ and its monoidal unit by $\pu$.
\end{definition}

\begin{remark}\label{rk:smcat}
In the presence of the symmetry axiom in \cref{symmoncatsymhexagon}, the two hexagon diagrams in \cref{hexagon-braided} are equivalent.  Thus, in \cref{symmoncatsymhexagon} the hexagon diagram may be replaced by the right hexagon in \cref{hexagon-braided}.  Moreover, a symmetric monoidal category is precisely a braided monoidal category that satisfies the symmetry axiom.  In particular, in each symmetric monoidal category, the unity diagrams in \cref{braidedunity} commute.
\end{remark}

\begin{definition}\label{def:closedcat}\index{closed category}\index{category!closed}\index{internal hom}
A symmetric monoidal category $(\C,\otimes,\tu)$ is \emph{closed} if, for each object $a \in \C$, the functor $- \otimes a \cn \C \to \C$ admits a right adjoint, which is called an \emph{internal hom}.  A right adjoint of $- \otimes a$ is denoted by $\Hom(a,-)$ or \label{notation:internal-hom}$[a,-]$.  If $\otimes$ is the Cartesian product $\times$, then we say $(\C, \times, [,])$ is \emph{Cartesian closed}.
\end{definition}

\begin{definition}[Diagrams]\label{def:diagramcat}\index{diagram category}\index{category!diagram}
For a small category $\cD$ and a category $\C$, the \emph{diagram category} \label{not:DC}$\DC$ has 
\begin{itemize}
\item functors $\cD \to \C$ as objects and
\item natural transformations between such functors as morphisms.
\end{itemize}
Identities and composition are given by, respectively, identity natural transformations and vertical composition.
\end{definition}

\begin{example}[Small Categories]\label{ex:cat}\index{category!of small categories}
The category $\Cat$\label{not:cat} of small categories and functors is Cartesian closed.
\begin{itemize}
\item The monoidal product is the Cartesian product $\times$.  
\item The monoidal unit $\boldone$ is the terminal category with only one object and its identity morphism.  To be exact, we fix the unique object of $\boldone$ to be $0$, but we sometimes denote it by $*$.  
\item The closed structure $[,]$ is given by diagram categories (\cref{def:diagramcat}). 
\end{itemize} 
Moreover, the category $\Cat$ is complete and cocomplete; see \cite[Section 1.4]{yau-involutive} for an elementary proof.
\end{example}

\subsection*{Monoids and Modules}

\begin{definition}[Monoids]\label{def:monoid}
Suppose $(\C, \otimes, \tu, \al, \la, \rho)$ is a monoidal category.
\begin{enumerate}
\item\label{def:monoid_i} A \emph{monoid}\index{monoid} in $\C$ is a triple\label{notation:monoid} $(a,\mu,\eta)$ consisting of 
\begin{itemize}
\item an object $a \in \C$,
\item a morphism $\mu \cn a \otimes a \to a$, called the \index{multiplication}\emph{multiplication}, and
\item a morphism $\eta \cn \tu \to a$, called the \index{unit}\emph{unit},
\end{itemize}
such that the following associativity and unity diagrams commute.
\begin{equation}\label{monoid_axioms}
\begin{tikzcd}[column sep=normal]
(a \otimes a) \otimes a \arrow{dd}[swap]{\mu \otimes 1_a} \rar{\alpha} 
& a \otimes (a \otimes a) \dar{1_a \otimes \mu}\\ 
& a \otimes a \dar{\mu}\\  
a  \otimes a \arrow{r}{\mu} & a
\end{tikzcd}\qquad
\begin{tikzcd}[column sep=normal]
\tu \otimes a \ar{d}[swap]{\eta \otimes 1_a} \ar{r}{\lambda_a} & a \ar[equal]{d}\\ 
a \otimes a \ar{r}{\mu} & a \ar[equal]{d}\\
a \otimes \tu \ar{u}{1_a \otimes \eta} \ar{r}{\rho_a} & a
\end{tikzcd}
\end{equation}
\item\label{def:monoid_ii} Suppose $(\C,\xi)$ is a symmetric monoidal category.  A \emph{commutative monoid}\dindex{commutative}{monoid} in $\C$ is a monoid $(a,\mu,\eta)$ such that the following diagram commutes.
\begin{equation}\label{com_monoid}
\begin{tikzpicture}[vcenter]
\def\h{2.5}
\draw[0cell]
(0,0) node (a1) {a \otimes a}
(a1)++(\h,0) node (a2) {a \otimes a}
(a1)++(.5*\h,-1) node (a3) {a}
;
\draw[1cell=.9]
(a1) edge node {\xi} (a2)
(a2) edge node[pos=.4] {\mu} (a3)
(a1) edge node[swap,pos=.4] {\mu} (a3)
;
\end{tikzpicture}
\end{equation}
\item\label{def:monoid_iii} A \emph{morphism} of (commutative) monoids
\[\begin{tikzcd}[column sep=large]
(a,\mu^a,\eta^a) \ar{r}{f} & (b,\mu^b,\eta^b)
\end{tikzcd}\] 
is a morphism $f \cn a \to b$ in $\C$ such that the following diagrams commute.
\[\begin{tikzcd}[column sep=large]
a \otimes a \dar[swap]{\mu^a} \rar{f\otimes f} & b \otimes b \dar{\mu^b}\\ 
a \rar{f} & b\end{tikzcd} \qquad
\begin{tikzcd}[column sep=large]
\tu \dar[equal] \rar{\eta^a} & a \dar{f}\\ \tu \rar{\eta^b} & b\end{tikzcd}\]
\end{enumerate}
This finishes the definition.
\end{definition}

\begin{definition}[Modules]\label{def:modules}
Suppose $(a,\mu,\eta)$ is a monoid in a monoidal category $\C$.
\begin{enumerate}
\item\label{def:modules_i} A \emph{right $a$-module}\index{module} is a pair $(x,\umu)$ consisting of
\begin{itemize}
\item an object $x \in \C$ and
\item a morphism $\umu \cn x \otimes a \to x$, called the \emph{right $a$-action},
\end{itemize}
such that the following associativity and unity diagrams commute.
\begin{equation}\label{module_axioms}
\begin{tikzcd}[column sep=normal]
(x \otimes a) \otimes a \arrow{dd}[swap]{\umu \otimes 1_a} \rar{\alpha} 
& x \otimes (a \otimes a) \dar{1_x \otimes \mu}\\ 
& x \otimes a \dar{\umu}\\  
x  \otimes a \arrow{r}{\umu} & x
\end{tikzcd}\qquad
\begin{tikzcd}[column sep=normal]
x \otimes a \ar{r}{\umu} & x \ar[equal]{d}\\
x \otimes \tu \ar{u}{1_x \otimes \eta} \ar{r}{\rho_x} & x
\end{tikzcd}
\end{equation}
\item\label{def:modules_ii} A \emph{left $a$-module} is defined analogously, with a \emph{left $a$-action} 
\[a \otimes x \fto{\umu} x\]
that makes the corresponding associativity and unity diagrams commute.
\item\label{def:modules_iii} If $(a,\mu,\eta)$ is a commutative monoid in a symmetric monoidal category $(\C,\xi)$, then each right $a$-module $(x,\umu)$ becomes a left $a$-module with left $a$-action given by the composite
\begin{equation}\label{xi_umu}
a \otimes x \fto[\iso]{\xi} x \otimes a \fto{\umu} x.
\end{equation}
\item\label{def:modules_iv} A \emph{morphism} of right $a$-modules
\[(x,\umu^x) \fto{f} (y,\umu^y)\]
is a morphism $f \cn x \to y$ in $\C$ such that the following diagram commutes.
\begin{equation}\label{module_morphism}
\begin{tikzpicture}[vcenter]
\def\v{-1.4}
\draw[0cell]
(0,0) node (a11) {x \otimes a}
(a11)++(2,0) node (a12) {x}
(a11)++(0,\v) node (a21) {y \otimes a}
(a12)++(0,\v) node (a22) {y}
;
\draw[1cell=.9]
(a11) edge node {\umu^x} (a12)
(a12) edge node {f} (a22)
(a11) edge node[swap] {f \otimes 1_a} (a21)
(a21) edge node {\umu^y} (a22)
;
\end{tikzpicture}
\end{equation}
\end{enumerate}
This finishes the definition.
\end{definition}

\subsection*{Monoidal Functors and Natural Transformations}

\begin{definition}\label{def:monoidalfunctor}\index{monoidal functor}\index{functor!monoidal}
For monoidal categories $\C$ and $\D$, a \emph{monoidal functor} 
\[(F, F^2, F^0) \cn \C \to \D\]
consists of
\begin{itemize}
\item a functor $F \cn \C \to \D$,
\item a morphism 
\begin{equation}\label{monoidal_unit}
\tu \fto{F^0} F\tu \inspace \D, 
\end{equation}
which is called the \index{unit constraint}\index{constraint!unit}\emph{unit constraint}, and
\item a natural transformation
\begin{equation}\label{monoidal_constraint}
\begin{tikzcd}[column sep=large]
Fa \otimes Fb \ar{r}{F^2_{a,b}} & F(a \otimes b)
\end{tikzcd}
\forspace a,b \in \C,
\end{equation}
which is called the \index{monoidal constraint}\index{constraint!monoidal}\emph{monoidal constraint}.
\end{itemize}
The above data are required to make the following unity and associativity diagrams commute for objects $a,b,c \in \C$.  
\begin{equation}\label{monoidalfunctorunity}
\begin{tikzcd}
\tu \otimes Fa \dar[swap]{F^0 \otimes 1_{Fa}} \rar{\lambda_{Fa}} & Fa \\ 
F\tu \otimes Fa \rar{F^2_{\tu,a}} & F(\tu \otimes a) \uar[swap]{F\lambda_a}
\end{tikzcd}
\qquad
\begin{tikzcd}
Fa \otimes \tu \dar[swap]{1_{Fa} \otimes F^0} \rar{\rho_{Fa}} & Fa \\ 
Fa \otimes F\tu \rar{F^2_{a,\tu}} & F(a \otimes \tu) \uar[swap]{F\rho_a}
\end{tikzcd}
\end{equation}
\begin{equation}\label{monoidalfunctorassoc}
\begin{tikzcd}[column sep=large]
\bigl(Fa \otimes Fb\bigr) \otimes Fc \rar{\alpha} \dar[swap]{F^2_{a,b} \otimes 1_{Fc}} 
& Fa \otimes \bigl(Fb \otimes Fc\bigr) \dar{1_{Fa} \otimes F^2_{b,c}}\\
F(a \otimes b) \otimes Fc \dar[swap]{F^2_{a \otimes b,c}} & Fa \otimes F(b \otimes c) \dar[d]{F^2_{a,b \otimes c}}\\
F\bigl((a \otimes b) \otimes c\bigr) \rar{F\alpha} & F\bigl(a \otimes (b \otimes c)\bigr)
\end{tikzcd}
\end{equation}
A monoidal functor is called
\begin{itemize}
\item \index{monoidal functor!strong}\index{strong monoidal functor}\emph{strong} if $F^0$ and $F^2$ are isomorphisms;
\item \index{monoidal functor!unital}\index{unital monoidal functor}\emph{unital} if $F^0$ is an isomorphism; 
\item \index{monoidal functor!strictly unital}\index{strictly unital monoidal functor}\emph{strictly unital} if $F^0$ is the identity morphism; and
\item \index{monoidal functor!strict}\index{strict monoidal functor}\emph{strict} if $F^0$ and $F^2$ are identities.
\end{itemize}
A monoidal functor $(F,F^2,F^0) \cn \C \to \D$ between braided monoidal categories is said to be \index{braided monoidal functor}\index{monoidal functor!braided}\emph{braided} if the following diagram commutes for $a,b \in \C$.  
\begin{equation}\label{monoidalfunctorbraiding}
\begin{tikzcd}[column sep=large]
Fa \otimes Fb \dar[swap]{F^2_{a,b}} \rar{\xi_{Fa,Fb}}[swap]{\cong} & Fb \otimes Fa \dar{F^2_{b,a}} \\ 
F(a \otimes b) \rar{F\xi_{a,b}}[swap]{\cong} & F(b \otimes a)
\end{tikzcd}
\end{equation}
A braided monoidal functor between symmetric monoidal categories is also called a \index{symmetric monoidal functor}\index{monoidal functor!symmetric}\emph{symmetric monoidal functor}.
\end{definition}

\begin{definition}\label{def:monoidalnattr}
For monoidal functors 
\[(F, F^2, F^0) \andspace (H,H^2,H^0) \cn \C \to \D\]
between monoidal categories $\C$ and $\D$, a \index{monoidal natural transformation}\index{natural transformation!monoidal}\emph{monoidal natural transformation} $\psi \cn F \to H$ is a natural transformation of the underlying functors such that the following unit constraint and monoidal constraint diagrams commute for $a,b \in \C$.
\begin{equation}\label{monnattr}
\begin{tikzpicture}[xscale=2,yscale=1.3,vcenter]
\draw[0cell=.9]
(0,0) node (x11) {\tu}
(x11)++(1,0) node (x12) {F\tu}
(x12)++(0,-1) node (x2) {H\tu}
(x12)++(1.2,0) node (y11) {Fa \otimes Fb}
(y11)++(1.5,0) node (y12) {Ha \otimes Hb}
(y11)++(0,-1) node (y21) {F(a \otimes b)}
(y12)++(0,-1) node (y22) {H(a \otimes b)}
;
\draw[1cell=.9]  
(x11) edge node {F^0} (x12)
(x11) edge node[swap] {H^0} (x2)
(x12) edge node {\psi_{\tu}} (x2)
(y11) edge node {\psi_a \otimes \psi_b} (y12)
(y12) edge node {H^2_{a,b}} (y22)
(y11) edge node[swap] {F^2_{a,b}} (y21)
(y21) edge node {\psi_{a \otimes b}} (y22)
;
\end{tikzpicture}
\end{equation}
A (monoidal) natural transformation is denoted diagrammatically using the following \index{2-cell!notation}2-cell notation.
\begin{equation}\label{twocellnotation}
\begin{tikzpicture}[xscale=2,yscale=1.7,baseline={(x1.base)}]
\draw[0cell=.9]
(0,0) node (x1) {\C}
(x1)++(1,0) node (x2) {\D}
;
\draw[1cell=.9]  
(x1) edge[bend left] node {F} (x2)
(x1) edge[bend right] node[swap] {H} (x2)
;
\draw[2cell]
node[between=x1 and x2 at .47, rotate=-90, 2label={above,\psi}] {\Rightarrow}
;
\end{tikzpicture}
\end{equation}
This finishes the definition.
\end{definition}

Discussion of pasting diagrams involving 2-cells can be found in \cite[Ch.\! 3]{johnson-yau}.

\section{Enriched Categories}
\label{sec:enrichedcat}

This section is a brief review of enriched category theory.  More detailed discussion can be found in \cite[Section 1.3]{johnson-yau}, \cite[Ch.\! 1]{cerberusIII}, and \cite{kelly-enriched}. 

\begin{definition}\label{def:enriched-category}
Suppose $(\V,\otimes,\tu,\alpha,\lambda,\rho)$ is a monoidal category (\cref{def:monoidalcategory}).  A \emph{$\V$-category} $(\C,\mcomp,\cone)$, which is also called a \index{category!enriched}\index{enriched category}\emph{category enriched in $\V$}, consists of
\begin{itemize}
\item a class \index{object!enriched category}$\Ob\C$ of \emph{objects};
\item for each pair of objects $a,b \in \C$, a \index{hom object}\index{object!hom}\emph{hom object} $\C(a,b) \in \V$, which is also called a \emph{hom $\V$-object} with domain $a$ and codomain $b$;
\item for objects $a,b,c \in \C$, a morphism\label{not:enrcomposition}
\begin{equation}\label{mcomp_abc}
\begin{tikzcd}[column sep=huge] 
\C(b,c) \otimes \C(a,b) \rar{\mcomp_{a,b,c}} & \C(a,c)
\end{tikzcd} \inspace \V
\end{equation}
called the \index{composition!enriched category}\emph{composition}; and
\item for each object $a \in \C$, a morphism\label{not:enridentity}
\begin{equation}\label{cone_a}
\begin{tikzcd}[column sep=large] 
\tu \rar{\cone_a} & \C(a,a)
\end{tikzcd} \inspace \V
\end{equation}
called the \emph{identity} of $a$, which is also denoted $\cone_a \cn a \to a$.
\end{itemize}
The following \index{associativity!enriched category}\emph{associativity diagram} and \index{unity!enriched category}\emph{unity diagram} are required to commute for objects $a,b,c,d \in \C$.
\begin{equation}\label{enriched-cat-associativity}
\begin{tikzcd}[cells={nodes={scale=.9}}]
  \big(\C(c,d) \otimes \C(b,c)\big) \otimes \C(a,b) \arrow{dd}[swap]{\mcomp \otimes 1}
  \rar{\alpha}
  & \C(c,d) \otimes \big(\C(b,c) \otimes \C(a,b)\big)  \dar{1 \otimes \mcomp} \\
& \C(c,d) \otimes \C(a,c) \dar{\mcomp}\\
\C(b,d) \otimes \C(a,b) \rar{\mcomp} & \C(a,d)  
\end{tikzcd}
\end{equation}
\begin{equation}\label{enriched-cat-unity}
\begin{tikzcd}[cells={nodes={scale=.9}}]
\tu \otimes \C(a,b) \dar[swap]{\cone_b \otimes 1} \rar{\lambda} & \C(a,b) \dar[equal] 
& \C(a,b) \otimes \tu \lar[swap]{\rho} \dar{1 \otimes \cone_a}\\
\C(b,b) \otimes \C(a,b) \rar{\mcomp} & \C(a,b) & \C(a,b) \otimes \C(a,a) \lar[swap]{\mcomp}
\end{tikzcd}
\end{equation}
This finishes the definition of a $\V$-category.
\end{definition}

\begin{definition}\label{def:enriched-functor}
For $\V$-categories $\C$ and $\D$, a \index{functor!enriched}\index{enriched!functor}\emph{$\V$-functor} $F \cn \C \to \D$ consists of
\begin{itemize}
\item an object assignment 
\[\Ob\C \fto{F} \Ob\D\] 
and
\item for each pair of objects $a,b \in \C$, a component morphism
\[\begin{tikzcd}[column sep=large]
\C(a,b) \rar{F_{a,b}} & \D\bigl(Fa,Fb\bigr)\end{tikzcd} \inspace \V,\]
which is usually abbreviated to $F$. 
\end{itemize}
The following two diagrams in $\V$ are required to commute for objects $a,b,c \in \C$.\index{composition!enriched functor}\index{identities!enriched functor} 
\begin{equation}\label{eq:enriched-composition}
\begin{tikzcd}[cells={nodes={scale=.9}}]
\C(b,c) \otimes \C(a,b) \rar{\mcomp} \dar[swap]{F \otimes F} & \C(a,c) \dar{F}\\
\D(Fb,Fc) \otimes \D(Fa,Fb) \rar{\mcomp} & \D(Fa,Fc)\end{tikzcd}
\qquad
\begin{tikzcd}[cells={nodes={scale=.9}}]
\tu \dar[equal] \rar{\cone_a} & \C(a,a) \dar{F}\\
\tu \rar{\cone_{Fa}} & \D(Fa,Fa)
\end{tikzcd}
\end{equation}
Identity and composition of $\V$-functors are defined separately on objects and component morphisms.
\end{definition}

\begin{definition}\label{def:enriched-natural-transformation}
Suppose $F,H \cn \C\to\D$ are $\V$-functors between $\V$-categories $\C$ and $\D$.  A \index{enriched!natural transformation}\index{natural transformation!enriched}\emph{$\V$-natural transformation} $\psi \cn F\to H$ consists of, for each object $a \in \C$, an $a$-component\index{component} $\V$-morphism
\[\begin{tikzcd}[column sep=large]
\tu \ar{r}{\psi_a} & \D(Fa,Ha),
\end{tikzcd}\]
which is also denoted $\psi_a \cn Fa \to Ha$, such that the following naturality diagram in $\V$ commutes for $a,b \in \C$.
\begin{equation}\label{enr_naturality}
\begin{tikzpicture}[x=40mm,y=12mm,vcenter]
\def\s{.85}
\draw[0cell=\s]
(0,0) node (a) {\C(a,b)}
(.25,1) node (b1) {\tu \otimes \C(a,b)}
(b1)++(1,0) node (c1) {\D(Fb,Hb) \otimes \D(Fa,Fb)}
(.25,-1) node (b2) {\C(a,b) \otimes \tu}
(b2)++(1,0) node (c2) {\D(Ha,Hb) \otimes \D(Fa,Ha)}
(c1)++(.25,-1) node (d) {\D(Fa,Hb)}
;
\draw[1cell=\s]
(a) edge node[pos=.3] {\la^\inv} node['] {\iso}(b1)
(a) edge node[swap,pos=.2] {\rho^\inv} node {\iso} (b2)
(b1) edge node {\psi_b \otimes F} (c1)
(b2) edge node {H \otimes \psi_a} (c2)
(c1) edge node[pos=.6] {\mcomp} (d)
(c2) edge node[swap,pos=.7] {\mcomp} (d)
;
\end{tikzpicture}
\end{equation}
Identity, vertical composition, and horizontal composition of $\V$-natural transformations are defined componentwise.
\end{definition}

\section{2-Categories}
\label{sec:twocategories}

This section is a brief review of 2-category theory.  Much more detailed discussion can be found in \cite{johnson-yau}.  

\begin{definition}\label{def:twocategory}\index{2-category}\index{category!2-}
A \emph{2-category} $\A$ consists of the following data and axioms.
\begin{description}
\item[Objects] It is equipped with a class $\A_0 = \Ob\A$ of \emph{objects}.
\item[1-cells] For each pair of objects $a,b \in\A_0$, it is equipped with a class $\A_1(a,b)$ of \index{1-cell}\emph{1-cells} from $a$ to $b$.  Such a 1-cell is denoted $a \to b$.
\item[2-cells] For 1-cells $f,f' \in \A_1(a,b)$, it is equipped with a set $\A_2(f,f')$ of \index{2-cell}\emph{2-cells} from $f$ to $f'$.  Such a 2-cell is denoted $f \to f'$.
\item[Identities]  It is equipped with
\begin{itemize}
\item an \emph{identity 1-cell} $1_a \in \A_1(a,a)$ for each object $a$ and
\item an \emph{identity 2-cell} $1_f \in \A_2(f,f)$ for each 1-cell $f \in \A_1(a,b)$.
\end{itemize}
\item[Compositions]
In the following compositions, $a$, $b$, and $c$ denote objects in $\A$.
\begin{itemize}
\item For 1-cells $f,f',f'' \in \A_1(a,b)$, it is equipped with the following \index{vertical composition!2-category}\emph{vertical composition} of 2-cells.\label{not:vcompiicell}
\[\begin{tikzcd}
\A_2(f',f'') \times \A_2(f,f') \rar{v} & \A_2(f,f'')
\end{tikzcd},\qquad v(\alpha',\alpha) = \alpha'\alpha.\] 
\item It is equipped with the following \index{horizontal composition!2-category}\emph{horizontal composition} of 1-cells.\label{not:hcompicell}
\[\begin{tikzcd}
\A_1(b,c) \times \A_1(a,b) \rar{h_1} & \A_1(a,c)
\end{tikzcd},\qquad h_1(g,f) = gf.\]
\item For 1-cells $f,f' \in \A_1(a,b)$ and $g,g' \in \A_1(b,c)$, it is equipped with the following \emph{horizontal composition} of 2-cells.\label{not:hcompiicell}
\[\begin{tikzcd}
\A_2(g,g') \times \A_2(f,f') \rar{h_2} & \B_2(gf,g'f')
\end{tikzcd},\quad h_2(\beta,\alpha) = \beta * \alpha.\]
\end{itemize}
\end{description}
The data above are required to satisfy axioms \pcref{twocat-i} through \pcref{twocat-iv} below.
\begin{enumerate}
\item\label{twocat-i} Vertical composition is associative and unital for identity 2-cells.
\item\label{twocat-ii} Horizontal composition preserves identity 2-cells and vertical composition.
\item\label{twocat-iii} Horizontal composition of 1-cells is associative and unital for identity 1-cells.
\item\label{twocat-iv} Horizontal composition of 2-cells is associative and unital for identity 2-cells of identity 1-cells.
\end{enumerate}
This finishes the definition of a 2-category.

Moreover, we define the following.
\begin{itemize}
\item For objects $a$ and $b$ in a 2-category $\A$, the \index{hom category!2-category}\emph{hom category} \label{not:homcat}$\A(a,b)$ is the category defined by the following data.
\begin{itemize}
\item Its objects are 1-cells $a \to b$.
\item Morphisms from a 1-cell $f \cn a \to b$ to a 1-cell $f' \cn a \to b$ are 2-cells $f \to f'$.
\item Its composition is vertical composition of 2-cells.
\item Its identities are identity 2-cells.
\end{itemize}
\item A 2-category is \index{locally small}\emph{locally small} if each hom category is a small category.
\item A 2-category is \emph{small} if it is locally small and has a set of objects.
\item The \index{underlying 1-category}\index{1-category!underlying}\emph{underlying 1-category} of a 2-category is the category defined by the following data.
\begin{itemize}
\item It has the same class of objects.
\item Its morphisms are 1-cells.
\item Its composition is horizontal composition of 1-cells.
\item Its identities are identity 1-cells.
\end{itemize}
\end{itemize}
We use the 2-cell notation in \cref{twocellnotation} for 2-cells in a 2-category.
\end{definition}

\begin{example}[Small Categories]\label{ex:catastwocategory}\index{2-category!of small categories}
$\Cat$ in \cref{ex:cat} is a 2-category, with natural transformations as 2-cells.  Horizontal and vertical compositions of 2-cells are those of natural transformations. 
\end{example}

\begin{example}[Small Enriched Categories]\label{ex:vcatastwocategory}\index{2-category!of small enriched categories}
Each monoidal category $\V$ has an associated 2-category $\Catv$ defined by the following data.
\begin{itemize}
\item Its objects are small $\V$-categories.
\item Its 1-cells are $\V$-functors.
\item Its 2-cells are $\V$-natural transformations.
\item Horizontal and vertical compositions are those of $\V$-functors and $\V$-natural transformations.
\end{itemize}
The 2-category $\Cat$ in \cref{ex:catastwocategory} is the specialization of $\Catv$ with $\V = (\Set, \times, *)$, the Cartesian closed category of sets and functions.
\end{example}

A $\Cat$-category is a 2-category.  The converse is almost true in the following sense.

\begin{proposition}\label{locallysmalltwocat}
A locally small 2-category is precisely a $\Cat$-category.
\end{proposition}

\begin{definition}\label{def:twofunctor}
For 2-categories $\A$ and $\B$, a \index{2-functor}\index{functor!2-}\emph{2-functor} $F \cn \A \to \B$ consists of
\begin{itemize}
\item an \emph{object assignment} 
\[\begin{tikzcd}[column sep=large]
\A_0 \ar{r}{F_0} & \B_0,
\end{tikzcd}\]
\item a \emph{1-cell assignment} 
\[\begin{tikzcd}[column sep=large]
\A_1(a,a') \ar{r}{F_1} & \B_1(F_0 a,F_0 a')
\end{tikzcd}\] 
for each pair of objects $a,a' \in \A$, and
\item a \emph{2-cell assignment}
\[\begin{tikzcd}[column sep=large]
\A_2(f,f') \ar{r}{F_2} & \B_2(F_1 f,F_1 f')
\end{tikzcd}\] 
for each pair of objects $a,a'$ and 1-cells $f,f' \in \A_1(a,a')$. 
\end{itemize}
The data above are required to satisfy axioms \cref{twofunctor-i,twofunctor-ii,twofunctor-iii} below, with $F_0$, $F_1$, and $F_2$ all abbreviated to $F$.
\begin{enumerate}
\item\label{twofunctor-i} The object and 1-cell assignments of $F$ constitute a functor between the underlying 1-categories of $\A$ and $\B$.
\item\label{twofunctor-ii} For each pair of objects $a,a' \in \A$, the 1-cell and 2-cell assignments of $F$ constitute a functor 
\[\begin{tikzcd}[column sep=large]
\A(a,a') \ar{r}{F} & \B(Fa,Fa')
\end{tikzcd}\] 
between hom categories.
\item\label{twofunctor-iii} $F$ preserves horizontal composition of 2-cells.
\end{enumerate}
This finishes the definition of a 2-functor.  Given another 2-functor $H \cn \B \to \C$, the \emph{composite 2-functor} 
\[\begin{tikzcd}[column sep=large]
\A \ar{r}{HF} & \C
\end{tikzcd}\]
is defined by separately composing the object, 1-cell, and 2-cell assignments.
\end{definition}

\begin{definition}\label{def:twonaturaltr}\index{2-natural!transformation}\index{natural transformation!2-}
For 2-functors $F,H \cn \A\to\B$ between 2-categories $\A$ and $\B$, a \emph{2-natural transformation} $\phi \cn F \to H$ consists of, for each object $a \in \A$, a \emph{component 1-cell} 
\[\begin{tikzcd}[column sep=large]
Fa \ar{r}{\phi_a} & Ha \inspace \B
\end{tikzcd}\]
such that the following two axioms hold:
\begin{description}
\item[1-cell naturality] For each 1-cell $f \cn a \to b$ in $\A$, the following two composite 1-cells in $\B(Fa,Gb)$ are equal.
\begin{equation}\label{onecellnaturality}
\begin{tikzcd}[column sep=large]
Fa \ar{d}[swap]{Ff} \ar{r}{\phi_a} & Ha \ar{d}{Hf}\\
Fb \ar{r}{\phi_b} & Hb
\end{tikzcd}
\end{equation}
\item[2-cell naturality] 
For each 2-cell $\theta \cn f \to g$ in $\A(a,b)$, the following two whiskered 2-cells are equal.
\begin{equation}\label{twocellnaturality}
\begin{tikzpicture}[xscale=3,yscale=1.4,vcenter]
\def\a{30} \def\s{.85}
\draw[0cell=\s]
(0,0) node (x11) {Fa}
(x11)++(1,0) node (x12) {Ha}
(x11)++(0,-1) node (x21) {Fb}
(x12)++(0,-1) node (x22) {Hb}
;
\draw[1cell=\s]  
(x11) edge node {\phi_a} (x12)
(x21) edge node {\phi_b} (x22)
(x11) edge[bend right] node[swap] {Ff} (x21)
(x11) edge[bend left] node {Fg} (x21)
(x12) edge[bend right] node[swap] {Hf} (x22)
(x12) edge[bend left] node {Hg} (x22)
;
\draw[2cell=.9]
node[between=x11 and x21 at .6, 2label={above,F\theta}] {\Rightarrow}
node[between=x12 and x22 at .6, 2label={above,H\theta}] {\Rightarrow}
;
\end{tikzpicture}
\end{equation}
This means that the following 2-cell equality holds in $\B(Fa,Gb)$.
\[G\theta * 1_{\phi_a} = 1_{\phi_b} * F\theta\]
\end{description}
This finishes the definition of a 2-natural transformation.  We use the 2-cell notation in \cref{twocellnotation} for 2-natural transformations.  A \index{2-natural!isomorphism}\index{natural isomorphism!2-}\emph{2-natural isomorphism} is a 2-natural transformation with each component 1-cell an isomorphism in the underlying 1-category.
\end{definition}

\begin{definition}\label{def:twonatcomposition}
Suppose $\phi \cn F \to H$ and $\varphi \cn H \to K$ are 2-natural transformations for 2-functors $F,H,K \cn \A \to \B$ between 2-categories $\A$ and $\B$.  The \emph{horizontal composite}\index{horizontal composition!2-natural transformation} 2-natural transformation
\[\begin{tikzcd}[column sep=large]
F \ar{r}{\varphi \phi} & K
\end{tikzcd}\]
is defined by, for each object $a \in \A$, the horizontal composite component 1-cell 
\[\begin{tikzcd}[column sep=large]
Fa \ar{r}{\phi_a} & Ha \ar{r}{\varphi_a} & Ka
\end{tikzcd}\]
in $\B(Fa,Ka)$.
\end{definition}

\begin{definition}\label{def:modification}
Suppose $\phi, \varphi \cn F \to H$ are 2-natural transformations for 2-functors $F,H \cn \A \to \B$ between 2-categories $\A$ and $\B$.  A \index{modification}\emph{modification} $\Phi \cn \phi \to \varphi$ consists of, for each object $a \in \A$, a \emph{component 2-cell}\label{not:componentiicell}
\[\begin{tikzcd}[column sep=large]
\phi_a \ar{r}{\Phi_a} & \varphi_a \inspace \B(Fa,Ha)
\end{tikzcd}\]
such that the following two whiskered 2-cells in $\B(Fa,Hb)$ are equal for each 1-cell $f \cn a \to b$ in $\A$.
\begin{equation}\label{modificationaxiom}
\begin{tikzpicture}[xscale=2.5,yscale=1.5,vcenter]
\def\a{35} \def\s{.85}
\draw[0cell=\s]
(0,0) node (x11) {Fa}
(x11)++(1,0) node (x12) {Ha}
(x11)++(0,-1) node (x21) {Fb}
(x12)++(0,-1) node (x22) {Hb}
;
\draw[1cell=\s]  
(x11) edge[bend left=\a] node[pos=.4] {\phi_a} (x12)
(x11) edge[bend right=\a] node[swap,pos=.6] {\varphi_a} (x12)
(x21) edge[bend left=\a] node[pos=.4] {\phi_b} (x22)
(x21) edge[bend right=\a] node[swap,pos=.6] {\varphi_b} (x22)
(x11) edge node[swap] {Ff} (x21)
(x12) edge node {Hf} (x22)
;
\draw[2cell=.9]
node[between=x11 and x12 at .44, rotate=-90, 2label={above,\Phi_a}] {\Rightarrow}
node[between=x21 and x22 at .44, rotate=-90, 2label={above,\Phi_b}] {\Rightarrow}
;
\end{tikzpicture}
\end{equation}
This means the equality
\[1_{Hf} * \Phi_a = \Phi_b * 1_{Ff}.\]
This is called the \index{modification!axiom}\emph{modification axiom}.
\end{definition}

\begin{definition}\label{def:modcomposition}
Suppose $F,H,K \cn \A \to \B$ are 2-functors between 2-categories $\A$ and $\B$.
\begin{itemize}
\item Suppose $\phi, \varphi, \psi \cn F \to H$ are 2-natural transformations, and $\Phi \cn \phi \to \varphi$ and $\Psi \cn \varphi \to \psi$ are modifications.  The \emph{vertical composite modification}\index{vertical composition!modification}
\begin{equation}\label{modificationvcomp}
\begin{tikzcd}[column sep=large]
\phi \ar{r}{\Psi \Phi} & \psi
\end{tikzcd}
\end{equation}
is defined by, for each object $a \in \A$, the vertical composite 2-cell
\[\begin{tikzcd}[column sep=large]
\phi_a \ar{r}{\Phi_a} & \varphi_a \ar{r}{\Psi_a} & \psi_a
\end{tikzcd}\]
in $\B(Fa,Ha)$.
\item Suppose $\Phi' \cn \phi' \to \varphi'$ is a modification for 2-natural transformations $\phi', \varphi' \cn H \to K$.  The \emph{horizontal composite modification}\index{horizontal composition!modification}
\begin{equation}\label{modificationhcomp}
\begin{tikzcd}[column sep=huge]
\phi' \phi \ar{r}{\Phi' * \Phi} & \varphi' \varphi
\end{tikzcd}
\end{equation}
is defined by, for each object $a \in \A$, the horizontal composite 2-cell
\[\begin{tikzcd}[column sep=huge]
\phi'_a \phi_a \ar{r}{\Phi'_a * \Phi_a} & \varphi'_a \varphi_a
\end{tikzcd}\]
in $\B(Fa,Ka)$.\defmark
\end{itemize}
\end{definition}

The following definition is from \cite[Definition 6.2.10]{johnson-yau}.

\begin{definition}\label{def:twoequivalence}
A 2-functor $F \cn \A \to \B$ is a \emph{2-equivalence}\index{2-equivalence} if there is a 2-functor $H \cn \B \to \A$ together with 2-natural isomorphisms
\[1_{\A} \iso HF \andspace FH \iso 1_{\B}.\]
In this case, we call $H$ a \emph{2-inverse} of $F$.
\end{definition}

There is a local characterization of an equivalence of categories as a functor that is (i) essentially surjective on objects and (ii) fully faithful (that is, bijective) on morphisms \cite[IV.4 Theorem 1]{maclane}.  The following 2-categorical analogue is proved in detail in \cite[Theorem 7.5.8]{johnson-yau}.  Another proof appears in \cite[Theorem 3.70]{yau-multgro}, which is obtained from a more general characterization of $\Cat$-multiequivalences.

\begin{theorem}\label{thm:twoequivalences}
A 2-functor $F \cn \A \to \B$ is a 2-equivalence if and only if $F$ is
\begin{itemize}
\item essentially surjective on objects,
\item fully faithful on 1-cells, and
\item fully faithful on 2-cells.
\end{itemize}
\end{theorem}

\section{Enriched Monoidal Categories}
\label{sec:enrmonoidalcat}

This section reviews enriched monoidal categories.  Detailed discussion of this material can be found in \cite[Chapters 1--3]{cerberusIII}.

\subsection*{Tensor Products of Enriched Categories}

To define enriched monoidal categories, first we need some preliminary definitions.  Suppose $(\V,\otimes,\tu,\alpha,\lambda,\rho,\xi)$ is a braided monoidal category (\cref{def:braidedmoncat}). 

\begin{definition}\label{definition:vtensor-0}
For $\V$-categories $\C$ and $\D$ (\cref{def:enriched-category}), the \emph{tensor product}\index{tensor product!enriched category}\index{enriched category!tensor product} $\C \otimes \D$ is the $\V$-category defined as follows.
\begin{description}
\item[Objects] Its class of objects is 
\[\Ob(\C \otimes \D) = \Ob\C \times \Ob\D.\]  
An object in $\C \otimes \D$ is denoted by $c \otimes d$ or $(c;d)$ for $c \in \C$ and $d \in \D$.
\item[Hom objects] For objects $c \otimes d$ and $c' \otimes d'$ in $\C \otimes \D$, the hom $\V$-object is given by
\[(\C \otimes \D)(c \otimes d, c' \otimes d') = \C(c,c') \otimes \D(d,d').\]
\item[Composition] For objects $c \otimes d$, $c'\otimes d'$, and $c''\otimes d''$ in $\C \otimes \D$, the composition is given by the following composite in $\V$, where $\xi$ is the unique coherence isomorphism that permutes the middle two factors using the associativity isomorphism and braiding in $\V$.
\[\begin{tikzpicture}[x=23mm,y=13mm]
\draw[0cell=.8]
(0,0) node (a) {\big(\C(c',c'')\otimes\D(d',d'')\big) \otimes \big(\C(c,c')\otimes\D(d,d')\big)}
(0,-1) node (b) {\big(\C(c',c'')\otimes\C(c,c')\big) \otimes \big(\D(d',d'')\otimes\D(d,d')\big)}
(2,0) node (c) {\C(c,c'') \otimes \D(d,d'')}
;
\draw[1cell=.8]
(a) edge node[swap] {\xi} (b)
(b) [rounded corners=2pt] -| node[pos=.25] {\mcomp \otimes \mcomp} (c);
\end{tikzpicture}\]
\item[Identities] The identity of an object $c \otimes d$ is the following composite.
\[\tu \fto[\iso]{\la^\inv} \tu \otimes \tu \fto{\cone_c \otimes \cone_d} \C(c,c) \otimes \D(d,d)\]
\end{description}
This finishes the definition of the $\V$-category $\C \otimes \D$.  Extending the tensor product $\otimes$ to $\V$-functors and $\V$-natural transformations componentwise, it becomes a 2-functor
\[\begin{tikzcd}[column sep=large]
\VCat \times \VCat \ar{r}{\otimes} & \VCat
\end{tikzcd}\]
on the 2-category $\VCat$ of small $\V$-categories (\cref{ex:vcatastwocategory}).
\end{definition}

\begin{definition}\label{definition:unit-vcat}
The \emph{unit $\V$-category}\index{unit enriched category}\index{enriched category!unit}\index{tensor product!enriched category!unit}\index{enriched!tensor product!unit} $\vtensorunit$ is the $\V$-category with one object $*$ and unique hom $\V$-object given by the monoidal unit $\tensorunit \in \V$.  The composition is given by the left unit isomorphism $\lambda_{\tu}$.  The identity of $*$ is given by the identity morphism $1_{\tensorunit}$.
\end{definition}

\begin{definition}\label{definition:vtensor-unitors}
We define the \emph{left unitor}\index{left unitor!enriched tensor product}\index{tensor product!enriched category!left unitor}\index{enriched category!tensor product!left unitor} $\ell^\otimes$ and the \emph{right unitor}\index{right unitor!enriched tensor product}\index{tensor product!enriched category!right unitor}\index{enriched category!tensor product!right unitor} $r^\otimes$ as the 2-natural isomorphisms (\cref{def:twonaturaltr})
  \[
  \begin{tikzpicture}[x=25mm,y=13mm,vcenter]
    \draw[0cell=.85]
    (0,1) node (a) {\VCat^2}
    (-.5,0) node (b) {\VCat}
    (.5,0) node (b') {\VCat}
    ;
    \draw[1cell=.9] 
    (b) edge node {\vtensorunit \times 1} (a)
    (a) edge node {\otimes} (b')
    (b) edge['] node (z) {1} (b')
    ;
    \draw[2cell] 
    node[between=a and z at .55, shift={(-.05,0)}, rotate=-90, 2label={above,\ell^\otimes}] {\Rightarrow}
    ;
  \end{tikzpicture}
  \qquad \text{and} \qquad
  \begin{tikzpicture}[x=25mm,y=13mm,vcenter]
    \draw[0cell=.85] 
    (0,1) node (a) {\VCat^2}
    (-.5,0) node (b) {\VCat}
    (.5,0) node (b') {\VCat}
    ;
    \draw[1cell=.9] 
    (b) edge node {1 \times \vtensorunit} (a)
    (a) edge node {\otimes} (b')
    (b) edge['] node (z) {1} (b')
    ;
    \draw[2cell] 
    node[between=a and z at .55, shift={(-.05,0)}, rotate=-90, 2label={above,r^\otimes}] {\Rightarrow}
    ;
  \end{tikzpicture}
  \]
whose components at a $\V$-category $\C$ are the $\V$-functors (\cref{def:enriched-functor})
  \[
  \vtensorunit \otimes \C \fto[\iso]{\ell^\otimes_{\C}} \C \xleftarrow[\iso]{r^\otimes_{\C}} \C
  \otimes \vtensorunit
  \]
given
\begin{itemize}
\item on objects by 
\[\ell^\otimes_{\C}(*,a) = a = r^\otimes_{\C}(a,*)\] 
and
\item on hom $\V$-objects by the left and right unit isomorphisms
  \[
  \tensorunit \otimes \C(a,b) \fto[\iso]{\la} \C(a,b) \xleftarrow[\iso]{\rho}
  \C(a,b) \otimes \tensorunit \inspace \V
  \]
for objects $a,b \in \C$.\defmark
\end{itemize} 
\end{definition}

\begin{definition}\label{definition:vtensor-assoc}
We define the 2-natural isomorphism $a^\otimes$, called the \emph{associator}\index{associator!enriched tensor product}\index{tensor product!enriched category!associator}\index{enriched!tensor product!associator},
  \[
  \begin{tikzpicture}[x=25mm,y=15mm]
    \draw[0cell=.85] 
    (0,0) node (a) {\VCat^3}
    (1,0) node (b) {\VCat^2}
    (0,-1) node (c) {\VCat^2}
    (1,-1) node (d) {\VCat}
    ;
    \draw[1cell=.9] 
    (a) edge node {\otimes \times 1} (b)
    (c) edge['] node {\otimes} (d)
    (a) edge['] node {1 \times \otimes} (c)
    (b) edge node {\otimes} (d)
    ;
    \draw[2cell] 
    node[between=a and d at .5, rotate=-90, 2label={above,a^\otimes}] {\Rightarrow}
    ;
  \end{tikzpicture}
  \]
whose component at small $\V$-categories $\C$, $\D$, and $\E$, is the $\V$-functor
  \[\begin{tikzcd}[column sep=huge]
  (\C \otimes \D) \otimes \E \ar{r}{a^\otimes_{\C,\D,\E}}[swap]{\iso} & \C \otimes (\D \otimes \E)
  \end{tikzcd}\]
given as follows.
\begin{itemize}
\item On objects it is given by 
\[a^\otimes_{\C,\D,\E}((c,d),e) = (c, (d,e))\]
for objects $c \in \C$, $d \in \D$, and $e \in \E$.
\item On hom $\V$-objects it is given by the associativity isomorphism
  \[\begin{tikzcd}
  \big(\C(c,c') \otimes \D(d,d')\big) \otimes \E(e,e') \ar{r}{\al}[swap]{\iso} & \C(c,c') \otimes \big(\D(d,d') \otimes \E(e,e')\big)
  \end{tikzcd}\]
in $\V$ for objects $c,c' \in \C$, $d,d' \in \D$, and $e,e'\in \E$.\defmark
\end{itemize} 
\end{definition}

In the previous definitions in this section, $\V$ is a braided monoidal category.  In the following definition, we need $\V$ to be a symmetric monoidal category (\cref{def:braidedmoncat}).

\begin{definition}\label{definition:vtensor-beta}
Suppose $(\V,\xi)$ is a symmetric monoidal category.  With $\twist$ denoting the 2-functor that permutes the two arguments, we define the 2-natural isomorphism $\beta^\otimes$, called the \emph{braiding}\index{braiding!enriched tensor product}\index{tensor product!enriched category!braiding}\index{enriched!tensor product!braiding},
  \[
  \begin{tikzpicture}[xscale=3, yscale=1.1]
    \def\v{-1} \def\h{1} \def\m{1} \def\q{15}
    \draw[0cell=.85] 
    (0,0) node (x11) {\VCat^2}
    ($(x11)+(\h,0)$) node (x12) {\VCat}
    ($(x11)+(\h/2,\v)$) node (x2) {\VCat^2}
    ;
    \draw[1cell] 
    (x11) edge node (i) {\otimes} (x12)
    (x11) edge[bend right=\q] node[swap,pos=.5] (a) {\twist} (x2)
    (x2) edge[bend right=\q] node[swap,pos=.5] {\otimes} (x12)
    ;
    \draw[2cell] 
    node[between=i and x2 at .5, shift={(-.1,0)}, rotate=-90, 2label={above,\beta^\otimes}] {\Rightarrow}
    ;
  \end{tikzpicture}
  \]
whose component at small $\V$-categories $\C$ and $\D$ is the $\V$-functor
\[\begin{tikzcd}[column sep=large]
\C \otimes \D \ar{r}{\beta^\otimes_{\C,\D}}[swap]{\iso} & \D \otimes \C
\end{tikzcd}\] 
defined as follows.  
\begin{itemize}
\item On objects it is given by
\[\beta^\otimes_{\C,\D}(c,d) = (d,c)\]
for objects $c \in \C$ and $d \in \D$.
\item On hom $\V$-objects it is given by the braiding in $\V$ 
\[\begin{tikzcd}[column sep=large]
\C(c,c') \otimes \D(d,d') \ar{r}{\xi}[swap]{\iso} & \D(d,d') \otimes \C(c,c')
\end{tikzcd}\]
for objects $c,c' \in \C$ and $d,d' \in \D$.\defmark
\end{itemize} 
\end{definition}

The following observation is proved in \cite[1.3.35]{cerberusIII}.

\begin{theorem}\label{theorem:vcat-mon}
For a braided monoidal category $\V$, the data\index{enriched category!2-category}\index{2-category!of small enriched categories} 
\[\big( \VCat,\otimes,\vtensorunit,a^\otimes,\ell^\otimes,r^\otimes \big)\]
constitute a monoidal category.\index{monoidal category!of small enriched categories}\index{enriched category!monoidal}\index{tensor product!enriched category!monoidal}  Moreover, if $\V$ is a symmetric monoidal category, then $(\VCat,\beta^\otimes)$ is a symmetric monoidal category.
\end{theorem}

\subsection*{Enriched Symmetric Monoidal Categories}

\begin{definition}\label{definition:monoidal-vcat}
For a braided monoidal category $\V$, a \emph{monoidal $\V$-category}\index{enriched category!monoidal}\index{enriched!monoidal - category}\index{category!enriched!monoidal}\index{monoidal category!enriched} is a sextuple
\[\big( \C,\vmtimes,\vmunitbox,a^\vmtimes,\ell^\vmtimes,r^\vmtimes \big)\]
consisting of the following data.
\begin{description}
\item[Base $\V$-category] $\C$ is a $\V$-category (\cref{def:enriched-category}), called the \emph{base $\V$-category}.\index{base enriched category!monoidal enriched category}
\item[Monoidal composition] It is the $\V$-functor\index{monoidal composition!monoidal enriched category}\index{composition!monoidal enriched category} (\cref{def:enriched-functor})
\begin{equation}\label{monvcat_moncomp}
\C \otimes \C \fto{\vmtimes} \C.
\end{equation}
\item[Monoidal identity] It is the $\V$-functor\index{monoidal identity!monoidal enriched category}
\begin{equation}\label{monvcat_monid}
\vtensorunit \fto{\vmunitbox} \C.
\end{equation}
The image of the object $* \in \vtensorunit$ is called the \emph{identity object}\index{identity object!monoidal enriched category} and also denoted $\vmunitbox$.
\item[Monoidal unitors] In the rest of this definition, we abbreviate the tensor product of $\V$-categories (\cref{definition:vtensor-0}) to juxtaposition, so $\C^2$ means $\C \otimes \C$.  The \emph{left monoidal unitor}\index{left monoidal unitor!monoidal enriched category}\index{unitor!monoidal enriched category} $\ell^\vmtimes$ and the \emph{right monoidal unitor}\index{right monoidal unitor!monoidal enriched category} $r^\vmtimes$ are $\V$-natural isomorphisms (\cref{def:enriched-natural-transformation})
    \begin{equation}\label{eq:monoidal-unitors}
    \begin{tikzpicture}[x=30mm,y=13mm,vcenter]
    \def\s{.9} \def\t{.8} \def\h{.1}
      \draw[0cell=\s] 
      (0,0) node (a) {\C}
      (\h,1) node (b) {\vtensorunit\C}
      (1-\h,1) node (c) {\C^2}
      (1,0) node (d) {\C}
      ;
      \draw[1cell=\t] 
      (a) edge node[pos=.25] {(\ell^\otimes)^{-1}} (b)
      (b) edge node (z) {\vmunitbox 1_{\C}} (c)
      (c) edge node[pos=.6] {\vmtimes} (d)
      (a) edge['] node (w) {1_{\C}} (d)
      ;
      \draw[2cell] 
      node[between=z and w at .5, rotate=-90, 2label={above,\ell^\vmtimes}] {\Rightarrow}
      ;
    \end{tikzpicture}
    \qquad
    \begin{tikzpicture}[x=30mm,y=13mm,vcenter]
    \def\s{.9} \def\t{.8} \def\h{.1}
      \draw[0cell=\s] 
      (0,0) node (a) {\C}
      (\h,1) node (b) {\C\vtensorunit}
      (1-\h,1) node (c) {\C^2}
      (1,0) node (d) {\C}
      ;
      \draw[1cell=\t] 
      (a) edge node[pos=.25] {(r^{\otimes})^{-1}} (b)
      (b) edge node (z) {1_{\C} \vmunitbox} (c)
      (c) edge node[pos=.6] {\vmtimes} (d)
      (a) edge['] node (w) {1_{\C}} (d)
      ;
      \draw[2cell] 
      node[between=z and w at .5, rotate=-90, 2label={above,\,r^\vmtimes}] {\Rightarrow}
      ;
    \end{tikzpicture}
    \end{equation}
with components at an object $c \in \C$ given by, respectively,
    \[
    \tensorunit \fto{\ell^\vmtimes_c} \C(\vmunitbox \vmtimes c, c) \andspace
    \tensorunit \fto{r^\vmtimes_c} \C(c \vmtimes \vmunitbox, c).
    \]
\item[Monoidal associator] It is the $\V$-natural isomorphism\index{associator!monoidal enriched category}\index{monoidal associator!monoidal enriched category} 
    \begin{equation}\label{eq:monoidal-assoc}
    \begin{tikzpicture}[x=20mm,y=20mm,vcenter]
    \def\s{.9}
      \draw[0cell=\s] 
      (0,0) node (a) {(\C^{2})\C}
      (225:.707) node (a') {\C(\C^{2})}
      (a)++(1,0) node (b) {\C^{2}}
      (a')++(0,-.6) node (c) {\C^{2}}
      (b)++(0,-1.1) node (d) {\C}
      ;
      \draw[1cell=\s] 
      (a) edge node {\vmtimes 1_\C} (b)
      (c) edge['] node {\vmtimes} (d)
      (a') edge['] node {1_\C \vmtimes} (c)
      (b) edge node {\vmtimes} (d)
      (a) edge['] node {a^\otimes} (a')
      ;
      \draw[2cell] 
      node[between=b and c at .5, rotate=225, 2labelalt={below,a^\vmtimes}] {\Rightarrow}
      ;
    \end{tikzpicture}
    \end{equation}
whose component at objects $x,y,z \in \C$ is a morphism in $\V$
    \[
    \tensorunit \fto{a^\vmtimes_{x,y,z}}
    \C\big((x \vmtimes y) \vmtimes z \scs x \vmtimes (y \vmtimes z)\big).
    \]
  \end{description}
The above data are subject to the following two axioms.
\begin{description}
\item[Unity axiom]\index{unity!monoidal enriched category}
The following two \index{unity!monoidal enriched category}\emph{middle
      unity pasting diagrams}\index{middle unity diagram!monoidal enriched category} have equal composites.
    \begin{equation}\label{vmonoidal-unit}
      \begin{tikzpicture}[x=20mm,y=15mm,baseline={(eq.base)}]
        \def\w{3} 
        \def\h{2} 
        \def\m{.5} 
        \def\s{.75}
        \draw[font=\Large] (0,0) node[rotate=90] (eq) {=}; 
        \newcommand{\boundary}{
          \draw[0cell=\s]
          (0,0) node (a) {(\C^2)\C}
          (1.25,0) node (b) {\C^2}
          (a)++(0,-1.5) node (b') {\C^2}
          (b')++(1.25,0) node (c) {\C}
          (a)++(-1,0) node (z) {(\C\vtensorunit)\C}
          (z)++(-1,-.75) node (w) {\C^2}
          ; 
          \draw[1cell=\s] 
          (a) edge node {\vmtimes 1} (b)
          (b) edge node {\vmtimes} (c)
          (b') edge node {\vmtimes} (c)
          (w) edge node {(r^{\otimes})^{-1} 1} (z)
          (z) edge node {(1\vmunitbox)1} (a)
          (w) edge[',out=-35,in=180] node {1^2} (b')
          ;
        }
        \begin{scope}[shift={(1.5,1.8)}]
          \boundary
          \draw[0cell=\s]
          (a)++(0,-.75) node (a') {\C(\C^2)}
          (a')++(-1,0) node (z') {\C(\vtensorunit\C)}
          ;
          \draw[1cell=\s] 
          (a) edge node {a^\otimes} (a')
          (a') edge node {1 \vmtimes} (b')          
          (w) edge['] node[pos=.6] {1 (\ell^{\otimes})^{-1}} (z')
          (z') edge node {1(\vmunitbox 1)} (a')
          (z) edge node {a^\otimes} (z')
          ;
          \draw[2cell] 
          node[between=b and b' at .5, rotate=225, 2labelalt={below,a^\vmtimes}] {\Rightarrow}
          node[between=z' and b' at .5, rotate=225, 2labelalt={below,1\ell^\vmtimes}] {\Rightarrow}
          (z') ++(135:.4) node {\divideontimes}
          ;
        \end{scope}
        \begin{scope}[shift={(1.5,-.3)}]
          \boundary
          \draw[1cell=\s] 
          (w) edge[',out=0,in=210] node (I) {1^2} (b)
          ;
          \draw[2cell] 
          node[between=I and a at .5, shift={(-.3,0)}, rotate=225, 2labelalt={below,r^\vmtimes 1}] {\Rightarrow}
          ;
        \end{scope}
      \end{tikzpicture}
    \end{equation}
    In the first pasting diagram in \cref{vmonoidal-unit}, the unlabeled region commutes by
    naturality of $a^\otimes$.  The
    region labeled $\divideontimes$ commutes by
    the middle unity for $\ell^\otimes$ and $r^\otimes$.
  \item[Pentagon axiom]\index{pentagon axiom!monoidal enriched category}
    The following two \emph{pentagon pasting
      diagrams} have equal composites.
    
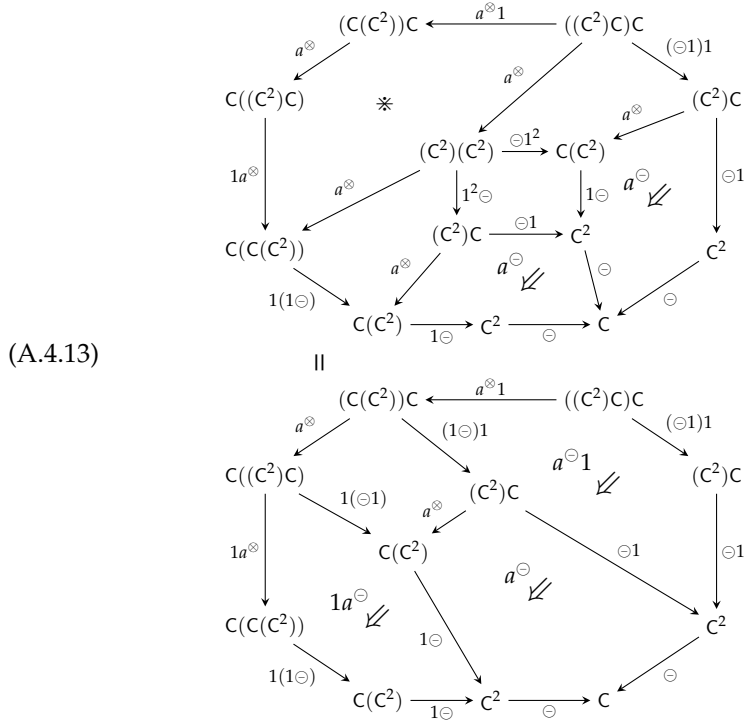
\begin{equation}\label{vmonoidal-pentagon-axiom}
      \begin{tikzpicture}[x=30mm,y=20mm,baseline={(eq.base)}]
        \def\w{2} 
        \def\hmarg{.25}
        \def\m{.75} 
        \def\s{.75}
        \draw[font=\Large] (0,0) node[rotate=90] (eq) {=}; 
        \newcommand{\boundary}{
          \draw[0cell=\s] 
          (.5,0) node (b) {(\C(\C^2))\C}
          (b)++(1,0) node (a) {((\C^2)\C)\C}
          (a)++(.5,-.5) node (z) {(\C^2)\C}
          (z)++(0,-1) node (y) {\C^2}
          (y)++(-.5,-.5) node (x) {\C}
          (0,-.5) node (c) {\C((\C^2)\C)}
          (c)++(0,-1) node (d) {\C(\C(\C^2))}
          (d)++(.5,-.5) node (e) {\C(\C^2)}
          (e)++(.5,0) node (f) {\C^2}
          ; 
          \draw[1cell=\s] 
          (a) edge['] node {a^\otimes1} (b)
          (b) edge['] node {a^\otimes} (c)
          (c) edge['] node {1a^\otimes} (d)
          (d) edge['] node {1(1\vmtimes)} (e)
          (e) edge['] node {1\vmtimes} (f)
          (f) edge['] node {\vmtimes} (x)
          (a) edge node {(\vmtimes 1) 1} (z)
          (z) edge node {\vmtimes 1} (y)
          (y) edge node {\vmtimes} (x)
          ;
        }
        \begin{scope}[shift={(-.25,-\hmarg)}]
          \boundary
          \draw[0cell=\s]
          (b)++(-50:.8) node (b') {(\C^2)\C}
          (c)++(-40:.8) node (c') {\C(\C^2)}
          ;
          \draw[1cell=\s] 
          (b) edge node {(1\vmtimes)1} (b')
          (c) edge node {1(\vmtimes 1)} (c')
          (b') edge['] node {a^\otimes} (c')
          (b') edge node {\vmtimes 1} (y)
          (c') edge['] node {1\vmtimes} (f)
          ;
          \draw[2cell] 
          node[between=z and b' at .5, rotate=225, 2label={below,a^\vmtimes 1}] {\Rightarrow}
          node[between=c' and y at .5, shift={(-.1,0)}, rotate=225, 2labelalt={below,a^\vmtimes}] {\Rightarrow}
          node[between=c' and e at .5, shift={(135:.1)}, rotate=225, 2labelalt={below,1 a^\vmtimes}] {\Rightarrow}

          ;
        \end{scope}
        \begin{scope}[shift={(-.25,2+\hmarg)}]
          \boundary
          \draw[0cell=\s]
          (-45:1.2) node (m) {(\C^2)(\C^2)}
          (m)++(0:.55) node (mr) {\C(\C^2)}
          (m)++(-90:.55) node (ml) {(\C^2)\C}
          (ml)++(0:.55) node (n) {\C^2}
          ;
          \draw[1cell=\s] 
          (m) edge node {1^2 \vmtimes} (ml)
          (m) edge node {\vmtimes 1^2} (mr)
          (ml) edge node {\vmtimes 1} (n)
          (mr) edge node {1 \vmtimes} (n)
          (ml) edge['] node {a^\otimes} (e)
          (z) edge['] node {a^\otimes} (mr)
          (a) edge['] node {a^\otimes} (m)
          (m) edge['] node {a^\otimes} (d)
          (n) edge node {\vmtimes} (x)
          ;
          \draw[2cell] 
          node[between=mr and y at .5, shift={(45:.05)}, rotate=225, 2labelalt={below,a^\vmtimes}] {\Rightarrow}
          node[between=ml and x at .5, rotate=225, 2labelalt={below,a^\vmtimes}] {\Rightarrow}
          (m)++(135:.45) node {\divideontimes}
          ;
        \end{scope}
      \end{tikzpicture}
    \end{equation}
    The central square in the first pasting diagram in \cref{vmonoidal-pentagon-axiom} commutes by
    2-functoriality of $\otimes$ in each variable.  The other unlabeled regions in the two pasting diagrams commute by 2-naturality
    of $a^\otimes$.  The pentagon labeled
    $\divideontimes$ commutes by the pentagon axiom for $a^\otimes$.
  \end{description}
  This finishes the definition of a monoidal $\V$-category.  
\end{definition}

The following notion of a mate is needed to define braided monoidal $\V$-categories.  This is similar to the mates of a pentagonator\index{pentagonator}\index{pentagonator!mate}\index{mate!pentagonator} in the bicategorical context; see \cite[12.1.4]{johnson-yau} for details.

\begin{definition}\label{eq:a-vmtimes-inv-mate}\index{mate!associator}\index{associator!mate}
For a braided monoidal category $\V$ and a monoidal $\V$-category $\C$ (\cref{definition:monoidal-vcat}), we denote by $a^\vmtimes_1$ the mate of $a^\vmtimes$ in \cref{eq:monoidal-assoc} given by the inverse of $a^\otimes$, as in the diagram below.
  
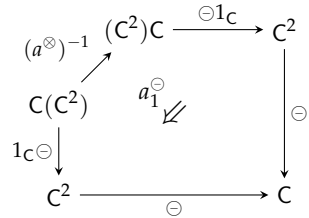
\[
  \begin{tikzpicture}[x=20mm,y=20mm]
  \def\s{.9}
    \draw[0cell=\s] 
    (0,0) node (a) {(\C^{2})\C}
    (225:.707) node (a') {\C(\C^{2})}
    (a)++(1,0) node (b) {\C^{2}}
    (a')++(0,-.6) node (c) {\C^{2}}
    (b)++(0,-1.1) node (d) {\C}
    ;
    \draw[1cell=\s] 
    (a) edge node {\vmtimes 1_\C} (b)
    (c) edge['] node {\vmtimes} (d)
      (a') edge['] node {1_\C \vmtimes} (c)
    (b) edge node {\vmtimes} (d)
    (a') edge node {(a^{\otimes})^{-1}} (a)
    ;
    \draw[2cell] 
    node[between=b and c at .5, rotate=225, 2labelalt={below,a^\vmtimes_1}] {\Rightarrow}
    ;
  \end{tikzpicture}
  \]
This finishes the definition.
\end{definition}

\begin{definition}\label{definition:braided-monoidal-vcat}
Suppose $\V$ is a symmetric monoidal category.  A \emph{braided
    monoidal $\V$-category}\index{braided monoidal category!enriched}\index{enriched!braided monoidal category}\index{enriched category!braided monoidal} is a pair $(\C,\beta^\vmtimes)$ consisting of
  \begin{itemize}
  \item a monoidal $\V$-category $\C$ (\cref{definition:monoidal-vcat}) and
  \item a $\V$-natural isomorphism $\beta^\vmtimes$, called the \emph{braiding}\index{braiding!braided monoidal enriched category} of $\C$, as follows.
\begin{equation}\label{bmvcat_braiding}
    \begin{tikzpicture}[xscale=3, yscale=1.1,vcenter]
      \def\v{-1} \def\h{1} \def\m{1} \def\q{15}
      \draw[0cell] 
      (0,0) node (x11) {\C^2}
      ($(x11)+(\h,0)$) node (x12) {\C}
      ($(x11)+(\h/2,\v)$) node (x2) {\C^2}
      ;
      \draw[1cell] 
      (x11) edge node (i) {\vmtimes} (x12)
      (x11) edge[bend right=\q] node[swap,pos=.5] (a) {\beta^\otimes} (x2)
      (x2) edge[bend right=\q] node[swap,pos=.5] {\vmtimes} (x12)
      ;
      \draw[2cell] 
      node[between=i and x2 at .5, shift={(-.1,0)}, rotate=-90, 2label={above,\beta^\vmtimes}] {\Rightarrow}
;
    \end{tikzpicture}
\end{equation}
\end{itemize}
The above data are subject to the following two axioms.
\begin{description}
\item[Left hexagon axiom] The following two
  \emph{left hexagon pasting diagrams} have equal composites.
  \begin{equation}\label{hexagon-bvmL}
    \begin{tikzpicture}[x=17mm,y=17mm,vcenter]
      \def\w{2.1} 
      \def\h{2} 
      \def\m{.5} 
      \def\s{.7}
      \draw[font=\Large] (\w+\m,0) node (eq) {=}; 
      \newcommand{\boundary}{
        \draw[0cell=\s] 
        (.5,0) node (b) {\C(\C^2)}
        (0,-.5) node (c) {(\C^2)\C}
        (b)++(0:1) node (a) {(\C^2)\C}
        (c)++(-90:1) node (d) {\C(\C^2)}
        (a)++(-45:1) node (y) {\C^2}
        (d)++(-45:1) node (e) {\C^2}
        (-45:3.12) node (z) {\C}
        ; 
        \draw[1cell=\s] 
        (a) edge['] node {a^\otimes} (b)
        (b) edge['] node {\beta^\otimes} (c)
        (c) edge['] node {a^\otimes} (d)
        (d) edge['] node {1\vmtimes} (e)
        (e) edge['] node {\vmtimes} (z)
        (a) edge node {\vmtimes 1} (y)
        (y) edge node {\vmtimes} (z)
        ;
      }
      \begin{scope}[shift={(0,.2*\h)}]
        \boundary
        \draw[0cell=\s]
        (a)++(-100:.8) node (b') {(\C^2)\C}
        (d)++(10:.8) node (c') {\C(\C^2)}
        ;
        \draw[1cell=\s] 
        (a) edge['] node {\beta^\otimes 1} (b')
        (c') edge['] node[pos=0] {1 \beta^\otimes} (d)
        (b') edge['] node {a^\otimes} (c')
        (b') edge['] node[pos=.2] {\vmtimes 1} (y)
        (c') edge node[pos=.2] {1 \vmtimes} (e)
        ;
        \draw[2cell=\s] 
        (b')++(25:.5) node[rotate=225, 2labelalt={below,\beta^\vmtimes 1}] {\Rightarrow}
        (c')++(242:.45) node[rotate=225, 2labelalt={below,1 \beta^\vmtimes}] {\Rightarrow}
        (z)++(135:1) node[rotate=225, 2labelalt={below,a^\vmtimes}] {\Rightarrow}
        node[between=c and b' at .5] {\divideontimes}
        ;
      \end{scope}
      \begin{scope}[shift={(\w+\m+\m,.2*\h)}]
        \boundary
        \draw[0cell=\s]
        (b)++(-50:.8) node (b'') {\C^2}
        (c)++(-40:.8) node (c'') {\C^2}
        ;
        \draw[1cell=\s] 
        (b) edge node {1 \vmtimes} (b'')
        (c) edge['] node {\vmtimes 1} (c'')
        (b'') edge['] node {\beta^\otimes} (c'')
        (b'') edge node {\vmtimes} (z)
        (c'') edge['] node {\vmtimes} (z)
        ;
        \draw[2cell=\s] 
        (b'')++(-90:.4) node[shift={(-45:.2)}, rotate=225, 2labelalt={below,\beta^\vmtimes}] {\Rightarrow}
        node[between=b'' and y at .5, rotate=225, 2labelalt={below,a^\vmtimes}] {\Rightarrow}
        node[between=c'' and e at .5, rotate=225, 2labelalt={below,a^\vmtimes}] {\Rightarrow}
        ;
      \end{scope}
    \end{tikzpicture}
  \end{equation}
  In the right pasting diagram above, the unlabeled region commutes by
  2-naturality of $\beta^\otimes$.  The
  hexagon labeled $\divideontimes$ commutes by the left hexagon axiom for
  $\beta^\otimes$.
\item[Right hexagon axiom] The following two
  \emph{right hexagon pasting diagrams}\index{hexagon diagram!braided monoidal enriched category} have equal composites.
  \begin{equation}\label{hexagon-bvmR}
    \begin{tikzpicture}[x=17mm,y=17mm,vcenter]
      \def\w{2.05} 
      \def\h{2} 
      \def\m{.5} 
      \def\s{.7}
      \draw[font=\Large] (\w+\m,0) node (eq) {=}; 
      \newcommand{\boundary}{
        \draw[0cell=\s] 
        (.5,0) node (b) {(\C^2)\C}
        (0,-.5) node (c) {\C(\C^2)}
        (b)++(0:1) node (a) {\C(\C^2)}
        (c)++(-90:1) node (d) {(\C^2)\C}
        (a)++(-45:1) node (y) {\C^2}
        (d)++(-45:1) node (e) {\C^2}
        (-45:3.12) node (z) {\C}
        ; 
        \draw[1cell=\s] 
        (a) edge['] node {(a^{\otimes})^{-1}} (b)
        (b) edge['] node {\beta^\otimes} (c)
        (c) edge['] node[pos=.7] {(a^{\otimes})^{-1}} (d)
        (d) edge['] node {\vmtimes 1} (e)
        (e) edge['] node {\vmtimes} (z)
        (a) edge node {1 \vmtimes} (y)
        (y) edge node {\vmtimes} (z)
        ;
      }
      \begin{scope}[shift={(0,.2*\h)}]
        \boundary
        \draw[0cell=\s]
        (a)++(-100:.8) node (b') {\C(\C^2)}
        (d)++(10:.8) node (c') {(\C^2)\C}
        ;
        \draw[1cell=\s] 
        (a) edge['] node {1 \beta^\otimes} (b')
        (c') edge['] node[pos=0] {\beta^\otimes 1} (d)
        (b') edge['] node {(a^{\otimes})^{-1}} (c')
        (b') edge['] node[pos=.2] {1 \vmtimes} (y)
        (c') edge node[pos=.2] {\vmtimes 1} (e)
        ;
        \draw[2cell=\s] 
        (b')++(25:.5) node[rotate=225, 2labelalt={below,1 \beta^\vmtimes}] {\Rightarrow}
        (c')++(242:.45) node[rotate=225, 2labelalt={below,\beta^\vmtimes 1}] {\Rightarrow}
        (z)++(135:.75) node[rotate=225, 2labelalt={below,(a^{\vmtimes}_1)^{-1}}] {\Rightarrow}
        node[between=c and b' at .5] {\divideontimes}
        ;
      \end{scope}
      \begin{scope}[shift={(\w+\m+\m,.2*\h)}]
        \boundary
        \draw[0cell=\s]
        (b)++(-50:.8) node (b'') {\C^2}
        (c)++(-40:.8) node (c'') {\C^2}
        ;
        \draw[1cell=\s] 
        (b) edge node {\vmtimes 1} (b'')
        (c) edge['] node {1 \vmtimes} (c'')
        (b'') edge['] node {\beta^\otimes} (c'')
        (b'') edge node {\vmtimes} (z)
        (c'') edge['] node {\vmtimes} (z)
        ;
        \draw[2cell=\s] 
        (b'')++(-90:.4) node[shift={(-45:.2)}, rotate=225, 2labelalt={below,\beta^\vmtimes}] {\Rightarrow}
        node[between=b'' and y at .6, rotate=225, 2labelalt={below,(a^{\vmtimes}_1)^{-1}}] {\Rightarrow}
        node[between=c'' and e at .7, shift={(0:.5)}, rotate=225, 2labelalt={below,(a^{\vmtimes}_1)^{-1}}] {\Rightarrow}
        ;
      \end{scope}
    \end{tikzpicture}
  \end{equation}
  In the right pasting diagram above, the unlabeled region commutes by
  2-naturality of $\beta^\otimes$.  The
  hexagon labeled $\divideontimes$ commutes by the right hexagon axiom for
  $\beta^\otimes$.  The 2-cell isomorphism
  $(a^{\vmtimes}_1)^{-1}$ is the inverse of $a^{\vmtimes}_1$ (\cref{eq:a-vmtimes-inv-mate}).
\end{description}
This finishes the definition of a braided monoidal $\V$-category.
\end{definition}

\begin{definition}\label{definition:symm-monoidal-vcat}
Suppose $\V$ is a symmetric monoidal category.  A \emph{symmetric
    monoidal $\V$-category}\index{symmetric monoidal category!enriched}\index{enriched!symmetric monoidal category}\index{enriched category!symmetric monoidal} is a braided monoidal $\V$-category $(\C,\beta^\vmtimes)$ that satisfies the following axiom.
  \begin{description}
  \item[Symmetry axiom]\index{symmetry axiom!symmetric monoidal enriched category} The following two
    \emph{symmetry pasting diagrams} have equal composites.
    \begin{equation}\label{vmonoidal-symmetry}
    \begin{tikzpicture}[x=16mm,y=14mm,vcenter]
      \def\w{2} 
      \def\h{2} 
      \def\m{.5} 
      \def\s{.75}
      \draw[font=\Large] (\w+\m,0) node (eq) {=}; 
      \newcommand{\boundary}{
        \draw[0cell=\s] 
        (0,0) node (a) {\C^2}
        (a)++(-60:2) node (c) {\C^2}
        (a)++(0:2) node (d) {\C}
        ; 
        \draw[1cell=\s] 
        (a) edge[',bend right=25] node {1} (c)
        (a) edge node {\vmtimes} (d)
        (c) edge[',bend right=25] node {\vmtimes} (d)
        ;
      }
      \begin{scope}[shift={(0,\h/2)}]
        \boundary
        \draw[0cell=\s]
        (c)++(90:2*.57735-.1) node (b) {\C^2}
        ;
        \draw[1cell=\s] 
        (a) edge['] node {\beta^\otimes} (b)
        (b) edge['] node {\beta^\otimes} (c)
        (b) edge['] node {\vmtimes} (d)
        ;
        \draw[2cell=\s] 
        (b)++(90:.57735/2+.07) node[rotate=-90, 2label={above,\beta^\vmtimes}] {\Rightarrow}
        (b)++(-30:.57735/2) node[rotate=-90, 2label={above,\beta^\vmtimes}] {\Rightarrow}
        (b)++(210:.57735/2+.07) node {\divideontimes}
        ;
      \end{scope}
      \begin{scope}[shift={(\w+\m+\m,\h/2)}]
        \boundary
      \end{scope}
    \end{tikzpicture}
    \end{equation}
  \end{description}
  The pasting diagram on the right is the identity $\V$-natural
  transformation of $\vmtimes$.  In the left pasting diagram, the region labeled
  $\divideontimes$ commutes by the symmetry axiom for $\beta^\otimes$.  We also call $\beta^\vmtimes$ the \emph{symmetry}\index{symmetry!symmetric monoidal enriched category} of $\C$.
\end{definition}

The following result, which is proved in \cite[1.5.5]{cerberusIII}, strengthens \cref{theorem:vcat-mon} to the $\Cat$-enriched context, where $(\Cat, \times, \boldone)$ is the Cartesian closed category in \cref{ex:cat}.

\begin{theorem}\label{theorem:vcat-cat-mon}
The following statements hold for a braided monoidal category $\V = (\V,\otimes,\xi)$.
\begin{enumerate}
\item There is a monoidal $\Cat$-category (\cref{definition:monoidal-vcat})
\[\big( \VCat,\otimes,\vtensorunit, a^\otimes, \ell^\otimes, r^\otimes \big).\]
\item If $\V$ is a symmetric monoidal category, then $(\VCat,\beta^\otimes)$ is a symmetric monoidal $\Cat$-category (\cref{definition:symm-monoidal-vcat}).
\end{enumerate}
\end{theorem}

\subsection*{Enriched Symmetric Monoidal Functors}

We need the following enriched variants of mates to define enriched monoidal functors.

\begin{definition}\label{eq:lr-vmtimes-mates}
Suppose $\V$ is a braided monoidal category and $\C$ is a monoidal $\V$-category.  We define $\ell^\vmtimes_1$\index{left unitor!mate}\index{mate!left unitor} and $r^\vmtimes_1$\index{right unitor!mate}\index{mate!right unitor} as the mates of $\ell^\vmtimes$ and $r^\vmtimes$ in \cref{eq:monoidal-unitors} given, respectively, by replacing $(\ell^\otimes)^{-1}$ and $(r^\otimes)^{-1}$ with $\ell^\otimes$ and $r^\otimes$, as follows.
    \[\begin{tikzpicture}[x=30mm,y=13mm,vcenter]
    \def\s{.8}
      \draw[0cell=\s] 
      (0,0) node (a) {\C}
      (.1,1) node (b) {\vtensorunit\C}
      (.9,1) node (c) {\C^2}
      (1,0) node (d) {\C}
      ;
      \draw[1cell=\s] 
      (b) edge['] node {\ell^{\otimes}} (a)
      (b) edge node (z) {\vmunitbox 1_{\C}} (c)
      (c) edge node {\vmtimes} (d)
      (a) edge['] node (w) {1_{\C}} (d)
      ;
      \draw[2cell] 
      node[between=z and w at .5, rotate=-90, 2label={above,\ell^\vmtimes_1}] {\Rightarrow}
      ;
    \end{tikzpicture}
    \qquad
    \begin{tikzpicture}[x=30mm,y=13mm,vcenter]
    \def\s{.8}
      \draw[0cell=\s] 
      (0,0) node (a) {\C}
      (.1,1) node (b) {\C\vtensorunit}
      (.9,1) node (c) {\C^2}
      (1,0) node (d) {\C}
      ;
      \draw[1cell=\s] 
      (b) edge['] node {r^{\otimes}} (a)
      (b) edge node (z) {1_{\C} \vmunitbox} (c)
      (c) edge node {\vmtimes} (d)
      (a) edge['] node (w) {1_{\C}} (d)
      ;
      \draw[2cell] 
      node[between=z and w at .5, rotate=-90, 2label={above,\,r^\vmtimes_1}] {\Rightarrow}
      ;
    \end{tikzpicture}\]
This finishes the definition.
\end{definition}

\begin{definition}\label{definition:monoidal-V-fun}
Suppose $\C$ and $\D$ are monoidal $\V$-categories with $\V$ a braided monoidal category.  A
  \emph{monoidal $\V$-functor}\index{monoidal functor!enriched}\index{enriched!monoidal functor}\index{functor!monoidal enriched} from $\C$ to $\D$ is a triple 
  \[\C \fto{(F,F^2,F^0)} \D\]
consisting of 
  \begin{itemize}
  \item a $\V$-functor $F\cn \C \to \D$ and
  \item $\V$-natural transformations $F^2$ and $F^0$, which are called, respectively, the \emph{monoidal constraint}\index{monoidal constraint!enriched} and the \index{unit constraint!enriched}\emph{unit constraint}, as displayed below.
  \begin{equation}\label{enr_constraints}
  \begin{tikzpicture}[x=20mm,y=15mm,vcenter]
  \def\s{.8}
    \draw[0cell=\s] 
    (0,0) node (k2) {\C^2}
    (1,0) node (l2) {\D^2}
    (0,-1) node (k) {\C}
    (1,-1) node (l) {\D}
    ;
    \draw[1cell=\s] 
    (k2) edge['] node {\vmtimes} (k)
    (l2) edge node {\vmtimes} (l)
    (k2) edge node {F \otimes F} (l2)
    (k) edge node {F} (l)
    ;
    \draw[2cell=\s] 
    node[between=l2 and k at .5, rotate=225, 2labelalt={below,F^2}] {\Rightarrow}
    ;
  \end{tikzpicture}
  \qquad\qquad
  \begin{tikzpicture}[x=20mm,y=17mm,vcenter]
  \def\s{.8}
    \draw[0cell=\s] 
    (0,0) node (i) {\vtensorunit}
    (240:1) node (k) {\C}
    (-60:1) node (l) {\D}
    ;
    \draw[1cell=\s] 
    (i) edge['] node {\vmunitbox} (k)
    (i) edge node {\vmunitbox} (l)
    (k) edge node {F} (l)
    ;
    \draw[2cell=\s] 
    (-90:.52) node[rotate=180, 2label={below,F^0}] {\Rightarrow}
    ;
  \end{tikzpicture}
  \end{equation}
  \end{itemize}
The above data are subject to the axioms \cref{enrmonfunctor-ass,enrmonfunctor-lunity,enrmonfunctor-runity}.
  \begin{description}
  \item[Associativity]\index{associativity!enriched monoidal functor} The following two
    \emph{associativity pasting diagrams} have equal composites.
    \begin{equation}\label{enrmonfunctor-ass}
    \begin{tikzpicture}[x=13mm,y=14mm,baseline={(eq.base)}]
      \def\wl{2.1} 
      \def\wr{.7} 
      \def\h{.5} 
      \def\m{.5} 
      \def\s{.7}
      \draw[font=\Large] (0,0) node (eq) {=}; 
      \newcommand{\boundary}{
        \draw[0cell=\s] 
        (0,0) node (a) {(\C^{2})\C}
        (225:.707) node (a') {\C(\C^{2})}
        (a')++(0,-1) node (c) {\C^{2}}
        (c)++(-20:1.5) node (d) {\C}
        (a)++(10:1.25) node (z) {(\D)^2\D}
        (z)++(-30:1) node (w) {\D^2}
        (w)++(0,-1.5) node (u) {\D}
        ; 
        \draw[1cell=\s] 
        (a) edge['] node {a^\otimes} (a')
        (a) edge node[pos=.7] {(FF)F} (z)
        (c) edge['] node {\vmtimes} (d)
        (a') edge['] node {1 \vmtimes} (c)
        (z) edge node {\vmtimes 1} (w)
        (w) edge node {\vmtimes} (u)
        (d) edge['] node {F} (u)
        ;
      }
      \begin{scope}[shift={(-\wl-\m,\h)}]
        \boundary
        \draw[0cell=\s]
        (d)++(90:1.5) node (b) {\C^{2}}
        ;
        \draw[1cell=\s] 
        (a) edge['] node {\vmtimes 1} (b)
        (b) edge node {\vmtimes} (d)
        (b) edge node {FF} (w)
        ;
        \draw[2cell=\s] 
        node[between=b and c at .5, rotate=225, 2labelalt={below,a^\vmtimes}] {\Rightarrow}
        node[between=z and b at .5, rotate=225, 2labelalt={below,F^2 1}] {\Rightarrow}
        node[between=w and d at .5, rotate=225, 2labelalt={below,F^2}] {\Rightarrow}
        ;
      \end{scope}
      \begin{scope}[shift={(\wr+\m,\h)}]
        \boundary
        \draw[0cell=\s]
        (z)++(225:.707) node (z') {\D(\D^2)}
        (z')++(270:1) node (v) {\D^2}
        ;
        \draw[1cell=\s] 
        (z) edge node {a^\otimes} (z')
        (z') edge node {1 \vmtimes} (v)
        (v) edge node {\vmtimes} (u)
        (a') edge['] node[pos=.3] {F(FF)} (z')
        (c) edge node {FF} (v)
        ;
        \draw[2cell=\s] 
        node[between=v and w at .6, rotate=225, 2labelalt={below,a^\vmtimes}] {\Rightarrow}
        node[between=v and d at .5, rotate=225, 2labelalt={below,F^2}] {\Rightarrow}
        node[between=a' and v at .5, rotate=225, 2labelalt={below,1 F^2}] {\Rightarrow}
        ;
      \end{scope}
    \end{tikzpicture}
    \end{equation}
    In the right pasting diagram, the unlabeled region
    commutes by naturality of $a^\otimes$. 
  \item[Left unity]\index{left unity!enriched monoidal functor} The following two \emph{left unity pasting
    diagrams} have equal composites.
    \begin{equation}\label{enrmonfunctor-lunity}
    \begin{tikzpicture}[x=22mm,y=20mm,baseline={(eq.base)}]
      \def\wl{1.6} 
      \def\wr{0} 
      \def\h{.25} 
      \def\m{.3} 
      \def\s{.7}
      \draw[font=\Large] (0,0) node (eq) {=}; 
      \newcommand{\boundary}{
        \draw[0cell=\s] 
        (0,0) node (a) {\vtensorunit\C}
        (1,0) node (b) {\D\C}
        (1.6,0) node (c) {\D^2}
        (0,-1) node (a') {\C}
        (1.6,-1) node (c') {\D}
        ; 
        \draw[1cell=\s] 
        (a) edge node {\vmunitbox 1} (b)
        (b) edge node {1 F} (c)
        (a') edge node {F} (c')
        (a) edge['] node {\ell^\otimes} (a')
        (c) edge node {\vmtimes} (c')
        ;
      }
      \begin{scope}[shift={(-\wl-\m,\h)}]
        \boundary
        \draw[0cell=\s]
        (a)++(-30:.707) node (k2) {\C^2}
        ;
        \draw[1cell=\s] 
        (k2) edge['] node {F1} (b)
        (k2) edge node {\vmtimes} (a')
        (a) edge['] node {\vmunitbox 1} (k2)
        ;
        \draw[2cell=\s] 
        (k2)++(0,.24) node[shift={(.07,0)}, rotate=-90, 2labelalt={below,F^0 1}] {\Rightarrow}
        (a')++(65:.45) node[rotate=180, 2labelalt={below,\hspace{6pt}\ell^\vmtimes_1}] {\Rightarrow}
        node[between=b and c' at .55, shift={(-.2,0)}, rotate=225, 2labelalt={below,F^2}] {\Rightarrow}
        ;
      \end{scope}
      \begin{scope}[shift={(\wr+\m,\h)}]
        \boundary
        \draw[0cell=\s]
        (a)++(-30:1) node (il) {\vtensorunit\D}
        ;
        \draw[1cell=\s] 
        (a) edge['] node {1 F} (il)
        (il) edge node[pos=.44] {\vmunitbox 1} (c)
        (il) edge['] node {\ell^\otimes} (c')
        ;
        \draw[2cell=\s] 
        (il)++(-5:.5) node[rotate=225, 2labelalt={below,\ell^\vmtimes_1}] {\Rightarrow}
        ;
      \end{scope}
    \end{tikzpicture}
    \end{equation}
    In the right pasting diagram, the lower unlabeled region
    commutes by naturality of $\ell^\otimes$.  The upper unlabeled
    triangle commutes by 2-functorality of $\otimes$. 
  \item[Right unity]\index{right unity!enriched monoidal functor} The following two \emph{right unity pasting
    diagrams} have equal composites.
    \begin{equation}\label{enrmonfunctor-runity}
    \begin{tikzpicture}[x=22mm,y=20mm,baseline={(eq.base)}]
      \def\wl{1.6} 
      \def\wr{0} 
      \def\h{.25} 
      \def\m{.3} 
      \def\s{.7}
      \draw[font=\Large] (0,0) node (eq) {=}; 
      \newcommand{\boundary}{
        \draw[0cell=\s] 
        (0,0) node (a) {\C\vtensorunit}
        (1,0) node (b) {\C\D}
        (1.6,0) node (c) {\D^2}
        (0,-1) node (a') {\C}
        (1.6,-1) node (c') {\D}
        ; 
        \draw[1cell=\s] 
        (a) edge node {1 \vmunitbox} (b)
        (b) edge node {F1} (c)
        (a') edge node {F} (c')
        (a) edge['] node {r^\otimes} (a')
        (c) edge node {\vmtimes} (c')
        ;
      }
      \begin{scope}[shift={(-\wl-\m,\h)}]
        \boundary
        \draw[0cell=\s]
        (a)++(-30:.707) node (k2) {\C^2}
        ;
        \draw[1cell=\s] 
        (k2) edge['] node {1F} (b)
        (k2) edge node {\vmtimes} (a')
        (a) edge['] node {1\vmunitbox} (k2)
        ;
        \draw[2cell=\s] 
        (k2)++(0,.24) node[shift={(.07,0)}, rotate=-90, 2labelalt={below,1F^0}] {\Rightarrow}
        (a')++(65:.45) node[rotate=180, 2labelalt={below,\hspace{6pt}r^\vmtimes_1}] {\Rightarrow}
        node[between=b and c' at .55, shift={(-.2,0)}, rotate=225, 2labelalt={below,F^2}] {\Rightarrow}
        ;
      \end{scope}
      \begin{scope}[shift={(\wr+\m,\h)}]
        \boundary
        \draw[0cell=\s]
        (a)++(-30:1) node (il) {\D\vtensorunit}
        ;
        \draw[1cell=\s] 
        (a) edge['] node {F1} (il)
        (il) edge node[pos=.44] {1\vmunitbox} (c)
        (il) edge['] node {r^\otimes} (c')
        ;
        \draw[2cell=\s] 
        (il)++(-5:.5) node[rotate=225, 2labelalt={below,r^\vmtimes_1}] {\Rightarrow}
        ;
      \end{scope}
    \end{tikzpicture}
    \end{equation}
    In the right pasting diagram, the lower unlabeled region
    commutes by naturality of $r^\otimes$.  The upper unlabeled
    triangle commutes by 2-functorality of $\otimes$.
  \end{description}
This finishes the definition of a monoidal $\V$-functor $(F,F^2,F^0)$.  Moreover, it  is
  \begin{itemize}
  \item \emph{unital}\index{unital monoidal!enriched functor} if $F^0$ is a $\V$-natural isomorphism and
  \item \emph{strictly unital}\index{strictly unital!monoidal enriched functor} if $F^0$ is the identity $\V$-natural transformation.
  \end{itemize}
Composition of monoidal $\V$-functors is defined by composing the $\V$-functors, pasting the monoidal constraints, and pasting the unit constraints.
\end{definition}

\begin{definition}\label{definition:braided-monoidal-vfunctor}
  Suppose $\C$ and $\D$ are braided monoidal $\V$-categories with
  $\V$ a symmetric monoidal category.  A \emph{braided monoidal}\index{braided monoidal functor!enriched}\index{enriched!braided monoidal functor}\index{functor!braided monoidal enriched}\index{monoidal functor!braided enriched} $\V$-functor
  \[
  \C \fto{(F,F^2,F^0)} \D
  \]
  is a monoidal $\V$-functor that satisfies the following axiom.
  \begin{description}
  \item[Braid axiom]\index{braid axiom!braided monoidal enriched functor} The following two \emph{braiding pasting diagrams}\index{braiding diagram!braided monoidal enriched functor} have equal composites.
    \begin{equation}\label{enrmonfunctor-braided}
    \begin{tikzpicture}[x=17mm,y=10mm,baseline={(eq.base)}]
      \def\wl{1} 
      \def\wr{1} 
      \def\h{1} 
      \def\m{.5} 
      \def\s{.8}
      \draw[font=\Large] (0,0) node (eq) {=}; 
      \newcommand{\boundary}{
        \draw[0cell=\s] 
        (0,0) node (a) {\C^2}
        (-1,-1) node (b) {\C^2}
        (0,-2) node (c) {\C}
        (1,0) node (a') {\D^2}
        (1,-2) node (c') {\D}
        ; 
        \draw[1cell=\s] 
        (a) [rounded corners=2pt] -| node[swap,pos=.7] {\beta^\otimes} (b)
        (a) edge node {FF} (a')
        (c) edge node {F} (c')
        (a') edge node {\vmtimes} (c')
        ;
		\draw[1cell=\s]
		(b) [rounded corners=2pt] |- node[swap,pos=.3] {\vmtimes} (c);
      }
      \begin{scope}[shift={(-\wl-\m,\h)}]
        \boundary
        \draw[0cell=\s]
        (0,-1) node (b') {\D^2}
        ;
        \draw[1cell=\s] 
        (a') edge['] node[pos=.75,inner sep=0pt] {\beta^\otimes} (b')
        (b') edge['] node[pos=.3] {\vmtimes} (c')
        (b) edge node {FF} (b')
        ;
        \draw[2cell=\s] 
        (b')++(.6,0) node[rotate=180, 2labelalt={below,\beta^\vmtimes}] {\Rightarrow}
        node[between=b' and c at .55, shift={(-.1,0)}, rotate=240, 2labelalt={below,F^2}] {\Rightarrow}
        ;
      \end{scope}
      \begin{scope}[shift={(\wr+\m,\h)}]
        \boundary
        \draw[1cell=\s] 
        (a) edge node {\vmtimes} (c)
        ;
        \draw[2cell=\s] 
        (b)++(.5,0) node[rotate=180, 2labelalt={below,\beta^\vmtimes}] {\Rightarrow}        
        node[between=c and a' at .55, rotate=225, 2labelalt={below,F^2}] {\Rightarrow}
        ;
      \end{scope}
    \end{tikzpicture}
    \end{equation}
  \end{description}
  In the left pasting diagram, the unlabeled region commutes
  by naturality of $\beta^\otimes$.  This finishes the definition of a
  braided monoidal $\V$-functor.  If $\C$ and $\D$ are symmetric monoidal $\V$-categories, then we
  call $F$ a \emph{symmetric monoidal $\V$-functor}.\index{symmetric monoidal functor!enriched}\index{enriched!symmetric monoidal - functor}\index{functor!symmetric monoidal enriched}\index{monoidal functor!symmetric enriched}
\end{definition}

\subsection*{Symmetric Monoidal Closed Categories}

\begin{definition}\label{definition:eval}
Suppose $(\V, \otimes,[,])$ is a symmetric monoidal closed category (\cref{def:closedcat}).  For an object $a \in \V$, the \index{evaluation!at $a$}\emph{evaluation at $a$} is the counit
\begin{equation}\label{evaluation}
[a,-] \otimes a \fto{\ev_{a,-}} 1_\V
\end{equation}
of the adjunction 
\[- \otimes a \cn \V \lradj \V \cn [a,-]\]
that is part of the closed structure.
\end{definition}

\begin{definition}\label{definition:canonical-v-enrichment}
Suppose $(\V,\otimes,\tu,\alpha,\lambda,\rho,\xi,[,])$ is a symmetric monoidal closed category.  The \index{canonical self-enrichment}\index{self-enrichment}\index{symmetric monoidal category!closed!self-enrichment}\index{category!symmetric monoidal!closed!self-enrichment}\emph{canonical self-enrichment} of $\V$, denoted $\Vse$, is the $\V$-category defined as follows.
\begin{description}
\item[Objects] $\Ob(\Vse) = \Ob(\V)$.
\item[Hom objects] The hom $\V$-object for objects $a,b \in \Vse$ is the object  
\[\Vse(a,b) = [a,b].\]
\item[Composition] For objects $a,b,c \in \Vse$, the composition $\V$-morphism
\[[b,c] \otimes [a,b] \fto{\mcomp_{a,b,c}} [a,c]\]
is the adjoint of the following composite $\V$-morphism.
\begin{equation}\label{eq:m-adj}
\begin{tikzpicture}[x=45mm,y=14mm,vcenter]
\draw[0cell]
(0,0) node (a) {([b,c] \otimes [a,b]) \otimes a}
(.1,-1) node (b) {[b,c] \otimes ([a,b] \otimes a)}
(.9,-1) node (c) {[b,c] \otimes b}
(1,0) node (d) {c}
;
\draw[1cell]
(a) edge node[swap,pos=.3] {\al} node[pos=.8] {\iso} (b)
(b) edge node {1 \otimes \ev} (c)
(c) edge['] node[pos=.6] {\ev} (d)
;
\end{tikzpicture}
\end{equation}
\item[Identities] For an object $a \in \Vse$, the identity is the $\V$-morphism
\[\tu \fto{\cone_a} [a,a]\]
adjoint to the left unit isomorphism 
\begin{equation}\label{i-adjoint}
\tu \otimes a \fto[\iso]{\la} a \inspace \V.
\end{equation}
\end{description}
This finishes the definition of the $\V$-category $\Vse$.  To simplify the notation, we also denote $\Vse$ by $\V$.
\end{definition}

A detailed proof of the following result appears in \cite[3.3.2]{cerberusIII}.  It states that a symmetric monoidal closed category is also symmetric monoidal in the self-enriched sense.

\begin{theorem}\label{theorem:v-closed-v-sm}\index{canonical self-enrichment!symmetric monoidal}\index{symmetric monoidal category!closed!self-enrichment}
For a symmetric monoidal closed category $\V$, the symmetric monoidal structure extends to $\Vse$ such that $\Vse$ is a symmetric monoidal $\V$-category.
\end{theorem}

The following statements about change of enrichment are proved in \cite[2.4.10 and 3.3.4]{cerberusIII}.

\begin{theorem}\label{thm:change-enrichment}
Suppose $f \cn \V \to \W$ is a symmetric monoidal functor between symmetric monoidal categories.  Then the following statements hold.
\begin{enumerate}
\item\label{change-enr-i} Changing enrichment\index{change of enrichment} along $f$ sends each symmetric monoidal $\V$-category $\C$ to a symmetric monoidal $\W$-category $\C_f$ with
\begin{itemize}
\item the same objects as $\C$ and
\item hom objects
\[\C_f(x,y) = f\C(x,y) \in \W \forspace x,y \in \C.\]
\end{itemize} 
Thus, if $\V$ is also closed, then $\Vse_f$ is a symmetric monoidal $\W$-category.
\item\label{change-enr-i'} Changing enrichment along $f$ sends each symmetric monoidal $\V$-functor $F \cn \C \to \D$ between symmetric monoidal $\V$-categories to a symmetric monoidal $\W$-functor
\[\C_f \fto{F_f} \D_f\]
with 
\begin{itemize}
\item the same object assignment as $F$ and
\item component $\W$-morphisms
\[\C_f(x,y) = f\C(x,y) \fto{(F_f)_{x,y} = f F_{x,y}} \D_f(Fx,Fy) = f\D(Fx,Fy)\]
for $x,y \in \C$.
\end{itemize}
\item\label{change-enr-ii} Suppose $\V$ and $\W$ are both closed.  Then $f$ induces a symmetric monoidal $\W$-functor
\[\Vse_f \fto{f'} \Wse\]
with the same object assignment as $f$.
\end{enumerate}
\end{theorem}

\chapter{Multicategories}
\label{ch:prelim_multicat}
This appendix is a brief review of enriched multicategories.  References are provided at the beginning of each section.  Below we list the sections in this appendix, along with the main concepts in each section.

\appsecname{sec:enrmulticat}
\begin{itemize}
\item \appentry{Profiles and $\V$-multicategories}{def:profile,def:enr-multicategory}
\item \appentry{$\V$-multifunctors and $\M$-algebras}{def:enr-multicategory-functor}
\item \appentry{Endomorphism $\V$-multicategories}{definition:EndK,proposition:monoidal-v-cat-v-multicat}
\item \appentry{Endomorphism $\V$-multifunctors}{def:EndF,EndF_multi}
\end{itemize}

\appsecname{sec:multi_iicat}
\begin{itemize}
\item \appentry{$\V$-multinatural transformations}{def:enr-multicat-natural-transformation,def:enr-multinatural-composition,v-multicat-2cat}
\item \appentry{$\Cat$-multiequivalences}{def:catmultiequivalence}
\end{itemize}

\appsecname{sec:mult_change_enr}
\begin{itemize}
\item \appentry{Change of enrichment along symmetric monoidal functors}{mult_change_enr,expl:mult_change_enr}
\item \appentry{Commutation with endomorphism}{End_dashf}
\end{itemize}

\appsecname{sec:pointed-diagrams}
\begin{itemize}
\item \appentry{Symmetric monoidal closed category of pointed objects}{def:pointed-objects,def:wedge-smash-phom,theorem:pC-sm-closed}
\item \appentry{Null objects and pointed unitary enrichment}{definition:zero,def:unitary-enrichment}
\item \appentry{Pointed functors and pointed Day convolution}{def:pointed-category,definition:Dgm-pV-convolution-hom}
\item \appentry{Symmetric monoidal closed categories and enriched multicategories of pointed diagrams}{thm:Dgm-pv-convolution-hom,theorem:diagram-omnibus}
\end{itemize}

\appsecname{sec:dstarv-multicat}
\begin{itemize}
\item \appentry{Multimorphism objects of $\DstarV$}{pmapsmafh}
\item \appentry{Symmetric group action of $\DstarV$}{dstarvsigmaaction}
\item \appentry{Multicategorical composition of $\DstarV$}{dstarv-multicomp}
\end{itemize}

Throughout this appendix, 
\[\big(\V,\otimes,\tu,\alpha,\lambda,\rho,\xi\big)\]
denotes a symmetric monoidal category (\cref{def:symmoncat}), such as the Cartesian closed categories $\Set$ or $\Cat$.

\section{Enriched Multicategories}
\label{sec:enrmulticat}

In this section, we recall the definitions of enriched multicategories (\cref{def:enr-multicategory}) and enriched multifunctors (\cref{def:enr-multicategory-functor}).  Then we review the endomorphism construction that associates to each enriched symmetric monoidal category an enriched multicategory (\cref{proposition:monoidal-v-cat-v-multicat,EndF_multi}).  More detailed discussion of this material can be found in \cite[Chapter 6]{cerberusIII} and \cite{yau-operad}.  We use the following notation to denote finite sequence of objects.

\begin{definition}\label{def:profile}
Suppose $C$\label{notation:s-class} is a class.  The class of finite, possibly empty tuples in $C$ is denoted by\index{profile}\label{notation:profs}
\[\Prof(C) = \coprodover{m \geq 0}\ C^{m}.\] 
\begin{itemize}
\item An element in $\Prof(C)$ is called a \emph{$C$-profile}.  A $C$-profile of length\index{length!of a profile} $n=\len\angc$ is denoted by\label{notation:us}  
\[\angc = \ang{c_j}_{j=1}^n = (c_1, \ldots, c_n) \in C^{n}.\]
The empty $C$-profile\index{empty profile} is denoted by $\ang{}$ or $\emptyset$.
\item An element in $\Prof(C)\times C$ is denoted by\label{notation:duc} $\smscmap{\angc; c'}$ with $\angc\in\Prof(c)$ and $c'\in C$.
\item The \emph{concatenation}\index{concatenation} of two $C$-profiles $\angc = \ang{c_i}_{i=1}^m$ and $\angd = \ang{d_j}_{j=1}^n$ is the $C$-profile\label{not:concat}
\begin{equation}\label{concatprof}
\angc \oplus \angd = \big(c_1, \ldots , c_m, d_1, \ldots , d_n\big).
\end{equation}
Concatenation is associative with the empty tuple $\ang{}$ as the strict unit.\defmark
\end{itemize}
\end{definition}

The \index{symmetric group}symmetric group on $n$ letters is denoted $\Sigma_n$.

\begin{definition}\label{def:enr-multicategory}
A \emph{$\V$-multicategory}\index{enriched multicategory}\index{category!enriched multi-}\index{multicategory!enriched}, also called a \emph{multicategory in $\V$}, is a triple\label{notation:enr-multicategory}  
\[(\M, \gamma, \operadunit)\]
consisting of the following data.
\begin{description}
\item[Objects] It is equipped with a class $\ObM$ of \index{object!enriched multicategory} \emph{objects}.  We write $\ProfM$ for $\Prof(\ObM)$.
\item[Multimorphisms] For each $\smscmap{\angx;x'} \in \ProfMM$ with $\angx = \ang{x_j}_{j=1}^n$, it is equipped with an \index{multimorphism}\emph{$n$-ary multimorphism object} in $\V$\label{notation:enr-cduc}
\begin{equation}\label{multimorphism_object}
\M\scmap{\angx;x'} = \M\mmap{x'; x_1,\ldots,x_n}
\end{equation}
with \emph{input profile}\index{input profile} $\angx$ and \emph{output}\index{output} $x'$.

\item[Symmetric group action]
For each $\smscmap{\angx;x'} \in \ProfMM$ with $n=\len\angx$  and permutation $\sigma \in
\Sigma_n$, $\M$ is equipped with a $\V$-isomorphism called the \emph{right $\sigma$-action}\index{right action} or \index{symmetric group!action}\emph{symmetric group action},
\begin{equation}\label{rightsigmaaction}
\begin{tikzcd}[column sep=large]
\M\scmap{\angx;x'} \rar{\sigma}[swap]{\cong} & \M\scmap{\angx\sigma; x'},\end{tikzcd}
\end{equation}
where\label{enr-notation:c-sigma}
\[\angx\sigma = \ang{x_{\sigma(j)}}_{j=1}^n 
= \big(x_{\sigma(1)}, \ldots, x_{\sigma(n)}\big) 
\in \ProfM\]
is the right permutation\index{right permutation} of $\angx$ by $\sigma$.  Set-theoretically, for $a \in \M\scmap{\angx;x'}$, we denote the right $\sigma$-action on $a$ by $a\sigma$, $a \cdot \sigma$, or $a^\sigma$.
\item[Units] Each object $x \in \M$ is equipped with a $\V$-morphism\label{notation:enr-unit-c} called the \index{colored unit}\emph{$x$-colored unit},
\begin{equation}\label{ccoloredunit}
\begin{tikzcd}[column sep=large]
\tu \ar{r}{\operadunit_x} & \M\scmap{x;x}.
\end{tikzcd}
\end{equation}
\item[Composition] For each
\begin{itemize}
\item $\scmap{\ang{x'};x''} \in \ProfMM$ with $\ang{x'} = \ang{x'_j}_{j=1}^n \in \ProfM$ and 
\item $\ang{x_j} = \ang{x_j^i}_{i=1}^{k_j} \in \ProfM$ for each $j\in\{1,\ldots,n\}$ with $\angx = \oplus_{j=1}^n \ang{x_j}$,
\end{itemize}
$\M$ is equipped with a $\V$-morphism called the \index{multicategory!composition}\index{composition!multicategory}\emph{composition} or \emph{multicategorical composition}\label{notation:enr-multicategory-composition},
\begin{equation}\label{eq:enr-defn-gamma}
\begin{tikzcd}[column sep=large]
\M\mmap{x'';\ang{x'}} \otimes
\txotimes_{j=1}^n \M\mmap{x_j';\ang{x_j}} \rar{\gamma} &
\M\mmap{x'';\ang{x}}.
\end{tikzcd}
\end{equation}
\end{description}
The data above are subject to the following axioms.
\begin{description}
\item[Unity]
Suppose $\scmap{\angx;x'} \in \ProfMM$ with $\angx = \ang{x_j}_{j=1}^n \in \ProfM$.\index{unity!enriched multicategory}
\begin{itemize}
\item The following \emph{right unity diagram}\index{right unity!enriched multicategory} commutes if $n \geq 1$.
\begin{equation}\label{enr-multicategory-right-unity}
\begin{tikzcd} 
\M\scmap{\angx;x'} \otimes \tensorunit^{\otimes n} \dar[swap]{1 \otimes (\otimes_j \operadunit_{x_j})} \ar[bend left=15]{dr}{\rho} &\\
\M\scmap{\angx;x'} \otimes \txotimes_{j=1}^n \M\scmap{x_j;x_j} \rar{\gamma} & \M\scmap{\angx;x'}
\end{tikzcd}
\end{equation}
\item
The \index{unity!enriched multicategory}\index{left unity!enriched multicategory}\emph{left unity diagram} below commutes.
\begin{equation}\label{enr-multicategory-left-unity}
\begin{tikzcd}
\tensorunit \otimes \M\scmap{\angx;x'} \dar[swap]{\operadunit_{x'} \otimes 1} \ar[bend left=15]{dr}{\lambda} &\\
\M\mmap{x';x'} \otimes \M\scmap{\angx;x'} \rar{\gamma} & \M\scmap{\angx;x'}
\end{tikzcd}
\end{equation}
\end{itemize}
\item[Associativity]
Suppose given
\begin{itemize}
\item $\scmap{\ang{x''};x'''} \in \ProfMM$ with $\ang{x''} = \ang{x''_j}_{j=1}^n \in \ProfM$,
\item $\ang{x_j'} = \ang{x'_{j,i}}_{i=1}^{k_{j}} \in \ProfM$ for each $j \in \{1,\ldots,n\}$ with $\ang{x'} = \oplus_{j=1}^n \ang{x_j'}$ and $k_j > 0$ for at least one $j$, and
\item $\ang{x_{j,i}} = \ang{x_{j,i,p}}_{p=1}^{\ell_{j,i}} \in \ProfM$ for each $j\in\{1,\ldots,n\}$ and each $i \in \{1,\ldots,k_j\}$ with $\ang{x_j} = \oplus_{i=1}^{k_j}{\ang{x_{j,i}}}$ and $\ang{x} =
\oplus_{j=1}^n{\ang{x_{j}}}$.
\end{itemize}
Then the following \index{associativity!enriched multicategory}\emph{associativity diagram} in $\V$ commutes.
\begin{equation}\label{enr-multicategory-associativity}
\begin{tikzpicture}[x=40mm,y=15mm,vcenter]
  \draw[0cell=.85] 
  (0,0) node (a) {\textstyle
    \M\mmap{x''';\ang{x''}}
    \otimes
    \biggl[\bigotimes\limits_{j=1}^n \M\mmap{x''_j;\ang{x'_{j}}}\biggr]
    \otimes
    \bigotimes\limits_{j=1}^n \biggl[\bigotimes\limits_{i=1}^{k_j} \M\mmap{x'_{j,i};\ang{x_{j,i}}}\biggr] 
  }
  (1,.8) node (b) {\textstyle
    \M\mmap{x''';\ang{x'}}
    \otimes
    \bigotimes\limits_{j=1}^{n} \biggl[\bigotimes\limits_{i=1}^{k_j} \M\mmap{x'_{j,i};\ang{x_{j,i}}}\biggr]
  }
  (0,-1.2) node (a') {\textstyle
    \M\mmap{x''';\ang{x''}} \otimes
    \bigotimes\limits_{j=1}^n \biggl[\M\mmap{x_j'';\ang{x_j'}} \otimes \bigotimes\limits_{i=1}^{k_j} \M\mmap{x'_{j,i};\ang{x_{j,i}}}\biggr]
  }
  (1,-2) node (b') {\textstyle
    \M\mmap{x''';\ang{x''}} \otimes \bigotimes\limits_{j=1}^n \M\mmap{x_j'';\ang{x_{j}}}
  }
  (1.2,-.6) node (c) {\textstyle
    \M\mmap{x''';\ang{x}}
  }
  ;
  \draw[1cell=.85]
  (a) edge[shorten <=-1ex,shorten >=-1ex] node {\iso} node['] {\mathrm{permute}} (a')
  (a) edge[shorten >=-4ex,shorten <=-1ex, transform canvas={xshift=-2.5em}] node[pos=.6] {(\ga,1)} (b)
  (b) edge node {\ga} (c)
  (a') edge[shorten >=-3ex, shorten <=-1ex,transform canvas={xshift=-2.5em}] node[swap,pos=.6] {(1,\textstyle\bigotimes_j \ga)} (b')
  (b') edge['] node {\ga} (c)
  ;
\end{tikzpicture}
\end{equation}
\item[Symmetric group action]
The identity permutation $\id_n \in \Sigma_n$ acts as the identity morphism on $\M\scmap{\angx;x'}$ whenever $\angx$ has length $n$.  Moreover, for permutations $\sigma,\tau\in\Sigma_n$, the following diagram in $\V$ commutes. 
\begin{equation}\label{enr-multicategory-symmetry}
\begin{tikzcd}
\M\scmap{\angx;x'} \arrow{rd}[swap]{\sigma\tau} \rar{\sigma} 
& \M\scmap{\angx\sigma;x'} \dar{\tau}\\
& \M\mmap{x';\angx\sigma\tau}
\end{tikzcd}
\end{equation}
\item[Equivariance]
Suppose $\len\ang{x_j} = k_j \geq 0$ in \eqref{eq:enr-defn-gamma}.\index{equivariance!enriched multicategory}
\begin{itemize}
\item For each permutation $\sigma \in \Sigma_n$, the following \index{top equivariance!enriched multicategory}\emph{top equivariance diagram} commutes.
\begin{equation}\label{enr-operadic-eq-1}
\begin{tikzcd}[column sep=large,cells={nodes={scale=.8}},
every label/.append style={scale=.8}]
\M\mmap{x'';\ang{x'}} \otimes \txotimes_{j=1}^n \M\mmap{x'_j;\ang{x_j}} 
\dar[swap]{\gamma} \rar{(\sigma, \sigma^{-1})}
& \M\mmap{x'';\ang{x'}\sigma} \otimes \txotimes_{j=1}^n \M\mmap{x'_{\sigma(j)};\ang{x_{\sigma(j)}}} \dar{\gamma}\\
\M\mmap{x'';\ang{x_1},\ldots,\ang{x_n}} \rar{\sigma\langle k_{\sigma(1)}, \ldots , k_{\sigma(n)}\rangle}
& \M\mmap{x'';\ang{x_{\sigma(1)}},\ldots,\ang{x_{\sigma(n)}}}
\end{tikzcd}
\end{equation}
In this diagram, the block permutation\index{block!permutation}\label{notation:enr-block-permutation}  
\begin{equation}\label{blockpermutation}
\sigma\langle k_{\sigma(1)}, \ldots , k_{\sigma(n)} \rangle \in \Sigma_{k_1+\cdots+k_n}
\end{equation}
permutes $n$ consecutive blocks of lengths $k_{\sigma(1)}, \ldots, k_{\sigma(n)}$, respectively, as $\sigma$ permutes $\{1,\ldots,n\}$, leaving the relative order within each block unchanged.
\item
Given permutations $\tau_j \in \Sigma_{k_j}$ for $1 \leq j \leq n$,
the following \index{bottom equivariance!enriched multicategory}\emph{bottom equivariance diagram} commutes.
\begin{equation}\label{enr-operadic-eq-2}
\begin{tikzcd}[cells={nodes={scale=.85}},
every label/.append style={scale=.85}]
\M\mmap{x'';\ang{x'}} \otimes \bigotimes_{j=1}^n \M\mmap{x'_j;\ang{x_j}}
\dar[swap]{\gamma} \rar{(1, \otimes_j \tau_j)} & 
\M\mmap{x'';\ang{x'}} \otimes \bigotimes_{j=1}^n \M\mmap{x'_j;\ang{x_j}\tau_j}\dar{\gamma} \\
\M\mmap{x'';\ang{x_1},\ldots,\ang{x_n}} \rar{\tau_1 \times \cdots \times \tau_n}
& \M\mmap{x'';\ang{x_1}\tau_1,\ldots,\ang{x_n}\tau_n}
\end{tikzcd}
\end{equation}
In this diagram, the block sum\index{block!sum}\label{notation:enr-block-sum} 
\begin{equation}\label{blocksum}
\tau_1 \times\cdots \times\tau_n \in \Sigma_{k_1+\cdots+k_n}
\end{equation} 
is the image of the $n$-tuple $(\tau_1, \ldots, \tau_n)$ under the canonical inclusion 
\[\Sigma_{k_1} \times \cdots \times \Sigma_{k_n} \to \Sigma_{k_1 + \cdots + k_n}.\]
\end{itemize}
\end{description}
This finishes the definition of a $\V$-multicategory.  A $\V$-multicategory is \emph{small}\index{multicategory!enriched!small}\index{small!enriched multicategory} if it has a set of objects.  

Moreover, we have the following variants, terminology, and notation.
\begin{itemize}
\item A \index{multicategory}\emph{multicategory} is a $\Set$-multicategory, where $(\Set,\times,*)$ is the Cartesian closed category of sets and functions.
\item For $n \in \{0,1,2\}$ we also say \emph{nullary} for 0-ary, \emph{unary} for 1-ary, and \emph{binary} for 2-ary.
\item If $\V = \Set$ and $\len\angx = n$, then an element in $\M\scmap{\angx;x'}$ is called an \emph{$n$-ary multimorphism} and denoted $\angx \to x'$.  For multimorphisms $f \cn \angxp \to x''$ and $f_j \cn \ang{x_j} \to x_j'$ for $1 \leq j \leq n$, their composition is also denoted by
\[\angx = \big(\ang{x_1}, \ldots, \ang{x_n}\big) \fto{(f_1,\ldots,f_n)} \angxp \fto{f} x''.\]
\item If $\V = \Cat$, then each $\M\smscmap{\angx;x'}$ is called an \emph{$n$-ary multimorphism category}.  Its objects are called \emph{$n$-ary 1-cells} and denoted $\angx \to x'$.  Its morphisms are called \emph{$n$-ary 2-cells} and denoted using the following 2-cell notation. 
\begin{equation}\label{narytwocell}
\begin{tikzpicture}[baseline={(x1.base)}]
\def\a{25}
\draw[0cell]
(0,0) node (x1) {\angx}
(x1)++(1.7,0) node (x2) {\phantom{\angx}}
(x2)++(-.1,.04) node (x2') {x'}
;
\draw[1cell=.9]  
(x1) edge[bend left=\a] node {} (x2)
(x1) edge[bend right=\a] node[swap] {} (x2)
;
\draw[2cell]
node[between=x1 and x2 at .5, rotate=-90] {\Rightarrow}
;
\end{tikzpicture}
\end{equation}
\item A \index{enriched operad}\index{enriched!operad}\index{operad!enriched}\emph{$\V$-operad}, also called an \emph{operad in $\V$}, is a $\V$-multicategory with one object, usually denoted $*$.  For a $\V$-operad $\M$, its $n$-ary multimorphism object is denoted by \label{not:nthobject}$\M_n$ or $\M(n)$, and the unit $\opu$ in \cref{ccoloredunit} and the composition $\ga$ in \cref{eq:enr-defn-gamma} are called, respectively, the \emph{operadic unit} and the \emph{operadic composition}.  A $\V$-operad $\M$ is \emph{reduced}\index{operad!reduced}\index{reduced operad} if $\M(0)$ is a terminal object in $\V$.  An \emph{operad}\index{multicategory!one object} is a $\Set$-operad.  
\end{itemize}   
This finishes the definition.
\end{definition}

\begin{example}[Underlying $\V$-Categories]\label{ex:unarycategory}\index{enriched multicategory!underlying enriched category}\index{enriched category!underlying - of an enriched multicategory}
Each $\V$-multicategory $(\M,\ga,\operadunit)$ has an underlying $\V$-category (\cref{def:enriched-category}) with
\begin{itemize}
\item the same class of objects as $\M$,
\item hom $\V$-objects given by $\M\smscmap{x;y}$ for objects $x,y \in \M$,
\item identities given by the colored units in $\M$, and
\item composition given by the restriction of $\ga$ to unary multimorphism objects. 
\end{itemize}  
For example, if $\V = (\Cat,\times,\boldone)$ (\cref{ex:cat}), then each $\Cat$-multicategory has an underlying $\Cat$-category, which is the same thing as a locally small 2-category (\cref{locallysmalltwocat}).
\end{example}

\begin{definition}\label{def:enr-multicategory-functor}
Suppose $(\M,\ga,\opu)$ and $(\N,\ga,\opu)$ are $\V$-multicategories.  A \emph{$\V$-multifunctor}\index{multifunctor!enriched}\index{enriched!multifunctor}\index{functor!multi-} 
\[\M \fto{F} \N\]
consists of 
\begin{itemize}
\item an \emph{object assignment}
\[\ObM \fto{F} \ObN\]
and
\item for each $\smscmap{\angx;y} \in \ProfMM$ with $\angx= \ang{x_j}_{j=1}^n$,  a \emph{component $\V$-morphism}
\begin{equation}\label{multifunctor_component}
\begin{tikzcd}[column sep=large]
\M\mmap{y;\ang{x}} \ar{r}{F_{(\angx \sscs y)}} & \N\mmap{Fy;F\ang{x}}
\end{tikzcd}
\end{equation}
where $F\angx = \ang{Fx_j}_{j=1}^n$.  We usually abbreviate $F_{(\angx \sscs y)}$ to $F_n$ or $F$.  When $\V=\Cat$, we call $F_{(\angx \sscs y)}$ an \emph{$n$-ary multimorphism functor} of $F$.
\end{itemize}
The data above are subject to the axioms \cref{enr-multifunctor-equivariance,enr-multifunctor-unit,v-multifunctor-composition} below.
\begin{description}
\item[Unity] 
The following diagram in $\V$ commutes for each $x \in \ObM$.
\begin{equation}\label{enr-multifunctor-unit}
\begin{tikzpicture}[x=25mm,y=15mm,vcenter]
  \draw[0cell] 
  (0,0) node (a) {\tu}
  (1,.5) node (b) {\M\mmap{x;x}}
  (1,-.5) node (b') {\N\mmap{Fx;Fx}}
  ;
  \draw[1cell] 
  (a) edge node {\operadunit_x} (b)
  (a) edge node[swap,pos=.6] {\operadunit_{Fx}} (b')
  (b) edge node {F} (b')
  ;
\end{tikzpicture}
\end{equation} 
\item[Composition] 
For objects $x''$, $\ang{x'}$, $\ang{x_j}$, and $\ang{x} = \oplus_{j=1}^n\ang{x_j}$ as in \eqref{eq:enr-defn-gamma}, the following diagram in $\V$ commutes.
\begin{equation}\label{v-multifunctor-composition}
\begin{tikzcd}[column sep=large,cells={nodes={scale=.8}},
every label/.append style={scale=.9}]
\M\mmap{x'';\ang{x'}} \otimes \bigotimes_{j=1}^n \M\mmap{x'_j;\ang{x_j}} \dar[swap]{F \,\otimes\, \bigotimes_{j=1}^n F} \ar{r}{\gamma} & \M\mmap{x'';\ang{x}} \ar{d}{F} \\  
\N\mmap{Fx'';F\ang{x'}} \otimes \bigotimes_{j=1}^n \N\mmap{Fx'_j;F\ang{x_j}} \ar{r}{\gamma}
 & \N\mmap{Fx'';F\ang{x}}
\end{tikzcd}
\end{equation}
\item[Symmetric group action] 
For each $\smscmap{\angx;y} \in \ProfMM$ and permutation $\sigma \in \Sigma_n$, the following diagram\index{equivariance!enriched multifunctor} in $\V$ commutes.
\begin{equation}\label{enr-multifunctor-equivariance}
\begin{tikzcd}[column sep=large]
\M\mmap{y;\ang{x}} \ar{d}[swap]{\sigma} \ar{r}{F} 
& \N\mmap{Fy;F\ang{x}} \ar{d}{\sigma}\\
\M\mmap{y;\ang{x}\sigma} \ar{r}{F} & \N\mmap{y;F\ang{x}\sigma}\end{tikzcd}
\end{equation}
\end{description}
This finishes the definition of a $\V$-multifunctor.  

Moreover, we define the following.
\begin{itemize}
\item A $\V$-multifunctor $F \cn \M \to \N$ is also called an \emph{$\M$-algebra} or an \emph{algebra over $\M$} in $\N$.\index{algebra}\index{operad!algebra}\index{multicategory!algebra}
\item For a $\V$-multifunctor $H \cn \N\to\P$, the \index{composition!enriched multifunctor}\emph{composition} 
\begin{equation}\label{multifunctors_compose}
\M\fto{HF}\P
\end{equation}
is the $\V$-multifunctor with object assignment given by the composite function 
\[\begin{tikzcd} \ObM \ar{r}{F} & \ObN \ar{r}{H} & \ObP
\end{tikzcd}\]
and component morphisms given by the composites
\begin{equation}\label{composite-multifunctor}
\begin{tikzcd}
\M\mmap{y;\ang{x}} \ar{r}{F} & \N\mmap{Fy;F\ang{x}} \ar{r}{H} & \P\mmap{GFy;GF\ang{x}}.
\end{tikzcd}
\end{equation}
\item The \index{identity!enriched multifunctor}\emph{identity $\V$-multifunctor} $1_{\M} \cn \M\to\M$ is defined by the identity object assignment and identity component morphisms.
\item A \emph{$\V$-operad morphism}\index{operad!morphism}\index{enriched operad!morphism}\index{morphism!enriched operad} is a $\V$-multifunctor between $\V$-multicategories with one object.
\item A \emph{nonsymmetric $\V$-multifunctor}\index{nonsymmetric!multifunctor}\index{multifunctor!nonsymmetric} has the same data as a $\V$-multifunctor, and it preserves units \cref{enr-multifunctor-unit} and composition \cref{v-multifunctor-composition}.  The symmetric group action axiom \cref{enr-multifunctor-equivariance} is not required.  Composition and identities are defined as above.
\item A \emph{multifunctor} is a $\Set$-multifunctor.\defmark
\end{itemize}
\end{definition}

\begin{example}[Underlying $\V$-Functors]\label{ex:un-v-functor}
Each $\V$-multifunctor restricts to a $\V$-functor between underlying $\V$-categories (\cref{ex:unarycategory}).  In particular, for $\V = (\Cat,\times,\boldone)$, each $\Cat$-multifunctor between $\Cat$-multicategories restricts to a 2-functor between underlying 2-categories.
\end{example}

\subsection*{Enriched Multicategories from Enriched Symmetric Monoidal Categories}

Suppose $(\V,\otimes,\tu,\xi)$ is a symmetric monoidal category.  To review the $\V$-multicategory associated to a symmetric monoidal $\V$-category, we use the following notation.

\begin{convention}\label{convention:norm}\index{convention!left normalized product}
For a monoidal $\V$-category $(\C,\vmtimes)$ (\cref{definition:monoidal-vcat}) and an $n$-tuple $\ang{c} = \ang{c_j}_{j=1}^n$ of objects in $\C$ for some $n \geq 0$, we define the \index{left normalized!product}\index{normalized!left - product}\emph{left normalized product} as the object
\[\bigvmtimes \ang{c} = \bigvmtimes_{j=1}^n c_j = \big(\cdots \, \big((c_1 \vmtimes c_2) \vmtimes c_3 \big) \, \cdots \big) \vmtimes c_n,\]
which is the identity object of $\C$ if $n=0$.
\end{convention}

\begin{definition}\label{definition:EndK}
Suppose $\V$ is a symmetric monoidal category and $(\C,\beta^\vmtimes)$ is a symmetric monoidal $\V$-category (\cref{definition:monoidal-vcat,definition:braided-monoidal-vcat,definition:symm-monoidal-vcat}).  The \index{endomorphism!enriched multicategory}\index{enriched multicategory!endomorphism}\emph{endomorphism $\V$-multicategory} of $\C $, denoted $\End(\C)$, is defined as follows.
\begin{description}
\item[Objects] $\Ob(\End(\C)) = \Ob\C$.
\item[Multimorphism objects] 
For each tuple $(\ang{c}; c') \in \Prof(\C) \times \C$, using \cref{convention:norm}, we define the $\V$-object
\[\End(\C)(\ang{c}; c') = \C\big(\! \bigvmtimes \ang{c}, c'\big).\]
\item[Symmetric group action]
For a permutation $\sigma \in \Si_n$, the right $\sigma$-action is defined as the following composite in $\V$.
\begin{equation}\label{EndC_symmetry}
      \ \qquad
      \begin{tikzpicture}[x=50mm,y=15mm,vcenter]
      \def\s{.9}
        \draw[0cell=\s] 
        (0,0) node (a) {
          \C \big(\! \bigvmtimes \ang{c}, c' \big)
        }
        (0,-1) node (b) {
          \C \big(\! \bigvmtimes \ang{c}, c' \big)
          \otimes
          \tensorunit
        }
        (1,-1) node (c) {
          \C \big(\! \bigvmtimes \ang{c}, c' \big)
          \otimes
          \C \big(\! \bigvmtimes \ang{c}\si, \bigvmtimes \ang{c} \big)
        }
        (1,0) node (d) {
          \C \big(\! \bigvmtimes \ang{c}\si, c' \big)
        }
        ;
        \draw[1cell=\s] 
        (a) edge['] node {\rho^\inv} (b)
        (b) edge node {1 \otimes \beta^\vmtimes_\si} (c)
        (c) edge['] node {\mcomp} (d)
        (a) edge node {\si} (d)
        ;
      \end{tikzpicture}
\end{equation}
    In this diagram, $\beta^\vmtimes_\si$ denotes the $\V$-natural
    isomorphism that permutes coordinates according to $\si$.
\item[Units] The $c$-colored unit of an object $c \in \End(\C)$ is the identity
\begin{equation}\label{EndC_unit}
\begin{tikzcd}[column sep=large]
\tensorunit \ar{r}{\cone_c} & \C(c,c) = \End(\C)\mmap{c;c}.
\end{tikzcd}
\end{equation}
\item[Composition] The composition $\ga$ in $\End(\C)$ is defined as the following composite in $\V$ for tuples of objects
\[c'', \quad \ang{c'} = \ang{c'_j}_{j=1}^n \scs \andspace
\ang{c_j} = \ang{c_{j,i}}_{i=1}^{k_j}\] 
with $j\in\{1,\ldots,n\}$ and $\ang{c} = \big(\ang{c_1}, \ldots, \ang{c_n}\big)$. 
\begin{equation}\label{EndC_composition}
\begin{tikzpicture}[x=50mm,y=15mm,vcenter]
\def\s{.9}
        \draw[0cell=\s] 
        (0,0) node (a) {
          \C \big(\! \bigvmtimes\ang{c'},c'' \big) \otimes
          {\textstyle\bigotimes\limits_{j=1}^n} \C \Big( \bigvmtimes_{i=1}^{k_j} c_{j,i} \, , \, c_j' \Big)
        }
        (0,-1) node (b) {
          \C \big(\! \bigvmtimes\ang{c'},c'' \big) \otimes
          \C \Big( \bigvmtimes_{j=1}^n\bigvmtimes_{i=1}^{k_j} c_{j,i} \, , \,
          \bigvmtimes \ang{c'} \Big)
        }
        (1,-1) node (c) {
          \C \Big( \bigvmtimes_{j=1}^n\bigvmtimes_{i=1}^{k_j} c_{j,i} \, , \, c'' \Big)
        }
        (1,0) node (d) {
          \C \big(\! \bigvmtimes \ang{c},c'' \big)
        }
        ;
        \draw[1cell=\s] 
        (a) edge[transform canvas={xshift={-2ex}}, shorten <=-1ex, shorten >=-1ex] node[swap] {1 \otimes \vmtimes^{n-1}} (b)
        (b) edge node {\mcomp} (c)
        (c) edge node['] {\iso} (d)
        (a) edge node {\ga} (d)
        ;
\end{tikzpicture}
\end{equation}
\end{description}
$\End(\C)$ is usually abbreviated to $\End\C$ or $\C$.
\end{definition}

The following result is proved in \cite[6.3.6]{cerberusIII}.

\begin{proposition}\label{proposition:monoidal-v-cat-v-multicat}
In the context of \cref{definition:EndK}, $\End(\C)$ is a $\V$-multicategory.
\end{proposition}

Recall $\V$-multifunctors in \cref{def:enr-multicategory-functor}.

\begin{definition}\label{def:EndF}
In the context of \cref{definition:EndK}, suppose 
\[(F,F^2,F^0) \cn \C \to \D\]
is a symmetric monoidal $\V$-functor between symmetric monoidal $\V$-categories \pcref{definition:monoidal-V-fun,definition:braided-monoidal-vfunctor}.  The \emph{endomorphism $\V$-multifunctor}\index{endomorphism!enriched multifunctor}\index{enriched multifunctor!endomorphism} of $F$, denoted 
\[\End(\C) \fto{\End F} \End(\D),\]
is defined as follows.
\begin{itemize}
\item $\End F$ has the same object assignment as $F$.
\item For each $(n+1)$-tuple $(\ang{c}; c') \in  \Prof(\C) \times \C$, the component $\V$-morphism of $\End F$ is the following composite.
\begin{equation}\label{EndF_component}
\begin{tikzpicture}[vcenter]
\def\h{4.3} \def\v{-1.4} \def\s{.8}
\draw[0cell=\s]
(0,0) node (a11) {\C(\bigvmtimes \ang{c}, c')}
(a11)++(\h,0) node (a12) {\D(\bigvmtimes \ang{Fc}, Fc')}
(a11)++(0,\v) node (a21) {\D(F(\bigvmtimes \ang{c}), Fc')}
(a21)++(0,\v) node (a31) {\D\big(F(\bigvmtimes \ang{c}), Fc'\big) \otimes \tu}
(a31)++(\h,-\v) node (a4) {\D\big(F(\bigvmtimes \ang{c}), Fc'\big) \otimes \D\big(\!\bigvmtimes \ang{Fc}, F(\bigvmtimes \ang{c})\big)}
;
\draw[1cell=\s]
(a11) edge node {(\End F)_{(\ang{c}; c')}} (a12)
(a11) edge node[swap] {F_{\vmtimes \ang{c}, c'}} (a21)
(a21) edge node[swap] {\rhoinv} (a31)
(a4) edge node[swap] {\mcomp} (a12)
;
\draw[1cell=\s]
(a31) [rounded corners=2pt] -| node[pos=.2] {1 \otimes F^n} (a4)
;
\end{tikzpicture}
\end{equation}
In the diagram \cref{EndF_component}, $F^n$ denotes
\begin{itemize}
\item the unique component of the unit constraint $F^0$ if $n=0$,
\item a component of the identity $\V$-natural transformation $1_F$ if $n=1$, and
\item a component of an iterate of the monoidal constraint $F^2$ if $n>1$.
\end{itemize}
\end{itemize}
$\End F$ is usually abbreviated to $F$.
\end{definition}

The following result is proved in \cite[6.3.10]{cerberusIII}.

\begin{proposition}\label{EndF_multi}
In the context of \cref{def:EndF}, $\End F$ is a $\V$-multifunctor.
\end{proposition}

\section{Enriched Multinatural Transformations}
\label{sec:multi_iicat}

This section reviews
\begin{itemize}
\item $\V$-multinatural transformations (\cref{def:enr-multicat-natural-transformation});
\item vertical and horizontal compositions of $\V$-multinatural transformations (\cref{def:enr-multinatural-composition});
\item the 2-category of small $\V$-multicategories, $\V$-multifunctors, and $\V$-multinatural transformations (\cref{v-multicat-2cat});
\item $\Cat$-multinatural transformations (\cref{expl:catmultitransformation}); and
\item $\Cat$-multiequivalences (\cref{def:catmultiequivalence}).
\end{itemize}
References for this section are \cite[Chapter 3]{yau-multgro} and \cite[Chapter 6]{cerberusIII}.

\begin{definition}\label{def:enr-multicat-natural-transformation}
Suppose $F,K \cn \M\to\N$ are $\V$-multifunctors between $\V$-multicategories.  A \emph{$\V$-multinatural transformation}\index{enriched!multinatural transformation}\index{natural transformation!enriched multi-}\index{multinatural transformation!enriched} $\theta \cn F\to K$ consists of, for each object $x \in \M$, a component morphism in $\V$
\[\begin{tikzcd}[column sep=large]
\tu \ar{r}{\theta_x} & \N\mmap{Kx;Fx}
\end{tikzcd}\]
such that the following \emph{$\V$-naturality diagram} in $\V$
commutes for each $\smscmap{\ang{x};y} \in \ProfMM$ with
$\angx= \ang{x_j}_{j=1}^n$.
\begin{equation}\label{enr-multinat}
\begin{tikzpicture}[x=25mm,y=12mm,vcenter]
  \draw[0cell=.8]
  (.5,0) node (a) {\M\mmap{y;\ang{x}}}
  (.7,1) node (b) {\tu \otimes \M\mmap{y;\ang{x}}}
  (3,1) node (c) {\N\mmap{Ky;Fy} \otimes \N\mmap{Fy;F\ang{x}}}
  (3.2,0) node (d) {\N\mmap{Ky;F\ang{x}}}
  (.7,-1) node (b') {\M\mmap{y;\ang{x}}\otimes \txotimes_{j=1}^n \tu}
  (3,-1) node (c') {\N\mmap{Ky;K\ang{x}} \otimes \txotimes_{j=1}^n \N\mmap{Kx_j;Fx_j}}
  ;
  \draw[1cell=.75] 
  (a) edge node[pos=.3] {\la^\inv} (b)
  (a) edge node[swap,pos=.2] {\rho^\inv} (b')
  (b) edge node {\theta_{y} \otimes F} (c)
  (c) edge node[pos=.7] {\ga} (d)
  (b') edge node {K \otimes \txotimes_{j=1}^n \theta_{x_j}} (c')
  (c') edge['] node[pos=.7] {\ga} (d)
  ;
\end{tikzpicture}
\end{equation}
When $n=0$, the bottom composite in \cref{enr-multinat} is interpreted as
\[\M(\ang{}; y) \fto{K} \N(\ang{}; Ky).\]
Moreover, the \emph{identity $\V$-multinatural transformation}\index{multinatural transformation!identity} $1_F \cn F\to F$ is defined by the $(Fx)$-colored units
\[(1_F)_x = \operadunit_{Fx}\]
for $x\in\ObM$.
\end{definition}

\begin{definition}[Compositions]\label{def:enr-multinatural-composition}
Suppose $\theta \cn F \to K$ is a $\V$-multinatural transformation between $\V$-multifunctors $F,K \cn \M \to \N$.
\begin{description}
\item[Vertical] Suppose $H \cn \M \to \N$ is a $\V$-multifunctor, and $\psi \cn K \to H$ is a $\V$-multinatural transformation.  The \emph{vertical composition}\index{vertical composition!enriched multinatural transformation}\label{notation:enr-operad-vcomp}
\begin{equation}\label{multinatvcomp}
\begin{tikzcd}[column sep=large]
F \ar{r}{\psi\theta} & H
\end{tikzcd}  
\end{equation} 
is the $\V$-multinatural transformation with, for each $x \in \ObM$, $x$-component given by the following composite in $\V$.
  \[
  \begin{tikzpicture}[x=45mm,y=12mm]
    \draw[0cell=.9] 
    (0,0) node (a) {\tu}
    (0,-1) node (b) {\tu \otimes \tu}
    (1,-1) node (c) {\N\mmap{Hx;Kx} \otimes \N\mmap{Kx;Fx}}
    (1,0) node (d) {\N\mmap{Hx;Fx}}
    ;
    \draw[1cell=.9] 
    (a) edge node['] {\la^\inv} (b)
    (b) edge node {\psi_x \otimes \theta_x} (c)
    (c) edge['] node {\ga} (d)
    (a) edge node {(\psi\theta)_x} (d)
    ;
  \end{tikzpicture}
  \]
\item[Horizontal] Suppose $F', K' \cn \N \to \P$ are $\V$-multifunctors, and $\vartheta \cn F' \to K'$ is a $\V$-multinatural transformation.  The \emph{horizontal composition}\index{horizontal composition!enriched multinatural transformation}\label{notation:enr-operad-hcomp}
\begin{equation}\label{multinathcomp}
\begin{tikzcd}[column sep=large]
F'F \ar{r}{\vartheta \ast \theta} & K'K
\end{tikzcd}
\end{equation} 
is the $\V$-multinatural transformation with, for each $x \in \ObM$, $x$-component given by the following composite in $\V$.
\[
\begin{tikzpicture}[x=50mm,y=12mm]
  \draw[0cell=.9] 
  (0,0) node (a) {\tu}
  (0,-2) node (b) {\tu \otimes \tu}
  (1,-2) node (c) {\P\mmap{K'Kx;F'Kx} \otimes \N\mmap{Kx;Fx}}
  (1,-1) node (d) {\P\mmap{K'Kx;F'Kx} \otimes \P\mmap{F'Kx;F'Fx}}
  (1,0) node (e) {\P\mmap{K'Kx;F'Fx}}
  ;
  \draw[1cell=.9] 
  (a) edge node {(\vartheta * \theta)_x} (e)
  (a) edge['] node {\la^\inv} (b)
  (b) edge node {\vartheta_{Kx} \otimes \theta_x} (c)
  (c) edge['] node {1 \otimes F'} (d)
  (d) edge['] node {\ga} (e)
  ;
\end{tikzpicture}
\]
\end{description}
This finishes the definition.
\end{definition}

\begin{theorem}\label{v-multicat-2cat}
For each symmetric monoidal category $\V$, there is a 2-category\index{2-category!of small enriched multicategories}\index{enriched multicategory!2-category}\index{multicategory!enriched!2-category} 
\[\VMulticat\]
consisting of the following data.
\begin{itemize}
\item The objects are small $\V$-multicategories.
\item For small $\V$-multicategories $\M$ and $\N$, the hom category 
\[\VMulticat(\M,\N)\]
is given by the following data.
\begin{itemize}
\item The objects are $\V$-multifunctors $\M \to \N$.
\item The morphisms are $\V$-multinatural transformations.
\item Composition is given by vertical composition of $\V$-multinatural transformations.
\item Identity morphisms are identity $\V$-multinatural transformations.
\end{itemize}
\item Identity 1-cells are identity $\V$-multifunctors.
\item Horizontal composition of 1-cells is the composition of $\V$-multifunctors.
\item Horizontal composition of 2-cells is the horizontal composition of $\V$-multinatural transformations.
\end{itemize}
\end{theorem}

In \cref{v-multicat-2cat}, the smallness of $\V$-multicategories are only needed to ensure that, given two $\V$-multifunctors $\M \to \N$, there is a set of $\V$-multinatural transformations between them.

\subsection*{Enrichment in $\Cat$}
\begin{explanation}[$\Cat$-Multinatural Transformations]\label{expl:catmultitransformation}
Suppose $\M$ and $\N$ are $\Cat$-multicategories, and $F,H \cn \M \to \N$ are $\Cat$-multifunctors.  Interpreting \cref{def:enr-multicat-natural-transformation} in the case $\V = \Cat$, a \index{multinatural transformation!Cat-@$\Cat$-}\index{Cat-multinatural@$\Cat$-multinatural!transformation}\emph{$\Cat$-multinatural transformation} 
\[F \fto{\theta} H\]
consists of, for each $x \in \Ob\M$, a component 1-ary 1-cell in $\N$
\begin{equation}\label{thetaccomponent}
\begin{tikzcd}[column sep=normal]
Fx \ar{r}{\theta_x} & Hx
\end{tikzcd}
\end{equation}
such that the following two \index{Cat-naturality conditions@$\Cat$-naturality conditions}naturality conditions hold for $k \geq 0$.
\begin{description}
\item[1-cell naturality] For each $k$-ary 1-cell $p \cn \angx \to y$ in $\M$ with
$\angx= \ang{x_j}_{j=1}^k$, we denote by
\begin{equation}\label{Fangcthetaangc}
\left\{\begin{split}
F\angx &= \ang{F x_j}_{j=1}^k \in (\Ob\N)^k\\
\theta_{\angx} &= \ang{\theta_{x_j}}_{j=1}^k \in \txprod_{j=1}^k \N\mmap{H x_j; F x_j}.
\end{split}\right.
\end{equation}
Then the following $k$-ary 1-cell equality holds, with composition taken in the $\Cat$-multicategory $\N$:
\begin{equation}\label{catmultinaturality}
\gamma\scmap{Hp; \theta_{\angx}} =
\gamma\scmap{\theta_{y}; Fp} \inspace \N\mmap{Hy;F\angx}.
\end{equation}
\item[2-cell naturality] For each $k$-ary 2-cell $f \cn p \to q$ in $\M\mmap{y;\angx}$, the following $k$-ary 2-cell equality holds, where $1_{\theta_{\angx}} = \ang{1_{\theta_{x_j}}}_{j=1}^k$:
\begin{equation}\label{catmultinaturalityiicell}
\gamma\scmap{Hf; 1_{\theta_{\angx}}} =
\gamma\scmap{1_{\theta_{y}}; Ff} \inspace \N\mmap{Hy;F\angx}.
\end{equation}
\end{description}
The 1-cell naturality and 2-cell naturality conditions, \cref{catmultinaturality,catmultinaturalityiicell}, together constitute the $\V$-naturality condition \cref{enr-multinat} when $\V$ is $\Cat$.  In the 0-ary case, with $k=0$, the left-hand sides of \cref{catmultinaturality,catmultinaturalityiicell} are interpreted as, respectively, $Hp$ and $Hf$.

A \emph{$\Cat$-multinatural isomorphism}\index{multinatural isomorphism!Cat-@$\Cat$-}\index{Cat-multinatural@$\Cat$-multinatural!isomorphism} is an invertible $\Cat$-multinatural transformation.  More explicitly, a $\Cat$-multinatural transformation $\theta$ is a $\Cat$-multinatural isomorphism if and only if each component 1-ary 1-cell $\theta_x$ \cref{thetaccomponent} is an isomorphism in the underlying 1-category of $\N$ (\cref{ex:unarycategory}).
\end{explanation}

The notions of equivalences, adjunctions, and adjoint equivalences make sense in any 2-category \cite[Chapter 6]{johnson-yau}.  In particular, these notions exist for the 2-category of small $\Cat$-multicategories, $\Cat$-multifunctors, and $\Cat$-multinatural transformations (\cref{v-multicat-2cat}).  The following definition makes these concepts explicit and removes the smallness requirement at the same time.

\begin{definition}\label{def:catmultiequivalence}
A $\Cat$-multifunctor $L \cn \M \to \N$ between $\Cat$-multicategories is a \emph{$\Cat$-multiequivalence} if there exist
\begin{itemize}
\item a $\Cat$-multifunctor $R \cn \N \to \M$ and
\item $\Cat$-multinatural isomorphisms (\cref{expl:catmultitransformation})
\begin{equation}\label{catmultiequnit}
\begin{tikzcd}[column sep=normal]
1_{\M} \ar{r}{\sfu}[swap]{\iso} & RL
\end{tikzcd} \andspace 
\begin{tikzcd}[column sep=normal]
LR \ar{r}{\sfv}[swap]{\iso} & 1_{\N}.
\end{tikzcd}
\end{equation}
\end{itemize}
Moreover, the quadruple $(L,R,\sfu,\sfv)$ is called an \emph{adjoint $\Cat$-multiequivalence} if it restricts to an adjunction between the underlying 1-categories of $\M$ and $\N$.  In this case, we call $L$ the \emph{left adjoint}, $R$ the \emph{right adjoint}, $\sfu$ the \emph{unit}, and $\sfv$ the \emph{counit}.  
\end{definition}

\begin{explanation}[Adjoint $\Cat$-Multiequivalences]\label{expl:catmultieq}
For the 2-category of small $\Cat$-multicategories, $\Cat$-multifunctors, and $\Cat$-multinatural transformations, each of the two triangle identities for an adjunction \cite[Definition 6.1.1]{johnson-yau} asserts the equality between two $\Cat$-multinatural transformations.  Each $\Cat$-multinatural transformation is determined by its 1-ary 1-cells \cref{thetaccomponent}.  Thus, given $\Cat$-multifunctors
\[\begin{tikzpicture}
\def\t{20}
\draw[0cell]
(0,0) node (a1) {\M}
(a1)++(1.7,0) node (a2) {\N}
;
\draw[1cell=.9]
(a1) edge[bend left=\t] node {L} (a2)
(a2) edge[bend left=\t] node {R} (a1)
;
\end{tikzpicture}\]
and $\Cat$-multinatural transformations
\[1_\M \fto{\sfu} RL \andspace LR \fto{\sfv} 1_\N,\]
the triangle identities defining an adjunction assert the commutativity of the following diagrams for objects $x \in \M$ and $y \in \N$.
\begin{equation}\label{triangleid}
\begin{tikzpicture}[vcenter]
\def\h{2} \def\v{-1.2}
\draw[0cell]
(0,0) node (a1) {Lx}
(a1)++(\h,0) node (a2) {LRLx}
(a2)++(0,\v) node (a3) {Lx}
(a2)++(\h,0) node (b1) {Ry}
(b1)++(\h,0) node (b2) {RLRy}
(b2)++(0,\v) node (b3) {Ry}
;
\draw[1cell=.9]
(a1) edge node {L\sfu_x} (a2)
(a2) edge node {\sfv_{Lx}} (a3)
(a1) edge node[swap] {1_{Lx}} (a3)
(b1) edge node {\sfu_{Ry}} (b2)
(b2) edge[shorten <=-.5ex] node {R\sfv_{y}} (b3)
(b1) edge node[swap] {1_{Ry}} (b3)
;
\end{tikzpicture}
\end{equation}
The left and right diagrams in \cref{triangleid} happen in the underlying 1-categories of, respectively, $\N$ and $\M$.  For an adjoint $\Cat$-multiequivalence, it is, furthermore, required that the components of $\sfu$ and $\sfv$ be isomorphisms in the underlying 1-categories of, respectively, $\M$ and $\N$. 
\end{explanation}

\section{Change of Enrichment}
\label{sec:mult_change_enr}

This section reviews change of enrichment of enriched multicategories \pcref{mult_change_enr}.   \cref{End_dashf} observes that change of enrichment commutes with the endomorphism construction.  Recall the 2-category $\VMulticat$ in \cref{v-multicat-2cat}.  The following observation is proved in \cite[6.2.9]{cerberusIII}; see also \cite[11.5.1]{yau-operad}.

\begin{proposition}\label{mult_change_enr}
For each symmetric monoidal functor 
\[\V \fto{(f,f^2,f^0)} \W\]
between symmetric monoidal categories, there is a change-of-enrichment\index{change of enrichment} 2-functor
\[\VMulticat \fto{\dashf} \WMulticat.\]
\end{proposition}

If there is no danger of confusion, we sometimes suppress $\dashf$ from the notation.

\begin{explanation}[Unraveling $\dashf$]\label{expl:mult_change_enr}
In the context of \cref{mult_change_enr}, given a $\V$-multicategory $(\M,\ga,\opu)$ \pcref{def:enr-multicategory}, the $\W$-multicategory $(\M_f,\ga_f,\opu_f)$ is defined as follows.
\begin{itemize}
\item $\M_f$ has the same objects as $\M$.
\item The $n$-ary multimorphism $\W$-objects are given by
\[\M_f(\angx;x') = f\M(\angx; x').\]
\item The symmetric group action on $\M_f$ is given by the image under $f$ of the symmetric group action on $\M$.
\item For each object $x \in \M_f$, the $x$-colored unit is the composite
\begin{equation}\label{Mf_unit}
\tu \fto{f^0} f\tu \fto{f\opu_x} f\M(x;x)
\end{equation}
of
\begin{itemize}
\item the unit constraint $f^0$ of $f$ \cref{monoidal_unit} and
\item the image under $f$ of the $x$-colored unit of $x \in \M$.
\end{itemize}
\item Using the notation in \cref{eq:enr-defn-gamma}, the composition in $\M_f$ is given by the following composite in $\W$.
\begin{equation}\label{Mf_gamma}
\begin{tikzpicture}[vcenter]
\draw[0cell=.9]
(0,0) node (a1) {f\M(\angxp ; x'') \otimes \txotimes_{j=1}^n f\M(\ang{x_j}; x_j')}
(a1)++(0,-1.4) node (a2) {f\big[ \M(\angxp ; x'') \otimes \txotimes_{j=1}^n \M(\ang{x_j}; x_j') \big]}
(a1)++(4.5,0) node (a3) {\phantom{f\M(\angx; x'')}}
(a3)++(0,.04) node (a3') {f\M(\angx; x'')}
;
\draw[1cell=.85]
(a1) edge node {\ga_f} (a3)
(a1) edge node[swap] {f^{n+1}} (a2)
(a2) [rounded corners=2pt] -| node[pos=.2] {f\ga} (a3')
;
\end{tikzpicture}
\end{equation}
In the above composite, $f^{n+1}$ denotes a component of an iterate of the monoidal constraint of $f$ \cref{monoidal_constraint}.
\end{itemize}

For a $\V$-multifunctor $F \cn \M \to \N$ \pcref{def:enr-multicategory-functor}, the $\W$-multifunctor
\[\M_f \fto{F_f} \N_f\]
has the same object assignment as $F$.  Its component $\W$-morphisms are given by
\[f\M(\angx; y) \fto{fF_{(\angx; y)}} f\N(F\angx; Fy).\]
For a $\V$-multinatural transformation $\theta \cn F \to K$ \pcref{def:enr-multicat-natural-transformation} between $\V$-multifunctors $F,K \cn \M \to \N$, the $\W$-multinatural transformation
\[F_f \fto{\theta_f} K_f\]
has, for each object $x \in \M$, $x$-component $\W$-morphism given by the composite
\[\tu \fto{f^0} f\tu \fto{f\theta_x} f\N(Fx; Kx)\]
of $f^0$ with the image of $\theta_x$ under $f$.
\end{explanation}

\subsection*{Change of Enrichment Commutes with $\End$}
Recall symmetric monoidal $\V$-categories \pcref{definition:monoidal-vcat,definition:braided-monoidal-vcat,definition:symm-monoidal-vcat} and their change of enrichment \pcref{thm:change-enrichment}.  The following observation says that change of enrichment for enriched symmetric monoidal categories and for enriched multicategories are compatible via the endomorphism construction \pcref{definition:EndK,def:EndF}.

\begin{proposition}\label{End_dashf}
Suppose 
\[\V \fto{(f,f^2,f^0)} \W\]
is a symmetric monoidal functor between symmetric monoidal categories.
\begin{enumerate}
\item\label{End_dashf_i} For each symmetric monoidal $\V$-category $\C$, there is an equality
\begin{equation}\label{EndC_f}
(\End\C)_f = \End(\C_f)
\end{equation}
of $\W$-multicategories.
\item\label{End_dashf_ii} For each symmetric monoidal $\V$-functor $F \cn \C \to \D$,  there is an equality
\begin{equation}\label{EndF_f}
(\End F)_f = \End(F_f)
\end{equation}
of $\W$-multifunctors.
\end{enumerate}
\end{proposition}

\begin{proof}
\parhead{Assertion \pcref{End_dashf_i}}. 
The equality \cref{EndC_f} of $\W$-multicategories is given explicitly as follows for each symmetric monoidal $\V$-category $(\C,\vmtimes)$.
\begin{description}
\item[Left-hand side] $\End\C$ is the endomorphism $\V$-multicategory of $\C$ \pcref{definition:EndK}, and $(\End\C)_f$ is the $\W$-multicategory obtained from $\End\C$ by changing enrichment along $f$ \pcref{mult_change_enr}.
\item[Right-hand side] $\C_f$ is the symmetric monoidal $\W$-category obtained from $\C$ by changing enrichment along $f$
(\cref{thm:change-enrichment} \pcref{change-enr-i}), and $\End(\C_f)$ is the endomorphism $\W$-multicategory of $\C_f$ \pcref{definition:EndK}.
\item[Objects] Each side of \cref{EndC_f} has the same objects as $\C$.
\item[Multimorphisms] For each $(n+1)$-tuple $(\angx; y) \in \Prof(\C) \times \C$, each side of \cref{EndC_f} has an $n$-ary multimorphism $\W$-object given by
\[f\C(\bigvmtimes \angx, y).\]
\item[Units] By \cref{EndC_unit,Mf_unit}, for each object $x \in \C$, its $x$-colored unit on either side of \cref{EndC_f} is the composite
\[\tu \fto{f^0} f\tu \fto{f\cone_x} f\C(x,x),\]
where $\cone_x$ is the identity of $x \in \C$ \cref{cone_a}.
\item[Symmetry] Using the notation in \cref{EndC_symmetry} and the functoriality of $f$ and $\otimes$ in $\W$, the right $\si$-actions of the two sides of \cref{EndC_f} are given by the boundary composites in the following diagram in $\W$.
\[\begin{tikzpicture}[vcenter]
\def\h{4.4} \def\v{-1.4} \def\t{1ex} \def\s{.8}
\draw[0cell=.8]
(0,0) node (a11) {f\C(\bigvmtimes \angc, c')}
(a11)++(.5*\h,-.5*\v) node (a12) {f\C(\bigvmtimes \angc, c') \otimes \tu}
(a11)++(0,\v) node (a21) {f\big( \C(\bigvmtimes \angc, c') \otimes \tu\big)}
(a12)++(.5*\h,.5*\v) node (a22) {f\C(\bigvmtimes \angc, c') \otimes f\tu}
(a21)++(0,\v) node (a31) {f\big[\C(\bigvmtimes \angc, c') \otimes \C(\bigvmtimes \angc\si, \bigvmtimes \angc) \big]}
(a22)++(0,\v) node (a32) {f\C(\bigvmtimes \angc, c') \otimes f\C(\bigvmtimes \angc\si, \bigvmtimes \angc)}
(a32)++(0,\v) node (a4) {f\C(\bigvmtimes \angc \si, c')}
;
\draw[1cell=.8]
(a11) edge node[swap] {f\rhoinv} (a21)
(a21) edge node[swap] {f(1 \otimes \be^{\vmtimes}_\si)} (a31)
(a31) edge node {f\mcomp} (a4)
;
\draw[1cell=\s]
(a11) [rounded corners=2pt] |- node[pos=.7] {\rhoinv} (a12)
;
\draw[1cell=\s]
(a12) [rounded corners=2pt] -|node[pos=.3] {1 \otimes f^0} (a22)
;
\draw[1cell=\s]
(a22) edge node[swap] {f^2} (a21)
(a22) edge node {f1 \otimes f\be^{\vmtimes}_\si} (a32)
(a32) edge node[swap] {f^2} (a31)
;
\end{tikzpicture}\]
In the diagram above, the top quadrilateral commutes by the right unity axiom for $(f,f^2,f^0)$ \cref{monoidalfunctorunity}.  The lower quadrilateral commutes by the naturality of $f^2$
\item[Composition] Using the notation in \cref{EndC_composition,Mf_gamma}, the functoriality of $f$, and the functoriality of $\otimes$ in $\W$, the composition of the two sides of \cref{EndC_f} are given by the boundary composites in the following diagram in $\W$, where the indices $j$ and $i$ run through, respectively, $\ufs{n}$ and $\ufs{k}_j$.
\[\begin{tikzpicture}[vcenter]
\def\h{4} \def\v{-1.5} \def\t{-2em} \def\s{.8}
\draw[0cell=\s]
(0,0) node (a11) {f\C(\bigvmtimes \angcp, c'') \otimes \txotimes_j f\C(\txvmt_i c_{j,i} \spc c_j')}
(a11)++(\h,.5*\v) node (a12) {f\C(\bigvmtimes \angcp, c'') \otimes f\big[ \txotimes_j \C(\txvmt_i c_{j,i} \spc c_j')\big]}
(a11)++(0,\v) node (a21) {f\big[ \C(\bigvmtimes \angcp, c'') \otimes \txotimes_j \C(\txvmt_i c_{j,i} \spc c_j') \big]}
(a12)++(0,\v) node (a22) {f\C(\bigvmtimes \angcp, c'') \otimes f\C(\txvmt_{j,i} c_{j,i} \spc \txvmt \angcp)}
(a21)++(0,\v) node (a31) {f\big[\C(\bigvmtimes \angcp, c'') \otimes \C(\txvmt_{j,i} c_{j,i} \spc \bigvmtimes \angcp) \big]}
(a22)++(0,\v) node (a32) {f\C(\txvmt_{j,i} c_{j,i} \spc c'')}
;
\draw[1cell=\s]
(a11) edge node[swap] {f^{n+1}} (a21)
(a21) edge node[swap] {f(1 \otimes \vmtimes^{n-1})} (a31)
;
\draw[1cell=\s]
(a31) [rounded corners=2pt] |- node[pos=.8] {f\mcomp} (a32)
;
\draw[1cell=\s]
(a11) [rounded corners=2pt] -| node[pos=.2] {1 \otimes f^{n}} (a12)
;
\draw[1cell=\s]
(a12) edge[transform canvas={xshift=\t}, shorten <=-2ex] node[swap] {f^2} (a21)
(a12) edge node {1 \otimes f(\vmtimes^{n-1})} (a22)
;
\draw[1cell=\s]
(a22) [rounded corners=2pt] |- node[swap,pos=.8] {f^2} (a31)
;
\end{tikzpicture}\]
In the diagram above, the top region commutes by the associativity axiom for $(f,f^2)$ \cref{monoidalfunctorassoc}.  The middle region commutes by the naturality of $f^2$.
\end{description}
This finishes the proof of the equality \cref{EndC_f} of $\W$-multicategories.

\parhead{Assertion \pcref{End_dashf_ii}}. 
Each of $(\End F)_f$ and $\End(F_f)$ has the same object assignment as $F$.  Using the notation in \cref{EndF_component} and the functoriality of $f$ and $\otimes$, the component $\W$-morphisms of $(\End F)_f$ and $\End(F_f)$ are the two boundary composites of the following diagram in $\W$. 
\[\begin{tikzpicture}[vcenter]
\def\h{4.4} \def\v{-1.4} \def\t{1ex} \def\s{.75}
\draw[0cell=.8]
(0,0) node (a11) {f\D(F(\bigvmtimes \angc), Fc')}
(a11)++(-1,1) node (a0) {f\C(\bigvmtimes \angc, c')}
(a11)++(.5*\h,-.5*\v) node (a12) {f\D(F(\bigvmtimes \angc), Fc') \otimes \tu}
(a11)++(0,\v) node (a21) {f\big[ \D(F(\bigvmtimes \angc), Fc') \otimes \tu\big]}
(a12)++(.5*\h,.5*\v) node (a22) {f\D(F(\bigvmtimes \angc), Fc') \otimes f\tu}
(a21)++(0,\v) node (a31) {\phantom{f\big[ \D(F(\bigvmtimes \angc), Fc') \otimes \tu\big]}}
(a31)++(1.2,0) node (a31') {f\big[\D(F(\bigvmtimes \angc), Fc') \otimes \D(\bigvmtimes \ang{Fc}, F(\bigvmtimes \angc)) \big]}
(a22)++(0,\v) node (a32) {f\D(F(\bigvmtimes \angc), Fc') \otimes f\D(\bigvmtimes \ang{Fc}, F(\bigvmtimes \angc))}
(a31')++(\h,0) node (a4) {f\D(\bigvmtimes \ang{Fc}, Fc')}
;
\draw[1cell=\s]
(a0) edge node[swap,pos=.2,inner sep=0pt] {f F_{\vmtimes \ang{c}, c'}} (a11)
(a11) edge node[swap] {f\rhoinv} (a21)
(a21) edge node[swap] {f(1 \otimes F^n)} (a31)
(a31') edge node {f\mcomp} (a4)
;
\draw[1cell=\s]
(a11) [rounded corners=2pt] |- node[pos=.8] {\rhoinv} (a12)
;
\draw[1cell=\s]
(a12) [rounded corners=2pt] -|node[pos=.3] {1 \otimes f^0} (a22)
;
\draw[1cell=\s]
(a22) edge node[swap] {f^2} (a21)
(a22) edge node {f1 \otimes f F^n} (a32)
(a32) edge node[swap] {f^2} (a31)
;
\end{tikzpicture}\]
In the diagram above, the top quadrilateral commutes by the right unity axiom for $(f,f^2,f^0)$ \cref{monoidalfunctorunity}.  The lower quadrilateral commutes by the naturality of $f^2$.  This finishes the proof of the equality \cref{EndF_f} of $\W$-multifunctors.
\end{proof}

\section{Pointed Diagram Categories}
\label{sec:pointed-diagrams}

This section reviews pointed diagram categories. 
\begin{itemize}
\item \cref{thm:Dgm-pv-convolution-hom} states that the category of pointed functors and pointed natural transformations is complete, cocomplete, and symmetric monoidal closed, with the monoidal structure given by the pointed Day convolution.
\item \cref{theorem:diagram-omnibus} states that this pointed diagram category is an enriched symmetric monoidal category and also an enriched multicategory.  This enriched multicategory of pointed diagrams is explained in more detail in \cref{sec:dstarv-multicat}.
\end{itemize}
Detailed discussion of this material can be found in \cite[Chapters 1--4]{cerberusIII}.

\subsection*{Symmetric Monoidal Closed Category of Pointed Objects}

\begin{definition}\label{def:pointed-objects}
Suppose $\A$ is a category and $\term$ is a terminal object in $\A$.  We define the category $\pA$ of pointed objects as the category under $\term$.  
\begin{itemize}
\item More explicitly, an object in $\pA$, called a \index{object!pointed}\index{pointed!object}\emph{pointed object}, is a pair $(a,i^a)$\label{not:aia} consisting of an object $a \in \A$ and a morphism $i^a \cn \term \to a$ in $\A$, called the \index{basepoint}\emph{basepoint} or the \index{pointed!structure}\emph{pointed structure}.  We abbreviate the pointed object $(\term, 1_{\term})$ to $\term$.
\item A morphism between pointed objects
\[(a,i^a) \fto{h} (b,i^b) \inspace \pA\]
is a morphism $h \cn a \to b$ in $\A$ such that the following diagram in $\A$ commutes.
\begin{equation}\label{ptmorphism-diagram}
\begin{tikzpicture}[xscale=1,yscale=1,vcenter]
\draw[0cell]
(0,0) node (t) {\term}
(t)++(2,.7) node (a) {a}
(t)++(2,-.7) node (b) {b}
;
\draw[1cell]
(t) edge node[pos=.6] {i^a} (a)
(t) edge node[swap,pos=.6] {i^b} (b)
(a) edge node {h} (b)
;
\end{tikzpicture}
\end{equation} 
A morphism in $\pA$ is called a \index{morphism!pointed}\index{pointed!morphism}\emph{pointed morphism}.  Composition and identity pointed morphisms are defined in $\A$.
\end{itemize} 
This finishes the definition of $\pA$.
\end{definition}

\begin{explanation}\label{expl:pAbicomplete}
If $\A$ is complete and cocomplete, then so is the category $\pA$ of pointed objects, with $\term$ as a terminal object.  We emphasize that $\term$ is \emph{not} required to be an initial object in $\A$.  For example, if $\A = \Cat$, then we can choose $\term$ to be the terminal category $\boldone$.  If $\A = \Set$, then we can choose $\term$ to be the one-element set $\{*\}$.
\end{explanation}

\begin{definition}\label{def:wedge-smash-phom}
Suppose $(\C,\otimes,\tu,[,])$ is a complete and cocomplete symmetric monoidal closed category (\cref{def:closedcat}) together with a chosen terminal object $\term$.  Suppose $(a,i^a)$ and $(b,i^b)$ are pointed objects.  We define the following pointed objects.
\begin{itemize}
\item The \index{wedge}\emph{wedge} $a \wed b$ is the pushout in $\C$ of the diagram\label{not:awedgeb}
\[a \xleftarrow{i^a} \term \fto{i^b} b.\]
Its pointed structure is given by the composite 
\[\term \fto{i^a} a \to a \wed b.\]
\item The \index{product!smash}\index{smash product}\emph{smash product} $a \sma b$ is defined as the following pushout in $\C$.
\begin{equation}\label{eq:smash}
\begin{tikzpicture}[x=50mm,y=15mm,vcenter]
      \draw[0cell] 
      (0,0) node (a) {(a \otimes \term) \bincoprod (\term \otimes b)}
      (1,0) node (b) {a \otimes b}
      (0,-1) node (c) {\term}
      (1,-1) node (d) {a \sma b}
      ;
      \draw[1cell=.9] 
      (a) edge node {(1_a \otimes i^b) \bincoprod (i^a \otimes 1_b)} (b)
      (c) edge node {i^{a \sma b}} (d)
      (a) edge node {} (c)
      (b) edge node {\pjt_{a,b}} (d)
      ;
\end{tikzpicture}
\end{equation}
\item The \index{smash unit}\emph{smash unit} $\stu$ is defined as the coproduct in $\C$
\begin{equation}\label{smash-unit-object}
\stu = \tu \bincoprod \term.
\end{equation}
Its pointed structure is given by the inclusion of the $\term$ summand.
\item The \index{hom!pointed}\index{pointed!hom}\emph{pointed hom} $[a,b]_*$ is defined as the following pullback in $\C$.
  \begin{equation}\label{eq:pHom-pullback}
  \begin{tikzpicture}[x=30mm,y=15mm,vcenter]
    \draw[0cell] 
    (0,0) node (a) {[a,b]_*}
    (1,0) node (b) {\term}
    (0,-1) node (c) {[a,b]}
    (1,-1) node (d) {[\term,b]}
    ;
    \draw[1cell] 
    (a) edge node {} (b)
    (c) edge node {[i^a,b]} (d)
    (a) edge node {} (c)
    (b) edge node {} (d)
    ;
  \end{tikzpicture}
  \end{equation}
In the previous diagram, the right vertical arrow is the composite
\[\term \iso [a,\term] \fto{[a,i^b]} [a,b] \fto{[i^a,b]} [\term,b].\]
This induces the pointed structure of $[a,b]_*$.\defmark
\end{itemize}
\end{definition}

A detailed proof of the following observation from \cite[4.20]{elmendorf-mandell-perm} in given in \cite[4.2.3]{cerberusIII}.

\begin{theorem}\label{theorem:pC-sm-closed}
Suppose $\C$ is a complete and cocomplete symmetric monoidal closed category with a chosen terminal object.  Then the quadruple in \cref{def:pointed-objects,def:wedge-smash-phom}
\[\left(\pC,\sma,\stu,[,]_*\right)\]
is a complete and cocomplete symmetric monoidal closed category.
\end{theorem}


\subsection*{Pointed Unitary Enrichment}

We use the following definition when the category in question is the indexing category of some diagram category.  Thus, we denote it by $\Dgm$.

\begin{definition}\label{definition:zero}
A \index{object!zero}\index{zero!object}\emph{zero object} in a category $\Dgm$ is an initial and terminal object.  For a given zero object $\zob \in \Dgm$, we define the following.
\begin{itemize}
\item A \index{morphism!zero}\index{zero!morphism}\emph{zero morphism} in $\Dgm$ is a morphism that factors through $\zob$.  A \emph{nonzero morphism} is a morphism that does not factor through $\zob$.
\item For objects $a,b \in \Dgm$, we denote the set of nonzero morphisms $a \to b$ by $\Dpunc(a,b)$.
\item Suppose $(\Dgm,\Dtimes,\tu)$ is a symmetric monoidal category.  A \index{null object}\index{object!null}\emph{null object} is a zero object $\zob \in \Dgm$ such that there are natural isomorphisms
\[a \Dtimes \zob \iso \zob \iso \zob \Dtimes a\]
for objects $a \in \Dgm$.\defmark
\end{itemize}
\end{definition}


\begin{definition}\label{def:unitary-enrichment}
Suppose we are given the following data.
\begin{itemize}
\item $(\V,\otimes,\tu,[,])$ is a complete and cocomplete symmetric monoidal closed category, and $\termv$ is a terminal object in $\V$, which we use to define the category $\pV$ of pointed objects (\cref{def:pointed-objects}).
\item $(\Dgm,\Dtimes,\tu)$ is a small permutative category, and $\zob$ is a null object in $\Dgm$.
\end{itemize} 
The \index{pointed!unitary enrichment}\index{unitary enrichment!pointed}\emph{pointed unitary enrichment} of $\Dgm$ over $\pV$, denoted $\Dhat$, is the $\pV$-category (\cref{def:enriched-category}) defined as follows.
\begin{itemize}
\item $\Ob\Dhat = \Ob\Dgm$.
\item For objects $a,b \in \Dhat$, the hom $\pV$-object is defined as the wedge
\begin{equation}\label{eq:Dhat-ptd-unitary-enr}
\Dhat(a,b) = \bigvee_{\Dpunc(a,b)} \stu
\end{equation}
of copies of the smash unit $\stu$ in \cref{smash-unit-object}, indexed by the set of nonzero morphisms $\Dpunc(a,b)$.  An empty wedge means $\termv \in \pV$.  
\item The identity of an object $a \in \Dhat$ is given by the copy of $\stu$ corresponding to the identity of $a \in \Dgm$.
\item Composition is induced by the isomorphism $\stu \sma \stu \iso \stu$ and the composition in $\Dgm$.\defmark
\end{itemize}
\end{definition}

By \cite[2.4.10]{cerberusIII}, $\Dhat$ is a symmetric monoidal $\pV$-category (\cref{definition:symm-monoidal-vcat}).

\subsection*{Pointed Day Convolution}

\begin{definition}\label{def:pointed-category}
A \index{category!pointed}\index{pointed!category}\emph{pointed category} is a pair $(\D,\bpt)$ consisting of a category $\D$ and a chosen object $\bpt \in \D$, called the \index{basepoint}\emph{basepoint}.  We emphasize that the basepoint $\bpt$ is not required to be a terminal object or a zero object.
\begin{itemize}
\item A \index{functor!pointed}\index{pointed!functor}\emph{pointed functor}
\[(\D,\bpt) \fto{F} (\C,\apt)\]
between pointed categories is a functor $F \cn \D \to \C$ such that $F\bpt = \apt$.  
\item A \index{natural transformation!pointed}\index{pointed!natural transformation}\emph{pointed natural transformation}
\[F \fto{\theta} H\]
between pointed functors $F,H \cn (\D,\bpt) \to (\C,\apt)$ is a natural transformation such that
\[F\bpt = \apt \fto{\theta_{\bpt} = 1_{\apt}} H\bpt = \apt \inspace \D.\]
\item Suppose $\D$ is small.  The \index{diagram!category!pointed}\index{pointed!diagram category}\emph{pointed diagram category} $\DstarC$\label{not:DstarC} is the category with
\begin{itemize}
\item pointed functors $(\D,\bpt) \to (\C,\apt)$ as objects,
\item pointed natural transformations as morphisms, and
\item composition and identities defined as those of natural transformations.
\end{itemize}
\end{itemize}
Suppose $\D$ is small with $\bpt \in \D$ an initial object, and $\apt \in \C$ is a terminal object, which is used to define the category $\pC$ of pointed objects \pcref{def:pointed-objects}.
\begin{itemize}
\item For a pointed functor $F \cn (\D,\bpt) \to (\C,\apt)$ and an object $x \in \D$, the object $Fx \in \C$ is equipped with the \emph{canonical basepoint} 
\[\apt = F\bpt \to Fx\]
given by the $F$-image of the unique morphism $\bpt \to x$.  By the functoriality of $F$, the $F$-image of each morphism in $\D$ becomes a pointed morphism in $\C$.  We regard $F$ equivalently as a pointed functor
\[(\D,\bpt) \fto{F} (\pC,\apt).\]
\item For pointed functors $F, H \cn (\D,\bpt) \to (\C,\apt)$ and a pointed natural transformation $\theta \cn F \to H$, the pointed naturality of $\theta$ implies that each component of $\theta$ is a pointed morphism.  We regard $\theta$ equivalently as a pointed natural transformation
\[\begin{tikzpicture}[vcenter]
\def\a{26}
\draw[0cell] 
(0,0) node (a1) {\phantom{X}}
(a1)++(1.8,0) node (a2) {\phantom{X}}
(a1)++(-.3,0) node (a1') {(\D,\bpt)}
(a2)++(.4,0) node (a2') {(\pC,\apt).}
;
\draw[1cell=.9]
(a1) edge[bend left=\a] node {F} (a2) 
(a1) edge[bend right=\a] node[swap] {H} (a2) 
;
\draw[2cell]
node[between=a1 and a2 at .43, rotate=-90, 2label={above,\theta}] {\Rightarrow}
;
\end{tikzpicture}\]
\item The pointed diagram category $\DstarC$ is equivalently regarded as the category with
\begin{itemize}
\item pointed functors $(\D,\bpt) \to (\pC,\apt)$ as objects and
\item pointed natural transformations between such pointed functors as morphisms.\defmark
\end{itemize}
\end{itemize}
\end{definition}

In the following definition, an empty wedge means the chosen terminal object.

\begin{definition}\label{definition:Dgm-pV-convolution-hom}
In the context of \cref{def:unitary-enrichment,def:pointed-category}, we define the following pointed functors from $(\Dgm,\zob)$ to $(\pV,\termv)$.
\begin{itemize}
\item The \index{pointed!diagrams!unit diagram}\index{unit diagram!pointed}\emph{monoidal unit diagram} is the pointed functor
\begin{equation}\label{ptdayunit}
\du = \bigvee_{\Dpunc(\tu,-)} \stu \cn (\Dgm,\zob) \to (\pV,\termv)
\end{equation}
with $\stu = \tu \bincoprod \termv$ denoting the smash unit in \cref{smash-unit-object}.  The wedge is indexed by the set $\Dpunc(\tu, -)$ of nonzero morphisms from the monoidal unit $\tu \in \Dgm$.
\item The \index{pointed!diagrams!Day convolution}\index{Day convolution!pointed diagrams}\emph{pointed Day convolution} of pointed functors
\[F, H \cn (\Dgm,\zob) \to (\pV,\termv)\]
is the $\pV$-coend
\begin{equation}\label{ptdayconv}
F \sma H = \ecint^{(c,d) \in \Dhat \otimes \Dhat}
\bigvee_{\Dpunc(c \Dtimes d, -)} (Fc \sma Hd)
\end{equation}
with $\Dpunc(c \Dtimes d, -)$ denoting the set of nonzero morphisms from $c \Dtimes d$.  The value of $F \sma H$ at $\zob$ is declared to be the chosen terminal object $\termv \in \pV$.
\item The \index{pointed!diagrams!hom diagram}\index{hom diagram!pointed diagrams}\emph{pointed hom diagram} is the $\pV$-end
\begin{equation}\label{ptdayhom}
[F,H]_* = \ecint_{d \in \Dhat}\, \big[Fd \scs H(- \Dtimes d)\big]_*
\end{equation}
where $[,]_*$ on the right-hand side is the pointed hom \cref{eq:pHom-pullback} in $\pV$.  The value of this $\pV$-end at $\zob$ is declared to be $\termv$.
\item The \index{mapping object!pointed diagrams}\index{pointed!diagrams!mapping object}\emph{pointed mapping object} is the $\pV$-end
\begin{equation}\label{ptdaymap}
\pMap(F,H) = \ecint_{d \in \Dhat}\, [Fd, Hd]_* \iso [F,H]_*(\tu)
\end{equation}
with $\tu$ denoting the monoidal unit in $\Dgm$.
\end{itemize}
We extend each of \cref{ptdayconv,ptdayhom,ptdaymap} componentwise to pointed natural transformations.
\end{definition}

\begin{explanation}[Pointed Day Convolution]\label{expl:pointedday}\index{Ga-category@$\Ga$-category!pointed Day convolution}\index{pointed!Day convolution}\index{Day convolution!pointed}
Unraveling the $\pV$-coend, the pointed Day convolution $F \sma H$ in \cref{ptdayconv} can be characterized externally by the following universal property, without any mention of coends.  For pointed functors 
\[F, H, K \cn (\Dgm,\zob) \to (\pV,\termv),\] 
a pointed natural transformation
\[\begin{tikzcd}[column sep=large]
F \sma H \ar{r}{\zeta} & K
\end{tikzcd}\]
consists of component $\pV$-morphisms
\[\begin{tikzcd}[column sep=huge]
Fc \sma Hd \ar{r}{\zeta_{c,d}} & K(c \Dtimes d)
\end{tikzcd} \forspace (c,d) \in \Dgm^2\]
such that, for each pair of $\Dgm$-morphisms
\[\begin{tikzcd}
c \ar{r}{f} & c'
\end{tikzcd} \andspace 
\begin{tikzcd}
d \ar{r}{g} & d'
\end{tikzcd}\]
the diagram 
\begin{equation}\label{ptdaynaturality}
\begin{tikzpicture}[xscale=3,yscale=1.2,vcenter]
\draw[0cell=.9]
(0,0) node (x11) {Fc \sma Hd}
(x11)++(1,0) node (x12) {K(c \Dtimes d)}
(x11)++(0,-1) node (x21) {Fc' \sma Hd'}
(x12)++(0,-1) node (x22) {K(c' \Dtimes d')}
;
\draw[1cell=.9] 
(x11) edge node {\zeta_{c,d}} (x12)
(x21) edge node {\zeta_{c',d'}} (x22)
(x11) edge node[swap] {Ff \sma Hg} (x21)
(x12) edge node {K(f \Dtimes g)} (x22)
;
\end{tikzpicture}
\end{equation}
in $\pV$ commutes.
\end{explanation}

\begin{explanation}[Pointed Mapping Objects]\label{expl:Fpointedmap}\index{pointed!mapping object}\index{mapping object!pointed}
For pointed functors 
\[F, H \cn (\Dgm,\zob) \to (\pV,\termv),\] 
by unraveling the defining $\pV$-end, the pointed mapping object $\pMap(F,H)$ in \cref{ptdaymap} can be characterized as follows.  With $\stu = \tu \bincoprod \termv$ denoting the smash unit in $\pV$ \cref{smash-unit-object}, a $\pV$-morphism
\[\stu \fto{\theta} \pMap(F,H) = \ecint_{d \in \Dhat}\, [Fd, Hd]_*\]
is precisely a pointed natural transformation (\cref{def:pointed-category})
\[F \fto{\theta} H.\]
Moreover, the pointed condition for $\theta$, which states
\[\theta_{\zob} = 1_{\termv} \cn F\zob = \termv \to H\zob = \termv,\]
is automatic because $\termv \in \V$ is a terminal object.  Therefore, $\theta$ is equivalently a natural transformation $F \to H$ between the underlying functors.
\end{explanation}

A detailed proof of the following observation, which is originally from \cite{day-convolution}, is given in \cite[3.7.22 and 4.3.37]{cerberusIII}.

\begin{theorem}\label{thm:Dgm-pv-convolution-hom}
In the context of \cref{def:pointed-category,def:unitary-enrichment,definition:Dgm-pV-convolution-hom}, the quadruple
\[\left( \DstarV, \sma, \du, [,]_* \right)\]
is a \index{pointed!diagram category!complete and cocomplete}complete and cocomplete \index{pointed!diagram category!symmetric monoidal closed}symmetric monoidal closed category.
\end{theorem}

The following result states that, in the context of \cref{thm:Dgm-pv-convolution-hom}, each terminal object-preserving symmetric monoidal functor induces a change-of-base symmetric monoidal functor.  The proof is given in \cite[3.8.4 and 4.3.34]{cerberusIII}.

\begin{theorem}\label{thm:Dgm_f}
Consider the following context.
\begin{itemize}
\item $(\V,\termv)$ and $(\W,\termw)$ are complete and cocomplete symmetric monoidal closed categories with chosen terminal objects $\termv$ and $\termw$, which are used to define the categories $\pV$ and $\pW$ of pointed objects \pcref{def:pointed-objects}.
\item $(\Dgm,\zob)$ is a small permutative category with a null object $\zob$, which is used to define the pointed diagram categories $\DstarV$ and $\DstarW$ \pcref{def:pointed-category}.
\item $f \cn \V \to \W$ is a symmetric monoidal functor such that $f(\termv) = \termw$.
\end{itemize}
Then post-composing and post-whiskering with $f$ induce a symmetric monoidal functor
\[\DstarV \fto{f_*} \DstarW\]
between pointed diagram categories.
\end{theorem}

By \cref{thm:Dgm-pv-convolution-hom,theorem:v-closed-v-sm}, the category $\DstarV$ of pointed functors and pointed natural transformations is a symmetric monoidal $(\DstarV)$-category.  Evaluation at the monoidal unit $\tu \in \Dgm$ defines a symmetric monoidal functor\index{evaluation!at $\pu$}
\begin{equation}\label{evtu}
\DstarV \fto{\ev_{\tu}} \pV.
\end{equation}  
It admits a strong symmetric monoidal left adjoint \index{evaluation!at $\pu$!left adjoint}$\ladj$.  A detailed proof of the following change-of-enrichment result is given in \cite[3.8.1, 3.9.15, and 4.3.37]{cerberusIII}.

\begin{theorem}\label{theorem:diagram-omnibus}
In the context of \cref{thm:Dgm-pv-convolution-hom}, the adjunction 
\begin{equation}\label{eq:pV-DstarV-adj}
\begin{tikzpicture}[baseline={(x.base)}]
\def\a{20} \def\h{1.8}
\draw[0cell] 
(0,0) node (x) {\pV}
(x)++(\h,0) node (y) {\phantom{\pV}}
(y)++(.2,0) node (y') {\DstarV}
;
\draw[0cell=.8]
(x)++(\h/2,0) node (b) {\bot}
;
\draw[1cell=.9]
(x) edge[bend left=\a] node {\ladj} (y) 
(y) edge[bend left=\a] node {\ev_{\tu}} (x) 
;
\end{tikzpicture}
\end{equation}
makes $\DstarV$ \index{pointed!diagram category!enriched}enriched, \index{pointed!diagram category!tensored and cotensored}tensored, and cotensored over $\pV$, with mapping objects given by $\pMap$ in \cref{ptdaymap}.  Changing enrichment along $\ev_{\tu}$, $\DstarV$ becomes a symmetric monoidal $\pV$-category by \cref{thm:change-enrichment} \pcref{change-enr-i}, and a $\pV$-multicategory by \cref{proposition:monoidal-v-cat-v-multicat}.
\end{theorem}

The forgetful functor $U \cn \pV \to \V$ that forgets the basepoint of each pointed object is symmetric monoidal.
\begin{itemize}
\item The unit constraint of $U$ is the $\V$-morphism
\[\tu \to \stu = \tu \bincoprod \term\]
given by the inclusion of the coproduct summand $\tu$ in the smash unit \cref{smash-unit-object}.
\item The monoidal constraint of $U$ at two pointed objects $(a,i^a$) and $(b,i^b)$ is the $\V$-morphism
\[a \otimes b \fto{\pjt_{a,b}} a \sma b\] 
in the definition \cref{eq:smash} of the smash product.
\end{itemize}
The composite symmetric monoidal functor
\begin{equation}\label{evtu_unpt}
\DstarV \fto{\ev_{\tu}} \pV \fto{U} \V
\end{equation}
is also denoted by $\ev_{\tu}$.  Changing enrichment along this $\ev_{\tu}$, $\DstarV$ becomes a symmetric monoidal $\V$-category by \cref{thm:change-enrichment} \pcref{change-enr-i}, and a $\V$-multicategory by \cref{proposition:monoidal-v-cat-v-multicat}.

\section{Enriched Multicategory of Pointed Diagrams}\label{sec:dstarv-multicat}

In this section, we unravel the $\pV$-multicategory $\DstarV$ in \cref{theorem:diagram-omnibus}, where $\V$ and $\Dgm$ are specified in \cref{def:unitary-enrichment}. 

\subsection*{Objects} 

The objects in the $\pV$-multicategory $\DstarV$ are pointed functors 
\[(\Dgm,\zob) \to (\pV,\termv)\]
as in \cref{def:pointed-category}.   These are functors $\Dgm \to \pV$ that send the null object $\zob \in \Dgm$ to the terminal object $\termv \in \pV$.

\subsection*{Multimorphism $\pV$-Objects}

For $n \geq 0$ and pointed functors
\[F_1, \ldots, F_n, H \cn (\Dgm,\zob) \to (\pV,\termv),\]
the multimorphism $\pV$-object with input profile $\angF = \ang{F_j}_{j=1}^n$ and output object $H$ is the pointed mapping object \cref{ptdaymap}
\begin{equation}\label{pmapsmafh}
\begin{split}
(\DstarV)(\angF \sscs H) 
&= \pMap\big(\!\wedge_{j=1}^n F_j \scs H\big) \\
&= \ecint_{d \in \Dhat}\, \big\lbrack\big(\!\wedge_{j=1}^n F_j\big)d \scs Hd\big\rbrack_*.
\end{split}
\end{equation}
By \cref{expl:Fpointedmap}, a $\pV$-morphism from the smash unit
\[\stu \fto{\theta} \pMap\big(\!\wedge_{j=1}^n F_j \scs H\big)\]
is precisely a pointed natural transformation
\begin{equation}\label{theta-smafjh}
\wedge_{j=1}^n F_j \fto{\theta} H.
\end{equation}
By \cref{expl:pointedday}, such a pointed natural transformation $\theta$ consists of component $\pV$-morphisms
\[\begin{tikzcd}[column sep=huge]
\wedge_{j=1}^n F_j c_j  \ar{r}{\theta_{c_1,\ldots,c_n}} & H(c_1 \Dtimes \cdots \Dtimes c_n)
\end{tikzcd} \forspace (c_1,\ldots,c_n) \in \Dgm^n\]
such that, for any $n$-tuple of $\Dgm$-morphisms
\[\begin{tikzcd}
c_j \ar{r}{f_j} & c'_j
\end{tikzcd} 
\forspace j \in \{1,\ldots,n\},\]
the following diagram in $\pV$ commutes. 
\begin{equation}\label{ptdaynat-n}
\begin{tikzpicture}[xscale=4.5,yscale=1.4,vcenter]
\draw[0cell=.9]
(0,0) node (x11) {\wedge_{j=1}^n F_j c_j}
(x11)++(1,0) node (x12) {H(c_1 \Dtimes \cdots \Dtimes c_n)}
(x11)++(0,-1) node (x21) {\wedge_{j=1}^n F_j c_j'}
(x12)++(0,-1) node (x22) {H(c_1' \Dtimes \cdots \Dtimes c_n')}
;
\draw[1cell=.9] 
(x11) edge node {\theta_{c_1,\ldots,c_n}} (x12)
(x21) edge node {\theta_{c_1',\ldots,c_n'}} (x22)
(x11) edge node[swap] {\sma_j F_j f_j} (x21)
(x12) edge[transform canvas={xshift=-2em}] node {H(f_1 \Dtimes \cdots \Dtimes f_n)} (x22)
;
\end{tikzpicture}
\end{equation}

If $n=0$, then $\wedge_{j=1}^n F_j$ is interpreted as the monoidal unit diagram $\du$ in \cref{ptdayunit}.  Thus, the 0-ary multimorphism $\pV$-object with output object $H$ is the pointed mapping object
\begin{equation}\label{pmapjh}
(\DstarV)(\du;H) = \pMap\left(\du, H\right) = \ecint_{d \in \Dhat}\, \left\lbrack \du d, Hd \right\rbrack_*.
\end{equation}

\subsection*{Units}

The $H$-colored unit
\begin{equation}\label{dstarvcoloredunit}
\stu \fto{\opu_H} \pMap(H,H)
\end{equation}
is given by the identity natural transformation $1_H \cn H \to H$.

\subsection*{Symmetric Group Action}

For a permutation $\sigma \in \Sigma_n$, the right $\sigma$-action
\begin{equation}\label{dstarvsigmaaction}
\begin{tikzcd}[column sep=large]
\pMap\left(\wedge_{j=1}^n F_j \scs H\right) \ar{r}{\sigma}[swap]{\iso} &
\pMap\left(\wedge_{j=1}^n F_{\sigma(j)} \scs H\right)
\end{tikzcd}
\end{equation}
is the $\pV$-isomorphism given as follows.  Considering a pointed natural transformation $\theta \cn \sma_{j=1}^n F_j \to H$ as in \cref{theta-smafjh}, the right $\sigma$-action sends $\theta$ to the pointed natural transformation
\begin{equation}\label{thetasigma-smafjh}
\wedge_{j=1}^n F_{\sigma(j)} \fto{\theta^\sigma} H
\end{equation}
with, for each $n$-tuple $\ang{c_j}_{j=1}^n \in \Dgm^n$ of objects, a component $\pV$-morphism given by the following composite.
\begin{equation}\label{thetasigmacomponent}
\begin{tikzpicture}[xscale=1,yscale=1.4,vcenter]
\draw[0cell=.9]
(0,0) node (x11) {\txsma_{j=1}^n F_{\sigma(j)} c_j}
(x11)++(5.5,0) node (x12) {H\left(c_1 \Dtimes \cdots \Dtimes c_n\right)}
(x11)++(0,-1) node (x21) {\txsma_{j=1}^n F_j c_{\sigmainv(j)}}
(x12)++(0,-1) node (x22) {H\left(c_{\sigmainv(1)} \Dtimes \cdots \Dtimes c_{\sigmainv(n)}\right)}
;
\draw[1cell=.9]  
(x11) edge node {\theta^\sigma_{c_1,\ldots,c_n}} (x12)
(x11) edge node {\iso} node[swap] {\sigma} (x21)
(x21) edge node {\theta_{c_{\sigmainv(1)}, \ldots, c_{\sigmainv(n)}}} (x22)
(x22) edge[shorten <=-.4ex] node {\iso} node [swap] {H(\sigmainv)} (x12)
;
\end{tikzpicture}
\end{equation}
The left, bottom, and right arrows in \cref{thetasigmacomponent} are given as follows.
\begin{itemize}
\item The left vertical arrow is the coherence isomorphism in $\pV$ that permutes the smash factors according to $\sigma$.
\item The bottom horizontal arrow is the component of $\theta$ for the $n$-tuple of objects $\ang{c_{\sigmainv(j)}}_{j=1}^n$ in $\Dgm$.
\item The right vertical arrow is $H$ applied to the coherence isomorphism in $\Dgm$ that permutes the $\Dtimes$-factors according to $\sigmainv$. 
\end{itemize}

\subsection*{Multicategorical Composition}

For each $j \in \{1,\ldots,n\}$ with $n \geq 1$, suppose given $k_j \geq 0$ pointed functors
\[F_j^1, \ldots, F_j^{k_j} \cn (\Dgm,\zob) \to (\pV,\termv).\]
The composition $\pV$-morphism 
\begin{equation}\label{dstarv-multicomp}
\begin{tikzpicture}[xscale=1,yscale=1.4,vcenter]
\draw[0cell=.85]
(0,0) node (x11) {\pMap\left(\sma_{j=1}^n F_j \scs H \right) \sma \bigwedge_{j=1}^n \pMap\left(\sma_{i=1}^{k_j} F_j^i \scs F_j\right)}
(x11)++(5.6,0) node (x12) {\pMap\left(\sma_{j=1}^n \sma_{i=1}^{k_j} F_j^i \scs H \right)}
;
\draw[1cell=.9]  
(x11) edge node {\ga} (x12)
;
\end{tikzpicture}
\end{equation}
sends pointed natural transformations \cref{theta-smafjh}
\[\sma_{j=1}^n F_j \fto{\theta} H \andspace \sma_{i=1}^{k_j} F_j^i \fto{\theta_j} F_j\]
to the pointed natural transformation
\[\ga\left(\theta; \ang{\theta_j}_{j=1}^n\right) \cn \sma_{j=1}^n \sma_{i=1}^{k_j} F_j^i \to H\]
with the following components.  For objects $c_j^i \in \Dgm$ with $j \in \{1,\ldots,n\}$ and $i \in \{1,\ldots,k_j$\}, we use the notation below.
\[\begin{aligned}
c_j &= c_j^1 \Dtimes \cdots \Dtimes c_j^{k_j} \in \Dgm & \ang{c_j} &= \left(c_j^1, \ldots, c_j^{k_j}\right) \in \Dgm^{k_j}\\
\angc &= \left(\ang{c_1}, \ldots, \ang{c_n}\right) \in \Dgm^{k_1+\cdots+k_n} &&\\
\end{aligned}\]
Then the $\angc$-component $\pV$-morphism is the following composite. 
\begin{equation}\label{dstarv-ga-component}
\begin{tikzpicture}[xscale=1,yscale=1.5,vcenter]
\draw[0cell=.9]
(0,0) node (x11) {\sma_{j=1}^n \sma_{i=1}^{k_j} F_j^i c_j^i}
(x11)++(5,0) node (x12) {H\left(\Dtimes_{j=1}^n c_j\right)}
(x11)++(0,-1) node (x21) {\sma_{j=1}^n \left(\sma_{i=1}^{k_j} F_j^i c_j^i\right)}
(x12)++(0,-1) node (x22) {\sma_{j=1}^n F_j c_j}
;
\draw[1cell=.9]  
(x11) edge node {\ga(\theta; \ang{\theta_j}_{j=1}^n)_{\angc}} (x12)
(x11) edge node[swap] {\iso} (x21)
(x21) edge node {\sma_{j=1}^n (\theta_j)_{\ang{c_j}}} (x22)
(x22) edge[shorten <=-.3ex] node [swap] {\theta_{c_1,\ldots,c_n}} (x12)
;
\end{tikzpicture}
\end{equation}
This finishes the description of the $\pV$-multicategory $\DstarV$.

\chapter{Classifying Space}
\label{ch:nerve}
This appendix reviews simplicial sets, the nerve functor \cref{nerve}, the geometric realization functor \cref{realization}, and the classifying space functor \cref{classifying_space}.  Much more detailed discussion of these topics can be found in \cite[Ch.\! 7]{cerberusIII} and \cite{gabriel_zisman,goerss-jardine,may_simplicial}.

\subsection*{Simplicial Sets} 
$\Delta$\label{not:Delta_cat} denotes the category
\begin{itemize}
\item whose objects are the ordered sets 
\[\ord{n} = \{0 < 1 < \cdots < n\} \forspace n \geq 0\] 
and
\item whose morphisms are weakly order-preserving functions.  
\end{itemize}
Morphisms of $\Delta$ are generated by\label{not:face_degen}
\begin{align*}
&\ord{n-1} \fto{d^i} \ord{n} \forspace 0 \leq i \leq n \andspace\\
&\ord{n+1} \fto{s^i} \ord{n} \forspace 0 \leq i \leq n,
\end{align*}
called the \index{coface}\index{face!co-}\emph{coface} and \index{degeneracy!co-}\index{codegeneracy}\emph{codegeneracy} morphisms, respectively. 
\begin{itemize}
\item 
The coface $d^i$ is the unique order-preserving injection whose image does not contain $i \in \ord{n}$. 
\item The codegeneracy $s^i$ is the unique order-preserving surjection such that the preimage of $i \in \ord{n}$ is the subset $\{i,i+1\} \subseteq \ord{n+1}$.
\end{itemize}
These morphisms satisfy the following \index{cosimplicial identities}\index{identities!cosimplicial}\index{simplicial identities!co-}\emph{cosimplicial identities}.
\[\left\{
\begin{aligned}
    d^jd^i & = d^i d^{j-1} \ \ifspace i < j\\
    s^j d^i & = d^i s^{j-1} \ \ifspace i < j\\
    s^j d^j & = 1 = s^j d^{j+1} \\
    s^j d^i & = d^{i-1} s^j \ \ifspace i > j+1\\
    s^j s^i & = s^i s^{j+1} \ \ifspace i \leq j\dqed
\end{aligned}
\right.\]
The coface morphisms, codegeneracy morphisms, and cosimplicial identities form a generating set of morphisms and relations for $\Delta$.

$\SSet$\label{not:SSet} denotes the category of functors $\Deltaop \to \Set$ and natural transformations between them. 
\begin{itemize}
\item An object in $\SSet$ is called a \index{simplicial set}\emph{simplicial set}.  
\item For a simplicial set $X \cn \Deltaop \to \Set$, we denote
\begin{itemize}
\item $X\ord{n}$\label{not:n_simplex} by $X_n$, called the set of \index{simplex}\emph{$n$-simplices},
\item $Xd^i$ by $d_i$, called the \index{face}\emph{$i$-th face morphism}, and
\item $Xs^i$ by $s_i$, called the \index{degeneracy}\emph{$i$-th degeneracy morphism}.
\end{itemize}
\item For each $n \geq 0$, the \index{standard $n$-simplex}\index{simplicial set!standard $n$-simplex}\emph{standard $n$-simplex}, denoted $\Delta^n \in \SSet$, is the simplicial set defined by the functor $\Delta(-,\ord{n}) \cn \Deltaop \to \Set$.  
\end{itemize}

\subsection*{Nerve}
Recall that $\Cat$ is the category of small categories and functors.  The \index{nerve}\emph{nerve functor} 
\begin{equation}\label{nerve}
\Cat \fto{\Ner} \SSet
\end{equation}
sends a small category $\C$ to the simplicial set $\Ner\C$ defined as follows.
\begin{itemize}
\item A 0-simplex in $\Ner\C$ is an object in $\C$. 
\item For $n>0$, an $n$-simplex is a diagram of $n$ composable morphisms in $\C$, 
\[c_0 \fto{f_1} c_1 \fto{f_2} c_2 \to \cdots \fto{f_{n-1}} c_{n-1} \fto{f_n} c_n,\]
also denoted by $(f_1,\ldots, f_n)$. 
\item The face morphisms $(\Ner\C)_1 \fto{d_i} (\Ner\C)_0$ are given by
\[\begin{split}
d_0(c_0 \to c_1) &= c_1 \andspace\\ 
d_1(c_0 \to c_1) &= c_0.
\end{split}\]
\item For $n>0$, the face morphisms $(\Ner\C)_n \fto{d_i} (\Ner\C)_{n-1}$ are defined as
\[d_i(f_1,\ldots,f_n) = \begin{cases}
(f_2,\ldots,f_n) & \text{if $i=0$,}\\
(f_1, \ldots, f_{i+1} f_i, \ldots, f_n) & \text{if $0 < i < n$, and}\\
(f_1,\ldots,f_{n-1}) & \text{if $i=n$}.
\end{cases}\]
\item The degeneracy morphisms are defined as
\[\begin{split}
s_0(c_0) &= 1_{c_0} \andspace\\
s_i(f_1,\ldots,f_n) &= (f_1, \ldots, f_i, 1_{c_i}, f_{i+1}, \ldots, f_n).
\end{split}\]
\end{itemize}
For a functor $F \cn \C \to \D$ between small categories, the morphism of simplicial sets
\[\Ner\C \fto{\Ner F} \Ner\D\]
is $\Ob F$ on 0-simplices.  For $n>0$, it sends an $n$-simplex $(f_1,\ldots,f_n)$ in $\Ner\C$ to the $n$-simplex $(Ff_1,\ldots,Ff_n)$ in $\Ner\D$.

\subsection*{Geometric Realization}
We denote by $\Top$ the category of compactly generated weak Hausdorff spaces, which are simply called \index{space}\emph{spaces} throughout this work, and continuous morphisms.  The \index{topological $n$-simplex}\emph{topological $n$-simplex} is the space\label{not:simpn}
\[\simp^n = \big\{(a_0,\ldots,a_n) \cn  0 \leq a_i \leq 1, \txsum_{i=0}^n a_i = 1 \big\} \subseteq \bR^{(n+1)}.\]
The \index{geometric realization}\emph{geometric realization functor} 
\begin{equation}\label{realization}
\SSet \fto{\Rea} \Top
\end{equation}
is defined as the coend
\[|X| = \ecint^{\ord{n} \in \Delta} X_n \times \simp^n.\]

\subsection*{Classifying Space}
The \index{classifying space}\emph{classifying space functor} $\cla$ is defined as the composite
\begin{equation}\label{classifying_space}
\begin{tikzpicture}[baseline={(a1.base)}]
\def\h{2.2}
\draw[0cell]
(0,0) node (a1) {\Cat}
(a1)++(\h,0) node (a2) {\phantom{\SSet}} 
(a2)++(0,.045) node (a2') {\SSet}
(a2)++(\h,0) node (a3) {\Top}
;
\draw[1cell=.9]
(a1) edge node {\Ner} (a2)
(a2) edge node {\Rea} (a3)
;
\draw[1cell=.9]
(a1) [rounded corners=2pt] |- ($(a2)+(-1,.6)$) -- node {\cla} ($(a2)+(1,.6)$) -| (a3)
 ;
\end{tikzpicture}
\end{equation}
of $\Ner$ in \cref{nerve} and $\Rea$ in \cref{realization}.

\begin{explanation}\label{expl:cla_space}
We mention a few properties of $\Ner$, $\Rea$, and $\cla$ here.
\begin{enumerate}
\item\label{expl:cla_space_i} Each of $\Ner$ and $\Rea$ preserves $G$-actions for a group $G$, and, therefore, so does $\cla$.
\item\label{expl:cla_space_ii} Each of $\Ner$ and $\Rea$ commutes, up to a natural isomorphism, with taking $G$-fixed subcategories/points for a group $G$, and, therefore, so does $\cla$.
\item\label{expl:cla_space_iii} A functor between small categories that admits either a left adjoint or a right adjoint, such as an equivalence of categories, is sent by $\Ner$ to a weak equivalence in $\SSet$.
\item\label{expl:cla_space_iv} Geometric realization preserves weak equivalences.  Thus, the classifying space functor $\cla$ sends each functor admitting an adjoint to a weak equivalence in $\Top$.
\item\label{expl:cla_space_v} $\Ner$ is a right adjoint, so it commutes with products.  Geometric realization commuets with finite products.  Thus, $\cla$ commutes with finite products.
\defmark
\end{enumerate}
\end{explanation}

\part*{Bibliography and Indices}

\backmatter
\bibliographystyle{sty/amsalpha3}
\bibliography{references}


\chapter*{List of Notations}
\def\s{.85} \def\v{-.13ex} \def\hor{\kern .1em}
\newcommand{\blob}{\>}

\newcommand{\entry}[3]{\scalebox{\s}{#1} \hor : \scalebox{\s}{#3}, \scalebox{\s}{\pageref{#2}} \blob\\[\v]}
\newcommand{\entryNoPage}[3]{\scalebox{\s}{#1} \hor : \scalebox{\s}{#3} \blob \\[\v]} 
\newcommand{\newchHeader}[1]{\blob\\ \textbf{\Cref{#1}}\blob\\}

\begin{tabbing}
\phantom{\textbf{Notation}} \= \kill

\textbf{Standard Notations}\blob\\ 
\entryNoPage{$\Ob(\C)$, $\Ob\C$}{not:objects}{objects in a category $\C$}
\entryNoPage{$\C(X,Y)$, $\C(X;Y)$}{not:morphisms}{set of morphisms $X \to Y$}
\entryNoPage{$1$}{not:idmorphism}{identity morphism, identity functor, identity natural transformation}
\entryNoPage{$\boldone$}{ex:terminal-category}{terminal category}
\entryNoPage{$h \circ f$, $hf$}{notation:morphism-composition}{composition of morphisms}
\entryNoPage{$\iso$, $\fto{\iso}$}{not:iso}{an isomorphism}
\entryNoPage{$F \cn \C \to \D$}{def:functors}{a functor}
\entryNoPage{$\Rightarrow$}{}{2-cell, natural transformation}
\entryNoPage{$\theta_?$}{thetax}{a component of a natural transformation $\theta$}
\entryNoPage{$\phi\theta$, $\theta' * \theta$}{not:vcomp}{vertical and horizontal compositions of natural transformations}
\entryNoPage{$L \dashv R$}{notation:adjunction}{an adjunction}
\entryNoPage{$\coprod$}{not:coprod}{a coproduct}
\entryNoPage{$\prod$, $\times$}{not:coprod}{a product}
\entryNoPage{$\vee$, $\sma$}{}{wedge and smash}
\entryNoPage{$\txint^?$, $\txint_?$}{}{coend and end}
\entryNoPage{$\Sigma_n$, $\id_n$}{not:Sigman}{symmetric group on $n$ letters, identity permutation}
\entryNoPage{$(\Set,\times,*)$}{notation:set}{category of sets and functions}
\entryNoPage{$(\Cat,\times,\boldone)$}{}{category of small categories and functors}

\newchHeader{ch:psalg}

\entry{$G$}{ch:psalg}{a group}
\entry{$(\Op,\ga,\opu$)}{def:intrinsic_pairing}{an operad}
\entry{$\intr$, $\intr_{j,k}$}{intr_jk}{intrinsic pairing}
\entry{$\Delta^j$, $\Delta$}{intr_jk}{diagonal}
\entry{$\As$}{ex:as_intrinsic}{associative operad}
\entry{$\phi\ang{k_1,\ldots,k_n}$}{as_gamma}{block permutation}
\entry{$\phi_1 \times \cdots \times \phi_n$}{as_gamma}{block sum}
\entry{$\ufs{n}$}{ufsn}{unpointed finite set $\{1,2,\ldots,n\}$}
\entry{$\txsum_{j \in \ufs{n}}$}{ufsn}{a sum $\txsum_{j=1}^n$}
\entry{$\twist_{j,k}$}{eq:transpose_perm}{transpose permutation}
\entry{$\lambda_{j,k}$}{lex_bijection}{lexicographic ordering}
\entry{$\Gcat$}{def:GCat}{2-category of small $G$-categories}
\entry{$g \cdot x$, $gx$}{gactioniso}{$g$-action}
\entry{$\Catg$}{def:Catg}{internal hom for $\Gcat$}
\entry{$(-)^G$}{def:fixedpoint}{$G$-fixed subcategory}
\entry{$\pcom$, $\pcom_{j,k}$}{pseudocom_isos}{pseudo-commutativity isomorphism}
\entry{$\twist$}{pseudocom_isos}{swap permutation}
\entry{$\tn S$}{def:translation_cat}{translation category of a set $S$}
\entry{$[b,a] , ! \cn a \to b$}{not:b_to_a}{unique morphism from $a$ to $b$}
\entry{$\BE$}{def:BE}{Barratt-Eccles operad}
\entry{$\GBE$}{def:GBE}{$G$-Barratt-Eccles operad}
\entry{$[-,-]$}{GBEn}{$G$-set of functions}
\entry{$\ga^G$}{GBE_gamma}{composition of $\GBE$}
\entry{$\ang{a_j}_{j \in \ufs{n}}$, $\ang{a_j}_{j=1}^n$}{angaj_ufsn}{an $n$-tuple $(a_1, a_2, \ldots , a_n)$}
\entry{$(\A,\gaA,\phiA)$}{def:pseudoalgebra}{an $\Op$-pseudoalgebra}
\entry{$\gaA_n$}{gaAn}{$n$-th $\Op$-action $G$-functor}
\entry{$\zero$, $\gaA_0(*)$}{pseudoalg_zero}{basepoint}
\entry{$\phiA$}{phiA}{associativity constraint}
\entry{$\dy_j$}{not:dyj}{insert $*$ or $\zero$}
\entry{$(f,\actf)$}{def:laxmorphism}{lax $\Op$-morphism with action constraint $\actf$}
\entry{$\AlglaxO$, $\AlgpspsO$, $\AlgstO$}{oalgps_twocat}{2-categories of $\Op$-pseudoalgebras}
\entry{$\Phi(\A,\gaA,\phiA)$}{naive_smgcat}{naive symmetric monoidal $G$-category of a $\BE$-pseudoalgebra}
\entry{$(f,f^2,f^0)$}{def:smGfunctor}{symmetric monoidal $G$-functor}
\entry{$\smgcat$, $\smgcatsg$, $\smgcatst$}{def:smGcat_twocat}{2-categories of naive symmetric monoidal $G$-category}
\entry{$\Phi$}{Phi_twofunctor}{2-functor $\AlglaxBE \to \smgcat$}
\entry{$\Psi(\A,\otimes,\zero,\alpha,\xi)$}{BEpseudo_from_smcat}{$\BE$-pseudoalgebra of a naive symmetric monoidal $G$-category}

\newchHeader{ch:multpso}

\entry{$\A_{[p,q]}$, $(\A_{[p,q]})^n$, $(\A^n)_{[p,q]}$}{Apq}{products of $G$-categories}
\entry{$\shuf_{n,i}$}{shufni}{$(n,i)$-strength}
\entry{$a^j$}{aj}{tuple of objects $(\ang{a_p}_{p=1}^{i-1} \scs a_{i,j}, \ang{a_p}_{p=i+1}^{k})$}
\entry{$(f,\laxf)$}{klax_Omorphism}{a $k$-lax $\Op$-morphism with action constraint $\laxf$}
\entry{$\fgr$}{ftilde}{functor $f( \ang{a_p}_{p=1}^{i-1}, - , \ang{a_p}_{p=i+1}^{k})$}
\entry{$\MultpsO(\ang{\A_i}_{i \in \ufs{k}}; \B)$}{multpso_k}{$G$-category of $k$-lax $\Op$-morphisms and $k$-ary $\Op$-transformations}
\entry{$\MultpspsO$, $\MultstO$}{multpspso_k}{pseudo and strict variants of $\MultpsO$}
\entry{$(f\sigma, \laxfsi)$}{klax_sigma}{$\si$-action on $(f,\laxf)$}
\entry{$\theta^\si$}{Otr_sigma}{$\si$-action on $\theta$}
\entry{$\MultvO$}{not:multpsvo_ab}{$\MultpsO$, $\MultpspsO$, or $\MultstO$ depending on $\va \in \{\sflax,\sfps,\sfst\}$}
\entry{$\gam$}{gam_fh}{multicategorical composition on $\MultpsO$}
\entry{$\MultvO^G$}{thm:multpso_gfixed}{$G$-fixed $\Cat$-multicategory of $\MultvO$}

\newchHeader{ch:ggcat}

\entry{$\ord{n}$}{ordn}{pointed finite set $\{0 < 1 < \cdots < n\}$}
\entry{$\Fsk$}{rk:Fsk}{category of $\ord{n}$ for $n \geq 0$ and pointed functions}
\entry{$(\sma, \ord{1}, \xi)$}{def:Fsk_permutative}{permutative structure on $\Fsk$}
\entry{$\Fsk^{(q)}$}{def:Fsk_smashpower}{$q$-th smash power of $\Fsk$}
\entry{$\vstar$}{not:Fsk_smash_zero}{initial-terminal basepoint of $\Fsk^{(0)}$}
\entry{$\ang{}$}{not:Fsk_smash_zero}{empty tuple}
\entry{$\ang{\ord{n}}$, $\ang{\ord{n}_i}_{i \in \ufs{q}}$}{angordn}{an object $(\ord{n}_1, \ldots, \ord{n}_q)$ in $\Fsk^{(q)}$}
\entry{$\ang{\psi}$}{angpsi}{a morphism $\ang{\psi_i}_{i \in \ufs{q}}$ in $\Fsk^{(q)}$}
\entry{$\Inj$}{def:injections}{category of $\ufs{n}$ for $n \geq 0$ and injections}
\entry{$h_*$}{def:injections}{a reindexing functor}
\entry{$\ord{n}_{\emptyset}$, $\psi_{\emptyset}$}{ordn_empty}{$\ord{1}$ and $1_{\ord{1}}$}
\entry{$\Gsk$}{def:Gsk}{category with finite tuples of pointed finite sets as objects}
\entry{$\vstar$}{Gsk_objects}{initial-terminal basepoint of $\Gsk$}
\entry{$\angordm$}{Gsk_morphisms}{an object in $\Gsk$}
\entry{$(f, \ang{\psi})$, $\upom$}{fangpsi}{a morphism in $\Gsk$}
\entry{$(\oplus, \ang{}, \xi)$}{def:Gsk_permutative}{permutative structure on $\Gsk$}
\entry{$\tau_{p,q}$}{Gsk_braiding}{block permutation interchanging $p$ and $q$ elements}
\entry{$\sma$}{def:smashFskGsk}{strict symmetric monoidal functor $\Gsk \to \Fsk$ defined by smash}
\entry{$\Gcatst$}{Gcatst}{2-category of small pointed $G$-categories}
\entry{$\Catgst$}{Gcatst_smc}{internal hom for $\Gcatst$}
\entry{$\bonep$}{gcatst_unit}{discrete trivial $G$-category with two objects}
\entry{$(\smag, \gu, \brkst)$}{GGCat_smc}{symmetric monoidal closed structure on $\GGCat$}
\entry{$\Gskpunc(-,-)$}{gu_angordm}{nonzero morphisms in $\Gsk$}
\entry{$\underline{\GGCat}$}{ggcat_smggcat}{symmetric monoidal self-enriched structure on $\GGCat$}
\entry{$\ev_{\ang{}}$}{evang_gcat}{evaluation at $\ang{} \in \Gsk$}
\entry{$\underline{\GGCat}_{\ev_{\ang{}}}$}{ggcat_smgcat}{symmetric monoidal $\Gcat$-category structure on $\GGCat$}
\entry{$\End\big(\underline{\GGCat}_{\ev_{\ang{}}}\big)$}{ggcat_gcatmulti}{$\Gcat$-multicategory structure on $\GGCat$}
\entry{$\GGCat(\ang{}; f)$}{ggcat_zero_gcat}{a 0-ary multimorphism $G$-category of $\GGCat$}
\entry{$\GGCat\big(\ang{f_i}_{i \in \ufs{k}}; f \big)$}{ggcat_kary_gcat}{a $k$-ary multimorphism $G$-category of $\GGCat$}
\entry{$\angordmdot$}{angordmdot}{an object in $\Gsk^k$}
\entry{$\si\angordmdot$}{si_angordmdot}{left $\si$-action on $\angordmdot$}

\newchHeader{ch:jemg}

\entry{$\compk$}{compk}{replacing the $k$-th entry}
\entry{$\coprod_{j \in \ufs{p}}\, S_j$}{partition}{a partition of a set}
\entry{$(a,\glu)$}{nsystem}{an $\angordn$-system}
\entry{$\ang{s}$}{marker}{a marker $\ang{s_j}_{j \in \ufs{q}} = \ang{s_j \subseteq \ufs{n}_j}_{j \in \ufs{q}}$}
\entry{$a_{\ang{s}}$}{a_angs}{$\angs$-component object}
\entry{$\glu_{x; \ang{s} \csp k, \ang{s_{k,i}}_{i \in \ufs{r}}}$}{gluing-morphism}{a gluing morphism}
\entry{$(\zero,1_\zero)$}{not:basesystem}{base $\angordn$-system}
\entry{$\theta_{\angs}$}{theta_angs}{$\angs$-component morphism}
\entry{$\Aangordn$}{not:Aangordn}{category of $\angordn$-systems in $\A$}
\entry{$\Aangordnsg$}{not:Aangordn}{full subcategory of $\Aangordn$ of strong $\angordn$-systems}
\entry{$\ftil$}{ftil_functor}{pointed $G$-functor $\Aangordm \to \Afangordm$}
\entry{$\left(\atil, \glutil\right)$}{ftil_aglu}{$f_*\angordm$-system $\ftil(a,\glu)$}
\entry{$\ftil_*\ang{s}$}{ftil_angs}{marker $\ang{s_{f(i)} \subseteq \ufs{m}_i}_{i \in \ufs{p}}$}
\entry{$\thatil$}{thatil_component}{morphism of $f_*\angordm$-systems $\ftil(\theta)$}
\entry{$\psitil$}{psitil_functor}{pointed $G$-functor $\Afangordm \to \Aangordn$}
\entry{$(a^{\psitil}, \glu^{\psitil})$}{angs_ordtun}{$\angordn$-system $\psitil(a,\glu)$}
\entry{$\psiinv\ang{s}$}{apsitil_angs}{marker $\ang{\psiinv_j s_j \subseteq \ufs{m}_{\finv(j)}}_{j \in \ufs{q}}$}
\entry{$\tha^{\psitil}$}{thapsitil_component}{morphism of $\angordn$-systems $\psitil(\theta)$}
\entry{$\Aupom$, $\Afpsi$}{AF}{pointed $G$-functor $\Aangordm \to \Aangordn$}
\entry{$\Aupomsg$}{AF_sg}{strong variant of $\Aupom$}
\entry{$\Adash$}{A_ptfunctor}{pointed functor $\Gsk \to \Gcatst$ defined by an $\Op$-pseudoalgebra $\A$}
\entry{$\Asgdash$}{A_ptfunctor}{strong variant of $\Adash$}
\entry{$\Jgo a$}{Jgoa}{object in $\GGCat(\ang{} \sscs \Adash)$ defined by $a \in \A$}
\entry{$\Jgosg a$}{Jgosga}{strong variant in $\GGCat(\ang{} \sscs \Asgdash)$}
\entry{$\Jgo f$}{Jgo_f}{object $\ang{\Aidash}_{i \in \ufs{k}} \to \Bdash$ defined by $(f,\laxf)$}
\entry{$a_{i,\, s_{i\crdot}}$}{Jgof_m_obj_comp}{$s_{i\crdot}$-component object of $(a_i,\glu^{a_i})$}
\entry{$\fgr$}{ftil_ApB}{pointed functor $f(\ang{a_{i,\, s_{i\crdot}}}_{i \in \ufs{k}}\, \compp - )$}
\entry{$\Jgo\theta$}{Jgotheta}{image of a $k$-ary $\Op$-transformation $\theta$ under $\Jgo$}
\entry{$\Jgo$, $\Jgosg$}{def:Jgo_multifunctor}{(strong) $J$-theory $\Gcat$-multifunctors}
\entry{$\Jgbe$, $\Jsgbe$, $\Jggbe$, $\Jsggbe$}{ex:JgBE}{multifunctorial (strong) $J$-theories of $\BE$ and $\GBE$}
\entry{$(\Gtop, \times, *, \Topg)$}{gtop_smclosed}{Cartesian closed category of $G$-spaces}
\entry{$(\cla,\clatwo,\clazero)$}{cla_gcat_gtop}{classifying space strong symmetric monoidal functor $\Gcat \to \Gtop$}
\entry{$\univ$}{def:g_universe}{complete $G$-universe}
\entry{$\Stein_V$, $\KU$}{ex:Einf_steiner}{Steiner operads}
\entry{$\LU$}{ex:Einf_linear}{linear isometries operad}

\newchHeader{ch:spectra}

\entry{$(\Gtopst, \sma, \stplus, \Topgst)$}{Gtopst_smc}{symmetric monoidal closed category of pointed $G$-spaces}
\entry{$\Topgst$}{topgst_gtopst_enr}{$\Gtopst$-category of pointed $G$-spaces}
\entry{$X_\splus$}{conv:disjoint_gbasept}{adjoining a disjoint $G$-fixed basepoint}
\entry{$\IU$}{def:indexing_gspace}{$G$-inner product spaces isomorphic to indexing $G$-spaces}
\entry{$S^V$}{def:indexing_gspace}{$V$-sphere}
\entry{$(\oplus,0,\xi)$}{def:iu_spaces_i}{symmetric monoidal structure on $\IU$}
\entry{$\IUsk$}{def:iu_spaces_iv}{small skeleton of $\IU$ of indexing $G$-spaces in $\univ$}
\entry{$\upphi_V$}{VV'}{chosen $G$-linear isometric isomorphism $V \fiso V'$ with $V' \in \IUsk$}
\entry{$\IU \to \Topgst$}{def:iu_space}{an $\IU$-space}
\entry{$X_V$}{expl:iu_space}{pointed $G$-space of an $\IU$-space $X$ and $V \in \IU$}
\entry{$X_f$}{iu_space_xf}{pointed homeomorphism of an $\IU$-space $X$ and $f \cn V \fiso W$}
\entry{$\IUT$}{IUT}{category of $\IU$-spaces and $\IU$-morphisms}
\entry{$\GIUT$}{GIUT}{category of $\IU$-spaces and $G$-equivariant $\IU$-morphisms}
\entry{$X \smau Y$}{x_smau_y}{smash product of $\IU$-spaces}
\entry{$\theta \smau \theta'$}{tha_smau_thap}{smash product of $\IU$-morphisms}
\entry{$(\smau, \iu, \au, \ellu, \ru, \beu)$}{def:IU_smgtop}{symmetric monoidal $\Gtop$-category structure on $\IUT$ or $\GIUT$}
\entry{$\AMod$}{amodule_mor}{$\Gtop$-category of $A$-modules and $A$-module morphisms}
\entry{$\umu_{V,W}$}{gmonoid_mod_vw}{$(V,W)$-component of $\umu \cn X \smau A \to X$}
\entry{$(\gsp,\mu,\eta)$}{def:g_sphere}{$G$-sphere}
\entry{$\GSp$, $\gspmod$}{def:gsp_module}{$\Gtop$-category of orthogonal $G$-spectra}
\entry{$X \smasg Y$}{def:gsp_sma}{smash product of orthogonal $G$-spectra}
\entry{$\psma$}{gsp_sma_coequal}{universal morphism from $X \smau Y$ to $X \smasg Y$}
\entry{$\upbe$}{gsp_sma_beta}{a coherence $\IU$-isomorphism}
\entry{$e'$}{xsmasgy_iden}{$e$ pre-composed with a unique coherence isomorphism in $(\IUsk,\oplus)$}
\entry{$\theta \smasg \theta'$}{def:gsp_mor_sma}{smash product of morphisms of orthogonal $G$-spectra}
\entry{$(\smasg, \gsp, \asg, \ellsg, \rsg, \bsg)$}{def:gsp_smgtop}{symmetric monoidal $\Gtop$-category structure on $\GSp$}

\newchHeader{ch:ggtop}

\entry{$(\GGTop, \smag, \gu, \brkst)$}{ggtop_smc}{symmetric monoidal closed category of $\Gskg$-spaces}
\entry{$\stplus$}{gu_angordm_ggtop}{smash unit in $\Gtopst$}
\entry{$\clast$}{clast}{symmetric monoidal functor $\GGCat \to \GGTop$ induced by $\cla$}
\entry{$(\clast')_{\evang}$}{clastpev}{symmetric monoidal $\Gtop$-functor induced by $\cla$}
\entry{$\End ((\clast')_{\evang})$}{End_clastpev}{$\Gtop$-multifunctor induced by $\cla$}
\entry{$\evang$}{evang_gtop}{evaluation at $\ang{} \in \Gsk$}
\entry{$\clastzero$}{clastzero}{unit constraint of $\clast$}
\entry{$\clasttwo$}{clasttwo}{monoidal constraint of $\clast$}
\entry{$\clatwobar$}{clasttwo_one}{morphism induced by $\clatwo$}
\entry{$(\smag, \gu, \ag, \ellg, \rg, \beg)$}{ggtop_smgtop}{symmetric monoidal $\Gtop$-category structure on $\GGTop$}

\newchHeader{ch:semg}

\entry{$(\Kg X, \umu)$}{Kg_object}{orthogonal $G$-spectrum associated to a $\Gskg$-space $X$}
\entry{$(S^V)^{\sma\angordn}$}{SVsman}{pointed $G$-space $\Topgst(\sma\angordn, S^V)$}
\entry{$\assm$, $\assm_{\angordn}$}{Kgx_action_vw}{pointed $G$-morphisms defining $\umu_{V,W}$ for $\Kg X$}
\entry{$\Kg$}{def:ggtop_gsp_mor}{$\Gtop$-functor $\GGTop \to \GSp$}
\entry{$\Kg\theta$}{Kg_theta}{image of $\theta \in \brk{X}{X'}$ under $\Kg$}
\entry{$\Kgzero$}{Kgzero}{unit constraint of $\Kg$}
\entry{$\Kgtwo$}{Kgtwo}{monoidal constraint of $\Kg$}
\entry{$\ktwo$}{smau_kgx}{$G$-equivariant $\IU$-morphism that induces $\Kgtwo_{X,Y}$}
\entry{$\smas$}{smash_s}{pointed $G$-morphism that induces $\ktwo_U$}
\entry{$(\Kg,\Kgtwo,\Kgzero)$}{Kg_smfunctor}{unital symmetric monoidal $\Gtop$-functor $\GGTop \to \GSp$}
\entry{$\End(\Kg)$}{EndKg}{$\Gtop$-multifunctor $\GGTop \to \GSp$}
\entry{$\Kgo$}{thm:Kgo_multi_i}{$\Gtop$-multifunctor $\MultpsO \to \GSp$}
\entry{$\Kgosg$}{thm:Kgo_multi_ii}{strong variant of $\Kgo$}

\newchHeader{ch:prelim}

\entry{$(\C,\otimes,\tu,\alpha,\lambda,\rho)$}{def:monoidalcategory}{a monoidal category}
\entry{$\xi$}{braiding_bmc}{braiding}
\entry{$\Hom$, $[-,-]$}{notation:internal-hom}{internal hom}
\entry{$\DC$}{not:DC}{diagram category of functors $\cD \to \C$ and natural transformations}
\entry{$(a,\mu,\eta)$}{notation:monoid}{a monoid}
\entry{$(x,\umu)$}{def:modules_i}{a right module}
\entry{$(F, F^2, F^0)$}{def:monoidalfunctor}{a monoidal functor}
\entry{$(\C,\mcomp,\cone)$}{def:enriched-category}{a $\V$-category with composition $\mcomp$ and identity $\cone$}
\entry{$\Catv$}{ex:vcatastwocategory}{2-category of small $\V$-categories}
\entry{$\C \otimes \D$}{definition:vtensor-0}{tensor product of $\V$-categories}
\entry{$\vtensorunit$}{definition:unit-vcat}{unit $\V$-category}
\entry{$\ell^\otimes$, $r^\otimes$}{definition:vtensor-unitors}{left and right unitors for $\otimes$}
\entry{$a^\otimes$}{definition:vtensor-assoc}{associator for $\otimes$}
\entry{$\beta^\otimes$}{definition:vtensor-beta}{braiding for $\otimes$}
\entry{$(\C,\vmtimes,\vmunitbox,a^\vmtimes,\ell^\vmtimes,r^\vmtimes)$}{definition:monoidal-vcat}{a monoidal $\V$-category}
\entry{$a^\vmtimes_1$}{eq:a-vmtimes-inv-mate}{mate of $a^\vmtimes$}
\entry{$\beta^\vmtimes$}{bmvcat_braiding}{braiding of a braided monoidal $\V$-category}
\entry{$\ell^\vmtimes_1$, $r^\vmtimes_1$}{eq:lr-vmtimes-mates}{mates of $\ell^\vmtimes$ and $r^\vmtimes$}
\entry{$(F,F^2,F^0)$}{definition:monoidal-V-fun}{a monoidal $\V$-functor}
\entry{$\ev$}{evaluation}{evaluation, counit of tensor-hom adjunction}
\entry{$\Vse$}{definition:canonical-v-enrichment}{canonical self-enrichment of a symmetric monoidal closed category $\V$}
\entry{$\C_f$}{thm:change-enrichment}{change of enrichment of $\C$ along $f$}

\newchHeader{ch:prelim_multicat}

\entry{$\Prof(C)$}{notation:profs}{class of $C$-profiles, $\coprod_{m \geq 0} C^{m}$}
\entry{$\angc$, $\ang{c_j}_{j=1}^n$}{notation:us}{a profile $(c_1, \ldots, c_n)$ of length $n$}
\entry{$\ang{}$, $\emptyset$}{notation:us}{empty profile}
\entry{$\angc \oplus \angd$}{not:concat}{concatenation of profiles}
\entry{$(\M, \gamma, \operadunit)$}{notation:enr-multicategory}{a $\V$-multicategory with composition $\gamma$ and unit $\operadunit$}
\entry{$\smscmap{\angx;x'}$}{notation:enr-cduc}{an element in $\ProfMM$}
\entry{$\M\scmap{\angx;x'}$}{multimorphism_object}{an $n$-ary multimorphism $\V$-object}
\entry{$\angx\sigma$}{enr-notation:c-sigma}{$\ang{x_{\sigma(j)}}_{j=1}^n$, right $\si$-action on $\angx$}
\entry{$a\sigma$, $a \cdot \sigma$, $a^\sigma$}{enr-notation:c-sigma}{right $\si$-action}
\entry{$\sigma\langle k_{\sigma(1)}, \ldots , k_{\sigma(n)} \rangle$}{blockpermutation}{block permutation induced by $\si$}
\entry{$\tau_1 \times\cdots \times\tau_n$}{blocksum}{block sum}
\entry{$\M_n$, $\M(n)$}{not:nthobject}{$n$-th object of a $\V$-operad}
\entry{$\End(\C)$}{definition:EndK}{endomorphism $\V$-multicategory of a symmetric monoidal $\V$-category $\C$}
\entry{$\End F$}{def:EndF}{endomorphism $\V$-multifunctor of a symmetric monoidal $\V$-functor $F$}
\entry{$\VMulticat$}{v-multicat-2cat}{2-category of small $\V$-multicategories}
\entry{$\dashf$}{mult_change_enr}{change-of-enrichment 2-functor $\VMulticat \to \WMulticat$}
\entry{$\pA$}{def:pointed-objects}{category of pointed objects in $\A$ with a terminal object $\term$}
\entry{$(a,i^a)$}{not:aia}{a pointed object}
\entry{$\wed$}{not:awedgeb}{wedge}
\entry{$\sma$}{eq:smash}{smash product}
\entry{$\stu$}{smash-unit-object}{smash unit $\tu \bincoprod \term$}
\entry{$[,]_*$}{eq:pHom-pullback}{pointed hom}
\entry{$\zob$}{definition:zero}{a zero object}
\entry{$\Dpunc(a,b)$}{definition:zero}{set of nonzero morphisms $a \to b$}
\entry{$\Dhat$}{eq:Dhat-ptd-unitary-enr}{pointed unitary enrichment}
\entry{$(\D,\bpt)$}{def:pointed-category}{a pointed category with basepoint $\bpt$}
\entry{$\DstarC$}{not:DstarC}{pointed diagram category}
\entry{$\du$}{ptdayunit}{monoidal unit diagram}
\entry{$\sma$}{ptdayconv}{pointed Day convolution}
\entry{$[,]_*$}{ptdayhom}{pointed hom diagram}
\entry{$\pMap$}{ptdaymap}{pointed mapping object}
\entry{$\left( \DstarV, \sma, \du, [,]_* \right)$}{thm:Dgm-pv-convolution-hom}{symmetric monoidal closed pointed diagram category}

\newchHeader{ch:nerve}

\entry{$\Delta$}{not:Delta_cat}{category of $\ordn$ for $n \geq 0$ and weakly order-preserving functions}
\entry{$d^i$, $s^i$}{not:face_degen}{coface and codegeneracy}
\entry{$\SSet$}{not:SSet}{category of simplicial sets}
\entry{$X_n$, $d_i$, $s_i$}{not:n_simplex}{set of $n$-simplices, $i$-th face, and $i$-th degeneracy}
\entry{$\Ner$}{nerve}{nerve functor $\Cat \to \SSet$}
\entry{$\simp^n$}{not:simpn}{topological $n$-simplex}
\entry{$\Rea$}{realization}{geometric realization functor $\SSet \to \Top$}
\entry{$\cla$}{classifying_space}{classifying space functor $\Cat \to \Top$}

\end{tabbing}

\printindex
\end{document}